\documentclass[english,a4paper,leqno]{smfbook}
\usepackage{hyperref}
\usepackage[all]{xy}

\newcommand{\vep}{\varepsilon}
\newcommand{\ds}{\displaystyle}
\newcommand{\Sz}{\mathrm{Sz}}
\newcommand{\Cz}{\mathrm{Cz}}
\newcommand{\Kz}{\mathrm{Kz}}

\newcommand{\bE}{\ensuremath{\mathbb E}}

\newcommand{\bL}{\ensuremath{\mathbb L}}
\newcommand{\bN}{\ensuremath{\mathbb N}}
\newcommand{\bM}{\ensuremath{\mathbb M}}
\newcommand{\bP}{\ensuremath{\mathbb P}}
\newcommand{\bQ}{\ensuremath{\mathbb Q}}
\newcommand{\bR}{\ensuremath{\mathbb R}}

\newcommand{\bZ}{\ensuremath{\mathbb Z}}

\newcommand{\sA}{\ensuremath{\mathsf A}}
\newcommand{\sB}{\ensuremath{\mathsf B}}

\newcommand{\sD}{\ensuremath{\mathsf D}}

\newcommand{\sH}{\ensuremath{\mathsf H}}
\newcommand{\sI}{\ensuremath{\mathsf I}}
\newcommand{\sJ}{\ensuremath{\mathsf J}}
\newcommand{\sK}{\ensuremath{\mathsf K}}
\newcommand{\sN}{\ensuremath{\mathsf N}}

\newcommand{\sP}{\ensuremath{\mathsf P}}

\newcommand{\sT}{\ensuremath{\mathsf T}}

\newcommand{\cA}{\ensuremath{\mathcal A}}
\newcommand{\cB}{\ensuremath{\mathcal B}}
\newcommand{\cC}{\ensuremath{\mathcal C}}

\newcommand{\cF}{\ensuremath{\mathcal F}}
\newcommand{\cG}{\ensuremath{\mathcal G}}
\newcommand{\cH}{\ensuremath{\mathcal H}}

\newcommand{\cL}{\ensuremath{\mathcal L}}
\newcommand{\cM}{\ensuremath{\mathcal M}}
\newcommand{\cN}{\ensuremath{\mathcal N}}

\newcommand{\cP}{\ensuremath{\mathcal P}}
\newcommand{\cQ}{\ensuremath{\mathcal Q}}

\newcommand{\cS}{\ensuremath{\mathcal S}}

\newcommand{\cU}{\ensuremath{\mathcal U}}
\newcommand{\cV}{\ensuremath{\mathcal V}}

\newcommand{\cZ}{\ensuremath{\mathcal Z}}

\newcommand{\alphabar}{\ensuremath{\bar{\alpha}}}
\newcommand{\betabar}{\ensuremath{\bar{\beta}}}
\newcommand{\abar}{\ensuremath{\bar{a}}}
\newcommand{\ebar}{\ensuremath{\bar{e}}}
\newcommand{\fbar}{\ensuremath{\bar{f}}}

\newcommand{\lbar}{\ensuremath{\bar{\ell}}}
\newcommand{\mbar}{\ensuremath{\bar{m}}}
\newcommand{\nbar}{\ensuremath{\bar{n}}}
\newcommand{\pbar}{\ensuremath{\bar{p}}}

\newcommand{\sbar}{\ensuremath{\bar{s}}}

\newcommand{\wbar}{\ensuremath{\bar{w}}}
\newcommand{\xn}{\ensuremath{(x_n)_{n=1}^\infty}}
\newcommand{\yn}{\ensuremath{(y_n)_{n=1}^\infty}}

\newcommand{\zn}{\ensuremath{(z_n)_{n=1}^\infty}}
\newcommand{\xnstar}{\ensuremath{(x_n^*)_{n=1}^\infty}}

\newcommand{\en}{\ensuremath{(e_n)_{n=1}^\infty}}

\newcommand{\ei}{\ensuremath{(e_i)_{i=1}^\infty}}

\newcommand{\coo}{\mathrm{c}_{00}}

\newcommand{\infin}[1]{{[#1]}^{\omega}}

\newcommand{\wtree}[2]{[{#1}]^{{#2}}}
\usepackage{smfthm}
\usepackage{enumerate}
\usepackage{amssymb}
\usepackage[english]{babel}
\usepackage{xcolor}
\usepackage{marginnote} % add margin notes

\theoremstyle{definition}
\newtheorem{exer}{Exercise}
\newtheorem{claim}{Claim}
\newtheorem{exam}{Example}

\newtheorem{prob}{Problem}

\newcommand{\bin}{\ensuremath{\mathsf{B}}}

\newcommand{\tree}{\ensuremath{\mathsf{T}}}

\newcommand{\dHk}{\ensuremath{\mathsf{d_{H_k}}}}
\newcommand{\wHk}{\ensuremath{\mathsf{\bar{w}H_k}}}
\newcommand{\dwHk}{\ensuremath{\mathsf{d_{\bar{w}H_k}}}}

\newcommand{\wTk}{\ensuremath{\mathsf{T^\omega_k}}}
\newcommand{\dwTk}{\ensuremath{\mathsf{d_{T^\omega_k}}}}
\newcommand{\wT}{\ensuremath{\mathsf{T^\omega_\infty}}}
\newcommand{\dwT}{\ensuremath{\mathsf{d_{T^\omega_\infty}}}}

\newcommand{\dJk}{\ensuremath{\mathsf{d_{J_k}}}}

\newcommand{\dIk}{\ensuremath{\mathsf{d_{I_k}}}}
\newcommand{\Ik}{\ensuremath{\mathsf{I_k}}}

\newcommand{\dIuk}{\ensuremath{\mathsf{d_{I_k^{\omega_1}}}}}
\newcommand{\Iuk}{\ensuremath{\mathsf{I_k^{\omega_1}}}}

\newcommand{\dia}{\ensuremath{\mathsf{D}}}
\newcommand{\diak}{\mathsf{D_k}}
\newcommand{\dwdiak}{\ensuremath{\mathsf{d_{D^{\omega}_{k}}}}}
\newcommand{\wdiak}{\mathsf{D^\omega_k}}

\newcommand{\sd}{\ensuremath{\mathsf{d}}}

\newcommand{\dsym}{\ensuremath{\mathsf{d_\triangle}}}

\newcommand{\cdist}[1]{\mathsf{c}_{#1}}
\newcommand{\sep}{\ensuremath{\mathrm{sep}}}
\newcommand{\US}{\mathbf{US}}
\newcommand{\UC}{\mathbf{UC}}
\newcommand{\AUS}{\mathbf{AUS}}
\newcommand{\AUC}{\mathbf{AUC}}
\newcommand{\AUF}{\mathbf{AUF}}
\newcommand{\REF}{\mathbf{REF}}
\newcommand{\BETA}{\mathbf{BETA}}
\newcommand{\ASO}{\mathbf{Asymp-c}_0}

\newcommand{\James}{\mathrm{J}}
\newcommand{\Tsi}{\mathrm{T}}
\newcommand{\eqd}{:=}
\newcommand{\xbar}{\bar{x}}

\newcommand{\sign}{\mathrm{sign}}

\newcommand{\cal}{\mathcal}
\newcommand{\eps}{\varepsilon}
\newcommand{\bA}{\ensuremath{\mathbb A}}
\newcommand{\bB}{\ensuremath{\mathbb B}}

\newcommand{\Mdb}{\mathbb M}

\newcommand{\Rdb}{\mathbb R}

\newcommand{\C}{\mathbb C}
\newcommand{\M}{\mathbb M}
\newcommand{\N}{\mathbb N}

\newcommand{\R}{\mathbb R}

\newcommand{\Z}{\mathbb Z}

\def\H{{\mathbb H}}

\newcommand{\Nk}{[\N]^k}

\newcommand{\car}{\textrm{\bbold 1}}
\font\bbold=bbold12

\newcommand{\cof}{\text{\rm cof}}
\newcommand{\spa}{\text{\rm span}}
\newcommand{\cspa}{\overline{\text{\rm span}}}
\newcommand{\supp}{\text{\rm supp}}
\newcommand{\ie}{{\it i.e.\ }}
\newcommand{\co}{\mathrm{c}_0}
\newcommand{\lip}{\text{\rm lip}}
\newcommand{\Lip}{\text{\rm Lip}}
\newcommand{\conv}{\text{\rm conv}}
\newcommand{\coLip}{\text{co-Lip}}

\newcommand{\n}{\overline n}
\newcommand{\m}{\overline m}

\newcommand{\wtoo}{\stackrel{w}{\longrightarrow}}
\newcommand{\wstoo}{\stackrel{w^*}{\longrightarrow}}

\newcommand{\bib}{\bibitem}

\newcommand{\abs}[1]{\left\vert#1\right\vert} %absolute value, use: $\abs{x}$
\newcommand{\norm}[1]{\left\Vert#1\right\Vert} %\norm, use: $\norm{x}$ or $\norm{x}_{X}$
\newcommand{\bnorm}[1]{\big\Vert#1\big\Vert} %\norm, use: $\bnorm{x}$ or $\bnorm{x}_{X}$
\newcommand{\set}[1]{\left\{#1\right\}} %set, use: $\set{1,2,\ldots}$
\newcommand{\diam}{\mathop{\mathrm{diam}}\nolimits} %diameter, use: $\diam(A)$
\newcommand{\dist}{\mathop{\mathrm{dist}}\nolimits} %the distance function or distortion, use: $\dist(x,A)$ or $\dist(f)$
\newcommand{\closedball}[1]{B_{#1}} %the unit ball or any ball of $X$, use: $\closedball{X}$ or $\closedball{X}(x,r)$
\newcommand{\sphere}[1]{S_{#1}} %the unit sphere, use: $\sphere{X}$
\newcommand{\linspan}{\mathrm{span}} %linear span, use: $\linspan\{x_n\}$ or $\linspan_\Rational \set{x_n}$

\newcommand{\restricted}{\upharpoonright} %the symbol of restriction, use: $f\restricted_A$

\newcommand{\Mid}{\mathrm{Mid}}%approximate midpoints, use: $\Mid(x,y,\delta)$
\newcommand{\indicator}[1]{{\mathbf 1}_{{#1}}}

\renewcommand{\le}{\leqslant}
\renewcommand{\ge}{\geqslant}

\usepackage[only,llbracket,rrbracket]{stmaryrd}

\begin{document}
	
	\title{Asymptotic and nonlinear geometries of Banach spaces and their interactions}
	
	\author[F. P. Baudier]{Florent P. Baudier}
	\address[F. P. Baudier]{Department of Mathematics, Texas A\&M University, College Station, TX 77843, USA}
	\email{florent@tamu.edu}

	\author[G. Lancien]{Gilles Lancien}
	\address[G. Lancien]{ Universit\'e Marie et Louis Pasteur, CNRS, LmB (UMR 6623), F-25000 Besançon, France}
	\email{gilles.lancien@univ-fcomte.fr}

	\maketitle 
	
	%\hskip 6cm \`A mes parents Chantal et Pierre.
	
	%\vskip 1cm
	
	%\hskip 4cm \`A mes parents André Lancien et Annick Garlantézec.
	
	\hskip 2cm\`A nos parents, 
	
	\hskip 6cm Andr\'e Lancien et Annick Garlant\'ezec, 
	
	\hskip 5cm Chantal Baudier (n\'ee Benoit) et Pierre Baudier.
	
	\setcounter{tocdepth}{5}
	\tableofcontents

	\chapter*{Introduction}
	%\addcontentsline{toc}{chapter}{Introduction}
	
	Banach spaces are fundamental examples of metric spaces and, in turn, of topological spaces. Thus, they can be considered from different perspectives whether one would like, or not, to take all the available structures into consideration, just a selection of them, or the bare minimum. Two quintessential results of the nonlinear geometry of Banach spaces are the Mazur-Ulam rigidity theorem and its extreme opposite, namely Kadets' theorem. The Mazur-Ulam theorem tells us that every surjective isometry fixing the origins between two Banach spaces is automatically a \emph{linear} isometry. Therefore, one can recover the linear structure of a Banach space simply from the knowledge with absolute precision, of its metric structure. At the other end of the spectrum, Kadets showed that two infinite-dimensional separable Banach spaces are always homeomorphic and hence it is impossible to distinguish these Banach spaces from topological data only. From a classification standpoint, these two results imply that in the category of metric spaces and isometries, every Banach space is a unique member of its equivalence class, whereas in the category of topological spaces and homeomorphisms, all the infinite-dimensional separable Banach spaces belong to the same equivalence class. The nonlinear geometry of Banach spaces originated, for the most part, in the desire to understand the classification problem for categories that lie between the two mentioned above. The dominating point of view was to consider Banach spaces in the category of metric spaces with morphisms whose degree of faithfulness lies between isometry and homeomorphism. For a long time, the most studied notions have arguably been Lipschitz isomorphisms and uniform isomorphisms, until coarse isomorphisms became also part of the picture more recently. Despite the efforts of many mathematicians, the Lipschitz classification problem remains open in the separable setting. This problem asks whether two separable Banach spaces that are Lipschitz isomorphic must already be linearly isomorphic. To this date, the authoritative monograph ``Geometric Nonlinear Functional Analysis'' by Y. Benyamini and J. Lindenstrauss \cite{BenyaminiLindenstrauss2000}, which appeared in 2000, is still the most thorough exposition of the beautiful theory that developed around nonlinear classification problems and beyond. A more concise exposition with which an outsider should probably start to get acquainted with the nonlinear geometry of Banach spaces is Chapter 14 of the second edition of the book ``Topics in Banach Space Theory'' by F. Albiac and N.~J. Kalton \cite{AlbiacKalton2016}. We also refer the reader to the book ``Nigel Kalton’s lectures in nonlinear functional analysis'' by A. Bowers \cite{Bowers2024}. One of the key approaches to tackle nonlinear classification problems was to understand to what extent a nonlinear map can be adequately linearized so that the nonlinear problems can be reduced to linear ones.
	
	An influential rigidity phenomenon discovered by M. Ribe in the mid-1970s opened up an entirely new research direction. Ribe's rigidity theorem states that uniformly homeomorphic Banach spaces are crudely finitely representable in each other. In looser terms, Banach spaces that look similar in a certain metric sense have essentially the same finite-dimensional subspaces.
	The (usually daunting) task of reformulating local properties of Banach spaces in purely metric terms became known as the Ribe program. This program was promoted early on by J. Lindenstrauss and later put forth by J. Bourgain. We refer to surveys by A. Naor \cite{Naor2012} and K. Ball \cite{Ball2013} for more details on the history, the goals and the stellar achievements in this program. Over the past 40 years, the explosion of work motivated by Ribe's rigidity theorem and which established a connection with theoretical computer science and geometric group theory, certainly deserves a book (or two!) on its own. The book ``Metric Embeddings'' by M.I. Ostrovskii \cite{Ostrovskii_book13} discusses, amongst many other topics, some aspects of the Ribe program. 
	
	The material presented in our book is not an exposition of the Ribe program but it is nevertheless strongly influenced by it. From its infancy in the mid-1980s, the metric reformulation of local properties of Banach spaces has been tied to the geometry of certain sequences of finite graphs or infinite but locally finite graphs. In the mid-2000s Nigel Kalton was probably the first to establish striking connections between the geometry of infinite non-locally finite graphs and infinite-dimensional or asymptotic properties of Banach spaces. For a first overview of his results, we refer to the survey article \cite{GLZ2014}.  Eventually, the task of finding purely metric characterizations of asymptotic properties of Banach spaces was referred to as the asymptotic Ribe program (or similar designations). It is worth pointing out that while several asymptotic properties have been metrically characterized, there is still currently no known asymptotic analog of Ribe's rigidity theorem! To acknowledge the tremendous influence of N. Kalton, beyond his untimely passing, on this research area and to clear any confusion that previous terminologies incurred, this program is now called the Kalton program.  
	
	In some sense, this book picks up where Benyamini and Lindenstraus' book left off. One of our goals was to cover some of the most striking advances in the nonlinear geometry of Banach spaces that have appeared over the past 25 years. As we said, we essentially left out the fascinating developments directly motivated by the Ribe program and their applications. Instead, we decided to focus on how the nonlinear geometry of a Banach space and its asymptotic geometry interact. By asymptotic geometry, we do not mean the geometric structure of finite-dimensional subspaces as the dimension grows to infinity as in the book ``Asymptotic Theory of finite-dimensional Normed Spaces'' by Milman and Schechtman \cite{MilmanSchechtman1986}. Asymptotic geometry, in our case, is concerned with the behavior of a Banach space that is governed by its asymptotic subspaces. An asymptotic subspace is a finite-dimensional subspace that can be realized in a specific asymptotic fashion. There are various ways to define what an asymptotic subspace could be. Each of them comes with its advantages and drawbacks and one aspect of this book is to describe those when it comes to studying nonlinear problems. Various aspects of the local theory of Banach spaces have already received monograph treatments by B. Beauzamy in his ``Introduction to Banach Spaces and their Geometry'' \cite{Beauzamy1985}, N. Tomczak-Jaegermann's ``Banach-Mazur Distances and Finite-Dimensional Operator Ideals'' \cite{Tomczak1989}, G. Pisier's ``The Volume of Convex Bodies and Banach Space Geometry'' and ``Absolutely Summing Operators'' by J. Diestel, H. Jarchow and A. Tonge \cite{DiestelJarchowTonge1995}. The same cannot be said about the asymptotic theory except for Beauzamy and Laprest\'e's book ``Mod\`eles \'etal\'es des espaces de Banach'' \cite{BeauzamyLapreste1984} with its thorough treatment of spreading models. We are not aware of a detailed treatment in a book form of the asymptotic geometry of Banach spaces, in particular of its two central notions: asymptotic uniform convexity and asymptotic uniform smoothness and their connections with the Szlenk index. We hope that the extended study of these linear notions provided in Chapters \ref{chapter:asymptotic-moduli} and \ref{chapter:Szlenk} of this book will help them gain visibility and traction outside of the Banach space theory community. In the first three chapters, we introduce the basic notions (Chapter \ref{chapter:basics}), we recall fundamental results about metrically universal spaces (Chapter \ref{chapter:universal}) and the linearization of Lipschitz maps machinery (Chapter \ref{chapter:linear-reductions}). There is hardly any originality added to these classical topics already covered in Benyamini-Lindenstrauss or Albiac-Kalton, but it is a natural place to start the conversation. In Chapter \ref{chapter:Lipschitz-free}, we present the Lipschitz-free space technique of Godefroy and Kalton for the isometric embedding rigidity problem. In Chapter \ref{chapter:asymptotic-moduli}, the notions of asymptotic uniform smoothness, asymptotic uniform convexity and Rolewicz's property $(\beta)$ are introduced and discussed. Their renorming theories in connection with the Szlenk index are detailed in Chapter \ref{chapter:Szlenk}. The approximate midpoint principle and the Gorelik principle circle of ideas and their applications to nonlinear classification problems form the content of Chapter \ref{chapter:AMP_I} and Chapter \ref{chapter:Gorelik}, respectively. From Chapter \ref{chapter:Johnson} on is what can be considered as the essence of this book: the interaction between the asymptotic geometry and the nonlinear geometry of Banach spaces. The ideas and vision of Nigel Kalton permeate virtually every corner of the theory presented there. If the reader is eager to understand what the Kalton program is about, this is the place to go!  
	
	We have tried to keep the book to a reasonable length. Since our primary focus was on asymptotic properties of Banach spaces, we had to make difficult choices when writing this book, and, unfortunately, a few major results about the metric characterization of merely infinite-dimensional properties are not covered. For instance, we did not include Ostrovskii's work on metric characterizations of the Radon-Nikod\'{y}m property. In \cite{Ostrovskii2014}, M.I. Ostrovskii proved that a dual Banach space does not have the Radon-Nikod\'{y}m property if and only if it admits a bi-Lipschitz embedding of an infinite diamond and in \cite{Ostrovskii2014b}, a characterization of the Radon-Nikod\'{y}m property in terms of the existence of bi-Lipschitz copies of thick families of geodesics is given. The reader interested in this topic is advised to consult \cite{Ostrovskii_survey} (and also \cite{Gartland}). In \cite{MS17}, P. Motakis and Th. Schlumprecht gave a metric characterization of the reflexivity of a Banach space in terms of the geometry of the Schreier family $(\cal S_\alpha)_{\alpha <\omega_1}$. The Motakis-Schlumprecht result is a technical tour de force that incorporates many ideas and tools presented in this book and we refer to the Notes section of Chapter \ref{chapter:trees} for more on this topic. 
	We also did not delve into operator versions of spatial properties and thus the work of Causey and Dilworth in \cite{CauseyDilworth18} is not presented here. Finally and regrettably, we have chosen to leave aside the important interactions between the nonlinear geometry of Banach spaces and descriptive set theory. Any other omission would certainly be due to the lack of expertise of the authors of this book!
	
	\medskip\noindent
	{\bf Acknowledgments.}
	First of all, we wish to express our gratitude to Sophie Grivaux, who initiated this project by kindly suggesting that we write a book on this subject for the series ``Cours Spécialisés'' of the ``Société Mathématique de France''. We also thank her for her infinite patience and support, even though this work took much longer than we initially expected. We also owe special thanks to Tony Proch\'azka, who considered pursuing this project with us. His initial enthusiasm and our multiple discussions on the subject have been very important to us. We are especially grateful to Gilles Godefroy, from whom we learned a large part of what is in this book. We also want to thank him for his constant and caring support. We thank all the friends and colleagues who read preliminary versions of the book, answered our questions (whether they knew it was related to the book or not!) and sent us comments or corrections that helped us to improve it. We particularly want to thank Ram\'on Aliaga, Estelle Basset, Bruno de Mendonça Braga, Ryan Causey, Daniel Li, Audrey Fovelle, Chris Gartland, Gilles Godefroy, Guillaume Grelier, Petr H\'{a}jek, Bill Johnson, Florence Lancien, Pavlos Motakis, François N\'etillard, Alexandre Nou, Yo\"el Perreau, Eva Perneck\'a, Colin Petitjean, Tony Proch\'azka, Andr\'es  Quilis,  Mat\'{\i}as Raja, Thomas Schlumprecht, and Quanhua Xu. We are extremely grateful to the four referees for their careful reading of the first submitted version of this book and their encouragements. Their sharp and constructive comments have been very helpful in improving our presentation. 
	
	\smallskip
	
	The first author is grateful for support from the National Science Foundation (NSF) through grants DMS-1800322, DMS-2055604, and DMS-2453662, and from the Department of Mathematics at Texas A\&M University, the College of Arts and Sciences, and the Association of Former Students of Texas A\&M for support for a sabbatical semester spent at Universit\'e Marie et Louis Pasteur during Fall 2024 when part of this book was written. He also wishes to thank the wonderful staff and faculty at the Laboratoire de Math\'ematiques de Besan\c con for the stellar working conditions and a welcoming and stimulating environment provided when spending a sabbatical semester there during Fall 2024, and to the Centre National de la Recherche Scientifique (CNRS) for his support via a Poste Rouge during this period. The second author wishes to thank the Department of Mathematics at Texas A\&M University for the invitation, the warm hospitality and great working conditions during Spring 2019, when we started this project. He also acknowledges partial support from the French ANR projects No. ANR-20-CE40-0006 and ANR-24-CE40-0892-01.
	
	\smallskip
	
	Without the unwavering support of Bill Johnson, Thomas Schlumprecht, and Gilles Pisier, his dear colleagues at Texas A\&M University, the first author would certainly have never been in a position to co-author this book. Thank y'all for everything!
	
	\smallskip The second author wishes to take this opportunity to thank all his former and present PhD students: Pierre Portal, Florent Baudier, Eva Perneck\'a, Aude Dalet, Colin Petitjean, Fran\c cois N\'etillard, Yoël Perreau, Audrey Fovelle, and Estelle Basset. By their enthusiasm, sharp questions, and great ideas, they made him love this subject even more, which lead to this book. Advising them has been a driving force and a privilege. 
	
	\smallskip
	
	And of course, this book would never have come to life without our wonderful and caring families, in particular our spouses ... especially our spouses ... and their endless patience dealing with us while we were working on this project! Their unquantifiable support kept us motivated and moving forward in ways that they probably can't imagine. Florence and Maria, you can rightfully claim ownership of this book. 
	
	\smallskip
	
	We have both been blessed to collaborate with Nigel Kalton during our careers. Working with him, benefiting from his generosity and humility, allowing us to grasp a bit of his mathematical genius, and becoming friends was an incredible gift. His enlightening vision of the nonlinear geometry of Banach spaces has inspired so much of almost every chapter that Nigel has been a companion throughout the preparation of this volume. This book is dedicated to his living memory.

	%%%%%%%%%%%%%%%%%%%%%%%%%%%%%%%%%%%%%%%%%%%%%%%%%%%%%%%%%%%%%%%%%%%%%%%%%%%%%%%%%%%%%%%%%%%%%%%%%%%%%%%%%%%%%%%%%%%%%%%%%%%%%%%%%%%%

	\chapter*{Description chapter by chapter}
	
	We now give a more precise description of the contents of this book. 
	
	\medskip
	
	In Chapter \ref{chapter:basics}, we briefly introduce some standard notation and definitions for Banach spaces and metric spaces. We also define the various notions of nonlinear embeddings and equivalences that will be studied in this book. Namely the bi-Lipschitz, uniform and coarse-Lipschitz embeddings and equivalences. We include an elementary proposition on the many characterizations of coarse-Lipschitz equivalences between Banach spaces. 
	
	In Chapter \ref{chapter:universal}, we review the classical universal spaces for the class of separable metric spaces. We start with the universality of $\ell_\infty$ for isometric embeddings and of the space $C([0,1])$ for isometric linear embeddings and separable Banach spaces. We conclude this chapter with a proof of Aharoni's embedding theorem: every separable metric space bi-Lipschitzly embeds into $\co$. In this book, we take the opportunity to present the proof, taken from \cite{KaltonLancien2008}, giving the optimal distortion in Aharoni's embedding theorem, a question left open in \cite{BenyaminiLindenstrauss2000}.
	
	Chapter \ref{chapter:linear-reductions} is devoted to the classical linearization of Lipschitz maps between Banach spaces. We refer the reader to the book \cite{BenyaminiLindenstrauss2000} for a much more thorough exposition of these tools. We have chosen to take what we believed to be the shortest path to the bi-Lipschitz classification results we wanted to present. To start, we have taken, as the definition of the Radon-Nikod\'{y}m property (RNP in short) of a Banach space $X$, one of its many characterizations: the validity of the Lebesgue differentiability theorem for Lipschitz maps from $\bR$ to $X$. This is a very arguable choice, but it is enough together with the construction of the Bochner integral, recalled in Appendix \ref{appendix:Bochner}, to get nontrivial examples of RNP spaces such as separable duals and reflexive Banach spaces and to deduce the finite-dimensional and infinite-dimensional Rademacher theorems on the Gateaux or weak$^*$ Gateaux differentiability of Lipschitz maps. With these results in hand, we can already prove the stability under Lipschitz equivalences of RNP, reflexivity and the class of spaces isomorphic to $\ell_2$. Then, we present the classical combination of Lipschitz retracts with invariant means techniques to obtain complementation results for Lipschitzly equivalent Banach spaces. This allows us to present the solution of the Lipschitz classification of $\ell_p$ and $L_p$  for $p\in (1,\infty)$. Then, we combine the Rademacher theorems with ultrapower techniques to prove Ribe's rigidity theorem. We conclude this chapter with the application, taken from \cite{BJLPS1999}, of a different kind of differentiablity argument, by showing the stability of Asplundness under Lipschitz quotients. 
	
	In Chapter \ref{chapter:Lipschitz-free}, we present a few applications of the study of the linear properties of Lipschitz-free spaces to the nonlinear classification of Banach spaces. We limit ourselves to the consequences of the seminal works by G. Godefroy and N. Kalton in \cite{GodefroyKalton2003} and \cite{Kalton2004} from the early 2000s. This theory has since grown tremendously and now deserves a book of its own. The Lipschitz-free space $\cF(M)$ of a metric space $M$ (also known as Arens-Eells or transportation cost space) is a canonical predual of the space of real-valued Lipschitz maps on $M$ with the fundamental property that a Lipschitz map defined on $M$ with values in a Banach space can be extended to $\cF(M)$ into a bounded linear map. We first prove the Godefroy-Kalton lifting theorem and its stunning application to the isometric rigidity of separable Banach spaces. Then, we show how they use Lipschitz-free spaces to prove that the bounded approximation property is stable under Lipschitz equivalences. This chapter also contains examples of separable Lipschitz-free spaces that are universal for separable metric spaces and for coarse or uniform (but not bi-Lipschitz) embeddings. Finally, we explain how Kalton used Lipschitz-free techniques to prove that any stable metric space coarsely and uniformly embeds into a reflexive Banach space. 
	
	The main objective of Chapter \ref{chapter:asymptotic-moduli} is to present the notions of asymptotic uniform convexity (AUC) and asymptotic uniform smoothness (AUS) of a Banach space and their quantitative versions: $q$-asymptotic uniform convexity for $q\ge 1$ ($q$-AUC),  $p$-asymptotic uniform smoothness for $p>1$ ($p$-AUS) and asymptotic uniform flatness (AUF). We introduce the corresponding moduli and detail the duality between the asymptotic uniform smoothness of a Banach space and the weak$^*$ asymptotic uniform convexity of its dual. This chapter contains a few classical renorming results: a separable Banach space admits an equivalent AUF norm if and only if it is isomorphic to a subspace of $\co$ (\cite{GKL2000}); for $p\in (1,\infty)$ a separable reflexive Banach space admits an equivalent $p$-AUS norm and an equivalent $p$-AUC norm if and only if it is isomorphic to a subspace of an $\ell_p$-sum of finite-dimensional spaces (\cite{JLPS2002}). We conclude with a description of the basic features of the property $(\beta)$ of Rolewicz.  
	
	In Chapter \ref{chapter:Szlenk}, we study in detail the renorming theory for the asymptotic properties defined in the previous chapter and provide isomorphic characterizations for when such renormings exist, especially for AUS and $p$-AUS renormings. We believe that this has not yet been exposed in any textbook and this explains the length of this chapter. This chapter is exclusively about the linear theory of Banach spaces. One characterization of AUS renormability of a Banach space $X$ can be read on its dual $X^*$ and is expressed in terms of the Szlenk index of $X$, an ordinal index introduced by W. Szlenk \cite{Szlenk1968} in the late 1960s. We provide a self-contained presentation of this fundamental object and its basic properties. Another characterization is obtained in terms of the behavior of weakly null trees in $X$. The use of weakly null trees in this context goes back to the works of Godefroy, Kalton and Lancien \cite{GKL2001} and Odell and Schlumprecht \cite{OdellSchlumprecht2002}. We present the most recent and precise results due to R.M. Causey (\cite{Causey2018c}, \cite{Causey3.5}), as they will be important to obtain some very precise nonlinear invariants, especially in Chapter \ref{chapter:Gorelik}. In particular, we describe the classes $\sT_p$, $\sA_p$ and $\sN_p$, introduced by Causey, which are defined in terms of various $p$-upper estimates for weakly null trees and their links with $p$-AUS renormability. 
	
	In the short Chapter \ref{chapter:AMP_I}, we present the approximate midpoint principle. The Mazur-Ulam theorem states that an onto isometry between normed spaces must be affine. One key ingredient is that isometries preserve the midpoint set of a pair of points. The notion of approximate midpoint set, inspired by this idea, was first introduced by P. Enflo in an unpublished paper where he showed that $L_1$ and $\ell_1$ are not uniformly homeomorphic. We introduce this important tool in this chapter. In the spirit of this book, we give a lower estimate for the Kuratowski measure of noncompactness of approximate midpoint sets in AUS spaces. On the other hand, we provide an upper estimate in AUC spaces. As an application, we show that a Banach space that is too asymptotically uniformly smooth does not admit a coarse-Lipschitz embedding into a Banach space that is too asymptotically uniformly convex. These ideas were already displayed in the book \cite{BenyaminiLindenstrauss2000} and the version we present is taken from N.L. Randrianarivony's thesis \cite{RandrianarivonyThesis}. The use of approximate midpoints has been pushed later to a much more sophisticated level by N. Kalton \cite{Kalton2013b}. This will be detailed in the important Chapter \ref{chapter:AMP_II}.
	
	Chapter \ref{chapter:Gorelik} concentrates on the Gorelik principle and its applications. The Gorelik principle was formulated by Johnson, Lindenstrauss and Schechtman in \cite{JLS1996} and used in their proof of the uniform rigidity of $\ell_p$ for $p\in(1,\infty)$. It was named after E. Gorelik who devised it to prove that for $p\in (2,\infty)$, $\ell_p$ is not uniformly homeomorphic to $L_p$. We present the original proof from \cite{JLS1996} of the uniform (actually even coarse-Lipschitz) rigidity of $\ell_p$, for $p\in (1,\infty)$. For $p\in (1,2)$, it relies on the approximate midpoint principle and for $p\in (2,\infty)$ on the Gorelik principle. It is important to mention that for $p\in (2,\infty)$, one can instead use the tools from Chapters \ref{chapter:Johnson} or \ref{chapter:Hamming}. For this reason, we focus in this chapter on the results where the Gorelik principle is still needed. We start with a version of the Gorelik principle for Lipschitz equivalences. We apply it to show that the class of $p$-AUS renormable Banach spaces is stable under Lipschitz equivalences and that a space Lipschitzly equivalent to $\co$ must be linearly isomorphic to $\co$ \cite{GKL2001}. This was already written in the textbook \cite{AlbiacKalton2016}. Then, we devise a version of the Gorelik principle for coarse-Lipschitz equivalences. We apply it to deduce the coarse-Lipschitz rigidity of AUS renormings. We actually prove much more precise results on the stability under coarse-Lipschitz equivalences of classes of spaces satisfying various upper-$p$-estimates for weakly null trees: the class of asymptotic-$\co$ spaces and the classes $\textsf{A}_p$ and $\textsf{N}_p$ introduced by R.M. Causey to provide a fine hierarchy of AUS renormable Banach spaces. 
	
	In Chapter \ref{chapter:Johnson}, we make our first encounter with the beautiful connection between the geometry of sequences of countably infinite graphs, here the Johnson graphs and the asymptotic geometry of Banach spaces. This different approach is better suited than the Gorelik principle for addressing questions on the stability of asymptotic uniform smoothness under bi-Lipschitz or coarse-Lipschitz \emph{embeddings}. A pivotal tool is the study of Ramsey-type concentration inequalities for Banach-valued Lipschitz maps defined on the Johnson graphs initiated by Kalton and Randrianarivony in \cite{KaltonRandrianarivony2008}. The stability of these concentration inequalities under nonlinear embeddings has many important applications: the coarse-Lipschitz rigidity of the class of subspaces of $\ell_p$ for $p\in(2,\infty)$, a coarse version of Tsirelson's theorem and the optimality of snowflake and compression exponents between $\ell_p$ spaces. In particular, Hilbert space does not coarsely embed into Tsirelson's original space $\Tsi^*$, making $\Tsi^*$ the first example of an infinite-dimensional Banach space with this property.
	A recurring theme of this chapter is the connection between the Johnson graphs and the theory of spreading models recalled in Appendix \ref{appendix:asymptotic}. 
	
	In Chapter \ref{chapter:Hamming}, we investigate the geometry of the Hamming graphs, a countably branching version of the Hamming cubes and its connection with the theory of asymptotic models of Halbeisen and Odell recalled in Appendix \ref{appendix:asymptotic}. All the results from Chapter \ref{chapter:Johnson} remain true if instead of the Johnson graphs, we consider the Hamming graphs. The main result of this chapter is a purely metric characterization of the class of reflexive and asymptotic-$\co$ Banach spaces in terms of a Ramsey-type concentration inequality for Lipschitz maps on the Hamming graphs. This result is a prototypical result in the Kalton program and is up to this day, the only metric characterization of this type. A consequence of this metric characterization is that this class is stable under coarse embeddings. We also provide a metric characterization in terms of metric space preclusion of the class of asymptotic-$\co$ spaces \emph{within the class of separable and reflexive Banach spaces}. This chapter also contains a fundamental result about the preservation of upper estimates for asymptotic models under coarse-Lipschitz embeddings. This extension of a similar result obtained by Kalton and Randrianarivony in \cite{KaltonRandrianarivony2008} for spreading models has never appeared elsewhere. 
	
	Chapter \ref{chapter:interlaced-graphs} is devoted to the foundational work of Kalton about the coarse universality problem for Banach spaces. It is specifically to solve this longstanding open problem that Kalton initiated the investigation of concentration inequalities for Banach-valued Lipschitz maps defined on certain nonlocally finite graphs. One of the main results of this chapter states that if a Banach space is coarsely universal for separable metric spaces, then one of its iterated duals must be nonseparable. Consequently, a reflexive Banach space cannot be coarsely universal. This chapter contains a detailed analysis of the geometry of the interlaced graphs and its close relationship with the summing basis of $\co$, as well as Kalton's property $Q$ and its recent reformulation in terms of upper stability (a natural bi-Lipschitz generalization of Krivine-Maurey isometric notion of stability). In particular, this chapter contains a purely metric characterization, within the class of Banach spaces with the alternating Banach-Saks property (or equivalently without $\ell_1$-spreading models), of the class of reflexive Banach spaces in terms of property $Q$. This metric characterization could arguably be considered as the birth certificate of the Kalton program. This chapter also contains a brief study of an uncountable version of property $Q$ and its applications to coarse embeddability of nonseparable spaces into $\ell_\infty$.
	
	In Chapter \ref{chapter:AMP_II}, we revisit the approximate midpoint principle. It is based on a paper by N. Kalton \cite{Kalton2013b} that was edited and published after its author suddenly passed away in August 2010. It contains the most sophisticated results on the preservation of asymptotic uniform convexity under coarse-Lipschitz embeddings, as well as striking applications. We try in this chapter to advertise and make accessible some of the most important of these results, which we believe, have been overlooked for too long. The first general result is the preservation under coarse-Lipschitz embeddings of a natural quantitative lower bound on the fundamental function of spreading models generated by weakly null sequences in relation to the modulus of asymptotic midpoint uniform convexity. This chapter culminates with the proof of the stability under coarse-Lipschitz embeddings of the following classes: subspaces, quotients and subspaces of quotients of $\ell_p$, for $p \in (1,\infty)$. We also include the proof of the uniqueness of the coarse structure of $(\sum_{n=1}^\infty \ell_r^n)_{\ell_p}$ for $1<p<r\le 2$ or $2\le r<p<\infty$. Let us also mention that this chapter contains some original linear results, interesting on their own, obtained in \cite{Kalton2013b} for these nonlinear purposes. 
	
	Chapter \ref{chapter:Counterexamples} is devoted to the presentation of the main counterexamples in the theory of the nonlinear classification of Banach spaces. We start with what we call the elementary lifting method. This construction of nonlinear equivalences is based on the existence of linear quotient maps that admit a nice (Lipschitz or uniformly continuous) lifting. The first example, due to Aharoni and Lindenstrauss \cite{AharoniLindenstrauss1978}, of a pair of nonseparable Banach spaces that are Lipschitz equivalent, but not linearly isomorphic, is based on this method. Then, we show how this method can be combined with the results on Lipschitz-free spaces from Chapter \ref{chapter:Lipschitz-free} to produce an example of a pair of separable Banach spaces that are uniformly homeomorphic but not linearly isomorphic. Another method, inspired by the first example of such a pair built by M. Ribe in \cite{Ribe1984}, has been formalized and generalized by N. Kalton in \cite{Kalton2012} and \cite{Kalton2013}. We describe this machinery and its most striking application: the class of separable Banach spaces with an equivalent $p$-AUS norm is not stable under uniform homeomorphisms. We conclude this chapter with the description of a stunning example, again due to Kalton \cite{Kalton2012}, of a pair of coarsely equivalent but not uniformly equivalent separable Banach spaces. 
	
	In Chapter \ref{chapter:trees}, we investigate the geometry of countably branching trees and their applications to various metric characterizations and nonlinear rigidity problems. The original application is the metric characterization, within the class of reflexive spaces, of the class of Banach spaces that are asymptotically uniformly convexifiable and asymptotically uniformly smoothable in terms of the Lipschitz geometry of the countably branching trees. To some extent, this metric characterization, which is an asymptotic analog of Bourgain's metric characterization of superreflexivity, which kick-started the Ribe program, motivated the systematic investigation of the metric reformulation of asymptotic properties and was often referred to as the initial step in an asymptotic variation of the Ribe program. The second application is a metric characterization of Banach spaces admitting an equivalent norm with  property $(\beta)$. The last application is the stability under coarse or uniform quotients of property $(\beta)$. The chapter starts with some elementary bi-Lipschitz embedding results into the Banach spaces $\ell_1$ and $\co$ that are further refined and extended. In particular, the connection between coarse embeddings of trees into Banach spaces and the spreading models or asymptotic models of the latter is investigated and a purely metric characterization of infinite dimensionality of a Banach space is expressed in terms of the coarse geometry of countably branching trees. The natural relationship between the Szlenk index of Banach spaces and equi-bi-Lipschitz embeddings of trees is explained. Then, we show how distortion lower bounds and compression function upper bounds can be derived from a Poincar\'e-type inequality, called beta $p$-convexity, for Lipschitz maps on the countably branching trees taking values in Banach spaces with property $(\beta_p)$. The notion of beta $p$-convexity, which is stable under certain nonlinear quotients, has not explicitly appeared in the work of Baudier and Gartland in \cite{BaudierGartland24} but only in preliminary unpublished work and we thank C. Gartland for allowing us to present this invariant in this book. All these results are finally combined together to prove the three applications mentioned above. 
	
	The last chapter, Chapter \ref{chapter:diamonds}, is dedicated to proving a metric characterization, within the class of reflexive spaces with an unconditional asymptotic structure, of the class of asymptotically uniformly convexifiable spaces. This characterization is reminiscent of the Johnson-Schechtman metric characterization of the class of uniformly convexifiable spaces in terms of the Lipschitz geometry of the $2$-branching diamond graphs. The most delicate aspect here revolves around the embeddability properties of the sequence of countably branching diamond graphs. Therefore, we first carefully analyze the graph metric on the diamond graphs and prove a bi-Lipschitz embedding result into $L_1$. The non bi-Lipschitz embeddability of the countably branching diamond graphs into asymptotically midpoint uniformly convex spaces is taken care of via a self-improvement argument. Then, we proceed to construct bi-Lipschitz embeddings with distortion at most $3$, of the countably branching diamond graphs into $\co$. These embeddings have finite support and can be used to obtain sharp upper bounds for the $\ell_p$-distortion of the countably branching diamond graphs. Finally, we explain how this particular $\co$-embedding can be implemented in reflexive Banach spaces with an asymptotic unconditional structure and that are not asymptotically uniformly convexifiable. The Szlenk index plays, yet again, a prominent role in this final chapter.
	
	We conclude with a few appendices containing some background material that is necessary for this book. We have postponed the presentation of these subjects either because they are very classical, or in an attempt to avoid perturbing the flow of the main text. They include the rudiments of Banach space theory, basics of the theory of the Bochner integral, infinite Ramsey theory and spreading models. They also contain some less classical material on asymptotic models, joint asymptotic models and asymptotic structure in Banach spaces. Only a few proofs are given but pointers to the relevant literature are provided.

	%%%%%%%%%%%%%%%%%%%%%%%%%%%%%%%%%%%%%%%%%%%%%%%%%%%%%%%%%%%%%%%%%%%%%%%%%%%%%%%%%%%%%%%%%%%%%%%%%%%%%
	
	\chapter[Basic notation and definitions]{Basic notation and definitions}
	\label{chapter:basics}
	
	\section{Banach spaces}
	
	In this book, we will work with real Banach spaces. They will mostly be denoted by the letters $X,Y,Z$. For a Banach space $(X,\norm{\cdot}_X)$, we denote by $B_X$ its closed unit ball and by $S_X$ its unit sphere. Thus, $B_X=\{x\in X \colon \norm{x}_X \le 1\}$  and $S_X=\{x\in X \colon \norm{x}_X = 1\}$. If $X$ and $Y$ are Banach spaces, $B(X,Y)$ denotes the space of all bounded linear maps (or operators) from $X$ to $Y$ endowed with its natural norm:
	$$\forall T \in B(X,Y),\ \ \norm{T}_{B(X,Y)}=\sup_{x\in B_X}\norm{Tx}_Y.$$ 
	In the particular case when $Y=\R$, $B(X,\R)$ is the topological dual of $X$, denoted by $X^*$ and for its elements, we will most of the time use the notation $x^*,y^*,z^*$. For $x^* \in X^*$, $\norm{x^*}_{X^*}=\sup_{x \in B_X}|x^*(x)|$. For $x\in X$ and $x^* \in X^*$, we will use interchangeably $x^*(x)$, $\langle x,x^*\rangle$, or $\langle x^*,x\rangle$ for the action of an element $x^* \in X^*$ on an element $x\in X$. In this spirit, we will regard $X$ as a subspace of its bidual $X^{**}$ through its canonical embedding. 
	
	An invertible and surjective operator between Banach spaces will be called a \emph{linear isomorphism} or just an \emph{isomorphism} if the context is clear. If there is a linear isomorphism between the Banach spaces $X$ and $Y$, we will naturally say that $X$ and $Y$ are \emph{linearly isomorphic} (or just isomorphic) and we write $X \simeq Y$. A linear isometry from $X$ into $Y$ is a linear map $T\colon X \to Y$ such that $\norm{Tx}_Y=\norm{x}_X$ for all $x \in X$. If moreover, this map is onto, we say that $X$ and $Y$ are isometric and write $X\equiv Y$. 
	
	If $X$ and $Y$ are two isomorphic Banach spaces, $d_{BM}(X,Y)$ denotes the Banach-Mazur distance between these spaces, defined to be the infimum of $\norm{T}\,\|T^{-1}\|$, as T ranges over all isomorphisms from $X$ onto $Y$. If $X$ and $Y$ are not isomorphic, we use the convention $d_{BM}(X,Y)=\infty$. 
	
	We call \emph{bounded projection}, or just projection, any idempotent operator. A subspace $Z$ of a Banach space $X$ is said to be \emph{complemented} if there is a projection from $X$ onto $Z$. This is the case if and only if $Z$ is closed and there is a
	closed subspace $W$ of $X$ so that $W \cap Z = \{0\}$ and $X = W + Z$; we then write $X = W \oplus Z$ and say that X is the \emph{direct sum} of $W$ and $Z$. 
	
	We now recall the notation for the most classical Banach spaces.
	\begin{itemize}
		\item If $K$ is a compact Hausdorff space, then $C(K)$ is the space of all real-valued continuous functions on $K$ equipped with the supremum norm, $\norm{f}_\infty=\sup_{x\in K}\abs{f(x)}$.
		\medskip
		\item[] Let $(\Omega,\Sigma,\mu)$ be a measure space.
		\item For $p\in [1,\infty)$, $L_p(\Omega,\Sigma,\mu)$, or just $L_p(\mu)$, is the space of all real-valued measurable functions such that $\int_\Omega |f|^p\,d\mu <\infty$, normed with $\norm{f}_p=(\int_\Omega |f|^p\,d\mu)^{1/p}$.
		\item $L_\infty(\Omega,\Sigma,\mu)$, or just $L_\infty(\mu)$, is the space of all real-valued measurable functions that are essentially bounded, normed by the essential supremum of $|f|$, i.e.    $\norm{f}_\infty=\inf \{C\ge 0 \colon \mu\{\abs{f}>C\}\}=0$.
		\item For $p\in [1,\infty]$, $\ell_p(\bN)$, or just $\ell_p$, is the space $L_p(\mu)$, when $\mu$ is the counting measure on $\bN$.
		\item Finally, we denote by $c_0(\bN)$, or just $\co$, the space of all real-valued sequences that are converging to $0$, equipped with the supremum norm. 
	\end{itemize}
	Assume that $I$ is a set and $\cU$ is a nonprincipal ultrafilter on $I$; assume
	also that for all $i\in I$, $X_i$ is a Banach space. Denote by $(\sum_{i\in I} X_i)_{\ell_\infty}$ the space of all families $x=(x_i)_{i\in I}$ such that $\norm{x}_\infty=\sup_{i\in I}\norm{x_i}_{X_i}<\infty$. For $x=(x_i)_{i\in I} \in (\sum_{i\in I} X_i)_{\ell_\infty}$, we define the seminorm $p(x)=\lim_{\cU}\norm{x_i}_{X_i}$. The quotient of $(\sum_{i\in I} X_i)_{\ell_\infty}$ by the closed subspace of all $x$ with $p(x)=0$ equipped with its natural norm is a Banach space called the \emph{ultraproduct} of the $X_i$ with respect to $\cU$ and is denoted by $(\prod_{i \in I} X_i)^{\cU}$. When all the $X_i$ are the same space $X$ we call the space thus obtained the \emph{ultrapower} of $X$ with respect to $\cU$ and write $X^{\cU}$.

	\section{Lipschitz, coarse-Lipschitz and coarse maps between metric spaces}
	\label{sec:nonlinear-maps}
	
	For metric spaces, we will mostly use the letters $M$ and $N$ to denote them and use $(M,d)$ or $(M,d_M)$ to indicate the associated distance. Let $(M,d)$ be a metric space. Then, for all $x\in M$ and all $r\ge 0$, $B_M(x,r)=\{y\in M,\ d(x,y)\le r\}$ denotes the closed ball of center $x$ and radius $r$. 
	
	\begin{defi} 
		Let $(M,d_M)$ and $(N,d_N)$ be two metric spaces. For a given map $f\colon M \to N$, one defines
		\begin{equation*}
			\rho_f(t)=\inf\big\{d_N(f(x),f(y)) : d_M(x,y)\ge  t\big\},
		\end{equation*}
		and
		\begin{equation*}
			\omega_f(t)={\rm sup}\{d_N(f(x),f(y)) : d_M(x,y)\le  t\}.
		\end{equation*}
		We adopt the convention $\sup(\emptyset)=0$ and $\inf(\emptyset)=+\infty$. The moduli $\rho_f$ and $\omega_f$ will be called the \emph{compression modulus} and the \emph{expansion modulus} of the map $f$, respectively.
	\end{defi}
	
	Note that a map $f\colon M \to N$ is uniformly continuous if and only if $\lim_{t \to 0} \omega_f(t)=0$.
	
	\begin{defi}
		\label{def:CLmaps}
		Let $(M,d_M)$ and $(N,d_N)$ be two metric spaces and $f\colon M \to N$.
		\begin{enumerate}[(a)]
			\item We say that $f$ is \emph{Lipschitz} if there exists $C\ge 0$ such that $\omega_f(t)\le Ct$ for all $t\ge 0$. In that case, the Lipschitz constant of $f$ is 
			\begin{equation*}
				\Lip(f)=\sup\Big\{\frac{d_N(f(x),f(y))}{d_M(x,y)} \colon x\neq y \ \text{in}\ M\Big\}.
			\end{equation*}
			\item We say that $f$ is \emph{coarse-Lipschitz} if there exist $A$ and $B$ in $[0,\infty)$ such that $\omega_f(t)\le At+B$ for all $t\ge 0$.
			\item We say that $f$ is \emph{coarse} if  $\omega_f(t)<\infty$ for all $t\ge 0$.
		\end{enumerate}
	\end{defi}
	
	%\begin{rema} We did not include it as a definition, but note that $f$ is uniformly continuous if and only if $\lim_{t \to 0} \omega_f(t)=0$.
	%\end{rema}
	
	The following elementary result is known as the Corson-Klee Lemma.  
	\begin{prop}
		\label{prop:CorsonKlee}
		Let $X$ be a Banach space, $(N,d)$ a metric space and $f\colon X \to N$. Assume that $\omega_f(t_0)<\infty$, for some $t_0>0$. Then, $f$ is coarse-Lipschitz. In particular, if $f$ is coarse or uniformly continuous, then it is coarse-Lipschitz. 
	\end{prop}
	
	\begin{proof}
		Consider $x\neq y \in X$. Then, the segment $[x,y]$ can be divided into $k$ segments of length at most $t_0$ with $k\le \frac{1}{t_0}\norm{x-y}_X+1$. It follows that
		$$d(f(x),f(y))\le k\omega_f(t_0) \le \frac{\omega_f(t_0)}{t_0}\norm{x-y}_X+\omega_f(t_0).$$
		This finishes the proof.
	\end{proof}
	
	\begin{rema} 
		The above argument can be extended to $f\colon M \to N$, where $M$ is a \emph{metrically convex} metric space. This means that for any $x\neq y \in M$ and any $\lambda \in [0,1]$, there exists $z \in M$ such that $d_M(x,z)=\lambda d_M(x,y)$ and $d_M(z,y)=(1-\lambda)d_M(x,y)$. 
	\end{rema}
	
	We will also need to introduce the following quantities. Let $(M,d_M)$ and $(N,d_N)$ be two metric
	spaces and let $f\colon M\to N$ be a map. If $(M,d_M)$ is unbounded, we define for all $s>0$, 
	\begin{equation*}
		\Lip_s(f)=\sup\Big\{\frac{d_N(f(x),f(y))}{d_M(x,y)} \colon d_M(x,y)\ge s\Big\}
	\end{equation*}
	and
	\begin{equation*}
		\Lip_\infty(f)=\inf_{s>0}\Lip_s(f).
	\end{equation*}
	
	The following characterizations of coarse-Lipschitz maps, which we state for Banach spaces, are actually valid for metrically convex metric spaces and easy to verify.
	
	\begin{prop}
		\label{prop:CL}
		Let $X$ and $Y$ be two Banach spaces and $f\colon X\to Y$ be a mapping. The following assertions are equivalent.
		\begin{enumerate}[(i)]
			\item The map $f$ is coarse-Lipschitz.
			\item There exist $A$ and $\theta$ in $[0,+\infty)$ such that for all $x,y\in X$,
			\begin{equation*}
				\norm{x-y}_X\ge \theta \Rightarrow \norm{f(x)-f(y)}_Y\le A\norm{x-y}_X.
			\end{equation*}
			\item The map $f$ is such that $\Lip_\infty(f)<\infty$.
		\end{enumerate}
	\end{prop}
	
	Note that in the above statement, $\Lip_\infty(f)$ coincides with the infimum of all $A\ge 0$ such that $(ii)$ is satisfied for some $\theta\ge 0$ and also with the infimum of all $A\ge 0$ such that Definition \ref{def:CLmaps} (b) is satisfied for some $B\ge 0$.

	\section{Bi-Lipschitz, coarse-Lipschitz and coarse embeddings}
	\label{sec:nonlinear-embeddings}
	
	In this section, we introduce various types of nonlinear embeddings that play a central role in this book. We start with the classical notion of bi-Lipschitz embeddability and its quantification, namely, the bi-Lipschitz distortion. 
	
	\begin{defi}
		\label{def:Lip-dist}
		If $(M,d_M)$ and $(N,d_N)$ are two metric spaces, the \emph{$N$-distortion} of $M$, denoted by $c_N(M)$, is defined as the infimum of those $D\in[1,\infty)$ for which there exist $s\in(0,\infty)$ and a map $f\colon M\to N$ so that for all $x,y\in M$,
		\begin{equation}
			\label{eq:distortion}
			s\,d_{M}(x,y)\le  d_{N}\big(f(x),f(y)\big)\le 
			sD\,d_{M}(x,y).
		\end{equation}
		When \eqref{eq:distortion} holds we say that $M$ bi-Lipschitzly embeds into $N$ with distortion $D$.
	\end{defi}
	
	The basic requirement of a bi-Lipschitz embedding is that the geometry of the embedded space is preserved at \emph{all scales} when representing it in the host space.
	
	\medskip When families of metric spaces are involved, we are usually concerned with families of embeddings with some universal control and we now introduce some convenient terminology and notation.
	
	\begin{defi}
		\label{def:rho-omega}
		Let $(M,d_M)$ and $(N,d_N)$ be metric spaces. Let $\rho,\omega\colon [0,\infty)\to \R \cup \{\infty\}$. We say that $M$ $(\rho,\omega)$-embeds into $N$ if there exists $f\colon M\to N$ such that for all $x,y\in M$ we have
		\begin{equation*}
			\rho(d_M(x,y))\le d_N(f(x),f(y))\le \omega(d_M(x,y)).
		\end{equation*}
		Let $\{M_i\}_{i\in I}$ be a collection of metric spaces. We say that $\{M_i\}_{i\in I}$ $(\rho,\omega)$-embeds into $N$ if for every $i\in I$, $M_i$ $(\rho,\omega)$-embeds into $N$.
	\end{defi}
	
	Note that for any map $f\colon M \to N$, $M$ $(\rho_f,\omega_f)$-embeds into $N$ and $(\rho_f,\omega_f)$ can be viewed as the best pair of functions such that this is true.
	
	Using the terminology from Definition \ref{def:rho-omega}, a family $\{M_i\}_{i\in I}$ of metric spaces \emph{equi-uniformly embeds} into $N$ if there exist nondecreasing functions $\rho,\omega\colon [0,\infty)\to[0,\infty)$ such that $\lim_{t\to 0}\omega(t)=0$, $\rho(t)>0$ for all $t>0$ and  $\{M_i\}_{i\in I}$ $(\rho,\omega)$-embeds into $N$. In the special case where $\{M_i\}_{i\in I}$ reduces to the singleton $\{M\}$, the notion boils down to the classical notion of uniform embedding, i.e. $M$ uniformly embeds into $N$ if there is a map from $M$ into $N$ that is injective, uniformly continuous and with a uniformly continuous inverse. 
	
	\begin{defi}
		\label{def:equi-Lip} 
		Let $\{M_i\}_{i\in I}$ be a collection of metric spaces and $N$ be a metric space.
		We say that $\{M_i\}_{i\in I}$ admits an \emph{equi-bi-Lipschitz embedding} into $N$ if there are nonzero linear maps $\rho,\omega\colon [0,\infty)\to [0,\infty)$ such that $\{M_i\}_{i\in I}$ $(\rho,\omega)$-embeds into $N$.  
		%		$$\{M_i\}_{i\in I}   \buildrel {E-L}\over {\hookrightarrow} N\ \ \text{and}\ \ M \buildrel {L}\over {\hookrightarrow} N,\  \text{respectively.}$$
	\end{defi}
	
	Note that if $\{M_i\}_{i\in I}=\{M\}$, then we recover the notion of bi-Lipschitz embeddability from Definition \ref{def:Lip-dist}. 
	
	\begin{rema}
		Equi-bi-Lipschitz embeddability is a stronger condition than merely assuming that $\sup_{i\in I}c_N(M_i)<\infty$ since it does not allow for arbitrarily large or arbitrarily small scaling factors in \eqref{eq:distortion}. However, if $N$ is a Banach space, rescaling is possible and the two notions coincide.    
	\end{rema}
	
	Bi-Lipschitz embeddability is a rather stringent property that can be relaxed in various ways. One example is uniform embeddability, which only takes into consideration the small-scale geometry. The Lipschitz and uniform geometry of Banach spaces has been extensively studied in \cite{BenyaminiLindenstrauss2000}, and around the time of publication of \cite{BenyaminiLindenstrauss2000}, other relaxations of bi-Lipschitz embeddability gained traction thanks to fundamental applications in geometric group theory, topology, or noncommutative geometry. These applications will not be discussed here, and we refer the interested reader to \cite{NowakYu}.
	
	The notion of coarse-Lipschitz embeddability is another relaxation that erases the small-scale geometry but preserves in a Lipschitz way the large-scale geometry. This notion is the same as the notion of quasi-isometric embeddability introduced by Gromov in \cite{Gromov1987} (see also the book \cite{GhysDelaHarpe1990} by Ghys and de la Harpe).
	
	\begin{defi}
		\label{def:CLembeddings} 
		Let $\{M_i\}_{i\in I}$ be a collection of metric spaces and $N$ be a metric space.
		We say that $\{M_i\}_{i\in I}$ admits an \emph{equi-coarse-Lipschitz embedding} into $N$ if  $\{M_i\}_{i\in I}$ $(\rho,\omega)$-embeds into $N$ for some (strictly) increasing and affine maps $\rho,\omega \colon [0,\infty)\to [0,\infty)$. 
		When $\{M_i\}_{i\in I}=\{M\}$, we simply say that $M$ admits a \emph{coarse-Lipschitz embedding} into $N$ and in this case there are constants $A>0$, $B\ge 0$ and a map $f\colon M\to N$ such that for all $x,y\in M$ we have
		\begin{equation*}
			\frac{1}{A}d_M(x,y) - B\le d_N(f(x),f(y))\le Ad_M(x,y) + B.
		\end{equation*}
		%We denote		$$\{M_i\}_{i\in I}   \buildrel {E-CL}\over {\hookrightarrow} N\ \ \text{and}\ \ M \buildrel {CL}\over {\hookrightarrow} N,\  \text{respectively.}$$	
	\end{defi}
	
	If one is willing to pass to ultrapowers, we can eliminate the additive noise. Indeed, a coarse-Lipschitz embedding between Banach spaces induces a bi-Lipschitz embedding between their ultrapowers, as we prove in the next proposition.
	
	\begin{prop}
		\label{prop:CL-ultrapower-Lip}
		Let $X$ and $Y$ be two Banach spaces and $\cU$ be a nonprincipal ultrafilter on $\bN$. If $X$ admits a coarse-Lipschitz embedding into $Y$, then $X^{\cU}$ admits a bi-Lipschitz embedding into $Y^{\cU}$.
	\end{prop}
	
	\begin{proof} 
		Let $f\colon X \to Y$ be a coarse-Lipschitz embedding. Define $f^{\cU}\colon X^{\cU} \to Y^{\cU}$ as follows. For a bounded sequence $(x_n)_{n=1}^\infty$ in $X$, denote by $x$ its equivalence class in $X^{\cU}$ and define $f^{\cU}(x) := y$, where $y$ is the equivalence class of $(\frac1nf(nx_n))_{n=1}^\infty$ in $Y^{\cU}$. Checking that $f^{\cU}$ is a bi-Lipschitz embedding is then immediate.
	\end{proof}
	
	If one is content with a coarser control of the large-scale structure, one arrives at the notion of coarse embeddability.
	
	\begin{defi}
		\label{def:coarse-embeddings} 
		Let $\{M_i\}_{i\in I}$ be a collection of metric spaces and $N$ be a metric space.
		We say that $\{M_i\}_{i\in I}$ \emph{equi-coarsely embeds} into $N$ if  $\{M_i\}_{i\in I}$ $(\rho,\omega)$-embeds into $N$ for some nondecreasing maps $\rho,\omega \colon [0,\infty)\to [0,\infty)$ satisfying $\lim_{t\to\infty}\rho(t)=\infty$. 
		When $\{M_i\}_{i\in I}=\{M\}$, we simply say that $M$ \emph{coarsely embeds} into $N$ and in this case there are nondecreasing maps $\rho,\omega \colon [0,\infty)\to [0,\infty)$ satisfying 
		\begin{equation*}
			\lim_{t\to\infty}\rho(t)=\infty,
		\end{equation*}
		and a map $f\colon M\to N$ such that for all $x,y\in M$ we have
		\begin{equation*}
			\rho(d_M(x,y))\le d_N(f(x),f(y))\le \omega(d_M(x,y)).
		\end{equation*}
		%We denote		$$\{M_i\}_{i\in I}   \buildrel {E-C}\over {\hookrightarrow} N\ \ \text{and}\ \ M \buildrel {C}\over {\hookrightarrow} N,\  \text{respectively.}$$
	\end{defi}
	
	\begin{rema}\,
		\begin{enumerate}
			\item If $f\colon M\to N$ is a coarse embedding, then it is a coarse map by definition.
			\item Unlike bi-Lipschitz or uniform embeddings that establish a bijection between the embedded space and its image, coarse-Lipschitz embeddings do not as injectivity is usually lost.
		\end{enumerate} 
	\end{rema}
	
	It is immediate that a bi-Lipschitz embedding is a coarse-Lipschitz embedding, which is itself a coarse embedding.
	It is also clear that coarse-Lipschitz or coarse embeddability are uninteresting notions if the spaces to be embedded have bounded diameters. In fact, every bounded metric space admits a coarse-Lipschitz embedding into any one-point metric space.
	
	It is quite remarkable that we do not know of examples in the Banach space setting that distinguish between coarse embeddability and uniform embeddability! 
	
	\begin{prob}
		\label{pb:coarse-vs-uniform}
		Let $X$ and $Y$ be Banach spaces. Is it true that $X$ uniformly embeds into $Y$ if and only if $X$ coarsely embeds into $Y$?
	\end{prob}
	
	\section{Lipschitz, coarse-Lipschitz, coarse and net equivalences}
	\label{sec:nonlinear-equivalences}
	
	In nonlinear classification problems for Banach spaces, one tries to understand the equivalence class of a Banach space according to some equivalence relation that retains some of its nonlinear geometry. We will introduce essentially three equivalence relations for which the classification problem is interesting. These nonlinear equivalence relations correspond to the nonlinear embeddings introduced in the previous section, and we gather their definitions below. 
	
	\begin{defi}
		\label{def:nonlinear-equivalences} 
		Let $M$ and $N$ be two metric spaces.
		\begin{enumerate}[(a)]
			\item We say that $M$ and $N$ are \emph{Lipschitz equivalent}, or \emph{Lipschitz homeomorphic}, if there exists $f\colon M \to N$ that is Lipschitz, bijective and such that $f^{-1}$ is Lipschitz.
			%		We denote $M \buildrel {L}\over \sim N$.
			\item We say that $M$ and $N$ are \emph{uniformly equivalent}, or \emph{uniformly homeomorphic}, if there exists $f\colon M \to N$ which is uniformly continuous, bijective and such that $f^{-1}$ is uniformly continuous.
			%		We denote $M \buildrel {UH}\over \sim N$.
			\item We say that $M$ and $N$ are \emph{coarse-Lipschitz equivalent} if there exist coarse-Lipschitz maps $f\colon M\to N$ and $g\colon N\to M$ such that
			$$\sup_{x\in M} d_M(g\circ f(x),x)<\infty,\ \text{and}\ \ \ \sup_{y\in N} d_N(f\circ g(y),y)<\infty.$$
			%		We denote $M \buildrel {CL}\over \sim N$.
			\item We say that $M$ and $N$ are \emph{coarsely equivalent} if there exist coarse maps $f\colon M\to N$ and $g\colon N\to M$  such that
			$$\sup_{x\in M} d_M(g\circ f(x),x)<\infty,\ \text{and}\ \ \ \sup_{y\in N} d_N(f\circ g(y),y)<\infty.$$
			%		We denote $M \buildrel {C}\over \sim N$.		
		\end{enumerate}
		
	\end{defi}
	
	\begin{rema} Let $M$ and $N$ be two metric spaces. It is obvious that if $f\colon M \to N$ is a bi-Lipschitz (respectively uniform) embedding, then $f\colon M \to f(M)$ is a Lipschitz (respectively uniform) equivalence. With a bit of work, one can see that the same is true for coarse-Lipschitz and coarse embeddings (see Exercise \ref{ex:embedding-equivalence}). 
	\end{rema}
	
	In the Banach space setting, Proposition \ref{prop:CorsonKlee} has several immediate consequences regarding the inter-relations between the various nonlinear equivalence relations above. The proof of the following elementary observation is left as an exercise.
	
	\begin{prop}
		\label{prop:CL-equiv->Lip-equiv-ultra}
		Let $\cU$ be a nonprincipal ultrafilter on $\bN$. If two Banach spaces $X$ and $Y$ are coarse-Lipschitz equivalent, then $X^{\cU}$ and $Y^{\cU}$ are Lipschitz equivalent.
	\end{prop}
	
	It follows from Proposition \ref{prop:CorsonKlee} that the notions of coarse-Lipschitz and coarse equivalences between Banach spaces are the same. The following statement is also a direct consequence of Proposition \ref{prop:CorsonKlee}. 
	
	\begin{prop}
		\label{prop:UH-equiv->CL-equiv}
		A uniform equivalence between two Banach spaces is also a coarse-Lipschitz equivalence. 
	\end{prop}
	
	The converse of this statement is clearly false. Indeed, if $X$ is Banach space, define $f\colon X \to X$ by $f(x)=x$ if $x \notin B_X$ and $f(x)=0$ if $x\in B_X$. Then, it is easy to check that $f$ is a coarse-Lipschitz equivalence but not a uniform homeomorphism. 
	
	It is also worth mentioning that the coarse or uniform embeddability of a Banach space $X$ into a Banach space $Y$ does not imply its coarse-Lipschitz embeddability. Indeed, as it will appear in the course of this book, $\ell_1$ coarsely and uniformly embeds into $\ell_2$,  but does not coarse-Lipschitz embed into $\ell_2$. 
	
	For a long time, it remained unknown whether two coarse-Lipschitz equivalent Banach spaces are always uniformly equivalent until a remarkable counterexample was discovered by N. Kalton \cite{Kalton2012} (see Section \ref{subsec:CEnotUH} for the proof).
	
	\begin{theo} 
		%There exist separable Banach spaces $X$ and $Y$ that are coarse-Lipschitz equivalent but not uniformly homeomorphic.
		There exists a pair of separable Banach spaces that are coarse-Lipschitz equivalent but not uniformly equivalent.
	\end{theo}

	Let us now concentrate on the notion of coarse-Lipschitz equivalence and relate it to the notion of net equivalence that can be found in \cite{BenyaminiLindenstrauss2000}.
	
	\begin{defi} 
		Let $0<a\le b$. An \emph{$(a,b)$-net} in the metric space $(M,d)$ is a subset $\cM$ of $M$ such that for every $z\neq z'$ in $\cM$,  $d(z,z')\ge a$ and for
		every $x$ in $M$, $d(x,\cM)< b$.
		
		Then, a subset $\cM$ of $M$ is a \emph{net} in $M$ if it is an $(a,b)$-net for some $0<a\le b$.
	\end{defi}
	
	\begin{rema}
		Let $b>0$. If $\cM$ is a maximal $b$-separated subset of $M$, then $\cM$ is a $(b,b)$-net in $M$.     
	\end{rema}
	
	We refer the reader to Proposition 10.22 in \cite{BenyaminiLindenstrauss2000} for the proof of the following important result.
	
	\begin{prop}\label{prop:equivalence of nets}
		Two nets in the same infinite-dimensional Banach space are always Lipschitz equivalent. 
	\end{prop}
	
	This yields the following natural definition.
	
	\begin{defi} 
		\label{def:net-equivalence}
		Two infinite-dimensional Banach spaces $X$ and $Y$ are said to be \emph{net
			equivalent} if there exist a net $\cM$ in $X$ and a net $\cN$ in $Y$ such
		that $\cM$ and $\cN$ are Lipschitz equivalent.
		% and we denote $X \buildrel {N}\over {\sim} Y$,
	\end{defi}
	
	In the next proposition, we will show that two infinite-dimensional Banach spaces are net-equivalent if and only if they are coarse-Lipschitz equivalent. Moreover, we will prove, using partitions of unity, that in this case, the maps $f$ and $g$ in the definition of a coarse-Lipschitz equivalence can be taken to be continuous. This will be crucial for the proof of a general version of the Gorelik principle, which we will present in Chapter \ref{chapter:Gorelik}. 
	
	\begin{prop}
		\label{prop:CL-equivalences} 
		Let $X$ and $Y$ be two Banach spaces and $f\colon X\to Y$. The following assertions are equivalent.
		\begin{enumerate}[(i)]
			\item  The map $f$ is a coarse-Lipschitz equivalence.
			\item There exist $\theta_0>0$ and $K\ge 1$ such that for all $\theta\ge \theta_0$ and all maximal $\theta$-separated subsets $\cM$ of $X$, $\cN=f(\cM)$ is a net in $Y$ and for all $x,x'\in \cM$, 
			\begin{equation*}
				\frac{1}{K}\norm{x-x'}_X\le \norm{f(x)-f(x')}_Y\le K\norm{x-x'}_X.
			\end{equation*}
			\item There exist two \underline{continuous} coarse-Lipschitz maps $\varphi \colon X\to Y$ and $\psi \colon Y\to X$ and a constant $C\ge 0$ such that  
			\begin{equation*}
				\sup_{x\in X} \norm{(\psi\circ \varphi)(x)-x}_Y\le C,\ \sup_{y\in Y} \norm{(\varphi\circ \psi)(y)-y}_X\le C,
			\end{equation*}
			and 
			\begin{equation*}
				\sup_{x\in X} \norm{\varphi(x)-f(x)}_Y\le C.
			\end{equation*}
		\end{enumerate}
		In particular, if $X$ and $Y$ are infinite-dimensional Banach spaces, then they are net-equivalent if and only if they are coarse-Lipschitz equivalent.
	\end{prop}
	
	\begin{proof} 
		$(iii)\Rightarrow (i)$ is trivial.
		
		$(i)\Rightarrow (ii)$: Assume that there exist $g\colon Y\to X$ and constants $C,D,M>0$ such that
		$$\forall x\in X\ \ \norm{(g\circ f)(x)-x}\le C,\ \ \ \ \forall y\in Y\ \ \norm{(f\circ g)(y)-y}\le C$$
		and
		$$\forall x,x'\in X\ \ \norm{f(x)-f(x')}\le D+M\norm{x-x'},$$
		$$\forall y,y'\in Y\ \ \norm{g(y)-g(y')}\le D+M\norm{y-y'}.$$
		Let $\theta_0= (2C+D)(M+1)$, $\theta\ge \theta_0$ and $\mathcal M$ be a maximal $\theta$-separated subset of $X$. Let now $x\neq x' \in \mathcal M$, $y=f(x)$ and $y'=f(x')$. Then, $$\norm{f(x)-f(x')}\le D+M\norm{x-x'}\le \theta+M\norm{x-x'}\le (M+1)\norm{x-x'}.$$ On the other hand, $\norm{g(y)-x}\le C$ and $\norm{g(y')-x'}\le C$, which implies that
		$\norm{g(y)-g(y')}\ge \norm{x-x'}-2C$ and therefore
		$$\norm{x-x'}\le 2C+D+M\norm{y-y'}\le \frac{\theta}{M+1}+M\norm{y-y'}\le \frac{\norm{x-x'}}{M+1}+M\norm{y-y'}.$$
		It follows that $\norm{x-x'}\le (M+1)\,\norm{y-y'}$. So, $f$ is a Lipschitz isomorphism from $\mathcal M$ onto $\mathcal N=f(\mathcal M)$ and $K=M+1$ satisfies the required inequalities. In particular, $\mathcal N$ is $a$-separated with $a=\theta(M+1)^{-1}$.\\ 
		Finally, let $z\in Y$. There exists $x\in \mathcal M$ such that $\norm{x-g(z)}\le \theta$. Then, $$\norm{f(x)-z}\le \norm{f(x)-f(g(z))}+C\le D+M\theta+C=b.$$
		This shows that $\mathcal N$ is an $(a,b)$-net in $Y$.

		$(ii)\Rightarrow (iii)$: For $\theta\ge \theta_0$, we pick $(x_i)_{i\in I}$ a maximal $\theta$-separated subset of $X$. For $i\in I$, let $y_i=f(x_i)$. Then, by assumption, $(y_i)_{i\in I}$ is an $(a,b)$-net in $Y$, for some $0<a\le b$ and, for some $K\ge 1$, we have
		$$\forall i,j\in I\ \ \ \frac1K \norm{x_i-x_j}\le \norm{y_i-y_j} \le K\,\norm{x_i-x_j}.$$
		We can find  a locally finite continuous partition of unity $(f_i)_{i\in I}$ subordinated to the open cover $(U_i)_{i\in I}$ of $X$ and a locally finite continuous partition of unity $(g_i)_{i\in I}$ subordinated to the open cover $(V_i)_{i\in I}$ of $Y$, where $U_i$ is the open ball in $X$ of center $x_i$ and radius $\theta$ and $V_i$ is the open ball in Y of center $y_i$ and radius $b$ (see for instance \cite[Appendix B]{BenyaminiLindenstrauss2000}). Then, we set
		$$\forall x\in X\ \ \varphi(x)=\sum_{i\in I}f_i(x)\,y_i\ \ {\rm and}\ \ \forall y \in Y\ \ \psi(y)=\sum_{i\in I}g_i(y)\,x_i.$$
		Note first that $\varphi$ and $\psi$ are continuous.\\
		Let $x\in X$ and pick $i\in I$ such that $\norm{x-x_i}\le \theta$. Now, if $f_j(x)\neq 0$, then $\norm{x-x_j}\le \theta$, $\norm{x_i-x_j}\le 2\theta$ and $\norm{y_i-y_j}\le 2\theta K$. It follows that
		$$\norm{\varphi(x)-y_i}=\Big\|\sum_{j\colon f_j(x)\neq 0} f_j(x)\,(y_j-y_i)\Big\|\le 2\theta K.$$
		Let now $x'\in X$ and $j\in I$ be so that $\norm{x'-x_j}\le \theta$. Then, we have
		$$\norm{\varphi(x)-\varphi(x')}\le 4\theta K+\norm{y_i-y_j}\le 4\theta K+K\norm{x_i-x_j}\le 6\theta K+K\norm{x-x'}.$$
		This shows that $\varphi$ is coarse-Lipschitz and $\Lip_\infty(\varphi)\le K$. A similar proof shows that the same is true for $\psi$.\\
		For $x\in X$, pick again $i\in I$ such that $\norm{x-x_i}\le \theta$. If $g_j(\varphi(x))\neq 0$, then $\norm{\varphi(x)-y_j}\le b$ and $\norm{y_i-y_j}\le \norm{\varphi(x)-y_i}+\norm{\varphi(x)-y_j}\le 2\theta K+b$. Therefore,
		$$\norm{\psi(\varphi(x))-x_i}=\Big\|\sum_{j\colon g_j(\varphi(x))\neq 0} g_j(\varphi(x))\,(x_j-x_i)\Big\|\le K(2\theta K+b).$$
		We deduce that
		$$\norm{\psi(\varphi(x))-x}\le \norm{\psi(\varphi(x))-x_i}+\norm{x-x_i}\le K(2\theta K+b)+\theta=C_1.$$
		Similarly, we get that there exists $C_2\ge 0$ such that for all $y\in Y$, $\norm{\varphi(\psi(y))-y}\le C_2$.\\
		Finally, recall that $f$ is coarse-Lipschitz. So, there exist $D,E\ge 0$ such that for all $x,x'\in X$, $\norm{f(x)-f(x')}\le D\norm{x-x'}+E$. Since
		$$\forall x\in X,\ \ \ \varphi(x)-f(x)=\sum_{j\colon f_j(x)\neq 0} f_j(x)\,(f(x_j)-f(x)),$$
		and $\norm{x_j-x}\le \theta$ whenever $f_j(x)\neq 0$, we obtain that
		$$\forall x\in X,\ \ \norm{\varphi(x)-f(x)}\le D\theta+E=C_3.$$
		We conclude the proof of this implication by taking $C=\max\{C_1,C_2,C_3\}$.
		
	\end{proof}

	\begin{rema}
		The equivalence between $(i)$ and $(ii)$ in Proposition \ref{prop:CL-equivalences} holds in the metric setting with the same proof.
	\end{rema}

	\section{Notes}
	
	The notion of coarse embedding is due to Gromov, but it appears difficult to track its exact origin precisely. In \cite[Chapter 4]{GromovRTG}, Gromov discusses the notion of uniformly bounded map, UB-equivalence and placing (corresponding to coarse map, coarse equivalence and coarse embedding in the modern terminology). In fact, Gromov uses the terminology uniformly equivalent on the large-scale in \cite[Section 0.2.
	D]{GromovGGT} and refers to the placings in \cite{GromovRTG} as uniform embeddings (see the Summary of Section 1.H. and Section 7.E. in \cite{GromovGGT}). Following \cite{GromovGGT}, the terminology uniform embeddings was retained in the geometers and geometric group theory community. For instance in G. Yu's breakthrough paper \cite{Yu00} entitled ``The coarse Baum–Connes conjecture for spaces which admit a uniform embedding into Hilbert space" and in \cite{DGLY02}, \cite{DadarlatGuentner03}, \cite{GuentnerKaminker04}, or \cite{BrownGuentner05} for instance, the notion of uniform embedding is what is now (almost universally) called a coarse embedding. A gradual shift in terminology seemed to have occurred around the mid-2000s. In \cite{KasparovYu06}, Kasparov and Yu proved ``the coarse geometric Novikov conjecture for spaces which admit a (coarse) uniform embedding into a uniformly convex Banach space". The fact that Yu's work establishes a striking link between the coarse geometric Novikov conjecture and the geometry of Banach spaces led Banach space geometers to become interested in these nonlinear embeddings preserving the large-scale structure. In the Banach space geometry community, the terminology uniform embedding already had a very precise meaning since the beginning of the 20th century, namely an injective and uniformly continuous map with a uniformly continuous inverse. Since in the language of coarse geometry (see \cite{Roe_book03} for instance), a uniform embedding $f\colon X \to Y$ in the sense of Gromov is a coarse equivalence between $X$ and $f(X)$, the clash of terminology was quickly resolved by using the now standard terminology coarse embedding. For the profound and beautiful applications of coarse embeddings to topology, geometric group theory and noncommutative geometry, we refer the interested reader to \cite{NowakYu} and \cite{DrutuKapovich} and the references therein. 
	
	\medskip
	
	The most advanced results around Problem \ref{pb:coarse-vs-uniform} are due to C. Rosendal \cite{Rosendal}. Let us  describe them briefly. We say that a map $f\colon X \to E$ is \emph{uncollapsed} if there exists $r>0$ such that $\inf\{\norm{f(x)-f(y)} \colon  \norm{x-y}\ge r\}>0$. Rosendal shows that if there exists a map $f\colon X \to E$ which is uncollapsed and uniformly continuous, then $X$ admits a simultaneously uniform and coarse embedding into $\ell_p(E)$, for $1\le p<\infty$. An embedding that is simultaneously uniform and coarse was called a \emph{strong embedding} by Kalton. As an immediate corollary, one gets that if $X$ is uniformly embeddable into $\ell_p$, $L_p$ (for some $1\le p<\infty$), a reflexive, super-reflexive, stable, super-stable space, or a space with nontrivial type or cotype, then $X$ admits a strong embedding into a space of the same kind.
	
	In the other direction, the following can be deduced from the above results. Assume that $X$ and $E$ are Banach spaces such that for every net $N$ in $X$ and every Lipschitz map $f\colon N \to E$, $f$ extends to a uniformly continuous map $g\colon X \to E$. Then, if $X$ coarsely embeds into $E$, then $X$ strongly embeds into $\ell_p(E)$, for $1 \le p <\infty$. A remark has to be made about the assumptions of this theorem. Indeed, A. Naor \cite{Naor2015} produced an example of Banach spaces $X$ and $E$, a net $N$ in $X$ and a Lipschitz map $f\colon N \to E$ so that, for any uniformly continuous $g\colon X \to E$, $\sup\{\norm{f(x)-g(x)} \colon \ x\in N\}=+\infty$. 
	
	\medskip 
	
	Around this question, we also wish to mention a result due to A. Swift who proves in \cite{Swift2018a} that for any given index set $\Gamma$, coarse, uniform and bi-Lipschitz embeddability into  $c_0(\Gamma)$ are equivalent notions for normed linear
	spaces.
	
	\medskip

	It is also important to recall that Proposition \ref{prop:equivalence of nets} is not true in finite dimensions greater than one. Indeed, there exists a net in $\R^2$ which is not Lipschitz equivalent to $\Z^2$. We refer the reader to \cite[Theorem 10.26]{BenyaminiLindenstrauss2000} for its proof. This result was obtained independently by McMullen \cite{McMullen} and Burago and Kleiner \cite{BuragoKleiner}.

	\section{Exercises}
	
	\begin{exer}
		Prove Proposition \ref{prop:CL}.   
	\end{exer}
	
	\begin{exer}\,
		\begin{enumerate}
			\item Show that $\bR$ and $\bZ$ are coarse-Lipschitz equivalent.
			\item Show that $\bZ$ does not coarsely embed into $\bN$.
		\end{enumerate}
	\end{exer}
	
	\begin{exer}
		Assume that $n>m \in \bN$. Show that $\bR^n$ does not coarsely embed into $\bR^m$. 
	\end{exer}
	
	\begin{proof}[Hint]
		Show that there is no map $f\colon \bR^n \to \bR^m$ such that $\omega_f(1)<\infty$ and $\rho_f(1)>0$ by comparing the maximal number of elements of a $1$-separated family in large balls of $(\bR^n,\norm{\cdot}_\infty)$ and $(\bR^m,\norm{\cdot}_\infty)$.
	\end{proof}
	
	\begin{exer}\,
		\begin{enumerate}
			\item Give an example of two unbounded metric spaces $M$ and $N$ and a map $f\colon M\to N$ such that $\omega_f(t_0)<\infty$ for some $t_0>0$, but $f$ is not coarse-Lipschitz.
			\item Give an example of  a map $f\colon M\to N$ between unbounded metric spaces such that $\Lip_\infty(f)<\infty$  but $f$ is not coarse-Lipschitz.
		\end{enumerate}
	\end{exer}
	
	\begin{exer}
		\label{exer:Swift}
		The density character of a metric space $M$, denoted by $\mathrm{dens}(M)$, is the smallest cardinality of a dense subset. 
		\begin{enumerate}
			\item Show that for every net $\cM$ in a metric space $M$ one has $\abs{\cM}\le \mathrm{dens}(M)$.
			\item Show that for every net $\cM$ in a normed space $X$ one has $\abs{\cM} = \mathrm{dens}(X)$.
			\item Show that if there exists a map $f$ from a normed space $X$ into a metric space $M$ such that $\lim_{t\to \infty}\rho_f(t)=\infty$, then $\mathrm{dens}(X)\le \mathrm{dens}(M)$. Deduce that separability is preserved under coarse embeddings of normed spaces.
		\end{enumerate}
	\end{exer}
	
	\begin{exer}
		\label{exo:CL-equiv->Lip-equiv-ultra}
		Let $\cU$ be a nonprincipal ultrafilter on $\bN$. Show that if two Banach spaces $X$ and $Y$ are coarse-Lipschitz equivalent, then $X^{\cU}$ and $Y^{\cU}$ are Lipschitz equivalent.
	\end{exer}
	
	\begin{exer}\label{ex:embedding-equivalence}
		Let $M$ and $N$ be two metric spaces and let $f\colon M\to N$.
		\begin{enumerate}
			\item Assume that $f$ is a coarse Lischitz embedding. Show that $f$ is a coarse-Lipschitz equivalence from $M$ onto $f(M)$. 
			\item Assume that $f$ is a coarse embedding. Show that $f$ is a coarse equivalence from $M$ onto $f(M)$.
		\end{enumerate}
		
	\end{exer}
	
	%%%%%%%%%%%%%%%%%%%%%%%%%%%%%%%%%%%%%%%%%%%%%%%%%%%%%%%%%%%%%%%%%%%%%%%%%%%%%%%%%%%%%%%%%%%%
	
	\chapter{ Metrically universal Banach spaces}
	\label{chapter:universal}
	
	In this very short chapter, we will describe the fundamental examples of Banach spaces that are universal for the class of separable metric or Banach spaces and for isometric or bi-Lipschitz embeddings. By a universal Banach space, for a given class $\cC$ of spaces and a given family $\cF$ of embeddings, we mean a space $U$ such that for any $X \in \cC$, there exists an embedding $f\in \cF$ from $X$ into $U$. When $U$ can be chosen in $\cC$, it is usually an important result. In this chapter, we only consider isometric and bi-Lipschitz embeddings and the class of embedded spaces will be the class of all separable metric or Banach spaces.
	
	\section{\texorpdfstring{The space $\ell_\infty$ is isometrically universal for separable metric spaces}{The space is isometrically universal for separable metric spaces}}
	
	We start with the applications of the most fundamental isometric embedding, the so-called \emph{Fr\'echet embedding}. We recall that, for a set $\Gamma$, $\ell_\infty(\Gamma)$ is the space of all bounded maps from $\Gamma$ to $\bR$ equipped with the supremum norm:
	$$\forall x=(x_\gamma)_{\gamma \in \Gamma},\ \ \norm{x}_\infty= \sup_{\gamma \in \Gamma}\abs{x_\gamma}.$$
	We most often simply write $\ell_\infty$ instead of $\ell_\infty(\bN)$.
	
	\begin{theo}
		\label{thm:Frechet} 
		Let $(M,d)$ be a metric space. Then, $(M,d)$ isometrically embeds into $\ell_\infty(M)$. In particular, if $M$ is finite, it isometrically embeds into $\ell_\infty^n$, where $n$ is the cardinality of $M$.
	\end{theo} 
	
	\begin{proof} 
		Fix $x_1 \in M$  and consider the map $f\colon M \to \ell_\infty(M)$ defined by 
		$$f(x)=\big(d(x,y)-d(y,x_1)\big)_{y\in M}.$$
		Note first that for all $x\in M$, $\|f(x)\|_\infty \le d(x,x_1)<\infty$. Then, 
		it is easily checked that for all $x,x' \in M$, $\|f(x)-f(x')\|_\infty=d(x,x')$.
	\end{proof}
	
	The target space can be somewhat reduced by using a dense family in the embedded space. For instance, we have:
	
	\begin{coro}[Fr\'echet embedding, 1906]
		\label{cor:univlinfini} 
		%Let $(M,d)$ be a separable metric space. Then, $(M,d)$ isometrically embeds into $\ell_\infty$.
		Every separable metric space isometrically embeds into $\ell_\infty$
	\end{coro}
	
	\begin{proof}
		Let $(M,d)$ be a separable metric space, $(x_n)_{n=1}^\infty$ be a dense sequence in $M$ and let $N:=\{x_n\colon n\in \bN\}$. From Theorem \ref{thm:Frechet}, we deduce that $(N,d)$ isometrically embeds into $\ell_\infty$. By completeness of $\ell_\infty$, this embedding extends to an isometric embedding from $(M,d)$ into $\ell_\infty$. In fact, this extension can be directly given by
		$$\forall x\in M,\ \ f(x)=\big(d(x,x_n)-d(x_n,x_1)\big)_{n\in \bN}.$$
	\end{proof}
	
	The space $\ell_\infty$ is thus easily seen to be isometrically universal for separable metric spaces. Its main disadvantage is that it is not itself separable. This problem can be overcome by using the classical Banach-Mazur theorem, which we prove in the next section. 
	
	\section{\texorpdfstring{The space $C([0,1])$ is linearly isometrically universal for separable Banach spaces}{The space is linearly isometrically universal for separable Banach spaces}}
	
	Recall that for a compact space $K$, $C(K)$ is the space of all real-valued continuous functions defined on $K$, equipped with the supremum norm. The \emph{Cantor space}, denoted by $\Delta$, is the space $\{0,1\}^\bN$ equipped with the product topology of the discrete topology on $\{0,1\}$.
	
	\begin{theo}[The Banach-Mazur theorem - 1932]
		\label{thm:BanachMazur} 
		%Let $(X,\norm{\cdot}_X)$ be a separable Banach space. Then, $X$ is linearly isometric to a subspace of $C(\Delta)$ and linearly isometric to a subspace of $C([0,1])$.
		Every separable Banach space is linearly isometric to a subspace of $C(\Delta)$ and linearly isometric to a subspace of $C([0,1])$.
	\end{theo}
	
	\begin{proof} 
		We only outline the main steps of the beautiful proof of this result and we refer the reader, for instance, to the textbook \cite{AlbiacKalton2016} for the details. 
		
		Let $(X,\norm{\cdot}_X)$ be a separable Banach space. First, since $X$ is separable, $K:=(B_{X^*},w^*)$ is a compact metrizable space, where $w^*$ denotes the weak$^*$ topology induced by $X$ on $X^*$.  Then, it follows from the Hahn-Banach theorem that the map $J\colon (X,\norm{\cdot}_X)\to (C(K),\norm{\cdot}_\infty)$ defined by $J(x)(x^*):=x^*(x)$ for $x^*\in K$, is a linear isometry. Next, we use the well-known fact that $K$, like any nonempty compact metrizable space, is the continuous image of $\Delta$. So, if $\varphi\colon\Delta \to K$ is a continuous surjection, the map $f \mapsto f\circ \varphi$ is a linear isometry from $C(K)$ into $C(\Delta)$. This shows that $X$ is linearly isometric to a subspace of $C(\Delta)$. 
		
		Now, if $T$ denotes the classical triadic Cantor subset of $[0,1]$, we have that $T$ is homeomorphic to $\Delta$ and therefore (by composition with the homeomorphism) that $C(\Delta)$ is linearly isometric to $C(T)$. The last step consists of proving that $C(T)$ is linearly isometric to a subspace of $C([0,1])$. The isometry is given by the following operation: extend any continuous map from $T$ to $\bR$ in such a way that it is affine on any connected component of $[0,1]\setminus T$.
	\end{proof}
	
	\begin{rema} 
		Note that since $\Delta$ and $[0,1]$ are metrizable, both spaces $C(\Delta)$ and $C([0,1])$ are separable. This is a consequence, for instance, of the Stone-Weierstrass theorem.
	\end{rema}
	
	We now use these universal spaces in the metric setting.
	
	\begin{coro} 
		%Let $(M,d)$ be a separable metric space. Then, $(M,d)$ isometrically embeds into $C(\Delta)$ and  into $C([0,1])$.
		Every separable metric space isometrically embeds into $C(\Delta)$ and into $C([0,1])$.
	\end{coro} 
	
	\begin{proof} 
		Let $(M,d)$ be a separable metric space and $f\colon M \to \ell_\infty$ be an isometric embedding given by Corollary \ref{cor:univlinfini}. Let $X$ be the closed linear span of $f(M)$. Since $M$ is separable, $X$ is a closed separable subspace of $\ell_\infty$, which, by Theorem \ref{thm:BanachMazur}, isometrically linearly embeds into $C(\Delta)$ and into $C([0,1])$.
	\end{proof}

	\section{\texorpdfstring{The space $\co$ is Lipschitz universal for separable metric spaces}{The space is Lipschitz universal for separable metric spaces}}
	
	In 1974, I. Aharoni \cite{Aharoni1974} proved that a much ``smaller'' Banach space could serve as a universal space for separable metric spaces and bi-Lipschitz embeddings. He showed that every separable metric space is bi-Lipschitz isomorphic to a subset of $\co$ and that this could be done with distortion $6+\eps$, for any $\eps>0$. He also showed that it could not be done with distortion less than $2$ in general. In fact, the embedding built by Aharoni took values in $c_0^+$, the positive cone of $\co$. The distortion of the embedding into $c_0^+$ was later improved by Assouad \cite{Assouad1978} and finally by Pelant \cite{Pelant1994}, who showed that every separable metric space embeds into $c_0^+$ with distortion $3$ and that the distortion $3$ is optimal in general. In this section, we will focus on the optimal distortion for bi-Lipschitz embeddings of separable metric spaces into $\co$, which turns out to be equal to $2$. We will start with the simple proof, due to Aharoni  \cite{Aharoni1974}, that $\ell_1$ cannot be embedded into $\co$ with distortion less than $2$. Then, we will prove that every separable metric space embeds into $\co$ with distortion $2$. This result is due to Kalton and Lancien \cite{KaltonLancien2008}. We will follow the presentation given by Baudier and Deville in \cite{BaudierDeville}.
	
	\begin{prop}
		\label{prop:c_0l_1distorsion}
		The Banach space $\ell_1$ does not admit a bi-Lipschitz embedding into $\co$ with distortion less than $2$.
	\end{prop}
	
	\begin{proof}
		Let $\vep \in (0,1]$ and assume that there exists $f\colon \ell_1 \to c_0$ such that for all $x,y \in \ell_1$,
		\begin{equation*}
			\norm{x-y}_1 \le \norm{f(x)-f(y)}_\infty \le (2-\vep)\norm{x-y}_1.
		\end{equation*}
		We denote by $(e_n)_{n=1}^\infty$ the canonical basis of $\ell_1$ and, for $x\in \ell_1$, we let $f(x)=\big(f_k(x)\big)_{k=1}^\infty \in c_0$.\\
		For $2<n\ne m$ in $\bN$, we consider
		$$A_{n,m}:=\big\{k\in \bN\colon |f_k(e_1+e_n)-f_k(e_2+e_m)|\ge 4\big\}.$$
		Since $f(e_1+e_n)-f(e_2+e_m) \in c_0$ and $\norm{f(e_1+e_n)-f(e_2+e_m)}_\infty\ge 4$, we have that $A_{n,m}$ is not empty. Next, we have that for all $k\in A_{n,m}$,
		\begin{align*}
			\abs{f_k(e_1)-f_k(e_2)}\ge & \abs{f_k(e_1+e_n)-f_k(e_2+e_m)}-\abs{f_k(e_1+e_n)-f_k(e_1)}\\
			&-\abs{f_k(e_2)-f_k(e_2+e_m)}\ge 4-2(2-\eps)=2\eps.   
		\end{align*}
		We now use the fact that $f(e_1)-f(e_2)\in c_0$ to deduce that there exists a finite set $A$ such that $A_{n,m}$ is included in $A$ for all $(n,m)$. Now, for all $k\in \bN$, the sequence $(f_k(e_n))_n$ is bounded by $\norm{f(0)}_\infty+(2-\vep)$. By an easy diagonal argument, we obtain that there exists an increasing sequence of integers $(n_i)_i$ so that for all $k\in \N$, the sequence $(f_k(e_{n_i}))_i$ is converging. Since $A$ is finite, we deduce that 
		$$\exists i_0\in \N\ \ \forall k\in A\ \ \forall j>i\ge i_0\ \ \ |f_k(e_{n_j})-f_k(e_{n_i})|<2\eps.$$
		But, for all $k\in A_{n,m}$,
		\begin{align*}
			\abs{f_k(e_n)-f_k(e_m)}\ge &\abs{f_k(e_1+e_n)-f_k(e_2+e_m)}-\abs{f_k(e_1+e_n)-f_k(e_n)}\\
			&-\abs{f_k(e_m)-f_k(e_2+e_m)}\ge 4-2(2-\eps)=2\eps.   
		\end{align*}
		This is clearly a contradiction.
	\end{proof}
	
	We now turn to the sharpest version of Aharoni's embedding theorem. Our precise statement is the following.
	
	\begin{theo}
		\label{thm:Aharoni}
		Let $(M,d)$ be a separable metric space. Then, there exists $f\colon M\to c_0$ such that for all $x,y\in M$ with $x\neq y$,
		\begin{equation*}
			\frac12 d(x,y) < \norm{f(x)-f(y)}_\infty \le d(x,y).
		\end{equation*}
	\end{theo}
	
	We first need some preparatory notation and lemmas. None of the lemmas requires a separability assumption, and hence for a moment $(M,d)$ is merely a metric space. Given $E$ a subset of $M\times M$, we let 
	$$\delta(E):=\inf\{d(x,y),\ (x,y)\in E\}\ \ \text{and}\ \ D(E):=\sup\{d(x,y),\ (x,y)\in E\}.$$
	as well as
	$$P_1(E):=\{x\in M\colon \exists y\in M\ (x,y)\in E\}\ \ \text{and}\ \ P_2(E):=\{y\in M\colon  \exists x\in M\ (x,y)\in E\}.$$
	Finally, $R(E):=P_1(E)\times P_2(E)$ is the ``rectangle envelope'' of $E$, or the smallest product of subsets of $M$ containing $E$. We can now state our first lemma.
	
	\begin{lemm}
		\label{lem:separate+vanish} Let $(M, d)$ be a metric space, $U,V,F$ be three nonempty subsets of $M$ and $\eps \ge 0$. Then, there  exists a $1$-Lipschitz map $f\colon M \to \bR$  such that :
		\begin{enumerate}[(i)]
			\item For all $x\in F$, $\abs{f(x)}\le \eps$.
			\item For all $(x,y)\in U\times V$, $f(x)-f(y)=\min\{\delta(U\times V), \delta(U\times F)+\delta(V\times F)+2\eps$\}.
		\end{enumerate}
	\end{lemm}
	
	\begin{proof} 
		We can find $s,t$ such that 
		$$-\delta(V\times F)-\eps \le s\le 0 \le t\le \delta(U\times F)+\eps$$
		and
		$$t-s = \min\{\delta(U\times V ), \delta(U\times F) + \delta(V\times F) + 2\eps\}.$$
		Then we set 
		$$f(x) := \min \{d(x, U) + t, d(x, V ) + s, d(x, F)+ \eps\},\ \ x\in M.$$
		The function f is $1$-Lipschitz as the infimum of 1-Lipschitz functions.\\
		If $x\in U$, then $f(x)=\min\{t,d(x, V ) + s, d(x, F)+ \eps\}=t$, since $d(x, V )+s \ge \delta(U\times V)+s \ge t$
		and $d(x, F) + \eps \ge \delta(U\times F) + \eps \ge t$. On the other hand, if $y \in V$, then  $f(y) = \min \{d(y, U) + t, s, d(y, F) + \eps\}= s$, because $s\le 0$ and $t\ge 0$. Therefore, if $x \in U$ and $y \in  V$, then $f(x)-f(y) = t-s$, which insures $(ii)$.\\
		Finally, if $x \in F$, then $f(x) = \min \{d(x, U) + t, d(x, V ) + s, \eps\}\le \eps$, while $d(x, U) + t \ge 0$ and $d(x, V ) + s \ge \delta(V\times F) + s \ge -\eps$, so $f(x) \ge -\eps$. This proves $(i)$.
	\end{proof}
	
	\begin{lemm}
		\label{lem:2balls}
		Let $(M,d)$ be a metric space and $B_1,B_2$ be two closed balls of $M$ with respective radii $r_1$ and $r_2$. Assume that $E$ is a subset of $B_1 \times B_2$ such that $\delta(E)>2(r_1+r_2)$. Then, there exists a partition $\{E_1,\ldots,E_N\}$ of $E$ such that for all $n\in \{1,\ldots,N\}$, 
		\begin{equation*}
			D(E_n)<2\delta(R(E_n)).
		\end{equation*}
	\end{lemm}
	
	\begin{proof} 
		We assume, as we may, that $r_1\ge r_2\ge 0$. Set $\eps:=\delta(E)-2(r_1+r_2)>0$. Pick $a_0:=\delta(E)<a_1<\dots<a_{N-1}<D(E)<a_N$ so that $a_n-a_{n-1}< \eps$ for $n=1,\dots,N$. Our partition of $E$ will be given by 
		$$E_n:=\{(x,y)\in E,\ a_{n-1}\le d(x,y)<a_n\},\ \text{for}\ n=1,\dots,N.$$
		It is clear from our construction that 
		$$D(E_n)<\delta(E_n)+\eps=\delta(E_n)+\delta(E)-2(r_1+r_2)\le 2\delta(E_n)-4r_2.$$
		Fix now $(u,v)\in R(E_n)$. Note that $v\in B_2$ and that there is $v'\in B_2$ such that $(u,v')\in E_n$. Then, we have
		$$D(E_n)<2\delta(E_n)-4r_2\le 2d(u,v')-4r_2\le 2d(u,v')-2d(v,v')\le 2d(u,v).$$
		Taking the infimum over $(u,v)\in R(E_n)$ gives that $D(E_n)<2\delta(R(E_n)).$
	\end{proof}
	
	For a finite subset $F$ of a metric space $(M,d)$ and $\alpha>0$, we let 
	$$A(F,\alpha):=\Big\{(x,y)\in M\times M\colon d(x,F)+d(y,F)+\alpha \le \frac12 d(x,y)\Big\}.$$
	
	\begin{lemm}
		\label{lem:levels}
		Let $(M,d)$ be a metric space. Assume that $F\subset G$ are two finite subsets of $M$ and $0<\alpha<\beta$. Then, there exists a partition $\{E_1,\dots,E_N\}$ of $\Delta:=A(G,\alpha)\setminus A(F,\beta)$ such that for all $n\in \{1,\dots,N\}$,
		$$D(E_n)<2\min\big\{\delta(R(E_n)),\delta(P_1(E_n)\times F)+\delta(F\times P_2(E_n))+2\beta\big\}.$$
	\end{lemm}
	
	\begin{proof} 
		Let $(x,y)\in A(G,\alpha)$. Then, 
		$$2d(x,G)+2d(y,G)\le d(x,y) \le d(x,G)+d(y,G)+\diam (G).$$ 
		This implies that $d(x,G)+d(y,G)\le \diam (G)$. Since $G$ is finite, we deduce that there exists a bounded subset $B$ of $M$ such that $\Delta \subset B\times B$. Using that $G$ is finite again, we find a partition $\{B_1,\dots,B_m\}$ of $B$ so that $|d(x,a)-d(x',a)|\le \frac{\alpha}{5}$ for all $j\in\{1,\dots,m\}$, $x,x' \in B_j$ and $a\in G$. It follows that 
		\begin{equation}
			\label{eq:levels}
			d(x,a)<d(B_j,a)+\frac{\alpha}{4}
		\end{equation}
		for all $j\in\{1,\dots,m\}$, $x\in B_j$ and $a \in G$. Since $G$ is finite, we can find $a_j \in G$ such that $d(B_j,a_j)=\delta(B_j\times G)$. Thus,  
		$$B_j\subset B(a_j,r_j)\ \ \text{with}\ \ r_j:=\delta(B_j\times G)+\frac{\alpha}{4}.$$ 
		For $k,j \le m$, we now set $C_{j,k}:=\Delta \cap (B_j \times B_k)$. Then, $(C_{j,k})_{j,k}$ is a partition of $\Delta$ such that $C_{j,k} \subset B(a_j,r_j)\times B_k$. We have that for all $(x,y)\in C_{j,k}$,
		$$d(x,y)\ge 2(d(x,G)+d(y,G)+\alpha)\ge 2(\delta(B_j\times G)+\delta(B_k\times G)+\alpha)=2(r_j+r_k)+\alpha.$$
		Thus $\delta(C_{j,k})>2(r_j+r_k)$ and we can apply Lemma \ref{lem:2balls} to each $C_{j,k}$ to deduce the existence of a partition $\{E_1,\ldots,E_N\}$ of $\Delta$ such that for all $n\in\{1,\dots,N\}$, $D(E_n)<2\delta(R(E_n))$. Then, only one inequality remains to be proved. So, consider $n,j,k$ so that $E_n\subset C_{j,k}$ and fix $(x,y)\in E_n$. Since $(x,y)\notin A(F,\beta)$, we have that $d(x,y)< 2(d(x,F)+d(y,F)+\beta)$. Taking the infimum over $a\in F$ in \eqref{eq:levels}, we deduce that 
		$$d(x,y)<2(\delta(B_j\times F)+\delta(B_k\times F)+\frac{\alpha}{2}+\beta).$$
		Recall that $P_1(E_n)\subset B_j$ and $P_2(E_n) \subset B_k$, so 
		$$d(x,y)< 2(\delta(P_1(E_n)\times F)+\delta(F\times P_2(E_n))+\frac{\alpha}{2}+\beta).$$
		Therefore, $D(E_n)<2(\delta(P_1(E_n)\times F)+\delta(F\times P_2(E_n))+2\beta)$.
	\end{proof}
	
	We can now build our embedding.
	
	\begin{proof}[Proof of Theorem \ref{thm:Aharoni}] 
		We will construct a sequence $(f_n)_{n=1}^\infty$ of $1$-Lipschitz maps from $M$ to $\bR$ such that for all $x\in M$, $(f_n(x))_n$ tends to $0$ and a partition $\{E_n\colon n\in \bN\}$ of $U:=\{(x,y)\in M\times M,\ x\neq y\}$ so that for any $n\in \bN$, there exists $c_n\in (0,\infty)$ with the property that $f_n(x)-f_n(y)=c_n$ for all $(x,y)\in E_n$ and $D(E_n)<2c_n$. It will then be clear that $f=(f_n)_{n=1}^\infty$ satisfies the desired properties for our embedding.
		
		We use now, at last, the fact that $M$ is separable. So, we pick a sequence $(x_k)_{k=1}^\infty$ which is dense in $M$, set $F_k:=\{x_1,\dots,x_k\}$ and fix a decreasing sequence $(\alpha_k)_{k=1}^\infty$  of positive reals that is converging to $0$. Finally, for $k\ge 1$, we set $\Delta_k:=A(F_{k+1},\alpha_{k+1})\setminus A(F_{k},\alpha_{k})$. First, we note that $\{\Delta_k\colon k\in \bN\}$ is a partition of $U$. Indeed, given $(x,y)\in U$, set $\sigma_k:=2(d(x,F_k)+d(y,F_k)+\alpha_k)$. We have that $0<d(x,y)<\sigma_1$, the sequence $(\sigma_k)_{k=1}^\infty$ is strictly decreasing and, by density of the $x_n$, tending to $0$. Therefore, there exists a unique $k\in \bN$ such that $\sigma_{k+1}\le d(x,y)<\sigma_k$.\\
		Then, by Lemma \ref{lem:levels}, there exist $0=n_0<n_1<\dots <n_k<\cdots$ and subsets $(E_n)_{n=1}^\infty$ of $M\times M$ such that for all $k\ge 1$, $\{E_n,\ n_k<n\le n_{k+1}\}$ is a partition of $\Delta_k$ with the property that $D(E_n)<2c_n$, where
		$$c_n:=\min\big\{\delta(R(E_n)),\delta(P_1(E_n)\times F_k)+\delta(F_k\times P_2(E_n))+2\alpha_k\big\}.$$
		For $n_k<n\le n_{k+1}$, we now use Lemma \ref{lem:separate+vanish} to produce a $1$-Lipschitz map $f_n\colon M\to \bR$ such that $\abs{f_n}\le \alpha_k$ on $F_k$ and $f(x)-f(y)=c_n$ for all $(x,y)\in E_n$.\\
		It remains to show that $(f_n(x))_{n=1}^\infty \in c_0$ for all $x\in M$. Fix $x\in M$ and $\eps>0$ and pick $j\in \bN$ such that $d(x,x_j)<\frac{\eps}{2}$ and $k\ge j$ so that $\alpha_k<\frac{\eps}{2}$. Then, since all the $f_n$ are $1$-Lipschitz, we get that for all $n>n_k$
		$$\abs{f_n(x)}\le \frac{\eps}{2}+\abs{f_n(x_j)}\le \frac{\eps}{2}+\alpha_k <\eps.$$
		We have shown that $\lim_n f_n(x)=0$. This finishes our proof.
	\end{proof}
	
	\begin{coro}[Aharoni's embedding theorem, 1974]
		Every separable metric space admits a bi-Lipschitz embedding into $\co$.
	\end{coro}
	
	\section{Notes}
	
	It is an important open question to know whether $c_0$ is a minimal separable Banach space that is universal for bi-Lipschitz embeddings and separable metric spaces. This can be precisely restated as follows.
	
	\begin{prob}
		\label{pb:Lipschitz-embed-rigidity-c0} 
		Let $X$ be a Banach space and assume that $c_0$ bi-Lipschitzly embeds into $X$. Does it imply that $c_0$ linearly embeds into $X$?
	\end{prob}
	
	If $K$ is an infinite compact metric space, then $C(K)$ contains an isometric linear copy of $c_0$. One can wonder for which $K$, $C(K)$ is a better universal space than $c_0$. If $K$ is uncountable, then it contains a closed subset homeomorphic to the Cantor space, so $C(K)$ is isometrically universal. Therefore, the interesting case is when $K$ is infinite and countable. This question has been solved by Proch\'azka and S\'anchez-Gonz\'alez \cite{ProchazkaSanchez} who constructed a separable metric space $M$ such that $M$ does not bi-Lipschitzly embed with distortion less than $2$ in any space $C(K)$ for $K$ compact and countable. However, it is not clear whether their
	example embeds into $\ell_1$ isometrically (or with distortion less than 2). The study of the minimal distortion of $\ell_1$ into $C(K)$ spaces has been addressed in \cite{BaudierFreemanSchlumprechtZsak}, where, among other things, it is proved that, for $k\in \bN$, this distortion is at least $\frac{k+1}{k}$ if the Cantor-Bendixon index of $K$ is less than or equal to $k+1$. 
	
	\section{Exercises}
	
	\begin{exer} 
		Show that every separable Banach space $X$ is linearly isometric to a subspace of $\ell_\infty$.
	\end{exer}
	\begin{proof}[Hint]
		Use a dense countable subset of $X$ and the associated norming functionals. 
	\end{proof}
	
	\begin{exer}
		\label{ex:co-sum-fds}
		\,
		\begin{enumerate}
			\item Show that for every $\vep>0$ and every finite-dimensional normed space $E$, there exists $N\in \N$ such that $E$ linearly embeds into $\ell_\infty^N$ with distortion at most $1+\eps$. 
			\item Let $(F_n)_{n=1}^\infty$ be a sequence of finite-dimensional Banach spaces. Show that $(\sum_{k=1}^\infty F_{k})_{c_0}$ linearly embeds into $\co$.
		\end{enumerate}
	\end{exer}
	\begin{proof}[Hint]
		For $1.$, consider a finite net in the unit sphere of $E$ and the associated norming functionals. 
	\end{proof}
	
	\begin{exer}
		Let $K$ be an infinite compact metric space. Show that every separable metric space admits a bi-Lipschitz embedding into $C(K)$ with distortion at most $2$.     
	\end{exer}
	
	\begin{exer} 
		Show that the distortion of any bi-Lipschitz embedding from $c$ (the space of converging sequences equipped with the supremum norm) into $\co$ is at least $2$.
	\end{exer}
	
	\begin{exer} 
		Let $p\in [1,\infty)$. Show that the distortion of any bi-Lipschitz embedding from $\ell_p$ into $\co$ is at least $2^{1/p}$.
		
		Note: it is shown in \cite{KaltonLancien2008} that this is optimal: $\ell_p$ bi-Lipschitzly embeds into $\co$ with distortion $2^{1/p}$.
	\end{exer}
	
	\begin{exer} A metric space $(M,d)$ is said to be \emph{uniformly discrete} if there exists $\theta >0$ such that $d(x,y)\ge \theta$ for all $x\neq y \in M$. We denote by $G_\infty$ the integer grid of $c_0$ that is the set $\Z^{\N} \cap c_0$, equipped with the metric induced by the norm of $c_0$. Prove that for any $\eps>0$ and any countable uniformly discrete metric space $(M,d)$, we have that $M$ bi-Lipschitz embeds into $G_\infty$ with distortion less than $2+\eps$. 
		
		\begin{proof}[Hint]
			Show that $G_\infty$ is a $(1,1)$-net in $c_0$, dilate the metric $d$ if needed, and use Theorem \ref{thm:Aharoni}.
		\end{proof}
		
	\end{exer}

	%{\bf suggestions} Autres espaces universels ? Urysohn ? Gurarii ? meilleure constante de plongement dans $C(K)$ pour $K$ d\'enombrable ?  Distorsion des $C(K)$ les uns dans les autres ?  {\bf Theorem 2.2 in paper with Beata!} Distance de Banach Mazur entre $c$ et $\co$ (preuve d'Audrey).

	%%%%%%%%%%%%%%%%%%%%%%%%%%%%%%%%%%%%%%%%%%%%%%%%%%%%%%%%%%%%%%%%%%%%%%%%%%%%%%%%%%%%%%%%%%%%%

	\chapter[Linear reductions of Lipschitz classification problems]{Linear reductions of Lipschitz classification problems}
	\label{chapter:linear-reductions}
	
	In this chapter, we survey the basic tools for linearizing Lipschitz maps and their first applications to the Lipschitz classification of some classical Banach spaces: essentially the $L_p$-spaces for $p\in (1,\infty)$. The first method we describe is naturally by differentiating and for that purpose we briefly introduce the Radon-Nikod\'{y}m property for Banach spaces. Then, we turn to the linearization of Lipschitz retracts through the use of invariant means. Finally, we describe a few more sophisticated results based on the use of Fréchet differentiability of real-valued Lipschitz maps defined on Asplund spaces. The material presented here is fairly standard and has been thoroughly treated in \cite{BenyaminiLindenstrauss2000} or \cite{AlbiacKalton2016}. Our swift presentation might omit some details at times but the reader can consult the two references just mentioned should the need be felt.

	\section{The Radon-Nikod\'{y}m property}
	
	In this section, we shall need some basic elements of the theory of the Bochner measurable and integrable functions. The reader can find the necessary background in Appendix \ref{appendix:Bochner} as well as pointers to the literature. 
	
	The Radon-Nikod\'{y}m property of a Banach space $X$ is most naturally defined in terms of the validity of the Radon-Nikod\'{y}m theorem for $X$-valued vector measures. Over the years, many characterizations have been discovered. Some of the most important are in terms of representability of operators from $L_1$ to $X$, convergence of bounded $X$-valued martingales, or dentability of bounded sets in $X$. We refer the reader to the monograph by Diestel and Uhl \cite{DiestelUhl1977} for a beautiful presentation of this subject. Our choice in this chapter is to adopt as our definition for the Radon-Nikod\'{y}m property its characterization in terms of the differentiability of $X$-valued Lipschitz functions defined on $\bR$, which is best suited for our purpose. Then, we shall have to work to exhibit nontrivial examples of spaces with the Radon-Nikod\'{y}m property and to derive infinite-dimensional differentiability results from our definition, which we now give.
	
	\begin{defi}
		\label{def:RNP}
		A real Banach space $X$ is said to have the \emph{Radon-Nikod\'{y}m property} (in short RNP) if any Lipschitz map from $\bR$ to $X$ is differentiable almost everywhere.
	\end{defi}
	
	The first example is given by Lebesgue's classical result on the almost everywhere differentiability of Lipschitz maps from $\bR$ to $\bR$ (see \cite{Lebesgue} for the original reference).
	
	\begin{theo}
		\label{thm:Lebesgue} 
		Let $f\colon \bR \to \bR$ be a Lipschitz function. Then, $f$ is almost everywhere differentiable and for all $x\in \bR$,
		\begin{equation*}
			f(x)=f(0)+\int_0^xf'(t)\,dt.
		\end{equation*}
		In particular, $(\bR,\abs{\cdot})$ has the Radon-Nikod\'{y}m property.
	\end{theo}
	
	The following corollary is an immediate consequence.
	
	\begin{coro}
		\label{cor:fdRNP}
		Every finite-dimensional Banach space has the Radon-Nikod\'{y}m property.
	\end{coro}
	
	The next proposition is elementary.
	
	\begin{prop}
		\label{pro:basic-RNP}\
		\begin{enumerate}[(i)]
			\item The Radon-Nikod\'{y}m property is stable under linear isomorphisms.
			\item If $Y$ is a closed subspace of a Banach space with RNP, then $Y$ has RNP.
			\item A Banach space $X$ has RNP if and only if all its closed separable subspaces have RNP.
		\end{enumerate}
	\end{prop}
	
	\begin{proof} 
		Statements $(i)$ and $(ii)$ are clear. So, assume that all the closed separable subspaces of a Banach space $X$ have RNP and let $f\colon\bR \to X$ be a Lipschitz map. Denote by $Y$ the closed linear span of $f(\bR)$. Then, $Y$ is a closed separable subspace of $X$. By assumption, $f\colon\bR \to Y$ is differentiable almost everywhere. This finishes the proof of the nontrivial implication in $(iii)$.
	\end{proof}
	
	\begin{rema}
		\begin{enumerate}
			\item When a Banach space property holds if and only if it holds for all closed separable subspaces, we say that the property is separably determined. Thus, the Radon-Nikod\'{y}m property is separably determined.
			\item The Radon-Nikod\'{y}m property does not pass to quotients since we will see shortly that $\ell_1$ has RNP, while $c_0$ fails RNP and, like every separable Banach space, is isometric to a quotient of $\ell_1$.
		\end{enumerate}
	\end{rema}
	The next statement, due to Dunford and Pettis \cite{DunfordPettis1940}, will provide many examples of infinite-dimensional spaces with the Radon-Nikod\'{y}m property.
	
	\begin{theo}
		\label{thm:dual-RNP}
		Let $X$ be a Banach space such that $X^*$ is separable. Then, $X^*$ has the Radon-Nikod\'{y}m property.
	\end{theo}
	
	\begin{proof} 
		Let $f\colon\bR \to X^*$ be a Lipschitz map. For $x\in X$ and $t\in \bR$, we define $f_x(t):=\langle x,f(t)\rangle$. Since $f_x$ is $\norm{x}\Lip(f)$-Lipschitz, it follows from Theorem \ref{thm:Lebesgue} that $f_x$ is differentiable almost everywhere. Since $X^*$ is separable, $X$ is separable. So, let us pick a dense sequence $(x_n)_{n=1}^\infty$ in $X$. Then, there is a measurable subset $N$ of $\bR$ of Lebesgue measure 0 so that for every $t\in \bR \setminus N$ and every $n\in \bN$, $f_{x_n}$ is differentiable at $t$. We now fix $t\in \bR \setminus N$ and $x\in X$. There is an increasing sequence $(n_k)_{k=1}^\infty$ in $\bN$ such that $\norm{x_{n_k}-x}\to 0$. Since $f$ is Lipschitz, the sequence $\big(f_{x_{n_k}}'(t)\big)_{k=1}^\infty$ is bounded in $\bR$, so, passing to a further subsequence, we may assume that it converges to some $l$ in $\bR$. It is now easy to check, using the Lipschitz property of $f$, that $f_x'(t)=l$. We have shown that for every $t\in \bR \setminus N$ and every $x\in X$, $f_x$ is differentiable at $t$. It also follows from the fact that $f$ is Lipschitz that the linear map $x\mapsto f'_x(t)$ is continuous on $X$. Therefore, for any $t\in \bR \setminus N$, there exists $g(t) \in X^*$ such that for all $x\in X$,
		\begin{equation*}
			f_x'(t)=\langle x,g(t)\rangle.
		\end{equation*}
		Let us extend $g$ as we may by setting $g(t)=0$, for $t\in N$. Everything until now was a routine argument, and we have not yet fully used the assumption that $X^*$ is separable. We now need to show that $g(t)=f'(t)$ almost everywhere. The key step is to prove that $g$ is measurable. Let $(y_k)_{k=1}^\infty$ be a dense sequence in $S_X$ and $x^* \in X^*$. Then, for all $t\in \bR \setminus N$,
		$$\norm{g(t)-x^*} = \sup_{k\in \bN}\langle y_k,g(t)-x^*\rangle = \sup_{k\in \bN}(f_{y_k}'(t)-\langle y_k,x^*\rangle).$$
		Since for all $k\in \bN$, the real-valued function vanishing on $N$ and equal to $f'_{y_k}$ on $\bR \setminus N$ is measurable, we deduce that $\norm{g-x^*}$ is measurable. It follows that for any open ball $B$ in $X^*$, $g^{-1}(B)$ is measurable. Finally, the separability of $X^*$ implies that any open set $U$ in $X^*$ is a countable union of open balls, so $g^{-1}(U)$ is measurable. So, we can apply Proposition \ref{prop:measurable} to deduce that $g$ is an $X^*$-measurable function.
		
		Note that $\norm{g}_{X^*}\le \Lip(f)$ on $\bR$. In particular, $g$ is an $X^*$-valued locally integrable function defined on $\bR$. Now, for any $x\in X$, the function $f_x$ is Lipschitz and $\bR$-valued and therefore, by Theorem \ref{thm:Lebesgue}, satisfies:
		$$\forall t\in \bR\ \ \langle x,f(t)\rangle= f_x(t)=f_x(0)+\int_0^t f_x'(s)\,ds=\langle x,f(0)\rangle + \int_0^t \langle x,g(s)\rangle \,ds.$$
		It follows that for all $t\in \bR$,
		\begin{equation*}
			f(t)=f(0)+\int_0^t g(s)\,ds.
		\end{equation*}
		Let us emphasize the fact that proving the measurability of $g$ was crucial to justify the existence of the above Bochner integral. It is now an immediate consequence of Theorem \ref{thm:Lebesgue-point} that for almost every $t\in \bR$, $f$ is differentiable at $t$ and $f'(t)=g(t)$.
	\end{proof}
	
	As a consequence, we have:
	
	\begin{coro}
		\label{reflexiveRNP} 
		Every reflexive Banach space has the Radon-Nikod\'{y}m property.
	\end{coro}
	
	\begin{proof} 
		Let $X$ be a reflexive Banach space and $Y$ be a closed separable subspace of $X$. Then, $Y$ is reflexive and thus isometric to $Y^{**}$ which is a separable dual space. It follows from Theorem \ref{thm:dual-RNP} that $Y$ has RNP and then from Proposition \ref{pro:basic-RNP} that $X$ has RNP.
	\end{proof}
	
	It is worth mentioning that there is a precise converse to Theorem \ref{thm:dual-RNP}. 
	
	\begin{theo}
		\label{thm:dual-RNP+} 
		Let $X$ be a Banach space. Then, $X^*$ has RNP if and only if every separable subspace of $X$ has a separable dual.
	\end{theo}
	
	\begin{proof} One direction is straightforward. Indeed, any separable subspace of $X^*$ isometrically embeds into $Y^*$, for some separable subspace $Y$ of $X$. So, if every separable subspace of $X$ has a separable dual, then every separable subspace of $X^*$ has RNP and hence so does $X^*$ by Proposition \ref{pro:basic-RNP}. 
		
		We refer again the reader to \cite{DiestelUhl1977} for the proof of the difficult direction and for a complete list of associated references.
	\end{proof}
	
	Let us now list classical examples and counterexamples.
	
	\begin{prop}
		\label{prop:examples-RNP}\,
		\begin{enumerate}[(i)]
			\item Let $(\Omega,\Sigma,\mu)$ be a measure space. For any $p\in (1,\infty)$, $L_p(\mu)$ has RNP.
			\item The space $\ell_1$ has RNP.
			\item The spaces $\co$, $L_1$, $C([0,1])$, $\ell_\infty$ and $L_\infty$ do not have RNP.
		\end{enumerate}
	\end{prop}
	
	\begin{proof} 
		Statement $(i)$ is an application of Corollary \ref{reflexiveRNP}, as $L_p(\mu)$ is reflexive if $p\in (1,\infty)$ and $(ii)$ an application of Theorem \ref{thm:dual-RNP} as $\ell_1$ is separable and isometric to the dual of $\co$.
		
		In order to prove $(iii)$, one can check (see Exercise \ref{ex:no-RNP}) that both maps $f\colon \bR \to c_0$ and $g\colon \bR \to L_1(\bR)$ defined by
		$$f(t):=\Big(\frac{\sin(nt)}{n}\Big)_{n=1}^\infty\ \ \text{and}\ \ g(t):=\car_{[0,|t|]},\ \ \text{for}\ t\in \bR$$
		are Lipschitz and nowhere differentiable. So, $\co$ and $L_1$ fail RNP. The other counterexamples follow as they all contain subspaces isometric to $\co$ and $L_1$.
	\end{proof}
	
	\begin{rema} 
		The first example of a separable Banach space with RNP that does not linearly embed into a separable dual was provided by McCartney and O'Brien in \cite{McCartneyOBrien1980}. 
	\end{rema}

	\section{Finite-dimensional Rademacher theorems}
	
	\subsection{Gateaux and Fr\'echet Differentiability}
	
	The first aim of this section is to show that the proof of the classical differentiability theorem of Rademacher for Lipschitz maps between finite-dimensional Banach spaces can be easily adapted to the setting of Lipschitz maps from $\bR^n$ to a Banach space with the Radon-Nikod\'{y}m property. For $n\in \bN$, let us denote by $\lambda_n$ the Lebesgue measure on $\bR^n$, $\Sigma_n$ the completion of the Borel sigma-algebra on $\bR^n$ and $B_n(u,r)$ the closed euclidean ball of center $u$ and radius $r$ in $\bR^n$.
	
	Recall that a map $f$ from an open set $U$ of a normed vector space $X$ to a normed vector space $Y$ is said to be \emph{Gateaux differentiable} at $x\in U$ if there exists $T \in B(X,Y)$  such that for all $u \in X$,
	\begin{equation}
		\label{eq:G-diff}
		\lim_{t \to 0}\frac{f(x+tu)-f(x)}{t}=T(u).
	\end{equation}
	
	We say that $f$ is \emph{Fr\'echet differentiable} at $x\in U$ if there exists $T \in B(X,Y)$  such that
	\begin{equation}
		\label{eq:F-diff}
		\lim_{\norm{h}_X \to 0}\frac{\norm{f(x+h)-f(x) - T(h)}_Y}{\norm{h}_X}=0.
	\end{equation}
	
	The crucial difference between \eqref{eq:G-diff} and \eqref{eq:F-diff} is that the first limit asserts the existence of directional derivatives for every direction, while in the second limit, the convergence is uniform over all directions. It is clear that being Fr\'echet differentiable is stronger than being Gateaux differentiable.
	If an operator satisfying \eqref{eq:G-diff} exists, then it is unique. We denote this operator by $D_f(x)$ and call it the \emph{Gateaux derivative} of $f$ at the point $x$. If \eqref{eq:F-diff} holds, then we call it the \emph{Fr\'echet derivative} of $f$ at the point $x$.
	
	\smallskip 
	The next theorem is an extension of Rademacher's differentiability theorem to RNP-target spaces. Its proof is identical to the usual proof fo the classical Rademacher Theorem, where the target space is finite-dimensional. 
	
	\begin{theo}
		\label{thm:Rademacher} 
		Let $U$ be an open subset of $\bR^n$, $Y$ be a Banach space with the Radon-Nikod\'{y}m property and $f\colon U\to Y$ be a Lipschitz map. Then, $f$ is almost everywhere Fr\'echet differentiable.
	\end{theo}
	
	Since $\bR^n$ is finite-dimensional and $f$ is Lipschitz, a classical compactness argument shows that $f$ is Fr\'echet differentiable at $x$ if and only if $f$ is Gateaux differentiable at $x$ (see Exercise \ref{ex:Fréchet=Gâteaux}). So, using the separability of $\bR^n$, we only need to show that $f$ is Gateaux differentiable at almost every interior point of any closed ball included in $U$. Now, a Lipschitz map defined on a closed ball of $\bR^n$ extends, via radial extension, to a Lipschitz map defined on $\bR^n$, so we may assume that $U=\bR^n$. We start with an easy lemma. Its proof is left as an exercise (see Exercise \ref{ex:additive}). 
	
	\begin{lemm}
		\label{lem:additive}
		Let $f\colon\bR^n \to Y$ be a Lipschitz map and $x\in \bR^n$. For $u\in \bR^n$, define; if it exists,
		\begin{equation*}
			h_u(x):=\lim_{t\to 0}\frac{f(x+tu)-f(x)}{t}.
		\end{equation*}
		Then, $f$ is Gateaux differentiable at $x$ if and only if for all $u\in \bQ^n$, $h_u(x)$ exists and the map $u\mapsto h_u(x)$ is additive on $\bQ^n$.
	\end{lemm}
	
	\begin{proof}[Proof of Theorem \ref{thm:Rademacher}] 
		The first step is to show that for any $u\in \bR^n$, $h_u(x)$ exists almost everywhere (note that, since $\Sigma_n$ is complete, it implies that the set $\{x\in \R^n:\ h_u(x)\ \text{exists}\}$ is measurable). So, fix $u\in \bR^n$ and assume on the contrary that there exists $A\in \Sigma_n$ such that $\lambda_n(A)>0$ and for all $x\in A$, $h_u(x)$ does not exist. Then, $u\neq 0$ and for $\xi\in u^\perp$, set  $A_\xi:=\{s\in \bR\colon \xi+su \in A\}$. It follows from Fubini's theorem that there exists $\xi\in u^\perp$ such that $\lambda_1(A_\xi)>0$. Then, for any $s\in A_\xi$, the map $t\mapsto f(\xi+tu)$ is not differentiable at $s$, which is a contradiction since $Y$ has RNP.
		%Since $\bQ^n$ is countable, we deduce immediately that there exists $B\in \Sigma_n$, so that $\lambda_n(\bR^n \setminus B)=0$ and for all $x\in B$ and all $u\in \bQ^n$, $h_u(x)$ exists. It remains to show that for almost every $x\in B$, the map $u\mapsto h_(x)$ is additive on $\bQ^n$.\\
		Consider now $\varphi\colon \bR^n\to [0,\infty)$ a $C^\infty$ map with support in $B_n(0,1)$ so that $\int_{\bR^n}\varphi\,d\lambda_n=1$ and for $k\in \bN$ let $\varphi_k(x) := k^n\varphi(kx)$. Now, for any $k\in \bN$ and $x\in \bR^n$ let
		\begin{equation*}
			g_k(x) := (f\ast \varphi_k)(x)=\int_{\bR^n} f(y)\varphi_k(x-y)\,d\lambda_n(y).
		\end{equation*}
		The map $g_k$ is $C^\infty$ on $\bR^n$. On the other hand, we can write 
		$$g_k(x)=\int_{\bR^n} \varphi_k(y)f(x-y)\,d\lambda_n(y).$$
		A simple application of the Lebesgue dominated convergence theorem, combined with the first step of our proof and the fact that $h_u \in L_\infty(\R^n)$, yields that for all $x\in \bR^n$ and all $u \in \bR^n$,
		\begin{equation*}
			Dg_k(x).u=(\varphi_k \ast h_u)(x).
		\end{equation*} 
		In particular, since $Dg_k(x)$ is linear, we get that for all $u,v \in \bR^n$,
		\begin{equation*}
			\varphi_k \ast (h_{u+v}-h_u-h_v)=0.
		\end{equation*}
		Fix now $u,v \in \bR^n$ and let $h:=h_{u+v}-h_u-h_v$. Note that $h$ is well defined almost everywhere and bounded. Since $\bR^n$ is separable and $f$ is Lipschitz, we may assume that $Y$ is separable. Moreover, it is easy to see that $h$ is weakly measurable. Thus, it follows from Proposition \ref{prop:Pettis} that $h$ is measurable. So, $h$ is locally Bochner integrable on $\bR^n$. We can therefore apply Theorem \ref{thm:Lebesgue-point} to deduce that for almost every $x\in \bR^n$,
		$$\lim_{r\to 0^+}\frac{1}{\lambda_n(B_n(x,r))}
		\int_{B_n(x,r)}\|h(y)-h(x)\|\,d\lambda_n(y)=0.$$
		Then, usual estimates on convolution by approximations of unity imply that for almost every $x\in \bR^n$, $\lim_{k\to \infty}(\varphi_k \ast h)(x)=h(x)$. It follows that $h=0$ almost everywhere.
		Finally, we use the countability of $\bQ^n$ to deduce that the assumptions of Lemma \ref{lem:additive} are satisfied almost everywhere and therefore that $f$ is Gateaux differentiable almost everywhere. 
	\end{proof}
	
	\begin{rema}
		\label{rem:differential-embedding}
		It is easy to see that if $f\colon X\to Y$ is a bi-Lipschitz embedding between the Banach spaces $X$ and $Y$ and if $f$ is Gateaux differentiable at $x\in X$, then $Df(x)$ is a linear embedding from $X$ into $Y$ with $\norm{Df(x)}\le \Lip(f)$ and $\|Df(x)^{-1}\|\le \Lip(f^{-1})$. 
	\end{rema}
	
	\subsection{\texorpdfstring{Weak$^*$ Gateaux differentiability}{Weak Gateaux differentiability}}

	We now turn to results on the weak$^*$ Gateaux differentiability of Lipschitz functions with values in a dual of a separable Banach space.
	
	\begin{defi}
		Let $U$ be an open subset of a Banach space $X$, $Y:=Z^*$ a dual Banach space, $f\colon U\to Z^*$ and $x\in U$. Then, we say that $f$ is \emph{weak$^*$ Gateaux differentiable at $x$} if there exists $T\in B(X,Z^*)$ such that for all $u\in X$ and all $z\in Z$,
		\begin{equation*}
			\lim_{t\to 0}\frac{\langle z,f(x+tu)-f(x)\rangle}{t}=\langle z,T(u)\rangle.
		\end{equation*}
		In that case $T$ is called the \emph{weak$^*$ Gateaux derivative} of $f$ at $x$ and denoted by $D_f^*(x)$.
	\end{defi}
	
	We can now state and prove a weak$^*$ analog of the Rademacher theorem. The main point here is that the dual Banach space does not have to be separable: if it is, then we would have a genuine Gateaux derivative, but if not, then we at least have a weak$^*$ Gateaux derivative. 
	
	\begin{theo}
		\label{thm:weak*-Rad}
		Let $U$ be an open subset of $\bR^n$ and $f\colon U\to Z^*$ be a Lipschitz map, where $Z$ is a \underline{separable} Banach space. Then, $f$ is almost everywhere weak$^*$ Gateaux differentiable.\\
		Moreover, if $f$ is a bi-Lipschitz embedding, then for almost every $x$ in $U$, $D^*_f(x)$ exists and is a linear embedding satisfying $\norm{D^*_f(x)}\le \Lip(f)$ and $\|D^*_f(x)^{-1}\|\le \Lip(f^{-1})$.
	\end{theo}
	
	\begin{proof} 
		Let $(z_k)_{k=1}^\infty$ be a dense sequence in $Z$. It follows from Theorem \ref{thm:Rademacher}  that the set $\Omega$ of all $x$ in $U$ such that for all $k\in \bN$, the map $y\mapsto \langle f(y),z_k\rangle$ is Gâteaux (and Fréchet) differentiable at $x$ is measurable and has full Lebesgue measure in $U$. Using the density of $(z_k)_{k=1}^\infty$ and the Lipschitz property of $f$, we easily see that $\Omega$ is the set of points of weak$^*$ Gateaux differentiability of $f$. It then follows from the weak$^*$ lower semi-continuity of $\norm{\cdot}_{Z^*}$ that for any $x\in \Omega$, $\norm{D^*f(x)}\le \Lip(f)$. 
		
		Assuming now that $f$ is a bi-Lipschitz embedding, we have to find points $x\in \Omega$ such that $D^*f(x)$ is a linear embedding. Assume as we may that $\Lip(f^{-1})=1$ and $\Lip(f)=b\ge 1$. Since the unit sphere of $\bR^n$ is separable, it is enough to show that for any $u$ in the unit sphere of $\bR^n$, the set $\{x\in \Omega\colon \norm{D^*f(x)(u)}\ge 1\}$, which is easily seen to be measurable,  is of full Lebesgue measure in $U$. So, assume this is not the case for some $u$ of norm one. Then, there exists $\delta >0$ such that the set $N:=\{x\in \Omega\colon \norm{D^*f(x)(u)}\le 1-\delta\}$ has positive Lebesgue measure. So, by Fubini's theorem, there exists $x_0\in u^\perp$ so that $\lambda_1(I)>0$ where $I:=\{t\in \bR\colon x_0+tu \in N\}$. It follows from Theorem \ref{thm:Lebesgue-point} applied to the indicator function of $I$ that there exist $x\in U$ and $h>0$ such that $\lambda_1(A)>h(1-\frac{\delta}{2b})$ where $A:=\{t\in (0,h)\colon x+tu \in N\}$. 
		
		Recall now that $\norm{f(x+hu)-f(x)}\ge h$. So, there exists $z\in S_Z$ such that
		$$\varphi(h)-\varphi(0)> h\Big(1-\frac{\delta}{2}\Big),$$
		where $\varphi(t):=\langle f(x+tu),z\rangle$. Note that the map $\varphi$ is $b$-Lipschitz and that $\varphi'(t)\le 1-\delta$ for all $t\in A$. It follows that
		\begin{align*}
			\varphi(h)-\varphi(0) & = \int_{[0,h]} \varphi'(t)\,dt=\int_A \varphi'(t)\,dt + \int_{[0,h]\setminus A} \varphi'(t)\,dt \\
			& \le h(1-\delta)+h\frac{\delta}{2}=h\Big(1-\frac{\delta}{2}\Big).
		\end{align*}
		This contradiction concludes the proof.
		
		%It follows from Theorem \ref{thm:Lebesgue-point} applied to the indicator function of $I$, that there exist a Lebesgue point $t_0$ and $h>0$ small enough such that $\lambda_1(I\cap(t_0-h/2,t_0+h/2))>h(1-\frac{\delta}{2b})$. Recall now that $\norm{f(x_0+(t_0+h/2)u)-f(x_0+(t_0-h/2)u)}\ge h$. So, there exists $z\in S_Z$ such that $$\varphi(t_0+h/2)-\varphi(t_0-h/2)> h\Big(1-\frac{\delta}{2}\Big),$$ where $\varphi(t):=\langle f(x_0+tu),z\rangle$. Note that the map $\varphi$ is $b$-Lipschitz and that $\varphi'(t)\le 1-\delta$ for all $t\in I\cap(t_0-h/2,t_0+h/2)$. It follows that 
		%\begin{align*}
		%		\varphi(t_0+h/2) - \varphi(t_0-h/2) & = \int_{[t_0-h/2,t_0+h/2]} \varphi'(t)\,dt=\int_{[t_0-h/2,t_0+h/2]\cap I} \varphi'(t)\,dt + \int_{[t_0-h/2,t_0+h/2]\setminus I} \varphi'(t)\,dt \\
		%                      	   & \le h(1-\delta)+h\frac{\delta}{2}=h\Big(1-\frac{\delta}{2}\Big).
		%\end{align*}
		%This contradiction concludes the proof.
	\end{proof}

	\section{Infinite-dimensional Rademacher theorems}
	\label{sec:inf-dim-Rademacher}
	
	The main issue encountered in trying to prove Rademacher-type theorems when the domain space is an infinite-dimensional Banach space stems from the absence of a Lebesgue measure on separable infinite-dimensional Banach spaces. Indeed, in Exercise \ref{ex:no-Lebesgue-measure} you are asked to show that any translation-invariant Borel measure on an infinite-dimensional separable Banach space must be either infinite for all sets or zero for all sets! Therefore, one needs to come up with a notion of being small for Borel subsets which behaves similarly to being small in the Lebesgue sense.
	
	\subsection{Haar null sets}
	
	Throughout this section, $X$ will be a \underline{separable} Banach space and $\mathrm{Bor}(X)$ will denote the $\sigma$-algebra of Borel subsets of $X$.
	
	\begin{defi} 
		Let $X$ be a separable Banach space and $A\in \mathrm{Bor}(X)$. We say that $A$ is \emph{Haar null} if there exists a probability measure $\lambda$ on $(X,\mathrm{Bor}(X))$ such that $\lambda(x+A)=0$ for every $x\in X$.
	\end{defi}
	
	It will be convenient to use the following equivalent formulations.
	
	\begin{prop}
		\label{prop:equivalent-Haar-null}
		Let $A \in \mathrm{Bor}(X)$ and $\eps>0$. The following assertions are equivalent.
		\begin{enumerate}[(i)]
			\item The Borel set $A$ is Haar null.
			\item There exists a nonzero $\sigma$-finite measure $\mu$ on $(X,\mathrm{Bor}(X))$ such that $\mu(x+A)=0$ for every $x\in X$.
			\item There exists a probability space $(\Omega,\Sigma,\mu)$ and a measurable function $\xi\colon (\Omega,\Sigma)\to (X,\mathrm{Bor}(X))$ such that:
			\begin{itemize}
				\item $\norm{\xi(\omega)}\le \eps$ for all $\omega \in \Omega$,
				\item $\int_\Omega \car_A(x+\xi(\omega))\,d\mu(\omega)=0$ for all $x\in X$.    
			\end{itemize}		
		\end{enumerate}
	\end{prop}
	
	\begin{proof} 
		Assertion $(i)$ trivially implies $(ii)$ and it is clear that $(iii)$ implies $(i)$. Indeed, assume that $(\Omega,\Sigma,\mu)$ and $\xi$ satisfy the conditions of $(iii)$. Then, for $B\in \mathrm{Bor}(X)$, set $\lambda(B) := \mu(\xi^{-1}(B))$. The probability measure $\lambda$ on $(X,\mathrm{Bor}(X))$ testifies that $A$ is Haar null.
		
		It remains to show that $(ii)$ implies $(iii)$. So, assume that $\mu$ is a nonzero $\sigma$-finite measure on $(X,\mathrm{Bor}(X))$ such that $\mu(x+A)=0$ for every $x\in X$. We write $X=\bigcup_{k=1}^\infty X_k$ with $\mu(X_k)<\infty$. Fix $(x_n)_{n=1}^\infty$ a dense sequence in $X$. For $n\in \bN$, let $A_{n,k}^\eps := B_X(x_n,\eps) \cap X_k$ (here the ball is closed). Elementary properties of measures imply that there exist $n,k\in \bN$ so that $0<\mu(A_{n,k}^\eps)<\infty$. For $E\in \mathrm{Bor}(X)$, define
		\begin{equation*}
			\nu(E) := \frac{\mu(E\cap A_{n,k}^\eps)}{\mu(A_{n,k}^\eps)}.
		\end{equation*}
		Then, $\nu$ is a probability measure on $(X,\mathrm{Bor}(X))$, $\xi(\omega) =(\omega-x_n)\car_{A_{n,k}^\eps}(\omega)$ is a measurable function from $(X,\mathrm{Bor}(X))$ to $(X,\mathrm{Bor}(X))$, $\norm{\xi}_\infty\le \eps$ and for all $x\in X$:
		\begin{align*}
			\int_X \car_A(x+\xi(\omega))\,d\nu(\omega)
			&=\frac{1}{\mu(A_{n,k}^\eps)}\int_{A_{n,k}^\eps}\car_A(x-x_n+\omega)\,d\mu(\omega)\\
			&\le \frac{1}{\mu(A_{n,k}^\eps)}\int_{X}\car_A(x-x_n+\omega)\,d\mu(\omega)\\
			&=\frac{1}{\mu(A_{n,k}^\eps)}\mu(x_n-x+A)=0.
		\end{align*}
	\end{proof}
	
	There are many sets that are not Haar null.
	
	\begin{prop}
		A nonempty open subset of a separable Banach space $X$ is never Haar null. 
	\end{prop} 
	
	\begin{proof}
		Assume that there exists a nonempty open subset of $X$ which is Haar null. Then, there exists $x\in X$ and $r>0$ such that $B_X(x,r)$ is Haar null. So, let $\lambda$ be a probability measure on $(X,\mathrm{Bor}(X))$ such that $\lambda(y+B_X(x,r))=0$ for every $y\in X$. It follows that $\lambda(B_X(z,r))=0$ for all $z\in X$. Since $X$ is separable, there exists a sequence $(z_n)_{n=1}^\infty$ in $X$ such that $\bigcup_{n=1}^\infty B(z_n,r)=X$. This implies that $\lambda(X)=0$; a contradiction.
	\end{proof}
	
	Next, we verify that Haar null sets coincide with Lebesgue null sets in finite-dimensional spaces. 
	
	\begin{lemm}
		Let $X$ be a separable Banach space.
		\begin{enumerate}[(i)]
			\item  Let $A\in \mathrm{Bor}(X)$ and assume that there exists a nontrivial finite-dimensional subspace $E$ of $X$ such that for any $x\in X$, $(A+x)\cap E$ is Lebesgue null, then $A$ is Haar null.
			\item  Let $A\in \mathrm{Bor}(\bR^n)$ where $n\in \bN$. Then, $A$ is Lebesgue null if and only if $A$ is Haar null.
		\end{enumerate}
	\end{lemm}
	
	\begin{proof} 
		$(i)$ Let $\lambda_E$ be a Lebesgue measure on $E$. For $B\in \mathrm{Bor}(X)$, define $\mu(B) :=\lambda_E(B\cap E)$. Then, $\mu$ clearly satisfies $(ii)$ of Proposition \ref{prop:equivalent-Haar-null}. 
		
		$(ii)$ Assume $A\in \mathrm{Bor}(\bR^n)$ is Lebesgue null, then for all $x\in \bR^n$, $A+x$ is Lebesgue null and we can apply $(i)$ to deduce that $A$ is Haar null.\\
		Assume now that $A\in \mathrm{Bor}(\bR^n)$ is Haar null and let $\lambda$ be a probability measure on $(\bR^n,\mathrm{Bor}(\bR^n))$ such that $\lambda(A+x)=0$ for all $x\in X$. By the translation invariance of the Lebesgue measure $\lambda_n$ on $\bR^n$, we get that for all $\omega \in \bR^n$, 
		\begin{equation*}
			\lambda_n(A) = \lambda_n(A-w) = \int_{\bR^n}\car_A(x+\omega)\,d\lambda_n(x).
		\end{equation*}
		Integrating this equality with respect to the probability measure $\lambda$ and using Fubini's theorem, we deduce
		$$\lambda_n(A)=\int_{\bR^n}\Big(\int_{\bR^n} \car_A(x+\omega)\,d\lambda(\omega)\Big)\,d\lambda_n(x)=\int_{\bR^n}\lambda(A-x)\,d\lambda_n(x)=0.$$
	\end{proof}
	
	Finally, we prove a crucial fact: the property of being Haar null is stable under countable unions. 
	
	\begin{prop}
		\label{prop:countable-union-Haar-null}
		Assume that $(A_n)_{n=1}^\infty$ is a sequence of Haar null subsets of $X$. Then, $\bigcup_{n=1}^\infty  A_n$ is Haar null.
	\end{prop}
	
	\begin{proof} 
		Here, we shall use $(iii)$ from Proposition \ref{prop:equivalent-Haar-null}. So, for any $n\in \bN$, there exist a probability space $(\Omega_n,\Sigma_n,\mu_n)$ and a measurable function $\xi_n:(\Omega_n,\Sigma_n)\to (X,\mathrm{Bor}(X))$ such that $\norm{\xi_n}_\infty\le 2^{-n}$ and for all $x\in X$, 
		\begin{equation*}
			\int_{\Omega_n} \car_{A_n}(x+\xi_n(\omega_n))\,d\mu_n(\omega_n)=0.
		\end{equation*}
		Consider now the product probability $\mu := \otimes_{n=1}^\infty \mu_n$ on $\Omega := \prod_{n=1}^\infty \Omega_n$ and $\xi(\omega) := \sum_{n=1}^\infty \xi_n(\omega_n)$, for $\omega=(\omega_n)_{n=1}^\infty \in \Omega$. Then, an easy application of Fubini's theorem yields that for all $x\in X$ and all $n\in \bN$,
		\begin{equation*}
			\int_{\Omega} \car_{A_n}(x+\xi(\omega))\,d\mu(\omega)=0.
		\end{equation*}
		Therefore, for all $x\in X$,
		\begin{equation*}
			\int_{\Omega} \sum_{n=1}^\infty\car_{A_n}(x+\xi(\omega))\,d\mu(\omega)=0.
		\end{equation*}
		Since $0\le \car_{\cup_{n=1}^\infty  A_n} \le \sum_{n=1}^\infty\car_{A_n}$, we deduce that for all $x\in X$,
		\begin{equation*}
			\int_{\Omega} \car_{\cup_{n=1}^\infty A_n}(x+\xi(\omega))\,d\mu(\omega)=0.
		\end{equation*}
		This finishes our proof.
	\end{proof}
	
	\subsection{Gateaux differentiability results}
	
	We are now ready to prove an infinite-dimensional version of the Rademacher theorem. This result is due to Aronszajn \cite{Aronszajn}, Christensen \cite{Christensen} and Mankiewicz \cite{Mankiewicz} independently. 
	
	\begin{theo}
		\label{thm:infinite-Rad}
		Let $X$ be a separable Banach space, $Y$ a Banach space with the Radon-Nikod\'{y}m property, and assume that $f\colon X\to Y$ is a Lipschitz map. Then, the set of points in $X$ at which $f$ is not Gateaux differentiable is Haar null.
	\end{theo}
	
	\begin{proof} 
		By separability of $X$ one can pick $(E_n)_{n=1}^\infty$ an increasing sequence of finite-dimensional subspaces of $X$ such that $\bigcup_{n=1}^\infty E_n$ is dense in $X$. Then, let us denote by $G_n$ the set of all $x$ in $X$ such that the map $v\mapsto f(x+v)$ is Gateaux differentiable at $0$ when considered as a map from $E_n$ to $Y$ and observe that $f$ is Gateaux differentiable at every $x\in\bigcap_{n=1}^\infty G_n$. Indeed, for any $n\in \bN$, there exists $T_n \in B(E_n,Y)$ such that for all $v\in E_n$,
		\begin{equation*}
			\lim_{t\to 0}\frac{f(x+tv)-f(x)}{t}=T_n(v).
		\end{equation*}
		The operator $T_{n+1}$ clearly extends $T_n$ and hence we can define for $v\in \bigcup_{n=1}^\infty E_n$, $T(v) := T_n(v)$ for all $n$ such that $v\in E_n$. Since $f$ is Lipschitz, $T$ is bounded with $\norm{T}\le \Lip(f)$ and extends to a bounded operator from $X$ into $Y$ that we still denote by $T$. It is not difficult to check, using the density of $\bigcup_{n=1}^\infty E_n$ in $X$ and the Lipschitz property of $f$, that $T$ is the Gateaux derivative of $f$ at $x$. On the other hand, $f$ is clearly not Gateaux differentiable at a point in $X\setminus(\bigcap_{n=1}^\infty G_n)$. To finish the proof it remains to show that $ X\setminus(\bigcap_{n=1}^\infty G_n) = \bigcup_{n=1}^\infty (X\setminus G_n)$ is Haar null. In order to see that, observe that it follows from Lemma \ref{lem:additive} that $A_n:=X\setminus G_n\in \mathrm{Bor}(X)$. If we denote by $\lambda_{E_n}$ a Lebesgue measure on $E_n$, we deduce from Theorem \ref{thm:Rademacher} that for all $x\in X$, $\lambda_{E_n}((A_n + x)\cap E_n)=0$. Then, we apply $(ii)$ from Proposition \ref{prop:equivalent-Haar-null} to deduce that $A_n$ is Haar null and by Proposition \ref{prop:countable-union-Haar-null} that $\bigcup_{n=1}^\infty A_n$ is Haar null.
	\end{proof}
	
	As an immediate corollary, we have:
	
	\begin{coro}
		\label{thm:linearization-RNP} 
		Let $X$ be a separable Banach space and $Y$ be a Banach space with the Radon-Nikod\'{y}m property and assume that $f\colon X\to Y$ is a bi-Lipschitz embedding. Then, $X$ is linearly isomorphic to a subspace of $Y$.
	\end{coro}
	
	\begin{proof} 
		By Theorem \ref{thm:infinite-Rad}, $f$ admits a point of Gateaux differentiability $x_0\in X$. Then, by Remark \ref{rem:differential-embedding}, $Df(x_0)$ is a linear embedding with $\norm{Df(x_0)}\le \Lip(f)$ and $\|Df(x_0)^{-1}\|\le \Lip(f^{-1})$.
	\end{proof}
	
	For instance, it is well known that $L_1$ does not linearly embed into $\ell_1$. Corollary  \ref{thm:linearization-RNP} allows us to upgrade this linear result and to conclude that $L_1$ does not even admit a bi-Lipschitz embedding into $\ell_1$. This is just one example amongst many others of an extension to the nonlinear setting of a linear embeddability obstruction.

	\subsection{\texorpdfstring{Weak$^*$ Gateaux differentiability results}{Weak Gateaux differentiability results}}
	
	There is also an analog of the above results for the weak$^*$ Gateaux differentiability of Lipschitz functions with values in a dual space. The proof is a combination of the arguments in the proofs of Theorem \ref{thm:weak*-Rad} and Theorem \ref{thm:infinite-Rad}. The results of this subsection, as well as Theorem \ref{thm:weak*-Rad} are due to Heinrich and Mankiewicz \cite{HeinrichMankiewicz1982}. 
	
	\begin{theo}
		\label{thm:w*-linearization-dual} 
		Let $X$ and $Z$ be two separable Banach spaces and $f\colon X\to Z^*$ be a Lipschitz map. Then, the following hold 
		\begin{enumerate}[(i)]
			\item There exists a Haar null set $A$ in $X$ such that $f$ is weak$^*$ Gateaux differentiable at every point in $X\setminus A$.
			\item For every point $x$ of weak$^*$ Gateaux differentiability of $f$, $\|D^*_f(x)\|\le \Lip(f)$.
			\item If in addition, $f$ is a bi-Lipschitz embedding, then there exist a Haar null set $B$ containing $A$ and such that for every $x\in X\setminus B$, $D^*_f(x)$ is a linear embedding from $X$ into $Z^*$, satisfying $\|(D^*_f(x))^{-1}\|\le \Lip(f^{-1})$.
		\end{enumerate}
	\end{theo}
	
	\begin{proof} 
		Denote by $A$ the set of points in $X$ at which $f$ is not weak$^*$ Gateaux differentiable. Since $Z$ is separable and $f$ is Lipschitz, we easily get that $A\in \mathrm{Bor}(X)$. Then, following the lines of the proof of Theorem \ref{thm:infinite-Rad} and using the first part of Theorem \ref{thm:weak*-Rad} instead of Theorem \ref{thm:Rademacher}, we get that $A$ is Haar null. It follows immediately from the weak$^*$ lower semi-continuity of $\norm{\cdot}_{Z^*}$ that for any $x\in \Omega := X\setminus A$, $\norm{D_f^*(x)}\le \Lip(f)$.
		
		Assume now that $f$ is a bi-Lipschitz embedding with $\Lip(f^{-1})=a$. Then, using the separability of $Z$, we infer that $\Omega' := \{x\in \Omega\colon \forall u\in X\  \norm{D_f^*(x)(u)}\ge a\norm{u}\} \in \mathrm{Bor}(X)$. Again, following the proof of Theorem \ref{thm:infinite-Rad} in combination with the second part of Theorem \ref{thm:weak*-Rad}, we deduce that $B:= X\setminus \Omega'$ is Haar null.
	\end{proof}
	
	Next, it will be important to note that in the above statement on bi-Lipschitz embeddings, we do not need to assume that the target space is the dual of a separable Banach space. The proof relies on the following result by Heinrich and Mankiewicz \cite{HeinrichMankiewicz1982}. We shall omit the proof of this purely linear result and refer the reader to \cite{HeinrichMankiewicz1982} or to \cite{SimsYost1989} (see also \cite[Chapter 14]{AlbiacKalton2016}) for a proof avoiding the use of model theory.
	
	\begin{theo}
		\label{thm:model-HM} 
		Let $Z$ be a Banach space and $Y$ be a separable subspace of $Z^*$. Then, there exists a separable Banach space $E$ such that $Y$ is linearly isometric to a subspace of $E^*$ and $E^*$ is linearly isometric to a subspace of $Z^*$. Moreover, the isometric embedding of $E^*$ into $Z^*$ is a Hahn-Banach extension map. 
	\end{theo}
	
	Then we deduce
	\begin{theo}
		\label{thm:embed-in-dual+}
		Let $X$ be a separable Banach space and $Z$ be any Banach space. Assume that there exists a bi-Lipschitz embedding $f$ from $X$ into $Z^*$. Then, there exists a linear embedding $T:X\to Z^*$ so that $\norm{T}\,\|T^{-1}\|\le \dist(f) :=\Lip(f)\Lip(f^{-1})$.
		%In particular if a separable Banach space $X$ Lipschitz embeds into a Banach space $Y$, then $X$ linearly embeds into $Y^{**}$.
	\end{theo}
	
	\begin{proof} 
		Assume $f\colon X\to Z^*$ is a bi-Lipschitz embedding and let $Y$ be the closed linear span of $f(X)$ in $Z^*$. Since $X$ is separable and $f$ is Lipschitz, $Y$ is also separable. By Theorem \ref{thm:model-HM}, there exists a separable Banach space $E$ such that $Y$ embeds isometrically into $E^*$ and $E^*$ embeds isometrically into $Z^*$. Now, it follows from Theorem \ref{thm:w*-linearization-dual} that $X$ linearly embeds into $E^*$ and therefore into $Z^*$ with distortion at most $\dist(f)$.
	\end{proof}

	\subsection{Lipschitz rigidity: first results and general problem}
	
	Equipped with the differentiability results obtained in the previous sections, we can already prove some fundamental Lipschitz rigidity results.
	\begin{theo}
		\label{thm:linearization-applications}
		\
		\begin{enumerate}[(i)]
			\item If a separable Banach space $X$ admits a bi-Lipschitz embedding into a Banach space $Y$ with the Radon-Nikod\'{y}m property, then $X$ linearly embeds into $Y$.
			\item Bi-Lipschitz embeddings between Banach spaces preserve the Radon-Nikod\'{y}m property and reflexivity.
			\item If an infinite-dimensional Banach space $X$ admits a bi-Lipschitz embedding into $\ell_2$, then it is linearly isomorphic to $\ell_2$. In particular, if $X$ is Lipschitz equivalent to $\ell_2$, then it is automatically isomorphic to $\ell_2$.
			\item If a separable Banach space $X$ admits a bi-Lipschitz embedding into a Banach space $Y$, then $X$ linearly embeds into $Y^{**}$.
		\end{enumerate}
	\end{theo}
	
	\begin{proof} 
		Statement $(i)$ follows from Corollary \ref{thm:linearization-RNP}.
		
		Statement $(ii)$ is a direct consequence of the combination of $(i)$ together with the facts that the Radon-Nikod\'{y}m property and reflexivity pass to subspaces and are separably determined.
		
		$(iii)$ is also a consequence of $(i)$ since $\ell_2$ has RNP and any infinite-dimensional subspace of $\ell_2$ is isometric to $\ell_2$.
		
		Finally, $(iv)$ is a direct consequence of Theorem \ref{thm:embed-in-dual+}.
	\end{proof}
	
	Assertion $(ii)$ implies that reflexivity and the Radon-Nikod\'{y}m property are \emph{Lipschitz invariants}, meaning that they are preserved under Lipschitz equivalences. Lipschitz invariants are useful when trying to distinguish or classify Banach spaces in the Lipschitz category. Theorem \ref{thm:linearization-applications} becomes particularly helpful when one tries to rule out bi-Lipschitz embeddings into Banach spaces for which the structure of their subspaces is elusive.
	
	The phenomenon in assertion $(iii)$ is what we will call \emph{Lipschitz rigidity}. In layman's terms, it says that if a Banach space has a metric structure that looks similar to that of Hilbert space, then it must be isomorphic to a Hilbert space. Therefore, one can recover the linear geometry of a Hilbert space if one merely knows its Lipschitz structure and one says that Hilbert space is \emph{Lipschitz rigid} (or has a unique Lipschitz structure). Whether the same phenomenon holds for other Banach spaces is a central problem in the nonlinear geometry of Banach spaces.
	
	\begin{prob}
		\label{prob:Lipschitz-rigidity}
		Is every separable Banach space Lipschitz rigid?   
	\end{prob}
	
	Separability cannot be dropped in Problem \ref{prob:Lipschitz-rigidity}, and we will present nonseparable counterexamples in Chapter \ref{chapter:Counterexamples}.
	
	Another interesting application relies on the principle of local reflexivity. Recall that $\cdist{N}(M)$ denotes the $N$-distortion of $M$, and whether $\cdist{N}(M)$ is finite is clearly a Lipschitz invariant. Given Banach spaces $X$ and $Y$ we can define a linear analog as follows:  
	\begin{equation*}
		d_Y^{BM}(X):=\inf\{\norm{T}\|T^{-1}\|\colon T\colon X\to Y \textrm{ isomorphic embedding}\}.
	\end{equation*}
	The following theorem states that for a sequence of finite-dimensional spaces $(E_k)_{k=1}^\infty$ the isomorphic invariant $\sup_{k\ge 1}d_X^{BM}(E_k)<\infty$ can be reformulated in a purely metric fashion.
	
	\begin{theo}
		If $X$ is a Banach space and $(E_k)_{k=1}^\infty$ a sequence of finite-dimensional spaces, then 
		$\sup_{k\ge 1} \cdist{X}(E_k)<\infty$ if and only if $\sup_{k\ge 1}d_X^{BM}(E_k)<\infty$.
	\end{theo}
	
	\begin{proof}
		For the nontrivial implication, observe that it follows from Theorem \ref{thm:linearization-applications} $(iv)$ that $\sup_{k\ge 1}d_{X^{**}}^{BM}(E_k)<\infty$. Invoking the principle of local reflexivity we can conclude that $\sup_{k\ge 1}d_{X}^{BM}(E_k)<\infty$.
	\end{proof}

	\section{Extension of Lipschitz maps and Lipschitz retracts}
	
	Whether Lipschitz maps can be extended while controlling the Lipschitz constant of the extension is a fundamental problem in metric geometry. In this section, we gather some elementary facts on the possibility of extending Lipschitz maps, which will be useful to us in the sequel. The first elementary extension result relies on the classical but fundamental ``Inf-Convolution'' formula, also referred to as McShane-Whitney’s extension theorem. 
	
	\begin{prop}
		\label{pro:Inf-Conv}
		Let $(M,d)$ be a metric space, $A$ a nonempty subset of $M$ and $f\colon A \to \R$ be a $C$-Lipschitz map. Then, the function $F\colon M\to \R$ defined by:
		$$F(x):=\inf\{f(y)+Cd(x,y)\colon y\in A\}\ \text{for all}\ x\in M$$
		is $C$-Lipschitz and coincides with $f$ on $A$.
	\end{prop}
	
	\begin{proof} 
		Fix $y_0\in A$. Then, it follows from the triangle inequality that for all $x\in M$ and all $y\in A$:
		$$f(y)+Cd(x,y)\ge f(y_0)-Cd(y,y_0)+Cd(x,y) \ge f(y_0)-Cd(x,y_0).$$
		This shows that $F$ is well defined on $M$. It is clear that $F=f$ on $A$. \\
		Consider now $x,x' \in M$ and $y\in A$. Then, $$f(y)+Cd(x,y)-F(x')\ge f(y)+Cd(x,y)-f(y)-Cd(x',y)\ge -Cd(x,x').$$
		Taking the infimum over all $y\in A$ in the above inequality, we deduce that $F(x)-F(x')\ge -Cd(x,x')$ and by symmetry that $F$ is $C$-Lipschitz.
	\end{proof}
	
	\begin{rema}\, 
		\label{rem:inf-conv}
		\begin{enumerate}
			\item The formula $G(x):=\sup\{f(y)-Cd(x,y)\colon y\in A\}$ provides another $C$-Lipschitz extension of $f$. In fact, for any $C$-Lipschitz extension $H$ of $f$, we have that $G\le H \le F$ on $M$.
			\item The terminology inf-convolution comes from the fact that if the domain space is a normed vector space $(X,\norm{\cdot})$, then
			$F(x)=\inf\{f(y)+g(z)\colon y+z=x,\ y\in A, z\in X\}$ where $g(z)=C\norm{z}$ and this formula is reminiscent of the classical infimum-convolution.
		\end{enumerate}
	\end{rema}
	
	By applying Proposition \ref{pro:Inf-Conv} for each coordinate, we deduce immediately the following corollary.
	
	\begin{coro}
		\label{1-retract}
		Let $\Gamma$ be a nonempty set, $(M,d)$  a metric space, $A$ a nonempty subset of $M$ and $f\colon A \to \ell_\infty(\Gamma)$ be a $C$-Lipschitz map. Then, $f$ admits a $C$-Lipschitz extension $F\colon M\to \ell_\infty(\Gamma)$. 
	\end{coro}
	
	We recall the links between the possibility of extending Lipschitz maps and the notion of absolute Lipschitz retract. We refer the reader to Chapter 1 in  \cite{BenyaminiLindenstrauss2000} for a thorough exposition of the subject.
	
	\begin{defi}\ 
		\begin{enumerate}[(a)]
			\item A nonempty subset $A$ of a metric space $M$ is a \emph{Lipschitz retract} of $M$ if there exists a Lipschitz map $r\colon M\to A$ whose restriction to $A$ is the identity on $A$.\\
			If $r$ is $C$-Lipschitz, for $C\ge 1$, then $A$ is  a \emph{$C$-Lipschitz retract} of $M$.
			\item Let $C\ge 1$. A metric space $M$ is called an \emph{absolute $C$-Lipschitz retract} if it is a $C$-Lipschitz retract of any metric space containing (an isometric copy of) it.
		\end{enumerate}
	\end{defi}
	
	We now give the classical characterizations of absolute Lipschitz retracts. 
	
	\begin{prop}
		\label{prop:absoluteretract}
		Let $M$ be a metric space and $C\ge 1$. The following assertions are equivalent.
		\begin{enumerate}[(i)]
			\item $M$ is an absolute $C$-Lipschitz retract.
			\item For every metric space $N$ and every nonempty subset $B$ of $N$, every Lipschitz map $f\colon B\to M$ can be extended to a Lipschitz map $F\colon N\to M$ with $\Lip(F)\le C\Lip(f)$.
			\item For every metric space $N$ containing $M$ and every metric space $X$, every Lipschitz map $f\colon M\to X$ can be extended to a Lipschitz map $F\colon N\to X$ with $\Lip(F)\le C\Lip(f)$.
		\end{enumerate}
	\end{prop}
	
	\begin{proof} 
		Let $N$ be a metric space containing $M$. Applying $(ii)$ to $B=M$ and $f=Id\colon M\to M$, we see that  $(ii)$ implies $(i)$. Then, applying $(iii)$ to $X=M$ and $f=Id\colon M\to M$, it is also immediate that $(iii)$ implies $(i)$.
		
		Assume that $M$ is an absolute $C$-Lipschitz retract, let $N$ be a metric space containing $M$ and let $X$ be a metric space. There is a $C$-Lipschitz retraction $r$ from $N$ onto $M$ and for any Lipschitz map $f\colon M\to X$, $F=f\circ r :N\to X$ is an extension of $f$ to $N$ with $\Lip(F)\le C\Lip(f)$. We have proved that $(i)$ implies $(iii)$.
		
		Let us finally show that $(i)$ implies $(ii)$ and assume again that $M$ is an absolute $C$-Lipschitz retract. By Theorem \ref{thm:Frechet}, there exists an isometric embedding $f\colon M\to \ell_\infty(M)$.  Thus, $M$ can be seen as a subset of $\ell_\infty(M)$. Therefore, there exists a $C$-Lipschitz retraction $r:\ell_\infty(M)\to M$. Let now $N$ be a metric space, $B$ a subset of $N$ and $f\colon B\to M$ Lipschitz. Considering $f$ as a function from $B$ to $\ell_\infty(M)$, by Corollary \ref{1-retract} we may extend it into a Lipschitz function $g\colon N\to \ell_\infty(M)$ without affecting its Lipschitz constant. Then, $F=r \circ g$ is the desired extension of $f$, and this concludes our proof.
	\end{proof}
	
	\begin{rema} 
		Note that Corollary \ref{1-retract} asserts that $\ell_\infty(\Gamma)$ is an absolute $1$-Lipschitz retract. 
	\end{rema}
	
	We now turn to the case of $\co$.
	
	\begin{prop}
		\label{prop:c0-abs-2-Lip-retract}
		There exists a $2$-Lipschitz retraction from $\ell_\infty$ onto $\co$ and $\co$ is an absolute $2$-Lipschitz retract. 
	\end{prop}
	
	\begin{proof} 
		For $x:=(x_n)_{n=1}^\infty \in \ell_\infty$, observe that $d(x,c_0)=\limsup_{n\to \infty} \abs{x_n}$. Then, set $r(x)_n=0$ if $\abs{x_n}<d(x,c_0)$ and $r(x)_n=(\abs{x_n}-d(x,c_0))\sign(x_n)$ if $\abs{x_n}\ge d(x,c_0)$. We leave it to the reader to check that $r$ is a $2$-Lipschitz retraction from $\ell_\infty$ onto $\co$. 
		
		For the second part of the statement, we will use item $(ii)$ in Proposition \ref{prop:absoluteretract}. So, let $N$ be a metric space, $B\subset N$ and $f\colon B\to c_0$ a Lipschitz map. By Corollary \ref{1-retract}, there exits $F\colon N\to \ell_\infty$ an extension of $f$ with $\Lip(F)=\Lip(f)$. Then, $r\circ F\colon N \to c_0$, where $r\colon \ell_\infty \to c_0$ is a $2$-Lipschitz retract, is an extension of $f$ with $\Lip(r \circ F)\le 2\Lip(f)$.
	\end{proof}
	
	\begin{rema} 
		Proposition \ref{prop:c0-abs-2-Lip-retract} does not have a linear analog. There is no bounded linear projection from $\ell_\infty$ onto $\co$. This result is due to Phillips \cite{Phillips1940}. We also refer the reader to the textbook \cite{AlbiacKalton2016} for the proof and more references. See also Exercise \ref{ex:c0-uncomplemented}.
	\end{rema}

	\section{Linearization of Lipschitz retracts via invariant means}
	\label{sec:invariant-means}
	
	We first introduce the notion of invariant mean on a semigroup. It will be enough for our purpose to consider only the case of abelian semigroups.
	
	\begin{defi} 
		let $(G,+)$ be an abelian semigroup. An \emph{invariant mean} $M$ on $G$ is a  linear functional on $\ell_\infty(G)$ such that:
		\begin{enumerate}[(a)]
			\item $M(\car_G)=1$.
			\item $M$ is positive, i.e. for any $f\in \ell_\infty(G)$ such that $f\ge 0$, $M(f)\ge0$.
			\item $M$ is translation invariant, i.e. for any $f\in \ell_\infty(G)$ and any $g\in G$, $M(\tau_gf)=M(f)$, where $\tau_gf(h):=f(g+h)$ for all $h\in G$.
		\end{enumerate}
	\end{defi}
	
	Note that an invariant mean $M$ on $G$ always belongs to $\ell_\infty(G)^*$ and that its norm is equal to $1$.
	
	\begin{theo} 
		%Let $(G,+)$ be an abelian semigroup. Then, there exists an invariant mean on $G$.
		There exists an invariant mean on every abelian semigroup.
	\end{theo}
	
	\begin{proof} 
		Consider $K:=\{\varphi \in \ell_\infty(G)^*\colon \varphi(\car_G)=1\ \text{and}\ \varphi\ \text{is\ positive}\}.$ It is clear that $K$ is a weak$^*$ closed convex subset of the unit ball of $\ell_\infty(G)^*$ and therefore weak$^*$ compact. For $g\in G$, define $T_g\colon \ell_\infty(G)^* \to \ell_\infty(G)^*$ by:
		$$\forall \varphi \in \ell_\infty(G)^*,\ \forall f\in \ell_\infty(G),\ \ T_g\varphi(f):=\varphi(\tau_g f).$$
		For any $g\in G$, $T_g(K)\subset K$ and since $G$ is abelian, $(T_g)_{g\in G}$ is a commuting family of bounded and weak$^*$ continuous operators on $\ell_\infty(G)^*$. Therefore, it follows from the Markov-Kakutani fixed point theorem that the family $(T_g)_{g\in G}$ has a common fixed point in $K$. This concludes the proof.
	\end{proof}
	
	We now turn to the definition of vector-valued invariant means. Let $(G,+)$ be an abelian semigroup and $M$ be an invariant mean on $G$. For a Banach space $X$, we denote by $\ell_\infty(G,X^*)$ the space of bounded functions from $G$ to $X^*$ equipped with the norm $\norm{f}:=\sup_{g\in G}\norm{f(g)}_{X^*}$. Now, for $f\in \ell_\infty(G,X^*)$, we define $\widehat M (f)\in X^*$ as follows:
	$$\forall x \in X,\ \ \langle x, \widehat M (f)\rangle := M\big(\langle x,f\rangle\big).$$
	It is then easy to verify that $\widehat M\colon \ell_\infty(G,X^*) \to X^*$ is a norm one linear map satisfying the following properties.
	\begin{enumerate}[(i)]
		\item For any $x^*\in X^*$, $\widehat M(\car_G \otimes x^*)=x^*$.
		\item $\widehat M$ is translation invariant. 
	\end{enumerate}
	
	\begin{rema}
		\label{rem:mean-limit}
		In the particular case of an invariant mean $M$ on $(\bN,+)$, it is important to note that if $u:=(u_n)_{n=1}^\infty \subset \bR$ converges to $l\in \bR$, then $M(u)=l$ (see Exercise \ref{ex:invariant-mean}).    
	\end{rema} 
	
	We have seen in the previous section that it is not always possible to linearize a Lipschitz retraction since there exists a Lipschitz retraction from $\ell_\infty$ onto $\co$, but $\co$ is not complemented in $\ell_\infty$. However, the following general result, due to J. Lindenstrauss \cite{Lindenstrauss1964}, will provide sufficient conditions for a Lipschitz retract to be complemented. The main tool of the proof will be the existence of invariant means.
	
	\begin{theo}
		\label{thm:linear-extension}
		Let $X$ and $Y$ be two Banach spaces and $X_0$ be a closed subspace of $X$. Assume that $f\colon X\to Y$ is Lipschitz and that $f_{\restriction_{X_0}}=S \in B(X_0,Y)$. Then, there exists $T\in B(X,Y^{**})$ such that $T_{\restriction_{X_0}}=S$.
	\end{theo}

	\begin{proof}
		Let $M\in \ell_\infty(X_0)^*$ and $N\in \ell_\infty(X)^*$ be invariant means on the abelian groups $(X_0,+)$ and $(X,+)$ respectively. Let $\widehat M\colon \ell_\infty(X_0,Y^{**})\to Y^{**}$ and $\widehat N\colon \ell_\infty(X,Y^{**})\to Y^{**}$ be the associated vector-valued invariant means. Then, for $z\in X$ and $x\in X_0$, we define $f_z(x):=f(z+x)-f(x)$ and we consider $f_z$ as an element of $\ell_\infty(X_0,Y) \subset \ell_\infty(X_0,Y^{**})$, so that we can set $g(z):=\widehat M(f_z)$, for $z\in X$. Thus, $g$ is a map from $X$ to $Y^{**}$. It is easily checked that $g$ is Lipschitz and $\Lip (g)\le \Lip (f)$.
		
		Since $f_{\restriction_{X_0}}$ is linear, we have that $f_z$ is constant equal to $f(z)$ on $X_0$. It follows that for all $z\in X_0$, $g(z)=\widehat M(f_z)=f(z)$. Thus,  $g_{\restriction_{X_0}}=f_{\restriction_{X_0}}=S$.\\
		Since $\widehat M$ is translation invariant, we get that for all $x\in X$ and all $z\in X_0$,
		$$g(x+z)=\widehat M(f_{x+z}-f_z)+g(z)=\widehat M(f_x)+g(z)=g(x)+g(z).$$
		
		Finally, we use $\widehat N$ to perform one last linearization as follows. For $z,x\in X$, set $g_z(x):=g(x+z)-g(x) \in Y^{**}$ and define, as we may, $T(z):=\widehat N(g_z)$. Again, we have that $\Lip(T)\le \Lip (g)\le \Lip (f)$. Also, we have shown that for any $z\in X_0$, $g_z$ is constant equal to $g(z)$ on $X$, so we get that $T_{\restriction_{X_0}}=g_{\restriction_{X_0}}=S$. Finally, using the translation invariance of $\widehat N$, as above, we get that $T$ is additive on $X$. Since $T$ is Lipschitz, it follows now easily that $T$ is linear.
	\end{proof}
	
	We deduce that Lipschitz retracts can be linearized when the target space is complemented in its bidual (and this last condition cannot be dropped).  
	
	\begin{coro}
		\label{cor:linearization-Lip-retract}
		Let $X$ be a Banach space and $Y$ be a closed subspace of $X$ such that $Y$ is complemented in $Y^{**}$. If there exists a Lipschitz retraction from $X$ onto $Y$, then $Y$ is complemented in $X$.
	\end{coro}
	
	\begin{proof}
		It follows from Theorem \ref{thm:linear-extension} that there exists $T\in B(X,Y^{**})$ such that the restriction of $T$ to $Y$ is the identity on $Y$. Letting $P$ be a bounded projection from $Y^{**}$ onto $Y$, one can verify that $PT$ is a bounded projection from $X$ onto $Y$.
	\end{proof}
	
	\section[Lipschitz rigidity of $L_p$-spaces]{\texorpdfstring{Lipschitz rigidity of $L_p$-spaces}{Lipschitz rigidity of  spaces}}
	
	Equipped with the crucial result about linearization of Lipschitz retracts from the previous section, we can improve the basic differentiability results from Section \ref{sec:inf-dim-Rademacher}. 
	We start with a general result due to Heinrich and Mankiewicz \cite{HeinrichMankiewicz1982}.
	
	\begin{theo}
		\label{thm:Lipschitz-isomorphisms}
		\
		\begin{enumerate}[(i)]
			\item Let $X$ be a Banach space that is complemented in its bidual. If $X$ is Lipschitz equivalent to a Banach space $Y$ with RNP, then every complemented separable subspace of $X$ is isomorphic to a complemented (separable) subspace of $Y$.
			%Let now $X_0$ be a complemented separable subspace of $X$. Then, $X_0$ is isomorphic to a complemented subspace of $Y$.
			\item Let $X$ be a Banach space that is Lipschitz equivalent to a separable reflexive Banach space $Y$. Then, $X$ is isomorphic to a complemented subspace of $Y$ and $Y$ is isomorphic to a complemented subspace of $X$.
			\item Assume that $X$ is a separable dual space that is Lipschitz equivalent to a Banach space $Y$. Then, $X$ is isomorphic to a complemented subspace of $Y$.
		\end{enumerate}
	\end{theo}
	
	\begin{proof} 
		Assertion $(i)$ is the heart of the matter. So, let $X_0$ be a complemented separable subspace of $X$, $f\colon X\to Y$ be a Lipschitz equivalence, and denote by $f_0$ its restriction to $X_0$. Since $X_0$ is separable and $Y$ has RNP, Theorem \ref{thm:infinite-Rad} ensures the existence of $x_0\in X_0$ such that $f_0$ is Gateaux differentiable at $x_0$. Let us assume, as we may, that $x_0=0$ and $f(0)=0$. For $n\in \bN$, define now $g_n\colon Y\to X$ by $g_n(y) := nf^{-1}(y/n)$. Note that $\norm{g_n(y)}_X \le C\norm{y}_Y$, where $C:=\Lip (f^{-1})$ and that $\Lip (g_n)\le C$. Let $M\in \ell_\infty(\bN)^*$ be an invariant mean on $(\bN,+)$ and $\widehat M\colon \ell_\infty(\bN,X^{**}) \to X^{**}$ be the associated vector-valued invariant mean and define $h\colon Y\to X^{**}$ by $h(y) := \widehat M\big((g_n(y))_{n\in \bN}\big)$. Since $\widehat M$ has norm one, we have that $\Lip (h)\le C$.\\
		For $x\in X_0$, let $y_n := nf(x/n)$. Since $\lim_{n\to \infty}\|y_n-D_{f_0}(0)x\|=0$, $g_n(y_n)=x$ and all the $g_n$ are $C$-Lipschitz, we get that
		$$\lim_{n\to \infty}\norm{g_n(D_{f_0}(0)x)-x}=0.$$
		It follows from Remark \ref{rem:mean-limit} that for all $x\in X_0$, $$h\big(D_{f_0}(0)x\big)=x.$$
		Since $X_0$ is complemented in $X$ and $X$ is complemented in $X^{**}$, there exists a bounded projection $P$ from $X^{**}$ onto $X_0$. If we let $\varphi := D_{f_0}(0)\circ P \circ h$, then it is easy to verify that $\varphi$ is a Lipschitz retraction from $Y$ onto $Z := D_{f_0}(0)(X_0)$. Note that, since $X$ is complemented in $X^{**}$ and $X_0$ is complemented in $X$, $X_0$ is also complemented in $X_0^{**}$. Since $Z$ is isomorphic to $X_0$, it is complemented in its bidual, and it follows from Corollary \ref{cor:linearization-Lip-retract} that $Z$ is complemented in $Y$ and therefore that $X_0$ is isomorphic to a complemented subspace of $Y$.
		
		$(ii)$ Theorem \ref{thm:linearization-applications} insures that $X$ is also reflexive. The conclusion then immediately follows from $(i)$ and the fact that reflexive Banach spaces have RNP.
		
		$(iii)$ Since $X$ is a separable dual, it has RNP and therefore so does $Y$. Then, the conclusion follows again directly from $(i)$ since the dual $E^*$ of a Banach space $E$ is always complemented in its bidual $E^{***}$.
	\end{proof}
	
	We can now solve the Lipschitz rigidity problem for $L_p$-spaces with $p\in (1,\infty)$.
	
	\begin{theo}
		\label{thm:Lipschitz-rigidity-Lp}
		For every $1<p<\infty$ $\ell_p$ and $L_p$ are Lipschitz rigid.
	\end{theo}
	
	\begin{proof} 
		The case of $\ell_p$ is the easiest since a classical theorem of Pe{\l}czy\'{n}ski \cite{Pelczynski1960} insures that, for $p\in [1,\infty)$, any infinite-dimensional complemented subspace of $\ell_p$ is isomorphic to $\ell_p$. Since by Theorem \ref{thm:Lipschitz-isomorphisms} $(ii)$, any Banach space that is Lipschitz equivalent to $\ell_p$ with $p\in(1,\infty)$ must be isomorphic to a complemented subspace of $\ell_p$ we are done.
		The description of the complemented subspaces of $L_p$ is not as simple, and a little more work is needed. It follows again from Theorem \ref{thm:Lipschitz-isomorphisms} $(ii)$ that there exist Banach spaces $X_1$ and $Y_1$ such that $L_p\simeq X \oplus_p X_1$ and $X\simeq L_p \oplus_p Y_1$. We can conclude with a decomposition argument where we will use freely the elementary fact that $L_p\oplus_p L_p$ and $\ell_p(L_p)$ are isometric to $L_p$. First observe that 
		\begin{equation*}
			L_p\oplus_p X\simeq L_p\oplus_p L_p \oplus_p Y_1 \equiv L_p \oplus_p Y_1 \simeq X,
		\end{equation*}
		and then note that
		\begin{align*}
			L_p\oplus_p X\simeq \ell_p(L_p) \oplus_p X \simeq \ell_p(X_1\oplus_p X) \oplus_p X & \equiv \ell_p(X_1)\oplus \ell_p(X)\oplus_p X  \\
			& \equiv \ell_p(X_1)\oplus \ell_p(X) \\
			& \simeq  \ell_p(X_1\oplus_p X)  \\
			& \simeq \ell_p(L_p) \equiv L_p.
		\end{align*}
		This finishes the proof.
	\end{proof}
	
	For $p=1$ in the sequence space situation, we have a partial result. 
	
	\begin{theo}
		\label{thm:Lipschitz-rigidity-l1}
		If $X$ is complemented in its bidual (in particular if $X$ is isomorphic to a dual space) and Lipschitz equivalent to $\ell_1$, then $X$ is (linearly) isomorphic to $\ell_1$.
	\end{theo}
	
	\begin{proof}
		Since $\ell_1$ has RNP, $X$ must be isomorphic to a complemented subspace of $\ell_1$ by Theorem \ref{thm:Lipschitz-isomorphisms} $(i)$ and hence it is isomorphic to $\ell_1$.
	\end{proof}
	
	Nevertheless, the general problem remains open.
	
	\begin{prob}
		\label{prob:Lipschitz-rigidity-l1}
		Is $\ell_1$ Lipschitz rigid?   
	\end{prob}
	
	The case $p=1$ is completely open for function spaces.
	
	\begin{prob}
		\label{prob:Lipschitz-rigidity-L1}
		Is $L_1$ Lipschitz rigid?  
	\end{prob}

	\section{Ultrapowers techniques and Ribe's rigidity theorem}
	
	In this section, we exploit the weak$^*$ Gateaux differentiability results combined with ultrapowers techniques. As a first application, we state and prove Ribe's rigidity theorem on the stability of local linear properties of Banach spaces under coarse-Lipschitz embeddings. We invite the reader to read the original proof by Ribe \cite{Ribe1976}, where there is no use of ultrapowers nor of differentiability, and which certainly provides more insight. 
	
	Recall that a Banach space $X$ is said to be \emph{crudely finitely representable into} a Banach space $Y$ if there exists a constant $C\ge 1$ such that for any finite-dimensional subspace $E$ of $X$, there exist a finite-dimensional subspace $F$ of $Y$ and an isomorphism $T$ from $E$ onto $F$ with $\norm{T}\,\|T^{-1}\|\le C$. If $C$ can be taken to be as close as $1$ as one wishes, then one talks about \emph{finite representability}.  
	
	\begin{theo}[Ribe's rigidity theorem, 1976]
		\label{thm:Ribe-rigidity}  
		Let $X$ and $Y$ be two Banach spaces. If $X$ coarse-Lipschitz embeds into $Y$, then $X$ is crudely finitely representable in $Y$. In particular, if $X$ and $Y$ are coarsely or uniformly equivalent, then $X$ and $Y$ are crudely finitely representable in each other. 
	\end{theo}
	
	\begin{proof} 
		Let $\cU$ be a nonprincipal ultrafilter on $\bN$ and assume that $X$ coarse-Lipschitz embeds into $Y$. By Proposition \ref{prop:CL-ultrapower-Lip}, there exists a bi-Lipschitz embedding $f$ from $X^{\cU}$ into $Y^{\cU}$. Let $D\ge 1$ be the distortion of $f$ and $E$ be a finite-dimensional subspace of $X$. It follows from Theorem \ref{thm:embed-in-dual+} that $E$ linearly embeds into $(Y^{\cU})^{**}$ with distortion $D$. Now, given $C>D$, it follows from the principle of local reflexivity and the finite representability of $Y^{\cU}$ into $Y$ that $E$ linearly embeds into $Y$ with distortion at most $C$. 
		The last statement follows from the fact observed in Section \ref{sec:nonlinear-equivalences} that coarse and uniform equivalences are coarse-Lipschitz equivalences.
	\end{proof}
	
	Recall that a Banach space is super-reflexive if each of its ultrapowers is reflexive. It is not trivial to show that super-reflexivity is an isomorphic invariant. This can also be deduced from the next theorem, which contains the fact that super-reflexivity is a coarse-Lipschitz invariant. 
	
	\begin{theo}
		\label{thm:CL-stability-SR}
		If a Banach space $X$ coarse-Lipschitz embeds into a super-reflexive Banach space $Y$, then $X$ is super-reflexive. Moreover, if $X$ is separable, then $X$ is isomorphic to a subspace of any ultrapower of $Y$.
	\end{theo}
	
	\begin{proof} 
		If there is a coarse-Lipschitz embedding from $X$ into $Y$, then it follows from Proposition \ref{prop:CL-ultrapower-Lip} that there is a bi-Lipschitz embedding from $X^{\cU}$ into $Y^{\cU}$. Since $X^{\cU}$ is usually not separable, we cannot necessarily differentiate, but thanks to the separable determination of reflexivity, we can nevertheless invoke Theorem \ref{thm:linearization-RNP} to conclude that $X^{\cU}$ is reflexive and hence $X$ is super-reflexive. However, if $X$ is separable, we can differentiate the restriction to $X$ of the bi-Lipschitz embedding from $X^{\cU}$ into $Y^{\cU}$ and conclude that $X$ isomorphically embeds into $Y^{\cU}$.
	\end{proof}
	
	Finally, for coarse-Lipschitz equivalences between separable super-reflexive spaces, we prove the following.  
	%result whose proof combines the previous technique with a separable saturation argument. 
	This result will be used crucially when establishing nonlinear rigidity results for quotients of $\ell_p$. 
	
	\begin{theo}
		\label{thm:CL-equiv-complemented-SR}
		Assume that a Banach space $X$ is coarse-Lipschitz equivalent to a separable and super-reflexive Banach space $Y$. Then, for any nonprincipal ultrafilter $\cU$ on $\bN$, $X$ is isomorphic to a  separable complemented subspace of $Y^{\cU}$.
	\end{theo}
	
	\begin{proof}
		Theorem \ref{thm:CL-stability-SR} already ensures that $X$ is super-reflexive and thanks to Proposition \ref{prop:CL-equiv->Lip-equiv-ultra} there exists a Lipschitz equivalence $f^{\cU}\colon X^{\cU}\to Y^{\cU}$. Now, since $X$ is reflexive, it is $1$-complemented in $X^{\cU}$. Indeed, the map $P\colon X^{\cU} \to X$ defined by $P((x_n)_{n=1}^\infty)= w-\lim_{n,\cU}x_n$ is a norm one projection from $X^{\cU}$ onto $X$. Since $Y^{\cU}$ has RNP and $X^{\cU}$ is reflexive, thus complemented in its bidual, we can directly apply item $(i)$ of Theorem \ref{thm:Lipschitz-isomorphisms} to deduce that $X$ is isomorphic to a separable complemented subspace of $Y^{\cU}$. 
	\end{proof}

	\section[Rigidity of Asplund spaces under Lipschitz quotients]{Rigidity of Asplund spaces under Lipschitz quotients.}
	\label{sec:Asplundrigidity}

	The study of Asplund spaces is a vast subject. We shall only collect the information needed for this section here. We refer the reader to the books \cite{DiestelUhl1977} and \cite{DGZ1993} for a thorough study of this notion. Recall that a \emph{$\cG_\delta$ set} is a countable intersection of open sets.
	
	\begin{defi} 
		A Banach space $X$ is called an \emph{Asplund space} if for any open convex subset $C$ of $X$ and any $f\colon C \to  \bR$ convex and continuous, $f$ is Fr\'echet differentiable on a dense $\cG_\delta$ subset of $C$.
	\end{defi}
	
	Like for the Radon-Nikod\'{y}m property, there are many characterizations of Asplund spaces (see \cite{DGZ1993} and references therein for an overview). We will only give the following, which is partly contained in Theorem \ref{thm:dual-RNP+}.
	
	\begin{theo}
		\label{thm:Asplund}
		\,
		Let $X$ be a Banach space. The following assertions are equivalent.
		\begin{enumerate}[(i)]
			\item The space $X$ is an Asplund space.
			\item Every separable subspace of $X$ has a separable dual.
			\item The space $X^*$ has the Radon-Nikod\'{y}m property. 
		\end{enumerate}
		In particular, a separable Banach space is an Asplund space if and only if its dual is separable.
	\end{theo}
	
	Note that being an Asplund space passes to subspaces and quotients. A crucial tool for this section is the following major result due to D. Preiss \cite{Preiss1990}.
	
	\begin{theo}
		\label{thm:Preiss} 
		%Let $U$ be an open subset of an Asplund  space $X$ and $f\colon  U \to \R$ be a Lipschitz map. Then, $f$ admits a point of Fr\'echet differentiability.
		Every map defined on an open subset of an Asplund space that is real-valued and Lipschitz admits a point of Fr\'echet differentiability.
	\end{theo}
	
	We have chosen not to give here the long and difficult proof of this beautiful result. Let us just mention the following important open question.
	
	\begin{prob} 
		Let $X$ be an Asplund space and $f\colon  X \to \R^2$ be a Lipschitz function. Does $f$ admit a point of Fr\'echet differentiability?
	\end{prob}
	
	One of the main results of this section, which is taken from \cite{BJLPS1999} is that the class of Asplund spaces is stable under Lipschitz equivalences. In fact it is even stable under Lipschitz quotients. We shall detail the proof of the more general result but we encourage the reader to rewrite the proof for Lipschitz equivalences (see Exercise \ref{ex:Asplundrigidity}).
	
	\begin{defi}
		\label{def:Lipchitz-quotient}
		Let $M$ and $N$ be two metric spaces. A map $f\colon  M \to N$ is said to be \emph{co-Lipschitz} if there exists $C>0$ such that for every $r>0$ and every $x\in M$, $B_N(f(x),\frac{r}{C}) \subset f(B_M(x,r))$. The infimum of all such constants $C$ is the co-Lipschitz constant of $f$ and is denoted by $\coLip(f)$.
		Then, we say that a map $f\colon  M \to N$ is a \emph{Lipschitz quotient mapping} if it is Lipschitz and co-Lipschitz.
		Finally, a metric space $N$ is said to be a \emph{Lipschitz quotient} of a metric space $M$ provided there is a Lipschitz quotient mapping from $M$ \underline{onto} $N$.
	\end{defi}
	
	It is clear that a Lipschitz isomorphism is a Lipschitz quotient. On the other hand, the open mapping theorem asserts that a linear quotient (namely, a continuous surjective linear map between two Banach spaces) is a Lipschitz quotient. When restricting a co-Lipschitz map to a susbet one might loose the co-Lipschitz property. However, this drawback can be circumvented in some favorable situation. 
	
	\begin{lemm}
		\label{lem:separable-saturation-quotient}
		Let $f$ be a continuous, co-Lipschitz mapping from a complete metric space $M$ onto a separable metric space $N$ and let $M_0$ be a separable subset of $M$. Then, there exists a separable closed subset $M_1$ of $M$ that contains $M_0$ and such that $f_{\restriction_{M_1}}$, the restriction of $f$ to $M_1$, is co-Lipschitz from $M_1$ onto $N$ with $\coLip(f_{\restriction_{M_1}})\le \coLip(f)$. If $M$ is a Banach space $X$, then $M_1$ can be taken to be a subspace of $X$.    
	\end{lemm}

	\begin{proof} 
		By rescaling the metric of $M$, we may assume that $\coLip(f)=1$. So, we have that for all $x\in M$ and all $r>0$, $f(B_M(x,r))$ contains $B_N(f(x),s)$ for all $s\in (0,r)$. Since $M_0$ and $N$ are separable, we can pick $D_0$ a countable subset of $M$ such that $M_0 \subset \overline{D_0}$ and $f(D_0)$ is dense in $N$. Then, using again the separability of $N$, we  build inductively $D_0\subset D_1 \subset \dots \subset D_n \dots$, countable subsets of $M$ so that for each $n\in \bN$, each $x\in D_n$ and all $s<r \in \bQ_+$, $B_N(f(x),s)$ is included in the closure of $f(B_M(x,r)\cap D_{n+1})$. If $M$ is a Banach space, we can also ensure that $D_n$ is closed under rational linear combinations. 
		
		We now denote by $M_1$, the closure of  $\cup_{n\ge 0} D_n$. It follows from the inductive construction that for all $x\in M_1$ and all $0<s<r$:
		$$B_N(f(x),s)\subset \overline{f(B_{M_1}(x,r))}.$$
		We deduce, using the continuity of $f$ and a classical ``successive approximations'' type of argument that for all $x\in M_1$ and all $0<s<r$:
		$$B_N(f(x),s)\subset f(B_{M_1}(x,r)).$$
		Indeed, assume that $x\in M_1$ and $d_N(y,f(x))\le s$ and fix $t\in (s,r)$ and $\eta\in (0,1)$ (to be precised later). Then, there exists $x_1 \in M_1$ such that $d_M(x,x_1)\le t$ and $d_N(y,f(x_1))\le \eta$. Next, we build inductively $(x_n)_{n=1}^\infty \subset M_1$ such that for all $n\in \bN$, $d_M(x_n,x_{n+1})\le 2\eta^n$ and $d_N(y,f(x_n))\le \eta^n$. Since $M$ is complete and $M_1$ is closed in $M$, the sequence $(x_n)_n$ converges to some $z$ in $M_1$ and $f(z)=y$ by the continuity of $f$. On the other hand, $d_M(x,z)\le t+2\sum_{n=1}^\infty \eta^n \le r$, if $\eta$ was carefully chosen beforehand to be sufficiently small. 
	\end{proof}
	
	The next proposition will allow us to reduce the upcoming Lipschitz rigidity proof to the separable setting.
	
	\begin{prop}
		\label{prop:separable-saturation-quotient}
		Let $f$ be a continuous, co-Lipschitz mapping from a Banach space $X$ onto a Banach space $Y$. Let $X_0$ be a separable subspace of $X$ and $Y_0$ be a separable subspace of $Y$. Then, there exist a separable subspace $X_1$ of $X$ and a separable closed subspace $Y_1$ of $Y$ so that $X_0\subset X_1$, $Y_0 \subset Y_1$, $f(X_1)=Y_1$ and $f_{\restriction_{X_1}}$ is co-Lipschitz from $X_1$ onto $Y_1$ with $\coLip(f_{\restriction_{X_1}})\le \coLip(f)$.
	\end{prop}
	
	\begin{proof}
		Note first that for any subset $Z$ of $Y$, the restriction of $f$ to $f^{-1}(Z)$ is co-Lipschitz onto $Z$ with co-Lipschitz constant at most $\coLip(f)$. Combining this with Lemma \ref{lem:separable-saturation-quotient}, we build separable closed subsets of $X$: $X_0\subset W_0\subset \dots \subset W_n\subset \dots$ such that for each $n$, $f_{\restriction_{W_n}}$ is co-Lipschitz from $W_n$ onto $f(W_n)$ with $\coLip(f_{\restriction_{W_n}})\le \coLip(f)$, $f(W_0)\supset Y_0$, $W_{n+1}$ contains the linear span of $W_n$ and $f(W_{n+1})$ contains the linear span of $f(W_n)$. Let $E_1:=\cup_{n\ge 0}W_n$ and $F_1:=f(E_1)=\cup_{n\ge 0}f(W_n)$. It is now clear that $E_1$ is a subspace of $X$, $F_1$ is a subspace of $Y$ and $f$ is co-Lipschitz from $E_1$ onto $F_1$ with co-Lipschitz constant at most $\coLip(f)$. Finally, using successive approximations as at the end of the previous proof, we deduce that $f$ is co-Lipschitz from $X_1$ onto $Y_1$ with co-Lipschitz constant at most $\coLip(f)$, where $X_1$ is the closure of $E_1$ and $Y_1$ is the closure of $F_1$.
	\end{proof}

	We can now state and prove the main result of this section, which is due to S. Bates, W.B. Johnson, J. Lindenstrauss, D. Preiss and G. Schechtman \cite[Theorem 3.18]{BJLPS1999}. 
	
	\begin{theo}
		\label{thm:Asplund-rigidity}
		%Let $X$ and $Y$ be Banach spaces. Assume that $f\colon X\to Y$ is a Lipschitz quotient onto $Y$ and that $X$ is an Asplund space. Then, $Y$ is an Asplund space.
		A Banach space that is a Lipschitz quotient of an Asplund space is automatically an Asplund space.
		In particular, being Asplund is a Lipschitz invariant.
	\end{theo}
	
	\begin{proof}
		It follows from Theorem \ref{thm:Asplund} and Proposition \ref{prop:separable-saturation-quotient} that we may assume that $X$ is separable (and thus $X^*$ is separable, because $X$ is Asplund) and that we need to prove that $Y^*$ is separable. Assume also that $\Lip(f)=1$ and that $\coLip(f)=C>0$. Suppose, aiming for a contradiction, that $Y^*$ is not separable. Then, there exists a family $(y^*_\gamma)_{\gamma \in \Gamma}$ in $S_{Y^*}$ with $\Gamma$ uncountable such that for all $\gamma \neq \gamma' \in \Gamma$,
		\begin{equation*}
			\norm{y^*_\gamma-y^*_{\gamma'}}>\frac12.
		\end{equation*}
		Let $f_\gamma := y^*_\gamma \circ f$ and note that $\Lip(f_\gamma)\le 1$. By Theorem \ref{thm:Preiss}, for any $\gamma \in \Gamma$, there exists $x_\gamma \in B_X$ such that $f_\gamma$ is Fr\'echet differentiable at $x_\gamma$. Then, letting $x^*_\gamma := Df_\gamma(x_\gamma)$, we deduce that for any $\gamma \in \Gamma$ there exists $\delta_\gamma >0$ such that for all $z \in X$,
		\begin{equation*}
			\norm{z}\le \delta_\gamma \Rightarrow \abs{f_\gamma(x_\gamma+z) - f_\gamma(x_\gamma)-x^*_\gamma(z)}\le \frac{\norm{z}}{10C}.
		\end{equation*}
		By passing to an uncountable subset of $\Gamma$, we may as well assume that $\delta := \inf_{\gamma \in \Gamma}\delta_\gamma >0$ and, since $X$ and $X^*$ are separable that for all $\gamma,\gamma' \in \Gamma$:
		\begin{equation*}
			\norm{x_\gamma - x_{\gamma'}} < \frac{\delta}{10C},\ \norm{x^*_\gamma - x^*_{\gamma'}}<\frac{1}{10C},
		\end{equation*}
		and
		\begin{equation*}
			\abs{f_\gamma(x_\gamma) - f_{\gamma'}(x_{\gamma'})} < \frac{\delta}{10C}.
		\end{equation*}
		Therefore, for $\norm{z}\le \delta$, we have
		\begin{align*}
			|f_\gamma(x_\gamma+z) - & f_{\gamma'}(x_{\gamma}+z)|\\
			\le & \abs{ f_\gamma(x_\gamma+z) - f_\gamma(x_\gamma)-x^*_\gamma(z)} + \abs{ f_\gamma(x_\gamma) - f_{\gamma'}(x_{\gamma'})} +
			\abs{x^*_\gamma(z) - x^*_{\gamma'}(z)}\\
			& + \abs{f_{\gamma'}(x_{\gamma'}+z) - f_{\gamma'}(x_{\gamma'}) - x^*_{\gamma'}(z)} + \abs{f_{\gamma'}(x_{\gamma}+z)-f_{\gamma'}(x_{\gamma'}+z)}\\
			\le & \frac{1}{10C}(\norm{z} + \delta + \norm{z} + \norm{z} + \delta) \le \frac{\delta}{2C}.
		\end{align*}
		On the other hand, remembering that $f(B_X(x_\gamma,\delta))\supset B_Y(f(x_\gamma), \frac{\delta}{C})$, we have
		\begin{align*}
			\sup_{z\in \delta B_X} |f_\gamma(x_\gamma+z) - & f_{\gamma'}(x_{\gamma}+z)| = \sup_{z\in \delta B_X} \abs{ (y^*_\gamma - y^*_{\gamma'})(f(x_\gamma+z))}\\
			& \ge \sup_{y\in \frac{\delta}{C}B_Y} \abs{ (y^*_\gamma - y^*_{\gamma'})(f(x_\gamma)+y)} \ge \frac{\delta}{C} \norm{ y^*_\gamma - y^*_{\gamma'}}>\frac{\delta}{2C}.
		\end{align*}
		This yields the desired contradiction and concludes the proof.
	\end{proof}
	
	\begin{rema} Recall that a map $f\colon (M,d_M)\to (N,d_N)$ between metric spaces is called a \emph{uniform quotient map} and $N$ is simply said to be a \emph{uniform quotient} of $M$, if $f$ is surjective, uniformly continuous and \emph{co-uniformly continuous}, i.e. for every $r>0$ there exists $\delta(r)>0$ such that for all $x\in M$, one has
		\begin{equation*}
			B_{N}\left(f(x),\delta(r)\right)\subset f(B_M(x,r)).
		\end{equation*}
		
		Note that uniform equivalences are examples of uniform quotient maps. As we will see in Chapter \ref{chapter:Counterexamples}, the analog of Theorem \ref{thm:Asplund-rigidity} is not true for uniform quotients. Indeed, Ribe exhibited in \cite{Ribe1984} a reflexive space $X$ that is uniformly equivalent to $X \oplus \ell_1$.
	\end{rema}

	\section{Notes}
	
	Recall that one of the main problems in the nonlinear geometry of Banach spaces asks whether any separable Banach space is Lipschitz rigid (see Problem \ref{prob:Lipschitz-rigidity}). Therefore, enlarging the list of classical Banach spaces that can be shown to be Lipschitz rigid is a natural program. Using the techniques developed in this chapter, one can identify a few more examples of spaces that are Lipschitz rigid: $\ell_p(\ell_q)$ and $\ell_p \oplus \ell_q$ for $1<p,q<\infty$ (see Exercise \ref{exe:Lipschitzrigidityl_p(l_q)}), or the James space $\James$ (see Exercise \ref{exe:LipschitzrigidityJames}). Another important example of a Lipschitz rigid space is $c_0$. However, the method to prove it is completely different and will be explained in Chapter \ref{chapter:Gorelik}. 
	
	Recall that a Banach space is \emph{prime} if every one of its infinite-dimensional complemented subspaces is isomorphic to the whole space. According to Theorem \ref{thm:Lipschitz-isomorphisms}, any separable reflexive Banach space that is prime will be Lipschitz rigid. Examples of prime spaces are $\ell_p$ ($1\le p\le \infty$) and $\co$. Since the original Tsirelson space, now denoted by $\Tsi^*$ and its dual $\Tsi$, are not prime \cite{CasazzaOdell}, the following problem has remained open. We refer the reader to Section \ref{sec:asymptotic-c_0} for the definition of $\Tsi^*$ and its main properties.  
	
	\begin{prob}
		\label{pb:Lipschitz-rigidity-co-Tsirelson}
		Is Tsirelson space $\Tsi^*$, or its dual $\Tsi$, Lipschitz rigid?   
	\end{prob}
	
	A Banach space $X$ is called \emph{(complementably) minimal} if every infinite-dimensional subspace of $X$ contains a (complemented) subspace isomorphic to $X$.
	It was shown in \cite{CJT84} that $\Tsi^*$ was the first example of a minimal space, other than $\ell_p$ or $\co$, and that $\Tsi^*$ is not complementably minimal. P. Casazza asked whether $\ell_p$ ($1\le p\le \infty$) and $\co$ are the only complementably minimal spaces, but it was shown in \cite{AndroulakisSchlumprecht03} that Schlumprecht space $\mathrm{S}$, which does not contain $\ell_p$ ($1\le p< \infty$) or $\co$, is complementably minimal. Consequently, $S$ is either prime or there exists a complemented subspace $Z$ of $S$ such that $S$ and $Z$ are examples of nonisomorphic  Banach spaces (with unconditional bases) that are isomorphic to complemented subspaces of each other (a refinement of the Schroeder-Bernstein problem for Banach spaces). Whether $S$ is prime is still an open problem, and its Lipschitz rigidity status has yet to be determined.
	
	\begin{prob}
		\label{pb:Lipschitz-rigidity-Schlumprecht}
		Is Schlumprecht space $\mathrm{S}$ Lipschitz rigid?   
	\end{prob}
	
	It seems to be a good place to recall Problem \ref{pb:Lipschitz-embed-rigidity-c0}: assuming that $c_0$ bi-Lipschitzly embeds into a Banach space $X$, does it imply that $c_0$ isomorphically embeds into $X$? It follows from Theorem \ref{thm:infinite-Rad} that such an $X$ cannot have RNP. N. Kalton proved that Problem \ref{pb:Lipschitz-embed-rigidity-c0} has a positive answer when restricted to Banach lattices, i.e., if $c_0$ bi-Lipschitzly embeds into a Banach lattice $X$, then $c_0$ isomorphically embeds into $X$. In order to see why this statement holds, let us anticipate some important results of this book (see Chapter \ref{chapter:interlaced-graphs} for all relevant definitions and results). Kalton proved in \cite{Kalton2007} that a Banach space whose unit ball uniformly embeds into a reflexive Banach space has Property $\cQ$. It follows that if $X$ equi-coarsely contains the interlacing graphs, in particular if $c_0$ bi-Lipschitzly embeds into $X$, then $B_X$ does not uniformly embed into a reflexive Banach space. But Kalton also proved in \cite{Kalton2007} that for a separable Banach lattice $X$, $B_X$ uniformly embeds into a reflexive Banach space if and only if $X$ does not contain any subspace isomorphic to $\co$. The conclusion follows. A more detailed discussion of this remark can be found in Section 5.2 of \cite{BLMS2020b}.

	\medskip 
	We refer the reader to \cite{BJLPS1999} for many other important results about uniform quotients. Let us just mention the following question, which we believe to be open. 
	
	\begin{prob}
		Are separable Lipschitz quotients isomorphic to linear quotients?
	\end{prob}

	\section{Exercises}
	
	\begin{exer}
		\label{ex:no-RNP}\
		\begin{enumerate}
			\item  Prove that $\co$ fails the Radon-Nikod\'{y}m property.
			\item Prove that $L^1(\R)$ fails the Radon-Nikod\'{y}m property.
		\end{enumerate}
	\end{exer}
	
	\begin{proof}[Hint:]
		$1.$ Consider the function $f\colon \R \to c_0(\bN)$ defined by
		$$f(t)=\Big(\frac{\sin nt}{n}\Big)_{n=1}^\infty,$$
		or, note that $c_0(\bN,\R)$ is $\R$-isomorphic to $c_0(\bN,\C)$ and use $$g(t)=\Big(\frac{e^{int}}{n}\Big)_{n=1}^\infty.$$
		$2.$ consider the function $f\colon \R \to L^1(\R)$ defined by $f(t)=\car_{[0,|t|]}$.
	\end{proof}
	
	\begin{exer}\label{ex:RNPl_psums}
		Let $1\le p<\infty$ and $(X_k)_{k=1}^\infty$ be a sequence of Banach spaces with the Radon-Nikod\'{y}m property. Show that $(\sum_{k=1}^\infty X_k)_{\ell_p}$ has the Radon-Nikod\'{y}m property.
	\end{exer}

	\begin{exer}
		\label{ex:additive}
		Prove Lemma \ref{lem:additive}   
	\end{exer}
	
	\begin{exer}
		\label{ex:Fréchet=Gâteaux}
		Let $E$ be a Banach space and $n\in\bN$. Prove that if $f\colon \R^n \to E$ is Lipschitz and Gateaux differentiable at $x\in \R^n$, then it is Fr\'{e}chet differentiable at $x$.  
	\end{exer}
	
	\begin{exer}
		\label{ex:no-Lebesgue-measure}
		Show that any translation-invariant Borel measure on an infinite-dimensional separable Banach space must be either infinite for all sets or zero for all sets
	\end{exer}
	
	\begin{proof}[Hint:]
		You can use Riesz lemma together with the fact that separable metric spaces are Lindel\"of spaces.
	\end{proof}
	
	\begin{exer}
		Let $X$ and $Y$ be two Banach spaces. Assume that $f\colon B_X\to Y$ is Lipschitz. Show that $f$ admits a Lipschitz extension $\tilde f\colon X\to Y$.
	\end{exer} 
	
	\begin{proof}[Hint:]
		Extend the restriction of $f$ to $S_X$ homogeneously outside $B_X$
	\end{proof}
	
	\begin{exer}
		Let $(M,d)$ be a metric space, $A$ a nonempty subset of $M$ and $f\colon A \to \bR$ be a $C$-Lipschitz map.
		\begin{enumerate}
			\item Show that the function defined by $G(x):=\sup\{f(y)-Cd(x,y)\colon y\in A\}$ provides $C$-Lipschitz extension of $f$. 
			\item Show that for any $C$-Lipschitz extension $H$ of $f$, we have that $G\le H \le F$ on $M$, where $F(x):=\inf\{f(y)+Cd(x,y)\colon y\in A\}$.
		\end{enumerate}
	\end{exer}
	
	\begin{exer} 
		Let $X$ and $Y$ be two Banach spaces. Assume that $f\colon X\to Y$ is a Lipschitz isomorphism and that $f$ is Fr\'echet differentiable at $x\in X$. Show that $Df(x)$ is a linear isomorphism from $X$ onto $Y$. 
	\end{exer}

	\begin{exer}\ 
		\begin{enumerate}
			\item Construct a bi-Lipschitz embedding from $\bR$ into a dual Banach space $X^*$, whose weak$^*$ derivative at $0$ is null. 
			\item Is this map differentiable at $0$?
		\end{enumerate}
	\end{exer}
	
	\begin{proof}[Hint:]
		Consider your favorite normalized weak$^*$ null sequence $(x_n^*)_{n=1}^\infty$ in a dual Banach space. Define $f$ such that $f(2^{-n})=2^{-n}x^*_n$ and extend it to $\bR$ in an affine way.
	\end{proof}
	
	\begin{exer}  
		We denote by $c$ the space of all real-valued converging sequences, equipped with the supremum norm. Assume that $C\in [1,2)$. Show that there is no $C$-Lipschitz retraction from $c$ onto $\co$.
	\end{exer}

	\begin{exer}
		\label{ex:c0-uncomplemented} 
		The goal of this exercise is to show that $\co$ is not complemented in $\ell_\infty$.
		\begin{enumerate}
			
			\item Let $\bQ=\{q_n\colon n\in \bN\}$ be an enumeration of $\bQ$. To each $\lambda \in \bR$ we assign a sequence of pairwise distinct rational numbers that converges to $\lambda$. We denote by $E_\lambda$ the set of values of this sequence and $N_\lambda:=\{n\in \bN\colon q_n\in E_\lambda\}$. Show that for any $\lambda \in \bR$, $N_\lambda$ is an infinite subset of $\bN$ and that for all $\lambda \neq \mu \in \bR$, $N_\lambda \cap N_\mu$ is finite.
			
			\item Let $f \in \ell_\infty^*$ such that $c_0\subset \mathrm{Ker}(f)$. Let us denote by $\car_\lambda \in \ell_\infty$ the indicator function of $N_\lambda$ and let $A_n :=\{\lambda \in \bR\colon |f(\car_\lambda)|> 1/n\}$. Show that the cardinality of $A_n$ is at most $n\norm{f}_{\ell_\infty^*}$. Deduce that $\{\lambda \in \bR\colon f(\car_\lambda)\neq 0\}$ is countable.
			
			\item Assume now that $P$ is a bounded projection from $\ell_\infty$ onto $\co$. Show that $\co=\bigcap_{n=0}^\infty \mathrm{Ker}(f_n)$, where $f_n := \phi_n\circ(I-P)$ and $\phi_n$ is the n$^{\rm th}$ coordinate functional defined on $\ell_\infty$.\\ Deduce that there exists $\lambda \in \bR$ such that $\car_\lambda \in \co$ and conclude.
			
		\end{enumerate}
	\end{exer}
	
	\begin{exer}
		\label{ex:invariant-mean} 
		Let $M$ be an invariant mean on $(\bN,+)$.\\
		Show that for any converging sequence $x=(x_n)_{n=0}^\infty$ in $\bR$, $M(x)=\lim_{n\to \infty}x_n$. 
	\end{exer}
	
	\begin{exer}\label{exe:Lipschitzrigidityl_p(l_q)}
		Let $(X,\norm{\cdot}_X)$ be a Banach space and $1\le  p<\infty$. We recall that $$\ell_p(X)=\Big\{x=(x_n)_{n=0}^\infty \in X^\bN: \ \sum_{n=0}^\infty \norm{x_n}_X^p<\infty\Big\},$$  
		equipped with $$\norm{x} := \Big(\sum_{n=0}^\infty \norm{x_n}_X^p\Big)^{1/p}$$ is a Banach space.
		The goal of this exercise is to prove some Lipschitz rigidity results for sums of classical sequence spaces.
		\begin{enumerate}
			\item Prove that if $X$ is separable, then $\ell_p(X)$ is separable.
			
			\item Let $1<p<\infty$ and $p'$ be the conjugate exponent of $p$. Show that  $\ell_p(X)^*$ is linearly isometric to $\ell_{p'}(X^*)$.
			
			\item Assume that $1<p,q<\infty$. Prove that $\ell_p(\ell_q)$ is reflexive.
			
			\item Show that $\ell_p(\ell_q)$ is Lipschitz rigid whenever $1<p,q<\infty$.
			%Let again $1<p,q<\infty$ and assume that $X$ is a Banach space Lipschitz equivalent to $\ell_p(\ell_q)$. Prove that $X$ is isomorphic to $\ell_p(\ell_q)$.
			
			\item Let $Z := \ell_1(\ell_p)$ or $Z := \ell_p(\ell_1)$ with $1<p<\infty$. Assume that $X$ is a dual Banach space that is Lipschitz equivalent to  $Z$. Prove that $X$ is isomorphic to $Z$.
			
			\item Show that $\ell_p \oplus \ell_q$ is Lipschitz rigid whenever $1<p, q <\infty$.
			%We assume  that $1<p\neq q <\infty$. Show that if a Banach space $X$ is Lipschitz equivalent to $\ell_p \oplus \ell_q$, then $X$ is isomorphic $\ell_p \oplus \ell_q$.\\
			
		\end{enumerate}
		
	\end{exer}\label{exe:LipschitzrigidityJames}
	
	\begin{proof}[Hint:]
		For 6. observe that the space $Y=\ell_p\oplus \ell_q$ has a $1$-unconditional basis and that $Y(Y)$ and $Y\oplus Y$ are isomorphic to $Y$. Alternatively, you can use the fact that $\ell_p$ and $\ell_q$ are totally incomparable and a theorem of {\`{E}}del'\v{s}te{\u\i}n and Wojtaszczyk stated in Appendix \ref{appendix:Banach}.
	\end{proof}

	\begin{exer}
		In this exercise, you will prove that James' space is Lipschitz rigid, and you are allowed to use any result about James' space from Appendix \ref{sec:James-space}.
		\begin{enumerate}
			\item Show that $\James$ has the Radon-Nikod\'{y}m property.
			\item Show that $\James$ is Lipschitz rigid.
		\end{enumerate}
	\end{exer}
	
	\begin{exer}
		\label{ex:Asplundrigidity}
		Let $X$ and $Y$ be Banach spaces and $f\colon X\to Y$ be a Lipschitz equivalence.
		\begin{enumerate}
			\item Give a simple proof of the following particular case of Proposition \ref{prop:separable-saturation-quotient}: for any separable subspace $Y_0$ of $Y$, there exists a separable subspace $X_1$ of $X$ and a separable subspace $Y_1$ of $Y$ such that $f_{\restriction_{X_1}}$ is a Lipschitz equivalence from $X_1$ onto $Y_1$.
			\item Rewrite also the proof of Theorem \ref{thm:Asplund-rigidity} in this simpler case.
		\end{enumerate}
		
	\end{exer}

	%\bigskip\noindent{\bf Exercise 7.} Consider the map $r:\ell_\infty \to c_0$ defined, for $x=(x_n)_n\in \ell_\infty$ by $r(x)_n= (|x_n|-\limsup_k{|x_k|})\text{sgn}(x_n)$ if $|x_n|>\limsup_k{|x_k|}$ and $r(x)_n=0$ otherwise.
	%\begin{enumerate}
	%\item Show that $r$ is a 2-Lipschitz retract from $\ell_\infty$ onto $\co$.

	%\item  Let $A$ be a nonempty subset of a metric space $(M,d)$ and $f\colon (A,d) \to \bR$ a $C$-Lipschitz map. Prove that the map $g$ defined by
	%$$\forall y \in M,\ \   g(y) = \inf\{f(x) +Cd(x, y),\ x\in A\},$$
	%is a $C$-Lipschitz extension of $f$ to $M$.
	%\item  Show that any $C$-Lipschitz map from $A$ to $\ell_\infty$ can be extended into a $C$-Lipschitz map from $M$ to $\ell_\infty$.
	%\item Deduce that any $C$-Lipschitz map from a subset of a metric space $M$ to $\co$ can be extended into a $2C$-Lipschitz map from  $M$ to $\co$.
	%\end{enumerate}

	%\bigskip Other ideas. $X$ nonseparable such that some open balls in $X$ are  Haar-null????

	%%%%%%%%%%%%%%%%%%%%%%%%%%%%%%%%%%%%%%%%%%%%%%%%%%%%%%%%%%%%%%%%%%%%%%%%%%%%%%%%%%%%%%%%

	\chapter{Lipschitz-free techniques in nonlinear geometry}\label{chapter:Lipschitz-free}
	
	In this chapter, we exploit yet another linearization procedure for Lipschitz maps, which is, in some sense, canonical. We explain how every metric space admits an isometric embedding into a normed vector space, whose image consists of linearly independent vectors, and such that every Lipschitz map defined on the metric space with values in a Banach space can be extended into a bounded linear map. For this reason, we call this normed vector space the Lipschitz-free space over the metric space. In different settings, this space is rightfully called, Arens-Eells, transportation cost, or Kantorovitch space. After describing the construction of the Lipschitz-free space, we prove the fundamental Godefroy-Kalton lifting theorem for separable Banach spaces and its striking application to the isometric rigidity of separable Banach spaces. In Section \ref{sec:free-BAP}, we see how Lipschitz-free spaces are instrumental in proving that the bounded approximation property is a Lipschitz invariant. In Section \ref{sec:structure-Lip-free}, we prove some structural results for Lipschitz-free spaces and apply them to the construction of important universal spaces for uniform or coarse-Lipschitz embeddings. Finally, we use Lipschitz-free techniques to show some embeddability results for stable metric spaces into reflexive Banach spaces.
	
	\section{Lipschitz-free spaces}
	
	Throughout this chapter, $(M,d)$ will be a \emph{pointed} metric space; that is, a metric space equipped with a distinguished point denoted by $0$. To this metric space we associate the space $\Lip_0(M)$ of all real-valued Lipschitz functions on $M$ that vanish at $0$, equipped with the Lipschitz norm:
	\begin{equation*}
		\norm{f}_L := \sup\Big\{\frac{\abs{f(x) - f(y)} }{d(x,y)}\colon x,y \in M,\ x\neq y\Big\}.
	\end{equation*}
	It is easy to check that $(\Lip_0(M),\norm{\cdot}_L)$ is a Banach space. Clearly, this construction does not fundamentally depend on the choice of the distinguished point. Indeed, for any $e\in M$ the map $f \mapsto f-f(e)$ is a linear isometry from $\Lip_0(M)$ onto $\Lip_e(M)$.
	
	For $x\in M$, we denote by $\delta_M(x) \in \Lip_0(M)^*$, the evaluation map at $x$, i.e. $\langle f,\delta_M(x)\rangle := f(x)$ for all $f\in \Lip_0(M)$. We refer to the map $\delta_M$ as the \emph{evaluation map}. 
	
	\begin{prop} 
		The evaluation map $\delta_M$ is an isometry from $M$ into $\Lip_0(M)^*$.
	\end{prop}
	
	\begin{proof} 
		For all $x,y \in M$ and all $f\in \Lip_0(M)$, $\abs{f(y) - f(x)}\le \norm{f}_Ld(x,y)$, so $\norm{\delta_M(y) - \delta_M(x)}_{\Lip_0(M)^*}\le d(x,y)$. On the other hand, the function $z \mapsto d(z,x)-d(0,x)$ is $1$-Lipschitz, vanishes at $0$ and satisfies $\langle f,\delta_M(y)-\delta_M(x)\rangle = d(x,y)$. So, $\norm{ \delta_M(y) - \delta_M(x) }_{\Lip_0(M)^*}\ge d(x,y)$.
	\end{proof}
	
	The Lipschitz-free space is formally defined as the closed linear span (in the dual of the Lipschitz function space) of the image of the evaluation map.
	
	\begin{defi}
		\label{def:Lipschitz-free}
		The \emph{Lipschitz-free space over} $M$ is the closed linear span of $\{\delta_M(x)\colon x\in M\}$ in $\Lip_0(M)^*$ and it is denoted by $\cF(M)$.
	\end{defi}
	
	The fundamental factorization property of $\cF(M)$ is that any Lipschitz map from $M$ to a Banach space $X$ extends, via the evaluation map, to a bounded linear map from $\cF(M)$ to $X$. 
	
	For a Banach space $X$ we denote by $\Lip_0(M,X)$ the space of Lipschitz maps from $M$ to $X$ (vanishing at $0$) equipped with the norm $\norm{f}_L := \Lip(f)$.
	
	\begin{theo}
		\label{thm:factorization-free-Lip} 
		Let $X$ be a Banach space, $M$ a pointed metric space and $f\colon M\to X$ be a Lipschitz map such that $f(0)=0$. Then, there exists a unique bounded linear map $\tilde{f}:\cF(M) \to X$ such that $f=\tilde{f}\circ \delta_M$. That is, we have the following commutative diagram:  
		$$ \xymatrix{
			M \ar[r]^f \ar[d]_{\delta_M}  & X.  \\
			{\mathcal F(M)} \ar[ur]_{\tilde{f}} 
		} 
		$$
		Moreover, the map $f\in \Lip_0(M,X) \mapsto \tilde{f}\in B(\cF(M),X)$ is a linear bijective isometry.
	\end{theo}
	
	\begin{proof} 
		We define first $\tilde{f}(\delta_M(x))=f(x)$ and extend it linearly to the linear span of the, clearly linearly independent, family $\{\delta_M(x)\colon x\in M\setminus \{0\}\}$ (note that $\delta_M(0)=0$). For $x_1,\dots,x_n \in M$ and $a_1,\dots,a_n \in \bR$, we have
		\begin{equation*}
			\Big\|\sum_{i=1}^n a_if(x_i)\Big\|_X = \sup_{x^* \in B_{X^*}} \sum_{i=1}^n a_i(x^*\circ f)(x_i) = \sup_{x^* \in B_{X^*}} \langle x^*\circ f,\sum_{i=1}^n a_i\delta_M(x_i)\rangle.
		\end{equation*}
		Since, for $x^*\in B_{X^*}$, $x^* \circ f \in \Lip_0(M)$ with $\norm{x^* \circ f}_L\le \Lip(f)$, we deduce that 
		\begin{equation*}
			\Big\|\sum_{i=1}^n a_if(x_i)\Big\|_X \le \Lip(f)\Big\|\sum_{i=1}^n a_i\delta_M(x_i)\Big\|_{\Lip_0(M)^*}.
		\end{equation*}
		So, $\|\tilde{f}\|\le \Lip(f)$. On the other hand, for all $x,y\in M$, $\|\tilde{f}(\delta_M(y)) - \tilde{f}(\delta_M(x))\|_X = \norm{f(y) - f(x)}_X$. This, combined with the fact that $\delta_M$ is an isometry, implies that $\|\tilde{f}\|\ge \Lip(f)$. Then, by density, $\tilde{f}$ uniquely extends to a bounded linear map, still denoted by $\tilde{f}\colon \cF(M) \to X$ with $\|\tilde{f}\| =  \Lip(f)$. So, $f\mapsto \tilde{f}$ is an isometry from $\Lip_0(M,X)$ into $B(\cF(M),X)$, which is obviously linear. Finally, if $T \in \cal B(\cF(M),X)$, we clearly have that $T=\tilde{f}$, where $f\colon x\mapsto T(\delta_M(x)) \in \Lip_0(M,X)$.
	\end{proof}
	
	Applying the previous statement with $X=\bR$, we immediately deduce a representation of  $\cF(M)$ as a natural predual of $\Lip_0(M)$. More precisely, we have:
	
	\begin{coro}
		\label{cor:dual-free=Lip}
		For $f\in \Lip_0(M)$, let $\Lambda(f) \in \cF(M)^*$ be the unique bounded linear map from $\cF(M)$ to $\bR$ such that $\Lambda(f) \circ \delta_M=f$. Then, $\Lambda$ is a linear isometry from $\Lip_0(M)$ onto $\cF(M)^*$.
	\end{coro}
	
	\begin{rema} 
		One key consequence of Corollary \ref{cor:dual-free=Lip} is that the topology of pointwise convergence on $M$ and the weak$^*$ topology induced by $\cF(M)$ on $\Lip_0(M)$ coincide on bounded subsets of $\Lip_0(M)$. 
	\end{rema}
	
	\noindent {\bf Examples.} As Exercise \ref{ex:Lipschitzfreeexamples} we leave to show that $\cF(\N)$ is linearly isometric to $\ell_1$ and $\cF(\R)$ is linearly isometric to $L_1$. 
	
	\medskip 
	The following corollary shows that the Lipschitz-free space over a subset of $M$ can be canonically seen as a subspace of $\cF(M)$. More precisely, we have
	
	\begin{coro}
		\label{cor:free-subset} 
		If $S$ is a subset of $M$ containing $0$, then the map $I_S$ defined, for all $x_1,\dots,x_n\in S$ and $a_1,\dots,a_n \in \bR$, by
		\begin{equation*}
			I_S\Big(\sum_{i=1}^n a_i\delta_S(x_i)\Big)=\sum_{i=1}^n a_i\delta_M(x_i),
		\end{equation*}
		extends to a linear isometry from $\cF(S)$ onto $\cF_M(S)$, the closed linear span of $\delta_M(S)$ in $\cF(M)$. 
	\end{coro}
	
	\begin{proof} 
		Let $x_1,\dots,x_n\in S$ and $a_1,\dots,a_n \in \bR$. The fact that 
		\begin{equation*}
			\Big\| \sum_{i=1}^n a_i\delta_S(x_i)\Big\|_{\cF(S)} \ge \Big\|\sum_{i=1}^n a_i\delta_M(x_i)\Big\|_{\cF(M)}
		\end{equation*}
		is clear, since restrictions to $S$ of $1$-Lipschitz functions on $M$ are $1$-Lipschitz on $S$. For the other inequality, pick $f \in \Lip_0(S)$ with $\Lip(f)=1$ such that 
		\begin{equation*}
			\Big\|\sum_{i=1}^n a_i\delta_S(x_i)\Big\|_{\cF(S)}=\sum_{i=1}^n a_if(x_i).
		\end{equation*}
		Then, by Proposition \ref{pro:Inf-Conv}, there exists $g \in \Lip_0(M)$, an extension of $f$ with $\Lip(g)=1$. Thus,
		\begin{equation*}
			\Big\|\sum_{i=1}^n a_i\delta_M(x_i)\Big\|_{\cF(M)}\ge \sum_{i=1}^n a_ig(x_i)=\Big\|\sum_{i=1}^n a_i\delta_S(x_i)\Big\|_{\cF(S)}.
		\end{equation*}
		
	\end{proof}
	
	The next corollary is concerned with the factorization of Lipschitz maps between metric spaces.
	
	\begin{coro}
		\label{cor:factorization-Lipschitz-maps} Let $M$ and $N$ be two pointed metric spaces and $f\colon M \to N$ be a Lipschitz map such that $f(0_M)=0_N$. 
		\begin{enumerate}[(i)]
			\item There exists a unique bounded linear map $\widehat{f}\colon \cF(M) \to \cF(N)$ such that 
			$$\delta_N \circ f = \widehat{f} \circ \delta_M.$$
			In other words, we have the following commutative diagram:
			$$\xymatrix{
				M \ar[r]^f \ar[d]_{\delta_M}  & N \ar[d]^{\delta_N} \\
				\mathcal F(M) \ar[r]_{\widehat{f}} & \mathcal F(N).
			}
			$$
			Furthermore, $\|\widehat{f}\|=\Lip(f)$.
			\item If $f$ is a Lipschitz equivalence, then $\widehat{f}$ is a linear isomorphism.
			\item If $f$ is a bi-Lipschitz embedding, then $\widehat{f}$ is a linear isomorphic embedding.
		\end{enumerate}
	\end{coro}
	
	\begin{proof} 
		$(i)$ This is immediate by applying Theorem \ref{thm:factorization-free-Lip} to the map $\delta_N \circ f\colon M\to \cF(N)$. 
		
		$(ii)$ Apply $(i)$ to $f$ and $f^{-1}$ and check that $\widehat{f^{-1}}=(\widehat{f})^{-1}$. 
		
		$(iii)$ Let $S:=f(M)$. By $(ii)$, $\widehat{f}$ is a linear isomorphism from $\cF(M)$ onto $\cF(S)$ seen as a subspace of $\cF(N)$ (Corollary \ref{cor:free-subset}).
	\end{proof}
	
	\begin{rema}
		\label{rem:composition-operator} Note that the adjoint operator of $\widehat{f}$ is identified with the composition operator $g\in \Lip_0(N) \mapsto g \circ f \in \Lip_0(M)$.
	\end{rema}
	
	We now turn to the particular case when the metric space is a Banach space $X$ and we naturally choose the distinguished element $0$ to be the origin of $X$. In that situation, there is another fundamental object denoted by $\beta_X$ and called the \emph{barycenter map}. 
	
	\begin{prop} 
		There exists a unique bounded linear map $\beta_X\colon \cF(X) \to X$ such that $\beta_X \circ \delta_X = Id_X$. Furthermore, $\norm{\beta_X}=1$. 
	\end{prop}
	
	\begin{proof} 
		Just apply Theorem \ref{thm:factorization-free-Lip} to the map $f=Id_X$. 
	\end{proof}
	
	The terminology ``barycenter'' comes from the fact that $\beta_X$ maps any finitely supported probability measure on $X$ to the barycenter of the points in $X$ as follows:
	$$\beta_X\Big(\sum_{i=1}^k p_i \delta_X(x_i)\Big) = \sum_{i=1}^k p_i x_i.$$
	
	\begin{rema}
		\label{rem:section-free}
		The linear map $\beta_X$ is an exact quotient map from $\cF(X)$ onto $X$, and $\delta_X$ is an isometric nonlinear lifting of $\beta_X$. Note also that there is a natural Lipschitz equivalence between $\cF(X)$ and $X \times \mathrm{Ker}(\beta_X)$ (see Exercise \ref{exe:section-free} for more details). This last observation will be used later to construct Banach spaces that are Lipschitz equivalent.
	\end{rema}
	
	\section{The Godefroy-Kalton lifting Theorem}
	
	The main result of this section is a fundamental lifting theorem due to G. Godefroy and N. Kalton \cite{GodefroyKalton2003}. Let us first try to explain the idea behind their result. We have already observed that the evaluation map $\delta_X$ is an isometric nonlinear lifting of the barycenter map $\beta_X$. It is natural to wonder whether there exists an isometric and \emph{linear} lifting. Observing that for all $x,y\in X$ and all $t \in (0,\infty)$,
	\begin{equation*}
		\beta_X\Big(\frac{\delta_X(y+tx)-\delta_X(y)}{t}\Big)=x, 
	\end{equation*}
	it is clear that if there was a point $y\in X$ where we could differentiate the evaluation map $\delta_X$, the Gateaux derivative at this point would be an isometric and linear lifting of $\beta_X$. Unfortunately, the evaluation map is \emph{nowhere} Gateaux differentiable (see Exercise \ref{ex:G-diff-delta}). In the case where the Banach space is separable, Godefroy and Kalton were able to overcome this roadblock by using a very elegant averaging argument.
	%of the map $y\mapsto \delta_X(y+x)-\delta_X(y)$. 
	
	\begin{theo}
		\label{thm:GK-lifting-free}
		Let $X$ be a \underline{separable} Banach space. Then, there exists a linear isometry $T$ from $X$ into $\cF(X)$ such that $\beta_X T = Id_X$ and $T(X)$ is $1$-complemented in $\cF(X)$. 
	\end{theo}
	
	\begin{proof}
		So assume that $(x_n)_{n=1}^\infty$ is  linearly independent and total in $X$ with $\norm{x_n}_X=2^{-n}$. We denote by $V$ the linear span of $(x_n)_{n=1}^\infty$ and let $H :=[0,1]^{\bN}$ and $\lambda$ be the product measure on $H$ of the Lebesgue measure of $[0,1]$. Then, for $n\in \bN$, we let $\bN_n := \bN\setminus \{n\}$, $H_n := [0,1]^{\bN_n}$ and $\lambda_n$ be the product measure on $H_n$ of the Lebesgue measure of $[0,1]$. For all $t \in H$, let
		\begin{equation*}
			\bar{x}(t) := \sum_{k=1}^\infty t_kx_k    
		\end{equation*}
		and, for all $t \in H_n$,
		\begin{equation*}
			\bar{x}_n(t) := \sum_{k\neq n} t_k x_k.
		\end{equation*}
		Finally, for $n\in \bN$, let 
		\begin{equation*}
			\mu_n :=\int_{H_n} \big(\delta_X(x_n+\bar{x}_n(t))-\delta_X(\bar{x}_n(t))\big)\,d\lambda_n(t).
		\end{equation*}
		It is easily checked that this $\cF(X)$-valued Bochner integral is well defined and, since it commutes with the bounded operator $\beta_X$ that $\beta_X(\mu_n) = x_n$. Then, there exists a unique linear map $T$ from $V$ to $\cF(X)$ such that $T(x_n) = \mu_n$ for all $n\in \bN$. Note that $\beta_XT = Id$ on $V$. For all $v\in V$, we will now show that $\norm{Tv}_{\cF(X)}\le \norm{v}_X$, i.e. for all $f\in \Lip_0(X)$, 
		\begin{equation}
			\label{eq:GKL-lifting-free}
			\big|\langle f,Tv\rangle\big|\le \norm{f}_{L}\norm{v}_X.
		\end{equation}
		Once \eqref{eq:GKL-lifting-free} is proved, we can immediately conclude the proof since $T$ can be extended to a bounded linear map from $X$ to $\cF(X)$ such that $\beta_XT = Id_X$ and $\norm{T}\le 1$. Moreover, since $\norm{\beta_X}=1$, $T$ is actually a linear isometry, and it automatically follows that $T\beta_X$ is a norm one projection from $\cF(X)$ onto $T(X)$.
		
		In order to prove \eqref{eq:GKL-lifting-free}, we consider first the case of a function $f\in \Lip_0(X)$ that is everywhere Gateaux differentiable (and we denote by $D_f(y)$ its Gateaux derivative at $y$). In this case, it follows from the fundamental theorem of Calculus and Fubini's theorem that
		\begin{align*} 
			\langle f, Tx_n\rangle & =\int_{H_n} \big(f(x_n+\bar{x}_n(t))-f(\bar{x}_n(t))\big)\,d\lambda_n(t)\\
			&= \int_{H_n}\int_0^1\langle D_f(\bar{x}_n(t)+ sx_n),x_n\rangle\,ds\,d\lambda_n(t) = \int_H \langle D_f(\bar{x}(t)),x_n\rangle \,d\lambda(t).
		\end{align*}
		We deduce that for all $v\in V$
		\begin{equation*}
			\langle f,Tv\rangle = \int_H \langle D_f(\bar{x}(t)), v\rangle \,d\lambda(t),
		\end{equation*}
		and hence $\big|\langle f,Tv\rangle\big|\le \norm{f}_{L}\norm{v}_X$ whenever $f$ is everywhere Gateaux differentiable.
		Note that, for a fixed $v :=\sum_{i=1}^n a_ix_i \in V$, the above reasoning applies to $f\in \Lip_0(X)$ such that for all $y\in X$, the map $z\mapsto f(y+z)$ is differentiable everywhere on $V_n$, the linear span of $\{x_1,\dots,x_n\}$. So, fix $v :=\sum_{i=1}^n a_ix_i \in V_n$, $f\in \Lip_0(X)$. Rademacher's theorem (Theorem \ref{thm:Rademacher}) ensures that for any $y\in X$, the map $z\mapsto f(y+z)$ is differentiable almost everywhere on $V_n$. The point is now to find a sequence $(g_k)_k$ of maps that are in $\Lip_0(X)$, differentiable everywhere on $V_n$ and such that $(g_k)_k$ converges to $f$ for the weak$^*$ topology induced by $\cF(X)$ on $\Lip_0(X)$. Indeed, for such maps we deduce from the observation above that $\big|\langle g_k,Tv\rangle\big|\le \norm{f}_{L}\norm{v}_X$ and it follows from weak$^*$ convergence that $\big|\langle f,Tv\rangle\big|\le \norm{f}_{L}\norm{v}_X$. 
		The $g_k$ can be constructed by smoothing $f$ as follows: for $k\in \bN$ and $y\in X$ let
		\begin{equation*}
			g_k(y) := \int_{[0,1]^n} f\Big(y+\frac{1}{k}\sum_{i=1}^n t_ix_i \Big)\,dt_1\cdots dt_n - \int_{[0,1]^n} f\Big(\frac{1}{k}\sum_{i=1}^n t_ix_i \Big)\,dt_1\cdots dt_n.
		\end{equation*}
		It is clear that $g_k \in \Lip_0(X)$ with $\norm{g_k}_L \le \norm{f}_L$ and $\lim_k\norm{g_k-f}_\infty=0$. The fact that, for any $y\in X$, the map  $z\mapsto g_k(y+z)$ is differentiable everywhere on $V_n$ is a consequence of the dominated convergence theorem. Finally, since the sequence $(g_k)_k$ is bounded in $\Lip_0(X)$ and pointwise converging to $f$ on $X$, it is converging to $f$ for the weak$^*$ topology induced by $\cF(X)$ on $\Lip_0(X)$.  
	\end{proof}
	
	The above lifting property can be seen as a universal lifting property in the separable setting. The following corollary illustrates this.
	
	\begin{coro}
		\label{cor:GK-lifting}
		Let $X$ be a Banach space, $Y$ be a closed subspace of $X$ so that $X/Y$ is separable and $Q\colon X \to X/Y$ be the quotient map. If there exists $\varphi\colon X/Y \to X$ Lipschitz so that $Q\circ \varphi = Id_{X/Y}$, then there exists $\Phi\colon  X/Y \to X$ linear and bounded such that $Q\circ \Phi = Id_{X/Y}$ and $\norm{\Phi}\le \Lip(\varphi)$.
	\end{coro}
	
	\begin{proof} 
		We may assume that $\varphi(0)=0$ and according to Theorem \ref{thm:factorization-free-Lip} there exists a bounded linear map $\tilde{\varphi}$ from $\cF(X/Y)$ to $X$ such that $\tilde{\varphi}\circ\delta_{X/Y}=\varphi$ and $\norm{\tilde{\varphi}} = \Lip(\varphi)$. It follows that $Q\circ \tilde{\varphi} \circ \delta_{X/Y} = Q\circ \varphi = Id_{X/Y} = \beta_{X/Y} \circ \delta_{X/Y}$. Since $Q\tilde{\varphi}$ and $\beta_{X/Y}$ are linear, bounded and coincide on a total subset of $\cF(X/Y)$, we deduce that  $Q\tilde{\varphi}=\beta_{X/Y}$. Theorem \ref{thm:GK-lifting-free} ensures the existence of a linear isometry $T\colon X/Y \to \cF(X/Y)$ such that $\beta_{X/Y}T = Id_{X/Y}$. Thus, we can conclude our proof by letting $\Phi := \tilde{\varphi}T$.
	\end{proof}
	
	The fact that, in the separable setting, Lipschitz sections can be linearized rules out a natural approach to construct separable Banach spaces that are Lipschitz equivalent but not linearly isomorphic. The lifting method will be discussed in Section \ref{sec:lifting-method}.
	
	\begin{rema} 
		The proof of Corollary \ref{cor:GK-lifting} can be carried out without using the theory of Lipschitz-free spaces. This direct approach is detailed in the textbook \cite{AlbiacKalton2016} (Theorem~14.3.3). 
	\end{rema}

	The separability assumption is crucial in Theorem \ref{thm:GK-lifting-free}. Nevertheless, we can state a local version of it, which we shall use in the next section. 
	
	\begin{prop}
		\label{prop:GK-local-lifting} 
		Let $X$ be a Banach space and suppose that $E$ is a finite-dimensional subspace of $X$. Then, there exists a linear isometry $T\colon X\to \cF(X)^{**}$ such that $T(E)\subset \cF(E)$ and $\beta_X^{**}T = Id_X$.
	\end{prop}
	
	\begin{proof}
		Let $(e_1,\dots,e_n)$ be a basis of $E$. We define $U:E \to \cF(E)\subset \cF(X)$ to be linear and such that for $1\le k \le n$,
		\begin{equation*}
			U(e_k) := \int_{\prod_{j\neq k}[0,1]}\Big(\delta_X\Big(e_k+\sum_{j\neq k}t_je_j\Big)-\delta_X\Big(\sum_{j\neq k}t_je_j\Big)\Big) \prod_{j\neq k}dt_j.
		\end{equation*}
		Following the proof of Theorem \ref{thm:GK-lifting-free}, we show that $U$ is a linear isometry from $E$ into $\cF(E)$ such that $\beta_XU = Id_E$. 
		
		Let us fix $\eps>0$ and a finite-dimensional subspace $G$ of $X$ containing $E$. We reproduce the above construction as follows. First we complete the basis $(e_1,\dots,e_n)$ of $E$ into a basis $(g_1,\dots,g_m)$ of $G$ so that $\sum_{i=n+1}^m \norm{g_i}<\eps$ and we define $T_{G,\eps}\colon G\to \cF(G)\subset \cF(X)$ to be linear and such that for $1\le k\le m$,
		\begin{equation*}
			T_{G,\eps}(g_k) := \int_{\prod_{j\neq k}[0,1]}\Big(\delta_X\Big(g_k+\sum_{j\neq k}t_jg_j\Big)-\delta_X\Big(\sum_{j\neq k}t_jg_j\Big)\Big) \prod_{j\neq k}dt_j.
		\end{equation*}
		Again, $T_{G,\eps}$ is a linear isometry from $G$ into $\cF(X)$ satisfying $\beta_XT_{G,\eps} = Id_G$. If $x\in X \setminus G$, then set $T_{G,\eps}(x) := 0$ and note that for all $1\le j\le n$,
		\begin{equation*}
			\norm{T_{G,\eps}(e_j)-U(e_j)}\le 2\eps.
		\end{equation*}
		It follows from Tychonoff's theorem that there exists a subnet $(T_\alpha)_{\alpha \in A}$ of the net $\{T_{G,\eps}\colon E\subset G,\ \eps>0\}$ (directed in the classical way) such that for any $x\in X$, $T_\alpha(x)$ is converging to an element $T(x) \in \cF(X)^{**}$, in the weak$^*$ topology of $\cF(X)^{**}$.  It is now easy to check that $T$ is linear, of norm at most one, extends $U$ (so $T(E)\subset \cF(E)$), satisfies $\beta_X^{**}T=Id_X$, and therefore $T$ is actually an isometry. 
	\end{proof}

	\section{Isometric rigidity of separable Banach spaces}
	\label{sec:isometric-rigidity}
	
	The ancestor of all nonlinear rigidity results is the celebrated Mazur-Ulam Theorem \cite{MazurUlam1932}. Its original proof relies in part on the notion of metric midpoints. In this section, we will recall this proof. Let us start with a definition.
	
	\begin{defi} 
		Let $(M,d)$ be a metric space. For $x,y \in M$, the set of \emph{metric midpoints} of $x$ and $y$ is the set, denoted by $\Mid(x,y)$ and defined as follows:
		\begin{equation*}
			\Mid(x,y) := \Big\{z\in M\colon d(x,z)=d(y,z)=\frac12 d(x,y)\Big\}
		\end{equation*}
	\end{defi}
	
	The metric midpoint set can be empty, but for normed spaces, we always have $\frac{x+y}{2} \in \Mid(x,y)$. We now state and prove the Mazur-Ulam Theorem.
	
	\begin{theo}[Mazur-Ulam Theorem, 1932]
		\label{MazurUlam}
		Let $(X,\norm{\cdot}_X)$ and $(Y,\norm{\cdot}_Y)$ be two normed spaces. If $f$ is an isometry from $X$ \underline{onto} $Y$ such that $f(0)=0$, then $f$ is linear.
		In particular, two normed spaces that are isometrically equivalent are automatically linearly isometric.
	\end{theo}
	
	\begin{proof} 
		The density of dyadic numbers in $[0,1]$ and a continuity argument ensure that it is enough to prove that midpoints are sent to midpoints, i.e. for all $x\neq x'$ in $X$, $f(\frac{x+x'}{2})=\frac{f(x)+f(x')}{2}$.
		
		So, let us fix $x\neq x' \in X$. Let $K_1 :=\Mid(x,x')$ and define inductively for all $n\in \bN$:
		\begin{equation*}
			K_{n+1} := \Big\{u \in K_n\colon \forall v\in K_n,\ \norm{u-v}_X\le \frac{\diam(K_n)}{2}\Big\}.
		\end{equation*}
		We first want to prove that $\bigcap_{n=1}^\infty K_n=\{\frac{x+x'}{2}\}$. Since the diameter of $K_n$ is tending to $0$, it is enough to show that for all $n\in \bN$, $\frac{x+x'}{2}\in K_n$. Assuming, as we may, that $x'=-x$ we will show that $0\in K_n$. In fact, we will show by induction that $K_n$ contains $0$ and is symmetric. This is clearly true for $K_1 =\{z\colon \norm{z-x}_X = \norm{z+x}_X = \norm{x}_X\}$. Now, assume that $K_n$ is symmetric and $0\in K_n$. For any $u \in K_n$, $-u \in K_n$ and hence $\norm{2u}_X \le \diam(K_n)$, which implies that $0\in K_{n+1}$. Consider now $u\in K_{n+1}$. We have that $u\in K_n$ and for any $v\in K_n$, $\norm{u-v}_X\le \frac12 \diam(K_n)$. Since $K_n$ is symmetric, this implies that $-u\in K_n$ and for any $v\in K_n$, $\norm{u+v}_X\le \frac12 \diam(K_n)$ and therefore that $-u \in K_{n+1}$. This concludes our induction.
		
		In $Y$, we similarly define $L_1 := \Mid(f(x),f(x'))$ and by induction
		\begin{equation*}
			L_{n+1}=\Big\{s \in L_n \colon \forall t\in L_n,\ \norm{s-t}_Y\le \frac{\diam(L_n)}{2}\Big\}.
		\end{equation*}
		We deduce from the previous proof that $\bigcap_{n=1}^\infty L_n=\{\frac{f(x)+f(x')}{2}\}$. It now follows easily from the fact that $f$ is a bijective isometry that for all $n\in \bN$, $f(K_n)=L_n$. Taking the intersections over $n$ gives that $f(\frac{x+x'}{2})=\frac{f(x)+f(x')}{2}$.
	\end{proof}
	
	\begin{rema} 
		The assumption that the image of $f$ is linear is crucial in this result. Indeed, it is very easy to construct nonlinear isometries. If you pick your favorite $1$-Lipschitz and nonlinear function $h\colon \bR \to \bR$ (e.g. the sine function), then the map $f\colon \bR \to \ell_\infty^2$ defined by $f(t) := (t,h(t))$ is clearly a nonlinear isometry.
	\end{rema}

	We continue this section with an extension of the Mazur-Ulam theorem. It is a description of the into isometries between Banach spaces due to T. Figiel \cite{Figiel1968}. 
	
	\begin{theo}
		\label{thm:Figiel} 
		Let $X$ and $Y$ be two Banach spaces and $f\colon X\to Y$ a (nonlinear) isometric embedding such that $f(0)=0$. Then, for all $x_1,\dots,x_n \in X$ and all $a_1,\dots,a_n \in \bR$,
		$$\Big\|\sum_{i=1}^n a_i f(x_i)\Big\|_Y \ge \Big\|\sum_{i=1}^n a_ix_i\Big\|_X.$$
		If $Z$ is the closed linear span of $f(X)$, then there exists a unique linear operator $Q\colon Z\to X$ such that $Q\circ f = Id_X$. Moreover, $\norm{Q}=1$. 
	\end{theo}
	
	\begin{proof} 
		The second part of the statement is immediate. Indeed, it follows from the inequality that the map $Q\colon\sum_{i=1}^n a_i f(x_i) \mapsto \sum_{i=1}^n a_ix_i$ is well defined and of norm one on the linear span of $f(X)$ and therefore uniquely extends to a bounded operator from $Z$ to $X$ that we still denote by $Q$. 
		
		\medskip For the proof of the inequality, we start with a one-dimensional lemma.
		
		\begin{lemm} 
			Let $x\in S_X$. Then, there exists $y^*_x\in S_{Y^*}$ such that for all $t \in \bR$,
			\begin{equation*}
				\langle y^*_x, f(tx)\rangle=t.
			\end{equation*}
		\end{lemm}
		
		\begin{proof} 
			Since $f$ is an isometry, for all $k\in \bN$ we can find $y^*_k\in S_{Y^*}$ such that $\langle y^*_k,f(kx)-f(-kx)\rangle=2k$. Since $f(0)=0$ and $y^*_k \circ f$ is $1$-Lipschitz, we necessarily have that $\langle y^*_k, f(tx)\rangle=t$ for all $t\in [-k,k]$. Finally, we conclude the proof by picking $y^*_x$ to be a weak$^*$ accumulation point of $(y^*_k)_k$ in $B_{Y^*}$, but actually necessarily in $S_{Y^*}$.
		\end{proof}
		
		We now turn to the finite-dimensional version of the above lemma.
		
		\begin{lemm}
			Let $x_1,\dots,x_n \in X$ and $E$ be the linear span of $\{x_1,\dots,x_n\}$. Let $x\in S_E$ be such that $\norm{\cdot}_E$ is differentiable at $x$. Then, there exists a unique $v^*_x\in S_{E^*}$ such that $v^*_x(x)=1$. Furthermore, we have that $(y^*_x \circ f)(y)=v^*_x(y)$ for all $y\in E$.
		\end{lemm}
		
		\begin{proof} 
			The existence of $v^*_x$ is just the Hahn-Banach theorem. Assume that $v^*,w^* \in S_{E^*}$ are such that $v^*(x)=w^*(x)=1$. Using the differentiability of $\norm{\cdot}_E$ at $x$, we can write, for $h\in E$: 
			\begin{equation*}
				2+(v^*-w^*)(th)=v^*(x+th)+w^*(x-th)\le \norm{x+th} + \norm{x-th} = 2+o(|t|).
			\end{equation*}
			Letting $t$ tend to $0$, we deduce that $v^*(h)=w^*(h)$ and thus the uniqueness of $v^*_x$ follows. Note that this property is satisfied by the derivative of the norm at $x$, which therefore coincides with $v^*_x$.
			
			Letting $\phi :=y^*_x \circ f$, we have that $\phi$ is $1$-Lipschitz and, by the previous lemma, $\phi(tx)=t$ for all $t\in \bR$. So, for $t\in \bR,\ t\neq 0$ and $y\in E$ we can write 
			\begin{align*}
				\abs{t-\phi(y)} = \abs{\phi(tx)-\phi(y)} \le \norm{tx-y} = \abs{t}\Big(1-v^*_x\big(\frac{y}{t}\big)+o\big(\frac{1}{\abs{t}}\big)\Big)
			\end{align*}
			We get $\phi(y)\ge v^*_x(y)$ by letting $t$ tend to $+\infty$ and $\phi(y)\le v^*_x(y)$ by letting $t$ tend to $-\infty$. This finishes the proof of this lemma.
			
		\end{proof}
		
		We can now conclude the proof of Figiel's theorem. Fix $x_1,\dots,x_n \in X$, denote by $E$ the linear span of $\{x_1,\dots,x_n\}$ and by $\cal S$ the set of all $x\in S_E$ such that $\norm{\cdot}_E$ is differentiable at $x$. It follows from the Rademacher Theorem \ref{thm:Rademacher} and an homogeneity argument that $\cal S$ is dense in $S_E$ and thus that the family $(v^*_x)_{x\in \cal S}$ is $1$-norming for $E$. We can deduce that 
		\begin{align*}
			\Big\|\sum_{i=1}^n a_i f(x_i)\Big\|_Y \ge \sup_{x\in \cal S}\Big|y^*_x\Big(\sum_{i=1}^n a_i f(x_i)\Big)\Big| = \sup_{x\in \cal S}\Big|v^*_x\Big(\sum_{i=1}^n a_i x_i\Big)\Big|=\Big\|\sum_{i=1}^n a_i x_i\Big\|_X.
		\end{align*}
		
	\end{proof}
	
	In \cite{GodefroyKalton2003}, Godefroy and Kalton combined their lifting theorem with Figiel's theorem to obtain the following striking result, which says, roughly speaking that any isometry between separable Banach spaces induces a linear isometry.
	
	\begin{theo}
		\label{thm:isometric-GK}
		%Let $X$ and $Y$ be Banach spaces. Assume that $X$ is separable and isometric to a subset of $Y$. Then, $X$ is linearly isometric to a subspace of $Y$.  
		A separable Banach space that isometrically embeds into another Banach space automatically admits a linear isometric embedding into it.
	\end{theo}

	\begin{proof}
		Let $X$ and $Y$ be Banach spaces and assume that $X$ is separable. Let $f\colon X \to Y$ be an isometry. We may assume without loss of generality that $f(0)=0$ and that $Y$ is the closed linear span of $f(X)$. By Figiel's theorem (Theorem \ref{thm:Figiel}), there exists a quotient map $Q \colon Y \to X$ such that $Q \circ f =Id_X$ and $\norm{Q}\le 1$ (in other words, $Q$ is an exact quotient map). It now follows from Corollary \ref{cor:GK-lifting} and the separability of $X$ that there exists $T\colon X \to Y$ linear and bounded such that $\norm{T}\le 1$ and $QT=Id_X$. Since $\norm{Q}\le 1$ and $\norm{T}\le 1$, we deduce that $T$ is a linear isometry. 
	\end{proof}
	
	\begin{rema} 
		Let us mention that the above proof contains more information: the nonlinear isometric embedding $f$ can be decomposed into the ``direct sum'' of the linear isometry $J$ and a nonlinear map. Indeed, the map $P:=JQ$ is a norm one projection from $Y$ onto $J(X)$ and $Pf=J$.
	\end{rema}
	
	\begin{rema} 
		The conclusion of Theorem \ref{thm:isometric-GK} is false if $X$ is not separable. Indeed, it follows from Proposition 4.1 in \cite{GodefroyKalton2003} that if $H$ is a nonseparable Hilbert space, then $H$ does not linearly embed into $\cF(H)$.  
	\end{rema}
	
	\section{Lipschitz rigidity of the bounded approximation property}
	\label{sec:free-BAP}
	
	A Banach space $X$ has the \emph{approximation property} (AP in short) if the identity operator can be approximated arbitrarily well on any compact subset by a bounded linear operator of finite rank.
	%, i.e. for any compact subset $K$ of $X$ and any $\eps >0$, there exists $T \in B(X)$ such that  $T$ has a finite rank and $\norm{Tx-x}\le \eps$, for all $x\in K$. 
	The bounded approximation property is a quantification of the approximation property.
	
	\begin{defi} 
		\label{def:BAP}
		Let $X$ be a Banach space and $\lambda \in [1,\infty)$. We say that $X$ has the \emph{$\lambda$-bounded approximation property} ($\lambda$-BAP for short) if for any compact subset $K$ of $X$ and any $\eps >0$, there exists $T \in B(X)$ such that $\norm{T}\le \lambda$, $T$ has  finite rank and $\norm{Tx-x}\le \eps$ for all $x\in K$.
		
		We say that $X$ has the \emph{bounded approximation property} (BAP for short), if it has the $\lambda$-BAP for some $\lambda \in [1,\infty)$. The $1$-BAP is called the \emph{metric approximation property} and is abbreviated as MAP.
	\end{defi}
	
	It is a well-known and elementary fact (see Exercise \ref{ex:BAP_lemma}) that the compact sets in Definition \ref{def:BAP} can be replaced by finite sets, or even finite subsets of a given set with dense linear span,  without altering the concept. It is also easy to see that a Banach space with a Schauder basis has the bounded approximation property by considering the basis projections. None of the approximation properties above passes to subspaces. However, the approximation property and its bounded version are inherited by \emph{complemented} subspaces. Therefore, it follows immediately from Theorem \ref{thm:GK-lifting-free} that if $\cF(X)$ has the bounded approximation property and $X$ is separable, then $X$ has the bounded approximation property. The separability assumption can be lifted with the help of Proposition \ref{prop:GK-local-lifting}. The following result, like all the results in this section about the BAP for Lipschitz-free spaces, is due to Godefroy and Kalton \cite{GodefroyKalton2003}.
	
	\begin{theo}
		\label{thm:BAP-LF->space}
		If the Lipschitz-free space of a Banach space has the $\lambda$-bounded approximation property, then the Banach space itself has the $\lambda$-bounded approximation property.    
	\end{theo}
	
	\begin{proof}
		Let $X$ be a Banach space such that $\cF(X)$ has the $\lambda$-BAP. Fix $x_1,\dots,x_n \in X$ and $\eps>0$ and denote by $E$ the linear span of $\{x_1,\dots,x_n\}$. By Proposition \ref{prop:GK-local-lifting}, we can find an isometry $T\colon X \to \cF (X)^{**}$ so that $T(E)\subset \cF(E)$ and $\beta_X^{**}T = Id_X$. Since $\cF(X)$ has the $\lambda$-BAP, we can find a finite-rank operator $R\colon \cF(X)\to \cF(X)$ such that $\norm{R}\le \lambda$ and $\norm{RTx_j-Tx_j}\le \eps$, for $1\le j \le n$. Then, consider $S := \beta_X^{**}R^{**}T$. Since $R$ has a finite rank, $R^{**}$ actually maps $\cF(X)^{**}$ to $\cF(X)$ and has finite rank. So, $S\colon X\to X$ is a finite-rank operator with $\norm{S}\le \lambda$ and $\norm{Sx_j-x_j} = \norm{\beta_X^{**}R^{**}Tx_j - \beta_X^{**}Tx_j} = \norm{\beta_X^{**}(RTx_j-Tx_j)}\le \eps$, for $1\le j \le n$. This finishes the proof.
	\end{proof}
	
	The fact that the bounded approximation property passes from a Banach space to its Lipschitz-free space is much more delicate and requires significantly more work. In fact, a nonlinear weakening of the bounded approximation property suffices to guarantee the bounded approximation property of the Lipschitz-free space. 
	
	\begin{defi} 
		\label{def:LBAP}
		Let $X$ be a Banach space and $\lambda \ge 1$. We say that $X$ has the \emph{$\lambda$-Lipschitz bounded approximation property} ($\lambda$-LBAP) if for every compact subset $K$ of $X$ and every $\eps>0$  there  exists a Lipschitz  map $f\colon X \to X$ with finite-dimensional range such that $\Lip(F)\le \lambda$ and $\norm{f(x)-x}\le \eps$ for all $x\in K$. 
	\end{defi}
	
	Note that in the above definition, we may equivalently ask that, moreover, $f(0)=0$. It is plain that the bounded approximation property implies the Lipschitz bounded approximation property. The following theorem is due to Godefroy and Kalton \cite{GodefroyKalton2003}.
	
	\begin{theo}
		\label{thm:LBAP->LF-BAP}
		If a Banach space has the $\lambda$-Lipschitz bounded approximation property, then its Lipschitz-free space has the $\lambda$-bounded approximation property.
	\end{theo}
	
	Interestingly, an immediate consequence of Theorem \ref{thm:LBAP->LF-BAP} and Theorem \ref{thm:BAP-LF->space} is that the bounded approximation property and its Lipschitz version coincide.
	
	\begin{rema}
		It might have been more natural in the definition of the LBAP to avoid any reference to the linear structure. For instance, one could have replaced the finite-dimensional range condition with the condition that the image of the unit ball is relatively compact. This delicate point is discussed by Godefroy in \cite{Godefroy2020}, where it is shown that the finite-dimensional range condition is necessary for Theorem \ref{thm:LBAP->LF-BAP} and where the notion of so-called Lipschitz approximable Banach spaces is explained and explored. 
	\end{rema}

	The main step in the proof of Theorem \ref{thm:LBAP->LF-BAP} is to treat the case of Lipschitz-free spaces over finite-dimensional Banach spaces. For the finite-dimensional case, it suffices to globally approximate the evaluation map by Lipschitz maps with finite-dimensional ranges. This is done by first finding a weak approximation on a ball and then by using a convex combination argument together with a proper rescaling. The weak approximation argument is the content of the next lemma.
	
	\begin{lemm}
		\label{lem:weak-approx-Fejer}
		Let $E$ be a finite-dimensional Banach space. There exists $r_0>0$ such that for any $\eps>0$, there is a sequence of Lipschitz functions $L_n \colon E \to \cF(E)$ with $L_n(0)=0$ and such that:
		\begin{enumerate}[(i)]
			\item the range of $L_n$ is included in a finite-dimensional space,
			\item $\limsup_{n \to \infty}\Lip(L_n)\le 1+\eps$,
			\item for all $x\in r_0B_E$ and all $f\in \Lip_0(E)$, $$\lim_{n\to \infty}\langle L_n(x),f\rangle = f(x).$$
		\end{enumerate}  
	\end{lemm}
	
	\begin{proof} 
		The idea of the proof is to ``truncate'' the evaluation map by convoluting it with Fej\'er kernels on well-chosen cubes. Therefore, our first order of business is to define these kernels. We start with the one-dimensional Dirichlet and Fej\'er kernels. For $n\in \bN\cup\{0\}$ and $t\in \bR$, let 
		\begin{equation*}
			D_n(t) := \sum_{k=-n}^n e^{ikt}
		\end{equation*}
		and
		\begin{equation*}
			F_n := \frac{1}{n+1}\big(D_0+\cdots+D_{n}\big).
		\end{equation*}
		The basic properties of the sequence $(F_n)_n$ of Fej\'er kernels are well known:
		\begin{enumerate}
			\item $F_n$ is nonnegative and even,
			\item $\int_{-\pi}^\pi F_n(t)\,\frac{dt}{2\pi}=1$,
			\item for all $s \in (0,\pi)$, $\lim_{n\to \infty}\int_{s\le |t|\le \pi}F_n(t)\,\frac{dt}{2\pi}=0$.
		\end{enumerate} 
		From now on, we fix $\eps>0$ and identify $E$ with $\bR^N$ where $N$ is the dimension of $E$. For simplicity, we denote by $\norm{\cdot}$ the norm of $E$ and we assume, as we may that $\norm{x}\le \norm{x}_\infty \le \beta\norm{x}$ for some constant $\beta$ where $\norm{\cdot}_\infty$ stands for the standard $\ell_\infty$-norm on $\bR^N$ and whose unit ball is denoted by $B_\infty$. Next, we extend the one-dimensional Fej\'er kernels to $\bR^N$. For that purpose, let us fix $a>1$ to be carefully chosen later and $n\in \bN$. For $t := (t_1,\dots,t_N)\in \bR^N$, consider the following rescaling of $F_n$ 
		\begin{equation*}
			G_n(t_k) := \frac{1}{8a}F_n\big(\frac{\pi t_k}{4a}\big), 
		\end{equation*}
		and let  
		\begin{equation*}
			H_n(t) := \prod_{k=1}^N G_n(t_k).
		\end{equation*}
		The basic properties of the Fej\'er kernels are easily seen to induce similar properties for the sequence of kernels $(H_n)_n$. In particular,
		\begin{equation*}
			\int_{4aB_\infty}H_n(t)\,dt=1,
		\end{equation*}
		and for all $r\in (0,4a)$,
		\begin{equation*}
			\lim_{n\to \infty}\int_{rB_\infty}H_n(t)\,dt=1.
		\end{equation*}
		For $x\in \bR^N$, we now define $\rho_n(x) \in \cF(E)$ as the following Bochner integral:
		\begin{equation*}
			\rho_n(x) := \int_{4aB_\infty} (H_n(x-t)-H_n(t))\delta_E(t)\,dt.
		\end{equation*}
		The key observation is that the maps $t\mapsto (H_n(x-t)-H_n(t))\delta_E(t)$, for $x\in \bR^N$, belong to a common finite-dimensional space of functions. Thus, there exists a finite-dimensional subspace $X_n$ of $\cF(E)$ such that $\rho_n$ maps $\bR^N$ into $X_n$. 
		
		One can estimate the Lipschitz constant of $\rho_n$ considered as a map on $aB_\infty$ as follows. For $x,y \in aB_\infty$ we let $u := \frac12(x-y)$ and $v := \frac12(x+y)$. We shall use the following auxiliary elements of $\cF(E)$:
		\begin{equation*}
			\mu := \int_{u+4aB_\infty}H_n(x-t)\delta_E(t)\,dt -\int_{4aB_\infty}H_n(t)\delta_E(t)\,dt
		\end{equation*}
		and 
		\begin{equation*}
			\nu := \int_{-u+4aB_\infty}H_n(y-t)\delta_E(t)\,dt -\int_{4aB_\infty}H_n(t)\delta_E(t)\,dt.
		\end{equation*}
		An elementary change of variable tells us that  
		\begin{equation*}
			\mu-\nu = \int_{4aB_\infty}H_n(v-w)\big(\delta_E(w+u)-\delta_E(w-u)\big)\,dw,
		\end{equation*}
		and hence 
		\begin{equation*}
			\norm{\mu-\nu}_{\cF(E)} \le 2 \norm{u}\int_{4aB_\infty}H_n(v-w)\,dw = \norm{x-y}.
		\end{equation*} 
		Next, we estimate $\|\rho_n(x)-\mu\|_{\cF(E)}$. From the nonnegativity of $H_n$ and the triangle inequality we have that 
		\begin{equation*}
			\norm{\rho_n(x)-\mu}_{\cF(E)} \le \int_{(4a+\norm{u}_\infty) B_\infty \setminus 4aB_\infty} H_n(x-t)\norm{\delta_E(t)}_{\cF(E)}\,dt.
		\end{equation*}
		Since $\norm{u}_\infty \le a$, $\norm{\delta_E(t)}_{\cF(E)}\le 5a$ on this domain of integration. Therefore,
		\begin{equation}
			\label{eq1:weak-approx-Fejer}
			\norm{\rho_n(x)-\mu}_{\cF(E)} \le 5a\int_{(4a+\norm{u}_\infty) B_\infty \setminus 4aB_\infty}H_n(x-t)\,dt.
		\end{equation}
		Remembering that $x:=(x_1,\dots,x_N)\in aB_{\infty}$, we have $\abs{x_k}\le a$ for any $k\le N$. So, if $4a\le \abs{t_k}\le 4a + \norm{u}_\infty$, then $3a\le \abs{x_k-t_k}\le 6a$. An elementary study of the one-dimensional Fej\'er kernels shows the existence of a constant $c(a,n)$ so that $|G_n(s)|\le c(a,n)$ whenever $3a\le |s|\le 6a$ with the property that $\lim_{n\to \infty}c(a,n)=0$ and thus
		\begin{equation}
			\label{eq2:weak-approx-Fejer}
			\int_{-4a-\norm{u}_\infty}^{4a+\norm{u}_\infty} G_n(x_k-t_k)\,dt_k \le 1+2c(a,n) \norm{u}_\infty \le 1+2c(a,n)M\norm{u}.
		\end{equation}
		Combining \eqref{eq1:weak-approx-Fejer} and \eqref{eq2:weak-approx-Fejer} we get
		\begin{equation*}
			\norm{\rho_n(x)-\mu}_{\cF(E)} \le 5a\big((1+2c(a,n)\beta\norm{u})^N-1\big).
		\end{equation*}
		This yields the existence of a constant $C(a,N,\beta,n)$ so that $\lim_{n\to \infty}C(a,N,\beta,n)=0$ and 
		\begin{equation*}
			\norm{\rho_n(x)-\mu}_{\cF(E)}\le C(a,N,\beta,n)\norm{x-y}.
		\end{equation*}
		Similarly, one can show that $\norm{\rho_n(y)-\nu}_{\cF(E)}\le C(a,N,\beta,n)\norm{x-y}$ and hence for all $x,y \in aB_\infty$ one has
		\begin{equation}
			\label{eq3_Lipconstant}
			\norm{\rho_n(x)-\rho_n(y)}_{\cF(E)}\le (1+2C(a,N,\beta,n)) \norm{x-y}    
		\end{equation}
		The approximating map $L_n$ is obtained by slightly modifying $\rho_n$ as follows. For $x\in \bR^N$, let
		\begin{equation*}
			\varphi(x) := \max\Big(1-\frac{\max\{\log\norm{x}_\infty,0\}}{\log a},0\Big)\in [0,1],
		\end{equation*}
		and define
		\begin{equation*}
			L_n(x) := \varphi(x)\rho_n(x).
		\end{equation*}
		First, observe that condition $(i)$ holds since it is clear that the map $L_n$ still takes values in the finite-dimensional subspace $X_n$ of $\cF(E)$. 
		
		Let us now estimate the Lipschitz constant of $L_n$ on $\bR^N$. To simplify, we denote by $c_n$ the Lipschitz constant of $\rho_n$ considered as a map on $aB_\infty$. Note that \eqref{eq3_Lipconstant} means that $c_n\le (1+2C(a,N,\beta,n))$). Assume first that $\norm{x}_\infty \le \norm{y}_\infty \le a$. Then,
		\begin{align*}
			\norm{L_n(x)-L_n(y)}_{\cF(E)} & \le \abs{\varphi(y)}\,\norm{\rho_n(x)-\rho_n(y)}_{\cF(E)} + \abs{\varphi(y)-\varphi(x)}\,\norm{\rho_n(x)}_{\cF(E)}\\
			& \le c_n \norm{x-y} + \frac{1}{\log a}\log\frac{\norm{y}_\infty}{\norm{x}_\infty}c_n\norm{x}\\
			& \le c_n\left(\norm{x-y} + \frac{1}{\log a} \Big(\frac{\norm{y}_\infty}{\norm{x}_\infty}-1\Big) \norm{x}_\infty\right)\\
			& \le c_n(1 + \beta(\log a)^{-1})\norm{x-y}.
		\end{align*}
		It is easy to see that the same estimate is valid if $\norm{x}_\infty \le \norm{y}_\infty$ and $\norm{y}_\infty \ge a$. So, we deduce that 
		\begin{equation*}
			\Lip(L_n)\le (1+2C(a,N,\beta,n))(1+\beta(\log a)^{-1}).
		\end{equation*}
		If $a$ was initially chosen so that $\beta(\log a)^{-1} < {\eps}$, condition $(ii)$ follows, i.e. 
		\begin{equation*}
			\limsup_{n \to \infty}\Lip(L_n)<1+\eps.
		\end{equation*}
		Finally, if $r_0=\beta^{-1}$, then for $x\in r_0B_E$ we have $x\in B_\infty$ and thus $L_n(x)=\rho_n(x)$. It now follows from the standard results on Fourier series and Fej\'er kernels that for all $x\in r_0B_E$ and all $f\in \Lip_0(E)$,
		\begin{equation*}
			\lim_{n\to \infty}\langle L_n(x),f\rangle=\lim_{n\to \infty}\big((H_n \ast f)(x)-(H_n \ast f)(0)\big)=f(x).
		\end{equation*}
	\end{proof}
	
	With the help of Lemma \ref{lem:weak-approx-Fejer} we can now prove that Lipschitz-free spaces over finite-dimensional Banach spaces have the bounded approximation property. In fact, we have the following more precise result.
	
	\begin{theo}
		\label{thm:fd-LF-MAP}
		%Let $E$ be a finite-dimensional Banach space. Then, $\cF(E)$ has the metric approximation property (MAP).
		The Lipschitz-free space over a finite-dimensional Banach space has the metric approximation property.
	\end{theo}
	
	\begin{proof}
		Let $r_0$ and $L_n$ be as in Lemma \ref{lem:weak-approx-Fejer} and let $x_1,\dots, x_k$ be fixed vectors in $r_0B_E$. Then, it follows from Mazur's theorem that there exists a sequence $(f_n)_n$ of convex combinations of the $L_n$ such that $\lim_n\norm{f_n(x_i) - \delta_E(x_i)}=0$ for all $1\le i\le k$. Therefore, using the compactness of $r_0B_E$, we can build a sequence of maps from $E$ to $\cF(E)$, still denoted by $(f_n)_n$, such that $f_n(0)=0$, the range of $f_n$ is included in a finite-dimensional space, $\Lip(f_n)\le 1+\frac{1}{n}$ and $\norm{f_n(x)-\delta_E(x)}\le \frac{1}{n^2}$, for all $x\in r_0B_E$. Using adequate rescaling, we can modify the map $f_n$ in a way that the new map still approximates the evaluation map but on a much larger ball, namely the ball of radius $nr_0$. To do so, consider $S_n \colon \cF(E)\to \cF(E)$, the unique bounded linear map so that for all $x\in E$
		\begin{equation*}
			S_n \circ \delta_E(x) = \delta_E(nx).
		\end{equation*} 
		This map is guaranteed to exist by Corollary \ref{cor:GK-lifting} and $\norm{S_n}=n$. Now, considering the map $g_n \colon x\in E \mapsto (S_n \circ f_n)(\frac{x}{n})$, we have that $\Lip(g_n)\le 1+\frac{1}{n}$. Moreover, for all $x\in nr_0B_E$,
		\begin{equation*}
			\norm{g_n(x) - \delta_E(x)} = \Big\| S_n\Big(g_n\big(\frac{x}{n}\big) - \delta_E\big(\frac{x}{n}\big)\Big)\Big\|\le \frac{1}{n}.
		\end{equation*}
		Invoking Theorem \ref{thm:GK-lifting-free}, there is a unique bounded linear map $\tilde{g}_n \colon \cF(E)\to \cF(E)$ such that $\tilde{g}_n\circ \delta_E = g_n$ and $\|\tilde{g}_n\|\le 1+\frac{1}{n}$. The range of $\tilde{g}_n$ is included in the closed linear span of the range of $g_n$ and is therefore finite-dimensional. Finally, we have that for all $x\in nr_0B_E$, $\|\tilde{g}_n(\delta_E(x))-\delta_E(x)\|\le \frac
		{1}{n}$. For a fixed compact subset $K$ of $\cF(E)$ and a fixed $\eps>0$, this clearly implies (modulo an approximation argument) that we can find $n$ large enough, so that $\|\tilde{g}_n(\mu) - \mu\|\le \eps$ for all $\mu \in K$, i.e. we have proved that $\cF(E)$ has the MAP. 
		
	\end{proof}
	
	The proof of Theorem \ref{thm:LBAP->LF-BAP} follows from Theorem \ref{thm:fd-LF-MAP} and a density argument.
	
	\begin{proof}[Proof of Theorem \ref{thm:LBAP->LF-BAP}]
		Assume that a Banach space $X$ has the $\lambda$-LBAP. By density, it is enough to show that for any $\eps>0$ and any $x_1,\ldots,x_n \in X$, there exits $T\colon \cF(X)\to \cF(X)$ a finite-rank linear operator such that $\norm{T}\le \lambda$ and $\norm{T\delta_X(x_i)-\delta_X(x_i)}\le \eps$, for $1\le i\le n$. By assumption, there exists a Lipschitz map $f\colon X\to X$ such that $E$, the linear span of $f(X)$ is finite-dimensional, $\Lip(f)\le \lambda$, $f(0)=0$ and $\norm{f(x_i)-x_i} = \norm{\delta_X(x_i)-\delta_X(f(x_i))}\le \frac{\eps}{2}$, for $1\le i \le n$. Let $\widehat{f}\colon \cF(X)\to \cF(X)$ be the unique linear operator so that $\widehat{f}\circ \delta_X = \delta_X\circ f$. We have that $\|\widehat{f}\|\le \lambda$ and the range of $\widehat{f}$ is included in $\cF(E)\subset \cF(X)$. Thus, we can apply Theorem \ref{thm:fd-LF-MAP} to find a finite-rank operator $R\colon\cF(E)\to \cF(E)$ such that $\norm{R}\le 1$ and $\|R\widehat{f}(\delta_X(x_i))-\widehat{f}(\delta_X(x_i))\|\le \frac{\eps}{2}$, for $1\le i\le n$. It is then enough to consider the finite-rank operator $T := R\widehat{f}$ to conclude.
	\end{proof}
	
	The next corollary collects the results and observations from this section.
	
	\begin{coro}
		\label{cor:LBAP}
		Let $X$ be an arbitrary Banach space and $\lambda\ge 1$. Then, the following conditions are equivalent:
		\begin{enumerate}[(i)]
			\item $X$ has the $\lambda$-bounded approximation property.
			\item $\cF(X)$ has the $\lambda$-bounded approximation property.
			\item $X$ has the $\lambda$-Lipschitz bounded approximation property.
		\end{enumerate}
	\end{coro}
	
	As another immediate corollary, we have the following important rigidity result.
	
	\begin{coro}
		\label{cor:BAP-Lip-invariant}
		The bounded approximation property is a Lipschitz invariant, i.e., if $X$ is a Banach space with the bounded approximation property and $Y$ is a Banach space Lipschitz equivalent to $X$, then $Y$ has the bounded approximation property.
	\end{coro}
	
	\begin{proof}
		%Let $X$ be a Banach space with the BAP. If a Banach space $Y$ is Lipschitz equivalent to $X$, then $\cF(X)$ and $\cF(Y)$ are linearly isomorphic. It follows from Corollary \ref{cor:LBAP} that $\cF(X)$ has the BAP, so $\cF(Y)$ also has the BAP and another application of Corollary \ref{cor:LBAP} implies that $Y$ has the BAP.
		Simply remember that two Banach spaces that are Lipschitz equivalent have linearly isomorphic Lipschitz-free spaces. The conclusion follows from Corollary \ref{cor:LBAP} and the fact that the BAP is an isomorphic invariant.
	\end{proof}

	\section{Structure of Lipschitz-free spaces}
	\label{sec:structure-Lip-free}
	
	In this section, we delve deeper into the structural properties of Lipschitz-free spaces. The main focus is on establishing criteria that guarantee that a Lipschitz-free space has classical Banach space properties such as the approximation property, the Radon-Nikod\'ym property and most importantly the Schur property. The two main tools to establish these structural results are explained in Section \ref{sec:Kalton-toolbox}. The first one is a general, yet powerful, result regarding the decomposition of Lipschitz-free spaces. The second one is an approximation result about weakly null sequences in Lipschitz-free spaces. In Section \ref{sec:free-spaces-Schur}, we give three criteria, of rather different natures, implying that large classes of Lipschitz-free spaces have the Schur property. These structural results almost immediately lead to the existence of Banach spaces that are universal for separable metric spaces and coarse-Lipschitz or uniform embeddings, but do not contain a bi-Lipschitz copy of $\co$.
	
	\subsection{Kalton's toolbox for Lipschitz-free spaces}
	\label{sec:Kalton-toolbox}
	
	We first detail a fundamental decomposition scheme for Lipschitz-free spaces due to N. Kalton \cite{Kalton2004}. So, consider an arbitrary metric space $(M,d)$, pointed at $0\in M$ and for $k\in \bZ$, let 
	\begin{equation*}
		M_k :=\{x\in M \colon \ d(x,0)\le 2^k\}.
	\end{equation*}
	The first basic idea is, roughly speaking, to decompose $\cF(M)$ into a direct sum of Lipschitz-free spaces over disjoint annuli of the form $M_{s_k}\setminus M_{r_k}$, for suitably chosen sequences $(r_k)_k$ and $(s_k)_k$ of integers. The first result is a simple lower $\ell_1$-estimate for sums of measures belonging to skipped blockings of this type.
	
	\begin{lemm}
		\label{lem:skipped-dyadic}
		Let $(M,d)$ be a pointed metric space.
		Suppose that $r_1<s_1<\cdots<r_n<s_n$ belong to $\bZ$ and $\gamma_k \in \cF(M_{s_k}\setminus M_{r_k})$ for $1\le k \le n$, then
		\begin{equation*}
			\Big\|\sum_{k=1}^n \gamma_k\Big\|_{\cF(M)}\ge \frac{2^\theta-1}{2^\theta+1}\sum_{k=1}^n \|\gamma_k\|_{\cF(M)},
		\end{equation*}
		where $\theta :=\min_{1\le k< n}(r_{k+1}-s_k)$.
	\end{lemm}
	
	\begin{proof}
		For $k\in \bZ$, pick $f_k\in \Lip_0(M)$ so that $\|f_k\|_L=1$ and $\langle \gamma_k,f_k\rangle=\|\gamma_k\|_{\cF(M)}$. Then, consider the map $g \colon \{0\}\cup \bigcup_{k=1}^n (M_{s_k}\setminus M_{r_k}) \to \bR$ defined by $g(0) := 0$ and $g(x) := f_k(x)$ if $x\in M_{s_k}\setminus M_{r_k}$. The map $g$ is Lipschitz. Indeed, if $x\in M_{s_k}\setminus M_{r_k}$ and $y\in M_{s_j}\setminus M_{r_j}$, for $j<k$ in $\bZ$, then it is easy to verify that $2^\theta d(y,0) \le d(x,0) \le \frac{d(x,y)}{1-2^{-\theta}}$ and hence
		\begin{equation*}
			\abs{g(x)-g(y)}\le d(x,0)+d(y,0)\le \frac{1+2^{-\theta}}{1-2^{-\theta}}d(x,y).
		\end{equation*}
		Since $\frac{1+2^{-\theta}}{1-2^{-\theta}} = \frac{2^\theta+1}{2^\theta-1} \ge 1$, it follows that $g$ is $\frac{2^\theta+1}{2^\theta-1}$-Lipschitz. Therefore, we can extend $g$ to $f\in \Lip_0(M)$ with $\norm{f}_L\le \frac{2^\theta+1}{2^\theta-1}$. Then,
		\begin{equation*}
			\Big\|\sum_{k=1}^n \gamma_k\Big\|_{\cF(M)}\ge \frac{2^\theta-1}{2^\theta+1}\sum_{k=1}^n\langle \gamma_k,f\rangle = \frac{2^\theta-1}{2^\theta+1}\sum_{k=1}^n\langle \gamma_k,f_k\rangle=\frac{2^\theta-1}{2^\theta+1}\sum_{k=1}^n \norm{\gamma_k}_{\cF(M)}.
		\end{equation*}
		
	\end{proof}
	
	Next, we describe a ``wavelet-like'' decomposition of $\cF(M)$. For that purpose, we introduce, for $n\in \bZ$, the linear operator $T_n\colon \cF(M)\to \cF(M_{n+1}\setminus M_{n-1})$ such that $T_n(\delta_M(x))=\psi_n(x)\delta_M(x)$, where $\psi_n\colon M \to [0,1]$ is defined as follows:
	$$
	\psi_n(x) := \begin{cases}
		0 \textrm{ if } x\in M_{n-1} \textrm{ or } x\notin M_{n+1}, \\
		\log_2d(x,0)-(n-1) \textrm{ if } x\in M_n\setminus M_{n-1},\\
		n+1-\log_2d(x,0) \textrm{ if } x\in M_{n+1}\setminus M_{n}.
	\end{cases} 
	$$ 
	The operator $T_n$ can be seen as a multiplier on the ``basis'' $(\delta_M(x))_{x\in M\setminus\{0\}}$ with symbol $\psi_n$, although $(\delta_M(x))_{x\in M\setminus\{0\}}$ is not a topological basis. Let us first detail how $T_n$ extends to $\cF(M)$. Observing that $1+\ln(t)\le t$ whenever $t\ge 1$, it follows that for all $x,y\neq 0$ with $d(y,0)\ge d(x,0)$, one has
	\begin{equation*}
		\abs{\psi_n(x) - \psi_n(y)}\le \log_2\frac{d(y,0)}{d(x,0)} \le \frac{1}{\ln 2}\frac{d(y,0)-d(x,0)}{d(x,0)}\le 2\frac{d(x,y)}{d(x,0)}.
	\end{equation*}
	Hence, since $0\le \psi_n\le 1$, we get
	\begin{align*}
		\norm{ T_n(\delta_M(x)) - T_n(\delta_M(y)) }_{\cF(M)} & \le \norm{ \psi_n(y)(\delta_M(x) - \delta_M(y))}_{\cF(M)}\\ 
		&+ \norm{(\psi_n(x) - \psi_n(y))\delta_M(x)}_{\cF(M)}\\
		& \le  d(x,y)+2\frac{d(x,y)}{d(x,0)}d(x,0)\\
		& \le 3 d(x,y).
	\end{align*}
	Then, it readily follows from Theorem \ref{thm:factorization-free-Lip} that the $3$-Lipschitz map $T_n\circ \delta_M\colon M \to \cF(M)$ extends to a  bounded linear operator, still denoted by $T_n$, from $\cF(M)$ to $\cF (M)$ with $\|T_n\|\le 3$. The sequence $(T_n)_{n\in \bZ}$ of operators can be used to decompose every element of a Lipschitz-free space as a sum of elements in the Lipschitz-free spaces over its dyadic annuli.
	
	\begin{prop}
		\label{pro:decomposition-free}
		Let $M$ be a metric space and $(T_n)_{n\in \bZ}$ defined as above. Then, for every $\gamma \in \cF(M)$, 
		\begin{enumerate}[(i)]
			\item $\sum_{n\in \bZ}\norm{T_n \gamma}_{\cF(M)}\le 108\norm{\gamma}_{\cF(M)}$,
			\item[] and
			\item $\gamma=\sum_{n\in \bZ}T_n\gamma.$
		\end{enumerate}
	\end{prop}
	
	\begin{proof}
		For $(ii)$, note that for all $x\in M$, $\sum_{n\in \bZ}T_n(\delta_M(x))=\delta_M(x)$ (this is a finite sum with only two terms) and it follows from $(i)$ and a density argument that $\sum_{n\in \bZ}T_n\gamma=\gamma$ for all $\gamma \in \cF(M)$.  Thus, it remains to prove $(i)$. 
		For a finitely nonzero family $(a_n)_{n\in \bZ}$ in $[-1,1]$, we consider the operator $S :=\sum_{n\in \bZ} a_nT_n$ on $\cF(M)$. For fixed $x,y \in M$, the quantity $S(\delta_M(x))-S(\delta_M(y))$ involves at most 4 values of $n$. It follows that $\norm{S(\delta_M(x))-S(\delta_M(y))}_{\cF(M)}\le 12d(x,y)$ and, as before, that $\norm{S}\le 12$. Now, for any $\gamma \in \linspan(\delta_M(M))$, $\sum_{n\in \bZ}T_n\gamma$ is a finite sum and hence trivially summable. It then follows from our estimate on $\norm{S}$, an approximation argument and Cauchy's criterion for summability that, for any $\gamma \in \cF(M)$, $\sum_{n\in \bZ}T_n\gamma$ is summable (or equivalently unconditionally convergent). 
		Now, for $r\in\{0,1,2\}$ and $\gamma \in \cF(M)$, observe that $T_{3n+r}(\gamma)\in \cF(M_{3n+r+1}\setminus M_{3n+r-1})$. Therefore, by Lemma~\ref{lem:skipped-dyadic} (with $\theta=1$) and our estimate on $\norm{S}$, we deduce that,
		\begin{equation*}
			\sum_{n\in \bZ}\norm{ T_{3n+r}\gamma}_{\cF(M)}\le 3\big\|\sum_{n\in \bZ}T_{3n+r}\gamma\big\|_{\cF(M)} \le 36\norm{\gamma}_{\cF(M)}.
		\end{equation*}
		The final estimate easily follows. 
	\end{proof}
	
	If we consider the map $T\colon \cF(M)\to \big(\sum_{k\in \bZ} \cF(M_k)\big)_{\ell_1}$ given by $T(\gamma) := (T_{k-1}(\gamma))_{k\in \bZ}$, then it follows immediately from Proposition \ref{pro:decomposition-free} that $T$ is an isomorphic embedding with distortion at most $108$. The distortion can certainly be slightly improved if instead of partitioning $\bZ$ with multiples of $3$ we use a partition with multiples of $m\ge 3$. Indeed, a simple verification shows that the modified embedding, still denoted by $T$, will then satisfy 
	\begin{equation*}
		\norm{\gamma}\le \norm{T\gamma}\le 12m\frac{2^{m-2}+1}{2^{m-2}-1}\norm{\gamma}.
	\end{equation*}
	Note that here we would use $m$ times the triangle inequality and a gap $\theta=m-2$. The best distortion we can get like this is about $77$ with $m=5$. In order to lower the distortion further, we need to be craftier. Using a partition of $\bZ$, we can not do much better than what we explained above, and we need to resort to some covering of $\bZ$ with $m$ sequences of subsets $(A_{m,n})_{n\in\bZ}$ that overlap in a certain fashion. By carefully picking this covering, we can get an almost isometric embedding. 
	
	\begin{theo}
		\label{thm:almost-isometric-embed-free}
		For any metric space $M$ and $\eps >0$, the space $\cF(M)$ is $(1+\eps)$-isomorphic to a subspace of $\big(\sum_{k\in \bZ} \cF(M_k)\big)_{\ell_1}$.
	\end{theo}
	
	\begin{proof}
		The salient point of the proof is to identify a specific covering of $\bZ$ with $m$ sequences of subsets consisting of finite disjoint subsets with large gaps between them. Given $m\ge 3$ and $s\ge 1$ to be chosen large enough later, we consider the $m$ sequences of subsets $(A_{0,n})_{n\in\bZ},\dots, (A_{m-1,n})_{n\in\bZ}$ where for $0\le r\le m-1$ and $n\in \bZ$,
		\begin{equation*}
			A_{r,n} := \{(mn+r)s+1,\dots,(m(n+1)+r-1)s\}. 
		\end{equation*} 
		To better grasp the properties of these collections of sets, note that in the particular case when $m=3$ and $s=2$, we cover $\bZ$ with $(A_{0,n})_{n\in\bZ} = (\{6n+1,\dots,6n+5\})_{n\in\bZ}$, $(A_{1,n})_{n\in\bZ} = (\{6n+3,\dots,6n+6\})_{n\in\bZ}$ and 
		and $(A_{2,n})_{n\in\bZ} = (\{6n+5,\dots,6n+8\})_{n\in\bZ}$. These collections have two crucial properties. Every integer is covered by exactly $2$ of these sequences of sets and the complements of these sequences are pairwise disjoint. These observations hold true in general, as it can be easily verified. We gather the elementary properties of such coverings in the following claim.
		
		\begin{claim}
			\label{clai:almost-isometric-embed-free}
			For all $m\ge 3, r\in \{0,\dots, m-1\}, s\ge 1$ and $n\in\bZ$, we have
			\begin{enumerate}[(i)]
				%\item $\abs{A_r,n} = (m-1)s$,
				\item $\min A_{r,n+1}-\max A_{r,n}= s+1$.
				\item[]For all $0\le r\le m-1$, if we let $A_r:=\cup_{n\in\bZ}A_{r,n}$, then
				\item $\bZ=\cup_{r=0}^{m-1} A_r$ and every $k\in \bZ$ belongs to exactly $m-1$ of the sets $A_0,\dots, A_{m-1}$,
				\item the collection $\{\bZ\setminus A_r\colon 0\le r\le m-1\}$ is a collection of pairwise disjoint sets.
			\end{enumerate}
		\end{claim}
		Now, we slightly rework the computation from Proposition \ref{pro:decomposition-free} by replacing $T_{3n+r}$ by $\sum_{i \in A_{r,n}} T_i$ where $0\le r\le m-1$ and $n\in \bZ$. Observe first that $\sum_{i \in A_{r,n}} T_i$ maps $\cF(M)$ into $\cF(M_{\max A_{r,n}+1}\setminus M_{\min A_{r,n}-1})$. 
		Therefore, it follows from assertion $(i)$ in Claim \ref{clai:almost-isometric-embed-free}, Lemma \ref{lem:skipped-dyadic} and Proposition \ref{pro:decomposition-free} that for all $\gamma \in \cF(M)$ and $r\in\{0,\dots,m-1\}$, 
		\begin{equation}
			\label{eq1:almost-ismoteric-emb-free}
			\sum_{n\in \bZ} \Big\|\sum_{i \in A_{r,n}} T_i\gamma\Big\| \le \frac{2^{s-1}+1}{2^{s-1}-1}\Big\|\sum_{i \in A_{r}} T_i\gamma\Big\| \le \frac{2^{s-1}+1}{2^{s-1}-1}\big(\norm{\gamma} + \sum_{i\in \bZ\setminus A_r}\norm{T_i \gamma}\Big).
		\end{equation}
		Summing over $r\in\{0,\dots,m-1\}$ in \eqref{eq1:almost-ismoteric-emb-free} and because the $\bZ\setminus A_r$ are pairwise disjoint, it follows from Proposition \ref{pro:decomposition-free} that
		\begin{equation}
			\label{eq2:almost-ismoteric-emb-free}
			\sum_{r=0}^{m-1}\sum_{n\in \bZ}\Big\|\sum_{i \in A_{r,n}}  T_i\gamma\Big\| \le \frac{2^{s-1}+1}{2^{s-1}-1}(m+108)\norm{\gamma}.
		\end{equation}
		On the other hand, since each integer in $\bZ$ belongs to exactly $(m-1)$ of the sets $A_0,\dots, A_{m-1}$ that cover $\bZ$, we have that for all $\gamma \in \cF(M)$,
		\begin{equation}
			\label{eq3:almost-isomteric-embed-free}
			\sum_{r=0}^{m-1}\sum_{n\in \bZ}\Big\|\sum_{i \in A_{r,n}} T_i\gamma\Big\|\ge \Big\|\sum_{r=0}^{m-1}\sum_{i \in A_{r}} T_i\gamma\Big\| = \Big\|(m-1) \sum_{k\in\bZ}T_k\gamma\Big\| = (m-1)\norm{\gamma}.
		\end{equation}
		Now, for all $0\le r\le m-1$, define $S_r\colon \cF(M) \to \big(\sum_{k\in \bZ} \cF(M_k)\big)_{\ell_1}$ by setting for all $k\in \bZ$,
		$$
		(S_r \gamma)(k) := \begin{cases}
			\sum_{i \in A_{r,n}} T_i\gamma \textrm{ if } k=\min A_{r,n} \textrm{ for some } n\in \bZ,\\
			0 \textrm{ otherwise},
		\end{cases}
		$$
		and $T\colon \cF(M)\to \ell_1^m\big(\sum_{k\in \bZ} \cF(M_k)\big)_{\ell_1}$ by letting $T\gamma := \frac{1}{m-1}(S_r\gamma)_{r=0}^{m-1}$.
		Then, it follows from \eqref{eq2:almost-ismoteric-emb-free} and \eqref{eq3:almost-isomteric-embed-free} that 
		\begin{equation*}
			\norm{\gamma}\le \norm{T\gamma}\le \frac{(2^{s-1}+1)(m+108)}{(2^{s-1}-1)(m-1)}\norm{\gamma}.
		\end{equation*}
		To conclude, it remains to observe that given $\vep>0$ we can certainly choose $m\ge 3$  and $s\ge 1$ large enough such that $\frac{(2^{s-1}+1)(m+108)}{(2^{s-1}-1)(m-1)}\le 1+\vep$ and also that since $\cF(M_k)$ isometrically embeds into any $\cF(M_l)$ whenever $l\ge k$, then $\ell_1^m\big(\sum_{k\in \bZ} \cF(M_k)\big)_{\ell_1}$ isometrically embeds into $\big(\sum_{k\in \bZ} \cF(M_k)\big)_{\ell_1}$.
	\end{proof}
	
	Intuitively, it should follow from the embedding described in Theorem \ref{thm:almost-isometric-embed-free} that the structure of a Lipschitz-free space over a metric space whose dyadic balls are relatively simple objects should be easier to describe.
	
	\medskip 
	We now turn to the general properties of weakly null sequences in free spaces. For a subset $S$ of $M$, recall that $\cF(S)$ can be seen as a subspace of $\cF(M)$. For $\gamma \in \cF(M)$, the norm of the equivalence class of $\gamma$ in $\cF(M)/\cF(S)$ is $$\norm{\gamma}_{\cF(M)/\cF(S)}=d(\gamma,\cF(S))=\inf\{\norm{\gamma-\mu} \colon \ \mu \in \cF(S)\}.$$
	Observing that $\cF(S)^\perp = \{f\in \Lip_0(M)\colon\ f=0\text{ on }S\}$, we deduce from a classical duality argument that 
	\begin{equation}
		\label{eq:quotient-norm-free}
		\norm{\gamma}_{\cF(M)/\cF(S)} = \sup\{\langle f,\gamma\rangle \colon f\in B_{\Lip_0(M)}\ \text{and}\ f=0\text{ on } S\}.
	\end{equation}
	The definition of $\cF(M)$ implies that for all $\gamma \in \cF(M)$, 
	\begin{equation*}
		\inf\{\norm{\gamma}_{\cF(M)/\cF(S)}\colon \abs{S}<\infty\}=0.
	\end{equation*} 
	It then follows easily that for any $K$ relatively compact subset in $\cF(M)$, we have
	\begin{equation*}
		\inf_{\abs{S}<\infty}\sup_{\gamma \in K}\norm{\gamma}_{\cF(M)/\cF(S)}=0.
	\end{equation*}
	Our first lemma establishes a similar behavior when $K$ is replaced by a weakly null sequence. The principle behind this lemma proved to be extremely useful beyond the results from \cite{Kalton2004} and can be called ``Kalton's gliding hump principle for free spaces''.  We first need to introduce, for $\delta>0$ and $S\subset M$ the $\delta$-neighborhood of $S$:
	\begin{equation*}
		[S]_\delta := \{x\in M\colon d(x,S)\le \delta\}.
	\end{equation*}
	
	\begin{lemm}
		\label{lem:weakly-null-free} 
		Suppose that $M$ is a bounded metric space and that $(\gamma_n)_{n=1}^\infty$ is a weakly null sequence in $\cF(M)$. Then, for all $\delta >0$, 
		\begin{equation*}
			\inf_{\abs{S}<\infty}\sup_{n\in \bN} \norm{\gamma_n}_{\cF(M)/\cF([S]_\delta)}=0.
		\end{equation*}
	\end{lemm}
	
	\begin{proof}
		A simple approximation argument shows that we may assume that each $\gamma_n$ is of finite support. Then, we assume, seeking a contradiction,  that there exist $\delta,\eps>0$ so that
		\begin{equation*}
			\inf_{\abs{S}<\infty}\sup_{n\in \bN} \norm{\gamma_n}_{\cF(M)/\cF([S]_\delta)}>\vep.
		\end{equation*}
		Our goal is to show that the sequence $(\gamma_n)_{n\in \bN}$ cannot be weakly convergent.
		Since $M$ is bounded, let $R>0$ be such that $d(x,0)\le R$ for all $x\in M$ and assume, as we may, that $\delta<R$. We can then construct a subsequence $(\mu_n)_{n=1}^\infty$ of $(\gamma_n)_{n=1}^\infty$ and an increasing sequence $(S_n)_{n=1}^\infty$ of finite subsets of $M$ with $S_0 := \{0\}$ such that $\norm{\mu_n}_{\cF(M)/\cF([S_{n-1}]_\delta)}>\vep$, for all $n\in \bN$ and $\supp(\mu_n)\subset S_n$. It follows from \eqref{eq:quotient-norm-free} that there exists $f_n\in \Lip_0(M)$ such that $\norm{f_n}_L\le 1$, $f_n=0$ on $[S_{n-1}]_\delta$ and $\langle \mu_n,f_n \rangle >\eps$. Then, $\langle \mu_n,f_n^+ \rangle -\langle \mu_n,f_n^- \rangle>\eps$. So, letting $g_n := f_n^+$ or $f_n^-$, we have  that $\norm{g_n}_L\le 1$, $g_n\ge 0$, $g_n =0 $ on $[S_{n-1}]_\delta$ and $\abs{\langle \mu_n,g_n \rangle}>\frac{\eps}{2}$. Next, we extend (by sup-convolution) the restriction of $g_n$ on $\supp(\mu_n)$ to all of $M$. So, define, for $x\in M$,
		\begin{equation*}
			\varphi_n(x) := \sup_{y \in \supp(\mu_n)}(g_n(y)-R\delta^{-1}d(x,y)).
			%\text{ and } h_n(x) :=\max\{0,\varphi_n(x)\}.
		\end{equation*}
		Since $\norm{g_n}_L \le  1< R\delta^{-1} $, it follows from Remark \ref{rem:inf-conv} after Proposition \ref{pro:Inf-Conv} that $\varphi_n=g_n$ on $\supp(\mu_n)$ and $\norm{\varphi_n}_L\le R\delta^{-1}$. Since $g_n\ge 0$, $\varphi_n^+=g_n$ on $\supp(\mu_n)$ and $\norm{\varphi_n^+}_L\le R\delta^{-1}$. For all $y\in M$, $g_n(y)\le d(y,0)\le R$, so $\varphi_n^+(x)=0$ whenever $d(x,\supp(\mu_n))\ge \delta$. On the other hand, using again that $\delta \le R$, we get that for all $y\in \supp(\mu_n)$, 
		$$g_n(y)-R\delta^{-1}d(x,y)\le g_n(x)+(1-R\delta^{-1})d(x,y) \le g_n(x).$$ Thus $0\le \varphi_n^+\le g_n$ on $M$ and in particular $\varphi_n^+=0$ on $[S_{n-1}]_\delta$. 
		It follows that the sets $(\{x\in M\colon \varphi_n^+(x)>0\})_{n=1}^\infty$ are pairwise disjoint. Therefore, $h:=\sup_{k\in \bN} \varphi_k^+ =\sum_{k=1}^\infty \varphi_n^+ \in \Lip_0(M)$ and $\norm{h}_L\le R\delta^{-1}$. By construction, we have that for all $k>n$, $\langle \mu_n, \varphi_k^+\rangle =0$. Since $(\mu_n)_n$ is weakly null, we may also assume, after passing to a further subsequence, that for all $k<n$, $\abs{\langle \mu_n, \varphi_k^+\rangle} \le \frac{\eps}{4\cdot2^k}$. It follows that 
		\begin{equation*}
			\abs{\langle \mu_n,h \rangle} \ge \abs{\langle \mu_n, \varphi_n^+ \rangle} - \frac{\eps}{4}.
		\end{equation*}
		However, for all $n\in \bN$, $\varphi_n^+$ and $g_n$ coincide on the support of $\mu_n$, so 
		$\abs{\langle \mu_n,\varphi_n^+ \rangle} = \abs{\langle \mu_n,g_n \rangle} >\frac{\eps}{2}$ and $\abs{\langle \mu_n,h \rangle} \ge \frac{\eps}{4}$. This contradicts the fact that the sequence $(\gamma_n)_n$ is weakly null.
	\end{proof}
	
	Note that Lemma \ref{lem:weakly-null-free} says that weakly null sequences in Lipschitz-free spaces are essentially supported on a small set consisting of finitely many arbitrarily small balls.
	
	\subsection{Lipschitz-free spaces with the Schur property}\label{sec:free-spaces-Schur}
	
	Recall that a Banach space $X$ is said to have the \emph{Schur property} (or to be a Schur space) if every weakly converging sequence in $X$ is norm converging. The Banach space $\ell_1(\Gamma)$, where $\Gamma$ is any infinite set, is the simplest example of an infinite-dimensional Schur space. In fact, it follows from Rosenthal's $\ell_1$-theorem that every infinite-dimensional Banach space with the Schur property must necessarily contain an isomorphic copy of $\ell_1$. The converse is not true since it is well known that $L_1$ is not a Schur space. The Schur property is stable under isomorphic embeddings and taking $\ell_1$-sums, but not under passing to quotients (see Exercise \ref{ex:RNP-Schur-l1-sums}). In this section, we give three sufficient conditions for a Lipschitz-free space to have the Schur property. Thanks to a recent work by Aliaga, Gartland, Petitjean and Proch\'azka \cite{AGPP2021}, such Lipschitz-free spaces are now perfectly understood, and we shall comment on them in the Notes section. For now, we will just observe that $\cF(\bZ)\equiv \ell_1$ is a Schur space but $\cF(\bR)\equiv L_1$ is not and these criteria have all in common making some assumption impacting the small scale structure of the metric space, which in particular prevents the Lipschitz-free space to contain an isomorphic copy of $L_1$ (and hence the metric space cannot contain a bi-Lipschitz copy of an interval of positive measure). 
	
	The first criterion is simple, yet has interesting consequences. It makes use of the embedding theorem from Section \ref{sec:Kalton-toolbox} and is based on the observation that a Lipschitz-free space over a bounded metric space that is uniformly discrete, i.e., for which $\inf_{x\neq y\in M}d(x,y)>0$, must be isomorphic to an $\ell_1$-space.
	
	\begin{lemm}
		\label{lem:Lip-free-bounded-unif-discrete}
		If $(M,d)$ is a bounded and uniformly discrete metric space (pointed at $0\in M$), then $\cF(M)$ is isomorphic to $\ell_1(M\setminus\{0\})$.
	\end{lemm}
	\begin{proof}
		If $(M,d)$ has at least two distinct points, let $\theta:=\inf_{x\neq y\in M}d(x,y)$ and $R:=\sup_{x\in M}d(0,x)>0$ and observe that boundedness and uniform discreteness of $M$ implies that $\Lip_0(M)$ and the set of bounded real-valued functions vanishing at $0$ coincide. Indeed, for all $f\in \Lip_0(M)$, we have  
		\begin{equation*}
			R^{-1}\norm{f}_\infty \le \norm{f}_L\le 2\theta^{-1}\norm{f}_\infty.
		\end{equation*}
		It is also not difficult to check that the identity map from $\Lip_0(M)$ to $\ell_\infty(M \setminus \{0\})$ is a weak$^*$ to weak$^*$ continuous isomorphism and thus it is the dual operator of an isomorphism between $\ell_1(M\setminus\{0\})$ and $\cF(M)$.
	\end{proof}
	
	A consequence of Lemma \ref{lem:Lip-free-bounded-unif-discrete} and the embedding theorem from the previous section is that a Lipschitz-free space over a countable and uniformly discrete metric space has every isomorphic property of $\ell_1$ that passes to $\ell_1$-sums and isomorphic embeddings. Moreover, it also has the approximation property.
	
	\begin{theo}
		\label{thm:Lip-free-over-countable-uniformly-discrete}
		%Let $M$ be a countable uniformly discrete metric space. Then, $\cF(M)$ is a Schur space with the Radon-Nikod\'{y}m property and the approximation property.  
		The Lipschitz-free space over a uniformly discrete metric space has the Schur property, the Radon-Nikod\'{y}m property and the approximation property. 
	\end{theo}
	
	\begin{proof}
		Let $(M,d)$ be a uniformly discrete metric space and let $M_k$ be  $B(0,2^k)$, the dyadic ball of radius $2^k$ centered at $0\in M$. It follows from Lemma \ref{lem:Lip-free-bounded-unif-discrete} that $\cF(M_k)$ is isomorphic to an $\ell_1$-space and therefore has the Radon-Nikod\'{y}m and the Schur properties. These properties being stable under taking $\ell_1$-sums and under isomorphic embeddings, we deduce from Theorem \ref{thm:almost-isometric-embed-free} (in fact Proposition \ref{pro:decomposition-free} would suffice) that $\cF(M)$ has the Radon-Nikod\'{y}m and the Schur properties.
		
		For the approximation property, we argue as follows. It follows from  Proposition \ref{pro:decomposition-free} that, given a compact subset $K$ of $\cF(M)$ and $\eps>0$, there exists $N\in \bN$ such that $\norm{S_N\gamma-\gamma}\le \frac{\eps}{2}$ for all $\gamma\in K$, where $S_N := \sum_{n=-N}^N T_n$ is a bounded operator with values in $\cF(M_{N+1})$. Then, since $\cF(M_{N+1})$ is isomorphic to an $\ell_1$-space and thus has the approximation property, there exists a finite-rank operator $R\colon \cF(M_{N+1}) \to \cF(M_{N+1})$ such that $\norm{R\mu-\mu}\le \frac{\eps}{2}$ for all $\mu$ in the compact set $S_N(K)$. Finally, $\norm{RS_N\gamma-\gamma}\le \eps$ for all $\gamma \in K$. This proves that $\cF(M)$ has the approximation property.
	\end{proof}
	
	There is an interesting application of Theorem \ref{thm:Lip-free-over-countable-uniformly-discrete} to Problem \ref{pb:Lipschitz-embed-rigidity-c0} regarding the Lipschitz rigidity of $\co$. As we are about to explain, there is a Banach space that contains a coarse-Lipschitz copy, but no bi-Lipschitz copy, of $\co$ and thus $\co$ is not coarse-Lipschitzly rigid. 
	
	We will denote by $[\bZ]_\infty^{(\bN)}$ the $\ell_\infty$-integer grid, i.e., the set of all finitely nonzero and $\bZ$-valued sequences equipped with the distance induced by the norm of $\co$.
	
	\begin{coro}
		The Banach space $\cF([\bZ]_\infty^{(\bN)})$ is universal for the class of separable metric spaces and coarse-Lipschitz embeddings. However, it is not universal for bi-Lipschitz embeddings and separable metric spaces.
	\end{coro}
	
	\begin{proof} 
		It is clear that $[\bZ]_\infty^{(\bN)}$ is a $(1,1)$-net of $\co$, meaning that it is $1$-separated and for any $x\in c_0$ there exists $y\in [\bZ]_\infty^{(\bN)}$ such that $\norm{x-y}_\infty\le 1$. A coarse-Lipschitz embedding from $\co$ into $[\bZ]_\infty^{(\bN)}$ is given by any map $f\colon c_0\to [\bZ]_\infty^{(\bN)}$ such that $\norm{f(x)-x}_\infty \le 1$ for all $x\in c_0$. It follows that $\co$ coarse-Lipschitz embeds into $\cF([\bZ]_\infty^{(\bN)})$ and, by Aharoni's theorem (Theorem \ref{thm:Aharoni}) that $\cF([\bZ]_\infty^{(\bN)})$ is a universal Banach space for the class of separable metric spaces and coarse-Lipschitz embeddings. 
		Since $[\bZ]_\infty^{(\bN)}$ is uniformly discrete, Theorem \ref{thm:Lip-free-over-countable-uniformly-discrete}, implies that $\cF([\bZ]_\infty^{(\bN)})$ has RNP. It now follows from the results from Chapter \ref{chapter:linear-reductions} that $\cF([\bZ]_\infty^{(\bN)})$ cannot contain a bi-Lipschitz copy of $\co$.
	\end{proof}
	
	For the second criterion, we will exploit the special properties of weakly null sequences in Lipschitz-free spaces. This criterion is due to C. Petitjean \cite{Petitjean2017}. As it is mentioned in \cite{Petitjean2017}, some of the ideas, as well as their applications, come from \cite{Kalton2004}. This second criterion relies on the notions of uniformly locally flat Lipschitz functions and norming subspaces.
	
	\begin{defi}
		Let $(M,d)$ be a metric space. A Lipschitz map $f\colon M \to \bR$ is said to be \emph{uniformly locally flat} if 
		\begin{equation*}
			\lim_{\eps \to 0}\ \sup \Big\{\frac{\abs{f(x)-f(y)}}{d(x,y)}\colon 0<d(x,y)\le \eps\Big\}=0.
		\end{equation*}
		We denote by $\lip_0(M)$ the space of all functions in $\Lip_0(M)$ that are uniformly locally flat on $M$.
	\end{defi}
	
	It is clear that $\lip_0(M)$ is a closed subspace of $\Lip_0(M)$ and it is often called the \emph{little Lipschitz space} over $M$. It is also obvious that every constant function is uniformly locally flat. An important observation is that when $(M,d)=(\bR, \abs{\cdot})$ the converse holds. Indeed, a locally flat function on $\bR$ must be differentiable almost everywhere, but since its derivative must vanish, it must be constant. Therefore, $\lip_0(\bR)=\{0\}$. On the other hand, $\lip_0(\bZ)=\Lip_0(\bZ)\equiv \ell_\infty$. 
	
	As for norming subspaces, recall that when $X$ is a Banach space and $c\in(0,1]$, a subspace $E$ of $X^*$ is called \emph{$c$-norming} for $X$ if $\frac{1}{c}\norm{x} \le  \sup_{x^*\in B_E}x^*(x)$ for all $x\in X$. 
	
	\begin{rema}
		\label{rem:1norming}
		It is worth pointing out that $E$ is $1$-norming if and only if $B_E$ is weak$^*$ dense in $B_{X^*}$ (for the nontrivial direction, use the Hahn-Banach theorem for the weak$^*$ topology).
	\end{rema}
	
	Subspaces that are $1$-norming in the setting of Lipchitz-free spaces admit the following characterization. 
	
	\begin{prop}
		\label{pro:1norming}
		Let $(M,d)$ be a metric space pointed at $0\in M$. For a subspace $E$ of $\Lip_0(M)=\cF(M)^*$ to be $1$-norming for $\cF(M)$, it is necessary and sufficient that for any finite subset $S$ of $M$ containing $0$, any $\eps>0$ and any $f\in \Lip_0(M)$, there exists $g\in E$ with $g(x)=f(x)$ for $x\in S$ and $\norm{g}_L\le (1+\eps)\norm{f}_L$.
	\end{prop}
	
	\begin{proof}
		Assume first that for any finite subset $S$ of $M$ containing $0$, any $\eps>0$ and any $f\in \Lip_0(M)$, there exists $g\in E$ with $g(x)=f(x)$ for $x\in S$ and $\norm{g}_L\le (1+\eps)\norm{f}_L$. Let $f\in B_{\Lip_0(M)}$ and $V$ be a weak$^*$ neighborhood of $f$. Recalling that on bounded sets of $\Lip_0(M)$ the weak$^*$ topology coincides with the topology of pointwise convergence on  $M$, we deduce the existence of a finite subset $S$ of $M$ containing $0$ and $\eps>0$ such that $W\cap B_{\Lip_0(M)} \subset V\cap B_{\Lip_0(M)}$, where $W := \{h \in \Lip_0(M)\colon \abs{f(x)-h(x)}<\eps \text{ for all }x\in S\}$. Since $S$ is finite, we can certainly pick $\eta>0$ such that $\eta \abs{f(x)}<\eps$ for all $x\in S$. By assumption, there exits $g\in E$ with $g(x)=f(x)$ for all $x\in S$ and $\norm{g}_L\le (1+\eta)$. Let $h :=g$ if $\norm{g}_L\le 1$ and $h :=\frac{g}{\norm{g}_L}$ otherwise. Then, the way $\eta$ was chosen and $h$ defined guarantees that $h\in W\cap B_{\Lip_0(M)}$. Since $E$ is a vector space, $h\in E$ and we have shown that $B_E$ is weak$^*$ dense in $B_{X^*}$. It follows from Remark \ref{rem:1norming} that $E$ is $1$-norming.
		
		The converse is a direct consequence of Helly's theorem. Let us, however, indicate an argument in this setting. So, assume that $E$ is $1$-norming and let $S :=\{0,x_1,\dots,x_n\}$ a finite subset of $M$, $f \in \Lip_0(M)$ and $H := \{g\in E\colon g(x_i)=f(x_i) \text{ for all }i\le n\}$. We may clearly assume that $f_{\restriction_S}\neq 0$ and since $H$ is closed we have $d(0,H)>0$. Then, we fix $b\in (0,d(0,H))$ and define $\Phi\colon \Lip_0(M) \to \bR^n$, by $\Phi(g) :=(g(x_i))_{i=1}^n$. Since $b<a$ and $\Phi$ is linear, $\Phi(f)$ does not belong to the $n$-dimensional convex set $C :=\Phi(bB_E)$. Moreover, the closure of $C$ in $\bR^n$ is $\Phi(bB_{X^*})$ since $\Phi$ is weak$^*$ to norm continuous, $B_E$ is weak$^*$ dense in $B_{X^*}$ and $B_{X^*}$ is weak$^*$ closed. We deduce, by applying the Hahn-Banach theorem in $\bR^n$ that there exist $a_1,\dots,a_n \in \bR^n$ such that for all $g \in bB_{X^*}$,
		\begin{equation}
			\label{eq:1norming}
			\sum_{i=1}^na_ig(x_i)\le \sum_{i=1}^na_if(x_i).
		\end{equation}
		Letting $\mu :=\sum_{i=1}^n a_i\delta_{x_i}\in \cF(M)$, we have  
		\begin{equation*}
			\norm{\mu}_{\cF(M)} = \sup_{g\in B_{X^*}}\abs{\langle g,\mu \rangle} \stackrel{\eqref{eq:1norming}}{\le} \frac{1}{b}\abs{\langle f,\mu \rangle}\le \frac{1}{b}\norm{\mu}_{\cF(M)}\norm{f}_L.
		\end{equation*}
		Since $b<d(0,H)$ was arbitrary, we have shown that $d(0,H)\le \norm{f}_L$. Since $d(0,H)=\inf\{\norm{g}_L\colon g\in E \text{ with }g(x_i)=f(x_i) \text{ for all }i\le n\}$ this finishes the proof. 
	\end{proof}
	
	As we observed above, $\lip_0(\bR)=\{0\}$ is not $1$-norming for the space $\cF(\bR)\equiv L_1$ which lacks the Schur property, while $\lip_0(\bZ)= \Lip_0(\bZ)$ is trivially $1$-norming for the Schur space $\cF(\bZ)\equiv \ell_1$. In fact, every Lipschitz-free space whose little Lipschitz space is large enough to be $1$-norming is a Schur space. 
	
	\begin{theo}
		\label{thm:1norming-Schur-free}
		Let $(M,d)$ be a metric space. If $\lip_0(M)$ is $1$-norming for $\cF(M)$, then $\cF(M)$ has the Schur property.
	\end{theo}
	
	\begin{proof}
		Since the Schur property is stable under $\ell_1$-sums and isomorphic embeddings, it follows from Theorem \ref{thm:almost-isometric-embed-free} that we may assume that $M$ is bounded. Assume now that $M$ is bounded and that $\cF(M)$ does not have the Schur property. Then, there exists a normalized weakly null sequence $(\gamma_n)_{n=1}^\infty$ in $\cF(M)$ and the main step of the proof will be to show that for all $\gamma \in \cF(M)$,
		\begin{equation}
			\label{eq:1norming-Schur-free}
			\liminf_{n\to \infty}\norm{\gamma+\gamma_n}\ge\norm{\gamma}+\frac12.
		\end{equation}
		Indeed, the validity of \ref{eq:1norming-Schur-free}) implies the existence of a subsequence $(\gamma_n')_n$ of $(\gamma_n)_n$ that is equivalent to the canonical basis of $\ell_1$; this is a blatant contradiction with the fact that $(\gamma_n)_n$ is weakly null and that $\ell_1$ has the Schur property. To see this, fix a strictly decreasing sequence $(c_N)_N$ in $(\frac13,\frac12)$ and build inductively $\gamma'_1,\ldots,\gamma'_N$ so that for all $a \in \bR^N$,
		\begin{equation*}
			\Big\|\sum_{k=1}^N a_i\gamma_k'\Big\| \ge c_N\sum_{k=1}^N \abs{a_k}.
		\end{equation*}
		Let us briefly indicate this classical recursive argument. Set $\gamma'_1 :=\gamma_1$ and assuming that $\gamma'_1,\dots,\gamma'_N$ have been chosen, observe that, by homogeneity, for all $a \in \bR^N$ and $t\in \bR$,
		\begin{equation*}
			\liminf_{n \to \infty}\Big\|\sum_{k=1}^N a_k\gamma'_k + t\gamma_n\Big\|\ge \Big\|\sum_{k=1}^N a_k\gamma'_k\Big\| + \frac12 \abs{t}\ge  c_N\sum_{k=1}^N \abs{a_k} + \frac12 \abs{t}.
		\end{equation*}
		Applying this, for each $(a,t)$ in a finite $\eta$-net of the unit ball of $\ell_1^{N+1}$, for $\eta>0$ small enough, combined with another homogeneity argument implies that for $n$ large enough we can pick $\gamma'_{N+1}=\gamma_n$.
		
		Returning to the proof of \eqref{eq:1norming-Schur-free}, fix $\gamma \in \cF(M)$ and $\eps>0$. Since $\lip_0(M)$ is $1$-norming, we can find $f\in \lip_0(M)$ such that $\norm{f}_L=1$ and $\langle \gamma,f \rangle\ge \norm{\gamma}-\eps$. Since $f\in \lip_0(M)$, there exists $\theta >0$ such that $\abs{f(x)-f(y)}\le \eps d(x,y)$ whenever $d(x,y)\le \theta$. We now fix $\delta \in (0,\frac{\eps\theta}{2})$ and it follows from Lemma \ref{lem:weakly-null-free} that 
		\begin{equation*}
			\inf_{|S|<\infty}\sup_{n\in \bN}\norm{\gamma_n}_{\cF(M)/\cF([S]_\delta)}=0.
		\end{equation*}  
		We can thus find a finite set $S$ containing $\{0\}$ so that there exists $\mu \in \cF(S)$ with $\norm{\mu-\gamma}\le \eps$ and for each $n\in \bN$, there exists $\mu_n \in \cF([S]_\delta)$ with $\norm{\mu_n-\gamma_n}\le \eps$. 
		
		Since $(\gamma_n)_n$ is normalized and weakly null, we have, by a standard gliding hump argument,  that for any finite-dimensional subspace $G$ of $\cF(M)$, $\liminf_n d(\gamma_n,G)= \liminf_n \norm{\gamma}_{\cF(M)/G}\ge \frac12$. In particular, $S$ being finite, $\liminf_{n} \norm{\gamma_n}_{\cF(M)/\cF(S)}\ge \frac12$. Hence, by the Hahn-Banach theorem, we can find a sequence $(f_n)_n$ in $\Lip_0(M)$ such that
		\begin{equation*}
			f_n =0 \text{ on S, } \norm{f_n}_L=1, \text{ and } \liminf_{n\to \infty} \langle \gamma_n,f_n \rangle\ge \frac12.
		\end{equation*}
		Denote now by $g_n$ the restriction of $f+f_n$ to  $[S]_\delta$ and consider $x,y \in [S]_\delta$. If $d(x,y)\le \theta$, then $\abs{g_n(x)-g_n(y)}\le (1+\eps)d(x,y)$. If $d(x,y)>\theta$, then we can find $u,v \in S$ with $d(x,u)\le \delta$ and $d(y,v)\le \delta$, so that 
		\begin{align*}
			\abs{g_n(x)-g_n(y)} & \le \abs{f(x)-f(y)} + \abs{f_n(x)-f_n(u)} + \abs{f_n(y)-f_n(v)}\\
			&\le d(x,y)+d(x,u)+d(y,v) \\
			& \le (1+\frac{2\delta}{\theta})d(x,y)\le (1+\eps)d(x,y).   
		\end{align*}
		So $g_n$ is $(1+\eps)$-Lipschitz on $[S]_\delta$. Then, since $\mu+\mu_n \in \cF([S]_\delta)$ and $\langle \mu, f_n\rangle =0$,
		$$\norm{\mu + \mu_n} \ge \frac{1}{1+\eps}\langle \mu+\mu_n,g_n\rangle = \frac{1}{1+\eps}\big(\langle \mu, f\rangle+\langle \mu_n, f\rangle+\langle \mu_n, f_n\rangle\big).$$
		Now $\langle \mu, f\rangle \ge \norm{\gamma} - 2\eps$ and $\liminf_n\langle \mu_n, f_n\rangle \ge \frac12 -\eps$. Moreover, since $\lim_n\langle \gamma_n, f\rangle=0$, we have $\limsup_n \abs{\langle \mu_n, f\rangle}\le \eps$. Summing up, we get 
		\begin{equation*}
			\liminf_{n\to \infty}\norm{\mu+\mu_n} \ge (1+\eps)^{-1}(\norm{\gamma} +\frac12-4\eps)
		\end{equation*} 
		and so 
		\begin{equation*}
			\liminf_{n\to \infty} \norm{\gamma+\gamma_n}\ge (1+\eps)^{-1}(\norm{\gamma} + \frac12-4\eps)-2\eps.
		\end{equation*}
		Since $\eps>0$ was arbitrary, this finishes the proof of (\ref{eq:1norming-Schur-free}).
	\end{proof}
	
	\begin{rema}
		Theorem \ref{thm:1norming-Schur-free} reveals again that if $M$ is uniformly discrete, then $\cF(M)$ has the Schur property. In Theorem \ref{thm:1norming-Schur-free}, the reason we need $\lip_0(M)$ to be $1$-norming is to control small distances. 
	\end{rema}
	
	The third and last criterion applies to Lipschitz-free spaces over certain metric transforms. Assume $\omega\colon [0,\infty) \to [0,\infty)$ is nondecreasing, subadditive and satisfies $\omega(t)=0$ if and only if $t=0$. Then, such a function induces a new metric on a metric space $(M,d)$ by letting $d_\omega :=\omega \circ d$ on $M$. The metric $d_\omega$ is usually called a \emph{metric transform} of $d$. It is not difficult to see that such properties are satisfied by any map $\omega\colon [0,\infty) \to [0,\infty)$ that is increasing, concave and satisfies $\omega(0)=0$. In the sequel, we consider metric transforms induced by gauges.
	
	%If a nondecreasing function $\omega\colon [0,\infty) \to [0,\infty)$ is concave and satisfies $\omega(0)=0$, it is not too difficult to see that $\omega$ is subadditive. Therefore, such function $\omega$ induces a new metric on a metric space $(M,d)$ by letting $d_\omega :=\omega \circ d$ on $M$. The metric $d_\omega$ is usually called a \emph{metric transform} of $d$ and, when $0<\alpha<1$, the metric space $(M,d^{\alpha})$ is usually referred to as the $\alpha$-snowflaking of $M$. 

	\begin{defi}
		A \emph{gauge} is a continuous, nondecreasing and subadditive function $\omega\colon [0,\infty) \to [0,\infty)$ that
		satisfies 
		\begin{enumerate}[(i)]
			\item $\omega(0)=0$,
			\item $\omega(t)\ge t$ for all $t\in [0,1]$
		\end{enumerate}
		We say that $\omega$ is a \emph{nontrivial gauge} if $\lim_{t\to 0}\frac{\omega(t)}{t}=\infty$. 
	\end{defi}
	The most natural nontrivial gauge is given by $t \mapsto t^\alpha$, where $0<\alpha<1$. For $\alpha \in (0,1)$, we will also need the map  $t\mapsto \omega_\alpha(t):=\max\{t,t^\alpha\}$, which is a nontrivial gauge, but not concave. We will use the convenient notation $\cF_\omega(M):=\cF(M,\omega \circ d)$ and $\cF_\alpha(M) :=\cF((M,d^\alpha))$. The dual of $\cF_\omega(M)$ is the space $\Lip_0(M,\omega \circ d)$. Lipschitz-free spaces of the form $\cF_\omega(M)$ for a nontrivial gauge have little Lipschitz spaces that are $1$-norming.
	
	\begin{theo}
		\label{thm:1norming-Lip-free} 
		Suppose that $(M,d)$ is any metric space and $\omega$ is a nontrivial gauge. Then, $\lip_0(M,\omega \circ d)$ is $1$-norming for $\cF_\omega(M)$.
	\end{theo}
	
	Before starting the proof, we first construct an approximating sequence for any gauge. We define, for $n\in \bN$:
	\begin{equation*}
		\omega_n(t) :=\inf\{\omega(s)+n(t-s)\colon 0\le s\le t\},\ \ t\ge 0.
	\end{equation*}
	We list the main properties of these maps in the following lemma. 
	\begin{lemm}
		\label{lem:approx-gauge} The sequence of maps $(\omega_n)_n$ satisfies the following properties.
		\begin{enumerate}[(i)]
			\item The sequence $(\omega_n)_{n=1}^\infty$ is nondecreasing and pointwise converging to $\omega$.
			\item For all $0\le t\le t'$, $\omega_n(t)\le \omega_n(t')\le \omega_n(t)+n(t'-t)$. In particular, $\omega_n$ is $n$-Lipschitz for all $n\in \bN$.
			\item For all $n\in \bN$, $\omega_n$ is subadditive.
			\item For all $0\le t\le t'$ and $n\in \bN$, $\omega_n(t')-\omega_n(t)\le \omega(t'-t)$.
		\end{enumerate}
	\end{lemm}
	
	\begin{proof}
		%For $n\in \bN$ and $0\le s\le t$, denote $\Phi_{t}^n(s)=\omega(s)+n(t-s)$.
		\begin{enumerate}[$(i)$]
			\item Clearly, for all $t\ge 0$, $\omega_n(t) \le \omega(t)$ and $(\omega_n)_{n=1}^\infty$ is nondecreasing. Fix now $t>0$ and $\lambda<\omega(t)$. Since $\omega$ is continuous, there exists $0\le t_0<t$ such that $\omega(s)+n(t-s) \ge \omega(s)>\lambda$ for all $s \in [t_0,t]$ and all $n\in \bN$. On the other hand, for $n$ large enough, $\omega(s)+n(t-s)>\lambda$ for all $s\in [0,t_0]$. It follows that $\lim_n\omega_n(t)=\omega(t)$. The case $t=0$ is trivial. 
			\item Let now $0\le t\le t'$. Then, for $t\le s\le t'$, $\omega_n(t) \le \omega(t) \le \omega(s) \le \omega(s)+n(t'-s)$,
			and hence 
			\begin{equation}
				\label{eq1:approx-gauge}
				\omega_n(t)\le \inf_{s\in[t,t']} [\omega(s)+n(t'-s)].
			\end{equation}
			On the other hand, for $0\le s \le t \le t'$,  
			\begin{equation*}
				\omega(s)+n(t-s) \le \omega(s)+n(t'-s) \le \omega(s)+n(t-s) + n(t' - t),
			\end{equation*} 
			so 
			\begin{equation}
				\label{eq2:approx-gauge}
				\omega_n(t) \le \inf\limits_{s \in [0, t]} [\omega(s)+n(t'-s)] \le \omega_n(t) + n(t' - t).
			\end{equation}
			It follows from \eqref{eq1:approx-gauge} and \eqref{eq2:approx-gauge} that $\omega_n(t)\le \omega_n(t')\le \omega_n(t)+n(t'-t)$.
			\item Let $n\in \bN$ and $t,t' \ge 0$. Since $\omega$ is continuous, there exist $s \in [0, t]$ and $s' \in [0, t']$ such that $\omega_n(t) = \omega(s)+n(t-s)$ and $\omega_n(t') = \omega(s')+n(t'-s')$. Since $\omega$ is subadditive, we have 
			\begin{align*}
				\omega_n(t) + \omega_n(t') & = \omega(s) + n(t - s) + \omega(s') + n(t' - s') \\
				& \ge \omega(s + s') + n(t + t' - (s + s')),
			\end{align*} 
			with $s + s' \in [0, t + t']$. Thus, $\omega_n(t) + \omega_n(t') \ge \omega_n(t + t')$.
			\item Let $0\le t \le t'$. Since $\omega_n$ is subadditive, we have $\omega_n(t')-\omega_n(t)\le \omega_n(t'-t)\le \omega(t'-t)$.
		\end{enumerate}    
	\end{proof}
	
	\begin{proof}[Proof of Theorem \ref{thm:1norming-Lip-free}]
		%In view of Theorem \ref{thm:1norming-Schur-free}, it is enough to show that $\lip_0(M,\omega \circ d)$ is $1$-norming for $\cF_\omega(M)$. 
		First, we fix $y\in M$ and let 
		\begin{equation*}
			h_n(x) := \omega_n(d(x,y)) - \omega_n(d(y,0)).
		\end{equation*}
		It follows from Lemma \ref{lem:approx-gauge} that  $h_n(0)=0$, $h_n$ is $1$-Lipschitz on $(M,d_\omega)$ and $n$-Lipschitz on $(M,d)$. Therefore, because $\omega$ is nontrivial, $h_n$ is uniformly locally flat on $(M,d_\omega)$. Also, for any $x\in M$, $h_n(x)-h_n(y)$ tends to $d_\omega(x,y)$. It follows that for any $\eps>0$ and any $x,y \in M$ there exists $h_{x,y}\in \lip_0(M,d_\omega)$ which is $(1+\eps)$-Lipschitz on $(M,d_\omega)$ and such that $\abs{h_{x,y}(x) - h_{x,y}(y)}=d_\omega(x,y)$.
		
		We are now ready to apply Proposition \ref{pro:1norming}. So, let $S$ be a finite subset of $M$ containing $0$, $\eps>0$ and $f\in \Lip_0(M,\omega \circ d)$. Then, for any $x,y\in S$, adding a constant to $h_{x,y}$ if necessary, we can find $g_{x,y} \in \lip(M,\omega \circ d)$ such that $g_{x,y}$ is $(1+\eps)$-Lipschitz on $(M,d_\omega)$, $g_{x,y}(x)=f(x)$ and $g_{x,y}(y)=f(y)$ (possibly $g_{x,y}(0)\neq 0)$. Then, we use the fact that $\lip(M,\omega \circ d)$ is a lattice to deduce that 
		\begin{equation*}
			g :=\max_{x\in S} \min_{y\in S\setminus\{x\}} g_{x,y}
		\end{equation*}
		belongs to $\lip(M,\omega \circ d)$. Since $g$ coincides with $f$ on $S$, it follows that $g(0)=0$ and hence $g\in \lip_0(M,\omega \circ d)$). Noting that $g$ is $(1+\eps)$-Lipschitz on $(M,d_\omega)$, it follows from Proposition \ref{pro:1norming} that $\lip_0(M,\omega \circ d)$ is $1$-norming for $\cF_\omega(M)$.
	\end{proof}
	
	The following corollary follows immediately from Theorem \ref{thm:1norming-Schur-free}.
	
	\begin{coro}
		\label{cor:Schur-Lip-gauge-free}
		Suppose that $M$ is any metric space and $\omega$ is a nontrivial gauge. Then, $\cF_\omega(M)$ has the Schur property.
	\end{coro}
	
	This last result is taken from \cite{Kalton2004}, where the question of whether $\cF_\omega(M)$ has the Radon-Nikod\'ym property whenever $\omega$ is a nontrivial gauge was left open. It is important to mention here that one of the main results from \cite{AGPP2021} is that for a metric space $M$, $\cF(M)$ has the Schur property if and only if it has the Radon-Nikod\'ym property if and only if the metric space $M$ is purely $1$-unrectifiable (i.e. it does not contain any bi-Lipschitz copy of a subset of $\bR$ of positive Lebesgue measure). So, $\cF_\omega(M)$ indeed has the Radon-Nikod\'ym property whenever $\omega$ is a nontrivial gauge.
	Let us also mention that there exists a uniformly discrete metric space $M$ such that $\cF(M)$ is a Schur space (since for a uniformly discrete metric space $\Lip_0(M)=\lip_0(M)$), has the Radon-Nikod\'ym property but is not isometrically a dual space (see \cite{GPPR2018}).
	Finally, let us mention a result due to A. Dalet \cite{Dalet2015} that preceded \cite{AGPP2021}: if $M$ is a countable compact metric space, then $\cF(M)$ is isometric to $\lip_0(M)^*$ (through the usual duality). In particular, being a separable dual space, it has the Radon-Nikod\'ym property. A. Dalet also showed in \cite{Dalet2015} that if $M$ is a countable compact metric space, then $\cF(M)$ has the metric approximation property. However, the following problem remains open.
	
	\begin{prob}
		\label{pb:Lipschitz-free-RNP}
		Does there exist a separable metric space $M$ such that $\cF(M)$ has the Radon-Nikod\'ym property but is not isomorphic to a subspace of a separable dual Banach space?   
	\end{prob}
	
	We conclude this section with another important universal space.
	
	\begin{coro}
		For $\alpha \in (0,1)$ and $t\mapsto \omega_\alpha(t):=\max\{t,t^\alpha\}$, $\cF_{\omega_\alpha}(c_0)$ provides an example of a Banach space which has the Schur property, the Radon-Nikod\'ym property and is universal for separable metric spaces and embeddings that are at the same time uniform and coarse-Lipschitz embeddings. However, it is not universal for separable metric spaces and bi-Lipschitz embeddings. 
	\end{coro}
	
	\begin{proof}
		The identity from $(c_0,\norm{\cdot}_\infty)$ to $(c_0,\omega_\alpha\circ \norm{\cdot}_\infty)$ is both a uniform and coarse-Lipschitz embedding (and even an equivalence). The universality of $\cF_{\omega_\alpha}(c_0)$ then follows from Aharoni's theorem and the isometric embeddability of $(c_0,\omega_\alpha\circ \norm{\cdot}_\infty)$ into $\cF_\alpha(c_0)$. Finally, since $\omega_\alpha$ is clearly a nontrivial gauge, $\cF_{\omega_\alpha}(c_0)$ has RNP and hence it cannot contain a bi-Lipschitz copy of $\co$.
	\end{proof}

	\section{Embeddability of stable metric spaces into reflexive spaces}
	\label{sec:stable-into-reflexive}
	
	The notion of stability was introduced, originally for (separable) Banach spaces, in the work of Krivine and Maurey  \cite{KrivineMaurey1981}. Stability will be discussed in detail in Section \ref{sec:stability}. For metric spaces, it takes the following form.
	
	\begin{defi} 
		A metric space $(M,d)$ is said to be \emph{stable} if for every nonprincipal ultrafilters $\cU$ and $\cV$ over $\bN$ and every bounded sequences $(x_n)_{n=1}^\infty$ and $(y_n)_{n=1}^\infty$ in $M$, we have
		\begin{equation*}
			\lim_{n\in \cU}\lim_{m \in \cal \cV} d(x_n,y_m) = \lim_{m\in \cV}\lim_{n \in \cU} d(x_n,y_m).
		\end{equation*}
	\end{defi}
	
	\begin{rema}\label{rem:iteratedlimits} This definition is equivalent to the property that for every bounded sequences $(x_n)_{n=1}^\infty$ and $(y_n)_{n=1}^\infty$ in $M$ such that both iterated limits $\lim_{n\to \infty}\lim_{m \to \infty} d(x_n,y_m)$ and $\lim_{m\to \infty}\lim_{n \to \infty} d(x_n,y_m)$ exist, they are equal. We shall only need the trivial implication. It is also easy to see that if $(x_n)_{n=1}^\infty$ and $(y_n)_{n=1}^\infty$ are unbounded sequences in $M$ such that both iterated limits $\lim_{n\to \infty}\lim_{m \to \infty} d(x_n,y_m)$ and $\lim_{m\to \infty}\lim_{n \to \infty} d(x_n,y_m)$ exist in $[0,\infty]$, then they are both equal to $\infty$. 
	\end{rema}
	
	It is worth pointing out that any $L_p$-space is stable for $p\in[1,\infty)$. The main result of this section is due to Kalton \cite{Kalton2007}, who showed that any stable metric space can be coarsely and uniformly embedded into a reflexive Banach space. In fact, Kalton's construction can even be performed so that the embedding is isometric for distances away from $0$ and $\infty$. The following improved statement is taken from \cite{BaudierLancien2015}, and we shall detail this more precise version. As the reader will see, the techniques used in the proof are inspired by Lipschitz-free spaces.
	
	\begin{theo}
		\label{thm:Kalton-stable-reflexive}
		Let $(M,d)$ be a stable metric space. Assume that $\rho,\omega :[0,\infty)\to [0,\infty)$ are continuous, nondecreasing and such that $t\mapsto \frac{\rho(t)}{t}$ and $t\mapsto \frac{\omega(t)}{t}$ are nonincreasing. Assume also that
		\begin{enumerate}[(i)]
			\item $t\le \omega(t)$ for $t\in [0,1]$ and $\omega(t)=t$ for $t\in[1,\infty)$,
			\item $\omega(0)=0$ and $\lim_{t\to 0}\frac{\omega(t)}{t}=+\infty$,
			\item $\rho(t)=t$ for $t\in [0,1]$ and $\rho(t)\le t$ for $t\in[1,\infty)$,
			\item $\lim_{t\to +\infty}\frac{\rho(t)}{t}=0$,
		\end{enumerate}
		Then, there exists a reflexive Banach space $(Y,\norm{\cdot})$ and a map $f\colon M\to Y$ such that for all $x,y \in M$, $$\rho(d(x,y))\le \|f(x)-f(y)\| \le \omega(d(x,y)).$$
	\end{theo}
	
	\begin{rema} 
		\label{rem:embedstable}
		This statement applies, for instance, to the following functions $\omega$ and $\rho$. Fix $0<A\le 1 \le B<\infty$ and set $\omega(t):=\sqrt{At}$  if $t\in [0,A]$, $\omega(t):=t$ if $t\ge A$, $\rho(t):=t$ if $t\in [0,B]$ and $\rho(t):=\sqrt{tB}$ if $t\ge B$. Note that, for these functions, Theorem \ref{thm:Kalton-stable-reflexive} provides an embedding which is isometric for distances in $[A,B]$.
	\end{rema}

	\begin{proof}[Proof of Theorem \ref{thm:Kalton-stable-reflexive}] 
		Fix an arbitrary point (denoted by $0$) in $M$ and define
		\begin{equation*}
			\text{Lip}_0^\omega(M):=\Big\{g\colon M\to M\ \mid \ g(0)=0\text{ and }\sup_{x\neq y\in M}\frac{|g(y)-g(x)|}{\omega(d(x,y))}<+\infty\Big\},
		\end{equation*}
		and, for all $g\in \text{Lip}_0^\omega(M)$, let 
		\begin{equation*}
			N_\omega(g) :=\sup_{x\neq y\in M}\frac{\abs{g(y)-g(x)}}{\omega(d(x,y))}.
		\end{equation*}
		At this point it is important to note that we have not assumed that $\omega$ is concave or subadditive. So, $\omega \circ d$ is not necessarily a metric. However, it is easily checked that $(\text{Lip}_0^\omega(M),N_\omega)$ is a Banach space. Define now the map $\delta\colon M\to \text{Lip}_0^\omega(M)^*$ by letting $\delta(x)(g)=g(x)$ for all $x\in M$. By construction, $N_\omega^*(\delta(x)-\delta(y))\le \omega(d(x,y))$ for all $x,y\in M$, where $N_\omega^*$ is  the dual norm of $N_\omega$. Again, we shall not speak of the Lipschitz-free space $\cF_\omega(M)$, as $\omega \circ d$ is not a necessarily a metric.
		
		Then, for $(p,q) \in \tilde M :=\{(p,q)\in M^2 \colon\ p\neq q\}$ and $x\in M$, we define 
		\begin{equation*}
			g_{p,q}(x) :=\max\{d(p,q)-d(q,x);0\}-\max\{d(p,q)-d(q,0);0\}.
		\end{equation*}
		For all $x,y\in M$ and $(p,q) \in \tilde M$, we have
		\begin{align}
			\abs{g_{p,q}(x)-g_{p,q}(y)}\le \min\{d(p,q);d(x,y)\}.
		\end{align}
		Now let
		\begin{equation*}
			h_{p,q} :=\frac{\rho(d(p,q))}{d(p,q)}g_{p,q}.
		\end{equation*}
		Note that $N_\omega(h_{p,q})\le 1$ and for all $(p,q) \in \tilde M$
		\begin{equation*}
			\abs{\langle h_{p,q},\delta(p)-\delta(q)\rangle}=\rho(d(p,q)).
		\end{equation*}
		We will need the following lemma.
		\begin{lemm}
			\label{lem:limits}
			For $(p,q) \in \tilde M$ and $(x,y)\in \tilde M$, consider 
			\begin{equation*}
				R_{p,q}(x,y) := \frac{\abs{h_{p,q}(x)-h_{p,q}(y)}}{\omega(d(x,y))}.
			\end{equation*}
			Then,
			\begin{enumerate}[(i)]
				\item $R_{p,q}(x,y)\to 0$ as $d(p,q)\to +\infty$ or $d(p,q)\to 0$, uniformly in $(x,y)\in \tilde M$,
				\item $R_{p,q}(x,y)\to 0$ as $d(x,y)\to +\infty$ or $d(x,y)\to 0$, uniformly in $(p,q)\in \tilde M$.
			\end{enumerate}
		\end{lemm}
		
		\begin{proof}[Proof of Lemma \ref{lem:limits}] 
			For all $(x,y)\in \tilde M$, we have
			\begin{equation*}
				R_{p,q}(x,y)\le \frac{\rho(d(p,q))}{d(p,q)}\frac{\min\{d(p,q);d(x,y)\}}{\omega(d(x,y))}.
			\end{equation*}
			$(i)$ First, for all $(x,y)\in \tilde M$,
			\begin{equation*}
				R_{p,q}(x,y)\le \frac{\rho(d(p,q))}{d(p,q)}\frac{d(x,y)}{\omega(d(x,y))}\le \frac{\rho(d(p,q))}{d(p,q)},
			\end{equation*}
			and thus $R_{p,q}(x,y)\to 0$ as $d(p,q)\to +\infty$  uniformly in $(x,y)\in \tilde M$.
			
			Assume now that $d(p,q)\le 1$. If $d(x,y)\le d(p,q)$ then
			\begin{equation*}
				R_{p,q}(x,y)\le\frac{d(x,y)}{\omega(d(x,y))}\le \frac{d(p,q)}{\omega(d(p,q))},
			\end{equation*} 
			because $t\mapsto \frac{t}{\omega(t)}$ is nondecreasing. Finally, if $d(x,y)\ge d(p,q)$ then
			\begin{equation*}
				R_{p,q}(x,y)\le\frac{d(p,q)}{\omega(d(x,y))}
				\le\frac{d(p,q)}{\omega(d(p,q))},
			\end{equation*}
			because $\omega$ is nondecreasing. Therefore, $R_{p,q}(x,y)\to 0$ as $d(p,q)\to 0$,  uniformly in $(x,y)\in \tilde M$.
			
			\smallskip
			
			$(ii)$ Assume first that $d(x,y)>1$. If $d(p,q)\le d(x,y)$ then
			\begin{equation*}
				R_{p,q}(x,y)\le\frac{\rho(d(p,q))}{d(p,q)}\frac{d(p,q)}{d(x,y)}\le \frac{\rho(d(x,y))}{d(x,y)},
			\end{equation*}
			because $\rho$ is nondecreasing. If $d(p,q)\ge d(x,y)$, then
			\begin{equation*}
				R_{p,q}(x,y)\le\frac{\rho(d(p,q))}{d(p,q)}\le \frac{\rho(d(x,y))}{d(x,y)},
			\end{equation*}
			because $t\to\frac{\rho(t)}{t}$ is nonincreasing. Therefore, $R_{p,q}(x,y)\to 0$ as $d(x,y)\to +\infty$  uniformly in $(p,q)\in \tilde M$.
			
			Assume now that $d(x,y)\le 1$. Then,
			\begin{equation*}
				R_{p,q}(x,y)\le\frac{\rho(d(p,q))}{d(p,q)}\frac{d(x,y)}{\omega(d(x,y))}\le \frac{d(x,y)}{\omega(d(x,y))}.
			\end{equation*}
			Therefore, $R_{p,q}(x,y)\to 0$ as $d(x,y)\to 0$  uniformly in $(p,q)\in \tilde M$.\\
		\end{proof}
		
		The next lemma, which is a key ingredient of the proof, is a consequence of Lemma \ref{lem:limits} and Eberlein-\v{S}mulian theorem.
		
		\begin{lemm}
			\label{lem:weak-compact} 
			The set $W:=\{h_{p,q}\colon (p,q)\in \tilde M\}$ is weakly relatively compact in $\Lip_0^\omega(M)$.
		\end{lemm}
		
		\begin{proof}[Proof of Lemma \ref{lem:weak-compact}] 
			It follows from the Eberlein-\v{S}mulian Theorem that it is enough to show that any sequence $(h_n)_{n=1}^\infty=(h_{p_n,q_n})_{n=1}^\infty$ with $(p_n,q_n)\in \tilde M$ admits a weakly convergent subsequence. The closed linear span of $\{h_n\colon n\ge 1\}$ in $\Lip_0^\omega(M)$, denoted by $[h_n\colon n\ge 1]$ is separable. So, there exists a countable subset $M_0$ of $M$ containing $0$ and all $p_n,q_n$ for $n\ge 1$ such that for all $g\in [h_n\colon n\ge 1]$
			\begin{equation*}
				\norm{g}:=\norm{g}_{\Lip_0^\omega(M)} = \norm{g_{\restriction_{M_0}}}_{\Lip_0^\omega(M_0)}.
			\end{equation*}
			Clearly, for each $x\in M$, the sequence $(h_n(x))_n$ is bounded in $\bR$. Thus, using a diagonal argument, we may assume that $(h_n)_{n=1}^\infty$ converges pointwise on $M_0$ to a function $h$. We may also assume that $\lim_{n\to +\infty}d(p_n,q_n)=r\in [0,+\infty]$ and that for all $x\in M_0$ $\lim_{n\to +\infty}d(x,q_n)$ exists in $[0,+\infty]$. In particular, we let $s:=\lim_{n\to +\infty}d(0,q_n)$. Now we define $V\colon [h_n\colon n\ge 1]\to \ell_\infty(\tilde{M_0})$, where $\tilde{M_0} :=\{(p,q)\in M_0^2\colon p\neq q\}$ by letting for all $g\in [h_n\colon n\ge 1]$,
			\begin{equation*}
				Vg := \Big(\frac{g(x)-g(y)}{\omega(d(x,y))}\Big)_{(x,y)\in \tilde{M_0}}.
			\end{equation*}
			By construction, $V$ is a linear isometry from $[h_n\colon n\ge 1]$ into $\ell_\infty(\tilde{M_0})$. So, it is enough for us to show that there exists a subsequence of $(Vh_n)_{n=1}^\infty$ which is weakly convergent to $Vh$. We denote by $A$ the closed subalgebra of $\ell_\infty(\tilde{M_0})$, generated by the constant functions, $V([h_n;\ n\ge 1])$, the maps $(x,y)\mapsto \frac{h(x)-h(y)}{\omega(d(x,y))}$, $(x,y)\mapsto \arctan(d(x,u))$ and $(x,y)\mapsto \arctan(d(y,u))$, for $u\in M_0$ and the map $(x,y)\mapsto \arctan(d(x,y))$. Since $A$ is separable, there exists a metrizable compactification $K$ of $\tilde{M_0}$ such that every $f$ in $A$ admits a unique continuous extension to $K$. With this description of $A$, let us first show that for any $\xi \in K$, there exists a subsequence of $(h_n)_n$ such that, $\lim_{n\to\infty}Vh_n(\xi)=Vh(\xi)$.
			
			%it follows from the dominated convergence theorem and a diagonal argument, that we only need to prove that for any $\xi \in K$, there exists a subsequence of $(h_n)_n$ such that, $\lim_{n\to\infty}Vh_n(\xi)=Vh(\xi)$. \\
			If $r=0$ or $r=\infty$, it follows from $(i)$ in Lemma \ref{lem:limits} that $\lim_{n\to \infty}\norm{Vh_n}_\infty=0$. Therefore, $\lim_{n\to \infty}\norm{h_n}=0$ and $h\equiv0$. Thus, we may assume that $0<r<\infty$.\\
			We now fix $\xi\in K$. Then, pick $(x_m,y_m)\in \tilde{M_0}$ such that $(x_m,y_m)\to \xi$ in $K$ and let $t=\lim_{m\to \infty}d(x_m,y_m)$. First note that if $t=0$ or $t=\infty$, it follows from $(ii)$ in Lemma \ref{lem:limits} that for all $n\ge 1$, $\lim_{n\to \infty}Vh_n(\xi)=Vh(\xi)=0$.\\
			Thus, for the sequel, we may and will assume that $r,t \in (0,\infty)$. Note that we already assumed that $\lim_nd(q_n,x_m)$. By taking a subsequence of $(x_m)_m$, we may also assume that $\lim_m\lim_n d(q_n,x_m)$ exists in $[0,+\infty]$. Similarly, $\lim_md(q_n,x_m)$ already exists because of the inclusion of the functions $(x,y)\mapsto \arctan(d(x,u))$ in $A$. Then, we can pass to a subsequence of $(q_n)_n$, depending on $\xi$, such that $\lim_n\lim_m d(q_n,x_m)$ exists in $[0,+\infty]$. Then, it follows from the stability of $M$ and Remark \ref{rem:iteratedlimits} that $\lim_m\lim_n d(q_n,x_m)=\lim_n\lim_m d(q_n,x_m)$ and therefore that $\lim_m\lim_n h_n(x_m)=\lim_n\lim_m h_n(x_m)$.\\
			Since $Vh\in C(K)$ and $(h_n)_{n\ge 1}$ converges pointwise to $h$ on $M_0$, we have
			\begin{equation*}
				Vh(\xi)=\lim_{m\to \infty}Vh(x_m,y_m)=\lim_{m\to \infty}\frac{h(x_m)-h(y_m)}{\omega(d(x_m,y_m))}=\lim_{m\to\infty}\lim_{n\to \infty}\frac{h_n(x_m)-h_n(y_m)}{\omega(d(x_m,y_m))}.
			\end{equation*}
			Finally, since $M$ is stable  and $Vh_n \in C(K)$ we obtain
			\begin{equation*}
				Vh(\xi)=\lim_{n\to\infty}\lim_{m\to \infty}\frac{h_n(x_m)-h_n(y_m)}{\omega(d(x_m,y_m))}=
				\lim_{n\to\infty}Vh_n(\xi).
			\end{equation*}
			Until now,  the subsequence $(h_n)_n$ depends on $\xi$, but the compact $K$ being metrizable, by using a diagonal argument to pass to a further subsequence, we may assume that $\lim_{n\to\infty}Vh_n(\xi)=Vh(\xi)$, for all $\xi \in K$. Finally, it follows from the dominated convergence theorem that $(Vh_n)_n$ is weakly converging to $Vh$ in $C(K)$. This concludes the proof of this lemma.

			%If $L$ is the spectrum of $B=\ell_\infty(\tilde{M_0})$, it is easily seen that $\{\xi\restrict{A},\ \xi\in L\}$ is closed in $(K,\sigma(A,A^*))$. It also norming for $A=C(K)$. Then, it follows that $K=\{\xi\restrict{A},\ \xi\in L\}$. Now, it known that $L$ is the Stone-Cech compactification of $\tilde{M_0}$.
			
		\end{proof}
		
		We now proceed with the construction of the embedding $f$. Consider the operator 
		\begin{align*}
			S  & \colon\ell_1(W) \to (\Lip_0^\omega(M),N_\omega)\\
			& \xi=(\xi_h)_{h\in W} \mapsto  S(\xi):=\sum_{h\in W}\xi_h h.
		\end{align*}
		Since every $h \in W$ is in the unit ball of $N_\omega$, we have that $\norm{S}\le 1$. Moreover, it follows from Lemma \ref{lem:weak-compact} and the Krein-\v{S}mulyan Theorem that $S$ is a weakly compact operator. Then, the isometric version, due to Lima, Nygaard and Oja \cite{LimaNygaardOja2000}, of a fundamental factorization theorem of Davis, Figiel, Johnson and Pe\l czy\'nski \cite{DavisFJP74} yields the existence of a reflexive Banach space $X$ and of linear maps $T\colon \ell_1(W)\to X$ and $U \colon X\to  \Lip_0^\omega(M)$ such that $\norm{T}\le 1$, $\norm{U}\le 1$ and $S=UT$. Then, we define $f\colon M\to Y=X^*$ by letting for all $x\in M$,
		\begin{equation*}
			f(x) := U^*(\delta(x)).
		\end{equation*}
		First, we clearly have that for all $x,y\in M$
		\begin{equation*}
			\norm{f(x)-f(y)}\le \norm{U^*}\,N_\omega^*\big(\delta(x)-\delta(y)\big)\le \omega(d(x,y)).
		\end{equation*}
		On the other hand, since $S^*=T^*U^*$ and $\norm{T^*}\le 1$, we have that for all $(x,y)\in \tilde M$,
		\begin{equation*}
			\norm{f(x)-f(y)}\ge \norm{S^*(\delta(x)-\delta(y))}_{\ell_\infty(W)}\ge |\langle h_{x,y},\delta(x)-\delta(y)\rangle |=\rho(d(x,y)).
		\end{equation*}
		This finishes the proof, as $f$ provides the desired embedding into the reflexive Banach space $Y=X^*$.
	\end{proof}
	
	The faithfulness of the embedding can be improved if the stable space has additional structural properties.
	
	\begin{coro}\, 
		\begin{enumerate}[(i)]
			\item If $M$ is a stable, bounded and uniformly discrete metric space, then $M$ embeds isometrically into a reflexive Banach space.
			\item If $M$ is a stable and uniformly discrete metric space, then $M$ embeds bi-Lipschitzly into a reflexive Banach space.
			\item Every stable metric space is coarse-Lipschitz embeddable into a reflexive Banach space.
		\end{enumerate}
	\end{coro}
	
	\begin{proof} 
		$(i)$ Let $0<A \le 1\le B<\infty$ such that $A\le \sep(M)\le \diam\,(M)\le B$. We just have to apply Theorem \ref{thm:Kalton-stable-reflexive} with the functions $\omega$ and $\rho$ given in Remark \ref{rem:embedstable}.
		
		$(ii)$ Fix $x_0 \in M$ and, for $k\in \bN$, denote by $B_k$ the ball with center $x_0$ and radius $2^k$. Then, by $(i)$, for all $k \in \bN$ there exists a reflexive Banach space $X_k$ such that $B_k$ isometrically embeds into $X_k$. Now, we can apply a barycentric gluing technique, due to the first author, to deduce that there exists a bi-Lipschitz embedding from $M$ into the reflexive space $X=(\sum_{k=1}^\infty X_k)_{\ell_2}$. We refer the reader to \cite{Baudier2022} for a survey on the applications of this technique and, more specifically, to Theorem 2 in \cite{Baudier2022} for our purpose.
		
		$(iii)$ Of course, we may assume that $(M,d)$ is unbounded. Fix a maximal $1$-separated subset $N$ of $M$. By $(ii)$, there exists a reflexive Banach space $X$ and a bi-Lipschitz embedding $g$ from $(N,d)$ into $X$. We can also define a map $h\colon M \to N$ such that for all $x\in M$, $d(x,h(x))\le 1$. Then, it is easy to check that $g \circ h$ is a coarse-Lipschitz embedding from $(M,d)$ into $X$.
		
	\end{proof}
	
	\section{\texorpdfstring{Lipschitz rigidity of $\ell_1$ and a Lipschitz-free space problem}{Lipschitz rigidity of and a Lipschitz-free space problem}}
	
	Recall that Problem \ref{prob:Lipschitz-rigidity-l1} asks whether a Banach space Lipschitz isomorphic to $\ell_1$ is necessarily linearly isomorphic to $\ell_1$. The following question, due to G. Godefroy, is closely related.
	
	\begin{prob}
		\label{pb:free-space -l1-comp-bidual}
		Is $\cF(\ell_1)$ complemented in its bidual? 
	\end{prob}
	
	Let us now explain why a positive answer to Problem \ref{pb:free-space -l1-comp-bidual} would imply a positive answer to Problem  \ref{prob:Lipschitz-rigidity-l1}. So, assume that $\cF(\ell_1)$ is complemented in $\cF(\ell_1)^{**}$ and that $X$ is Lipschitz equivalent to $\ell_1$. Then, by Corollary \ref{cor:factorization-Lipschitz-maps} $(ii)$, $\cF(X)$ is isomorphic to $\cF(\ell_1)$ and therefore complemented in its bidual $\cF(X)^{**}$. Since $X$ is separable, by the Godefroy-Kalton lifting theorem (Theorem \ref{thm:GK-lifting-free}), it is linearly isometric to a complemented subspace $Z$ of $\cF(X)$. Hence, there exists a bounded projection from $\cF(X)^{**}$ onto $Z$. Restricting this projection to $Z^{\perp\perp}$ (seeing $Z$ as a subspace of $\cF(X)$ and $Z^{\perp\perp}$ as a subspace of $\cF(X)^{**}$ and as the bidual of $Z$) shows that $Z$ and therefore $X$ is complemented in its bidual. Finally, we apply Theorem \ref{thm:Lipschitz-rigidity-l1} to deduce that $X$ is linearly isomorphic to $\ell_1$.
	
	\section{Notes}
	
	A few years after the seminal papers by Godefroy and Kalton \cite {GodefroyKalton2003} and Kalton \cite{Kalton2004}, the study of the linear structure of Lipschitz-free spaces became a very important subject of research. It developed in many different directions and has provided a wealth of beautiful and deep results. In these notes, we will only mention a few of them that are directly related to the themes of this chapter. 
	
	Let us first address the subject of approximation properties for Lipschitz-free spaces. The harmonic analysis approach presented in Section \ref{sec:free-BAP} has been pushed further by Doucha and Kaufmann in \cite{DouchaKaufmann}, where they prove that if a metric compact group $G$ is equipped with an arbitrary left-invariant metric, then the corresponding Lipschitz-free space $\cF(G)$ has the metric approximation property.  A different simple approach is quite naturally to relate the bounded approximation property for the Lipschitz-free space over a metric space $M$ to the existence of linear extension operators from $\Lip_0(N)$ to $\Lip_0(M)$, for  $N\subset M$. In fact, in \cite{Godefroy2015}, Godefroy provides, for a compact metric space $M$, a nice characterization of the bounded approximation property for $\cF(M)$ in terms of the existence of linear near-extension operators. Linear extension operators are used in \cite{LancienPernecka} to prove that if $M$ is a doubling metric space, then $\cF(M)$ has the $\lambda$-bounded approximation property, where $\lambda$ depends on the doubling constant of $M$. This last result applies to the subsets of $\R^n$. The main open question in this direction is probably the following:
	\begin{prob} 
		Does there exist a universal constant $\lambda \ge 1$ such that for any $n\in \N$ and any $M\subset \R^n$, $\cF(M)$ has the $\lambda$-bounded approximation property? Can $\lambda$ be taken equal to $1$?
	\end{prob}
	Another question is to study, in the setting of Lipschitz-free spaces, improved versions of the bounded approximation property, such as the commuting bounded approximation property, the $\pi$-property, the existence of a Schauder finite-dimensional decomposition, or even of a Schauder basis. In particular, do they pass from a Banach space to its Lipschitz-free space? For all the properties that we have listed, this question is open.  Let us mention that H\'ajek and Perneck\'a proved in \cite{HajekPernecka} that $\cF(\ell_1)$ admits a Schauder basis. At this moment, $\ell_1$ is the only known example of an infinite-dimensional Banach space such that its Lipschitz-free space admits a Schauder basis. We just underline the following problem.
	
	\begin{prob}
		\label{pb:Lip-free-Schauder-basis}
		Let $X$ be a Banach space with a Schauder basis. Does $\cF(X)$ admit a Schauder basis?
	\end{prob}
	
	We now comment on the main recent developments related to Section \ref{sec:structure-Lip-free}. We will concentrate on the article \cite{AGPP2021} by Aliaga, Gartland, Petitjean and Proch\'azka. Let us first give a definition. 
	
	\begin{defi} 
		Let $(M,d)$ be a metric space. We say that $M$ is \emph{1-rectifiable} if
		it equals the union of countably many Lipschitz images of subsets of $\bR$, up to a
		set of null Hausdorff 1-measure (denoted by $\cH^1$). We say that $M$ is \emph{purely 1-unrectifiable} if it contains no 1-rectifiable subset of positive Hausdorff 1-measure; equivalently, if $\mathcal H^1(\gamma(A)) = 0$ for every
		$A \subset \R$ and every Lipschitz map $\gamma\colon A \to M$.
	\end{defi}
	
	Note that the class of purely 1-unrectifiable metric spaces includes the complete countable metric spaces and the complete discrete metric spaces. The ultimate result proved in \cite{AGPP2021} establishes, amongst other things, a characterization of Lipschitz-free spaces containing no isomorphic copy of $L_1$. 
	\begin{theo} 
		Let $M$ be a metric space. Then, the following are equivalent:
		\begin{enumerate}[(i)]
			\item The completion of $M$ is purely 1-unrectifiable,
			\item  $\cF(M)$ has the Radon-Nikod\'{y}m property,
			\item $\cF(M)$ has the Krein-Milman property,
			\item $\cF(M)$ has the Schur property,
			\item $\cF(M)$ contains no isomorphic copy of $L_1$.
		\end{enumerate}
	\end{theo}
	
	A Lipschitz operator version of this result has been very recently obtained in \cite{FJLPPQ}. Let $f$ be a 0-preserving  Lipschitz map between complete pointed metric spaces $M$ and $N$. Recall that $\widehat{f}\colon \cF(M)\to \cF(N)$ is its canonical linearization. Then, the following are equivalent:
	\begin{enumerate}[(i)]
		\item $f$ is curve-flat, that is, for every $K \subset \bR$ compact and every Lipschitz map $\gamma \colon K \to M$,
		\begin{equation*}
			\lim_{\substack{y \to x \\ y\in K}}\frac{d(f(\gamma(x)),f(\gamma(y)))}{|x-y|} = 0 \quad \text{for $\lambda$-almost every $x \in K$, }
		\end{equation*}
		\item  $\widehat{f}$ is a Radon-Nikodým operator,
		\item $\widehat{f}$ is a Dunford-Pettis operator,
		\item $\widehat{f}$ does not fix any copy of $L_1$.
	\end{enumerate}

	\section{Exercises}
	
	\begin{exer}
		\label{ex:origin-change-free-spaces}
		%Show that changing the origin of a metric space just yields the construction of a free-space that is linearly isometric to the initial one.
		Show that Lipschitz-free spaces over a fixed metric space constructed by taking different origins are linearly isometric. 
	\end{exer}
	
	\begin{exer}
		Let $(M,d)$ be a metric space. Show that $\{\delta_M(x)\colon x\in M\}$ is a linearly independent family in $\Lip_0(M)^*$.
	\end{exer}
	
	\begin{exer}
		Let $M$ and $N$ be two pointed metric spaces and $f\colon M \to N$ be a Lipschitz map such that $f(0_M)=0_N$. Recall that there exists a unique bounded linear map $\widehat{f}\colon \cF(M) \to \cF(N)$ such that 
		\begin{equation*}
			\delta_N \circ f = \widehat{f} \circ \delta_M.
		\end{equation*}
		Show that the adjoint operator of $\widehat{f}$ can be identified with the composition operator $g\in \Lip_0(N) \mapsto g \circ f \in \Lip_0(M)$.
		%Prove Remark \ref{rem:composition-operator}.
	\end{exer}
	
	\begin{exer} 
		\label{exe:section-free}
		Let $X$ be a Banach space.
		\begin{enumerate}
			\item Show that the barycenter map $\beta_X$ is a quotient map from $\cF(X)$ onto $X$.
			\item Show that $\cF(X)$ is Lipschitz equivalent to $X \times \mathrm{Ker}(\beta_X)$.
		\end{enumerate}
	\end{exer}
	
	\begin{proof}[Hint]
		For 2. consider the map $\phi:=(\beta_X,Id_{\cF(X)}-\delta_X\circ\beta_X)$.
		%induces a Lipschitz equivalence from $\cF(X)$ onto $X \times \mathrm{Ker}(\beta_X)$, whose inverse is $\phi^{-1}(x,\mu)=\mu+\delta_X(x)$, where $(x,\mu)\in X \times \mathrm{Ker}(\beta_X)$.
	\end{proof}
	
	\begin{exer}
		\label{ex:G-diff-delta} 
		Let $X$ be a nontrivial Banach space. Show that the map $\delta_X \colon X\to \cF(X)$ is nowhere Gateaux differentiable. 
	\end{exer}
	
	\begin{exer}\label{ex:Lipschitzfreeexamples} \ 
		\begin{enumerate}
			\item Show that $\cF(\bN)$ is isometric to $\ell_1$;
			\item Show that $\cF(\bR)$ is isometric to $L_1$. 
		\end{enumerate}
	\end{exer}
	
	The next two exercises discuss alternative proofs of the Mazur-Ulam rigidity theorem.
	
	\begin{exer}
		Prove the Mazur-Ulam rigidity theorem using Figiel's linear left inverse theorem.     
	\end{exer}
	
	The following alternative proof of the Mazur-Ulam rigidity theorem is taken from \cite{Nica2012}.
	
	\begin{exer} 
		Let $X$ and $Y$ be normed vector spaces and assume that $f\colon X\to Y$ is a surjective isometry. Fix $x_1\neq x_2 \in X$ and let
		\begin{equation*}
			\Delta(f) := \norm{f\big(\frac{x_1+x_2}{2}\big)-\frac12\big(f(x_1)+f(x_2)\big)}.
		\end{equation*}
		\begin{enumerate}
			\item Show that $\Delta(f)\le \frac12 \norm{x_1-x_2}$.
			\item Denote by $r$ the reflection about $\frac12(f(x_1)+f(x_2))$ and let $g := f^{-1}\circ r\circ f$. Prove that $\Delta(g)=2\Delta(f)$.
			\item Deduce that $\Delta(f)=0$ and that $f$ is affine.
			\item Deduce the Mazur-Ulam rigidity theorem.
		\end{enumerate}
	\end{exer}
	
	\begin{exer}
		\label{ex:BAP_lemma} 
		Let $X$ be a Banach space and $\lambda\ge 1$ and $A$ be subset of $X$ with dense linear span in $X$. Prove that $X$ has the $\lambda$-BAP if and only if for any $x_1,\dots,x_n \in A$ and any $\eps >0$,  there exits $T \in B(X)$ such that $\norm{T}\le \lambda$, $T$ has  finite rank and $\norm{Tx_i-x_i}\le \eps$ for all $i\le n$.
	\end{exer}
	
	\begin{exer}
		\label{ex:RNP-Schur-l1-sums}\,
		\begin{enumerate}
			\item Show that the Schur property is stable under isomorphic embeddings.
			\item Show that the Schur property is stable under taking $\ell_1$-sums.
		\end{enumerate}
	\end{exer}
	
	\begin{exer} 
		Let $(M,d)$ be a metric space, $(f_n)_n$ a sequence in $\Lip_0(M)$ and $f\in \Lip_0(M)$. Show that the sequence $(f_n)_n$ is converging to $f$ for the weak$^*$ topology induced by $\cF(M)$ on $\Lip_0(M)$ if and only if it is bounded in $\Lip_0(M)$ and pointwise converging to $f$ on $M$.
	\end{exer}
	
	The goal of the next exercise is to show a result due to Godefroy and Ozawa \cite{GodefroyOzawa2014} that asserts the existence of a compact metric space $K$ such that $\cF(K)$ does not have the approximation property. For this purpose, we will use a celebrated result due to Per Enflo \cite{Enflo1973} asserting the existence of a separable Banach space $X$ without the approximation property.
	\begin{exer}
		\label{ex:Godefroy-Ozawa}\, 
		\begin{enumerate}
			\item We adopt the notation of the proof of Theorem \ref{thm:GK-lifting-free}. Show that there exists a compact subset $K$ of $X$ such that the linear map $T$ actually takes values in $\cF(K)$ (seen as a subspace of $\cF(X)$).
			\item Prove that $X$ is linearly isometric to a $1$-complemented subspace of $\cF(K)$.
			\item Deduce that $\cF(K)$ fails to have the approximation property. 
		\end{enumerate}
	\end{exer}
	
	\begin{exer} 
		Assume that $M$ is a countable uniformly discrete metric space. Show that 
		$c_{[\bZ]_\infty^{(\bN)}}(M) \le 2$. 
	\end{exer}
	\begin{proof}[Hint]
		Use Aharoni's embedding theorem and dilations. 
	\end{proof}
	
	\begin{exer}
		Prove Remark \ref{rem:iteratedlimits}.
	\end{exer}

	%%%%%%%%%%%%%%%%%%%%%%%%%%%%%%%%%%%%%%%%%%%%%%%%%%%%%%%%%%%%%%%%%%%%%%%%%%%%%%%%%%%%%%%%%

	\chapter{Asymptotic geometric moduli of Banach spaces}
	\label{chapter:asymptotic-moduli}
	
	There are many natural ways to study an infinite-dimensional Banach space from an asymptotic standpoint. Typically, an asymptotic notion forgets to some extent some local information about the structure of a Banach space. In the mid-1990s, the fundamental influence of certain asymptotic linear properties on the nonlinear geometry of Banach spaces was recognized by various groups of Banach space geometers. These asymptotic properties manifested themselves under various disguises but were all related to certain geometric moduli. The goal of this chapter is to present the definitions of the most important moduli and to study their basic isometric properties. 
	The geometric notions discussed in this chapter are of paramount importance to the topic of this book.

	%%%%%%%%%%%%%%%%%%%%%%%%%%%%%%%%%%%%%%%%%%%%%%%%%
	
	\section{Asymptotic uniform smoothness and convexity}
	\label{sec:auc-aus}
	
	Arguably, the most important asymptotic geometric moduli are the asymptotic uniform smoothness and convexity moduli. These two moduli are two of various moduli that have been first introduced by V.D. Milman in \cite{Milman1971} under a different terminology. The modern terminology was coined in \cite{JLPS2002} where their importance was highlighted. As their name indicates, they are asymptotic versions of the classical uniform convexity and smoothness moduli, and this will be further discussed at the end of this chapter. The definitions of uniform convexity and smoothness and of the corresponding moduli are recalled in Appendix \ref{appendix:us-uc}. Unless otherwise specified, all Banach spaces considered in this chapter will be infinite-dimensional. 
	
	\begin{defi}[AUC and AUS] 
		Let $(X,\norm{\cdot})$ be an infinite-dimensional Banach space. For $t\ge 0$, $x\in S_X$ the unit sphere of $X$ and $Y$ a closed linear subspace of $X$, let 
		\begin{equation*}
			\bar{\rho}_X(t,x,Y) := \sup_{y\in S_Y}\norm{x + ty}-1\ \ \ \ {\rm and}\ \ \
			\ \bar{\delta}_X(t,x,Y) := \inf_{y\in S_Y}\norm{x + ty}-1.
		\end{equation*}
		Let $\cof(X)$ be the set of all closed finite-codimensional subspaces of $X$, we set 
		\begin{equation*}
			\bar{\rho}_X(t,x) := \inf_{Y \in \cof(X)}\bar{\rho}(t,x,Y)\ \ \ \ {\rm and}\ \ \
			\bar{\delta}_X(t,x) := \sup_{Y \in \cof(X)}\bar{\delta}(t,x,Y).
		\end{equation*}
		We define the \emph{modulus of asymptotic smoothness} as
		\begin{equation}
			\bar{\rho}_X(t) :=\sup_{x\in S_X}\ \bar{\rho}_X(t,x) = \sup_{x\in S_X}\inf_{Y \in \cof(X)}\sup_{y\in S_Y}\norm{x + ty}-1,
		\end{equation}
		and the norm of $X$ is said to be \emph{asymptotically uniformly smooth} (in short AUS) if
		\begin{equation}
			\lim_{t \to 0}\frac{\bar{\rho}_X(t)}{t}=0.
		\end{equation}
		Finally, the \emph{modulus of asymptotic uniform convexity} is 
		\begin{equation}
			\bar{\delta}_X(t) := \inf_{x\in S_X}\ \bar{\delta}(t,x) = \inf_{x\in S_X}\sup_{Y \in \cof(X)}\inf_{y\in S_Y}\norm{x + ty}-1,
		\end{equation}
		and the norm of $X$ is said to be {\it asymptotically uniformly convex} (in short
		AUC) if 
		\begin{equation}
			\forall t>0, \quad \bar{\delta}_X(t)>0.
		\end{equation}
	\end{defi}
	
	\begin{rema}
		It is clear that we should only care about the behavior of the moduli when $t$ is close to $0$ since $\bar{\delta}_X(t)>0$ whenever $t>2$.
	\end{rema}
	
	In the following proposition, we record the basic properties of the moduli. We defer the proofs of these elementary properties to the exercises (see Exercises \ref{ex:moduli-basic} and \ref{ex:AUS-subspace-quotient}), as solving them provides an excellent way to familiarize oneself with the moduli.
	\begin{prop}
		\label{AUSsubspacequotient}
		Let $X$ be a Banach space and $Y$ be a closed subspace of $X$. Then, for all $t\in [0,\infty)$,
		\begin{enumerate}[(i)]
			\item $0\le \bar{\delta}_X(t)\le \bar{\rho}_X(t)\le t$,
			\item $\bar{\delta}_X$ is nondecreasing and $1$-Lipschitz,
			\item $\bar{\rho}_X$ is nondecreasing, $1$-Lipschitz and convex,
			\item $\bar{\rho}_Y(t)\le \bar{\rho}_X(t)$,
			\item $\bar{\delta}_Y(t)\ge \bar{\delta}_X(t)$,
			\item $\bar{\rho}_{X/Y}(t)\le \bar{\rho}_X(2t)$.
		\end{enumerate}
	\end{prop}
	
	The convexity of the AUS modulus already provides a hint that AUS behaves better than AUC. It follows from $(iv)$ and $(v)$ that AUS and AUC pass to subspaces and from $(vi)$ that AUS passes to quotients. As we will see, AUC does \emph{not} pass to quotients.
	
	There is significant flexibility in the definition of the moduli. For instance, we can optimize over the set of weak neighborhoods of $0$ instead of the finite-codimensional subspaces. This useful fact is recorded in the following identities for the pointed moduli. The proof is a standard functional analysis exercise and is the content of Exercise \ref{ex:moduli-weak-version}.
	
	\begin{prop}
		\label{prop:moduli-weak-opt}
		Let $(X,\norm{\cdot})$ be a Banach space and denote by $\cV_w(0)$ the set of weak neighborhoods of $0$. Then, for every $t>0$ and $x\in S_X$, 
		\begin{equation}
			\bar{\rho}_X(t,x)=\inf_{V \in \cV_w(0)}\sup_{v\in V\cap S_X} \norm{x+tv}-1,
		\end{equation}
		and
		\begin{equation}
			\bar{\delta}_X(t,x)=\sup_{V \in \cV_w(0)}\inf_{v\in V\cap S_X} \norm{x+tv}-1.    
		\end{equation}
		
	\end{prop}

	The fact that the pointed moduli govern the behavior of weakly null sequences, or more generally, nets, will be crucial in practice. This is recorded in the following proposition, which is an immediate consequence of Proposition \ref{prop:moduli-weak-opt}.
	%Note that by a convexity argument, we can optimize over vectors of norm at least one in the AUS case (and at most one in the AUC case) (cf Exercice \ref{ex}).
	
	\begin{prop}
		\label{prop:moduli-nets}
		Let $(X,\norm{\cdot})$ be a Banach space. Then, for every $t>0$ and $x\in S_X$, 
		\begin{equation}
			\bar{\rho}_X(t,x)=\sup\{\limsup_\alpha \norm{x+tx_\alpha}-1 \colon x_\alpha\wtoo_\alpha 0, \norm{x_\alpha}=1\}
		\end{equation}
		and
		\begin{equation}
			\bar{\delta}_X(t,x)=\inf \{ \liminf_\alpha \norm{x + tx_\alpha}-1 \colon x_\alpha\wtoo_\alpha 0, \norm{x_\alpha}=1\}.
		\end{equation}
		In particular, for every weakly null normalized sequence $(x_n)_n$, one has 
		$$1 + \bar{\delta}_X(t,x)\le \liminf_n\norm{x+t x_n}\le \limsup_n\norm{x+t x_n}\le 1 + \bar{\rho}_X(t,x).$$
	\end{prop}
	
	If the Banach space does not contain $\ell_1$, by calling on Rosenthal's $\ell_1$-theorem \cite{Rosenthal1974}, we can simply use sequences instead of nets. The proof is an immediate consequence of Exercise \ref{ex:moduli-rosenthal}.
	
	\begin{theo}\label{thm:moduli-no-l1}
		Let $(X,\|\cdot\|)$ be a Banach space that does not contain $\ell_1$.
		% and denote by $\cS_w(0)$ the set of weakly null normalized sequences in $X$.
		Then for every $t>0$ and $x\in S_X$, 
		\begin{equation}
			\bar{\rho}_X(t,x)=\sup\{ \limsup_n \norm{x+tx_n}-1\colon x_n\wtoo_n 0, \norm{x_n}=1\},
		\end{equation}
		and
		\begin{equation}
			\bar{\delta}_X(t,x)=\inf \{ \liminf_n \norm{x+tx_n}-1 \colon x_n\wtoo_n 0, \norm{x_n}=1\}.
		\end{equation}
	\end{theo}
	
	Of course, Theorem \ref{thm:moduli-no-l1} includes the case of Banach spaces that are reflexive and those that have a separable dual. In these two situations, the proofs do not rely on Rosenthal's theorem and are arguably simpler (cf Exercise \ref{ex:sequential-moduli}).
	
	\medskip
	With the help of Proposition \ref{prop:moduli-nets}, we can estimate the moduli for the classical $\ell_p$-spaces.
	The following lemma follows from an elementary gliding hump argument (cf Exercise \ref{ex:moduli-lp-sums}).
	\begin{lemm}
		Let $p\in[1,\infty)$ and $(F_n)_n$ be a sequence of finite-dimensional Banach spaces and consider the spaces $D_p := (\sum_{n=1}^\infty F_n)_{\ell_p}$ and $D_\infty := (\sum_{n=1}^\infty F_n)_{c_0}$. Then, for all $x\in S_X$ and every weakly null normalized net $(x_\alpha)_\alpha$ in $D_p$ we have
		$$\lim_\alpha \norm{x+t x_\alpha} = \begin{cases}
			\max\{1,\abs{t}\}, \qquad\textrm{ if } p=\infty\\
			(1+ |t|^p)^{1/p}, \qquad \textrm{ if } p\in[1,\infty)\\
		\end{cases}  
		$$
	\end{lemm}
	
	In particular, we deduce the following estimates whenever $t\in(0,1]$.
	\begin{equation}
		\bar{\rho}_{c_0}(t)=\bar{\delta}_{c_0}(t)=0,  
	\end{equation}
	\begin{equation}
		\bar{\rho}_{\ell_p}(t)=\bar{\delta}_{\ell_p}(t)=(1+t^p)^{\frac1p}-1, \text{ if } p\in[1,\infty).
	\end{equation}
	
	Based on these estimates, one can certainly say that $\co$ is the most asymptotically uniformly smooth space and $\ell_1$ the most asymptotically uniformly convex space.
	
	The special behavior of the modulus of asymptotic uniform smoothness of $\co$ deserves its own terminology.
	
	\begin{defi}\label{defi:AUF} 
		A Banach space $X$ is said to be \emph{asymptotically uniformly flat} (in short AUF) if there exists $t_0>0$ such that for all $t\in (0,t_0)$, $\bar{\rho}_X(t)=0.$
	\end{defi}
	
	\begin{rema} Since $\ell_\infty$ is linearly isometrically universal for separable Banach spaces, $\ell_\infty$ cannot be AUS nor AUC. More precisely (see Exercise \ref{ex:moduli-l-infty}), one can show that for all $t\ge 0$:
		\begin{equation}
			\bar{\rho}_{\ell_\infty}(t)=t \text{ and } \bar{\delta}_{\ell_\infty}(t)= \max\{t-1,0\}.
		\end{equation}
	\end{rema}
	
	\begin{defi}
		For $p\in (1,\infty)$, the modulus of asymptotic uniform smoothness is said to be \emph{of power type $p$} if there exists $C>0$ such that
		$$\forall t\in [0,\infty),\ \bar{\rho}_X(t)\le Ct^p.$$
		For $q\in [1,\infty)$, the modulus of asymptotic uniform convexity is said to be \emph{of power type $q$} if there exists $c>0$ such that
		$$\forall t\in [0,1),\ \bar{\delta}_X(t)\ge ct^q.$$
	\end{defi}
	
	As in the local setting, it is usually more convenient to work with inequalities than with the moduli. 
	
	\begin{defi}
		For $p\in (1,\infty)$, we say that a norm $\norm{\cdot}$ is \emph{$p$-asymptotically uniformly smooth} if there exists $C>0$ such that for every $x\in X$ and every bounded weakly null net $(x_\alpha)_\alpha$ in $X$,
		\begin{equation}
			\limsup_\alpha \norm{x+ x_\alpha}^p\le \norm{x}^p + C\limsup_\alpha\norm{x_\alpha}^p.
		\end{equation}
		Similarly, for $q\in[1,\infty)$, we say that a norm $\norm{\cdot}$ is \emph{$q$-asymptotically uniformly convex}
		if there exists $c>0$ such that for every $x\in X$ and every bounded weakly null net $(x_\alpha)_\alpha$ in $X$,
		\begin{equation}
			\liminf_\alpha \norm{x+ x_\alpha}^q \ge \norm{x}^q + c\liminf_\alpha\norm{x_\alpha}^q.
		\end{equation}
	\end{defi}
	
	The following lemma follows from Proposition \ref{prop:moduli-nets} (see Exercise \ref{ex:moduli-inequalities}).
	
	\begin{lemm}
		\label{lem:AUS-powertype}
		Let $(X, \norm{\cdot})$ be a Banach space, $p\in (1,\infty)$ and $q\in[1,\infty)$. Then, 
		\begin{enumerate}[(i)]
			\item $X$ has a modulus of asymptotic uniform smoothness of power type $p$ if and only if $\norm{\cdot}$ is $p$-asymptotically uniformly smooth.
			\item $X$ has a modulus of asymptotic uniform convexity of power type $q$ if and only if $\norm{\cdot}$ is $q$-asymptotically uniformly convex.
		\end{enumerate}
		In case $(i)$ we shall simply say that $X$ is $p$-AUS and in case $(ii)$ that $X$ is $q$-AUC. 
	\end{lemm}

	\begin{rema}
		The space $(\sum_{n=1}^\infty F_n)_{\ell_p}$ is the prototypical space that has both asymptotic uniformly smooth and convex moduli of power type $p$. A remarkable result from \cite{JLPS2002}, which we will prove in this chapter, states that conversely, if a separable reflexive Banach space has an equivalent norm with modulus of asymptotic uniform smoothness of power type $p$ and an equivalent norm with modulus of asymptotic uniform convexity also of power type $p$, for the same $p\in (1,\infty)$, then it is isomorphic to a subspace of an $\ell_p$-sum of finite-dimensional spaces.
	\end{rema}

	\begin{rema}
		We already saw that asymptotic uniform convexity and asymptotic uniform smoothness pass to subspaces and that asymptotic uniform smoothness passes to quotients, but we left the quotient question temporarily open for asymptotic uniform convexity. Since the moduli estimates above show that $\co$ is not asymptotically uniformly convex, but is (like any separable Banach space) isometric to a quotient of $\ell_1$, asymptotic uniform convexity does \emph{not} pass to quotients.   
	\end{rema}
	
	We conclude this section with the definition of a third geometric modulus, the modulus of asymptotic midpoint uniform convexity, which will only be used in Chapter \ref{chapter:AMP_II}. 
	
	\begin{defi} 
		Let $(X,\norm{\cdot})$ be a Banach space. For $t\ge 0$, the \emph{modulus of asymptotic midpoint uniform convexity} is defined as
		\begin{equation*}
			\hat{\delta}_{\norm{\cdot}}(t) := \inf_{x \in S_X}\ \sup_{Y\in \cof(X)}\ \inf_{y\in Y, \norm{y}\ge 1}\  \frac{\norm{x+ty}+\norm{x-ty}}{2}-1.    
		\end{equation*}
		Banach spaces which satisfy $\hat{\delta}_{\norm{\cdot}}(t)>0$ for all $t\in (0,1)$ are called \emph{asymptotically midpoint uniformly convex}\index{asymptotically midpoint uniformly convex} (in short AMUC).  
	\end{defi}
	
	It is clear from the definitions that $\bar{\delta}_X\le \hat{\delta}_{\norm{\cdot}}$ and therefore if a space is AUC, then it is AMUC. The converse is false, as it has been proved in \cite{DKRRZ2016}. Whether the converse is true up to renorming is an open question.
	
	\begin{prob} 
		\label{pb:AMUC->AUCable}
		Does every asymptotically midpoint uniformly convex Banach space admit an equivalent asymptotically uniformly convex norm?
	\end{prob}

	%%%%%%%%%%%%%%%%%%%%%%%%%%%%%%%%%%%%%%%%%%%%%%%%%
	
	\section{\texorpdfstring{Weak$^*$ AUC and the duality between AUS and AUC}{Weak AUC and the duality between AUS and AUC}}
	The duality between the classical notions of uniform convexity and smoothness and their related moduli is well known and understood. For the asymptotic analogs, the situation is slightly more delicate. The analog duality theory holds for reflexive spaces, but in general, there is only a clean duality between asymptotic uniform smoothness and a weak$^*$ version of asymptotic uniform convexity. This motivates the introduction of the following modulus.
	
	\begin{defi}[weak$^*$ AUC] 
		Let $(X,\norm{\cdot})$ be a Banach space. Denote by $\cof^*(X^*)$ the set of all weak$^*$ closed subspaces of
		$X^*$ of finite codimension. For all $t \ge 0$ and all $x^* \in S_{X^*}$, define
		\begin{equation*}
			\bar{\delta}_X^*(t,x^*) := \sup_{E\in \cof^*(X^*)}\inf_{y^*\in S_E}(\norm{x^* + ty^*}-1)\ \ \text{and}\ \ \bar{\delta}_X^*(t) := \inf_{x^*\in S_{X^*}}\bar{\delta}_X^*(t,x^*).
		\end{equation*}
		The norm of $X^*$ is said to be {\it weak$^*$ asymptotically uniformly convex} (in short AUC$^*$) if
		\begin{equation*}
			\forall t>0,\ \bar{\delta}_X^*(t)>0.
		\end{equation*}
		For $q\in [1,\infty)$, the weak$^*$ asymptotic uniformly convex modulus is said to be \emph{of power type $q$} if there exists $c>0$ such that
		$$\forall t\in [0,1),\ \bar{\delta}_X^*(t)\ge ct^q.$$
	\end{defi}
	
	The basic and elementary properties of the AUC$^*$ modulus are recorded below (see Exercise \ref{ex:AUC*-modulus}).
	
	\begin{prop}
		\label{prop:AUC*-modulus}
		Let $X$ be a Banach space. Then, for all $t\in (0,1]$,
		\begin{enumerate}[(i)]
			\item $0\le \bar{\delta}^*_X(t)\le \bar{\delta}_{X^*}(t)$,
			\item $\bar{\delta}^*_X$ is nondecreasing and $1$-Lipschitz.
		\end{enumerate}
	\end{prop}

	Similarly to the previously defined moduli, we can optimize this time over the set of weak$^*$ neighborhoods of $0$ instead of the weak$^*$ closed finite-codimensional subspaces. This useful fact is recorded in the following identities for the pointed moduli whose proof is deferred to Exercise \ref{ex:moduli-weak*-opt}.
	
	\begin{prop}
		\label{prop:moduli-weak*-opt}
		Let $(X,\norm{\cdot})$ be a Banach space and denote by $\cV_{w^*}(0)$ the set of weak$^*$ neighborhoods of $0$. Then, for every $t>0$ and $x^*\in S_{X^*}$, 
		\begin{equation*}
			\bar{\delta}^*_X(t,x^*)=\sup_{V \in \cV_{w^*}(0)}\inf_{v^*\in V\cap S_{X^*}} \norm{x^* + tv^*}-1.    
		\end{equation*}
	\end{prop}
	
	\begin{rema} 
		Obviously, when $X$ is reflexive, $\cof^*(X^*)=\cof(X^*)$ and $\bar{\delta}^*_X(t) = \bar{\delta}_{X^*}(t)$. The space $\ell_1$ is AUC and AUC$^*$ for the weak$^*$ topology coming from $\co$. However, it can be shown that $\ell_1$ admits isometric preduals such that it is not AUC$^*$ for the corresponding weak$^*$ topology (see Exercise \ref{ex:Szlenk index of C(K)} at the end of the next chapter for this and more). This justifies the notation $\bar{\delta}^*_X$, which underlines the fact that this modulus depends on $X$ rather than on $X^*$. Another example is given by the James tree space $\mathrm{JT}$. Since it is separable with a nonseparable dual, it will follow from results in the next chapter that $JT$ does not admit an equivalent norm with a dual AUC$^*$ norm. However, M. Girardi \cite{Girardi2001} showed that $JT^*$ is AUC.
	\end{rema}
	
	The characterization of the modulus in terms of nets reads as follows.
	
	\begin{prop}
		\label{prop:w*-modulus-nets}
		Let $(X,\norm{\cdot})$ be a Banach space. Then, for every $t>0$ and $x^*\in S_{X^*}$, 
		\begin{equation*}
			\bar{\delta}^*_X(t,x^*)=\inf \{ \liminf_\alpha \norm{x^* + tx^*_\alpha}-1 \colon x^*_\alpha\wstoo_\alpha 0, \norm{x^*_\alpha}=1\}.  
		\end{equation*}
		In particular, for every weak$^*$ null normalized sequence $(x^*_n)_n$ one has 
		\begin{equation*}
			1 + \bar{\delta}^*_X(t,x^*)\le \liminf_n\norm{x^*+t x^*_n}.
		\end{equation*}
	\end{prop}
	
	When $X$ is a separable Banach space, $\bar{\delta}^*_X$ can be described by the behavior of weak$^*$ null sequences in $X^*$ (see Exercise \ref{ex:w*-sequential-moduli}).
	
	\begin{prop}
		\label{prop:w*-modulus-sep}
		Let $(X,\norm{\cdot})$ be a separable Banach space. Then, for every $t>0$ and $x^*\in S_{X^*}$, 
		\begin{equation*}
			\bar{\delta}^*_X(t,x^*)=\inf \{ \liminf_n \norm{x^* + tx^*_n}-1 \colon x^*_n\wstoo_n 0, \norm{x^*_n}=1\}.
		\end{equation*}
	\end{prop}
	
	Our next statement contains all the information on the duality between asymptotic uniform smoothness and weak$^*$ asymptotic uniform convexity. It is a straightforward extension to the nonseparable setting of Proposition 2.6 in \cite{GKL2001}. The proof is taken from \cite{DKLR2017}.
	
	\begin{prop}
		\label{prop:Young} 
		Let $X$ be a Banach space and $0<\sigma,\tau<1$.
		\begin{enumerate}[(a)]
			\item If $\bar{\rho}_X(\sigma)<\sigma\tau$, then  $\bar{\delta}_X^*(6\tau)> \sigma\tau$.
			\smallskip
			\item If $\bar{\delta}_X^*(\tau)> \sigma\tau$, then $\bar{\rho}_X(\sigma)<\sigma\tau$.
		\end{enumerate}
	\end{prop}
	
	The proof relies on the following elementary lemma, which says that every finite-codimensional subspace in a Banach space essentially norms some weak$^*$ closed finite-codimensional subspace in the dual.
	
	\begin{lemm}
		\label{lem:codim} 
		Let $(X,\norm{\cdot})$ be a Banach space, $\eps>0$ and $Y \in \cof(X)$. Then, there exists $Z\in \cof^*(X^*)$ such that 
		\begin{equation}
			\label{eq:lem-codim}
			\forall z^*\in Z,\ \sup_{y\in B_Y}\abs{z^*(y)}\ge \frac{1}{2+\eps}\norm{z^*}.
		\end{equation}
	\end{lemm}
	
	\begin{proof} 
		There are $x^*_1, \dots, x^*_n\in X^*$ such that $Y=\bigcap_{i=1}^n \ker x_i^* $. Denote by $F$ the linear span of $\{x_1^*,\dots,x_n^*\}$ and pick an $\eta$-net $\{u_1^*,\dots,u_k^*\}$ of $S_F$ where $\eta>0$ is such that $(1-2\eta)^{-1}\le 1+\eps$. Now pick $u_i\in S_X$ so that $u_i^*(u_i)\ge 1-\eta$ and consider the (clearly weak$^*$ closed) subspace $Z := \{u_1,\ldots,u_k\}^\perp$ of $X^*$. In order to show \eqref{eq:lem-codim}, observe that if follows from a classical algebraic fact that $Y^\perp=F$ and since $Y^*$ is canonically isometric to $X^*/Y^\perp$, we have that 
		\begin{equation*}
			\forall x^* \in X^*,\ \sup_{y\in B_Y}\abs{x^*(y)} = \norm{x^*}_{X^*/Y^\perp} = d(x^*,F).
		\end{equation*}
		So, fix $z^*\in Z$ and consider $u^*\in F$. Since there exists $1\le i\le k$ such that $\norm{u^* - \norm{u^*}u_i^*}\le \eta\norm{u^*}$, it follows that $\norm{z^*-u^*}\ge \abs{(z^*-u^*)(u_i)}\ge (1-2\eta)\norm{u^*}$. Therefore, by our choice of $\eta$,
		\begin{equation*}
			\norm{z^*}\le \norm{z^*-u^*} + \norm{u^*} \le (1+(1-2\eta)^{-1})\norm{z^*-u^*}\le (2+\eps)\norm{z^*-u^*},
		\end{equation*}
		and \eqref{eq:lem-codim} follows since $\sup_{y\in B_Y}\abs{z^*(y)} = d(z^*,F) = \inf_{u^*\in F}\norm{z^*-u^*}$.
	\end{proof}
	
	\begin{proof}[Proof of Proposition \ref{prop:Young}] 
		
		$(a)$ Assume that $\bar{\rho}_X(\sigma)<\sigma\tau$. Let $x^*\in S_{X^*}$, $\eps>0$ and choose $x\in S_X$ such that $x^*(x)\ge 1-\eps$. By our assumption on $\bar{\rho}_X$, there exists a finite-codimensional subspace $Y$ of $X$ such that
		\begin{equation}
			\label{eq:23}
			\forall y\in S_Y\ \ \norm{x+\sigma y}<1+\sigma\tau.
		\end{equation}
		By considering $Y\cap \ker(x^*)$ if needed, we may also assume that $Y\subset \ker(x^*)$. Then, it follows from Lemma \ref{lem:codim} that there exists $Z \in \cof^*(X^*)$ such that
		\begin{equation}
			\label{eq:24}
			\forall z^*\in Z\ \ \sup_{y\in B_Y}\abs{z^*(y)}\ge \frac{1}{2+\eps}\norm{z^*}
		\end{equation}
		and of course, we can also assume that $z^*(x)=0$, for $z^*\in Z^*$. Therefore, for all $y\in S_Y$ and all $z^*\in S_Z$ 
		\begin{equation}
			\label{eq:25}
			\norm{x^*+6\tau z^*}\stackrel{\eqref{eq:23}}{\ge} \langle x^*+6\tau z^* , \frac{x+\sigma y}{1+\sigma\tau}\rangle\ge \frac{1}{1+\sigma\tau}\big(1-\vep+6\tau\sigma z^*(y)\big),
		\end{equation}
		and after taking the supremum over $y\in S_Y$ in \eqref{eq:25} and invoking \eqref{eq:24}, we get that
		\begin{equation*}
			\forall z^*\in S_Z\ \ \norm{x^*+6\tau z^*} \ge \frac{1}{1+\sigma\tau}\big(1-\eps+\frac{6\sigma\tau}{2+\eps}\big).
		\end{equation*}
		Finally, letting $\eps$ tend to $0$, we conclude that
		\begin{equation*}
			\bar{\delta}_X^*(6\tau)\ge \frac{1+3\sigma\tau}{1+\sigma\tau}-1 > \sigma\tau.
		\end{equation*}
		
		$(b)$ Assume that $\bar{\delta}_X^*(\tau)> \sigma\tau$. We shall prove by contradiction that $\bar{\rho}_X(\sigma)\le \sigma\tau$ (the strict inequality will follow formally). So, assume that for some $x\in S_X$, $\bar{\rho}_X(\sigma,x)>\sigma\tau$. We can certainly pick $\rho$ such that $\sigma\tau <\rho < \bar{\rho}_X(\sigma,x)$ and $\bar{\delta}_X^*(\tau)> \rho$ and thus 
		\begin{equation*}
			\forall Y\in \cof(X)\ \ \exists x_Y\in S_Y\ \ \norm{x+\sigma x_Y}>1+\rho.
		\end{equation*}
		If we preorder $\cof(X)$ by reverse inclusion, i.e., for $Y,Z\in \cof(X)$, we declare that $Y\preceq Z$ if and only if $Z\subseteq Y$, then $(\cof(X),\preceq)$ becomes a directed set and $(x_Y)_{Y\in \cof(X)}$ a net in $S_{X}$ that is eventually in any finite-codimensional subspace and, in particular, weakly null. Now, for any $Y\in \cof(X)$, we can pick a norming functional $y^*_Y\in S_{X^*}$ such that 
		\begin{equation}
			\label{eq:26}
			\langle y^*_Y, x+\sigma x_Y\rangle>1+\rho.
		\end{equation} 
		By weak$^*$ compactness of $B_{X^*}$, we can find a subnet $(y^*_\alpha)_{\alpha \in A}$ of the net $(y^*_Y)_{Y\in \cof(X)}$, as well as $y^*\in B_{X^*}$ and $c\ge0$ such that $(y^*_\alpha)_{\alpha \in A}$ is weak$^*$ converging to $y^*$ and $(\norm{y^*_\alpha-y^*})_{\alpha \in A}$ converges to $c$. There are two cases to consider:
		
		\underline{Case $c<\tau$:}
		
		\noindent Fix $\eps>0$. Since $A$ is a directed set that is cofinal in $\cof(X)$, there exists $\alpha_0\in A$ such that for all $\alpha\ge \alpha_0$,
		\begin{equation*}
			x_\alpha\in \ker(y^*),\ \abs{\langle y^*_\alpha-y^*, x\rangle}<\eps,\ {\rm and}\ \norm{y^*_\alpha-y^*}<c+\eps.
		\end{equation*}
		Therefore, setting $x^*_\alpha :=y^*_\alpha-y^*$, for any $\alpha \in A$, it follows that for all $\alpha\ge \alpha_0$
		\begin{equation*}
			\langle y^*+x^*_\alpha , x+\sigma x_\alpha \rangle = \langle y^* , x \rangle + \langle x^*_\alpha , x \rangle + \sigma 
			\langle x^*_\alpha , x_\alpha \rangle\le 1+\eps+\sigma(c+\eps).
		\end{equation*}
		Since $c<\tau$ and $\sigma\tau <\rho$, if $\eps$ was initially chosen small enough we would have for all $\alpha\ge \alpha_0$
		\begin{equation*}
			\langle y^*+x^*_\alpha , x+\sigma x_\alpha \rangle = \langle y^*_\alpha , x+\sigma x_\alpha \rangle <1+\rho,
		\end{equation*}
		contradicting \eqref{eq:26}.
		
		\underline{Case $c\ge \tau$:}
		
		\noindent Fix $\eps>0$ to be chosen small enough later. If we can show that $\norm{y^*}\le 1-\sigma c$, then picking $\alpha_0\in A$ as already done above so that for all $\alpha\in A$, $\alpha\ge \alpha_0$, we have $x_\alpha\in \ker(y^*)$, $\abs{\langle y^*_\alpha - y^* , x\rangle}<\eps$ and $\norm{y^*_\alpha - y^*}\le c+\eps$. Then, for all $\alpha\ge \alpha_0$,
		\begin{equation*}
			\langle y^*+x^*_\alpha , x+\sigma x_\alpha \rangle = \langle y^* , x \rangle + \langle x^*_\alpha , x \rangle + \sigma \langle x^*_\alpha , x_\alpha \rangle \le 1-\sigma c +\eps+\sigma (c+\eps) = 1 +\vep(1+\sigma),
		\end{equation*}
		which is again a contradiction with \eqref{eq:26} for $\vep <\rho/(1+\sigma)$. 
		
		It remains to prove that $\norm{y^*}\le 1-\sigma c$ and we may assume that $y^*\neq 0$. Recall that $\bar{\delta}_X^*(\tau)> \rho$ and hence there exists  $Z \in \cof^*(X^*)$ such that
		\begin{equation*}
			\forall z^*\in S_Z\ \  \Big\|y^* + \tau\norm{y^*}\,z^*\Big\|\ge (1+\rho)\norm{y^*}.
		\end{equation*}
		Since $(x^*_\alpha)_{\alpha \in A}$ is weak$^*$ converging to $0$ and $(\norm{x^*_\alpha})_{\alpha \in A}$ converges to $c$, we have that $d(c^{-1}x^*_\alpha,S_Z)$ tends to $0$. Therefore, we deduce
		\begin{equation*}
			\exists \alpha_1\in A\ \ \forall \alpha\ge \alpha_1\ \ \Big\|y^*+\tau c^{-1}\norm{y^*}\,x_\alpha^*\Big\|\ge (1+\rho)\norm{y^*}-\eps.
		\end{equation*}
		Note that $\lambda=\tau c^{-1}\norm{y^*}\in [0,1]$ and thus we can write
		\begin{equation*}
			y^*+\lambda x^*_\alpha=\lambda(y^*+x^*_\alpha)+(1-\lambda)y^*.
		\end{equation*}
		Using the convexity of the norm, we deduce that
		\begin{equation*}
			\norm{y^*+\lambda x^*_\alpha}\le \lambda \norm{y^*+x^*_\alpha} + (1-\lambda)\norm{y^*} = \lambda+(1-\lambda)\norm{y^*}.
		\end{equation*}
		Therefore,
		\begin{equation*}
			(1+\rho)\norm{y^*} - \eps \le \tau c^{-1}\norm{y^*} + (1-\tau c^{-1}\norm{y^*})\norm{y^*}.
		\end{equation*}
		Letting $\rho$ tend to $\sigma\tau$ and $\eps$ tend to $0$, we get that
		\begin{equation*}
			(1+\sigma\tau) \norm{y^*}\le \tau c^{-1}\norm{y^*}+(1-\tau c^{-1}\norm{y^*})\norm{y^*}.
		\end{equation*}
		Dividing by $\norm{y^*}$, we obtain $\sigma\tau \le \tau c^{-1}(1-\norm{y^*})$ and as claimed $\norm{y^*}\le 1-\sigma c$.
		
	\end{proof}
	
	\begin{rema}\label{rem:dualityAUF} 
		Note that it follows directly from Proposition \ref{prop:Young} that a Banach space $X$ is AUF if and only if there exists $c>0$ such that $\bar{\delta}^*_{X}(t)\ge ct$ for all $t\ge 0$ (in other words if and only if the norm of $X^*$ is $1$-AUC$^*$).
	\end{rema}
	
	Proposition \ref{prop:Young} can be rephrased in terms of Young's duality between $\bar{\rho}_X$ and $\bar{\delta}^*_{X}$. We recall some basic definitions.
	\begin{defi}[Dual Young function]
		Let $f$ be a continuous and nondecreasing function on $[0,1]$ such that $f(0)=0$. Its \emph{dual Young function} is the function $f^*\colon [0,1]\to [0,\infty)$ defined by
		\begin{equation}
			f^*(t) := \sup\{st-f(s):\ 0\le s\le 1\}.
		\end{equation}
	\end{defi}
	Let $C\ge 1$ and $f,g$ be continuous and nondecreasing functions on $[0,1]$ satisfying $f(0)=g(0)=0.$
	We will say that $f,g$ are \emph{$C$-equivalent} if $f(t)\ge g(t/C)$ and $g(t)\ge f(t/C)$ for all $t\in [0,1]$. The duality between AUC$^*$ and AUS takes the following form in terms of dual Young functions.

	\begin{coro}
		\label{cor:Young2} 
		Let $X$ be a Banach space. Then, for all $s\in [0,1]$
		\begin{equation}
			\label{eq:Young1/2}
			(\bar{\delta}^*_{X})^*(s)\ge \bar{\rho}_X\big(\frac{s}{2}\big),
		\end{equation}
		and
		\begin{equation}
			\label{eq:Young1/6}
			(\bar{\delta}^*_{X})^*\big(\frac{s}{6}\big)\le \bar{\rho}_X(s).
		\end{equation}
		In particular, $\bar{\rho}_X$ and $(\bar{\delta}^*_{X})^*$ are $6$-equivalent.
	\end{coro}
	
	\begin{proof} 
		Consider first $t=\frac2s\bar{\rho}_X(\frac{s}{2})\in [0,1]$. Then, $\bar{\rho}_X(\frac{s}{2})=\frac{s}{2}t$. So, it follows from Proposition \ref{prop:Young} (b)  that $\bar{\delta}^*_{X}(t)\le \frac{s}{2}t$. Therefore, $(\bar{\delta}^*_{X})^*(s)\ge ts-\bar{\delta}^*_{X}(t)\ge \frac{s}{2}t=\bar{\rho}_X(\frac{s}{2}).$
		
		Assume now that $(\bar{\delta}^*_{X})^*(\frac{s}{6}) > \bar{\rho}_X(s).$ Then, there exists $t\in [0,1]$ such that $\frac{s}{6}t-\bar{\delta}^*_{X}(t)>\bar{\rho}_X(s)$. Thus, $\bar{\delta}^*_{X}(t)<\frac{s}{6}t-\bar{\rho}_X(s)\le\frac{s}{6}t$. It now follows from Proposition \ref{prop:Young}~(a) that $\bar{\rho}_X(s)\ge \frac{s}{6}t$. But this implies that $\bar{\delta}^*_{X}(t)<0$, which is impossible.
	\end{proof}
	
	The following result is now an immediate consequence of the discussion above (with the convention that $\infty$-AUS means AUF).
	
	\begin{coro}
		\label{cor:duality}
		Let $(X, \norm{\cdot}_X)$ be a Banach space and $p,q \in [1,\infty]$ be two conjugate exponents. Then,
		\begin{enumerate}[(1)]
			\item $\norm{\cdot}_X$ is AUS if and only if $\norm{\cdot}_{X^*}$ is AUC$^*$ and
			\item $\norm{\cdot}_X$ is $p$-AUS if and only if $\norm{\cdot}_{X^*}$ is $q$-AUC$^*$.
			\item[] Moreover, if $X$ is reflexive, then 
			\item $\norm{\cdot}_X$ is AUS if and only if $\norm{\cdot}_{X^*}$ is AUC and
			\item $\norm{\cdot}_X$ is $p$-AUS if and only if $\norm{\cdot}_{X^*}$ is $q$-AUC.
		\end{enumerate}
	\end{coro}

	\section{Orlicz functions and Kalton's iterated norms}
	
	We have already observed that the modulus of asymptotic uniform smoothness is a function from $[0,\infty)$ to $[0,\infty)$ that is nondecreasing, continuous and convex. Moreover, $\bar{\rho}_X(0)=0$ and $\lim_{t\to\infty} \bar{\rho}_X(t)=\infty$. Functions with these properties are known as Orlicz functions.
	
	\begin{defi}[Orlicz functions] 
		A map $F \colon [0,\infty)\to [0,\infty)$ is said to be an \emph{Orlicz function} if it is nondecreasing, convex, continuous and satisfies $F(0)=0$ and $\lim_{t\to \infty} F(t)= \infty$. If moreover, $F(t)>0$ for all $t>0$, then $F$ is said to be \emph{nondegenerate}. 
	\end{defi}
	
	A typical and important example of an Orlicz function is the function $t \mapsto t^p$, for $p\in [1,\infty)$ and we can define a class of sequence spaces that generalizes the classical $\ell_p$-spaces.
	
	\begin{defi}[Orlicz sequence spaces]
		Let $F$ be an Orlicz function and let
		\begin{equation}
			\ell_{F} :=\Big\{x=(x_n)_{n=1}^\infty \in \bR^\bN\colon \exists r>0\ \
			\sum_{n=1}^\infty F\Big(\frac{\abs{x_n}}{r}\Big)<\infty\Big\}.
		\end{equation}
		When equipped with the \emph{Luxemburg norm}
		\begin{equation}
			\norm{x}_{\ell_{F}} := \inf\Big\{r>0,\ \sum_{n=1}^\infty F\Big(\frac{\abs{x_n}}{r}\Big)\le 1\Big\},
		\end{equation}
		the vector space $\ell_F$ becomes a Banach space, called the \emph{Orlicz sequence space} associated with $F$.
	\end{defi}
	
	\begin{rema}
		\begin{enumerate}
			\item Without additional assumption on the Orlicz function, an Orlicz space is in general nonseparable. We refer to \cite{LindenstraussTzafriri1977} for a detailed study of Orlicz spaces.
			\item If $F_p(t) := t^p$, for $p\in [1,\infty)$, then $\ell_{F_p}$ is a separable Banach space that coincides with $\ell_p$ and $\norm{\cdot}_{\ell_{F_p}}$ is equal to $\norm{\cdot}_p$.
		\end{enumerate}
	\end{rema} 
	
	If an Orlicz function is merely dominated by $F_p$, or dominates $F_p$ near the origin, the corresponding norms are still somewhat comparable in the following sense.  
	
	\begin{prop}\label{prop:Orlicz-pnorm} Let $F$ be an Orlicz function and $p\in [1,\infty)$.
		\begin{enumerate}[(i)]
			\item Assume now that there exists $t_0>0$ such that $F(t)=0$ for $t \in [0,t_0]$ (in other words $F$ is degenerate). Then, $\norm{\cdot}_{\ell_F}$ is equivalent to $\norm{\cdot}_\infty$ on $c_{00}$.
			\item Assume that there exists a constant $c>0$ such that $F(t)\le ct^p$ for all $t\ge 0$. Then, for all $x\in c_{00}$, $\norm{x}_{\ell_F}\le c^{1/p}\norm{x}_p$.
			\item Assume that there exist $b,c>0$ such that $F(t)\ge ct^p$ for all $t\in [0,b]$. Then, there exists $\gamma>0$ such that for all $x\in c_{00}$, $\norm{x}_{\ell_F}\ge \gamma\norm{x}_p$.
		\end{enumerate} 
	\end{prop} 
	
	\begin{proof} 
		For $(i)$ and $(ii)$, the arguments are straightforward. For $(iii)$, the conclusion follows easily from the observation that there must exist a constant $\gamma>0$ such that $F(t)\ge \gamma^p t^p$ whenever $F(t)\le 1$.
	\end{proof}
	
	\begin{rema}
		More generally, if for two Orlicz functions $F_1$ and $F_2$, we have $F_1 = F_2$ on a neighborhood of the origin, then $\ell_{F_1}$ and $\ell_{F_2}$ coincide as sets and the corresponding norms are equivalent.
	\end{rema}
	
	There is a close connection between Orlicz functions and absolute norms that will become very handy when studying the behavior of weakly null sequences in Banach spaces. 
	
	\begin{defi}[Absolute norms] 
		We say that a norm $N$ on $\bR^2$ is an \emph{absolute norm} if $N(a,b)=N(\abs{a},\abs{b})$ for all $(a,b)\in \bR^2$.
		An absolute norm is \emph{normalized} if $N(0,1)=N(1,0)=1$.
	\end{defi}
	
	Some of the basic properties of absolute norms are recorded in the following lemma, whose proof constitutes the goal of Exercise 
	\ref{ex:absolute-norm}.
	
	\begin{lemm}
		\label{lem:absolute-norm}
		Let $N$ be a norm on $\bR^2$. 
		\begin{enumerate}
			\item If $N$ is absolute and normalized, then $\norm{(a,b)}_\infty \le N(a,b)\le \norm{(a,b)}_1$ for all $a,b\in \bR$.
			\item $N$ is absolute if and only if $N(a_1,b_1)\le N(a_2,b_2)$ for all $a_1,a_2,b_1,b_2\in \bR$ satisfying $\abs{a_1}\le \abs{a_2}$ and $\abs{b_1}\le \abs{b_2}$.
			% see \cite{BauerStoerWitzgall66} or use Krein-Milman?
			%    \item $N(a,b)=N(b,a)$.
		\end{enumerate}
	\end{lemm}

	Given an absolute norm on $\bR^2$, one can associate with it a sequence space via an inductive procedure. These iterated norms have been instrumental in providing deep insights in the asymptotic geometry of Banach spaces. It seems that their first use is due to Kalton in \cite{Kalton1993}.
	
	\begin{defi}[Kalton's iterated norms] 
		Let $N$ be a normalized absolute norm on $\bR^2$ and denote by $(e_n)_{n=1}^\infty$ the canonical basis of $c_{00}$. We define inductively a norm on $c_{00}$, denoted by $\norm{\cdot}_{\Lambda_N}$, as follows: 
		\begin{itemize}
			\item $\norm{a_1 e_1}_{\Lambda_N} := N(0,\abs{a_1})$ for all $a_1\in \bR$,
			\item $\big\|\sum_{j=1}^n a_je_j\big\|_{\Lambda_N} := N\big(\big\|\sum_{j=1}^{n-1} a_je_j\big\|_{\Lambda_N}, \abs{a_n}\big)$ for all $(a_1,\dots,a_n)\in \bR^n,$ if $n\ge 2$.
		\end{itemize}
		The sequence space $\Lambda_N$ is defined as the completion of $c_{00}$ under the norm $\norm{\cdot}_{\Lambda_N}$.
	\end{defi}
	
	Iterated norms enjoy some very useful properties that can be easily derived from the absolute property and/or the inductive definition.
	The proof of the following lemma is the content of Exercise \ref{ex:iterated-unc-shift}.
	
	\begin{lemm}
		\label{lem:iterated-unc-shift} 
		Let $N$ be a normalized absolute norm.
		\begin{enumerate}
			\item The canonical basis of $\Lambda_N$ is normalized and $1$-unconditional.
			\item For any $n\in \bN$ and any $(a_1,\dots,a_n)\in \bR^n$, we have:
			\begin{equation}
				\big\|\sum_{i=1}^na_ie_{i+1}\big\|_{\Lambda_N} = \big\|\sum_{i=1}^na_ie_{i}\big\|_{\Lambda_N}.
			\end{equation}
		\end{enumerate}
	\end{lemm}
	
	The next lemma describes how normalized absolute norms generate quite regular Orlicz functions and vice versa.
	
	\begin{lemm}
		\label{lem:absolute-Orlicz}\,
		\begin{enumerate}
			\item Let $N$ be a normalized absolute norm on $\bR^2$. The map $F_N$ defined by $F_N(t) := N(1,t)-1$ is a $1$-Lipschitz Orlicz function so that $\lim_{t\to \infty}\frac{F_N(t)}{t}=1$.
			\item Let $F$ be an Orlicz function that is Lipschitz and satisfies $\lim_{t\to \infty}\frac{F(t)}{t}=1$ and $F(t)\ge t-1$ for all $t\ge 0$. Then, the map $N_F$ defined by $$N_F(a,b) := \begin{cases}
				\abs{b} \text{ if } a=0,\\
				\abs{a}\big(1+F(\frac{\abs{b}}{\abs{a}})\big) \text{ if } a\neq 0,
			\end{cases}$$
			is a normalized absolute norm on $\bR^2$.
			\item $F_{N_F}=F$, $N_{F_N}=N$ and in particular $F(\cdot) = N_F(1, \cdot) - 1$.
		\end{enumerate} 
	\end{lemm}
	
	The proof of Lemma \ref{lem:absolute-Orlicz} is a routine argument that is postponed to Exercise \ref{ex:absolute-Orlicz}.
	
	Our next statement shows that in the above situation, the spaces $\ell_F$ and $\Lambda_{N_F}$ coincide. More precisely,  we have.
	
	\begin{prop}
		\label{pro:Orlicz-iterated}
		Let $N$ be a normalized absolute norm on $\bR^2$ and set $F_N(t) := N(1,t)-1$. Then, for all $x\in c_{00}$, \begin{equation*}
			\frac12\norm{x}_{\ell_{F_N}}\le \norm{x}_{\Lambda_N} \le e\norm{x}_{\ell_{F_N}}.
		\end{equation*}
	\end{prop}
	
	\begin{proof}
		We start with the leftmost inequality. So, assume that $x=\sum_{i=1}^n a_i e_i$ and without loss of generality $\big\|\sum_{i=1}^na_ie_{i}\big\|_{\Lambda_N}=1$. Then, it follows from the previous lemma that  
		\begin{equation*}
			2=1+\big\|\sum_{i=1}^na_ie_{i}\big\|_{\Lambda_N} = \norm{e_1}_{\Lambda_N} + \big\|\sum_{i=1}^na_ie_{i+1}\big\|_{\Lambda_N} \ge  \big\|e_1+\sum_{i=1}^na_ie_{i+1}\big\|_{\Lambda_N}
		\end{equation*}
		Using the definition of $\norm{\cdot}_{\Lambda_N}$ and the unconditionality of $(e_i)_i$ with respect to $\norm{\cdot}_{\Lambda_N}$, we deduce that 
		\begin{align*}
			\big\| e_1 + \sum_{i=1}^n a_i e_{i+1} \big \|_{\Lambda_N} & = N(\| e_1 + \sum_{i=1}^{n-1} a_i e_{i+1} \|_{\Lambda_N}, \abs{a_n}) \\
			& = \| e_1 + \sum_{i=1}^{n-1} a_i e_{i+1} \|_{\Lambda_N} \Big(1+ F_N \Big( \frac{\abs{a_n} }{ \|e_1 + \sum_{i=1}^{n-1} a_i e_{i+1}\|_{\Lambda_N}} \Big) \Big)\\
			& \ge \big\|e_1+\sum_{i=1}^{n-1}a_ie_{i+1}\big\|_{\Lambda_N}\Big(1 + F_N \Big(\frac{\abs{a_n}}{2}\Big)\Big).
		\end{align*}
		Iterating and using the elementary fact that $\Pi_{i=1}^n (1 + z_i)\ge 1 +\sum_{i=1}^n z_i$ we get 
		\begin{equation*}
			2 \ge \big\| e_1 + \sum_{i=1}^n a_i e_{i+1} \big \|_{\Lambda_N} \ge \prod_{i=1}^n \Big(1+F_N\big(\frac{\abs{a_i}}{2}\big)\Big)\ge 1+\sum_{i=1}^n F_N\big(\frac{\abs{a_i}}{2}\big),
		\end{equation*}
		from which it readily follows that  $\big\|\sum_{i=1}^na_ie_{i}\big\|_{\ell_{F_N}}\le 2$. 
		
		We now turn to the rightmost inequality and assume that $\big\|\sum_{i=1}^na_ie_{i}\big\|_{\ell_{F_N}}=1$. 
		We wish to show that $\big\|\sum_{i=1}^na_ie_{i}\big\|_{\Lambda_N}\le e$. So, let us assume, as we may, that $\big\|\sum_{i=1}^na_ie_{i}\big\|_{\Lambda_N} >1$. Denote by $r$ the smallest integer in $\{1,\dots ,n\}$ such that $\big\|\sum_{i=1}^ra_ie_{i}\big\|_{\Lambda_N} >1$. Then, for all $j>r$ (if any):
		\begin{equation*}
			\big\|\sum_{i=1}^ja_ie_{i}\big\|_{\Lambda_N} \le \big\|\sum_{i=1}^{j-1}a_ie_{i}\big\|_{\Lambda_N}(1 + F_N(|a_j|)),
		\end{equation*}
		and iterating, we get
		\begin{equation*}
			\big\|\sum_{i=1}^ja_ie_{i}\big\|_{\Lambda_N} \le \big\|\sum_{i=1}^{r}a_ie_{i}\big\|_{\Lambda_N}\prod_{k= r+1}^j(1 + F_N(|a_k|)).
		\end{equation*}
		If $r>1$, then $\big\|\sum_{i=1}^{r-1} a_ie_{i}\big\|_{\Lambda_N}\le 1$ and by absolute monotonicity of $N$ we have
		\begin{equation*}
			\big\|\sum_{i=1}^ra_ie_{i}\big\|_{\Lambda_N} \le N(1, \abs{a_r}) = 1 + F_N(\abs{a_r}).
		\end{equation*} 
		Similarly, if $r=1$, then  
		\begin{equation*}
			\|a_1e_1\|_{\Lambda_N} = N(0,|a_1|)\le N(1,|a_1|) = 1 + F_N(|a_1|).    
		\end{equation*}
		Since $F_N\ge0$, using the elementary inequality $\prod_{i=1}^n(1+z_i)\le e^{\sum_{i=1}^n z_i}$, in both cases we have
		\begin{equation*}
			\big\|\sum_{i=1}^na_ie_{i}\big\|_{\Lambda_N}\le \prod_{i=1}^n (1+F_N(|a_i|))\le e^{\sum_{i=1}^n F_N(|a_i|)}\le e,
		\end{equation*}
		where the last inequality holds since we assumed that $\big\|\sum_{i=1}^na_ie_{i}\big\|_{\ell_{F_N}}=1$.
	\end{proof}
	
	We can immediately deduce from Proposition \ref{prop:Orlicz-pnorm} the following statement that will be of special interest.
	\begin{coro}
		\label{cor:iterated_pnorm} 
		Let $N$ be a normalized absolute norm on $\bR^2$ and set $F_N(t) := N(1,t)-1$. Let $p\in [1,\infty)$.
		\begin{enumerate}[(i)]
			\item Assume that there exists a constant $c>0$ such that $F_N(t)\le ct^p$ for all $t\ge 0$. Then, there exists $\gamma >0$ such that for all $x\in c_{00}$, $\norm{x}_{\Lambda_N}\le \gamma\norm{x}_p$.
			\item Assume that there exist $b,c>0$ such that $F_N(t)\ge ct^p$ for all $t\in [0,b]$. Then, there exists $\gamma>0$ such that for all $x\in c_{00}$, $\norm{x}_{\Lambda_N}\ge \gamma \norm{x}_p$.
		\end{enumerate} 
	\end{coro} 
	
	The fact that the Orlicz space associated with the modulus of asymptotic uniform smoothness, or some convexified version of the modulus of asymptotic uniform convexity, is isomorphic to an iterated norm sequence will be extremely useful when dealing with weakly null trees instead of merely weakly null sequences or nets. As an illustration, we describe how Kalton's iterated norms can be used to construct finite-dimensional Schauder decompositions (FDD) satisfying nontrivial upper estimates. We refer to the classical textbook \cite{LindenstraussTzafriri1977} for all the necessary background on this subject, which will also be required for other arguments later.
	
	A Banach space $X$ admits a FDD if there is a sequence of finite-dimensional and nonzero subspaces $(E_n)_{n=1}^\infty$ such that every $x\in X$ can be decomposed as $x=\sum_{n=1}^\infty x_n$ for some unique sequence $(x_n)_{n=1}^\infty \in \Pi_{n=1}^\infty E_n$. If $I$ is an interval of $\bN$ and $x$ is decomposed as above, we define $P_I(x) := \sum_{n\in I}x_n$ and for $n\in \bN$, $S_n := P_{[1,n]}$. Recall that the FDD (or projection) constant of $(E_n)_{n=1}^\infty$ is the finite number $\sup_n\norm{S_n}$ and that the space can always be renormed such that $\sup_{n\le k}\|P_{[n,k]}\|\le 1$ in which case we say that the FDD is bimonotone. If $\sup_n\norm{S_n}\le K$ we say that $(E_n)_{n=1}^\infty$ is a $K$-FDD.
	
	Let $p\in [1,\infty]$. We say that a FDD $(E_n)_{n=1}^\infty$ satisfies an \emph{upper-$p$-estimate with constant $C\in [1,\infty)$}, or a $(C,p)$-upper estimate, if 
	\begin{equation}
		\forall n\in \bN\ \ \forall (x_1,\ldots,x_n)\in E_1\times \cdots \times E_n,\ \ \Big\|\sum_{i=1}^n x_i\Big\|\le C\big(\sum_{i=1}^n\|x_i\|^p\big)^{\frac{1}{p}},
	\end{equation}
	and a \emph{lower-$p$-estimate with constant $c>0$}, or a $(c,p)$-lower estimate, if 
	\begin{equation}
		\forall n\in \bN\ \ \forall (x_1,\ldots,x_n)\in E_1\times \cdots \times E_n,\ \ \Big\|\sum_{i=1}^n x_i\Big\|\ge c\big(\sum_{i=1}^n\|x_i\|^p\big)^{\frac{1}{p}},   
	\end{equation}
	with the classical convention when $p=\infty$.
	%The FDD is said to be \emph{shrinking} if any bounded block sequence with respect to $(E_n)_{n=1}^\infty$ is weakly null.  
	
	\begin{rema}
		Of course every FDD satisfies an upper-$1$-estimate with constant $C=1$ and a lower-$\infty$-estimate with constant $c= \frac{1}{2K}$ where $K := \sup_n\norm{P_n}$ is the FDD constant.     
	\end{rema}
	
	Given a sequence of vector subspaces $(E_n)_{n=1}^\infty$ of a Banach space $X$, we will say that $(F_k)_{k=1}^\infty$ is a \emph{blocking} of $(E_n)_{n=1}^\infty$ if there exist $1=m_1<\dots<m_i<\dots$ so that $F_k=\sum_{m_{k}\le n < m_{k+1}}E_n$, for $k\ge 1$ and a \emph{skipped blocking} if there exist $1\le m_1<\dots<m_i<\dots$ so that $F_k=\sum_{m_{2k-1}\le n < m_{2k}}E_n$, for $k\ge 1$. 
	
	A system $(x_n,x_k^*)_{(n,k)\in \bN^2}$ in $X\times X^*$ is called a \emph{biorthogonal system} if $x_n^*(x_k)=\delta_{n,k}$ for all $n,k\in \bN$. A biorthogonal system $(x_n,x_k^*)_{(n,k)\in \bN^2}$ in $X\times X^*$ is called \emph{fundamental} if the linear span of $(x_n)_{n=1}^\infty$ is norm dense in $X$ and \emph{total} if the linear span of $(x_n^*)_{n=1}^\infty$ is weak$^*$ dense in $X^*$. A fundamental and total biorthogonal system is called a \emph{Markukushevich basis} or \emph{M-basis}. If in addition the linear span of $(x_n^*)_{n=1}^\infty$ is norm dense in $X^*$, the biorthogonal system is called \emph{shrinking}. Every separable Banach space admits a Markushevich basis and every Banach space with separable dual admits a shrinking biorthogonal system. We refer the reader to the textbook \cite{HMVZ2008} for all references and proofs on the subject. 
	
	Note that a Banach space with a shrinking FDD admits a shrinking biorthogonal system.
	
	The following proposition is a fundamental and rather general, blocking argument whose proof shows the versatility of Kalton's iterated norms. This argument could be found for $p$-AUS Banach spaces in \cite{JLPS2002}, but this general version has been written for the purpose of this book. 
	
	\begin{prop}
		\label{prop:pAUS-blocking}
		Let $X$ be a Banach space admitting a shrinking biorthogonal system $(x_n,x^*_k)_{(n,k) \in \bN}$ in $X\times X^*$. Then, there exists  a sequence $1=m_1<\dots<m_n<\dots$ so that if $(E_n)_{n=1}^\infty$ is given by $E_n := \spa\{x_j\}_{j=m_n}^{m_{n+1}-1}$, then every skipped blocking $(F_k)_{k=1}^\infty$ of $(E_n)_{n=1}^\infty$ is a $2$-FDD of its closed linear span that satisfies for all $k\in \bN$ and all $(z_1,\dots,z_k) \in F_1 \times \dots \times F_k$,
		\begin{equation}
			\Big\|\sum_{i=1}^k z_i\Big\|\le 2e\big\|(\norm{z_i})_{i=1}^k\big\|_{\ell_{\bar{\rho}_X}}.
		\end{equation}
	\end{prop}
	
	\begin{proof}
		Let $N:=N_{\bar{\rho}_X}$ be the absolute norm on $\bR^2$ induced by the modulus of asymptotic uniform smoothness. In view of Proposition \ref{pro:Orlicz-iterated}, it will be enough to ensure that every skipped blocking $(F_k)_{k=1}^\infty$ of $(E_n)_{n=1}^\infty$ is a $2$-FDD of its closed linear span such that for all $k\in \bN$ and all $z_1 \in F_1, \dots, z_k\in F_k$,
		\begin{equation}
			\Big\|\sum_{i=1}^k z_i\Big\|\le 2\big\|(\norm{z_i})_{i=1}^k\big\|_{\Lambda_N},
		\end{equation}
		where $\norm{\cdot}_{\Lambda_N}$ is the iterated norm induced by $N$.
		
		We fix, as we may, a sequence $(\eps_k)_{k=1}^\infty$ in $(0,1)$ so that $\prod_{k=1}^\infty (1+4\eps_k)<2$. The key technical point is the following claim.
		\begin{claim}
			\label{claim:biortho-blocking}
			There are $1=m_1<m_2<\dots<m_k<\dots$ such that for all $k \in \bN$, all $x \in X_k := \spa(x_i)_{i< m_{k}}$ and all $z \in Z_k := \cspa (x_i)_{i\ge m_{k+1}}$,
			\begin{equation}
				\label{eq:upper-estimate-condition}
				\norm{x+z}\le N(\norm{x},\norm{z})(1+4\eps_k),
			\end{equation}
			and
			\begin{equation}
				\label{eq:FDD-condition}
				\norm{x}\le (1+2\eps_k)\norm{x+z}.
			\end{equation}
		\end{claim}
		
		\begin{proof}[Proof of Claim \ref{claim:biortho-blocking}]
			We choose $(m_k)_{k=1}^\infty$ inductively as follows. First, set $m_1 := 1$ and assume that for $k\in \bN$,  $m_1<m_2<\dots<m_{k}$ have been chosen. It follows from the definitions of $\bar{\rho}_X$ and $N$ that for any $t\ge 0$, any $x\in S_X$ and any $\vep_k>0$, there exists $Y\in \cof(X)$ such that for all $y\in S_Y$, $\norm{x+ ty}\le N(1,t)(1+\eps_k)$. We claim that there exists $Z\in \cof(X)$ such that $\norm{x+ tz}\le N(1,t)(1+2\eps_k),$ whenever $x\in S_{X_{k}}$, $z\in S_Z$ and $t\ge 0$. Indeed, since $\|x+tz\|\le 1+t$ and $N(1,t)\ge t$, we may assume that $t\le \frac{1}{2\eps_k}$ and the claimed statement is now an easy consequence of the compactness of $S_{X_{k}}$ and $[0,\frac{1}{2\eps_k}]$. Then, an homogeneity argument gives that for all $x\in X_{k}$ and $z\in Z$ we have
			\begin{equation*}
				\norm{x+z}\le N(\norm{x},\norm{z})(1+2\eps_k).
			\end{equation*}
			Next, we apply the standard Mazur technique for constructing basic sequences to argue that, by taking a smaller $Z\in \cof(X)$ if necessary, we can also get that for all $x\in X_{k}$ and $z\in Z$,
			\begin{equation*}
				\norm{x}\le (1+\eps_k)\norm{x+z}.
			\end{equation*}
			Finally, we use the fact that the linear span of $(x_n^*)_n$ is norm dense in $X^*$ in order to show that, up to a small perturbation (changing $\eps_k$ into $2\eps_k$), the space $Z$ can be taken to be of the form $Z_{k}:= \cspa (x_i)_{i\ge m_{k+1}}$ for some $m_{k+1}>m_{k}$. Let us give a word of explanation on this last classical approximation argument. So, let $Z\in \cof(X)$ be the space just obtained. There exist $y_1^*,\dots,y^*_r \in S_{X^*}$ such that $Z=\cap_{i=1}^r \ker y^*_i$. Let $F$ be the linear span of $\{y^*_1,\dots,y^*_r\}$ and note that $Z^\perp=F$. Fix $\eta>0$ and pick a finite $\eta$-net $(z^*_s)_{s\in S}$ of $S_F$. Since the linear span of $(x^*_n)_n$ is norm dense in $X^*$, there exist $m_{k+1}>m_k$ and $(v^*_s)_{s\in S}$ in the linear span of $\{x^*_1,\dots,x^*_{m_{k+1}-1}\}$ such that $\norm{v^*_s-z^*_s}\le  \eta$ for all $s \in S$. Note that it follows from biorthogonality that for $v\in Z_k$, $v^*_s(v)=0$ for all $s\in S$. 
			Now, identifying $(X/Z)^*$ with $Z^\perp=F$, we get that for $v\in S_{Z_k}$, $d(v,Z)=\sup_{x^*\in S_F}x^*(v)\le 2\eta$. It follows that there exists $z\in S_Z$ so that $\norm{v-z}\le 6\eta$. Note that for all $x\in X_k$, $1=\norm{z}\le (2+\eps_k)\norm{x+z}$. Thus, $\norm{x+v}\ge \frac{1}{2+\eps_k}-6\eta$ and $\norm{x}\le (1+\eps_k)\norm{x+z}\le (1+\eps_k)\norm{x+v}+12\eta$. So, if $\eta$ was initially chosen small enough, we obtain that $\norm{x}\le (1+2\eps_k)\norm{x+v}$ for all $x\in X_k$ and all $v\in S_{Z_k}$. Then, \eqref{eq:FDD-condition} follows by homogeneity. The proof of \eqref{eq:upper-estimate-condition} is similar. 
		\end{proof}
		To finish the proof, simply set $E_n := \spa\{x_j\}_{j=m_n}^{m_{n+1}-1}$ and consider a skipped blocking $(F_k)_k$ of $(E_n)_n$. It follows from \eqref{eq:FDD-condition} and the choice of $(\eps_k)_{k=1}^\infty$ that $(F_k)_k$ is a $2$-FDD. The upper estimate follows from \eqref{eq:upper-estimate-condition} since for all $z_1\in F_1, \dots, z_{k+1}\in F_{k+1}$, 
		\begin{align*}
			\Big\|\sum_{i=1}^{k+1} z_i\Big\| & \le N\Big(\Big\|\sum_{i=1}^k z_i\Big\|, \norm{z_{k+1}} \Big)(1+4\eps_k)\\
			& \le N\Big(N\Big(\Big\|\sum_{i=1}^{k-1} z_i\Big\|, \norm{z_{k}} \Big)(1+4\eps_{k-1}), \norm{z_{k+1}} \Big)(1+4\eps_k)\\
			& \le N\Big(N\Big(\Big\|\sum_{i=1}^{k-1} z_i\Big\|, \norm{z_{k}} \Big), \norm{z_{k+1}} \Big)(1+4\eps_{k-1})(1+4\eps_k)\\
			& \quad \vdots\\
			& \le \prod_{i=1}^k(1+4\eps_i)\big\|(\norm{z_i})_{i=1}^{k+1}\big\|_{\Lambda_N} \le 2\big\|(\norm{z_i})_{i=1}^{k+1}\big\|_{\Lambda_N},
		\end{align*}
		where we have used the monotonicity properties of absolute norms as well as the definition of Kalton's iterated norm.
	\end{proof}

	\section{\texorpdfstring{Asymptotic uniform flatness and subspaces of $\co$}{Asymptotic uniform flatness and subspaces of}}
	
	Recall that $\co$ is the prototypical example of a Banach space that is asymptotically uniformly flat. It is relatively simple to check that this is also true for any subspace of $\co$, and we leave the proof of the next proposition to the reader's diligence as Exercise \ref{ex:AUF}.
	
	\begin{prop}
		\label{prop:sub-c0-AUF}
		Every subspace of $\co$ is asymptotically uniformly flat.
	\end{prop}
	
	The goal of this section is to show that Proposition \ref{prop:sub-c0-AUF} admits the following converse, which will be crucial for the proof of the Lipschitz rigidity of $\co$ and of its subspaces in Section \ref{sec:Lip-rigidity-c_0}.
	
	\begin{theo}
		\label{thm:AUF->subspace-c_0}
		Let $X$ be a separable Banach space. If $X$ admits an equivalent norm that is asymptotically uniformly flat, then $X$ is isomorphic to a subspace of $\co$.
	\end{theo}
	
	Theorem \ref{thm:AUF->subspace-c_0} appeared in \cite{GKL2000} but could be found between the lines in \cite{KaltonWerner1995}. We will give here a simpler proof taken from \cite{JLPS2002}. This proof can be seen as a warm-up for a similar result about reflexive spaces that are both $p$-AUS and $p$-AUC. These results are based on some more or less delicate blocking arguments for finite-dimensional decompositions. An elementary observation is that in order to embed a space with an FDD into $\co$ it is sufficient to find a blocking satisfying an upper-$\infty$-estimate.
	
	\begin{prop}\label{prop:c0-blocking}
		If a Banach space $X$ with an FDD admits a blocking with an upper-$\infty$-estimate, then $X$ is isomorphic to a subspace of $\co$.
	\end{prop}
	
	\begin{proof}
		Assume that there is a blocking $(F_k)_k$ of an FDD $(E_n)_n$ of $X$, associated with the sequence $1=m_1<\dots<m_i<\dots$ that satisfies an upper-$\infty$-estimate with constant $C\ge 1$. If $x=(x_n)_n$ is the decomposition of $x$ with respect to $(E_n)_n$, then we let $P_{F_k}(x) := \sum_{n=m_k}^{m_{k+1}-1} x_n$. Clearly $\sup_k \norm{P_{F_k}}\le 2K$ where $K$ is the FDD constant of $(E_n)_n$. Now define 
		$T\colon X=\oplus_{n=1}^\infty E_n \to (\sum_{k=1}^\infty F_k)_{\co}$ by $T(x) = (P_{F_k}(x))_k$.
		Observe that $T$ is clearly linear and well defined since $\lim_k\norm{P_{F_k}(x)}\to 0$ for all $x\in X$, as $(m_i)_i$ is strictly increasing. Moreover, it follows from the upper-$\infty$-estimate that $T$ is a linear embedding. Indeed,  
		\begin{equation*}
			\norm{T((x_n)_n)} = \sup_k \norm{\sum_{n=m_k}^{m_{k+1}-1} x_n} \le 2K \norm{\sum_{k=1}^\infty \sum_{n=m_k}^{m_{k+1}-1} x_n} = 2K\norm{(x_n)_n},
		\end{equation*} 
		while 
		\begin{equation*}
			\norm{T((x_n)_n)} \ge  \frac{1}{C} \norm{\sum_{k=1}^\infty \sum_{n=m_k}^{m_{k+1}-1} x_n} = \frac{1}{C}\norm{(x_n)_n}.
		\end{equation*}
		Since $T$ clearly has dense range, we deduce that it is an isomorphism. Finally, the conclusion follows by observing that $(\sum_{k=1}^\infty F_k)_{\co}$ embeds into $\co$ (see Exercise \ref{ex:co-sum-fds}).
	\end{proof}

	Another ingredient in the proof of Theorem \ref{thm:AUF->subspace-c_0} is the following observation, due to Lindenstrauss, which follows from the separable injectivity of $\co$.
	
	\begin{lemm}\label{lem:3-space-property-c0}
		Let $Y$ be a closed subspace of a separable Banach space $X$. If $Y$ and $X/Y$ embed isomorphically into $\co$, then $X$ embeds isomorphically into $\co$.
	\end{lemm}
	
	\begin{proof}
		Let $U$ be a linear embedding from $Y$ into $\co$ and $V$ be a linear embedding from $X/Y$ into $\co$. Since $X$ is separable, by Sobczyk's Theorem (see Theorem 2.f.5 in  \cite{LindenstraussTzafriri1977}), there exists a bounded linear map $S: X \to c_0$ which extends $U$. It is now easy to check that the map $T$ defined by $Tx := (Sx,VQx)$, where $Q$ is the quotient map from $X$ onto $X/Y$, is a linear embedding from $X$ into $c_0 \oplus_\infty c_0 \equiv \co$. Let us detail this classical argument. So, let $c>0$ such that $\norm{Uy}\ge c\norm{y}$ for all $y\in Y$, $\norm{Vz}\ge c\norm{z}$ for all $z\in X/Y$ and fix $\eta \in (0,1)$. Let $x\in X$ and observe that if $d(x,Y)\ge \eta \norm{x}$, then $\norm{VQx}\ge c\norm{Qx}\ge c\eta \norm{x}$. If $d(x,Y)< \eta \norm{x}$, then there exists $y\in Y$ such that $\norm{x-y}\le \eta\norm{x}$. Then, $\norm{y}\ge (1-\eta)\norm{x}$, $\norm{Sy}=\norm{Uy}\ge c(1-\eta)\norm{x}$ and $\norm{Sx}\ge \big(c(1-\eta)-\norm{S}\eta\big)\norm{x}$. It is now clear that an appropriate initial choice of $\eta$ yields that $T$ is a linear embedding. 
	\end{proof}
	
	As we will see shortly, the proof of Theorem \ref{thm:AUF->subspace-c_0} reduces to the case of asymptotically uniformly flat Banach spaces admitting a special type of FDD, namely a shrinking FDD. Recall that a FDD $(E_n)_{n=1}^\infty$ of a Banach space $X$ is \emph{shrinking} if and only if for all $f\in X^*$, $\lim_n \norm{f_{\restricted \spa(\cup_{k\ge n} E_k)}}=0$. Another way to say this is the following. Let $P_n$ be the canonical projection from $X$ onto $E_n$, associated with the FDD $(E_n)_{n=1}^\infty$ and consider $P_n$ as a bounded operator on $X$. Then, set $F_n := P_n^*(X^*)$. If $(E_n)_{n=1}^\infty$ is a shrinking FDD of $X$, then $(F_n)_{n=1}^\infty$ is an FDD of $X^*$ with associated projections given by  $(P_n^*)_n$. Also, $F_n$ is canonically isomorphic to $E_n^*$ (through $P_n^*$, if $P_n$ is viewed as an operator from $X$ onto $E_n$). In that case,  note that it is easy to build a shrinking biorthogonal system $(x_n,x^*_k)_{n,k}$ in $X \times X^*$ by pasting together finite biorthogonal systems in $E_n\times F_n$. The core of the argument is to show that when $X$ is asymptotically uniformly flat and admits a shrinking FDD $(E_n)_{n=1}^\infty$, then $X$ is isomorphic to a $\co$-sum of a certain blocking of the FDD.
	
	\begin{proof}[Proof of Theorem \ref{thm:AUF->subspace-c_0}] 
		Assume first that $X$ has a shrinking FDD $(E_n)_{n=1}^\infty$. Since $X$ is AUF there is $t_0\in(0,1]$ such that $\bar{\rho}_X(t)=0$ for all $t\le t_0$. Since $X$ admits a shrinking biorthogonal system, we can apply Proposition \ref{prop:pAUS-blocking} and $(i)$ of Proposition \ref{prop:Orlicz-pnorm} to deduce the existence of a constant $C\ge 1$ and of a sequence $1=m_1<\dots<m_n<\dots$ so that the associated blocking $(G_n)_{n}$ of $(E_i)_i$ is such that for every skipped blocking $(F_k)_k$ of $(G_n)_{n}$ and for all $k\in \bN$ and all $(x_1,\dots,x_k) \in F_1 \times \cdots \times F_k$,
		\begin{equation*}
			\Big\|\sum_{i=1}^k x_i\Big\|\le C\big\|(\norm{x_i})_{i=1}^k\big\|_{\infty}.
		\end{equation*}
		We need to add a small warning here. When applying Proposition \ref{prop:pAUS-blocking} to the biorthogonal system associated with the FDD $(E_n)_n$, we need to make sure that the blocking we obtain is also a blocking of the FDD. We leave it to the reader to convince himself or herself that this is possible. 
		
		Considering now the skipped blockings  $(G_{2n})_{n=1}^\infty$ and $(G_{2n-1})_{n=1}^\infty$ of $(G_n)_{n=1}^\infty$ and applying the triangle inequality we get that for all $k\in \bN$ and all $(x_1,\dots,x_k) \in G_1 \times \dots \times G_k$,
		\begin{equation*}
			\Big\|\sum_{i=1}^k x_i\Big\|\le 2C\big\|(\norm{x_i})_{i=1}^k\big\|_{\infty}.
		\end{equation*}
		The conclusion of the proof then follows from Proposition \ref{prop:c0-blocking}.

		We now turn to the general case. For that, it is convenient to anticipate some results from Chapter \ref{chapter:Szlenk}. Since $X$ is separable and AUF and in particular AUS, we deduce from Theorem \ref{thm:Szlenk-Asplund} and Proposition \ref{prop:Sz-omega} that $X^*$ is separable. We can now appeal to a result of Johnson and Rosenthal \cite{JohnsonRosenthal1972}, ensuring the existence of a subspace $Y$ of $X$ such that both $Y$ and $X/Y$ admit a shrinking FDD (see also \cite[Theorem 1.g.2]{LindenstraussTzafriri1977} or \cite[Theorem 4.34]{FHHMZ2011}). Since being AUF passes to subspaces and quotients, we deduce from the first part of the proof that both $Y$ and $X/Y$ linearly embed into $\co$. The conclusion follows from Lemma \ref{lem:3-space-property-c0}.
	\end{proof}
	
	Combining the results in this section, we obtain the following remarkable characterization of subspaces of $\co$ in terms of the existence of an asymptotically uniformly flat renorming.
	
	\begin{coro}
		\label{cor:AUF-c0}
		%    Let $X$ be a separable Banach space. Then, $X$ admits an equivalent norm that is asymptotically uniformly flat if and only if $X$ is isomorphic to a subspace of $\co$.
		A Banach space $X$  is isomorphic to a subspace of $\co$ if and only if it is separable and admits an equivalent norm that is asymptotically uniformly flat.
	\end{coro}

	\section{\texorpdfstring{Subspaces of $\ell_p$-sums of finite-dimensional spaces}{Subspaces of -sums of finite-dimensional spaces}}
	
	We have already seen that a subspace of an $\ell_p$-sum of finite-dimensional spaces is $p$-AUS and $p$-AUC. The following theorem, also taken from \cite{JLPS2002}, is a converse statement.
	
	\begin{theo}
		\label{thm:subspaces-of-lp-sums} 
		Let $p\in (1,\infty)$. If a reflexive and separable Banach space admits an equivalent $p$-AUS norm and an equivalent $p$-AUC norm, then it is isomorphic to a subspace of an $\ell_p$-sum of finite-dimensional spaces. 
	\end{theo}
	
	The proof of Theorem \ref{thm:subspaces-of-lp-sums} relies on similar ideas already discussed in the proof of the characterization of asymptotically uniformly flat spaces. However, the technical details are slightly more complicated and require some additional knowledge of Banach space theory.
	
	It is well known that a separable Banach space may fail to have a Schauder basis, but in any case it admits a Markushevich basis (see \cite[Section 4.8]{FHHMZ2011} or \cite{HMVZ2008} and references therein). A consequence of this fact is that any separable reflexive Banach space admits a shrinking biorthogonal system (cf. Exercise \ref{ex:reflexive-biorthogonal}). 
	
	\begin{proof}[Proof of Theorem \ref{thm:subspaces-of-lp-sums}] 
		We may and do assume that the norm of $X$ is $p$-AUS. Let $(x_n,x^*_k)_{n,k \in\bN}$ be a shrinking biorthogonal system in $X\times X^*$. Then, consider $(m_n)_{n=1}^\infty$ and $(E_n)_{n=1}^\infty$ provided by Proposition \ref{prop:pAUS-blocking}. We now use the fact that $X$ also admits an equivalent $p$-AUC norm and by duality, since $X$ is reflexive, $X^*$ admits an equivalent $q$-AUS norm, where $q$ is the conjugate exponent of $p$. By passing to a further subsequence of $(m_n)_n$ and adjusting $C$, we may therefore assume that if $F_n := \spa\{x_j^*\}_{j=m_n}^{m_{n+1}-1}$, then every skipped blocking of $(F_n)_n$ is a $2$-FDD of its closed linear span satisfying an upper-$q$-estimate with constant $C$. 
		
		Let $n \ge 3$, $x\in E_{n}$ and $y\in Y_{n-2} := E_1+\dots+E_{n-2}$. Then, $\norm{x}\le \norm{x-y}+\norm{y}\le 3\norm{x-y}$. Therefore, for all $x\in E_n$,
		\begin{equation*}
			\norm{x}\le 3\,d(x,Y_{n-2})=3\sup_{x^*\in B_{Y_{n-2}^\perp}} x^*(x).
		\end{equation*}
		It is easily checked that $Y_{n-2}^\perp$ is the closed linear span of the $F_i$ for $i\ge n-1$. Then, using the compactness of $B_{E_n}$ and a standard approximation argument, we could ensure, by passing to a further blocking if needed, that for all $x\in E_n$,
		\begin{equation*}
			\norm{x}\le 4\sup_{x^*\in B_{H_n}} x^*(x),
		\end{equation*}
		where $H_n := F_{n-1}+F_n+F_{n+1}$. Therefore, we may assume that $H_n$ is 4-norming for $E_n$. 
		
		Consider now the map $U\colon X\to (\sum_{n=1}^\infty H_{4n})_q^* = (\sum_{n=1}^\infty H_{4n}^*)_p $ defined by $Ux((y^*_n)_{n=1}^\infty)=\sum_{n=1}^\infty y^*_n(x)$. Since $(H_{4n})_n$ is a skipped blocking of $(F_n)_n$ it satisfies upper-$q$-estimates. Therefore, $U$ is a well-defined, bounded linear map. We recall that $(E_{4n})_n$ is a $2$-FDD of its closed linear span, which we denote by $Y$. So, for $x\in Y$, we can write $x=\sum_{n=1}^\infty x_n$ with $x_n \in E_{4n}$. Since $H_{4n}$ is finite-dimensional and $4$-norming for $E_{4n}$, we can pick a sequence $(y_n^*)_{n=1}^\infty \in \prod_{n=1}^\infty H_{4n}$ such that $\sum_{n=1}^\infty \norm{y_n^*}^q\le 1$ and
		$\sum_{n=1}y_n^*(x_n)\ge \frac14 \big(\sum_{n=1}^\infty \norm{x_n}^p\big)^{1/p}$. Since $E_{4n}$ satisfies an upper-$p$-estimate with constant $C$, we deduce that for all $x\in Y$, so that $x=\sum_{n=1}^\infty x_n$ with $x_n \in E_{4n}$, we have
		\begin{equation*}
			\norm{Ux}\ge \Big|\sum_{n=1}^\infty y^*_n(x)\Big|=\Big|\sum_{n=1}^\infty y^*_n(x_n)\Big|\ge \frac14 \big(\sum_{n=1}^\infty \norm{x_n}^p\big)^{1/p} \ge \frac{1}{4C} \norm{x}.
		\end{equation*}
		Thus $U$ is an isomorphism from $Y$ onto a subspace of $(\sum_{n=1}^\infty H_{4n})_q^*$.
		
		Observe that $Y^\perp$ is the closed linear span of the $F_k$ for $k\notin 4\bN$ and let $G_n := F_{4n-3}+F_{4n-2}+F_{4n-1}$. Since $(G_n)_{n=1}^\infty$ is a skipped blocking of $(F_n)_{n=1}^\infty$ it is a $2$-FDD of its closed linear span $Y^\perp$, satisfying an upper-$q$-estimate. It follows from the reflexivity of $X$ that $X/Y=(X/Y)^{**}= (Y^{\perp})^*=(\oplus_n G_n)^*=\oplus_n G_n^*$ and thus $(G_n^*)_{n=1}^\infty$ can be viewed as a FDD of $X/Y$ with a lower-$p$-estimate. Since being $p$-AUS and reflexive passes to quotients, it follows that $X/Y$ is a reflexive Banach space which is $p$-AUS. So, we can apply again Proposition \ref{prop:pAUS-blocking} to find a blocking $(B_k^*)_k$ of $(G_n^*)_n$ such that every skipped blocking of $(B_k^*)_k$ satisfies both an upper-$p$-estimate and a lower-$p$-estimate. Since the FDD $(B_k^*)_k$ is the union of two of its skipped blockings, we have that $(B_k^*)_k$ itself satisfies both an upper-$p$-estimate and a lower-$p$-estimate. Let us detail this last argument. First, the fact that $(B_k^*)_k$ satisfies an upper-$p$-estimate is just a consequence of the triangle inequality. For the lower-$p$-estimate, we have to use duality. Since our setting is reflexive, the FDD $(B_k^*)_k$ is both shrinking and boundedly complete and, by passing to a further blocking if necessary, its dual FDD $(C_k)_k$ is such that every skipped blocking of $(C_k)_k$ satisfies an upper-$q$-estimate, where $q$ is the conjugate exponent of $p$ (by the existence of an equivalent $q$-AUS norm on its linear span and Proposition  \ref{prop:pAUS-blocking}). Therefore, $(C_k)_k$ itself satisfies an upper-$q$-estimate and, by duality (and reflexivity) again, $(B_k^*)_k$ satisfies a lower-$p$-estimate. Therefore, there is an isomorphism $T \colon X/Y \to (\sum_k B_k^*)_p$.
		
		Finally if $Q$ is the quotient map from $X$ onto $X/Y$, we define $S\colon X \to Z=(\sum_{n=1}^\infty H_{4n})_q^* \oplus_p (\sum_{k=1}^\infty B_k^*)_p$ by $Sx := (Ux,TQx)$. Since $U$ and $T$ are linear embeddings, we deduce that $S$ is a linear embedding (the argument is identical to the end of the proof of Lemma \ref{lem:3-space-property-c0}).   
	\end{proof}
	
	\begin{rema} It is important to mention that the reflexivity assumption is important for Theorem \ref{thm:subspaces-of-lp-sums}. Indeed, the classical James space $\James$ (see Appendix \ref{sec:James-space}) is not reflexive but admits an equivalent $2$-AUS norm and an equivalent $2$-AUC norm \cite{LancienBesac}. 
		
	\end{rema}
	
	\section{\texorpdfstring{The property $(\beta)$ of Rolewicz}{The property of Rolewicz}}\label{sec:Rolewicz}
	
	This section is devoted to the study of a property of norms introduced by S. Rolewicz in \cite{Rolewicz1987} and now called the property $(\beta)$ of Rolewicz. For its definition, we shall use a characterization due to D. Kutzarova \cite{Kutzarova1990}. 
	
	\begin{defi} 
		An infinite-dimensional Banach space $X$ is said to have \emph{property $(\beta)$} if for any $t\in (0,1]$,  there exists $\delta>0$ such that for any $x$ in $B_X$ and any $t$-separated sequence $(x_n)_{n=1}^\infty$ in $B_X$, there exists $n\ge 1$ so that
		$${\norm{\frac{x+x_n}{2}}}\le 1-\delta.$$
		For a given $t\in (0,1]$, we denote by $\bar{\beta}_X(t)$ the supremum of all $\delta\ge 0$ so that the above property is satisfied.
	\end{defi}
	
	\begin{rema}
		\begin{enumerate}
			\item Because of the symmetry of $B_X$, $\norm{\frac{x+x_n}{2}}$ can be equivalently replaced by $\norm{\frac{x-x_n}{2}}$ in this definition.
			\item Note that if $X$ has property $(\beta)$, then for all $t\in(0,1]$ we have $\bar{\beta}_X(t)\in(0,1)$. 
		\end{enumerate}
		
	\end{rema}
	
	Property $(\beta)$ is a natural extension of uniform convexity from a particular asymptotic perspective (see Exercise \ref{ex:uc-beta}). Since we are considering infinite separated sequences rather than a pair of separated vectors, and because every infinite separated sequence must escape any finite-dimensional subspace, property $(\beta)$ does not capture local finite-dimensional phenomena. Similarly to uniform convexity, property $(\beta)$ implies reflexivity. This fact is due to Rolewicz \cite{Rolewicz1987} and is an easy application of the existence of James' sequences in nonreflexive Banach spaces.
	
	\begin{prop}
		\label{prop:beta->reflexive} 
		Every Banach space with property $(\beta)$ is reflexive.
	\end{prop}
	
	\begin{proof} 
		Let $X$ be a nonreflexive Banach space. Then, it follows from \cite{James1963/64} that for any $\theta \in [\frac12,1)$, there exist sequences $(x_n)_{n=1}^\infty$ in $S_X$ and $(x_n^*)_{n=1}^\infty$ in $S_{X^*}$ such that $x_n^*(x_i)>\theta$ for all $i \ge n$ and $x_n^*(x_i)=0$ for all $ i < n$. This clearly implies that $(x_n)_{n=1}^\infty$ is $\frac12$-separated, while $\|x_1+x_n\|\ge 2\theta$ for all $n\in \bN$. Since $\theta$ can be chosen arbitrarily close to $1$, this implies that $X$ does not have property $(\beta)$.
	\end{proof}
	
	Asymptotic uniform convexity and property $(\beta)$ are two different ways of relaxing the notion of uniform convexity. The main purpose of this section is to explore in detail the links between property $(\beta)$ and the asymptotic uniform smoothness or convexity of a Banach space. In particular, we will explain that property $(\beta)$ also implies some smoothness properties. We start with two isometric results.  The proof of the first one is an adaptation of ideas that can be found between the lines in \cite{Huff1980} and explicitly in \cite[Lemma 4.0.2]{DKR2016}.
	
	\begin{theo}
		\label{thm:beta->AUC}
		Any Banach space with property $(\beta)$ is asymptotically uniformly convex. More precisely, $\bar{\delta}_X(t)\ge \bar{\beta}_X(\frac{t}{2})$ for all $t\in (0,1)$.
	\end{theo}
	
	\begin{proof} 
		Assume that $X$ is a Banach space with property $(\beta)$ and let $t\in (0,1)$. Fix $x\in S_X$ and assume without loss of generality that $\bar{\delta}_X(t,x)<\delta\le 1$. By definition, for all $Y \in \cof(X)$ there exists $x_Y \in S_Y$ such that $\norm{x + tx_Y}\le 1+\delta$. We pick $x^*\in S_{X^*}$ such that $x^*(x)=1$. Then, we build inductively a sequence  $(x_n^*)_{n=1}^\infty$ in $S_{X^*}$ such that for all $n\ge 1$, $x_n^*(x_{Y_n})=1$, where $Y_1 = \ker(x^*)$ and $Y_n = Y_1\cap \bigcap_{i=1}^{n-1} \ker(x_i^*)$ for $n\ge 2$.\\
		Next we set $y_n=(1+\delta)^{-1}(x + tx_{Y_n})$. It readily follows from our construction that $(y_n)_{n=1}^\infty$ is a $\frac{t}{2}$-separated sequence in $B_X$. We can now use the property $(\beta)$ of $X$ to deduce that, for any $\eta >0$,  there exists $n>1$ such that 
		$$\Big\|\frac{y_1+y_n}{2}\Big\|\le 1-\bar{\beta}_X(\frac{t}{2})+\eta.$$ 
		On the other hand,
		\begin{equation*}
			\frac{\norm{y_1+y_n}}{2}\ge x^*\Big(\frac{y_1+y_n}{2}\Big)=\frac{1}{1+\delta}.
		\end{equation*} 
		So we have shown that, for any $\eta >0$,  $\delta \ge \big(1-\bar{\beta}_X(\frac{t}{2})+\eta \big)^{-1}-1$. Since this true for any $x\in S_X$ and any $\delta>\bar{\delta}_X(x,t)$ and since $\bar{\beta}_X(\frac{t}{2})\in (0,1)$, we conclude that 
		\begin{equation*}
			\bar{\delta}_X(t)\ge \big(1-\bar{\beta}_X(\frac{t}{2})\big)^{-1}-1\ge \bar{\beta}_X(\frac{t}{2})>0.
		\end{equation*}
		Therefore, $X$ is AUC. 
	\end{proof}
	
	Next, we show that if $X$ is a reflexive Banach space whose norm is both AUS and AUC, then $X$ has property $(\beta)$. This result is due to D. Kutzarova \cite{Kutzarova1990} with a different terminology. We shall present here a sharper quantitative version of this result, which is taken from \cite{DKRRZ2013}.
	
	\begin{theo}
		\label{thm:beta->AUS&AUC} 
		For every reflexive Banach space $X$ whose norm is $p$-AUS and $q$-AUC for some $1<p\le q<\infty$, there exists $c>0$ such that for all $t \in (0,1]$,
		\begin{equation*}
			\bar{\beta}_{X}(t)\ge ct^{\frac{p(q-1)}{p-1}}.
		\end{equation*}
		In particular, $X$ has property $(\beta)$.
	\end{theo}
	
	\begin{proof}
		Recall that since $X$ is $p$-AUS and $q$-AUC, there exist $a,b>0$ such that for every weakly null sequence $(w_n)_n\subset X$ and $x\in X$, one has 
		\begin{equation}
			\label{eq:lower-q-est}
			(\norm{x}^q + b\liminf_n \norm{w_n}^q)^{\frac{1}{q}}\le \liminf_n \norm{x+w_n}
		\end{equation}
		and 
		\begin{equation}
			\label{eq:upper-p-est}
			\limsup_n \norm{x+w_n} \le (\norm{x}^p + a\limsup_n\norm{w_n}^p)^{\frac{1}{p}}.
		\end{equation}
		For convenience, we will assume, as we may, that $b\le 2^{-q}$. Let now $t\in (0,1]$ and fix $\delta >\bar{\beta}_{X}(t)$. Then, there exist a $t$-separated sequence $(x_n)_n$  in $B_X$ and $x\in B_X$ such that $\norm{x+x_n}\ge 2(1-\delta)$ for all $n\in \bN$. Since $X$ is reflexive, we may assume that $(x_n)_n$ is weakly converging to $u\in B_X$. Passing to a further subsequence, we may also assume that $\lim_{n \to \infty} \norm{x_n-u}=:s$. Note that clearly $s\le 2$ and since $(x_n)_n$ is $t$-separated in $B_X$ we also have $s\ge \frac{t}{2}$. Passing again to another subsequence if needed, we can assume that $\inf_n\norm{x_n-u}\ge \frac{s}{2}$. Observe now that there is $\alpha_q\in(0,1)$ such that for all $z\in(0,1]$, $(1+z)^{1/q}\ge 1+\alpha_q z$. Using the $q$-AUC property of $X$, we will show the following nontrivial upper bound on the norm of $u$: 
		\begin{equation}
			\label{eq:upper-bound-u}
			\norm{u}\le 1 - \frac{\alpha_qb}{2}s^q.
		\end{equation}
		Since $\alpha_qbs^q\le 1$, in order to prove \eqref{eq:upper-bound-u} we may assume that $\norm{u}\ge \frac12$. It follows that
		\begin{align*}
			1\ge \liminf_{n \to \infty}\norm{u+x_n-u} & \stackrel{\eqref{eq:lower-q-est}}{\ge} \Big( \norm{u}^q + b\liminf_n \norm{x_n-u}^q \Big)^{1/q}\\
			& \ge \norm{u}\Big(1+b\frac{s^q}{\norm{u}^q} \Big)^{1/q}\\
			& \ge \norm{u}(1+bs^q)^{1/q}\\
			& \ge \norm{u} +\norm{u}\alpha_qbs^q \quad (\textrm{since } bs^q\in(0,1]),
		\end{align*}   
		and we deduce from $\norm{u}\ge \frac12$ that $\norm{u}\le 1-\frac{\alpha_q b}{2}s^q$, which proves \eqref{eq:upper-bound-u}.
		
		We will now use the $p$-AUS property to find a midpoint of the form $\frac{x+x_{n_0}}{2}$ that is deep inside the unit ball. In order to achieve this, let $0<\lambda<1$ to be carefully chosen later, and note that it follows readily from the triangle inequality that 
		%$$\|(1-\lambda)x+\lambda u\|\le (1-\lambda)+\lambda(1-\frac{b}{2}s^q)=1-\frac{b}{2}\lambda s^q.$$
		$$\norm{(1-\lambda)x+\lambda u}\le 1-\frac{\alpha_q b}{2}\lambda s^q.$$
		%So $y=(1-\frac{b}{2}\lambda s^q)^{-1}((1-\lambda)x+\lambda u) \in B_X$. 
		We then have
		\begin{align*}
			\limsup_{n\to \infty}\norm{(1-\lambda)x + \lambda x_n} 
			& = \limsup_{n \to \infty}\Big\| (1-\lambda)x+\lambda u + \lambda(x_n-u)\Big\|\\
			& \stackrel{\eqref{eq:upper-p-est}}{\le} (\norm{(1-\lambda)x + \lambda u}^p + a\lambda^p\limsup_n \norm{x_n-u}^p)^{\frac{1}{p}}\\
			& \le \Big(\big(1-\frac{\alpha_q b\lambda}{2}s^q\big)^p + a\lambda^p\limsup_n \norm{x_n-u}^p\Big)^{\frac{1}{p}}\\
			&\le \Big(1-\frac{\alpha_q b\lambda}{2}s^q\Big)\Big(1+a\lambda^ps^p\big(1- \frac{\alpha_q b\lambda}{2} s^q\big)^{-p}\Big)^{\frac{1}{p}}\\
			&\le \Big(1-\frac{\alpha_q b\lambda}{2}s^q\Big)\Big(1+\frac{1}{p}2^pa\lambda^ps^p\Big).
		\end{align*}
		In particular, there exists $n_0 \in \bN$ such that 
		\begin{equation*}
			\norm{(1-\lambda)x+x_{n_0}} \le 1-\frac{\alpha_q b\lambda}{2}s^q + 2^{p}a\lambda^ps^p.
		\end{equation*}
		We now turn back to the initial choice of $\lambda$ that we pick in $(0,\frac12)$ such that $2^{p}a\lambda^p s^p \le \frac{\alpha_q b\lambda}{4} s^q$ (which we can since $p\in(1,\infty)$), i.e. such that $\lambda \le \Big(\frac{\alpha_q b}{2^p4a}\Big)^{\frac{1}{p-1}}s^{\frac{q-p}{p-1}}$. Thus, we can pick $\lambda \in (0,\frac12)$ satisfying this inequality and also such that 
		\begin{equation*}
			\norm{(1-\lambda)x + \lambda x_{n_0}} \le 1-\frac{\alpha_q b\lambda}{4}s^q\le 1-c_1s^\sigma,
		\end{equation*}
		for $\sigma := \frac{p(q-1)}{p-1}$ and some $c_1>0$ (depending on $a,b,p,q$). We now write 
		\begin{align*}
			\Big\|\frac{x+x_{n_0}}{2}\Big\| 
			& = \Big\|\frac{1}{2(1-\lambda)}\big((1-\lambda)x+\lambda x_{n_0}\big)+\frac{1-2\lambda}{2(1-\lambda)}x_{n_0}\Big\| \\
			& \le \frac{1}{2(1-\lambda)}(1-c_1s^\sigma)+\frac{1-2\lambda}{2(1-\lambda)}=1-\frac{c_1s^\sigma}{2(1-\lambda)} \le 1-\frac{c_1}{2}s^\sigma \\
			& \le 1-\frac{c_1}{2^{\sigma+1}}t^\sigma=1-ct^\sigma,\ \ \ \text{with}\ c := \frac{c_1}{2^{\sigma+1}}.
		\end{align*}
		We deduce that $1-\delta \le 1-ct^\sigma$ and since $\delta >\bar{\beta}_X(t)$ was arbitrary, we conclude that $\bar{\beta}_X(t)\ge ct^\sigma=ct^{\frac{p(q-1)}{p-1}}$.
	\end{proof}

	\begin{rema} 
		It is important to note that the original proof from \cite{Kutzarova1990} is also quantitative. Keeping track of the exponents in \cite{Kutzarova1990}, one would get the following slightly weaker estimate: $\bar{\beta}_{X}(t)\ge ct^{\frac{pq}{p-1}}$, for $t\le 1$.
	\end{rema}
	
	\begin{rema} 
		In chapter \ref{chapter:Szlenk} we shall prove an ``isomorphic'' converse of this theorem by showing that a Banach space admits an equivalent norm with property $(\beta)$ if and only if it is reflexive and it admits an equivalent AUS norm and an equivalent AUC norm.   
	\end{rema}

	\section{Notes}
	The asymptotic moduli studied in this chapter were introduced by V.D. Milman in a series of short papers written in Russian \cite{Milman67,Milman68,Milman69} where numerous results were announced without proofs. In the paper \cite{Milman1971}, written in English, proofs were provided and many moduli were studied in terms of $\beta$ and $\delta$-averages in the original terminology. Given a family $\cB$ of subspaces of a Banach space $X$ and a real-valued bounded function $f$ defined on the unit sphere of $X$, the $\beta$-average (resp. $\delta$-average) is defined as $\beta[f,\cB]:=\sup_{Y\in \cB}\inf_{y\in Y; \norm{y}=1}f(y)$ (resp. $\delta[f,\cB]:=\inf_{Y\in \cB}\sup_{y\in Y; \norm{y}=1}f(y)$). For instance, $\bar{\delta}_X(x,t)$ (resp. $\bar{\rho}_X(x,t)$) corresponds in Milman's terminology to $\beta^0(t;x,X):=\beta(t;x,\cB^0):=\beta_y[f_0,\cB^0]$ (resp. $\delta^0(t;x,X):=\delta(t;x,\cB^0):=\delta_y[f_0,\cB^0]$) (where $\cB^0$ is the set of finite-codimensional subspaces and $f_0(t;x,y):=\norm{x+ty}-1)$). Then, $\bar{\delta}_X(t)$ (resp. $\bar{\rho}_X(t)$) corresponds in Milman's terminology to $\beta^0(t;x,X):=\beta(t;x,\cB^0):=\beta_y[f_0,\cB^0]$ (resp. $\delta^0(t;x,X):=\delta(t;x,\cB^0):=\delta_y[f_0,\cB^0]$).
	The fundamental role played by the moduli $\inf_{x\in S_X}\beta^0(t;x,X)$ and $\sup_{x\in S_X}\delta^0(t;x,X)$ was recognized by Johnson, Lindenstrauss, Preiss and Schechtman in their influential paper \cite{JLPS2002} were they coined the modern terminology and notation. 
	
	Milman had proved (cf. \cite[Corollary 4.1, Theorem 5.5 c.]{Milman1971} that if $X$ is asymptotically uniformly flat (resp. $\bar{\delta}_X(t)\ge ct$), then $X$ contains a copy of $\co$ (resp. $X$ contains a copy of $\ell_1$). In fact, Milman also knew that these spaces are saturated with those copies. At the time of Milman's work, it was still an open question whether every infinite-dimensional Banach space contained a copy of $\ell_p$ ($1\le p<\infty$) or $\co$ (The Fundamental Structure Conjecture in \cite[page 145]{Milman1971}). The introduction of the moduli by Milman seems to have been greatly motivated by trying to answer positively the Fundamental Structure Conjecture. Indeed, Milman stated a conjecture (Conjecture $1$ in \cite[page 145]{Milman1971}) which, if true, would imply that every infinite-dimensional Banach space contains a copy of $\ell_p$ ($1\le p<\infty$) or $\co$. The derivation of the veracity of the Fundamental Structure Conjecture from Conjecture $1$ (\cite[Conditional Theorem 5.2]{Milman1971} made heavy use of the moduli. It is also interesting to note that Milman's Conjecture $1$ is closely connected to the Distortion Problem (cf. Conjecture $2$ in \cite[page 145]{Milman1971}), which was still open at the time. Both problems turned out to have a negative answer as shown to us by Tsirelson \cite{Tsirelson1974} and Odell and Schlumprecht \cite{OdellSchlumprecht1994}. Nevertheless, Milman's ideas ended up being extremely influential. One outgrowth of those is the notion of \emph{asymptotic structure}, which was formalized and systematically studied by Maurey, Milman and Tomczak-Jaegermann \cite{MMTJ} in the mid 1990s. The asymptotic structure is central to the topic of this book. Many linear isomorphic properties, which are asymptotic in the sense of asymptotic structure, can be shown to be invariant under various nonlinear maps such as Lipschitz, uniform, or coarse isomorphisms or embeddings.
	
	Kalton and Werner showed that if a separable Banach space $X$ does not contain a copy of $\ell_1$ and  satisfies property $m_p$ (i.e. $\limsup_n\norm{x+x_n}= \norm{(\norm{x},\limsup_n\norm{x_n})}_p$ for every $x\in X$ and weakly null sequence $(x_n)_n$), then $X$ (linearly) embeds almost isometrically into an $\ell_p$-sum of finite-dimensional spaces (when $1\le p<\infty$) or $\co$ (when $p=\infty$) (see \cite[Theorems 3.3 and 3.5]{KaltonWerner1995}). 
	Godefroy, Kalton and Lancien used a dual approach in \cite{GKL2000} to show that if $X$ is AUF, then $X$ isomorphically embeds into $\co$. The central notion in \cite{GKL2000} is a quantitative refinement of the weak$*$-Kadec-Klee property (see Exercise \ref{ex:LKK*-AUF}). The argument in \cite{GKL2000} is different from that of \cite{JLPS2002} presented here since it does not use the reduction to the case of Banach spaces with shrinking FDD. However, the isomorphism constant in \cite{GKL2000} is lower than that from \cite{JLPS2002}.
	
	In the remarkable paper \cite{LimaLova2012}, Lima and Randrianarivony noticed the relevance of Rolewicz's property $(\beta)$ in the study of nonlinear quotients of Banach spaces, and consequently, property $(\beta)$ became a central tool towards understanding the nonlinear geometry of Banach spaces.

	\section{Exercises}
	
	\begin{exer}
		\label{ex:moduli-basic}
		Let $X$ be a Banach space. Show the following properties of the moduli. 
		\begin{enumerate}[(i)]
			\item For all $t>0$, $0\le \bar{\delta}_X(t)\le \bar{\rho}_X(t)\le t$.
			\item $\bar{\delta}_X$ is nondecreasing and $1$-Lipschitz.
			\item $\bar{\rho}_X$ is nondecreasing, $1$-Lipschitz and convex.
		\end{enumerate}
	\end{exer}

	\begin{exer}
		\label{ex:AUS-subspace-quotient}  Let $X$ be a Banach space and $Y$ be a closed subspace of $X$.
		\begin{enumerate}
			\item Show that $\bar{\rho}_Y(t)\le \bar{\rho}_X(t)$ and $\bar{\delta}_Y(t)\ge \bar{\delta}_X(t)$.
			\item Let $X_0 \in \cof(X)$, $\eps >0$ and $Q\colon X \to X/Y$ be the quotient map. Show that there exists $Y_0 \in \cof(X/Y)$ such that $\frac{1}{2+\eps}B_{Y_0}\subset Q(B_{X_0})$.
			\item Show that $\bar{\rho}_{X/Y}(t)\le \bar{\rho}_X(2t)$.
		\end{enumerate}
		\begin{proof}[Hint]
			For 2. Show that, by quotienting all spaces by $X_0 \cap Y$, we may assume that $Y$ is finite-dimensional. Then, use $F_\perp$, where $F$ is  a finite-dimensional subspace of $X^*$ that is $(1+\eps)$-norming for $Y$. Or have a look at Lemma 2.7 in \cite{JLPS2002}.  
		\end{proof}
		
	\end{exer}
	
	\begin{exer}
		Let $X$ be a Banach space. Show that the map $t \mapsto t^{-1}\bar{\delta}_X(t)$ is increasing on $(0,1]$.
	\end{exer}
	
	\begin{exer}
		\label{ex:moduli-weak-version}
		Let $(X,\norm{\cdot})$ be a Banach space and denote by $\cV_w(0)$ the set of weak neighborhoods of $0$. Show that for every $t>0$ and $x\in S_X$, 
		\begin{enumerate}
			\item $\bar{\rho}_X(t,x)=\inf_{V \in \cV_w(0)}\sup_{v\in V\cap S_X} \norm{x+tv}-1$.
			\item $\bar{\delta}_X(t,x)=\sup_{V \in \cV_w(0)}\inf_{v\in V\cap S_X} \norm{x+tv}-1$.
		\end{enumerate}   
	\end{exer}
	
	\begin{exer}
		\label{ex:moduli-unit-ball}
		Let $(X,\norm{\cdot})$ be a Banach space. Show that for every $t>0$ and $x\in S_X$, 
		\begin{enumerate}
			\item $\bar{\rho}_X(t,x)=\inf_{Y \in \cof(X)}\sup_{y\in Y, \norm{y}\le 1} \norm{x+ty}-1$
			\item $\bar{\delta}_X(t,x)=\sup_{Y \in \cof(X)}\inf_{y\in Y, \norm{y}\ge 1} \norm{x+ty}-1$
		\end{enumerate}
	\end{exer}
	
	\begin{exer}
		\label{ex:sequential-moduli}
		Let $X$ be a Banach space with separable dual. Show that for every $t>0$ and $x\in S_X$,
		\begin{enumerate}
			\item $\bar{\rho}_X(t,x)=\sup\big\{\limsup_{n\to \infty}(\norm{x+x_n}-1)\colon  x_n \wtoo 0,\ \norm{x_n}\le t\big\},$
			\item $\bar{\delta}_X(t,x)=\inf\big\{\liminf_{n\to \infty}(\norm{x+x_n}-1)\colon x_n \wtoo 0,\ \norm{x_n}\ge t\big\}.$
		\end{enumerate}
	\end{exer}
	
	\begin{exer}
		\label{ex:moduli-rosenthal}
		Let $X$ be a Banach space that does not contain $\ell_1$. Show that the formulas from Exercise \ref{ex:sequential-moduli} still hold.  
	\end{exer}
	
	\begin{proof}[Hint:]
		Use that every bounded subset of $X$ is weakly sequentially dense in its weak closure (cf. Appendix \ref{appendix:Banach})
	\end{proof}

	\begin{exer}
		\label{ex:moduli-lp-sums}
		Let $p\in[1,\infty]$ and $(F_n)_n$ be a sequence of finite-dimensional Banach spaces and consider the space $X:=(\sum_{n=1}^\infty F_n)_{\ell_p}$. Then, for all $x\in X$ and all weakly null net $(x_\alpha)_\alpha$ in $X$ we have
		\begin{equation*}
			\lim_\alpha \norm{x+t x_\alpha}^p = \norm{x}^p + t^p\lim_\alpha\norm{x_\alpha}^p,
		\end{equation*}
		with the usual interpretation when $p=\infty$.
	\end{exer} 
	
	\begin{exer}[$\ell_\infty$ is neither AUC nor AUC]
		\label{ex:moduli-l-infty}
		Show that:
		\begin{enumerate}
			\item $\bar{\rho}_{\ell_\infty}(t)= t$.
			\item $\bar{\delta}_{\ell_\infty}(t)= \max\{t-1,0\}$.
		\end{enumerate}
	\end{exer}
	
	\begin{proof}[Hint:]
		Use that $\ell_\infty$ contains isometric copies of $\ell_1$ and $\co$. 
	\end{proof}
	
	\begin{exer}[$L_1$ is not AUC]
		\label{ex:L1-not-AUC}
		Let $(r_n)_n$ be the Rademacher sequence on $[0,1]$, i.e. $r_n(t):=\sign \sin (2^n\pi t)$.
		\begin{enumerate}
			\item Show that for all $f\in L_1[0,1]$, $\lim_n\int_0^1 f(t)r_n(t)dt=0$ and deduce that the sequence $(fr_n)_n$ is weakly null.
			\item Show that $\norm{f(1+r_n)}_1\to_{n\to \infty} \norm{f}_1$.
			\item Show that for all $t\in[0,1)$, $\bar{\delta}_{L_1[0,1]}(t)=0$.
		\end{enumerate}
	\end{exer}
	
	\begin{proof}[Hint:]
		For $(1)$ remember that $(r_n)_n$ is an orthogonal sequence in $L_2[0,1]$ and that $L_2$ is dense in $L_1$.
	\end{proof}
	
	\begin{exer}
		\label{ex:moduli-inequalities} 
		Prove Lemma \ref{lem:AUS-powertype}
	\end{exer}

	\begin{exer}
		\label{ex:AUC*-modulus}
		Let $X$ be a Banach space. Show that for all $t\in (0,1]$,
		\begin{enumerate}[(i)]
			\item $0\le \bar{\delta}^*_X(t)\le \bar{\delta}_{X^*}(t)$,
			\item $\bar{\delta}^*_X$ is nondecreasing and $1$-Lipschitz.
		\end{enumerate}
	\end{exer}
	
	\begin{exer}
		\label{ex:moduli-weak*-opt}
		Let $(X,\norm{\cdot})$ be a Banach space and denote by $\cV_{w^*}(0)$ the set of weak$^*$ neighborhoods of $0$. Show that for every $t>0$ and $x^*\in S_{X^*}$, 
		\begin{equation*}
			\bar{\delta}^*_X(t,x^*)=\sup_{V \in \cV_{w^*}(0)}\inf_{v^*\in V\cap S_{X^*}} \norm{x^* + tv^*}-1.    
		\end{equation*}
	\end{exer}
	
	%\begin{hint}
	%Recall that if $Z$ is a weak$^*$-closed subspace of $X^*$ and $Z^\perp$ denote the pre-orthogonal of $Z$ in $X$, then $(X^*/Z)^*=Z^\perp$.
	%\end{hint}
	
	\begin{exer}\label{ex:w*-sequential-moduli}
		Let $X$ be a separable Banach space. Show that for every $t>0$ and $x^*\in S_{X^*}$
		\begin{equation*}
			\bar{\delta}^*_X(t,x^*)=\inf\big\{\liminf_{n\to \infty}(\norm{x^*+x^*_n}-1)\colon x^*_n \wstoo 0,\ \norm{x^*_n}\ge t\big\}.
		\end{equation*}
	\end{exer}

	\begin{exer}
		Let $X$ be a Banach space. Show that if $X^*$ is asymptotically uniformly smooth, then $X$ is reflexive.
	\end{exer}
	
	\begin{proof}[Hint:] 
		Assume that $X$ is not reflexive and use the fact that $\norm{\cdot}_{X^{**}}$ is AUC$^*$ and Goldstine's theorem to get a contradiction. 
	\end{proof}
	
	\begin{exer}
		\label{ex:absolute-norm} 
		Let $N$ be norm on $\bR^2$. We say that $N$ is \emph{absolutely monotonic} if $N(a_1,b_1)\le N(a_2,b_2)$ for all $a_1,a_2,b_1,b_2\in \bR$ with $\abs{a_1}\le \abs{a_2}$ and $\abs{b_1}\le \abs{b_2}$. Show the following statements.
		\begin{enumerate}
			\item If $N$ is absolute and normalized, then $\norm{(a,b)}_\infty \le N(a,b)\le \norm{(a,b)}_1$ for all $a,b\in \bR$.
			\item Every absolute norm is absolutely monotonic.  
			\item Every absolutely monotonic norm is absolute.
		\end{enumerate}
	\end{exer}
	
	\begin{exer}
		\label{ex:iterated-unc-shift}
		Let $N$ be an absolute norm on $\bR^2$. Show that 
		\begin{enumerate}
			\item $\Lambda_N$ is a norm,
			\item $\norm{\sum_{i=1}^{n} a_i e_{i}}_{\Lambda_N} \le \norm{\sum_{i=1}^{n} b_i e_{i}}_{\Lambda_N}$ for all $n\ge 1$, $(a_i)_{i=1}^n, (b_i)_{i=1}^n\subset \bR$ such that $\abs{a_i}\le \abs{b_i}$ for $1\le i \le n$,
			\item $\norm{\sum_{i=1}^{n} \vep_i a_i e_{i}}_{\Lambda_N} = \norm{\sum_{i=1}^{n} a_i e_{i}}_{\Lambda_N}$ for all $n\ge 1$, $(\vep_i)_{i=1}^n\in \{\pm 1\}^n$, $(a_i)_{i=1}^n\in \bR^n$,
			\item $\norm{\sum_{i=1}^{n} a_ie_{i+1}}_{\Lambda_N} = \norm{\sum_{i=1}^{n}a_ie_{i}}_{\Lambda_N}$ for all $n\ge 1$, $(a_i)_{i=1}^n\in \bR^n$.
		\end{enumerate}
	\end{exer}
	
	%\begin{proof}[Solution of Exercise \ref{ex:iterated-unc-shift} 4.]
	%We prove by induction on $n\in \bN$.\\Note first that for all $a\in \bR$, $\|ae_2\|_{\Lambda_N} = \norm{0e_1+ a_1 e_2}_{\Lambda_N}=N(0,\abs{a})=\|ae_1\|_{\Lambda_N}$, so the statement is true for $n=1$. Assume now that it is true for $n\in \bN$ and consider $(a_1,\ldots,a_{n+1})\in \bR^{n+1}$. Then, $$\big\|\sum_{i=1}^{n+1}a_ie_{i+1}\big\|_{\Lambda_N}=N\big(\big\|\sum_{i=1}^na_ie_{i+1}\big\|_{\Lambda_N},|a_{n+1}|\big).$$ So, we can use the induction hypothesis to deduce that $$\big\|\sum_{i=1}^{n+1}a_ie_{i+1}\big\|_{\Lambda_N}= N\big(\big\|\sum_{i=1}^na_ie_{i}\big\|_{\Lambda_N},|a_{n+1}|\big)= \|\sum_{i=1}^{n+1}a_ie_{i}\big\|_{\Lambda_N}.$$ \end{proof}
	
	\begin{exer}\label{ex:absolute-Orlicz}\,
		\begin{enumerate}
			\item Let $N$ be a normalized absolute norm on $\bR^2$. Show that the map $F_N$ defined by $F_N(t) := N(1,t)-1$ is a $1$-Lipschitz Orlicz function so that $\lim_{t\to \infty}\frac{F_N(t)}{t}=1$.
			\item Let $F$ be an Orlicz function that is Lipschitz and such that $F(t)\ge t-1$ for all $t\ge 0$. 
			\begin{enumerate}
				\item Show that the map $t\mapsto \frac{F(t)}{t}$ is nondecreasing and that $\theta:=\lim_{t\to \infty}\frac{F(t)}{t}$ exists.
				\item Show that the map $N_F$ defined by $$N_F(a,b):= \begin{cases}
					\theta\abs{b} \text{ if } a=0,\\
					\abs{a}\big(1+F(\frac{\abs{b}}{\abs{a}})\big) \text{ if } a\neq 0,
				\end{cases}$$
				is an absolute norm on $\bR^2$.
			\end{enumerate} 
		\end{enumerate}
	\end{exer}

	\begin{exer}
		\label{ex:uc-beta}
		Let $X$ be a Banach space. Show the following statements.
		\begin{enumerate}
			\item For all $t\in(0,1)$, $\delta_X(t)\le \bar{\delta}_X(t)$, where $\delta_X$ is the modulus of uniform convexity.
			\item For all $t\in(0,1)$, $\rho(t) \ge \frac12 \bar{\rho}_X(t)$, where where $\rho_X$ is the modulus of uniform smoothness.
			\item A Banach space is uniformly convex if and only if for all $t>0$, there exists $\delta>0$ such that for all $x\in B_X$ and $x_1,x_2\in B_X$ with $\norm{x_1-x_2}\ge t$ we have \begin{equation*}
				\min_{i\in\{1,2\}} \norm{\frac{x+x_i}{2}}\le 1- \delta.
			\end{equation*}
			\item Let $\beta_X(t) := \inf\{1-\min_{i\in\{1,2\}} \norm{\frac{x+x_i}{2}} \colon x, x_1, x_2\in B_X,\ \norm{x_1-x_2}\ge t\}$. 
			Show that $\delta_X$ and $\beta_X$ are equivalent.
			\item Every uniformly convex Banach space has property $(\beta)$ and $\bar{\beta}_X\ge \beta_X$.
		\end{enumerate}
	\end{exer}
	
	\begin{exer}
		Show that if a Banach space $X$ has property $(\beta)$, then every subspace and every quotient of $X$ has property $(\beta)$.    
	\end{exer}
	
	%\begin{exer}
	%Show that a blocking of an FDD is an FDD of its closed linear span.
	%\end{exer}

	%\begin{exer}\label{exer:3-space-property}
	%Let $Z$ be Banach space and $Y$ be a closed subspace of a Banach space $X$ such that every bounded operator $T \colon Y\to Z$ can be extended to a bounded operator on $X$. If $X$ and $X/Y$ embed isomorphically into $Z$, then $X$ embeds isomorphically into $Z\oplus Z$.
	%\end{exer}
	
	\begin{exer}
		\label{ex:AUF}
		Let $X$ be an infinite-dimensional subspace of $\co$. Show that for any $t\in [0,1]$, $\bar{\rho}_X(t)=0$.    
	\end{exer}
	
	\begin{exer}
		\label{ex:reflexive-biorthogonal}
		Show that any separable reflexive Banach space admits a biorthogonal system $(x_n,x_k^*)_{n,k \in \bN}$ in $X\times X^*$ such that the linear span of $(x_n)_{n=1}^\infty$ is norm dense in  $X$ and the linear span of $(x_n^*)_{n=1}^\infty$ is norm dense in $X^*$. 
	\end{exer}

	\begin{exer}
		\label{ex:LKK*-AUF}
		Let $X$ be a separable Banach space. We say that $X$ is Lipschitz weak$^*$ Kadec–Klee (in short, LKK$^*$) if there exists $c\in (0,1]$ such that its dual norm satisfies the following property: for any $x^*$ in $X^*$ and any weak$^*$ null sequence $(x^*_n)_{n\ge 1}$ in $X^*$, $\limsup\norm{x^* + x^*_n}\ge  \norm{x^*} + c\limsup\norm{x^*_n}$. 
		\begin{enumerate}
			\item Show that $X$ is Lipschitz weak$^*$ Kadec-Klee if and only if $X$ is asymptotically uniformly flat.
			\item Show that $X$ is Lipschitz weak$^*$ Kadec-Klee with constant $c=1$ if and only if $\bar{\rho}_X(1)=0$.
		\end{enumerate}
	\end{exer}
	
	%\begin{exer}
	%Show that if $X$ has a shrinking FDD $(E_n)_{n=1}^\infty$ then for all $t\ge 0$ and $x\in S_X$, we have
	%$\bar{\rho}_X(t,x)= \lim_n \bar{\rho}_X(t,x, \spa(\cup_{k\ge n}E_k))$
	%\end{exer}
	%\begin{proof}[Hint]
	%Analysing the proof of Proposition \ref{prop:pAUS-blocking} should go a long way in solving this exercise.
	%\end{proof}

	%%%%%%%%%%%%%%%%%%%%%%%%%%%%%%%%%%%%%%%%%%%%%%%%%%%%%%%%%%%%%%%%%%%%%%%%%%%%%%%%%%%%%%%%%

	\chapter[Szlenk index and asymptotic renormings]{Szlenk index and asymptotic renormings}
	\label{chapter:Szlenk}
	%Introduction to the linear asymptotic geometry of Banach spaces
	Asymptotic uniform smoothness or convexity are not isomorphic properties. So, equivalent norms may disagree on that matter. In this chapter, we develop the renorming theory for these asymptotic properties and provide isomorphic characterizations for when such renormings exist.
	One characterization is intimately connected with an ordinal index introduced by Szlenk \cite{Szlenk1968} in the late 1960s. Another characterization is expressed in terms of upper or lower estimates for branches of weakly null trees in Banach spaces. The use of weakly null trees goes back at least to the mid-1990s. Weakly null trees appeared implicitly in the asymptotic structure notion introduced by Maurey, Milman and Tomczak-Jaegermann \cite{MMTJ}, and they were used explicitly by Godefroy, Kalton and Lancien \cite{GKL2001}. What seems to be the first explicit tree-like intrinsic characterization is due to Kalton  \cite{Kalton2001}. The systematic study of weakly null trees and their branches was then taken up by Odell and Schlumprecht \cite{OdellSchlumprecht2002} and others afterward. We will focus in particular on some more recent results by Causey from \cite{Causey2018c} and \cite{Causey3.5}.

	\section{Szlenk index: basic theory}
	
	One of the many advantages of the Szlenk index is that it provides the right isomorphic invariant characterizing the existence of an equivalent asymptotically uniformly smooth norm. It will also be instrumental in providing many \emph{nonlinear} invariants. The definition of the Szlenk index is based on a particular example of a Cantor-Bendixon derivation. Such derivation indices are defined using a transfinite recursive procedure. The first index of this type is the Cantor-Bendixon index; an index which plays a central role in the study of compact topological spaces. The Szlenk index was introduced with a slightly different definition by W. Szlenk in \cite{Szlenk1968}, in order to show that there is no separable reflexive Banach space that contains an isomorphic copy of all separable reflexive Banach spaces.
	
	\begin{defi}[Szlenk index]
		\label{defi:Szlenk-index}
		Let $X$ be a Banach space and $K$ be a bounded subset of $ X^*$. For each $\eps>0$, the Szlenk derivation operation is defined by 
		\begin{equation*}
			s_\eps(K) := K\setminus\{V\subset X^*\colon V\text{ weak$^*$ open and }\mathrm{diam}(V\cap K)\le \eps\}. 
		\end{equation*}
		Given an ordinal $\xi$, the derived set of order $\xi$, denoted by $s_\eps^\xi(K)$, is defined inductively by letting 
		\begin{itemize}
			\item $s_\eps^0(K) := K$,  
			\item $s_\eps^{\xi+1}(K) := s_\eps(s^\xi_\eps(K))$,
			\item $s_\eps^{\xi}(K) := \cap_{\zeta<\xi}s^\zeta_\eps(K)$ if $\xi$ is a limit ordinal. 
		\end{itemize} 
		We also use the obvious convention $s_0^\xi(K) := K$, for any ordinal $\xi$.
		
		Note that if $K$ is a weak$^*$ compact subset of $X^*$, then $s_\eps^{\xi}(K)$ is weak$^*$ compact for all $\eps>0$ and all ordinals $\xi$. 
		
		The ordinal $\Sz(X,\eps)$ is defined as the least ordinal $\xi$ so that $s_\eps^\xi(B_{X^*})=\emptyset$, if such ordinal exists and we set $\Sz(X,\eps)=\infty$ otherwise with the obvious convention that $\infty$ is larger than any ordinal. The \emph{Szlenk index of $X$} is the ordinal number (or $\infty$) defined by
		\begin{equation*}
			\Sz(X) := \sup_{\eps>0}\Sz(X,\eps).
		\end{equation*}
	\end{defi}
	
	The Szlenk index is based on an ``extrinsic procedure'' as it is defined in terms of a certain peeling process for the unit ball of its dual. However, the peeling process refers to the weak$^*$ topology on $X^*$ induced by $X$ and is thus governed by the properties of $X$. Since the weak$^*$ topology and the norm topology coincide if and only if $X$ is finite-dimensional, it follows quickly that $\Sz(X)=1$ if $X$ is finite-dimensional. On the other hand, $\Sz(X)\ge \omega$ whenever $X$ is infinite-dimensional (see Exercise \ref{ex:Szlenk-infinite-dimension} for a quantitative version of this statement). Therefore, $X$ is finite-dimensional if and only if $\Sz(X)<\omega$, in which case necessarily $\Sz(X)=1$. It is intuitively clear that, when possible, the speed at which we can peel the entire unit ball of the dual of $X$ will depend on the geometric shape of this dual ball. Exercise \ref{ex:w*UKK-implies-Szlenk-omega} is a particularly enlightening example. 
	
	\begin{rema}
		\label{rem:Szlenk-nets}
		If we denote by $\cV_{w^*}(x^*)$ the collection of all weak$^*$ open sets in $X^*$ containing $x^*\in X^*$, then $s_\vep(K)=\{x^*\in K\colon \forall V\in\cV_{w^*}(x^*) \diam(V\cap K) > \vep\}$.
		Thus, if there is a net $(x^*_\alpha)_{\alpha\in A}$ in $K$ that converges to $x^*$ for the weak$^*$ topology and such that $\norm{x^*_\alpha-x^*}> \vep$ for all $\alpha\in A$, it is clear that $x^*\in s_\vep(K)$. On the other hand, if $x^*\in s_\vep(K)$, then for every $V\in\cV_{w^*}(x^*)$ there is $x^*_V\in K$ such that $\norm{x^*_V-x^*}> \frac{\vep}{2}$. Since the neighborhood basis $\cV_{w^*}(x^*)$ can be directed by reverse inclusion, it follows that the net $(x^*_V)_{V\in \cV_{w^*}(x^*)}$ converges to $x^*$ for the weak$^*$ topology. Therefore, for every $x^*\in s_\vep(K)$, there is a net $(x^*_\alpha)_{\alpha\in A}$ in $K$ that converges to $x^*$ for the weak$^*$ topology and such that $\norm{x^*_\alpha-x^*}> \frac{\vep}{2}$ for all $\alpha\in A$. These useful observations will be used without further explanation in the sequel.
	\end{rema}
	
	\begin{rema}
		\label{rem:stabilization}
		When $X$ is separable, the transfinite sequence $(s_\eps^\xi(B_{X^*}))_\xi$ must stabilize (on the empty set or not) before $\omega_1$ (the first uncountable ordinal). Indeed, since $X$ is separable, the weak$^*$ topology on $B_{X^*}$ admits a countable basis of open sets $(V_n)_{n\in \bN}$. Assume that $(s_\eps^\xi(B_{X^*}))_{\xi<\gamma}$ is strictly decreasing. Then, for each $\xi<\gamma$, there exists $n\in \bN$ such that $V_n\cap s_\eps^\xi(B_{X^*})\neq \emptyset$ and $V_n\cap s_\eps^{\xi+1}(B_{X^*})= \emptyset$. This induces an injection from $[0,\gamma)$ into $\bN$ and implies that $\gamma<\omega_1$. 
	\end{rema}
	The next proposition describes the behavior of the Szlenk index under taking subspaces or quotients.
	
	\begin{prop}
		\label{Szsubspace}
		Let $X$ be a Banach space and $Y$ be a closed subspace of $X$.
		\begin{enumerate}[(i)]
			\item For any $\eps>0$, $\Sz(X/Y,\eps) \le \Sz(X,\eps)$.
			\item For any $\eps>0$, $\Sz(Y,\eps) \le \Sz(X,\frac{\eps}{2})$.
		\end{enumerate}
	\end{prop}
	
	\begin{proof} 
		The proof of $(i)$ is rather immediate. Indeed, the unit ball of $(X/Y)^*$ can be identified with $B_{Y^\perp}$, which is easily seen to be a weak$^*$ closed subset of $B_{X^*}$. Moreover, the weak$^*$ topology induced by $X/Y$ on $Y^\perp$ coincides with the restriction to $Y^\perp$ of the weak$^*$ topology induced by $X$ on $X^*$. Then, the conclusion follows from a straightforward transfinite induction argument showing that for any ordinal $\xi$, $s_\eps^\xi(B_{Y^\perp}) \subset s_\eps^\xi(B_{X^*})$.
		
		For $(ii)$, denote by $Q$ the quotient map from $X^*$ onto $X^*/Y^{\perp}$. Since the unit ball of $Y^*$ can be identified with $B_{X^*/Y^{\perp}}=Q(B_{X^*})$, it is enough to show that for any weak$^*$ closed subset $L$ of $B_{X^*/Y^{\perp}}$ and any weak$^*$ closed subset $K$ of $B_{X^*}$ such that $L\subset Q(K)$, we have that $s_\eps(L) \subset Q(s_{\frac{\eps}{2}}(K))$. A straightforward transfinite induction will then conclude the proof. So, let $z\in s_\eps(L)$ and $(z_\alpha)_{\alpha \in A}$ be a net in $L$ so that $(z_\alpha)_\alpha$ tends weak$^*$ to $z$ and $\norm{z-z_\alpha}> \frac{\eps}{2}$ for all $\alpha \in A$. For each $\alpha \in A$, there exists $u_\alpha \in K$ so that $Q(u_\alpha)=z_\alpha$. Since $K$ is weak$^*$ compact, we may assume, by taking a subnet, that $(u_\alpha)_{\alpha \in A}$ is weak$^*$ converging to some $u\in K$. Now, since $Q$ is weak$^*$ continuous, $Q(u)=z$. Finally,  $\norm{u_\alpha-u}> \frac{\eps}{2}$ for all $\alpha \in A$, because $\|Q\|\le 1$. This implies that $u\in s_{\frac{\eps}{2}}(K)$.
	\end{proof}
	
	It is sometimes more convenient to study the Szlenk index for separable Banach spaces. However, many results relying on Szlenk index arguments can be extended to nonseparable Banach spaces since the Szlenk index is separably determined. The statement below about the separable determination of the Szlenk index (up to the first uncountable ordinal) was proved in \cite{Lancien1996}. The proof given below is a variation of a similar proof from \cite{Lancien2006} for the dentability index and is an alternative to the original proof given in \cite{Lancien1996}. 
	
	\begin{theo}
		\label{thm:separable-determination-Szlenk} 
		Let $X$ be a Banach space.
		\begin{enumerate}[(i)]
			\item\label{item:a-sep-det-Szlenk} Let $\alpha<\omega_1$ and $\eps>0$. If $\Sz(X,\eps)>\alpha$, then there exists a closed separable subspace $Y$ of $X$ such that $\Sz(Y,\eps/2)>\alpha$.
			\item\label{item:b-sep-det-Szlenk} If $\Sz(X)<\omega_1$, then there exists a closed separable subspace $Y$ of $X$ such that $\Sz(Y)=\Sz(X)$.
		\end{enumerate}
	\end{theo}
	
	For the proof, we need to recall a few definitions about trees. We denote by $\bN^{<\omega}$ the set of finite sequences of natural numbers (including the empty sequence $\emptyset$), which we equip with its natural ordering: for $s,t \in \bN^{<\omega}$, $s\preceq t$ if $t$ extends $s$. For $s=(s_1,\dots,s_k)$ and  $t=(t_1,\dots,t_l)$ in $\bN^{<\omega}$ , we set $s \frown t := (s_1,\dots,s_k,t_1,\dots,t_l)$. A subset $T$ of $\bN^{<\omega}$ is a \emph{tree} if $s\in T$ whenever there exists $t\in T$ with $s \preceq t$. We can define a natural derivation on a tree $T$ by setting $T'$ to be the set of nonmaximal elements of $T$. Then, we can iterate this derivation and set 
	\begin{itemize}
		\item $T^{\xi+1}=(T^\xi)'$,
		\item $T^\xi=\cap_{\zeta<\xi} T^\zeta$ if $\xi$ is a limit ordinal.
	\end{itemize}
	A tree is \emph{well-founded} if there exists an ordinal $\xi$ such that $T^\xi=\emptyset$. Then, we define the height of $T$, denoted by $ht(T)$, as the infimum of all ordinals $\xi$ such that $T^\xi \subset \{\emptyset\}$.  Note that for cardinality reasons, the height of a well-founded tree on $\bN^{<\omega}$ must be countable and that a tree is well-founded if and only if it does not contain any infinite totally ordered subset. 
	
	We now introduce a special family  $(T_\alpha)_{\alpha < \omega_1}$ of trees on $\bN^{<\omega}$ defined inductively as follows:
	\begin{itemize}
		\item $T_0 = \{\emptyset\}$.
		\item $T_{\alpha+1} = \{\emptyset\} \cup \bigcup
		\limits_{n=1}^\infty (n)^\frown T_\alpha$, where $(n)^\frown
		T_\alpha := \{(n)^\frown s\colon s \in T_\alpha\}$.
		\item $T_\alpha = \{\emptyset\} \cup \bigcup
		\limits_{n=1}^\infty (n)^\frown T_{\alpha_n}$, if $\alpha$
		is a limit ordinal, where $(\alpha_n)_{n=1}^\infty$ is a fixed enumeration of the
		ordinals less than $\alpha$.
	\end{itemize}
	
	\begin{lemm}
		\label{lem:heightoftrees}\,
		\begin{enumerate}[(i)]
			\item For all $\alpha <\omega_1$, $ht(T_\alpha) = \alpha+1$.
			\item For $\alpha <\omega_1$ and $s\in T_\alpha$, let $T_\alpha (s) := \{t \in
			\bN^{< \omega}\colon  s^\frown t \in T_\alpha\}$ and  $h_\alpha (s) :=
			ht(T_\alpha(s))$. Then, $T_\alpha (s) = T_{{h_\alpha} (s)}$.
			\item For any $\alpha < \omega_1$, there exists a bijection $\phi_\alpha \colon \bN \cup\{0\} \to
			T_\alpha$ such that for all $s,s' \in T_\alpha$, $s \preceq s'$ implies
			$\phi_\alpha^{-1}(s) \le \phi_\alpha^{-1}(s')$.
		\end{enumerate}
	\end{lemm} 
	
	\begin{proof} 
		The proofs of assertions $(i)$ and $(ii)$ are simple transfinite inductions better left as Exercise \ref{ex:heightoftrees}. Let us say a word about assertion $(iii)$.  Let $\{{\cal B}_n\}_{n=1}^\infty$ be an enumeration of the
		branches of $T_\alpha$ (a branch of a well-founded tree is the set of elements less than or equal to some maximal element of the tree). Note that these branches are finite. In order to define $\phi_\alpha$ we enumerate
		successively ${\cal B}_1,\ {\cal B}_2\setminus {\cal B}_1,\dots,\ {\cal B}_{n+1}
		\setminus \cup_{k=1}^n {\cal B}_k,\dots$ (each enumeration of
		${\cal B}_{n+1} \setminus \cup_{k=1}^n {\cal B}_k$ following the
		natural partial order on $T_\alpha$).
	\end{proof}
	
	\begin{proof}[Proof of Theorem \ref{thm:separable-determination-Szlenk}] 
		Statement $(ii)$ is an easy consequence of $(i)$. So, let us fix $\eps>0$, assume that $\Sz(X,\eps)>\alpha$ and let $x^*\in s_\eps^\alpha(B_{X^*})$. We will show that there exist a separable subspace $Y$ of $X$ and a family $(x_s^*)_{s \in T_\alpha}$ in $B_{X^*}$ such that
		\begin{itemize}
			\item $x_\emptyset^* = x^*$,
			\item  $\forall s \in (T_\alpha)',\ \forall n \in \bN,\ 
			\big\|(x_{s^\frown (n)}^* - x_s^*)_{\restriction_Y}\big\|_{Y^*} > \eps/2$,
			\item  $\forall s \in (T_\alpha)',\  ((x_{s^\frown n}^* -
			x_s^*)_{\restriction_Y})_n$ is weak$^*$ null in $Y^*$. 
		\end{itemize}
		Then, an easy induction on $\beta$ will show that for all $\beta \le \alpha$ and all $s \in T_\alpha^\beta$, $(x_s^*)_{\restriction_Y} \in s_{\eps/2}^\beta(B_{Y^*})$ and in particular that $(x^*_s)_{\restriction_Y} \in s_{\eps/2}^\alpha(B_{Y^*})$ which will conclude the proof. In order to achieve this construction, we will build, by induction on $n \in \bN$, $
		(x_{{\phi_\alpha}(n)}^*)_{n=0}^\infty$ in $B_{X^*}$ and
		$(x_n)_{n=1}^\infty$ in $B_X$ so that
		\begin{enumerate}[(a)]
			\item $x_{{\phi_\alpha}(0)}^* = x_\emptyset^* = x^*$,
			\item for all $n \in \bN,\ x_{{\phi_\alpha}(n)}^* \in
			s_\eps^{h_\alpha(\phi_\alpha(n))}(B_{X^*})$,
			\item for all $n \in \bN,\ (x_{{\phi_\alpha}(n)}^* - x_{s_n}^*)
			(x_n) > \eps/2$, where $\phi_\alpha (n) =
			s_n^\frown (k_n)$ with $k_n \in \bN$.
			\item  For all $n \ge  2$ and all $1 \le  k \le 
			n-1,\ \big|(x_{{\phi_\alpha}(n)}^* - x_{s_n}^*)(x_k)\big| \le  2^{-n}$.
		\end{enumerate}
		Assume $x_{{\phi_\alpha}(k)}^*$ for $0 \le  k \le  n-1$ and $x_k$ for $1 \le  k \le  n-1$ have been constructed and
		satisfy conditions (a) to (d). By Lemma \ref{lem:heightoftrees}, there exists $i_n \le n-1$ such that $\phi_\alpha(n) =
		\phi_\alpha(i_n)^\frown (k_n)$ with $k_n \in \bN$. By induction hypothesis
		$x_{{\phi_\alpha}(i_n)}^* \in s_\eps^{h_\alpha(\phi_\alpha(i_n))}(B_{X^*}) \subset  s_\eps^{h_\alpha(\phi_\alpha(n))+1}(B_{X^*})$. So, for any weak*-neighborhood $V$ of $x_{{\phi_\alpha}(i_n)}^*,\ \diam(V \cap
		s_\eps^{h_\alpha(\phi_\alpha(n))}(B_{X^*})) > \eps$. In particular, there exists $x_{{\phi_\alpha}(n)}^* \in
		s_\eps^{h_\alpha(\phi_\alpha(n))}(B_{X^*})$ such that $\big\| x_{{\phi_\alpha}(n)}^* - x_{{\phi_\alpha}(i_n)}^*\big\| > \eps/2$ and $\big|{(x_{{\phi_\alpha}(n)}^* - x_{{\phi_\alpha}(i_n)}^*) (x_k)}\big| \le  2^{-n}$ for all $\le  k \le  n-1$. We conclude the induction by choosing $x_n$ in $B_X$ such
		that $(x_{{\phi_\alpha}(n)}^* - x_{{\phi_\alpha}(i_n)}^*)
		(x_n) > \eps/2$.
		
		Let now $Y$ be the closed linear span of $\{x_n\}_{n=1}^\infty$. The space $Y$ and the family $(x_s^*)_{s \in T_\alpha}$ constructed by induction satisfy the properties announced at the beginning of this proof. 
	\end{proof}
	
	In the context of separable spaces, the Szlenk index detects a few topological properties of $X$ that we recall here. The first property is ``Asplundness'' which was already defined in Section \ref{sec:Asplundrigidity}. Recall that a Banach space $X$ is \emph{Asplund} if and only if every separable subspace of $X$ has a separable dual. In particular, a separable Banach space is Asplund if and only if $X^*$ is separable. The second property is weak$^*$ fragmentability: a weak$^*$ compact subset $K$ of $X^*$ is \emph{weak$^*$ fragmentable} if, for any $\eps>0$ and any nonempty weak$^*$ closed subset $L$ of $K$, there exists a weak$^*$ open subset $U$ of $X^*$ such that $L\cap U \neq \emptyset$ and $\diam(L\cap U)\le \eps$. The last property is a Baire category notion: a function between two topological spaces is called a \emph{Baire class one function} if it is a pointwise limit of continuous functions. The Baire Great Theorem provides a characterization of Baire class one functions, namely, a function $f$ from a complete metric space $M$ into a normed linear space $Y$ is of Baire class one if and only if, for every nonempty closed subset $F$ of $M$, $f_{\restriction F}$ has a point of continuity. We refer the reader to Theorem I.4.1 in \cite{DGZ1993} for a very nice proof of this classical result.
	
	\begin{theo}
		\label{thm:Szlenk-Asplund} 
		Let $X$ be a separable Banach space. The following assertions are equivalent.
		\begin{enumerate}[(i)]
			\item $X^*$ is separable.
			\item $B_{X^*}$ is weak$^*$ fragmentable.
			\item $\Sz(X)<\omega_1$.
			\item $Id_{X^*} \colon (B_{X^*},w^*) \to (B_{X^*},\norm{\cdot})$ is a Baire class one function.
		\end{enumerate}
	\end{theo}
	
	\begin{proof} 
		The implication $(i) \Rightarrow (ii)$ can be proved directly using a Baire category argument (see Exercise \ref{ex:fragmented-separable-dual}). Here we give a proof of the contrapositive that does not rely on the Baire category theorem. Assume that $B_{X^*}$ is not weak$^*$ fragmentable. Then, there exists a weak$^*$ compact nonempty subset $K$ of $B_{X^*}$ and $\eps >0$ such that for any nonempty weak$^*$ open subset $U$ of $K$, $\diam(U)>\eps$. We will construct inductively an uncountable separated family in the unit ball of $X^*$. To this end, set $U_{\emptyset} := K$ and pick, as we may, $x^*_0,x^*_1$ in $U_{\emptyset}$ and $x_{\emptyset}\in S_X$ such that $(x^*_0-x^*_1)(x_\emptyset)>\eps$. Then, consider the sets 
		$$U_0 := \{x^* \in U_{\emptyset},\ x^*(x_\emptyset)> x^*_0(x_\emptyset)-\frac{\eps}{3}\}$$ 
		and 
		$$U_1 := \{x^* \in U_{\emptyset},\ x^*(x_\emptyset)< x^*_1(x_\emptyset)+\frac{\eps}{3}\}.$$ 
		It is clear that $U_0$ and $U_1$ are nonempty weak$^*$ open subsets of $K$ and hence their diameters are larger than $\eps$. Moreover, for all $x^* \in U_{0}$ and $y^*\in U_{1}$, $(x^*-y^*)(x_\emptyset) > \frac{\eps}{3}$. Denote by $2^{<\omega}$ the set of all finite sequences in $\{0,1\}$ and let $2^\omega := \{0,1\}^\bN$. We can then build inductively a collection $(U_s)_{s\in 2^{<\omega}}$ of nonempty weak$^*$ open subsets of $K$ together with $(x^*_s)_{s\in 2^{<\omega}}$ in $K$ and $(x_s)_{s\in 2^{<\omega}}$ in $B_{X^*}$ such that $U_t \subset U_s$ whenever $t$ extends $s$ (denoted by $s\prec t$) and
		\begin{equation}
			\label{separation}
			\forall s\in 2^{<\omega},\ \forall x^* \in U_{(s,0)},\ \forall y^*\in U_{(s,1)},\ \ (x^*-y^*)(x_s) > \frac{\eps}{3}.
		\end{equation}
		For an infinite sequence $\sigma \in 2^\omega$ we pick, as we may by weak$^*$ compactness, $x^*_\sigma$ in $\bigcap_{s\prec \sigma}F_s$, where $F_s$ is the weak$^*$ closure of $U_s$. Consider now $\sigma \neq \tau \in 2^\omega$ and let $s\in 2^{<\omega}$ their greatest common ancestor (or beginning). We have, for instance, that  $x^*_\sigma \in F_{(s,0)}$ and $x^*_\tau \in F_{(s,1)}$. It follows easily from (\ref{separation}) that $\|x^*_\sigma-x^*_\tau\|\ge \frac{\eps}{3}$. Since $2^\omega$ is uncountable, we deduce that $X^*$ is not separable.

		Let us prove $(ii) \Rightarrow (iii)$. So, assume that $B_{X^*}$ is weak$^*$ fragmentable. Then, for any $\eps >0$, the transfinite sequence $(s_\eps^\xi(B_{X^*}))_\xi$ can only stabilize on the empty set, and it follows from Remark \ref{rem:stabilization} that for any $\eps>0$, $\Sz(X,\eps)<\omega_1$. Thus, $\Sz(X)=\sup_n \Sz(X,2^{-n}) <\omega_1$, since a countable supremum of countable ordinals is a countable ordinal.

		For $(iii) \Rightarrow (ii)$ Assume that $B_{X^*}$ is not weak$^*$ fragmentable. So, there exists a nonempty weak$^*$ closed subset $K$ of $B_{X^*}$ and $\eps>0$ such that every weak$^*$ open subset of $K$ has diameter greater than $\eps$. It follows that $K\subset s_\eps^{\xi}(B_{X^*})$ for any ordinal $\xi$. Therefore, $\Sz(X)=\infty$ and in particular $\Sz(X)\ge \omega_1$.

		The implication $(ii) \Rightarrow (iv)$ follows from the Baire Great Theorem after observing that the weak$^*$ fragmentation property implies that the restriction of the identity to any weak$^*$ closed subset of $B_{X^*}$ has a point of continuity. Let us detail this argument. If there is a weak$^*$ closed subset $F$ of $B_{X^*}$ such that the restriction of the identity to $F$ does not have a point of continuity, then for all $x\in F$, $\textrm{osc}(\textrm{id})(x):=\inf\{\diam(V\cap F)\colon V\in \cV_{w^*}(x)\}>0$. Therefore, $F=\cup_{n} F_n$ where $F_n:=\{x\in F\colon \forall V\in \cV_{w^*}(x), \diam(V\cap F)\ge \frac1n\}$. It is well known that oscillation functions are upper-semi-continuous and thus $F_n$ is weak$^*$ closed for all $n\in \bN$, since $F_n=\textrm{osc}(\textrm{id})^{-1}([\frac1n,\infty))\cap F$. Since $F$ is weak$^*$ compact, the Baire category theorem tells us that there must be $n_0\in \bN$ such that $F_{n_0}$ has a nonempty interior in $F$ that we denote by $\textrm{int}(F_{n_0})$. We have that $\textrm{int}(F_{n_0})=V\cap F=V\cap F_{n_0}$ for some weak$^*$ open subset $V$ of $X^*$. If  $U$ is a weak$^*$ open set of $X^*$ such that $F_{n_0}\cap U\neq \emptyset$, then $U\cap \textrm{int}(F_{n_0})=U\cap V\cap F_{n_0}\neq \emptyset$. It follows that $\diam (U\cap V \cap F)\ge \frac{1}{n_0}$ and thus that $\diam(U\cap F_{n_0})\ge \frac{1}{n_0}$. This shows that $B_{X^*}$ is not weak$^*$ fragmentable.    
		
		$(iv) \Rightarrow (i)$. Assume that there exists a sequence $(f_n)_n$ of continuous functions from $(B_{X^*},w^*)$ to $(B_{X^*},\norm{\cdot})$ that converges pointwise on $B_{X^*}$ to the identity of $X^*$. Since $(B_{X^*},w^*)$ is a metrizable compact set (as $X$ is separable), $(B_{X^*},w^*)$ is separable and its continuous image under $f_n$ is norm-separable for any $n\in \bN$. Since $B_{X^*}$ is a subset of the norm-closure of the separable set $\bigcup_{n\in \bN}f_n(B_{X^*})$, we deduce that $B_{X^*}$ is separable, from which the separability of $X^*$ follows.
	\end{proof}
	
	\begin{rema}
		The argument for the implication $(ii) \Rightarrow (iii)$ is not specific to the Szlenk derivation but is a general fact about derivations on Polish spaces and is the content of Exercise \ref{ex:Lindelof}.
	\end{rema}
	Equipped with Theorem \ref{thm:Szlenk-Asplund} and Proposition \ref{Szsubspace}, we can compute the Szlenk index of some classical Banach spaces.
	
	\begin{coro}
		\label{cor:Szlenk-extreme-cases}
		$\Sz(\ell_1)=Sz(L_1)=\Sz(L_\infty)=\Sz(\ell_\infty)=\infty$.
	\end{coro}
	
	Note that none of the spaces in Corollary \ref{cor:Szlenk-extreme-cases} are asymptotically uniformly smooth. In fact, we will see that asymptotically uniformly smooth Banach spaces must have a very small Szlenk index, namely, it cannot be bigger than the first infinite countable ordinal $\omega$. To obtain this result, we initiate a quantitative study of the Szlenk index.

	\section{Szlenk index: quantitative theory}
	\label{sec:Szlenk-quantitative}
	First of all, observe that it follows from a weak$^*$ compactness argument that if $\vep>0$ and $\alpha$ is a limit ordinal such that $s_\vep^\beta(B_{X^*})\neq \emptyset$ for all $\beta<\alpha$, then $s_\vep^\alpha(B_{X^*})\neq \emptyset$. In other words, for $\vep>0$, $Sz(X,\vep)$ is always a successor ordinal or equal to $\infty$. We start with an elementary proposition (the proof is left as Exercise \ref{ex:Szlenk-infinite-dimension}).
	
	\begin{prop} 
		Let $X$ be a Banach space. Then, $X$ is finite-dimensional if and only if $\Sz(X)=1$.
	\end{prop}
	
	Another observation is that the Szlenk index of a Banach space can only take a very specific form.  
	
	\begin{prop}
		Let $X$ be a Banach space. Either $\Sz(X)=\infty$ or there exists an ordinal $\alpha$ such that $Sz(X)=\omega^\alpha$.    
	\end{prop}
	
	The proof of this fact is detailed in Exercise \ref{ex:value-for-Szlenk}. The original idea goes back to a paper by Sersouri \cite{Sersouri1989}. In the sequel, we will only be interested in infinite-dimensional Banach spaces, for which the conditions $\Sz(X)\le \omega$ and $\Sz(X)=\omega$ are therefore equivalent.
	
	\medskip At a finer level, given $\vep>0$, it is important to estimate $\Sz(X,\vep)$ whenever it is finite, in terms of $\vep$.
	
	\begin{defi}
		A Banach space $X$ has \emph{Szlenk power type $q$}, for some $q\in [1,+\infty)$, if there exists $C>0$ such that for all $\eps \in(0,1)$, 
		\begin{equation*}
			\Sz(X,\eps)\le  C\eps^{-q}.    
		\end{equation*}
	\end{defi}
	
	To understand when a Banach space has a Szlenk index of power type we need to make a fundamental observation.
	
	\begin{prop}
		\label{prop:submult} Let $X$ be a Banach space.
		The function $\vep\mapsto \Sz(X,\vep)$ is submultiplicative. More precisely, for all $\eps,\eps'\in (0,1)$ we have
		\begin{equation*}
			\Sz(X,\eps\eps')\le 
			\Sz(X,\eps)\cdot \Sz(X,\eps').
		\end{equation*}
	\end{prop}
	
	\begin{proof} 
		Let $\eps\in(0,1)$ and $\eps'\in(0,1)$. It is enough to show
		(omitting the obvious case when  $\Sz(X,\eps)=\infty$) that for any
		ordinal $\alpha$,
		$$ s_{\eps\eps'}^{\Sz(X,\eps)\cdot \alpha}(B_{X^*})\subset s_{\eps'}^\alpha(B_{X^*}).$$
		This will be achieved with a transfinite induction on $\alpha$. The statement
		is clearly true for $\alpha=0$ and it passes to limit ordinals since for all ordinals $\gamma$ and limit ordinals $\alpha$, $\gamma \cdot \alpha=\sup_{\beta<\alpha}\gamma \cdot \beta$. So, let us
		assume it is true for some ordinal $\alpha$. Let now $x^*\in B_{X^*}\setminus
		s_{\eps'}^{\alpha+1}(B_{X^*})$. We need to show that $x^*\notin
		s_{\eps\eps'}^{\Sz(X,\eps)\cdot (\alpha+1)}(B_{X^*})$, so we may assume that
		$x^*\in s_{\eps'}^\alpha(B_{X^*})$. Then, there is a weak$^*$ open set $V$
		containing $x^*$ such that $\diam(V\cap s_{\eps'}^\alpha(B_{X^*}))\le \eps'$ and
		by the induction hypothesis $\diam(V\cap s_{\eps\eps'}^{\Sz(X,\eps)\cdot \alpha}(B_{X^*}))\le
		\eps'$. But every set $C$ with diameter at most 
		$\eps'$ satisfies $s_{\eps\eps'}^{\Sz(X,\eps)}(C)=\emptyset$ (see Exercise \ref{ex:w*UKK-implies-Szlenk-omega}). So
		$x^*\notin s_{\eps\eps'}^{\Sz(X,\eps)}(s_{\eps\eps'}^{\Sz(X,\eps)\cdot \alpha}(B_{X^*})) = s_{\eps\eps'}^{\Sz(X,\eps)\cdot \alpha + \Sz(X,\eps)}(B_{X^*})=s_{\eps\eps'}^{\Sz(X,\eps)\cdot(\alpha + 1)}(B_{X^*})$, where for the last equality we use the fact that ordinal multiplication is distributive from the left over ordinal addition.
	\end{proof}
	
	An immediate corollary of the submultiplicativity of the Szlenk index follows from a classical fact for $\bN$-valued submultiplicative functions.
	\begin{coro}
		\label{coro:szlenk-power-type}
		If $X$ is a Banach space with $\Sz(X)=\omega$, then there exists $q\in [1,+\infty)$ such that $X$ has Szlenk power type $q$.
	\end{coro}
	
	It is clear that if $X$ has Szlenk index of power type $q$, then it has Szlenk index of power type $p$ for every $p>q$. If we let 
	\begin{equation*}
		\bar{q}_X := \inf\{q\in [1,\infty)\colon \exists C>0,\ \forall \eps\in (0,1),\ \Sz(X,\eps)\le  C\eps^{-q}\},
	\end{equation*}
	then it follows from Proposition \ref{prop:submult} that $\bar{q}_X<\infty$. The parameter $\bar{q}_X$ is closely related to the notion of $q$-summable Szlenk index.
	
	\begin{defi}
		Let $q\in [1,\infty)$. A Banach space $X$ is said to have a \emph{$q$-summable Szlenk index} if there exists $c>0$ so that for all $\eps_1,\dots,\eps_n>0$ the condition
		\begin{equation*}
			s_{\eps_1}\ldots s_{\eps_n}(B_{X^*}):= s_{\eps_1}(s_{\eps_{2}}(\dots s_{\eps_{n-1}}(s_{\eps_n}(B_{X^*}))\ldots ))\neq\emptyset
		\end{equation*}
		implies that $\eps_1^q+\dots+\eps_n^q\le  c$.\\
		Instead of $1$-summable Szlenk index, we shall simply say \emph{summable Szlenk index}. 
	\end{defi}
	
	Clearly, a Banach space with $q$-summable Szlenk index has Szlenk power type $q$. The next proposition is a partial converse. 
	
	\begin{prop}
		\label{prop:Szlenk-power-summable}
		Let $X$ be a Banach space such that $\Sz(X)= \omega$. Then, for any $s \in (\bar{q}_X,\infty)$, $X$ has an $s$-summable Szlenk index.
	\end{prop}
	
	\begin{proof} 
		Since $\Sz(X)\le \omega$ we have that $\bar{q}_X<\infty$. Let $s>\bar{q}_X$ and pick $q\in (\bar{q}_X,s)$. Thus, there exists $C>0$ such that $\Sz(X,\eps)\le  C\eps^{-q}$ for all $\eps \in (0,1)$. Assume now that $ s_{\eps_1}\ldots s_{\eps_n}(B_{X^*})\neq \emptyset$ with $\eps_1,\ldots,\eps_n \in (0,1)$. For each $k\in \bN$, let $A_k := \{j \colon 2^{-k}< \eps_j \le 2^{-k+1}\}$. Since $\Sz(X,2^{-k})\le  C2^{qk}$, we have that the cardinality of $A_k$ does not exceed $C2^{qk}$ and therefore that $\sum_{j \in A_k}\eps_j^s \le C2^s2^{(q-s)k}$. Summing over $k\in \bN$, we get that $\sum_{j=1}^n\eps_j^s \le C2^s(1-2^{q-s})^{-1}$. This concludes the proof. 
	\end{proof}

	In this book, the main applications of Szlenk index techniques come from the remarkable links between the existence of AUS renormings and the value of the Szlenk index. Recall that a Banach space is AUS if and only if its dual is weak$^*$ AUC. The first statement is elementary.
	
	\begin{prop}
		\label{prop:Sz-omega} Let $X$ be an AUS Banach space and $\eps_1,\dots,\eps_n \in (0,1)$. Let $\bar{\delta}^*_X$ be the AUC$^*$-modulus of the dual norm on $X^*$.
		\begin{enumerate}[(i)]
			\item Assume that 
			$$\prod_{i=1}^n\Big(1+\delta^*_X\big(\frac{\eps_i}{8}\big)\Big)^{-1} \le \frac12.$$
			Then, $s_{\eps_1}\dots s_{\eps_n}(B_{X^*})=\emptyset$.
			\item In particular, there exists a universal constant $C\ge 1$ such that for any AUS Banach space $X$ and any $\eps\in(0,1)$,  
			\begin{equation*}
				\Sz(X,\eps)\le C\Big(\bar{\delta}^*_X\big(\frac{\eps}{C}\big)\Big)^{-1},
			\end{equation*}
			\item If $X$ is AUS (and infinite-dimensional), then $\Sz(X)= \omega$.
		\end{enumerate}
	\end{prop}
	
	\begin{proof} 
		$(i)$ Let $\eps\in (0,1)$. It is easy to see (using  Remark \ref{rem:Szlenk-nets}, a rescaling argument and Proposition \ref{prop:w*-modulus-nets}) that  $s_{\eps}(B_{X^*})\subset \lambda B_{X^*}$, where $\lambda := (1+\bar{\delta}^*_X(\frac{\eps}{2}))^{-1}$. Then, an easy induction combined with a homogeneity argument shows that
		$$s_{\eps_1}\ldots s_{\eps_n}(B_{X^*})\subset \prod_{i=1}^n\Big(1+\delta^*_X\big(\frac{\eps_i}{2}\big)\Big)^{-1}B_X.$$
		So, assuming $\prod_{i=1}^n\Big(1+\delta^*_X\big(\frac{\eps_i}{8}\big)\Big)^{-1} \le \frac12$, we get that $s_{\frac{\eps_1}{4}}\ldots s_{\frac{\eps_n}{4}}(B_{X^*})\subset \frac12B_{X^*}$. Since $B_{X^*}\setminus \frac12 B_{X^*}$ contains a translate of $\frac14 B_{X^*}$, we deduce that $s_{\frac{\eps_1}{4}}\dots s_{\frac{\eps_n}{4}}(\frac14B_{X^*})=\emptyset$ and by homogeneity that $s_{\eps_1}\dots s_{\eps_n}(B_{X^*})=\emptyset$.
		
		$(ii)$ It follows from $(i)$ that if $n\in \N$ is such that $n\ge \ln 2 (\ln (1+\bar{\delta}^*_X(\frac{\eps}{8})))^{-1}$, then $s_{\eps}^n(B_{X^*})=\emptyset$. The conclusion now follows from elementary calculus.
		
		$(iii)$ is an immediate consequence of $(ii)$. 
		
	\end{proof}
	
	\begin{rema}
		\label{rem:q-AUC*->q-summable-Szlenk}
		It follows from item $(i)$ of the previous proposition that If $X$ is $p$-AUS, then $X$ has a $q$-summable Szlenk index, where $q$ is the conjugate exponent of $p$. We leave the details as Exercise \ref{ex:q-AUC*->q-summable-Szlenk}.
	\end{rema}
	
	As an application, we have the following estimate. We refer to Exercise \ref{ex:Szlenk-ell_p} for a more precise version. 
	\begin{prop}
		\label{prop:Sz-lp} 
		Let $p\in (1,\infty)$ and $q$ be its conjugate exponent. There exists $A\ge 1$ such that for any sequence $(F_n)_{n=1}^\infty$ of finite-dimensional normed spaces and any $\eps \in (0,1)$, 
		\begin{equation*}
			\eps^{-q}\le \Sz(X,\eps) \le A\eps^{-q},
		\end{equation*}
		where $X=(\sum_{n=1}^\infty F_n)_{\ell_p}$. 
	\end{prop}
	
	\begin{proof} 
		Since $X^*$ is $q$-AUC$^*$ with $\bar{\delta}^*_X(t)=(1+t^q)^{1/q}-1$,  the right-hand side inequality is a particular case of Proposition \ref{prop:Sz-omega}.\\
		We now turn to the other inequality. For each $n \in \bN$, we fix $e_n^*\in S_{F_n^*}$. Note that the sequence $(e^*_n)_n$ is weak$^*$ null in $X^*=(\sum_{n=1}^\infty F_n^*)_{\ell_q}$. Now, let $\eps \in (0,1)$ and $x^*\in X^*$ such that $\norm{x^*}^q<1-\eps^q$. So, $\norm{x^*}^q<1-\delta^q$, for some $\delta>\eps$. It is clear that $\limsup_n\|x^*+\delta e_n^*\|<1$ and it immediately follows from Remark \ref{rem:Szlenk-nets} that $x^*\in s_\eps(B_{X^*})$. Since $s_\eps(B_{X^*})$ is weak$^*$ closed, we deduce that $(1-\eps^q)^{1/q}B_{X^*}\subset s_\eps(B_{X^*})$. Iterating the argument, we infer that for all $n\in \bN$, $(1-n\eps^q)^{1/q}B_{X^*}\subset s_\eps^n(B_{X^*})$. Therefore, $\Sz(X,\eps)\ge \eps^{-q}$.
	\end{proof}
	
	The following corollary is now an immediate consequence of Proposition \ref{prop:Sz-omega}, Proposition \ref{prop:Sz-lp}, and the permanence properties of the Szlenk index.
	\begin{coro}
		Let $p\in(1,\infty)$ and $q$ be the conjugate exponent of $p$. We have 
		\begin{equation*}
			\bar{q}_{\ell_p} := \min\{s\in [1,\infty)\colon \exists C>0,\ \forall \eps\in (0,1),\ \Sz(\ell_p,\eps)\le C\eps^{-s}\}=q,
		\end{equation*}
		and 
		\begin{equation*}
			\bar{q}_{L_p} := \min\{s\in [1,\infty)\colon \exists C>0,\ \forall \eps\in (0,1),\ \Sz(L_p,\eps)\le C\eps^{-s}\}= \max\{2,q\}.
		\end{equation*}
		In particular, $\Sz(L_p)=\Sz(\ell_p)=\omega$.
	\end{coro}
	
	The fact that the converse of Proposition \ref{prop:Sz-omega} holds is an important theorem and will be a highlight of this chapter. It was first obtained by Knaust, Odell and Schlumprecht \cite{KOS1999} in the separable setting and later by Raja \cite{Raja2010} in the nonseparable case. The following precise quantitative version is taken from \cite{GKL2001}. We state it first here, but we will come back to its proof and give other related renorming results in the remainder of this chapter. Recall that $\bar{q}_X$ is the infimum of all $q\in[1,\infty)$ such that $X$ has Szlenk index of power type $q$.
	
	\begin{theo}
		\label{thm:Sz-omega->AUS} Let $X$ be a Banach space with $\Sz(X)\le \omega$.\\
		For any $p\in (1,\bar{p}_X)$, where $\bar{p}_X$ is the conjugate exponent of $\bar{q}_X$, $X$ admits an equivalent $p$-AUS norm.
	\end{theo}
	
	\section{\texorpdfstring{Summable Szlenk index and the $\sA$-game: a warmup.}{Summable Szlenk index and the -game: a warmup.}}
	\label{sec:summable-A}
	
	In this section, we begin the study of certain combinatorial games and their connections with the Szlenk index. For Banach spaces with separable duals, these games can be formulated in terms of properties of countably branching weakly null trees and their branches. However, in order to be able to handle arbitrary Banach spaces, we will need to consider different types of trees and adjust what we mean by the branches being weakly null. Thus, we start this section with the definition of general weakly null trees in a Banach space and the corresponding notation. 
	
	Let $D$ be an arbitrary set and let:
	\begin{itemize}
		\item $\emptyset$, be the empty sequence,
		\item $D^{\le n} := \{\emptyset\}\cup \cup_{i=1}^n D^i$, be the set of all finite sequences in $D$ of length at most $n$,
		\item $D^{<\omega} := \cup_{n=1}^\infty D^{\le n}$, be the set of all finite sequences in $D$,
		\item $D^\omega$, be the set of all infinite sequences in $D$,
		\item $D^{\le \omega} := D^{<\omega}\cup D^\omega$, be the set of all (finite or infinite) sequences in $D$.
	\end{itemize}
	For $s,t\in D^{<\omega}$, we let $s\frown t$ be the concatenation of $s$ with $t$ and we simply write $s\frown a$ instead of $s\frown (a)$ whenever $a\in D$. We denote by  $\abs{t}$ the length of $t$.  For $0\le i\le \abs{t}$, we let $t_{\restriction i}$ be the initial segment of $t$ having length $i$, where $t_{\restriction 0} := \emptyset$. If $s\in D^{<\omega}$, we denote by  $s\prec t$ the relation that $s$ is a proper initial segment of $t$, i.e. $s=t_{\restriction i}$ for some $0\le i<\abs{t}$. We write $s \preceq t$ if $s\prec t$ or $s=t$. For  $s=(s_1,\dots,s_n) \in  D^{<\omega}\setminus \{\emptyset\}$, we set $s^-$ the immediate predecessor of $s$ for $\preceq$, that is  $s^-:=(s_1,\dots,s_{n-1})$ if $n \ge 2$ and $s^-=\emptyset$ if $n=1$. It is well known that  $(D^{<\omega}, \prec)$ is a tree in the set-theoretic sense.
	
	If $(D, \le)$ is a directed set and $(x_t)_{t\in D^{<\omega}}\subset X$, we say that $(x_t)_{t\in D^{<\omega}}$ is a \emph{weakly null tree} in $X$, provided that for each $t\in D^{<\omega}$, $(x_{t\frown s})_{s\in D}$ is a weakly null net. A function $\varphi \colon D^{<\omega}\to D^{<\omega}$ is said to be a \emph{pruning} provided that 
	\begin{enumerate}[(i)]
		\item $\varphi$ preserves the tree ordering, i.e. if $s\prec t$, then $\varphi(s)\prec \varphi(t)$,
		\item $\varphi$ preserves the length, i.e. $\abs{\varphi(t)} = \abs{t}$ for all $t\in D^{<\omega}$,
		\item[] and moreover 
		\item[(ii')] if $\varphi((t_1, \dots, t_k)) = (s_1, \dots, s_k)$, then $t_i\le s_i$ for all $i\in \{1,\dots, k\}$. 
	\end{enumerate} 
	We define prunings $\varphi\colon D^{\le n}\to D^{\le n}$ similarly. As it can be seen in the next two lemmas, one reason one might want to prune a tree is to stabilize maps that are defined on the leaves of trees.  
	
	\begin{lemm}
		\label{lem:pruning1}
		Let $(D,\le)$ be a directed set, $F$ a finite set, $n\in \bN$ and $f \colon D^n \to F$ a function. Then, there exists a pruning $\varphi \colon D^{\le n}\to D^{\le n}$ such that $f \circ \varphi$ is constant on $D^n$.
	\end{lemm}
	
	\begin{proof}
		We work by induction on $n$. So, assume first that $n=1$. For $x\in F$, let $I_x := \{t\in D \colon f(t)=x\}$ and observe that $\cup_{x\in F}I_x = D$. Since $D$ is cofinal in itself and $F$ is finite, there must exist some $x\in F$ such that $I_x$ is cofinal in $D$. This means that for any $t\in D$, there exists $s_t\in I_x$ such that $t\le s_t$. Letting $\varphi(t) := s_t$ and $\varphi(\emptyset) := \emptyset$, we have that $\varphi \colon D^{\le 1}\to D^{\le 1}$ is a pruning and $f \circ \varphi=x$ on $D$ and that $f \circ \varphi$ is constant on $D^n$.
		
		Assume now that the result holds for $n\in \bN$ and consider $f\colon D^{n+1}\to F$. For each $t\in D$, define $f_t\colon D^n\to F$ by $f_t(t_1, \dots, t_n) := f(t, t_1, \dots, t_n)$.  By the inductive hypothesis, there exists a pruning $\varphi_t\colon D^{\le n}\to D^{\le n}$ and $x_t\in F$ such that $f_t\circ \varphi_t= x_t$ on $D^n$. Define $g\colon D^{1}\to F$ by $g((t))=x_t$. By the base case, there exists a pruning $\psi \colon D^{\le 1}\to D^{\le 1}$ such that $g \circ \psi$ is constant on $D$.  Define $\varphi \colon D^{\le n+1}\to D^{\le n+1}$ by, $\varphi(\emptyset) := \emptyset$, $\varphi(t) := \psi(t)$ and $\varphi(t, t_1, \dots, t_k) := \psi(t)\smallfrown \varphi_{\psi(t)}(t_1, \dots, t_k)$, for $1\le k\le n$. It is straightforward to verify that $\varphi$ is a pruning and that $f \circ \varphi$  is constant on $D^{n+1}$.  
	\end{proof}
	
	The following corollary is then immediate.
	
	\begin{coro}
		\label{cor:pruning2} 
		Let $(D,\le)$ be a directed set, $(K,d)$ a totally bounded metric space, $n\in \bN$ and $f\colon D^n \to K$ a function. Then, for any $\eps>0$, there exist a pruning $\varphi \colon D^{\le n}\to D^{\le n}$ and a subset $B$ of $K$ with $\diam(B)\le \vep$ such that $(f \circ \varphi)(t) \in B$ for all $t\in D^n$.
	\end{coro}
	
	As previously mentioned, if we restrict ourselves to Banach spaces with separable duals, then we can get by with the classical notion of weakly null trees where $(D,\le)$ is $(\bN,\le)$. For arbitrary Banach spaces, we will consider weakly null trees where $D = \cN_{w}(0)$ is the set of weak neighborhoods of $0$ in $X$ directed by reverse inclusion. We now define a particular type of two-player game. This game is an example of a mathematical game with perfect information where two players play alternately. We refer the reader, for instance, to the textbook by A. Kechris \cite{Kechris1995}, or \cite{Godefroy-Baire}, for more details on the theory of mathematical games. 
	For $c>0$, $\alpha\in\bN\cup\{\omega\}$ and a function $R \colon X^\alpha \to [0,\infty)$ called the \emph{referee function}, the $G(c, R,\alpha)$-game is a (finite if $\alpha\in\bN$ or infinite if $\alpha=\omega$) game where:
	\begin{itemize}
		\item Player A chooses $U_1\in \cN_{w}(0)$,
		\item Player B chooses $x_1\in U_1\cap B_X$,
		\item Player A chooses $U_2\in \cN_{w}(0)$,
		\item Player B chooses $x_2\in U_2\cap B_X$,
		\item ...
		\item and the game continues in this way until $(x_i)_{i=1}^\alpha$ has been chosen.
	\end{itemize}  
	At the end of the game, we compute the value of the $\alpha$-tuple $(x_i)_{i=1}^\alpha$ obtained by Player B under the referee function. Player A wins the $G(c, R, \alpha)$-game if $R((x_i)_{i=1}^\alpha) \le c$ and Player B is declared  the winner otherwise. We will consider the $G(c, R, \alpha)$ game for various choices of referee function $R$.   
	Let us detail the notion of strategies and winning strategies in these games. First, we define the notions of strategies and winning strategies for Player A. For a Banach space $X$, a \emph{strategy for Player A} is a function $\chi\colon B_X^{<\alpha}\to \cN_{w}(0)$. If $\chi$ is a strategy for Player A, then we say that $(x_i)_{i=1}^\alpha \subset B_X$ is $\chi$-\emph{admissible} if $x_{j+1}\in \chi((x_i)_{i=1}^{j})$ for all $0 \le j <\alpha$. We say the strategy is a \emph{winning strategy for Player A} if any sequence admissible with respect to it satisfies the winning condition of the game for Player A. 
	%Note that in theory, the choice of $U_{j+1}$ by Player A should depend on $(U_1,x_1,\ldots, U_{j},x_{j})$ in his strategy, but the outcome of the game only depends on the sequence $(x_i)_{i=1}^\alpha$. So, it allows us to consider strategies as function of $(x_1,\ldots,x_{j})$ only for the choice of $U_{j+1}$. 
	Winning strategies for Player B are defined similarly. 
	When one of the two players has a winning strategy, we say that the game is \emph{determined}. In the case of a game with finitely many rounds, i.e. in the case $\alpha\in\bN$, it follows from an elementary transfinite argument that the game $G(c, R, \alpha)$ is determined (see \cite[Proposition 6.2]{Godefroy-Baire}). We now consider a finite $G(c, R, \alpha)$-game for a particular choice of referee function $R$.
	
	\begin{defi}[The class $\sA_p$]\label{def:game-Ap}
		Let $p\in (1,\infty]$, $c>0$ and $n\in\bN$. The game $A(c,p,n)$ in a Banach space $X$ is the game $G(c, N_{q,n}, n)$ where $\frac1p +\frac1q = 1$ and
		\begin{equation}
			N_{q,n}((x_i)_{i=1}^n) := \inf\Big\{c\in (0,\infty] \colon \forall a=(a_i)_{i=1}^n \in \bR^n\ \Big\|\sum_{i=1}^n a_ix_i\Big\|\le c\norm{a}_p\Big\},
		\end{equation}
		We denote by $\textsf{a}_{p,n}(X)$ the infimum of all $c>0$ such that Player A has a winning strategy in the $A(c,p,n)$ game and we let $\textsf{a}_p(X) := \sup_n \textsf{a}_{p,n}(X)$. Note that $\textsf{a}_p(X)$ is the infimum of all $c>0$ such that for each $n\in\bN$, Player A has a winning strategy in the $A(c,p,n)$-game if such a $c$ exists and $\textsf{a}_p(X)=\infty$ otherwise. 
		Finally, we denote by $\sA_p$ the class of all Banach spaces $X$ such that $\textsf{a}_p(X)<\infty$. 
	\end{defi}
	
	Membership in the class $\sA_p$ can be conveniently reformulated in terms of norm estimates on branches of finite but arbitrarily large weakly null trees indexed by tuples of elements in a weak neighborhood basis of $0$. We omit the proof of the proposition below as we will give a complete proof of a similar result for an infinite game in the next section (cf. Proposition \ref{prop:Tp-game-trees}). 
	
	\begin{prop}
		\label{prop:Ap-game-trees}
		Let $1<p\le \infty$ and $X$ be a Banach space. The following assertions are equivalent:  
		\begin{enumerate}[(i)]
			\item $X\in \textsf{\emph{A}}_p$. 
			\item There exists a constant $c>0$ such that for any weak neighborhood basis $D$ at $0$ in $X$, any $n\in \bN$ and any weakly null tree $(x_t)_{t\in D^{\le n}}\subset B_X$, there exists $\tau\in D^n$ such that for all $(a_i)_{i=1}^n \in \bR^n$,
			\begin{equation*}
				\Big\| \sum_{i=1}^n a_ix_{ \tau_{\restriction i} } \Big\| \le c \norm{(a_i)_{i=1}^n}_p.
			\end{equation*}
			\item There exists a constant $c>0$ such that for any weak neighborhood basis $D$ at $0$ in $X$, any $n\in \bN$ and any weakly null tree $(x_t)_{t\in D^{\le n}}\subset S_X$, there exists $\tau\in D^n$ such that for all $(a_i)_{i=1}^n \in \bR^n$,
			\begin{equation*}
				\Big\| \sum_{i=1}^n a_ix_{ \tau_{\restriction i} } \Big\| \le c \norm{(a_i)_{i=1}^n}_p.
			\end{equation*}
		\end{enumerate}
	\end{prop}
	
	\begin{rema} Membership of a Banach space $X$ in the class $\sA_p$ can also be characterized in terms of properties of its asymptotic structure (see Appendix \ref{sec:asymptotic-structure} for the definitions). We shall not use it, but we just mention it, as it is easy to see that $X \in \sA_p$ if and only if there exists a constant $C\ge 1$ such that for all $k\in \N$ and all $E \in \{X\}_k$ with normalized basis $(e_i)_{i=1}^k$, the basis $(e_i)_{i=1}^k$ satisfies upper $\ell_p$-estimates with constant $C$.
	\end{rema}
	
	The next result describes a crucial link between the Szlenk index and the class $\sA_p$. This general and sharp version is due to R.M. Causey \cite{Causey3.5}.
	
	\begin{prop}
		\label{prop:Sz-q-summable->Ap}
		Let $p\in(1,\infty]$ and $q\in[1,\infty)$ be conjugate exponents. If a Banach space $X$ has a $q$-summable Szlenk index, then $X\in \sA_p$.
	\end{prop}
	
	\begin{proof}
		%Assume that $X$ has a $q$-summable Szlenk index and let $M>0$ be such that $\sum_{i=1}^n\eps_i^q \le M$, whenever $s_{\eps_{1}}\ldots s_{\eps_n}(B_{X^*})\neq \emptyset$. We will show that $\textsf{a}_p(X)\le M$, i.e. $X\in \sA_p$. 
		Assume that $X\notin \sA_p$. We will show that $X$ does not have a $q$-summable Szlenk index, i.e. for every $c>0$ we need to find an integer $n\ge 1$ and $\vep_1,\dots,\vep_n \ge 0$ such that $s_{\eps_{1}}\dots s_{\eps_n}(B_{X^*})\neq \emptyset$ and $\sum_{i=1}^n\vep_i^q\ge c$.
		
		So, let $c>0$. By Proposition \ref{prop:Ap-game-trees}, there exist a weak neighborhood basis $D$ at $0$ in $X$, $n\in \bN$ and a weakly null tree $(x_t)_{t\in D^{\le n}}\subset B_X$ such that for all $\tau \in D^n$, there exists $a^\tau\in B_{\ell_p^n}$ satisfying $\norm{ \sum_{i=1}^n a^\tau_i x_{\tau_{\restriction i}} } > c$. For $\tau \in D^n$, pick $y^*_\tau \in S_{X^*}$ such that 
		$$\Big\langle y^*_\tau,\sum_{i=1}^n a^\tau_i x_{\tau_{\restriction i}}\Big\rangle > c.$$
		We fix $\delta>0$ and define $f \colon D^n\to B_{\ell_\infty^n}$ by $f(\tau) := (y^*_{\tau}(x_{\tau_{\restriction i}}))_{i=1}^n$. By Corollary \ref{cor:pruning2} there exist a pruning $\varphi \colon D^{\le n}\to D^{\le n}$ and $b\in B_{\ell_\infty^n}$ such that for all $\tau \in D^n$ and $1\le i \le n$, 
		\begin{equation}\label{eq:pruning}
			\Big|y^*_{\varphi(\tau)}(x_{\varphi(\tau)_{\restriction i} }) - b_i\Big|<\delta.
		\end{equation}
		Now set $v^*_\tau := y^*_{\varphi(\tau)}$, for $\tau \in D^n$ and $u_t := x_{\varphi(t)}$, for $t\in D^{\le n}$. Note that $(u_t)_{t\in D^{\le n}}$ is still weakly null and that $u_{t_{\restriction i}} = x_{\varphi(t_{\restriction i})} = x_{\varphi(t)_{\restriction i} }$(since $\varphi$ preserves the tree ordering and the length). 
		
		\begin{claim}
			For $1\le i  \le n$, let $\eps_i :=\max\{0,\abs{b_i}-2\delta\}$. Then, for all $0\le j \le n$ and all $t\in D^{n-j}$, there exists $x^*_t \in s_{\eps_{n-j+1}}\dots s_{\eps_n}(B_{X^*})$ and if $j<n$, this $x^*_t$ can be chosen such that $\abs{x^*_t(u_{t_{\restriction i}}) - b_i}\le \delta$ for all $1\le i\le n-j$.    
		\end{claim}
		
		Assuming the claim has been proven, we can conclude as follows. Since $x^*_\emptyset \in s_{\eps_{1}}\dots s_{\eps_n}(B_{X^*})$, we have that $s_{\eps_{1}}\dots s_{\eps_n}(B_{X^*})\neq \emptyset$ and we need to estimate $\sum_{i=1}^n\eps_i^q$. Letting $I := \{1\le i\le n \colon \abs{b_i}>2\delta\}$ and observing that $\norm{a^\tau}_{\ell_1^n}\le n \norm{a^\tau}_{\ell_p^n} \le n$, we have that for any $\tau\in D^n$, 
		\begin{align*}
			c < \Big\langle v^*_\tau,\sum_{i=1}^n a_i^{\varphi(\tau)} u_{t_{\restriction i}}\Big\rangle & = \sum_{i=1}^n a_i^{\varphi(\tau)} ( v_\tau^*(u_{t_{\restriction i}}) - b_i) + \sum_{i=1}^n a_i^{\varphi(\tau)} b_i\\ 
			&\le \delta \sum_{i=1}^n \abs{a_i^{\varphi(\tau)}} +  \sum_{i=1}^n \abs{a_i^{\varphi(\tau)}} \abs{b_i} \\
			& \le 3\delta n + \sum_{i\in I}\abs{a_i^{\varphi(\tau)}}\abs{b_i} \le 5\delta n +\sum_{i\in I}\abs{a_i^{\varphi(\tau)}}\eps_i \\ & \le 5\delta n + \Big(\sum_{i\in I}\abs{a_i^{\varphi(\tau)}}^p\Big)^{1/p}\Big(\sum_{i\in I}\eps_i^q\Big)^{1/q}\\
			& \le 5\delta n + \Big(\sum_{i =1}^n\eps_i^q\Big)^{1/q}.
		\end{align*} 
		Since $\delta>0$ was arbitrary, we conclude that $\sum_{i =1}^n\eps_i^q \ge c^q$. Since $c>0$ was also arbitrary, we have shown that  $X$ does not have a $q$-summable Szlenk index. 
		
		\smallskip
		It remains to prove the claim by induction on $0\le j \le n$. The base case is settled by the convention $s_{\eps_{n+1}}s_{\eps_n}(B_{X^*})=B_{X^*}$ and by letting $x^*_\tau := v^*_\tau$, $\tau \in D^n$. Assume now the claim is proved for some $0\le j<n$. By the inductive hypothesis for each $t\in D^{n-j-1}$ and $U\in D$, there exists $x^*_{t\frown (U)}\in s_{\eps_{n-j+1}}\ldots s_{\eps_n}(B_{X^*})$ such that for each $1\le i< n-j$, $\abs{x^*_{t\frown (U)}(u_{t_{\restriction i}})-b_i}\le \delta$ and $\abs{x^*_{t\frown (U)}(u_{t\smallfrown (U)})-b_{n-j}} \le \delta$. If $\eps_{n-j}=0$, then let $x^*_t := x^*_{t\frown (U)}$ for some $U\in D$ and note that the conclusions are satisfied, since $s_{\eps_{n-j}}s_{\eps_{n-j+1}}\dots s_{\eps_n}(B_{X^*}) = s_{\eps_{n-j+1}}\dots s_{\eps_n}(B_{X^*})$. Assume now that $\eps_{n-j}>0$ and choose $x^*_t$ to be a weak$^*$ cluster point of $(x^*_{t\frown U})_{U\in D}$. Then, $\abs{x^*_t(u_{t_{\restriction i}})-b_i} \le \delta$ for each $1\le i< n-j$. Fix now a weak$^*$ neighborhood $V$ of $x^*_t$ in $X^*$. Since $x^*_t$ is a weak$^*$ cluster point of $(x^*_{t\frown (U)})_{U\in D}$ and $(u_{t\smallfrown (U)})_{U\in D}$ is weakly null, there exists $U\in D$ such that $x^*_{t\frown (U)}\in V$ and $\abs{x^*_t(u_{t\smallfrown (U)})}<\delta$. Recall that by the induction hypothesis, $x^*_{t\frown (U)}\in s_{\eps_{n-j+1}}\dots s_{\eps_n}(B_{X^*})$,  and since the derived sets are weak$^*$ closed, $x^*_t\in s_{\eps_{n-j+1}}\dots s_{\eps_n}(B_{X^*})$. Thus, 
		\[\text{diam}(V\cap s_{\eps_{n-j+1}}\ldots s_{\eps_n}(B_{X^*})) \ge |(x^*_{t\frown (U)}-x^*_t)(u_{t\smallfrown (U)}) | > |b_{n-j}|-2\delta=\eps_{n-j}.\]    
		Therefore, $x^*_t\in s_{\eps_{n-j}}\dots s_{\eps_n}(B_{X^*})$.  This finishes the inductive proof of our claim. 
	\end{proof}
	
	The following corollary follows easily from the link between the Szlenk index and asymptotic uniform smoothness described in Remark \ref{rem:q-AUC*->q-summable-Szlenk}.
	
	\begin{coro}
		Let $p\in(1,\infty)$. If a Banach space is $p$-AUS, then $X\in \sA_p$.
	\end{coro}
	
	As we will explain in Section \ref{sec:A-N-T} the converse of Proposition \ref{prop:Sz-q-summable->Ap} is true: see Theorem \ref{thm:A-theorem-full}. In Section \ref{sec:asymptotic-c_0}, we will study in detail the case $p=\infty$.

	\section{\texorpdfstring{Asymptotically uniformly smoothable spaces and the $\sT$-game}{Asymptotically uniformly smoothable spaces and the -game}}
	\label{sec:AUS-T}

	For convenience, we will denote by $\langle\AUS\rangle$ (resp. $\langle p$-$\AUS\rangle$) the class of Banach spaces that admit an equivalent norm that is AUS (resp. $p$-AUS).
	In this section, we detail characterizations of the class $\langle p$-$\AUS\rangle$ in terms of games and in terms of upper-$p$-estimates for weakly null trees. 
	
	\subsection{\texorpdfstring{The class $\sT_p$}{The class}}
	
	The class $\sT_p$ is defined similarly to the class $\sA_p$, but this time the two-player game on the Banach space $X$ has countably infinitely many rounds. We still denote by $\cN_{w}(0)$ the set of weak neighborhoods of $0$ in $X$.
	\begin{defi}[The class $\sT_p$]
		\label{def:game-Tp}
		Let $p\in (1,\infty]$. The game $T(c,p)$ in a Banach space $X$ is the game $G(c,N_{q,\omega}, \omega)$ where $\frac1p +\frac1q=1$ and
		\begin{equation}
			N_{q,\omega}((x_i)_{i=1}^\infty) := \inf\Big\{c\in (0,\infty] \colon \forall a=(a_i)_{i=1}^\infty \in c_{00},\ \Big\|\sum_{i=1}^\infty a_ix_i \Big\| \le c\norm{a}_p\Big\},
		\end{equation}
		We let $\textsf{t}_{p}(X)$ be the infimum of all $c>0$ such that Player A has a winning strategy in the $T(c,p)$ game, provided such a $c$ exists and we let $\textsf{t}_p(X)=\infty$ otherwise.
		Finally, we denote by $\sT_p$ the class of all Banach spaces $X$ such that $\textsf{t}_p(X)<\infty$. 
	\end{defi}
	
	\begin{rema}
		\label{rem:T_pHahnBanach}
		The above notation can be justified by observing that 
		$$N_{q,\omega}((x_i)_{i=1}^\infty) = \sup\{\norm{(x^*(x_i))_{i=1}^\infty}_{\ell_q}\colon\ x^*\in B_{X^*}\}.$$    
	\end{rema}
	
	The fact that the game $T(c,p)$ is determined is more delicate due to the infinite nature of every play. However, note that Player B wins if and only if there exists $a\in c_{00}$ such that $\norm{\sum_{i=1}^\infty a_ix_i}> c\norm{a}_p$. It follows that the game is over after a finite number of rounds whenever Player B wins.
	
	\begin{prop}
		\label{prop:Tp-determination}
		Let $p\in (1,\infty]$ and $c>0$. Then, the $T(c,p)$ game is determined.  
	\end{prop}
	
	\begin{proof}[Sketch of proof]
		The $T(c,p)$ game is a game with perfect information and without draw, but it is infinite.  The fact that such games are determined is called the Gale–Stewart theorem and is due to the fact that Player A wins if and only if $(x_i)_{i=1}^\infty\in N_{q,\omega}^{-1}([0,c])$ which is a closed set in the product of the discrete topology on $X^\bN$. The proof relies on the following argument: if Player B does not have a winning strategy, then any ``defensive'' strategy for Player A (meaning a strategy that does not allow a winning strategy for Player B starting at the next round; an easy induction shows that such a defensive strategy must exist) is a winning strategy for Player A.
	\end{proof}
	
	\begin{rema}
		Understanding the topological restrictions under which an infinite two-player game with perfect information is determined is a topic on its own. Using Bernstein sets one can define games that are not determined. On the other hand, Martin's determinacy theorem \cite{Martin1975} states that if the outcome of the game is a Borel set, then the game is determined. 
	\end{rema}
	
	The next proposition is a reformulation of the class $\sT_p$ in terms of upper-$p$-estimates on branches of certain trees.
	
	\begin{prop}
		\label{prop:Tp-game-trees}
		Let $1<p\le \infty$ and $X$ be a Banach space. The following assertions are equivalent:  
		\begin{enumerate}[(i)]
			\item $X\in \sT_p$. 
			\item There exists a constant $c>0$ such that for any weak neighborhood basis $D$ at $0$ in $X$ and any weakly null tree $(x_t)_{t\in D^{<\omega}}\subset B_X$, there exists $\tau\in D^\omega$ such that for all $(a_i)_{i=1}^\infty \in c_{00}$,
			\begin{equation*}
				\Big\| \sum_{i=1}^\infty a_ix_{ \tau_{\restriction i} } \Big\| \le c \norm{(a_i)_{i=1}^\infty}_p.
			\end{equation*}
			\item There exists a constant $c>0$ such that for any weak neighborhood basis $D$ at $0$ in $X$ and any weakly null tree $(x_t)_{t\in D^{<\omega}}\subset S_X$, there exists $\tau\in D^\omega$ such that for all $(a_i)_{i=1}^\infty \in c_{00}$,
			\begin{equation*}
				\Big\| \sum_{i=1}^\infty a_ix_{ \tau_{\restriction i} } \Big\| \le c \norm{(a_i)_{i=1}^\infty}_p.
			\end{equation*}
		\end{enumerate}
	\end{prop}
	
	\begin{proof} 
		We first show that $(i)$ implies $(ii)$. Assume that $(ii)$ is not satisfied. For a given $c>0$, we will show that there is a winning strategy for Player B in the $T(c,p)$ game. By assumption, there exists a weak neighborhood basis $D$ at $0$ in $X$ and a weakly null tree $(x_t)_{t\in D^{<\omega}}\subset B_X$ such that for all $\tau\in D^\omega$, $N_{q,\omega}( (x_{\tau_{\restriction_i}})_{i=1}^\infty ) > c$. When Player A chooses $U_1\in \cN_{w}(0)$, then Player B chooses $V_1\in D$ such that $x_1 := x_{(V_1)} \in U_1$ (this is possible because the net $(x_{(U)})_{U \in D}$ is weakly null). Then, Player A chooses $U_2\in \cN_{w}(0)$, to which Player B's response is $x_2 := x_{(V_1,V_2)}\in U_2$ (again possible, since $(x_{(V_1,V)})_{V\in D}$ is weakly null). The play continues in this way and the result is $(x_{\tau_{\restriction_i} })_{i=1}^\infty$ for some $\tau\in D^\omega$, which satisfies $N_{q,\omega}( (x_{\tau_{\restriction_i} })_{i=1}^\infty ) > c$.
		
		We now prove that $(ii)$ implies $(i)$. Assume now that $(i)$ is not satisfied and fix $c>0.$ Since the $T(c,p)$ game is determined, Player B must have a winning strategy. We define by induction on $k$ and for $U_1,\dots, U_k \in \cN_{w}(0)$, the vector $x_{(U_1,\dots, U_k)}$ to be  Player B's response according to his winning strategy following the choices $U_1,x_{(U_1)},\dots, U_k$ in the $T(c,p)$ game. Clearly $(x_t)_{t\in \cN_{w}(0)^{<\omega}}\subset B_X$ is a weakly null tree and $N_{q,\omega}( (x_{\tau_{\restriction_i}})_{i=1}^\infty ) > c$ for all $\tau\in \cN_{w}(0)^\omega$.
		
		Obviously $(ii)$ implies $(iii)$. So, let us assume $(iii)$ and prove that $(ii)$ holds. First, we fix, as we may, a weakly null tree $(y_t)_{t\in D^{<\omega}}\subset S_X$ (for any $U\in D$, pick $y_U\in U\cap S_X$, which is possible as $X$ is infinite-dimensional and let $y_{(U_1,\dots, U_n)} := y_{U_n}$, for $(U_1,\dots, U_n)\in D^{<\omega}$). Let $(x_t)_{t\in D^{<\omega}}$ be a weakly null tree in $B_X$. We also fix $\eps_i>0$, for $i\ge 1$. We let 
		$$
		z_t := \begin{cases}
			y_t \textrm{ if }\abs{t}=i \textrm{ and }\norm{x_t}<\eps_i,\\ 
			\frac{x_t}{\norm{x_t}} \textrm{ otherwise.}
		\end{cases}
		$$
		It is easily checked that  $(z_t)_{t\in D^{<\omega}}$ is a weakly null tree in $S_X$. Therefore, there exists $\tau\in D^\omega$ such that for all $(a_i)_{i=1}^\infty \in c_{00}$,
		\begin{equation*}
			\Big\|\sum_{i=1}^\infty a_iz_{ \tau_{\restriction i} } \Big\| \le c \norm{a}_p.
		\end{equation*}
		If $A := \{i\in \bN\colon \norm{x_{ \tau_{\restriction i} } }\ge \eps_i\}$ and $a\in c_{00}$, we let 
		$$
		b_i := \begin{cases}
			a_i\norm{x_{ \tau_{\restriction i} }} \textrm{ if } i\in A,\\
			0 \textrm{ otherwise.}
		\end{cases}
		$$
		Then, if $q$ is the conjugate of $p$, we have 
		\begin{align*}
			\Big\|\sum_{i=1}^\infty a_ix_{ \tau_{\restriction i} } \Big\| &\le \Big\|\sum_{i
				\in A} a_ix_{ \tau_{\restriction i} }\Big\|+\norm{(\eps_i)_{i=1}^\infty}_q\norm{a}_p=\Big\|\sum_{i=1}^\infty b_iz_{ \tau_{\restriction i} }\Big\|+\norm{(\eps_i)_{i=1}^\infty}_q\norm{a}_p\\
			&\le c\norm{b}_p+\norm{(\eps_i)_{i=1}^\infty}_q\norm{a}_p \le (c+\norm{(\eps_i)_{i=1}^\infty}_q)\norm{a}_p.
		\end{align*}
		Since $\norm{(\eps_i)_{i=1}^\infty}_q$ can be chosen as small as we wish, this finishes the proof. 
	\end{proof}

	\subsection{Asymptotically uniformly smooth renormings}
	\label{sec:Tp=p-AUSs}
	
	The main result of this section is a characterization of the class of $p$-asymptotically uniformly smoothable spaces in terms of the class $\sT_p$. Here $p\in(1,\infty)$ and we delay the discussion of the case $p=\infty$ to Section \ref{sec:asymptotic-c_0}. The result is due to Odell and Schlumprecht \cite{OdellSchlumprecht2006RACSAM} for separable Banach spaces and the general case was obtained by R.M. Causey \cite{Causey2018}.
	The proof we present here is due to Causey and uses a completely different approach than the one from Odell and Schlumprecht. 
	
	\begin{theo}
		\label{thm:Tp=p-AUS} 
		Let $X$ be a Banach space and $p\in(1,\infty)$. Then,
		$$ 
		X\in \sT_p \textrm{ if and only if } X \in\langle p\textrm{-}\AUS\rangle.
		$$
	\end{theo}
	
	\begin{proof}
		Let us prove first the easier implication, namely that if $X\in\langle p$-$\AUS\rangle$, then Player A has a winning strategy for the $\sT_p$-game. So, assume, as we may, that $X$ is $p$-AUS and therefore that $\bar{\rho}_X(t)\le ct^p$, for some $c >0$. Let $ N_{\bar{\rho}_X}$ be the  absolute norm  associated with $\bar{\rho}_X$ and recall that, by Corollary \ref{cor:iterated_pnorm}, there exists $\gamma>0$ such that $\norm{\cdot}_{\Lambda_{N_{\bar{\rho}_X}}}\le \gamma \norm{\cdot}_p$. It follows from a, by now, classical compactness argument (cf the proof of Claim \ref{claim:biortho-blocking} together with Proposition \ref{prop:moduli-weak-opt}) that for any $\eps>0$ and any finite-dimensional subspace $F$ of $X$ there exists $U \in \cN_{w}(0)$ such that for all $x \in F$, $y \in U\cap B_X$ and $t\ge 0$,
		\begin{equation*}
			\norm{x+ty}\le (1+\eps)N_{\bar{\rho}_X}(\norm{x}, t).
		\end{equation*}
		Fix a sequence $(\eps_n)_{n=1}^\infty$ in $(0,1)$ such that $\prod_{n=1}^\infty (1+\eps_n)\le 2$. We can describe a winning strategy for Player A in the $T(2\gamma,p)$ game as follows. Player A's initial choice $U_1$ is arbitrary. Once $U_1, x_1, \dots, U_n, x_n$ have been chosen, let $F$ be the linear span of $\{x_1,\dots,x_n\}$. Then, we pick $U_{n+1}\in \cN_{w}(0)$ such that for all $x \in F$, $y \in U_{n+1}\cap B_X$ and $t\ge 0$,
		\begin{equation*}
			\norm{x+ty}\le (1+\eps_n)N_{\bar{\rho}_X}(\norm{x},t).
		\end{equation*}
		Then, an easy induction yields that at the end of the game:
		$$\forall a \in c_{00},\ \ \Big\|\sum_{n=1}^\infty a_n x_n\Big\|\le \prod_{n=1}^\infty (1+\eps_n)\norm{a}_{\Lambda_{N_{\bar{\rho}_X}}}\le 2\gamma\norm{a}_p.$$
		We have proved that this strategy is indeed a winning strategy for Player A in the $T(2\gamma,p)$ game. Therefore, $X \in \textsf{T}_p$.
		
		Showing that every Banach space in $\sT_p$ admits an equivalent norm that is $p$-AUS is arguably the more difficult implication. The construction of the equivalent norm we reproduce here is due to R.M. Causey \cite{Causey2018} and is an asymptotic analog of Pisier's $p$-smooth renorming of Banach spaces with martingale type $p$ \cite{Pisier1975}.  For the rest of the proof, we fix a weak neighborhood basis $D$ at $0$ in $X$. If $X$ has $\textsf{T}_p$, then by Proposition \ref{prop:Tp-game-trees} we can find $c\ge 1$ such that for any weakly null tree $(x_t)_{t\in D^{<\omega}}\subset B_X$, there exists $\tau \in D^\omega$ such that for all scalar sequences $a=(a_i)_{i=1}^\infty \in c_{00}$, 
		\begin{equation}
			\Bigl\|\sum_{i=1}^\infty a_ix_{\tau_{\restriction i}}\Bigr\|^p\le c^p\sum_{i=1}^\infty \abs{a_i}^p.
		\end{equation}   
		We first note that for any $x\in X$ and any  weakly null tree $(x_t)_{t\in D^{<\omega}}\subset B_X$, there exists $\tau\in D^\omega$ such that for all scalar sequences $a \in c_{00}$, 
		\begin{equation}
			\label{eq:f1}
			\Bigl\|x+\sum_{i=1}^\infty a_ix_{\tau_{\restriction _i}}\Bigr\|^p\le (2c)^p\Bigl[\norm{x}^p+ \sum_{i=1}^\infty \abs{a_i}^p\Bigr].
		\end{equation}
		Indeed, for an appropriate branch $\tau$, it holds that 
		\begin{align*} 
			\Bigl\|x+\sum_{i=1}^\infty a_ix_{\tau_{\restriction _i}}\Bigr\|^p & \le 2^{p-1}\Bigl( \norm{x}^p + \Bigl\|\sum_{i=1}^\infty a_ix_{\tau_{\restriction _i}}\Bigr\|^p\Bigr) 
			%    & \le (2c)^p\max\Bigl\{\norm{x}^p,\sum_{i=1}^\infty |a_i|^p\Bigr\} \\
			\le 2^{p-1}c^p\Bigl[\norm{x}^p+ \sum_{i=1}^\infty |a_i|^p\Bigr].
		\end{align*}
		Simple convexity arguments of this nature will not be detailed in the sequel.
		Now, let $A:=2c$ and define 
		\begin{equation*}
			f(x) := \Bigl[ \underset{(x_t)}{\ \sup\ }\underset{\tau}{\ \inf\ }\underset{(a_i)}{\ \sup\ } \frac{1}{A^p}\Bigl\|x+\sum_{i=1}^\infty a_ix_{\tau_{\restriction _i}}\Bigr\|^p-\sum_{i=1}^\infty \abs{a_i}^p\Bigr]^{1/p},
		\end{equation*}
		where the outer supremum is taken over all weakly null trees $(x_t)_{t\in D^{<\omega}}$ in  $B_X$, the infimum is taken over $\tau\in D^\omega$ and the inner supremum is taken over all scalar sequences $a\in c_{00}$. It follows from taking $x_t=0$ for all $t$ and $(a_i)_i=0$ that $f(x)\ge slant \frac{\norm{x}}{A}$ for all $x\in X$. On the other hand, it follows from  (\ref{eq:f1}) that $f(x)\le \norm{x}$. We also have that  $f(\lambda x)=|\lambda|f(x)$, for each $x\in X$ and each scalar $\lambda$. Let us detail this last fact. It is clear that $f(0)=0$ since $f$ is equivalent to $\norm{\cdot}$. Now, if $x\in X$ and $\lambda\neq 0$ we have 
		\begin{align*}
			f(\lambda x) & =\Bigl[ \underset{(x_t)}{\ \sup\ }\underset{\tau}{\ \inf\ }\underset{(a_i)}{\ \sup\ } \frac{1}{A^p}\Bigl\|\lambda x+\sum_{i=1}^\infty a_ix_{\tau_{\restriction _i}}\Bigr\|^p-\sum_{i=1}^\infty \abs{a_i}^p\Bigr]^{1/p}\\
			& =\abs{\lambda}\Bigl[ \underset{(x_t)}{\ \sup\ }\underset{\tau}{\ \inf\ }\underset{(a_i)}{\ \sup\ } \frac{1}{A^p}\Bigl\|x+\sum_{i=1}^\infty \frac{a_i}{\lambda}x_{\tau_{\restriction _i}}\Bigr\|^p-\sum_{i=1}^\infty \abs{\frac{a_i}{\lambda}}^p\Bigr]^{1/p}\\
			& =\abs{\lambda}\Bigl[ \underset{(x_t)}{\ \sup\ }\underset{\tau}{\ \inf\ }\underset{(b_i)}{\ \sup\ } \frac{1}{A^p}\Bigl\|x+\sum_{i=1}^\infty b_ix_{\tau_{\restriction _i}}\Bigr\|^p-\sum_{i=1}^\infty \abs{b_i}^p\Bigr]^{1/p}\\
			&= \abs{\lambda}f(x).
		\end{align*}
		
		The issue at this point is that there is no reason why $f$ should satisfy the triangle inequality. A classical way to remedy this problem is to ``convexify'' the function $f$. For that purpose, we set 
		\begin{equation*}
			\hat{f}(x) := \inf\Bigl\{\sum_{i=1}^n f(x_i): n\in\bN, x=\sum_{i=1}^n x_i\Bigr\}.
		\end{equation*}
		It is elementary to verify that $\hat{f}\colon X\to [0,\infty)$, satisfies $\hat{f}(\lambda x)=\abs{\lambda}\hat{f}(x)$ for all $\lambda\in\bR$ and most importantly that $\hat{f}$ satisfies the triangle inequality. 
		%In fact, since $f$ is positively homogeneous it is well-known (and also elementary to verify) that $\hat{f}$ is convex and hence it coincides with the Minkowski functional of the convex balanced set $\{x\in X \colon \hat{f}(x)\le 1\}$; in particular it is a norm. 
		It also follows straightforwardly from the definition of $\hat{f}$ that $\frac{\norm{x}}{A}\le \hat{f}(x) \le f(x) \le  \norm{x}$ for all $x\in X$. Therefore, $\hat{f}$ defines an equivalent norm on $X$ and with a bit more work one can see that $ \{ x\in X \colon \hat{f}(x)\le 1\}$ is the closed convex hull of $\{x\in X\colon f(x)<1\}$. It remains to prove that $(X,\hat{f})$ is $p$-AUS. According to Proposition \ref{prop:moduli-nets} and Lemma \ref{lem:AUS-powertype} it is sufficient to show that there is $K>0$ such that for every $x\in X$ with $\hat{f}(x)\le 1$, every weakly null net $(x_U)_{U\in D}\subset B_{(X,\hat{f})}$ and every $t>0$, 
		\begin{equation}
			\label{eq:f}
			\underset{U}{\limsup} \hat{f}(x+t x_U) \le 1 + Kt^p.
		\end{equation}
		So let $x\in X$ with $\hat{f}(x)\le 1$, $(x_U)_{U\in D}\subset B_{(X,\hat{f})}$ be a weakly null net and $t>0$. As noted above, $B_{(X,\hat{f})}$ is the closed convex hull of $\{y\in X \colon f(y)<1\}$. Therefore, for each $\eps>0$, we can find $y_1, \ldots, y_k\in X$ with $f(y_i)<1$ and convex coefficients $w_1, \ldots, w_k$ such that $\hat{f}(x-\sum_{i=1}^k w_iy_i)<\eps$.  Then, 
		\begin{align} 
			\notag \limsup_{U} \hat{f}(x+ t x_U) -1 & \le \eps + \limsup_U \sum_{i=1}^k w_i (\hat{f}(y_i + tx_U)-1)\\ 
			\label{eq:fhat}& \le \eps + \sum_{i=1}^k w_i (\limsup_U\hat{f}(y_i + tx_U)-1).
		\end{align} 
		The key technical step in order to achieve \eqref{eq:f} will be to show the following claim.
		\begin{claim}
			\label{claim:Tp->pAUS}
			For every $y\in X$, every weakly null net $(y_U)_{U\in D}\subset B_{(X,\norm{\cdot})}$ and each $t>0$, 
			\begin{equation}
				\label{eq:f2}
				\underset{U}{\lim\sup} f(y+t y_U)^p \le f(y)^p+t^p.
			\end{equation}
		\end{claim}
		Indeed, assuming Claim \ref{claim:Tp->pAUS}, we can conclude as follows. Recall that $\hat{f}(x_U)\le 1$, for $U\in D$ and $f(y_i)<1$, for $1\le i\le n$. Therefore, $\norm{y_i}\le A$, for $1\le i\le n$ and if we let $y_U := A^{-1}x_U$, for $U\in D$, then $(y_U)_{U\in D}$ is a weakly null net in $B_{(X,\norm{\cdot})}$. Applying \eqref{eq:f2} to $y_i$, for $1\le i\le n$, we get
		\begin{align*}
			\limsup_U \hat{f}(y_i+t x_U)^p & \le \limsup_U f(y_i + tA y_U)^p\\ 
			& \stackrel{\eqref{eq:f2}}{\le} f(y_i)^p + t^pA^p < 1+t^pA^p.
		\end{align*}
		%Therefore, by concavity of the function  $t\mapsto (1+t)^{1/p}$,
		%\begin{equation*}
		%\limsup_U \hat{f}(y_i + t x_U)-1 \le (1+ t^p A^p)^{1/p}-1 \le \frac{A^pt^p}{p}.
		%\end{equation*}
		Plugging this directly in $\eqref{eq:fhat}$ we get
		\begin{align*}
			\limsup_{U} \hat{f}(x+ t x_U) -1 & \le \eps+ \sum_{i=1}^k w_i \frac{A^pt^p}{p}= \eps+ \frac{A^pt^p}{p},
		\end{align*}
		but since $\eps>0$ was arbitrary, this finishes the proof of (\ref{eq:f}). 
		
		To complete the proof of the theorem, it remains to prove Claim \ref{claim:Tp->pAUS}. Let $\vep>0$ and observe that after passing to a subnet and relabeling, we can assume that for all $U$ 
		\begin{equation*}
			\limsup_U f(y+ t y_U)^p -\vep < f(y+ t y_U)^p.
		\end{equation*}
		For each $U$, by definition of $f(y+ty_U)$, we find a weakly null tree $(x^U_t)_{t\in D^{<\omega}}\subset B_X$ such that for each $\tau\in D^\omega$, there exists $a\in c_{00}$ satisfying
		\begin{equation}
			\label{eq:f5}
			\frac{1}{A^p}\Bigl\| y + t y_U + \sum_{i=1}^\infty a_i x^U_{\tau_{\restriction i}} \Bigr\|^p - \sum_{i=1}^\infty \abs{a_i}^p > \limsup_U f(y+ t y_U)^p -\vep .
		\end{equation}
		We define the weakly null tree $(x_t)_{t\in D^{<\omega}}\subset B_X$ by letting $x_{(U)} := y_U$ and $x_{(U, U_1, \dots, U_k)} := x^U_{(U_1, \ldots, U_k)}$ for $1\le k< \infty$. By definition of $f(y)$, there exists $\sigma\in D^\omega$ such that for all $b\in c_{00}$, 
		\begin{equation}
			\label{eq:f6}
			\frac{1}{A^p}\Bigl\| y + \sum_{i=1}^{\infty} b_i x_{\sigma_{\restriction_i}}\Bigr\|^p - \sum_{i=1}^{\infty}\abs{b_i}^p \le f(y)^p+\eps.
		\end{equation}   
		Write $\sigma=(U, U_1, \dots, U_n, \dots)$, $\tau=(U_1, \dots, U_n, \dots)$ and let $a := a(U,\tau)\in c_{00}$ satisfying \eqref{eq:f5}.  By taking $b_1=t$ and  $b_{i+1}=a_i$ for $i\ge 1$, it follows that 
		\begin{align*} 
			\limsup_U f(y+ t y_U)^p -\vep  & \stackrel{\eqref{eq:f5}}{<} \frac{1}{A^p}\Bigl\| y + ty_U + \sum_{i=1}^\infty a_i x^U_{\tau_{\restriction i}}\Bigr\|^p -\sum_{i=1}^\infty \abs{a_i}^p \\ 
			& = \frac{1}{A^p}\Bigl\| y + \sum_{i=1}^{\infty}b_ix_{\sigma_{\restriction i}}\Bigr\|^p - \sum_{i=1}^{\infty}\abs{b_i}^p + t^p \\
			& \stackrel{\eqref{eq:f6}}{\le} f(y)^p + \eps + t^p. 
		\end{align*}
		Therefore, $\limsup_U f(y+ t y_U)^p  < f(y)^p + 2\eps + t^p$, but since $\eps>0$ was arbitrary, we have proved \eqref{eq:f2}. 
		
	\end{proof}
	
	\section{\texorpdfstring{The class $\sN_p$ and the convex Szlenk index}{The class  and the convex Szlenk index}}
	In this section, we introduce and briefly discuss a variant of the $\sA$-game and a variant of the Szlenk index. 
	The $\sN$-game is very similar to the $\sA$-game, albeit that we only consider constant coefficients in the referee function. 
	
	\begin{defi}[The class $\sN_p$]
		\label{def:game-Np}
		Let $p\in (1,\infty]$, $c>0$ and $n\in\bN$. The game $N(c,p,n)$ in a Banach space $X$ is the game $G(c, R_{q,n}, n)$ where $\frac1p +\frac1q=1$ and
		\begin{equation}
			R_{q,n}((x_i)_{i=1}^n) := \inf\Big\{c\in (0,\infty] \colon \Big\|\sum_{i=1}^n x_i\Big\|\le c n^{\frac1p}\Big\},
		\end{equation}
		We let $\textsf{n}_{p,n}(X)$ be the infimum of all $c>0$ such that Player A has a winning strategy in the $N(c,p,n)$ game and we let $\textsf{n}_p(X) := \sup_n \textsf{n}_{p,n}(X)$. We note that $\textsf{n}_p(X)$ is the infimum of all $c>0$ such that for each $n\in\bN$, Player A has a winning strategy in the $N(c,p,n)$-game if such a $c$ exists and $\textsf{n}_p(X)=\infty$ otherwise. 
		Finally, we denote by $\sN_p$ the class of all Banach spaces $X$ such that $\textsf{n}_p(X)<\infty$. 
	\end{defi}
	
	We record the tree reformulation of the $\sN_p$-game and omit its proof as the argument is a straightforward modification of the proof of the analogous result for the $\sT_p$-game.
	
	\begin{prop}
		\label{prop:Np-game-trees}
		Let $1<p\le \infty$ and $X$ be a Banach space. The following assertions are equivalent.  
		\begin{enumerate}[(i)]
			\item $X\in \sN_p$. 
			\item There exists a constant $c>0$ such that for any weak neighborhood basis $D$ at $0$ in $X$, any $n\in \bN$ and any weakly null tree $(x_t)_{t\in D^{\le n}}\subset B_X$, there exists $\tau\in D^n$ such that $\norm{\sum_{i=1}^n x_{\tau_{\restriction_i}}}\le cn^{1/p}$. 
		\end{enumerate}
	\end{prop}
	
	The class $\sN_p$ also admits a characterization in terms of a renorming property. This renorming theorem is obtained via a derivation index that is somewhat slower than the Szlenk derivation, and which is a convexified version of it.
	
	\begin{defi} Let $X$ be  Banach space and $K$ be a weak$^*$ compact subset of $ X^*$. For each $\eps>0$, define $c_\eps(K)$ to be the weak$^*$ closed convex hull of $s_\eps(K)$. 
		Given an ordinal $\xi$, the derived set of order $\xi$, denoted by $c_\eps^\xi(K)$, is defined inductively by letting 
		\begin{itemize}
			\item $c_\eps^0(K) := K$,  
			\item $c_\eps^{\xi+1}(K) := c_\eps(c^\xi_\eps(K))$,
			\item $c_\eps^{\xi}(K) := \cap_{\zeta<\xi}c^\zeta_\eps(K)$ if $\xi$ is a limit ordinal.
		\end{itemize}
		As usual, for $\eps=0$, we set $c_0^{\xi}(K) := K$ for all ordinal $\xi$. We then define $\Cz(X,\eps)$ as the least ordinal $\xi$ so that $c_\eps^\xi(B_{X^*})=\emptyset$, if such ordinal exists and we set $\Cz(X,\eps)=\infty$ otherwise. The \emph{convex Szlenk index of $X$} is then defined as
		\begin{equation*}
			\Cz(X) := \sup_{\eps>0}\Cz(X,\eps).
		\end{equation*}  
	\end{defi}
	
	See Exercise \ref{ex:Cz-Kz} or \cite{HajekLancien} for another description of the convex Szlenk index.
	
	\begin{rema} Let $X$ be a Banach space. It is obvious that $\Cz(X)\ge \Sz(X)$. It follows from Theorem \ref{thm:T-theorem-full} that if $\Sz(X)=\omega$, then $\Cz(X)=\omega$ (see next section for details). The fact that for separable Banach spaces the convex Szlenk index and the Szlenk index are always equal is due to Lancien, Proch\'azka and Raja \cite{LPR}, but will not be used in this book. 
	\end{rema}
	
	A more general version of the following result is proved in \cite{Causey3.5}. 
	
	\begin{theo}
		\label{thm:N-theorem} 
		Fix $p\in(1, \infty)$ and let $q$ be conjugate to $p$. Let $X$ be a Banach space. The following assertions are equivalent. 
		\begin{enumerate}[(i)]
			\item $X\in \sN_p$.  
			\item There exist a constant $M\ge 1$ and a constant $c>0$ such that for each $t_0\in (0,1]$, there exists a norm $\abs{\cdot}$ on $X$ such that $M^{-1}\abs{x}\le \norm{x}_X\le M \abs{x}$ for all $x\in X$ and $\bar{\rho}_{\abs{\cdot}}(t_0)\le c t_0^p$.
			\item There exists a constant $C>0$ such that $\Cz(X,\eps) \le C\eps^{-q}$ for all $\eps\in (0,1)$.
		\end{enumerate}
	\end{theo}
	
	\begin{proof}[Sketch of Proof] 
		The argument for $(i) \Rightarrow (ii)$ is similar to the proof showing that if $X\in \sT_p$, then $X\in \langle p$-$\AUS\rangle$. This was essentially proved in \cite[Theorem 4.2]{GKL2001} for separable spaces, and the proof can be adjusted to handle nonseparable spaces (for more details see \cite{GKL2001} or \cite{Causey3.5}). The implication $(ii) \Rightarrow (iii)$ is a consequence of the AUS-AUC$^*$  duality and the proof of Proposition \ref{prop:Sz-omega}, since $(ii)$ implies that, in this situation, the derived sets for the Szlenk derivation are included in balls which are convex sets. The proof of $(iii)\Rightarrow (i)$ goes along the same lines as the proof showing that $q$-summable Szlenk index implies membership in $\sA_p$ (cf Proposition \ref{prop:Sz-q-summable->Ap}), and we refer to \cite{Causey3.5} for more details. 
		Note that an alternative proof of the implication $(ii) \Rightarrow (i)$ relying on iterated norms to prove upper tree estimates is explained in \cite{CauseyFovelleLancien2023}.
	\end{proof}

	\section{\texorpdfstring{An ordinal index characterization of the class $\langle \AUS \rangle$}{An ordinal index characterization of the class}}
	\label{sec:A-N-T}
	
	The goal of this section is to provide a characterization of the class of asymptotically uniformly smoothable spaces in terms of the Szlenk index. To this end, we develop our understanding of the relationship between the various classes $\sN_p$, $\sA_p$ and $\sT_p$.
	Let us first recall the following characterization of $\sT_p = \langle p$-$\AUS\rangle$, which is a direct corollary of Theorem \ref{thm:Tp=p-AUS},  Proposition \ref{prop:Tp-game-trees} and Corollary \ref{cor:duality}. 
	
	\begin{theo}
		\label{thm:T-theorem-full}
		Fix $p\in(1, \infty)$ and let $q$ be conjugate to $p$. Let $X$ be a Banach space. The following assertions are equivalent. 
		\begin{enumerate}[(i)]
			\item $X\in \sT_p$.
			\item There exists a constant $c>0$ such that for any weak neighborhood basis $D$ at $0$ in $X$ and any weakly null tree $(x_t)_{t\in D^{<\omega}}\subset B_X$, there exists $\tau\in D^\omega$ such that for all $(a_i)_{i=1}^\infty \in c_{00}$,
			\begin{equation*}
				\Big\| \sum_{i=1}^\infty a_ix_{ \tau_{\restriction i} } \Big\| \le c \norm{(a_i)_{i=1}^\infty}_p.
			\end{equation*}
			\item There exist a constant $M\ge 1$, a constant $C>0$ and a norm $\abs{\cdot}$ on $X$ such that $M^{-1}\norm{x}_X\le \abs{x}\le M \norm{x}_X$ for all $x\in X$ and $\bar{\rho}_{\abs{\cdot}}(t)\le C t^p$ for all $t\in(0,1]$.
			\item $X$ admits an equivalent norm whose dual norm is $q$-AUC$^*$.
		\end{enumerate}
	\end{theo}
	
	\begin{rema}
		It is interesting to point out that all the assertions above imply that the Szlenk index of $X$ is $q$-summable, but there is currently no characterization of the class $\sT_p$ in terms of a derivation index. It is in fact unlikely to find one. Indeed, it is proved in \cite{CauseyLancien2023} that the class of separable Banach spaces in $\sT_p$ is analytic non-Borel in the class of separable Banach spaces in the sense of Bossard \cite{Bossard}, while these derivations are usually Borel operations, still in the sense of Bossard.      
	\end{rema}
	
	In this section, we study further the class $\sA_p$ and give two additional characterizations of $\sA_p$: one in terms of a renorming property and another in terms of the behavior of the Szlenk derivation. 
	
	\begin{theo}
		\label{thm:A-theorem-full} 
		Fix $p\in(1, \infty)$ and let $q$ be conjugate to $p$. Let $X$ be a Banach space. The following assertions are equivalent.  
		\begin{enumerate}[(i)]
			\item\label{it:Ti} $X\in \sA_p$. 
			\item\label{it:Tii} There exists a constant $c>0$ such that for any weak neighborhood basis $D$ at $0$ in $X$, any $n\in \bN$ and any weakly null tree $(x_t)_{t\in D^{\le n}}\subset B_X$, there exists $\tau\in D^n$ such that for all $(a_i)_{i=1}^n \in \bR^n$,
			\begin{equation*}
				\Big\| \sum_{i=1}^n a_ix_{ \tau_{\restriction i} } \Big\| \le c \norm{(a_i)_{i=1}^n}_p.
			\end{equation*}%\item There exists a constant $c$ such that for any $n\in\nn$ and any $(e_i)_{i=1}^n \in \{X\}_n$, $\|(e_i)_{i=1}^n\|_q^w\le c$. 
			\item\label{it:Tiii} There exist a constant $M\ge 1$ and a constant $C>0$ such that for any $t_0\in (0,1]$ there exists a norm $\abs{\cdot}$ on $X$ satisfying $M^{-1}\norm{x}_X\le \abs{x}\le M \norm{x}_X$ for all $x\in X$ and
			$$\forall t\ge t_0,\ \ \bar{\rho}_{\abs{\cdot}}(t)\le Ct^p.$$
			\item\label{it:Tiv} $X$ has $q$-summable Szlenk index.   
		\end{enumerate}
	\end{theo}
	\begin{proof}[Sketch of Proof] The equivalence between $(i)$ and $(ii)$ is Proposition \ref{prop:Ap-game-trees}.
		
		The proof of $(ii) \Rightarrow (iii)$ is a modification of the proof of the similar implication in Theorem \ref{thm:T-theorem-full} (cf. Theorem \ref{thm:Tp=p-AUS}), and the details can be found in \cite{CauseyFovelleLancien2023}. We describe the main ingredients. Again, we fix 
		a weak neighborhood basis $D$ at $0$ in $X$. Let us assume that $(ii)$ is satisfied. Similarly to the proof of Theorem \ref{thm:Tp=p-AUS} we let $A := 2c$ and for $n \in \bN$ and $x\in X$, we set:   
		\begin{equation*}
			f_n(x) := \Bigl[ \underset{(x_t)}{\ \sup\ }\underset{\tau}{\ \inf\ }\underset{(a_i)}{\ \sup\ } \frac{1}{A^p}\Bigl\|x+\sum_{i=1}^n a_ix_{\tau_{\restriction_i}}\Bigr\|^p-\sum_{i=1}^n \abs{a_i}^p\Bigr]^{1/p},
		\end{equation*} 
		where the outer supremum is taken over all weakly null collections $(x_t)_{t\in D^{\le n}}$ in  $B_X$, the infimum is taken over $\tau \in D^n$, and the inner supremum is taken over all scalar sequences $(a_i)_{i=1}^n$. As in the proof of Theorem \ref{thm:Tp=p-AUS}, we have that $f_n$ is positively homogeneous and satisfies $\frac{\norm{x}}{A}\le f_n(x) \le \norm{x}$ for all $x\in X$. The key step will be to show the following analogue of Claim \ref{claim:Tp->pAUS}.
		\begin{claim}
			\label{claim:Ap->pAUS}
			For every $y\in X$, every weakly null net $(y_U)_{U\in D}\subset B_{(X,\norm{\cdot})}$, every $n\in \bN$ and every $t>0$, 
			\begin{equation}
				\label{eq:fn}
				\underset{U}{\lim\sup}\, f_n(y+t y_U)^p \le f_{n+1}(y)^p+t^p.
			\end{equation}
		\end{claim}
		The proof follows exactly the same lines as for Claim \ref{claim:Tp->pAUS}. Let us just stress that the only difference is that we are not dealing with trees of infinite height anymore. From weakly null trees of height $n$ rooted at the points $y+t y_U$, $U\in D$, we can only build a weakly null tree of height $(n+1)$ rooted at $y$. This should explain our statement. 
		
		In this situation, an additional step is to take a $p$-average of the $f_n$. So, fix $N\in \bN$ and define 
		\begin{equation*}
			g_N(x) := \Big(\frac{1}{N}\sum_{n=1}^{N}f_n(x)^p\Big)^{1/p}.
		\end{equation*}   
		Clearly, we still have that for all $x\in X$ and $N\in\bN$, $\frac{\norm{x}}{A} \le g_N(x) \le \norm{x}$ and $g_N$ is positively homogeneous. Then, applying Claim \ref{claim:Ap->pAUS} for each $n\in \{1,\dots,N\}$, we obtain that for any weakly null net $(y_U)_{U\in D}\subset B_{(X,\norm{\cdot})}$, any $N\in \bN$, any $y\in AB_{(X,\norm{\cdot})}$ and any $t>0$, 
		\begin{equation}
			\label{eq:f3}
			\underset{U}{\lim\sup}\, g_N(y+t y_U)^p \le g_N(y)^p+t^p + \frac{A^p}{N}.
		\end{equation}
		
		The last stage of the proof is again to convexify our function $g_N$. So, we set 
		\begin{equation*}
			\abs{x}_N = \inf\Bigl\{\sum_{i=1}^n g_N(x_i)\colon n\in\bN, x=\sum_{i=1}^n x_i\Bigr\},
		\end{equation*} 
		which defines an equivalent norm on $X$ satisfying $\frac{\norm{x}}{A}\le \abs{x}_N \le \norm{x}$. Moreover, $B_{(X,\abs{\cdot}_N)}$ is the closed, convex hull of $\{x\in X\colon g_N(x)<1\}$. We shall now prove that 
		\begin{equation}
			\label{ex:f4}
			\forall t>0,\ \ \bar{\rho}_{(X, \abs{\cdot}_N)}(t) \le \frac{A^p}{p}\Big(t^p + \frac{1}{N}\Big).
		\end{equation}
		First, we fix $y\in X$ such that $g_N(y)<1$ and from this it follows that $\norm{y}\le A$. Then, fix $t>0$ and a weakly null net $(y_U)_{U\in D}\subset B_{(X,\abs{\cdot}_N)}$ and define $x_U := A^{-1}y_U$ and observe that $(x_U)_{U\in D}\subset B_{(X,\norm{\cdot})}$ is weakly null. Finally, apply (\ref{eq:f3}) to get
		\begin{align*}
			\underset{U}{\lim\sup}\,\abs{y + t y_U}^p_N & \le \underset{U}{\lim\sup}\, g_N(y+t A x_U)^p\\ 
			&\le  g_N(y)^p+t^pA^p +\frac{A^p}{N}<1+t^pA^p +\frac{A^p}{N}.
		\end{align*}
		Therefore, by concavity of the function  $h(r) :=(1+r)^{1/p}$,
		\begin{equation*}
			\underset{U}{\lim\sup}\, \abs{y+t y_U}_N-1 \le  \Big(1+ t^p A^p + \frac{A^p}{N}\Big)^{1/p}-1 \le \frac{A^p}{p}\Big(t^p+\frac1N\Big).
		\end{equation*}
		Next, fix $x\in B_{(X, \abs{\cdot}_N)}$. As noted above, $B_{(X, \abs{\cdot}_N)}$ is the closed convex hull of $\{y\in X \colon g_N(y)<1\}$. Therefore, for each $\eta>0$, we can find $y_1, \dots, y_k\in X$ with $g_N(y_i)<1$ and convex coefficients $w_1, \dots, w_k$ such that $|x-\sum_{i=1}^k w_iy_i|_N<\eta$. Then,
		\begin{align*} 
			\underset{U}{\lim\sup}\, \abs{x+ t y_U}_N -1 & \le \eta + \underset{U}{\lim\sup}\, \sum_{i=1}^k w_i (\abs{y_i+ty_U}_N-1)\\ 
			&\le \eta + \sum_{i=1}^k w_i \Big(\frac{A^p}{p}\big(t^p+\frac1N\big)\Big) = \eta + \frac{A^p}{p}\big(t^p+\frac1N\big).
		\end{align*} 
		Since $\eta>0$ was arbitrary, this finishes the proof of (\ref{ex:f4}). 
		
		To conclude, it is clear by taking $N$ large enough in (\ref{ex:f4}) that for any $t_0 \in (0,1]$, there exists an equivalent norm $\abs{\cdot}$ on $X$ such that $\frac{\norm{x}}{A}\le \abs{x}\le \norm{x}$ and for any $t \ge t_0$, $\bar{\rho}_{(X,\abs{\cdot})}(t) \le \frac{(A+1)^p}{p}t^p$. We have thus proved that $X$ satisfies $(iii)$.
		
		The implication $(iii) \Rightarrow (iv)$ is a consequence of the AUS-AUC$^*$ duality and Proposition \ref{prop:Sz-omega}. We indicate its proof. So, assume $(iii)$ is satisfied and let $q$ be the conjugate of $p$. Then, it follows from Corollary \ref{cor:Young2} that there exists $\gamma \in (0,1]$ so that for any $t_0\in (0,1]$ there exists a norm $\abs{\cdot}$ on $X$ satisfying 
		$$\forall x\in X,\ M^{-1}\norm{x}_X \le \abs{x} \le M\norm{x}_X\ \ \text{and}\ \ \forall t\in [t_0,1],\ \bar{\delta}^*_{\abs{\cdot}}(t)\ge \gamma t^q.$$ 
		Fix now $\eps_1,\ldots,\eps_n \in (0,1]$ and pick an equivalent norm $\abs{\cdot}$ as above for $t_0 := \min\{\frac{\eps_1}{4M^2},\dots,\frac{\eps_n}{4M^2}\}$. Assuming that $s_{\eps_1}\dots s_{\eps_n}(B_{X^*})$ is nonempty, it follows immediately that $s_{\eps_1}\dots s_{\eps_n}(MB_{(X^*,\abs{\cdot}^*)})$ is nonempty and by homogeneity so is $s_{\frac{\eps_1}{M}}\dots s_{\frac{\eps_n}{M}}(B_{(X^*,\abs{\cdot}^*)})$. Thus, if we denote temporarily by $\sigma_\eps$ the Szlenk derivation on $X^*$ where the diameter is taken with respect to the norm $\abs{\cdot}^*$, we have that  $\sigma_{\frac{\eps_1}{M^2}}\dots \sigma_{\frac{\eps_n}{M^2}}B_{(X^*,\abs{\cdot}^*)}$ is nonempty. Then, classical manipulations on the Szlenk derivation imply that  
		\begin{equation*}
			\frac12B_{(X^*,\abs{\cdot}^*)} \subset \sigma_{\frac{\eps_1}{4M^2}}\dots \sigma_{\frac{\eps_n}{4M^2}}B_{(X^*,\abs{\cdot}^*)} \subset \prod_{k=1}^n\big(1+\frac{\gamma\eps_k^q}{8^qM^{2q}}\big)^{-1} B_{(X^*,\abs{\cdot}^*)}.
		\end{equation*}
		The argument for the first inclusion can be found in Exercise \ref{ex:value-for-Szlenk}  (see also \cite[proof of Proposition 3.3]{Lancien2006}) and the second inclusion follows from Proposition \ref{prop:Sz-omega}. Then, the $q$-summability of the Szlenk index clearly follows. 
		
		Finally, $(iv) \Rightarrow (i)$ is Proposition \ref{prop:Sz-q-summable->Ap}. 
	\end{proof}
	
	We have chosen to give the details of the renorming characterization of $\sA_p$, because it will be crucial for proving, in Chapter \ref{chapter:Gorelik}, its stability under coarse-Lipschitz equivalences. We shall also prove the stability of $\sN_p$ under coarse-Lipschitz equivalences, similarily based on its renorming characterization that we now reproduce (without proof), to allow the comparison of the various characterizations of the classes $\sA_p$, $\sN_p$ and $\sT_p$ in terms of upper estimates on branches of trees, renorming properties and (when applicable) derivation indices.  
	
	\begin{theo}
		\label{thm:N-theorem-full}
		Fix $p\in(1, \infty)$ and let $q$ be conjugate to $p$. Let $X$ be a Banach space. The following assertions are equivalent. 
		\begin{enumerate}[(i)]
			\item $X\in \sN_p$.
			\item There exists a constant $c>0$ such that for any weak neighborhood basis $D$ at $0$ in $X$, any $n\in \bN$ and any weakly null tree $(x_t)_{t\in D^{\le n}}\subset B_X$, there exists $\tau\in D^n$ such that $\norm{\sum_{i=1}^n x_{\tau_{\restriction_i}}}\le cn^{1/p}$.
			\item There exist a constant $M\ge 1$ and a constant $C>0$ such that for each $t_0\in (0,1]$, there exists a norm $\abs{\cdot}$ on $X$ such that $M^{-1}\norm{x}_X\le \abs{x}\le M \norm{x}_X$ for all $x\in X$ and $\bar{\rho}_{\abs{\cdot}}(t_0)\le C t_0^p$.
			\item The convex Szlenk index of $X$ is of power type $q$.
		\end{enumerate}
	\end{theo}
	
	All the properties in Theorem \ref{thm:T-theorem-full} imply the corresponding properties in Theorem \ref{thm:A-theorem-full}, which in turn imply the corresponding properties in Theorem \ref{thm:N-theorem-full}. This is readily visible at the level of the renorming properties or for the upper estimates on branches of trees. In particular, we have $\sT_p\subseteq \sA_p \subseteq \sN_p$. These inclusions follow as easily from the definitions in terms of games. Indeed, any winning strategy for Player A in the $T(c,p)$-game can be used to create a winning strategy for Player A in the $A(c,p,n)$-game for each $n\in\bN$, by simply playing according to the strategy for $n$ turns. Similarly, any winning strategy for Player A in the $A(c,p,n)$-game is also a winning strategy for Player A in the $N(c,p,n)$ game.
	
	\smallskip
	
	A basic approximation principle is crucial in proving the next proposition, which establishes a quantitative link between the finite $\sN$-game and the infinite $\sT$-game.
	
	\begin{prop}
		\label{prop:upper-finite-trees->upper-infinite-trees}
		Let $X$ be a Banach space and $1<p\le \infty$. If $X\in \sN_p$, then $X\in \sT_r$ for every $1<r<p$.
	\end{prop}
	
	\begin{proof} Clearly, it is enough to show the statement for  $1<p<\infty$. 
		So, assume $p\in (1,\infty)$ and $X\in \sN_p$ and fix $c>\textsf{n}_p(X) \ge 1$ (excluding the trivial case of finite-dimensional spaces).  For each $n\in \bN$, let $\phi_n$ be a winning strategy for Player A in the $N(c,p,n)$ game. Assume also, as we may, that $\phi_n$ takes values in the symmetric weak neighborhoods of $0$. We will use these strategies to choose a winning strategy $\phi$ for Player A in the $T(K,r)$ game, for some $K$ depending on the parameters $c$, $p$ and $r$. We first define $\phi(\emptyset) := \phi_1(\emptyset)$. Then, for $x_1,\dots,x_n \in B_X$, we let 
		\begin{equation*}
			\phi(x_1, \dots, x_n) := \bigcap_{i=1}^n \bigcap_{\sigma \in \Sigma_i(x_1, \dots, x_n)}\phi_i(\sigma),
		\end{equation*} 
		where $\Sigma_i(x_1, \dots, x_n)$ denotes the set of all sequences of the form $(\eps_1 x_{k_1}, \dots,\eps_{s} x_{k_s})$ with $s\le i-1$, $k_1<k_2<\dots<k_s$,  and $\eps_1,\dots,\eps_s \in\{-1,1\}$ together with the empty sequence. This completes the recursive construction of $\phi$. At this point, the introduction of sign choices in the construction of the strategy $\phi$ might seem a bit mysterious, but it will be explained when the time comes. 
		
		Assume now that Player A plays according to $\phi$ and thus forces Player B to choose $x_1,\ldots,x_n,\dots$ so that $x_{n+1}\in \phi(x_1, \ldots, x_n)\cap B_X$ for all $n\ge 0$. In other words, $(x_i)_{i=1}^\infty$ is admissible for $\phi$ and it follows from the definition of $\phi$ and the symmetry of the neighborhoods that for any natural numbers  $m\le k_1<\dots < k_m$ and any $\eps_1,\dots,\eps_m \in \{-1,1\}$, $(\eps_ix_{k_i})_{i=1}^m$ is an admissible sequence for the strategy $\phi_m$ for Player A in the $N(c,p,m)$ game and hence 
		\begin{equation}
			\label{eq:slice}
			\Bigl\|\sum_{i=1}^m\eps_ix_{k_i}\Bigr\|\le cm^{1/p}.
		\end{equation}
		Observe now that, by splitting $\{1,\dots,n\}$ in dyadic intervals, \eqref{eq:slice} holds (up to a multiplicative factor) for \emph{any} $n$-tuple of vectors played by Player B. We detail this observation. Consider $1\le k_1<\dots<k_n$ and $\eps_1,\dots,\eps_n \in \{-1,1\}$. Let $N\in \bN$ be such that $2^{N-1}\le n<2^N$. Then,
		\begin{align*}
			\Bigl\|\sum_{i=1}^n \eps_ix_{k_i}\Bigr\| & \le \sum_{j=0}^{N-2}\Bigl\|\sum_{i=2^{j}}^{2^{j+1}-1} \eps_ix_{k_i}\Bigr\| + \Bigl\| \sum_{i=2^{N-1}}^{n} \eps_ix_{k_i} \Bigr\|\\
			& \stackrel{\eqref{eq:slice}}{\le} \sum_{j=0}^{N-2} c 2^{\frac{j}{p}} + c (n-2^{N-1}+1)^{\frac{1}{p}}\\
			& \le \sum_{j=0}^{N-1} c 2^{\frac{j}{p}} \le c\Big(\frac{2^{N/p}-1}{2^{1/p}-1}\Big)\le c\frac{2^{1/p}}{2^{1/p}-1} n^{1/p}.
		\end{align*}
		Therefore, for all $1\le k_1<\dots<k_n$ and $\eps_1,\dots,\eps_n \in \{-1,1\}$, one has
		\begin{equation}
			\label{eq:slice2}
			\Bigl\|\sum_{i=1}^n \eps_ix_{k_i}\Bigr\| \le d  n^{\frac{1}{p}},
		\end{equation}
		where $d := c2^{1/p}(2^{1/p}-1)^{-1}$.
		If $r<p$, let $q$ be the conjugate exponent of $p$, $s$ be the conjugate exponent of $r$ and fix $x^*\in B_{X^*}$. According to Remark \ref{rem:T_pHahnBanach}, in order to prove that $\phi$ is a winning strategy for Player A in the $T(K,r)$-game (we have not forgotten that we still have to say what is $K$), it is sufficient to show that $\sum_{i=1}^\infty \abs{x^*(x_i)}^s \le K^s$. Since $c>1$, note that for all $i\in \bN$, $\abs{x^*(x_i)}\le 1 \le  d$. For $j\in \bN$, let $I_j := \{i\in \bN \colon \abs{x^*(x_i)}\in (d2^{-j},d2^{-j+1}]\}$. It follows from \eqref{eq:slice2} that for all choices of signs $(\eps_i)_{i\in I_j}$, $\norm{\sum_{i\in I_j} \eps_i x_i}\le d\abs{I_j}^{1/p}$. It is worth pointing out that this last inequality holds regardless of what the $I_j$ are. The reason to pick the $I_j$ in this particular form and the fact that we can pick a choice of signs allows us to derive a crucial upper bound on their size. Indeed, there is always a choice of signs such that 
		\begin{equation*}
			\sum_{i\in I_j}\abs{x^*(x_i)} = x^*\Big(\sum_{i\in I_j} \vep_i x_i \Big) \le \Big\|\sum_{i\in I_j} \vep_i x_i \Big\|\le d\abs{I_j}^{1/p},
		\end{equation*}
		and it follows from the definition of $I_j$ and $q$ that $\abs{I_j}\le 2^{jq}$. 
		Therefore, since $q<s$ we have
		\begin{equation*}
			\sum_{i=1}^\infty \abs{x^*(x_i)}^s \le \sum_{j=1}^\infty 2^{j(q-s)}2^sd^s =: K^s,
		\end{equation*}
		for some constant $K$ depending on the parameters $d$, $q$ and $s$ and therefore on $c$, $p$ and $r$ as anticipated. As mentioned above, this concludes the proof thanks to Remark \ref{rem:T_pHahnBanach}.
	\end{proof}
	
	\begin{rema}
		Another way to look at Proposition \ref{prop:upper-finite-trees->upper-infinite-trees} is to say that if we have upper-$p$-estimates on branches of trees of \emph{finite but arbitrarily large height with flat coefficients}, then we have, for every  $r<p$, upper-$r$-estimates on branches of trees of \emph{infinite height with general coefficients}.  
	\end{rema}
	
	\begin{coro}
		\label{cor:containements1}
		Let $1<p\le \infty$. Then, 
		\begin{equation*}
			\sT_p \subseteq \sA_p \subseteq \sN_p \subseteq \bigcap_{1<r<p}\sT_r.
		\end{equation*}
	\end{coro}
	
	We now have all the tools to obtain an ordinal index characterization of the class $\langle \AUS\rangle$ and Theorem \ref{thm:Sz-omega->AUS}, which was promised at the end of Section \ref{sec:Szlenk-quantitative}.
	
	\begin{theo}
		\label{thm:Sz-omega=AUS}
		Let $X$ be a Banach space. The following assertions are equivalent.
		\begin{enumerate}[(i)]
			\item\label{it:S-w=AUS-ii} $\Sz(X)\le \omega$
			\item\label{it:S-w=AUS-iii} $X$ admits an equivalent $p$-AUS norm, for some $p\in(1,\infty)$. 
			%$X\in \langle p$-$\AUS\rangle$, for some $p\in(1,\infty)$.
			\item\label{it:S-w=AUS-i} $X$ admits an equivalent AUS norm. 
			%$X\in \langle \AUS\rangle$
		\end{enumerate}
	\end{theo}
	
	\begin{proof}
		$(i) \Rightarrow (ii)$ is the fundamental implication. If $\Sz(X)\le \omega$, then $\vep\mapsto \Sz(X,\vep)$ is an integer-valued map and by submultiplicativity (Proposition \ref{prop:submult}) the Szlenk index of $X$ is of power type $q$, for some $q\in(1,\infty)$. It then follows from Proposition \ref{prop:Szlenk-power-summable} that $X$ has $r$-summable Szlenk index for every $r>q$ and hence $X\in \sA_{s}$ for every $s<p$ where $p$ is the conjugate of $q$ (Proposition \ref{prop:Sz-q-summable->Ap}). Noting that $\sA_{s}\subseteq \sN_{s}$, it follows from the approximation argument above (Proposition \ref{prop:upper-finite-trees->upper-infinite-trees}) that $X\in \sT_s$ for all $s<p$. Finally, for every $s<p$, we can renorm $X$ with a $s$-AUS norm according to Theorem \ref{thm:Tp=p-AUS}.\\
		$(ii) \Rightarrow (iii)$ is trivial.\\
		$(iii) \Rightarrow (i)$ follows from Proposition \ref{prop:Sz-omega}.
	\end{proof}

	\begin{rema}
		\label{rem:Szlenk/AUS-power-type}
		The proof of Theorem \ref{thm:Sz-omega=AUS} gives the following quantitative estimate. If $X$ has a Szlenk index of power type $q\in [1,\infty)$, then for every $s<p$, where $\frac1p+\frac1q=1$, $X$ admits an equivalent $s$-AUS norm.
	\end{rema}

	\section{\texorpdfstring{The classes $\sA_\infty$, $\sN_\infty$ and $\sT_\infty$.}{The classes ,  and .}}\label{sec:asymptotic-c_0}
	
	In this section, we discuss the extreme case $p=\infty$ for the classes $\sA_p$, $\sN_p$ and $\sT_p$. 
	In the next theorem, we gather the main characterizations of the class $\sT_\infty$ in the separable setting.
	\begin{theo}
		\label{thm:T-infty-theorem-full}
		Let $X$ be a separable Banach space. The following assertions are equivalent.
		\begin{enumerate}[(i)]
			\item\label{it:T-infty-i} $X\in \sT_\infty$.
			\item\label{it:T-infty-ii} There exists a constant $c>0$ such that for any weak neighborhood basis $D$ at $0$ in $X$ and any weakly null tree $(x_t)_{t\in D^{<\omega}}\subset B_X$, there exists $\tau\in D^\omega$ such that for all $(a_i)_{i=1}^\infty \in c_{00}$,
			\begin{equation*}
				\Big\| \sum_{i=1}^\infty a_ix_{ \tau_{\restriction i} } \Big\| \le c \sup_{i\in \bN}\abs{a_i}.
			\end{equation*}
			\item\label{it:T-infty-iii} $X \in \langle \AUF \rangle$.
			\item\label{it:T-infty-iv} $X$ is isomorphic to a subspace of $\co$.
		\end{enumerate}
		
	\end{theo}
	
	\begin{proof}[Sketch of Proof]
		The equivalence between $(i)$ and $(ii)$ is the case $p=\infty$ of Proposition \ref{prop:Tp-game-trees} and the equivalence between $(iii)$ and $(iv)$ is Corollary \ref{cor:AUF-c0}. The implication  $(iii)$ implies $(ii)$ can be done using the iterated norms technique. The fact that $(ii)$ implies  $(iii)$ is true in general. It was essentially shown in \cite{GKL2001} for separable spaces, and the general case is due to Causey  \cite{Causey2018}.    
	\end{proof}
	
	Concerning the classes $\sA_\infty$ and $\sN_\infty$, we first observe that a classical extreme point argument shows that they coincide.
	
	\begin{prop}
		\label{prop:Ainfty=Ninfty}
		$\sA_\infty = \sN_\infty$.
	\end{prop}
	
	\begin{proof}
		We only need to show that $\sN_\infty \subseteq \sA_\infty$. This argument is similar to and simpler than the argument from the proof of Proposition \ref{prop:upper-finite-trees->upper-infinite-trees}. So, assume that $X \in \sN_\infty$ and let $c\ge 1$ such that Player A has a winning strategy  $\phi_n$ in the $N(c,\infty,n)$-game. Assume again, as we may, that $\phi_n$ takes values in the symmetric weak neighborhoods of $0$. We now describe a winning strategy $\psi_n \colon B_X^{<n} \to \cN_{w}(0)$ in the $A(c,\infty,n)$ game. We first set $\psi_n(\emptyset) := \phi_n(\emptyset)$. Then, for $k<n$ and $x_1,\dots,x_k \in B_X$, we set 
		\begin{equation*}
			\psi_n(x_1,\dots,x_k) := \bigcap_{\eps \in \{-1,1\}^k}\phi_n(\eps_1 x_1,\dots,\eps_k x_k).
		\end{equation*}
		Assume Player A plays according to $\psi_n$ and thus forces Player B to choose $x_{k+1}\in \psi_n(x_1,\dots,x_k)\cap B_X$, for $0 \le k<n$. In other words, assume $(x_i)_{i=1}^n$ is admissible for the strategy $\psi_n$. We note that for any $\eps_1,\ldots,\eps_n \in \{-1,1\}$, $(\eps_ix_i)_{i=1}^n$ is admissible for the strategy $\phi_n$ of Player A in the $N(c,\infty,n)$ game and hence for all $\eps_1,\ldots,\eps_n \in \{-1,1\}$
		$$\Big\|\sum_{i=1}^n\eps_ix_i\Big\|\le c.$$
		Since the unit ball of $\ell_\infty^n$ is the convex hull of $\{-1,1\}^n$, we deduce that for all $a \in \ell_\infty^n$,
		$$ \Big\|\sum_{i=1}^na_ix_i\Big\|\le c\|a\|_\infty.$$
		This shows that $\psi_n$ is a winning strategy for Player A in the $A(c,\infty,n)$ game.
	\end{proof}

	The class $\sA_\infty$ has a simple and elegant reformulation in terms of the notion of asymptotic structure of Maurey, Milman and Tomczak-Jaegermann \cite{MMTJ}. 
	
	\begin{defi} 
		\label{def:asymptotic-c0}
		Let $X$ be a Banach space and $C\ge 1$. We recall that $\cof(X)$ is the set of closed finite-codimensional subspaces of $X$. We say that $X$ is \emph{$C$-asymptotic-$\co$} if for any $n\in \bN$,
		\begin{align*}
			\exists X_1 \in\cof(X),&\ \forall x_1\in B_{X_1},\ \exists X_2\in\cof(X),\ \forall x_2\in B_{X_2},\dots  ,\exists X_n\in\cof(X),\ \forall x_n\in B_{X_n},\\
			&\forall (a_1,\dots,a_n)\in \bR^n,\  \Big\|\sum_{i=1}^na_ix_i\Big\|\le C\max_{1\le i\le n}\abs{a_i}.
		\end{align*}
		Then $X$ is {\em asymptotic-$\co$} if it is $C$-asymptotic-$\co$ for some $C\ge 1$.
	\end{defi}
	
	\begin{prop}
		$X$ is asymptotic-$\co$ if and only if $X \in \sA_\infty$
	\end{prop}
	
	\begin{proof}[Sketch of Proof]
		Recall that $\cN_{w}(0)$ is the set of weak neighborhoods of $0$ in $X$. A standard approximation argument allows us to see that $X$ is asymptotic-$\co$ if and only if there exists $C\ge 1$ such that for any $n\in \bN$,
		\begin{align*}
			\exists U_1 \in \cN_{w}(0),&\ \forall x_1\in B_{X}\cap U_1,\dots,  \exists U_n\in \cN_{w}(0),\ \forall x_n\in B_{X}\cap U_n,\\
			&\forall (a_1,\ldots,a_n)\in \bR^n,\  \Big\|\sum_{i=1}^na_ix_i\Big\|\le C\max_{1\le i\le n}\abs{a_i}.
		\end{align*}
		
		In other words, $X$ is asymptotic-$\co$ if and only if and only if there exists $C\ge 1$ such that Player A has a winning strategy in the $A(C,\infty,n)$-game for all $n\in \bN$, which in turn is equivalent to $X$ being a member of $\textsf{A}_\infty$.    
	\end{proof}
	
	We have already seen that every Banach space with a summable Szlenk index must be in the class $\sA_\infty$. The summability of the Szlenk index is actually equivalent to being in $\sA_\infty$ and also, as we will explain, to the following property of the convex Szlenk index.
	
	\begin{defi}
		\label{def:Szlenk-power-type}
		A Banach space $X$ has \emph{proportional convex Szlenk index} if there exists $K>0$ such that for all $\eps \in(0,1)$, $\Cz(X,\eps)\le K\eps^{-1}$.    
	\end{defi}
	
	Since the convex Szlenk index is a slower derivation index than the Szlenk index, the next proposition, taken from \cite[Lemma 2.4]{GKL2001}, is not that surprising.
	
	\begin{prop}
		\label{prop:Szlenkconvex-proportional}
		If a Banach space has proportional convex Szlenk index, then it has summable Szlenk index. 
	\end{prop}
	
	We start with a lemma regarding the Szlenk derivation that establishes a simple relation between $s_{2n\eps}(K)$ and $s_\eps^n(K)$ where $K$ is a convex weak$^*$ compact subset of $X^*$.
	
	\begin{lemm}
		\label{lem:iteratingCz}
		Let $X$ be a Banach space and $K$ be a convex weak$^*$ compact subset of $X^*$. Then, for any $n\in \bN\cup\{0\}$ and any $\eps>0$, $s_{2n\eps}(K) \subset s_\eps^n(K)$ and also $c_{2n\eps}(K) \subset c_\eps^n(K)$.
	\end{lemm}
	
	\begin{proof} 
		The case $n=0$ is trivial, so let $n\ge 1$ and $x^*\in s_{2n\eps}(K)$. 
		Then, there exists a net $(x^*_\alpha)_{\alpha \in A} \subset K$ weak$^*$ converging to $x^*$ and such that $\norm{x^*-x^*_\alpha}>n\eps$ for all $\alpha \in A$. Observe that for any given sequence $\alpha_1<\dots <\alpha_n$ in $A$, we have, since $K$ is convex, that $\frac{1}{n}(x^*_{\alpha_1}+\dots +x^*_{\alpha_n})\in K$. Since $\inf_{\alpha\in A} \frac1n\norm{x^*-x^*_\alpha}\ge \vep$ and $w^*-\lim_{\alpha_n\in A}\frac{1}{n}(x^*_{\alpha_1}+\dots +x^*_{\alpha_n}) = \frac{1}{n}(x^*_{\alpha_1}+\dots +x^*_{\alpha_{n-1}}+x^*)$, we get that $\frac{1}{n}(x^*_{\alpha_1}+\dots +x^*_{\alpha_{n-1}}+x^*) \in s_\eps(K)$. Iterating this argument we deduce that $x^*\in s_\eps^n(K)$.
		The statement about the convex Szlenk index follows immediately since by definition $c_{2n\eps}(K)$ is the weak$^*$ closed convex hull of $s_{2n\vep}(K)$ and $s_{\vep}^n(K)\subset c_\vep^n(K)$ with the latter set being convex. 
	\end{proof}
	
	\begin{proof}[Proof of Proposition \ref{prop:Szlenkconvex-proportional}]
		Let $X$ be a Banach space and assume that there exists $K>0$ such that $\Cz(X,\eps)\le K\eps^{-1}$ for all $\eps >0$. Let $\eps_1,\dots,\eps_n >0$ such that $s_{\eps_1}\dots s_{\eps_n}(B_{X^*}) \neq \emptyset$ and hence $c_{\eps_1}\dots c_{\eps_n}(B_{X^*})\neq \emptyset$. In order to show that $s := \sum_{i=1}^n \eps_i$ is bounded above by some universal constant, we will use Lemma \ref{lem:iteratingCz} and we need to find an appropriate $\vep>0$. Let us set $\vep := \frac{s}{4n}$ and pick for all $1\le i\le n$, $k_i\in \bN\cup\{0\}$ such that $2k_i\eps \le \eps_i<2(k_i+1)\eps$. Since obviously, $c_{\eps_1}\dots c_{\eps_n}(B_{X^*}) \subset c_{2k_1\eps}\dots c_{2k_n\eps}(B_{X^*})$, we deduce from Lemma \ref{lem:iteratingCz} that $c_{\eps_1}\dots c_{\eps_n}(B_{X^*}) \subset c_\eps^{\sum_{i=1}^n k_i}(B_{X^*})$. Because $\Cz(X,\eps)\le K\eps^{-1}$, we must have $\sum_{i=1}^n k_i \le K\eps^{-1}$. Now recall that the right-hand side inequality in the definition of the $k_i$ tells us that $\sum_{i=1}^n k_i\ge \frac{s}{2\vep} -n$. It follows from our choice of $\eps$ that $n\le \frac{4n}{s}K$ and reorganizing we have $s=\sum_{i=1}^n\vep_i \le 4K$.

		%It follows from  Denote $s=\sum_{i=1}^n \eps_i$ and $\eps= \frac{s}{4n}$. For any $i\in \{1,\ldots,n\}$, there exists $k_i \in \bN\cup\{0\}$ such that $2k_i\eps \le \eps_i<2(k_i+1)\eps$. Note that, summing over $i$ the left hand side inequality, we obtain that $s\le 2\eps\sum_{i=1}^n k_i+2n\eps$, so $s\le 4\eps\sum_{i=1}^n k_i$. Obviously, $c_{\eps_1}\ldots c_{\eps_n}(B_{X^*}) \subset c_{2k_1\eps}\ldots c_{2k_n\eps}(B_{X^*})$. We can now apply Lemma \ref{lem:iteratingCz} to deduce that $c_{\eps_1}\ldots c_{\eps_n}(B_{X^*}) \subset c_\eps^{\sum_{i=1}^n k_i}(B_{X^*})$. It follows that $c_\eps^{\sum_{i=1}^n k_i}(B_{X^*}) \neq \emptyset$ and therefore that $\sum_{i=1}^n k_i \le K\eps^{-1}$ and $s \le 4K$. We have shown that $X$ has a summable Szlenk index. 
	\end{proof}
	
	We are now in a position to compile the main characterizations for the class $\sA_\infty=\sN_\infty$, which includes the remarkable fact that summable Szlenk index implies proportional convex Szlenk index, originally due to Godefroy, Kalton and Lancien \cite{GKL2001}. 
	\begin{theo}
		\label{thm:N-infty-theorem-full} 
		Let $X$ be a Banach space. The following assertions are equivalent.
		\begin{enumerate}[(i)]
			\item\label{it:N-infty-i} $X\in \sN_\infty=\sA_\infty$.
			\item\label{it:N-infty-ii} $X$ is asymptotic-$\co$
			\item\label{it:N-infty-iii} There exists a constant $C>0$ such that for any weak neighborhood basis $D$ at $0$ in $X$, any $n\in \bN$ and any weakly null tree $(x_t)_{t\in D^{\le n}}\subset B_X$, there exists $\tau\in D^n$ such that $\norm{\sum_{i=1}^n x_{\tau_{\restriction_i}}}\le C$.
			\item\label{it:N-infty-iv} There exist a constant $M\ge 1$ and $t_0\in (0,1]$ such that for all $\vep>0$ there exists a norm $\abs{\cdot}$ on $X$ satisfying $M^{-1}\norm{x}_X\le \abs{x}\le M \norm{x}_X$ for all $x\in X$ and $\bar{\rho}_{\abs{\cdot}}(t_0)\le \vep$.
			\item\label{it:N-infty-v} $X$ has proportional convex Szlenk index.
			\item\label{it:N-infty-vi} $X$ has summable Szlenk index.
		\end{enumerate}
	\end{theo}
	
	\begin{proof}[Sketch of Proof]
		The equivalence between $(i)$ and $(ii)$ was proven above and $(i)$ is equivalent to $(iii)$ is by now standard. The renorming implication $(iii) \Rightarrow (iv)$ is Theorem 4.2 in  \cite{GKL2001} for the separable case and extended in \cite{Causey3.5} for the general case. In \cite{Causey3.5}, the renorming technique is similar to the one used in the proof of Theorem \ref{thm:A-theorem-full}. The implication  $(iv) \Rightarrow (iii)$ can be obtained, as usual, by an iterated norm argument. The implication $(iv) \Rightarrow(v)$ is a consequence of the AUS-AUC$^*$ duality and an adaptation of the proof of Proposition \ref{prop:Sz-omega}, while $(v)\Rightarrow (vi)$ is Proposition \ref{prop:Szlenkconvex-proportional}. Finally, $(vi) \Rightarrow (i)$ is the case $p=\infty$ of Proposition \ref{prop:Sz-q-summable->Ap}. 
	\end{proof}
	
	Comparing Theorem \ref{thm:T-infty-theorem-full} and Theorem \ref{thm:N-infty-theorem-full}, it is expected that $\sT_\infty$ does not coincide with $\sA_\infty=\sN_\infty$. 
	This is indeed the case as there exist reflexive infinite-dimensional asymptotic-$\co$ spaces. The first example of such a space that we will denote by $\Tsi^*$, was originally constructed by Tsirelson in \cite{Tsirelson1974}. It was the first example of a Banach space that does not contain any isomorphic copies of $\ell_p$ or $\co$. Soon after, in \cite{FigielJohnson1974}, it became clear that the more natural space to define is $\Tsi$, the dual of $\Tsi^*$, because the norm of this space is more conveniently described. It has since become common to refer to $\Tsi$ as Tsirelson space instead of $\Tsi^*$. It is arguably not an understatement that the appearance of Tsirelson-type spaces disrupted the most optimistic beliefs about the structure theory of Banach spaces, but, on the other hand, it certainly magnified the incredible richness of Banach space theory.
	
	Figiel and Johnson in \cite{FigielJohnson1974} gave an implicit formula that describes the norm of $\Tsi$ as follows. For $E,F\in [\bN]^{<\omega}$ and $n\in \bN$, when one writes $n\le E$ we mean that $n\le \min E$ and $E<F$ means that $\max(E)< \min(F)$. We call a sequence $(E_j)_{j=1}^n$ of finite subsets of $\bN$ \emph{admissible} if $n\le E_1<E_2<\dots<E_n$. Let $(e_j)_{j=1}^\infty$ be the canonical basis of $c_{00}$.  For $x =\sum_{j=1}^\infty \lambda_j e_j\in c_{00}$ and $E\in [\bN]^{<\omega}$ we write $Ex := \sum_{j\in E} \lambda_j e_j$. As it was observed in \cite{FigielJohnson1974}, there is a unique norm $\|\cdot\|_T$ on $c_{00}$ such that for every $x\in c_{00}$:
	\begin{equation}
		\label{eq:Tnorm_V0}
		\norm{x}_T = \max\Big\{\norm{x}_\infty, \frac12 \sup\sum_{j=1}^n \norm{E_jx}_T\Big\},
	\end{equation}
	where the supremum is taken over all $n\in\bN$ and all admissible sequences $(E_j)_{j=1}^n$. The space $\Tsi$ is the completion of $c_{00}$ with this norm and the unit vector basis is a 1-unconditional basis of $\Tsi$. Then, it was proven in \cite{FigielJohnson1974} that $\Tsi$ does not contain a subspace isomorphic to $\ell_1$, which together with the easy observation that $\Tsi$ certainly does not contain a subspace isomorphic to $\co$, yields by James' Theorem  \cite[Theorem 2]{James1950} that $\Tsi$ must be reflexive. The following property of $\Tsi^*$ (see \cite[Lemma 4]{Tsirelson1974}) is essential:
	
	If $(x_j)_{j=1}^n$ is a block sequence with $n\le \supp(x_1)$, then
	\begin{equation}
		\label{eq:dual-of-T}
		\Big\| \sum_{j=1}^n x_j\Big\|_{\Tsi^*}\le 2\max_{1\le j\le n}\norm{x_j}_{\Tsi^*}.
	\end{equation}
	The fact that $\Tsi^*$ is $2$-asymptotic-$\co$ is then a consequence of the above estimate (see Exercise \ref{ex:Tsirelson}). This provides us with an example of a reflexive Banach space that is in $\sA_\infty=\sN_\infty$ but not in $\sT_\infty$.
	
	\begin{coro}
		$\sT_\infty \subsetneq \sA_\infty=\sN_\infty$.
	\end{coro}
	
	\section{\texorpdfstring{Separable determination of the classes $\sA_p$, $\sN_p$ and $\sT_p$.}{Separable determination of the classes ,  and .}}
	\label{sec:sep-det-A-N-T}
	
	In this section, we take the opportunity to provide a unified proof, taken from \cite{CauseyFovelleLancien2023}, of the separable determination of all the properties considered in this chapter. Before stating the main theorem of this section, let us recall that we have already seen (Theorem \ref{thm:separable-determination-Szlenk}) that the Szlenk index, when countable, is separably determined. The fact that having summable Szlenk index and having power type Szlenk index are separably determined is due to Draga and Kochanek in \cite{DragaKochanek2016}. We start with a simple but fundamental statement about selecting weakly null sequences from weakly null nets in AUS-able Banach spaces.
	
	\begin{prop}
		\label{prop:baba} 
		Let $X$ be a Banach space with $\Sz(X)\le \omega$. Let $D$ be a weak neighborhood base at $0$ in $X$. For any weakly null net $(x_U)_{U\in D}\subset B_X$, there exists a function $f\colon \bN \to D$ such that $(x_{f(n)})_{n=1}^\infty$ is a weakly null sequence. 
	\end{prop}
	
	\begin{proof}
		Since $\Sz(X)\le \omega$, Theorem \ref{thm:Sz-omega=AUS} tells us that $X\in \textsf{T}_r$ for some $1<r<\infty$. Let $c>\textsf{t}_r(X)$ and $\psi \colon B_X^{<\omega} \to \cN_{w}(0)$ be a winning strategy for Player A in the $T(c,r)$ game. Let $V_1 := \psi(\emptyset)$ and $U_1 \in D$ such that $x_{U_1}\in V_1 \cap B_X$. We define inductively $V_{n+1} := \psi(x_{U_1},\dots,x_{U_n})$ and $U_{n+1} \in D$ such that $x_{U_{n+1}}\in V_{n+1}\cap B_X$.  Define $f(n) :=U_n$ and note that $N_{s,\omega}( (x_{f(n)})_{n=1}^\infty ) = N_{s,\omega}( (x_{U_n})_{n=1}^\infty ) \le c< \infty$ where $1/r+1/s=1$. This implies that for any $x^* \in X^*$, $(x^*(x_{f(n)}))_{n=1}^\infty \in \ell_s$. In particular, $(x_{f(n)})_{n=1}^\infty$ is weakly null.   
	\end{proof}
	
	\begin{theo}
		\label{thm:separable-determination} 
		Let $X$ be a Banach space with $\Sz(X)\le \omega$ and $p\in (1,\infty]$. Then, 
		\begin{equation*}
			\textsf{\emph{t}}_p(X) := \sup \{\textsf{\emph{t}}_p(E)\colon E\text{\ is a separable subspace of }X\},
		\end{equation*}
		and this supremum is attained, although possibly infinite. 
		
		In particular, if $X$ is a Banach space all of whose separable subspaces lie in $\sT_p$, then $X$ lies in $\sT_p$.  
		
		The same is true of $\textsf{\emph{a}}_p(X)$ and $\textsf{\emph{n}}_p(X)$ and the same conclusion holds for the classes $\sA_p$ and $\sN_p$.
	\end{theo}
	
	\begin{proof} 
		It is clear that $\textsf{t}_p(X)\ge \sup \{\textsf{t}_p(E)\colon E\text{\ is a separable subspace of }X\}$. 
		
		If $c < \textsf{t}_p(X)$, then there exist a weak neighborhood basis $D$ at $0$ and a weakly null tree $(x_t)_{t\in D^{<\omega}}$ such that for each $\tau\in D^\omega$, $N_{q,\omega}( (x_{\tau_{\restriction_i}})_{i=1}^\infty ) >c$. 
		First, we build $\varphi\colon \bN^{<\omega}\to D^{<\omega}$ which preserves lengths and immediate predecessors such that $(x_{\varphi(t)})_{t\in \bN^{<\omega}}$ is weakly null. We define $\varphi(t)$ by induction on $\abs{t}$. By Proposition \ref{prop:baba} applied to $(x_{(U)})_{U\in D}$, there exists $f\colon\bN\to D$ such that $(x_{(f(n))})_{n=1}^\infty$ is weakly null.  Define $\varphi((n)) := (f(n))$. Next, if $\varphi(t)$ has been defined, apply Proposition \ref{prop:baba} to $(x_{\varphi(t)\smallfrown (U)})_{U\in D}$ to select $g\colon \bN\to D$ such that $(x_{\varphi(t)\smallfrown (g(n))})_{n=1}^\infty$ is weakly null and let $\varphi(t\smallfrown(n)) := \varphi(t)\smallfrown (g(n))$. This completes the construction. Next, define $y_t := x_{\varphi(t)}$, for $t \in \bN^{<\omega}$. It follows from our construction that for any $\sigma\in \bN^\omega$, there exists a unique $\tau\in D^\omega$ such that $\varphi(\sigma_{\restriction_i})=\tau_{\restriction_i}$ for all $i\in\bN$, so that 
		\begin{equation*}
			N_{q,\omega}( (y_{\sigma_{\restriction_i}})_{i=1}^\infty ) = N_{q,\omega}( (x_{\tau_{\restriction_i}})_{i=1}^\infty ) > c.
		\end{equation*}  
		Therefore, if $F$ is the closed linear span of $(y_t)_{t\in \bN^{<\omega}}$, then $F$ is separable and $\textsf{t}_p(F)>c$. This shows that $\textsf{t}_p(X)\le \sup \{\textsf{t}_p(E)\colon E\text{\ is a separable subspace of }X\}$. 
		
		Next, let $R$ be the set of rational numbers $r$ such that $\textsf{t}_p(X)>r$. For each $r\in R$, let $E_r$ be a separable subspace of $X$ such that $\textsf{t}_p(E_r)>r$ and let $E$ be the closed linear span of $\cup_{r\in R}E_r$. Then, $E$ is separable and $\textsf{t}_p(E)=\textsf{t}_p(X)$, so the supremum is attained. 
		
		If $X$ is a Banach space all of whose separable subspaces lie in $\textsf{T}_p$, then $\textsf{t}_p(X) = \sup \{\textsf{t}_p(E): E\text{\ is a separable subspace of }X\}$ must be finite. Indeed, if the supremum were infinite, then since it is attained, there would exist some separable subspace $E\subset  X$ such that $\textsf{t}_p(E)=\infty$ and $E$ does not belong to $\textsf{T}_p$.  
		
		The arguments for $\textsf{a}_p(X)$, $\textsf{n}_p(X)$, $\textsf{A}_p$ and $\textsf{N}_p$ are similar.
		
	\end{proof}
	
	\begin{rema} 
		It is worth noticing that this provides an alternative proof of the fact that $\Sz(X)\le \omega$ is separably determined. Indeed, assume that $\Sz(X) >\omega$. Then, for any  $p\in \bQ \cap (1,\infty)$, $X$ does not belong to $\textsf{T}_p$.  So for any $p\in \bQ \cap (1,\infty)$, there exists a separable subspace $E_p$ of $X$ so that  $E_p$ is not in $\textsf{T}_p$. Then, the closed linear span $E$ of these $E_p$ is a separable subspace of $X$ which does not belong to $\cup_{1<p<\infty}\textsf{T}_p$. So, by Theorem \ref{thm:Sz-omega=AUS}, $\Sz(E)>\omega$.
	\end{rema}

	\section{\texorpdfstring{Property $(\beta)$ and simultaneously AUC and AUS renormings.}{Property and simultaneously AUC and AUS renormings.}}
	
	We saw in Section \ref{sec:Rolewicz} that a Banach space with property $(\beta)$ is automatically reflexive and asymptotically uniformly convex. However, there are examples of Banach spaces that are reflexive and asymptotically uniformly convex (or with property $(\beta)$ or even uniformly convex!) but that are not asymptotically uniformly smooth (cf. Example 5 in \cite{Kutzarova1990}). The first result of this section states that a Banach space with property $(\beta)$ admits an equivalent norm that is asymptotically uniformly smooth.
	
	\begin{theo}
		\label{thm:beta->AUS}
		Let $X$ be a Banach space with property $(\beta)$. Then, $X$ admits an equivalent asymptotically uniformly smooth norm. 
		%(or equivalently $\Sz(X)\le \omega$).
	\end{theo}
	
	This result goes back to a paper by D. Kutzarova \cite{Kutzarova1990}, where it was proved that a Banach space with a Schauder basis and property $(\beta)$ has an equivalent AUS norm. The general case (separable or not) follows, for instance, from the fact that a reflexive Banach space $X$ satisfies $\Sz(X)\le \omega$ if and only if all its subspaces with a Schauder basis do (see \cite{DKLR2017} for a proof). 
	In \cite{Kutzarova1990}, the focus was on the notion of \emph{nearly uniformly smooth} spaces (NUS in short) and it was actually shown that a Banach space with a Schauder basis and property $(\beta)$ admits an equivalent NUS norm. The NUS property, which we shall not define here, is due to S. Prus \cite{Prus1989} and it turns out that a Banach space is NUS if and only if it is reflexive and AUS.
	
	The proof of Theorem \ref{thm:beta->AUS} that we are about to give does not rely on the fact that having Szlenk index at most $\omega$ is determined by Schauder basic sequences. Instead, we have chosen to give a ``coordinate-free'' version of Kutzarova's argument and we prove that a Banach space with property $(\beta)$ belongs to $\sN_p$ for some $p>1$. The conclusion will then follow from the results in Section \ref{sec:A-N-T}. 
	
	%Recall for the convenience of the reader that Mazur's lemma tells us that for any finite-dimensional space $E$ in an infinite-dimensional Banach space $X$ and any $\vep>0$, there is a weakly open neighborhood $U$ of $0$ such that for all $\lambda\in \bR$, $e\in E$ and $z\in U$, we have $(1+\vep)\norm{e+\lambda z}\ge \norm{e}$. 

	\begin{rema}
		\label{rem:Mazur}
		It is a well-known application of Mazur's lemma that for every $\vep>0$ and every normalized weakly null sequence in an infinite-dimensional Banach space, one can extract a subsequence that is $(1+\vep)$-basic. A similar argument together with Lemma \ref{lem:pruning1}, insures that if $(x_t)_{t\in D^{\le n}}$ is a normalized weakly null tree (where $D$ is a weak neighborhood basis of $0$), then there exists a pruning $\varphi \colon D^{\le n} \to D^{\le n}$ such that $(x_{\varphi(\tau)_{\restriction 1}},\dots,x_{\varphi(\tau)_{\restriction n}})$ is $(1+\vep)$-basic for all $\tau \in D^n$.
	\end{rema}
	
	A key ingredient leading up to Theorem \ref{thm:beta->AUS} is the following asymptotic analog of a result of Gurarii-Gurarii \cite{GurariiGurarii1971}. 
	
	\begin{lemm}
		\label{lem:beta->AUS} 
		Let $X$ be a Banach space and assume that $\bar{\beta}_X(\vep)>0$ for some $\vep\in(0,1)$. Then, there exists $\delta \in (0,1)$ such that for any weak neighborhood basis $D$ of $0$ in $X$, any $n\ge 2$ and any weakly null tree $(x_t)_{t\in D^{\le n}}$ in $S_X$, there exists a pruning $\psi\colon D^{\le n} \to D^{\le n}$ such that for all $\tau \in D^n$, all $1\le i< n$, the following vectors are well defined:
		\begin{equation*}
			u^i_{\psi(\tau)} := \sum_{j=1}^ix_{\psi(\tau)_{\restriction j}},\ 
			y^i_{\psi(\tau)} := \frac{u^i_{\psi(\tau)}}{\|u^i_{\psi(\tau)}\|},\ v^i_{\psi(\tau)} := \sum_{j=i+1}^nx_{\psi(\tau)_{\restriction j}},\ z^i_{\psi(\tau)} := \frac{v^i_{\psi(\tau)}}{\|v^i_{\psi(\tau)}\|}
		\end{equation*}
		and
		\begin{equation*}
			\|y^i_{\psi(\tau)}+z^i_{\psi(\tau)}\|\le 2-2\delta.
		\end{equation*}
	\end{lemm}
	
	\begin{proof}
		Since $\bar{\beta}_X(\vep)>0$, there is $\delta\in(0,1)$ such that for all $x, z_1,z_2,\dots$ in $S_X$, if $\inf_{i\neq j} \norm{z_i-z_j}\ge \vep$, then $\norm{x+z_{i_0}}\le 2-2\delta$ for some $i_0\in \bN$.
		
		Let $D$ be a weak neighborhood basis of $0$ in $X$, $n\ge 2$ and $(x_t)_{t\in D^{\le n}}$ a weakly null tree in $S_X$. It is clearly enough to find a pruning satisfying the conclusion of the lemma for a given $i\in \{1,\ldots,n-1\}$. So, let us fix $1\le i <n$. After a first pruning, we may assume by  Remark \ref{rem:Mazur} that $(x_{\tau_{\restriction 1}},\dots,x_{\tau_{\restriction n}})$ is $2$-basic for all $\tau \in D^n$. In particular, $\norm{u_\tau^i}\ge \frac12$ and $\norm{v_\tau^i}\ge \frac12$ for all $\tau \in D^n$, so $y^i_{\psi(\tau)}$ and $z^i_{\psi(\tau)}$ are well defined. Before we proceed to the heart of the proof, we recall that, by Lemma \ref{lem:pruning1}, it is enough to find one branch $\tau \in D^n$ such that $\norm{y^i_{\tau}+z^i_{\tau}}\le 2-2\delta.$
		
		Let us fix $\sigma \in D^i$. We will build inductively $\tau_k := \sigma \smallfrown \sigma_k \in D^n$, for $k\in \bN$ such that the sequence $(z_{\tau_k}^i)_{k=1}^\infty$ is $\eps$-separated. Pick any $\sigma_1 \in D^{n-i}$ and let $\tau_1 :=\sigma \smallfrown \sigma_1$. Assume that $\sigma_1,\dots,\sigma_k \in D^{n-i}$ have been constructed. Since $\eps\in (0,1)$, we can find $\eta>0$ so that $1-2n\eta>\eps$. For each $l\in \{1,\dots,k\}$, pick $x^*_l\in S_{X^*}$ such that $x^*_l(z_{\tau_l}^i)=1$. Let now $V := \{x\in X\colon \abs{x^*_l(x)}<\eta\ \text{for all}\ 1\le l \le k\}$. Since $V$ is a weak neighborhood of $0$ and our tree is weakly null, we can pick recursively $U_{i+1},\dots, U_n \in D$ so that for all $j\ge i+1$ $x_{\sigma\smallfrown (U_{i+1},\dots, U_j)}\in V$ and we let $\tau_{k+1} := \sigma\smallfrown (U_{i+1},\dots, U_n)$. Then, for all $1\le l \le k$, $|x^*_l(v^i_{\tau_{k+1}})|\le n\eta$ and $|x^*_l(z^i_{\tau_{k+1}})|\le 2n\eta$. It follows that for all $1\le l \le k$, $\|z^i_{\tau_{k+1}}-z^i_{\tau_{l}}\|>\eps$. This finishes the inductive construction of $(\tau_k)_{k=1}^\infty$. We now note that for all $k\in \bN$, $y_{\tau_k}^i=\big(\sum_{j=1}^i x_{\sigma_{\restriction j}}\big)\|\sum_{j=1}^i x_{\sigma_{\restriction j}}\|^{-1}:=x$. Then, it follows from the definition of the $(\beta)$-modulus recalled at the beginning of the proof that there exists $k_0\in \bN$ such that $\|x+z_{\tau_{k_0}}^i\|=\|y_{\tau_{k_0}}^i+z_{\tau_{k_0}}^i\|\le 2-2\delta$. This concludes the proof of this lemma.
	\end{proof}

	We shall also need the following elementary fact.
	\begin{lemm}
		\label{lem:beta+calculus} 
		Let $\delta \in (0,1)$. Then, there exists $p>1$ and $\nu \in (0,\frac12)$ such that whenever $y,z \in S_X$ satisfy $\norm{y+z}\le 2-2\delta$, then, $\norm{y+tz}^p\le 1+t^p$ for all $t\in [1-\nu,1+\nu]$
	\end{lemm}
	
	\begin{proof} 
		There exists $p>1$ such that $(2-\delta)^p<2$. Then, the functions $f\colon t\mapsto (1+t-\delta)^p$ and $g\colon t\mapsto 1+t^p$ are continuous on $[0,\infty)$ and satisfy $f(1)<g(1)$. So, there exists $\nu \in (0,\frac12)$ such $f(t)\le g(t)$ for all $t\in [1-\nu,1+\nu]$. Now let $y,z \in S_X$ be such that $\norm{y+z}\le 2-2\delta$. Then, for all $t\in [1-\nu,1+\nu]$:
		\begin{equation*}
			\norm{y + tz}^p \le \Big(\frac12\norm{y} + \big(t-\frac12\big)\norm{z} + \frac12 \norm{ y + z}\Big)^p \le (1+t-\delta)^p\le 1+t^p.
		\end{equation*}
	\end{proof}
	
	We are now ready to prove the theorem below which immediately implies Theorem \ref{thm:beta->AUS} via Proposition \ref{prop:upper-finite-trees->upper-infinite-trees} and Theorem \ref{thm:T-theorem-full}. The result we present now is taken from \cite{BaudierLancien2025}.
	
	\begin{theo}
		Let $X$ be a Banach space. If $\bar{\beta}_X(\vep)>0$ for some $\vep\in(0,1)$, then there exists $p\in(1,\infty)$ such that $X\in \sN_p$.
	\end{theo}
	
	\begin{proof} 
		Let $\delta \in (0,1)$ be given by Lemma \ref{lem:beta->AUS}. Let $p>1$ and $\nu \in (0,\frac12)$ be the constants associated with $\delta$ through Lemma \ref{lem:beta+calculus}. We set $C:=\frac{3}{\nu}$. We will show by induction on $n$ that for any weak neighborhood basis $D$ at $0$ in $X$, any $n\in \bN$ and any weakly null tree $(x_t)_t\in D^{\le n}$  in $S_X$, there exists $\tau \in D^n$ such that  $\norm{\sum_{j=1}^n x_{\tau_{\restriction j}}}\le Cn^{1/p}$. Since $C\ge 3$, this is clearly true for $n\in \{1,2,3\}$. So, assume $n>3$ and that our induction hypothesis is true for values smaller than $n$. Consider a weak neighborhood basis $D$ of $0$ in $X$ and a weakly null tree $(x_t)_{t\in D^{\le n}}$  in $S_X$. For $1\le i<n$ and $\tau \in D^n$, we adopt the notation $u_\tau^i,y_\tau^i,v_\tau^i,z_\tau^i$ from Lemma \ref{lem:beta->AUS} that we complete with $u^0_\tau:=v^n_\tau:=0$ and $v^0_\tau:=u^n_\tau:=\sum_{j=1}^nx_{\tau_{\restriction j}}$.
		
		It clearly follows from our induction hypothesis and the pruning Lemma \ref{lem:pruning1} that, after pruning, we may assume that for all $1\le i<n$ and all $\tau \in D^n$, $\|u_\tau^i\|\le Ci^{1/p}$ and $\|v_\tau^i\|\le C(n-i)^{1/p}$. Next, we use Lemma \ref{lem:beta->AUS} to justify that, again after pruning, we may assume that for all $1\le i<n$ and all $\tau \in D^n$, $\|y_\tau^i+z_\tau^i\|\le 2-2\delta$. 
		
		Fix now $\tau \in D^n$. The proof will be complete if we show that there exists $1\le i<n$ such that $\|u_\tau^i+v_\tau^i\|\le Cn^{1/p}$. Note that $\|u^0_\tau\|<\|v^0_\tau\|$, $\|u^n_\tau\|>\|v^n_\tau\|$ and for all $0\le i<n$, $\big|\|u^{i+1}_\tau\|-\|u^{i}_\tau\|\big|\le 1$ and $\big|\|v^{i+1}_\tau\|-\|v^{i}_\tau\|\big|\le 1$. It is an easy exercise to check that this implies the existence of $0\le i_0\le n$ such that  $\big|\|u^{i_0}_\tau\|-\|v^{i_0}_\tau\|\big|\le 1$. Assume first that $\|u^{i_0}_\tau\|\le \frac{C}{3}$ or $\|v^{i_0}_\tau\|\le \frac{C}{3}$. Then, $\|u_\tau^{i_0}+v_\tau^{i_0}\|\le \frac{2C}{3}+1 \le C \le Cn^{1/p}$. So, we can assume that $\|u^{i_0}_\tau\|> \frac{C}{3}$ and $\|v^{i_0}_\tau\|> \frac{C}{3}$, which implies that $1\le i_0<n$ and also that $(1-\nu)\|u^{i_0}_\tau\|\le \|v^{i_0}_\tau\| \le (1+\nu)\|u^{i_0}_\tau\|$. It now follows from the fact that $\|y_\tau^{i_0}+z_\tau^{i_0}\|\le 2-2\delta$, Lemma \ref{lem:beta+calculus} and an homogeneity argument that 
		$$\|u_\tau^{i_0}+v_\tau^{i_0}\|^p\le \|u_\tau^{i_0}\|^p+\|v_\tau^{i_0}\|^p \le C^p\big(i_0+(n-i_0)\big)=C^pn.$$
		This finishes our inductive proof. 
	\end{proof}

	\begin{rema} 
		We have chosen to present a proof using the general pruning Lemma \ref{lem:pruning1} and trees indexed by weak neighborhood bases directed by reverse inclusion. Alternatively, we could have used the separable reduction results in order to work with countably branching trees and use the classical Ramsey Theorem for colorings in $[\bN]^k$ in place of the general pruning lemma. 
	\end{rema}
	
	\begin{rema} Note that it follows from Proposition \ref{prop:beta->reflexive}, Theorem \ref{thm:beta->AUC}, Theorem \ref{thm:beta->AUS}, the AUS-AUC$^*$ duality and Proposition \ref{prop:Sz-omega} that a Banach space $X$ with property $(\beta)$ satisfies $\max\{\Sz(X), \Sz(X^*)\}\le \omega$. Therefore, if a Banach space has an equivalent norm with property $(\beta)$, then  $\Sz(X)\le \omega$ and $\Sz(X^*)\le \omega$. In fact, this result can be derived from a different argument based on metric considerations and the geometry of graphs. These results appear in Chapter \ref{chapter:trees}, to which we refer for the relevant definitions. The remark is then an immediate consequence of the following two statements. First, if a Banach space $X$ is such that $\Sz(X) >\omega$ or $\Sz(X^*)> \omega$, then the $\omega$-regular tree of infinite height $T_\infty^\omega$, equipped with its graph metric, bi-Lipschitzly embeds into $X$. Second, if a Banach space $X$ has an (equivalent) norm with property $(\beta)$, then $T_\infty^\omega$ does not bi-Lipschitzly embed into $X$. In fact, Chapter \ref{chapter:trees} contains more and the bi-Lipschitz embeddability of $T_\infty^\omega$ characterizes Banach spaces with property $(\beta)$ among reflexive spaces. 
	\end{rema}
	
	A natural question in renorming theory is whether a Banach space admitting an equivalent AUS norm and an equivalent AUC norm can be renormed with an equivalent norm that is simultaneously AUS and AUC. If the space is reflexive, this problem can be recast as a Baire category problem and the Baire category Theorem can be invoked to answer positively this question (see Exercise \ref{ex:Baire-AUS-AUC}). We can also ask a more quantitative question.
	%namely, whether $X\in \langle p$-$\AUS\rangle \cap \langle q$-$\AUC\rangle$ implies that $X\in \langle p$-$\AUS~\&~q$-$\AUC \rangle$.
	
	\begin{prob}
		\label{pb:AUS-AUC-sim-renormings}
		Let $1<p \le q <\infty$. 
		If a Banach space has an equivalent $q$-asymptotically uniformly convex norm and an equivalent $p$-asymptotically uniformly smooth norm, does it admit an equivalent norm that is both $q$-asymptotically uniformly convex and $p$-asymptotically uniformly smooth?
	\end{prob}
	
	It follows from Theorem \ref{thm:subspaces-of-lp-sums} that this quantitative problem has a positive solution for separable reflexive spaces when $p=q\in(1,\infty)$. The general quantitative problem, i.e. when $1<p \le q< \infty$, also has a positive solution for reflexive Banach space with an FDD (see \cite[Proposition 2.10]{JLPS2002} and the discussion leading to it where this result is essentially credited to Prus \cite{Prus1989}). The FDD assumption was eventually removed in \cite{OdellSchlumprecht2002} by reducing the problem to the FDD case using, among other things, ideas reminiscent of Zippin's embedding theorem. 
	
	The renorming in the next proposition is quantitatively not optimal, at least in the reflexive case, but the proof only uses elementary tools which are inspired by the Asplund averaging method initiated by E. Asplund in \cite{Asplund1967}. Moreover, the reflexivity assumption is not needed. We believe this had not been noticed before. 
	
	\begin{prop}
		\label{prop:sim-AUS-AUC}
		Let $X$ be a  Banach space and $1<p\le q <\infty$. If $X$ admits an equivalent $p$-AUS norm and an equivalent $q$-AUC norm, then $X$ admits an equivalent norm which is simultaneously $p$-AUS and $q(p+2)$-AUC.
	\end{prop}
	
	\begin{proof} 
		We shall adapt the proof of Proposition IV.5.2 in \cite{DGZ1993}, which is due to John and Zizler \cite{JohnZizler1979}. The technique is a variant of the so-called Asplund averaging method initiated by E. Asplund in \cite{Asplund1967}. So, assume that there exist $a,b>0$ such that $X$ admits an equivalent norm $N$ satisfying $\bar{\delta}_N(t)\ge at^q$ for all $t\in (0,1)$ and an equivalent norm $M$ satisfying $\bar{\rho}_M(t)\le bt^p$ for all $t\in (0,\infty)$. Assume also, as we may, that for all $x\in X$, $N(x)\le M(x)\le CN(x)$ with $C\ge 1$. Note that the AUS-AUC$^*$ duality  ensures the existence of $b_1>0$ so that $\bar{\delta}_{M}^*(t)\ge b_1t^{p'}$, where $p'$ is the conjugate exponent of $p$. We now define a sequence of equivalent norms on $X^*$ as follows. For all $n\ge 1$ and all $x^*\in X^*$, let 
		\begin{equation*}
			\norm{x^*}_{*,n} := N^*(x^*)+n^{-1}M^*(x^*),
		\end{equation*}
		where $M^*$ and $N^*$ are the dual norms of $M$ and $N$, respectively. Since $M^*$ and $N^*$ are weak$^*$ lower semi-continuous, so is $\norm{\cdot}_{*,n}$. It is therefore the dual norm of an equivalent norm that we denote by $\norm{\cdot}_n$ on $X$. 
		It is elementary to verify that for all $n\ge 1$ and all $x^*\in X^*$,
		\begin{equation}
			\label{eq:equiv-norm}
			\max\{ (1+n^{-1})M^*(x^*), N^*(x^*) \} \le \norm{x^*}_{*,n}\le \min\{(C + \frac1n)M^*(x^*), (1+\frac1n)N^*(x^*)\}.
		\end{equation}
		The AUC property is clearly preserved under sums of norms. More precisely, if $\|\ \|_1$ and $\|\ \|_2$ are equivalent norms and $\|\ \|_2$ is AUC, then $\|\ \|_1+\|\ \|_2$ is AUC. On the quantitative level, since $\norm{\cdot}_n^*=N^*+\frac1nM^*$, it is easy to check that there exists $b_2>0$ such that for all $n\ge 1$ and all $t\in (0,1)$,  $\bar{\delta}^*_{\norm{\cdot}_n}(t)\ge b_2 n^{-1}t^{p'}$ (we leave the details to the reader). Then, using again the AUS-AUC$^*$ duality, we deduce the  existence of $b_3>0$ such that $\bar{\rho}_{\norm{\cdot}_n}(t)\le n^{p-1} b_3 t^p$ for all $n\in \bN$ and all $t\in (0,\infty)$. Finally, for all $x\in X$, we set
		\begin{equation}
			\label{eq:def-equiv-norm}
			\abs{x} := \sum_{n=1}^\infty n^{-(p+1)} \norm{x}_n.
		\end{equation}
		Letting $s := \sum_{n=1}^\infty n^{-(p+1)}$, it follows from \eqref{eq:equiv-norm} and the definition of $\norm{\cdot}_n$ that for all $x\in X$,
		\begin{equation*}
			\frac{s}{2}N(x) \le \abs{x} \le sN(x),
		\end{equation*}
		and hence, for all $n\ge 1$, we have
		\begin{equation}
			\label{eq:equiv-norm-2}
			\frac{s}{2}\norm{x}_n \le \abs{x} \le 2s \norm{x}_n.
		\end{equation}
		In particular, $\abs{\cdot}$ is an equivalent norm on $X$ and we will now verify the smoothness and convexity properties of $\abs{\cdot}$.
		\begin{itemize}
			\item $\abs{\cdot}$ is $p$-AUS: Let $x\in X$ so that $|x|=1$ and $t\in (0,1)$. For every $n\ge 1$, since $\bar{\rho}_{\norm{\cdot}_n}(t)\le n^{p-1} b_3 t^p$, there is $Y_n\in \cof(X)$ such that for all $y\in Y_n$, 
			\begin{align}
				\notag \norm{x+y}_n &= \norm{x}_n\Big\|\frac{x}{\norm{x}_n}+\frac{y}{\norm{x}_n}\Big\|_n\\
				\notag            &\le \norm{x}_n(1+n^{p-1}b_3 \norm{x}_n^{-p}\norm{y}_n^p)\\
				\label{eq:aux-n}       & \le \norm{x}_n + n^{p-1}b_3 \norm{x}_n^{1-p}\norm{y}_n^p.
			\end{align}
			Pick $n_0\ge 1$ such that $4s^{-1}\sum_{n>n_0}n^{-(p+1)}\le Kt^p$ with $K := b_3 4^ps^{-1}\sum_{n=1}^{\infty} n^{-2}\in(0,\infty)$. Then, for all $y\in Y := \cap_{n=1}^{n_0} Y_n \in \cof(X)$ with $\abs{y}\le t$, we have
			\begin{align*}
				\abs{x+y} & = \sum_{n=1}^\infty n^{-(p+1)}\norm{x+y}_n \\
				& \stackrel{\eqref{eq:aux-n}}{\le}  \sum_{n=1}^{n_0} n^{-(p+1)}\norm{x}_n + \sum_{n=1}^{n_0} n^{-2} b_3 \norm{x}_n^{1-p}\norm{y}_n^p + \sum_{n=n_0+1}^\infty n^{-(p+1)}\norm{x+y}_n \\
				& \stackrel{\eqref{eq:def-equiv-norm}\land \eqref{eq:equiv-norm-2}}{\le} \abs{x} + \sum_{n=1}^{n_0} n^{-2} b_3 (2s)^{p-1} (2s^{-1}t)^p + \sum_{n=n_0+1}^\infty n^{-(p+1)} 4s^{-1} \\
				& \le 1 +2Kt^p,
			\end{align*}
			where the last inequality follows from the definition of $K$ and the choice of $n_0$.
			
			\medskip
			\item $\abs{\cdot}$ is $4q$-AUC: It follows from \eqref{eq:equiv-norm} that for all $n\ge 1$ and all $x\in X$,
			\begin{equation}
				\label{eq:equiv-N-norm-n}
				\Big(1-\frac1n\Big)N(x)\le \Big(1+\frac1n\Big)^{-1} N(x)\le \norm{x}_n \le N(x).
			\end{equation}
			Consider $x\in X$ so that $|x|=1$ and $t\in (0,\frac12)$. Since $\bar{\delta}_N(t)\ge at^q$, there is $Y\in \cof(X)$ such that for all $y\in Y$ with $N(y)\le N(x)$, we have:
			\begin{align}
				\label{eq:N-q-AUC}
				N(x+y) \ge N(x)(1+ a N(x)^{-q}N(y)^q).
			\end{align}
			Given $\alpha\in(0,\infty)$ a large enough number (depending only on $s$, $a$ and $q$) to be chosen later, pick $n_0\ge 2$ so that $n_0-1\le \alpha t^{-q}\le n_0$. Then, for all $y\in Y$ such that $\abs{y}=t$, we have that $N(y)\le 2s^{-1}\abs{y}\le s^{-1}\abs{x}\le N(x)$, so that for all $n\ge n_0$,
			\begin{align*}
				\norm{x+y}_{n}  & \stackrel{\eqref{eq:equiv-N-norm-n}}{\ge} \Big(1-\frac1n\Big)N(x+y) \\
				& \stackrel{\eqref{eq:N-q-AUC}}{\ge} \Big(1-\frac1n\Big)N(x)(1+ a N(x)^{-q}N(y)^q)\\
				& \stackrel{\eqref{eq:equiv-N-norm-n}}{\ge} \Big(1-\frac1n\Big)\norm{x}_n (1 + a (2\norm{x}_n)^{-q}\norm{y}_n^q)\\         
				& \stackrel{\eqref{eq:equiv-norm-2}}{\ge} \Big(1-\frac1n\Big)\norm{x}_n (1 + a 4^{-q}s^q(2s)^{-q}t^q)\\  
				& \ge \Big(1-\frac{t^q}{\alpha}\Big)\big(1+ a 8^{-q}t^q\big)\norm{x}_{n},
			\end{align*}
			where in the last inequality, we use the definition of $n_0$ and the fact that $n\ge n_0$.
			Therefore, assuming, as we may, that $\alpha\ge \max\{2a^{-1} 8^q, 2\}$ and letting $D := \frac{a}{2} 8^{-q}(2s)^{-1}$ we have that for all $n\ge n_0$
			\begin{equation}
				\label{eq:eq3}
				\norm{x+y}_{n} \ge \norm{x}_n + Dt^q.
			\end{equation}
			Now, given $\vep_1, \dots, \vep_{n_0-1}\in (0,1)$, the standard Mazur technique applied to the equivalent renormings $(X,\norm{\cdot}_1), \dots, (X,\norm{\cdot}_{n_0-1})$, allows us to construct $Z\in \cof(X)$ such that $Z\subset Y$ and for all $n\in \{1, \dots, n_0-1\}$ and all $z\in Z$,
			\begin{equation}
				\label{eq:mazur}
				\norm{x+z}_n\ge (1-\eps_n)\norm{x}_n.
			\end{equation}
			Therefore, for all $z\in Z$ with $\abs{z}=t$ we have
			\begin{align*}
				\abs{x+z} & = \sum_{n=1}^\infty n^{-(p+1)}\norm{x+z}_n \\
				& = \underbrace{\sum_{n=1}^{n_0-1} n^{-(p+1)}\norm{x+z}_n}_{A} + n_0^{-(p+1)}\norm{x + z}_{n_0} + \sum_{n=n_0+1}^{\infty} n^{-(p+1)}\norm{x+z}_n.
			\end{align*}
			Let us estimate $A$ separately. It follows from \eqref{eq:mazur}, \eqref{eq:equiv-norm-2} and the fact that $\abs{x}=1$ and $\abs{y}=t\le \frac12$ that 
			\begin{equation}
				A \ge \sum_{n=1}^{n_0-1} n^{-(p+1)}(\norm{x}_n -2s^{-1}\vep_n),
			\end{equation}
			and hence 
			\begin{equation}
				\label{eq:alpha}
				A \ge \sum_{n=1}^{n_0-1} n^{-(p+1)}\norm{x}_n - \frac{D}{2n_0^{p+1}}t^q.
			\end{equation}
			Assuming that we had chosen the $\vep_1, \dots,\vep_{n_0-1}$ small enough, which we can do since the choices depend only on $s$, $a$, $q$ and $t$.
			Consequently,
			\begin{align*}
				\abs{x+z} & \stackrel{\eqref{eq:alpha}\land \eqref{eq:eq3}}{\ge}  \sum_{n=1}^{n_0-1} n^{-(p+1)}\norm{x}_n - \frac{D}{2n_0^{p+1}}t^q+ n_0^{-(p+1)} \norm{x}_{n_0}\\ 
				&+\frac{D}{n_0^{p+1}}t^q+  \sum_{n=n_0+1}^{\infty} n^{-(p+1)}\norm{x}_n\\
				& \ge \abs{x} + \frac{D}{2n_0^{p+1}}t^q\\
				& \ge 1 + \frac{D}{2(\alpha+1)^{p+1}}t^{q(p+2)},
			\end{align*}
			where in the last inequality we used the fact that $n_0\le \alpha t^{-q}+1 \le t^{-q}(\alpha+1)$.
		\end{itemize}
	\end{proof}

	\begin{rema} 
		In the reflexive case, by duality, we could deduce from Proposition \ref{prop:sim-AUS-AUC} the existence of an equivalent norm which is $q$-AUC and AUS with a power type in $(1,p)$. Still in the reflexive setting, a Baire category argument also allows to obtain a simultaneous AUC and AUS renorming (see Exercise \ref{ex:Baire-AUS-AUC}). But the very nature of the Baire category Theorem makes it quite hopeless to get a quantitative result such as Proposition \ref{prop:sim-AUS-AUC} with this kind of tool.
	\end{rema}
	
	The next statement summarizes the main linear characterizations of the class of Banach spaces with an equivalent norm with property $(\beta)$ (in short $\langle \mathrm{BETA}\rangle$) with some extra information on the power type of the modulus.  
	
	\begin{theo}
		\label{thm:beta-renorming} Let $X$ be a Banach space. The following assertions are equivalent.
		\begin{enumerate}[(i)]
			\item\label{it:beta-i} $X$ admits an equivalent norm with property $(\beta)$.
			\item\label{it:beta-ii} $X$ is reflexive and admits an equivalent AUS norm and an equivalent AUC norm.
			\item\label{it:beta-iii} $X$ is a reflexive space with $\Sz(X)\le \omega$ and $\Sz(X^*)\le \omega$.
			\item\label{it:beta-iv} $X$ is reflexive and there exists $p\in (1,\infty)$ such that $X$ admits an equivalent $p$-AUS norm and there exists $q\in (1,\infty)$ such that $X$ admits an equivalent $q$-AUC norm.
			\item\label{it:beta-v} $X$ is reflexive and there exist $p,q\in (1,\infty)$ such that $X$ admits an equivalent norm which is simultaneously $p$-AUS and $q$-AUC.
			\item\label{it:beta-vi} There exists $r\in (1,\infty)$ such that $X$ admits an equivalent norm with property $(\beta)$ with power type $r$.
		\end{enumerate}
	\end{theo}
	
	\begin{proof}
		$(i) \Rightarrow (ii)$ follows from Proposition \ref{prop:beta->reflexive}, Theorem \ref{thm:beta->AUC} and Theorem \ref{thm:beta->AUS}.\\
		$(ii) \Rightarrow (iii)$ follows from the reflexivity of $X$ together with the AUS-AUC$^*$-duality and Proposition \ref{prop:Sz-omega}.\\
		$(iii) \Rightarrow (iv)$ follows from the reflexivity of $X$ together with the AUS-AUC$^*$-duality and Theorem \ref{thm:Sz-omega=AUS}.\\
		$(iv) \Rightarrow (v)$ is Proposition \ref{prop:sim-AUS-AUC}.\\
		$(v) \Rightarrow (vi)$ is Theorem \ref{thm:beta->AUS&AUC}.\\
		$(vi) \Rightarrow (i)$ is trivial.
	\end{proof}

	\section{Notes}
	
	For separable Banach spaces, Theorem \ref{thm:Sz-omega=AUS} is originally due to Knaust, Odell and Schlumprecht \cite{KOS1999} using blocking techniques for spaces with a boundedly complete FDD and reducing the separable case to it. A different approach was taken by Godefroy, Kalton and Lancien in \cite{GKL2001} that led to better quantitative results. Yet another proof was given by Raja in \cite{Raja2010}, which also applied to nonseparable Banach spaces. The results in \cite{KOS1999}, \cite{GKL2001} and \cite{Raja2010}, are essentially stated in terms of the weak$^*$ UKK property on the dual space, but it is not too difficult to show that this property is equivalent to AUC$^*$ (see Exercise \ref{ex:w*UKK-implies-Szlenk-omega}). 
	The proof we have presented here follows Causey's reworking of the renorming problem which builds on the approach from \cite{GKL2001} and makes the similarity to Pisier's renorming technique more transparent. In particular, the classes $\sN_p$, $\sA_p$ and $\sT_p$ were introduced in \cite{Causey3.5}.
	
	More can be said about the various inclusions between these classes. We will denote by $\textsf{Sz}_\omega$ the class of all Banach spaces with Szlenk index equal to $\omega$. And if $1\le q< \infty$, we denote by $\textsf{Sz}_q$ the class of all Banach spaces $X$ such that $X$ has Szlenk power type $q$.
	\begin{theo}
		\label{thm:containments2} 
		Let $p\in (1,\infty]$ and $p'$ be its conjugate exponent. 
		\begin{enumerate}[(i)]
			\item $\textsf{\emph{T}}_p\subset \textsf{\emph{A}}_p\subset \textsf{\emph{N}}_p \subset \bigcap_{1<r<p}\textsf{\emph{T}}_r$. 
			\item $\bigcap_{1<r<p}\textsf{\emph{T}}_r=\bigcap_{1<r<p}\textsf{\emph{A}}_r=\bigcap_{1<r<p}\textsf{\emph{N}}_r=\bigcap_{q>p'}\textsf{\emph{Sz}}_q.$
			\item $\textsf{\emph{Sz}}_\omega=\bigcup_{1\le q< \infty}\textsf{\emph{Sz}}_q.$
		\end{enumerate}
	\end{theo}
	
	\begin{proof} 
		Statement $(i)$ is just Proposition \ref{prop:upper-finite-trees->upper-infinite-trees} and $(iii)$ is Proposition \ref{prop:submult}. The fact that $\bigcap_{1<r<p}\sT_r=\bigcap_{1<r<p}\sA_r=\bigcap_{1<r<p}\sN_r$ is a consequence of $(i)$. We deduce from Theorem \ref{thm:Tp=p-AUS} that $\sT_r \subset \textsf{{Sz}}_{r'}$ for all $r \in (1,\infty]$. Finally, it follows from Propositions \ref{prop:Szlenk-power-summable} and \ref{prop:Sz-q-summable->Ap} that $\textsf{{Sz}}_{r'} \subset \textsf{{A}}_{s}$ for all $1<s<r\le \infty$. This concludes the proof of $(ii)$. 
	\end{proof}
	
	It is important to mention that R. Causey also proved in \cite{Causey3.5} that the following inclusions are strict. 
	
	\begin{theo}\ 
		\begin{enumerate}[(i)]
			\item For $1<p<\infty$, $\textsf{\emph{T}}_p\subsetneq \textsf{\emph{A}}_p\subsetneq \textsf{\emph{N}}_p\subsetneq \bigcap_{1<r<p}\textsf{\emph{T}}_r$. 
			\item $\textsf{\emph{T}}_\infty\subsetneq \textsf{\emph{A}}_\infty= \textsf{\emph{N}}_\infty\subsetneq \bigcap_{1<r<\infty}\textsf{\emph{T}}_r$.
		\end{enumerate}   
	\end{theo}
	
	To the best of our knowledge, it has been quite unnoticed that the class $\textsf{{Sz}}_1$ coincides with the class of Banach spaces with summable Szlenk index or equivalently the class of asymptotic-$c_0$ Banach spaces. It is essentially the contents of Proposition 6.9 in \cite{KOS1999} for the case of Banach spaces with a shrinking FDD. The general case can be deduced from arguments on the determination by quotients with a basis as in \cite{Lancien1996}. A direct proof using weak$^*$ null trees and prunings is given in \cite{CauseyFovelleLancien_tocome}. 
	
	\medskip
	Let us also mention that there is a version of Theorem \ref{thm:Sz-omega=AUS} for higher values of the Szlenk index. Indeed, Lancien, Proch\'azka and Raja \cite{LPR} proved that for $\alpha \in [0,\omega_1)$, a Banach space $X$  satisfies $\Sz(X)\le \omega^{\alpha+1}$ if and only if it admits an equivalent norm whose dual norm in $\omega^\alpha$-AUC$^*$, which means that for every $\eps>0$, there exists $\delta>0$ such that $s_\eps^{\omega^\alpha}(B_{X^*})\subset (1-\delta)B_{X^*}$. 
	
	\medskip
	The origins of Theorem \ref{thm:beta->AUS} can be traced back to the late 1980s and works of Prus \cite{Prus1989} and Kutzarova \cite{Kutzarova1990}. Prus (re)introduced the notion of nearly uniform smoothness (NUS) in \cite{Prus1989} where he characterized Banach spaces with a basis admitting an equivalent NUS norm in terms of the existence of certain upper estimates for basic sequences. Prus mentioned without details that a similar proof will work for separable spaces using biorthogonal systems.
	
	In \cite{Kutzarova1990}, Kutzarova showed that in a separable Banach space with property $(\beta)$, every basic sequence has a blocking satisfying upper $q$-estimates for some $q>1$. Thus, invoking Prus' result, it follows that every Banach space with a basis and property $(\beta)$ admits an equivalent norm that is NUS. To obtain Theorem \ref{thm:beta->AUS} for Banach spaces with a basis, it remains to observe that a Banach space has an equivalent NUS norm if and only if it is reflexive and has an equivalent AUS norm. 
	This can be done directly, or we could invoke Prus' duality result between NUS and Huff's notion of nearly uniform convexity (NUC) \cite{Huff1980}. Indeed, Huff showed that a Banach space is NUC if and only if it is reflexive and has the uniform Kadec-Klee property (UKK).
	To prove Theorem \ref{thm:beta->AUS} in the separable case, one can either observe that Kutzarova's proof can be adjusted to deal with shrinking Markushevich bases and remember Prus' remark that his proof works for separable spaces using biorthogonal systems. Another possibility would be to argue that property $(\beta)$ passes to quotients and that the property $\Sz(X)\le \omega$ is determined by quotients with a basis as shown in \cite{Lancien1996}. 
	The proof we give here is arguably more direct and works without any separability assumption.
	
	To the best of our knowledge, the following local analog of Problem \ref{pb:AUS-AUC-sim-renormings} remains open.
	
	\begin{prob}
		\label{pb:US-UC-sim-renormings}
		Let $1<p \le 2 \le  q <\infty$. 
		If a Banach space has an equivalent $q$-uniformly convex norm and an equivalent $p$-uniformly smooth norm, does it admit an equivalent norm that is both $q$-uniformly convex and $p$-uniformly smooth?
	\end{prob}
	
	Note that the case $p=q=2$ is settled positively since it is a theorem of Kwapien that a Banach space with type $2$ and cotype $2$ must be isomorphic to Hilbert space.
	
	\section{Exercises}
	
	\begin{exer}
		\label{ex:heightoftrees}
		Prove assertions $(i)$ and $(ii)$ in Lemma \ref{lem:heightoftrees}.
	\end{exer}
	
	\begin{exer}
		\label{ex:Szlenk-infinite-dimension} 
		Let $X$ be an infinite-dimensional Banach space. 
		\begin{enumerate}
			\item Show that for any $\eps \in (0,1)$ and $n\in \bN$ so that $n<\frac{1}{\eps}$, $(1-n\eps)B_{X^*} \subset s_\eps^n(B_{X^*})$.
			\item Deduce that $\Sz(X)\ge \omega$.
			\item Show that a Banach space $X$ is finite-dimensional if and only if $\Sz(X)=1$.
		\end{enumerate}  
	\end{exer}
	
	\begin{exer}
		\label{ex:w*UKK-implies-Szlenk-omega}
		Let $X$ be a Banach space  and $K$ be a bounded subset of $X^*$.
		\begin{enumerate}
			\item Show that for all $x^*\in X^*$,  all $\lambda,\vep>0$ and all ordinals $\alpha$:  $s_{\lambda\vep}^\alpha(x^*+\lambda K)= x^*+\lambda s_\vep^\alpha(K)$.
			\item Show that for any weak$^*$ open subset $V$ of $X^*$ and any ordinal $\alpha$: $s_\eps^\alpha(V \cap K)=V\cap s_\eps^\alpha(K)$.    \item The canonical dual norm $\norm{\cdot}_{X^*}$ on $X^*$ is said to be \emph{weak$^*$ uniformly Kadets-Klee} (weak$^*$ UKK in short) if for any $\vep>0$ there exists $\theta>0$ such that $\norm{x^*}_{X^*}\le 1-\theta$ for all $x^*\in B_{X^*}$ for which every weak$^*$ neighborhood $V$ of $x^*$ in $B_{X^*}$ satisfies $\diam(V)\ge \vep$. Show that the canonical dual norm of $X^*$ is weak$^*$ UKK if and only if it is AUC$^*$.
			\item Show that if the canonical dual norm of $X^*$ is weak$^*$ UKK and $X$ is infinite-dimensional, then  $\Sz(X)=\omega$.
		\end{enumerate}
	\end{exer}
	
	\begin{exer}\label{ex:Szlenk index of C(K)}
		Let $K$ be a countable compact metric space.
		
		\begin{enumerate}
			\item Show that $C(K)^*$ is isometric to $\ell_1(K)$ and therefore to $\ell_1$.
			\item Let $\alpha$ be a countable ordinal and assume that $x\in K^{(\alpha)}$ (the Cantor-Bendixon derived set of $K$ of order $\alpha$). Show that the Dirac measure $\delta_x$ belongs to $s_1(B_{\ell_1(K)})$, where $\ell_1(K)$ is equipped with the weak$^*$  topology $\sigma(\ell_1(K),C(K))$. 
			\item Give an example of a weak$^*$ topology on $\ell_1$, for which $\|\ \|_1$ is not AUC$^*$. 
			\item Show that there exists an uncountable family of Banach spaces that are mutually nonisomorphic but so that their duals are all isometric to $\ell_1$. 
		\end{enumerate}
	\end{exer}
	
	\begin{exer}
		\label{ex:value-for-Szlenk}
		Let $X$ be a Banach space and $\alpha$ an ordinal.
		\begin{enumerate}
			\item Show that $\frac12 s_\vep^\alpha(B_{X^*}) + \frac12 B_{X^*} \subseteq s^\alpha_{\vep/2}(B_{X^*})$.
			\item Show that if $s_\vep^\alpha(B_{X^*})\neq \emptyset$, then $0\in s_{\vep/2}^\alpha(B_{X^*})$ and for all $n\ge 1$ $0\in s_{\vep/2^{n+1}}^{\alpha\cdot 2^n}(B_{X^*})$.
			\item If $\Sz(X)\neq \infty$, show that $\Sz(X)=\omega^\alpha$ for some ordinal $\alpha$.
		\end{enumerate}
	\end{exer}
	
	\begin{exer} 
		\label{ex:Lindelof}\, 
		\begin{enumerate}
			\item Recall that a topological space $X$ is \emph{Lindel\"of} if every open cover has a countable subcover and \emph{hereditarily Lindel\"of} if every subspace of $X$ is Lindel\"of. Show that if $(F_\alpha)_{\alpha <\omega_1}$ is a transfinite sequence of closed subsets of a hereditarily Lindel\"of space such that $F_\beta \subseteq F_\alpha$ whenever $\alpha<\beta$, then there is $\gamma <\omega_1$ such that $F_\alpha=F_\gamma$ for all $\alpha \ge \gamma$. 
			\item Let $P$ be a Polish space, i.e. a separable and completely metrizable topological space and denote by $\cF(P)$ the set of closed subsets of $P$. Let $D\colon \cF(P)\to \cF(P)$ be a map satisfying:
			\begin{itemize}
				\item $D(F)\subseteq F$ for all $F\in \cF(P)$,
				\item $D(F)\subsetneq F$ for all $\emptyset\neq F\in \cF(P)$,
				\item $D(F_1)\subseteq D(F_2)$ whenever $F_1\subseteq F_2$ and $F_1,F_2 \in \cF(P)$.
			\end{itemize}
			For every ordinal $\xi$, define 
			\begin{itemize}
				\item $D^0(F) := F$,
				\item $D^{\xi +1}(F):= D(D^\xi(F))$,
				\item $D^\xi(F) := \cap_{\zeta<\xi} D^{\zeta}(F)$ if $\xi$ is a limit ordinal.
			\end{itemize}
			Show that there exists a countable ordinal $\alpha$ such that $D^{\alpha}(F)=\emptyset$ for all $F\in\cF(P)$. 
		\end{enumerate}
	\end{exer}

	\begin{exer}
		\label{ex:fragmented-separable-dual}
		Show directly using a Baire category argument that if $X^*$ is separable, then $(B_{X^*},w^*)$ is weak$^*$ fragmentable.
	\end{exer} 
	
	The goal of the next exercise is to provide a shorter (but less precise) argument for the proof of $(i)$ implies $(ii)$ in Proposition \ref{thm:Szlenk-Asplund}.
	
	\begin{exer}
		Assume that $X$ is a Banach space with separable dual and $K$ is a weak$^*$ compact subset of $X^*$. Consider $(x_k^*)_k$ a $\norm{\cdot}$-dense sequence in $X^*$ and for $n\in \bN$, denote by $B_{n,k}$ the closed ball of center $x_k^*$ and radius $\frac{1}{n}$ in $X^*$ and let $F_{n,k} := K\cap B_{n,k}$. 
		\begin{enumerate}
			\item Show that $F_{n,k}$ is weak$^*$ closed in $K$.
			\item Use Baire's Lemma to conclude. 
		\end{enumerate}
		Note: the proof in the text is more precise in the sense that it shows the existence of a subset of $K$ which is homeomorphic to the Cantor space when it is equipped with the weak$^*$ topology, but uniformly separated in norm.     
	\end{exer}
	
	The goal of the next exercise is to give a more precise version of Proposition \ref{prop:Sz-lp}.
	
	\begin{exer}
		\label{ex:Szlenk-ell_p}
		Let $p\in (1,\infty)$, $q$ its conjugate exponent and $(F_n)_{n=1}^\infty$ be a sequence of finite-dimensional normed spaces. 
		\\ If $X := (\sum_{n=1}^\infty F_n)_{\ell_p}$, show that for all $\eps \in (0,1)$,
		\begin{equation*}
			\Sz(X,\eps) = \lfloor 2^q\eps^{-q} \rfloor +1.
		\end{equation*} 
	\end{exer}
	
	In the next exercise, you are asked to prove Remark \ref{rem:T_pHahnBanach}.
	\begin{exer}
		Let $p\in (1,\infty]$ and $q\in [1,\infty)$ such that $\frac1p +\frac1q=1$. Show that for any sequence $(x_i)_{i=1}^\infty$ in a Banach space $X$ it holds 
		\begin{align*}
			\inf\Big\{c\in (0,\infty] \colon \forall a=(a_i)_{i=1}^\infty \in c_{00}, & \ \Big\|\sum_{i=1}^\infty a_ix_i \Big\| \le c\norm{a}_p\Big\} \\
			& = \sup\{\norm{(x^*(x_i))_{i=1}^\infty}_{\ell_q}\colon\ x^*\in B_{X^*}\}.
		\end{align*}
	\end{exer}
	
	\begin{exer}
		\label{ex:q-AUC*->q-summable-Szlenk}
		Let $q\in[1,\infty)$.
		Show that if $X$ is $q$-AUC$^*$, then $X$ has a $q$-summable Szlenk index.
	\end{exer}
	
	\begin{exer}
		\label{ex:Cz-Kz}
		Let $X$ be a Banach space and $K$ be a weak$^*$ compact subset of $ X^*$. For each $\eps>0$, we introduce the following derivation. 
		\begin{equation*}
			k_\eps(K) := K\setminus\{S\subset X^*\colon S\text{ weak$^*$ open slice and }\alpha(V\cap K)\le \eps\}, 
		\end{equation*}
		where $\alpha(V\cap K)$ is the Kuratowski measure of noncompactness of $V \cap K$ (see Definition \ref{def:Kuratowski} in Section \ref{sec:Kuratowski}). Given an ordinal $\xi$, the derived set of order $\xi$, denoted by $k_\eps^\xi(K)$, is defined inductively by letting:
		\begin{itemize}
			\item $k_\eps^0(K) := K$,  
			\item $k_\eps^{\xi+1}(K) := k_\eps(k^\xi_\eps(K))$,
			\item $k_\eps^{\xi}(K) := \cap_{\zeta<\xi}k^\zeta_\eps(K)$ if $\xi$ is a limit ordinal.
			\item We also use the obvious convention $k_0^\xi(K) := K$, for any ordinal $\xi$. 
		\end{itemize} 
		The ordinal $\Kz(X,\eps)$ is defined as the least ordinal $\xi$ so that $k_\eps^\xi(B_{X^*})=\emptyset$, if such ordinal exists and  we set $\Kz(X,\eps)=\infty$ otherwise. Finally
		\begin{equation*}
			\Kz(X) := \sup_{\eps>0}\Kz(X,\eps).
		\end{equation*}
		Show that $\Kz(X)=\Cz(X)$.
	\end{exer}
	\begin{proof}[Hint]
		The reader can have look at \cite{HajekLancien}, if needed.
	\end{proof}
	
	\begin{exer}
		\label{ex:Tsirelson} 
		Show that for any $c>2$ and any $n\in \bN$, Player A has a winning strategy in the $A(c,\infty,n)$-game in Tsirelson space $\Tsi^*$.
	\end{exer}

	\begin{exer}
		\label{ex:Baire-AUS-AUC}
		Let $(X,\norm{\cdot})$ be a Banach space. Let $S$ be the set of continuous semi-norms on $(X,\norm{\cdot})$ and $P$ be the set of norms on $X$ that are equivalent to $\norm{\cdot}$. Denote by $B$ the closed unit ball of $\norm{\cdot}$. For $M,N \in S$, we set 
		\begin{equation*}
			d(M,N) := \sup\{\abs{M(x)-N(x)}\colon x\in B\}.
		\end{equation*}
		\begin{enumerate}
			\item Show that $(S,d)$ is a complete metric space and that $P$ is open in $(S,d)$ and deduce that $(P,d)$ is a Baire space. 
			\item Show that the set $\mathcal C$ of equivalent AUC norms on $X$ is a $\mathcal G_\delta$ subset of $P$.
			\item Show that if $\mathcal C$ is not empty, then it is dense in $P$.
			\item Assume that $X$ is reflexive and equip the set of equivalent norms on $X^*$ with a similar metric and show that it is homeomorphic to $(P,d)$. 
			\item We still assume that $X$ is reflexive and denote by $\mathcal S$ the set of equivalent AUS norms on $X$. Show that if $\mathcal C$ and  $\mathcal S$ are not empty, then $\mathcal C \cap \mathcal S$ is a dense $\mathcal G_\delta$ subset of $P$.
		\end{enumerate}    
	\end{exer}
	
	\begin{exer} 
		Let $X$ be a Banach space. Show that $X$ admits an equivalent norm with property $(\beta)$ if and only if any separable subspace of $X$ does. 
	\end{exer}
	
	\begin{exer}
		\label{ex:Szlenk of sums}\,
		\begin{enumerate}
			\item Let $p\in (1,\infty)$ and $(X_n)_{n=1}^\infty$ be a sequence of Banach spaces such that for all $n\in \bN$, $\Sz(X_n)= \omega$. Show that if $X := (\sum_{n=1}^\infty X_n)_{\ell_p}$, then $\Sz(X)\le \omega^2$.
			\item Let $(p_n)_{n=1}^\infty \subset (1,\infty)$ be a  strictly decreasing sequence tending to $1$.\\ If $Y := (\sum_{n=1}^\infty \ell_{p_n})_{\ell_2}$, show that $Y$ is reflexive and $\Sz(Y)= \omega^2$. 
			\item Show that $\Sz(Y^*)=\omega$. 
		\end{enumerate} 
	\end{exer}
	
	%%%%%%%%%%%%%%%%%%%%%%%%%%%%%%%%%%%%%%%%%%%%%%%%%%%%%%%%%%%%%%%%%%%%%%%%%%%%%%%%%%%%%%%%%%%%%%%%%%%%%%%%%
	
	\chapter[The basic approximate midpoint principle]{The basic approximate midpoint principle}
	\label{chapter:AMP_I}
	
	The proof of the Mazur-Ulam theorem relies partially on the fact that isometries preserve the (exact) midpoint set of a pair of points. In this chapter, we develop an elegant technique that leads to a coarse-Lipschitz embeddability obstruction between metric spaces whose approximate midpoint sets are of significantly different sizes. The use of approximate metric midpoints in nonlinear Banach space geometry was first introduced by Per Enflo in an unpublished paper where he showed that $L_1$ and $\ell_1$ are not uniformly homeomorphic. 
	%(see Exercise \ref{exer:EnfloFirstUseOfMidpoint}).
	The main result of this section is a comparison of the modulus of asymptotic uniform convexity with the modulus of asymptotic uniform smoothness under a coarse-Lipschitz embedding assumption. This result is obtained by combining quantitative estimates for approximate midpoint sets in AUS or AUC spaces together with a basic approximate midpoint principle.
	Quantitative estimates for approximate midpoint sets can be conveniently formulated in terms of the Kuratowski measure of noncompactness which we introduce in this chapter. Asymptotically uniformly smooth spaces tend to have large approximate midpoint sets, while asymptotically uniformly convex spaces have small approximate midpoint sets. Therefore, the basic approximate midpoint principle is very useful in showing that a Banach space that is sufficiently asymptotically uniformly smooth does not admit a coarse-Lipschitz embedding into a Banach space that is sufficiently asymptotically uniformly convex. 
	
	\section{The approximate midpoint principle}
	
	We start with the definition of approximate midpoints. 
	
	\begin{defi}
		\label{def:ApproximateMidpoints}
		Let $(M,d)$ be a metric space.
		For $\delta \in [0,1)$ and $x,y \in M$, the set of $\delta$-\emph{approximate metric midpoints}\index{approximate metric midpoint}\index{metric midpoint!approximate} of $x$ and $y$ is defined as 
		\begin{align*}
			\Mid(x,y,\delta) & :=\closedball{M}\Big(x,\frac{1+\delta}{2}d(x,y)\Big) \cap \closedball{M}\Big(y,\frac{1+\delta}{2}d(x,y)\Big)\\
			& = \Big\{ z\in M\colon \max\{d(x,z), d(y,z)\}\le \frac{1+\delta}{2}d(x,y) \Big\}.
		\end{align*}
	\end{defi}
	
	The notion of approximate midpoint set is a natural extension of the notion of (exact) midpoint set since it follows from the definitions that $\Mid(x,y,0)=\Mid(x,y)$. 
	
	In the next lemma, we record several useful facts about approximate midpoint sets in normed spaces.
	
	\begin{lemm}
		\label{lem:aprox-mid-basic}
		Let $X$ be a normed space, $x,y \in X$ and $\delta \in (0,1)$. Then,
		\begin{enumerate}[(i)]
			\item $\Mid(x,y,\delta)$ is centrally symmetric around $\frac{x+y}2$.
			\item $\Mid(x,y,\delta) \subset \closedball{X}(\frac{x+y}2,\frac{1+\delta}{2}\norm{x-y})$.
			\item $\closedball{X}(\frac{x+y}2,\frac\delta2\norm{x-y}) \subset \Mid(x,y,\delta)$.
		\end{enumerate}
	\end{lemm}
	
	We defer the elementary justifications to Exercise \ref{ex:first-midpoints}.
	
	\begin{rema}
		Lemma \ref{lem:aprox-mid-basic} says that if $(X\norm{\cdot})$ is a normed space and $x,y\in X$, then $\Mid(x,y,\delta)$ is the ball centered at $\frac{x+y}2$ of an equivalent norm on $X$.
	\end{rema}
	
	In the proof of the Mazur-Ulam Theorem (Theorem \ref{MazurUlam}), we used the observation that if $f$ is an isometric embedding, then $f(\Mid(x,y))\subset \Mid(f(x),f(y))$. This fact might no longer be true if one merely considers bi-Lipschitz embeddings. However, when $f \colon M \to N$ is a Lipschitz map and $x,y \in M$ form a pair of points that saturates the Lipschitz condition, i.e. $d(f(x),f(y))=\Lip(f) d(x,y)$, then it follows from the triangle inequality that the inclusion $f(\Mid(x,y))\subset \Mid(f(x),f(y))$ still holds. Up to passing to approximate midpoint sets, this midpoint principle can be extended to accommodate coarse-Lipschitz maps. The next lemma is a well-known version of the approximate midpoint principle for coarse-Lipschitz maps (cf. \cite[Lemma 10.11]{BenyaminiLindenstrauss2000}). 
	
	%When $\Lip(f)$ is only nearly attained, metric midpoints become approximate metric midpoints and, once the epsilons are unleashed, approximate metric midpoints are mapped to approximate metric midpoints with different error. A precise statement is below.
	
	%\begin{lemm}\label{lem:LipschitzApproxMid}
	%Let $X,Y$ be Banach spaces and $f\colon X \to Y$ a Lipschitz map.
	%If $x,y \in X$ and $\vep>0$ are such that $\norm{f(x)-f(y)}\ge \Lip(f)(1+\vep)^{-1}\norm{x-y}$ then for every $\delta\in(0,1)$ we have
	%\[
	%f(\Mid(x,y,\delta))\subset \Mid(f(x),f(y),(1+\delta)(1+\vep)-1).
	%\]
	%\end{lemm}
	
	%\begin{proof}
	%Indeed, if $z\in \Mid(x,y,\delta)$, we have
	%$$\norm{f(x)-f(z)}\le  \Lip(f)\frac{1+\delta}2\norm{x-y}\le  %\frac{(1+\delta)(1+\vep)}2 \norm{f(x)-f(y)}$$
	%and similarly for $\norm{f(z)-f(y)}$.
	%\end{proof}

	\begin{lemm}
		\label{lem:CLApproxMid} 
		Let $(X,\norm{\cdot})$ be a normed space, $(N,d_N)$ be a metric space and  $f \colon (X,\norm{\cdot}) \to (N,d_N)$ be a coarse-Lipschitz map such that $\Lip_\infty(f)>0$.
		Then, for every $\vep>0$ and $t>0$ there are $x,y \in X$ with $\norm{x-y}\ge t$ such that for every $\delta \in (0,\frac12)$ we have
		\begin{equation*}
			f(\Mid(x,y,\delta)) \subset \Mid(f(x),f(y),\delta(1+\vep)+\vep).    
		\end{equation*}
		In particular, for every $t>0$ and every $\delta \in (0,\frac12)$
		there are $x,y \in X$ with $\norm{x-y}\ge t$ and such that
		\begin{equation}
			f(\Mid(x,y,\delta)) \subset \Mid(f(x),f(y),2\delta).   
		\end{equation}
	\end{lemm}
	
	\begin{proof} Fix $t>0$ and $\eps>0$. Let $\eta \in (0,1)$ be such that $(1+\eta)(1-\eta)^{-1}<1+\vep$. Since $f$ is coarse-Lipschitz, we have that $\Lip_\infty(f) := \inf_{s>0}\Lip_s(f)<\infty$. We also assumed that $\Lip_\infty(f)>0$, so there exists $s>0$ be such that $\Lip_s(f)<(1+\eta)\Lip_\infty(f)$ and recalling that 
		\begin{equation*}
			\Lip_s(f) := \sup\Big\{\frac{d_N(f(x),f(y))}{\norm{x-y}}\colon \norm{x-y}\ge s\Big\},
		\end{equation*}
		we can assume without loss of generality that $s>t$.
		Since $\Lip_{4s}(f)\le \Lip_s(f)<\infty$, we can pick $x_1,x_2 \in X$ such that 
		$$\norm{x_1-x_2}\ge 4s\ \ \text{and}\ \  d_N(f(x_1),f(x_2))\ge  (1-\eta)\Lip_{4s}(f)\norm{x_1-x_2}.$$
		Let $\delta \in (0,\frac12)$ and let $z \in \Mid(x_1,x_2,\delta)$.
		Then, for $\{i,j\}=\{1,2\}$ 
		\begin{equation*}
			\norm{x_i-z}\ge \norm{x_i-x_j}-\frac{1+\delta}2\norm{x_i-x_j}=\frac{1-\delta}2\norm{x_i-x_j}\ge \frac14\norm{x_i-x_j}\ge s.
		\end{equation*}
		Thus, for $i\in\{1,2\}$ we have
		\begin{align*}
			d_N(f(x_i),f(z))&\le \Lip_s(f)\|x_i-z\| \le  (1+\eta)\Lip_\infty(f) \norm{x_i-z}\\
			&\le  (1+\eta)\Lip_{4s}(f)\frac{1+\delta}2\norm{x_1-x_2}\\
			&\le  \frac{1+\eta}{1-\eta}\frac{1+\delta}2 d_N(f(x_1),f(x_2))
			\le \frac{(1+\eps)(1+\delta)}{2}d_N(f(x_1),f(x_2)).
		\end{align*}
		%The same computation is true for $\norm{f(y)-f(z)}$ and the result follows.
		The first claim follows by letting $(x,y) := (x_1, x_2)$.
		By picking $\vep>0$ so small that $\delta+\vep+\delta \vep<2\delta$, we get the second claim. 
		%(Notice that with this order of quantifiers, we can do it for $\delta \in (0,1)$.)
	\end{proof}
	
	%\begin{rema}
	%In the proof of the approximate midpoint principle (Lemma \ref{lem:CLApproxMid}), we only used that the domain space is a normed space to guarantee that $\Mid(x_1,x_2,\delta)$ is not empty and the proof would go through verbatim if, for instance, the domain space was a geodesic metric space.  
	%\end{rema}
	
	In order to use the approximate midpoint principle for nonembeddability purposes, we need to estimate the size of approximate midpoint sets. In \cite{BenyaminiLindenstrauss2000}, upper and lower estimates for the approximate midpoint sets that were implicit in the works of Enflo, Bourgain and Gorelik on $\ell_p$-spaces, were formalized in spaces with lower/upper-$p$-estimates. Ultimately, in Randrianarivony's Ph.D. thesis \cite{RandrianarivonyThesis}, these estimates were generalized for arbitrary Banach spaces in terms of the moduli of asymptotic uniform smoothness and convexity. This is the level of generality that we choose to present here.
	
	On the one hand, the size of approximate midpoint sets is controlled from below by the modulus of asymptotic uniform smoothness.
	\begin{lemm}
		\label{lem:midpointsAUS}
		Let $X$ be an infinite-dimensional Banach space.
		For every $x,y \in X$, every $\delta \in (0,1)$ and every $\eta>0$ such that $\bar{\rho}_X(\eta)<\delta$, there exists $Z \in \cof(X)$ such that
		\begin{equation}
			\frac{x+y}{2}+\eta \frac{\norm{x-y}}{2}B_Z \subset \Mid(x,y,\delta).
		\end{equation}
	\end{lemm}
	
	\begin{proof}
		Up to translation and rescaling, it is enough to study the sets of the form $\Mid(-x,x,\delta)$ for $x \in \sphere{X}$. So, let $x\in \sphere{X}$, $\delta \in (0,1)$ and $\eta>0$ so that $\bar{\rho}_X(\eta)<\delta$.
		By definition of $\bar{\rho}_X(\eta)$, there exists $Z \in \cof(X)$ such that
		$x+\eta B_Z \subset (1+\delta)B_X$ and thus, by symmetry, $-x+\eta B_Z \subset (1+\delta)B_X$. So, we have that $\eta B_Z \subset \Mid(-x,x,\delta)$.
	\end{proof}
	
	On the other hand, the size of approximate midpoint sets is controlled from above by the modulus of asymptotic uniform convexity. The argument in this situation is less straightforward and it is based on a comparison of the unit ball of a space with the unit balls of its finite-codimensional subspaces.
	
	\begin{lemm}
		\label{lem:compact-pert}
		Let $X$ be an infinite-dimensional Banach space and $Z$ a finite-codimensional subspace of $X$. Then, there is a compact subset $K$ of $X$ such that 
		\begin{equation}
			B_X\subset 3B_Z+K.
		\end{equation}
	\end{lemm}
	
	\begin{proof}
		Let $Q\colon X\to X/Z$ be the canonical quotient map and $\eps \in (0,1)$ such that $(2+2\eps)(1-\eps)^{-1}\le 3$. Pick $(u_{i})_{i=1}^k$ an $\frac{\vep}{2}$-net in $B_{X/Z}$. Now, for all $i\in\{1,\dots,k\}$, pick $x_i\in (1+\eps)B_X$ such that $u_i=Q(x_i)$. Let $F:=\linspan\{x_1, \dots,x_k\}$ and $x\in B_X$. Since $Q$ is of norm $1$, we can find $i_0\in \{1,\dots,k\}$ such that $\norm{Q(x)-u_{i_0}}_{X/Z}<\vep$. Thus, there is $z\in Z$ such that $\norm{x-x_{i_0}-z}<\vep$. Note that $\norm{z}\le 2+2\eps$. We have shown that $B_X \subset (2+2\vep)B_Z + (1+\eps)B_F + \vep B_X.$ After iterating $n$-times we get 
		\begin{equation*}
			B_X\subset (1+\vep+\dots+\vep^n)\big((2+2\vep)B_Z + (1+\eps)B_F\big) + \vep^{n+1} B_X.
		\end{equation*}
		Since $B_F$ is compact and $B_Z$ is closed, $(2+2\vep)B_Z + (1+\eps)B_F$ is closed and we easily deduce that 
		$$B_X \subset \frac{1}{1-\vep}\big((2+2\vep)B_Z + (1+\eps)B_F\big).$$
		The conclusion follows by taking $K :=\frac{1+\eps}{1-\vep}B_F$.
	\end{proof}
	
	%\begin{rema}
	%The compact perturbation in Lemma \ref{lem:compact-pert} is necessary even for hyperplanes since there are Banach spaces that are not isomorphic to their hyperplanes \cite{Gowers}.
	%\end{rema}
	
	\begin{lemm}
		\label{lem:midpointsAUC}
		Let $X$ be an infinite-dimensional Banach space. Then,
		for every $x,y \in X$, every $\delta \in (0,1)$ and every $t>0$ such that $\bar{\delta}_X(t)>\delta$, there exist a compact subset $K$ of $X$ and $Z \in \cof(X)$ such that
		\begin{equation}
			\Mid(x,y,\delta) \subset \frac32t\norm{x-y}B_Z + K.
		\end{equation}
	\end{lemm}
	
	\begin{proof}
		Up to translation and rescaling, it is again enough to study the sets of the form $\Mid(-x,x,\delta)$ for $x \in \sphere{X}$.
		%\mnote{\tch{Does $t$ always exist?} Yes but not in $(0,1)$}
		So let $x\in \sphere{X}$, $\delta \in (0,1)$ and $t>0$ such that $\bar{\delta}_X(t)>\delta$. By continuity of $\bar{\delta}_X$ we may pick $t'<t$ such that $\bar{\delta}_X(t')>\delta$. By definition of $\bar{\delta}_X(t')$, there exists $Z \in \cof(X)$ such that for all $z \in S_Z$ we have $\norm{x+t'z}>1+\delta$ and thus $t'z \notin \Mid(-x,x,\delta)$. By convexity of the approximate midpoint set, it follows that $Z\cap \Mid(-x,x,\delta) \subset t'B_Z$.
		On the other hand, since $\Mid(-x,x,\delta)$ is the ball of an equivalent norm on $X$, it follows from Lemma \ref{lem:compact-pert}
		%we have by the Gorelik principle (Theorem~\ref{Gorelik}, case (2)), applied to $f=Id_X$ and $M=\frac{t}{t'}>1$,
		that there exists a compact $K \subset X$ such that
		\begin{equation}
			\Mid(-x,x,\delta) \subset 3(Z\cap \Mid(-x,x,\delta))+K \subset 3t'B_Z + K,    
		\end{equation}
		which gives the desired result.
	\end{proof}
	
	If $X$ is a subspace of $Y$, it is easy to see that $\bar{\delta}_Y\le \bar{\rho}_X$. The following theorem greatly generalizes this fact up to multiplicative constants.
	
	\begin{theo}
		\label{thm:CL-delta-rho}
		Let $X$ and $Y$ be infinite-dimensional Banach spaces and $f\colon X \to Y$ be a coarse-Lipschitz embedding. Then, there is $C>0$ such that for all $t\in(0,1)$,
		\begin{equation}
			\bar{\delta}_Y(t)\le 2\bar{\rho}_X(Ct).
		\end{equation}
	\end{theo}
	
	\begin{proof} Assume that for some $s>0$, there are $A,B>0$ such that 
		$$\forall x_1,x_2 \in X,\ \norm{x_1-x_2}\ge s \Rightarrow \frac{1}{A}\norm{x_1-x_2}\le \norm{f(x_1)-f(x_2)}\le B\norm{x_1-x_2}.$$
		Fix $t \in (0,1)$ and assume that $0<2\delta<\bar{\delta}_Y(t)$. Fix $r>s$ (to be adjusted later). By the approximate midpoint principle, there are $x_1,x_2\in X$ such that $\norm{x_1-x_2}\ge r$ and 
		$ f(\Mid(x_1,x_2,\delta)) \subset \Mid(f(x_1),f(x_2),2\delta)$ and using translations if needed, we can certainly assume that there is $x\in X$ such that $\norm{x}\ge \frac{r}{2}$, $f(-x)=-f(x)$ and $f(\Mid(-x,x,\delta)) \subset \Mid(-f(x),f(x),2\delta)$. By Lemma \ref{lem:midpointsAUC}, there exist a compact subset $K$ of $Y$ and $Z_1 \in \cof(Y)$ such that
		\begin{equation*}
			\Mid(-f(x), f(x),2\delta) \subset 3t\norm{f(x)}B_{Z_1} + K.    
		\end{equation*}
		Let now $C>0$ be such that $\bar{\rho}_X(Ct)< \delta$. Then, by Lemma \ref{lem:midpointsAUS}, one can find $Z_2\in \cof(X)$ such that  
		\begin{equation*}
			Ct\norm{x} B_{Z_2} \subset \Mid(-x,x,\delta),
		\end{equation*}
		and hence 
		\begin{equation*}
			f(Ct\norm{x}B_{Z_2})\subset f(\Mid(-x,x,\delta))\subset \Mid(-f(x),f(x),2\delta) \subset 3t\|f(x)\|B_{Z_1} + K.
		\end{equation*}
		Since $Z_2$ is infinite-dimensional, there is a $Ct\norm{x}$-separated sequence $(z_n)_n$ in $Ct\norm{x}B_{Z_2}$. If, as we could, $r$ was chosen so that $Ctr \ge 2s$, we get that $(f(z_n))_n$ is $\frac{C}{A}t\norm{x}$-separated. Then, we can write $f(z_n)=y_n+k_n$ with $y_n \in 3t\norm{f(x)}B_{Z_1}$ and $k_n\in K$. Since $K$ is compact, $\sup_{n\neq m}\norm{y_n-y_m}\ge \frac{C}{A}t\norm{x}$. Since $B_{Z_1}$ has diameter $2$, we deduce that $\frac{C}{A}t\norm{x}\le 6t\norm{f(x)}\le 6Bt\norm{x}$ (the last inequality follows from the fact that $\norm{2x}\ge r>s$ and $f(-x)=-f(x)$). Summarizing, we have shown that if 
		$2\bar{\rho}_X(Ct)<\bar{\delta}_Y(t)$, then $Ct\le 6ABt$. It follows from the continuity of $\bar{\rho}_X$ that $\bar{\delta}_Y(t)\le 2\bar{\rho}_X(6ABt)$. 
	\end{proof}

	\begin{rema}
		\label{rema:AMUCMidpoints}\, 
		\begin{enumerate}
			\item The constant $C$ in Theorem \ref{thm:CL-delta-rho} can be taken to be $C=6\Lip_\infty(f)\cdot \textrm{co}\Lip_\infty(f)$.
			\item It also follows from Theorem \ref{thm:CL-delta-rho} that if $X$ isomorphically embeds into $Y$, then $\bar{\delta}_Y(t)\le 2\bar{\rho}_X(6\textrm{d}_{\textrm{BM}}(X,Y)t)$ for all $t\in (0,1)$.  
			\item The proofs above would have also worked for the AMUC modulus:
			\begin{equation*}
				\hat{\delta}_X(t) := \inf_{x \in S_X} \sup_{Y\in \cof(X)} \inf_{y\in S_Y} \frac12(\norm{x+ty}+\norm{x-ty})-1.    
			\end{equation*}
			We shall use this modulus in Chapter \ref{chapter:AMP_II} where we will present a more sophisticated application of the approximate midpoint principle.
		\end{enumerate}
	\end{rema}
	
	The following immediate corollary appeared first in this form in the Ph.D. thesis of F. N\'etillard \cite{NetillardThesis}, but the main ideas can be found in \cite{JLS1996}. 
	
	\begin{coro}
		\label{cor:CLembAUC}
		Let $1\le p<q \le \infty$.
		Let $X$ be $q$-AUS and $Y$ be $p$-AUC Banach spaces.
		Then, there is no coarse-Lipschitz embedding of $X$ into $Y$. In particular, there is no coarse-Lipschitz embedding of an $\ell_q$-sum of finite-dimensional spaces into an $\ell_p$-sum of finite-dimensional spaces.
	\end{coro}
	
	\begin{rema}
		\label{rem:Bourgain-Enflo}
		It follows from Corollary \ref{cor:CLembAUC} that $\ell_2$ does not coarse-Lipschitz embed into $\ell_p$ if $p<2$ and hence $L_p$ cannot be uniformly homeomorphic to $\ell_p$ when $p\in[1,2)$. This recovers the results from Enflo (unpublished for $p=1$) and Bourgain \cite{Bourgain1987} (for $p\in(1,2)$).
	\end{rema} 
	
	\section[Kuratowski measure of noncompactness]{Approximate midpoint sets and Kuratowski measure of noncompactness}
	\label{sec:Kuratowski}
	
	\begin{defi}
		\label{def:Kuratowski}
		For a metric space $(M,d)$, its \emph{Kuratowski measure of noncompactness}\index{Kuratowski measure of noncompactness} $\alpha(M)$ is defined as the infimum of all $\vep>0$ such that $M$ can be covered by finitely many sets of diameter at most  $\vep$.
	\end{defi}
	Recall that a subset $A$ of $M$ is \emph{$\eps$-separated} if $d(x,y)\ge \eps$ for all $x,y\in A$ so that $x\neq y$.
	It will also be convenient to use another measure of noncompactness, denoted by $\beta(M)$ and defined as the infimum of all $\eps>0$ such that $M$ has no infinite $\eps$-separated  subset. 
	
	We start with an elementary lemma where we record the basic properties of the Kuratowski measure of noncompactness.
	\begin{lemm}
		\label{lem:Kuratowski}\
		\begin{enumerate}[(i)]
			\item Let $(M,d)$ be a metric space. Then, $\beta(M) \le \alpha(M)\le 2 \beta(M)$.
			
			\item If $(K,d)$ is compact, then $\alpha(K)=0$.
			
			\item If $A,B$ are subsets of a normed space $X$, then $\alpha(A+B)\le \alpha(A)+\alpha(B)$.
			
			\item If $X$ is an infinite-dimensional normed space, then $1\le \alpha(B_X)\le 2$.
			
			% \item Let $X$ be a Banach space and $K$ be a compact subset of $X$. Then, $$\alpha(K+rB_X)\le  2r.$$
			
			\item Let $M,N$ be metric spaces, $S\subset M$ and $f\colon S\to N$ be a map such that $\rho(d_M(x,y))\le  d_N(f(x),f(y))$ for all $x,y\in S$ and some increasing continuous function $\rho\colon [0,\infty)\to [0,\infty)$. Then,
			\begin{equation*}
				\rho\Big(\frac12\alpha(S)\Big)\le  \alpha(f(S)).
			\end{equation*}
			%If $f\colon M \to Y$ satisfies $A\norm{x-y}\le  \norm{f(x)-f(y)}$ for all $x,y\in M$ such that $\norm{x-y}\ge C$, then $\alpha(M)>C$ implies $\alpha(f(M))\ge A\alpha(M)$.
			\item If $(X,\norm{\cdot})$ is an infinite-dimensional normed space, then  
			\begin{equation*}
				\frac{\delta}{2}\norm{x-y}\le \alpha(\Mid(x,y,\delta))\le (1+\delta)\norm{x-y}.
			\end{equation*}
		\end{enumerate}
	\end{lemm}
	
	\begin{proof} 
		$(i)$ Let $\eps>0$. Assume that $A$ is an infinite $\eps$-separated subset of $M$ and consider a cover of $M$ with sets of diameters at most $\eps'$, for some $\eps'\in (0,\eps)$. Then, as any element of this cover can contain at most one element of $A$, it is clear that this cover is infinite. This proves that $\beta(M)\le \alpha(M)$.
		
		For a given $\eps>0$, consider a maximal $\eps$-separated subset $A$ of $M$ and assume that $A$ is finite. Then, $M$ can be covered by the open balls of radius $\eps$ (and diameter at most $2\eps$) centered at the elements of $A$. This is enough to conclude that $\alpha(M)\le 2\beta(M)$.
		
		%$(b)$ It is a consequence of the following obvious facts: $\alpha(K)=0$, $\alpha(B_X)\le 2$ and for all $A,B$ subsets of $X$, $\alpha(A+B)\le \alpha (A)+ \alpha(B)$.
		
		$(ii)-(iii)-(iv)$ are obvious.
		
		$(v)$ Without loss of generality, assume that $\rho(\frac12\alpha(S))>0$ and pick $\eps$ such that $0<\eps<\rho(\frac12\alpha(S))$. Since $\rho$ is continuous, there is $t<\alpha(S)$ such that $\eps<\rho(\frac12 t)$. By point $(i)$, there is a $\frac12 t$-separated sequence $(x_n)_{n\ge 1}$ in $S$. Then, by our assumption, the sequence $(f(x_n))$ is $\rho(\frac12 t)$-separated in $f(S)$. Using again point $(i)$, we get that $\alpha(f(S)) \ge \beta(f(S))\ge \rho(\frac12 t)>\eps$. Since $\eps<\rho(\frac12\alpha(S))$ was arbitrary, we deduce that $\alpha(f(S)) \ge \rho(\frac12\alpha(S))$.
		
		$(vi)$ follows immediately from $(iv)$ and Lemma \ref{lem:aprox-mid-basic}.
	\end{proof}
	
	We now state a general nonembedding result that can be proved using only the combination of the approximate midpoint principle and the Kuratowski measure of noncompactness.
	\begin{prop}
		\label{prop:NoCLembeddingGeneral}
		Let $(X,\norm{\cdot})$ be a normed space and $(N,d_N)$ be a metric space such that for some $C_X,C_N>0$,  $p,q\in [1,\infty]$ and for all $\delta \in (0,1)$, we have
		\begin{enumerate}
			\item $\alpha(\Mid(y,y',\delta))\le C_N\delta^{\frac1p}d_N(y,y')$ for all $y,y' \in N$,
			\item[] and
			\item $\alpha(\Mid(x,x',\delta))\ge C_X\delta^{\frac1q} \norm{x-x'}$ for all $x,x' \in X$.
		\end{enumerate}
		If $1\le p<q\le \infty$, then there is no coarse-Lipschitz embedding of $X$ into $N$.
	\end{prop}
	
	\begin{proof}
		Assume that $f \colon X \to N$ is a coarse-Lipschitz embedding and $A,B>0$ and $C\ge 1$ are such that for all $x,y \in X$:
		\begin{equation*}
			A \norm{x-y}-C\le d_N(f(x),f(y))\le B\norm{x-y}+C.    
		\end{equation*}
		Let us now fix $\delta \in (0,\frac12)$ and $t>0$, to be specified later. By Lemma \ref{lem:CLApproxMid}, there exist $x,y \in X$ with $\norm{x-y}>t$ such that $f(\Mid(x,y,\delta)) \subset \Mid(f(x),f(y),2\delta)$. Then, by Lemma \ref{lem:Kuratowski}~$(v)$,
		\begin{equation*}
			\alpha(f(\Mid(x,y,\delta))) \ge \frac{A}{2}\alpha(\Mid(x,y,\delta))-C\ge \frac{AC_X}{2}\delta^{1/q}\norm{x-y}-C,
		\end{equation*}
		while
		\begin{equation*}
			\alpha(\Mid(f(x),f(y),2\delta))\le C_N(2\delta)^{1/p}d_N(f(x),f(y))\le C_N(2\delta)^{1/p}(B\norm{x-y}+C).
		\end{equation*}
		Therefore, since $\delta<\frac12$, 
		\begin{equation*}
			\frac{AC_X}{2}\delta^{1/q}\norm{x-y}\le C_N(2\delta)^{1/p}B\norm{x-y}+(C_N+1)C
		\end{equation*}
		and 
		\begin{equation*}
			1\le \frac{2C_NB}{AC_X}2^{1/p}\delta^{1/p-1/q}+\frac{2C(C_N+1)}{AC_X\delta^{1/q}t}.
		\end{equation*}
		Since $1/p-1/q>0$, this is impossible if $\delta$ was first chosen small enough and then $t$ large enough. This concludes the proof. 
	\end{proof} 
	
	\begin{rema} 
		Clearly, in the above statement the maps $\delta \mapsto \delta^{1/p}$ and $\delta \mapsto \delta^{1/q}$ can be replaced by more general functions $\phi,\psi$ satisfying $\lim_{\delta\to 0^+}\frac{\phi(\delta)}{\psi(\delta)}=0$. 
	\end{rema}
	
	Now, we can reformulate the results from the previous section in terms of the Kuratowski measure of noncompactness.
	
	\begin{lemm}
		\label{lem:midpointsAUS-K}
		Let $X$ be an infinite-dimensional Banach space. 
		For every $x,y \in X$ and $\delta \in (0,1)$, we have
		\begin{equation*}
			\alpha(\Mid(x,y,\delta))\ge \sup\set{\eta\colon \bar{\rho}_X(\eta)<\delta}\frac{\norm{x-y}}2.     
		\end{equation*}
		In particular, if $X$ is $p$-AUS, for some $p \in (1,\infty)$, then there exists $C_X>0$ such that for all $x,y \in X$ and $\delta\in (0,1)$,
		\begin{equation*}
			\alpha(\Mid(x,y,\delta))\ge C_X\delta^{1/p}\norm{x-y}.
		\end{equation*}
	\end{lemm}
	
	\begin{proof} 
		Since $X$ is infinite-dimensional, $\alpha(B_X)\ge 1$ and the lower estimate of $\alpha(\Mid(x,y,\delta))$ is then an immediate consequence of Lemma \ref{lem:midpointsAUS}.
		In the case when $\bar{\rho}_X(t)\le Ct^p$, setting $C_X := \frac12\left(\frac1C\right)^{1/p}$ yields the required estimate.
	\end{proof}

	\begin{lemm}
		\label{lem:midpointsAUC-K}
		Let $X$ be an infinite-dimensional Banach space. For every $x,y \in X$ and $\delta \in (0,1)$, we have
		\begin{equation*}
			\alpha(\Mid(x,y,\delta))\le 3 \inf\set{t \colon \bar{\delta}_X(t)>\delta}\norm{x-y}.     
		\end{equation*}
		In particular, if $X$ is $p$-AUC, then there exists $C_X>0$ such that for all $x,y \in X$ and $\delta\in (0,1)$,
		\begin{equation*}
			\alpha(\Mid(x,y,\delta))\le C_X\delta^{1 /p}\norm{x-y}.
		\end{equation*}
	\end{lemm}
	
	\begin{proof}
		The estimate of $\alpha(\Mid(x,y,\delta))$ follows from Lemma \ref{lem:midpointsAUC} and Lemma \ref{lem:Kuratowski}~$(ii)$-$(iii)$-$(iv)$. Assume now that $X$ is $p$-AUC and therefore that  $\bar{\delta}_X(t)\ge Ct^p$, for some $C\in (0,1)$ and all $t\in (0,1)$. Then, setting $C_X := 3\left(\frac1C\right)^{1/p}$ yields the required estimate.
	\end{proof}

	\section{Notes}
	Another application of the approximate midpoint principle is due to Weston in \cite{Weston93}, where he shows that $L_p$ and $\ell_q$ are not uniformly equivalent when $0<p,q\le 1$. In this range, we are not exclusively dealing with Banach spaces anymore but merely with metric spaces that are not necessarily metrically convex (but are quasi-Banach spaces). 
	
	We refer the reader to Chapter \ref{chapter:AMP_II} and Chapter \ref{chapter:diamonds} for more sophisticated developments of the applications of the approximate midpoint principle and the asymptotic (midpoint) uniform convexity to the nonlinear geometry of Banach spaces. 
	
	The Kuratowski measure of noncompactness was introduced by Kuratowski in \cite{Kuratowski1930}. We refer the reader to \cite{BanasGoebel1980} or \cite{ATDBLA1997} for more information on various measures of noncompactness and their applications.
	
	\section{Exercises}
	
	\begin{exer}
		\label{ex:Mazur-Ulam-strictly-convex}
		Let $X$ and $Y$ be two normed spaces. Assume that $\norm{\cdot}_Y$ is strictly convex and that $f$ is an isometry from $X$ into $Y$ such that $f(0)=0$. Prove that $f$ is linear.
	\end{exer}
	
	\begin{exer}
		\label{ex:first-midpoints}
		Let $X$ be a Banach space, $x,y \in X$ and $\delta \in (0,1)$.
		Prove that
		\begin{enumerate}
			\item $\Mid(x,y,\delta)$ is centrally symmetric around $\frac{x+y}2$.
			\item $\Mid(x,y,\delta) \subset \closedball{X}(\frac{x+y}2,\frac{1+\delta}{2}\norm{x-y})$.
			\item $\closedball{X}(\frac{x+y}2,\frac\delta2\norm{x-y}) \subset \Mid(x,y,\delta)$.
			\item $\frac{\delta}{2}\norm{x-y}\le \alpha(\Mid(x,y,\delta))\le (1+\delta)\norm{x-y}$
		\end{enumerate}
	\end{exer}
	
	%\begin{proof}
	% We have $\frac{x+y}2+z-y=-(\frac{x+y}2-z-x)$ and $\frac{x+y}2+z-x=-(\frac{x+y}2-z-y)$ which implies the symmetry.
	
	% Since $\Mid(x,y,\delta)$ is symmetric, we have for every $\frac{x+y}2+z \in \Mid(x,y,\delta)$ that also $\frac{x+y}2-z \in \Mid(x,y,\delta)$. It follows that $\norm{z}\le  \frac12(\norm{z+\frac{x+y}2}+\norm{z-\frac{x+y}2})\le  \frac{1+\delta}2\norm{x-y}$.
	
	% Let $\norm{z}\le  \frac\delta2$. By the triangle inequality we get $\norm{x-\frac{x+y}2-z}\le  \frac{1+\delta}2\norm{x-y}$ and similarly for $y$.
	%\end{proof}
	
	The goal of the next exercise is to give a proof of Enflo's unpublished result that does not rely on the fact that $L_1$ contains an isomorphic copy of $\ell_2$ (it will use Rademacher functions though, but not Khintchine's inequalities).
	\begin{exer} 
		Let $Y$ and $Z$ be Banach spaces and $X :=Y \oplus_1 Z$ with corresponding projections $P\colon X\to Y$ and $Q\colon X\to Z$.
		\begin{enumerate}
			\item Let $x,y \in X$ be such that $\norm{P(x-y)}\ge (1-\vep)\norm{x-y}$ for some $\eps \in (0,1)$.
			Prove that
			\begin{equation*}
				\Mid(x,y,\delta) \subset \frac{x+y}2+\frac{1+\delta}2\norm{x-y}B_Y+\frac{\delta+\vep}2\norm{x-y}B_Z.
			\end{equation*}
			\item Show that for every $x,y\in \ell_1$,  $\alpha(\Mid(x,y,\delta))\le \delta\norm{x-y}_1$.
			\item Let $f\neq g \in L_1 := L_1[0,1]$. Show that there is a sequence $(h_n)_n \subset\Mid(f,g,0)$ such that for all $m\neq n$, 
			\begin{equation*}
				\norm{h_n-h_m}_1=\frac12 \norm{f-g}_1.
			\end{equation*}
			
			\item Deduce that $\alpha(\Mid(f,g,0))\ge \frac12 \norm{f-g}_1$.
			\item Show that there is no coarse-Lipschitz embedding from $L_1$ into  $\ell_1$.
		\end{enumerate}
	\end{exer}
	
	\begin{proof}[Hint:]
		For 2., consider the linear span $Y$  of the first $N$ elements of the canonical basis of $\ell_1$.
		For 3., assume that $g=0$ and $\norm{f}_1=1$ and consider $h_n := f\frac{1+r_n}{2}$, where $(r_n)_n$ is the Rademacher sequence.
	\end{proof}
	%\begin{proof}
	%Let $\frac{x+y}2+z \in \Mid(x,y,\delta)$. Then,
	%\[
	%\begin{aligned}
	%\frac{1+\delta}2\norm{x-y}&\ge \Big\|{\frac{x-y}2\pm z}\Big\|\\
	%&=\Big\|{P\Big(\frac{x-y}2\Big)\pm P(z)}\Big\|+\Big\|{Q\Big(\frac{x-y}2\Big)\pm Qz}\Big\|.
	%\end{aligned}
	%\]
	%Thus by summing, applying the triangle inequality in both terms and dividing by 2, we get
	%\[
	%\frac{1+\delta}2 \norm{x-y}-\Big\|{P\Big(\frac{x-y}2\Big)}\Big\|\ge \norm{Qz}.
	%\]
	%Therefore,  $(\delta+\vep)\norm{\frac{x-y}2}\ge \norm{Qz}$.\\
	
	%On the other hand, by Exercise \ref{ex:first-midpoints}~(b), the assumption $\frac{x+y}2+z \in \Mid(x,y,\delta)$ also implies that $\norm{Pz}\le  \norm{z}\le  \frac{1+\delta}2\norm{x-y}$.
	
	%By translation, we may assume that $g=0$. Let us define a measure $\mu$ on $[0,1]$ by $\mu(A)=\frac{\int_A \abs{f}d\lambda}{\norm{f}_1}$. Then, $T:L_1(\mu) \to L_1$ given by $T\varphi = \varphi f$ is a bijective isometry such that $T(\norm{f}_1\indicator{[0,1]})=f$. Let $(r_n)$ be a Rademacher sequence on $([0,1],\mu)$. Then, setting $h_n=T(\norm{f}_1\frac{r_n+1}{2})$ will yield the conclusion.
	%\end{proof}
	
	In the next exercise, we propose a different proof that $\co$ does not admit a coarse-Lipschitz embedding into $\ell_p$ when $p\in(1,\infty)$.
	
	\begin{exer}
		\label{ex:c0-midpoints}
		Let $1<p<\infty$, $Y$ and $Z$ be Banach spaces and $X :=Y \oplus_p Z$ with corresponding projections $P\colon X\to Y$ and $Q\colon X\to Z$.
		\begin{enumerate}
			\item Let $x,y\in X$ be such that $\norm{P(x-y)}\ge (1-\vep)\norm{x-y}$ for some $\eps \in (0,1)$. 
			Prove that
			\begin{equation*}
				\Mid(x,y,\delta) \subset \frac{x+y}2+\frac{1+\delta}2\norm{x-y}B_Y+\frac{((1+\delta)^p-(1-\vep)^p)^{1/p}}2\norm{x-y}B_Z.
			\end{equation*}
			\item Show that  $\alpha(\Mid(x,y,\delta))\le 4\delta^{1/p}\norm{x-y}_p$ for all $x,y \in \ell_p$.
			\item Let $x,y \in c_0$. Show that for every $\delta\in (0,1)$ and every $\eta\in (0,\delta)$, there is $N \in \bN$ such that
			\begin{equation*}
				\frac{x+y}2+\frac{1+\eta}{2}\norm{x-y}_\infty B_{E_N} \subset \Mid(x,y,\delta),    
			\end{equation*}
			where $E_N:=\set{z \in c_0\colon z(n)=0 \textrm{ for } n<N}$.
			\item Show that $\beta(B_{c_0})=\alpha(B_{c_0})=2$ and deduce that $\alpha(\Mid(x,y,\delta))\ge \norm{x-y}_\infty$.
			\item Show that there is no coarse-Lipschitz embedding from $\co$ into  $\ell_p$ for $1< p<\infty$.
		\end{enumerate}
	\end{exer}
	
	%\begin{proof}
	%Indeed, it is enough to take $N$ such that for every $n\ge N$ we have $\max\set{\abs{x(n)},\abs{y(n)}}<\frac{\delta}2\norm{x-y}$.
	
	%Use that $(1+\delta)^p-1\le  2^p\delta$ for $0\le  \delta \le  1$,
	
	%Assume that $f\colon X\to \ell_1$ satisfies $A\norm{x-y}\le  \norm{f(x)-f(y)}\le  B\norm{x-y}$.
	%Let $\vep>0$ be such that $\vep B<A$. Let $x,y \in X$ be such that $\norm{f(x)-f(y)}\ge \frac{\norm{f}_L}{1+\vep} \norm{x-y}$ and let $\delta>0$ be arbitrary. Exercises~\ref{e:L1midpoints} and \ref{e:c0midpoints} yield that $A \norm{x-y}\le  \alpha(f(\Mid(x,y))$.Lemma~\ref{l:l1sum} implies on the other hand that $\alpha(\Mid(f(x),f(y),(1+\delta)(1+\vep)-1))\le  \vep B\norm{x-y}$ if $\delta$ is small enough. Now using Lemma~\ref{l:LipschitzApproxMid} together with the monotonicity (with respect to inclusion) of the Kuratowski measure we obtain $A \norm{x-y}\le  \vep B\norm{x-y}$ which is impossible by the choice of $\vep$.
	%\end{proof}

	%%%%%%%%%%%%%%%%%%%%%%%%%%%%%%%%%%%%%%%%%%%%%%%%%%%%%%%%%%%%%%%%%%%%%%%%%%%%%%%%%%%%%%%%%%
	
	\chapter[The Gorelik principle and applications]{The Gorelik principle and rigidity of asymptotic renormings}
	\label{chapter:Gorelik}
	
	%The goal of this chapter is to study the stability of the classes $\textsf{T}_p$, $\textsf{A}_p$ and $\textsf{N}_p$ under Lipschitz equivalences, uniform homeomorphisms and coarse-Lipschitz equivalences. In the first section, we detail the elementary properties of coarse-Lipschitz equivalences. In Section 2, we prove the Gorelik principle, which will be the main nonlinear tool of this chapter. We apply it in Section 3, to prove that $\textsf{T}_p$ is stable under Lipschitz equivalences. In Section 4, we show how this can be used to deduce the Lipschitz rigidity of $\co$ and its subspaces. In Section 5, we detail a refinement of the techniques from Section 3 that will lead to the coarse-Lipschitz rigidity of $\textsf{A}_p$ and $\textsf{N}_p$. We conclude, with a few comments on the class $\textsf{A}_\infty=\textsf{N}_\infty$ of asymptotic-$\co$ spaces.
	
	In this chapter, we study several versions of a principle originally due to Gorelik. This principle was formulated by Johnson, Lindenstrauss and Schechtman and used in their proof of the uniform rigidity of $\ell_p$ for $p\in(1,\infty)$. There are now alternative proofs of the uniform rigidity of $\ell_p$ in the reflexive range that do not require the Gorelik principle. However, various versions of the Gorelik principle are still needed to show finer results such as the nonlinear rigidity of certain asymptotic renormings. These are the applications of the Gorelik principle that we will discuss in this chapter.
	Of course, these finer results also imply the nonlinear rigidity of some classes of Banach spaces, and for historical purposes, we will provide a proof of the uniform rigidity of $\ell_p$ for the range $p\in(2,\infty)$ that relies on the Gorelik principle.
	In Section \ref{sec:strategy}, we discuss a strategy that leads to the uniform rigidity of $\ell_p$. In Section \ref{sec:Gorelik} we state and prove several variants of the Gorelik principle. We apply the Lipschitz variant in Section \ref{sec:Lip-rigidity-AUS} to show that the existence of an equivalent asymptotically uniformly smooth norm with power-type $p\in (1,\infty]$ is stable under Lipschitz equivalences. From this, we deduce the Lipschitz stability of subspaces of $\co$ and the Lipschitz rigidity of $\co$ in Section \ref{sec:Lip-rigidity-c_0}. In Section \ref{sec:CL-rigidity-AUS}, we use a coarse-Lipschitz variant of the Gorelik principle to show how asymptotically uniformly smooth renormings are stable under coarse-Lipschitz equivalences. In Section \ref{sec:CL-rigidity-classes}, we discuss the coarse-Lipschitz rigidity of $\ell_p$ and of the classes $\textsf{A}_p$, $\textsf{N}_p$ and $\textsf{A}_\infty=\textsf{N}_\infty$ the class of asymptotic-$\co$ spaces.
	
	\section{\texorpdfstring{Uniform rigidity of $\ell_p$: the strategy}{Uniform rigidity of : the strategy}}
	\label{sec:strategy}
	
	To show that every Banach space that is uniformly homeomorphic to $\ell_p$ must be isomorphic to $\ell_p$, we will use several important and beautiful results from the structure theory of Banach spaces dating back to the early 1970s together with Lipschitz linearization techniques from the early 1980s and finally results about the nonlinear and asymptotic structure of Banach spaces from the late 1990s. Since some of the arguments require $p$ to be in the reflexive range, some fundamental open problems will remain for the extreme values of $p$. The strategy we discuss below to attack uniform rigidity problems can be found in a 1996 groundbreaking paper by Johnson, Lindenstrauss and Schechtman \cite{JLS1996}.
	
	So, assume that $X$ is a Banach space such that $X$ is uniformly homeomorphic to $\ell_p$. At this point, we can already assert that $X$ must be separable. The first major step is to realize that under this nonlinear assumption, $X$ must be isomorphic to a subspace of $L_p[0,1]$. This can be done in essentially two different ways.  One could invoke Ribe's rigidity theorem \cite{Ribe1976}, which says that two separable uniformly homeomorphic Banach spaces must be crudely finitely representable into each other. In particular, $X$ must be infinite-dimensional. It then follows from the classical relationship between finite representability and ultraproducts that $X$ must be isomorphic to a subspace of an ultraproduct $\ell_p^\cU$ of $\ell_p$, where $\cU$ is a nonprincipal ultrafilter on some set $I$. 
	
	When $p=2$,  this is already enough to conclude the uniform rigidity of $\ell_2$. Indeed, any ultraproduct $\ell_2^\cU$ is a (nonseparable) Hilbert space. Then, since $X$ is separable and infinite-dimensional, we obtain that $X$ is isomorphic to $\ell_2$. 
	
	\begin{theo}[Uniform rigidity of $\ell_2$]
		\label{thm:uniform-rigidity-ell_2}
		If $X$ is uniformly homeomorphic to $\ell_2$, then $X$ is linearly isomorphic to $\ell_2$.
	\end{theo}
	
	Alternatively, one could have used differentiability techniques instead of Ribe's rigidity theorem, but then we must restrict $p$ in the reflexive range $(1,\infty)$. The argument would go as follows. We already know that uniformly homeomorphic Banach spaces must be coarse-Lipschitzly equivalent, and this automatically implies that $X$ must be infinite-dimensional. In this process, we gained metric accuracy but lost topological data. To regain control of the topological structure, we can pass to an ultraproduct since then $X^\cU$ is Lipschitz homeomorphic to $\ell_p^\cU$ and in turn $X$ bi-Lipschitzly embeds into $\ell_p^\cU$. Any ultraproduct $\ell_p^\cU$ is isometrically isomorphic to $L_p(\mu)$ for some measure space $(\Omega,\Sigma,\mu)$. Since $L_p(\mu)$ has RNP when $1<p<\infty$ and since $X$ is separable, we can deduce by differentiating (Theorem \ref{thm:infinite-Rad}) that $X$ is linearly isomorphic to a subspace of $L_p(\mu)$. When $p=2$, we can then conclude as above and obtain Theorem \ref{thm:uniform-rigidity-ell_2}. It is worth pointing out that the first proof of Theorem \ref{thm:uniform-rigidity-ell_2} is due to Enflo \cite{Enflo1970} and predates the full development of differentiability techniques as well as Ribe's rigidity theorem.
	
	The above argument for $p\neq 2$ breaks down since $(\ell_p)^{\cU}$ is an $L_p$-space which has a much richer subspace structure. Nevertheless, in the reflexive range $p\in(1,\infty)$, $X$ must be reflexive, and by considering the map sending a bounded sequence in $X$ to its weak limit along $\cU$, this has the nice consequence of making it complemented in $X^{\cU}$. We can now apply statement $(i)$ of Theorem \ref{thm:Lipschitz-isomorphisms} to deduce that $X$ is isomorphic to a complemented subspace of $L_p(\mu)$ for some measure space $(\Omega,\Sigma,\mu)$. Here again, since $X$ is separable, a classical argument from measure theory yields that $X$ is isomorphic to a complemented subspace of a separable $L_p$-space, which we can take without loss of generality to be $L_p[0,1]$ (see \cite[p. 9]{Wojtaszczyk1996} for instance for more details). 
	Now, the issue is that besides $\ell_p$, $L_p[0,1]$ has many more complemented subspaces; $\ell_2$ being one of them. 
	However, it follows from the results of Johnson and Odell \cite{JohnsonOdell1974} that any infinite-dimensional complemented subspace of $L_p[0,1]$ that does not contain an isomorphic copy of $\ell_2$ is isomorphic to $\ell_p$. We recall how to derive this last claim. An infinite-dimensional complemented subspace $X$ of $L_p[0,1]$ with $p\in(1,\infty)$, is either a $\cL_p$-space (see definition in Exercise \ref{ex:script}) or is isomorphic to a Hilbert space. The second alternative is ruled out if $X$ does not contain $\ell_2$. But since Johnson and Odell proved that any separable $\cL_p$-space ($1<p\neq 2 <\infty$) that does not contain $\ell_2$ must be isomorphic to $\ell_p$, the claim follows. 
	From the discussion above we can extract the following proposition.
	%Note that the preceding arguments would have gone through with the weaker assumption that $X$ coarse-Lipschitz embeds into $\ell_p$ albeit assuming upfront that $X$ is infinite-dimensional. 
	
	\begin{prop}
		\label{prop:rigidity-ell_p}
		Let $p\in(1,\infty)$.
		If $X$ is uniformly homeomorphic to $\ell_p$ and $X$ does not contain an isomorphic copy of $\ell_2$, then $X$ is isomorphic to $\ell_p$.
	\end{prop}
	
	According to Proposition \ref{prop:rigidity-ell_p}, we have two paths leading to the uniform rigidity of $\ell_p$. Either we can try to show that if $X$ is uniformly homeomorphic to $\ell_p$, then $X$ must retain enough of the geometry of $\ell_p$ so that it will prevent $X$ from containing a copy of $\ell_2$, or we can try to show that $\ell_2$ does not coarse-Lipschitz embeds into $\ell_p$ for any $p\in (1,\infty)$. As we will see, these two approaches turned out to be successful.
	
	It is worth pointing out that a desirable structural property of $\ell_p$ serving as an obstruction to containing $\ell_2$ cannot be of a local nature. Indeed, by Dvoretzky's theorem $\ell_2$ is finitely representable in any infinite-dimensional Banach space, and loosely speaking no Banach space can have better local properties than the Hilbert space. This observation will be expanded upon in Chapter \ref{chapter:Hamming}. It was observed in \cite{JLS1996} that Bourgain's extension of Enflo's approximate midpoint argument used to show that $L_p$ is not uniformly homeomorphic to $\ell_p$ when $p\in(1,2)$, actually gives that $\ell_2$ does not coarse-Lipschitz embeds into $\ell_p$ when $p\in(1,2)$ (cf. Remark \ref{rem:Bourgain-Enflo}), thus resolving positively the uniform rigidity of $\ell_p$ in this range. To settle the situation when $p\in(2,\infty)$, Johnson, Lindenstrauss and Schechtman formalized an idea of Gorelik used in his proof showing that $L_p$ is not uniformly homeomorphic to $\ell_p$ when $p\in(2,\infty)$; they called it the \emph{Gorelik principle}.
	Using the Gorelik principle and the approximate midpoint principle, they showed that a Banach space that is uniformly homeomorphic to $\ell_p$ cannot have a quotient with an unconditional basis with a lower $r$-estimate for $r<p$. Then, they went on to show that this fact, when $p\in(2,\infty)$ together with additional classical results from Banach space theory, prevents the containment of $\ell_2$.
	
	\begin{theo}[Uniform rigidity of $\ell_p$]
		\label{thm:uniform-rigidity-ell_p}
		Let $p\in (1,\infty)$. Any Banach space that is uniformly homeomorphic to $\ell_p$ must be linearly isomorphic to $\ell_p$.
	\end{theo}
	
	It is worth pointing out that the arguments above are valid under a weaker coarse-Lipschitz equivalence assumption and we will provide the missing details of the proof of Theorem \ref{thm:uniform-rigidity-ell_p} in the proof of the more general Theorem \ref{thm:CL-rigidity-ell_p}.
	
	\begin{rema}
		\label{rem:uniform-rigidity-ell_p}\ 
		\begin{enumerate}
			\item For $p\in(1,\infty)$, the fact that $X$ is uniformly homeomorphic to $\ell_p$ implies that $X$ is isomorphic to a \emph{complemented} subspace of $L_p[0,1]$ can also be obtained using Ribe's approach. But this time we need the strongest version of Ribe's rigidity theorem from \cite{Ribe1978}.
			\item As we will see in Chapter \ref{chapter:Hamming}, Kalton and Randrianarivony studied concentration inequalities on certain nonlocally finite graphs \cite{KaltonRandrianarivony2008} and a by-product of their work is that $\ell_2$ does not coarsely embed into $\ell_q$ whenever $q\in(2,\infty)$, thereby giving a second proof of the uniform rigidity of $\ell_p$ in this range. This proof does not rely on the Gorelik principle.
			\item In \cite{JLS1996} it was also proved using the Gorelik principle that $\ell_{p_1}\oplus \dots \oplus \ell_{p_n}$ is uniformly rigid for $1<p_1<\dots< p_n<2$ or $2<p_1<\dots<p_n<\infty$. We postpone the proof of this result to the next chapter, where we shall replace the Gorelik principle with the use of concentration inequalities on Hamming graphs. As we will see, this method, due to N. Kalton and L. Randrianarivony \cite{KaltonRandrianarivony2008}, will allow us to also show the uniqueness of the uniform rigidity of $\ell_{p_1}\oplus \dots \oplus \ell_{p_n}$ for $1<p_1<\dots<p_n<\infty$, all different from $2$.
		\end{enumerate}
	\end{rema}

	\section{The Gorelik Principle}
	\label{sec:Gorelik}
	
	The tool that we shall now describe is the Gorelik principle. It was initially devised by Gorelik in \cite{Gorelik1994} to prove that for $p\in (2,\infty)$, $\ell_p$ is not uniformly
	homeomorphic to $L_p$. The Gorelik principle was explicitly formulated by Johnson, Lindenstrauss and Schechtman in \cite{JLS1996} as follows (cf Exercise \ref{exe:Gorelik} for a precise quantitative statement):
	
	``\textbf{The Gorelik Principle.} A uniform homeomorphism between Banach spaces cannot take a large ball of a finite-codimensional subspace into a small neighborhood of a subspace of infinite codimension."
	
	The Gorelik principle appeared naturally in \cite{Gorelik1994} when trying to estimate the size of the approximate midpoint sets. In Chapter \ref{chapter:AMP_I}, an upper estimate on the size of the approximate midpoint sets relies on the fact that the unit ball of an infinite-dimensional Banach space is included in a multiple of the unit ball of any of its finite-codimensional subspace up to a compact perturbation (see Lemma \ref{lem:compact-pert}). One way to look at the Gorelik principle is as a nonlinear extension of this elementary observation.
	
	\begin{theo}[The Gorelik principle, Lipschitz version]\label{thm:Gorelik-Lip}
		Let $X$ and $Y$ be two Banach spaces. Assume that $f\colon X\to Y$ is a homeomorphism such that $\Lip(f^{-1})<D$, then for any $s>0$ and any finite-codimensional subspace $X_0$ of $X$, there is a compact subset $K$ of $Y$ so that
		\begin{equation}
			sB_Y \subset f(2DsB_{X_0}) + K.
		\end{equation}
	\end{theo}
	
	\begin{rema}
		By taking $f=id \colon X \to Y=X$ in Theorem \ref{thm:Gorelik-Lip} and $s=1$ one recovers Lemma \ref{lem:compact-pert}.
	\end{rema}
	
	For the finer applications, we will in fact need a coarse-Lipschitz variant of the Gorelik principle (Theorem \ref{thm:Gorelik-CL}), from which Theorem \ref{thm:Gorelik-Lip} is an immediate consequence. A crucial ingredient leading up to the proof of the Gorelik principle relies on Brouwer's fixed point theorem and the existence of Bartle-Graves continuous selectors (see Lemma VII.3.2 in \cite{DGZ1993}).
	
	\begin{prop}
		\label{prop:Brouwer} Let $X_0$ be a finite-codimensional subspace
		of a Banach space $X$ and let $0<c<d$. Then,  there exists a compact subset $K$ of
		$dB_X$ such that for every continuous map $\phi\colon K\to X$ satisfying $\norm{\phi(x)-x}\le c$ for
		all $x\in K$, we have that $\phi(K)\cap X_0\neq \emptyset$.
	\end{prop}
	
	\begin{proof} 
		Denote by $Q$ the quotient map from $X$ to $X/X_0$. The Bartle-Grave selection theorem ensures the existence of a (not necessarily linear) continuous lifting, i.e. a continuous map $\psi\colon X/X_0 \to X$ such that $Q\circ \psi=Id_{X/X_0}$ and $K:=\psi(cB_{X/X_0})$ is included in $dB_X$. Note that, since $X/X_0$ is finite-dimensional and $\psi$ continuous, $K$ is a compact subset of $dB_X$. Assume now that $\phi\colon K \to X$ is continuous and such that $\norm{\phi(x)-x}\le c$ for all $x\in K$. Then, set $g(y)=y-(Q\circ \phi \circ \psi)(y)$, for $y\in cB_{X/X_0}$. Then, for all $y\in cB_{X/X_0}$
		\begin{equation*}
			\norm{g(y)} = \norm{ Q\big(\psi(y)-\phi(\psi(y))\big)}\le \norm{\psi(y)-\phi(\psi(y))}\le c.
		\end{equation*}
		Since $cB_{X/X_0}$ is clearly nonempty, compact and convex, we are in a position to apply Brouwer's fixed point theorem, and the existence of a fixed point for $g$ in $cB_{X/X_0}$ yields the desired conclusion.
	\end{proof}
	
	We will now state a general version of the Gorelik principle that can be used to study coarse-Lipschitz equivalent Banach spaces. 
	%It also contains the extra information needed to study Lipschitz isomorphisms.
	
	\begin{theo}[The Gorelik principle, coarse-Lipschitz version]
		\label{thm:Gorelik-CL}
		Let $X$ and $Y$ be two Banach spaces. Assume that $f\colon X\to Y$ and $g\colon Y\to X$ are continuous and  that there exist constants $A,B,C > 0$ such that for all $y,y'\in Y$
		\begin{equation*}
			\norm{g(y)-g(y')}\le A\norm{y-y'} + B,
		\end{equation*}
		\begin{equation*}
			\sup_{x\in X} \norm{(g\circ f)(x)-x}\le C \textrm{ and }
			\sup_{y\in Y} \norm{(f\circ g)(y)-y}\le C.
		\end{equation*}
		Let $\lambda\in(0,1)$. Then, for all $\alpha \ge \frac{2(B+C)}{A(1-\lambda)}$ and all finite-codimensional subspace $X_0$ of $X$, there is a compact subset $K$ of $Y$ so that
		\begin{equation*}
			\lambda\alpha B_Y\subset f(2A\alpha B_{X_0}) + CB_Y + K.
		\end{equation*}
	\end{theo}
	
	\begin{proof}
		Let $\mu:=\frac{1+\lambda}{2}\in(0,1)$, $\alpha_0:=\frac{2(B+C)}{A(1-\lambda)}=\frac{B+C}{A(\mu-\lambda)}$ and $\alpha\ge \alpha_0$. Let also $X_0$ be a finite-codimensional subspace of $X$. It follows from Proposition \ref{prop:Brouwer} that there exists a compact subset $K_0$ of
		$A\alpha B_X$ such that for every continuous map $\phi\colon K_0\to X$ satisfying $\norm{\phi(x)-x}\le \mu A\alpha$ for all $x\in K_0$, we have that $\phi(K_0)\cap X_0\neq \emptyset$.
		
		Consider now $y\in \lambda\alpha B_Y$ and define $\phi\colon K_0\to X$ by
		$\phi(x)=g(y+f(x))$. Then, $\phi$ is clearly continuous and for all $x\in K_0$,
		\begin{align*}
			\norm{\phi(x)-x}&\le \norm{g(y+f(x))-g(f(x))} + C\\
			&\le A\norm{y} + B + C\le A\lambda \alpha + B + C\le \mu A\alpha.
		\end{align*}
		Hence, there exists $x_0\in K_0$ so that $\phi(x_0)\in X_0$. Since $\norm{x_0}\le A\alpha$ and $\norm{\phi(x_0)-x_0}\le \mu A\alpha$, we have that $\phi(x_0)\in 2A\alpha B_{X_0}$.
		Finally, the decomposition 
		\begin{equation*}
			y = f(\phi(x_0)) - (f\circ g)(y + f(x_0)) + (y + f(x_0))   - f(x_0)
		\end{equation*}
		together with the fact that $\norm{(f\circ g)(y+f(x_0))-(y+f(x_0))}\le C$ implies that $y\in f(2A\alpha B_{X_0}) + CB_Y + K$, where $K:=-f(K_0)$ is a compact subset of $Y$.
	\end{proof}
	
	\begin{rema} 
		Theorem \ref{thm:Gorelik-Lip} is quite miraculous. It states in a quantitative way that for any finite-codimensional subspace $Z$ of $X$, the image of a dilation of the unit ball in $Z$ is so large that up to a compact perturbation (plus a uniformly bounded perturbation for Theorem \ref{thm:Gorelik-CL}), it contains a ball of $Y$ of proportional size. This is similar to the behavior of linear isomorphisms, and it indicates how this will provide information on the nonlinear classification of Banach spaces. We leave it to the reader to check that Theorem \ref{thm:Gorelik-Lip} is a consequence of Theorem \ref{thm:Gorelik-CL}, and we invite her or him to write down the details of the simpler proof of Theorem \ref{thm:Gorelik-Lip}, where the ideas may appear less mysterious.
	\end{rema}

	\section[Lipschitz rigidity of AUS norms]{Lipschitz rigidity of asymptotically uniformly smooth renormings}
	\label{sec:Lip-rigidity-AUS}
	
	The main result of this section, which appeared first in \cite{GKL2000} in a separable setting, states that the existence of an equivalent asymptotically uniformly smooth norm is stable under Lipschitz isomorphisms. A proof can also be found in the second edition of Albiac-Kalton textbook \cite[paragraph 14.6]{AlbiacKalton2016}. The general case can be deduced by arguments of separable saturation and separable determination of the moduli. However, we shall detail here a direct proof in the general case, taken from  \cite{DaletLancien2017}. 
	%The only modification is that we deal directly with the definition of the asymptotic moduli instead of using weak$^*$-null or weakly null sequences.
	
	\begin{theo}
		\label{thm:Lip-rigidity-AUS}
		Let $X$ and $Y$ be two Banach spaces that are Lipschitzly isomorphic and let $f\colon X\to Y$ be a Lipschitz equivalence. For every $D>\dist(f)$, there exists an equivalent norm $\abs{\cdot}_Y$ on $Y$ such that $\norm{\cdot}_Y\le \abs{\cdot}_Y\le D\norm{\cdot}_Y$ and for all $t\in [0,1]$,
		\begin{equation*}
			\bar{\rho}_{\abs{\cdot}_Y}\Big(\frac{t}{96 D}\Big)\le \bar{\rho}_{\norm{\cdot}_X}(t).
		\end{equation*}
	\end{theo}
	
	It is quite remarkable that the quantitative loss in Theorem \ref{thm:Lip-rigidity-AUS} is only of the order of a multiplicative factor of the distortion. We can thus deduce the Lipschitz rigidity of several classes of Banach spaces that we have studied in Chapter \ref{chapter:Szlenk}.
	
	%It then follows clearly that for any $p\in (1,\infty]$, admitting an equivalent $p$-AUS norm is stable under Lipschitz equivalences. As we have seen in Chapter \ref{chapter:Szlenk}, for $p\in (1,\infty]$, the class $\textsf{T}_p$ coincides with the class of Banach spaces admitting an equivalent $p$-AUS norm. Also, the class $\textsf{\emph{Sz}}_\omega$ of all Banach space with Szlenk index $\omega$ coincides with $\bigcup_{1<p\le \infty}\textsf{T}_p$. So, we can state:
	
	\begin{coro}
		\label{cor:Lip-rigidity-AUS}
		The following classes are stable under Lipschitz equivalences:
		\begin{enumerate}
			\item For any $p\in (1,\infty]$, the class of Banach spaces with an equivalent norm that is asymptotically uniformly smooth with power type $p$ (i.e. $\langle p$-$\AUS\rangle$ for $p\in (1,\infty)$ or $\langle \AUF \rangle$ for $p=\infty$), or equivalently the class $\sT_p$.
			\item The class $\langle \AUS \rangle$ of all Banach spaces admitting an equivalent asymptotically uniformly smooth norm, or equivalently the class of Banach spaces with Szlenk index equal to $\omega$, namely $\mathbf{Sz}_\omega$.
		\end{enumerate}
	\end{coro}
	
	There is a relatively natural way to construct an equivalent norm on a Banach space $Y$ given a Lipschitz equivalence between $X$ and $Y$: one can simply take the Minkowski gauge of the closed convex hull of the image of the Lipschitz equivalence under the elementary molecules of the vector-valued Lipschitz-free space $\cF(X,Y)$, i.e. 
	\begin{equation*}
		\abs{y}_Y:= \inf\left\{\lambda>0 \colon \frac{1}{\lambda} y\in \overline{\rm conv}\,\Big\{\frac{\delta_x-\delta_x'}{\norm{x-x'}_X}(f),\ x\neq x' \in X\Big\}\right\}.
	\end{equation*}
	Working with this norm is not very convenient, and it will be easier to deal with its dual norm. Therefore, we will turn our attention to the behavior of weak$^*$ asymptotically uniformly convex norms under Lipschitz equivalences. 
	Theorem \ref{thm:Lip-rigidity-AUS} will be an immediate corollary of the usual AUS-AUC$^*$ duality from Corollary \ref{cor:Young2}
	and the stability under Lipschitz equivalences of the existence of an equivalent weak$^*$ asymptotically uniformly convex norm. 
	
	\begin{theo}
		\label{thm:Lip-rigidity-AUC*}
		Let $X$ and $Y$ be two Banach spaces that are Lipschitzly isomorphic and let $f\colon X\to Y$ be a Lipschitz equivalence. For every $D>\dist(f)$, there exists an equivalent norm $\abs{\cdot}_Y$ on $Y$ such that $\norm{\cdot}_Y\le \abs{\cdot}_Y\le D\norm{\cdot}_Y$ and for all $t\in [0,1]$,
		\begin{equation*}
			\bar{\delta}^*_{\abs{\cdot}_Y}(t)\ge \bar{\delta}^*_{\norm{\cdot}_X}\Big(\frac{t}{8D}\Big).
		\end{equation*}
	\end{theo}
	
	\begin{proof}
		Assume as we may that $\Lip(f)\le 1$ and let $D>\Lip(f^{-1})$. The equivalent norm on $Y$ will be defined as the Minkowski gauge of a well-chosen closed convex symmetric set. Let
		\begin{equation*}
			C:=\overline{\rm conv}\,\Big\{\frac{f(x)-f(x')}{\norm{x-x'}_X},\ x\neq x' \in X\Big\}.
		\end{equation*}
		Clearly, $C$ is closed convex symmetric and $C\subset B_Y$. We claim that $\frac{1}{D}B_Y\subset C$. To see this, let $y\in Y$ such that $\norm{y}_Y=\frac{1}{D}$. For $t\in [0,+\infty)$, denote by $x_t=f^{-1}(ty)$. We have that $\norm{x_1-x_0}_X\le 1$ and $\norm{x_D-x_0}_X\ge 1$. So, there exists $t_0\in [1,D]$ such that $\norm{x_{t_0}-x_0}_X=1$. It follows that $t_0y\in C$. Since $C$ is convex and symmetric, we deduce that $\frac1D B_Y\subset C$. So, if we let $\abs{\cdot}_Y$ be the Minkowski functional of $C$, we have that $\abs{\cdot}_Y$ is an equivalent norm on $Y$ such that $\norm{\cdot}_Y\le \abs{\cdot}_Y\le D\norm{\cdot}_Y$. It then follows from the definition of $C$ that the dual norm of an element $y^*\in Y^*$ with respect to the new norm on $Y$ is given by 
		\begin{equation}\label{eq:AUS-rigidity-explicit-formula}
			\abs{y^*}_{Y^*}=\sup\Big\{\frac{\langle y^*,f(x)-f(x')\rangle}{\norm{x-x'}_X},\ x\neq x'\in X\Big\}.
		\end{equation}
		Let $t\in (0,1]$. Assume as we may that $\bar{\delta}^*_{\norm{\cdot}_X}\big(\frac{t}{8D}\big)>0$ and pick $0<\delta <\bar{\delta}_{\norm{\cdot}_X}^*(\frac{t}{8D})$. Given $y^*\in Y^*$ such that $\abs{y^*}_{Y^*}=1$ we need to find a weak$^*$ closed finite-codimensional subspace $Z$ of $Y^*$ so that $\abs{y^*+z^*}_{Y^*}-1\ge \delta$ for all $z^*\in Z$ with $\abs{z^*}_{Y^*}=t$. So, let $\eta >0$ and pick $x \neq x'\in X$ such that 
		\begin{equation*}
			\langle y^*,f(x)-f(x')\rangle \ge (1-\eta)\norm{x-x'}_X.
		\end{equation*}
		We may assume, using translations if needed, that $x'=-x$ and $f(x')=-f(x)$, so that we have
		\begin{equation}\label{eq:AUS-rigidity-eq1}
			\langle y^*,f(x)\rangle \ge (1-\eta)\norm{x}_X.
		\end{equation}
		To make the connection with the geometry of $X$ we use statement (b) in Proposition \ref{prop:Young} to assert that $\bar{\rho}_{\norm{\cdot}_X}(\frac{8D\delta}{t})<\delta$. Consequently, there exists a finite-codimensional subspace $X_0$ of $X$ such that for all $z \in \frac{8D\delta \norm{x}_X}{t}\,B_{X_0}$,
		\begin{equation}\label{eq:AUS-rigidity-eq2}
			\norm{x+z}_X\le (1+\delta)\norm{x}_X.
		\end{equation}
		Now comes the key application of the Gorelik principle. Pick $s<\frac{4\delta \norm{x}_X}{t}$. It follows from the Gorelik principle (Theorem \ref{thm:Gorelik-Lip}) that there exists a compact subset $K$ of $Y$ such that
		\begin{equation}\label{eq:AUS-rigidity-eq3}
			sB_Y\subset f\Big(\frac{8D\delta \norm{x}_X}{t}B_{X_0}\Big) + K.
		\end{equation}
		Consider a finite $\eps$-net $F$ of $K$, for some $\vep>0$ to be chosen small enough later (only depending on $t$ and $\eta$) and let $E$ be the finite-dimensional subspace of $Y$ spanned by $F \cup \{f(x)\}$. Note that $E^\perp$ is a weak$^*$ closed finite-codimensional subspace of $Y^*$.  For any $z^*\in Z:=E^\perp$ such that $\abs{z^*}_{Y^*}=t$, we have $\norm{z^*}_{Y^*}\ge t$ and, if $\eps>0$ was initially chosen small enough, we deduce from \eqref{eq:AUS-rigidity-eq3} that there exists $z_0 \in \frac{8D\delta \norm{x}_X}{t}B_{X_0}$ such that
		\begin{equation}\label{eq:AUS-rigidity-eq4}
			\langle z^*, f(z_0)\rangle \ge (s-\eta)t.
		\end{equation}
		It remains to verify that $\abs{y^*+z^*}_{Y^*}\ge 1+\delta$. To this end we will test $y^*+z^*$ against the vector $\frac{f(z_0)-f(x')}{\norm{z_0-x'}_X} \in C$. Since $\langle z^*,f(x)\rangle =0$ we have
		\begin{align*}
			\langle y^* + z^*, f(z_0) - f(x')\rangle & = \langle y^*, f(z_0) + f(x) \rangle + \langle z^*,f(z_0)\rangle\\
			& \stackrel{\eqref{eq:AUS-rigidity-eq4}}{\ge} \langle y^*, f(z_0) - f(x) \rangle + 2\langle y^*, f(x)\rangle + (s-\eta)t \\
			& \stackrel{\eqref{eq:AUS-rigidity-eq1}}{\ge}  2(1-\eta)\norm{x}_X  + (s-\eta)t - \langle y^*, \frac{f(x) - f(z_0)}{\norm{x -z_0}}\rangle\norm{x-z_0}\\
			& \stackrel{\eqref{eq:AUS-rigidity-eq2}}{\ge} 2(1-\eta)\norm{x}_X  + (s-\eta)t - (1+\delta)\norm{x}_X,
		\end{align*}
		where in the last inequality we also used the fact that $\abs{y^*}_{Y^*}=1$. Note that since $\eta$ can be taken arbitrarily small, we can assume without loss of generality that $z_0\neq x'$. Using again the definition of $\abs{\cdot}_Y$ and the fact that $\norm{z_0 - x'}_X \le (1+\delta)\norm{x}_X$ thanks to \eqref{eq:AUS-rigidity-eq2} we get
		\begin{equation*}
			\abs{y^*+z^*}_{Y^*} \ge \Big((1-\eta)\norm{x}_X-(\delta+\eta)\norm{x}_X + (s-\eta)t\Big)\Big((1+\delta)\norm{x}_X\Big)^{-1},
		\end{equation*}
		and letting $s$ tend to $\frac{4\delta \norm{x}_X}{t}$ and $\eta$ tend to $0$, we deduce that $\abs{y^*+z^*}_{Y^*} \ge \frac{1+3\delta}{1+\delta}$. Therefore, since $\delta\in(0,1)$ we have $\abs{y^*+z^*}_{Y^*}>1+\delta$ and 
		\begin{equation*}
			\bar{\delta}^*_{\abs{\cdot}_Y}(y^*,t) \ge \delta.
		\end{equation*}
		Finally, after letting $\delta$ tend to $\bar{\delta}^*_X(\frac{t}{8D})$ in the above estimate that does not depend on $y^*$ in the unit sphere of $\abs{\cdot}_Y$, we obtain $\bar{\delta}^*_{\abs{\cdot}_Y}(t)\ge \bar{\delta}^*_{\norm{\cdot}_X}\big(\frac{t}{8D}\big)$. 
	\end{proof}
	
	\begin{rema} As noted in Remark XI.4.9 in \cite{Godefroy-Baire}, this renorming is quite canonical. Indeed, let $\bar{f}\colon \cF(X) \to Y$ be the bounded linear operator such that $\bar{f} \circ \delta_X= f$. Then, $\bar{f}$ is a quotient map and $|\ |_Y$ is the corresponding quotient norm. See Exercise \ref{ex:Gorelikandfreespaces} for the details. 
	\end{rema}
	
	\section{\texorpdfstring{Nonlinear rigidity of $\co$}{Nonlinear rigidity of }}
	\label{sec:Lip-rigidity-c_0}
	In this section, we draw a few more important consequences of Theorem \ref{thm:Lip-rigidity-AUS}. 
	We have seen (Theorem \ref{thm:AUF->subspace-c_0}) that a separable infinite-dimensional Banach space admits an equivalent asymptotically uniformly flat norm if and only if it is isomorphic to an infinite-dimensional subspace of $\co$. An immediate consequence of Theorem \ref{thm:Lip-rigidity-AUS} (or Corollary \ref{cor:Lip-rigidity-AUS}) is thus the stability under Lipschitz equivalences of the class of Banach spaces that are isomorphic to a subspace of $\co$. 
	
	\begin{coro}
		\label{cor:Lip-rigidity-subspaces-c_0}
		The class of Banach spaces that are isomorphic to a subspace of $\co$ is stable under Lipschitz equivalences.
	\end{coro}
	
	In particular, if a Banach space is Lipschitz equivalent to $\co$, then it is isomorphic to a closed infinite-dimensional subspace of $\co$, but we do not know at this point if this subspace is isomorphic to $\co$. Thankfully, Johnson and Zippin \cite{JohnsonZippin1972} told us that checking whether a closed infinite-dimensional subspace of $\co$ is isomorphic to $\co$ is equivalent to verifying whether it is a $\cL_\infty$-space. It is not difficult to see that $\co$ is a $\cL_\infty$-space (see Exercise \ref{ex:script}). Since Henrich and Mankiewicz \cite{HeinrichMankiewicz1982} showed that being a $\cL_\infty$-space is stable under Lipschitz equivalences (in fact under coarse-Lipschitz equivalences), it follows from this discussion and Corollary \ref{cor:Lip-rigidity-subspaces-c_0} that the space $\co$ is Lipschitzly rigid. 
	
	\begin{theo}[Lipschitz rigidity of $\co$]
		\label{thm:Lip-rigidity-c_0}
		A Banach space is Lipschitzly isomorphic to $\co$ if and only if it is linearly isomorphic to $\co$.
	\end{theo}
	
	\begin{rema}
		For a proof of the nonlinear stability of the class of $\cL_\infty$-spaces, we refer the reader to the original paper \cite{HeinrichMankiewicz1982}, or the presentation given in \cite[Theorem 10.5 and Appendix F]{BenyaminiLindenstrauss2000}.    
	\end{rema}
	
	The only apparent limitation to settle the uniform, or coarse-Lipschitz, rigidity of $\co$ stems from the fact that Theorem \ref{thm:Lip-rigidity-AUS} is valid for Lipschitz equivalences only. 
	
	\begin{prob}[Uniform or coarse-Lipschitz rigidity of $\co$]
		\label{prob:uniform-rigidity-c_0} 
		Is a Banach space uniformly (or coarse-Lipschitzly) equivalent to $\co$ linearly isomorphic to $\co$?
	\end{prob} 
	
	A positive answer to Problem \ref{prob:uniform-rigidity-c_0} would follow from a positive answer to the following problem.
	
	\begin{prob}[Uniform or coarse-Lipschitz stability of $\langle \AUF\rangle=\sT_\infty$]
		\label{prob:uniform-rigidity-AUF}
		Is the class of Banach spaces admitting an equivalent asymptotically uniformly flat norm, or equivalently, the class $\sT_\infty$, stable under uniform or coarse-Lipschitz equivalences?
	\end{prob}
	
	Problems \ref{prob:uniform-rigidity-c_0} and \ref{prob:uniform-rigidity-AUF} are two important open problems in the nonlinear geometry of Banach spaces. A partial result indicating that Problem \ref{prob:uniform-rigidity-c_0} could have a positive solution will be given later once we push the application of the Gorelik principle a little further (see Theorem \ref{thm:CL-rigidity-c_0} in Section \ref{sec:CL-rigidity-classes}).

	\section[Coarse-Lipschitz rigidity of AUS norms]{Coarse-Lipschitz rigidity of asymptotically uniformly smooth renormings}
	\label{sec:CL-rigidity-AUS}
	Theorem \ref{thm:Lip-rigidity-AUS} is essentially optimal in the sense that if we downgrade the assumption of Lipschitz equivalence to a uniform equivalence, then the fine quantitative conclusion does not have to hold anymore. Indeed, N. Kalton showed in \cite{Kalton2013} that for any $p\in (1,\infty)$, there exist two uniformly homeomorphic separable Banach spaces $X$ and $Y$ so that $X$ is $p$-AUS but $Y$ does not admit any equivalent $p$-AUS norm. In other words, for every $p\in(1,\infty)$, the class $\sT_p=\langle p$-$\AUS\rangle$ is not stable under uniform equivalences. We will discuss Kalton's construction in Chapter \ref{chapter:Counterexamples}. However, if we are willing to give up on the preservation of power-type estimates, what remains true is that the existence of an equivalent norm that is asymptotically uniformly smooth is stable under coarse-Lipschitz equivalences. Thus, the main result of this section is the stability under coarse-Lipschitz equivalence of the class $\langle \AUS \rangle$. This was shown in \cite[Theorem 5.3]{GKL2001} in the separable case. 
	
	In this section, we give a proof that does not require the separability assumption. The idea was indicated in \cite{GLZ2014} and detailed in \cite{DaletLancien2017}. In fact, we will present it as a consequence of a more precise result from \cite{CauseyFovelleLancien2023}: for $p\in (1,\infty]$, the classes $\sA_p$ and $\sN_p$ are stable under coarse-Lipschitz equivalences. As in Section \ref{sec:Lip-rigidity-AUS}, we deduce this result from a dual statement about weak$^*$ asymptotically uniformly convex renormings.

	%As we have already seen, a uniform homeomorphism between Banach spaces is easily seen to be a coarse-Lipschitz  equivalence, while N. Kalton showed in \cite{Kalton2012} that there exist two Banach spaces that are coarse-Lipschitz equivalent but not uniformly homeomorphic. 

	\begin{theo}
		\label{thm:CL-rigidity-AUC*}
		Let $X$ and $Y$ be two Banach spaces and $D>1$. Assume that $f\colon X\to Y$ and $g\colon Y\to X$ are {\bf continuous} with $\Lip_\infty(f)\le 1$, $\Lip_\infty(g)< D$ and that there exists a constant $C\ge 0$ such that $\sup_{x\in X} \norm{(g\circ f)(x)-x}\le C$ and $\sup_{y\in Y}\norm{(f\circ g)(y)-y}\le C$.
		
		Then, for any $\eps$ in $(0,1)$, there exists an equivalent norm $\abs{\cdot}_\vep$ on $Y$ such that $\frac{1}{1+\eps}\norm{\cdot} _Y\le \abs{\cdot}_\vep\le D\norm{\cdot}_Y$ and for all $t\in [0,1]$, 
		\begin{equation}
			\label{eq:CL-rigidity-AUC*-eq1}
			\bar{\delta}^*_{\abs{\cdot}_\vep}(t)\ge \bar{\delta}^*_{\norm{\cdot}_X}\big(\frac{t}{48D^2}\big)-\eps.
		\end{equation}
	\end{theo}
	
	\begin{proof}
		Fix $\vep\in(0,1)$ and for $k\in \bN$ define
		\begin{equation*}
			\label{eq:CL-rigidity-AUC*-eq2}
			C_k:=\overline{\rm conv}\,\Big\{\frac{f(x)-f(x')}{\norm{x-x'}_X}\colon  x,x'\in X,\ \norm{x-x'}_X\ge 2^k\Big\}.
		\end{equation*}
		Then, $(C_k)_{k=1}^\infty$ is a decreasing sequence of closed convex and symmetric subsets of $Y$. Since $\Lip_\infty(f)\le 1$, it follows that $\Lip_{2^k}(f)\le 1+\vep_k$ for all $k\ge 1$, where $(\eps_k)_{k=1}^\infty$ is a sequence of positive numbers tending to 0. In terms of sets, this means that $C_k\subset (1+\eps_k)B_Y$. In particular, there exists $k_0\in \bN$ such that for all $k\ge k_0$,
		\begin{equation*}
			\label{eq:CL-rigidity-AUC*-eq3}
			C_k \subset (1+\eps)B_Y\subset 2B_Y.
		\end{equation*}
		Now, fix $k\in \bN$, $y\in S_Y$ and let $y_0:=f(0)$. It follows easily from our assumptions  that $\lim_{t\to \infty} \norm{g(ty)}_Y=\infty$. Recall also that $f(g(ty))=ty +u_t$ with $\norm{u_t}_Y\le C$. So, for $t$ large enough
		\begin{equation}
			\label{eq:CL-rigidity-AUC*-eq4}
			\frac{f(g(ty))-y_0}{\norm{g(ty)}_Y}=\frac{ty}{\norm{g(ty)}_Y}+\frac{u_t-y_0}{\norm{g(ty)}_Y}\in C_k.
		\end{equation}
		In particular $\big(\frac{ty}{\norm{g(ty)}}\big)_t$ is bounded for $t$ large enough. Then, it follows from the assumption $\Lip_\infty(g)<D$ that there exist $\alpha\ge \frac1D$ and a sequence $(t_n)_n$ tending to $+\infty$ such that $\big(\frac{t_n}{\norm{g(t_ny)}_X}\big)_n$ tends to $\alpha$. Since $C_k$ is closed, we obtain from \eqref{eq:CL-rigidity-AUC*-eq4} that $\alpha y\in C_k$. Finally, we use the fact that $C_k$ is convex and symmetric to deduce that $\frac1D y\in C_k$ and thus, $y\in S_Y$ being arbitrary, we get  that $\frac1D B_Y \subset C_k$.\\
		So, if we let $\abs{\cdot}_k$ be the Minkowski functional of $C_k$, we have that for all $k\ge k_0$, $\abs{\cdot}_k$ is an equivalent norm on $Y$ such that 
		\begin{equation}
			\label{eq:CL-rigidity-AUC*-eq4.5}
			\frac{1}{1+\vep}\norm{\cdot}_Y\le \abs{\cdot}_k\le D\norm{\cdot}_Y.
		\end{equation} 
		It will be useful to describe the dual norm of $\abs{\cdot}_k$, also denoted by $\abs{\cdot}_k$, as follows
		\begin{equation*}
			\label{eq:CL-rigidity-AUC*-eq5}
			y^*\in Y^*\mapsto  \abs{y^*}_k =\sup \left\{ \frac{\langle y^*,f(x)-f(x')\rangle}{\norm{x-x'}_X}\colon 
			x,x'\in X,\ \norm{x-x'}_X\ge 2^k\right\}.
		\end{equation*}
		Note that our assumptions also imply the existence of $B\ge 0$ such that for all $y,y'\in Y$
		\begin{equation}
			\label{eq:CL-rigidity-AUC*-eq6}
			\norm{g(y)-g(y')}_Y\le D\norm{y-y'}_Y + B.
		\end{equation}
		This will enable us to apply the Gorelik principle as it is stated in Theorem \ref{thm:Gorelik-CL}.
		
		\medskip The key lemma is the following.
		
		\begin{lemm}
			\label{lem:CL-rigidity-AUC*} 
			Let $t\in (0,1)$ and assume that $\bar{\delta}^*_{\norm{\cdot}_X}\big(\frac{t}{48D^2}\big)>0$. Let $\eps>0$, $y^*\in Y^*$ such that $ \norm{y^*}_{Y^*}\le D$ and $k_1\in \bN$ such that $k_1\ge k_0$ and
			\begin{equation*}
				\label{eq:CL-rigidity-AUC*-eq7}
				2^{-k_1} <\min\Big\{ \frac{\eps}{2(1+CD)},  \frac{6D^2}{(B+C)}t^{-1}\bar{\delta}^*_{\norm{\cdot}_X}\big(\frac{t}{48D^2}\big) \Big\}.
			\end{equation*}
			%\begin{equation*}
			%  \bar{\delta}^*_{\norm{\cdot}_X}\big(\frac{t}{48D^2}\big)> \frac{(B+C)}{6D^2}t2^{-k_1},
			%\end{equation*}
			%and 
			%\begin{equation*}
			%    2^{-k_1}(CD+1)\le \frac{\eps}{4}.
			%\end{equation*}
			Then, for all $k\ge k_1$, there exists a finite-dimensional subspace $E$ of $Y$ so that for all $z^*\in E^\perp$ with $\frac{t}{2}\le \norm{z^*}_{Y^*}\le tD$, we have
			\begin{equation}
				\label{eq:CL-rigidity-AUC*-eq8}
				\abs{y^*+z^*}_k\ge 2\abs{y^*}_{k+1}-\abs{y^*}_k+\bar{\delta}^*_{\norm{\cdot}_X}\big(\frac{t}{48D^2}\big)-\frac{\eps}{2}.
			\end{equation}
		\end{lemm}
		
		Assuming that the key lemma has been proved, we show how to conclude the proof of the theorem.
		Note that a simple convexity argument (see Exercise \ref{ex:AUC*-modulus}) shows that for any space $Z$, the function $t\mapsto t^{-1}\bar{\delta}^*_Z(t)$ is increasing on $(0,1]$.
		
		Assume first that $\bar{\delta}^*_{\norm{\cdot}_X}\big(\frac{1}{48D^2}\big)\le\frac{\eps}{2}.$
		Then, for any $t \in (0,1]$ we trivially have that
		\begin{equation*}
			\label{eq:CL-rigidity-AUC*-eq9}
			\bar{\delta}^*_{\norm{\cdot}_Y}(t)\ge 0> \bar{\delta}^*_{\norm{\cdot}_X}\big(\frac{t}{48D^2}\big)-\eps,
		\end{equation*}
		and the original norm on $Y$ works.
		
		Assume now that $\bar{\delta}^*_{\norm{\cdot}_X}\big(\frac{1}{48D^2}\big)>\frac{\eps}{2}.$ Since $\bar{\delta}^*_{\norm{\cdot}_X}$ is continuous, there exists $t_0\in (0,1)$ so that $\bar{\delta}^*_{\norm{\cdot}_X}\big(\frac{t_0}{48D^2}\big)=\frac{\eps}{2}.$ As above, we easily have that for any equivalent norm $N$ on $Y$ and any $t\in (0,t_0]$, $\bar{\delta}^*_N(t)\ge 0\ge \bar{\delta}^*_{\norm{\cdot}_X}\big(\frac{t}{48D^2}\big)-\eps$. So, we only have to treat the problem for $t\in [t_0,1]$. Let us pick $k_1\in \bN$ satisfying the assumptions of Lemma \ref{lem:CL-rigidity-AUC*} for $t_0$. It then follows from the monotonicity of $t\mapsto t^{-1}\bar{\delta}^*_{\norm{\cdot}_X}(t)$ that the conclusion of Lemma \ref{lem:CL-rigidity-AUC*} applies for any $t\in [t_0,1]$ and any $k\ge k_1$.\\
		Pick now $N\in \bN$ such that $\frac{4D}{N}<\frac{\eps}{2}$ and define
		\begin{equation}
			\label{eq:CL-rigidity-AUC*-eq10}
			\abs{y^*}_\vep:=\frac1N \sum_{k=k_1+1}^{k_1+N}\abs{y^*}_k
		\end{equation}
		which is the dual norm of a norm, still denoted by $\abs{\cdot}_\vep$, on $Y$. Note that it follows from \eqref{eq:CL-rigidity-AUC*-eq4.5} and \eqref{eq:CL-rigidity-AUC*-eq10} that 
		\begin{equation*}
			\label{eq:CL-rigidity-AUC*-eq11}
			D^{-1}\norm{y^*}_{Y^*}\le \abs{y^*}_\vep \le (1+\eps)\norm{y^*}_{Y^*}\le 2 \norm{y^*}_{Y^*}.
		\end{equation*}
		Let $y^*\in Y^*$ with $\abs{y^*}_\vep=1$. It follows from Lemma \ref{lem:CL-rigidity-AUC*} that for any $t\in [t_0,1]$, there exists a finite-dimensional subspace $E$ of $Y$ so that for all $k\in [k_1,k_1+N]$ and all $z^*\in E^\perp$ such that $\abs{z^*}_\vep=t$, we have
		\begin{equation}
			\label{eq:CL-rigidity-AUC*-eq12}
			\abs{y^*+z^*}_k\ge 2\abs{y^*}_{k+1}-\abs{y^*}_k+\bar{\delta}^*_{\norm{\cdot}_X}\big(\frac{t}{48D^2}\big)-\frac{\eps}{2},
		\end{equation}
		which implies, after summing \eqref{eq:CL-rigidity-AUC*-eq12} for $k_1+1\le k \le k_1+N$ that
		\begin{equation*}
			\label{eq:CL-rigidity-AUC*-eq13}
			\abs{y^*+z^*}_\vep \ge \abs{y^*}_\vep + \frac{2}{N}\big(\abs{y^*}_{k_1+N+1} - \abs{y^*}_{k_1+1}\big)+
			\bar{\delta}^*_{\norm{\cdot}_X}\big(\frac{t}{48D^2}\big)-\frac{\eps}{2}.
		\end{equation*}
		Since $\abs{y^*}_\vep=1$, we have that $\norm{y^*}_{Y^*}\le D$ and $\abs{y^*}_{k_1+1}\le 2D$. So
		\begin{equation*}
			\label{eq:CL-rigidity-AUC*-eq14}
			\abs{y^*+z^*}_\vep \ge 1 + \bar{\delta}^*_{\norm{\cdot}_X}\big(\frac{t}{48D^2}\big)-\frac{\eps}{2}-\frac{4D}{N} > 1 + \bar{\delta}^*_{\norm{\cdot}_X}\big(\frac{t}{48D^2}\big)-\eps.
		\end{equation*}
		This shows that for all $t\in [t_0,1]$, $\bar{\delta}^*_{\abs{\cdot}_\vep}(t)\ge \bar{\delta}^*_{\norm{\cdot}_X}\big(\frac{t}{48D^2}\big)-\eps$, completing our proof.
		
		\medskip 
		It remains to prove the key lemma.
		
		\begin{proof}[Proof of Lemma \ref{lem:CL-rigidity-AUC*}]
			Let $y^*\in Y^*$ with $\norm{y^*}_{Y^*}\le D$ and
			pick, as we may, $0<\delta <\bar{\delta}^*_{\norm{\cdot}_X}(\frac{t}{48D^2})$ such that
			\begin{equation}
				\label{eq:CL-rigidity-AUC*-eq15}
				6D^2\delta2^{k_1}>(B+C)t.
			\end{equation}
			Let $k\ge k_1$, $\eta\in (0,1)$ and choose $x\neq x'\in X$ such that $\norm{x-x'}_X\ge 2^{k+1}$ and
			\begin{equation*}
				\label{eq:CL-rigidity-AUC*-eq16}
				\langle y^*,f(x)-f(x')\rangle \ge (1-\eta)\abs{y^*}_{k+1}\norm{x-x'}_X.
			\end{equation*}
			We may assume that $x'=-x$ and $f(x')=-f(x)$, so that we have
			\begin{equation}
				\label{eq:CL-rigidity-AUC*-eq17}
				\langle y^*,f(x)\rangle \ge (1-\eta)\abs{y^*}_{k+1}\norm{x}_X.
			\end{equation}
			Since $0<\delta <\bar{\delta}^*_{\norm{\cdot}_X}(\frac{t}{48D^2})$, it follows from statement (b) in Proposition \ref{prop:Young} that $\bar{\rho}_{\norm{\cdot}_X}(\frac{48D^2\delta}{t})<\delta$. So, there exists a finite-codimensional subspace $X_0$ of $X$ such that for all $z \in \frac{48D^2\delta \norm{x}_X}{t}\,B_{X_0}$,
			\begin{equation}
				\label{eq:CL-rigidity-AUC*-eq18}
				\norm{x+z}_X\le (1+\delta)\norm{x}_X\ \ {\rm and}\ \ \norm{x+z}_X\ge \norm{x}_X\ge 2^k.
			\end{equation}
			Note that the second inequality follows from the fact that $X_0$ may be chosen to be included in the kernel of norming functional for $x$. 
			
			From \eqref{eq:CL-rigidity-AUC*-eq15}, \eqref{eq:CL-rigidity-AUC*-eq6} and Theorem \ref{thm:Gorelik-CL}, applied with $\lambda=\frac12$ and $\alpha=t^{-1}24D\delta\norm{x}_X$, we infer the existence of a compact subset $K$ of $Y$ such that
			\begin{equation*}
				\label{eq:CL-rigidity-AUC*-eq19}
				\frac{12D\delta\norm{x}_X}{t}B_Y\subset f\Big(\frac{48D^2\delta\norm{x}_X}{t}B_{X_0}\Big) + CB_Y + K.
			\end{equation*}
			
			As in the proof in the bi-Lipschitz case, fix $\eta'>0$, pick a finite $\eta'$-net $F$ of $K$ and let $E$ be the linear span of $F\cup\{f(x)\}$. Let now $z^*\in E^\perp$ such that $\frac{t}{2} \le \norm{z^*}_{Y^*}\le tD$. Then, there exists $z_0\in B_Y$ such that 
			\begin{equation*}
				\label{eq:CL-rigidity-AUC*-eq20}
				(1-\eta)6D\delta\norm{x}_X \le  \Big\langle z^*, \frac{12D\delta\norm{x}_X}{t}z_0\Big\rangle = \langle z^*, f(x_0)\rangle + \langle z^*,y_0\rangle + \langle z^*, k_0\rangle ,
			\end{equation*} with $x_0 \in \frac{48D^2\delta\norm{x}_X}{t}B_{X_0}$, $\norm{y_0}_Y\le C$ and $k_0\in K$. Therefore,
			
			\begin{equation}
				\label{eq:CL-rigidity-AUC*-eq21}
				\langle z^*,f(x_0)\rangle \ge (1-\eta)6D\delta \norm{x}_X -(CD+1),
			\end{equation}
			if $\eta'$ was initially chosen small enough.\\
			Now,
			\begin{align*}
				\langle y^* &+ z^*  , f(x_0) - f(x')\rangle = \langle y^*, f(x_0) + f(x) \rangle + \langle z^*,f(x_0)\rangle\\
				& \stackrel{\eqref{eq:CL-rigidity-AUC*-eq21}}{\ge} \langle y^*, f(x_0) - f(x) \rangle + 2\langle y^*, f(x)\rangle + (1-\eta)6D\delta \norm{x}_X -(CD+1) \\
				& \stackrel{\eqref{eq:CL-rigidity-AUC*-eq17}}{\ge}  2(1-\eta)\abs{y^*}_{k+1}\norm{x}_X + (1-\eta)6D\delta \norm{x}_X -(CD+1)\\ 
				&-\langle y^*, \frac{f(x) - f(x_0)}{\norm{x -x_0}}\rangle\norm{x-x_0}\\
				& \stackrel{\eqref{eq:CL-rigidity-AUC*-eq18}}{\ge} 2(1-\eta)\abs{y^*}_{k+1}\norm{x}_X  + (1-\eta)6D\delta \norm{x}_X -(CD+1) - \abs{y^*}_{k}(1+\delta)\norm{x}_X,
			\end{align*}
			It follows again from \eqref{eq:CL-rigidity-AUC*-eq18} that 
			\begin{equation*}
				\label{eq:CL-rigidity-AUC*-eq22}
				\abs{y^*+z^*}_k\ge \frac{1}{1+\delta}\Big(2\abs{y^*}_{k+1} - (1+\delta)\abs{y^*}_k  - 2\eta\abs{y^*}_{k+1}  + (1-\eta)6D\delta -
				2^{-k_1}(CD+1)\Big).
			\end{equation*}
			So, it follows from our initial choice of $k_1$ and after letting $\eta \to 0$ that
			\begin{align*}
				\label{eq:CL-rigidity-AUC*-eq23}
				\abs{y^*+z^*}_k &\ge \frac{2}{1+\delta}\abs{y^*}_{k+1} - \abs{y^*}_k + \frac{\delta 6D}{1+\delta}-\frac{\eps}{2} \\
				&= 2\abs{y^*}_{k+1} - \abs{y^*}_k + \frac{\delta}{1+\delta} (6D -2\abs{y^*}_{k+1})-\frac{\eps}{2} \ge  2\abs{y^*}_{k+1} - \abs{y^*}_k + \delta -\frac{\vep}{2},
			\end{align*}
			where in the last inequality we use the fact that $\abs{y^*}_{k+1}\le 2D$ and $\delta\le 1$.
			
			Letting  $\delta$ tend to $\bar{\delta}^*_{\norm{\cdot}_X}(\frac{t}{48D^2})$, we obtain
			\begin{equation*}
				\label{eq:CL-rigidity-AUC*-eq24}
				\abs{y^*+z^*}_k\ge 2\abs{y^*}_{k+1}-\abs{y^*}_k+\bar{\delta}^*_{\norm{\cdot}_X}\big(\frac{t}{48D^2}\big)-\frac{\eps}{2}.
			\end{equation*}
		\end{proof} 
		This finishes the proof of the lemma and thus of the theorem.
	\end{proof}

	\begin{coro}
		\label{cor:CL-rigidity-AUS} 
		Let $X$ and $Y$ be two Banach spaces that are coarse-Lipschitz equivalent. There exists a constant $M>1$ such that for any $\eps\in (0,1)$, there exists an equivalent norm $\abs{\cdot}_\vep$ on $Y$ such that $\frac{1}{1+\eps}\norm{\cdot} _Y\le \abs{\cdot}_\vep\le D\norm{\cdot}_Y$ and for all $t\in [0,1]$,
		\begin{equation}
			\label{eq:CL-rigidity-AUS}
			\bar{\rho}_{\abs{\cdot}_\vep}\big(\frac{t}{M}\big)\le
			\bar{\rho}_{\norm{\cdot}_X}(t)+\eps.
		\end{equation}
	\end{coro}
	
	\begin{proof} 
		%Let $\varphi,\psi$ be continuous, monotone and non-decreasing functions on $[0,1]$ with $\varphi(0)=\psi(0)=0$. If there exist $\alpha\ge 1$ and $\eps >0$ such that for all $t\in [0,1]$, $\varphi(t)\ge \psi(t/\alpha)-\eps$, then it is clear that for all $t\in [0,1]$, $\varphi^*(t/\alpha)\le \psi^*(t)+\eps$. 
		It follows from Proposition \ref{prop:CL-equivalences}  that the coarse-Lipschitz equivalence can be realized by continuous maps and hence the assumptions of Theorem \ref{thm:CL-rigidity-AUC*} are fulfilled. Then, we can apply the $\AUS$-$\AUC^*$ duality ( Corollary \ref{cor:Young2}) to get that if $\abs{\cdot}_\vep$ is the norm given by Theorem \ref{thm:CL-rigidity-AUC*}, then for all $t\in [0,1]$:
		\begin{equation*}
			\bar{\rho}_{\abs{\cdot}_\vep}\big(\frac{t}{576D^2}\big)\stackrel{\eqref{eq:Young1/2}}{\le}
			(\bar{\delta}^*_{\abs{\cdot}_\vep})^*\big(\frac{t}{288D^2}\big)\stackrel{\eqref{eq:CL-rigidity-AUC*-eq1}}{\le}
			(\bar{\delta}^*_{\norm{\cdot}_X})^*\big(\frac{t}{6}\big)+\eps \stackrel{\eqref{eq:Young1/6}}{\le}
			\bar{\rho}_{\norm{\cdot}_X}(t)+\eps.
		\end{equation*}
	\end{proof}
	
	\begin{rema} A close look at the main proof of this section as well as at the proof of Proposition \ref{prop:CL-equivalences} shows that the constant $M$ is still at most $576D^2$ where $D>1$ is such that there are maps $f\colon X\to Y$ and $g\colon Y\to X$ with $\Lip_\infty(f)\le 1$, $\Lip_\infty(g)< D$ and a constant $C\ge 0$ such that $\sup_{x\in X} \norm{(g\circ f)(x)-x}\le C$ and $\sup_{y\in Y}\norm{(f\circ g)(y)-y}\le C$. This without assuming the continuity of $f$ and $g$. 
	\end{rema}
	
	It would be very tempting to take the limit when $k\to \infty $ of $|\cdot|_k$. But passing to the limit in Lemma \ref{lem:CL-rigidity-AUC*} is not allowed, because of the critical dependence on $k$ of the finite-dimensional subspace $E$ of $Y$. This is confirmed by the counterexample of Kalton mentioned in the introduction of this section.
	Quantitative results such as Corollary \ref{cor:CL-rigidity-AUS} will only lead to the stability of classes of asymptotic renormings that do not rely on very small values of the parameter of their moduli. Two renorming results from Chapter \ref{chapter:Szlenk}, namely Theorem \ref{thm:A-theorem-full} for the class $\sA_p$ and Theorem \ref{thm:N-theorem-full} for the class $\sN_p$, are now expressing their subtle flavor and versatility. Recall that for $p\in(1,\infty)$, we have the chain of strict inclusions $\sT_p\subsetneq \sA_p\subsetneq \sN_p\subsetneq \bigcap_{1<r<p}\sT_r$ and  
	$\sT_\infty\subsetneq \sA_\infty= \sN_\infty\subsetneq \bigcap_{1<r<\infty}\sT_r$ and following Causey we will set $\sP_p:=\bigcap_{1<r<p}\sT_r$. 
	
	\begin{coro}
		\label{cor:CL-rigidity-classes} 
		For every $p\in (1,\infty]$, the classes $\sA_p$, $\sN_p$ and $\sP_p$ are stable under coarse-Lipschitz equivalences.
	\end{coro}
	
	\begin{proof}
		For the class $\sA_p$, resp. $\sN_p$, this is a direct consequence of Corollary \ref{cor:CL-rigidity-AUS} and Theorem \ref{thm:A-theorem-full}, resp. Theorem \ref{thm:N-theorem-full}. 
		For the class $\sP_p$ it follows from the simple observation that $\sP_p=\bigcap_{1<r<p} \sA_r$ and the first part of this corollary.
	\end{proof}
	
	Corollary \ref{cor:CL-rigidity-classes} appeared in \cite{CauseyFovelleLancien2023} but the case $p=\infty$ and $\sP_p$ with $p\in(1,\infty)$ was originally from \cite{GKL2001}. 
	
	Since we know from Chapter \ref{chapter:Szlenk} that $\langle \AUS \rangle = \bigcup_{1<p\le \infty} \sT_p = \bigcup_{1<p\le \infty}\sN_p$ and hence a union of coarse-Lipschitz stable classes, we obtain the main result of this section.
	\begin{coro}[Coarse-Lipschitz rigidity of $\langle \AUS \rangle = \mathbf{Sz}_\omega$]
		\label{cor:CL-rigidity-<AUS>}  
		The class of Banach spaces admitting an equivalent asymptotically uniformly smooth norm, i.e., the class $\langle \AUS \rangle = \mathbf{Sz}_\omega$, is stable under coarse-Lipschitz equivalences. 
	\end{coro}
	
	Recall that for  $1\le q< \infty$, we denote by $\mathbf{Sz_q}$ the class of Banach spaces $X$ with Szlenk index of power-type $q$, i.e. there exists $C>0$ such that for all $\eps\in (0,1)$, $\Sz(X,\eps)\le  C\eps^{-q}$.
	Also  
	\begin{equation*}
		\bar{q}_X:=\inf\{q\in [1,\infty) \colon X\in \mathbf{Sz_q}\}.
	\end{equation*} 
	Corollary \ref{coro:szlenk-power-type} tells us that the parameter $\bar{q}_X$ is finite for every space Banach space $X\in \mathbf{Sz_\omega}$. Another nice consequence of Corollary \ref{cor:CL-rigidity-classes} is that this parameter is a coarse-Lipschitz invariant.
	
	\begin{coro}[Coarse-Lipschitz invariance of the infimal Szlenk power-type]
		Assume that $X$ and $Y$ are coarse-Lipschitz equivalent Banach spaces in $\mathbf{Sz}_\omega$. Then, $\bar{q}_X=\bar{q}_Y$.    
	\end{coro}
	
	\begin{proof}
		By symmetry, it is sufficient to show that $\bar{q}_X\le \bar{q}_Y$. So, assume that $Y\in \mathbf{Sz_q}$. Then, it follows from Remark \ref{rem:Szlenk/AUS-power-type} that $Y\in \sP_p$ where $p$ is conjugate to $q$. By Corollary \ref{cor:CL-rigidity-classes} one has that $X\in \sP_p$ which means that $X\in \cap_{1<s<p} \sT_s$ and hence $X\in \mathbf{Sz_r}$ for all $r>q$. Therefore, $\bar{q}_X\le q$ and it follows that $\bar{q}_X\le \bar{q}_Y$. 
	\end{proof}
	
	We conclude this section by observing that while the class $\langle p$-$\AUS\rangle$ is not stable under coarse-Lipschitz equivalence we nevertheless have the following consequence of Corollary \ref{cor:CL-rigidity-classes} and the inclusion $\sT_p\subset \sN_p \subset \cap_{1<s<p} \sT_s$.
	
	\begin{coro}\label{cor:ausstableCL}
		Let $p\in (1,\infty]$ and assume that $X$ and $Y$ are coarse-Lipschitz equivalent Banach spaces. If $X$ is $p$-asymptotically uniformly smooth, then for any $s\in (1,p)$, $Y \in \langle s$-$\AUS\rangle$.
		%In other words, if $X \in \textsf{\emph{T}}_p$, then $Y\in \bigcap_{1<s<p}\textsf{\emph{T}}_s=\textsf{\emph{P}}_p$.
	\end{coro}
	
	Let us mention that this last statement, from which the stability of $\sP_p$ also follows, was deduced in \cite{GKL2001} from Theorem \ref{thm:CL-rigidity-AUC*}, by considering the following dual equivalent norm:
	$$|y|^*=\sum_{k=1}^\infty \frac{1}{k^2}|y|^*_{2^{-k}}.$$

	\section{More on coarse-Lipschitz rigidity problems}
	\label{sec:CL-rigidity-classes}
	%$\textsf{T}_\infty$ and $\textsf{A}_\infty$
	
	It is now time to pay our debts. First, we provide details for the coarse-Lipschitz rigidity of $\ell_p$ for $p\in(1,\infty)$. The uniform rigidity of $\ell_p$ is due to Johnson, Lindenstrauss and Schechtman \cite{JLS1996}. The coarse-Lipschitz rigidity (or under net equivalences) of $\ell_p$ can be found in \cite{BenyaminiLindenstrauss2000} (Theorem 10.21), where an ad hoc Gorelik principle is devised for a net equivalence from $\ell_p$ to another Banach space. The proof we give is based on our more general version of the Gorelik principle.

	\begin{theo}[Coarse-Lipschitz rigidity of $\ell_p$]
		\label{thm:CL-rigidity-ell_p}
		Let $p\in (1,\infty)$. Any Banach space that is coarse-Lipschitzly equivalent to $\ell_p$ must be linearly isomorphic to $\ell_p$.
	\end{theo}
	
	\begin{proof} 
		The case $p\in(1,2)$ is the easiest. Indeed, it follows from Corollary \ref{cor:CLembAUC} that $\ell_2$ does not coarse-Lipschitz embed into $\ell_p$, and we conclude with Proposition \ref{prop:rigidity-ell_p}, which is valid under a coarse-Lipschitz assumption. 
		If $p\in(2,\infty)$ and $2<s<p$, Corollary \ref{cor:ausstableCL} insures that $X$ admits an equivalent norm which is $s$-AUS. In particular, $X$ cannot have a subspace isomorphic to $\ell_2$. Indeed, it follows that $X\in \sT_s$ but clearly $\ell_2\notin \sT_s$. 
		
		%Indeed, it follows from Propositions \ref{prop:Sz-omega} and \ref{prop:Sz-lp} that there exists $A>0$ such that for all $\eps \in (0,1)$, Sz$(X,\eps)\le A\eps^{-q}$, where $q<2$ is the conjugate exponent of $s$, while Sz$(\ell_2,\eps)\ge \eps^{-2}$. This would be in contradiction with Proposition \ref{prop:Sz-subspace}. 
	\end{proof}
	
	Next, we prove a partial result for the coarse-Lipschitz rigidity problem of $\co$.
	\begin{theo}
		\label{thm:CL-rigidity-c_0}
		Assume that $X$ is a Banach space that is coarse-Lipschitz equivalent to $\co$. Then, its dual $X^*$ is linearly isomorphic to $\ell_1$.
	\end{theo}
	
	\begin{proof} 
		We have already observed that under this assumption, $X$ is a $\cL_\infty$-space. It follows from Theorem \ref{cor:ausstableCL} and Theorem \ref{thm:Szlenk-Asplund} that $X^*$ is separable. To conclude, we use a result of Lewis and Stegall \cite{LewisStegall1973} which states that the dual of every $\cL_\infty$-space with a separable dual is isomorphic to $\ell_1$.
	\end{proof}

	\begin{rema} In relation with Problem \ref{prob:uniform-rigidity-c_0}, we recall that $\textsf{A}_\infty$ (the class of asymptotic $\co$ spaces) is stable under coarse-Lipschitz equivalences. So, at this point, it is only known that a Banach space coarse-Lipschitz equivalent to $\co$ is $\textsf{A}_\infty$ and $\mathcal L_\infty$. Another hope was to show that a separable Banach space which is $\textsf{A}_\infty$ and $\mathcal L_\infty$ is necessarily $\textsf{T}_\infty$ (see conjecture after Theorem 5.6 in \cite{GKL2001}). Let us mention here that this question has been solved negatively by Argyros, Gasparis and Motakis in \cite{AGM2016}, who showed the existence of a separable Banach space $X$ which is $\textsf{A}_\infty$ and $\mathcal L_\infty$ but so that every infinite-dimensional subspace of $X$ contains an infinite-dimensional reflexive subspace.
	\end{rema}

	%\noindent {\bf Remark.} Let us point out that we have tried to keep the constant close to optimal in the estimation of the moduli of weak$^*$ asymptotic uniform convexity. However, for the estimates on the moduli of asymptotic uniform smoothness, we have brutally used the duality inequalities and certainly lost much on the optimality.

	\begin{rema} We recall that the original Tsirelson space \cite{Tsirelson1974}, denoted by $\Tsi^*$, is an example of a separable reflexive space in $\textsf{A}_\infty$. It is worth mentioning that the class of reflexive $\textsf{A}_\infty$ spaces satisfies an even stronger rigidity property, as it is stable under coarse or uniform embeddings \cite{BLMS2020}. In this last result, which will be proved in Chapter \ref{chapter:Hamming}, the role of reflexivity is crucial. Indeed, remember that Aharoni's theorem (Theorem \ref{thm:Aharoni}) states that any separable metric space Lipschitz embeds into $\co$.   
	\end{rema}
	
	\section{Notes}
	
	It is worth pointing out that Bessaga \cite{Bessaga1965} showed that if a locally convex space $X$ is uniformly homeomorphic to a Banach space, then $X$ is itself a Banach space, while Enflo \cite{Enflo1970} obtained the same conclusion for locally bounded spaces uniformly homeomorphic to a Banach space with nontrivial roundness.  
	One of the main results of this chapter is that the class $\mathbf{Sz}_\omega$ is stable under coarse-Lipschitz equivalences. This result cannot be extended to higher ordinals. Indeed, the class of Asplund spaces is not stable under uniform homeomorphisms. Ribe proved in \cite{Ribe1984} that if $(p_n)_{n=1}^\infty \subset (1,\infty)$ is a  strictly decreasing sequence tending to $1$, then the space $X=(\sum_{n=1}^\infty \ell_{p_n})_{\ell_2}$ is uniformly homeomorphic to $X\oplus \ell_1$. Note that $X$ is reflexive with $\Sz(X)= \omega^2$ (see Exercise \ref{ex:Szlenk of sums}), while $X\oplus \ell_1$ is not Asplund. This example will be detailed in Chapter \ref{chapter:Counterexamples}. However, we have seen (Theorem \ref{thm:Asplund-rigidity}) that being Asplund is stable under Lipschitz isomorphisms and even under Lipschitz quotients. Using tools from descriptive theory, Y. Dutrieux gave in \cite{Dutrieux2001} the following ``uniform'' version of this result: there exists a universal function $\psi\colon (0,\omega_1) \to (0,\omega_1)$ such that $\Sz(Y)\le \psi(\Sz(X))$ whenever $X$ is a separable Asplund space and $Y$ is a Lipschitz quotient of $X$. In particular, there exists a universal function $\phi\colon (0,\omega_1) \to (0,\omega_1)$ such that $\Sz(Y)\le \phi(\Sz(X))$ whenever $Y$ is Lipschitz equivalent to $X$. We know that $\phi$ can be found with $\phi(\omega)=\omega$. Whether the Szlenk index is invariant under Lipschitz equivalences is still open.
	
	\begin{prob}[Lipschitz invariance of the Szlenk index]
		\label{prob:Lip-invariance-Szlenk}
		Assume that $X$ and $Y$ are two separable Banach spaces that are Lipschitz equivalent. Does this imply that $\Sz(Y)=\Sz(X)$?
	\end{prob} 
	
	It has been shown by C. Samuel \cite{Samuel1984} that the Szlenk index determines the isomorphism class of a $C(K)$ space, for $K$ a metrizable compact space. Thus, a positive answer to the previous question would imply a positive answer to the following important problem.  
	
	\begin{prob}[Lipschitz rigidity of $C(K)$-spaces]
		\label{prob:Lip-rigidity-C(K)}
		Let $K$ and $L$ be two metrizable compact spaces and assume that $C(K)$ is Lipschitz equivalent to $C(L)$. Does this imply that $C(K)$ is linearly isomorphic to $C(L)$?
	\end{prob}
	
	It follows from a work by S. Grivaux \cite{Grivaux2003} that two normed spaces of countable algebraic dimension are linearly isomorphic if and only if their completions are linearly isomorphic (we also refer the reader to Corollary X.1.20 in \cite{Godefroy-Baire}). The following statement then easily follows from Theorem \ref{thm:Lip-rigidity-c_0}.
	
	\begin{coro}
		Assume that $X$ is a normed space of countable algebraic dimension and that $X$ is Lipschitz isomorphic to $(c_{00},\norm{\cdot}_\infty)$. Then, $X$ is linearly isomorphic to $(c_{00},\norm{\cdot}_\infty)$.
	\end{coro}

	\section{Exercises} 
	
	\begin{exer}
		\label{exe:Gorelik}
		\textbf{The original Gorelik principle.} Let $X$ and $Y$ be infinite-dimensional Banach spaces and $f\colon X \to Y$ be a homeomorphism with uniformly continuous inverse $f^{-1}$. Suppose that $\alpha,\beta>0$ are such that there exist a finite-codimensional subspace $X_0$ of $X$ and an infinite-codimensional subspace $Y_0$ of $Y$ for which 
		\begin{equation*}
			f(\alpha B_{X_0}) \subseteq Y_0 + \beta B_Y.
		\end{equation*}
		Show that $\omega_{f^{-1}}(2\beta)\ge \frac{\alpha}{4}$.
	\end{exer}

	\begin{exer}
		\label{ex:Gorelikcodim1} 
		\textbf{Gorelik principle, a simpler proof for one-codimensional subspaces.}
		Assume that $X_0$ is a one-codimensional subspace of the Banach space $X$. Explain how the proof of the Gorelik principle can be simplified in this case by using a bounded linear selector (instead of a Bartle-Graves selector) and the Intermediate Value Principle. 
	\end{exer}
	
	\begin{exer}
		\label{ex:Gorelikforw*null}
		\textbf{Sequential variants of the Gorelik principle.}
		Let $X$ and $Y$ be two Banach spaces. Assume that $f\colon X\to Y$ is a homeomorphism such that $\Lip(f^{-1})<D$, $s>0$ and $X_0$  is a finite-codimensional subspace of $X$. 
		\begin{enumerate}
			\item Show that  for any weak$^*$ null sequence $(y_k^*)_k$ in $Y^*$ such that $\norm{y_k^*}\ge \eps >0$ for all $k$, there exists a sequence $(x_k)_k$ in $2Ds\,B_{X_0}$ such that
			\begin{equation*}
				\liminf_{k\to \infty}\langle f(x_k),y^*_k\rangle \ge s\eps.
			\end{equation*}
			\item Assume moreover that $X^*$ is separable. Show that  for any weak$^*$ null sequence $(y_k^*)_k$ in $Y^*$ such that $\norm{y_k^*}\ge \eps >0$ for all $k$, there exists a subsequence $(y^*_{k_i})_i$ of $(y^*_k)_k$ and a weakly null sequence $(x_i)_i$ in $2Ds\, B_X$ such that for all $i \in \bN$,
			\begin{equation*}
				\liminf_{i\to \infty}\langle f(x_i),y^*_{k_i}\rangle \ge s\eps.
			\end{equation*}
		\end{enumerate} 
	\end{exer}
	
	\begin{exer}
		\label{ex:Gorelikinc_0}
		\textbf{A simpler proof of the Gorelik principle for Lipschitz isomorphisms with $\co$.}
		Let $(e_n)_{n=1}^\infty$ be the canonical basis of $\co$ and, for $n\in \bN$, $F_n$ be the linear span of $\{e_1,\ldots,e_n\}$, $X_n$ the closed linear span of $\{e_k,\ k>n\}$ and $\Pi_n$ the canonical projection from $X=c_0$ onto $F_n$. 
		\begin{enumerate}
			\item Let $\phi:B_{F_n}\to c_0$ be a continuous map such that $\norm{\phi(x)-x}\le 1$ for all $x\in B_{F_n}$. Show that there exists $x\in B_{F_n}$ such that $\phi(x)\in B_{X_n}$. Hint: apply Brouwer's theorem to the map $F\colon x\mapsto x-\Pi_n(\phi(x))$.
			\item Assume that $Y$ is a Banach space and $f\colon c_0 \to Y$ is a Lipschitz isomorphism with $\Lip(f^{-1})\le D$. Show that for any $n\in \bN$, there exists a compact subset $K$ of $Y$ such that
			\begin{equation*}
				\frac{1}{D}B_Y \subset f(B_{X_n}) + K.
			\end{equation*}
			\item Show that this version of the Gorelik principle is enough to prove that $Y$ admits an equivalent asymptotically uniformly flat norm.  
		\end{enumerate}  
	\end{exer}

	\begin{exer}
		\label{ex:Gorelikandfreespaces} 
		Let $X$ and $Y$ be two Banach spaces and $f\colon X \to Y$ be a Lipschitz map. Let $\bar{f}\colon \cF(X) \to Y$ be the bounded linear operator such that $\bar{f} \circ \delta_X= f$. 
		\begin{enumerate}
			\item Assume that $f$ is a Lipschitz isomorphism. Show that the equivalent norm built in the proof of Theorem \ref{thm:Lip-rigidity-AUS} is given by 
			\begin{equation*}
				y^*\in Y^* \mapsto \abs{y^*}_{Y*}=\norm{(\bar{f})^*(y^*)}_{\Lip_0(X)}.
			\end{equation*} 
			\item Show that $\bar{f}:\cF(X)\to Y$ is a linear quotient map and that $|\ |_Y$ is the corresponding quotient norm.
			\item Assume that $f$ is a Lipschitz quotient. Show that the above formula still defines an equivalent norm on $Y^*$. 
		\end{enumerate}
	\end{exer}
	
	\begin{exer}
		\label{ex:script}
		Recall that for $p\in [1,\infty]$, a Banach space $X$ is a $\cL_p$-space if there exists $\lambda \ge 1$ such that for any finite-dimensional subspace $F$ of $X$ there exists a finite-dimensional subspace $G$ of $X$ satisfying $F \subset G$ and $d_{BM}(G,\ell_p^{\dim(G)})\le \lambda$.
		\begin{enumerate}
			\item Show that the spaces $\co$ and $C(K)$ for $K$ compact are $\cL_\infty$-spaces
			\item Let $p\in [1,\infty)$. Show that $\ell_p$ and $L_p$ are $\cL_p$-spaces.
		\end{enumerate}   
	\end{exer}

	%\begin{exer}
	%Tony's direct proof of $X^*$ separable and $X^*$ weak$^*$-fragmentable??    
	%\end{exer}

	%%%%%%%%%%%%%%%%%%%%%%%%%%%%%%%%%%%%%%%%%%%%%%%%%%%%%%%%%%%%%%%%%%%%%%%%%%%%
	
	\chapter[Geometry of the Johnson graphs and applications]{Geometry of the Johnson graphs and applications}
	\label{chapter:Johnson}
	
	In Chapter \ref{chapter:Gorelik} we have seen how the Gorelik principle could be used to prove the stability of asymptotic uniform smoothness properties under Lipschitz or coarse-Lipschitz \emph{equivalences} between Banach spaces. It is clear from its statement and its proof that the Gorelik principle is not suited for addressing questions on the stability of asymptotic uniform smoothness under Lipschitz or coarse-Lipschitz \emph{embeddings}. An incredibly powerful and elegant tool to handle nonlinear embedding problems was introduced by Kalton in \cite{Kalton2007} and further developed by Kalton and Randrianarivony in \cite{KaltonRandrianarivony2008}.
	The key idea is to study concentration properties for Lipschitz maps defined on infinite graphs with countable degree, which are generated by certain bases in classical Banach spaces. These graphs pick up the nonlinear and asymptotic geometry of the Banach spaces with the relevant bases. The goal is to prove a scale of concentration inequalities using geometric features of the target space and to show that these inequalities are stable under a corresponding scale of nonlinear embeddings. Intuitively, the strongest possible concentration inequality will serve as an obstruction to the weakest notion of nonlinear embeddings.
	In Section \ref{sec:Johnson}, we introduce the Johnson graphs and explain their connection with symmetric difference metrics. In Section \ref{sec:J-concentration} we describe and prove several concentration inequalities for Lipschitz maps on the Johnson graphs which take values in reflexive asymptotically uniformly smooth spaces. We then discuss the stability of these concentration inequalities under nonlinear embeddings. In Section \ref{sec:J-concentration-consequences} we detail a few immediate consequences of the concentration inequalities: the coarse-Lipschitz rigidity of the class of subspaces of $\ell_p$ for $p\in(2,\infty)$, a coarse version of Tsirelson's theorem, the optimality of snowflake and compression exponents between $\ell_p$ spaces (see the notes in Section \ref{sec:Notes-Johnson}). The coarse-Lipschitz rigidity of finite sums of $\ell_p$-spaces is dealt with in Section \ref{sec:CL-rigidity-sums-of-lp} while in Section \ref{sec:CL-rigidity-l2+lp}, we discuss the problem of the coarse-Lipschitz rigidity of $\ell_2\oplus \ell_p$.
	In this chapter, we also shine a light on several nonlinear embedding problems, some of which are among the most important open problems in the nonlinear geometry of Banach spaces.

	\section[Johnson graphs: basic properties and embeddings]{Johnson graphs: basic properties and embeddings into Banach spaces}
	\label{sec:Johnson}
	
	If we denote by $\bar{e}:=\ei$ the canonical basis of $\ell_1$, one can induce, for every $k\in \bN$ and every infinite subset $\bM$ of $\bN$, a metric on $[\bM]^{k}:=\{ S\subset \bM \colon |S|=k\}$, denoted by $\sd_{\ebar}$ for simplicity, as follows:
	\begin{equation}
		\label{eq:sym-diff}
		\sd_{\ebar}(S_1,S_2):=\Big\|{\sum_{i\in S_1} e_i -\sum_{i\in S_2} e_i}_1\Big\|.
	\end{equation}
	It is immediate that $\sd_{\ebar}(S_1,S_2) = \abs{S_1\triangle S_2}$
	where $S_1\triangle S_2:=(S_1\setminus S_2) \cup (S_2\setminus S_1)$ is the symmetric difference between $S_1$ and $S_2$ and hence $\sd_{\ebar}$ coincides with the symmetric difference metric $\dsym$. Therefore, another interpretation of $\eqref{eq:sym-diff}$ is that the embedding $S\in [\bM]^k \mapsto \sum_{i\in S} e_i$ is an isometric embedding of $([\bM]^k,\dsym)$ into $\ell_1$.
	The symmetric difference metric can be realized (up to a dilation) as a graph distance. Indeed, if we consider $[\bM]^k$ as the vertex set of a graph $\sJ_k(\bM)$ where two distinct vertices are adjacent if and only if their intersection has $k-1$ elements, then one can show that the graph metric on $\sJ_k(\bM)$, denoted by $\dJk$, satisfies $\dJk=\frac12 \dsym$. It is also easy to see that the Johnson graphs on two infinite subsets of $\bN$ are graph-isomorphic (see Exercise \ref{ex:J-graph}) and hence isometric as metric spaces. Therefore, we will simply use the notation $\sJ_k$ for the $k$-th Johnson graph if it does not lead to an ambiguity.  Unless otherwise specify, we denote a vertex in $\sJ_k$ by $\mbar:=\{m_1,m_2, \dots, m_k\}$ where $m_1<m_2<\dots<m_k$ are all elements in the underlying set.
	
	As we already explained, $(\sJ_k,\dJk)$  admits a (scaled) isometric embedding into $\ell_1$. If we consider $\ei$ as the canonical basis of $\ell_p$, then it is clear that the natural embedding $\mbar\in \sJ_k \stackrel{h_p}{\mapsto} \sum_{i=1}^k e_{m_i}\in \ell_p$ satisfies $\norm{h_p(\mbar) - h_p(\nbar) }_p = 2^{\frac{1}{p}}\dJk(\mbar,\nbar)^{\frac{1}{p}}$. This means that the $\frac1p$-snowflake of $(\sJ_k,\dJk)$ admits a (scaled) isometric embedding into $\ell_p$. As we will see this is the best snowflake exponent that can be achieved.
	
	The $k$-th Johnson graph $\sJ_k$ is an infinite graph with countable degree and diameter $\diam(\sJ_k)= k$, which prevents good embeddability into finite-dimensional Banach spaces. How well the Johnson graphs can be embedded into infinite-dimensional Banach spaces is governed by the properties of their weakly null sequences. In Proposition \ref{prop:J-graph-embedding} below we show that if $Y$ is an infinite-dimensional Banach space admitting a spreading model $E$, generated by a weakly null sequence, then there always exists a $1$-Lipschitz map from $\sJ_k$ into $Y$ whose compression modulus is essentially bounded from below by the fundamental function of the spreading model. Recall that the \emph{fundamental function of a spreading model $S$ of a Banach space} is given for all $r\in \bN$ by $\varphi_S(r):=\norm{ \sum_{i=1}^r e_i}_S$ where $\ei$ is the fundamental sequence of the spreading model $S$. We refer to Section \ref{sec:spreading-models} of Appendix \ref{appendix:asymptotic}  for the basics on spreading models. 
	
	\begin{prop}
		\label{prop:J-graph-embedding}
		Let $Y$ be an infinite-dimensional Banach space admitting a spreading model $S$, generated by a normalized weakly null sequence, then for every $\nu>0$ and $k\ge 1$ there exists a map  $h_{k,\nu} \colon \sJ_k \to Y$ such that for all $\mbar,\nbar \in \sJ_k$,
		\begin{equation}
			\label{eq:J-graph-embedding-eq1}
			\frac{1}{4(1+\nu)}\varphi_S\left(\dJk(\mbar,\nbar)\right)\le \norm{h_{k,\nu}(\mbar)-h_{k,\nu}(\nbar)}_Y\le \dJk(\mbar,\nbar).
		\end{equation}
	\end{prop}
	
	\begin{proof}
		By Proposition \ref{prop:weakly-null-spreading-model}  there is a weakly null normalized basic sequence $\yn$ with basis constant  not larger than $(1+\nu)$,
		and thus the bimonotonicity constant is at most $2(1+\nu)$, generating a spreading model $S$ such that for all $r\ge1$, for all $r\le n_1< n_2< \dots < n_{r}$ and for all $(\vep_i)_{i=1}^{r}\in\{-1,1\}^r$ one has
		\begin{equation}
			\label{eq:J-graph-embedding-eq2}
			\Big\|\sum_{i=1}^{r} \vep_i y_{n_i}\Big\| \ge \frac{1}{1+\nu}\Big\|\sum_{i=1}^{r} \vep_i e_i\Big\| \ge \frac{1}{2(1+\nu)}\varphi_S(r).
		\end{equation}
		
		Let $h:=h_{k,\nu}\colon \sJ_k\to Y$  be defined by $h(\mbar):= \frac12\ds\sum_{i=1}^k y_{km_i}$ for all $\mbar\in \sJ_k$. The map $h$ is clearly $1$-Lipschitz since after cancellations it follows that for all $\mbar,\nbar\in \sJ_k$
		\begin{equation*}
			\label{eq:J-graph-embedding-eq3}
			\Big\| h(\mbar) - h(\nbar)\Big\|= \frac12 \Big\| \sum_{i=1}^{\dsym(\mbar,\nbar)} \vep_i y_{kq_i}\Big\|,
		\end{equation*}
		for some $q_1<q_2<\dots<q_{\dsym(\mbar,\nbar)}$ and $(\vep_i)_{i=1}^{\dsym(\mbar,\nbar)}\in\{-1,1\}^{\dsym(\mbar,\nbar)}$.
		It follows from \eqref{eq:J-graph-embedding-eq2} that
		\begin{align*}
			\norm{ h(\mbar) - h(\nbar) } & \ge \frac{1}{4(1+\nu)}\varphi_S(\dsym(\mbar,\nbar)) = \frac{1}{4(1+\nu)}\varphi_S(2\dJk(\mbar,\nbar))\\
			& \ge \frac{1}{4(1+\nu)}\varphi_S(\dJk(\mbar,\nbar)).
		\end{align*}
		The last inequality follows from the fact that $(e_i)_i$ is $1$-suppression unconditional and therefore $\varphi_S$ is non-decreasing.     
	\end{proof}
	
	It is clear that the sequence $\{\sJ_k\}_{k=1}^\infty$ of Johnson graphs equi-bi-Lipschitzly embeds into any Banach space that contains an isomorphic copy of $\ell_1$. It follows from Proposition \ref{prop:J-graph-embedding} that any Banach space with a spreading model isomorphic to $\ell_1$ will also equi-bi-Lipschitzly contain the Johnson graphs. This is a significant improvement, as having a spreading model isomorphic to $\ell_1$ is a much weaker requirement than containing $\ell_1$. A concrete example distinguishing these two properties is given by the reflexive asymptotic-$\ell_1$ Tsirelson space $\Tsi$, introduced in Section \ref{sec:asymptotic-c_0}.
	
	Proposition \ref{prop:J-graph-embedding} does not provide a nontrivial compression function if the fundamental function of every spreading model is bounded. This happens when all the spreading models are isomorphic to $\co$. However, Proposition \ref{P:4} in Section \ref{sec:spreading-models} ensures that if a Banach space admits at least one spreading model generated by a weakly null sequence that is not isomorphic to $\co$, then the associated fundamental function is unbounded. In this case, Proposition \ref{prop:J-graph-embedding} yields an equi-coarse embedding of the sequence of Johnson graphs. 
	%The following corollary makes this statement more precise. 
	
	\begin{coro}
		\label{cor:J-graph-embeddings}
		If a Banach space $Y$ has a spreading model $S$, generated by a normalized weakly null sequence that is not isomorphic to $\co$, then 
		the sequence of Johnson graphs $\{\sJ_k,\dJk\}_{k\in \bN}$ equi-coarsely embeds into $X$.
		%for any $\nu>0$ and $k\ge 1$ there is a coarse embedding $h_{\nu,k}\colon \sJ_k \to Y$  such that $h_{\nu,k}$ is $1$-Lipschitz and the modulus of compression of $h_{\nu,k}$ satisfies $\rho_{h_{\nu,k}}(t)\ge \frac{1}{8(1+\nu)}\varphi_S(t)$, for $t>0$, where $\varphi_S$ is the fundamental function of the spreading model $S$.
	\end{coro}
	
	%%%%%%%%%%%%%%%%%%%%%%%%%%%%%%%%%%%%%%%%%%%%%%%%%%%%%%%%%%%%%%%%%%%%%%%%%%%%%%%%%%%%%%

	\section[Concentration inequalities on Johnson graphs]{Concentration inequalities for Lipschitz maps on Johnson graphs} 
	\label{sec:J-concentration}
	
	The Johnson graphs are $\ell_1$-objects in nature. They capture to some extent the extreme asymptotically uniformly convex behavior of $\ell_1$ and can be used to provide embedding obstructions into Banach spaces that are asymptotically uniformly smooth. Of course, it follows from Aharoni's embedding theorem that asymptotic uniform smoothness (or even flatness) alone is never an obstruction to nonlinear embedding, even as good as bi-Lipschitz ones. To obtain nonlinear embedding obstructions, asymptotic uniform smoothness together with some degree of reflexivity are usually required. As we will see, the degree of faithfulness of any embedding of the Johnson graphs into a reflexive Banach space will deteriorate as the target space becomes more and more asymptotically uniformly smooth. This is another occurrence of the paradigm that one cannot embed a space that is too convex into a reflexive space that is too smooth.
	
	A chief reason why reflexivity plays a crucial role is that maps on $[\bN]^k$ and taking values into a reflexive space can be, in some sense, linearized, thereby establishing a connection between the nonlinear world to the linear one. This fundamental principle does not have much to do with the metric $d_k
	$ with which $[\bN]^k$ is equipped, but for the applications, additional properties of the metric $d_k$ will matter. Recall that for an infinite subset $\M$ of $\bN$, a subset $(x_{\mbar})_{\mbar\in[\bM]^{\le k}}$ of a Banach space $X$ is a \emph{weakly null tree} if for all $\mbar \in [\bM]^{\le k-1}$, the sequence $(x_{(\mbar \smallfrown l)})_{l\in \bM}$ is weakly null. 
	We now state and prove the fundamental linearization principle. We need a convenient piece of notation. Given $\mbar \in [\bM]^{\le k}$ of length $l$, $1\le j\le l$ and $m\in \bM$, we define a replacement operator $\sigma_{j,m}$, by $\sigma_{j, m}(\mbar):=(m_1,\dots,m_{j-1}, m_, m_{j+1},\dots, m_{l})$. 
	
	\begin{lemm}
		\label{lem:linearization} 
		Let $X$ be a reflexive Banach space and $d_k$ be a metric on $[\bN]^k$. Fix $k\in\bN$ and $f \colon ([\bN]^k, d_k)\to X$ a map with bounded image. There exist an infinite subset $\bM$ of $\bN$ and a weakly null (rooted) tree $(x_{\mbar})_{\mbar\in[\bM]^{\le k}}$ in $X$ so that:
		\begin{enumerate}[(i)]
			\item for all $\mbar \in [\bM]^k$,  
			\begin{equation}
				\label{eq:linearization-eq1}
				f(\mbar) = x_\emptyset + \sum_{i=1}^k x_{(m_1,\dots, m_i)},
			\end{equation}
			\item for all $\mbar\in [\bM]^{\le  k}\setminus\{\emptyset\}$ of length $i$, $1\le i \le k$:
			\begin{equation}
				\norm{x_{\mbar}} \le  \Lip(f)  \sup\{ d_k(\nbar,\sigma_{i,m}(\nbar)):\ \nbar\in [\bM]^{k}, m\in \bM \text{ s.t. }\sigma_{i,m}(\nbar)\in [\bM]^k\}.
			\end{equation}
		\end{enumerate}
	\end{lemm}
	
	\begin{proof} 
		We prove the first part of the statement by induction on $k\in\bN$, which only uses the weak compactness in $X$. If $k=1$, since $X$ is reflexive, we can take a subsequence   $\big(f(n_i)\big)_{i\in\bN}$ of $\big(f(n)\big)_{n\in\bN}$ which converges weakly to some $x_\emptyset\in X$. Then, take  $\bM=\{n_i \colon i\in \bN\}$ and $x_{(n_i)}=f(n_i)-x_\emptyset$. 
		
		Assume our claim to be true for $k$ with  $k\in\bN$ and let $f\colon [\bN]^{k+1}\to X$ have a bounded image. Since $X$ is reflexive and $[\bN]^k$ is countable, we can  use weak compactness and a diagonal argument to find $\bL \in [\bN]^\omega$ such that for all $\mbar \in [\bN]^k$, there exists $g(\mbar)\in X$ such that 
		\begin{equation*}
			w-\lim_{l\to\infty, l\in \bL} f(\mbar\frown l)=g(\mbar). 
		\end{equation*}
		Finally, for $\mbar \in [\bL]^{k+1}$, we set $x_{\mbar} = f(\mbar) - g(m_1,\dots,m_k)$ and we apply the induction hypothesis to $g\colon [\bL]^{k}\to X$, which provides us with an infinite set $\bM\subset \bL$ and a weakly null tree $(x_{\mbar})_{\mbar\in[\bM]^{\le k}}$ so that $g(\mbar) = x_\emptyset + \sum_{i=1}^k x_{(m_1,m_2,\dots, m_i)}$ for all $\mbar=(m_1,m_2,\dots, m_{k})\in[\bM]^{k}$ and thus for all $\mbar=(m_1,\dots, m_{k+1})\in[\bM]^{k+1}$,
		\begin{equation*}
			f(\mbar) = x_{\mbar}+g(m_1,\dots, m_{k})=x_\emptyset + \sum_{i=1}^{k+1} x_{(m_1,\dots, m_i)}.
		\end{equation*}
		
		To prove the second part of the statement, assume after relabeling that $\bM=\bN$. Fix $i\in \{1,\dots,k\}$ and let $\mbar = (m_1,m_2,\dots, m_i)$ in $[\bN]^{\le k}\setminus \{\emptyset\}$ and write $\mbar_{\restriction i-1}:=(m_1,m_2,\dots, m_{i-1})$ and observe from \eqref{eq:linearization-eq1} we can deduce that  
		\begin{equation}
			x_{\mbar} = w\text{\;-}\lim_{m_i'\to\infty}\lim_{m_{i+1}\to\infty}\!\!\!\!\!\cdots\!\!\! \lim_{m_k\to\infty}\big( f(\mbar,m_{i+1},\dots,m_k) - f(\mbar_{\restriction i-1},m_i',m_{i+1},\dots,m_k)\big)
		\end{equation}
		It now follows from the lower semicontinuity of the norm with respect to the weak topology that we have
		\begin{align*}
			\norm{x_{\mbar}} &\le \lim_{m_i'\to\infty}\lim_{m_{i+1}\to\infty}\!\!\!\!\!\cdots\!\!\! \lim_{m_k\to\infty}\norm{f(\mbar,m_{i+1},\dots,m_k) - f(\mbar_{\restriction i-1},m_i',m_{i+1},\dots,m_k)}\\
			&\le \Lip(f)  \sup_{\nbar\in [\bM]^{k}, m\in \bM\colon \sigma_{i,m}(\nbar)\in [\bM]^k} d_k(\nbar,\sigma_{i,m}(\nbar))
			%\lim_{m_i'\to\infty}\lim_{m_{i+1}\to\infty}\!\!\!\!\!\cdots\!\!\! \lim_{m_k\to\infty} d_{k}\big(\m\frown (n_{i+1},\ldots, n_k),\m'\frown (n_i,n_{i+1},\ldots, n_k)\big)\\
		\end{align*}
		
	\end{proof}
	
	\begin{rema}
		\label{rem:linearization-J}
		When $d_k=\dJk$ we have $$\sup\{ \dJk(\nbar,\sigma_{i,m}(\nbar))\colon i\in \{1,\dots,k\}, \nbar\in [\bM]^{k}, m\in \bM \text{ s.t. }\sigma_{i,m}(\nbar)\in [\bM]^k\}=1,$$ and hence the norms of the nodes of the weakly null tree are controlled by the Lipschitz constant of the map. 
	\end{rema}
	
	Note that for every map $f \colon ([\bN]^k, d_k)\to (X,d_X)$ we have the trivial inequality
	\begin{equation*}
		\diam(f([\bN]^k)):=\sup_{\mbar,\nbar\in[\bN]^k}d_X(f(\mbar),f(\nbar))\le \Lip(f)\diam(([\bN]^k,d_k))
	\end{equation*}
	Informally, we will say that $X$-valued Lipschitz maps on $([\bN]^k, d_k)$ concentrate if this trivial inequality can be improved by eventually passing to an infinite subset of $\bN$.
	
	The next theorem encapsulates the following concentration phenomenon: Lipschitz maps on the Johnson graphs will concentrate if they take their values in reflexive and asymptotically uniformly smooth spaces. Moreover, the smoother the space is, the stronger the concentration.
	
	\begin{theo}
		\label{thm:J-KR-BLMS}
		Let $X$ be a reflexive Banach space.
		\begin{enumerate}[(i)]
			\item  If $X$ has property $\sN_p$ for some $p\in (1,\infty)$, then for all  $k\in\bN$ and all $f\colon ([\bN]^k, \dJk)\to X$, there exists an infinite subset $\bM$ of $\bN$ so that
			\begin{equation}
				\sup_{\mbar,\nbar\in[\bM]^k} \norm{f(\mbar)-f(\nbar))} \le 3\textsf{n}_p(X)k^{1/p}\Lip(f).
			\end{equation}
			\item If $X$ has property $\sN_\infty=\sA_\infty$, then for all  $k\in\bN$ and all $f\colon ([\bN]^k, \dJk)\to X$, there exists an infinite subset $\bM$ of $\bN$ so that
			\begin{equation}
				\sup_{\mbar,\nbar\in[\bM]^k} \norm{f(\mbar)-f(\nbar))} \le 3\textsf{n}_\infty(X)\Lip(f).
			\end{equation}
		\end{enumerate}
	\end{theo}
	
	\begin{proof} Let $k\in \bN$ and $f\colon ([\bN]^k, \dJk)\to X$ and assume, as we may, that $\Lip(f)=1$. Since the closed linear span of the image of $f$ is separable, we may also assume that $X$ is separable and, since it is reflexive, that it has a countable neighborhood basis at $0$. By Lemma \ref{lem:linearization} and Remark \ref{rem:linearization-J} following it, there exist $\bL \in [\bN]^\omega$ and a weakly null tree $(x_{\mbar})_{\mbar\in[\bL]^{\le k}}$ in $B_X$ so that for all $\nbar \in [\bL]^k$, 
		\begin{equation*}
			f(\nbar) = x_\emptyset + \sum_{i=1}^k x_{(n_1,\dots,n_i)} = \sum_{\mbar \preceq \nbar} x_{\mbar}.
		\end{equation*}
		Since $X$ has $\sN_p$, $\textsf{n}_p(X)<\infty$. Thus, it follows from Proposition \ref{prop:Np-game-trees} (and the separability of $X$) that for any $\bM \in [\bL]^\omega$, there exists $\nbar \in [\bM]^k$ such that $\norm{\sum_{i=1}^k x_{(n_1,\dots,n_i)}}\le \frac32\textsf{n}_p(X)k^{1/p}$. We can therefore use Ramsey's theorem (see Theorem \ref{thm:infinite_Ramsey} in Appendix \ref{appendix:Ramsey}) to deduce the existence of $\bM \in [\bL]^\omega$ such that $\norm{\sum_{i=1}^k x_{(n_1,\dots,n_i)}}\le \frac32\textsf{n}_p(X)k^{1/p}$ for all $\nbar \in [\bM]^k$. This immediately implies that $\diam (f([\M]^k))\le 3\textsf{n}_p(X)k^{1/p}$ and assertion $(i)$ follows.
		
		The proof of assertion $(ii)$ is identical.
	\end{proof}
	
	It will be convenient to introduce the following terminology. 
	Given a function $\gamma\colon (0,\infty)\to [1,\infty)$ and a metric space $(X,d_X)$, we say that the $X$-valued Lipschitz maps on $([\bN]^k,d_k)$ satisfy a \emph{$O(\gamma(\diam(d_k)))$-concentration inequality} if there exists a constant $C\ge 1$ such that for all  $k\in\bN$  and all $f\colon ([\bN]^k, d_k)\to X$ Lipschitz, there exists $\bM \in [\bN]^\omega$ so that 
	\begin{equation}
		\sup_{\mbar,\nbar\in[\bM]^k} d_X(f(\mbar),f(\nbar)) \le C\gamma(\diam(d_k))\Lip(f).
	\end{equation}
	
	If $C=1$ above we simply say that we have a $\gamma(\diam(d_k))$-concentration inequality.
	
	Of course, we are interested in concentration functions $\gamma$ such that $\lim_{t\to \infty}\gamma(t)/t=0$ and in metrics $d_k$ such that $\diam([\bM]^k,d_k)\approx \diam([\bN]^k,d_k)$ for every infinite subset $\bM$ in $\bN$. This last requirement on the metrics is clearly true for the Johnson metrics and will always be satisfied in the other examples that we will cover.
	Since every $p$-asymptotically uniformly smooth space is in $\sN_p$ (Theorem \ref{thm:N-theorem-full}) and the class $\sA_\infty=\sN_\infty$ coincides with the class of asymptotic-$\co$ spaces (Theorem \ref{thm:N-infty-theorem-full}) we can rephrase Theorem \ref{thm:J-KR-BLMS} according to the newly introduced terminology. 
	
	\begin{coro}\, 
		\label{cor:J-KR-BLMS}
		\begin{enumerate}
			\item Let $p\in(1,\infty)$ and $X$ be a reflexive and $p$-asymptotically uniformly smooth Banach. Then, the $X$-valued Lipschitz maps on $\sJ_k$ satisfy a $O(k^{\frac1p})$-concentration inequality.
			\item Let $X$ be a reflexive asymptotic-$\co$ Banach space. Then, the $X$-valued Lipschitz maps on $\sJ_k$ satisfy a $O(1)$-concentration inequality.
		\end{enumerate}
	\end{coro}
	
	The first assertion in Corollary \ref{cor:J-KR-BLMS} is originally due to Kalton and Randrianarivony \cite{KaltonRandrianarivony2008}, where concentration inequalities on the Johnson graphs were first studied. Their original proof is essentially self-contained and only relies on a basic induction and Ramsey-type arguments. The second assertion in Corollary \ref{cor:J-KR-BLMS} is from \cite{BLS2018} for the special case of the asymptotic-$\co$ Tsirelson space $\Tsi^*$ and its original proof was inspired by \cite{KaltonRandrianarivony2008}. The use of the well-known fundamental linearization result to reduce the problem of obtaining concentration inequalities from upper estimates on weakly null trees is from \cite{BLMS2020}. In this section, we have decided to follow the latter conceptual approach as it unifies the proofs of the two concentration inequalities, which are now immediate consequences of the tools developed in Chapter \ref{chapter:Szlenk}.
	
	Concentration inequalities, as above, are instrumental in the study of nonlinear embeddings as they provide obstructions to nonlinear embeddings.
	
	\begin{prop}
		\label{prop:J-obstruction}
		Let $(X,d_X)$ be a metric space.
		\begin{enumerate}[(i)]
			\item If the $X$-valued Lipschitz maps on $\sJ_k$ satisfy a $\gamma(k)$-concentration inequality with $\lim_{t\to \infty}\gamma(t)/t=0$, then $(\sJ_k)_{k\in \bN}$ does not admit an equi-coarse-Lipschitz embedding into $X$.
			\item If the $X$-valued Lipschitz maps on $\sJ_k$ satisfy a $O(1)$-concentration inequality, then $(\sJ_k)_{k\in \bN}$ does not admit an equi-coarse embedding into $X$.
		\end{enumerate}
	\end{prop}
	
	\begin{proof}
		The proofs are elementary.
		
		$(i)$ Assume that the $X$-valued Lipschitz maps on $\sJ_k$ satisfy a $\gamma(k)$-concentration inequality and for the sake of a contradiction that we have a sequence of maps $\{f_k\colon \sJ_k\to X\colon k\in \bN\}$ witnessing the equi-coarse-Lipschitz embeddability of $(\sJ_k)_{k\in \bN}$, i.e. there are $(A,B,C,D)\in(0,\infty)^2\times[0,\infty)^2$ so that for all $k\ge 1$, $\rho_{f_k}(t)\ge At-B$ and $\omega_{f_k}(t)\le Ct+D$. Since by assumption and the fact that $\dJk$ is a graph metric, there exists $\bM \in [\bN]^\omega$ such that $\diam (f_k([\bM]^k))\le (C+D) \gamma(k)$, it follows that $A\diam (\sJ_k)-B\le (C+D)\gamma(k)$ and hence $A-\frac{B}{k} \le (C+D)\gamma(k)/k$ for all $k\in \bN$. This is clearly impossible if $\lim_{t\to \infty}\gamma(t)/t=0$.
		
		$(ii)$ Now, let us assume that the $X$-valued Lipschitz maps on $\sJ_k$ satisfy a $\lambda$-concentration inequality, for some $\lambda>0$. Should there exist a sequence of maps $\{f_k\colon \sJ_k\to X\colon k\in \bN\}$ witnessing the equi-coarse embeddability of $(\sJ_k)_{k\in \bN}$ and since $\dJk$ is a graph metric, we would then have $\rho_{f_k}(k)\le \lambda \omega_{f_k}(1)$ leading again to a contradiction whenever for all $k\ge 1$, $\rho\le \rho_{f_k}\le \omega_{f_k}\le \omega$ with $\lim_{t\to\infty}\rho(t)=\infty$ and $\omega(1)<\infty$.
	\end{proof}
	
	%The attentive reader will have noticed that the sequence $(\sJ_k)_{k\in \bN}$ cannot admit an equi-coarse-Lipschitz embedding into a space $X$ for which every $X$-valued Lipschitz map on $\sJ_k$ satisfies a $O(k^{\frac1p})$-concentration.  Similarly, a $O(1)$-concentration provides an obstruction to  equi-coarse embeddability. Based on these observations it is now clear why concentration inequalities are instrumental in the study of nonlinear embeddings. We postpone the elementary details of these claims to the next section and 
	
	We now turn to the permanence properties of the concentration inequalities.
	
	\begin{prop}
		\label{prop:JFC-invariant} Let $X$ and $Y$ be Banach spaces.
		\begin{enumerate}[(a)]
			\item Assume that the $Y$-valued Lipschitz maps on $\sJ_k$ satisfy a $O(\gamma(k))$-concentration inequality. If $X$ coarse-Lipschitz embeds into $Y$, then the $X$-valued Lipschitz maps on $\sJ_k$ satisfy a $O(\gamma(k))$-concentration inequality.
			\item Assume that the $Y$-valued Lipschitz maps on $\sJ_k$ satisfy a $O(1)$-concentration inequality. If $X$ coarsely embeds into $Y$, then the $X$-valued Lipschitz maps on $\sJ_k$ satisfy a $O(1)$-concentration inequality.
		\end{enumerate}
	\end{prop}
	
	\begin{proof} $(a)$ Assume that the $Y$-valued Lipschitz maps on $\sJ_k$ satisfy a $\lambda \gamma(k)$-concentration inequality, for some $\lambda \ge 1$ and that $g\colon X\to Y$ is a coarse-Lipschitz embedding. Let $(A,B,C,D)\in(0,\infty)^2\times[0,\infty)^2$ so that $\rho_g(t)\ge At-B$ and $\omega_g(t)\le Ct+D$. Let $k\in \bN$ and $f\colon ([\bN]^k,\dJk) \to X$ be a  Lipschitz function. Then, for all $\mbar,\nbar \in [\bN]^k$, 
		\begin{equation*}
			\norm{(g\circ f)(\mbar) - (g\circ f)(\nbar)}_{Y}\le C\Lip(f)\dJk(\mbar,\nbar)+D.
		\end{equation*} 
		Since $\dJk$ is a graph metric, we deduce that $\Lip(g\circ f)\le C\Lip(f)+D$. It follows that there exists $\bM \in [\bN]^\omega$ such that $\diam ((g\circ f)([\bM]^k))\le \lambda(C\Lip(f)+D)\gamma(k)$. Therefore, $A\diam (f([\M]^k))-B\le \lambda(C\Lip(f)+D)\gamma(k)$ and
		\begin{equation*}
			\diam (f([\M]^k))\le \frac{\lambda(C\Lip(f)+D)+B}{A}\gamma(k).
		\end{equation*}
		Since the target space is a Banach space, we could have assumed without loss of generality,  by rescaling if needed, that $\Lip(f) = 1$ and hence we have proved that the $X$-valued Lipschitz maps on $\sJ_k$ satisfy a $\mu \gamma(k)$-concentration inequality where 
		$\mu=\frac{\lambda(C+D)+B}{A}$.
		
		\medskip $(b)$ Assume that the $Y$-valued Lipschitz maps on $\sJ_k$ satisfy a $\lambda$-concentration inequality, for some $\lambda \ge 1$ and that $g\colon X\to Y$ is a coarse embedding with compression and expansion moduli $\rho_g,\omega_g\colon [0,\infty)\to[0,\infty)$. We will show that the $X$-valued Lipschitz maps on $\sJ_k$ satisfy a $\mu $-concentration inequality for any $\mu$ such that $\rho(\mu)>\lambda \omega_g(1)$. Note that such a $\mu$ exists, since $\omega_g(1)<\infty$ and $\lim_{t\to \infty} \rho_g(t)=\infty$. Let $f\colon [\bN]^k\to X$ be a Lipschitz map and assume, as we may, that $\Lip(f)=1$. 
		Because $\dJk$ is a graph metric it follows that $\Lip(g\circ f)=\omega_{g \circ f}(1)\le \omega_g(1)$. It follows that there exists $\M \in [\bN]^\omega$ such that $\diam ((g\circ f)([\M]^k))\le \lambda \omega_g(1)$. This implies that $\rho_g(\diam f([\M]^k))\le \lambda \omega_g(1)$, which concludes the proof. 
	\end{proof}
	
	\begin{rema}
		In the proofs of Proposition \ref{prop:JFC-invariant} we do not use that the domain space is a Banach space, and the proofs remain valid if the domain space is a metric space. Also, if the target space is merely a metric space the proofs show that Johnson-concentration inequalities are stable under bi-Lipschitz embeddings.
	\end{rema}

	\section{First consequences of the concentration inequalities on the Johnson graphs}
	\label{sec:J-concentration-consequences}
	
	Given any normalized sequence $\xbar:=\xn$ in a Banach space $X$, the map $h_{\xbar}\colon ([\bN]^k,\dJk)$ given by $h_{\xbar}(\nbar)=\sum_{i=1}^k x_{n_i}$ is easily seen to be $2$-Lipschitz. If $\xbar=\en$ the canonical basis of $\ell_p$, then we set $h_p:=h_{\xbar}$ and observe that for any infinite subset $\bM$ of $\bN$ and $\mbar\prec \nbar \in [\bM]^k$ we have $\norm{h_p(\mbar)-h_p(\nbar)}= (2k)^{1/p}$. Therefore, for every $p\in[1,\infty)$ there is an $\ell_p$-valued map on the Johnson graphs that does not $O(k^{1/q})$-concentrate for any $q\in(p,\infty)$. Consequently, the following corollary immediately follows from the coarse-Lipschitz stability result in Proposition \ref{prop:JFC-invariant}.
	
	\begin{coro}
		\label{cor:lp-not-CL-in-lq}
		Let $1\le p<q<\infty$.
		Then, $\ell_p$ does not coarse-Lipschitz embed into $\ell_q$.
	\end{coro}
	
	Corollary \ref{cor:lp-not-CL-in-lq} can be used to complete the proof of the uniform rigidity of $\ell_p$ for $p\in(2,\infty)$ via Proposition \ref{prop:rigidity-ell_p}. Note that this argument does not involve a Gorelik principle.
	
	Another striking consequence of the concentration inequality in the power-type regime for Lipschitz maps on the Johnson graphs, via Corollary \ref{cor:lp-not-CL-in-lq}, is the coarse-Lipschitz rigidity of subspaces of $\ell_p$ in the case $p\in(2,\infty)$. The situation for the other values of $p$, as well as the case of quotients or subspaces of quotients of $\ell_p$, is more difficult and will be treated in Chapter \ref {chapter:AMP_II}. 
	
	\begin{theo}
		\label{thm:CL-rigidity-subspaces-lp>2}
		Let $p\in[2,\infty)$. If $X$ is a Banach space that coarse-Lipschitz embeds into $\ell_p$, then $X$ is isomorphic to a subspace of $\ell_p$.
	\end{theo}
	
	\begin{proof} 
		Let $\cU$ be a nontrivial ultrafilter on $\bN$. As usual, we deduce from our assumption that $X^{\cU}$ Lipschitz embeds into $(\ell_p)^{\cU}$, which is reflexive, and a now classical differentiation argument implies that $X$ is isomorphic to a subspace of $L_p$. The conclusion is then clear if $p=2$. So, let us assume that $2<p<\infty$ and, aiming for a contradiction, that $X$ is not isomorphic to a subspace of $\ell_p$. It then follows from a result by W.B Johnson and E. Odell \cite{JohnsonOdell1974} that $\ell_2$ linearly embeds into $X$ and therefore coarse-Lipschitz embeds into $\ell_p$. This is in contradiction with Corollary \ref{cor:lp-not-CL-in-lq}.
	\end{proof}
	It is worth pointing out that Theorem \ref{thm:CL-rigidity-subspaces-lp>2} does not seem achievable via a Gorelik principle, as any Gorelik principle seems to require a nonlinear \emph{equivalence}.
	
	\medskip The Tsirelson spaces $\Tsi$ and $\Tsi^*$ were the first examples of Banach spaces that do not contain any isomorphic copies of $\co$ or $\ell_p$ for any $p\in[1,\infty)$.
	For the asymptotic-$\co$ Tsirelson space $\Tsi^*$, the noncontainment property can be drastically improved.
	
	\begin{coro}
		\label{cor:T*-BLS}
		$\Tsi^*$ does not contain any coarse copies of $\co$ or $\ell_p$ for any $p\in[1,\infty)$.
	\end{coro}
	
	\begin{proof}
		It suffices to show that there is always a Lipschitz map $h$ on $\sJ_k$ into $\ell_p$ or $\co$ that does not $O(1)$-concentrate.
		We have already seen above that $h_p$ is such a map for $\ell_p$. For $\co$ consider the map  $h_0(\nbar)=\sum_{i=1}^k s_{n_i}$, where $s_n=\sum_{i=1}^n e_n$ is the $n$-th vector of the summing basis of $\co$. It remains to observe that $h_0$ is $1$-Lipschitz as a map from $\sJ_k$ to $\co$ and that for all $\M \in [\bN]^\omega$, $\diam(h_0(\sJ_k(\M)))=k$. It is instructive for the reader to try to provide the elementary proofs for these last two claims. We will discuss the summing basis in Chapter \ref{chapter:interlaced-graphs}, where the proofs of these facts can be found. Note that we could also have used Aharoni's theorem to assert that the $\sJ_k$ equi-bi-Lipschitzly embed into $c_0$. 
	\end{proof}

	\section{\texorpdfstring{Coarse-Lipschitz rigidity of $\ell_p\oplus \ell_q$ for $1<p,q<\infty$}{Coarse-Lipschitz rigidity of}}
	\label{sec:CL-rigidity-sums-of-lp}
	
	The main purpose of the work by Kalton and Randrianarivony \cite{KaltonRandrianarivony2008} was to show the uniqueness of the coarse structure of $\ell_p \oplus \ell_q$, for $1<p<2<q<\infty$.  It also implicitly contained a ``Gorelik free'' argument for the cases $1<p<q<2$ and $2<p<q<\infty$ that we will include. 
	
	%We leave it to the reader to check that the use of the Gorelik principle can also be replaced by Corollary \ref{CLembAUS} to show the uniqueness of the coarse structure of $\ell_p$, for $1<p<\infty$ (see Exercise \ref{ex:Hammingl_p}).
	
	The crucial step in \cite{KaltonRandrianarivony2008} was to show the following nonlinear statement. The argument is a clever mix of the use of the approximate midpoints techniques from Chapter \ref{chapter:AMP_I} and of Johnson graphs techniques from this chapter.
	
	\begin{prop}
		\label{prop:lr-notCL-in-lp+lq}
		Let $1\le p<q<\infty$ and $r\in [1,\infty)\setminus \{p,q\}$. Then, $\ell_r$ does not coarse-Lipschitz embed into $\ell_p\oplus \ell_q$.
	\end{prop}
	
	\begin{proof}
		We distinguish three cases.
		\begin{description}
			\item[$1\le p <q < r$] In this case, $\ell_r$ is $r$-AUS and $\ell_p \oplus \ell_q$ is $q$-AUC. Therefore, the conclusion is a consequence of Corollary \ref{cor:CLembAUC}.
			\item[$1\le  r<p<q$] In this case, $\ell_p \oplus \ell_q$ is reflexive and $p$-AUS and the conclusion follows from Corollary \ref{cor:J-KR-BLMS}, Proposition \ref{prop:JFC-invariant} and the fact that there is an $\ell_r$-valued map $h_r$ on the Johnson graphs that does not $O(k^{1/p})$-concentrate.
			\item[$1\le p<r<q<\infty$] This is the most interesting and delicate case, as this is the situation where the Johnson graphs technique will provide the answer that the Gorelik principle, or the approximate midpoint principle alone, does not. 
		\end{description}
		
		So let $f\colon \ell_r \to \ell_p\oplus_\infty \ell_q$ be a coarse-Lipschitz embedding. Without loss of generality we can assume that there exists $C\ge 1$ such that for all $x,y\in \ell_r$ with $\norm{x-y}_r\ge 1$ we have 
		\begin{equation*}
			\norm{x-y}_r\le \norm{f(x)-f(y)}\le C\norm{x-y}_r.
		\end{equation*}
		Write $f=(f_p,f_q)$ and fix $k\in \bN$ and $\eps \in (0,\frac13)$. Let  $(e_n)_n$ be the canonical basis of $\ell_r$ and $E_N$ be the closed linear span of $(e_n)_{n=N}^\infty$. 
		
		We start with the application of the approximate midpoint technique to the coarse-Lipschitz map $f_p\colon \ell_r\to \ell_p$. We fix $\delta \in (0,\frac12)$, which we will make precise later. We recall that Lemma \ref{lem:CLApproxMid} ensures that for any $t\ge 1$, there exist $x,y \in \ell_r$ so that $\norm{x-y}_r\ge t$ and
		\begin{equation*}
			f_p(\Mid(x,y,\delta)) \subset \Mid(f_p(x), f_p(y),2\delta).
		\end{equation*}
		By Lemmas \ref{lem:midpointsAUS} and \ref{lem:midpointsAUC} and because $\Lip_1(f_p)\le C$, we get that there exist $Y \in \cof(\ell_r)$ and $K$ a compact subset of $\ell_p$ such that
		\begin{equation*}
			f_p(u + c_r\delta^{1/r}\norm{x-y}_rB_Y)\subset c_pC\delta^{1/p}\norm{x-y}_rB_{\ell_p} + K,
		\end{equation*}
		where $u:=\frac{x+y}{2}$ and $c_r,c_p$ are positive constants depending respectively on $r$ and $p$. By approximating, we may as well assume that $Y=E_N$, for some large enough $N\in \bN$. If $\delta>0$ was initially chosen such that $c_r^{-1}c_pC\delta^{1/p - 1/r}<\eps$ (this is possible because $p<r$) and if the approximate midpoint principle was applied with $t\ge 1$ such that $c_r\delta^{1/r}t\ge k$, we can deduce that there exist $\tau>k$, $u\in
		\ell_r$, $N\in \bN$ and $K$ a compact subset of $\ell_p$ such that
		\begin{equation*}
			f_p(u+\tau B_{E_N}) \subset \eps\tau B_{\ell_p} + K.
		\end{equation*}
		Let $\bM=\{n\in \bN \colon n>N\}$ and define $\varphi\colon [\bM]^k \to \ell_r$ by 
		$$\varphi(\mbar) := u+\tau k^{-1/r}(e_{m_1}+\dots+e_{m_k}),\ \mbar \in [\bM]^k.$$
		It is clear that $\varphi([\bM]^k)\subset u+\tau B_{E_N}$ and thus we have that
		$(f_p\circ \varphi)([\Mdb]^k)\subset \eps\tau B_{\ell_p} + K$. Thus, by Ramsey's theorem (Corollary \ref{cor:Ramsey_concentration}), 
		there is an infinite subset $\bM_1$ of $\bM$ such that $\diam (f_p\circ \varphi)([\bM_1]^k)\le 3\eps\tau.$
		
		Despite $\diam \varphi([\bM]^k)\ge \tau$ we cannot conclude yet since $f_p$ is not necessarily a coarse-Lipschitz embedding. Therefore, we need to turn our attention to the map $f_q$. We have
		$\norm{\varphi(\mbar) - \varphi(\nbar)}_r\ge \tau k^{-1/r}\ge 1$ whenever $\mbar \neq \nbar$ and $\norm{\varphi(\mbar) - \varphi(\nbar)}_r \le \tau k^{-1/r} (2\dJk(\mbar,\nbar))^{1/r} $. So, $\Lip(f_q\circ
		\varphi)\le \Lip(f\circ \varphi)\le 2\tau Ck^{-1/r}$. Then, since $\ell_q$ is reflexive and $q$-AUS, by Corollary \ref{cor:J-KR-BLMS}, there exists an infinite subset $\bM_2$ of $\bM_1$ such that
		$\diam (f_q\circ \varphi)([\bM_2]^k)\le C_qC\tau k^{1/q-1/r}$, for some constant $C_q>0$ that depends only on $q$. So, if $k$ was initially chosen large enough and since $\frac1q - \frac1r<0$, we get that  $\diam (f_q\circ \varphi)([\bM_2]^k)\le \eps\tau$.
		
		Finally we have that $\diam (f\circ \varphi)([\bM_2]^k)\le
		3\eps\tau$, while $\diam \varphi([\bM_2]^k)\ge \tau,$ which is impossible because $\eps<\frac13$.
	\end{proof}
	
	With the help of Proposition \ref{prop:lr-notCL-in-lp+lq}, we can prove the following coarse and uniform rigidity result.
	
	\begin{theo}
		\label{thm:rigidity-lp+lq} 
		Let $1<p<q<\infty$ such that $2\notin\{p,q\}$. Assume that $X$ is Banach
		space that is coarse-Lipschitz equivalent to $\ell_p\oplus \ell_q$. Then, $X$ is linearly isomorphic to $\ell_p\oplus\ell_q$.
		In particular, for these values of $p$ and $q$, $\ell_p\oplus \ell_q$ is uniformly and coarsely rigid.
	\end{theo}
	
	\begin{proof}
		Yet again, the key point is to show that $X$ does not contain any isomorphic copy of
		$\ell_2$, but this follows clearly from the Proposition \ref{prop:lr-notCL-in-lp+lq}. The fact that this is enough to imply the conclusion is a bit more delicate. The argument below gives us a golden opportunity to showcase some gems of the classical linear theory of Banach spaces.  
		
		The proof starts pretty much as we did for the uniform rigidity of $\ell_p$. Let $\cU$ be a nonprincipal ultrafilter on $\bN$ and observe that as in the proof of Theorem \ref{thm:CL-rigidity-ell_p}, we have that $X^{\cU}$ is Lipschitz isomorphic to $(\ell_p\oplus \ell_q)^{\cU}$. For identical reasons, we deduce that $X$ is isomorphic to a complemented subspace of $L_p \oplus L_q$. The key linear weapon is a result of W.B. Johnson \cite{Johnson1976} stating that if $r>2$ and $X$ is a Banach space not containing an isomorphic copy of $\ell_2$, then any bounded operator $T$ from $X$ to $L_r$ factors through $\ell_r$. This means that there exist $A\in B(X,\ell_r)$ and $B\in B(\ell_r,L_r)$ such that $T=BA$.
		
		Therefore, if we assume that $2<p<q$, then combining Johnson's theorem with the fact that $X$ is isomorphic to a complemented subspace of $L_p \oplus L_q$, we deduce easily that $X$ is isomorphic to a complemented subspace of $\ell_p\oplus\ell_q$. The spaces $\ell_p$ and $\ell_q$ are totally incomparable, meaning that they have no isomorphic infinite-dimensional subspaces. We can now use a theorem
		of {\`{E}}del\v{s}te{\u\i}n and  Wojtaszczyk \cite{EW1976} which states that a complemented subspace of the direct sum of two totally incomparable Banach spaces is isomorphic to a direct sum of complemented subspaces of the summands. Therefore, $X$ is isomorphic to $F\oplus G$, where $F$ is a complemented subspace of $\ell_p$ and $G$ is a complemented subspace of $\ell_q$. Then, it follows from Pe{\l}czy{\'n}ski's classical description of the complemented subspaces of $\ell_r$ \cite{Pelczynski1960} that $F$ is isomorphic to $\ell_p$ or finite-dimensional and $G$ is isomorphic to $\ell_q$ or finite-dimensional. Summarizing, since $X$ must be infinite-dimensional, we have that $X$ is isomorphic to $\ell_p$, $\ell_q$ or $\ell_p\oplus
		\ell_q$. But we already know that $\ell_p$ and $\ell_q$ have a unique coarse-Lipschitz structure. So, $X$ cannot be at the same time isomorphic to $\ell_p$ or $\ell_q$ and coarse-Lipschitz equivalent to $\ell_p \oplus \ell_q$. Therefore, $X$ is isomorphic to $\ell_p\oplus \ell_q$.
		
		Assume now that $1<p<q<2$ and note that $X^*$ is isomorphic to a complemented subspace of $L_{p'}\oplus L_{q'}$, where $2<q'<p'$ with $p'$ and $q'$ the conjugate exponents of $p$ and $q$ respectively. Since $L_{p'}\oplus L_{q'}$ and therefore $X^*$ are of type 2, we deduce from a theorem of B. Maurey (see Theorem 7.4.8 in \cite{AlbiacKalton2016}) that if $X^*$ contains an isomorphic copy of $\ell_2$, then $\ell_2$ must be complemented in $X^*$ and thus $\ell_2$ linearly embeds into $X^{**}=X$, which is not. Thus, $X^*$ does contain an isomorphic copy of $\ell_2$ and we deduce, again from Johnson's theorem, that $X^*$ is isomorphic to a complemented subspace of $\ell_{p'}\oplus\ell_{q'}$. The reflexivity of $X$ implies that $X$ is isomorphic to a complemented subspace of $\ell_p\oplus\ell_q$. The end of the argument is now the same as in the first case above.
		
		Assume finally that $p<2<q$. The application of Johnson's theorem yields that $X$ is isomorphic to a complemented subspace of $L_p\oplus\ell_q$. In this situation, the spaces $L_p$ and $\ell_q$ are also totally incomparable, so by the {\`{E}}del\v{s}te{\u\i}n-Wojtaszczyk Theorem, $X$ is isomorphic to $F\oplus
		G$, where $F$ is a complemented subspace of $L_p$ and $G$ is a complemented subspace of $\ell_q$. As we have already seen, $G$ is isomorphic to $\ell_q$ or finite-dimensional. We use again the fact that $F$ cannot contain an isomorphic copy of $\ell_2$, to apply the Johnson-Odell theorem from \cite{JohnsonOdell1974} and infer that $F$ is isomorphic to $\ell_p$ or finite-dimensional. The conclusion follows easily, as in the two previous cases above.
		
	\end{proof}
	
	\begin{rema} 
		As we already said, the cases $1<p<q<2$ and $2<p<q$, were originally settled in \cite{JLS1996} and we leave it to the reader to replace the use of Corollary \ref{cor:J-KR-BLMS} by the application of the Gorelik principle in the range $2<p<q$ (see Exercise \ref{ex:Gorelikl_p+l_q}). The case $1<p<2<q$ was solved in \cite{KaltonRandrianarivony2008} and the use of Johnson graphs was first devised on this occasion. It is very unlikely that it could be replaced by the Gorelik principle. The reason is that the Gorelik principle only applies to coarse-Lipschitz equivalences, whereas the concentration phenomenon for Johnson graphs applies to coarse-Lipschitz maps. This was crucial in this situation, as the component $f_q$ of a coarse-Lipschitz equivalence $f=(f_p,f_q)$ has no reason to be a coarse-Lipschitz equivalence.
	\end{rema}
	
	\begin{rema}
		\label{rema:fintesum} Theorem \ref{thm:rigidity-lp+lq} extends to finite sums of $\ell_p$ spaces. More precisely,
		if $1<p_1<\dots<p_n<\infty$ are all different from 2, then $\ell_{p_1}\oplus \dots \oplus
		\ell_{p_n}$ has a unique coarse structure. This statement is proved in \cite{KaltonRandrianarivony2008}. See Exercise \ref{ex:finitesum}.
	\end{rema}

	\section{\texorpdfstring{Towards the coarse-Lipschitz rigidity of $\ell_2\oplus\ell_p$}{Towards the coarse-Lipschitz rigidity of}}
	\label{sec:CL-rigidity-l2+lp}
	
	It is clear that the linear arguments in the proof of Theorem \ref{thm:rigidity-lp+lq} cannot be used if $p$ or $q$ is equal to $2$. 
	This leaves open the following rigidity problem.
	\begin{prob}
		\label{prob:rigidity-lp+l2}
		Let $p\in(1,\infty)\setminus\{2\}$.
		Is $\ell_2\oplus\ell_p$ uniformly, or coarsely, rigid?
	\end{prob}
	
	Problem \ref{prob:rigidity-lp+l2} seems to be a first step towards the following fundamental problem. 
	
	\begin{prob}
		\label{prob:rigidity-Lp}
		Let $p\in(1,\infty)$.
		Is $L_p$ uniformly, or coarsely, rigid?
	\end{prob}

	The paper \cite{KaltonRandrianarivony2008} contains some important partial information on the coarse structure of $\ell_2\oplus\ell_p$ that we shall now explain. The main new statement is the following.
	
	\begin{theo}
		\label{thm:lp(l2)-notCL-in-l2+lp} 
		Let $p\in (2,\infty)$. Then, $\ell_p(\ell_2)$ does not coarse-Lipschitz embed into $\ell_2\oplus \ell_p$.
	\end{theo}
	
	To prove Theorem \ref{thm:lp(l2)-notCL-in-l2+lp} we need the following variant of Lemma \ref{lem:midpointsAUS} whose proof is deferred to Exercise \ref{ex:ampl_p(l_2)}
	
	\begin{lemm}
		\label{lem:midpointslpl2}
		Let $p\in (2,\infty)$. For $i \in \bN$, we denote by $(e^{(i)}_{j})_{j=1}^\infty$ the canonical basis of the $i^{th}$ copy of $\ell_2$ in $\ell_p(\ell_2)$ and, for $N\in \bN$, we denote by $E_N$ the closed linear span of $\{e^{(i)}_{j}\colon j\in \bN,\ i\ge N\}$. Then, there exists $c_1>0$ such that for every $x,y \in \ell_p(\ell_2)$, every $\delta \in (0,1)$, there exists $N\in \bN$ so that
		\begin{equation*}
			\frac{x+y}{2}+c_1\delta^{1/p}\norm{x-y}B_{E_N} \subset \Mid(x,y,\delta).
		\end{equation*}
	\end{lemm}
	
	\begin{proof}[Proof of Theorem  \ref{thm:lp(l2)-notCL-in-l2+lp}] 
		Assume that $\ell_p(\ell_2)$ coarse-Lipschitz embeds into $\ell_2\oplus \ell_p$. Then, there is a map $f\colon X:=\ell_p(\ell_2) \to \ell_2\oplus_\infty \ell_p$ and a constant $C\ge 1$ such that for all $x,y\in X$ with $\norm{x-y}\ge 1$ we have 
		\begin{equation*}
			\norm{x-y}\le \norm{f(x)-f(y)}\le C\norm{x-y}.
		\end{equation*}
		Write $f=(f_2,f_p)$ and fix $k\in \bN$ and $\eps \in (0,\frac13)$.
		
		First, we apply the approximate midpoint principle to the map $f_2$. Since $p>2$, if we combine Lemmas \ref{lem:CLApproxMid}, \ref{lem:midpointslpl2} and \ref{lem:midpointsAUC}, then we get as in the proof of Proposition \ref{prop:lr-notCL-in-lp+lq} that there exist $\tau >k$, $x\in X$, $N\in \bN$ and a compact subset $K$ of $\ell_2$ so that
		\begin{equation*}
			f_2(x+\tau B_{E_N}) \subset  \eps\tau B_{\ell_2} + K.
		\end{equation*}
		We now define a map $\phi \colon [\bN]^k \to \ell_p(\ell_2)$, by 
		$$\phi(\nbar):= x + \tau k^{-1/2}\sum_{i=1}^k e^{(N)}_{n_i},\ \nbar=(n_1,\dots,n_k)\in [\bN]^k.$$ 
		Clearly $\phi([\bN]^k) \subset x+\tau B_{E_N}$ and therefore $(f_2 \circ \phi)([\bN]^k) \subset \eps\tau B_{\ell_2} + K.$ Then, we can apply Ramsey's theorem (Corollary \ref{cor:Ramsey_concentration}) to obtain that there
		is an infinite subset $\bM_1$ of $\bN$ such that $\diam (f_2\circ \phi)([\bM_1]^k)\le 3\eps\tau.$
		
		Note that $\norm{\phi(\mbar) - \phi(\nbar)}\le \sqrt{2}\tau k^{-1/2}\sqrt{\dJk(\mbar,\nbar)}$ and we deduce that $\phi$ is $\sqrt{2}\tau k^{-1/2}$-Lipschitz. On the other hand, $\norm{\phi(\mbar) - \phi(\nbar)}\ge \tau k^{-1/2}\ge 1$ whenever $\m \neq \n$. Therefore, $f_p \circ \phi$ is $C\sqrt{2}\tau k^{-1/2}$-Lipschitz. Since $\ell_p$ is $p$-AUS, Corollary \ref{cor:J-KR-BLMS} now implies that there is an infinite subset $\bM_2$ of $\bM_1$ such that $\diam (f_p\circ \phi)([\bM_2]^k)\le c_2C\tau k^{1/p-1/2}$, for some universal constant $c_2>0$. Since $p>2$, if $k$ was initially chosen large enough, we get that $\diam (f_p\circ \phi)([\bM_2]^k)\le 3\eps\tau$ and thus $\diam (f\circ \phi)([\bM_2]^k)\le 3\eps\tau <\tau$. But $\diam\, \phi([\bM_2]^k) = \sqrt{2}\tau$, which implies that $\diam (f\circ \phi)([\bM_2]^k)\ge \sqrt{2}\tau$.
		This contradiction concludes our proof.
	\end{proof}
	
	As a corollary, we have.
	
	\begin{coro} 
		\label{cor:Lp-notCL-in-l2+lp}
		Let $p\in [1,\infty)\setminus \{2\}$. Then, $L_p$ does not coarse-Lipschitz embed into $\ell_2 \oplus \ell_p$.
	\end{coro}
	
	\begin{proof} 
		Assume first that $p\in (2,\infty)$. Since it is well known that $\ell_p(\ell_2)$ is isometric to a subspace of $L_p$, this is a direct consequence of Theorem \ref{thm:lp(l2)-notCL-in-l2+lp}.
		
		Assume now that $p\in [1,2)$. Then, for any $r\in (p,2)$, it follows from the existence of $r$-stable random variables (see for instance Theorem D.8 in \cite{BenyaminiLindenstrauss2000}) that $\ell_r$ is isometric to a subspace of $L_p$. The statement then follows from Proposition \ref{prop:lr-notCL-in-lp+lq}.
	\end{proof}
	
	To the best of our knowledge, the following problem is open.
	
	\begin{prob}
		\label{pb:lp(l2)-CL-in-l2+lp}
		Let $p\in[1,2)$. Does $\ell_p(\ell_2)$ coarse-Lipschitz embed into $\ell_2 \oplus \ell_p$?
	\end{prob}

	\section[Johnson-concentration, reflexivity and spreading models]{Johnson-concentration inequalities, reflexivity and spreading models}\label{sec:J-reflexivity-spreading}
	
	In this section, we discuss the impact on the geometry of $X$ of the existence of concentration inequalities for $X$-valued maps on the Johnson graphs. The first interesting fact is that a nontrivial concentration inequality implies reflexivity.
	
	\begin{prop}
		\label{prop:J-reflexivity}
		If the $X$-valued Lipschitz maps on $\sJ_k$ satisfy a $O(\gamma(k))$-concentration inequality with $\lim_{t\to \infty}\gamma(t)/t=0$, then $X$ is reflexive. 
	\end{prop}
	
	\begin{proof}
		Assume that $X$ is nonreflexive. By one of James' characterization of reflexive spaces \cite{James1972}, there exists a sequence $(x_n)_{n=1}^\infty \subset B_X$ such that for all $k\in \bN$, all $x\in \conv(\{x_1,\dots,x_k\})$ and all $y\in \conv(\{x_{i} \colon i\ge k+1\})$, one has $\norm{ x-y }\ge \frac12$. In particular, for all $k\in \bN$ and all $(m_1,m_2,\dots, m_{2k})\in[\bN]^{2k}$,
		\begin{equation}
			\label{eq:James-reflexivity}
			\Big\|\sum_{i=1}^k x_{m_i}-\sum_{i=k+1}^{2k} x_{m_i}\Big\|\ge \frac{k}{2}.
		\end{equation}
		For $k\in \bN$, consider $h(\nbar)=\sum_{i=1}^k x_{n_i}$, for $\nbar=(n_1,\dots, n_k)$ in $[\bN]^k$. Since $(x_n)_{n=1}^\infty \subset B_X$, $h$ is $1$-Lipschitz for $d_{\sJ_k}$. On the other hand, it follows from (\ref{eq:James-reflexivity}) that for  any $\bM\in [\bN]^\omega$, $\diam(h([\bM]^k))\ge \frac{k}{2}$. This shows that there is always an $X$-valued Lipschitz map on $\sJ_k$ that cannot satisfy a $O(\gamma(k))$-concentration inequality with $\lim_{t\to \infty}\gamma(t)/t=0$ whenever $X$ is not reflexive.
	\end{proof}
	
	The fact that $\ell_1$ admits a coarse and uniform embedding into $\ell_2$ shows that, in general, reflexivity is not preserved under nonlinear embeddings. However, it follows immediately from Proposition \ref{prop:J-reflexivity}, Proposition \ref{prop:JFC-invariant} and Corollary \ref{cor:J-KR-BLMS} that if we restrict ourselves to smooth spaces, then reflexivity is preserved under nonlinear embeddings. In fact, the nonlinear embeddings can be less faithful if the spaces are smoother.
	
	\begin{coro}
		Let $X$ be a Banach space and $Y$ be a reflexive Banach space.
		\begin{enumerate}
			\item If $X$ coarse-Lipschitz embeds into $Y$ and $Y$ is asymptotically uniformly smooth, then $X$ is reflexive.
			\item If $X$ coarsely embeds into $Y$ and $Y$ is asymptotic-$\co$, then $X$ is reflexive.
		\end{enumerate}
	\end{coro}
	
	We also know that under the assumption of item (1) (resp. (2)) the $X$-valued Lipschitz maps on $\sJ_k$ satisfy a $O(k^{1/p})$-concentration inequality for some $p\in(1,\infty)$ (resp. a $O(1)$-concentration inequality) but we do not know if it implies that $X$ must be asymptotically uniformly smooth (resp. asymptotic-$\co$).
	However, one can deduce information on the spreading models of the space $X$ from concentration inequalities for $X$-valued Lipschitz maps on $\sJ_k$.
	
	\begin{prop}
		\label{prop:J-spreading-model}
		If the $X$-valued Lipschitz maps on $\sJ_k$ satisfy a $\gamma(k)$-concentration inequality, then for every spreading model $S=[\en]$ of $X$ that is generated by a normalized weakly null sequence, one has that for all $k\ge 1$,
		\begin{equation*}
			\Big\|\sum_{i=1}^k e_i \Big\|_S\le 2\gamma(k).
		\end{equation*} 
		In particular, if the $X$-valued Lipschitz maps on $\sJ_k$ satisfy a $O(1)$-concentration inequality, then every spreading model of $X$ generated by a normalized weakly null sequence is isomorphic to $\co$.
	\end{prop}
	
	\begin{proof}
		Let $\xn$ be a normalized weakly null sequence in $X$ that generates a spreading model $S=[\en]$. Fixing $k\in\bN$ and $\delta>0$ and passing to an appropriate subsequence we may assume that for any $k\le j_1<\cdots<j_k$ and any $a_1,\ldots,a_k$ in $[-1,1]$ we have
		\begin{equation}
			\label{E:4.4.1}
			\abs{\Big\|\sum_{i=1}^k a_i x_{j_i}\Big\|_X - \Big\| \sum_{i=1}^k a_i e_i \Big\|_S } < \delta.
		\end{equation}
		
		Define $h \colon[\bN]^k\to X$ as usual via the formula $h(\nbar) = \sum_{i=1}^k x_{n_i}$ for $\nbar=\{n_1,n_2,\dots,n_k\}\in [\bN]^k$. 
		Then, $h$ is $1$-Lipschitz with respect to the Johnson graph metric and it follows from the concentration property that there exists $\bM\in[\bN]^\omega$ such that for all $\mbar\prec\nbar\in [\bM]^k$ one has
		\begin{equation}
			\norm{ h(\mbar) - h(\nbar) }_X = \Big\|\sum_{i=1}^k x_{m_i} - \sum_{i=1}^k x_{n_i} \Big\|_X \le \gamma(k).
		\end{equation}
		
		Fix $\vep>0$. Since the sequence is weakly null, for all $\mbar \in[\bM]^k$, we can find $x^*\in S_{X^*}$ and $\nbar\in[\bM]^k$ such that $\mbar\prec\nbar$ and 
		
		\begin{equation*}
			x^*\Big(\sum_{i=1}^k x_{m_i} - \sum_{i=1}^k x_{n_i}\Big)\ge (1-\vep)\Big\|\sum_{i=1}^k x_{m_i}\Big\|_X - \vep,
		\end{equation*}
		
		and thus 
		$$\Big\|\sum_{i=1}^k x_{m_i} \Big\|_X\le \frac{1}{1-
			\vep}(\gamma(k)+\vep).$$ 
		Therefore, if one chooses $\delta,\vep$ small enough, we have that for all $k\ge 1$, $\norm{ \sum_{i=1}^k e_i }_S\le 2\gamma(k)$.
		The last statement follows from the unconditionally of spreading models and a classical extreme point argument.
	\end{proof}
	
	Proposition \ref{prop:J-spreading-model} shows the limits of applicability of Johnson graphs techniques. Concentration inequalities for $X$-valued Lipschitz maps on Johnson graphs seem to only contain partial information about the spreading models of $X$.
	In order to get more information about the asymptotic geometry of $X$, we need to consider another sequence of graphs: the Hamming graphs. The study of the geometry of the Hamming graphs is detailed in Chapter \ref{chapter:Hamming}.

	\section{Notes}\label{sec:Notes-Johnson}

	It is well known that $L_p$ admits a coarse embedding into $\ell_p$ when $p\in[1,2)$. Indeed, it goes back to the work of Bretagnolle, Dacunha-Castelle and Krivine \cite{BDCK} (based on earlier work on Schoenberg on negative-type kernels \cite{Schoenberg1938}) that the $p/q$-snowflaking of $L_p$ embeds isometrically into $L_q$ whenever $1\le p\le q \le  2$. In particular, $(L_p, \norm{\cdot}_p^{p/2})$ admits an isometric embedding into $L_2$. The limitation on the range of $p$ is inherent in negative-type kernel techniques. Thus, $L_p$ coarsely embeds into $L_2$ when $p\in[1,2)$. Since Nowak \cite{Nowak2006} had shown that $\ell_2$ coarsely embeds into $\ell_p$ for every $p\in[1,\infty)$, it follows that we can coarsely embed $L_p$ into $\ell_q$ whenever $p\in[1,2]$ and $q\in[1,\infty)$. It would be very interesting to construct directly a coarse embedding of $L_p$ into $\ell_p$ without going through Hilbert space. For more on this problem and a finer quantitative analysis of the coarse embedding between $L_p(\mu)$-spaces we refer to \cite{Baudier2016}. 
	
	The situation for $p\in(2,\infty)$ has yet to be elucidated. What we know from Corollary \ref{cor:lp-not-CL-in-lq} or Corollary \ref{cor:Lp-notCL-in-l2+lp} is that $L_p$ does not coarse-Lipschitz embed into $\ell_p$ when $p\in(2,\infty)$. It is quite remarkable that the following problem remains open.
	\begin{prob}
		\label{pb:Lp-coarse-lp}
		Does $L_p$ coarsely embed into $\ell_p$ when $p\in(2,\infty)$?
	\end{prob}
	Problem \ref{pb:Lp-coarse-lp} is arguably one of the most fundamental open problems in the coarse geometry of Banach spaces.
	The results in the chapter give some partial information on Problem \ref{pb:Lp-coarse-lp}.
	Mendel and Naor provided in \cite{MendelNaor2008} a simple explicit formula for an isometric embedding of the $p/q$-snowflaking of $L_p$ into $L_q$ for all $1\le p<q<\infty$. Given $1\le p<q<\infty$, identifying the optimal $s\in(0,1)$ such that the $s$-snowflaking of $L_p$ embeds into $L_q$ was a longstanding open problem. The work of Naor and Schechtman on the metric $X_p$-inequalities \cite{NaorSchechtman2016}, \cite{Naor2016} resolved this problem when $2\le p<q$ as they show that $s=p/q$ is the best that can be achieved.
	It is also known that the $p/q$-snowflaking of $\ell_p$ admits a bi-Lipschitz embedding into $\ell_q$ for all $1\le p<q<\infty$. This can be done in various ways depending on the values of $p$ and $q$ (see for instance \cite{AlbiacBaudier2015} for an approach using wavelets in the spirit of Assouad's embedding theorem, or Ostrovskii \cite{Ostrovskii2009} for a discretization of the Mendel-Naor embedding).
	A nice consequence of the work of Kalton and Randrianarivony on the geometry of the Johnson graphs is that $s=p/q$ is the optimal  snowflake exponent $s$ such that the $s$-snowflaking of $\ell_p$ can bi-Lipschitzly embed into $\ell_q$ when $1\le p<q<\infty$.
	The next theorem contains a quantitative refinement of Corollary \ref{cor:lp-not-CL-in-lq}.
	\begin{theo}
		\label{thm:snowflake-lp-lq}
		Let $1\le p<q<\infty$. Then,
		\begin{equation}
			\mathsf{s}_{\ell_q}(\ell_p)=\alpha_{\ell_q}(\ell_p)=\frac{p}{q},
		\end{equation}
		where $\mathsf{s}_{\ell_q}(\ell_p):=\sup\{s\in (0,1],\ (\ell_p,\|\ \|_p^s)\ \text{bi-Lipschitz embeds into}\  (\ell_q,\|\ \|_q)\}$ and 
		$\mathsf{\alpha}_{\ell_q}(\ell_p):=\sup\{s\in (0,1],\ (\ell_p,\|\ \|_p^s)\ \text{coarse-Lipschitz embeds into}\  (\ell_q,\|\ \|_q)\}$.
	\end{theo}

	It would be very surprising that there are different answers to Problem \ref{pb:Lp-coarse-lp} for different values of $p\in(2,\infty)$.
	It follows from Theorem \ref{thm:snowflake-lp-lq} and the fact that $L_p$ contains an isometric copy of $\ell_2$ that 
	$\alpha_{\ell_p}(L_p)\le \frac{2}{p}$ and thus $\lim_{p\to \infty}\alpha_{\ell_p}(L_p)=0$.
	This means that if there is a coarse embedding of $L_p$ into $\ell_p$ for any $p\in(2,\infty)$, the faithfulness of this embedding has to deteriorate as $p$ gets bigger. We can think of this as some (very partial) indication that the answer to Problem \ref{pb:Lp-coarse-lp} might be negative.
	
	Corollary \ref{cor:T*-BLS} can be seen as a coarse version of Tsirelson's result \cite{Tsirelson1974}, who constructed the first example of a Banach space that does not contain any linear copy of $\ell_p$ for $1\le p<\infty$ or $\co$. It also answers negatively a question that was open for a while: does any infinite-dimensional Banach space contain a coarse copy of $\ell_2$? This problem was solved in \cite{BLS2018}. Note also that if an infinite-dimensional Banach space $X$ coarsely embeds into $\Tsi^*$, then it must uniformly contain isomorphic copies of the $\ell_\infty^n$ (see Exercise \ref{ex:trivialcotype}) and therefore has trivial cotype (in the next chapter, we will actually show that $X$ must be asymptotic-$\co$). Thus, it follows from Mendel-Naor metric cotype obstruction \cite{MendelNaor2008} that $X$ is not coarsely embeddable into $\ell_2$. As a consequence, there is no minimal infinite-dimensional Banach space for coarse embeddings. 
	
	\section{Exercises}
	
	\begin{exer}\label{ex:J-graph}
		Let $\M_1$ and $\M_2$ be two infinite subsets of $\bN$. Show that $\sJ_k(\M_1)$ and $\sJ_k(\M_2)$ are graph-isomorphic.    
	\end{exer}
	
	\begin{exer}
		\label{ex:finitesum}
		Show Remark \ref{rema:fintesum}. Hint: use an induction on $n\in \bN$.  
	\end{exer}

	\begin{exer}
		\label{ex:Gorelikl_p+l_q}
		Assume that $2<p<q<\infty$. Prove that $\ell_p\oplus \ell_q$ has a unique coarse-Lipschitz structure without using Corollary \ref{cor:J-KR-BLMS}. Use the Gorelik principle instead.  
	\end{exer}
	
	\begin{exer}
		\label{ex:ampl_p(l_2)}
		Show Lemma \ref{lem:midpointslpl2}.   
	\end{exer}
	
	\begin{exer}\label{ex:trivialcotype} 
		Assume that $X$ is an infinite-dimensional Banach space that coarsely embeds into $\Tsi^*$. Show that $X$ uniformly contains isomorphic copies of the $\ell_\infty^k$. Hint: use the results of Section \ref{sec:J-concentration-consequences}. 
		
		For a more elementary approach: use that the $X$-valued Lipschitz maps satisfy a $O(1)$-concentration inequality to show that for any $k\in \bN$, there exists a basic sequence $(x_n)_n$ in $S_X$ such that for some $C\ge 1$, we have that $\big\|\sum_{i=1}^k\eps_ix_{n_i}\big\|\le C$ for all $(\eps_1,\ldots,\eps_k)\in \{-1,1\}^k$.
		
	\end{exer}

	%%%%%%%%%%%%%%%%%%%%%%%%%%%%%%%%%%%%%%%%%%%%%%%%%%%%%%%%%%%%%%%%%%%%%%%%%%%%%%%%%
	
	%%%%%%%%%%%%%%%%%%%%%%%%%%%%%%%%%%%%%%%%%%%%%%%%%%%%%%%%%%%%%%%%%%%%%%%%%%%%%%%%%

	\chapter[Geometry of the Hamming graphs and applications]{Geometry of the Hamming graphs and applications}
	\label{chapter:Hamming}

	In this chapter, we delve into the geometry of the Hamming graphs. All the results from Chapter \ref{chapter:Johnson} remain true if instead of the Johnson graphs, we consider the Hamming graphs. Therefore, we will only state these results and leave the slight modifications of the proofs to the reader. We will detail only those results for which Johnson graphs techniques are insufficient and the Hamming graphs are truly needed. The highlight of this chapter is a purely metric characterization of the class of reflexive and asymptotic-$\co$ Banach spaces in terms of a concentration inequality for Lipschitz maps on the Hamming graphs. A consequence of this metric characterization is that this class is stable under coarse embeddings. We also provide a metric characterization of the class of asymptotic-$\co$ spaces \emph{within the class of separable and reflexive Banach spaces}. This metric characterization is in terms of metric space preclusion. A recurrent theme of this chapter is the connection between Hamming graphs and asymptotic models of Banach spaces. It culminates in the last section, where the preservation of upper estimates for asymptotic models under coarse-Lipschitz embeddings is discussed.
	
	\section[Hamming graphs: basic properties and embeddings]{Hamming graphs: basic properties and embeddings into Banach spaces}
	\label{sec:Hamming}
	
	Consider a sequence of sequences of vectors (we call this an array) $[a]:=\{(a^{(j)}_i)_{i\in \bN} \colon j\in \bN\}$ whose vectors live in a Banach space $X$. This array induces, for every $k\in \bN$ and every infinite subset $\bM$ of $\bN$, a semi-metric on $[\bM]^{k}:=\{ S\subset \bM \colon |S|=k\}$, denoted by $\sd_{[a]}$ for simplicity, as follows:
	\begin{equation*}
		\sd_{[a]}(\mbar, \nbar):=\Big\|\sum_{i=1}^k a^{(i)}_{m_i} -\sum_{i=1}^k a^{(i)}_{n_i}\Big\|_X.
	\end{equation*}
	This semi-metric is a genuine metric if the array consists of linearly independent vectors. In particular, if the array is given by the canonical basis of $\ell_1(\ell_1)$, i.e. one takes $a^{(j)}_i:=e^{(j)}_i$ for all $i,j\ge 1$ where $e^{(j)}_i$ is the $i$-th vector of the canonical basis of the $j$-th copy of $\ell_1$ in $\ell_1(\ell_1)$, then $\sd_{[e]}$ is a metric and 
	\begin{equation}\label{eq:Hamming}
		\sd_{[e]}(\mbar, \nbar) = \Big\|\sum_{i=1}^k e^{(i)}_{m_i} -\sum_{i=1}^k e^{(i)}_{n_i}\Big\|_{\ell_1(\ell_1)} = \sum_{i=1}^k \|e^{(i)}_{m_i} - e^{(i)}_{n_i} \|_{1} = 2\abs{ \{ 1\le i \le k \colon m_i \neq n_i\} }
	\end{equation}
	The metric $\sd_{[e]}$ can also be realized, up to a dilation, as a graph metric. To see this, consider $[\bM]^k$ as the vertex set of a graph $\sH_k(\bM)$ where two distinct vertices are adjacent if and only if they differ in exactly one coordinate. This graph structure is well known on the hypercube $\{0,1\}^k$ and leads to what has been called the Hamming cubes. Note that $\sH_k(\bM)$ is an infinite graph with countable degree and diameter $\diam(\sH_k(\bM))= k$ and it would be natural to call those graphs the countably branching Hamming cubes. For the sake of simplicity, we will refer to them as the \emph{Hamming graphs}. One can show that the graph metric on $\sH_k(\bM)$, denoted by $\dHk$, satisfies $\dHk=\frac12 \sd_{[e]}$. It is also easy to see that the Hamming graphs on two infinite subsets of $\bN$ are graph-isomorphic (see Exercise \ref{exe:H-graph}) and hence isometric as metric spaces. Therefore, we will simply use the notation $\sH_k$ for the $k$-th Hamming graph if it does not lead to ambiguity. The identity in $\eqref{eq:Hamming}$ tells us that the embedding $\mbar\in [\bM]^k \mapsto \sum_{i=1}^k e^{(i)}_{m_i}$ is a scaled isometric embedding of $([\bM]^k,\dHk)$ into $\ell_1(\ell_1)\equiv\ell_1$.
	
	Similarly to the situation with the Johnson graphs, if we consider the array $[e]$ as the canonical basis of $\ell_p$, then it is clear that the natural embedding 
	$$\mbar\in \sH_k \stackrel{h_p}{\mapsto} \sum_{i=1}^k e^{(i)}_{m_i}\in \ell_p(\ell_p)\equiv \ell_p
	$$
	satisfies 
	$$\norm{h_p(\mbar) - h_p(\nbar) }_p = 2^{\frac{1}{p}}\dHk(\mbar,\nbar)^{\frac{1}{p}}.$$ This means that the $\frac1p$-snowflake of $(\sH_k,\dHk)$ admits a (scaled) isometric embedding into $\ell_p$.
	
	While how well the Johnson graphs can be embedded into infinite-dimensional Banach spaces is governed by the properties of their weakly null sequences, the embeddability of the Hamming graphs is closely related to the properties of weakly null arrays. It is not too difficult to show that Proposition \ref{prop:J-graph-embedding} remains true for the Hamming graphs (see Exercise \ref{exe:H-graph-embedding-spreading-models}). In particular, the sequence $\{\sH_k\}_{k=1}^\infty$ of Hamming graphs equi-bi-Lipschitzly embeds into any Banach space that admits a spreading model isomorphic to $\ell_1$. Proposition \ref{prop:H-graph-embedding} below, where we make our first encounter with asymptotic models (see Section \ref{sec:asymptotic-models} of Appendix \ref{appendix:asymptotic} for the basic definitions and properties), gives a much stronger conclusion. The proof is similar to the one of Proposition \ref{prop:J-graph-embedding} with a higher level of difficulty. The essential property of spreading models that was needed there was the suppression unconditionally of spreading models; the spreading property did not play any significant role. Asymptotic models are also suppression unconditional but there is an additional subtlety here. Given a normalized array $[a]$ in a Banach space $X$ the canonical map $h_k\colon \sH_k\to X$ given by $h_{k}(\mbar):= \frac12\sum_{i=1}^k a^{(i)}_{m_i}$ is clearly $1$-Lipschitz. A slight issue arises here since expressions of the form $h_k(\mbar)- h_k(\nbar)$ have potentially two vectors on the same row of the array, but in order to use the suppression unconditionality property of asymptotic models we can only allow one vector per row. Note that this problem did not arise in the context of Johnson graphs and spreading models. Luckily, this technical annoyance can be taken care of using joint spreading models techniques (see Proposition \ref{prop:jointspreadingmodels} in Section \ref{sec:joint-spreading-models}). 
	
	\begin{prop}
		\label{prop:H-graph-embedding}
		Let $X$ be an infinite-dimensional Banach space admitting an asymptotic model $A=[\en]$, generated by a normalized weakly null array, then for every $\nu>0$ and $k\ge 1$ there exists a map  $h_{k,\nu} \colon \sH_k \to X$ such that for all $\mbar,\nbar \in \sH_k$,
		\begin{equation}
			\label{eq:H-graph-embedding-eq1}
			\frac{1}{2(1+\nu)}\Big\|\sum_{i\colon m_i\neq n_i}e_i \Big\|_A\le \norm{h_{k,\nu}(\mbar)-h_{k,\nu}(\nbar)}_X\le \dHk(\mbar,\nbar).
		\end{equation}
	\end{prop}
	
	\begin{proof}
		Let $\{x^{(i)}_j)_{j\in \bN} \colon i\in \bN\}$ be a normalized weakly null array in $X$ that generates an asymptotic model $A=[(e_j)_{j=1}^\infty]$. Fixing $k\in\bN$ and $\delta>0$ and passing to appropriate subsequences of the array we may assume that for any $j_1<\dots<j_k$ and any $a_1,\dots,a_k$ in $[-1,1]$ we have
		\begin{equation*}
			\abs{ \Big\|\sum_{i=1}^k a_i x^{(i)}_{j_i} \Big\|_X - \Big\|\sum_{i=1}^k a_i e_i \Big\|_A } < \delta.
		\end{equation*}
		In addition,  by applying Proposition \ref{prop:jointspreadingmodels} we may also assume that for any $i_1,\dots,i_{2k}$ in $\{1,\dots,k\}$ and any pairwise different $l_1,\dots,l_{2k}$ in $\bN$ the sequence $(x^{(i_j)}_{l_j})_{j=1}^{2k}$ is $(1+\delta)$-suppression unconditional.
		
		Let $h:=h_{k,\nu}\colon \sH_k\to Y$  be defined by $h_{k,\nu}(\mbar):= \ds\frac12\ds\sum_{i=1}^k x^{(i)}_{km_i+i}$ for all $\mbar\in \sH_k$. Observe first that for $m_1 <\cdots<m_k$ we have $k< km_1 + 1< km_2+2<\cdots<km_k+k$. Then, if $\mbar = \{m_1,\dots,m_k\}$,
		$\nbar = \{n_1,\dots,n_k\}$ and $F = \{i:m_i\neq n_i\}$ we have  
		\begin{equation*}
			\norm{ h(\mbar) - h(\nbar) }_X = \frac12\Big\|\sum_{i\in F} (x^{(i)}_{km_i+i} - x^{(i)}_{kn_i+i}) \Big\|_X \le  \dHk(\mbar,\nbar),
		\end{equation*}
		and hence the map $h$ is $1$-Lipschitz.
		
		Also, note that $km_i + i = kn_{i'} + i'$ if and only if $i = i'$ and $m_i = n_{i'}$. We deduce that the sequence $(x^{(i)}_{km_i+i})_{i\in F}\cup(x^{(i)}_{kn_i+i})_{i\in F}$ is $(1+\delta)$-suppression unconditional. Therefore, we have
		\begin{equation*}
			\norm{ h(\mbar) - h(\nbar) }_X \ge \frac{1}{2(1+\delta)} \Big\| \sum_{i\in F} x^{(i)}_{km_i+i} \Big\|_X \ge \frac{1}{2(1+\delta)}\Big(\Big\|\sum_{i\in F}e_i \Big\|_A - \delta\Big),    
		\end{equation*}
		and choosing $\delta>0$ small enough gives the result.
	\end{proof}
	
	It follows from Proposition \ref{prop:H-graph-embedding} that any Banach space with an asymptotic model isomorphic to $\ell_1$ will equi-bi-Lipschitzly contain the Hamming graphs. It is worth pointing out that there are Banach spaces that do not have any spreading model isomorphic to $\ell_1$ but that have an asymptotic model isomorphic to $\ell_1$ (\cite{HalbeisenOdell2004}).
	
	\begin{rema}
		Proposition \ref{prop:H-graph-embedding} does not provide a nontrivial compression function when all asymptotic models are isomorphic to $\co$. Given an asymptotic model $A=[\en]$, let 
		$$\varphi_{A,k}(r):= \inf_{1\le r\le k ;  \abs{F}=r}\Big\|\sum_{i\in F} e_i\Big\|_A\ \ \text{and}\ \  \varphi_{A}(r):=\inf_{\abs{F}=r}\Big\|\sum_{i\in F} e_i\Big\|_A.$$ 
		Then, Proposition \ref{prop:H-graph-embedding} says that if a Banach space $X$ has an asymptotic model $A$, generated by a normalized weakly null array $[x]$, then for any $k\ge 1$ there is a map $h_{k}\colon \sH_k \to X$ of the form $h_{k}(\mbar):= \sum_{i=1}^k x^{(i)}_{m_i}$ whose modulus of compression satisfies $\rho_{h_{k}}(r)\gtrsim \varphi_{A,k}(r) \ge \varphi_{A}(r)$. The issue here is that $\varphi_{A}$ might not necessary tend to $\infty$ when $r\to\infty$. In fact, the space $\Tsi^*(\Tsi^*)$ is an example of a reflexive nonasymptotic$-\co$ Banach space admitting an asymptotic model $A$ that is not isomorphic to $\co$ with such property. Indeed, if it were not the case then, by Proposition \ref{prop:H-graph-embedding}, $\Tsi^*(\Tsi^*)$ would contain equi-coarse copies of the Hamming graphs via an embedding of the canonical form above, but this is not possible as shown in \cite[Proposition 4.6]{BLMS2020}.
		
	\end{rema}

	\section[Hamming-concentration and metric characterizations]{Hamming-concentration and applications to metric characterizations}
	\label{sec:H-concentration}
	
	The Johnson-concentration inequalities were immediately derived from the linearization lemma (Lemma \ref{lem:linearization}).
	The only thing that we needed in that case was that for all $1\le i \le k$, $\nbar\in [\bM]^{k}$ and $m\in \bM$ such that $\sigma_{i,m}(\nbar)\in [\bM]^k$, $\dJk(\nbar,\sigma_{i,m}(\nbar))\colon=1$, where we recall that $\sigma_{i,m}(\nbar):=(n_1,\dots,n_{i-1}, m, n_{i+1},\dots, n_k)$. Since in this situation $\dHk(\nbar,\sigma_{i,m}(\nbar))\colon=1$ as well, Theorem \ref{thm:J-KR-BLMS} remains true for the Hamming graphs in place of the Johnson graphs.
	
	We only state the most useful consequences.
	
	\begin{theo}\, 
		\label{thm:H-KR-BLMS}
		\begin{enumerate}
			\item Let $p\in(1,\infty)$ and $X$ be a reflexive Banach space with property $\sN_p$ (in particular $p$-asymptotically uniformly smooth). Then, the  $X$-valued Lipschitz maps on $\sH_k$ satisfy a $O(k^{\frac1p})$-concentration inequality.
			\item\label{item2:H-KR-BLMS} Let $X$ be a reflexive asymptotic-$\co$ Banach space. Then, the $X$-valued Lipschitz maps on $\sH_k$ satisfy a $O(1)$-concentration inequality.
		\end{enumerate}
	\end{theo}
	
	What is truly remarkable is that the converse of item $\eqref{item2:H-KR-BLMS}$ is true! Before we explain why, we collect for further reference the basic properties of the Hamming-concentration inequalities. The proofs are essentially the same as in the Johnson graph case and we invite the reader to verify this.
	The first proposition gives embedding obstructions for the Hamming graphs.
	
	\begin{prop}
		\label{prop:H-obstruction}
		Let $(X,d_X)$ be a metric space.
		\begin{enumerate}[(i)]
			\item If the $X$-valued Lipschitz maps on $\sH_k$ satisfy a $\gamma(k)$-concentration inequality with $\lim_{t\to \infty}\gamma(t)/t=0$, then $(\sH_k)_{k\in \bN}$ does not admit an equi-coarse-Lipschitz embedding into $X$.
			\item If the $X$-valued Lipschitz maps on $\sH_k$ satisfy a $O(1)$-concentration inequality, then $(\sH_k)_{k\in \bN}$ does not admit an equi-coarse embedding into $X$.
		\end{enumerate}
	\end{prop}
	
	The second proposition states the permanence properties of the Hamming-concentration inequalities.
	
	\begin{prop}
		\label{prop:HFC-invariant} Let $X$ and $Y$ be Banach spaces.
		\begin{enumerate}[(a)]
			\item Assume that the $Y$-valued Lipschitz maps on $\sJ_k$ satisfy a $O(\gamma(k))$-concentration inequality. If $X$ coarse-Lipschitz embeds into $Y$, then the $X$-valued Lipschitz maps on $\sJ_k$ satisfy a $O(\gamma(k))$-concentration inequality.
			\item Assume that the $Y$-valued Lipschitz maps on $\sJ_k$ satisfy a $O(1)$-concentration inequality. If $X$ coarsely embeds into $Y$, then the $X$-valued Lipschitz maps on $\sJ_k$ satisfy a $O(1)$-concentration inequality.
		\end{enumerate}
	\end{prop}
	
	Since every edge in $\sH_k$ is an edge in $\sJ_k$, it follows that $\dJk\le \dHk$ and hence any map $f\colon [\bM]^k \to X$ that is $C$-Lipschitz for $\dJk$ is $C$-Lipschitz for $\dHk$. 
	%Note also that if $\mbar\prec\nbar$, then $\dJk(\mbar,\nbar)=\dHk(\mbar,\nbar)=k$.   
	It follows from this observation that Proposition \ref{prop:J-reflexivity} and Proposition \ref{prop:J-spreading-model} remain true for the Hamming graphs. We only record the analog of the former proposition as we will be able to crucially improve the latter. 
	
	\begin{prop}
		\label{prop:H-reflexivity}
		If the $X$-valued Lipschitz maps on $\sH_k$ satisfy a $O(\gamma(k))$-concentration inequality with $\lim_{t\to \infty}\gamma(t)/t=0$, then $X$ is reflexive. 
	\end{prop}
	
	Hamming graphs techniques give us natural access to information about asymptotic models.
	
	\begin{prop}
		\label{prop:H-asymptotic-model}
		If the $X$-valued Lipschitz maps on $\sH_k$ satisfy a $\gamma(k)$-concentration inequality, then for every asymptotic model $A=[\en]$ of $X$ that is generated by a weakly null normalized array, one has that for all $k\ge 1$,
		\begin{equation*}
			\Big\|\sum_{i=1}^k e_i \Big\|_A\le 2\gamma(k).
		\end{equation*} 
		
		In particular, if the $X$-valued Lipschitz maps on $\sH_k$ satisfy an $O(1)$ concentration inequality, then every asymptotic model of $X$ generated by a normalized weakly null array is isomorphic to $\co$.
	\end{prop}
	
	\begin{proof}
		Let $[x]$ be a normalized weakly null array in $X$ that generates an asymptotic model $A=[\en]$. Fixing $k\in\bN$ and $\delta>0$ and passing to an appropriate subarray we may assume that for any $j_k>\dots >j_1\ge k$ and any $a_1,\ldots,a_k$ in $[-1,1]$ we have
		\begin{equation*}
			\label{eq1:H-asymptotic-model}
			\abs{\Big\|\sum_{i=1}^k a_i x^{(i)}_{j_i} \Big\|_X - \Big\|\sum_{i=1}^k a_i e_i \Big\|_A } < \delta.
		\end{equation*}
		Define $h \colon[\bN]^k\to X$ by $h(\nbar) = \sum_{i=1}^k x^{(i)}_{n_i}$ for $\nbar=\{n_1,n_2,\dots,n_k\}\in [\bN]^k$. 
		Then, $h$ is $1$-Lipschitz with respect to the Hamming graph metric and it follows from the concentration property that there exists $\bM\in[\bN]^\omega$ such that for all $\mbar\prec\nbar\in [\bM]^k$ one has
		\begin{equation*}
			\norm{ h(\mbar) - h(\nbar) }_X = \Big\| \sum_{i=1}^k x^{(i)}_{m_i} - \sum_{i=1}^k x^{(i)}_{n_i} \Big\|_X \le \gamma(k).
		\end{equation*}
		Since the array is weakly null, if follows from the weak lower-semi-continuity of the norm that for all $\mbar \in[\bM]^k$, we have $\|\sum_{i=1}^k x^{(i)}_{m_i}\|\le \gamma(k)$. Therefore, if one chooses $\delta$ small enough, we have that for all $k\ge 1$, $\| \sum_{i=1}^k e_i\|_A\le 2\gamma(k)$.
		The last statement follows from the unconditionality of asymptotic models and a classical extreme point argument.
	\end{proof}
	
	A theorem of Freeman, Odell, Sari and Zheng \cite{FOSZ2018} states that if every asymptotic model generated by a weakly null array of a separable Banach space $X$ that does not contain an isomorphic copy of $\ell_1$ is isomorphic to $\co$, then $X$ must be asymptotic-$\co$. This theorem is remarkable as it says that if every weakly null tree of infinite height that is induced by an infinite weakly null array has upper $\ell_\infty$-estimates, then every weakly null tree of finite but arbitrarily large height has upper $\ell_\infty$-estimates. The fact that we can deduce information about arbitrary trees from very structured trees is quite unexpected (at least to the authors!).
	
	\begin{theo}
		\label{thm:BLMS}
		Let $X$ be a Banach space. Then,
		$X$ is reflexive and asymptotic-$\co$ if and only if the $X$-valued Lipschitz maps on the Hamming graphs satisfy a $O(1)$ concentration inequality, i.e.
		there exists a constant $C\ge 1$ such that for all $k\in\bN$ and all $f\colon ([\bN]^k, \dHk)\to X$ Lipschitz, there exists $\bM \in [\bN]^\omega$ so that 
		\begin{equation}
			\sup_{\mbar,\nbar\in[\bM]^k} \norm{ f(\mbar) - f(\nbar)}_X \le C\Lip(f).
		\end{equation}
	\end{theo}
	
	\begin{proof}
		The necessary condition follows from Theorem \ref{thm:H-KR-BLMS}. For the sufficient condition, it follows from Proposition \ref{prop:H-reflexivity} that $X$ must be reflexive and from Proposition \ref{prop:H-asymptotic-model} that every asymptotic model of $X$, generated by a weakly null array, is isomorphic to $\co$. Therefore, thanks to the Odell-Freeman-Sar\'i-Zheng theorem, $X$ must be asymptotic-$\co$.
	\end{proof}
	
	Theorem \ref{thm:BLMS} is one of a few clean results from the Kalton program. Indeed, it provides a purely metric characterization of a class of Banach spaces defined in terms of an asymptotic property. From the coarse stability of $O(1)$-concentration inequalities (Proposition \ref{prop:HFC-invariant}) we can immediately decide the coarse rigidity of the class of reflexive and asymptotic$-\co$ Banach spaces. At the time this book was written, finding a class of metric/Banach spaces that is coarsely rigid was still a rare event.
	
	\begin{coro}[Coarse rigidity of $\REF\cap \ASO$]
		\label{cor:coarse-rigidity-RAS0}
		If a Banach space $X$ coarsely embeds into a Banach space that is reflexive and asymptotic-$\co$, then $X$ is itself reflexive and asymptotic-$\co$
	\end{coro}
	
	The metric characterization in Corollary \ref{cor:coarse-rigidity-RAS0} is in terms of a concentration inequality which is a certain type of a Poincar\'e-type inequality. In the remainder of this section, we explain how using ideas of the previous sections, one can derive a purely metric characterization, in terms of metric spaces preclusion, of the class of asymptotic-$\co$ spaces \emph{within the class of reflexive spaces}. This type of characterizations has applications to the study of the descriptive set-theoretic complexity of classes of Banach spaces. 
	
	Given an arbitrary normalized $1$-suppression unconditional basis $\ebar = (e_j)_{j\in\bN}$ of a Banach space $X$ and an infinite subset $\bM$ of $\bN$, we define for every $k\in\bN$ a map $\sd_{\ebar} \colon [\bM]^k\times[\bM]^k\to[0,\infty)$ such that for every $\mbar=\{m_1,m_2,\dots,m_k\}$ and $\nbar=\{n_1,n_2,\dots,n_k\}$ in $[\bM]^k$ if we let $F := \{j\colon m_j\neq n_j\}$, then
	\begin{equation}
		\sd_{\ebar}(\mbar,\nbar) = \Big\| \sum_{j\in F} e_j\Big\|.
	\end{equation}
	It is elementary to verify that $\sd_{\ebar}$ satisfies the definitions and symmetry axioms of a metric. The triangle inequality follows from the unconditionality condition. Indeed, if $\mbar=\{m_1,\dots,m_k\}$, $\nbar=\{n_1,\dots,n_k\}$ and $\lbar = \{l_1,\dots,l_k\}$, set $F := \{j\colon m_j\neq n_j\}$, $G := \{j\colon m_j\neq l_j\}$ and $H := \{j \colon n_j\neq l_j\}$.
	Since the set $F\subset G\cup H$  we have
	\begin{equation*}
		F = F\cap(G\cup H) = (F\cap G)\cup ((F\setminus G)\cap H).
	\end{equation*}
	It follows from $1$-suppression unconditionality that
	\begin{align*}
		\sd_{\ebar}(\mbar,\nbar) &= \Big\| \sum_{j\in F} e_j\Big\| \le   \Big\|\sum_{j\in F\cap G} e_j\Big\| +  \Big\| \sum_{j\in (F\setminus G)\cap H} e_j \Big\| \\
		&\le   \Big\| \sum_{j\in G} e_j \Big\| +  \Big\| \sum_{j\in H} e_j \Big\|
		= \sd_{\ebar}(\mbar,\lbar) + \sd_{\ebar}(\lbar,\nbar).
	\end{align*}
	
	The metric $\sd_{\ebar}$ is similar to the Hamming metric in the sense that for $\mbar = \{m_1,\dots,m_k\}$ and  $\nbar = \{n_1,\dots,n_k\}$ the distance $\sd_{\ebar}(\mbar,\nbar)$ is determined by the set $F\subset \{1,2,\dots,k\}$ of coordinates $i$ on which $m_i$ and $n_i$ differ. The following important features directly follow from the definition of the metric and classical Banach space theory.
	\begin{lemm}
		\label{lem:Hamming-type-metrics}
		Let $\ebar = (e_j)_{j\in \bN}$ be a normalized 1-suppression unconditional basis of a Banach space.
		\begin{enumerate}[(i)]
			\item If $\ebar = (e_j)_{j\in\bN}$ is the unit vector basis of $\ell_1$, then $\sd_{\ebar}$ is the Hamming distance $\dHk$ on $[\bN]^k$. Hence, for any normalized 1-suppression unconditional basic sequence $\ebar = (e_j)_{j\in\bN}$ and any $\mbar$, $\nbar$ in $[\bN]^k$ we have $\sd_{\ebar}(\mbar,\nbar) \le  \dHk(\mbar,\nbar)$.
			\item For every $k\in\bN$ and every $\bM\in[\bN]^\omega$ we have 
			\begin{equation*}
				\diam([\bM]^k, \sd_{\ebar})= \Big\| \sum_{j=1}^{k} e_j \Big\|.
			\end{equation*}
			In particular, $\lim_{k\to\infty} \diam \big([\bM]^k, \sd_{\ebar}\big) = \infty$ if and only if $\ebar = \en$ is not equivalent to the unit vector basis of $\co$.
		\end{enumerate}
	\end{lemm}
	
	Item $(ii)$ says that if $\bar e = \en$ is not equivalent to the unit vector basis of $\co$, then the sequence of metric spaces $\big([\bN]^k, \sd_{\ebar}\big)_k$ is hereditarily unbounded, in the following sense:
	\begin{equation*}
		\lim_{k\to\infty}\inf_{\bM\in[\bN]^\omega}\diam\big([\bM]^k, \sd_{\ebar}\big) = \lim_{k\to\infty} \Big\| \sum_{i=1}^k e_i \Big\| = \infty.
	\end{equation*}
	The domination of the metric $\sd_{\ebar}$ by the Hamming metric allows us to use the concentration inequalities in Theorem \ref{thm:H-KR-BLMS} to obtain nonembeddability obstructions for the sequence $\big([\bN]^k, \sd_{\ebar}\big)_k$. For instance, assume that $X$ is asymptotic-$\co$ and reflexive and let $\ebar = \en$ be a normalized $1$-suppression unconditional sequence. 
	%such that $\lim_{k\to \infty} \diam\big([\bM]^k,\sd_{\ebar})\big)=\infty$.
	The crucial observation here is that the domination property in Lemma \ref{lem:Hamming-type-metrics} (i), can be equivalently restated by saying that the identity maps from  $([\bN]^k, \dHk)$ to $([\bN]^k, \sd_{\ebar}))$ are $1$-Lipschitz and a straightforward application of Theorem \ref{thm:H-KR-BLMS} shows that there exists $C\in [1,\infty)$ so that for every $1$-suppression unconditional basis $\ebar=\en$, every $k\in\bN$ and every $1$-Lipschitz map $f\colon \big([\bN]^k, \sd_{\ebar}\big)\to X$ there exists $\bM\in[\bN]^\omega$ so that
	\begin{equation}
		\label{eq:Hamming-type-metrics}
		\diam\big( f([\bM]^k)\big) \le  C.
	\end{equation}
	
	It follows now from inequality \eqref{eq:Hamming-type-metrics} and $(ii)$ of Lemma \ref{lem:Hamming-type-metrics} that assuming $\lim_{k\to \infty}\inf_{\bM\in[\bN]^\omega}\diam\big([\bM]^k,\sd_{\ebar}\big)=\infty$ clearly prevents the equi-bi-Lipschitz embeddability of the sequence $\big([\bN]^k, \sd_{\ebar}\big)_{k\in \bN}$ into an asymptotic-$\co$ reflexive Banach space $X$, or in other words implies $\sup_{k\in\bN}\cdist{X}\big([\bN]^k, \sd_{\ebar}\big)=\infty$. We thus proved:
	
	\begin{theo}
		\label{thm:H-type-concentration}
		Let $X$ be a separable asymptotic-$\co$ reflexive Banach space. Then, for every $1$-suppression unconditional
		basis $\ebar = \en$ such that
		$\lim_{k\to \infty}\inf_{\bM\in[\bN]^\omega}\diam\big([\bM]^k,\sd_{\ebar}\big)=\infty$ one has $\sup_{k\in\bN}c_X\big([\bN]^k, \sd_{\ebar}\big)=\infty.$
	\end{theo}
	
	It follows from the definition of $\sd_{\ebar}$, the domination property $\sd_{\ebar}\le \dHk$ and the $1$-suppression unconditionality of asymptotic models that Proposition \ref{prop:H-graph-embedding} can be restated as the following bi-Lipschitz embedding result.
	
	\begin{prop}
		\label{prop:H-type-embedding}
		Let $X$ be a Banach space and $A= [\en]$ be an asymptotic model generated by a normalized weakly null array in $X$. Then, for any $k\in \bN$ and $\vep>0$, the metric space $([\bN]^k, \sd_{\ebar})$ bi-Lipschitzly embeds into $X$ with distortion at most $(2+\vep)$.
	\end{prop}
	
	We are now in a position to obtain the following metric characterization.  
	
	\begin{coro}
		\label{cor:asymp-co-metric-charac}
		Let $X$ be a separable and reflexive Banach space. Then, $X$ is asymptotic-$\co$ if and only for every $1$-suppression unconditional
		basis $\ebar = \en$ such that
		$\lim_{k\to \infty}\inf_{\bM\in[\bN]^\omega}\diam\big([\bM]^k,\sd_{\ebar}\big)=\infty$ one has $\sup_{k\in\bN}c_X\big([\bN]^k, \sd_{\ebar}\big)=\infty.$
	\end{coro}
	
	\begin{proof}
		One implication is Theorem \ref{thm:H-type-concentration}. For the other implication, assume that $X$ is a separable reflexive Banach space that is not asymptotic-$\co$. The Odell-Freeman-Sar\'i-Zheng theorem tells us that there must be an asymptotic model of $X$ that is not equivalent to $\co$ and we conclude using Proposition \ref{prop:H-type-embedding}.
	\end{proof}

	\section[Weighted Hamming graphs and asymptotic models]{Weighted Hamming graphs and upper estimates on asymptotic models}
	\label{sec:weighted-Hamming}
	
	It follows from Proposition \ref{prop:H-asymptotic-model} and the stability of concentration inequalities that upper estimates \emph{for constant coefficients} on asymptotic models are preserved under nonlinear embeddings.
	In this section, the goal is to obtain similar results but for \emph{arbitrary coefficients}. With this goal in mind, we consider as in \cite{KaltonRandrianarivony2008} a weighted version of the Hamming metric, which is defined as follows. For a given set of positive weights $\wbar=(w_1,\dots,w_k)\in (0,\infty]^k$, define for all $\mbar, \nbar \in [\bM]^k$, 
	\begin{equation*}
		\dwHk(\mbar,\nbar) := \sum_{i \colon m_i \neq n_i} w_i.
	\end{equation*}
	
	We will simply denote by $\wHk$ any metric space of the form $([\bM]^k,\dwHk)$ for some infinite subset $\bM$ of $\bN$.
	
	The concentration inequality below is due to Kalton and Randrianarivony \cite{KaltonRandrianarivony2008}. Recall that $\ell_{\bar{\rho}_Y}$ is the Orlicz sequence space associated with the Orlicz function $\bar{\rho}_Y$
	
	\begin{theo}
		Let $X$ be a reflexive Banach space, $k\in \bN$ and $\wbar=(w_1,\dots,w_k)\in (0,\infty]^k$. Then, for all $f\colon ([\bN]^k, \dwHk)\to X$ Lipschitz, there exists $\bM \in [\bN]^\omega$ so that 
		\begin{equation}
			\sup_{\mbar,\nbar\in[\bM]^k} \norm{ f(\mbar) - f(\nbar)}_X \le 2e \norm{\wbar}_{\ell_{\bar{\rho}_X}}\Lip(f).
		\end{equation}
	\end{theo}
	
	\begin{proof}
		Observe first that it follows from Proposition \ref{pro:Orlicz-iterated} that  for any weakly null tree $(x_{\nbar})_{\nbar\in [\bN]^{\le k}\setminus \{\emptyset\} }$ in $X$, there exists $\nbar \in [\bN]^k$ such that  
		$$\Big\|\sum_{i=1}^k x_{n_1,\dots,n_k}\Big\|_X\le 2e\Big\| \big(\norm{x_{n_1,\dots,n_i}}_X\big)_{i=1}^k \Big\|_{ \ell_{\bar{\rho}_X} }.$$ 
		Then, the conclusion  simply follows from Lemma \ref{lem:linearization}, an application of Ramsey's theorem and the fact that for all $1\le i \le k$, $\mbar\in [\bN]^{k}$ and $m_i'\in \bN$ such that $\sigma_{m_i'}(\mbar)\in [\bN]^k$, $\dwHk(\mbar,\sigma_{m_i'}(\mbar))=w_i$.    
	\end{proof}

	The following result was first stated in \cite{Kalton2013b} for spreading models but the extension to asymptotic models is straightforward.  
	
	\begin{theo}
		\label{thm:upper-estimates-asymp-models} 
		Let $Y$ be a reflexive Banach space and assume that $X$ is a Banach space that coarse-Lipschitz embeds into $Y$. Then, there exists a constant $C>0$ such that for any asymptotic model $ A :=[\en]$ generated by a normalized weakly null array in $X$, all $k\in \bN$ and $a_1,\dots, a_k\in \bR$, 
		\begin{equation*}
			\Big\|\sum_{i=1}^k a_i e_i \Big\|_{A} \le C \norm{ (a_i)_{i=1}^k }_{\ell_{\bar{\rho}_Y}}.
		\end{equation*}
	\end{theo}
	
	\begin{proof}
		Let $f\colon X\to Y$ be a coarse-Lipschitz embedding. Without loss of generality one can assume that there are constants $B,D\in(0,\infty)$ such that for all $x_1,x_2\in X$ 
		\begin{equation}
			\label{eq1bis:H-asymptotic-model}
			\frac1D \norm{x_1-x_2}_X -B \le \norm{f(x_1) -f(x_2) }_Y \le D \norm{ x_1 -x_2 }_X +B.
		\end{equation}
		Let $[x]$ be a normalized weakly null array in $X$ that generates an asymptotic model $A=[\en]$. Fixing $k\in\bN$ and $\delta>0$ and passing to an appropriate subarray we may assume that for any $j_k>\dots >j_1\ge k$ and any $a_1,\ldots,a_k$ in $[-1,1]$ we have
		\begin{equation}
			\label{eq2:H-asymptotic-model}
			\abs{\Big\|\sum_{i=1}^k a_i x^{(i)}_{j_i}\Big\|_X - \Big\|\sum_{i=1}^k a_i e_i \Big\|_A } < \delta.
		\end{equation}
		Given $(a_1,\ldots,a_k)\in [-1,1]^k$ and $\lambda>0$, define $h_\lambda \colon[\bN]^k\to Y$ by 
		$$h_\lambda(\nbar) = f\Big(\frac{\lambda}{2} \sum_{i=1}^k a_i x^{(i)}_{n_i}\Big)\ \  \text{for}\ \  \nbar=\{n_1,n_2,\dots,n_k\}\in [\bN]^k.$$ 
		Then, $h_\lambda$ is $(D\lambda + B(\min_i\abs{a_i})^{-1})$-Lipschitz with respect to $\dwHk$ with $\wbar=(\abs{a_1},\dots,\abs{a_k})$ and it follows from the concentration property that there exists $\bM\in[\bN]^\omega$ such that for all $\mbar\prec\nbar\in [\bM]^k$ one has
		\begin{equation*}
			\norm{ h_\lambda(\mbar) - h_\lambda(\nbar) }_X \le
			4e\norm{(a_i)_{i=1}^k }_{\ell_{\bar{\rho}_Y}}(D\lambda + B(\min_i\abs{a_i})^{-1}).
		\end{equation*}
		Hence, using the left-hand side of \eqref{eq1bis:H-asymptotic-model} we deduce that
		\begin{equation*}
			\Big\|\sum_{i=1}^k a_i x^{(i)}_{m_i} - \sum_{i=1}^k a_i x^{(i)}_{n_i}\Big\|\le \frac{8eD}{\lambda}\norm{(a_i)_{i=1}^k }_{\ell_{\bar{\rho}_Y}}(D\lambda + B(\min_i\abs{a_i})^{-1})+\frac{2BD}{\lambda}. 
		\end{equation*}
		If we let $n_1, n_2,\dots,n_k \to \infty$ and $\lambda\to \infty$ we have 
		\begin{equation*}
			\Big\| \sum_{i=1}^k a_i x^{(i)}_{m_i}\Big\|\le  8D^2e\norm{(a_i)_{i=1}^k }_{\ell_{\bar{\rho}_Y}},  
		\end{equation*}
		and the conclusion follows from \eqref{eq2:H-asymptotic-model}.
	\end{proof}

	\section{Notes}
	In view of the above results, we introduce the following definition that has been introduced by A. Fovelle \cite{Fovelle2022}.
	
	\begin{defi} Let $X$ be a Banach space.
		\begin{enumerate}
			\item Let $p\in (1,\infty)$ and $\lambda\ge 1$. We say that $X$ has the \emph{$p$-Hamming Full Concentration Property} with constant $\lambda$, denoted by $\lambda$-HFC$_p$, if for all $f\colon  ([\bN]^k,d^{(k)}_\H)\to X$ Lipschitz, there exists $\M \in [\bN]^\omega$ satisfying $\diam (f([\M]^k))\le \lambda k^{1/p}\Lip(f)$. We say that $X$ has the \emph{$p$-Hamming Full Concentration Property}, denoted by HFC$_p$, if it has $\lambda$-HFC$_p$ for some $\lambda \ge 1$.
			\item We say that $X$ has the \emph{$\infty$-Hamming Full Concentration Property}, denoted by HFC$_\infty$ with constant $\lambda$, if for all $f\colon  ([\bN]^k,d^{(k)}_\H)\to X$ Lipschitz, there exists $\M \in [\bN]^\omega$ satisfying $\diam (f([\M]^k))\le \lambda\Lip(f).$ We say that $X$ has the \emph{$\infty$-Hamming Full Concentration Property}, denoted by HFC$_\infty$, if it has $\lambda$-HFC$_\infty$ for some $\lambda \ge 1$.
		\end{enumerate}
	\end{defi}
	
	\begin{rema}
		With this terminology, we can summarize Theorem \ref{thm:H-KR-BLMS} in one sentence: if a reflexive Banach space belongs to $\textsf{N}_p$, for $p\in (1,\infty]$, then it has HFC$_p$.   
	\end{rema}
	
	All these results leave open the following important problem.
	
	\begin{prob}[Coarse-Lipschitz rigidity of $\langle \REF \& \AUS\rangle$] 
		Assume that a Banach space $X$ coarse-Lipschitz embeds into a reflexive asymptotically uniformly smooth Banach space. Does it follow that $X$ admits an equivalent norm that is asymptotically uniformly smooth?    
	\end{prob}
	
	It is important to mention here that A. Fovelle recently proved in \cite{Fovelle2024} that, for $p\in (1,\infty)$, HFC$_p$ does not imply AUS renormability. She even built separable spaces with HFC$_p$ and arbitrarily high Szlenk index. 
	
	As we have seen, one of the features of HFC$_p$ is that it implies reflexivity. A weaker form of this concentration property can be used to address the nonreflexive setting. 
	
	\begin{defi} 
		\label{def:HFCp}
		Let $X$ be a Banach space. Let $p\in (1,\infty]$ and $\lambda\ge 1$. We say that $X$ has the \emph{$p$-Hamming Concentration Property} with constant $\lambda$, denoted by $\lambda$-HC$_p$, if for all $f\colon  ([\bN]^k,d^{(k)}_\H)\to X$ Lipschitz, there exist $\m,\n \in [\bN]^k$ such that $d^{(k)}_\H(\m,\n)=k$ and $\|f(\m)-f(\n)\|\le \lambda k^{1/p}\Lip(f)$.  We say that $X$ has the \emph{$p$-Hamming Concentration Property}, denoted by HC$_p$, if it has $\lambda$-HC$_p$ for some $\lambda \ge 1$.
	\end{defi}
	
	In \cite{LancienRaja2018}, G. Lancien and M. Raja prove that a quasi-reflexive Banach space which belongs to $\textsf{N}_p$ has HC$_p$. In \cite{Fovelle2022} and \cite{Fovelle2024}, A. Fovelle studies the stability of HC$_p$ and HFC$_p$ under $\ell_p$-sums. This allows her to provide examples of non-quasi-reflexive Banach spaces with HC$_p$ and, as we already mentioned, of Banach spaces with HFC$_p$ and arbitrarily high Szlenk index (in particular not AUS renormable). As it is explained in \cite{BLMS2020} (Section 5), there exists a quasi-reflexive nonreflexive Banach space $X$, which is asymptotic-$\co$. Therefore, $X$ has HC$_\infty$ but not HFC$_\infty$. Since  HC$_\infty$ clearly implies the non-equi-coarse embeddability of the Hamming graphs, this proves that HFC$_\infty$ is not equivalent to the non-equi-coarse embeddability of the Hamming graphs. It leaves open the following. 
	
	\begin{prob}
		Assume that the Hamming graphs do not equi-coarsely embed into a Banach space $X$. Does this imply that $X$ has HC$_\infty$? 
	\end{prob}
	
	The property of not being asymptotic-$\co$ cannot be characterized by the equi-bi-Lipschitz embeddability of $([\bN]^k, \sd_{\ebar})$, $k\in\bN$, for some $\ebar$, where $\ebar$ only comes out of a countable subset of hereditarily $1$-suppresion unconditional basis. Indeed, we have the following proposition.
	
	\begin{prop} 
		\label{P:3.10}
		Let
		$$D\subset \left\{ (d_k)_{k\in\bN}\colon  \begin{matrix}d_k \text{ is a metric on $[\bN]^k$, which is dominated by $\dHk$ and }\\ \limsup_{k\to\infty}\inf_{\bM\in[\bN]^\omega} \diam([\bM]^k,d_k) =\infty\end{matrix}\right\}$$
		be countable.
		Then, there exists a reflexive Banach space $X$,  which is not asymptotic-$\co$, so that for all $ (d_k)_{k\in\bN}\in D$ and
		for all sequences $(\Psi_k)_{k\in\bN}$, where $\Psi_k\colon ([\bN]^k, d_k )\to X$ is $1$-Lipschitz it follows that
		$$\lim_{k\to\infty} \inf_{\bM\in[\bN]^{\omega}} \frac{\diam (\Psi_k([\bM]^k, d_k ))}{\diam ([\bM]^k, d_k)}=0.$$
		In particular, the $\Psi_k$ cannot be equi-bi-Lipschitz embeddings.
	\end{prop}
	
	The space $X$ above is the dual of Schlumprecht space constructed in \cite{Schlumprecht91}.

	\section{Exercises}
	
	\begin{exer}\label{exe:H-graph}
		Let $\M_1$ and $\M_2$ be two infinite subsets of $\bN$. Show that $\sH_k(\M_1)$ and $\sH_k(\M_2)$ are graph-isomorphic.   
	\end{exer}
	
	\begin{exer}
		\label{exe:H-graph-embedding-spreading-models}
		Let $Y$ be an infinite-dimensional Banach space admitting a spreading model $E$, generated by a normalized weakly null sequence. Show that for every $\nu>0$ and $k\ge 1$ there exists a map $h_{k,\nu} \colon \sH_k \to Y$ such that for all $\mbar,\nbar \in \sH_k$,
		\begin{equation*}
			\frac{1}{8(1+\nu)}\varphi_E\left(\dHk(\mbar,\nbar)\right)\le \norm{h_{k,\nu}(\mbar)-h_{k,\nu}(\nbar)}_Y\le \dHk(\mbar,\nbar).
		\end{equation*}   
	\end{exer}
	
	\begin{proof}[Hint]
		Consider a bijection between $\bN\times \bN$ and $\bN$.
	\end{proof}
	
	%\begin{exer}\label{ex:trivialcotype} Let $X$ be a Banach space satisfying the conclusion of Corollary \ref{coro:basisasymptoticc0}. Show that $X$ uniformly contains isomorphic copies of the $\ell_\infty^k$'s. Indication: prove that for any $k\in \bN$, there exists a basic sequence $(x_n)_n$ in $S_X$ such that for some $C\ge 1$, we have that $\big\|\sum_{i=1}^k\eps_ix_{n_i}\big\|\le C$, for all $(\eps_1,\ldots,\eps_k)\in \{-1,1\}^k$.
	%\end{exer}

	\begin{exer}\,
		\begin{enumerate}
			\item Let $p\in (1,\infty)$. Show that HC$_p$ is stable under coarse-Lipschitz embeddings;
			\item Prove that HC$_\infty$ is stable under coarse embeddings.
		\end{enumerate}
	\end{exer}
	
	\begin{exer}
		Let $p\in (1,\infty)$. Show that HC$_p$ is equivalent to the following property (denoted by HIC$_p$ in \cite{Fovelle2022} for $p$ Hamming Interlacing Concentration property):\\
		There exists $\lambda \ge 1$ such that for all $f\colon  ([\bN]^k,d^{(k)}_\H)\to X$ Lipschitz, there exists $\M \in [\bN]^\omega$ with the property that for all $\m=(m_1,\ldots,m_k),\n=(n_1,\ldots,n_k) \in [\M]^k$ such that $m_1<n_1<\ldots<m_k<n_k$, we have $\|f(\m)-f(\n)\|\le \lambda k^{1/p}\Lip(f)$.
		\begin{proof}[Hint]
			Make a repeated and clever use of the Ramsey Theorem, or have a look at \cite{Fovelle2022}. 
		\end{proof} 
	\end{exer}
	
	\begin{exer}
		Let $p\in [1,\infty)$. 
		Show that $\{(\sH_k,\dHk) \}_{k\in \bN}$ equi-coarsely embeds into $(\sum_{k=1}^\infty \ell_p^k(\Tsi^*))_{\Tsi^*}$ when $p>1$ and equi-bi-Lipschitzly when $p=1$.
		What can you say about the compression rates when $p>1$? 
	\end{exer}

	%%%%%%%%%%%%%%%%%%%%%%%%%%%%%%%%%%%%%%%%%%%%%%%%%%%%%%%%%%%%%%%%%%%%%%%%%%%%%%%%%%%%%%

	\chapter{Geometry of the Interlaced Graphs and Applications}
	\label{chapter:interlaced-graphs}
	
	The highlight of this chapter is a result of Kalton about the coarse universality of Banach spaces. Kalton showed that if a Banach space is coarsely universal for separable metric spaces, one of its iterated duals must be nonseparable. Consequently, a reflexive Banach space cannot be coarsely universal, answering a long-standing open problem. To prove this beautiful theorem, Kalton introduced influential new techniques and ideas, namely the investigation of concentration inequalities for Lipschitz maps defined on certain nonlocally finite graphs and taking values into Banach spaces. These ideas were already in full display in Chapters \ref{chapter:Johnson} and \ref{chapter:Hamming} but their origins are found in Kalton's approach to the coarse universality problem that is presented here. To tackle this problem, the geometry of a sequence of graphs introduced by Kalton, the interlaced graphs, plays a predominant role due to its close relationship with the summing basis of $\co$. The basic properties of Kalton's interlaced graphs are presented in Section \ref{subsec:interlaced-graphs}. In Section \ref{sec:stability}, we recall and discuss Krivine-Maurey's isometric notion of stability and a natural bi-Lipschitz version of it called infrasup-stability.
	Kalton's property $Q$ is the focal point of Section \ref{sec:Q} and Section \ref{sec:coarse-universality}. It is shown that property $Q$ is an obstruction to the equi-coarse embeddability of the interlaced graphs and that it coincides with the notion of infrasup-stability. Application of property $Q$ to the coarse universality problem is studied in Section \ref{sec:coarse-universality}, where the fact that reflexive spaces have property $Q$ is proved. Section \ref{sec:Q-applications} presents applications of property $Q$ to coarse rigidity and metric characterization problems. In Section \ref{sec:Q-Szlenk}, connections between property $Q$ and the Szlenk index are established. And finally, in Section \ref{sec:Q-omega-1} an uncountable version of property $Q$ is described with a few applications to coarse embeddability of nonseparable spaces into $\ell_\infty$.
	
	\section{Kalton's interlaced graphs}
	\label{sec:intelaced-graphs}
	
	Let $s_n:=\sum_{i=1}^n e_i$ be the $n$-th unit vector of $\sbar:=(s_n)_{n=1}^\infty$ the summing basis of $\co$. This basis induces a metric on $[\bN]^{<\omega}$ via the formula 
	\begin{equation}
		\label{eq:interlaced-summing}
		\sd_{\sbar}(\mbar,\nbar):= \Big\|\sum_{i=1}^{\abs{\mbar}} s_{m_i} - \sum_{i=1}^{\abs{\nbar}} s_{n_i}\Big\|_\infty.
	\end{equation}
	Kalton observed that this metric has a graph-theoretical interpretation that is relatively simple to describe when one considers the metric on the subset $[\bN]^k$. 
	
	Following Kalton \cite{Kalton2007}, for $\bM \in [\bN]^\omega$, we equip $[\bM]^k$ with a graph structure by declaring $\mbar\neq\nbar\in [\bM]^k$ adjacent if and only if they interlace, i.e.
	$$
	n_1 \le m_1 \le n_2 \le  \dots \le n_k \le m_k\ \  \text{or}\ \  m_1 \le n_1 \le m_2\le   \ldots \le m_k \le n_k.
	$$
	We denote by $\Ik(\bM)$, or simply $\Ik$, the interlacing graph on $[\bM]^k$.
	The shortest path distance on $\Ik$ is denoted by $\dIk$ and called the interlaced (or interlacing) metric. 
	
	\begin{rema}
		We do not make any reference to $\bM$ in our notation $\dIk$, because the distance $\dIk$ is independent of the set $\bM$. By this, we mean that if $\bM$ and $\bL$ are two infinite subsets of $\bN$ and $\mbar,\nbar$ both belong to $[\bM]^k$ and $[\bL]^k$, then the shortest paths from $\mbar$ to $\nbar$ in $[\bM]^k$ and in $[\bL]^k$ have the same lengths. In particular, $[\bM]^k$ is a metric subspace of $[\bL]^k$ whenever $\bM \in [\bL]^\omega$.
	\end{rema}
	
	In the next proposition, we give two closed formulas for the interlaced metric that were given in \cite{LPP2020}. These formulas come in handy when providing rigorous arguments to back up statements, as in the remark above, which are intuitively clear but need a justification for pairs of vertices that are not adjacent.
	
	\begin{prop}
		\label{pro:interlaced-closed-formulas}
		Let $\bM \in [\bN]^\omega$. Then, for all $\nbar,\mbar \in [\bM]^k$,
		\begin{align}
			\label{eq1:interlaced-closed-formulas}    \dIk(\mbar,\nbar) & = \max \{\big| \abs{S\cap \mbar} - \abs{S\cap \nbar} \big| \colon S \text{ interval of } \bM \} \\
			\label{eq2:interlaced-closed-formulas}                    & = \ds\max_{i\ge 0} (\abs{\llbracket 1,i \rrbracket \cap \mbar} - \abs{\llbracket 1,i \rrbracket \cap \nbar}) - \ds\min_{i \ge 0} (\abs{\llbracket 1,i \rrbracket \cap\mbar} - \abs{\llbracket 1,i \rrbracket \cap \nbar}).
		\end{align}
	\end{prop}
	
	\begin{proof}
		For all $\nbar,\mbar \in [\bM]^k$, let $d_k(\mbar,\nbar):=\sup \{\big| \abs{S\cap \mbar} - \abs{S\cap \nbar} \big| \colon S \text{ interval of } \bM \} $. It is easily seen that $d_k$ is a metric on $[\bM]^k$ (in fact, it is the restriction of a metric on $[\bM]^{<\omega}$, see Exercise \ref{exer:interlaced-closed-formula}).
		Let us first show that equality \eqref{eq2:interlaced-closed-formulas} holds, i.e. 
		$$ d_k(\mbar,\nbar) = \ds\max_{i\ge 0} F_{\mbar,\nbar}(i) - \ds\min_{i\ge 0} F_{\mbar,\nbar}(i),$$
		where $F_{\mbar,\nbar}(i):=\sum_{j=1}^i \big(\car_{\mbar}(j)-\car_{\nbar}(j)\big)$ for $ i\ge 1$ and $F_{\mbar,\nbar}(0):= 0$.
		Indeed, for any interval $S=\llbracket a,b \rrbracket\subset \bM$,
		$$
		\abs{S\cap \mbar}-\abs{S\cap\nbar} = \sum_{j=a}^b\car_{\mbar}(j)- \sum_{j=a}^b \car_{\nbar}(j) = F_{\mbar,\nbar}(b) - F_{\mbar,\nbar}(a-1).
		$$
		In particular, for any interval $S$ of $\bM$, 
		$$\big| \abs{S\cap \mbar} - \abs{S\cap\nbar} \big|\le  \max_{i\ge 0} F_{\mbar,\nbar}(i) - \min_{i\ge 0} F_{\mbar,\nbar}(i).$$
		For the reverse inequality, observe that for all $i\ge 0$, $F_{\mbar,\nbar}(i)\in \llbracket -\abs{\nbar},  \abs{\mbar}\rrbracket$ and let $a\le b$ such that $$\{F_{\mbar,\nbar}(a), F_{\mbar,\nbar}(b)\}=\{\max_{i\ge 0} F_{\mbar,\nbar}(i),\min_{i\ge 0} F_{\mbar,\nbar}(i)\}.$$ 
		If $a=b$, there is nothing to prove so assume that $a<b$. Considering $S=\llbracket a+1,b \rrbracket\subset \bM$ it is immediate that 
		$$\big|\abs{S\cap \mbar} - \abs{S\cap\nbar}\big| = \big| F_{\mbar,\nbar}(b) - F_{\mbar,\nbar}(a) \big| = \max_{i\ge 0} F_{\mbar,\nbar}(i) - \min_{i\ge 0} F_{\mbar,\nbar}(i)$$
		which finishes the proof of equality \eqref{eq2:interlaced-closed-formulas}.
		
		To prove equality \eqref{eq1:interlaced-closed-formulas}, observe first that $\dIk(\mbar,\nbar)=1$ if and only if $\mbar$ and $\nbar$ interlace if and only if $d_k(\mbar,\nbar)=1$ and thus, by applying the triangle inequality for $d_k$, we deduce that $d_k\le \dIk$. To show the reverse inequality, it is enough to prove that given $\mbar\neq \nbar$ in $[\M]^k$, there is a path of length $d_k(\mbar,\nbar)$ in $\sI_k(\bM)$ that connects $\mbar$ and $\nbar$. Obviously, we only need to consider a pair of points such that $d_k(\mbar,\nbar)\ge 2$ and if we can show that there is $\lbar\in [\M]^k\setminus\set{\mbar,\nbar}$ such that $d_k(\mbar,\nbar)=d_k(\mbar,\lbar) + d_k(\lbar,\nbar)$, then the conclusion follows from an elementary induction. 
		
		Since $F_{\mbar,\nbar}(i)= -F_{\nbar,\mbar}(i)$, we can assume without loss of generality that $\max_{i\ge 0} F_{\mbar,\nbar}(i)>0$. Since $d_k(\mbar,\nbar)>0$, the sets $\arg\max(F_{\mbar,\nbar})$ and $\arg\min(F_{\mbar,\nbar})$ are necessarily disjoint. We select inductively $\{a_1<\dots<a_p\}\subset \arg\max(F_{\mbar,\nbar})$ and $\{b_1<\dots<b_q\}\subset \arg\min(F_{\mbar,\nbar})$ (with $p\ge 1$ and $q\ge 0$) with the property that
		
		\begin{itemize}
			\item $a_1=\min \arg\max(F_{\mbar,\nbar})$,
			\item For $i\ge 1$,\ $b_i=\min \left(\{n>a_i\} \cap \arg\min(F_{\mbar,\nbar})\right)$, if this is not empty.
			\item $a_{i+1}=\min\left( \{n>b_i\} \cap \arg\max(F_{\mbar,\nbar})\right)$, if this set is not empty.
		\end{itemize}
		Note that since $\max_{i\ge 0} F_{\mbar,\nbar}(i)>0$, then $\{a_1,\dots,a_p\} \subset \mbar \setminus \nbar$ and $\{b_1,\dots,b_q\} \subset \nbar\setminus \mbar$.
		Note also that either $p=q$ or $p=q+1$.
		In the latter case, we define $b_p:=r$ for some $r$ such that $r>a_p$ and $F_{\mbar,\nbar}(r-1)>F_{\mbar,\nbar}(r)$.
		Such $r$ must exist since $\max_{i\ge 0} F_{\mbar,\nbar}(i)>0$ and $F_{\mbar,\nbar}(\max\{n_k,m_k\})=0$ and we can take $r \in \nbar \setminus \mbar$.
		We will set
		$$
		\lbar=\mbar \cup \{b_1,\dots,b_p\} \setminus \{a_1,\dots,a_p\}.
		$$
		It is clear that $\lbar \in [\M]^k$.
		We also have $\max_{i\ge 0} F_{\lbar,\nbar}(i) = \max_{i\ge 0} F_{\mbar,\nbar}(i)-1$ and $\min_{i\ge 0} F_{\lbar,\nbar}(i) = \min_{i\ge 0} F_{\mbar,\nbar}(i)$.
		Indeed, the point $\lbar$ is constructed in such a way that when $F_{\mbar,\nbar}$ attains its maximum for the first time  (going from the left), $F_{\lbar,\nbar}$ is reduced by $1$ after that point and stays reduced by $1$ until the next time the minimum of $F_{\mbar,\nbar}$ is attained (or until the point $r$) where this reduction is corrected back; and so on. But since the correction occurs at a point where the minimum is attained and the difference between the minimum and the maximum is at least $2$, the maximum is effectively decreased by $1$ and the minimum stays unchanged. Thus, $d_k(\lbar,\nbar)=d_k(\mbar,\nbar)-1$.
		Also, observe that by construction $F_{\lbar,\nbar}\le F_{\mbar,\nbar}$ and since the sets $\set{a_1,\dots,a_p}$ and $\set{b_1,\dots,b_p}$ are interlaced we also have $F_{\mbar,\nbar}-1\le F_{\lbar,\nbar}$.
		Therefore, since $F_{\mbar,\nbar} = F_{\mbar,\lbar} + F_{\lbar,\nbar}$, we have that $0\le F_{\mbar,\lbar}\le 1$ and hence $d_k(\mbar,\lbar)=1$, since by construction $\mbar \neq \lbar$. Finally, we have proved that $d_k(\mbar,\nbar)=d_k(\mbar,\lbar) + d_k(\lbar,\nbar)$.
	\end{proof}

	\begin{prop}
		\label{pro:interlaced-co-embedding} 
		For all $k\in \bN$ and $\mbar,\nbar\in [\bN]^k$,
		\begin{equation}
			\label{eq:interlaced-co-embedding}
			\frac12 \dIk(\mbar,\nbar) \le \sd_{\sbar}(\mbar,\nbar)\le \dIk(\mbar,\nbar)
		\end{equation}
		%Moreover, $\cdist{\co}(\sI_k) \le 2$ for all $k\in \bN$.
		%define  $h_k\colon([\bN]^k, \dIk)\to c_0$ by $h_k(\nbar)=\sum_{i=1}^k s_{n_i}$. Then, for all $\mbar,\nbar \in [\bN]^k$ we have 
		%\begin{equation}
		%\frac12 \dIk(\mbar,\nbar)\le \norm{ h_k(\mbar) - h_k(\nbar)}_\infty\le \dIk(\mbar,\nbar).
		%\end{equation}
	\end{prop}
	
	\begin{proof}
		It is straightforward to verify that if $\mbar$ and $\nbar$ interlace, then $\sd_{\sbar}(\mbar,\nbar)=1$, from which is follows immediately that $\sd_{\sbar}\le \dIk$. For the other inequality, observe that $k-\sum_{j=1}^i \car_{\mbar}(j) = e^*_{i+1}(\sum_{j=1}^k s_{m_j})$ for all $i\ge 0$ and $\mbar\in [\bN]^k$ and thus 
		$$\sum_{j=1}^i \big(\car_{\mbar}(j)-\car_{\nbar}(j)\big) = e^*_{i+1}\Big(\sum_{j=1}^k s_{n_j}\Big) - e^*_{i+1}\Big(\sum_{j=1}^k s_{m_j}\Big) = e^*_{i+1}(h_k(\nbar) - h_k(\mbar)),$$
		where  $h_k\colon([\bN]^k, \dIk)\to c_0$ is defined by by $h_k(\nbar)=\sum_{i=1}^k s_{n_i}$. Obviously, $e^*_{0}(h_k(\nbar) - h_k(\mbar))=0$. It then follows from Proposition \ref{pro:interlaced-closed-formulas} that for all $\mbar,\nbar \in [\bN]^k$
		$$\dIk(\mbar,\nbar) = \max_{i\in \bN} e^*_{i}(h_k(\mbar) - h_k(\nbar)) - \min_{i\in\bN} e^*_{i}(h_k(\mbar) - h_k(\nbar)).$$
		Since for all $\mbar,\nbar \in [\bN]^k$, $e^*_{\max\{m_k,n_k\}+1}(h_k(\mbar) - h_k(\nbar)) = 0$ we have 
		$$\min_{i\in \bN} e^*_{i}(h_k(\mbar) - h_k(\nbar))\le 0\le \max_{i\in\bN} e^*_{i}(h_k(\mbar) - h_k(\nbar)).$$ 
		The lower bound in \eqref{eq:interlaced-co-embedding} easily follows.
	\end{proof}
	
	\begin{rema}
		The map $h_k$ provides a concrete embedding from $\sI_k$  into $\co$ with distortion at most $2$. It follows that $\cdist{\co}(\sI_k) \le 2$, which was also insured by the sharp version of Aharoni's theorem (Theorem \ref{thm:Aharoni}). In fact $\cdist{\co}(\sI_k) = 2$ (see Exercise \ref{exer:I_kdistortion}).
	\end{rema}
	
	%This upper bound on the bi-Lipschitz distortion is easily seen to be tight by considering 

	\begin{rema}
		If in \eqref{eq:interlaced-summing} we replace the sup-norm with the bimonotone version of the summing norm, then the interlaced metric coincides exactly with this metric and the interlaced graphs admit isometric embeddings into $\co$ equipped with this equivalent norm (cf. Exercise \ref{exer:interlaced-isometric}). 
	\end{rema}
	
	%\subsection{Interlaced graphs and the summing basis of $\co$}
	\label{subsec:interlaced-graphs}
	
	%\subsection{Embeddings of the interlaced graphs and spreading models}
	\label{subsec:interlaced-graphs-embeddings}
	
	\section{Krivine-Maurey stability}
	\label{sec:stability}
	
	In this section, we make a (not so) quick detour to study the fascinating notion of stability introduced by Krivine and Maurey in \cite{KrivineMaurey1981}. The results gathered in this section will be crucial later for the geometry of the interlaced graphs and their applications to rigidity results. 
	
	We have already made a brief encounter with the notion of stability in Section \ref{sec:stable-into-reflexive}, where a connection between stability and reflexivity was established via nonlinear embeddings. Recall that a metric space $(M,d)$ is said to \emph{stable} if for every nonprincipal ultrafilters $\cU$ and $\cV$ over $\bN$ and every bounded sequences $(x_n)_{n=1}^\infty$ and $(y_n)_{n=1}^\infty$ in $M$, we have
	\begin{equation}
		\label{eq:stability}
		\lim_{n\in \cal V}\lim_{m \in \cal U} d(x_n,y_m)=\lim_{m\in \cal U}\lim_{n \in \cal V} d(x_n,y_m).
	\end{equation}
	
	A Banach space is stable if its canonical metric induced by its norm is stable. Stability is an isometric property and is inherited by subsets. Despite condition \eqref{eq:stability} seeming quite restrictive, the class of stable metrics is rather rich and contains many interesting metrics. 
	
	\begin{exam}
		\label{exa:trivial_stable}
		Proper metric spaces, and in particular finite, compact, bounded geometry, or locally finite metric spaces, are stable.
	\end{exam}
	For all the metric spaces in Example \ref{exa:trivial_stable}  the closed balls are either finite or compact and if $\xn$, $\yn$ are bounded sequences, then given any nonprincipal ultrafilters $\cU,\cV$ on $\bN$, there exist $x$ and $y$ such that $\lim_{n\in\cU}x_n=x$ and $\lim_{n\in\cV}y_n=y$. By continuity properties of the distance function we thus have, 
	\begin{equation*}
		\lim_{m\in\cU}\lim_{n\in\cV} d(x_m,y_n) = d(x,y) = \lim_{n\in\cV}\lim_{m\in\cU}d(x_m,y_n).
	\end{equation*}
	In particular, finite-dimensional Banach spaces, finitely generated groups equipped with their canonical word metric and compactly generated groups equipped with their canonical proper metric are stable metric spaces. By the classical Riesz Theorem,  a Banach space is proper if and only if it is finite-dimensional. 
	
	We now describe infinite-dimensional Banach spaces that are stable, since it will provide us with examples of stable metric spaces that do not belong to the list of (trivially) stable spaces from Example \ref{exa:trivial_stable}.
	\begin{exam}
		\label{exa:Hilbert-stable}
		Hilbert space is stable.
	\end{exam}
	It is fairly easy to show that Hilbert space is stable using the classical representation of Hilbert norm in terms of the scalar product. Assume that $\xn$ and $\yn$ are two bounded sequences in Hilbert space. For any two ultrafilters $\cU$ and $\cV$, since Hilbert space is reflexive, the two sequences are weakly convergent along $\cU$ and $\cV$ to, say, $x$ and $y$, respectively. Let $a :=\lim_{n\in\cU}\norm{x_n}_2^2$ and $b :=\lim_{n\in\cV}\norm{y_n}_2^2$. Then, 
	\begin{align*}
		\lim_{m\in\cU}\lim_{n\in\cV}\norm{x_m-y_n}_2^2&=\lim_{m\in\cU}\lim_{n\in\cV}\left(\norm{x_m}_2^2+\norm{y_n}_2^2-2\langle x_m,y_n\rangle\right)\\
		&=\lim_{m\in\cU}(\|x_m\|_2^2+b-2\langle x_m,y\rangle\\
		&=a+b-2\langle x,y\rangle=\lim_{n\in\cV}(\|y_n\|_2^2+a-2\langle x,y_n\rangle\\
		&=\lim_{n\in\cV}\lim_{m\in\cU}\left(\norm{x_m}_2^2+\norm{y_n}_2^2-2\langle x_m,y_n\rangle\right)\\
		&=\lim_{n\in\cV}\lim_{m\in\cU}\norm{x_m-y_n}.
	\end{align*}
	This completes the proof.
	
	\begin{exam}
		\label{exa:lp-stable}
		The Banach space $\ell_p$, for $p\in[1,\infty)$, is stable.
	\end{exam}
	
	There are several different approaches to prove the statement in Example \ref{exa:lp-stable}. It is an elementary, but somewhat tedious, exercise to show that for $p\in[1,\infty)$, the $\ell_p$-sum of stable Banach spaces is stable (see Exercise \ref{exe:lp-sum-of-stable}) and one can then argue that $\ell_p$ is by definition the $\ell_p$-sum of countably many copies of $(\bR,\abs{\cdot})$, which is a proper and hence stable, space. This argument builds on the following elementary lemma whose proof is left as Exercise \ref{exe:asymp-lp}. 
	
	\begin{lemm}
		\label{lem:asymp-lp}
		Let $\xn$ be a bounded sequence in $\ell_p$ and $\cU$ a nonprincipal
		ultrafilter on $\bN$. Suppose that $\lim_{n\in\cU}e_i^*(x_n)=0$ for all $i\in\bN$, i.e. $\xn$ converges coordinatewise to $0$ with respect to $\cU$. Then, for every $z\in\ell_p$,
		\begin{equation}
			\label{eq:asymp-lp}
			\lim_{n\in\cU}\norm{z+x_n}_p^p=\norm{z}_p^p+\lim_{n\in\cU}\norm{x_n}_p^p.
		\end{equation}
	\end{lemm}
	
	We indicate the idea of its proof. The conclusion of the lemma clearly holds if $z$ and the $x_n$ have disjoint supports. Since $\xn$ converges coordinate-wise to $0$ by reaching far out in the sequence, one can select $x_n$ such that the essential contributions to the norm of $x_n$ and of $z$ are supported on essentially disjoint supports, and a classical approximation and truncation argument gives the conclusion. 
	
	\medskip
	Then, for every bounded sequence $\yn$ and every nonprincipal ultrafilter $\cU$ on $\bN$, there exist $y\in \ell_p$ and $\mu\ge 0$ so that if $z\in\ell_p$, then $\lim_{n\in\cU}\norm{z-y_n}_p^p=\norm{z-y}^p+\mu^p$. Indeed, by weak or weak$^*$ compactness in $\ell_p$, there exists $y\in\ell_p$ such that for all $i\in\bN$ $\lim_{n\in\cU}e_i^*(y_n)=e_i^*(y)$ and the conclusion follows from Lemma \ref{lem:asymp-lp} with $\xn = (y-y_n)_{n=1}^\infty$ and $\mu=\lim_{n\in\cU}\norm{y-y_n}$. We are now in a position to prove that $\ell_p$ is stable. For a pair of bounded sequences $\xn$, $\yn$ and a pair of nonprincipal ultrafilters $\cU,\cV$ on $\bN$, there exist $x,y\in \ell_p$, $\mu,\nu\ge 0$ such that for all $z\in \ell_p$ the following equalities hold: 
	$$\lim_{n\in\cU}\norm{z-x_n}_p^p=\norm{z-x}_p^p+\mu^p$$
	and
	$$\lim_{n\in\cV}\norm{z-y_n}_p^p=\norm{z-y}_p^p+\nu^p.$$
	It follows that $$\lim_{m\in\cV}\lim_{n\in\cU}\norm{x_m-y_n}_p^p=\norm{x-y}_p^p+\mu^p+\nu^p=\lim_{n\in\cU}\lim_{m\in\cV}\norm{x_m-y_n}_p^p.$$
	
	\begin{exam}
		\label{exam:co}
		The Banach space $\co$ and more generally, any Banach space containing an isomorphic copy of $\co$, is not stable.
	\end{exam}
	
	The fact that $\co$ is not stable can easily be checked directly by considering the sequences $(-e_n)_{n\in\bN}$ and $(s_n)_{n\in \bN}$, where $\en$ and $(s_n)_{n\in \bN}$ are the canonical basis and the summing basis, respectively, of $\co$. Indeed, $\lim_{m\in\cV}\lim_{n\in\cU}\norm{e_n+s_m}_\infty=1$ while $\lim_{n\in\cU}\lim_{m\in\cV}\norm{e_n+s_m}_\infty=2$. In order to show the second part of Example \ref{exam:co} we need to invoke James' $\co$-distortion theorem \cite{James64} which says that if a Banach space $X$ contains an isomorphic copy of $\co$, then $X$ will contain for every $\epsilon>0$ a subspace that is $(1+\epsilon)$-isomorphic to $\co$. More examples of spaces that are not stable will be given shortly.
	
	The stability of the function space $L_p[0,1]$ is much more difficult to obtain than the stability of the sequence space $\ell_p$ since the validity of equality \eqref{eq:asymp-lp} for the sequence space does not hold in the function space case, and thus another argument is needed.
	
	\begin{exam}
		\label{exam:Lp-stable}
		The Banach space $L_p[0,1]$, for $p\in[1,\infty)$, is stable.
	\end{exam}
	
	It is clear that a snowflaking of a metric space is stable if and only if the original metric on the space is stable. For $p\in[1,2)$, $L_p[0,1]$ will be stable since it is well known that the $\frac{p}{2}$-snowflaking of $L_p[0,1]$ embeds isometrically into the Hilbert space $L_2[0,1]$. This argument fails for $p>2$ and the stability of $L_p[0,1]$ for $p>2$ is more difficult to prove and requires a fine understanding of the relationship between stability and reflexivity. Krivine and Maurey gave in \cite{KrivineMaurey1981} a representation of the norm of a stable Banach space which provides a direct relationship between stability and reflexivity. 
	
	\begin{theo}
		\label{thm:stable-representation}
		Let $X$ be a Banach space and fix $p\in[1,\infty)$. Then, $X$ is stable if and only if there exist a reflexive Banach space $Y$, a dense subset $B$ of the unit ball of $X$ and maps $g\colon B\to Y$, $h\colon B\to Y^*$ so that for all $x,y\in B$ we have $\norm{x-y}^p=\langle g(x),h(y)\rangle$ where $\langle \cdot,\cdot\rangle$ is the duality product between $Y$ and $Y^*$.
	\end{theo}
	
	The theorem below, which includes Example \ref{exam:Lp-stable}, follows from a similar but more delicate representation of stable norms that was also proved in \cite{KrivineMaurey1981}. For a Banach space $X$ we denote by $L_p(\Omega,\cB, \mu;X)$, or simply $L_p(\Omega;X)$,  the Banach space of Bochner equivalence classes of Bochner $p$-integrable and $X$-valued functions defined on the measured space $(\Omega, \cB, \mu)$. 
	
	\begin{theo}
		Let $p\in[1,\infty)$. If $X$ is a stable Banach space, then $L_p([0,1];X)$ is also stable.
	\end{theo}  
	
	Since in \eqref{eq:stability} we have equality, it is possible to extend the stability property to bounded maps on product spaces, as long as we are careful about which permutations of the coordinates are allowed. The following crucial observation is inspired by \cite[Th\'eor\`eme page 273]{KrivineMaurey1981} and \cite[proof of Proposition 5.1]{Raynaud83} (cf. \cite[Lemma 9.19]{BenyaminiLindenstrauss2000} and \cite[Lemma 6.17]{Ostrovskii_book13}.
	
	\begin{theo}
		\label{thm:stability}
		The following assertions are equivalent:
		\begin{enumerate}[(i)]
			
			\item $(M,d)$ is a stable metric space.
			
			\item For every $k\ge 2$, every $1\le \ell < k$, every permutation $\pi\colon [k]:=\{1,2,\dots,k\}\to [k]$ which preserves the order on $\{1,\dots,\ell\}$ and $\{\ell+1,\dots,k\}$, all nonprincipal ultrafilters $\cU_1, \dots ,\cU_k$ on $\bN$ and all bounded maps $f\colon \bN^\ell\to M$ and $g\colon \bN^{k-\ell}\to M$, the following identity holds:
			\begin{align}\label{eq:Rstab1}
				\nonumber \lim_{n_1\in\cU_1}\cdots\lim_{n_k\in\cU_k}& d\big( f(n_1,\dots,n_\ell),g(n_{\ell+1},\dots,n_k) \big) =\\
				&\lim_{n_{\pi^{-1}(1)}\in\cU_{\pi^{-1}(1)}}\cdots\lim_{n_{\pi^{-1}(k)}\in\cU_{\pi^{-1}(k)}}d\big( f(n_{1}, \dots,n_{\ell}),g(n_{\ell+1},\dots,n_{k}) \big).
			\end{align}
			
			\item For every $k\ge 2$, every $1\le \ell < k$, every permutation $\pi\colon [k]\to [k]$ which preserves the order on $\{1,\dots,\ell\}$ and $\{\ell+1,\dots,k\}$, every nonprincipal ultrafilter $\cU$ on $\bN$ and all bounded maps $f\colon \bN^\ell\to M$ and $g\colon \bN^{k-\ell}\to M$, the following identity holds:
			\begin{align}\label{eq:Rstab}
				\nonumber \lim_{n_1\in\cU}\cdots\lim_{n_k\in\cU}&d\big( f(n_1,\dots,n_\ell),g(n_{\ell+1},\dots,n_k) \big)=\\
				&\lim_{n_1\in\cU}\cdots\lim_{n_k\in\cU}d\big( f(n_{\pi(1)}, \dots,n_{\pi(\ell)}),g(n_{\pi(\ell+1)},\dots,n_{\pi(k)}) \big).
			\end{align}
			
			\item For every $k\ge 1$, every $1\le \ell < k$, every permutation $\pi\colon [k]\to [k]$ which preserves the order on $\{1,\dots,\ell\}$ and $\{\ell+1,\dots,k\}$ and all bounded maps $f\colon [\bN]^{\ell}\to M$ and $g\colon [\bN]^{k-\ell}\to M$ we have for every infinite subset $\bM$ of $\bN$,
			\begin{align}\label{eq:1-upper-stable}
				\nonumber \inf_{\nbar\in[\bM]^k}d(f(n_1,\dots,n_\ell),&g(n_{\ell+1},\dots,n_k))\le \\
				&\sup_{\nbar\in [\bM]^k}d(f(n_{\pi(1)},\dots,n_{\pi(\ell)}),g(n_{\pi(\ell+1)},\dots,n_{\pi(k)})).
			\end{align}
			
		\end{enumerate}
	\end{theo}

	\begin{proof}
		The proof of $(i) \Rightarrow (ii)$ is by induction on $k\ge 2$. If $k=2$ the only nontrivial permutation to consider is $\pi:= \left(\begin{smallmatrix}1 & 2 \\ 2 & 1\end{smallmatrix}
		\right)=\pi^{-1}.$ The identity 
		\begin{eqnarray*}
			\lim_{n_1\in\cU_1}\lim_{n_2\in\cU_2}d\big( f(n_1),g(n_2) \big) = \lim_{n_2\in\cU_2}\lim_{n_1\in\cU_1}d\big( f(n_1),g(n_2) \big)
		\end{eqnarray*}
		follows from the definition of stability.
		
		Assume now that $\eqref{eq:Rstab1}$ holds for some $k\ge 2$. Let $1\le \ell < k+1$ and $\pi\colon [k+1]\to [k+1]$ be a permutation which preserves the order on $\{1,\dots,\ell\}$ and $\{\ell+1,\dots, k+1\}$ We think of the inverse permutation $\pi^{-1}$ as a (reverse) ordering according to which the limits should be taken: first along the variable originally associated to $\cU_{\pi^{-1}(k-1)}$, then along the variable originally associated to $\cU_{\pi^{-1}(k)}$ and so on... The only requirement on the permutation $\pi$ means that, in the ordering given by $\pi^{-1}$, we can choose to alternate between the first $\ell$ original limits and the last original $k+1-\ell$ limits as long as we do not swap the order of the first $\ell$ original limits nor swap the order of the the last original $k+1-\ell$ limits. More formally, $[1\le \pi^{-1}(i)< \pi^{-1}(j) \le \ell ] \Rightarrow [i<j] \land [\ell+1\le \pi^{-1}(i)< \pi^{-1}(j) \le k+1] \Rightarrow [i<j]$. 
		
		We first completely detail the case when there is a pair of two consecutive limits amongst the first $\ell$ original limits which are still consecutive in reordering given by $\pi^{-1}$, i.e. if $2\le \pi^{-1}(i)+1 = \pi^{-1}(i+1) \le \ell$ for some $1\le i \le k$, we consider $\cU_{\pi^{-1}(i)} \otimes \cU_{\pi^{-1}(i+1)}$ the ultrafilter tensor product (or Fubini product) defined by 
		\begin{align*}
			&\cU_{\pi^{-1}(i)} \otimes \cU_{\pi^{-1}(i+1)} :=\\ 
			&\Big\{ A \subseteq \bN^2 \colon \big\{ m_i\in \bN \colon \{ m_{i+1}\in \bN \colon (m_i, m_{i+1}) \in A\} \in \cU_{\pi^{-1}(i+1)} \big\} \in \cU_{\pi^{-1}(i)} \Big\}.
		\end{align*}
		As $\cU_{\pi^{-1}(i)}$ and $\cU_{\pi^{-1}(i+1)}$ are nonprincipal ultrafilters, $\cU_{\pi^{-1}(i)}\otimes \cU_{\pi^{-1}(i+1)}$ is a nonprincipal ultrafilter on $\bN^2$ and moreover for any bounded function $h$ on $\bN^2$ it holds 
		$$\lim_{m_i\in\cU_{\pi^{-1}(i)}}\lim_{m_{i+1}\in\cU_{\pi^{-1}(i+1)}} h(m_i,m_{i+1}) = \lim_{(m_i,m_{i+1})\in \cU_{\pi^{-1}(i)}\otimes \cU_{\pi^{-1}(i+1)}} h(m_i,m_{i+1}).$$
		Now if one picks a bijection $\varphi \colon \bN^2 \to \bN$ and write $\varphi^{-1}=(\varphi_1^{-1},\varphi_2^{-1})$, then the pushforward-ultrafilter $\varphi_*(\cU_{\pi^{-1}(i)}\otimes \cU_{\pi^{-1}(i+1)})$ is a nonprincipal ultrafilter on $\bN$ satisfying 
		$$\lim_{m_i\in\cU_{\pi^{-1}(i)}}\lim_{m_{i+1}\in\cU_{\pi^{-1}(i+1)}} h(m_i,m_{i+1}) =  \lim_{n \in \varphi_*\big( \cU_{\pi^{-1}(i)}\otimes \cU_{\pi^{-1}(i+1)} \big)} h(\varphi_1^{-1}(n),\varphi_2^{-1}(n)).$$
		Based on this discussion if $\tilde{f} \colon \bN^{\ell-1} \to M$ is given by $$\tilde{f}(m_1, \dots, m_i , \dots, m_{\ell-1}) := f(m_1, \dots, m_{i-1}, \varphi_1^{-1}(m_i),\varphi_2^{-1}(m_i), m_{i+1},\dots, m_{\ell-1}),$$
		and
		$$\tilde{\cU}_{r}:= \begin{cases}
			\cU_r  \hskip 3cm\text{ if }1\le r \le \pi^{-1}(i)-1,\\
			\varphi_*(\cU_{\pi^{-1}(i)}\otimes \cU_{\pi^{-1}(i+1)})  \hskip .3cm \text{ if } r=\pi^{-1}(i),\\
			\cU_{r+1}  \hskip 2.7cm\text{ if } \pi^{-1}(i) + 1 \le r\le k,
		\end{cases}$$
		then we have 
		\begin{eqnarray}\label{eq:43}
			%\lim_{n_1\in\cU_1}\cdots \lim_{n_{i-1}\in\cU_{i-1}} \lim_{n_i\in\cU_i} \lim_{n_{i+1}\in\cU_{i+1}} \lim_{n_{i+2}\in\cU_{i+2}} \cdots \lim_{n_k\in\cU_k}d(f(n_1,\dots,  n_{i-1}, n_i, n_{i+1}, n_{i+2},n_\ell),g(n_{\ell+1},\dots,n_k))\\
			\nonumber \lim_{n_1\in\cU_1}\cdots\lim_{n_{k+1}\in\cU_{k+1}}d\big( f(n_1,\dots,n_\ell),g(n_{\ell+1},\dots,n_{k+1}) \big) = \\
			\lim_{m_1\to\tilde{\cU}_{1}} \cdots \lim_{m_{k}\to\tilde{\cU}_{k}}d\big( \tilde{f}(m_1, \dots, m_{\ell-1}), g(m_{\ell},\dots,m_{k})\big).
		\end{eqnarray}
		
		Also, if $$\tilde{\pi}(r) := \begin{cases}
			\pi(r)  \hskip 2cm\text{ if } 1\le r\le \pi^{-1}(i) \le \ell-1,\\
			\pi(r+1) \hskip 1.5cm \text{ if }  \pi^{-1}(i)+1 \le r \le k\text{ and } \pi(r+1)< i,\\
			\pi(r+1)-1 \hskip 1cm \text{ if } \pi^{-1}(i)+1 \le r \le k \text{ and } \pi(r+1)> i,
		\end{cases}$$
		then $\tilde{\pi}$ is a permutation of $[k]$ that preserves the order on $\{1,\dots, \ell-1\}$ and $\{\ell, \dots, k\}$ and 
		$$\tilde{\pi}^{-1}(s) = \begin{cases}
			\pi^{-1}(s) \hskip 2cm\text{ if } 1\le s \le i \text{ and } \pi^{-1}(s) \le \pi^{-1}(i),\\
			\pi^{-1}(s) -1 \hskip 1.5cm\text{ if } 1\le s < i \text{ and }\pi^{-1}(s) > \pi^{-1}(i),\\
			%	                          \pi^{-1}(r+1) \text{ if } i+1\le r \le k \text{ and }\pi^{-1}(r+1) < \pi^{-1}(i),\\
			\pi^{-1}(s+1)-1 \hskip 1cm \text{ if } i+1\le s \le k \text{ and }\pi^{-1}(s+1) > \pi^{-1}(i).
		\end{cases}$$
		It is not difficult, though a bit tedious, to verify that $\tilde{\pi}^{-1}\circ \tilde{\pi} = \tilde{\pi} \circ \tilde{\pi}^{-1} = id$.
		%\begin{proof}
		%If $1\le r\le \pi^{-1}(i) \le \ell-1$ then $\pi(r)\le i$ since $\pi$ preserves the order and thus $s:=\pi(r)\le i$. Also $\pi^{-1}(s)=r\le \pi^{-1}(i)$. Therefore,$\tilde{\pi}^{-1}\circ \tilde{\pi}(r) = \tilde{\pi}^{-1}(\pi(r)) = \pi^{-1}(\pi(r)) =r$. If $\pi^{-1}(i)+1 \le r \le k$ and $\pi(r+1)< i$, let $s:=\pi(r+1)$ and thus $s<i$.Also, $\pi^{-1}(s)=r+1\ge \pi^{-1}(i)+2>\pi^{-1}(i)$ and then $\tilde{\pi}^{-1}\circ \tilde{\pi}(r) = \tilde{\pi}^{-1}(\pi(r+1)) = \pi^{-1}(\pi(r+1))-1 =r$.Finally, If $\pi^{-1}(i)+1 \le r \le k$ and $\pi(r+1)> i$, let $s:=\pi(r+1)-1$ and thus $s\ge i$. In fact $s>i$ otherwise, $\pi(r+1)=i+1 =\pi (\pi^{-1}(i+1)$ and hence $\pi^{-1}(i+1) =r+1$ and $\pi^{-1}(i)=r$; a contradiction. Also, $\pi^{-1}(s+1) =r+1 \ge \pi^{-1}(i)+2 >\pi^{-1}(i)$. Therefore, $\tilde{\pi}^{-1}\circ \tilde{\pi}(r) = \tilde{\pi}^{-1}(\pi(r+1)-1) = \pi^{-1}(\pi(r+1)-1+1)-1 =r$. Since $\tilde{\pi}$ and $\tilde{\pi}^{-1}$ are bijections we also have $\tilde{\pi}\circ \tilde{\pi}^{-1}(s)=s$ for al $s$.  
		%\end{proof}
		To help understand these permutations, below is an illustrative example with $\ell=4$ and $k+1=7$. If 
		$$\pi = \left(\begin{smallmatrix}
			1 & 2 & 3 & 4 & 5 & 6 & 7 \\
			1 & 3 & 4 & 6 & 2 & 5 & 7 \\
		\end{smallmatrix}
		\right)
		\text{ and hence }
		\pi^{-1} = \left(\begin{smallmatrix}
			1 & 2 & 3 & 4 & 5 & 6 & 7 \\
			1 & 5 & 2 & 3 & 6 & 4 & 7 \\
		\end{smallmatrix}
		\right),
		$$
		then
		$$
		\tilde{\pi} = \left(\begin{smallmatrix}
			1 & 2 & 3 & 4 & 5 & 6  \\
			1 & 3 & 5 & 2 & 4 & 6  \\
		\end{smallmatrix}
		\right)
		\text{ and }
		\tilde{\pi}^{-1} = \left(\begin{smallmatrix}
			1 & 2 & 3 & 4 & 5 & 6  \\
			1 & 4 & 2 & 5 & 3 & 6  \\
		\end{smallmatrix}
		\right).$$
		In this situation $i=3$, $\pi^{-1}(2)= 2$ and $\pi^{-1}(4)= 3 = \pi^{-1}(2)+1$.\\
		The permutations $\tilde{\pi}^{-1}$ and $\tilde{\pi}$ are obtained by identifying the consecutive limits and shifting accordingly the other entries.
		Moreover,                          
		\begin{eqnarray}\label{eq:44}
			%\lim_{n_1\in\cU_{\pi^{-1}(1)}}\cdots \lim_{n_{\pi(i)}\in\cU_{i}} \lim_{n_{\pi(i+1)}\in\cU_{i+1}} \cdots \lim_{n_k\in\cU_{\pi^{-1}(k)}}d(f(n_{\pi(1)}, \dots, n_{\pi(i)}, n_{\pi(i+1)}, \cdots, n_{\pi(\ell)}),g(n_{\pi(\ell+1)},\dots,n_{\pi(k)}))\\
			\nonumber \lim_{n_{\pi^{-1}(1)}\in\cU_{\pi^{-1}(1)}}\cdots \lim_{n_{\pi^{-1}(k+1)}\in\cU_{\pi^{-1}(k+1)}} d\big( f(n_{1}, \dots,n_{\ell}),g(n_{\ell+1},\dots,n_{k+1}) \big)= \\
			\lim_{m_{\tilde{\pi}^{-1}(1)}\to \tilde{\cU}_{\tilde{\pi}^{-1}(1)}}\cdots \lim_{m_{\tilde{\pi}^{-1}(k)}\to \tilde{\cU}_{\tilde{\pi}^{-1}(k)}}d\big( \tilde{f}(m_{1}, \dots, m_{\ell-1}), g(m_{\ell},\dots, m_{k}) \big).
		\end{eqnarray}
		We can now use the induction hypothesis to conclude that all the quantities in \eqref{eq:43} and \eqref{eq:44} are equal.
		
		The case when there is a pair of consecutive limits amongst the last $k+1-\ell$ original limits which are still consecutive in the reordering given by $\pi^{-1}$, i.e. if $\ell+2\le \pi^{-1}(i)+1 = \pi^{-1}(i+1) \le k+1$ for some $1\le i \le k$, can be treated similarly.
		
		Now, if there is no pair of consecutive limits amongst the first $\ell$, nor the last $k+1-\ell$, original limits which are still consecutive in the reordering given by $\pi^{-1}$, i.e. for no $i\in \{1,\dots, \ell-1\} \cup \{\ell+1, \dots, k\}$ do we have $\pi^{-1}(i+1) = \pi^{-1}(i)+1$, then $\pi^{-1}$ is necessarily of one of the following two types: 
		
		(a) $\pi^{-1} = \left(\begin{smallmatrix}
			&  \cdots    &  & \cdots  &  		 k-1 & k 	& k+1\\
			&   \cdots & 	   &  \cdots	&  \ell-1 & k+1 & \ell\\
		\end{smallmatrix}
		\right)$
		
		(b)
		$\pi^{-1} = \left(\begin{smallmatrix}
			&  \cdots    &  & \cdots  &  	 k-1 & k 	& k+1\\
			&   \cdots & 	   &  \cdots	& k & \ell & k+1\\
		\end{smallmatrix}\right).$ 
		
		\noindent In case (a), by stability we have that for all $(m_1, \dots ,m_{k-1})\in \bN^{k-1}$
		\begin{eqnarray*}
			\lim_{n_{k+1}\in\cU_{k+1}} \lim_{n_{\ell}\in\cU_{\ell}} d\big( f(m_1,\dots,m_{\ell-1}, n_\ell),g(m_{\ell},\dots,m_{k-1}, n_{k+1}) \big)=\\
			\lim_{n_{\ell}\in\cU_{\ell}}\lim_{n_{k+1}\in\cU_{k+1}}d\big( f(m_1,\dots,m_{\ell-1}, n_\ell),g(m_{\ell},\dots,m_{k-1}, n_{k+1}) \big),
		\end{eqnarray*}
		and hence
		\begin{eqnarray*}
			\lim_{n_{\pi^{-1}(1)}\in\cU_{\pi^{-1}(1)}}\cdots \lim_{n_{\pi^{-1}(k+1)}\in\cU_{\pi^{-1}(k+1)}} d\big( f(n_{1}, \dots,n_{\ell}), g(n_{\ell+1},\dots,n_{k+1}) \big)= \\
			\lim_{n_{\pi^{-1}(1)}\in\cU_{\pi^{-1}(1)}}\cdots \lim_{n_{\ell-1}\in\cU_{\ell-1}} \lim_{n_{k+1}\in\cU_{k+1}} \lim_{n_\ell\in\cU_{\ell}}d\big( f(n_{1}, \dots, n_{\ell}),g(n_{\ell+1},\dots,n_{k+1}) \big)=\\
			\lim_{n_{\pi^{-1}(1)}\in\cU_{\pi^{-1}(1)}}\cdots \lim_{n_{\ell-1}\in\cU_{\ell-1}} \lim_{n_\ell\in\cU_{\ell}}\lim_{n_{k+1}\in\cU_{k+1}}d\big( f(n_{1}, \dots, n_{\ell}),g(n_{\ell+1},\dots,n_{k+1}) \big).
		\end{eqnarray*}
		We have thus reduced our problem to one of the previous situations already handled above. A similar reduction can be done in case (b).
		
		The implication $(ii)\Rightarrow (iii)$ is just formal. 
		
		To prove $(iii) \Rightarrow (iv)$ let $\bM$ be an infinite subset of $\bN$ and $\cU$ be an ultrafilter on $\bN$ containing $\bM$. Since for every $\nbar\in[\bM]^k$,
		\begin{multline*}
			d(f(n_{\pi(1)},\dots,n_{\pi(\ell)}),g(n_{\pi(\ell+1)},\dots,n_{\pi(k)})) \le \\
			\sup_{\nbar\in [\bM]^k}d(f(n_{\pi(1)},\dots,n_{\pi(\ell)}),g(n_{\pi(\ell+1)},\dots,n_{\pi(k)}))
		\end{multline*}
		and since the ultrafilter $\cU$ contains $\bM$ we have
		\begin{multline*}
			\lim_{n_1\in\cU}\cdots\lim_{n_k\in\cU}d(f(n_{\pi(1)},\dots,n_{\pi(\ell)}),g(n_{\pi(\ell+1)},\dots,n_{\pi(k)})) \le \\
			\sup_{\nbar\in [\bM]^k}d(f(n_{\pi(1)},\dots,n_{\pi(\ell)}),g(n_{\pi(\ell+1)},\dots,n_{\pi(k)})).
		\end{multline*}
		By $(iii)$,
		\begin{eqnarray*}
			\lim_{n_1\in\cU}\cdots\lim_{n_k\in\cU}d(f(n_1,\dots,n_\ell),g(n_{\ell+1},\dots,n_k))\le &\\
			\sup_{\nbar\in [\bM]^k}d(f(n_{\pi(1)},\dots,n_{\pi(\ell)}),g(n_{\pi(\ell+1)},\dots,n_{\pi(k)})).&
		\end{eqnarray*}
		Since $\cU$ contains $\bM$, \eqref{eq:1-upper-stable} immediately follows. 
		
		The implication $(iv)\Rightarrow (i)$ essentially follows from \cite[Th\'eor\`eme page 276]{KrivineMaurey1981}, where the argument is given for separable Banach spaces. We will show that if $(M,d)$ is not stable, then \eqref{eq:1-upper-stable} is violated for $k=2$ and $\pi$ the interlacing permutation, which is the only nontrivial permutation in this case. So, condition \eqref{eq:1-upper-stable}, in this case, is simply
		\begin{equation}
			\label{eq:stable-ineq}
			\inf_{n_1<n_2}d\big( f(n_1),g(n_2) \big)\le  \sup_{n_1<n_2}d\big( f(n_2), g(n_1) \big).
		\end{equation}
		Assume that there are ultrafilters $\cU_1$ and $\cU_2$ on $\bN$ and two bounded functions $f,g\colon \bN \to M$ such that 
		$$\lim_{n_1\in\cU_1}\lim_{n_2\in\cU_2}d\big( f(n_1),g(n_2) \big) \neq \lim_{n_2\in\cU_2}\lim_{n_1\in\cU_1}d\big( f(n_1),g(n_2) \big).$$
		We only treat the case $$\lim_{n_1\in\cU_1}\lim_{n_2\in\cU_2}d\big( f(n_1),g(n_2) \big) < \lim_{n_2\in\cU_2}\lim_{n_1\in\cU_1}d\big( f(n_1),g(n_2) \big),$$
		as the other case can be handled similarly.
		Let $\alpha, \beta>0$ such that $$\lim_{n_1\in\cU_1}\lim_{n_2\in\cU_2}d\big( g(n_2), f(n_1) \big) <\alpha <\beta < \lim_{n_2\in\cU_2}\lim_{n_1\in\cU_1}d\big( g(n_2),f(n_1) \big).$$
		There is $B\in \cU_2$ and $A\in \cU_1$ such that $\lim_{n_1\in\cU_1}d\big( g(n_2),f(n_1) \big)>\beta$ whenever $ n_2\in B$,  and $\lim_{n_2\in\cU_2}d\big( g(n_2),f(n_1) \big)<\alpha$ whenever $ n_1\in A $.
		We now construct recursively two sequences of integers $(a_j)_{j=1}^\infty$ and $(b_j)_{j=1}^\infty$ together with decreasing sequences of subsets $(A_j)_{j=1}^\infty$ in $\cU_1$ and $(B_j)_{j=1}^\infty$ in $\cU_2$ as follows.
		Let $B_1:=B$ and pick $b_1\in B_1$. Then, let $A_1:=A\cap \{ n_1\in \bN \colon d\big( g(b_1), f(n_1) \big)>\beta\}$ and observe that $A_1\in \cU_1$. Pick $a_1\in A_1\subset A$ and let $B_2:=B_1\cap\{n_2\in \bN \colon d\big( g(n_2),f(a_1) \big)<\alpha\}$ which is clearly in $\cU_2$. Assume now that $a_{k}$, $A_{k}$ and $B_{k+1}$ have been defined. First, pick $b_{k+1}\in B_{k+1}$, then let $A_{k+1}:=A_k\cap \{ n_1\in \bN \colon d\big( g(b_{k+1}), f(n_1) \big)>\beta\}$ and pick $a_{k+1}\in A_{k+1}$. Finally, we let $B_{k+2}:=B_{k+1}\cap \{n_2\in \bN \colon d\big( g(n_2),f(a_{k+1}) \big)<\alpha\}$.
		By construction, it is plain that for all $k\in \bN$, we have $a_k\in A_k$ and $b_k\in B_k$. Moreover, if $k\le \ell$, then $a_\ell\in A_\ell\subseteq A_k$, which implies that $d\big( g(b_{k}), f(a_\ell) \big)>\beta$ and in turn $\inf_{k< \ell}d\big( g(b_{k}), f(a_\ell) \big)>\beta$.
		On the other hand, when $k>\ell$ we have $b_k\in B_k\subseteq B_{\ell+1}$ and hence $d\big( g(b_k),f(a_{\ell}) \big)<\alpha$. Therefore, $\sup_{k>\ell}d\big( g(b_k),f(a_{\ell}) \big)<\alpha$ and 
		$$\sup_{\ell<k}d\big( g(b_k),f(a_{\ell}) \big)<\alpha<\beta <\inf_{k< \ell}d\big( g(b_{k}), f(a_\ell) \big),$$
		which is easily seen to contradict \eqref{eq:stable-ineq}.

	\end{proof}
	
	Statement $(iv)$ in Theorem \ref{thm:stability} suggests the following bi-Lipschitz relaxation of stability.
	
	\begin{defi}
		\label{def:infrasup-stable}
		A metric space $(M,d)$ is said to be \emph{$C$-infrasup-stable} if for every $k\ge 1$, every $1\le \ell < k$, every permutation $\pi\colon [k]\to [k]$ which preserves the order on $\{1,\dots,\ell\}$ and $\{\ell+1,\dots,k\}$ and all bounded maps $f\colon [\bN]^{\ell}\to M$ and $g\colon [\bN]^{k-\ell}\to M$ we have for every infinite subset $\bM$ of $\bN$,
		\begin{multline}
			\label{eq:infrasup-stable}
			\inf_{\nbar\in[\bM]^k}d(f(n_1,\dots,n_\ell),g(n_{\ell+1},\dots,n_k))\le \\
			C \sup_{\nbar\in [\bM]^k}d(f(n_{\pi(1)},\dots,n_{\pi(\ell)}),g(n_{\pi(\ell+1)},\dots,n_{\pi(k)})).
		\end{multline} 
		A metric space is \emph{infrasup-stable} if it is $C$-infrasup-stable for some $C\ge 1$.
	\end{defi}
	
	The observation below is immediate, and its proof is left to the reader. 
	
	\begin{prop}\label{prop:Lipstabilty infrasup}
		If $(M,d)$ embeds bi-Lipschitzly into an infrasup-stable space $(N,d)$, then $M$ is infrasup-stable. 
		Moreover, if we denote by $C_{iss}(M)$ the smallest constant $C$ such that \eqref{eq:infrasup-stable} holds, then $C_{iss}(M)\le \cdist{N}(M)C_{iss}(N)$.
	\end{prop}
	
	It follows from Theorem \ref{thm:stability} that the class of $1$-infrasup-stable spaces coincides with the class of stable spaces. As we have already noted, $\co$ is not stable, but it is not immediately evident whether $\co$ admits an equivalent norm that is stable or whether $\co$ admits a bi-Lipschitz embedding into a stable space. The fact that being infrasup-stable is a bi-Lipschitz invariant makes it a more convenient and versatile notion to work with when addressing this type of question. Moreover, many results for stable metrics that are of an isomorphic nature can be easily extended to the infrasup-stable metrics. 
	
	\begin{prop}
		\label{pro:co-not-infrasup-stable}
		$\co$ is not infrasup-stable. 
	\end{prop}
	
	\begin{proof}
		Denoting as always by $(s_n)_{n\in\bN}$ the summing basis of $\co$, it is easy to see that for all $\nbar\in [\bN]^{2k}$, 
		\begin{align*}
			\Big\|\sum_{i=1}^k s_{n_i}-\sum_{i=k+1}^{2k} s_{n_i}\Big\|_\infty= k& \qquad\textrm{and} & \Big\|\sum_{i=1}^k s_{n_{2i-1}}-\sum_{i=1}^k s_{n_{2i}}\Big\|_\infty= 1.
		\end{align*}
		Therefore, if we let $f(n_1,\dots, n_k) = g(n_1,\dots, n_k) = \sum_{i=1}^k s_{n_i}$ and take $\pi$ to be the interlacing permutation, then the smallest $C$ such that $\eqref{eq:infrasup-stable}$ holds satisfies $C\ge k$.
	\end{proof}
	
	\begin{rema}
		With an argument similar to the one in the proof of Proposition \ref{pro:co-not-infrasup-stable}, we can also show that $\sup_{k\in\bN}C_{iss}(\Ik)=\infty$. 
		In particular, it follows that $(\Ik)_{k\in \bN}$ does not embed equi-bi-Lipschitzly and thus $\co$ does not embed bi-Lipschitzly into a stable space. We will prove something much stronger below.
	\end{rema}
	
	Infrasup-stability is not only a bi-Lipschitz invariant, but it is also preserved under coarse embeddings under the assumption that the domain space is a Banach space.
	
	\begin{lemm}
		\label{lem:coarse-uniform-infrasup-stable}
		If a Banach space $X$ coarsely, or uniformly, embeds into an infrasup-stable metric space, then $X$ is infrasup-stable.
	\end{lemm}
	
	\begin{proof}
		We prove first the coarse statement.
		Since $(X,\norm{\cdot})$ coarsely embeds into an infrasup-stable metric space $(M,d)$, there exist $\varphi\colon X\to M$ and nondecreasing functions $\rho,\omega\colon [0,\infty)\to [0,\infty)$ such that 
		\begin{equation*}
			\rho(\norm{x-y})\le d(\varphi(x),\varphi(y))\le\omega(\norm{x-y}),
		\end{equation*}
		with $\lim_{t\to\infty}\rho(t)=\infty$. Let $1\le \ell < k$, $f\colon [\bN]^{\ell}\to X$ and $g\colon [\bN]^{k-\ell}\to X$ bounded maps and $\pi\colon [k]\to [k]$ a permutation that preserves the order on $\{1,\dots,\ell\}$ and $\{\ell+1,\dots,k\}$, be given. Let $\bM$ an infinite subset of $\bN$ and let $$\alpha\eqd\sup_{\nbar\in[\bM]^k} \norm{f(n_{\pi(1)},\dots,n_{\pi(\ell)})-g(n_{\pi(\ell+1)},\dots,n_{\pi(k)})}.$$
		Clearly, we may assume that $\alpha>0$. Then, 
		\begin{equation*}
			\sup_{\nbar\in [\bM]^k}d \Big(\varphi \big(\frac{1}{\alpha} f(n_{\pi(1)},\dots,n_{\pi(\ell)}) \big),\varphi \big(\frac{1}{\alpha}g(n_{\pi(\ell+1)},\dots,n_{\pi(k)}) \big) \Big)\le \omega(1). 
		\end{equation*}
		Let $\beta\eqd \frac{1}{\alpha}\inf_{\nbar\in[\bM]^k} \norm{f(n_1,\dots,n_\ell)-g(n_{\ell+1},\dots,n_k)}$ and observe that 
		\begin{equation*}
			\inf_{\nbar\in[\bM]^k}d \Big( \varphi \big( \frac{1}{\alpha}f(n_1,\dots,n_\ell) \big),\varphi \big( \frac{1}{\alpha}g(n_{\ell+1},\dots,n_k) \big) \Big)\ge \rho(\beta). 
		\end{equation*}
		The infrasup-stability assumption gives
		\begin{equation}\label{eq:contradiction-coarse}
			\rho(\beta) \le C_{iss}(M)\omega(1)<\infty.
		\end{equation}
		Let $C>0$ such that $\rho(C)> C_{iss}(M)\omega(1)$ (such a $C$ exists since by assumption $\lim_{t\to\infty}\rho(t)=\infty$). If  $\beta > C$, then whenever $\norm{x-y}\ge \beta$ we have $d \big( \varphi(x),\varphi(y) \big)\ge \rho \big( \norm{x-y} \big)\ge \rho(C)> C_{iss}(M)\omega(1)$, but this contradicts \eqref{eq:contradiction-coarse}. Therefore, $\beta\le C$ necessarily, and this completes the proof since it precisely means that 
		\begin{multline*}
			\inf_{\nbar\in[\bM]^k}\norm{ f(n_1,\dots,n_\ell) - g(n_{\ell+1},\dots,n_k) }\le \\
			C \sup_{\nbar\in [\bM]^k} \norm{ f(n_{\pi(1)},\dots,n_{\pi(\ell)}) - g(n_{\pi(\ell+1)},\dots,n_{\pi(k)}) }.
		\end{multline*}
		
		Now, for the uniform case we assume that $\lim_{t\to 0}\omega(t)=0$ and $\rho(t)>0$ for all $t>0$. Let $\beta\eqd \inf_{\nbar\in[\bM]^k} \norm{f(n_1,\dots,n_\ell)-g(n_{\ell+1},\dots,n_k)}$. Of course, we may assume that $\beta >0$. Then, 
		\begin{equation*}
			\inf_{n\in[\bM]^k}d \Big( \varphi \big( \frac{1}{\beta}f(n_1,\dots,n_\ell) \big),\varphi \big( \frac{1}{\beta}g(n_{\ell+1},\dots,n_k) \big) \Big)\ge \rho(1)
		\end{equation*}
		Let $$\alpha\eqd\frac{1}{\beta}\sup_{\nbar\in[\bM]^k} \norm{f(n_{\pi(1)},\dots,n_{\pi(\ell)})-g(n_{\pi(\ell+1)},\dots,n_{\pi(k)})},$$
		and observe that 
		\begin{equation*}
			\sup_{n\in[\bM]^k}d \Big( \varphi \big( \frac{1}{\beta}f(n_{\pi(1)},\dots,n_{\pi(\ell)}) \big), \varphi \big( \frac{1}{\beta}g(n_{\pi(\ell+1)},\dots,n_{\pi(k)})\big) \Big)\le \omega(\alpha). 
		\end{equation*}
		The infrasup-stability assumption then gives
		\begin{equation}
			\label{eq:contradiction-uniform}
			\rho(1)\le C_{iss}(M)\omega(\alpha)
		\end{equation}
		Let $C>0$ such that $\omega(C)< \frac{\rho(1)}{C_{iss}(M)}$ (such a $C$ exists since by assumption $\rho(1)>0$ and $\lim_{t\to 0}\omega(t)=0$). If  
		$\alpha < C$, then whenever $\norm{x-y}\le \alpha$ we have $d \big( \varphi(x),\varphi(y) \big)\le \omega \big( \norm{x-y} \big)\le \omega(C)< \frac{\rho(1)}{C_{iss}(M)}$, but this contradicts \eqref{eq:contradiction-uniform}. Therefore, $\alpha\ge C$ necessarily, and this completes the proof since it precisely means that 
		\begin{multline*}
			\inf_{\nbar\in[\bM]^k} \norm{f(n_1,\dots,n_\ell) - g(n_{\ell+1},\dots,n_k) } \le \\ \frac{1}{C} \sup_{\nbar\in [\bM]^k} \norm{ f(n_{\pi(1)},\dots,n_{\pi(\ell)}) - g(n_{\pi(\ell+1)},\dots,n_{\pi(k)})}.
		\end{multline*}
	\end{proof}
	
	\begin{rema}
		We do not know if Lemma \ref{lem:coarse-uniform-infrasup-stable} still holds if $X$ is merely a metric space.
	\end{rema}
	
	The next corollary emphasizes how far $\co$ is from being infrasup-stable.
	
	\begin{coro}
		$\co$ does not coarsely, nor uniformly, embed into an infrasup-stable metric space.
	\end{coro}

	We finish this section with an alternative view of infrasup-stability. There is a certain asymmetry in Theorem \ref{thm:stability} $(iv)$ above and
	hence in the definition of infrasup-stability, in the sense that the
	permutation is applied only to the right-hand side of the defining
	inequality. We will now introduce an a priori stronger version of
	infrasup-stability by ``symmetrizing'' this inequality. We will also
	introduce a different and often more convenient formulation in terms of cuts instead of
	permutations.
	
	Fix $k,l,m\in\bN$ with $k=l+m$. By a \emph{cut of $\{1,\dots,k\}$ of
		size $l$ }we mean simply a subset $P$ of $\{1,\dots,k\}$ of size
	$l$. If $P=\{p_1,\dots,p_l\}$ with $p_1<\dots<p_l$, then for
	$n_1<\dots<n_k$ in $\bN$, we write $(n_i\colon i\in P)$ for
	$(n_{p_1},\dots,n_{p_l})$. The proof of Theorem \ref{thm:stability} shows the following.
	
	\begin{prop}
		\label{prop:cut-stable-characterization}
		For a metric space $(M,d_M)$ and a constant $C\ge1$, the following
		are equivalent.
		\begin{enumerate}[(i)]
			\item
			For all $k,l,m\in\bN$ with $k=l+m$, every ultrafilter $\cU$,
			bounded functions $f\colon [\bN]^{l}\to M$ and
			$g\colon[\bN]^{m}\to M$ and cuts $P,Q$ of $\{1,\dots,k\}$ of
			size $l$, we have
			\begin{multline*}
				\lim_{n_1\in \cU}\cdots\lim_{n_k\in \cU}
				d\big(f(n_i\colon i\in P),g(n_i\colon i\in P^c)\big)\le \\
				C\lim_{n_1\in \cU}\dots\lim_{n_k\in \cU}
				d\big(f(n_i\colon i\in Q),g(n_i\colon i\in Q^c)\big)\ .
			\end{multline*}
			\item
			For all $k,l,m\in\bN$ with $k=l+m$, bounded functions
			$f\colon[\bN]^{l}\to M$ and $g\colon[\bN]^{m}\to M$,
			infinite sets $\bL\subset\bN$ and cuts $P,Q$ of $\{1,\dots,k\}$ of
			size $l$, we have
			\begin{multline*}
				\inf_{\nbar\in[\bL]^{k}}d\big(f(n_i\colon i\in P),g(n_i\colon i\in
				P^c)\big)\le \\
				C\cdot\sup_{\nbar\in[\bL]^{k}}d\big(f(n_i\colon i\in Q),g(n_i\colon i\in
				Q^c)\big)\ .
			\end{multline*}
		\end{enumerate}
	\end{prop}
	
	\begin{defi}
		We say that a metric space $(M,d)$ is \emph{$C$-cut-stable }if it
		satisfies one of the two equivalent conditions in
		Proposition~\ref{prop:cut-stable-characterization}. We say $M$ is
		\emph{cut-stable }if it is $C$-cut-stable for some $C\ge1$.
	\end{defi}
	
	\section{\texorpdfstring{Kalton's Property $Q$}{Kalton's Property}}
	\label{sec:Q}
	
	Infrasup-stability is a natural bi-Lipschitz relaxation of stability, and it can also be used as an obstruction to coarse embeddings for Banach spaces. But at this point, it is not clear whether the interlaced graphs could be embedded coarsely into an infrasup-stable space. In \cite{Kalton2007}, Kalton derived another coarse invariant from the notion of stability, which he called \emph{property $Q$}. In fact, Kalton's property $Q$ can be derived from the weaker notion of infrasup-stability (which was introduced later in \cite{Baudier2022}). Quite remarkably, we will see that property $Q$ and infra sup-stability coincide.
	
	\begin{defi} 
		\label{def:Q}
		A metric space $(M,d)$ has \emph{property $Q$} if there exists a constant $C\in(0,\infty)$ so that for every $k\in \bN$ and every Lipschitz map $f\colon (\Nk,\sd_{\sI_k})\to M$ there is an infinite subset $\bM$ of $\bN$ such that for all $\mbar,\nbar\in [\bM]^k$,
		\begin{equation}
			\label{eq:Q}
			d \big( f(\mbar),f(\nbar) \big)\le C\cdot\Lip(f).
		\end{equation}
		We define $Q_M$ to be the infimum of such constants $C$, when $M$ has property $Q$ and $Q_M=\infty$ otherwise\footnote{Note that it is the inverse of the constant defined in \cite{Kalton2007}}.
	\end{defi}
	
	Property $Q$ is thus a concentration property for maps defined on the interlaced graphs that is completely analogous to the concentration property for Hamming graphs we have seen in Section \ref{sec:Hamming}. Property $Q$ is easily seen to be an obstruction to coarse embeddings.
	
	\begin{prop}
		\label{prop:interlaced-notCE->Q}
		The sequence of interlaced graphs $(\sI_k)_{k\in\bN}$ does not equi-coarsely embed into a metric space with property $Q$.
	\end{prop}
	
	\begin{proof}
		Assume $M$ is a metric space with property $Q$ and let $C\ge 1$ be given by the definition of property $Q$. Aiming for a contradiction, assume that the family $(\Ik)_{k\in \bN}$ equi-coarsely embeds into $M$. That is, there are maps $f_k \colon ([\bN]^k,\dIk) \to X$ and two functions $\rho, \omega \colon [0,+\infty)\to[0,+\infty)$  such that $\lim_{t \to \infty} \rho (t)=\infty$ and
		$$\forall k\in \bN\ \ \forall t>0\ \ \rho(t) \le  \rho_{f_k}(t)\ \ \text{and}\ \ \omega_{f_k}(t) \le  \omega(t)<\infty.$$
		Thus, for every $k\in \bN$, there exists an infinite subset $\mathbb M_k$ of $\bN$ such that $\diam(f([\M_k]^k)))\le C\omega(1)$. Since $\diam([\M_k]^k)=k$, this implies that for all $k\in \bN$, $\rho(k) \le C\omega(1)$. This contradicts the fact that $\lim\limits_{t \to \infty} \rho(t)=\infty$.
	\end{proof}
	
	Kalton had observed that every stable space has property $Q$, but it is not much more difficult to show that infrasup-stable spaces also have property $Q$. We refer the reader to Proposition \ref{prop:Lipstabilty infrasup} for the definition of $C_{iss}(M)$. 
	
	\begin{prop}
		\label{prop:stable->Q}
		Every infrasup-stable metric space has property $Q$. More precisely, if $(M,d)$ is an infrasup-stable metric space, then $Q_M\le 4C_{iss}(M)$. 
	\end{prop}
	
	\begin{proof}
		Let $(M,d)$ be $C$-infrasup-stable and fix $k\ge 1$ and a Lipschitz map $f\colon (\Nk,\dIk)\to M$. Assume, as we may, that $\Lip(f)>0$ and define
		$$\cA \eqd \Big\{\nbar\in [\bN]^{2k}\colon d \big( f(n_1,\dots,n_k),f(n_{k+1},\dots, n_{2k}) \big)\le 2C\cdot\Lip(f) \Big\}.$$
		By Ramsey theorem there exists an infinite subset $\bM$ of $\bN$ such that either $[\bM]^{2k}\subset \cA$ or $[\bM]^{2k}\cap \cA=\emptyset$.
		If the first possibility happens, then for all $\nbar\in [\bM]^{2k}$ we have 
		$$d \big( f(n_1,\dots,n_k),f(n_{k+1},\dots, n_{2k}) \big)\le 2C\cdot\Lip(f).$$
		If $\nbar,\mbar\in[\bM]^k$ we can choose $\xbar\in[\bM]^k$ such that $x_1>\max\{n_k,m_k\}$ and thus
		$$d \big( f(\nbar),f(\mbar) \big)\le d \big( f(\nbar),f(\xbar)\big) + d\big( f(\xbar),f(\mbar) \big)\le 4C\cdot\Lip(f).$$
		It remains to show that the second possibility cannot happen. Consider the interlacing permutation $\pi\colon [2k]\to [2k]$ defined by 
		\begin{equation*}
			\pi(i)=\begin{cases}
				2i-1 \quad if \quad 1\le i\le k\\
				2(i-k) \quad if \quad k+1\le i\le 2k.
			\end{cases}
		\end{equation*}
		It is immediate to verify that $\pi$ preserves the order on $\{1,\dots,k\}$ and on $\{k+1,\dots,2k\}$ and if $\nbar\in[\bM]^{2k}$ then
		\begin{multline*}
			d(f(n_{\pi(1)},\dots,n_{\pi(k)}),f(n_{\pi(k+1)},\dots,n_{\pi(2k)}))= \\
			d(f(n_1,n_3,\dots,n_{2k-1}),f(n_2,n_4,\dots,n_{2k}))\le \Lip(f).
		\end{multline*}
		Therefore,
		\begin{equation*}
			\sup_{\nbar\in [\bM]^{2k}}d(f(n_{\pi(1)},\dots,n_{\pi(k)}),f(n_{\pi(k+1)},\dots,n_{\pi(2k)}))\le \Lip(f),
		\end{equation*}
		and by $C$-infrasup stability it holds
		\begin{equation*}
			\inf_{\nbar\in [\bM]^{2k}}d(f(n_{1},\dots,n_{k}),f(n_{k+1},\dots,n_{2k}))\le C\Lip(f).
		\end{equation*}
		In particular there exists $\nbar\in[\bM]^{2k}$ such that 
		\begin{equation*}
			d(f(n_{1},\dots,n_{k}),f(n_{k+1},\dots,n_{2k}))\le 2C\Lip(f),
		\end{equation*}
		but this implies that $[\bM]^{2k}\cap \cA\neq\emptyset$; a contradiction.
	\end{proof}
	
	%\begin{rema}
	%A careful analysis of the above proof reveals that for every infrasup-stable space $M$ one has $Q_M\le 2C_{iss}$.
	%\end{rema}
	
	The notion of infrasup-stability seems stronger than property $Q$ as infrasup-stability says something about all permutations that preserve the order, whereas Property $Q$ only seems to say something about the interlacing permutation. In the remainder of this section, we will show that a space with property $Q$ is in fact also infrasup-stable and thus property $Q$ and infrasup-stability are two faces of the same coin. 
	
	We will need the following stronger version of the interlacing adjacency relation. Given $k\in\bN$ and $\mbar,\nbar\in[\bN]^{k}$, we say that the pair $(\mbar,\nbar)$ is \emph{strongly interlacing }if either $m_1\le n_1<m_2\le n_2<\dots<m_k\le n_k$ or $n_1\le m_1<n_2\le m_2<\dots< n_k\le m_k$.
	
	\begin{lemm}
		\label{lem:strongly-interlacing}
		Given $k\in\bN$ and an interlacing pair $(\mbar,\nbar)$ in
		$[2\bN]^{k}$, there exists $\pbar\in[\bN]^{k}$ such that
		the pairs $(\mbar,\pbar)$ and $(\nbar,\pbar)$ are strongly interlacing.
	\end{lemm}
	
	\begin{proof}
		Without loss of generality, we may assume that
		\[
		m_1\le n_1\le m_2\le n_2\le\dots\le m_k\le n_k\ .
		\]
		Define $\pbar=\{p_1,\dots,p_k\}$ in $[\bN]^{k}$ by
		\[
		p_i=\begin{cases}%
			n_i-1 & \text{if }n_i=m_{i+1}\ ,\\
			n_i & \text{if }n_i<m_{i+1}\ .
		\end{cases}
		\]
		It is straightforward to verify that this works.
	\end{proof}
	
	We are now ready to prove that property $Q$ implies infrasup-stability. 
	
	\begin{prop}
		\label{prop:property-Q-implies-cut-stable}
		A metric space with property~$Q$ is cut-stable.
	\end{prop}
	
	\begin{proof}  
		Let $(M,d)$ be a metric space with property~$Q$. It will be more convenient to work with the ``cut version" of infrasup-stability and we will show that
		$M$ is $C$-cut-stable with $C=8Q_M+1$. Fix integers $1\le l<k$,
		bounded functions $f\colon[\bN]^{l}\to M$,
		$g\colon[\bN]^{k-l}\to M$ and cuts $P,Q$ of $\{1,\dots,k\}$ of
		size $l$. Set
		\begin{align*}
			A &= \inf_{\nbar\in[\bN]^{k}}
			d\big(f(n_i\colon i\in P),g(n_i\colon i\in P^c)\big)\\
			B &=\sup_{\nbar\in[\bN]^{k}}
			d\big(f(n_i\colon i\in Q),g(n_i\colon i\in Q^c)\big)\ .
		\end{align*}
		Put $\bL:=\{2nk\colon n\in\bN\}$ and let $L_f$ and $L_g$ be the Lipschitz
		constants of $f_{\restriction_{[\bL]^{l}}}$ and of
		$g_{\restriction_{[\bL]^{k-l}}}$, respectively. We first show that
		$L_f\le4B$ and $L_g\le4B$. To see this, set
		$\bL'=\{nk\colon n\in\bN\}$ and let $(\mbar,\nbar)$ be an interlacing
		pair in $[\bL]^{l}$. By Lemma~\ref{lem:strongly-interlacing}
		there exists $\pbar\in [\bL']^{l}$ such that the pairs
		$(\mbar,\pbar)$ and $(\pbar,\nbar)$ are strongly interlacing. We
		can then choose $m_1<\dots<m_k$ and $p_1<\dots<p_k$ in $\bL'$ such
		that $\mbar=\{m_i\colon i\in Q\}$, $\pbar=\{p_i \colon i\in Q\}$ and $m_i=p_i$
		for all $i\in Q^c$. It follows that
		\begin{multline*}
			d(f(\mbar),f(\pbar))\le
			d\big(f(m_i \colon i\in Q),g(m_i \colon i\in Q^c)\big)\\
			+
			d\big(f(p_i\colon i\in Q),g(p_i\colon i\in Q^c)\big)
			\le 2B.
		\end{multline*}
		Similarly, we get $d(f(\pbar),f(\nbar))\le2B$ and hence
		$d(f(\mbar),f(\nbar))\le 4B$, as required. A similar argument shows that
		$L_g\le 4B$.
		
		Now fix $C>Q_M$. Since $M$ has property~$Q$, there exists
		$\bA\in\infin{\bL}$ such that $d(f(\mbar),f(\nbar))\le CL_f$ for all
		$\mbar,\nbar\in [\bA]^{l}$ and $d(g(\mbar),g(\nbar))\le CL_g$ for all
		$\mbar,\nbar\in [\bA]^{k-l}$. Now fix any $n_1<\dots<n_k$ in $\bA$. We
		get
		\begin{eqnarray*}
			A &\le & d\big(f(n_i\colon i\in P),g(n_i\colon i\in P^c)\big)\\
			&\le& d\big(f(n_i\colon i\in P),f(n_i\colon i\in Q)\big)\\
			&& +
			d\big(f(n_i\colon i\in Q),g(n_i\colon i\in Q^c)\big)\\
			&& +
			d\big(g(n_i \colon i\in Q^c),g(n_i\colon i\in P^c)\big)
			\le (8C+1)B.
		\end{eqnarray*}
		The result follows.
	\end{proof}

	The next corollary now follows from Proposition \ref{prop:property-Q-implies-cut-stable}, Proposition \ref{prop:interlaced-notCE->Q} and the fact that cut-stability trivially implies infrasup-stability. 
	\begin{coro}
		\label{cor:Q=cut=infrasup-stable}
		Let $(M,d)$ be a metric space. The following assertions are equivalent.
		\begin{enumerate}[(i)]
			\item $M$ has property~$Q$.
			\item $M$ is cut-stable.
			\item $M$ is infrasup-stable.
		\end{enumerate}
	\end{coro}

	\section{Coarsely universal Banach spaces}
	\label{sec:coarse-universality}
	
	The first objective of Kalton in \cite{Kalton2007} was to show that a reflexive Banach space cannot contain a coarse copy of every separable metric space. Note that it follows from Aharoni's embedding theorem that a metric space is coarsely universal for separable metric spaces if and only if it contains a coarse copy of $\co$. 
	
	The quick way to show that $\co$ does not coarsely embed into a reflexive space is to prove that reflexive spaces have property $Q$. In order to do so, we fix once and for all a nonprincipal ultrafilter $\cU$ on $\bN$ and we define a ``derivation-like'' operation for maps on $[\bN]^k$ that we shall denote by $\partial:=\partial_\cU$. Let $X$ be a Banach space and $k\in \bN$. For a bounded function $f\colon [\bN]^k \to X$ we define $\partial f\colon [\bN]^{k-1}\to X^{**}$ by letting for all $\nbar\in [\bN]^{k-1}$,
	\begin{equation*}
		\partial f(\nbar) := w^{*}-\lim_{n_k \in \cU}f(n_1,\dots,n_{k-1},n_k),
	\end{equation*}
	where the limit is meant for the weak$^*$ topology induced by $X^*$ on $X^{**}$. Note that for $1\le i \le k$, $\partial^i f$ is a bounded map from $[\bN]^{k-i}$ into $X^{(2i)}$ (the dual of order $2i$ of $X$) and that $\partial^k f$ is an element of $X^{(2k)}$. We first need to prove a series of simple lemmas about this
	operation $\partial$.
	
	\begin{lemm}
		\label{lem:real-concentration}
		Let $h\colon [\bN]^k \to \bR$ be a bounded map and $\eps>0$. Then,
		there is an infinite subset $\bM$ of $\bN$ such that for all $\nbar \in [\bM]^k$, 
		\begin{equation*}
			\abs{h(\nbar)-\partial^k h}< \eps.
		\end{equation*}
	\end{lemm}
	
	\begin{proof}
		Let $(\vep_j)_{j\in \bN}\subset (0,\infty)$ such that $\sum_{j=1}^\infty \vep_j \le \vep$. The set $\bM=\{m_1,\dots,m_i,\dots\}$ is built by induction on $i$ so that for any subset $\nbar$ of $\{m_1,\dots,m_i\}$ with $1\le \abs{\nbar}\le \min(i,k)$, we have
		$\abs{\partial^{k-\abs{\nbar}}h(\nbar)-\partial^k h}<\sum_{j=1}^i \vep_j$.
		
		\noindent For $i=1$, we easily pick $m_1$ such that
		$\abs{\partial^{k-1}h(m_1)-\partial^k h}< \vep_1.$
		
		\noindent Assume now that $m_1,\dots,m_i$ have been constructed. Then, for every $\nbar \subset \{m_1,\dots,m_i\}$ with $\abs{\nbar}\le k-1$, there is $\bM_{\nbar} \in \cU$
		such that for all $m \in \bM_{\nbar}$ with $m>m_i$, 
		\begin{align*}
			\abs{\partial^{k-\abs{\nbar}-1}h(\nbar,m)-\partial^k h} & \le \abs{\partial^{k-\abs{\nbar}-1}h(\nbar,m)- \partial^{k-\abs{\nbar}}h(\nbar)} + \abs{\partial^{k-\abs{\nbar}}h(\nbar)-\partial^k h}\\
			&\le \vep_{i+1} +\sum_{j=1}^i \vep_i.    
		\end{align*}
		If $A:=\bigcap_{\nbar} \bM_{\nbar}$ where $\nbar$ runs through the (finite) set of subsets of $\{m_1,\dots ,m_i\}$ satisfying $\abs{\nbar}\le k-1$, then $A \neq \emptyset$ since $\cU$ is an ultrafilter and we can pick $m_{i+1}\in A$. The verification that the set $\{m_1,\dots,m_{i+1}\}$ satisfies the desired properties is immediate.
	\end{proof}
	
	In the next lemma, we record an elementary identity for the derivation operation.
	
	\begin{lemm}
		\label{lem:derivation}
		Let $X$ be a Banach space and let $f\colon [\bN]^k\to X$ and $g\colon [\bN]^k \to X^*$ be two bounded
		maps. Define $f \otimes g\colon [\bN]^{2k} \to \bR$ by
		\begin{equation*}
			(f\otimes g)(n_1,\dots,n_{2k}) := \langle f(n_2,n_4,\dots,n_{2k}),g(n_1,n_3,\dots,n_{2k-1})\rangle.
		\end{equation*}
		Then, $\partial^2(f\otimes g) = \partial f \otimes \partial g.$
	\end{lemm}
	
	\begin{proof}
		It simply follows from the canonical identification of a Banach space inside its bidual, the definition of the $w^*$-topology and the observation that the $w^{*}$-limit is the classical limit in $\bR$.
	\end{proof}
	
	The interlacing relation appears naturally in the definition of the tensor product defined in Lemma \ref{lem:derivation}, and combining Lemma \ref{lem:derivation} with Lemma \ref{lem:real-concentration}, we obtain the following inequality for Lipschitz maps on the interlaced graphs and taking values in an arbitrary Banach space.
	
	\begin{theo}
		\label{thm:general-concentration} 
		Let $X$ be a Banach space and $\vep>0$. Then, for every $k\in \bN$ and every Lipschitz map $f\colon ([\bN]^k, \dIk)\to X$ there is an infinite subset $\bM$ of $\bN$ such that for all $\nbar\in [\bM]^k$,
		\begin{equation}
			\norm{ f(\nbar)}\le \norm{\partial^k f} + \Lip(f)+\vep.
		\end{equation}
	\end{theo}
	
	\begin{proof}
		For all $\nbar\in [\bN]^k$, we can find $g(\nbar)\in S_{X^*}$ such
		that $\langle f(\nbar),g(\nbar)\rangle= \norm{f(\nbar)}$. Since dual norms are weak$^*$ lower semicontinuous, $\|\partial^k g\|\le 1$. Then, by an iterated application of the
		Lemma \ref{lem:derivation} we get that
		\begin{equation*}
			\abs{\partial^{2k}(f\otimes g)} = \abs{\langle \partial^k f,\partial^k g\rangle}\le
			\norm{\partial^k f}.
		\end{equation*}
		Then, by Lemma \ref{lem:real-concentration}, there is an
		infinite subset $\bM_0$ of $\bN$ such that for all
		$\pbar\in [\M_0]^{2k}$, 
		\begin{equation}
			\label{eq:general-concentration}
			\abs{(f\otimes g)(\pbar)}\le
			\norm{\partial^k f} + \vep.
		\end{equation} 
		We then write $\bM_0 = \bM\cup \bM'$ for $\bM:=\{m_1,m_2,\dots,m_i,\dots\}$ and $\bM':=\{m'_1,m'_2,\dots,m'_i,\dots\}$ with $m_1<m'_1<m_2<m'_2<\dots<m_i<m'_i<\dots$. For all $\nbar=(m_{i_1},\dots,m_{i_k})\in [\bM]^k$, since $m_{i_1}<m'_{i_1}<m_{i_2}<m'_{i_2}<\dots<m_{i_k}<m'_{i_k}$, we have
		\begin{align*}
			\norm{f(\nbar)} = \langle f(\nbar),g(\nbar)\rangle & \le \abs{\langle
				f(m_{i_1},\dots,m_{i_k}) - f(m'_{i_1},\dots,m'_{i_k}), g(m_{i_1},\dots,m_{i_k})\rangle}  \\ 
			& \qquad + \abs{\langle
				f(m'_{i_1},\dots,m'_{i_k}), g(m_{i_1},\dots,m_{i_k})\rangle} \\
			& \stackrel{\eqref{eq:general-concentration}}{\le} \Lip(f) + \norm{\partial^k f} + \vep.
		\end{align*}
	\end{proof}
	
	The deduction from Theorem \ref{thm:general-concentration} that reflexive spaces have property $Q$ is elementary.
	
	\begin{coro}
		\label{cor:reflexive-Q}
		Every reflexive Banach space has property $Q$.
	\end{coro}
	
	\begin{proof}
		Let $f\colon ([\N]^k,\dIk) \to X$  be a Lipschitz map. Since $X$ is reflexive, $\partial^k f$ is an element of $X$. Then, by applying Theorem \ref{thm:general-concentration} to the map $f-\partial^k f$, we get that for any $\vep>0$ there is an infinite subset $\bM$ of $\bN$ such that $\diam f([\bM]^k) \le 2\,\Lip(f) + \vep$, which concludes the proof.
	\end{proof}
	
	\begin{rema}
		It follows from the proof of Corollary \ref{cor:reflexive-Q} that for every reflexive space $X$ one has $Q_X\le 2$.
	\end{rema}
	
	The next corollary follows from Lemma \ref{lem:coarse-uniform-infrasup-stable}, Corollary \ref{cor:Q=cut=infrasup-stable} and Corollary \ref{cor:reflexive-Q}.
	
	\begin{coro}
		\label{cor:Q-CE-co-reflexive}\,
		\begin{enumerate}
			\item If a Banach space coarsely, or uniformly, embeds into a space with property $Q$, then it must have property $Q$.
			\item The Banach space $\co$ does not coarsely, or uniformly, embed into a reflexive space.
		\end{enumerate}
	\end{coro}
	
	Thus, a Banach space that is coarsely universal (for separable metric spaces) cannot be reflexive. With a little bit more care, the following stronger result was shown by Kalton in \cite{Kalton2007}. Intuitively, Theorem \ref{thm:coarse-universal} says that if a Banach space $X$ contains for all $k \in \bN$,  uncountably many well-separated $1$-Lipschitz images of $\Ik$ and if $X$ coarsely embeds into a Banach space $Y$, then $Y$ cannot have all its iterated duals separable.
	
	\begin{theo}
		\label{thm:coarse-universal}
		Let $X$ be a Banach space satisfying the following property:
		
		$(\star)$ There exist an uncountable set $J$ and for every $j\in J$ and $k\in \bN$ a $1$-Lipschitz map $f_j^k\colon \Ik \to X$ such that 
		\begin{equation*}
			\lim_{k\to \infty}\, \inf_{i\neq j\in J}\,\inf_{\bM\in [\bN]^\omega} \, \sup_{\nbar\in [\bM]^k} \norm{f_i^k(\nbar) -f_j^k(\nbar)}=\infty.
		\end{equation*}
		If $X$ coarsely embeds into a Banach space $Y$, then there is $r\in \bN$ such that $Y^{(2r)}$, the $2r$-th iterated dual of $Y$, is nonseparable.
	\end{theo}
	
	\begin{proof}
		Assume that $Y$ is a Banach space such that for all $k\in \bN$, $Y^{(2k)}$ is separable. Let $h\colon X\to Y$ be a coarse embedding and $\{(f_j^k\colon \Ik \to X)_{k\in\bN}\}_{j\in J}$ be an uncountable collection of sequences of $1$-Lipschitz maps such that 
		\begin{equation}
			\label{eq:coarse-universal}
			\lim_{k\to \infty}\,  \inf_{i\neq j\in J}\, \inf_{\bM\in [\bN]^\omega} \, \sup_{\nbar\in [\bM]^k} \norm{f_i^k(\nbar) -f_j^k(\nbar)}=\infty.
		\end{equation}
		Fix $k\in \N$ and let $g_j:= h\circ f_j^k$ for all $j\in J$. Since $h$ maps bounded sets to bounded subsets of $Y$, $(g_j)_{j\in J}$ is an uncountable family of bounded maps from $\Ik$ to $Y$. We claim that there are $i\neq j\in J$ and $\bM\in [\bN]^\omega$ such that for all $\nbar\in[\bM]^k$, 
		\begin{equation}
			\label{eq2:coarse-universal}
			\norm{g_i(\nbar) - g_j(\nbar)}_Y\le \Lip(g_i) + \Lip(g_j) + 1.
		\end{equation}
		Indeed, since $Y^{(2k)}$ is separable and $\partial^k g_j\in Y^{(2k)}$ for all $j$ in the uncountable set $J$, there exist $i\neq j \in J$ such that $\norm{\partial^k g_i - \partial^k g_j}_{Y^{(2k)}}<1$. Then, we can apply Theorem
		\ref{thm:general-concentration} to $(g_i-g_j)$ and \eqref{eq2:coarse-universal} follows easily. It remains to observe that for all $ \nbar\in [\bM]^k$, 
		\begin{equation*}
			\rho_h\Big(\norm{f^k_i(\nbar) - f^k_j(\nbar)}_X\Big) \le \norm{g_i(\nbar) - g_j(\nbar)}_Y \le \Lip(g_i)+\Lip(g_j)+1\le 2\omega_h(1) + 1.
		\end{equation*}
		Now, letting $k$ tend to $\infty$, we see the above inequality together with  \eqref{eq:coarse-universal} contradict the fact that $\lim_{t\to\infty} \rho_h(t)=\infty$ and $\omega_h(1)<\infty$. 
	\end{proof}

	\begin{coro}
		If $\co$ coarsely embeds into a Banach space $Y$, then one of the iterated duals of $Y$ must be nonseparable.   
	\end{coro}
	
	\begin{proof}
		We simply need to verify that $\co$ satisfies property $(\star)$. As usual, let $(e_k)_{k=1}^\infty$ be the canonical basis of $\co$. For an infinite subset $A$ of $\bN$ we now define for all $n\in \bN$,
		\begin{equation*}
			s_A(n) := \sum_{i\le n \colon i\in A} e_i
		\end{equation*}
		and, for fixed $k\in \bN$ and all $\nbar:=(n_1,\dots,n_k)\in [\bN]^k$,
		\begin{equation*}
			f_A(\nbar) :=\sum_{i=1}^k s_A(n_i).
		\end{equation*}
		Note that $f_A$ is $1$-Lipschitz on $([\N]^k,\dIk)$. It is also not difficult to see that if $A\neq B$ and if $q\in A\triangle B$, then for all $\bM\in [\bN]^\omega$, we can always choose $\nbar=(n_1,\dots,n_k)\in [\bM]^k$ with $n_1> q$ so that $\norm{f_A^k(\nbar)-f_B^k(\nbar)}_\infty = k$. Therefore, $(f_A)_{A\in [\bN]^\omega}$ is an uncountable family satisfying property $(\star)$ for $\co$.
	\end{proof}
	\begin{rema} 
		Note that the bidual of the prototypical example of a Banach space satisfying property $(\star)$, namely $\co$, is nonseparable. However, since $\ell_1$ coarsely embeds into $\ell_2$, we cannot hope to prove Theorem \ref{thm:coarse-universal} under the mere assumption that the embedded space has a nonseparable bidual. It is shown in \cite{BLMS2020b} that a certain presence of $\ell_1$ in the embedded space is essentially the only obstruction in the following sense. If  $X$ is a separable Banach space with nonseparable bidual, not containing any isomorphic copy of $\ell_1$ and such that no spreading model generated by a normalized weakly null sequence in $X$ is equivalent to the $\ell_1$-unit vector basis, then a Banach space containing a coarse copy of $X$ cannot have all its iterated duals separable. This applies to $\co$ of course, but also, for instance, to the James tree space $\mathrm{JT}$ or its predual $\mathrm{JT}_*$ (see Exercise \ref{exe:JT} for a direct proof that $\mathrm{JT}$ satisfies property $(\star)$). 
	\end{rema}
	
	The following natural problems regarding coarse universality are open. 
	
	\begin{prob}
		\label{prob:coarse-universality-bidual}
		Can a Banach space with separable bidual contain a coarse copy of $\co$? 
	\end{prob}
	
	\begin{prob}
		\label{prob:coarse-universality-dual}
		Can a separable dual space contain a coarse copy of $\co$?    
	\end{prob}
	
	The same problems are open for coarse-Lipschitz or uniform embeddings, but we want to emphasize one in particular.    
	
	\begin{prob}
		\label{prob:CL-universality-dual}
		Can a separable dual space contain a coarse-Lipschitz copy of $\co$?    
	\end{prob}
	
	The reason we wanted to emphasize Problem \ref{prob:CL-universality-dual} is the following. It is worth pointing out (see Chapter \ref{chapter:Lipschitz-free}) that Kalton showed in \cite{Kalton2004} that $\co$ uniformly and coarse-Lipschitzly embeds into a Banach space which has the Schur and RNP properties and therefore does not contain any isomorphic (even bi-Lipschitz) copy of $\co$. We also recall that this space is a Lipschitz-free space over a snowflaking of $\co$. Consequently, it would follow from a negative answer to Problem \ref{prob:CL-universality-dual} that this space would provide an example of a separable Lipschitz-free space with RNP that does not linearly embed into a separable dual, and this would answer Problem \ref{pb:Lipschitz-free-RNP}.

	\section[Restricted coarse rigidity of reflexivity]{Restricted coarse rigidity and metric characterization of reflexivity}
	\label{sec:Q-applications}
	
	We now turn our attention to two spectacular results obtained by Kalton in \cite{Kalton2007} as a by-product of his investigation of the coarse universality problem. Kalton obtained a purely metric characterization of reflexivity within the class of Banach spaces that do not have $\ell_1$-spreading models. From this characterization, one can deduce that reflexivity is preserved under coarse embedding under the assumption that the embedded space does not have $\ell_1$-spreading models. Such an assumption is obviously necessary since $\ell_1$ coarsely embeds into $\ell_2$, for instance.
	
	Let us start with a fundamental result that bounds from below alternating arithmetic means of (subsequences of) bounded sequences in Banach spaces with property $Q$ in terms of the distance between their weak$^*$ limit points (in $X^{**}$) to $X$. Of course, this lower bound is only nontrivial for nonreflexive Banach spaces.
	
	\begin{theo}
		\label{thm:reflexivity-ABS-Q}
		Let $(x_n)_{n=1}^\infty$ be a bounded sequence in a Banach space $X$.
		Assume $X$ has property $Q$ and $x^{**}\in X^{**}$ is weak$^*$ limit point of $(x_n)_{n=1}^\infty$. Then, for any  $\vep >0$, there exists a subsequence $(y_n)_{n=1}^\infty$ of $(x_n)_{n=1}^\infty$ such that for all $k\in \bN$ and all $\nbar \in [\bN]^{2k}$,
		\begin{equation*}
			d(x^{**},X) \le Q_X(1+\vep)\frac{1}{k}\Big\|\sum_{i=1}^{2k} (-1)^i y_{n_i}\Big\|.
		\end{equation*}
	\end{theo}
	
	\begin{proof}
		Let $\vep>0$, $\theta:=d(x^{**},X)$ and assume, as we may, that $\theta>0$. Fix some parameters $\lambda \in (0,1)$ and $\alpha\in(0,1)$.
		First, we claim that: 
		\begin{claim}
			\label{clai1:reflexivity-ABS-Q}
			There is a subsequence $(z_n)_{n=1}^\infty$ of $(x_n)_{n=1}^\infty$ such that for any $r\in \bN$, any $u\in \conv\{z_n \colon n\le r\}$ and any $v\in \conv\{z_n \colon n> r\}$, 
			\begin{equation}
				\label{eq1:reflexivity-ABS-Q}
				\norm{u - v}\ge \lambda \theta.    
			\end{equation}
		\end{claim}
		
		Claim \ref{clai1:reflexivity-ABS-Q} is an elementary application of the functional analytic version of Helly's theorem that we leave as Exercise \ref{exe:Helly}. 
		
		\medskip
		Then, we claim that:
		
		\begin{claim}
			\label{clai2:reflexivity-ABS-Q}
			There is a subsequence $(y_n)_{n=1}^\infty$ of $(z_n)_{n=1}^\infty$ such that for any $k\in \N$, there exists $b_k \in (0,\infty)$ so that 
			\begin{equation}
				\label{eq2:reflexivity-ABS-Q}
				b_k-\alpha \theta \le \Big\|\sum_{j=1}^{2k} (-1)^j y_{n_j}\Big\| \le b_k,
			\end{equation}
			whenever $\alpha k\le n_1<\cdots <n_{2k}$.
		\end{claim}
		
		Claim \ref{clai2:reflexivity-ABS-Q} follows from a standard combination of a Ramsey argument and a diagonal argument (Exercise \ref{exe:Ramsey-diagonal} invites you to provide the details). 
		
		\medskip
		Fix one more parameter $\eta >0$. Note that $(y_n-y_m)_{m<n}$ is semi-normalized and tends weakly to $0$ as $m<n$ tend to $\infty$. Thus, applying a standard gliding hump argument and passing to a further subsequence, we can assume that 
		
		\begin{equation}
			\label{eq2bis:reflexivity-ABS-Q}
			\Big\|\sum_{j=1}^{2r} (-1)^j y_{n_j}\Big\| \le (1+\eta)b_k,
		\end{equation}
		whenever $\alpha k\le n_1<\cdots <n_{2r}$ and $1\le r \le k$.
		
		The proof can now be concluded swiftly. Set $\bM:=\{n\in \bN \colon n\ge \alpha k\}$ and fix $(n_1,\dots,n_{2k})\in [\N]^{2k}$. Observe that for any $x\in (0,\infty)$, necessarily $n_{\lceil x\rceil}\ge x$ and thus there is an integer $s\le \lceil\alpha k -1\rceil$ such that if we let $m_j=n_{s+j}$ for all $j\in \{1,\ldots,2k-s\}$ and pick integers $m_{2k}>m_{2k-1}>\dots>m_{2k-s}$, then $(m_1,\ldots,m_{2k})\in [\bM]^{2k}$. Therefore, there exists $\beta\in \{-1,1\}$ such that if we let $C:=\sup_{n\in \bN} \norm{x_n}$, then
		\begin{equation}
			\label{eq3:reflexivity-ABS-Q}
			\Big\|\sum_{j=1}^{2k} (-1)^jy_{n_j} + \beta \sum_{j=1}^{2k} (-1)^jy_{m_j}\Big\| \le 2 s C \le 2\alpha C k
		\end{equation} 
		It follows that
		\begin{equation}
			\label{eq4:reflexivity-ABS-Q}
			\Big\|\sum_{j=1}^{2k} (-1)^jy_{n_j}\Big\| \stackrel{\eqref{eq2:reflexivity-ABS-Q}\land\eqref{eq3:reflexivity-ABS-Q}}{\ge} b_k-\alpha \theta-2\alpha Ck.
		\end{equation}
		So far, we have not used any specific property of $X$, but we will now use property $Q$ to estimate $b_k$ from below and conclude the proof.
		Let $f(\nbar) := \sum_{j=1}^k y_{n_j}$, for $\nbar \in [\bM]^k$. It follows from \eqref{eq1:reflexivity-ABS-Q} that whenever $\nbar,\mbar \in [\bM]^k$ are such that $m_k<n_1$ we have by \eqref{eq1:reflexivity-ABS-Q} that 
		\begin{equation}
			\label{eq5:reflexivity-ABS-Q}
			\norm{f(\mbar) - f(\nbar)}\ge  \lambda \theta k.
		\end{equation}
		It follows from \eqref{eq2bis:reflexivity-ABS-Q} that $f$ is $(1+\eta)b_k$-Lipschitz as a map from $([\bM]^k,\dIk)$ to $X$ and invoking property $Q$ we can find $\mbar,\nbar\in [\bM]^k$ with $m_k<n_1$ satisfying 
		\begin{equation}
			\label{eq6:reflexivity-ABS-Q}
			\norm{f(\mbar)-f(\nbar)}\le \frac{Q_X}{\lambda} (1+\eta)b_k.
		\end{equation}
		Comparing \eqref{eq5:reflexivity-ABS-Q} and \eqref{eq6:reflexivity-ABS-Q}, we deduce that \begin{equation}
			\label{eq7:reflexivity-ABS-Q}
			Q_X (1+\eta)b_k\ge \theta k \lambda^2.
		\end{equation}
		
		Note that $X$ is not finite-dimensional and thus $Q_X>0$ (see Exercise \ref{exe:Q-FD}). Therefore and it follows from \eqref{eq4:reflexivity-ABS-Q} and \eqref{eq7:reflexivity-ABS-Q} that 
		\begin{equation*}
			\label{eq8:reflexivity-ABS-Q}
			\Big\|\sum_{j=1}^{2k} (-1)^jy_{n_j}\Big\| \ge \frac{\theta k}{Q_X}\Big(\frac{\lambda^2}{1+\eta} - \frac{\alpha Q_X}{k}- \frac{2\alpha C Q_X}{\theta}\Big).
		\end{equation*}
		and assuming, as we may, that $\alpha\in(0,1)$, $1-\lambda>0$ and $\eta \in (0,1)$ were initially chosen small enough, we get 
		\begin{equation*}
			\frac{1}{k}\Big\|\sum_{j=1}^{2k} (-1)^jy_{n_j}\Big\| \ge \frac{d(x^{**},X)}{(1+\eps) Q_X}.
		\end{equation*}
		
	\end{proof}

	\begin{rema} The statement of Theorem \ref{thm:reflexivity-ABS-Q} remains valid if we only assume that $x^{**}$ is a weak$^*$ cluster point of the sequence $(x_n)_{n=1}^\infty$ (see Exercise \ref{exer:cluster}). As you will see, this generality is not needed for the applications that we will now detail. 
	\end{rema}
	
	\begin{defi} 
		A Banach space $X$ has  the \emph{alternating Banach-Saks property} if every bounded sequence $(x_n)_{n=1}^\infty$ in $X$ has a subsequence $(y_n)_{n=1}^\infty$ so that $\frac1k \sum_{j=1}^{k}(-1)^j y_j$ tends in norm to $0$.
	\end{defi}
	
	Evidently, $\ell_1$ does not have the alternating Banach-Saks property. It is not difficult to see that a Banach space with a bounded sequence generating an $\ell_1$-spreading model cannot have the alternating Banach-Saks property (see Section  \ref{sec:spreading-models} for more on spreading models). B. Beauzamy \cite{Beauzamy1979} proved this is the only obstruction, namely, having the alternating Banach-Saks property is equivalent to not having a spreading model of a bounded sequence in $X$ that is equivalent to the canonical basis of $\ell_1$. In particular, a Banach space failing to have the alternating Banach-Saks property has a trivial type (or equivalently uniformly contains the $\ell_1^n$). The following metric characterization of reflexivity within the class of Banach spaces without $\ell_1$-spreading models is a consequence of Theorem \ref{thm:reflexivity-ABS-Q}.
	
	\begin{theo}
		\label{thm:ABS->Q=reflexive}
		Let $X$ be a Banach space with the alternating Banach-Saks property. Then, $X$ has property $Q$ if and only if $X$ is reflexive.
	\end{theo}
	
	\begin{proof} 
		We have already seen that every reflexive space has property $Q$. So, assume that $X$ has property $Q$ and the alternating Banach-Saks property. By separable determination of reflexivity, we can further assume that $X$ is separable. We will show that $X$ is reflexive by observing that for every $x^{**} \in X^{**}$ one has $d(x^{**},X)=0$. So, given $x^{**}\in X^{**}$, since $X$ is separable and does contain a linear copy of $\ell_1$, it follows from Rosenthal's theorem that there is a bounded sequence $(x_n)_{n=1}^\infty$ in $X$ such that $x^{**}$ is a weak$^*$ limit point of $(x_n)_{n=1}^\infty$. By the alternating Banach-Saks property there exists a subsequence $(y_n)_{n=1}^\infty$ of $(x_n)_{n=1}^\infty$ such that $\frac1k \norm{\sum_{j=1}^{k}(-1)^j y_j}$ tends to $0$ as $k$ tends to $\infty$. Invoking Theorem \ref{thm:reflexivity-ABS-Q}, we deduce the existence of a further subsequence $(z_n)_{n=1}^\infty$ such that for all $k\in \bN$ and all $\nbar\in [\bN]^{2k}$
		\begin{equation}
			d(x^{**},X) \le 4Q_X\frac{1}{2k}\Big\| \sum_{i=1}^{2k} (-1)^i z_{n_i} \Big\|,
		\end{equation}
		from which it follows that $d(x^{**},X)=0$.
	\end{proof}
	
	The following striking coarse rigidity result, within the class of Banach spaces without $\ell_1$-spreading models, follows from Theorem \ref{thm:ABS->Q=reflexive} and Corollary \ref{cor:Q-CE-co-reflexive}.
	
	\begin{coro}
		Let $X$ be a Banach space without $\ell_1$-spreading models. If $X$ coarsely embeds into a reflexive space, then $X$ is reflexive. 
	\end{coro}
	
	It follows from Theorem \ref{thm:ABS->Q=reflexive} that a nonreflexive Banach space with property $Q$ must have trivial type. This observation provides new examples of spaces failing property $Q$.
	
	\begin{coro}
		\label{cor:NR+TT->NQ}
		Every nonreflexive space with nontrivial type fails to have property $Q$. 
	\end{coro}
	
	Examples of such spaces were constructed by James \cite{James1978} and later by Pisier and Xu \cite{PisierXu1987}. Pisier and Xu also gave quasi-reflexive examples. However, James' nonreflexive but quasi-reflexive space $\James$ and its dual $\James^*$ have trivial type (see \cite[Theorem 2.b.8, Theorem 2.i.5]{FetterGamboa}) and Corollary \ref{cor:NR+TT->NQ} does not apply. 
	Nevertheless, Kalton was able to use Theorem \ref{thm:reflexivity-ABS-Q} and explicit computations in  $\James$ and  $\James^*$ to show that they also fail property $Q$.
	
	\begin{theo} 
		\label{thm:J-J*-fail-Q}
		The James space $\James$ and its dual $\James^*$ fail property $Q$.
	\end{theo}
	
	\begin{proof} 
		Let $(e_n)_{n=1}^\infty$ be the canonical basis of $\James$ and let $s_n :=\sum_{i=1}^n e_i$, for $n\ge 1$. It is well known that the sequence $(s_n)_{n=1}^\infty$ is weak$^*$ converging to an element $x^{**}\in \James^{**}\setminus \James$, which can be identified with the sequence which is constant and equal to $1$. A simple computation gives that for all $n_1<\dots <n_{2k}$ in $\bN$,
		\begin{equation*}
			\Big\| \sum_{i=1}^{2k} (-1)^i s_{n_i} \Big\|_{\James} = (2k)^{1/2}.
		\end{equation*}
		Since $\James$ is not reflexive, Theorem \ref{thm:reflexivity-ABS-Q} clearly implies that $\James$ fails property $Q$.
		
		In $\James^*$, the dual basis $(e_n^*)_{n=1}^\infty$ of $(e_n)_{n=1}^\infty$ is a bounded sequence that is converging to an element $z\in \James^{***}\setminus \James^*$.  It follows easily from the definition of $\norm{\cdot}_{\James}$ and the Cauchy-Schwarz inequality,  that for all $n_1<\dots <n_{2k}$ in $\bN$,
		\begin{equation}
			\Big\| \sum_{i=1}^{2k} (-1)^i e^*_{n_i} \Big\|_{\James^*} \le k^{1/2}.
		\end{equation}
		Again, Theorem \ref{thm:reflexivity-ABS-Q} implies that $\James^*$ fails property $Q$.
	\end{proof}
	
	\begin{rema} 
		It is proved in \cite{LPP2020} that although they fail property $Q$, the spaces $\James$ and $\James^*$ do not contain equi-coarse copies of the interlaced graphs $(\Ik)_{k\in \bN}$. Thus, the converse of Proposition \ref{prop:interlaced-notCE->Q} is false, i.e., property $Q$ and the non-equi-coarse embeddability of the interlaced graphs are not equivalent properties.
	\end{rema}

	\section{\texorpdfstring{Property $Q$, Szlenk index and separable duals}{Property, Szlenk index and separable duals}}
	\label{sec:Q-Szlenk}
	
	As we have already mentioned, it is an open problem whether $\co$  coarsely, or coarse-Lipschitzly, embeds into a separable dual, or even into a Banach space with separable bidual.
	In this section, we present two results from \cite{BLPP2020} which provide some insights about these problems. One result says that if $X^*$ contains a coarse-Lipschitz copy of $\co$, then the Szlenk index of $X$ cannot be too small, while the other result states that the dual of an asymptotic-$\co$ space does not coarsely contain $\co$.
	
	\begin{prop}
		\label{pro:weak*-linearization} 
		Let $X$ be a separable Banach space,  $k\in\bN$,  $\vep>0$ and $f \colon ([\bN]^k, \dIk)\to X^*$ a Lipschitz map. Then, there exist $\bM \in [\bN]^\omega$, a weak$^*$ null (rooted) tree $(x^*_{\mbar})_{\mbar\in[\bM]^{\le k}}$ in $X^*$ and constants $K_1,\dots,K_k$ in $[0,\Lip(f)]$ such that
		\begin{enumerate}[(i)]
			\item For all $\nbar \in [\bM]^k$, $f(\nbar) = \sum_{\mbar \preceq \nbar} x^*_{\mbar}.$
			\item For all $\mbar\in [\bM]^{\le k}\setminus\{\emptyset\}$, $\norm{x^*_{\mbar}} \le  \Lip(f)$.
			\item For all $i\in\{1,\dots,k\}$ and all $(n_1,\dots,n_i) \in [\bM]^i$,
			$$K_i\le  \|x^*_{(n_1,\dots,n_i)}\|\le K_i+\vep.$$
			\item If we let $\bM :=\{l_1<\dots <l_n<\dots\}$, $y^*_{\emptyset} :=z^*_{\emptyset}=0$ and,\\
			for $\nbar := (n_1,\dots,n_i) \in [\bN]^{\le k}\setminus \{\emptyset\}$,
			\begin{equation*}
				y^*_{\nbar} := x^*_{(l_{2n_1}, \dots ,l_{2n_i})} - x^*_{(l_{2n_1+1}, \dots ,l_{2n_i+1})}
			\end{equation*}
			and
			\begin{equation*}
				z^*_{\nbar} := \sum_{\mbar \preceq \nbar} y^*_{\mbar},  
			\end{equation*}
			then $(z^*_{\nbar})_{\nbar\in [\bN]^{\le k}}$ is a weak$^*$ continuous (rooted) tree in $\Lip(f)B_{X^*}$ such that for every $i \in \{1, \dots , k\}$ and every $(n_1,\dots,n_i) \in [\bN]^i$,
			\begin{equation*}
				\|z^*_{(n_1,\dots,n_i)} - z^*_{(n_1,\dots,n_{i-1})}\| \ge \frac{K_i}{2}.
			\end{equation*}
		\end{enumerate}
	\end{prop}
	
	\begin{proof}
		The construction of a weak$^*$ null (rooted) tree $(x^*_{\mbar})_{\mbar\in[\bM]^{\le k}}$ in $X^*$ satisfying $(i)-(ii)$ goes verbatim as in the proof of Lemma \ref{lem:linearization} albeit we are using weak$^*$ compactness in the dual of a separable space instead of weak compactness in a reflexive space. The stabilization property in $(iii)$ is an immediate consequence (after passing to a further subtree) of Ramsey's theorem and the compactness of $[0,\Lip(f)]$.
		Using the weak$^*$ lower semi-continuity of $\norm{\cdot}_{X^*}$, we may also assume, by passing to a further subtree that for every $i \in \{1, \dots, k\}$ and all $(n_1,\dots,n_i) \in [\bN]^i$,
		\begin{equation*}
			\| x^*_{(l_{2n_1}, \dots ,l_{2n_i})} - x^*_{(l_{2n_1+1}, \dots ,l_{2n_i+1})}\| \ge \frac{K_i}{2}.
		\end{equation*}
		The last assertion $(iv)$ follows then from the observation that for every $\nbar \in [\bN]^i$, $(l_{2n_1}, \dots ,l_{2n_i})$ and $(l_{2n_1+1}, \dots ,l_{2n_i+1})$ interlace and  
		\begin{equation}
			\Big\| \sum_{j=1}^{i} y^*_{(n_1,\dots ,n_j)} \Big\| \stackrel{(i)}{=} \norm{ f(l_{2n_1}, \dots ,l_{2n_i}) - f(l_{2n_1+1}, \dots ,l_{2n_i+1})}.
		\end{equation}
	\end{proof}
	
	\begin{rema}
		It is easy to see that $(i)$-$(ii)$ relies only on a very specific property of the interlaced metric and hold for other metrics such as the Hamming metric, for instance (see Exercise \ref{exe:weak*-linearization}).
	\end{rema}
	
	The reason why item $(iv)$ in Proposition \ref{pro:weak*-linearization} is relevant to our purpose comes from the following lemma that can be easily derived from the definition of the Szlenk index.
	
	\begin{lemm}
		\label{lem:separated-trees-Szlenk}
		Let $X$ be a Banach space and assume that $(x^*_{\nbar})_{\nbar\in [\bN]^{\le k}}$ is a weak$^*$ continuous (rooted) tree in $B_{X^*}$ such that there exist $i_1<\dots<i_l$ in $\{0, \dots , k-1\}$ and $\vep_{i_1},\dots,\vep_{i_l}>0$ satisfying for all $s\in\{1,\dots,l\}$ and $\nbar\in[\bN]^{i_s}$,
		\begin{equation*}
			\limsup_{t\to \infty}\|x^*_{(\nbar,t)} - x^*_{\nbar}\| > \vep_{i_s}.
		\end{equation*}
		Then, $x^*_{\emptyset}\in s_{\vep_{i_l}}\dots s_{\vep_{i_1}}(B_{X^*})$.
	\end{lemm}
	
	We already know that $\ell_1$ has property $Q$ since it is a stable space. But, $\ell_1=\co^*$ is also the dual of a Banach space with summable Szlenk index. The next result, established in \cite{BLPP2020}, provides an alternative proof of the fact that $\ell_1$  has property $Q$.
	
	\begin{theo}
		\label{thm:summable-Szlenk->dual-Q}
		Let $X$ be a Banach space with summable Szlenk index. Then, $X^*$ has property $Q$.
	\end{theo}
	
	\begin{proof} 
		Let $f \colon [\bN]^k \to X^*$ be a map. If one wants to use Proposition \ref{pro:weak*-linearization}, we need to observe that there exists a separable subspace $Y$ of $X$ such that the closed linear span of $f([\bN]^k)$ isometrically embeds into $Y^*$. Indeed, since $[\bN]^k$ is countable, the closed linear span of $f([\bN]^k)$ is a separable subspace of $X^*$; let us call it $Z$. There exists a separable subspace $Y$ of $X$ such that for all $x^* \in Z$, $\norm{x^*}_{X^*} = \sup_{y\in B_Y} \abs{x^*(y)}$. So, $Z$ linearly and isometrically embeds into $Y^*$. Therefore, we may assume that $X$ is separable. Fix $\vep\in (0,\frac{1}{2k})$ and let $\bM\in [\bN]^\omega$ and $K_1,\dots,K_k$ be given by Proposition \ref{pro:weak*-linearization}. It follows from item $(iv)$ that
		\begin{equation}
			0 \in s_{\frac{K_1}{2}}\ldots s_{\frac{K_k}{2}} (B_{X^*}).
		\end{equation}
		Since the Szlenk index of $X$ is summable, we deduce that $\sum_{i=1}^k K_i \le 2C$, where $C$ is the ``summable Szlenk index constant'' of $X$.
		Then, we deduce from items $(ii)$ and $(iii)$ of Proposition \ref{pro:weak*-linearization} that
		\begin{equation}
			\diam (f([\bM]^k)) \le 2\sum_{i=1}^k K_i+2k\vep \le 4C+1.
		\end{equation}
	\end{proof}
	
	This can be restated as follows. 
	\begin{coro}
		\label{cor:asymp-co->dual-Q}
		The dual of an asymptotic-$\co$ space has property $Q$. In particular, it cannot contain a coarse copy of $c_0$. 
	\end{coro}
	
	\begin{proof} The first sentence is a consequence of Theorems \ref{thm:N-infty-theorem-full} and \ref{thm:summable-Szlenk->dual-Q}. Then, the second sentence follows from Aharoni's theorem. 
	\end{proof}
	
	\begin{rema}
		The space $\co$, as well as all its subspaces, are asymptotic-$\co$, but of course not every asymptotic-$\co$ Banach space is isomorphic to a subspace of $\co$. A typical example is $\Tsi^*$, the asymptotic-$\co$ Tsirelson space. Thus, Corollary \ref{cor:asymp-co->dual-Q} applies to the asymptotic-$\ell_1$ Tsirelson space $\Tsi$, since $\Tsi=(\Tsi^{*})^{*}$, providing an alternative proof that $T$ has property $Q$. It is worth pointing out that $J(\Tsi^*)$, the Jamesification of $\Tsi^*$ (for its canonical basis), is a nonreflexive quasi-reflexive Banach space that is asymptotic-$\co$ (see \cite{BLMS2020} and references therein). It can be shown that $J(\Tsi^*)^*$ is not stable and hence Kalton's techniques do not apply in this case, but nevertheless Corollary \ref{cor:asymp-co->dual-Q} applies to the dual of $J(\Tsi^*)$.    
	\end{rema}

	Recall that a separable Banach space has a separable dual if and only if its Szlenk index is countable. The following result asserts that if the interlaced graphs equi-bi-Lipschitzly embed into a dual space $X^*$, then the Szlenk index of $X$ is larger than the first infinite ordinal. 
	
	\begin{theo}
		\label{thm:EL-interlaced->Szlenk-large}
		Let $X$ be a Banach space. If $\sup_{k\in \bN}\cdist{X^*}(\Ik)<\infty$, then $\Sz(X) > \omega$.
	\end{theo}
	
	\begin{proof} 
		Without loss of generality, we may assume that $X$ is separable and that there exists $A\in (0,1]$ so that for any $k\in \bN$ there exists $f_k \colon ([\bN]^k, \dIk) \to X^*$ such that for all $\mbar,\nbar \in [\bN]^k$, 
		\begin{equation}
			A\dIk(\mbar,\nbar)\le \norm{f(\mbar) - f(\nbar)} \le \dIk(\mbar,\nbar).
		\end{equation}
		For a fixed $k\in \N$ and a given $\vep\in (0,\frac{A}{4})$, we consider $\bM\in [\bN]^\omega$ given by Proposition \ref{pro:weak*-linearization}. Let $\mbar$ and $\nbar$ in $[\bN]^k$ such that $m_k<n_1$. By triangle inequality, we have
		\begin{equation}
			\Big\|\sum_{\emptyset \prec \sbar \preceq \mbar} x^*_{\sbar}\Big\| + \Big\| \sum_{\emptyset \prec \sbar \preceq \nbar} x^*_{\sbar} \Big\| \ge \norm{ f(\mbar) - f(\nbar)} \ge Ak.
		\end{equation}
		It then follows from item $(iii)$ of Proposition \ref{pro:weak*-linearization} that
		$2\sum_{i=1}^{k} K_i \ge Ak-2k\vep$.
		Now, if we set $[k]:=\{1, \dots,k\}$, $J := \{ i\in [k]\colon K_i > \frac{A}{8}\} $ and  $N := \abs{J}$, we have that
		\begin{equation}
			\frac{Ak}{2}-k\vep \le \sum_{i=1}^{k} K_i = \sum_{i\in [k] \setminus J} K_i + \sum_{j\in J} K_j \le \frac{A}{8}k + N \Lip(f).
		\end{equation}
		From our choice of $\vep$ and since $\Lip(f)\le 1$, it follows that $N \ge \frac{Ak}{8}$.
		Finally, we deduce from item $(iv)$ in Proposition \ref{pro:weak*-linearization} and Lemma \ref{lem:separated-trees-Szlenk} that $0\in s_{\frac{A}{16}}^N(B_{X^*})$
		and therefore that $\Sz(X,\frac{A}{16}) \ge \frac{Ak}{8}$. Since $k$ was arbitrary and by weak$^*$ compactness $Sz(X,\frac{A}{16})$ is not a limit  ordinal, we deduce that $\Sz(X,\frac{A}{16})>\omega$. This concludes the proof. 
	\end{proof}
	
	\begin{rema} 
		In Theorem \ref{thm:EL-interlaced->Szlenk-large}, the assumption $\sup_{k\in \bN}\cdist{X^*}(\Ik)<\infty$ can be upgraded to the assumption that the interlaced graphs are equi-coarse-Lipschitzly embeddable into $X^*$ (see Exercise \ref{exe:EL=ECL-interlaced}). 
	\end{rema}
	
	\begin{rema} 
		Theorem \ref{thm:EL-interlaced->Szlenk-large} is optimal in the following sense. It was proved in \cite{BLPP2020} that there is a separable Banach space $X$ with $\Sz(X)=\omega^2$ and such that $\sup_{k\in \bN}\cdist{X^*}(\Ik)<\infty$. The Banach space $X$ is constructed as a Lipschitz-free space.
	\end{rema}

	\section{\texorpdfstring{Property $Q^{\omega_1}$ and coarse embeddings into $\ell_\infty$}{Property and coarse embeddings into}}
	\label{sec:Q-omega-1}
	
	In \cite{Kalton2012b}, Kalton considered an uncountable version of property $Q$ which involves an uncountably branching version of the interlaced graphs. Kalton's motivation was the study of the following important question raised by Benyamini and Lindenstrauss in \cite{BenyaminiLindenstrauss2000}: is every  Banach space a Lipschitz retract of its bidual? Kalton answered negatively by showing the existence of a nonseparable Banach space that is not a uniform retract of its bidual (Theorem 4.5 in \cite{Kalton2012b}). However, the following question remains open. 
	
	\begin{prob}
		\label{prob:Lipschitz-retraction}
		Let $X$ be a separable Banach space. Is there a Lipschitz retraction of $X^{**}$ onto $X$?
	\end{prob}  
	
	In this section, we shall not address this question, and we only introduce property $Q^{\omega_1}$ and give one application from \cite{Kalton2012b} to coarse embeddings into $\ell_\infty$. Of course, since $\ell_\infty$ is isometrically universal for separable metric spaces, the results will be meaningful for nonseparable metric spaces. Let us first define the relevant objects and notions.
	
	Let $\Omega_1:=[1,\omega_1)$ be the set of countable ordinals. For a subset $\Theta$ of $\Omega_1$ and for $k\in \bN$, we let
	\begin{equation}
		[\Theta]^k :=\{\alphabar=(\alpha_1,\dots,\alpha_k) \colon \alpha_1<\dots <\alpha_k \in \Theta\}.
	\end{equation}
	We will now equip $[\Omega_1]^k$ with a graph structure by declaring $\alphabar\neq\betabar$ to be adjacent if they interlace that is if
	\begin{equation*}
		\alpha_1\le\beta_1\le\cdots\le\alpha_n\le\beta_n\ \ {\rm or}\ \
		\beta_1\le\alpha_1\le\cdots\le \beta_n\le\alpha_n.
	\end{equation*}
	We denote by $\Iuk$ the vertex set $[\Omega_1]^k$ equipped with this graph structure and by $\dIuk$ the shortest path metric on $\Iuk$. Note that $\diam(\Iuk)=k$. For any $\alphabar,\betabar\in[\Omega_1]^k$, we write $\alphabar <\betabar$ if $\alpha_k<\beta_1$. Note that in that case $\dIuk(\alphabar,\betabar)=k$. %For $\alphabar :=(\alpha_1,\dots,\alpha_k)\in [\Omega_1]^k$ and $1\le i \le k$, $\abar_{|_i}=(\alpha_1,\ldots,\alpha_i)$ and $\abar_{|_0}=\emptyset$ is the empty sequence.
	
	The analogy with property $Q$ of the following definition is clear.
	
	\begin{defi}
		\label{def:Q-omega-1}Let $(M,d_M)$ be a metric space. We say that $M$ has \emph{property $Q^{\omega_1}$} if there exists $C\ge 1$ such that for every $k \in \bN$ and every Lipschitz map $f \colon ([\Omega_1]^k,\dIuk) \to M$, there exists an {\bf uncountable} subset $\Theta$ of $\Omega_1$ such that
		\begin{equation}
			\sup_{\alphabar,\betabar \in [\Theta]^k}  d_M(f(\alphabar),f(\betabar)) \le  C\Lip(f).
		\end{equation}
		We define $Q_M^{\omega_1}$ to be the infimum of such constants $C$, when $M$ has property $Q^{\omega_1}$ and $Q_M^{\omega_1}=\infty$ otherwise.
	\end{defi}
	
	The following proposition is elementary, and its proof is the same as in the countable setting and left to the reader as Exercise \ref{exe:basics-Q-omega-1}.
	\begin{prop}
		\label{pro:basics-Q-omega-1} \,
		\begin{enumerate}[(i)]
			\item If $(M,d_M)$ has property $Q^{\omega_1}$, then the sequence of graphs $(\Iuk)_{k\in \N}$ does not equi-coarsely embed into $M$.
			\item Assume that $X$ is a Banach space that coarsely embeds into $(M,d_M)$ and that $M$ has property $Q^{\omega_1}$. Then, $X$ has property $Q^{\omega_1}$.
		\end{enumerate}
	\end{prop}
	
	It is not difficult to see that the space of continuous real-valued functions on the compact space $[1,\omega_1]$ (for the order topology) fails property $Q^{\omega_1}$. Note that $C([1,\omega_1])$ is nonseparable. The argument is similar to the one showing that $\co$ fails $Q$ using the summing basis of $\co$.  
	
	\begin{prop}
		\label{pro:not-Q-omega-1}
		The space $C([1,\omega_1])$ fails property $Q^{\omega_1}$.
	\end{prop}
	
	\begin{proof} 
		Since $[1,\gamma]$ is a clopen set for every $\gamma<\omega_1$, the indicator function $\indicator{[1,\gamma]}$ is a continuous function on $[1,\omega_1]$. These indicator functions will play a role similar to the summing basis in $\co$.
		Let $k\in \bN$ and define $f_k \colon  [\Omega_1]^k \to C([1,\omega_1])$ by $f_k(\alphabar) := \sum_{i=1}^k \indicator{[1,\alpha_i]}$ whenever $\alphabar := (\alpha_1,\dots,\alpha_k) \in [\Omega_1]^k$.
		It is clear that $f_k$ is $1$-Lipschitz for the interlaced metric on $[\Omega_1]^k$. On the other hand, for any uncountable subset $\Theta$ of $\Omega_1$, we can find $\alphabar, \betabar \in [\Theta]^k$ with $\alpha_k <\beta_1$ and thus such that $\norm{ f_k(\alphabar) - f_k(\betabar) }_\infty = k$. This proves that $C([1,\omega_1])$ fails property $Q^{\omega_1}$. 
	\end{proof}
	
	The following corollary, which will be useful later, follows from the fact from \cite{Par1963} that  $C([1,\omega_1])$ is linearly isometric to a subspace of $\ell_\infty/c_0$. However, we indicate a direct proof in Exercise \ref{exe:l/c}.
	
	\begin{coro}
		\label{cor:not-Q-omega-1}
		The space $\ell_\infty/\co$ fails property $Q^{\omega_1}$. 
	\end{coro}
	
	The presence of the sequence $(\Ik)_{k\in \bN}$ in a Banach space implies the nonseparability of some of its iterated duals. The sequence $(\Iuk)_{k\in \bN}$ is made up of much bigger graphs, and one could expect a much stronger nonseparability conclusion for spaces containing it. To confirm this intuition, we need to introduce the notions of small, large and very large subsets of $[\Omega_1]^k$ and prove a Ramsey-type theorem for very large subsets. 
	
	For any subset $\Theta$ of $\Omega_1$ (even for $\Theta=\emptyset$), we agree that $[\Theta]^0=\{\emptyset\}$. For $k\in \bN$ and a subset $A$ of $[\Omega_1]^k$, we define the following derivation
	\begin{equation*}
		\partial A := \big\{(\alpha_1,\dots,\alpha_{k-1}) \in [\Omega_1]^{k-1} \colon \{\beta\colon (\alpha_1, \dots, \alpha_{n-1}, \beta)\in A\}\ \text{is
			uncountable}\big\},
	\end{equation*}
	and for $1\le i\le k-1$, $\partial^{i+1}A=\partial(\partial^iA)$. 
	Then, we say that $A\subset [\Omega_1]^k$ is \emph{large} if $\emptyset\in\partial^k A$. A subset $A\subset [\Omega_1]^k$ that is not large, i.e. $\partial^k A=\emptyset$, is called \emph{small} and we will say that $A\subset[\Omega_1]^k$ is \emph{very large} if its complement is small. 
	
	It is not immediately clear that a very large set is large, but this follows from the fact that being small is stable under taking countable unions.
	
	\begin{lemm}
		\label{lem:small-union} 
		Let $(A_n)_{n=1}^\infty$ be a sequence of small subsets of $[\Omega_1]^k$. Then, $\bigcup_{n=1}^\infty A_n$ is small. In particular, a very large set is large, and a countable intersection of very large sets is very large.
	\end{lemm}
	
	\begin{proof} 
		It is trivial that $\bigcup_{n=1}^\infty \partial A_n \subseteq \partial(\bigcup_{n=1}^\infty A_n)$ and it follows from a pigeonhole argument that we also have $\partial(\bigcup_{n=1}^\infty A_n)\subseteq\bigcup_{n=1}^\infty \partial A_n$. Therefore, $\partial(\bigcup_{n=1}^\infty A_n)=\bigcup_{n=1}^\infty \partial A_n$ and iterating gives that $\partial^k(\bigcup_{n=1}^\infty A_n) = \bigcup_{n=1}^\infty \partial^k A_n =\emptyset$ since the $A_n$ are all small.
		Obviously, $[\Omega_1]^k$ is large, and according to the previous observation, a set and its complement cannot be both small, thus easily yielding the second part of the statement.
	\end{proof}
	
	If $A$ is large, then it means that there is an $\omega_1$-branching subtree of height $k$ whose leaves are all in $A$, but it does not necessarily mean that there is an uncountable subset $\Theta$ of $\Omega_1$ such that $[\Theta]^k\subseteq A$. The Ramsey-type result below states that this stronger statement is true whenever $A$ is very large.
	
	\begin{lemm}
		\label{lem:very-large-Ramsey} 
		If $A$ is a very large subset of $[\Omega_1]^k$, then there is an
		uncountable subset $\Theta$ of $\Omega_1$ such that $[\Theta]^k\subseteq A$.
	\end{lemm}
	
	\begin{proof} 
		The argument is a simple (but clever) application of Zorn's lemma. Consider the set $\cC := \{S\subseteq \Omega_1 \colon \forall i\in \{0,\dots,k\}, [S]^{i}\cap \partial^{k-i}(\Omega_1 \setminus A)=\emptyset\}$ partially ordered by inclusion. For all $1\le i \le k$, $[\emptyset]^{i}=\emptyset$ and, since $A$ is very large $\partial^{k}(\Omega_1 \setminus A)=\emptyset$. So, $\emptyset \in \cC$ and hence $\cC$ is not empty. Moreover, it is easy to see that every chain in $\cC$ has an upper bound in $\cC$ (we are only considering finite tuples in $S$). Thus, by Zorn's Lemma, we can pick a maximal subset $\Theta$ of $\Omega_1$ such that for all $i\in \{0,\ldots,k\}$, $[\Theta]^i \cap \partial^{k-i}(\Omega_1 \setminus A)=\emptyset$ and in particular $[\Theta]^k\subseteq A$. It is then enough to show that $\Theta$ is uncountable. Observe first that by definition of $\Theta$, we have that for any $0\le i\le k-1$ and any $\alphabar \in [\Theta]^i$, the set $\{\beta \colon (\alphabar,\beta) \in \partial^{k-i-1}(\Omega_1 \setminus A)\}$ is at most countable and denote by $h_i(\alphabar)<\omega_1$ the supremum of this set (the supremum of the empty set being $0$). Assume, aiming for a contradiction, that $\Theta$ is countable, we can deduce that $\sigma$, the supremum of all $h_i(\alphabar)$ for $0\le i\le k-1$ and $\alphabar  \in [\Theta]^i$ is countable. It is straightforward to verify that $\sigma \ge \sup \Theta$ and that the set $\Theta \cup \{\sigma+1\}$ contradicts the maximality of $\Theta$.
	\end{proof}
	
	We can now prove the following theorem, which will surely not surprise the reader.
	
	\begin{theo}
		\label{thm:separable-Q-omega-1}
		Let $M$ be a separable metric space. Then, $M$ has property $\cal Q^{\omega_1}$.
	\end{theo}
	
	\begin{proof}
		Let $(M,d_M)$ be a separable metric space and $f \colon ([\Omega_1]^k,\dIuk) \to M$ be a Lipschitz map. According to Lemma \ref{lem:very-large-Ramsey} it is sufficient to find a subset $F$ of $M$ of diameter at most $\Lip(f)$ and such that $f^{-1}(F)$ is very large. Also,  note that we would be done if we could choose the set $F$ such that the following claim holds:
		\begin{claim}
			Given $x,y\in F$ and $\vep>0$ such that $d_M(x,y) > \Lip(f)+2\vep$, the sets $f^{-1}(B_M(x,\vep))$ and $f^{-1}(B_M(y,\vep))$ contain vertices $\alphabar$ and $\betabar$ (one in each) that interlace. 
		\end{claim}
		Indeed, assuming for the sake of a contradiction that $\diam(F)> \Lip(f)$ and letting $\vep>0$ and $x,y\in F$ such that $d_M(x,y) > \Lip(f)+2\vep$, we would then deduce that $d_M(f(\alphabar),f(\betabar)) \ge d_M(x,y) - d_M(f(\alphabar),x) - d_M(f(\alphabar),y)> \Lip(f)$, contradicting the definition of $\Lip(f)$. 
		
		\medskip
		Finding interlacing vertices in two subsets $A$ and $B$ of $[\Omega_1]^k$ is easy if the two subsets are large. Indeed, if $A$ and $B$ are two large subsets of $[\Omega_1]^k$, we first pick $\alpha_1 \in \partial^{k-1}A$. Since $\partial^{k-1}B$ is uncountable, we can pick $\beta_1 \in \partial^{k-1}B$ such that $\beta_1 >\alpha_1$. Now, for similar reasons, we can pick $\alpha_2>\beta_1$ and such that $(\alpha_1,\alpha_2)\in \partial^{k-2}A$ and we can continue this procedure to get $\alphabar \in A$ and $\betabar \in B$ so that $\alphabar$ and $\betabar$ interlace. Thus, it is enough that our construction of the set $F$ guarantees that $f^{-1}(B_M(x,\vep))$ is large for all $x \in F$ and all $\vep>0$.  This is where the separability of $M$ comes into play via the Lindel\"{o}f property of second countable spaces. Here are the details. Recall that by definition, $f^{-1}(F)$ is very large if $f^{-1}(F)^c = f^{-1}(F^c)$ is small.
		Let $\Sigma$ be the set of open subsets $V$ of $X$ such that $f^{-1}(V)$ is small and consider the open set $U := \bigcup\{V \colon V\in \Sigma \}$. Since $U$ is separable and $\Sigma$ is an open cover of $U$, there exists a sequence $(V_n)_{n=1}^\infty$ of sets in $\Sigma$ such that $U = \bigcup_{n=1}^\infty V_n$. Thus, by Lemma \ref{lem:small-union}, $f^{-1}(U)$ is small. Moreover, by construction, for any open subset $V$ of $X$, $f^{-1}(V)$ is small if and only if $V\subset U$. The set $F:=U^c$, which is closed and clearly nonempty, will do the job. First, $f^{-1}(F)=f^{-1}(U)^c$ is very large and for all $x\in F$ and $\vep>0$, $B_M(x,\vep)$ is an open set that is not contained in $U$, so its  inverse image by $f$ is large. 
	\end{proof}
	
	A coordinate-wise application of Theorem \ref{thm:separable-Q-omega-1} gives the following result.
	
	\begin{theo} 
		\label{thm:linfty-Q-omega-1}
		The Banach space $\ell_\infty$ has property $\cal Q^{\omega_1}$.
	\end{theo}
	
	\begin{proof} 
		Let $f \colon ([\Omega_1]^k,\dIuk) \to \ell_\infty$ be a Lipschitz map and for $\alphabar \in [\Omega_1]^k$ let $f(\alphabar) := ((f_n(\alphabar))_{n=1}^\infty$, where each function $f_n$ is $\Lip(f)$-Lipschitz. It follows from the separability of $\bR$ and Theorem \ref{thm:separable-Q-omega-1}, or more accurately its proof, that for every $n\in \bN$, there is a very large subset $B_n$ of $[\Omega_1]^k$ such that $\diam(f_n(B_n))\le \Lip(f)$. Therefore, there exists $\xi_n \in \bR$ such that $f_n(B_n)\subset [\xi_n-\frac{\Lip(f)}{2},\xi_n+\frac{\Lip(f)}{2}]$. From Lemma \ref{lem:small-union}, we deduce that $B:=\cap_{n=1}^\infty B_n$ is also very large. In particular, $B$ is not empty, which implies that $\xi :=(\xi_n)_{n=1}^\infty \in \ell_\infty$ and we have that $f(B)$ is included in $\xi+\frac{\Lip(f)}{2}B_{\ell_\infty}$. Since $B$ is very large, the conclusion follows from Lemma \ref{lem:very-large-Ramsey}.
	\end{proof}

	\begin{rema}
		If $M$ is a separable metric space or the space $\ell_\infty$, the proofs of Theorem \ref{thm:separable-Q-omega-1} and Theorem \ref{thm:linfty-Q-omega-1} give that $Q_M^{\omega_1}\le 1$.
	\end{rema}
	
	The following corollary is an immediate consequence of Proposition \ref{pro:basics-Q-omega-1}, Proposition \ref{pro:not-Q-omega-1}, Corollary \ref{cor:not-Q-omega-1} and Theorem \ref{thm:linfty-Q-omega-1}.
	
	\begin{coro}
		\label{cor:ellinftyquotient-not-coarse-into-ellinfty}
		$C([1,\omega_1])$ and $\ell_\infty/\co$ do not coarsely embed into $\ell_\infty$.
	\end{coro}

	\section{Notes}
	\label{sec:notes-Q}
	
	Stable norms were introduced, originally for separable Banach spaces, by Krivine and Maurey \cite{KrivineMaurey1981} to address a long-standing open problem in Banach space theory, namely, whether every infinite-dimensional Banach space contains an isomorphic copy of $\co$ or $\ell_p$, for some $p\in[1,\infty)$. In the early 1970s, Tsirelson \cite{Tsirelson1974} built an example of a space that does not have this property. However, Krivine and Maurey showed that for stable Banach spaces, such Tsirelson-type construction cannot happen since every stable Banach space contains an isomorphic copy (actually almost isometric) of $\ell_p$ for some $p\in[1,\infty)$. Krivine-Maurey's result is one of many pieces of evidence that stable Banach spaces have a much more regular structure than arbitrary Banach spaces. We point the reader to \cite{Guerre-book92} for an extensive account of the theory of stable Banach spaces.
	
	Krivine-Maurey stability theorem can also be used to provide examples of Banach spaces, other than $\co$, without any stable renormings. For instance, Tsirelson-like spaces, which do no contain any isomorphic copy of $\ell_p$ for any $p\in[1,\infty)$, do not admit any equivalent stable norm. Note also that if a Banach space contains an isomorphic copy of $\ell_p$ for some $p\in[1,\infty)$, then it will contain an unconditional basic sequence. Hence, Banach spaces with no unconditional basic sequences, in particular hereditarily indecomposable Banach spaces (see \cite{Maurey03} for a discussion of how to construct such spaces), also do not admit any equivalent norm that is stable. 
	
	In \cite{Raynaud83}, Raynaud showed that the conclusion of the Krivine-Maurey theorem is valid for any Banach space whose unit ball admits a uniformly equivalent stable metric. It was well known that spreading basic sequences in stable Banach spaces are unconditional, and Raynaud showed that this property still holds if we merely assume that the unit ball of the Banach space embeds uniformly into a stable space. 
	
	To find more examples of Banach spaces that are not stable, we need to resort to noncommutative $L_p$-spaces. We refer to \cite{PisierXu_HB03} for a thorough discussion of noncommutative $L_p$-spaces.
	
	It seems that the natural extension of the notion of stability to arbitrary metric spaces, as recalled above, was first studied by Garling \cite{Garling82}.
	
	The notion of infrasup stable was introduced by the first author in \cite{Baudier2022} under a different and more elusive terminology: upper stable. The equivalence between infrasup stability and property $\cQ$ is from \cite{BSZ2024} (where the terminology upper stable was still in effect). This paper also contains an alternative proof of the fact that reflexive spaces have property $\cQ$. It is shown there that every reflexive Banach space is infrasup stable by relating the latter property with the asymptotic structure of a Banach space.
	
	In Section \ref{sec:stable-into-reflexive}, we presented a result of Kalton which stipulates that every stable metric space admits an almost isometric embedding into a reflexive Banach space. In \cite{Kalton2007}, Kalton asked whether some kind of converse holds. Problem \ref{prob:Kalton} below has become known as one of Kalton's embedding problems. 
	
	\begin{prob}
		\label{prob:Kalton}
		Let $X$ be a separable reflexive Banach space.
		\begin{enumerate}[(i)]
			\item Does $X$ coarsely embed into a stable metric space? 
			\item Does $X$ uniformly embed into a stable metric space?
		\end{enumerate}
	\end{prob}
	
	Kalton also formulated the following variant of Problem \ref{prob:Kalton} (ii).
	
	\begin{prob}
		\label{prob:Kalton-uniform-ball}
		Does the unit ball of every separable reflexive Banach space uniformly embed into a stable metric space? 
	\end{prob}
	
	We know three types of spaces with property $Q$: stable spaces, reflexive spaces and duals of asymptotic-$c_0$ spaces. The first two kinds of spaces admit a coarse embedding into a reflexive Banach space. The following problem, which also asks whether the converse of the first assertion in Corollary \ref{cor:Q-CE-co-reflexive} is true, is open. 
	
	\begin{prob}
		\label{prob:Q->CE-reflexive}
		Does every metric space with property $Q$ coarsely embed into a reflexive Banach space?
	\end{prob}  
	
	The following subquestion is also open.
	
	\begin{prob}
		\label{prob:asymp-co->dual-CE-reflexive}
		Let $X$ be an asymptotic-$\co$ Banach space. Does $X^*$ coarsely embed into a reflexive Banach space?
	\end{prob}
	
	%The connection between stable functions and weak-compactness goes back at least to the work of Grothendieck \cite{Grothendieck52} (see BenYaccov, etc, stability in model theory....).

	\section{Exercises}
	\label{exer:Q}
	
	\begin{exer}
		\label{exer:interlaced-closed-formula}
		For all $\mbar,\nbar\in [\bN]^{<\infty}$, let 
		$$d(\mbar,\nbar):=\sup \{\big| \abs{S\cap \mbar} - \abs{S\cap \nbar} \big| \colon S \text{ interval of } \bN \}.$$
		Show that $d$ is a metric on $[\bN]^{<\infty}$. 
	\end{exer}
	
	\begin{exer}\label{exer:I_kdistortion}
		Show that $\cdist{\co}(\sI_k) = 2$.
	\end{exer}
	
	\begin{proof}[Hint]
		Adapt the proof of Proposition \ref{prop:c_0l_1distorsion}. 
	\end{proof}
	
	\begin{exer}
		\label{exer:interlaced-isometric}
		Let $(s_i)_{i\in \bN}$ be the summing basis of $\co$ and define the usual bimonotone version of the summing norm as follows:
		\begin{align}
			\label{eq:summing-norm}
			\Big\|\sum_{i=1}^\infty a_i s_i\Big\|_{\mathrm{sum}} = \sup\Big\{\Big|\sum_{j=r}^s a_j \Big| \colon r,s\in \bN,\,  r\le s\Big\}.
		\end{align}
		\begin{enumerate}
			\item Show that for all $x\in \co$, $\norm{x}_\infty\le \norm{x}_{\mathrm{sum}}\le 2\norm{x}_\infty$ and that the factor $2$ cannot be improved.
			\item[] \hskip -.8cm In the following questions for all $\mbar,\nbar \in [\bN]^{<\omega}$ we define the summing distance by $$\mathsf{d_{sum}}(\mbar,\nbar) := \Big\|\sum_{i=1}^{\abs{\mbar}} s_{m_i} - \sum_{i=1}^{\abs{\nbar}} s_{n_i}\Big\|_{\mathrm{sum}}.$$
			\item Show that in \eqref{eq:summing-norm} one only needs to consider intervals whose boundaries are sign-changes of the $a_i$. More precisely, for a sequence $(a_i)_{i=1}^\infty$ in $c_{00}$ let $0=m_0<m_1<\dots< m_l$ be chosen in $\bN$ so that for all $i\le s$ the signs of $a_j$ on $j\in [m_{i-1}+1, m_i]$ are the same
			(\ie all nonnegative or all nonpositive) then
			\begin{align}
				\label{eq2:summing-norm}
				\Big\|\sum_{i=1}^\infty a_i s_i\Big\|_{\mathrm{sum}} = \sup\Big\{\Big|\sum_{j=r}^s \sum_{i=m_{j-1}+1}^{m_j} a_i\Big|   \colon 1\le r\le s\le l\Big\}.
			\end{align}
			\item Show that for all $\mbar,\nbar \in [\N]^{<\omega}$ if we write $\mbar\triangle \nbar$ in increasing order as $\mbar\triangle \nbar =\{x_1,x_2,\dots,x_t\}$, then 
			\begin{align}
				\label{summing metric for sets}
				\mathsf{d_{sum}}(\mbar,\nbar) & = \max\left\{\big|\abs{S\cap \mbar} - \abs{S\cap \nbar}\big| \colon S \text{ is an interval of }\bN\right\}\\
				& = \max\left\{\big|\abs{\llbracket x_i, x_j \rrbracket\cap \mbar} - \abs{\llbracket x_i, x_j\rrbracket\cap \nbar} \colon 1\le i\le j\le t \right\}.\notag
			\end{align}
			\item We define the interlaced graph structure on the set of vertices $[\bN]^{<\omega}$. We declare that two vertices $\mbar$ and $\nbar$ in $[\N]^{<\omega}$
			are adjacent if and only if $\mbar \neq \nbar $ and one of the following interlacing relations holds
			\begin{enumerate}[(i)]
				\item $\abs{\mbar} = \abs{\nbar} +1 :=k+1$ and $m_1 \le n_1 \le m_2\le  \dots \le m_k \le n_k \le m_{k+1}$,
				\item $\abs{\nbar} = \abs{\mbar} +1 :=k+1$ and $n_1 \le m_1\le n_2\le  \dots \le n_k \le m_k \le n_{k+1}$,
				\item $\abs{\nbar} = \abs{\mbar} := k$, $m_1 \le n_1 \le m_2\le  \dots \le m_k \le n_k$, or
				\item $\abs{\nbar} = \abs{\mbar} := k$, $n_1 \le m_1\le n_2\le  \dots \le n_k \le m_k$.
			\end{enumerate}
			We also connect the empty set with all singletons. The graph metric for the graph structure is denoted by $\mathsf{d_I}$.
			
			\begin{enumerate}
				\item Show that $\mathsf{d_I}$ and $\mathsf{d_{sum}}$ coincide on $[\bN]^{<\omega}$.
				\item Show that for all $\mbar,\nbar\in[\bN]^k$ there is a path of length $\mathsf{d_I}(\mbar,\nbar)$ from $\mbar$ to $\nbar$ in the interlaced graph on $[\bN]^{<\omega}$ that stays in $[\N]^k$ and deduce that the restriction of $\mathsf{d_I}$ to $[\bN]^k$ is $\dIk$.
				\item Conclude that the interlaced graph $([\bN]^{<\omega},\mathsf{d_I})$ embeds isometrically into $\co$ equipped with the bimonotone version of the summing norm.
			\end{enumerate}
		\end{enumerate}
	\end{exer}
	
	\begin{exer}
		\label{exe:3-stable}
		Show that for all bounded sequences $\xn$ and $\yn$ in a metric space $(M,d)$ and every nonprincipal ultrafilters $\cU$ and $\cV$ on $\bN$ one has 
		\begin{equation*}
			\lim_{m\in \cU}\lim_{n\in \cV} d(x_m,y_n) \le 3 \lim_{n\in \cV}\lim_{m\in \cU} d(x_m,y_n)
		\end{equation*}
	\end{exer}
	
	%This is Proposition 3.6 in \cite{RO22}
	
	\begin{exer}
		\label{exe:diamond-not-stable}
		Give an example of a metric space that is not stable but admits a bi-Lipschitz embedding into a stable space.
	\end{exer}
	
	\begin{proof}[Hint]
		Show that $\dia_1^\omega$,  the countably branching diamond graph of height $1$, is not stable (see Chapter \ref{chapter:diamonds} for its definition).
		%Consider $(a_n)_n=(y_1,t,y_3,t,y_5,t,\dots, t,y_{2k+1},t,\dots)$, $(b_n)_n=(y_2,s,y_4,s,y_6,s,\dots, s,y_{2k},s,\dots)$, $\cU$ an ultrafilter that contains the even numbers and $\cV$ and ultrafilter that contains the odd number. Then, $\lim_{m\in \cU}\lim_{n\in \cV} d(a_m,b_n)=2$ while $\lim_{n\in \cV}\lim_{m\in \cU} d(a_m,b_n)=1$.
	\end{proof}

	\begin{exer}
		\label{exe:asymp-lp}
		Write a formal proof of  Lemma \ref{lem:asymp-lp}.
	\end{exer}
	
	\begin{exer}
		Prove Proposition   \ref{prop:Lipstabilty infrasup}.
	\end{exer}
	
	\begin{exer}
		\label{exe:lp-sum-of-stable} 
		Let $p\in[1,\infty)$. Show that if for all $n\in \bN$, $X_n$ is a stable Banach space, then $(\sum_{n=1}^\infty X_n)_{\ell_p}$ is stable.
	\end{exer}
	
	\begin{exer}
		\label{exe:Q-FD}
		Let $X$ be a Banach space. Show that the following assertions are equivalent
		\begin{enumerate}
			\item $X$ is finite-dimensional
			\item $Q_X=0$
			\item $Q_X<1$
		\end{enumerate}
	\end{exer}
	
	\begin{exer}
		\label{exe:weak*-linearization}
		Let $X$ be a separable Banach space and $d_k$ be a metric on $[\bN]^k$. Fix $k\in\bN$ and $f \colon ([\bN]^k, d_k)\to X^*$ a map with bounded image. Show that there exist an infinite subset $\bM$ of $\bN$ and a $(x^*_{\mbar})_{\mbar\in[\bM]^{\le k}}$ in $X^*$ so that:
		\begin{enumerate}[(i)]
			\item for all $\mbar \in [\bM]^k$,  
			\begin{equation*}
				\label{eq:weak*-linearization-eq1}
				f(\mbar) = x^*_\emptyset + \sum_{i=1}^k x^*_{(m_1,\dots, m_i)},
			\end{equation*}
			\item for all $\mbar\in [\bM]^{\le  k}\setminus\{\emptyset\}$ of length $l$, 
			\begin{equation*}
				\norm{x^*_{\mbar}} \le  \Lip(f)  \sup\{ d_k(\nbar,\sigma_{l,m}(\nbar))\colon \nbar\in [\bM]^{k}, m\in \bM \text{ s.t. }\sigma_{l,m}(\nbar)\in [\bM]^k\}.
			\end{equation*}
		\end{enumerate}    
	\end{exer}

	\begin{exer}
		\label{exe:Helly}
		Prove Claim \ref{clai1:reflexivity-ABS-Q}
	\end{exer}
	
	\begin{exer}
		\label{exe:Ramsey-diagonal}
		Prove Claim \ref{clai2:reflexivity-ABS-Q}
	\end{exer}
	
	\begin{exer}\label{exer:cluster} 
		Show that the conclusion of Theorem \ref{thm:reflexivity-ABS-Q} still holds if we only assume that $x^{**}$ is a weak$^*$ cluster point of the sequence $(x_n)_{n=1}^\infty$. 
	\end{exer}

	The argument in the next exercise is Kalton's original proof that a stable metric space has property $Q$.
	
	\begin{exer}
		Let $(M,d)$ be a stable metric space
		\begin{enumerate}
			\item Given a bounded map $f\colon([\N]^k, \dIk) \to M$ and $\eps>0$, show that there exist an infinite subset $\bM$ of $\bN$ and $u\in \bR$ such that for all $\nbar,\mbar \in [\bM]^k$ with $m_k<n_1$,
			\begin{equation*}
				\abs{u- d(f(\nbar),f(\mbar))} \le \frac{\eps}{3}.
			\end{equation*}
			\item Show that $M$ has property $Q$.
		\end{enumerate}
	\end{exer}
	
	%\begin{proof}[Solution] 
	%Assume that $X$ is a stable Banach space and $f\colon ([\N]^k,\dk) \to X$ is a bounded map and fix $\eps>0$. Note that, by Ramsey's theorem (Theorem \ref{Ramsey1}), we can assume by passing to an infinite subset of $\N$ and re-indexing that there exists $u\in \bR$ such that $$\forall \n,\m \in [\N]^k\ \text{such that}\ m_k<n_1,\ \ \big|u- \|f(\n)-f(\m)\|\big| \le \frac{\eps}{3}.$$ Fix now a non principal ultrafilter $\cal U$ on $\N$. It follows from the definition of the interlaced metric that $$\lim_{m_1\in \cal U}\lim_{n_1 \in \cal U}\ldots \lim_{m_k\in \cal U}\lim_{n_k \in \cal U} \|f(m_1,\ldots,m_k)-f(n_1,\ldots,n_k)\| \le \Lip (f).$$ Using the stability of $X$, we can exchange the limits in the above inequality and deduce that there exists an infinite subset $\M$ of $\N$ such that $$\forall \m,\n \in [\M]^k\ \text{such that}\ m_k<n_1,\ \ \|f(\m)-f(\n)\| \le \Lip (f)+\frac{\eps}{3}.$$ For any $\m,\n \in [\M]^k$, there exists $\lbar \in [\M]^k$ such that $l_1>n_k$ and $l_1>m_k$. This, together with the previous inequalities yields that $\diam\, f([\M]^k) \le 2\,\Lip (f)+\eps$.
	%\end{proof}

	\begin{exer}
		Show directly (without using the notion of infrasup stability) that if a Banach space $X$ coarsely embeds into a space with property $Q$, then $X$ has property $Q$.
	\end{exer}
	
	%\begin{proof}[Solution] 
	%Suppose that $g:Y\to X$ is a coarse embedding. Let $f\colon ([\N]^k,\dk) \to Y$  be a Lipschitz map that we may, by homogeneity, assume to be $1$-Lipschitz. Since $\dk$ is a graph metric, we deduce that $g \circ f$ is $\omega_g(1)$-Lipschitz. It now follows from (i) that there exists an infinite subset $\M$ of $\N$ such that $\diam (g \circ f)([\M]^k) \le 3\omega_g(1)$. Finally, since $\lim_{t \to \infty}\rho_g(t)=\infty$, $C=\inf\{t,\ \rho_g(t)>3\omega_g(1)\}<\infty$ and $\diam\, g([\M]^k) \le C$.
	%\end{proof}
	
	\begin{exer}
		\label{exe:JT}
		Show that the James Tree space $\mathrm{JT}$ satisfies property $(\star)$ from Theorem \ref{thm:coarse-universal}. 
		\begin{proof}[Hint]
			Use the uncountable set of infinite branches in $\{0,1\}^{<\omega}$.
		\end{proof}
	\end{exer}
	
	\begin{exer}
		\label{exe:EL=ECL-interlaced}
		Show that if $(\Ik)_{k\in\bN}$ admits an equi-coarse-Lipschitz embedding into a Banach space $X$, then $(\Ik)_{k\in\bN}$ admits an equi-bi-Lipschitz embedding into $X$.
	\end{exer}

	\begin{exer}
		\label{exe:basics-Q-omega-1}\,
		\begin{enumerate}[(i)]
			\item Show that If $X$ has property $Q^{\omega_1}$, then the sequence of graphs $(\Iuk)_{k\in \N}$ does not equi-coarsely embed into a Banach space with property $Q^{\omega_1}$.
			\item Assume that $X$ is a Banach space that coarsely embeds into $Y$ and that $Y$ has property $Q^{\omega_1}$. Show that $X$ has property $Q^{\omega_1}$.
		\end{enumerate}
	\end{exer}
	
	\begin{exer}[Direct proof of the failure of $Q^{\omega_1}$ for $\ell_\infty/\co$.]
		\label{exe:l/c}
		\begin{enumerate}
			\item Show that there is a family $(A_\mu)_{\mu<\omega_1}$ of infinite subsets of $\bN$ such that for all $\mu<\nu$, $A_\nu \setminus A_\mu$ is infinite and $A_\mu \setminus A_\nu$ is finite.
			\item Show that $\ell_\infty/\co$ fails property $Q^{\omega_1}$. 
		\end{enumerate}
	\end{exer}
	\begin{proof}[Hints]\,
		\begin{enumerate}
			\item Use transfinite induction.
			\item Let $x_\mu \in \ell_\infty/c_0$ be the image by the quotient map $Q\colon \ell_\infty \to \ell_\infty/c_0$ of the indicator function of $\N \setminus A_\mu$ and consider 
			$f(\overline{\alpha}) :=\sum_{i=1}^k x_{\alpha_i}$ for all $\alphabar :=(\alpha_1,\dots,\alpha_k) \in [\Omega_1]^k$.
		\end{enumerate}
	\end{proof}

	%%%%%%%%%%%%%%%%%%%%%%%%%%%%%%%%%%%%%%%%%%%%%%%%%%%%%%%%%%%%%%%%%%%%%%%%%%%%%%
	
	\chapter[The approximate midpoint principle revisited]{The approximate midpoint principle revisited}
	\label{chapter:AMP_II}
	
	This chapter is based on a paper by N. Kalton \cite{Kalton2013b} that was published after its author suddenly passed away in August 2010. This amazingly deep paper still contains the most advanced results on the preservation of asymptotic uniform convexity under coarse-Lipschitz embeddings. For various reasons, we believe that the numerous results from this paper have been unfairly overlooked. In this chapter, we attempted to make accessible to a wide audience some of the most important of these results. There are several mathematical gems in this chapter. 
	
	The first one is the preservation under coarse-Lipschitz embeddings of a natural quantitative lower bound on the fundamental function of spreading models generated by weakly null sequences. This result is a consequence of a more general result that establishes a very fine and delicate relationship between certain averages in the Banach space and its modulus of asymptotic midpoint uniform convexity. This is carried out in the first section that ends with a general consequence of the approximate midpoint principle for coarse-Lipschitz embeddings into asymptotically uniformly convex or even asymptotically midpoint uniformly convex Banach spaces.
	
	The second and third sections concern understanding conditions under which the weak $p$-co-Banach-Saks property (a sequential property) implies $p$-asymptotic uniform convexity (a tree property). Such a condition relies on the notion of a random $L_p$-norm.
	
	In the fourth section, we prove the coarse-Lipschitz rigidity of subspaces and quotients of $\ell_p$, for $p \in (1,\infty)$.
	
	In the last section, we prove the uniqueness of the coarse (or uniform) structure of $(\sum_{n=1}^\infty \ell_r^n)_{\ell_p}$ for $1<p<r\le 2$ or $2\le r<p<\infty$. This section contains some original linear results from \cite{Kalton2013b} about the subspace structure of $C_p$ that are crucial for applications to the nonlinear classification of Banach spaces, but are also of independent interest.
	
	\section[Spreading models and coarse-Lipschitz embeddings]{Properties of spreading models preserved under coarse-Lipschitz embeddings}
	
	In Section \ref{sec:weighted-Hamming}, we showed using asymptotic uniform smoothness arguments that for reflexive spaces, upper bounds on the behavior of spreading models are preserved under coarse-Lipschitz embeddings. In this section, we derive lower bounds on the behavior of spreading models. These lower bounds, which only apply to the fundamental function of spreading models, are obtained without a reflexivity assumption and are related to asymptotic (midpoint) uniform convexity.
	
	\subsection{An inequality in asymptotically midpoint uniformly convex Banach spaces}
	\label{sec:amuc}
	
	We first recall the definition of the modulus of asymptotic midpoint uniform convexity already introduced in Chapter \ref{chapter:asymptotic-moduli}, Section \ref{sec:auc-aus}. 
	
	\begin{defi}
		Let $X$ be a Banach space. For $t\ge 0$, let
		\begin{equation}
			\hat{\delta}_X(t) := \inf_{x \in S_X}\ \sup_{Y\in \cof(X)}\ \inf_{y\in Y, \norm{y}\ge 1}\  \frac{\norm{x+ty}+\norm{x-ty}}{2}-1.
		\end{equation}
		Banach spaces that satisfy $\hat{\delta}_X(t)>0$ for all $t\in (0,1)$ are called \emph{asymptotically midpoint uniformly convex}\index{asymptotically midpoint uniformly convex}. 
	\end{defi}
	
	It easily follows from the triangle inequality that the function $t \mapsto \frac{\hat{\delta}_X(t)}{t}$ is nondecreasing. Therefore, one can ``convexify'' the modulus as follows:
	\begin{equation}
		t\in[0,\infty) \mapsto \hat{\delta}^c_X(t) :=\int_0^t \frac{\hat{\delta}_X(s)}{s}\,ds.
	\end{equation}
	The map $\hat{\delta}^c_X$ is convex and using again that the map $t \mapsto \frac{\hat{\delta}_X(t)}{t}$ is nondecreasing we obtain that for all $t \ge 0$, 
	$$\hat{\delta}_X\big(\frac{t}{2}\big) \le \hat{\delta}^c_X(t) \le \hat{\delta}_X(t),$$
	so that $\hat{\delta}_X$ is equivalent to the convex function $\hat{\delta}^c_X$. Note also that $\hat{\delta}^c_X$ is $1$-Lipschitz and that $\hat{\delta}^c_X$ fulfills all the assumptions of the Orlicz functions in the second item of Lemma \ref{lem:absolute-Orlicz} (see Exercise \ref{exe:deltahat}). 
	
	\medskip Recall that the separation parameter of a sequence $\xn$ in a Banach space $(X,\norm{\cdot})$ is defined as $\sep(\xn):=\inf_{m\neq n} \norm{x_m-x_n}$. We say that $\xn$ is $\theta$-separated, for some $\theta>0$ whenever $\sep(\xn)\ge \theta$. The following inequality, which improves over the triangle inequality for asymptotically midpoint uniformly convex spaces, will be crucial in the sequel.
	
	\begin{prop}
		\label{pro:separation-AMUC}
		Let $X$ be a Banach space and $u\neq v$ in $X$. Assume that $(x_n)_{n=1}^\infty$ is a bounded and $\theta$-separated sequence in $X$. Then,
		\begin{equation}
			\liminf_{n\to \infty}(\norm{u-x_n} + \norm{v-x_n})\ge \norm{u-v} + \norm{u-v}\hat{\delta}_X\Big(\frac{\theta}{\norm{u-v}}\Big).
		\end{equation}
	\end{prop}
	
	\begin{proof} 
		By homogeneity, it is sufficient to prove the statement when $\norm{u-v}=1$. Fix $\eta\in(0,\theta)$ and pick a finite-codimensional subspace $Y$ of $X$ such that for all $z\in Y$ with $\norm{z}\ge \theta$ we have 
		\begin{equation*}
			\frac12\big(\norm{ u-v+z } + \norm{ u-v-z }\big)\ge 1+\hat{\delta}_X(\theta) - \eta.
		\end{equation*}
		Since $X/Y$ is finite-dimensional and $(x_n)_{n=1}^\infty$ is bounded, a compactness argument for $\norm{\cdot}_{X/Y}$ yields the existence of $m\neq n$ so that $\norm{x_m-x_n}_{X/Y}=d(x_n-x_m,Y)<\eta$. Choose $y\in Y$ so that $y\neq 0$ and $\norm{x_n-x_m-y}<\eta$ and let $z :=\norm{ x_n-x_m }\frac{y}{\norm{ y }}$. We have that $z\in Y$, $\norm{z}\ge \theta$ and 
		\begin{equation*}
			\norm{x_n-x_m-z}\le \norm{x_n-x_m-y} + \norm{y-z} < \eta + \big| \norm{y} - \norm{x_n-x_m} \big|< 2\eta.
		\end{equation*}
		We deduce that 
		\begin{equation*}
			\frac12\big(\norm{ u-v+(x_m-x_n)} + \norm{u-v-(x_m-x_n)}\big)\ge  1+\hat{\delta}_X(\theta) -3\eta,
		\end{equation*}
		and it follows from the triangle inequality that
		\begin{equation*}
			\frac12\big(\norm{u-x_n} + \norm{v-x_m} + \norm{u-x_m} + \norm{v-x_n}\big)\ge  1+\hat{\delta}_X(\theta) -3\eta,
		\end{equation*}
		from which we can easily deduce that for $j=m$ or $j=n$, we have 
		\begin{equation*}
			\norm{u-x_j} + \norm{v-x_j}\ge 1+\hat{\delta}_X(\theta) -3\eta.
		\end{equation*}
		Since $\eta\in(0,\theta)$ was arbitrary and the compactness argument could be applied to any subsequence of $\xn$, the statement follows.
	\end{proof}

	\subsection[Lower estimates on the fundamental functions of  spreading models]{Lower estimates on the fundamental functions of spreading models}
	
	Using an iterated norm argument, it is not too difficult to see (see Exercise \ref{exe:ff-sm}) that for any spreading model $S$ generated by a weakly null sequence in a Banach space $X$ (with separable dual) there is a constant $C>0$ such that for all $k\ge 1$ and $a_1,\dots,a_k\in \bR$, 
	\begin{equation}
		\label{eq:lower-upper-sm}
		\frac{1}{C}\Big\|\sum_{i=1}^k a_i e_i\Big\|_{\ell_{\bar{\delta}^c_X}} \le \Big\|\sum_{i=1}^k a_i e_i\Big\|_S \le C\Big\|\sum_{i=1}^k a_i e_i\Big\|_{\ell_{\bar{\rho}_X}}, 
	\end{equation}
	where, as always, $\en$ represents indifferently the canonical basis in the respective sequence spaces and $\bar{\delta}^c_X$ is the convexification of the (not necessarily convex) modulus of asymptotic convexity. In Chapter \ref{chapter:Hamming}, we have already discussed the stability of the upper estimate in \eqref{eq:lower-upper-sm} under coarse-Lipschitz embeddings. In this section, we are concerned with the preservation of the lower estimate in \eqref{eq:lower-upper-sm}
	under coarse-Lipschitz embeddings.
	
	Since $\bar{\delta}^c_X \le \hat{\delta}^c_X $ and spreading models generated by weakly null sequences are unconditional, the following proposition contains an improvement of the lower estimate in \eqref{eq:lower-upper-sm}. We only state it for a comparison with Theorem \ref{thm:CL-AMUC}, so we leave the proof as Exercise \ref{exe:lower-AMUC-sm}. 
	
	\begin{prop}
		\label{prop:lower-AMUC-sm}
		For any spreading model $S$ generated by a normalized sequence in a Banach space $X$, there is a constant $c>0$ such that for all $k\ge 1$ and $a_1,\dots,a_k\in \bR$, 
		\begin{equation}
			\label{eq:lower-AMUC-sm}
			c\Big\|\sum_{i=1}^k a_i e_i\Big\|_{\ell_{\hat{\delta}^c_X}} \le \bE \Big\|\sum_{i=1}^k \vep_i a_i e_i\Big\|_S,
		\end{equation}
		where $\bE$ refers to the average over all choices of signs $\vep:=(\vep_1,\dots,\vep_k)\in\{-1,1\}^k$.
	\end{prop}
	
	The main goal of this section is to show that an estimate similar to \eqref{eq:lower-AMUC-sm} is preserved under coarse-Lipschitz embeddings, albeit that this estimate can be proven for constant coefficients only. We start by establishing some notation and preparatory lemmas.
	
	Given a normalized sequence $\xn$ in a Banach space $X$, we will be interested in the following numerical quantities. Given $k\ge 1$, we let 
	\begin{equation}
		\sigma_k(\xn) := \sup \Big\{ \bE \bnorm{ \sum_{i=1}^k \eps_i x_{n_i} } \colon (n_1,\dots,n_k) \in [\bN]^k \Big\},
	\end{equation}
	and set $\sigma_0(\xn):=0$.
	
	In a first elementary lemma, we record some monotonicity properties of the sequence $(\sigma_k(\xn))_{k\ge 0}$.
	
	\begin{lemm}
		\label{lem:sigmak-monotonicity}
		Let $X$ be a Banach space, $\xn$ be a normalized sequence in $X$ and $\sigma_k:=\sigma_k(\xn)$ for $k\ge 0$. Then, the sequence $(\sigma_k)_{k=0}^\infty$ is nondecreasing and the sequence $\big(\frac{\sigma_k}{k}\big)_{k=1}^\infty$ is nonincreasing.
	\end{lemm}
	
	\begin{proof} 
		The fact that $(\sigma_k)_{k=0}^\infty$ is nondecreasing follows from an elementary convexity argument. Indeed, for a fixed $(n_1,\dots,n_{k+1}) \in [\bN]^{k+1}$, we have:
		\begin{align*}
			\bE\Big\|\sum_{i=1}^{k+1} \eps_i x_{n_i}\Big\| & =  \bE\Big( \frac12\Big\| \sum_{i=1}^{k}\eps_i x_{n_i} + x_{n_{k+1}} \Big\| + \frac12 \Big\|\sum_{i=1}^{k}\eps_i x_{n_i}-x_{n_{k+1}} \Big\| \Big)\\
			&\ge \bE \Big\| \sum_{i=1}^{k} \eps_i x_{n_i} \Big\|,
		\end{align*}
		and taking suprema leads to the conclusion.
		
		Now, for $k\in \bN$ and $(n_1,\ldots,n_{k+1})\in [\bN]^{k+1}$, we have 
		\begin{equation*}
			\bE\Big\| \sum_{i=1}^{k+1}\eps_i x_{n_i}\Big\|  = \frac{1}{k} \bE\Big\|\sum_{i=1}^{k+1} \sum_{j=1,j\neq i}^{k+1} \eps_j x_{n_j} \Big\| \le \frac{1}{k} \sum_{i=1}^{k+1} \bE\Big\| \sum_{j=1,j\neq i}^{k+1} \eps_j x_{n_j} \Big\| \le \frac{k+1}{k}\sigma_k,    
		\end{equation*}
		from which it follows immediately that the sequence $\big(\frac{\sigma_k}{k}\big)_{k=1}^\infty$ is  nonincreasing.
	\end{proof}
	
	In the next lemma, we prove a rather fine inequality that improves over the basic triangle inequality for sequences where $\sigma_k(\xn) = o(k)$. 
	
	\begin{lemm}
		\label{lem:Orlicz-F_k}
		Let $X$ be a Banach space and $\xn$ be a normalized sequence in $X$. For $k\ge 1$, let $\sigma_k:=\sigma_k(\xn)$ and consider the Orlicz function $F_k$ defined by
		\begin{equation}
			\label{eq:F_k}
			F_k(t) :=\begin{cases} 
				\frac{\sigma_k}{k}t,           \qquad \qquad \text{ if } 0\le t\le \frac{1}{\sigma_k},\\
				t+\frac1k-\frac{1}{\sigma_k}, \quad  \text{ if } \frac{1}{\sigma_k}\le t<\infty.
			\end{cases}
		\end{equation}
		Then, for all $j,k\ge 1$ and all $a_1,\dots a_j\in \bR$, we have 
		\begin{equation}
			\bE\Big\|\sum_{i=1}^j \vep_i a_i x_i\Big\| \le 2\Big\|\sum_{i=1}^j a_i e_i\Big\|_{\ell_{F_k}}.
		\end{equation}
	\end{lemm}
	
	\begin{proof} 
		It follows from the (refined version) of Kahane's contraction principle that 
		\begin{equation}
			\label{eq:KCP}
			\bE\Big\|\sum_{i=1}^j \vep_i a_i x_i\Big\| \le \sum_{i=1}^{j}(a_i^* - a_{i+1}^*)\bE\Big\|\sum_{r=1}^i \vep_r x_r\Big\| \le \sum_{i=1}^{j}(a_i^* - a_{i+1}^*)\sigma_i,
		\end{equation}
		where $(a_i^*)_{i=1}^j$ is a decreasing rearrangement of $(|a_i|)_{i=1}^j$ and in particular $a^*_1\ge a^*_{2}\ge \dots a^*_j\ge a^*_{j+1}:=0$.
		
		Assume first that $j\ge k$. We can write 
		\begin{equation*}
			\sum_{i=1}^{j}(a_i^* - a_{i+1}^*)\sigma_i = \sum_{i=1}^{k}(a_i^* - a_{i+1}^*)\sigma_i + \sum_{i=k+1}^{j}(a_i^* - a_{i+1}^*)\sigma_i.
		\end{equation*} 
		Then remembering that $k\mapsto \frac{\sigma_k}{k}$ is nonincreasing, we have 
		\begin{equation}
			\label{eq1:Orlicz-F_k}
			\sum_{i=1}^{j}(a_i^* - a_{i+1}^*)\sigma_i\le \sum_{i=1}^{k}(a_i^* - a_{i+1}^*)\sigma_i + \frac{\sigma_k}{k}\sum_{i=k+1}^j i(a_i^*-a_{i+1}^*).
		\end{equation}
		Using simple Abel transforms, we have that 
		\begin{align}
			\label{eq2:Orlicz-F_k}
			\sum_{i=1}^{k}(a_i^* - a_{i+1}^*)\sigma_i = \sum_{i=1}^{k} a^*_i(\sigma_i - \sigma_{i-1}) + a_1^*\sigma_0 -a_{k+1}^*\sigma_k 
			= \sum_{i=1}^{k} a^*_i(\sigma_i - \sigma_{i-1}) -a_{k+1}^*\sigma_k
		\end{align}
		and
		\begin{equation}
			\label{eq3:Orlicz-F_k}
			\sum_{i=k+1}^j i(a_i^*-a_{i+1}^*) = \sum_{i=k+1}^j a_i^*(i-(i-1)) + ka^*_{k+1} - ja^*_{j+1} = \sum_{i=k+1}^j a_i^* + ka^*_{k+1}.
		\end{equation}
		A combination of $\eqref{eq:KCP}, \eqref{eq1:Orlicz-F_k},\eqref{eq2:Orlicz-F_k}$ and $\eqref{eq3:Orlicz-F_k}$ gives
		\begin{equation}
			\bE\Big\|\sum_{i=1}^j \vep_i a_i x_i\Big\| \le  \sum_{i=1}^{k} a^*_i(\sigma_i - \sigma_{i-1}) + \frac{\sigma_k}{k}\sum_{i=k+1}^j a_i^*.
		\end{equation}
		
		Now, if $\abar:=(a_1,\dots,a_j)\in \bR^j$ is such that 
		$$1=\norm{\abar}_{\ell_{F_k}}=\inf\Big\{s>0\colon \sum_{i=1}^j F_k\Big(\frac{\abs{a_i}}{s}\Big)\le 1\Big\},$$
		then by continuity of $F_k$ we have that $\sum_{i=1}^j F_k(a_i^*)\le 1$ and since $F_k$ is increasing it follows that $F_k(a_k^*) \le \frac{1}{k}$. By definition of $F_k$ we have that $a_k^*\le \frac{1}{\sigma_k}$ and thus  since $F_k$ is increasing, we have that $F_k(a_i^*) = \frac{\sigma_k}{k}a_i^*$ for all $i\ge k$.
		Moreover, $F_k(t)\ge t-\frac{1}{\sigma_k}$ for all $t\ge 0$ and hence
		\begin{align*}
			\bE\Big\|\sum_{i=1}^j \vep_i a_i x_i\Big\| & \le \sum_{i=1}^{k} (F_k(a_i^*) +\frac{1}{\sigma_k})(\sigma_i - \sigma_{i-1}) + \sum_{i=k+1}^j F_k(a_i^*)\\
			& = 1 + \sum_{i=1}^{k} F_k(a_i^*)(\sigma_i - \sigma_{i-1})+ \sum_{i=k+1}^j F_k(a_i^*).
		\end{align*}
		Observing that $\abs{\sigma_i - \sigma_{i-1}}\le 1$ for all $1\le i\le k$, we conclude that 
		\begin{equation}
			\bE\Big\|\sum_{i=1}^j \vep_i a_i x_i\Big\| \le  1 + \sum_{i=1}^{j} F_k(a_i^*) \le 2.
		\end{equation}
		The conclusion follows from homogeneity. 
		
		The remaining case $j<k$ is simpler. Indeed, in this case 
		\begin{align*}
			\bE\Big\|\sum_{i=1}^j \vep_i a_i x_i\Big\| & \le \sum_{i=1}^{j} a^*_i(\sigma_i - \sigma_{i-1}) \\
			& \le \sum_{i=1}^{j} \Big(F_k(a_i^*) +\frac{1}{\sigma_k}\Big)(\sigma_i - \sigma_{i-1})\\
			& \le \frac{\sigma_j}{\sigma_k} + \sum_{i=1}^{j} F_k(a_i^*),
		\end{align*}
		and if we have $\norm{\abar}_{\ell_{F_k}}=1$, this last quantity is at most $2$ since $(\sigma_i)_i$ is nondecreasing and we conclude again by homogeneity. 
	\end{proof}
	
	The theorem below is the main result of this section and has a wide range of applicability. It provides a powerful obstruction to coarse-Lipschitz embeddability that is based on asymptotic uniform convexity considerations.
	
	\begin{theo}
		\label{thm:CL-AMUC}
		Let $X$ and $Y$ be two Banach spaces so that $X$ coarse-Lipschitz embeds into $Y$. Then, there exists a constant $c>0$ such that for every $\theta>0$, every $\theta$-separated normalized sequence $(x_n)_{n=1}^\infty$ in $X$ and every $k\in \bN$, there exists $(n_1,\dots,n_k)\in [\bN]^k$ such that
		\begin{equation}
			\label{eq:CL-AMUC}
			c\theta\Big\|\sum_{i=1}^k e_i\Big\|_{\ell_{\hat{\delta}^c_Y}} \le \bE\Big\|\sum_{i=1}^k \eps_i x_{n_i}\Big\|_{X}.
		\end{equation}
	\end{theo}
	
	\begin{proof}
		In this proof, it will be convenient to equip $\{-1,1\}$ with the uniform probability measure, and we will consider $\Delta:=\{-1,1\}^\bN$ as the probability space equipped with the corresponding product measure. Also, $(\eps_n)_{n=1}^\infty$ will be any sequence of independent Rademacher random variables that we may and will view as the coordinate functionals on $\Delta$. If $x$ is a vector in a Banach space $X$ and $f\colon \Omega \to \bR$ is any function, then $f\otimes x\colon \Omega \to X$ will be the map $\omega \mapsto f(\omega)x$. 
		
		Let us fix a $\theta$-separated normalized sequence $(x_n)_{n=1}^\infty$ in $X$. For $k\in \bN$, we still denote by $F_k$ the function defined in \eqref{eq:F_k}. Since $F_k$ is an Orlicz function satisfying $\lim_{t\to \infty} \frac{F_k(t)}{t}=1$, we can define the corresponding absolute norm on $\bR^2$, denoted by $N_k$. Recall that $N_k(0,1)=1$ and for all $t\ge 0$, 
		\begin{equation*}
			N_k(1,t)=1+F_k(t)
		\end{equation*}
		It follows from Lemma \ref{lem:Orlicz-F_k} and Proposition \ref{pro:Orlicz-iterated} that for all $\xi \in c_{00}$, $\norm{T\xi}_{L_1(\Delta; X)}\le 4\norm{\xi}_{\Lambda_{N_k}}$, where $T \colon c_{00}\to L_1(\Delta; X)$ is the operator defined by 
		\begin{equation*}
			T\xi := \sum_{j=1}^{\infty}\xi_j\eps_j\otimes x_j.
		\end{equation*}
		On the other hand, the map $F_Y:=\hat{\delta}^c_Y$ is also an Orlicz function satisfying $\lim_{t\to \infty} \frac{F_Y(t)}{t}=1$. So, we similarly associate with it an absolute norm on $\bR^2$ that we denote by $N_Y$. Note that for all $t\ge 0$, 
		\begin{equation*}
			N_Y(1,t) = 1+\hat{\delta}^c_Y(t)\le 1+\hat{\delta}_Y(t).
		\end{equation*}
		Proposition \ref{pro:separation-AMUC} tells us that for any $y,z \in Y$ and any $\gamma$-separated bounded sequence $\yn$ in $Y$
		\begin{equation}
			\label{eq:separation-NF}
			\liminf_{n\to \infty}(\norm{y-y_n}_Y + \norm{z-y_n}_Y)\ge N_Y(\norm{y-z}_Y, \gamma).    
		\end{equation}
		Assume now that $X$ coarse-Lipschitz embeds into $Y$. Then, one can assume without loss of generality that there exists $f\colon X\to Y$ continuous and a constant $K\ge 1$ such that $f(0)=0$ and for all $ x,z\in X$,
		\begin{equation*}
			\norm{x-z}_X - 1 \le \norm{ f(x)-f(z) }_Y \le K\norm{x-z}_X + 1.
		\end{equation*}
		The key step of the proof will be to show that for all $k\in \bN$,
		\begin{equation}
			\label{eq1:CL-AMUC}
			N_Y\Big(1,\frac{\theta}{8\sigma_kK}\Big)\le 1+\frac{2}{k}.
		\end{equation}
		Indeed, assuming \eqref{eq1:CL-AMUC} holds, we have that $\hat{\delta}^c_Y\Big(\frac{\theta}{8\sigma_kK}\Big)\le\frac{2}{k}$ and, by convexity of $\hat{\delta}^c_Y$, $\hat{\delta}^c_Y\Big(\frac{\theta}{16\sigma_kK}\Big)\le \frac{1}{k}$. In particular, $\norm{ e_1 + \dots + e_k}_{\ell_{\hat{\delta}^c_Y}} \le 16K\theta^{-1}\sigma_k$ and from the definition of $\sigma_k$ we deduce that \eqref{eq:CL-AMUC} is satisfied with $c:=\frac{1}{32K}$.
		
		\medskip
		Inequality  \eqref{eq1:CL-AMUC} is derived from a clever application of the approximate midpoint principle, stated in Lemma \ref{lem:CLApproxMid}, for a well-chosen map between the normed space $(c_{00},\Lambda_{N_k})$ and the Banach space $L_1(\Delta;Y)$. So, fix $k\in \bN$ and consider $g \colon (c_{00},\Lambda_{N_k})\to L_1(\Delta;Y)$ defined by $g:=f\circ T$. Since $f$ is continuous and $T\xi$ simple for all $\xi \in c_{00}$, it is easy to see that $g$ is well defined, measurable, but also coarse-Lipschitz since for all $\xi,\eta \in c_{00}$,
		\begin{equation*}
			\norm{ g(\eta) - g(\xi) }_{L_1(\Delta; Y)} \le K \norm{ T(\eta-\xi) }_{L_1(\Delta; X)} + 1 \le 4K \norm{ \eta-\xi }_{\Lambda_{N_k}} + 1.
		\end{equation*}
		Note also that since $f(0)=0$, we have for all $s\ge 0$,
		\begin{equation*}
			\norm{ g(se_1) }_{L_1(\Delta; Y)} = \frac12 \big(\norm{ f(sx_1) }_Y + \norm{ f(-sx_1) }_Y\big)\ge s-1,
		\end{equation*}
		which shows that $\Lip_\infty(g)\ge 1$. The approximate midpoint principle asserts that given $t_0>0$, there exist $t>t_0$ and $\eta,\zeta \in c_{00}$ so that 
		\begin{equation*}
			\norm{\eta-\zeta}_{\Lambda_{N_k}} = 2t
		\end{equation*}
		and
		\begin{equation*}
			g\Big(\Mid\big(\eta,\zeta,\frac1k\big)\Big)\subset \Mid\Big(g(\eta),g(\zeta), \frac2k\Big).
		\end{equation*}
		Let $\xi :=\frac12(\eta+\zeta)$ and $m \in \bN$ such that $\eta,\zeta\in \spa\{e_1,\dots,e_{m-1}\}$. Then,  for any $j\ge m$, $\xi+t\sigma_k^{-1}e_j \in \Mid(\eta,\zeta,\frac1k)$. Indeed,
		\begin{align*}
			\norm{ \eta - (\xi + t\sigma_k^{-1}e_j) }_{\Lambda_{N_k}} & =  N_k\Big(\Big\| \frac{\eta-\xi}{2} \Big\|_{\Lambda_{N_k}}, \frac{ \abs{t} }{\sigma_k}\Big) = \norm{ \zeta - (\xi + t\sigma_k^{-1}e_j) }_{\Lambda_{N_k}}\\
			& = N_k(t,t\sigma_k^{-1}) = t(1+F_k(\sigma_k^{-1})) = \Big\|\frac{\eta-\zeta}{2}\Big\|_{\Lambda_{N_k}}\Big(1+\frac1k\Big).  
		\end{align*}
		Consequently, the collection of functions $(h_j)_{j\ge m}$, where 
		\begin{equation*}
			h_j := f\circ\Big(\sum_{i=1}^{m-1}\xi_i\eps_i\otimes x_i + t\sigma_k^{-1}\eps_j\otimes x_j\Big)
		\end{equation*}
		lies in $\Mid(g(\eta),g(\zeta), \frac2k)$. We would like to apply \eqref{eq:separation-NF} pointwise to $(h_j)_{j\ge m}$ but there is a separation issue. To remedy this issue, we slightly modify $(h_j)_{j\ge m}$ and set 
		\begin{equation*}
			\tilde{h}_j := f\circ \Big(\sum_{i=1}^{m-1}\xi_i\eps_i\otimes x_i + t\sigma_k^{-1}\eps_m\otimes x_j\Big),
		\end{equation*}
		for all $j\ge m$.
		Observe that $\tilde{h}_j \in L_1(\{-1,1\}^m;Y)$ and since $\xn$ is $\theta$-separated, we have for all $\omega \in \Delta$ and $i>j\ge m$,
		\begin{equation*}
			\norm{ \tilde{h}_i(\omega) - \tilde{h}_j(\omega)}_Y \ge \theta t\sigma_k^{-1}-1.
		\end{equation*} 
		If $t_0>\sigma_k\theta^{-1}$, then $\theta t\sigma_k^{-1}-1>0$ and thus for all $\omega \in \Delta$: 
		\begin{align*}
			\liminf_{j \to \infty}(\bnorm{ g(\eta)(\omega) - \tilde{h}_j(\omega)}_Y + & \bnorm{ g(\zeta)(\omega) - \tilde{h}_j(\omega)}_Y)\\
			& \stackrel{\eqref{eq:separation-NF}}{\ge} N_Y(\norm{ g(\eta)(\omega) - g(\zeta)(\omega)}_Y,\theta t\sigma_k^{-1}-1).
		\end{align*}
		It follows from integration over $\omega \in \{-1,1\}^m$ and convexity (Jensen's inequality to be more precise) that
		\begin{align}
			\label{eq2:CL-AMUC}
			\liminf_{j\to\infty} (\bnorm{ g(\eta) - \tilde{h}_j}_{L_1(Y)} + \bnorm{ g(\zeta) - \tilde{h}_j }_{L_1(Y)}) \ge N_Y(\norm{ g(\eta) - g(\zeta) }_{L_1(Y)}, \theta t\sigma_k^{-1}-1).
		\end{align}
		Since $g(\eta)$ and $g(\zeta)$ only depend on the first $(m-1)$ coordinates in $\Delta$, it is easy to see that for all $j\ge m$, 
		$\norm{ \tilde{h}_j - g(\eta)}_{L_1(Y)} = \norm{ h_j - g(\eta)}_{L_1(Y)}$ and $\norm{ \tilde{h}_j - g(\zeta)}_{L_1(Y)} = \norm{ h_j - g(\zeta)}_{L_1(Y)}$. This means that $\tilde{h}_j \in \Mid(g(\eta),g(\zeta), \frac2k)$ whenever $j\ge m$ and incorporating this piece of information in \eqref{eq2:CL-AMUC} one gets
		\begin{equation}
			\label{eq3:CL-AMUC}
			N_Y(\norm{ g(\eta) - g(\zeta) }_{L_1(Y)}, \theta t\sigma_k^{-1}-1) \le
			\Big(1+\frac{2}{k}\Big) \norm{ g(\eta) - g(\zeta)}_{L_1(Y)}.
		\end{equation}
		At this point, it is good to remember that $\norm{ g(\eta) - g(\zeta)}_{L_1(Y)} \le 4K \norm{ \eta-\zeta}_{\Lambda_k} + 1 = 8Kt + 1$. Since one can check that for any $u>0$ the map $s\mapsto N_Y(s,u)-s$ is nonincreasing, we deduce from \eqref{eq3:CL-AMUC} that
		\begin{equation*}
			N_Y(8Kt+1,\theta t\sigma_k^{-1}-1)-(8Kt+1) \le \frac{2(8Kt+1)}{k}.
		\end{equation*}
		and by positive homogeneity of $N_Y$, one gets
		\begin{equation*}
			N_Y\Big(1,\frac{\theta t\sigma_k^{-1}-1}{8Kt+1}\Big) \le 1+\frac2k,
		\end{equation*}
		Letting $t_0$ (and therefore $t$) tend to $\infty$, we obtain \eqref{eq1:CL-AMUC}, thereby completing the proof.
	\end{proof}
	
	Since spreading models generated by weakly null sequences are unconditional, Theorem \ref{thm:CL-AMUC} has an immediate corollary on the quantitative behavior of weakly null spreading models in the domain space in terms of the modulus of asymptotic midpoint uniform convexity of the target space. It is instructive to compare Corollary \ref{cor:lower-estimates-spreading-models} with Theorem \ref{thm:upper-estimates-asymp-models}.
	
	\begin{coro}
		\label{cor:lower-estimates-spreading-models}
		Let $X$ and $Y$ be two Banach spaces so that $X$ coarse-Lipschitz embeds into $Y$. Then, there exists a constant $C>0$ such that for any spreading model $S$ of a normalized weakly null sequence in $X$ and for all $k\in \bN$,
		\begin{equation*}
			\Big\|\sum_{i=1}^k e_i\Big\|_{\ell_{\hat{\delta}^c_Y}} \le C\Big\|\sum_{i=1}^k e_i\Big\|_S.
		\end{equation*}
	\end{coro}
	
	\section{Weak Banach-Saks properties}
	\label{sec:weak-BS}
	
	In Theorem \ref{thm:CL-rigidity-subspaces-lp>2} of Section \ref{sec:J-concentration-consequences}, we have shown that, for $p\in[2,\infty)$, if a Banach space $X$ admits a coarse-Lipschitz embedding into $\ell_p$, then $X$ is isomorphic to a subspace of $\ell_p$. In a nutshell, the proof goes by saying that $X$ must be isomorphic to a subspace of $L_p$ that does not contain a copy of $\ell_2$ and that such spaces are necessarily subspaces of $\ell_p$ as shown by Johnson and Odell. This argument relies on asymptotic smoothness considerations and breaks down when $p\in(1,2)$. Trying to salvage this argument for the range $p\in(1,2)$, we can use duality and claim that $X^*$ must be a quotient of $L_q$, where $q\in(2,\infty)$ is the conjugate exponent. If we knew whether $X^*$ had the weak $q$-Banach-Saks property, we could invoke a classical result of W.B. Johnson (see Theorem III.2 and Remark III.3 in \cite{Johnson1976}) which ensures that (under this additional assumption) $X^*$ must be isomorphic to a quotient of $\ell_q$ and by duality $X$ must be isomorphic to a subspace of $\ell_p$. This is just one example of the importance of the weak Banach-Saks property when establishing nonlinear rigidity results. In the next section, we delve deeper into Banach-Saks-like properties, their duality and their relationships with the asymptotic properties studied so far. These results will be crucial in obtaining several profound rigidity theorems. 
	
	The weak Banach-Saks property seems to have been considered first by Rosenthal and is sometimes referred to as the Banach-Saks-Rosenthal property. In \cite{Kalton2013b}, the adjective ``weak'' was omitted but to avoid confusion with the original Banach-Saks property we will not omit it in the sequel. More precisely, we are interested in a quantification of the weak Banach-Saks property.
	
	\begin{defi} 
		\label{def:wBS}
		Let $p \in (1,\infty)$. We say that a Banach space $X$ not containing $\ell_1$ has the \emph{weak $p$-Banach-Saks property} if there is a constant $C>0$ so that every normalized weakly null sequence
		$\xn$ has a subsequence $(x_{n_j})_{j=1}^{\infty}$ such that for all $k \in \bN$,  \begin{equation*}
			\Big\| \sum_{j=1}^k x_{n_j} \Big\| \le Ck^{1/p}.
		\end{equation*}
	\end{defi}
	
	In the reflexive setting, the weak $p$-Banach-Saks property has a dual notion whose definition is given below. 
	
	\begin{defi} 
		\label{def:wcoBS}
		Let $p \in (1,\infty)$. We say that a Banach space $X$ not containing $\ell_1$ has the \emph{weak $p$-co-Banach-Saks property} if there is a constant $c>0$ so that every normalized weakly null sequence $\xn$ has a subsequence $(x_{n_j})_{j=1}^{\infty}$ such that for all $k \in \bN$,  
		\begin{equation*}
			\Big\| \sum_{j=1}^k x_{n_j} \Big\| \ge ck^{1/p}.
		\end{equation*}
	\end{defi}
	
	\begin{rema}
		In terms of spreading models, a Banach space $X$ not containing $\ell_1$ has the weak $p$-Banach-Saks property if and only if there is a constant $C>0$ so that for every spreading model $S :=[(e_i)_{i=1}^{\infty}]$ generated by a normalized weakly null sequence in $X$ and every $k\in \bN$, we have $\big\|\sum_{i=1}^k e_i\big\|_S\le Ck^{1/p}$ and similarly, $X$ has the weak $p$-co-Banach-Saks property if and only if there is a constant $c>0$ so that for every spreading model $S :=[(e_i)_{i=1}^{\infty}]$ generated by a normalized weakly null sequence in $X$ and every $k\in \bN$, we have $\big\|\sum_{i=1}^k e_i\big\|_S\ge ck^{1/p}$.
	\end{rema}

	\begin{prop}
		\label{pro:wBS-duality}
		Let $p \in (1,\infty)$, $q$ its conjugate exponent and $X$ be a reflexive Banach space. If $X$ has the weak $p$-Banach-Saks property, then $X^*$ has the weak $q$-co-Banach-Saks property. 
	\end{prop}
	
	\begin{proof}  
		Assume that $X$ has the weak $p$-Banach-Saks property with constant $C>0$. Let $(x_n^*)_{n=1}^{\infty}$ be a normalized weakly null sequence in $X^*$. We may pick a normalized sequence $\xn$ in $X$ with $x_n^*(x_n)=1$.
		Since $X$ is reflexive, we can assume, by passing to a subsequence, that $\xn$ converges weakly to $x\in X$. Then, $\norm{ x_n-x }\le 2$ and so passing to a further subsequence, we can assume that for all $k\in \bN$,
		\begin{equation*}
			\lim_{n_1\to\infty}\cdots\lim_{n_k\to\infty} \norm{ x_{n_1} + \dots + x_{n_k} - kx}\le 2Ck^{1/p}.
		\end{equation*}  
		However
		\begin{equation*}
			\lim_{n_1\to\infty}\cdots\lim_{n_k\to\infty}\langle x_{n_1}+\dots+x_{n_k}-kx,x_{n_1}^*+\dots+x_{n_k}^*\rangle=k,
		\end{equation*}  
		which implies that for any spreading model $(e_j)_{j=1}^n$ of $(x^*_j)_{j=1}^{\infty}$ we must have
		\begin{equation*}
			\Big\|\sum_{j=1}^k e_j \Big\|_S\ge \frac{1}{2C}k^{1/q}.
		\end{equation*}
	\end{proof}
	
	We describe below some simple links between these Banach-Saks properties and asymptotic uniform smoothness and convexity. Since this result will not be needed, we leave the proof as Exercise \ref{exo:AUS-BS}.
	\begin{prop}\label{prop:AUS-BS} Let $X$ be a Banach space with separable dual and $p\in (1,\infty)$.
		\begin{enumerate}[(i)]
			\item If $X$ is $p$-asymptotically uniformly smoothness, then $X$ has the weak $p$-Banach-Saks property. 
			\item If $X$ is $p$-asymptotically uniformly convex, then $X$ has the weak $p$-co-Banach-Saks property. 
		\end{enumerate}
	\end{prop}
	
	With these definitions in hand, we deduce from Theorem \ref{thm:upper-estimates-asymp-models} the following statement.
	
	\begin{theo}
		\label{thm:CL-into-ref+pAUS->ref+weak-pBS}
		Let $p\in (1,\infty)$. If a Banach space $X$ coarse-Lipschitz embeds into a reflexive $p$-asymptotically uniformly smooth Banach space, then $X$ is reflexive and has the weak $p$-Banach-Saks property.
	\end{theo}
	
	We can also rephrase Corollary \ref{cor:lower-estimates-spreading-models}.
	
	\begin{theo}
		\label{thm:CL-into-pAMUC->weak-pcoBS} 
		Let $p\in (1,\infty)$. If a Banach space $X$ coarse-Lipschitz embeds into a $p$-asymptotically midpoint uniformly convex Banach space, then $X$ has the weak $p$-co-Banach-Saks property.  
		%Let $X$ and $Y$ be two Banach spaces so that $X$ coarse-Lipschitz embeds into $Y$ and let $p \in (1,\infty)$. Assume that $Y$ is $p$-AMUC. Then, $X$ has the weak $p$-co-Banach-Saks property. 
	\end{theo}
	
	Understanding under which conditions the weak $p$-co-Banach-Saks property automatically implies the existence of a $p$-asymptotically (midpoint) uniformly convex norm will be of critical importance in order to obtain new coarse-Lipschitz rigidity results. One such sufficient condition is discussed in the next section.  
	
	\section{\texorpdfstring{Random $L_p$-norms}{Random norms}}
	
	The notion of random $L_p$-norm first appeared in the works of Haydon, Raynaud and Levy on ultraproducts in \cite{LevyRaynaud1984}, \cite{HaydonLevyRaynaud1984} and \cite{HaydonLevyRaynaud1991}.
	
	\begin{defi} 
		\label{def:random-Lp-norm}
		Let $p\in (1,\infty)$. We say that a Banach space $X$ has a \emph{random $L_p$-norm} if there is a map $V\colon X\to L_p(\Omega)$ where $\Omega:=(\Omega, \Sigma, \mu)$ is a measure space such that
		\begin{enumerate}
			\item $V(x)\ge 0$ for all $x\in X$,
			\item $V(\lambda x)=\abs{\lambda}V(x)$, for all $x\in X$ and $\lambda\in\bR$,
			\item $V(x_1+x_2)\le V(x_1) + V(x_2)$ for all $x_1,x_2\in X$,
			\item $\norm{V(x)}_{L_p} = \norm{x}_X$ for all $x\in X$.
		\end{enumerate}
		The map $V$, which does not have to be linear, is said to be a random $L_p$-norm on $X$ or to realize $\norm{\cdot}_X$ as a random $L_p$-norm.
	\end{defi}
	
	\begin{rema}
		Note that it easily follows from the definition that for all $x_1,x_2 \in X$, $\abs{V(x_1)-V(x_2)}\le V(x_1-x_2)$ and we deduce from condition $4.$ that $V\colon X \to L_p(\Omega)$ is $1$-Lipschitz. Therefore, if $X$ is separable, so is $V(X)$ and we can replace $L_p(\Omega)$ by $L_p :=L_p[0,1]$.
	\end{rema} 
	
	Random $L_p$-norms arise naturally in the context of vector-valued integration. Given a Banach space $X$ and an $X$-valued Bochner space $L_p(\Omega; X)$, the map $V\colon L_p(\Omega; X) \to L_p(\Omega)$ defined by $V(f)=\|f\|_X$, for  $f\in L_p(\Omega; X)$, realizes $\norm{\cdot}_{L_p(\Omega;X)}$ as a random $L_p$-norm.

	Random $L_p$-norms whose associated maps $V$ satisfy a pointwise version of the $p$-uniform smoothness inequality will play an important role in the sequel. For a refresher about the standard notion of $p$-uniform smoothness and $q$-uniform convexity, we refer the reader to Appendix \ref{appendix:us-uc}.
	
	\begin{defi} 
		Let $p,r\in (1,\infty)$. A random $L_p$-norm $V$ on a Banach space $X$ is said to be of \emph{type $r$} if there is a constant $C\ge 1$ such that for all $x_1,x_2\in X$, 
		\begin{equation}
			\label{eq:random-Lp-norm-type-r}
			\frac{V(x_1+x_2)^r + V(x_1-x_2)^r}{2}\le V(x_1)^r + C^r V(x_2)^r.
		\end{equation}
	\end{defi}
	
	We start with elementary properties of random $L_p$-norms of type $r$. The first inequality implies that a random $L_p$-norm of type $r$ is a random $L_p$-norm of type $s$ for all $1\le s\le r$.
	
	\begin{prop}
		\label{prop:random-Lp-norm-inequalities} 
		Let $p,r\in(1,\infty)$ and assume that $V$ is a random $L_p$-norm of type $r$ on a Banach space $X$ for some constant $C\ge 1$. Then, 
		\begin{enumerate}[(i)]
			\item for all $x_1,x_2\in X$ and $1\le s\le r$,
			\begin{align*}
				\frac{V(x_1+x_2)^s + V(x_1-x_2)^s}{2}&\le \big(V(x_1)^r + C^r V(x_2)^r\big)^{s/r}\\
				&\le V(x_1)^s+C^sV(x_2)^s,
			\end{align*}
			\item for all $x_1,\dots,x_k\in X$ and $1\le s\le p\le  r$
			\begin{equation*}
				\bE_\vep \Big\|V\Big(\sum_{j=1}^k \vep_j x_j\Big)\Big\|_{L_s}^s \le C^s \Big\|\Big(\sum_{j=1}^k V(x_j)^r\Big)^{1/r}\Big\|_{L_s}^s,
			\end{equation*}
			where $\bE_\vep$ is the average over all choices of signs $(\vep_j)_{j=1}^k$, or equivalently the integration over independent Rademacher variables $(\vep_j)_{j=1}^k$.
		\end{enumerate}
	\end{prop}
	
	\begin{proof} 
		The proof of $(i)$ is immediate by using the concavity of the function $t\mapsto t^{s/r}$ on $[0,+\infty)$ and the fact that $C\ge 1$.
		
		The proof of $(ii)$ is essentially the same. We first use the concavity of $t\mapsto t^{s/r}$ to get that 
		\begin{equation*}
			\Big(\bE_\vep V\Big(\sum_{j=1}^k \vep_j y_j\Big)^s\Big)^{r/s} \le \bE_\vep V\Big(\sum_{j=1}^k \vep_j y_j\Big)^r.
		\end{equation*}
		Since $V$ is a random $L_p$-norm of type $r$ and constant $C\ge 1$, a straightforward induction on $k$ gives that
		\begin{equation*}
			\bE_\vep V\Big(\sum_{j=1}^k \vep_j x_j\Big)^r \le C^r \sum_{j=1}^k V(x_j)^r.
		\end{equation*}
		Therefore, 
		\begin{equation*}
			\bE_\vep V\Big(\sum_{j=1}^k \vep_j x_j\Big)^s \le C^s \Big(\sum_{j=1}^k V(x_j)^r\Big)^{s/r},
		\end{equation*}
		and the inequality follows by integrating over $\Omega$.
	\end{proof}
	
	Being a random $L_p$-norm of type $r$ passes to subspaces and ultrapowers. This second point will be crucial in applications.
	
	\begin{prop}
		\label{prop:random-norm-stability-subspaces-ultraproducts} 
		Let $1<p\le r<\infty$ and $(X,\norm{\cdot}_X)$ be a Banach space with a random $L_p$-norm of type $r$.
		\begin{enumerate}[(i)]
			\item If $Y$ is a subspace of $X$, then the induced norm on $Y$ is a random $L_p$-norm of type $r$.
			\item For every nonprincipal ultrafilter $\cU$ on $\bN$, $\norm{\cdot}_{X^{\cU}}$ is a random $L_p$-norm of type $r$ on $X^{\cU}$ the ultrapower of $X$ with respect to $\cU$.
		\end{enumerate}
	\end{prop}
	
	\begin{proof} 
		The first point is immediate by restricting to $Y$ the map $V\colon X\to L_p(\Omega)$ realizing the fact that the norm of $X$ is a random $L_p$-norm of type $r$.
		
		For the second point, recall first that for a nonprincipal ultrafilter $\cU$, the ultrapower $L_p(\Omega)^{\cU}$ is an $L_p$-space. Then, define $V^{\cU}\colon X^{\cU} \to L_p(\Omega)^{\cU}$, by letting
		\begin{equation*}
			V^{\cU}([(x_n)_{n}]^{\cU}) := [(V(x_n)_n)]^{\cU}
		\end{equation*} 
		Routine verifications show that the norm of $\norm{\cdot}_{X^{\cU}}$ is a random $L_p$-norm of type $r$ realized by $V^{\cU}$.
	\end{proof}
	
	The notion of random $L_p$-norm of type $r$ turns out to be a strengthening of the classical notion of $p$-uniform smoothness.
	
	\begin{prop}
		\label{prop:random-Lp-type->pUS} 
		Let $p\in(1,2]$. If $X$ is a Banach space whose norm is a random $L_p$-norm of type $r$, for some $r\in[p,2]$, then the norm of $X$ is $p$-uniformly smooth. 
	\end{prop}
	
	\begin{proof} 
		Assume that the norm of $X$ is a random $L_p$-norm of type $r$ with constant $C>0$ that is realized by a map $V\colon X\to L_p(\Omega,\mu)$. Since $\frac{p}{r} \in (0,1]$, it follows from assertion $(i)$ of Proposition \ref{prop:random-Lp-norm-inequalities} that the norm is a random $L_p$-norm of type $p$ and thus for all $x_1,x_2 \in X$, 
		\begin{align*} 
			\frac{\norm{x_1 + x_2}_X^p + \norm{x_1 - x_2}_Y^p}{2} & = \int_{\Omega}\frac{V(x_1+x_2)^p + V(x_1-x_2)^p}{2}\,d\mu\\    
			%                           &\le \int_{\Omega}\big(V(x_1)^r + C^rV(x_2)^r\big)^{p/r}\,d\mu \\
			& \le \int_{\Omega}\big(V(x_1)^p + C^pV(x_2)^p\big)d\mu = \norm{x_1}_X^p + C^p\norm{x_2}_X^p
		\end{align*}
		%Since $\alpha:=\frac{p}{r} \in (0,1]$, it follows from homogeneity and elementary calculus that $(s+t)^\alpha \le s^\alpha+t^\alpha$, for all $s,t\ge 0$.
		This implies that $X$ is $p$-uniformly smooth.  
	\end{proof}
	
	We now describe an important class of Banach spaces with random $L_p$-norms with type $r$. Recall that a family of Banach spaces $(X_i)_{i\in I}$ is said to be \emph{equi-$p$-uniformly smooth} if there exists $C>0$ such that for all $i\in I$ and all $t\ge 0$, $\rho_{X_i}(t)\le Ct^p$. 
	
	It is a classical fact that the $\ell_q$-sum of equi-$p$-uniformly smooth Banach spaces is $\min\{p,q\}$-uniformly smooth. In the light of Proposition \ref{prop:random-Lp-type->pUS}, the next proposition is a strengthening of this classical fact.
	
	\begin{prop}
		\label{prop:lp-sum-random-Lp-norm-type->random-Lp-norm-type} 
		Let $1<p\le r\le 2$ and $(X_n)_{n=1}^\infty$ be a sequence of equi-$r$-uniformly smooth Banach spaces. Then, the canonical norm of $(\sum_{n=1}^\infty X_n)_{\ell_p}$ is a random $L_p$-norm of type $r$.
	\end{prop}
	
	\begin{proof} 
		For $x :=(x_n)_{n=1}^\infty \in (\sum_{n=1}^\infty X_n)_{\ell_p}$, let $V(x) := (\norm{x_n}_{X_n})_{n=1}^\infty \in \ell_p$. As already discussed after Definition \ref{def:random-Lp-norm}, $V$ realizes the canonical norm on $(\sum_{n=1}^\infty X_n)_{\ell_p}$ as a random $L_p$-norm. Let us show that it is of type $r$. First, the $X_n$ being equi-$r$-uniformly smooth, we have that there exists a constant $D\ge 1$ such that for all $n$ and all $x_n,y_n\in X_n$:
		\begin{equation*}
			\frac{\norm{x_n+y_n}_{X_n}^r + \norm{x_n-y_n}_{X_n}^r}{2}\le \norm{x_n}_{X_n}^r + D^r\norm{y_n}_{X_n}^r.
		\end{equation*}
		The previous inequality is equivalent to the fact that for all $x :=(x_n)_{n=1}^\infty$ and $y :=(y_n)_{n=1}^\infty$ in $(\sum_{n=1}^\infty X_n)_{\ell_p}$,
		\begin{equation*}
			\frac{V(x+y)^r + V(x-y)^r}{2}\le V(x)^r + D^r V(y)^r,   
		\end{equation*}
		meaning that $V$ is a random $L_p$-norm of type $r$.
	\end{proof}
	
	The key technical result of this section is that $p$-asymptotic uniform convexity is implied by the weak $p$-co-Banach Saks property under an additional smoothness assumption. This result is truly remarkable in the sense that one deduces an $\ell_p$ lower estimate on branches of \emph{arbitrary} infinite weakly null trees only from  an $\ell_p$ lower estimate for constant coefficients on  subsequences of weakly null sequences.
	
	\begin{theo}
		\label{theo:coBS->AUC}
		Let $1<p < r\le 2$ and $Y$ be a separable Banach space with a random $L_p$-norm of type $r$. If $X$ is a quotient of $Y$ and has the weak $p$-co-Banach Saks property, then $X$ is $p$-asymptotically uniformly convex.
	\end{theo}
	
	\begin{proof} 
		Under the assumption of the theorem, we will show that $\bar{\delta}_{X}(t)\gtrsim t^p$. Since $X$ must be separable and reflexive, it is sufficient to show that there exists a universal constant $\gamma>0$ such that for every $t\in (0,1]$, every $x\in S_X$ and every weakly null sequence $(x_n)_n$ in $tS_X$, $\liminf_{n\to\infty}\norm{x+x_n}_X^p \ge 1 + \gamma^p t^p$. 
		
		%Let $Z$ be the subspace of $Y$ such that $X=Y/Z$ and $Q\colon Y\to X=Y/Z$ be the corresponding quotient map. 
		Since $Y$ is reflexive, it follows from a weak compactness argument and the weak lower semi-continuity of $\norm{\cdot}_Y$ that we can pick $v_n\in Y$ so that $\norm{v_n}_Y = \norm{x+x_n}_X$ and $Q(v_n) = x + x_n$, where $Q\colon Y\to X$ is the quotient map. Note that $\sup_{n\ge 1}\norm{v_n}_Y\le 2$ and we may also assume that $(v_n)_n$ converges weakly to $y\in Y$ with $\norm{y}_Y\le 2$. Letting $y_n:= v_n-y$, we have that $\norm{y_n}_Y\le 4$ and  since $Q$ is weak to weak continuous, it follows that $Q(y)=x$ and $Q(y_n)=x_n$ for all $n\ge 1$.
		
		Assume now that $\norm{\cdot}_Y$ is a random $L_p$-norm of type $r$ with constant $C>0$ that is realized by $V\colon Y\to L_p[0,1]$. Then, it follows from the observations above that
		\begin{equation}
			\label{eq0:random-norm}
			\liminf_{n\to\infty}\norm{x+x_n}_X^p = \liminf_{n\to\infty}\norm{y+y_n}_Y^p = \liminf_{n\to\infty}\norm{V(y+y_n)}^p_p.
		\end{equation}
		
		The heart of the matter is to prove the following claim.
		
		\begin{claim}
			\label{claim:random-norm}
			There is a constant $\gamma>0$, independent of $t$,  and a sequence of measurable sets $(G_n)_n$ in $[0,1]$ such that 
			$$\lim_{n\to \infty}\lambda(G_n)=0\ \ \text{and}\ \ \liminf_{n\to \infty} \norm{\car_{G_n}V(y_n)}_p \ge \gamma t.$$
		\end{claim}
		
		Assuming Claim \ref{claim:random-norm}, we can conclude as follows. First, observe that by disjointness we have 
		\begin{equation}
			\label{eq:random-norm}
			\liminf_{n\to\infty}\norm{V(y+y_n)}^p_p \ge \underbrace{\liminf_{n\to\infty}\norm{\car_{G_n}V(y+y_n)}^p_p}_{A} + \underbrace{\liminf_{n\to\infty}\norm{\car_{[0,1]\setminus G_n}V(y+y_n)}^p_p}_{B}.
		\end{equation}
		
		We can estimate $A$ from below simply by using the special properties of random $L_p$-norms. Indeed,
		\begin{equation*}
			\norm{\car_{G_n}V(y_n)}_p \le \norm{\car_{G_n}Vy}_p + \norm{\car_{G_n}V(y+y_n)}_p,
		\end{equation*}
		and it follows from Claim \ref{claim:random-norm} that 
		\begin{equation}
			\label{eq4:random-norm}
			\liminf_{n\to\infty}\norm{\car_{G_n}V(y+y_n)}_p\ge \gamma t.
		\end{equation}
		
		To estimate $B$ from below, Kalton introduced a very clever and useful semi-norm on $Y$. For all $z\in Y$, let 
		\begin{equation*}
			N_y(z) :=\int_0^1 V(y)(s)^{p-1}V(z)(s)\,ds.
		\end{equation*}
		H\"{o}lder's inequality implies that for all $z\in Y$,
		\begin{equation}
			\label{eq2:random-norm}
			N_y(z)\le \norm{y}_Y^{p-1}\norm{z}_Y.
		\end{equation}
		Note that it follows from $\eqref{eq2:random-norm}$ that $N_y$ is well-defined and that every bounded linear functional on $(Y,N_y)$ is in $Y^*$. Moreover,
		\begin{equation*}
			N_y(y)= \norm{Vy}_p^p = \norm{y}_Y^p.
		\end{equation*} 
		Therefore, by the Hahn-Banach Theorem, there is a linear functional $y^*\in (Y,N_y)^*$ with norm $1$ such that  
		\begin{equation*}
			y^*(y)= N_y(y) = \norm{y}_Y^p.
		\end{equation*}
		Since $(y_n)_{n\ge 1}$ is weakly null and $y^*\in Y^*$, we have 
		\begin{equation*}
			\norm{y}_Y^p=\lim_{n\to\infty}y^*(y+y_n)\le \liminf _{n\to\infty} N_y(y+y_n) = \liminf _{n\to\infty} \int_0^1 (Vy(s))^{p-1}V(y+y_n)(s)\,ds.
		\end{equation*}
		Applying H\"{o}lder's inequality again,
		\begin{equation*}
			\int_{G_n} V(y)(s)^{p-1}V(y+y_n)(s)\,ds \le
			\left(\int_{G_n}V(y)(s)^pds\right)^{1-1/p}\norm{y+y_n}_Y,
		\end{equation*}
		and, using the fact that $\lim_{n\to\infty}\abs{G_n}=0$, we deduce
		\begin{equation*}
			\lim_{n\to\infty}\int_{G_n} (Vy(s))^{p-1}V(y+y_n)(s)\,ds=0.
		\end{equation*}
		Therefore,
		\begin{equation*}
			\label{eq3:random-norm}
			\norm{y}_Y^p \le \liminf_{n\to\infty}\int_{[0,1]\setminus G_n}(Vy(s))^{p-1}V(y+y_n)(s)\,ds.
		\end{equation*}
		Invoking H\"{o}lder's inequality one more time, we have that
		\begin{align}
			\label{eq5:random-norm} \norm{y}_Y^p  & \le  \liminf_{n\to\infty}\norm{\car_{[0,1]\setminus G_n} V(y)}_p^{p-1}\norm{\car_{[0,1]\setminus G_n}V(y+y_n)}_p\\
			\notag  & \le \norm{y}_Y^{p-1}\liminf_{n\to\infty}\norm{\car_{[0,1]\setminus G_n}V(y+y_n)}_p,
		\end{align}
		where in the last inequality we remembered that $V(y)\ge 0$.
		Since $\norm{y}_Y \ge \norm{Qy}_X = \norm{x}_X = 1$, it follows from \eqref{eq5:random-norm} that  
		\begin{equation}
			\label{eq6:random-norm} 
			1 \le \liminf_{n\to\infty}\norm{\car_{[0,1]\setminus G_n}V(y+y_n)}_p.
		\end{equation}
		Combining \eqref{eq0:random-norm}-\eqref{eq:random-norm}-\eqref{eq4:random-norm}-\eqref{eq6:random-norm} we deduce that 
		\begin{equation*}
			1 + \gamma^pt^p \le \liminf_{n\to\infty}\norm{x+x_n}^p_X.
		\end{equation*}

		It remains to prove Claim \ref{claim:random-norm}. For this we state yet another intermediate claim.
		
		\begin{claim}\label{claim10} There exists a constant $\alpha>0$ such that for all $\theta \in (0,1)$, 
			\begin{equation}
				\label{eq7:random-norm} 
				\liminf_{n\to\infty}\norm{V(y_n)}_{p,\theta}\ge \alpha t,
			\end{equation}
			where for any $0<\theta<1$ and any $f\in L_p$, we let 
			\begin{equation}
				\label{eq8:random-norm} 
				\norm{f}_{p,\theta} := \sup_{\lambda(E)\le \theta}\norm{\car_E f}_p.
			\end{equation}
			Above and in the sequel, $\lambda(E)$ denotes the Lebesgue measure of a measurable subset $E$ of $[0,1]$.
		\end{claim}
		
		A first important point to make is that the infimum in \eqref{eq8:random-norm} is attained. This crucial observation follows from elementary considerations on the decreasing rearrangement of $\abs{f}$, which we detail below. 
		%We follow the notation from \cite[Section 2.a.]{LindenstraussTzafriri1979}. 
		Given any measurable function $\phi\colon [0,1] \to [0,\infty)$, the \emph{distribution function} of $\phi$ (with respect to the Lebesgue measure) is the non-decreasing and right-continuous map defined, for $t\in [0,\infty)$, by 
		\begin{equation*}
			\lambda_\phi(t) := \lambda (\phi >t) := \lambda(\{s\in [0,1]\colon \phi(s)>t\}).
		\end{equation*}
		Then, the \emph{decreasing rearrangement} of $\phi$ is given for all $s\in[0,1]$, by
		\begin{equation*}
			\phi^*(s) := \inf\{t>0 \colon \lambda_\phi(t)\le s\}.
		\end{equation*}
		The function $\phi^*$ is nonincreasing, right-continuous and has the same distribution function as $\phi$. Given $\theta \in (0,1)$, it follows from these definitions and basic measure theory that 
		\begin{equation*}
			\lambda(\{t\in[0,1]\colon \abs{f(t)} > \abs{f}^*(\theta)\})\le \theta
		\end{equation*} 
		and 
		\begin{equation*}
			\lambda(\{t\in [0,1]\colon \abs{f(t)}\ge \abs{f}^*(\theta)\}) \ge \theta.
		\end{equation*} 
		Therefore, there exists a measurable subset $E$ of $[0,1]$ with $\lambda(E) = \theta$ such that
		\begin{equation}
			\label{eq:maximalE}
			\{t\in [0,1]\colon \abs{f(t)} > \abs{f}^*(\theta)\} \subset E \subset  \{t\in[0,1] \colon \abs{f(t)} \ge \abs{f}^*(\theta)\}.
		\end{equation}
		Then, by comparing $f$ on $E\setminus F$ and $f$ on $F\setminus E$, it easily follows from \eqref{eq:maximalE} that for any measurable $F$ such that $\lambda(F)\le \theta$, $\norm{\car_F f}_p \le \norm{\car_E f}_p$.
		
		A second observation is that, given $f\in L_p$, the map $\psi\colon \theta\in[0,1] \mapsto \sup_{\lambda(E)\le \theta}\norm{\car_E f}_p$ is continuous thanks to the uniform integrability of $\abs{f}^p$. Consequently, if $0<a<\norm{f}_p$, there exists $\theta \in (0,1)$ such that $\psi(\theta)=a$. Consider now $\theta_a := \inf\{\theta \in (0,1)\colon \psi(\theta)=a\}$. We have that $\norm{f}_{p,\theta_a}=\psi(\theta_a)=a$ and if we pick $E_a$ measurable such that $\norm{\car_{E_a} f}_p = \norm{f}_{p,\theta_a}$, we clearly have that $E_a$ is a measurable subset of $[0,1]$ with minimal measure such that $\norm{\car_{E_a} f}_p=a$. 
		
		\smallskip 
		After these two crucial observations are made, we can quickly derive the conclusion of Claim \ref{claim:random-norm} from Claim \ref{claim10} as follows. There exists $n_\theta \in \N$ such that $\norm{V(y_n)}_{p} \ge \norm{V(y_n)}_{p,\theta} > \frac{\alpha}{2}t$ for all $n \ge n_\theta$. In particular, for all $n \ge n_1$, $\norm{V(y_n)}_{p} > \frac{\alpha}{2}t$. Then, as we explained above, for $n\ge n_1$, we can pick $G_n$ measurable and of minimal measure such that $\norm{\car_{G_n}Vy_n}_p=\frac{\alpha}{2}t$. By \eqref{eq:maximalE}, we have that for $n\ge n_\theta$ there exists $E_n$ measurable with $\lambda(E_n) = \theta$ such that $\norm{V(y_n)}_{p,\theta}=\norm{\car_{E_n}Vy_n}_p>\frac{\alpha}{2}t$. Moreover, by considering sets of the form $[0,t]\cap E_n$, we can find a measurable subset $F_n$ of $E_n$ (and thus satisfying $\lambda(F_n)\le \theta$) such that $\norm{\car_{F_n}Vy_n}_p = \frac{\alpha}{2}t$. But, the minimality of $\lambda(G_n)$ with respect to the latter property implies that $\lambda(G_n)\le \theta$ for all $n\ge n_\theta$. We have constructed $(G_n)_n$ so that  $\lambda(G_n)\le \theta$ for all $n\ge n_\theta$ and $\norm{\car_{G_n}Vy_n}_p=\frac{\alpha}{2}t$ for all $n\ge n_1$. This concludes the proof of Claim \ref{claim:random-norm}.

		\smallskip 
		We finally turn to the proof of Claim \ref{claim10}, for which the weak $p$-co-Banach-Saks property is needed. So, assume that $X$ has the weak $p$-co-Banach-Saks property with constant $c>0$. Fix $\theta \in (0,1)$. By passing to a first subsequence, we may assume that $\liminf_{n\to\infty}\norm{V(y_n)}_{p,\theta}=\lim_{n\to\infty}\norm{V(y_n)}_{p,\theta}$. Let $(e_i)_{i=1}^\infty$ be a spreading model of $(x_n)_{n=1}^\infty$. Since $(x_n)_{n=1}^\infty$ is weakly null, $(e_i)_{i=1}^\infty$ is $2$-unconditional and by the weak $p$-co-Banach-Saks property of $X$ and passing to a further subsequence (using Proposition \ref{prop:weakly-null-spreading-model} and a diagonal argument), we may assume that for every $k\in \bN$ there is $N_k \in \N$ such that for all natural numbers $N_k\le n_1<n_2<\dots<n_k$, 
		\begin{equation*}
			\bE_\vep\Big\|\sum_{j=1}^k\vep_j x_{n_j}\Big\|_X\ge \frac{ct}{3}k^{1/p},
		\end{equation*}
		and, since $Q(y_{n_j})=x_{n_j}$ and $\norm{Q}=1$, 
		\begin{equation*}
			\bE_\vep\Big\|\sum_{j=1}^k\vep_j y_{n_j}\Big\|_Y\ge \frac{ct}{3}k^{1/p}.
		\end{equation*}
		From Proposition \ref{prop:random-Lp-norm-inequalities}, we deduce that
		\begin{align*} 
			\Big(\bE_\vep\Big\|\sum_{j=1}^k\vep_j y_{n_j}\Big\|_Y^p\Big)^{1/p}
			&=\Big(\bE_\vep\Big\|V\Big(\sum_{j=1}^k\vep_j y_{n_j}\Big)\Big\|_p^p\Big)^{1/p} \\
			& \le C\Big\|\Big(\sum_{j=1}^k V(y_{n_j})^r\Big)^{1/r}\Big\|_p,
		\end{align*}
		where $C$ is the constant of the random $L_p$-norm of type $r$.
		Therefore, for all natural numbers $N_k\le n_1<n_2<\cdots<n_k$, 
		\begin{equation*}
			\Big\|\Big(\sum_{j=1}^k V(y_{n_j})^r\Big)^{1/r}\Big\|_p\ge \frac{ct}{3C}k^{1/p}.
		\end{equation*}
		Our next goal is to estimate $\big\|\big(\sum_{j=1}^k V(y_{n_j})^r\big)^{1/r}\big\|_p$ from above. In order to achieve this, we need to estimate the contribution of each $V(y_{n_j})$ on a carefully chosen partition of $[0,1]$ adapted to $V(y_{n_j})$. For the rest of this argument, fix $0<\theta<1$ and for each $n \in \bN$, let $E_n$, chosen as in \eqref{eq:maximalE}, be a measurable subset of $[0,1]$ with measure $\lambda(E_n)=\theta$ such that $\norm{V(y_n)}_{p,\theta}=\norm{\car_{E_n}V(y_n)}_p$. Remember that $\sup_n \norm{V(y_n)}_p\le 4$ and that we assumed that $(\norm{V(y_n)}_{p,\theta})_n$ converges.
		
		It follows from \eqref{eq:maximalE} that $\theta^{1/p}\abs{V(y_n)}^*(\theta)\le \norm{\car_{E_n} V(y_n)}_p$ and observing that $\abs{V(y_n)} \le \abs{V(y_n)}^*(\theta)$ on $[0,1]\setminus E_n$, we deduce that 
		\begin{equation}
			\label{eq:maxcomplement}
			\norm{\car_{[0,1]\setminus E_n}V(y_n)}_\infty \le \norm{\car_{E_n} V(y_n)}_p\theta^{-1/p} \le \norm{V(y_n)}_p\theta^{-1/p}\le 4\theta^{-1/p}.   
		\end{equation}
		Now, letting $E^c_n :=[0,1]\setminus E_n$, it follows from \eqref{eq:maxcomplement} that
		\begin{equation}
			\label{eq:equation1_coBSAUC}
			\norm{\car_{E^c_n}V(y_n)}_r\le \norm{\car_{E^c_n}V(y_n)}_\infty \lambda(E^c_n)^{1/r}\le 4\theta^{-1/p}(1-\theta)^{1/r}.    
		\end{equation}
		Using the triangle inequality in $\ell_r$, we obtain that for any natural numbers $N_k\le n_1<\dots<n_k$, we have the following pointwise inequality
		\begin{align*}
			\Big(\sum_{j=1}^k V(y_{n_j})^r\Big)^{1/r} \le \Big(\sum_{j=1}^k \car_{E_{n_j}}V(y_{n_j})^r\Big)^{1/r} +\Big(\sum_{j=1}^k \car_{E^c_{n_j}}V(y_{n_j})^r\Big)^{1/r},
		\end{align*}
		and it follows from the triangle inequality in $L_p$ that
		\begin{align*}
			\Big\|\Big(\sum_{j=1}^k V(y_{n_j})^r\Big)^{1/r}\Big\|_p &  \le \Big\| \Big(\sum_{j=1}^k \car_{E_{n_j}}V(y_{n_j})^r\Big)^{1/r}\Big\|_p + \Big\|\Big(\sum_{j=1}^k \car_{E^c_{n_j}}V(y_{n_j})^r\Big)^{1/r}\Big\|_p.
		\end{align*}
		To estimate the contribution on the $E^c_{n_j}$'s we use the fact that $p\le r$ to obtain
		\begin{align*}
			\Big\|\Big(\sum_{j=1}^k V(y_{n_j})^r\car_{
				E^c_{n_j}}\Big)^{1/r}\Big\|_p &\le \Big\|\Big(\sum_{j=1}^k V(y_{n_j})^r\car_{E^c_{n_j}}\Big)^{1/r}\Big\|_r =\Big\|\sum_{j=1}^k V(y_{n_j})^r\car_{E^c_{n_j}}\Big\|_1^{1/r}\\
			&\le \Big(\sum_{j=1}^k \big\|V(y_{n_j})^r\car_{E^c_{n_j}}\big\|_1\Big)^{1/r} = \Big(\sum_{j=1}^k \big\|V(y_{n_j})\car_{E^c_{n_j}}\big\|^r_r\Big)^{1/r} \\
			& \stackrel{\eqref{eq:equation1_coBSAUC}}{\le} 4\theta^{-1/p}(1-\theta)^{1/r}k^{1/r}.
		\end{align*}
		On the other hand, to estimate the contribution on the $E_{n_j}$'s, we use again $p\le r$ and the fact that $\norm{V(y_{n_j})}_{p,\theta}=\norm{\car_{E_{n_j}}V(y_{n_j})}_p$ to get 
		\begin{align*}
			\Big\|\Big(\sum_{j=1}^k V(y_{n_j})^r\car_{E_{n_j}}\Big)^{1/r}\Big\|_p &\le \Big\|\Big(\sum_{j=1}^k V(y_{n_j})^p\car_{ E_{n_j}}\Big)^{1/p}\Big\|_p  = \Big(\sum_{j=1}^k\norm{V(y_{n_j})}_{p,\theta}^p\Big)^{1/p}.
		\end{align*}
		Therefore, we have shown that for any $N_k\le n_1<\dots<n_k$,
		\begin{equation*}
			\frac{ct}{3C}k^{1/p} \le \Big(\sum_{j=1}^k\norm{V(y_{n_j})}_{p,\theta}^p\Big)^{1/p} + 4\theta^{-1/p}(1-\theta)^{1/r}k^{1/r}.
		\end{equation*}
		Since this is true for any $N_k\le n_1<\dots<n_k$ and $(\norm{V(y_{n})}_{p,\theta})_n$ converges, we have
		\begin{equation*}
			\frac{ct}{3C}k^{1/p} \le k^{1/p}\lim_{n\to\infty}\norm{V(y_{n})}_{p,\theta} + 4\theta^{-1/p}(1-\theta)^{1/r}k^{1/r}.
		\end{equation*}
		Since $\lim_{k\to\infty}k^{1/r-1/p}=0$ whenever $p<r$, we could have chosen $k$ large enough at the beginning to ensure that $\displaystyle\lim_{n\to \infty}\norm{V(y_n)}_{p,\theta}\ge \frac{c}{4C}t$. It is worth pointing out that this is the only place in the argument where we are using that $p<r$. Finally, it remains to recall that the first subsequence, still denoted $(y_n)_n$,  was chosen so that  $\liminf_{n\to \infty}\norm{V(y_n)}_{p,\theta}=\lim_{n\to \infty}\norm{V(y_n)}_{p,\theta}$. This concludes the proof.

	\end{proof}

		We can now deduce the following important corollary.
		
		\begin{coro}
			\label{cor:random-Lp-norm-type+pcoBS->Cp-subspace} 
			Let $1<p<r\le 2$ and $Y$ be a separable Banach space with a random $L_p$-norm of type $r$.
			If $X$ is a quotient of $Y$ and has the weak $p$-co-Banach Saks property, then $X$ is isomorphic to a closed subspace of an $\ell_p$-sum of finite-dimensional spaces.
		\end{coro}
		
		\begin{proof} 
			By Theorem \ref{theo:coBS->AUC} we already know that $X$ is $p$-AUC. The assumption on $Y$ together with Proposition \ref{prop:random-Lp-type->pUS} ensures that $Y$ is $p$-uniformly smooth and so is every quotient of $Y$ and in particular $X$. Since $X$ is $p$-uniformly smooth and $p$-AUC, it is in particular reflexive, $p$-AUS and $p$-AUC. The conclusion follows from Theorem \ref{thm:subspaces-of-lp-sums} after observing that $X$ must also be separable since it is a quotient of a separable space.
		\end{proof}
		
		\begin{rema}
			Theorem \ref{theo:coBS->AUC} and Corollary \ref{cor:random-Lp-norm-type+pcoBS->Cp-subspace} are not true without the smoothness assumption. Indeed, let $1<p<\infty$ and $\Tsi_p$ be the $p$-convexification of Tsirelson space $\Tsi$, which will be discussed in more detail in Chapter \ref{chapter:Counterexamples}. The space $\Tsi_p$ has the weak $p$-co-Banach-Saks property since it only has $\ell_p$-spreading models, but it does not even have an equivalent norm that is $p$-AUC since $\Tsi_p$ is certainly $p$-AUS (by strict $p$-convexity), reflexive, separable, but $\Tsi_p$ does not contain any $\ell_q$ subspace for $1<q<\infty$. 
		\end{rema}
		
		\section{\texorpdfstring{Coarse-Lipschitz stability of subspaces and quotients of $\ell_p$}{Coarse-Lipschitz stability of subspaces and quotients of}}\label{section:CLstabilitySQl_p}
		
		The results from the previous section provide the missing ingredients to complete the classification of the subspaces of $\ell_p$, up to coarse-Lipschitz equivalences, for any value of $p\in(1,\infty)$. To obtain these rigidity results, we need to understand under what conditions a coarse-Lipschitz embedding of a Banach space into an $\ell_p$-sum of finite-dimensional spaces can be upgraded to an isomorphic embedding into a similar space. We will see two situations where this is possible.
		
		\begin{theo}
			\label{theo:CL-embed-into-quotient-of-lp-sums-US} 
			Let $p\in (1,2]$ and $(F_n)_{n=1}^\infty$ be a sequence of finite-dimensional spaces such that $(F_n)_{n=1}^\infty$ is equi-$r$-uniformly smooth for some $r\in(p,2]$. 
			If a Banach space $X$ coarse-Lipschitz embeds into a quotient of $(\sum_{n=1}^\infty F_n)_{\ell_p}$, then $X$ linearly embeds into an $\ell_p$-sum of finite-dimensional spaces. 
		\end{theo}
		
		\begin{proof} 
			Let $Y:=(\sum_{n=1}^\infty F_n)_{\ell_p}$ and $Z$ be a quotient of $Y$. First, observe that the norm of $Y$ is a random $L_p$ norm of type $r$ (Proposition \ref{prop:lp-sum-random-Lp-norm-type->random-Lp-norm-type}) and hence it is $p$-uniformly smooth by Proposition \ref{prop:random-Lp-type->pUS}. In particular, $Y$ is super-reflexive, $p$-AUS and also $p$-AUC. It is by now standard that so is $Z$. Since $X$ coarse-Lipschitz embeds into a $q$-AUC space, it follows from Theorem \ref{thm:CL-into-pAMUC->weak-pcoBS} that $X$ has the weak $p$-co-Banach-Saks property. To reach the conclusion, it would be sufficient, according to Corollary \ref{cor:random-Lp-norm-type+pcoBS->Cp-subspace}, to show that $X$ is isomorphic to a quotient of a separable space with a random $L_p$-norm of type $r$. This can be done with a little bit of Banach space gymnastics involving ultrapowers and differentiability arguments, which we detail.
			%Fixing a non-principal ultrafilter $\cal U$ on $\bN$, we have, by Proposition \ref{prop:random-norm-stability-subspaces-ultraproducts}, that  $Y^{\cU}$ is also super-reflexive with a random $L_p$ norm of type $r$.
			First, fix a nonprincipal ultrafilter $\cU$ on $\bN$. Then, observe that $X$ must be separable and since $Z$ is super-reflexive, Theorem \ref{thm:CL-stability-SR} tells us that $X$ is linearly isomorphic to a separable subspace of $Z^{\cU}$. Clearly, $Z^{\cU}$ is isometric to a quotient of $Y^{\cU}$. Thus, $X$ is linearly isomorphic to a separable subspace $Z_1$ of a quotient of $Y^{\cU}$ and therefore to a quotient of a subspace $Y_1$ of $Y^{\cU}$. It is important to note that, since $Z_1$ is separable, $Y_1$ can also be taken separable. Indeed, if $Q$ is a quotient map from $Y_1$ onto $Z_1$ and $\phi$ is a continuous Bartle-Graves selector for $Q$, then the closed linear span of $\phi(Z_1)$ is separable and its image by $Q$ is equal to $Z_1$. By Proposition \ref{prop:random-norm-stability-subspaces-ultraproducts}, we have that $X$ is isomorphic to a quotient of a separable space with a random $L_p$-norm of type $r$ and this concludes the proof.
			
			%Since $E$ is separable with a random $L_p$ norm of type $r$ and $X$ has the weak $p$-co-Banach-Saks property, we can apply Corollary \ref{cor:random-Lp-norm-type+pcoBS->Cp-subspace} to deduce that $X$ has property $(\tilde{m_p})$.
		\end{proof}
		
		Before we prove a version of Theorem \ref{theo:CL-embed-into-quotient-of-lp-sums-UC} with a dual assumption, we make a convenient observation. In their study of $M$-ideals of compact operators in \cite{KaltonWerner1995}, Kalton and Werner introduced the following definition.
		
		\begin{defi} 
			\label{defi:mp}
			A Banach space $X$ has \emph{property $(m_p)$} if for every $x\in X$ and every weakly null sequence $(x_n)_{n=1}^{\infty}$,
			\begin{equation*}
				\limsup_{n\to\infty}\norm{x+x_n}^p = \norm{x}^p + \limsup_{n\to\infty}\norm{x_n}^p.
			\end{equation*}
		\end{defi}
		
		Note that every $\ell_p$-sum of finite-dimensional spaces has property $(m_p)$ and it is immediate that if $X$ isomorphically embeds into such an $\ell_p$-sum, then $X$ has an equivalent norm with property $(m_p)$.  The following proposition follows from these observations, Theorem \ref{thm:subspaces-of-lp-sums} and the AUS-AUC duality in the reflexive setting.
		
		\begin{prop}
			\label{prop:tilde-mp-equivalences}
			Let $p\in (1,\infty)$, $q$ be its conjugate exponent and $X$ be a separable reflexive Banach space. The following assertions are equivalent
			\begin{enumerate}[(i)]
				\item $X$ is isomorphic to a subspace of an $\ell_p$-sum of finite-dimensional spaces.
				\item $X$ has an equivalent norm with property $(m_p)$.
				\item $X$ admits an equivalent $p$-AUC norm and an equivalent $p$-AUS norm.
				\item $X^*$ is isomorphic to a subspace of an $\ell_q$-sum of finite-dimensional spaces.
			\end{enumerate}   
		\end{prop}
		
		Following Kalton and Werner, it will sometimes be convenient to say that a separable Banach space $X$ has \emph{property $(\tilde m_p)$}, for $p\in (1,\infty)$, if it is isomorphic to a subspace of an $\ell_p$-sum of finite-dimensional spaces. In this terminology, the equivalence $(i)\iff (iv)$ above reads as follows: a separable reflexive Banach space $X$ has property $(\tilde{m}_p)$ if and only if $X^*$ has property $(\tilde{m}_q)$. 
		%closed subspace of a space $(\sum_{n=1}^{\infty}F_n)_{\ell_p}$ where $(F_n)_{n=1}^{\infty}$ is  a sequence of finite-dimensional spaces, or equivalently, if it is isomorphic to a subspace of $C_p$. 
		On the isometric side, we send the reader to Exercise \ref{exe:duality-mp-mq} for the duality regarding property $(m_p)$ in the reflexive setting and to \cite{KaltonWerner1995} for the duality in a more general setting.
		
		%\begin{prop}\label{duality_mp} Let $X$ be a reflexive separable Banach space, $p\in (1,\infty)$ and $q$ the conjugate exponent of $p$. Then, $X$ has $(m_p)$ if and only if $X^*$ has $(m_q)$.\end{prop} We refer the reader to \cite{KaltonWerner1995} for the proof of an even more general statement. 
		
		\smallskip 
		Recall that a family of Banach spaces $(X_i)_{i\in I}$ is said to be \emph{equi-$q$-uniformly convex} if there exists $C>0$ such that for all $i\in I$ and all $t\in [0,1]$, $\delta_{X_i}(t)\ge Ct^q$.
		
		\begin{theo}
			\label{theo:CL-embed-into-quotient-of-lp-sums-UC} 
			Let $p\in [2,\infty)$ and $(F_n)_{n=1}^\infty$ be a sequence of finite-dimensional spaces such that $(F_n)_{n=1}^\infty$ is equi-$r$-uniformly convex for some $r\in[2,p)$.
			If a Banach space $X$ coarse-Lipschitz embeds into a quotient of $(\sum_{n=1}^\infty F_n)_{\ell_p}$, then $X$ linearly embeds into an $\ell_p$-sum of finite-dimensional spaces. 
		\end{theo}
		
		\begin{proof}
			Let $Y:=(\sum_{n=1}^\infty F_n)_{\ell_p}$ and $Z$ be a quotient of $Y$. Since $Y$ is reflexive and $p$-AUS, so is $Z$. Thus, it follows from Theorem \ref{thm:CL-into-ref+pAUS->ref+weak-pBS} that $X$ has the weak $p$-Banach-Saks property and by duality (Proposition \ref{pro:wBS-duality}) that $X^*$ has the weak $q$-co-Banach-Saks property, where $q$ is the conjugate exponent of $p$. After another round of Banach space gymnastics, we can show that $X^*$ is isomorphic to a quotient with a random $L_q$-norm with type $s$, where $s$ is the conjugate exponent of $r$. Indeed, by standard duality arguments, the space $Y^*$ is an $\ell_q$-sum of equi-$s$-uniformly smooth spaces. It follows that the norm of $Y^*$ is a random $L_q$-norm of type $s$ and hence $Y^*$ is $q$-uniformly smooth. In particular, $Y^*$ is super-reflexive and so are $Y$ and $X$. Now, arguing as in Theorem \ref{theo:CL-embed-into-quotient-of-lp-sums-US}, $X$ is isomorphic to a separable subspace of a quotient of $Y^{\cU}$ for a nonprincipal ultrafilter $\cU$ on $\bN$ and hence $X^*$ is isomorphic to a quotient of a subspace $Y_1$ of $(Y^{\cU})^*$. Note that since $X^*$ is separable, $Y_1$ can be taken to be separable as well. Since $Y$ is super-reflexive, $(Y^{\cU})^*$ can be identified with $(Y^*)^{\cU}$ that has a random $L_q$-norm of type $s$ (Proposition \ref{prop:random-norm-stability-subspaces-ultraproducts}). We have proved that $X^*$ has the weak $q$-co-Banach-Saks property and is isomorphic to a quotient of a separable space with a random $L_q$-norm of type $s$ with $1<q<s\le 2$. Invoking Corollary \ref{cor:random-Lp-norm-type+pcoBS->Cp-subspace}, we deduce that $X^*$ has property $(\tilde{m}_q)$ and therefore that $X$ has property $(\tilde{m}_p)$ by Proposition \ref{prop:tilde-mp-equivalences}.
		\end{proof}

		To shorten some arguments in the proofs to come, let us adopt the following notation. For a given Banach space $X$, we denote by:
		\begin{itemize}
			\item $\cS(X)$ the class of all Banach spaces isomorphic to a subspace of $X$, 
			\item $\cQ(X)$ the class of all Banach spaces isomorphic to a quotient of $X$, 
			\item $\cS\cQ(X)$ the class of all Banach spaces isomorphic to a subspace of a quotient of $X$,
			\item $\cQ\cS(X)$ the class of all Banach spaces isomorphic to a quotient of a subspace of $X$.
		\end{itemize} 
		
		Note that $\cS\cQ(X)\subset \cQ\cS(X)$ since a subspace of a quotient is always a quotient of a subspace. Then, in the reflexive setting, by elementary duality, we have that a quotient of a subspace is always a subspace of a quotient and in the reflexive setting $\cS\cQ(X) = \cQ\cS(X)$. We are now ready for the rigidity results of this section.
		
		\begin{theo} 
			Let $p\in (1,\infty)$ and let $X$ be a Banach space.
			\begin{enumerate}[(i)]
				\item If $X$ coarse-Lipschitz embeds into $\ell_p$, then $X$ linearly embeds into $\ell_p$. %\in \cal S(\ell_p)$.
				\item If $X$ is coarse-Lipschitz equivalent to a quotient of $\ell_p$, then $X$ is linearly isomorphic to a quotient of $\ell_p$. %\in \cal Q(\ell_p)$.
			\end{enumerate}
		\end{theo}
		
		\begin{proof} 
			If $p=2$, it follows from Theorem \ref{thm:CL-stability-SR} that in all cases $X$ is isomorphic to a Hilbert space. So, let us assume that $p\neq 2$. We note that Theorem \ref{theo:CL-embed-into-quotient-of-lp-sums-US}, or Theorem \ref{theo:CL-embed-into-quotient-of-lp-sums-UC}, implies that in all cases $X$ has property $(\tilde{m}_p)$. In particular, $X$ is separable, reflexive, has the weak $p$-Banach-Saks and the weak $p$-co-Banach-Saks properties and does not have any subspace isomorphic to $\ell_s$, for any $s\neq p$. For the remainder of this proof, we fix a nonprincipal ultrafilter $\cal U$ on $\bN$ and denote by $q$ the conjugate exponent of $p$.
			
			$(i)$ Since $X$ is separable, it follows from Theorem \ref{thm:CL-stability-SR} that $X$ is linearly isomorphic to a subspace of $(\ell_p)^{\cU}$, which is an $L_p$-space and classical measure theory implies that $X$ is moreover isomorphic to a subspace of a separable $L_p$-space and hence to a subspace of $L_p := L_p[0,1]$.
			
			Since $X\in \cS(L_p)$ while $\ell_2 \notin \cS(X)$, the case $2<p<\infty$ is an immediate consequence of a classical result of Johnson and Odell (\cite[Theorem 1]{JohnsonOdell1974}).
			
			We now suppose that $1<p<2$. By duality we have that $X^* \in \cQ(L_q)$. Since $X$ has property $(\tilde{m}_p)$, $X^*$ has property $(\tilde{m}_q)$ and, in particular, $X^*$ has the weak $q$-Banach-Saks property. Since $q>2$, we can now invoke another classical result of W.B. Johnson (see Theorem III.2 and Remark III.1 in \cite{Johnson1977}) which insures that $X^* \in \cQ(\ell_q)$ and thus $X\in \cS(\ell_p)$.
			
			$(ii)$ Assume that $X$ is coarse-Lipschitz equivalent to $Y\in \cal Q(\ell_p)$. We first show that $X$ must be a quotient of $L_p$. Indeed, since $Y$ is separable and super-reflexive, by Theorem \ref{thm:CL-equiv-complemented-SR}, there exist $X_1$ and $Y_1$ such that $X_1$ is a separable complemented subspace of $X^{\cU}$ containing $X$, $Y_1$ is a separable complemented subspace of $Y^{\cU}$  and $X_1$ and $Y_1$ are linearly isomorphic. Observe now that $Y^{\cU}$ is a quotient of $(\ell_p)^{\cU}=L_p(\Omega)$ and since $Y_1$ is complemented in $Y^{\cU}$, $Y_1$ is also a quotient of $L_p(\Omega)$. Denote by  $Q$ a quotient map from $L_p(\Omega)$ onto $Y_1$ and $\phi$ a continuous Bartle-Graves selector of $Q$. Then, classical measure theory implies that $\phi(Y_1)$ is included in a separable subspace of $L_p(\Omega)$, which is itself an $L_p$-space. So, $Y_1$ and therefore $X_1$ belong to $\cQ(L_p)$. Finally, since $X$ is reflexive, it is complemented in  $X^{\cU}$ and thus in $X_1$. So, $X \in \cQ(L_p)$.
			
			If $2<p<\infty$, Johnson's theorem \cite{Johnson1977} combined with the fact that $X$ has the weak $p$-Banach-Saks property implies that $X\in \cQ(\ell_p)$.
			
			If $1<p<2$, we have that $X^* \in \cS(L_q)$ and $X^*$ has $(\tilde{m}_q)$, so $\ell_2 \notin \cS(X^*)$. We use again the result of Johnson and Odell \cite{JohnsonOdell1974} to deduce that $X^* \in \cS(\ell_q)$ and thus $X\in \cQ(\ell_p)$.
		\end{proof}
		
		In order to handle coarse-Lipschitz embeddings into quotients of $\ell_p$ we need to call upon a deep result by Kalton and Werner from \cite{KaltonWerner1995}. 
		
		\begin{theo}
			\label{theo:Kalton-Werner}
			Let $p\in (1,\infty)$ and assume that $X$ is a separable Banach space that is linearly isomorphic to a subspace of an $\ell_p$-sum of finite-dimensional spaces. Then,
			\begin{enumerate}[(i)]
				\item $X$ is linearly isomorphic to a quotient of a space $(\sum_{n=1}^{\infty}E_n)_{\ell_p}$ where $(E_n)_{n=1}^{\infty}$ is a sequence of finite-dimensional subspaces of $X$ \underline{itself}.
				\item $X$ is linearly isomorphic to a subspace of a space $(\sum_{n=1}^{\infty}E_n)_{\ell_p}$
				where $(E_n)_{n=1}^{\infty}$ is a sequence of finite-dimensional quotients of
				$X$ \underline{itself}.
			\end{enumerate}
		\end{theo}
		
		We will only need assertion $(i)$ of Theorem \ref{theo:Kalton-Werner} to prove the last rigidity result of this section.
		
		\begin{theo} 
			Let $p\in (1,\infty)$ and let $X$ be a Banach space. If $X$ coarse-Lipschitz embeds into a quotient of $\ell_p$, then $X$ is linearly isomorphic to a subspace of a quotient of $\ell_p$. %\in \cal S \cal Q(\ell_p)$.
		\end{theo}
		
		\begin{proof}
			Assume that $X$ coarse-Lipschitz embeds into $Y \in \cQ(\ell_p)$. The usual arguments yield that $X$ is linearly isomorphic to a separable subspace of $Y^{\cU}$, which is a quotient of $(\ell_p)^{\cU}$ and since $X$ is separable, $X$ is isomorphic to a subspace $Z$ of a quotient $E$ of $L_p$. At this crucial point, we need to observe that $X$ has property $\tilde{m_p}$ and so does $Z$. Applying Theorem \ref{theo:Kalton-Werner} to $Z$, we have that $Z$ is linearly isomorphic to a quotient of $(\sum_{n=1}^\infty E_n)_{\ell_p}$, where $(E_n)_n$ is a sequence of finite-dimensional subspaces of $Z$ itself. Let us denote by $Q$ the quotient map from $L_p$ onto $E$. We now use the fact that $L_p$ is a $\cal L_p$-space to deduce for all $n\ge 1$, the existence of a finite-dimensional subspace $H_n$ of $L_p$ which is $2$-isomorphic to some $\ell_p^{m_n}$ and such that $F_n:=Q(H_n) \supseteq E_n$. It clearly follows that $(\sum_{n=1}^\infty F_n)_{\ell_p}$ is isomorphic to a quotient of $(\sum_{n=1}^\infty H_n)_{\ell_p}$, which is $2$-isomorphic to $\ell_p$. So, $(\sum_{n=1}^\infty F_n)_{\ell_p} \in \cQ(\ell_p)$ and remembering that $E_n\subseteq F_n$ we have that $(\sum_{n=1}^\infty E_n)_{\ell_p} \in \cS\cQ(\ell_p)$. Since we are dealing with reflexive spaces, we conclude that $Z$ and therefore $X$ belong to $\cal Q \cal S \cal Q(\ell_p)=\cal S \cal Q(\ell_p)$.
		\end{proof}
		
		\section{\texorpdfstring{Coarse-Lipschitz rigidity of $(\sum_{n=1}^\infty \ell_r^n)_{\ell_p}$}{Coarse rigidity of}}
		
		The main goal of the last section of this chapter is to show, under some restrictions on $p$ and $r$, that $(\sum_{n=1}^\infty \ell_r^n)_{\ell_p}$ has a unique coarse-Lipschitz structure. We first start with a simple observation.
		
		\begin{lemm}
			\label{lem:lp(X)=X}
			If $Y := (\sum_{n=1}^\infty \ell_r^n)_{\ell_p}$ for some $1\le p,r\le \infty$, then $Y$ is linearly isomorphic to $\ell_p(Y)$. 
		\end{lemm}
		
		\begin{proof}
			Choose a sequence $(n_k)_{k=1}^{\infty}$ in $\bN$ such that for all $k$, $n_k\le k$ and $\{k \colon n_k=j\}$ is infinite for each $j \in \bN$. The later condition implies that $\ell_p(Y)$ is isometric to  $(\sum_{k=1}^{\infty}\ell_r^{n_k})_{\ell_p}$ and the former  that it is complemented in $Y := (\sum_{k=1}^\infty \ell_r^k)_{\ell_p}$.  Hence, for some Banach space $Z$, we have 
			\begin{equation*}
				Y\simeq \ell_p(Y)\oplus Z\simeq \ell_p(Y)\oplus\ell_p(Y)\oplus Z\simeq \ell_p(Y) \oplus Y \simeq \ell_p(Y).
			\end{equation*}   
		\end{proof}
		
		Now, assume that $X$ is coarse-Lipschitz equivalent to $Y:=(\sum_{n=1}^\infty \ell_r^n)_{\ell_p}$. If we could show that under this assumption $Y$ is complemented in $X$ and that $X$ is complemented in $Y$, then a Pe{\l}czy{\'n}ski decomposition argument with, the help of Lemma \ref{lem:lp(X)=X}, would allow us to conclude that $X$ and $Y$ are isomorphic. It is quite remarkable that there are situations where we will be able to do this. The results of the previous section will be of tremendous help, since for some values of $p$ and $r$, $X$ must already be isomorphic to a subspace of an $\ell_p$-sum of finite-dimensional spaces. Therefore, the crux of the problem boils down to a better understanding of the linear structure of subspaces of $C_p$. This purely linear task is the focal point of the next section, where we present those linear results due to Kalton in \cite{Kalton2013}.
		
		\subsection{\texorpdfstring{Linear structure of subspaces of $C_p$}{Linear structure of subspaces of}}
		\label{sec:Cp}
		%N. Kalton also needed to produce new linear results that we shall describe in the next two subsections. The special properties of $\ell_p$-sums of finite-dimensional spaces will be crucial for the main nonlinear classification results of this chapter. 
		Recall that, for $p\in [1,\infty)$, $C_p:=(\sum_{n=1}^\infty G_n)_p$, where $(G_n)_{n=1}^\infty$ is a sequence which is dense for the Banach-Mazur distance in the class of finite-dimensional normed spaces. It is classical that, up to isomorphism, $C_p$ does not depend on the choice of the dense sequence $(G_n)_{n=1}^{\infty}$ (see Exercise \ref{ex:C_p}). Note also that a Banach space is isomorphic to a subspace of an $\ell_p$-sum of finite-dimensional spaces if and only if it is isomorphic to a subspace of $C_p$.
		
		\begin{rema}
			$C_p$ is isomorphic to a space also called $C_p$ that was first introduced by Johnson in \cite{Johnson1971} and further studied by Johnson and Zippin in \cite{JohnsonZippin1972,JohnsonZippin1973,JohnsonZippin1974}.
		\end{rema}

		The first result about the linear structure of subspaces of $C_p$ gives a simple criterion that is useful to construct complemented subspaces of subspaces of $C_p$. Recall that $d_{BM}$ is the Banach-Mazur distance between two isomorphic Banach spaces. 
		
		\begin{theo}
			\label{theo:complemented-sub-Cp} 
			Let $p\in(1,\infty)$, $X$ be a subspace of $C_p$ and $(E_n)_{n=1}^{\infty}$ be a sequence of finite-dimensional subspaces of $X$. If there exists some constant $\lambda\ge 1$ such that for every $m$ and $n$ in $\bN$, $\ell_p^m(E_n)$ is $\lambda$-isomorphic to a $\lambda$-complemented subspace in $X$, then $(\sum_{n=1}^{\infty}E_n)_{\ell_p}$ is isomorphic to a complemented subspace of $X$.
		\end{theo}
		
		\begin{proof}
			%[Proof of Theorem \ref{theo:complemented-sub-Cp}] 
			To prove that $(\sum_{n=1}^{\infty}E_n)_{\ell_p}$ is isomorphic to a complemented subspace of $X$, it is clearly sufficient to factor the identity on $(\sum_{n=1}^{\infty}E_n)_{\ell_p}$ via $X$. Indeed, if there are two bounded operators $A\colon (\sum_{n=1}^{\infty}E_n)_{\ell_p}\to X$ and $B\colon X\to (\sum_{n=1}^{\infty}E_n)_{\ell_p}$ such that $I_{(\sum_{n=1}^{\infty}E_n)_{\ell_p}} = BA$, then $A$ is an isomorphism from $(\sum_{n=1}^{\infty} E_n)_{\ell_p}$ onto a subspace $Y$ of $X$ and that $AB$ is a bounded projection from $X$ onto $Y$.
			Reducing the problem further, it will be sufficient to find a constant $\gamma>0$ and two sequences of operators $(A_n\colon E_n\to X)_{n=1}^\infty$ and $(B_n\colon X\to E_n)_{n=1}^\infty$ so that:
			\begin{enumerate}[(1)]
				\item $I_{E_n} = B_nA_n$, for all $n\ge 1$,
				\item \begin{equation*}
					\Big\|\sum_{n=1}^{\infty}A_n u_n\Big\|_X\le \gamma\Big(\sum_{n=1}^{\infty}\norm{u_n}^p\Big)^{1/p},\  \textrm{ whenever }  u_n\in E_n \textrm{ for } n\in\bN,
				\end{equation*} 
				\item \begin{equation*}
					\Big\|\sum_{n=1}^{\infty}B_n^*u_n^*\Big\|_{X^*}\le \gamma\Big(\sum_{n=1}^{\infty}\norm{u_n^*}^q\Big)^{1/q},\ \textrm{whenever } u_n^*\in E_n^* \textrm{ for } n\in \bN \textrm{ and } \frac1q+\frac1p =1,
				\end{equation*} 
				\item $A_n(E_n)$ is included in the orthogonal of $\bigcup_{k<n}B_k^*(E_k^*)$, 
				\item $B_n^*(E_n^*)$ is included in the orthogonal of $\bigcup_{k<n}A_k(E_k)$.
			\end{enumerate}
			Here and later, since we are working in a reflexive setting we will identify the orthogonal and the pre-orthogonal of a set. 
			
			Coming back to our claim, should we have constructed such factorizations of the identity operators on the $E_n$'s, we may define $A\colon (\sum_{n=1}^{\infty}E_n)_{\ell_p}\to X$ and $B\colon X\to (\sum_{n=1}^{\infty}E_n)_{\ell_p}$ by $A((u_n)_{n=1}^{\infty}) := \sum_{n=1}^{\infty}A_nu_n$ and $Bx := (B_nx)_{n=1}^\infty$. Clearly, $(2)$ implies that $A$ is well-defined and bounded. To see that $B$ is well-defined and bounded, simply observe that $X$ is reflexive and $B$ is the adjoint of $C\colon (\sum_{n=1}^{\infty}E_n^*)_{\ell_q}\to X^*$ given by $C((u_n^*)_{n=1}^{\infty}) := \sum_{n=1}^{\infty} B_n^* u_n^*$, which is well-defined and bounded by $(3)$. Finally, $(1)$, $(4)$ and $(5)$ are crucial when showing that $I_{(\sum_{n=1}^{\infty}E_n)_{\ell_p}} = BA$ once we realize that, for instance, $(4)$ implies that $A_n(E_n)\subset \ker(B_k)$ for all $k<n$. These last two facts are not terribly difficult, but the tedious verifications are better left as an exercise. 
			
			The operators $A_n\colon E_n\to X$ and $B_n\colon X\to E_n$ as above will be constructed simultaneously in an inductive fashion. This inductive procedure is based on two claims. The first claim below, whose elementary proof is deferred to Exercise \ref{exe:homogeneous_mp}, follows from the fact that $X$ is separable and reflexive and so we may assume, after renorming if needed, that $X$ has property $(m_p)$.
			
			\begin{claim}
				\label{claim:homogeneous-mp}
				For every finite-dimensional subspace $E$ of $X$ and every $\eta>0$, there exists a finite-dimensional subspace $H$ of $X^*$ such that for all $x\in E$ and $y\in H^\perp$,
				\begin{equation*}
					(1-\eta)(\norm{x}^p+\norm{y}^p)\le \norm{x+y}^p \le (1+\eta)(\norm{x}^p+\norm{y}^p).
				\end{equation*}
			\end{claim} 
			
			The second claim is the following factorization result regarding the identity map on a finite-dimensional subspace $E_n$ of $X$ that satisfies the assumption of the theorem.
			\begin{claim}
				\label{claim:factorization}
				Given any finite-dimensional subspaces $G\subseteq X$, $H\subseteq X^*$ and $n\in\bN$, there exist operators $A\colon E_n\to X$ and $B\colon X\to E_n$ such that 
				\begin{enumerate}[(i)]
					\item $I_{E_n} = BA$ with $\max\{\norm{A},\norm{B}\}\le 2\lambda$,
					\item $A(E_n)\subseteq H^{\perp}$,
					\item $B^*(E_n^*)\subseteq G^{\perp}$.
				\end{enumerate}   
			\end{claim}
			
			Before proving Claim \ref{claim:factorization}, we explain how with the help of Claim \ref{claim:homogeneous-mp} and Claim \ref{claim:factorization} we can proceed with the inductive construction of the two sequences of operators $(A_n\colon E_n\to X)_{n=1}^\infty$ and $(B_n\colon X\to E_n)_{n=1}^\infty$ satisfying conditions $(1)-(5)$ above and so that $I_{E_n} = B_nA_n$ for all $n\in \bN$.
			
			Fix a sequence $(\eta_k)_{k=1}^\infty \subset (0,1)$ such that $\prod_{k=1}^\infty (1+\eta_k)\le 2$. We will build inductively $A_n\colon E_n\to X$ and $B_n\colon X\to E_n$ as above with the following more precise estimates:
			\begin{equation*}
				\Big\|\sum_{k=1}^{n} A_k u_k\Big\|\le 2\lambda \prod_{k=1}^n(1+\eta_k)\Big(\sum_{k=1}^{n}\norm{u_k}^p\Big)^{1/p},\qquad u_k\in E_k,\ k=1,\dots,n
			\end{equation*} 
			and
			\begin{equation*}
				\Big\|\sum_{k=1}^{n} B_k^* u_k^*\Big\|\le 2\lambda \prod_{k=1}^n(1+\eta_k)\Big(\sum_{k=1}^{n}\norm{u_k^*}^q\Big)^{1/q},\qquad u_k^*\in E_k^*,\ k=1,\dots,n.
			\end{equation*}
			By Claim \ref{claim:factorization}, we can find $A_1\colon E_1 \to X$ and $B_1\colon X \to E_1$ so that $\max\{\norm{A_1},\norm{B_1}\}\le 2\lambda$ and $I_{E_1} = A_1B_1$. Assume now that $A_1,\dots, A_n$ and $B_1,\dots, B_n$ have been constructed. Let $E$ be the linear span of $\bigcup_{k\le n}A_k(E_k)$. By Claim \ref{claim:homogeneous-mp}, there exists a finite-dimensional subspace $H$ of $X^*$ such that for all $x\in E$ and $y\in H^\perp$, 
			\begin{equation*}
				\norm{x+y} \le (1+\eta_{n+1})(\norm{x}^p+\norm{y}^p)^{1/p}.
			\end{equation*}
			We may assume that $H$ contains the linear span of $\bigcup_{k\le n}B_k^*(E_k^*)$, which we denote by $F$. Since $X$ is reflexive with property $(m_p)$, $X^*$ has property $(m_q)$ by duality and by Claim \ref{claim:homogeneous-mp}, there exists a finite-dimensional subspace $G$ of $X$, containing the linear span of $\bigcup_{k\le n}A_k(E_k)$ such that for all $x^*\in F$ and $y^*\in G^\perp$, 
			\begin{equation*}
				\norm{x^* + y^*} \le (1+\eta_{n+1})(\norm{x^*}^q+\norm{y^*}^q)^{1/q}.
			\end{equation*}
			We now apply Claim \ref{claim:factorization} to find $A_{n+1}\colon E_{n+1} \to X$ and $B_{n+1}\colon X\to E_{n+1}$ such that $I_{E_{n+1}} = A_{n+1}B_{n+1}$, $\max\{\norm{A_{n+1}},\norm{B_{n+1}}\}\le 2 \lambda$, $A_{n+1}(E_{n+1})\subset H^\perp$ and $B_{n+1}^*(E_{n+1}^*)\subset G^\perp$. It follows that whenever $u_k \in E_k$, for  $k=1,\dots,n+1$, we have
			\begin{align*}
				\Big\|\sum_{k=1}^{n+1} A_k u_k\Big\|  & \le (1+\eta_{n+1})\Big(\Big\|\sum_{k=1}^{n} A_k u_k \Big\|^p + \norm{A_{n+1}u_{n+1}}^p\Big)^{1/p}\\
				&\le (1+\eta_{n+1})\Big((2\lambda)^p\prod_{k=1}^n(1+\eta_k)^p \sum_{k=1}^n\norm{u_k}^p + (2\lambda)^p\norm{u_{n+1}}^p\Big)^{1/p}\\
				&\le 2\lambda \prod_{k=1}^{n+1} (1+\eta_k)\Big(\sum_{k=1}^n\norm{u_k}^p\Big)^{1/p}.
			\end{align*}
			The proof of the second estimate is identical and this concludes the inductive construction. 
			
			\smallskip 
			We now turn to the proof of Claim \ref{claim:factorization}. Fix first $m\in \bN$ to be chosen large enough later. It follows from the hypothesis that the identity on $\ell_p^m(E_n)$ factorizes through $X$. More precisely, there are operators $S\colon \ell_p^m(E_n)\to X$ and $T\colon X\to \ell_p^m(E_n)$ with $I_{\ell_p^m(E_n)} = TS$ and $\max\{\norm{S},\norm{T}\}\le \lambda$. If we write $S((u_j)_{j=1}^m) = \sum_{j=1}^m S_j u_j$ and $Tx = (T_jx)_{j=1}^m$, then it is easy to see that , for all $1\le j\le m$, $\max\{\norm{S_j},\norm{T_j}\}\le \lambda$ and $T_j S_j = I_{E_n}$. Then, since $G$ and $H$ are finite-dimensional we can take a projection $P_G\colon X\to G$ with $\norm{P_G}\le d_G$ and a projection $P_H\colon X^*\to H$ with $\norm{P_H}\le d_H$, where $d_G$ and $d_H$ are the dimensions of $G$ and $H$ respectively. Since $X$ is reflexive, we can consider $Q=P_H^*\colon X \to X$, where $P_H$ is viewed as a map from $X^*$ to $X^*$, and of course $Q^*=P_H$. It is then elementary to verify that $B_j := T_j(I_X -P_G)$ and $A_j:=(I_X-Q)S_j$ satisfy $B_j^*(E_n^*)\subseteq G^\perp$ and $A_j(E_n)\subseteq H^\perp$. However, there are two issues here. The first one is that $\norm{B_j}\le \lambda(1+d_G)$ and $\norm{A_j}\le \lambda(1+d_H)$, and the norms are too big whenever $G$ and $H$ are at least two-dimensional.
			For the second one, note that 
			\begin{equation}
				B_jA_j = I_{E_n} + R_j,
			\end{equation}
			where $R_j:=T_j(P_GQ - Q -P_G)S_j$ does not necessarily vanish. To resolve both issues, the key observation here is that $S^*\colon X^* \to \ell_q^m(E_n^*)$ can also be decomposed as $S^*x^* = (S_j^*x^*)_{j=1}^m$ and for any $\vep>0$ if $m$ was initially chosen large enough (depending only on $\lambda$, $G$, $H$ and $\eps$), then there exists $j_0\le m$ such that 
			\begin{equation}
				\label{eq:pigeonhole}
				\norm{T_{j_0}}_{G\to E_n} + \norm{S^*_{j_0}}_{H \to E_n^*} \le \vep.
			\end{equation}
			The inequality above follows from the finite-dimensionality of $G$ and $H$, the uniform bound on the norms of the operators, and a pigeonhole argument that can be applied in a quite general situation (see Exercise \ref{exe:pigeonhole}). Note also that by chasing down the identifications $(S_j)^*=(S^*)_j$ and this justifies using the ambiguous notation $S_j^*$. 
			
			Inequality \eqref{eq:pigeonhole} is the remedy to all our troubles since for this $j_0$, we have (remembering that $Q^*=P_H$ and $\norm{P_H}\le d_H$) 
			\begin{align*}
				\norm{B_{j_0}A_{j_0} - I_{E_n}} & = \norm{T_{j_0}P_G QS_{j_0} - T_{j_0}QS_{j_0} -T_{j_0}P_GS_{j_0}}\\
				& \le \norm{T_{j_0}P_G}(\norm{QS_{j_0}}+\norm{S_{j_0}})+ \norm{T_{j_0}}\norm{QS_{j_0}}\\
				& \le \norm{T_{j_0}}_{G\to E_n}\norm{P_G}(\norm{S^*_{j_0}P_H} + \lambda) + \lambda\norm{S^*_{j_0}P_H}\\
				& \le \vep d_G(\norm{S^*_{j_0}}_{H\to E_n^*}\norm{P_H} +\lambda) + \lambda \norm{S^*_{j_0}}_{H\to E_n^*}\norm{P_H}\\
				& \le \vep^2d_Gd_H + \vep d_G\lambda + \lambda\vep d_H\\
				& \le (\lambda\vep(d_H+d_G))^2.
			\end{align*}
			Therefore, $\norm{B_{j_0}A_{j_0} - I_{E_n}}\le \frac14$ whenever $\vep\le \frac{1}{2\lambda(d_G+d_H)}$. 
			This means that $B_{j_0}A_{j_0}$ is invertible whose inverse $D_{j_0}\colon E_n\to E_n$ satisfies $\norm{D_{j_0}}\le \frac43$.
			If we let $B := B_{j_0}=T_{j_0}(I_X-P_G)\colon X\to E_n$ and $A := A_{j_0}D_{j_0}=(I_X-Q)S_{j_0}D_{j_0} \colon E_n\to X$, then clearly $BA=I_{E_n}$ while we still have $A(E_n) \subset H^{\perp}$ and $B^*(E_n^*)\subset G^\perp$. Finally, we estimate the norms of $A$ and $B$. Note that 
			\begin{equation*}
				\norm{B}\le \norm{T_{j_0}} + \norm{T_{j_0}}_{G\to E_n}\norm{P_G}\le \lambda +  \vep d_G.
			\end{equation*}
			On the other hand,
			\begin{equation*}
				\norm{A} \le \frac{4}{3}(\norm{S_{j_0}} + \norm{QS_{j_0}})\le \frac43(\lambda + \norm{S_{j_0}^*P_H}) \le \frac43(\lambda + \norm{S_{j_0}^*}_{H\to E_n^*}\norm{P_H})\le \frac43(\lambda +\vep d_H).
			\end{equation*}
			It follows that $\max\{\norm{A},\norm{B}\}\le 2\lambda$ whenever $\vep \le \frac{1}{2\lambda(d_G+d_H)}$ and this finishes the proof of Claim \ref{claim:factorization} and hence the proof of the theorem. 
		\end{proof}
		
		\begin{rema} 
			The use of properties $(m_p)$ and $(m_q)$ in the above proof, instead of properties $(\tilde{m}_p)$ and $(\tilde{m}_q)$, is merely a convenience, in order to simplify the control of the constants.
		\end{rema}
		
		It is clear why Theorem \ref{theo:complemented-sub-Cp} should be useful to show that spaces are complemented in $\ell_p$-sums of finite-dimensional spaces. What we also need to understand is when subspaces of $C_p$ that are complemented in a Banach space $Y$ are isomorphic to $\ell_p$-sums of finite-dimensional spaces that are themselves complemented in $Y$. This problem can be partially understood assuming the subspace of $C_p$ has the approximation property. We will need the following result regarding special resolutions of the identity for subspaces of $C_p$ that have the approximation property. We refer to Appendix \ref{appendix:approximation} for the results pertaining to approximation properties of Banach spaces that are needed in the sequel. 
		
		\begin{prop}
			\label{prop:resolution-identity-Cp}
			Let $p\in (1,\infty)$ and $X$ be a subspace of $C_p$ with the approximation property. Then, there exist a constant $C\ge 1$ and a sequence of finite-rank operators $A_n\colon X\to X$ such that
			\begin{enumerate}[(1)]
				\item $x=\sum_{n=1}^{\infty}A_nx$, for all  $x\in X$,
				\item $A_jA_k=0$ whenever $\abs{k-j}>1$,
				\item for all $x \in X$, \begin{equation*}
					\frac1C \norm{x}\le \Big(\sum_{n=1}^{\infty}\norm{A_n x}^p\Big)^{1/p}\le C\norm{x}.
				\end{equation*}
			\end{enumerate}  
		\end{prop}
		
		\begin{proof} 
			%We may assume that $X$ is a subspace of a space $Z=(\sum_{n=1}^{\infty}G_n)_{\ell_p}$, where $(G_n)_n$ is a sequence dense in all finite-dimensional spaces in the sense of the Banach-Mazur distance.
			Since $X$ is a reflexive space with the approximation property, $X$ has the metric approximation property. Moreover, since $X$ is separable, $X$ has the commuting metric approximation property and there is a sequence of finite-rank operators $R_n\colon X\to X$  such that 
			\begin{equation}
				\label{eq:approximation}
				x=\lim_{n\to\infty} R_n x, \quad\textrm{ for all }x\in X,
			\end{equation}
			and
			\begin{equation}
				\label{eq:commutation}
				R_mR_n = R_nR_m=R_n,\quad\textrm{whenever }m>n.
			\end{equation}  
			%(see Proposition 2.1 in \cite{CasazzaKalton1990} for a general proof in the nonreflexive setting and for historical references).
			Thanks to this commutativity property, it is now easy to produce finite-rank operators satisfying conditions (1) and (2). Indeed, one could simply take $A_n=C_n-C_{n-1}$ where $(C_n)_{n=1}^\infty$ is a block sequence of $(R_n)_{n=1}^\infty$ whose blocks are convex combinations (and $C_0=0$). Eying up condition (3), we let $P_n \colon C_p\to (\sum_{k=1}^{n}G_k)_{\ell_p}\subseteq C_p$ be the partial sum operator associated with the FDD $(G_n)_n$ of $C_p$. We will denote by $S_n$ the restriction of $P_n$ to $X$ and we set $S_0=0$. We now consider $S_n$ and $R_n$ as operators in $K(X,C_p)$, the space of compact operators from $X$ to $C_p$. Here comes the crucial observation. Since for all $x\in X$, $x=\lim_{n\to \infty} S_n(x)$, we certainly have that for all $x\in X$ and $z^*\in C_p^*$, $\lim_{n\to\infty} z^*((S_n-R_n)(x))=0$. This means that $(R_n-S_n)_{n=1}^\infty$ converges to $0$ in the weak operator topology, and hence also in the dual weak operator topology since $X$ is reflexive, and ultimately in the weak topology of $K(X,C_p)$ (see Appendix \ref{appendix:FA} for more details and references).
			
			Therefore it follows from Mazur's theorem that, given $(\vep_k)_{k=1}^{\infty}$ with $\vep_k>0$ and such that $\sum_{k=1}^{\infty}\vep_k<\frac18$, we can find an increasing sequence of integers $(m_k)_{k=0}^{\infty}$ with $m_0=0$ (we also agree that $m_{-1}=0$ and $\eps_0=0$) and a nonnegative sequence $(a_j)_{j=1}^{\infty}$ such that for all $k\in \bN$,
			\begin{equation}
				\label{eq:weaklynulloperators}
				\sum_{j=m_{k-1}+1}^{m_k}a_j=1\ \ \text{and}\ \  \Big\|\sum_{j=m_{k-1}+1}^{m_k} a_j(R_j-S_j)\Big\|_{X\to C_p} \le \vep_k.
			\end{equation}
			For $k\ge 1$, we let $T_k := \sum_{j=m_{k-1}+1}^{m_k} a_jR_j$ and $A_k := T_k-T_{k-1}$ with the convention that $T_0=0$. As we have already observed, it follows from \eqref{eq:approximation} that $x=\sum_{n=1}^\infty A_n x$, for all $x\in X$, and from \eqref{eq:commutation} we deduce that $A_jA_k=0$ whenever $\abs{k-j}>1$.
			
			Now, if we let $U_k := \sum_{j=m_{k-1}+1}^{m_k} a_jS_j$ and $B_k := U_k-U_{k-1}$, for $k\ge 1$, and $U_0=0$, then it follows from \eqref{eq:weaklynulloperators} that for all $k\in \bN$,
			\begin{equation*}
				\norm{A_k - B_k}\le \vep_k+\vep_{k-1}.
			\end{equation*}
			Hence, a simple application of the triangle inequality in $\ell_p$ gives that for all $x\in X$, 
			\begin{align}
				\label{eq:AB_p_estimates}
				\Big|\Big(\sum_{k=1}^{\infty}\norm{A_kx}^p\Big)^{1/p}-\Big(\sum_{k=1}^{\infty}\norm{B_kx}^p\Big)^{1/p}\Big|
				&\le \Big(\sum_{k=1}^\infty \norm{(A_k-B_k)x}^p\Big)^{1/p}\le \frac14\norm{x}.   
			\end{align}
			
			To understand why condition (3) holds, it will be convenient to see $B_k$ as the restriction to $X$ of a multiplier with respect to the FDD $(G_n)_n$ of $C_p$ with associated coefficients that we denote by $(b^k_n)_n$. More precisely, $B_kx=\sum_{n=1}^\infty b^k_n x_n$, where for every $x\in X$, $x := (x_n)_{n=1}^\infty\in X$ with $x_n:=(P_n-P_{n-1})x\in G_n$, and hence
			\begin{equation*}
				\sum_{k=1}^{\infty}\norm{B_kx}^p = \sum_{k=1}^{\infty}\sum_{n=1}^{\infty}(b^k_n)^p\norm{x_n}^p.
			\end{equation*}
			Writing an exact formula for the coefficient $b^k_n$ in terms of the $a_j$'s is rather irrelevant here. However, what is important to notice is that a careful study of these coefficients reveals that for each $k \in \bN$, $(b_n^k)_n \subset [0,1]$ and is supported in $[m_{k-2}+2,m_k]$ and that $b^k_n+b^{k+1}_n=1$ for all $n \in (m_{k-1},m_k]$. Then, for a fixed $n\in \bN$ there exists a unique $k_n\in \bN$ so that $n\in (m_{k_n-1},m_{k_n}]$ and $b^{k_n}_n + b^{k_n+1}_n = 1$, while $b^j_n=0$ for all $j\notin \{k_n,k_n+1\}$. From this we deduce that
			%we set $c_n := b^k_n$, $d_n := b^{k+1}_n$. So, we can write, for $x := (x_n)_{n=1}^\infty\in X$,
			\begin{equation*}
				\sum_{k=1}^{\infty}\sum_{n=1}^{\infty}(b^k_n)^p\norm{x_n}^p = \sum_{n=1}^{\infty}\Big(\sum_{k=1}^{\infty}(b_{n}^{k})^p \Big)\norm{x_n}^p =\sum_{n=1}^{\infty}\Big((b_{n}^{k_n})^p + (b_n^{k_n+1})^p)\Big)\norm{x_n}^p.
			\end{equation*}
			It follows that for all $x \in X$,
			\begin{equation*}
				\frac12 \norm{x}\le \frac{1}{2^{1-1/p}}\norm{x} \le \Big(\sum_{k=1}^{\infty}\norm{B_kx}^p\Big)^{1/p} \le \norm{x}.
			\end{equation*}
			Thus, by \eqref{eq:AB_p_estimates}, we get
			\begin{equation*}
				\frac{1}{4}\norm{x} \le \Big(\sum_{k=1}^{\infty}\norm{A_kx}^p\Big)^{1/p} \le \frac{5}{4}\norm{x}.
			\end{equation*} 
		\end{proof}
		
		\begin{rema}
			It follows from the proof that we can take $C\le 4$ and thus the constant is independent of $p$ and $X$.
		\end{rema}
		
		Proposition \ref{prop:resolution-identity-Cp} is instrumental in showing that if a subspace of $C_p$ with the approximation property is complemented in a space paved by a sequence of finite-dimensional spaces $(E_n)_{n=1}^\infty$, then it is isomorphic to a complemented subspace of an $\ell_p$-sum of a subsequence of $(E_n)_{n=1}^\infty$. 
		
		\begin{theo}
			\label{theo:complemented-sub-Cp-pavings}
			Let $p\in (1,\infty)$ and $X$ be a subspace of $C_p$ with the approximation property. Assume that  $X$ is complemented in a Banach space $Z$ admitting a directed family $(E_i)_{i\in I}$ of finite-dimensional subspaces of $Z$ such that $\bigcup_{i\in I}E_i$ is dense in $Z$. Then, there exists a sequence $(i_n)_{n=1}^{\infty}$ in $I$ such that $X$ is isomorphic to a complemented subspace of $(\sum_{n=1}^{\infty}E_{i_n})_{\ell_p}$.
		\end{theo}
		
		\begin{proof} 
			Let $(A_n)_n$ be the finite-rank operators given by Proposition \ref{prop:resolution-identity-Cp}. We claim that our assumption implies that for each $n\in \bN$, there are $i_n\in I$ and a finite-dimensional subspace $H_n$ of $Z$ such that $E_{i_n}$ is $2$-isomorphic to $H_n$ and $A_n(X) \subset H_n$. This claim follows from an elementary approximation argument, which we simply sketch. Pick an Auerbach basis $(e_1,\dots,e_k)$ of $A_n(X)$ and fix $\eta >0$. Then, there exist $i_n\in I$ and $f_1,\dots,f_k\in S_{E_{i_n}}$  so that $\norm{f_i-e_i}<\eta$ for all $i\le k$. Note that we have used here all the assumptions on the family $(E_i)_{i\in I}$. Then, we can use that the linear span of $\{f_1,\dots,f_k\}$ is $k$-complemented in $E_{i_n}$ (with complement $G$) and show that, if $\eta$ is small enough, the $f_i$'s are linearly independent and the linear map $T\colon E_{i_n}\to Z$ defined by $T(f_i) := e_i$ and $T(g) := g$ for $g\in G$ is a 2-isomorphism from $E_{i_n}$ onto a space $H_n$ containing $A_n(X)$.
			
			Observe now that, since each $H_n$ is $2$-isomorphic to $E_{i_n}$, $(\sum_{n=1}^{\infty} H_n)_{\ell_p}$ is $2$-isomorphic to $(\sum_{n=1}^{\infty} E_{i_n})_{\ell_p}$, and $X$ will be isomorphic to a complemented subspace of $(\sum_{n=1}^{\infty} E_{i_n})_{\ell_p}$ as soon as we show that $X$ is isomorphic to a complemented subspace of $(\sum_{n=1}^{\infty} H_n)_{\ell_p}$. To prove the latter, we define $J\colon X\to (\sum_{n=1}^{\infty} H_n)_{\ell_p}$ by $Jx :=(A_n x)_{n=1}^{\infty}$. It follows from the discussion above and condition (1) in Proposition \ref{prop:resolution-identity-Cp} that $J$ is a well-defined isomorphic embedding. It remains to construct a bounded projection from $(\sum_{n=1}^{\infty} H_n)_{\ell_p}$ onto $J(X)$. It is immediate to see that $JQ$ would be such a projection if we can define $Q\colon (\sum_{n=1}^{\infty}H_n)_{\ell_p}\to X$ linear and bounded such that $QJ=I_X$. Since $X$ is assumed to be complemented in $Z$, let $P\colon Z\to X$ be a bounded projection and define $Q\colon (\sum_{n=1}^{\infty}H_n)_{\ell_p}\to X$ by
			\begin{equation*}
				Q((h_n)_{n=1}^{\infty}) := \sum_{n=1}^{\infty} A_nP(h_{n-1} + h_{n} + h_{n+1}),
			\end{equation*}
			with the convention $h_0=0$. Before we properly justify why $Q$ is a well-defined bounded operator, we explain why $QJ = Id_X$. Indeed, using conditions (1) and (2) in Proposition \ref{prop:resolution-identity-Cp}, we see that, for all $x\in X$,
			\begin{equation*}
				QJ(x) = \sum_{n=1}^{\infty}A_nP\Big(\sum_{k=n-1}^{n+1}A_kx\Big) = \sum_{n=1}^{\infty}A_n\Big(\sum_{k=n-1}^{n+1}A_kx\Big) = \sum_{n=1}^{\infty} A_n\Big(\sum_{k=1}^\infty A_k x\Big) = x.
			\end{equation*}
			The well-definiteness of $Q$ is mostly a computational matter. Consider first  $(x_n)_{n=1}^\infty$ a finitely nonzero sequence with $x_n\in A_n(X)$. Since condition (2) of  Proposition \ref{prop:resolution-identity-Cp} guarantees that $A_jx_n=0$, unless $j-1\le n \le j+1$, we have, using the convention $x_0=0$, that 
			
			\begin{align*} 
				\Big\| \sum_{k=1}^{\infty} x_{k} \Big\|^p & \le C^p \sum_{n=1}^\infty \Big\| A_n\Big(\sum_{k=1}^{\infty} x_{k}\Big)\Big\|^p \\
				& = C^p\sum_{n=1}^{\infty} \norm{A_{n} x_{n-1} + A_{n} x_{n} + A_{n} x_{n+1}}^p\\
				& \le 3^{p-1}C^p\sum_{n=1}^\infty\norm{A_{n}x_{n-1}}^p+\norm{A_{n}x_n}^p+\norm{A_{n}x_{n+1}}^p
			\end{align*}
			Since condition (3) in Proposition \ref{prop:resolution-identity-Cp} implies that $\sup_{n}\norm{A_n} \le C$ we have 
			\begin{equation}
				\label{eq:complemented-sub-Cp-pavings}
				\Big\| \sum_{k=1}^{\infty} x_{k} \Big\|^p \le 3^pC^{2p}\sum_{n=1}^{\infty}\norm{x_{n}}^p,
			\end{equation}
			and hence for all $(x_n)_{n=1}^\infty$ a finitely nonzero sequence with $x_n\in A_n(X)$,
			\begin{equation*}
				\Big\|\sum_{n=1}^{\infty} x_n\Big\| \le 3C^2 \Big(\sum_{n=1}^{\infty}\norm{x_n}^p\Big)^{1/p}. 
			\end{equation*}
			Consider now $h :=(h_n)_{n=1}^{\infty} \in (\sum_{n=1}^{\infty} H_n)_{\ell_p}$  finitely nonzero, then $x_n:= A_nP(h_{n-1}+h_n+h_{n+1})\in A_n(X)$ and only finitely many of them are nonzero. Therefore, it follows from the above computation that  
			\begin{align*}
				\norm{Qh} & \le 3C^2\Big(\sum_{n=1}^{\infty}\Big\|A_n\sum_{k=n-1}^{n+1} Ph_k\Big\|^p\Big)^{1/p}\\
				& \le 3C^3\Big(\sum_{n=1}^{\infty}\Big\|\sum_{k=n-1}^{n+1} Ph_k\Big\|^p\Big)^{1/p}\\
				& \le 3^{1-1/p}3C^3\Big(\sum_{n=1}^{\infty}\sum_{k=n-1}^{n+1}\Big\| Ph_k\Big\|^p\Big)^{1/p}\\
				& \le 3^2C^3\Big(\sum_{n=1}^{\infty}\norm{P h_n}^p\Big)^{1/p}\\
				&  \le 3^2C^3\norm{P}\,\norm{h}.
			\end{align*}
			Thus $Q$ extends to a bounded operator from $(\sum_{k=1}^{\infty} H_k)_{\ell_p}$ to $X$.
		\end{proof}
		
		\subsection{\texorpdfstring{The coarse and uniform structure of $(\sum_{n=1}^\infty \ell_r^n)_{\ell_p}$}{The coarse and uniform of}}
		
		The previous sections of this chapter provided all the tools needed to prove the uniqueness of the coarse-Lipschitz structure of $(\sum_{n=1}^\infty \ell_r^n)_{\ell_p}$ for some values of $p$ and $r$. 
		
		\begin{theo}\label{thm:rigidity-l_p_l_r^n} 
			Assume that $1<p<r\le 2$ or $2\le r<p<\infty$. If a Banach space $X$ is coarse-Lipschitz equivalent to $(\sum_{n=1}^\infty \ell_r^n)_{\ell_p}$, then $X$ is linearly isomorphic to $(\sum_{n=1}^\infty \ell_r^n)_{\ell_p}$. 
			
			In particular, the space $(\sum_{n=1}^\infty \ell_r^n)_{\ell_p}$ has a unique coarse or uniform structure for $p$ and $r$ as above. 
		\end{theo}
		
		\begin{proof} 
			Assume that $X$ is coarse-Lipschitz equivalent to $Y := (\sum_{n=1}^\infty \ell_r^n)_{\ell_p}$. Note that by Theorem \ref{theo:CL-embed-into-quotient-of-lp-sums-US} for $1<p<r\le 2$, or Theorem \ref{theo:CL-embed-into-quotient-of-lp-sums-UC} for $2\le r<p<\infty$, we know that $X$ must be isomorphic to a subspace of $C_p$.
			
			The proof that $(\sum_{n=1}^\infty \ell_r^n)_{\ell_p}$ is isomorphic to a complemented subspace of $X$ goes as follows. Fix a nonprincipal ultrafilter $\cU$ on $\bN$ and observe that $(\sum_{n=1}^\infty \ell_r^n)_{\ell_p}$ is super-reflexive and so is $X$, thanks to Theorem \ref{thm:CL-stability-SR}. Since $X$ is necessarily separable, by Theorem \ref{thm:CL-equiv-complemented-SR} we know that $Y=(\sum_{n=1}^\infty \ell_r^n)_{\ell_p}$ must be isomorphic to a complemented subspace of $X^{\cU}$. In particular, for every $m$ and $n$ in $\bN$, $\ell_p^m(\ell_r^n)$ is uniformly isomorphic to a uniformly complemented subspace of $X^\cU$. Therefore, it follows from the principle of local reflexivity for ultrapowers (see Appendix \ref{appendix:ultraproducts} and Appendix \ref{appendix:complemented-subspaces} for precise statements) that the $\ell_p^m(\ell_r^n)$'s are uniformly isomorphic to uniformly complemented subspaces of $X$. Since $X$ is a subspace of $C_p$, we can invoke Theorem \ref{theo:complemented-sub-Cp} to conclude that $(\sum_{n=1}^\infty \ell_r^n)_{\ell_p}$ is isomorphic to a complemented subspace of $X$.
			%By Theorem \ref{thm:CL-equiv-complemented-SR}, there exists a complemented separable subspace $X_1$ of $X^{\cU}$, containing $X$ that is linearly isomorphic to a complemented subspace $Y_1$ of $Y^{\cU}$. Since $Y$ is complemented in $Y^{\cU}$ and therefore in $Y_1$, we deduce that $Y$ is isomorphic to a complemented subspace of $X^{\cU}$. Therefore, it follows from a result due to Heinrich \cite{Heinrich1980} (see also Proposition F.6 (v) in \cite{BenyaminiLindenstrauss2000}) that $\ell_p^m(\ell_r^n)$, $m,n \in \bN$ are uniformly isomorphic to uniformly complemented subspaces $(F_{m,n})$ of $X$. Since $X$ is a subspace of $C_p$, it follows from Theorem \ref{theo:complemented-sub-Cp} that $(\sum_{n=1}^\infty \ell_r^n)_{\ell_p}$ is isomorphic to a complemented subspace of $X$.
			
			Proving that $X$ is isomorphic to a complemented subspace of $Y=(\sum_{n=1}^\infty \ell_r^n)_{\ell_p}$ is a bit more delicate. Since we want to use Theorem \ref{theo:complemented-sub-Cp-pavings}, we first need to show that $X$ has the approximation property. Here again, since $Y$ is a separable and super-reflexive Banach space, we have that $X$ is isomorphic to a complemented subspace of $Y^{\cU}$. Since the bounded approximation property passes to complemented subspaces and is preserved under isomorphisms, we can safely focus on $Y^\cU$. Moreover, since $(\sum_{n=1}^\infty \ell_r^n)_{\ell_p}$ is evidently complemented in $\ell_p(\ell_r)$, we have that $Y^\cU$ is complemented in $(\ell_p(\ell_r))^{\cU}$. The subtle point here is that while the bounded approximation property does not necessarily pass to ultrapowers, a strengthening of it, introduced by Pe{\l }czy{\'n}ski and Rosenthal in \cite{PelczynskiRosenthal1974} and called the uniform approximation property, does. This was shown by Heinrich in \cite{Heinrich1980} where he also showed that $\ell_p(\ell_r)$ has the uniform approximation property (see Appendix \ref{appendix:approximation} for more details). The outcome of the discussion above is that $X$ inherits the bounded approximation property from $(\ell_p(\ell_r))^{\cU}$. 
			%Proving that $X$ is isomorphic to a complemented subspace of $(\sum_{n=1}^\infty \ell_r^n)_{\ell_p}$ is a bit more delicate. Since we want to use Theorem \ref{theo:complemented-sub-Cp-pavings}, we first need to show that $X$ has the bounded approximation property. It follows from the work of Pe{\l }czy{\'n}ski and Rosenthal \cite{PelczynskiRosenthal1974} and Heinrich \cite{Heinrich1980} that $(\ell_p(\ell_r))^{\cU}$ has the bounded approximation property. Indeed, it follows from Theorem 9.4 in \cite{Heinrich1980} that $\ell_p(\ell_r)$ has the uniform approximation property (introduced in \cite{PelczynskiRosenthal1974}) and Theorem 9.1 in \cite{Heinrich1980} implies that $(\ell_p(\ell_r))^{\cU}$ has the uniform approximation property and therefore the bounded approximation property. Since $Y$ is complemented in $\ell_p(\ell_r)$, $Y^{\cU}$ is complemented in $(\ell_p(\ell_r))^{\cU}$ and consequently $Y^{\cU}$ has the bounded approximation property. Since $X$ is complemented in $X^{\cU}$ and therefore in $X_1$, we deduce that $X$ is isomorphic to a complemented subspace $Z$ of $Y^{\cU}$ and hence inherits the bounded approximation property from $Y^{\cU}$. 
			%We can also apply Theorem \ref{theo:CL-embed-into-quotient-of-lp-sums} to deduce that $X$ and therefore $Z$ have property $(\tilde{m_p})$. 
			Since $X$ is isomorphic to a subspace of $C_p$ and has the bounded approximation property, we can apply Theorem \ref{theo:complemented-sub-Cp-pavings} with $Z=Y^{\cU}$ in the notation of Theorem \ref{theo:complemented-sub-Cp-pavings}. It remains to exhibit the appropriate directed family of finite-dimensional subspaces of $Y^{\cU}$. By a result of M. Levy and Y. Raynaud \cite{LevyRaynaud1984}, $(\ell_p(\ell_r))^{\cU}$ is isometric to a band of $L_p(L_r):=L_p(\mu)(L_r(\nu))$ for appropriate measures $\mu$ and $\nu$. Any band in the reflexive Banach lattice $L_p(L_r)$ is $1$-complemented (see Section 1.a in \cite{LindenstraussTzafriri1979} for the definitions and results). Thus $(\ell_p(\ell_r))^{\cU}$ is isometric to $1$-complemented subspace of $L_p(L_r)$. Then, it follows from a routine approximation argument using simple functions that, for any finite-dimensional subspace $E$ of $L_p(L_r)$ there exists $n\in \bN$ and a finite-dimensional subspace $F$ of $L_p(L_r)$ such that $E \subset F$ and $F$ is $2$-isomorphic to $\ell_p^n(\ell_r^n)$ which is $1$-complemented in $E_n:=(\sum_{k=1}^{n}\ell_r^k)_{\ell_p}$. As we already mentioned, $Y^{\cU}$ is $1$-complemented in $(\ell_p(\ell_r))^{\cU}$ and thus in $L_p(L_r)$. We deduce that for any finite-dimensional subspace $E$ of $Y^{\cU}$, there exists $n\in \bN$ and a finite-dimensional subspace $F$ of $Y^{\cU}$ such that $E \subset F$ and $F$ is $2$-isomorphic to a $2$-complemented subspace of $E_n$. In particular, the family of all subspaces of $Y^{\cU}$ that are $2$-isomorphic to a $2$-complemented subspace of some $E_n$ is a directed family, whose union is obviously dense in $Y^{\cU}$. Then, it follows from Theorem \ref{theo:complemented-sub-Cp-pavings} that there exists a sequence $(m_n)_{n=1}^\infty$ in $\bN$ such that $X$ is isomorphic to a complemented subspace of $(\sum_{n=1}^\infty E_{m_n})_{\ell_p}$. The conclusion follows from the fact that $(\sum_{n=1}^\infty E_{m_n})_{\ell_p}$ is $1$-complemented in $\ell_p(Y)\simeq Y$. 
			
			Finally, since $X$ is isomorphic to a complemented subspace of $Y$, $Y$ is isomorphic to a complemented subspace of $X$ and  $Y\simeq \ell_p(Y)$, we conclude that $X$ is isomorphic to $Y$ by playing a bit of Pe{\l }czy{\'n}ski's accordion.

		\end{proof}
		
		\begin{rema} 
			When $r=2$, the above results reduce to the uniqueness of the coarse structure of $\ell_p$. Indeed, it is a classical result of  Pe{\l}czy{\'n}ski \cite{Pelczynski1960} that $(\sum_{n=1}^\infty \ell_2^n)_{\ell_p}$ is linearly isomorphic to $\ell_p$ (see Exercise \ref{exe:accordion}). 
		\end{rema}
		
		\section{Notes} 
		
		The paper \cite{Kalton2013b} by N. Kalton contains other mathematical gems. Let us mention a few of them. First, the most general version of Theorem \ref{thm:rigidity-l_p_l_r^n} is the following (we refer again to Appendix \ref{appendix:approximation} for the definition of the uniform approximation property).
		
		\begin{theo}
			\label{theo:Kalton-general}\,
			\begin{enumerate}[(i)]
				\item Suppose $1<p<r\le 2$ and that $(E_n)_{n=1}^{\infty}$ is an increasing sequence of uniformly complemented finite-dimensional subspaces of an $r$-uniformly smooth  Banach space with the uniform approximation property. Then, $(\sum_{n=1}^{\infty}E_n)_{\ell_p}$ has a unique coarse (or uniform) structure.
				\item Suppose $2\le r<p<\infty$ and that $(E_n)_{n=1}^{\infty}$ is an increasing sequence of uniformly complemented finite-dimensional subspaces of an $r$-uniformly convex Banach space with the uniform approximation property.
				Then, $(\sum_{n=1}^{\infty}E_n)_{\ell_p}$ has a unique coarse (or uniform) structure.
			\end{enumerate}
		\end{theo}
		The details of the proof of Theorem \ref{theo:Kalton-general} can be found in the master's thesis of Audrey Fovelle \cite{Fovelle-master}.
		
		Recall that a Banach space $X$ is said to have the \emph{strong Schur property} if there is a constant
		$c>0$ so that if $(x_n)_{n=1}^{\infty}$ is a normalized sequence in $X$ with
		$\sep\{x_n\}_{n=1}^{\infty}:=\delta>0$, then there is a subsequence such that for all $k\ge 1$ and $(a_n)_n \subseteq \bR$ 
		\begin{equation*}
			\Big\|\sum_{j=1}^ka_jx_{n_j}\Big\|\ge c\delta\sum_{j=1}^k\abs{a_j}.
		\end{equation*} 
		This quantitative refinement of the classical Schur property was first considered (implicitly) by Johnson and Odell \cite{JohnsonOdell1974}, where it was observed that there are subspaces of $L_1$ with the Schur property but failing the strong Schur property. The strong Schur property seems to have been first explicitly considered in unpublished lecture notes by Rosenthal in \cite{Rosenthal1979}. The first explicit and published mention of the strong Schur property seems to appear in a paper of Bourgain and Rosenthal \cite{BourgainRosenthal1980} where they provided an example of an uncomplemented subspace of $L_1$ with the strong Schur property but that does not embed isomorphically into $\ell_1$. Bourgain and Rosenthal then asked in \cite[page 68]{BourgainRosenthal1980} whether a complemented subspace of $L_1$ with the strong Schur property must be isomorphic to $\ell_1$. This fundamental and long-standing open problem about the structure of Banach spaces is still open. An alternative formulation of the strong Schur property was later given by Kalton in \cite{Kalton2001}.
		
		In \cite[Theorem 6.4]{Kalton2013b}, Kalton provides two new characterizations of subspaces of $L_1$ with the strong Schur property. One in terms of the modulus of asymptotic midpoint uniform convexity and the other in terms of the anti-Banach-Saks property introduced by Kalton in \cite{Kalton2013b}. Following Kalton, let us say that a Banach space $X$ has the \emph{anti-Banach-Saks property} if there is a constant $c>0$ so that for every spreading model $(e_j)_{j=1}^\infty$ of a normalized sequence in $X$, we have for all $k\in \bN$,
		\begin{equation*}
			\Big\|\sum_{j=1}^ke_j\Big\|_S\ge ck.
		\end{equation*}
		It is worth pointing out that the strong Schur property implies the anti-Banach-Saks property. This nontrivial fact was proved in \cite{Kalton2013b}. 
		\begin{theo} 
			\label{theo:strong-Schur}
			Let $X$ be a closed subspace of $L_1$. The following conditions on $X$ are equivalent:
			\begin{enumerate}[(i)]
				\item $X$ has the anti-Banach-Saks property.
				\item $X$ has the strong Schur property.
				\item For some $c>0$ we have $\hat\delta_X(t)\ge ct$, for all $t>0$. 
			\end{enumerate}
		\end{theo}
		
		Then, Kalton deduced from Theorem \ref{thm:CL-AMUC} the following striking result.
		
		\begin{theo}
			The class of Banach spaces that isomorphically embed into $L_1$ and have the strong Schur property is stable under coarse-Lipschitz embeddings.    
		\end{theo}
		
		It follows from the theorem above that if a Banach space $X$ coarse-Lipschitz embeds into $\ell_1$, then $X$ is linearly isomorphic to a subspace of $L_1$, and the above theorem implies that $X$ has the strong Schur property, but we do not know if $X$ linearly embeds into $\ell_1$. Let us record this important question.
		
		\begin{prob} 
			If a Banach space coarse-Lipschitz embeds into $\ell_1$, does it isomorphically embed into $\ell_1$?
		\end{prob}
		
		Note that an important theorem of Maurey, Johnson, and Schechtman \cite{MaureyJohnsonSchechtman2009} states that the class of $\mathcal L_1$-spaces is stable under uniform homeomorphisms. Subsequently, Kalton \cite{Kalton2012} proved that this theorem is optimal by constructing two separable $\mathcal L_1$ spaces that are uniformly homeomorphic but not isomorphic. Unlike in the reflexive range, a $\mathcal L_1$-space is not necessarily isomorphic to a complemented subspace of an $L_1$-space.
		
		We refer the reader to \cite{KaltonWerner1995} and \cite {KaltonWerner1995b} for the original applications of the works of N. Kalton and D. Werner on $M$-ideals. Another application, in the spirit of this book, can be found in a paper by Dutta and Godard \cite{DuttaGodard}. They study, for a Banach space with property $(M)$ or $(M^*)$ in the sense of Kalton \cite{Kalton1993}, the links between its Szlenk index and its modulus of asymptotic uniform smoothness. They prove that a separable reflexive space with property $(M)$ embeds isomorphically into an $\ell_p$-sum of finite-dimensional spaces if and only if its Szlenk index has power type $q$, where $q$ is conjugate to $p$. Without property $(M)$, we know from Chapter \ref{chapter:Szlenk} that this is no longer true in general. They also show that a separable Banach space in $\sA_\infty$ with property $(M^*)$ linearly embeds into $c_0$. They also derive applications to the nonlinear classification of Orlicz sequence spaces. 
		
		\section{Exercises}
		
		\begin{exer}
			\label{exe:deltahat}
			Let $(X,\norm{\cdot})$ be a Banach space.
			\begin{enumerate}
				\item Show that $t \mapsto \frac{\hat{\delta}_X(t)}{t}$ is increasing. 
				\item Prove that for all $t\ge 0$, $\hat{\delta}_X\big(\frac{t}{2}\big) \le \hat{\delta}^c_X(t) \le \hat{\delta}_X(t)$. 
				\item Prove that $\hat{\delta}^c_X$ is $1$-Lipschitz.
				\item Prove that $\lim_{t\to \infty}\frac{\hat{\delta}^c_X(t)}{t}=1$.
			\end{enumerate}
		\end{exer}

		\begin{exer}
			\label{exe:ff-sm}
			Let $S$ be a spreading model of a Banach space $X$ that is generated by a weakly null sequence. Show that there is $C>0$ such that for all $k\in \bN$ and $a_i,\dots,a_k\in \bR$, 
			\begin{equation*}
				\frac{1}{C}\Big\|\sum_{i=1}^k a_i e_i\Big\|_{\ell_{\bar{\delta}^c_X}} \le \Big\|\sum_{i=1}^k a_i e_i\Big\|_S \le C\Big\|\sum_{i=1}^k a_i e_i\Big\|_{\ell_{\bar{\rho}_X}}, 
			\end{equation*}
		\end{exer}
		\begin{proof}[Hint]
			Consider the iterated norms $\Lambda_N$ for $N(1,t) := 1+\bar{\rho}_X(t)$ or $N(1,t) := 1+\bar{\delta}^c_X(t)$.
		\end{proof}
		
		\begin{exer}
			\label{exe:lower-AMUC-sm}
			%Prove Proposition \ref{prop:lower-AMUC-sm}.
			Show that for any spreading model $S$ generated by a normalized sequence in a Banach space $X$, there is a constant $c>0$ such that for all $k\ge 1$ and $a_1,\dots,a_k\in \bR$, 
			\begin{equation*}
				c\Big\|\sum_{i=1}^k a_i e_i\Big\|_{\ell_{\hat{\delta}^c_X}} \le \bE_\eps \Big\|\sum_{i=1}^k \vep_i a_i e_i\Big\|_S,
			\end{equation*}
			where $\bE_\eps$ refers to the average over all choices of signs $\vep:=(\vep_1,\dots,\vep_k)\in\{-1,1\}^k$.
		\end{exer}
		
		\begin{exer}
			\label{exe:UC-AUC}
			Let $(X,\norm{\cdot})$ be a Banach space and denote by $\rho_X$ its modulus of uniform smoothness and by $\delta_X$ its modulus of uniform convexity (see Appendix \ref{appendix:us-uc} for the definitions).
			\begin{enumerate}
				\item Show that for all $t\in(0,1)$, $\bar{\rho}_X(t)\le 2\rho_X(t)$.
				\item Show that for all $t\in(0,1)$, $\delta_X(t)\le \bar{\delta}_X(t)$.
			\end{enumerate}
		\end{exer}
		
		\begin{exer}
			\label{exo:AUS-BS} 
			Let $X$ be a Banach space with separable dual and $1<p<\infty$.
			\begin{enumerate}
				\item Show that if $X$ is $p$-asymptotically uniformly smooth, then $X$ has the weak $p$-Banach-Saks property. 
				\item Show that if $X$ is $p$-asymptotically uniformly convex, then $X$ has the weak $p$-co-Banach-Saks property. 
			\end{enumerate}
		\end{exer}
		
		\begin{exer}
			\label{exe:duality-mp-mq}
			Let $1<p<\infty$, $q$ be the conjugate exponent of $p$ and $X$ be a reflexive separable Banach space. Show that $X$ has property $(m_p)$ if and only if its dual $X^*$ has property $(m_q)$.    
		\end{exer}
		
		\begin{exer}
			\label{ex:C_p}
			\,
			\begin{enumerate}
				\item Show that there exists a sequence $(G_n)_{n=1}^\infty$ which is dense in the class of finite-dimensional normed spaces for the Banach-Mazur distance.
				\item Let $p\in [1,\infty)$. Show that $C_p$ does not depend, up to isomorphism, on the choice of the sequence $(G_n)_{n=1}^\infty$.
			\end{enumerate}
		\end{exer}

		\begin{exer}
			\label{exe:homogeneous_mp}
			%Prove Lemma \ref{lem:homogeneous_mp}.
			Let $p\in (1,\infty)$ and $X$ be a separable reflexive Banach space with property $(m_p)$. Show that for every finite-dimensional subspace $E$ of $X$ and every $\eta>0$, there exists a finite-dimensional subspace $H$ of $X^*$ such that for all $x\in E$ and $y\in H^\perp$,
			\begin{equation*}
				(1-\eta)(\norm{x}^p+\norm{y}^p)\le \norm{x+y}^p \le (1+\eta)(\norm{x}^p+\norm{y}^p).
			\end{equation*}
		\end{exer}
		\begin{proof}[Hint]
			Use the compactness of $B_E$, an approximation argument and homogeneity. 
		\end{proof}

		\begin{exer}
			\label{exe:accordion}
			Let $p\in[1,\infty]$. Show that $(\sum_{n=1}^\infty \ell_2^n)_{\ell_p}$ is linearly isomorphic to $\ell_p$.
		\end{exer}
		
		\begin{exer}
			\label{exe:pigeonhole}
			Let $\lambda\ge 1$. Assume that $E$ and $F$ are finite-dimensional Banach spaces, $X$ and $Y$ are Banach spaces, and that for all $m\ge 1$, $T^{(m)}\colon E\to \ell_p^m(X)$ and $R^{(m)}\colon F\to \ell_q^m(Y)$ are bounded linear operators such that 
			\begin{equation*}
				\sup_{m\ge 1}\max\Big\{\|T^{(m)}\|,\|R^{(m)}\|\Big\}\le \lambda.
			\end{equation*} 
			If we write $T^{(m)}=(T_1^{(m)},\dots,T_m^{(m)})$ and $R^{(m)}=(R_1^{(m)},\dots,R_m^{(m)})$, then for every $\vep>0$ there is $m_0:=m_0(\lambda,F,G,\vep)\ge 1$ and $j_0\le m_0$ such that 
			\begin{equation*}
				\norm{T^{(m_0)}_{j_0}}_{E\to X} + \norm{R^{(m_0)}_{j_0}}_{F\to Y} \le \vep.
			\end{equation*}
		\end{exer}

		%%%%%%%%%%%%%%%%%%%%%%%%%%%%%%%%%%%%%%%%%%%%%%%%%%%%%%%%%%%%%%%%%%%%%%%%%%%%%%

		\chapter{Counterexamples to some rigidity problems}
		\label{chapter:Counterexamples}
		
		In this chapter, we present the main counterexamples in the theory of the nonlinear classification of Banach spaces. We describe two methods. The first one that we call the lifting method, is based on the existence of ``smooth" (Lipschitz or uniformly continuous) selections for certain quotient maps. It allowed Aharoni and Lindenstrauss to produce the first example of a pair of Banach spaces that are Lipschitz equivalent, but not linearly isomorphic (see Section \ref{sec:Aharoni-Lindenstrauss-example}). These spaces are \emph{nonseparable} and the existence of such a counterexample in the separable setting is a famous open question that we already mentioned in Problem \ref{prob:Lipschitz-rigidity}. We show in Section \ref{sec:example-UH-not-Isomorphic} how this lifting method can be combined with the results on Lipschitz-free spaces from Chapter \ref{chapter:Lipschitz-free} to produce an example of a pair of separable Banach spaces that are uniformly homeomorphic but not linearly isomorphic. Historically, the first example of such a situation was built by Ribe with a completely different method. This approach, the second method alluded to above, has been analyzed, formalized and sharpened by Kalton in the articles \cite{Kalton2012} and \cite{Kalton2013}. We detail Kalton's presentation in Section \ref{sec:Ribe-Kalton-machinery} and deduce Ribe's example from it. Then, we apply this approach to show that, for $p\in (1,\infty)$, the class of separable Banach spaces with an equivalent $p$-asymptotically uniformly smooth norm is not stable under uniform homeomorphisms. Finally, we build the stunning example, due to Kalton, of a separable pair of coarsely equivalent but not uniformly equivalent Banach spaces. On our way we corrected a small mistake in the original paper \cite{Kalton2012} (see Remark \ref{rema:mistake}).

		\section{The lifting method}
		\label{sec:lifting-method}
		
		In this very short section, we describe a very simple method for constructing nonlinear equivalences between Banach spaces. Let $Y$ be a closed subspace of a Banach space $X$ and  $Q\colon X\to X/Y$ be the quotient map. A map $\psi \colon X/Y \to X$ satisfying $Q\circ \psi = Id_{X/Y}$ is called a \emph{lifting} (or section) of $Q$. It is well known that if $Q$ has a linear and continuous section, then $X$ is linearly isomorphic to $Y\oplus X/Y$. 
		Dealing with not necessarily linear sections, one only needs to change the conclusion accordingly.
		
		\begin{prop}
			\label{prop:section-method}
			Let $X$ be a Banach space and $Y$ be a closed subspace of $X$. If the quotient map $Q\colon X\to X/Y$ admits a lifting that is either:
			\begin{enumerate}[(i)]
				\item linear and continuous,
				\item Lipschitz,
				\item uniformly continuous,
				\item coarse-Lipschitz,
				\item continuous,
			\end{enumerate}
			then the Banach spaces $X$ and $Y \oplus X/Y$ will be, respectively,
			\begin{enumerate}[(i')]
				\item linearly isomorphic,
				\item Lipschitz equivalent,
				\item uniformly homeomorphic,
				\item coarse-Lipschitz equivalent,
				\item homeomorphic.
			\end{enumerate}
		\end{prop}
		
		\begin{proof} 
			We only treat the Lipschitz case and hence the implication $(ii)\implies (ii')$, as the other cases can be treated in an identical manner. Let $\psi \colon X/Y \to X$ be a Lipschitz map satisfying $Q\circ \psi = Id_{X/Y}$ and consider the map $U\colon X\to Y \oplus X/Y$ defined by $U(x) := (x-\psi \circ Q(x),Qx)$. It is immediate that $U$ is Lipschitz with values in $Y \oplus X/Y$. For $(y,z)\in Y \oplus X/Y$, define $T(y,z) := y+\psi(z)$. The map $T\colon Y \oplus X/Y \to X$ is clearly Lipschitz and it is easily checked that $T\circ U = Id_X$ and $U\circ T = Id_{Y \oplus X/Y}$.
		\end{proof}
		
		The Bartle-Graves Theorem ensures the existence of a continuous lifting $\psi$ of the quotient map $Q$. Therefore, Proposition \ref{prop:section-method} tells us that $X$ and $Y \oplus X/Y$ are always homeomorphic. This observation is a very particular case of a deep result of Kadec \cite{Kadec1966} (in the separable case) and Torunczyck \cite{Torunczyk1981} (in the general case), asserting that two infinite-dimensional Banach spaces with the same density character are homeomorphic. In fact, the continuous content of the above proposition is one of the arguments in the proofs of these difficult results. 
		
		Proposition \ref{prop:section-method} naturally guides us towards a clearly identified goal: to understand when a quotient map has a rather smooth nonlinear lifting. One such situation occurs in the context of Lipschitz-free spaces. Indeed, if $X$ is a Banach space, then the barycentric map $\beta_X \colon \cF(X) \to X$ is a quotient map and $\delta_X$ is an isometric lifting of $\beta_X$. Thus, it follows from the previous proposition that $\cF(X)$ is Lipschitz isomorphic to $X\oplus \mathrm{Ker}(\beta_X)$ (this was already mentioned in Remark \ref{rem:section-free}). Note that if $X$ is a separable Banach space, then the Godefroy-Kalton lifting theorem (Theorem \ref{thm:GK-lifting-free}) implies that  $\cF(X)$ is actually linearly isomorphic to $X\oplus \mathrm{Ker}(\beta_X)$. This means that in the context of Lipschitz-free space over separable Banach spaces, the lifting method cannot produce counterexamples to the Lipschitz  rigidity problem. However, as we will see in the next section, this is not the case in the nonseparable setting.

		\subsection{Aharoni-Lindenstrauss pair of Lipschitz equivalent but nonisomorphic nonseparable Banach spaces}
		\label{sec:Aharoni-Lindenstrauss-example}
		
		The goal of this section is to describe the first example, due to Aharoni and Lindenstrauss \cite{AharoniLindenstrauss1978}, of a pair of Banach spaces that are Lipschitz isomorphic but not linearly isomorphic. To motivate the Aharoni-Lindenstrauss construction, observe that it is an elementary exercise to show that if $T\colon X/Y \to X$ is a linear section of the quotient map $Q\colon X \to X/Y$, then $P:=I-TQ$ is a bounded projection from $X$ onto $Y$. Conversely, if $P\colon X \to Y$ is a bounded projection, then $T\colon X/Y \to X$ given by $T(Qx):= (I-P)x$ is well-defined and is a linear lifting of $Q$. In the nonlinear world, things are different. It is still true that if $\psi\colon X/Y \to X$ is a Lipschitz lifting of the quotient map $Q\colon X \to X/Y$, then $r:=I-\psi\circ Q$ is a Lipschitz retract from $X$ onto $Y$. However, the existence of a Lipschitz retract $r\colon X \to Y$ does not necessarily yield a Lipschitz lifting of the quotient map. If one wants to produce a pair of Lipschitz equivalent but not linearly isomorphic Banach spaces using the lifting method, one needs to look for Banach spaces admitting a Lipschitz retract onto one of their noncomplemented subspaces. Therefore, one is naturally led to consider the space $\ell_\infty$ together with its noncomplemented subspace $\co$ since $\co$ is a $2$-Lipschitz retract of $\ell_\infty$, as shown in Proposition \ref{prop:c0-abs-2-Lip-retract}. Unfortunately, trying to show that $\ell_\infty/\co\oplus \co$ is Lipschitz equivalent to $\ell_\infty$ is doomed since $\ell_\infty/\co\oplus \co$ does not even coarsely embed into $\ell_\infty$ (see Corollary \ref{cor:ellinftyquotient-not-coarse-into-ellinfty})! To produce their counterexample, Aharoni and Lindenstrauss consider $\co$ as sitting in a subspace, let us call it $X_{AL}$, of $\ell_\infty$ whose construction relies on the following classical and extremely useful combinatorial lemma.
		
		\begin{lemm}
			\label{lem:almostdisjoint} 
			There exists a family $(N_\lambda)_{\lambda \in \bR}$ of infinite subsets of $\bN$ such that $N_\lambda \cap N_\mu$  is finite for all $\lambda \neq \mu \in \bR$. 
		\end{lemm}
		
		The idea behind one possible proof of this lemma is simple: for each $\lambda \in \bR$ select a strictly increasing sequence of rationals converging to $\lambda$ and then identify this sequence with an infinite subset of $\bN$ by using a fixed bijection from $\bQ$ onto $\bN$. 
		
		The subspace $X_{AL}$ is simply the closed linear span in $\ell_\infty(\bN)$ of $c_0(\bN)\cup\{\car_{N_\lambda}\colon  \lambda \in \bR\}$, where $\car_{N_\lambda}$ is the indicator function of $N_\lambda$. As we will see, $X_{AL}$ is Lipschitz equivalent to $\co(\bR)$, where for a set $\Gamma$, we define $c_0(\Gamma)$ to be the set of all $x := (x_\gamma)_{\gamma \in \Gamma} \in \bR^\Gamma$ such that for all $\eps>0$ there exists a finite subset $F$ of $\Gamma$ so that $\abs{x_\gamma}\le \eps$ for all $\gamma \in \Gamma \setminus F$. We equip $c_0(\Gamma)$ with the supremum norm $\norm{\cdot}_\infty$ and we denote by $(e_\gamma)_{\gamma \in \Gamma}$ its canonical basis. As in the countable case, the dual of $c_0(\Gamma)$ is isometric to $\ell_1(\Gamma)$, the space of all summable families $y := (y_\gamma)_{\gamma \in \Gamma} \in \bR^\Gamma$ equipped with $\norm{y}_1 := \sum_{\gamma \in \Gamma}\abs{y_\gamma}$, through the duality $\langle x,y\rangle := \sum_{\gamma \in \Gamma}x_\gamma y_\gamma.$
		
		\begin{theo}
			\label{thm:Aharoni-Lindenstrauss-example}
			The (nonseparable) Banach space $X_{AL}$ is Lipschitz equivalent but not linearly isomorphic to $c_0(\bR)$.
		\end{theo}

		\begin{proof}
			We split the proof into several steps.
			
			\medskip
			
			\underline{Step 1:} $X_{AL}$ is not isomorphic to $c_0(\bR)$.
			
			\medskip
			
			\noindent This step is quite elementary and relies on the observation that the existence of a countable family of separating linear functionals for a space $X$ (i.e. there is $(x_n^*)_{n=1}^\infty\subset X^*$ such that for all $x\in X\setminus\{0\}$, there exists $n\in \bN$ such that $x_n^*(x)\neq 0$) is stable under isomorphisms (to see this use the adjoint). The restriction to $X_{AL}$ of the coordinate linear functionals on $\ell_\infty(\bN)$ provides a countable family of separating functionals on $X_{AL}$. On the other hand, if $(z_n^*)_{n=1}^\infty$ is a countable family in $\ell_1(\bR)=c_0(\bR)^*$, the union of the supports of the $z_n^*$ is countable, so there exists $\lambda \in \bR$ such that $z_n^*(e_\lambda)=0$ for all $n\in \bN$. Consequently, there is no countable family of separating functionals in $c_0(\bR)^*$, from which we deduce that $X_{AL}$ is not isomorphic to $c_0(\bR)$.
			
			\medskip 
			
			\underline{Step 2:} $X_{AL}/\co(\bN)\oplus_\infty \co(\bN)$ is linearly isometric to $\co(\bR)$
			
			\medskip
			
			\noindent Since $c_0(\bR)\oplus_\infty c_0(\bN)$ is easily seen to be linearly isometric to $c_0(\bR)$, it is sufficient to show that $X_{AL}/\co(\bN)$ is linearly isometric to $\co(\bR)$. If $(N_\lambda)_{\lambda \in \bR}$ is the uncountable collection of infinite subsets of $\bN$ such that $N_\lambda \cap N_\mu$ is finite for all $\lambda \neq \mu \in \bR$ that was used in the definition of $X_{AL}$, we claim that the map $J\colon c_{00}(\bR)\to X_{AL}/\co(\bN)$ defined by $J((x_\lambda)_{\lambda \in \bR}) := \sum_{\lambda \in \bR} x_{\lambda} Q(\car_{N_\lambda})$ is a linear isometry. To see that, fix $\lambda_1,\dots,\lambda_n$ distinct in $\bR$ and $x_{\lambda_1},\dots,x_{\lambda_n}$ in $\bR$. Since
			\begin{equation*}
				\Big\|\sum_{i=1}^n x_{\lambda_i} Q(\car_{N_{\lambda_i}})\Big\|_{X_{AL}/\co(\bN)} = \inf_{y\in \co(\bN)} \Big\| \sum_{i=1}^n x_{\lambda_i} Q(\car_{N_{\lambda_i}}) -y\Big\|,
			\end{equation*}
			a moment of thought reveals that 
			$$\Big\|\sum_{i=1}^n x_{\lambda_i} Q(\car_{N_{\lambda_i}})\Big\|_{X_{AL}/\co(\bN)} = \max_{1\le i\le n}\abs{x_{\lambda_i}}$$ 
			if the $N_{\lambda_i}$ formed a collection of infinite and pairwise disjoint subsets of $\bN$. But this does not have to be the case, so we proceed as follows. Remembering that for $i\neq j$, $N_{\lambda_i}\cap N_{\lambda_j}$ is finite, the set $M_i := N_{\lambda_i} \setminus (\cup_{j\neq i}N_{\lambda_i}\cap N_{\lambda_j})$ is infinite and disagree with $N_{\lambda_i}$ only on a finite set. By construction the $M_i$ form a collection of infinite pairwise disjoint subsets of $\bN$ such that $Q(\car_{N_{\lambda_i}}) = Q(\car_{M_i})$ and it follows from the above observation that $J$ is a linear isometry. This isometry extends to a linear isometry, still denoted by $J$, from $c_0(\bR)$ into $X_{AL}/\co(\bN)$. Finally, since $Q(\car_{N_\lambda})$ is in the image of $J$ for all $\lambda \in \bR$, the image of $J$ is dense in $X_{AL}/\co(\bN)$ and $J$ being isometric, we have that $J(c_0(\bR)) = X_{AL}/\co(\bN)$, completing the proof of Step 2.
			
			\medskip
			
			\underline{Step 3:} The quotient map $Q\colon X_{AL} \to X_{AL}/\co$ admits a $2$-Lipschitz lifting.
			
			\medskip
			
			\noindent According to Step 2, it is sufficient to find a $1$-Lipschitz map $\psi \colon c_0(\bR)^+ \to X_{AL}$ such that $Q\circ \psi = J$, where $c_0(\bR)^+$ is the positive cone of $c_0(\bR)$. Indeed, such a map $\psi$ can then be extended to $c_0(\bR)$, by setting $\psi(x):=\psi(x^+)-\psi(x^-)$, where $x^+$ and $x^-$ are the positive and negative parts of $x\in c_0(\bR)$. Since $x\mapsto x^+$ and $x\mapsto x^-$ are $1$-Lipschitz on $c_0(\bR)$, $\psi$ is $2$-Lipschitz on $c_0(\bR)$ and thanks to $Q$ and $J$ being linear, we have that $Q\circ \psi(x) = J(x)$ for all $x\in \co(\bR)$. 
			
			Now, observe that any $x := (x_\lambda)_{\lambda \in \bR} \in c_0(\bR)^+$, can be decomposed (not uniquely) as $x=\sum_{j=1}^\infty a_je_{\lambda_j}$, where $(a_j)_{j=1}^\infty$ is nonincreasing and tending to $0$ and define 
			\begin{equation*}
				\psi(x) := \sum_{i=1}^\infty a_j\car_{M_j},
			\end{equation*}
			where $M_1 := N_{\lambda_1}$ and for $j>1$, $M_j := N_{\lambda_j}\setminus \cup_{k<j}N_{\lambda_k}$.
			Note that $\psi(x)$ does not depend on the choice of the decomposition $x=\sum_{i=1}^\infty a_j e_{\lambda_j}$ with $(a_j)_{j=1}^\infty$ nonincreasing. Since $N_{\lambda_j}\setminus M_j$ is finite, $\psi(x)\in X_{AL}$ and 
			\begin{equation*}
				Q\circ\psi(x)=\sum_{i=1}^\infty a_jQ(\car_{M_j})=\sum_{i=1}^\infty a_jQ(\car_{N_{\lambda_j}})=J(x).
			\end{equation*}
			
			For $n\in \bN$, let $A_n := \cspa\{e_\lambda \colon \lambda\in \bR\ \text{such that}\ n\notin N_\lambda\}\subset c_0(\bR)$. We claim that $\psi(x)_n=d(x,A_n)$, which will insure that $\psi$ is $1$-Lipschitz on $c_0(\bR)^+$. 
			
			First, assume that $n\notin \cup_{j=1}^\infty N_{\lambda_j}=\cup_{j=1}^\infty M_j$. Then, $\psi(x)_n=0$. On the other hand, for all $j\in \bN$, $e_{\lambda_j}\in A_n$, so $x\in A_n$ and $d(x,A_n)=0$.
			
			Assume now that $n\in \cup_{j=1}^\infty N_{\lambda_j}$. Then, there exists a unique $j_0\in \bN$ such that $n\in M_{j_0} \subset N_{\lambda_{j_0}}$ and $\psi(x)_n = a_{j_0}$. Remark that $n\notin N_{\lambda_j}$ if $j<j_0$ and that $n\in N_{\lambda_{j_0}}$. This, combined with the fact that the sequence $(a_k)_{k=1}^\infty$ is positive and nonincreasing, implies that 
			\begin{equation*}
				d(x,A_n) = d\Big(\sum_{j=1}^\infty a_j e_{\lambda_j}, \cspa\{e_\lambda \colon \lambda\in\bR\ \text{such that}\ n\notin N_\lambda\}\Big) = a_{j_0} = \psi(x)_n.
			\end{equation*}
			
		\end{proof}
		
		\begin{rema} 
			As we already mentioned, it is still an open question to know whether two Lipschitz isomorphic \underline{separable} Banach spaces are always linearly isomorphic. At some point, it was hoped that the above argument could be adapted in the separable setting. The Godefroy-Kalton lifting Theorem (Corollary \ref{cor:GK-lifting}) ruined these hopes!
		\end{rema}
		
		\begin{rema} We refer the reader to the Notes of this chapter for an example of a compact nonmetrizable space $K$ such that $C(K)$ is Lipschitz equivalent but not linearly isomorphic to a space $c_0(\Gamma)$. 
		\end{rema}
		
		\subsection{A Lipschitz-free space construction of a pair of uniformly equivalent but nonisomorphic separable Banach spaces} 
		\label{sec:example-UH-not-Isomorphic}
		
		The first examples of pairs of nonisomorphic separable Banach spaces that are uniformly homeomorphic were provided by M. Ribe in \cite{Ribe1984}. Kalton's reworking of Ribe's delicate construction will be discussed in Section \ref{sec:Ribe-Kalton-machinery}. In this section, we discuss a Lipschitz-free space construction combined with the lifting method, due to Kalton \cite{Kalton2004}, to provide similar examples. 
		
		Let $(X,\norm{\cdot})$ be a Banach space and $\theta \in (0,1)$. For $x,y\in X$, let
		\begin{equation*}
			d_\theta(x,y) := \max\{\norm{x-y}^\theta,\norm{x-y}\}
		\end{equation*}
		and denote by $X^\theta$ the metric space $(X,d_\theta)$. Since $\norm{x-y}\le d_\theta(x,y)$, $Id_X$ is a $1$-Lipschitz map from $X^\theta$ to $X$. By Corollary \ref{cor:factorization-Lipschitz-maps}, there exists a linear map $\hat{I} \colon  \cF(X^\theta) \to \cF(X)$ such that for all $x\in X$, $\delta_X(x)=\hat{I} \circ \delta_{X^\theta}(x)$ and $\|\hat{I}\|\le 1$. The map $Q := \beta_X \hat{I}\colon \cF(X^\theta) \to X$ is bounded linear and onto and such that $Q \circ \delta_{X^\theta}(x)=x$ for all $x\in X$. Finally, the map $\delta_{X^\theta}$ is uniformly continuous from $X$ to $\cF(X^\theta)$. The existence of this uniformly continuous lifting of $Q$ implies, by Proposition \ref{prop:section-method} that $\cF(X^\theta)$ is uniformly homeomorphic to $X\oplus \mathrm{Ker}(Q)$. Therefore, each time one can find a Banach space $X$ failing an isomorphic property passing to subspaces that $\cF(X^\theta)$  possesses, we have a new pair of uniformly equivalent but nonisomorphic Banach spaces. Below is one such example. 
		
		\begin{theo}
			\label{thm:Kalton-Lip-free-counterexample}
			Let $X$ be a separable Banach space failing the Schur property (e.g. $\co$ or $\ell_p$ for $1<p<\infty$). Then, $\cF(X^\theta)$ is uniformly homeomorphic but not linearly isomorphic to $X\oplus \mathrm{Ker}(Q)$.
		\end{theo}
		
		\begin{proof} 
			By Corollary \ref{cor:Schur-Lip-gauge-free}, $\cF(X^\theta)$ has the Schur property, while $X$ and thus  $X\oplus \text{Ker}(Q)$ does not. 
		\end{proof}
		
		\begin{rema} 
			At the time Kalton proved Theorem \ref{thm:Kalton-Lip-free-counterexample}, it was unclear whether the spaces could be Lipschitz equivalent. But as we already mentioned, it follows from the work of Aliaga, Gartland, Petitjean and Proch\'azka \cite{AGPP2021} that $\cF(X^\theta)$ has the Radon-Nikod\'ym property. So, we can deduce from the classical differentiability results that, if $X$ is a separable Banach space failing the Schur or the Radon-Nikod\'ym property, then $X \oplus \text{Ker}(Q)$ does not even bi-Lipschitzly embed into $\cF(X^\theta)$.
		\end{rema}
		
		\section{The Ribe-Kalton machinery for counterexamples}
		\label{sec:Ribe-Kalton-machinery}
		
		In \cite{Ribe1984}, Ribe proved the following. 
		
		\begin{theo}\label{theo:RibeExample}
			Let $(p_n)_{n=1}^\infty$ be a strictly decreasing sequence in $(1,2)$ which is converging to $1$, then $(\sum_{n=1}^\infty \ell_{p_n})_{\ell_2}$ is uniformly homeomorphic to $(\sum_{n=1}^\infty \ell_{p_n})_{\ell_2}\oplus \ell_1$   
		\end{theo}
		A proof of Theorem \ref{theo:RibeExample} will be given in Section \ref{sec:Ribe-example}. Note that these spaces cannot be isomorphic since only one of them is reflexive. This pair of spaces was the first counterexample to the uniform rigidity problem. Ribe's original construction was quite delicate and was further improved, simplified and streamlined in \cite{AharoniLindenstrauss1985}, \cite{Benyamini1985} and \cite{JLS1996}. In \cite{Kalton2012} and \cite{Kalton2013}, Kalton overhauled Ribe's argument and turned it into a quite abstract tool but with a significantly wider scope. We start the study of this circle of ideas with some well-known generalities about the uniform structure of spaces, their balls and their spheres.
		
		Given $X$ and $Y$ Banach spaces and $f\colon S_X \to S_Y$. 
		The \emph{radial extension} of $f$ is given by 
		\begin{equation*}
			\fbar(x) :=\begin{cases}
				\norm{x}f(\frac{x}{\norm{x}}) \quad \textrm{ if } x\neq 0,\\
				0 \quad \textrm{ if } x=0.
			\end{cases} 
		\end{equation*}
		Note that a radial extension is always a \emph{positively homogeneous} map, i.e. $\fbar(\alpha x) = \alpha \fbar(x)$ for all $x\in X$ and $\alpha\ge 0$.                                        
		
		If $f$ were a Lipschitz homeomorphism between the spheres, then it is easy to see that $f$ induces, via its radial extension, a Lipschitz homeomorphism between $X$ and $Y$ and also between their unit balls.
		However, if $f$ is merely a uniform homeomorphism between the spheres, its radial extension is solely a bijection that is norm-preserving, positively homogeneous and such that $\fbar$ and its inverse are uniformly continuous when restricted to the respective balls. In Exercise \ref{ex:balls-to-spheres}, you are asked to prove this statement as well as its converse.
		
		In the uniform category, it is simply not true that a uniform homeomorphism between spheres implies a uniform homeomorphism between the spaces. Indeed, we have seen that, for $1\le p \neq q<\infty$,  $\ell_p$ is not uniformly homeomorphic to $\ell_q$ but it is well known that the Mazur map $M_{p,q}\colon \ell_p \to \ell_q$ given by 
		\begin{equation*}
			M_{p,q}((x_n)_{n=1}^\infty):= (\sign(x_n)\abs{x_n}^{p/q})_{n=1}^\infty
		\end{equation*}
		is a uniform homeomorphism (with inverse $M^{-1}_{p,q}=M_{q,p}$) between the corresponding spheres. 
		Note that $\norm{M_{p,q}(x)}_q = \norm{x}_p^{p/q}$ and that $M_{p,q}$ is not positively homogeneous. This can be remedied by taking the radial extension of the restriction of $M_{p,q}$ to the unit sphere, and with a slight abuse of notation, the map 
		\begin{equation*}
			x\mapsto \bar{M}_{p,q}(x):= \norm{x}_p^{1-p/q} M_{p,q}(x)
		\end{equation*}
		is a norm-preserving and positively homogeneous uniform homeomorphism between the corresponding spheres and also the corresponding balls.
		
		In \cite{Nahum2001}, Nahum showed that if two Banach spaces $X$ and $Y$ are uniformly equivalent, then $S_{X\oplus \bR}$ and $S_{Y\oplus \bR}$ are also uniformly equivalent. 
		The following theorem can be extracted from a careful investigation of the proof of Nahum's theorem (see Proposition 9.17 in \cite{BenyaminiLindenstrauss2000} for the proof).
		
		\begin{theo}
			\label{thm:Nahum} 
			Let $X$ and $Y$ be two Banach spaces. Assume that $\psi\colon X \to Y$ is a bijection such that $\omega_\psi(t)\le Lt+C$ and $\omega_{\psi^{-1}}(t)\le Lt+C$, for some $L>0$, $C>0$ and all $t\ge 0$. Then, there exists a constant $L_1:=L_1(L,X,Y)>0$ and a bijection $\varphi \colon S_{X\oplus_\infty \bR} \to S_{Y\oplus_\infty \bR}$ such that for all $t\ge 0$, 
			\begin{equation*}
				\omega_\varphi(t)\le L_1\omega_\psi(t) \quad  \mathrm{ and } \quad  \omega_{\varphi^{-1}}(t)\le L_1\omega_{\psi^{-1}}(t).
			\end{equation*}
			Moreover, if $\psi$ is a uniform homeomorphism from $X$ onto $Y$, then $\varphi$ and $\varphi^{-1}$ are uniformly continuous. 
		\end{theo}
		
		Kalton observed a simple but crucial consequence of Nahum's theorem as stated above.
		
		\begin{coro}
			\label{cor:Nahum}  
			If there is a bijection from $X$ onto $Y$ that is coarse-Lipschitz as well as its inverse, then for some $L>0$ and every $\vep>0$ there is a bijection $f_\vep \colon S_{X\oplus\bR} \to S_{Y\oplus \bR}$ such that for all $t\ge 0$, $\max\{\omega_{f^{-1}_{\vep}}(t), \omega_{f_{\vep}}(t)\}\le Lt+\vep$, i.e. for all $u,v\in S_{X\oplus\bR}$,
			\begin{equation*}
				\norm{f_{\vep}(u) - f_{\vep}(v)} \le L\norm{u-v}+\vep,
			\end{equation*}
			and for all $u,v\in S_{Y\oplus\bR}$, 
			\begin{equation*}
				\norm{f^{-1}_{\vep}(u) - f^{-1}_{\vep}(v)} \le L\norm{u-v}+\vep,
			\end{equation*}
			Moreover, if $X$ and $Y$ are uniformly homeomorphic, the maps $f_\vep$ and $f^{-1}_{\vep}$ can be taken to be uniformly continuous.
			%$X\oplus\mathbb R$ and $Y\oplus \mathbb R$ are close (respectivelyuniformly close). Hence, if $X$ and $Y$ are coarse-Lipschitz homeomorphic (respectively uniformly homeomorphic) Banach spaces that are linearly isomorphic to their hyperplanes, then $X$ and $Y$ are close (respectively uniformly close).
		\end{coro}
		
		\begin{proof}  
			As we will often do, we will only detail the more delicate uniform homeomorphism case.
			Let $\psi \colon X\to Y$ be a uniform homeomorphism.  Then, there exist constants $L,C>0$ such that for all $t>0$,
			\begin{equation*}
				\omega_\psi(t)\le Lt+C \quad  \mathrm{ and } \quad \omega_{\psi^{-1}}(t)\le Lt+C.
			\end{equation*} 
			For $n\in \bN$, let $\psi_n(x) := n^{-1}\psi(nx)$ and observe that
			\begin{equation*}
				\omega_{\psi_n}(t)\le Lt+\frac{C}{n} \quad  \mathrm{ and } \quad \omega_{\psi_{n}^{-1}}(t)\le Lt+\frac{C}{n}.
			\end{equation*}
			If $\varphi_n$ is the map associated to $\psi_n$ given by Theorem \ref{thm:Nahum} and $L_1$ is the positive constant given by the same theorem, then $\varphi_n \colon  S_{X\oplus_\infty \bR} \to S_{Y\oplus_\infty \bR}$ is a uniform homeomorphism such that
			\begin{equation*}
				\omega_{\varphi_n}(t)\le L_1Lt+\frac{L_1C}{n} \quad  \mathrm{ and } \quad \omega_{\varphi_{n}^{-1}}(t)\le L_1Lt+\frac{L_1C}{n}.
			\end{equation*}
			
			%Then we can use Proposition \ref{almostLipschitz} to conclude. 
		\end{proof}
		
		Kalton introduced the following general terminology. 
		
		\begin{defi}
			Let $(M,d_M)$ and $(N,d_N)$ be unbounded metric spaces. 
			\begin{itemize}
				\item A map $f\colon M\to N$ is of \emph{CL-type $(L,\eps)$} if for all $t\ge 0$,
				\begin{equation*}
					\omega_f(t)\le Lt+\eps.
				\end{equation*}
				\item A map $f\colon M\to N$ is a \emph{coarse-Lipschitz
					homeomorphism} if it is a \emph{bijection} and $f$ and $f^{-1}$ are coarse-Lipschitz maps. In this case, we say that $f$ is a \emph{CL-homeomorphism} of type $(L,\eps)$ if both $f$ and $f^{-1}$ are of CL-type $(L,\eps)$. 
				\item Finally, $M$ and $N$ are said to be \emph{almost Lipschitz isomorphic} if for some $L$ and every $\eps>0$ there is a CL-homeomorphism $f_\vep\colon M\to N$ of type
				$(L,\eps)$. Moreover, if $f_\vep$ and $f^{-1}_\vep$ are also uniformly continuous, then  $M$ and $N$ are said to be \emph{uniformly almost Lipschitz isomorphic}.
			\end{itemize}
		\end{defi}
		
		\begin{rema}
			Note that the notion of coarse-Lipschitz homeomorphism is formally stronger than the notion of coarse-Lipschitz equivalence.
		\end{rema} 
		
		Corollary \ref{cor:Nahum} rephrased using this new terminology states that:
		\begin{itemize}
			\item If $X$ and $Y$ are coarse-Lipschitz homeomorphic, then $S_{X\oplus \bR}$ and $S_{Y\oplus \bR}$ are almost Lipschitz isomorphic.
			\item If $X$ and $Y$ are uniformly homeomorphic, then $S_{X\oplus \bR}$ and $S_{Y\oplus \bR}$ are uniformly almost Lipschitz isomorphic.
		\end{itemize}

		\subsection{Close and uniformly close Banach spaces}
		\label{sec:closeness}
		
		In this section, we study the notion of \emph{closeness} between Banach spaces. This notion is central in Kalton's reworking of Ribe's circle of ideas. The notion of closeness appears naturally when one tries to push nontrivial information about spheres to global information for
		the spaces. With this goal in mind, the following classical and elementary lemma, which we state without proof, comes handy.
		
		\begin{lemm}
			\label{lem:polar} 
			Let $X$ be a Banach space and $x_1,x_2\in X$ with $\norm{x_1}\ge \norm{x_2}>0$. Then,
			\begin{equation}
				\label{eq1:polar}
				\Big\| \frac{x_1}{\norm{x_1}}-\frac{x_2}{\norm{x_2}}\Big\| \le 2\frac{\norm{x_1-x_2}}{\norm{x_1}},
			\end{equation} 
			and
			\begin{equation}
				\label{eq2:polar}
				\norm{x_1-x_2}\le \norm{x_1}-\norm{x_2} +\norm{x_2}\Big\| \frac{x_1}{\norm{x_1}}-\frac{x_2}{\norm{x_2}}\Big\|\le 3\norm{x_1-x_2}.
			\end{equation}
		\end{lemm}
		
		The next lemma contains such simple quantitative relation between the $CL$-type of the restriction to the unit sphere of a positively homogeneous map with global inequalities for the map itself.
		
		\begin{lemm}
			\label{lem:local-global-homogeneous}  
			Suppose $f \colon X \to Y$ is positively homogeneous and fix $L\ge 1$, $\eps \in (0,1]$ and $K\ge 0$. 
			\begin{enumerate}[(i)]
				\item If $f_{\restriction_{S_X}}$ is of CL-type $(L,\eps)$ and $\|
				f_{\restriction_{S_X}}(x)\|_Y\le K$ for all $x\in S_X$, then for all $x_1,x_2\in X$,
				\begin{equation*}
					\norm{f(x_1)-f(x_2)}_Y\le (4L+2K)\max\{\norm{x_1-x_2}_X,\eps\norm{x_1}_X,\eps\norm{x_2}_X\}.
				\end{equation*}
				%$\norm{f}_\eps \le 2K+4L.$
				\item If for all $x_1,x_2\in X$,
				\begin{equation*}
					\norm{f(x_1)-f(x_2)}_Y\le L\max\{\norm{x_1-x_2}_X,\eps\norm{x_1}_X,\eps\norm{x_2}_X\},
				\end{equation*}
				then $f_{\restriction_{S_X}}$ is of CL-type $(L,L\eps).$
			\end{enumerate}
		\end{lemm}
		
		\begin{proof}
			Let $\varphi:=f_{\restriction_{S_X}}$.

			$(i)$ Suppose, as we may, that $\norm{x_1}\ge \norm{x_2}>0.$  Then,
			\begin{align*} 
				\|f(x_1)-f(x_2)\| &=\Big\|\norm{x_1} \varphi\Big(\frac{x_1}{\norm{x_1}}\Big)-\norm{x_2}\varphi\Big(\frac{x_2}{\norm{x_2}}\Big)\Big\|\\
				&\le K\norm{x_1-x_2} +\norm{x_2}\Big\|\varphi\Big(\frac{x_1}{\norm{x_1}}\Big)-\varphi\Big(\frac{x_2}{\norm{x_2}}\Big)\Big\|\\
				&\le K\norm{x_1-x_2} +L\norm{x_2}\Big\| \frac{x_1}{\norm{x_1}}-\frac{x_2}{\norm{x_2}}\Big\| + \eps\norm{x_2}\\
				& \stackrel{\eqref{eq1:polar}}{\le} (K+2L) \norm{x_1-x_2}+\eps\norm{x_2}\\
				&\le (2K+4L)\max\{\norm{x_1-x_2},\eps\norm{x_2}\}.
			\end{align*}
			
			$(ii)$ If $x_1,x_2\in S_X$, then 
			\begin{equation*}
				\norm{\varphi(x_1)-\varphi(x_2)} \le L\max\{\norm{x_1-x_2},\eps\}\le L\norm{x_1-x_2}+L\eps.
			\end{equation*}
		\end{proof}
		
		\begin{rema}
			Note that the positive homogeneity is not used to derive the trivial estimate in Lemma \ref{lem:local-global-homogeneous} $(ii)$.
		\end{rema}
		
		The next elementary proposition about a local-global phenomenon is the last stop before the formal introduction of the notion of (uniform) closeness.
		
		\begin{prop}
			\label{prop:ALI-spheres->closeness}
			If two Banach spaces $X$ and $Y$ have almost Lipschitz isomorphic spheres, then there is a constant $L>0$ and for every $\vep>0$ a bijection $f_\vep \colon X \to Y$ such that 
			\begin{enumerate}[(i)]
				\item $f_\vep$ and $f^{-1}_\vep$ are positively homogeneous,
				\item $\sup_{\norm{x}_X\le 1} \norm{f_\vep(x)}_Y<\infty$ and $\sup_{\norm{y}_Y\le 1} \norm{f^{-1}_\vep(y)}_X<\infty$,
				\item for all $x_1,x_2\in X$,
				\begin{equation*}
					\norm{f_\vep(x_1)-f(x_2)}_Y\le L\max\{\norm{x_1-x_2}_X,\eps\norm{x_1}_X,\eps\norm{x_2}_X\},
				\end{equation*}
				and for all $y_1,y_2\in Y$,
				\begin{equation*}
					\norm{f^{-1}_\vep(y_1)-f^{-1}_\vep(y_2)}_X\le L\max\{\norm{y_1-y_2}_Y,\eps\norm{y_1}_Y,\eps\norm{y_2}_Y\},
				\end{equation*}
			\end{enumerate}
			Moreover, if $X$ and $Y$ have uniformly almost Lipschitz isomorphic spheres, then the restrictions of the maps $f_{\vep}$ and $f^{-1}_{\vep}$ to the corresponding unit balls are uniformly continuous.
		\end{prop}
		
		\begin{proof}
			By assumption, there is $D>0$ and for every $\vep>0$ a bijection $\varphi_\vep\colon S_X\to S_Y$ so that $\varphi_\vep$ and $\varphi^{-1}_\vep$ have CL-type $(D,\eps)$. Applying Lemma \ref{lem:local-global-homogeneous} to the radial extensions of $\varphi_\vep$ and $\varphi^{-1}_\vep$ provides the desired maps $f_\vep$ and $f^{-1}_\vep$ (which are uniform homeomorphisms when restricted to the unit balls under the additional assumption of the second part) and we get the conclusion with $L=4D+2$. 
		\end{proof}
		
		A map $f\colon X \to Y$ between Banach spaces such that $\sup_{\norm{x}_X\le 1} \norm{f(x)}_Y<\infty$ will be called \emph{bounded}. The conclusion of Proposition \ref{prop:ALI-spheres->closeness} led Kalton to the following definition.
		
		\begin{defi}
			\label{def:closeness}
			Given $L\ge 1$, two Banach spaces $X$ and $Y$ are said to be \emph{$L$-close}, if for every $\eps>0$ we can find a bijection $f_\vep\colon X\to Y$ such that 
			\begin{enumerate}[(i)]
				\item $f_\vep$ and $f^{-1}_\vep$ are positively homogeneous,
				\item $f_\vep$ and $f^{-1}_\vep$ are bounded,
				\item for all $x_1,x_2\in X$,
				\begin{equation*}
					\norm{f_\vep(x_1)-f(x_2)}_Y\le L\max\{\norm{x_1-x_2}_X,\eps\norm{x_1}_X,\eps\norm{x_2}_X\},
				\end{equation*}
				and for all $y_1,y_2\in Y$,
				\begin{equation*}
					\norm{f^{-1}_\vep(y_1)-f^{-1}_\vep(y_2)}_X\le L\max\{\norm{y_1-y_2}_Y,\eps\norm{y_1}_Y,\eps\norm{y_2}_Y\},
				\end{equation*}
			\end{enumerate}
			Furthermore, $X$ and $Y$ are said to be \emph{close} if they are $L$-close for some $L\ge 1$. 
			
			Moreover, if the restrictions of $f_\vep$ and $f^{-1}_\vep$ to the corresponding unit balls are uniformly continuous, we say that $X$ and $Y$
			are $L$-\emph{uniformly close} and, respectively, \emph{uniformly close}.
		\end{defi}
		
		Using the above terminology, observe that Proposition \ref{prop:ALI-spheres->closeness} shows that two Banach spaces with (uniformly) almost Lipschitz isomorphic spheres are (uniformly) close.
		
		An immediate consequence of Corollary \ref{cor:Nahum}, Proposition \ref{prop:ALI-spheres->closeness} and Definition \ref{def:closeness} is:
		
		\begin{coro}
			\label{cor:CLH->augmented-close}\, 
			\begin{enumerate}[(i)]
				\item If $X$ and $Y$ are coarse-Lipschitz homeomorphic Banach spaces, then $X\oplus\bR$ and $Y\oplus \bR$ are close.
				\item If $X$ and $Y$ are uniformly homeomorphic Banach spaces, then $X\oplus\bR$ and $Y\oplus \bR$ are uniformly close.
			\end{enumerate}
		\end{coro}
		
		\begin{rema}
			In Corollary \ref{cor:CLH->augmented-close}, if $X$ and $Y$ are isomorphic to their hyperplanes, then the conclusions hold with $X$ and $Y$ instead of $X\oplus\bR$ and $Y\oplus \bR$, respectively.  
		\end{rema}

		We now introduce a fair amount of convenient notation to simplify the statements and the proofs of the upcoming results.  Let $X$ and $Y$ be two Banach spaces.
		\begin{itemize}
			\item The space of all maps $f\colon X\to Y$ which are positively homogeneous and bounded will be denoted by $\mathcal H(X,Y)$. If we let $\norm{f}:=\sup_{\norm{x}_X\le 1} \norm{f(x)}_Y$, it is clear that $\norm{\cdot}$ is a norm on $\mathcal H(X,Y)$ with which it becomes a Banach space containing the space $B(X,Y)$ of all bounded linear operators.
			\item If $\eps \in (0,1]$, then for every $f\in \mathcal H(X,Y)$ we define
			$\norm{f}_\eps$ to be the least constant $L$ so that for all $x_1,x_2\in X$,
			\begin{equation*}
				\norm{f(x_1)-f(x_2)}_Y\le L\max\{\norm{x_1-x_2}_X,\eps\norm{x_1}_X,\eps\norm{x_2}_X\}.
			\end{equation*}
			\item We denote by $\mathcal G(X,Y)$ the subset of $\cH(X,Y)$ of all $f$ such  that $f$ is a bijection and $f^{-1}\in\cH(Y,X).$  
			\item For $f\in \mathcal G(X,Y)$, let
			\begin{equation*}
				[[f]] := \max(\norm{f},\norm{f^{-1}}),\qquad [[f]]_\eps := \max(\norm{f}_\eps,\norm{f^{-1}}_\eps).
			\end{equation*}
			\item We also define a metric $\Delta$ on $\cG(X,Y)$ by
			\begin{equation*}
				\Delta(f,g) := \max(\norm{f-g},\norm{f^{-1}-g^{-1}}).
			\end{equation*}
			\item We define the subspace $\mathcal {HU}(X,Y)$ as the set of $f \in \mathcal H(X,Y)$ such that the restriction of $f$ to $B_X$ (and hence to any bounded set) is uniformly continuous.
			\item Let $\mathcal {GU}(X,Y)$ be the subset of $\cG(X,Y)$ of all $f$ such that $f\in\mathcal {HU}(X,Y)$ and $f^{-1}\in\mathcal {HU}(Y,X).$ 
		\end{itemize}

		\begin{rema}
			According to the notation above, $X$ and $Y$ are \emph{$L$-close,} respectively \emph{uniformly $L$-close} for  some $L\ge1$, if for every $\eps>0$ we can find $f\in\cG(X,Y)$ (respectively
			$f\in\mathcal{GU}(X,Y)$) with $[[f]]_\eps\le L$ and $X$ and $Y$
			are \emph{(uniformly) close} if there exists $L$ such that they are
			(uniformly) $L$-close.
		\end{rema}
		
		In the next lemma, we gather some elementary properties of $\norm{\cdot}_\vep$ and leave the details of the proofs to the dutiful reader.
		\begin{lemm}
			Let $X$ and $Y$ be Banach spaces.
			\begin{enumerate}[(i)]
				\item For each $\eps \in (0,1]$,  $\norm{\cdot}_\eps$
				is a norm on $\mathcal H(X,Y)$ which is equivalent to the
				original norm. More precisely, for all $f\in \mathcal H(X,Y)$,
				\begin{equation*}
					\norm{f}\le \norm{f}_\eps\le \frac{2}{\eps}\norm{f}.
				\end{equation*}
				\item For all $f\in \cH(X,Y)$, $\vep \mapsto \norm{f}_\eps$ is decreasing.
				\item $\sup_{\eps\in (0,1]}\norm{f}_\eps<\infty$ if and only if $f$ is Lipschitz. In that case, $\Lip(f)=\sup_{\eps \in (0,1]}\norm{f}_\eps$.
			\end{enumerate}
		\end{lemm}
		
		Before we delve more into the stability properties of the notion of (uniform) closeness, we first show that (uniform) closeness between Banach spaces is in fact the same notion as (uniform) almost Lipschitz equivalence between their spheres. In order to do this, we first introduce a normalized version of $f\in\cG(X,Y)$:
		\begin{equation*}
			\hat f(x) := \begin{cases}
				\norm{x}\frac{f(x)}{\|f(x)\|} \qquad \textrm{ if }x\neq 0, \\
				0\qquad\qquad\qquad \textrm{ if } x=0.
			\end{cases}
		\end{equation*}
		It is clear that $\hat f\in \mathcal G(X,Y)$ with $(\hat f)^{-1}=\widehat{f^{-1}}$,  and is norm-preserving. Furthermore, if $f\in\mathcal {GU}(X,Y)$, then $\hat f\in\mathcal {GU}(X,Y)$. In the next lemma, we record elementary quantitative relations between the various numerical quantities introduced above for a map and its normalization. 
		\begin{lemm}
			\label{lem:quant}\, 
			\begin{enumerate}[(i)]
				\item If $f\in \cG(X,Y)$, then $\hat f\in \cG(X,Y)$ and for all $\eps \in (0,1]$,
				\begin{equation*}
					[[\hat f]]_\eps \le 2[[f]]_\eps^2+1.
				\end{equation*}
				\item If $f,g\in\cG(X,Y)$, then
				\begin{equation*}
					\Delta(\hat f,\hat g) \le 2[[f]]\Delta(f,g).
				\end{equation*}
			\end{enumerate}
		\end{lemm}
		
		\begin{proof}
			$(i)$ If $x_1,x_2\in X$ are both nonzero and $\norm{x_1}\ge \norm{x_2}$, then, using Lemma \ref{lem:polar},
			\begin{align*} 
				\norm{\hat f(x_1) - \hat f(x_2)} & \stackrel{\eqref{eq2:polar}}{\le} \norm{x_1} - \norm{x_2} + \norm{x_1}\Big\|\frac{f(x_1)}{\norm{f(x_1)}} - \frac{f(x_2)}{\norm{f(x_2)}}\Big\|\\
				& \stackrel{\eqref{eq1:polar}}{\le} \norm{x_1 - x_2} + 2\norm{x_1}\frac{\norm{f(x_1) - f(x_2)}}{\norm{f(x_1)}}\\
				& \le \norm{x_1 - x_2} + 2\norm{f^{-1}}_\epsilon\norm{f(x_1) - f(x_2)}\\
				& \le (1 + 2\norm{f^{-1}}_\eps\norm{f}_\vep)\max(\norm{x_1 - x_2},\eps\norm{x_1},\eps\norm{x_2}),
			\end{align*}
			and the conclusion follows form the fact that $[[f]]_\vep:=\max\{\norm{f}_\eps,\norm{f^{-1}}_\vep\}$ and by carrying a similar calculation with $(\hat{ f})^{-1} =\widehat{f^{-1}}.$
			
			$(ii)$ If $\norm{x}=1$ we have, again using Lemma \ref{lem:polar},
			\begin{equation*}
				\norm{\hat f(x) - \hat g(x)} \stackrel{\eqref{eq1:polar}}{\le} 2\frac{\norm{f(x) - g(x)}}{\max\{\norm{f(x)}, \norm{g(x)}\}} \le 2 \min\{\norm{f^{-1}}, \norm{g^{-1}}\}\norm{f(x) - g(x)}.
			\end{equation*}
			Combining with a similar estimate on the inverses gives the conclusion.
		\end{proof}
		
		As promised, we can rephrase the notion of (uniform) closeness in terms of (uniform) almost Lipschitz isomorphisms.
		
		\begin{prop}
			\label{prop:close<->ALE-spheres} 
			Let $X$ and $Y$ be two Banach spaces. Then,
			\begin{enumerate}[(i)]
				\item $X$ and $Y$ are close if and only if $S_X$ and $S_Y$ are almost Lipschitz isomorphic.
				\item $X$ and $Y$ are uniformly close if and only if $S_X$ and $S_Y$ are uniformly almost Lipschitz isomorphic.
			\end{enumerate}
		\end{prop}
		
		\begin{proof} 
			We only treat the uniform case $(ii)$. Suppose first that $X$ and $Y$ are uniformly close.
			Then, there exists $L\ge 1$ such that for all $\eps \in (0,1]$ there exists
			$f:=f_{\eps}\in\mathcal {GU}(X,Y),$ with $[[f]]_\eps\le L.$  Then by Lemma \ref{lem:quant}
			$[[\widehat f]]_\eps\le 2L^2+1$ and by Lemma \ref{lem:local-global-homogeneous}, $\widehat f_{\restriction_{S_X}}$ is a uniform homeomorphism onto $S_Y$ such that $\widehat f$ and $\widehat f^{-1}$ both have CL-type $(2L^2+1, (2L^2+1)\eps).$
			
			Conversely if $\varphi \colon S_X\to S_Y$ is a uniform homeomorphism so
			that $\varphi$ and $\varphi^{-1}$ have CL-type $(L,\eps)$ then, again by Lemma \ref{lem:local-global-homogeneous}, the homogeneous extension $f$ of $\varphi$ belongs to $\mathcal {GU}(X,Y),$ and $[[f]]_\eps\le 4L+2.$
		\end{proof}
		
		\subsection {\texorpdfstring{Stability of closeness and uniform closeness under $\ell_p$-sums}{Stability of closeness and uniform closeness under sums}}
		
		Let us first observe that the notion of (uniform) closeness induces a transitive relation on the class of Banach spaces.
		
		\begin{prop}
			\label{prop:composition}
			Let $X$ and $Y$ be Banach spaces. 
			\begin{enumerate}[(i)]
				\item If $X$ and $Y$ are close and $Y$ and $Z$
				are close, then $X$ and $Z$ are close.
				\item If $X$ and $Y$ are uniformly close and $Y$ and $Z$
				are uniformly close, then $X$ and $Z$ are uniformly close.
			\end{enumerate}
		\end{prop}
		
		Proposition \ref{prop:composition} follows immediately from the elementary Lemma \ref{lem:composition} below whose straightforward proof can be safely omitted.
		
		\begin{lemm}
			\label{lem:composition}
			Let $X,Y,Z$ be Banach spaces. Then, for all $\eps \in (0,1)$, $f\in\mathcal H(X,Y)$ and
			$g\in\mathcal H(Y,Z)$, we have
			\begin{equation*}
				\norm{g\circ f}_\eps \le \norm{g}_\eps\norm{f}_\eps.
			\end{equation*} 
		\end{lemm}
		
		The stability of (uniform) closeness under $\ell_p$-sums is more technical and will be crucial in the forthcoming construction of counterexamples. The statement is the following. 
		
		\begin{prop}
			\label{prop:closeness-ellp-sums}  
			Let $(X_n)_{n=1}^{\infty}$  and $(Y_n)_{n=1}$
			be two sequences of Banach spaces. If there exists $L>0$ such that $X_n$
			and $Y_n$ are $L$-close (resp. uniformly $L$-close) for all $n\in\bN$,  then $(\sum_{n=1}^{\infty}X_n)_{\ell_p}$ and
			$(\sum_{n=1}^{\infty}Y_n)_{\ell_p}$ are $3L$-close (resp.
			uniformly $3L$-close), for $1\le p\le \infty$, where the case $p=\infty$ is interpreted as a $\co$-sum. 
		\end{prop}
		
		\begin{proof} 
			We will detail the uniformly continuous case only and restrict to the case of $1\le p<\infty.$  The other cases are simpler. Let $L \ge 1$ such for all $n\in \bN$, $X_n$ and $Y_n$ are uniformly $L$-close. Let 
			$\eps \in (0,1)$. We pick a sequence $(\eps_n)$ in $(0,\eps]$ with
			$\lim_{n\to\infty}\eps_n=0$. For all $n\in \bN$, there exists  $f_n\in\mathcal {GU}(X_n,Y_n),$ with $[[f_n]]_{\eps_n}\le L.$ Define $f\colon (\sum_{n=1}^{\infty}X_n)_{\ell_p}\to (\sum_{n=1}^{\infty}Y_n)_{\ell_p}$ by
			$f((x_n)_{n=1}^{\infty}) :=(f_n(x_n))_{n=1}^{\infty}.$ Note that $f$ is a bijection and that $f^{-1}((y_n)_{n=1}^{\infty})=(f_n^{-1}(y_n))_{n=1}^{\infty}$, for $(y_n)_{n=1}^{\infty} \in (\sum_{n=1}^{\infty}Y_n)_{\ell_p}$. Let 
			$x :=(x_n)_{n=1}^{\infty}\in (\sum_{n=1}^{\infty}X_n)_{\ell_p}$ and
			$x':=(x'_n)_{n=1}^{\infty}\in(\sum_{n=1}^{\infty}X_n)_{\ell_p}.$
			Then \begin{align*} 
				\norm{f(x) - f(x')}^p & \le L^p\sum_{n=1}^{\infty}\max\{\norm{x_n - x'_n}^p, \eps^p\norm{x_n}^p, \eps^p\norm{x'_n}^p\}\\
				& \le L^p(\norm{x - x'}^p + \eps^p(\norm{x}^p + \norm{x'}^p)\\
				& \le 3L^p \max\{\norm{x - x'}^p,\eps^p\norm{x}^p,
				\eps^p\norm{x'}^p\}.
			\end{align*} 
			Thus $\norm{f}_\eps \le 3L$ and combining with a similar estimate for
			$f^{-1}$ we obtain that $[[f]]_\eps\le 3L.$
			
			To check the uniform continuity of $f$ and $f^{-1}$ on the respective unit balls, we note first that $\omega_{g_n}(t),\omega_{h_n}(t)\le Lt+\eps_n$, where $g_n$ is the restriction of $f_n$ to $B_{X_n}$ and $h_n$ is the restriction of $f_n^{-1}$ to $B_{Y_n}$. We now set: 
			\begin{equation*}
				\omega(t):=\sup_n\max\{\omega_{g_n}(t),\omega_{h_n}(t)\}.
			\end{equation*} 
			One can check that it follows from the uniform continuity of the maps $g_n$ and $h_n$ and the fact that $\lim_n \eps_n=0$ that $\lim_{t\to 0}\omega(t)=0.$
			
			It remains to verify that $f$ is uniformly continuous on the unit ball of $(\sum_{n=1}^{\infty} X_n)_{\ell_p}$. Suppose $\eta>0$ and pick $\nu \in (0,1]$ so that $3L\sqrt\nu + \omega(2\sqrt\nu)\le \eta$. For $x,x'\in B_{(\sum_{n=1}^{\infty}X_n)_{\ell_p}}$ such that $\norm{x - x'}<\nu$, we will show that $\norm{f(x) - f(x')}<\eta$. For any $\bA\subset \bN$, we have 
			\begin{align}
				\norm{f(x) -f(x')} & = \Big(\sum_{n\in\bA} \norm{f_n(x_n) - f_n(x'_n)}^p + \sum_{n\notin\bA} \norm{f_n(x_n) - f_n(x'_n)}^p\Big)^{1/p} \notag\\
				& \le \Big(\sum_{n\in\bA} \norm{f_n(x_n) - f_n(x'_n)}^p \Big)^{1/p}  + \Big(\sum_{n\notin\bA} \norm{f_n(x_n) - f_n(x'_n)}^p\Big)^{1/p} \notag \\
				& \le \Big(\sum_{n\in\bA} \norm{f_n(x_n)}^p\Big)^{1/p} + \Big(\sum_{n\in\bA} \norm{f_n(x'_n)}^p \Big)^{1/p} + \Big(\sum_{n\notin\bA} \norm{f_n(x_n) - f_n(x'_n)}^p\Big)^{1/p} \label{eq:closeness-ellp-sums}.
			\end{align}
			If we let
			\begin{equation*}
				\bA := \Big\{n\colon \norm{x_n - x'_n}>\sqrt\nu \max\{\norm{x_n}, \norm{x'_n}\}\Big\},
			\end{equation*}
			then
			\begin{equation*}
				\Big( \sum_{n\in\bA}\max\{\norm{x_n}, \norm{x'_n}\}^p\Big)^{1/p}\le \sqrt\nu
			\end{equation*}
			and since $\norm{f_n}\le \norm{f_n}_{\eps_n}\le L$, we have
			\begin{equation}
				\label{eq2:closeness-ellp-sums}
				\Big(\sum_{n\in\mathbb A}\|f_n(x_n)\|^p\Big)^{1/p} \le L\sqrt\nu\ \ \ \text{and}\ \ \ \Big(\sum_{n\in\mathbb A}\|f_n(x'_n)\|^p\Big)^{1/p} \le L\sqrt\nu.  
			\end{equation}
			To upper bound the last term in \eqref{eq:closeness-ellp-sums}, assume that $\norm{x_n}\ge \norm{x'_n}>0$ and write 
			\begin{equation*}
				f(x_n) - f(x'_n) = f\Big(\frac{x_n}{\norm{x_n}}\Big)\Big(\norm{x_n} - \norm{x_n'}\Big) + \norm{x'_n}\Big(f\Big(\frac{x_n}{\norm{x_n}}\Big) - f\Big(\frac{x'_n}{\norm{x'_n}}\Big)\Big)
			\end{equation*}
			%$\norm{f_n(x_n)}\ge \norm{f_n(x'_n)}>0$ and write 
			%\begin{equation*}
			%   f_n(x_n) - f_n(x'_n) =\frac{f_n(x_n)}{\norm{f_n(x_n)}}\Big(\norm{f_n(x_n)} - \norm{f_n(x_n')}\Big) + \norm{f_n(x'_n)}\Big(\frac{f_n(x_n)}{\norm{f_n(x_n)}} - \frac{f_n(x'_n)}{\norm{f_n(x'_n)}}\Big)
			%\end{equation*}
			and using Lemma \ref{lem:polar}, we get that for all $n\in \bN$ such that $\min\{\norm{x_n}, \norm{x'_n}\}>0$,
			\begin{equation}
				\label{eq4:closeness-ellp-sums}
				\norm{f_n(x_n) - f_n(x'_n)} \le L\norm{x_n - x'_n} + \min\{\norm{x_n}, \norm{x'_n}\}\,\omega\left(\frac{2\norm{x_n - x'_n}}{\max\{\norm{x_n}, \norm{x'_n}\}}\right).
			\end{equation}
			Thus, for all $n\notin \bA$ we have 
			\begin{equation*}
				\norm{f_n(x_n) - f_n(x'_n)} \le L\norm{x_n - x'_n} + \min\{\norm{x_n},\norm{x'_n}\}\,\omega(2\sqrt\nu).
			\end{equation*}
			Note that the above inequality was immediate when $\min\{\norm{x_n}, \norm{x'_n}\}=0$. Hence
			\begin{equation}
				\label{eq3:closeness-ellp-sums}
				\Big(\sum_{n\notin \bA} \norm{f_n(x_n) - f_n(x'_n)}^p\Big)^{1/p}\le
				L\nu +\omega(2\sqrt\nu).
			\end{equation}
			Combining \eqref{eq:closeness-ellp-sums}, \eqref{eq2:closeness-ellp-sums} and \eqref{eq3:closeness-ellp-sums}, we have
			\begin{equation*}
				\norm{f(x) - f(x')} \le  3L\sqrt\nu + \omega(2\sqrt\nu) <\eta.
			\end{equation*}
			This and a similar calculation for
			$f^{-1}$ give the conclusion that $f\in \mathcal {GU}(X,Y)$ and we
			already know that $[[f]]_{\eps}\le 3L$. Since $\eps \in (0,1]$ was arbitrary we have shown that $(\sum_{n=1}^{\infty}X_n)_{\ell_p}$ and $(\sum_{n=1}^{\infty}Y_n)_{\ell_p}$ are $3L$-uniformly close.
		\end{proof}
		
		For further reference, we record an immediate corollary.
		\begin{coro}
			\label{cor:closness-ellp-sums}
			If $X$ and $Y$ are close (resp. uniformly close) Banach spaces, then
			\begin{enumerate}[(i)]
				\item $\ell_p(X)$ and $\ell_p(Y)$ are close (resp. uniformly close) for $1\le p<\infty$,
				\item $c_0(X)$ and $c_0(Y)$ are close (resp. uniformly close).
			\end{enumerate}
		\end{coro}
		
		\subsection{Deriving uniform homeomorphisms from uniform closeness} 
		\label{sec:closenness->CLU-homeomorphism}
		
		In Section \ref{sec:closeness}, we saw that uniformly homeomorphic Banach spaces are uniformly close (up to a one-dimensional perturbation). 
		In this section, we investigate how we can deduce uniform (respectively coarse-Lipschitz) homeomorphisms from uniform closeness (respectively closeness). By definition, uniform closeness between $X$ and $Y$ provides the existence of a collection of maps $(f_{\vep})_{\vep\in(0,1]}$ in $\mathcal{GU}(X,Y)$. These maps are getting closer and closer to Lipschitz ones as $\vep$ gets smaller and smaller. In order to construct a uniform homeomorphism out of this collection of maps, one needs to be able to glue them together in a clever way.  
		This crucial ``continuous gluing procedure'' for families of homogeneous maps of CL-type $(L,\eps)$, for $\eps \in (0,1]$, is carried out in the proof of the next proposition. The gluing procedure requires an additional regularity condition on the collection of maps. For convenience, the collection of maps will be indexed by $t\in [0,\infty)$ and we will consider their $\norm{\cdot}_{e^{-2t}}$ norms.
		
		\begin{prop}
			\label{prop:base-construction-gluing}  
			Let $X$ and $Y$ be Banach spaces and suppose there is a collection of maps $(f_{t})_{t\in[0,\infty)}$ in $\cH(X,Y)$ with the property that for some constant $K$ we have:
			\begin{equation}
				\label{eq1:base-construction-gluing} 
				\sup_{t\in[0,\infty)} \norm{f_t}_{e^{-2t}}\le K
			\end{equation} 
			and
			\begin{equation}
				\label{eq2:base-construction-gluing} 
				\norm{f_t - f_s} \le K(\abs{t-s}+e^{-2t}+e^{-2s}), \qquad  t,s\ge 0.
			\end{equation}
			Then, the map $F\colon X\to Y$ defined by 
			\begin{equation*}
				F(x) := \begin{cases}
					f_{\log\norm{x}}(x) \quad \textrm{ if } \norm{x}\ge 1, \\
					f_{0}(x) \qquad  \quad\textrm{ if } 0\le \norm{x} \le 1,\\
				\end{cases} 
			\end{equation*}
			is coarse-Lipschitz.
			
			Moreover, if $(f_{t})_{t\in[0,\infty)}$ is in $\mathcal {HU}(X,Y)$ and the map $t\mapsto f_t$ is continuous, then $F$ is uniformly continuous.
		\end{prop}
		
		\begin{proof} 
			For convenience, let
			$$\log^+(t) := \begin{cases}
				\log(t) \quad \textrm{ if }t\ge 1, \\
				0 \qquad  \quad\textrm{ if } 0\le t \le 1,
			\end{cases} 
			$$
			and observe that $F(x) = f_{\log^+\norm{x}}(x)$. We first show that $F$ is coarse-Lipschitz.
			
			Recall that $\norm{f}\le \norm{f}_\eps$ for all $f \in \cH(X,Y)$ and all $\eps \in (0,1]$. So, we immediately have $\norm{F(x)}\le K\norm{x}$ for all $x\in X$, so that $F$ is continuous at the origin. \\
			Now, suppose $\norm{x}\ge\norm{z}>0$. Then,
			\begin{align}
				\norm{F(x) - F(z)} & = \norm{f_{\log^+\norm{x}}(x) - f_{\log^+\norm{z}}(z)} \notag\\
				& \le \norm{f_{\log^+\norm{x}}(x) - f_{\log^+\norm{x}}(z)} + \norm{f_{\log^+\norm{x}}(z) - f_{\log^+\norm{z}}(z)} \label{eq3:base-construction-gluing}.
			\end{align}
			To estimate the first term in \eqref{eq3:base-construction-gluing}, we distinguish two cases. If $\norm{x}\ge 1$, then
			\begin{align*}
				\norm{f_{\log^+\norm{x}}(x) - f_{\log^+\norm{x}}(z)} & = \norm{f_{\log\norm{x}}(x) - f_{\log\norm{x}}(z)}\\
				&\le \norm{(f_{\log\norm{x}})}_{\norm{x}^{-2}}\max\{\norm{x-z},\norm{x}^{-1}\}\\
				&\stackrel{\eqref{eq1:base-construction-gluing}}{\le} K\max\{\norm{x-z},\norm{x}^{-1}\}.
			\end{align*}
			On the other hand, if $\norm{x}\le 1$ we have
			\begin{equation*}
				\norm{f_{\log^+\norm{x}}(x) - f_{\log^+\norm{x}}(z)} = \norm{f_0(x) - f_0(z)}\stackrel{\eqref{eq1:base-construction-gluing}}{\le} K\max\{\norm{x-z},\norm{x}\}.
			\end{equation*} 
			Therefore, in general, we have
			\begin{equation}
				\label{eq4:base-construction-gluing}
				\norm{f_{\log^+\norm{x}}(x) - f_{\log^+\norm{x}}(z)}\le K\max\{\norm{x-z},\min\{\norm{x},\norm{x}^{-1}\}\}.
			\end{equation}
			To estimate the second term in \eqref{eq3:base-construction-gluing}, note that it follows from \eqref{eq2:base-construction-gluing} and an elementary cases analysis that
			\begin{equation*}
				\norm{f_{\log^+\norm{z}}(z) - f_{\log^+\norm{x}}(z)}\le K\norm{z}\left(\log\frac{\norm{x}}{\norm{z}}+\min\{1,\norm{x}^{-2}\} + \min\{1,\norm{z}^{-2}\}\right).
			\end{equation*}
			The concavity of the $\log$ function gives us 
			\begin{equation*}
				\norm{z}\log\frac{\norm{x}}{\norm{z}}\le \norm{x}-\norm{z}\le \norm{x-z}
			\end{equation*} 
			from which it follows easily that 
			\begin{equation}
				\label{eq5:base-construction-gluing}
				\norm{f_{\log^+\norm{z}}(z) - f_{\log^+\norm{x}}(z)} \le K\norm{x-z} + K\min\{\norm{x},\norm{x}^{-1}\} + K\min\{\norm{z},\norm{z}^{-1}\}.
			\end{equation}
			Plugging \eqref{eq4:base-construction-gluing} and \eqref{eq5:base-construction-gluing} in \eqref{eq3:base-construction-gluing} we have 
			\begin{equation}
				\label{eq6:base-construction-gluing}
				\norm{F(x) - F(z)}\le 2K\norm{x-z} + 2K\min\{\norm{x},\norm{x}^{-1}\} + K\min\{\norm{z},\norm{z}^{-1}\}.
			\end{equation} 
			Yet another elementary cases analysis reveals that
			\begin{equation*}
				\norm{F(x) - F(z)}\le 2K\norm{x-z} + 3K,
			\end{equation*} 
			thus showing that $F$ is coarse-Lipschitz.
			
			It remains to show that $F$ is uniformly continuous under the additional assumptions of the uniform continuity of each $f_{t}$ and the continuity of the map $t\mapsto f_t$ on $B_X$. Given $\eps>0$ we need to find a $\delta>0$ such that $\norm{F(x) - F(z)} \le \vep$ whenever $\norm{x-z}\le \delta$. Without loss of generality, we will assume in the sequel that $\norm{x}\ge \norm{z}$. For points that are far away from the origin but close to each other, we just need to call upon \eqref{eq6:base-construction-gluing}. Indeed, if $\norm{x}\ge \norm{z}\ge a$ and $\norm{x-z}<\delta_0$, then
			\begin{equation*}
				\norm{F(x) - F(z)}\stackrel{\eqref{eq6:base-construction-gluing}}{\le} 2K\delta_0 + 3Ka^{-1},
			\end{equation*}
			and it suffices to choose $\delta_0>0$ so that $2K\delta_0<\frac{\eps}{2}$ and $a>1$ so that $3\frac{K}{a}<\frac{\eps}{2}$ to conclude that $\norm{F(x) - F(z)}\le \vep$.
			It remains to take care of the case where $z$ is close to the origin. 
			
			Assume from now on that $\norm{z}\le a$ and let $b := \log(a+1)>1$. Noting that since $t\mapsto f_t$ is continuous, it is uniformly continuous on $[0,b]$, we can certainly pick an integer $N>b$ so large that, if $0\le \sigma,\tau\le b$ and $\abs{\sigma-\tau}\le \frac{b}N$, we have
			\begin{equation}
				\label{eq7:base-construction-gluing}
				\norm{f_{\sigma} - f_{\tau}}<\frac{\eps}{3e^b}.
			\end{equation}  
			Moreover, since each $f_t$ is uniformly continuous on $B_X$, we can pick $\delta_1\in(0,1)$ small enough, to be specified later, such that if $u,v\in B_X$ with $\norm{u - v}\le \delta_1$ we have for all $0\le k\le N$,
			\begin{equation}
				\label{eq8:base-construction-gluing}
				\norm{f_{kb/N}(u) - f_{kb/N}(v)}< \frac{\eps}{3e^b}.
			\end{equation}
			%with $\delta_1<\min\{\frac{b}{N},\delta_0\}$
			If $\norm{x-z}<\delta_2$ for some $\delta_2\in(0,1)$, to be chosen smaller later, then $\norm{x}\le a+\delta_2<a+1$ (since $\norm{z}\le a$), $\log^+\norm{x}\le b$ (by definition of $b$) and in particular $\norm{x}\le e^b$. We need to consider two cases. If $\norm{x}\le 1$, then
			\begin{equation*}
				\norm{F(x) - F(z)} = \norm{f_0(x) - f_0(z)} \le \frac{\eps}{3e^b}<\frac{\eps}{3}.
			\end{equation*}
			On the other hand, if $\norm{x}>1$, we write 
			\begin{align*}
				\norm{F(x) - F(z)} \le & \norm{f_{\log^+\norm{x}}(x) - f_{kb/N}(x)} + \norm{f_{kb/N}(x) - f_{kb/N}(z)} \\
				&+ \norm{f_{kb/N}(z) - f_{\log^+\norm{z}}(z)}.
			\end{align*}
			Since $\log^+$ is $1$-Lipschitz and we had chosen $\delta_2<\frac{b}{N}$, we can certainly find $0\le k\le N$ such that 
			\begin{equation*}
				\max\Big\{\Big|\log^+\norm{x}-\frac{kb}{N}\Big|, \Big|\log^+\norm{z}-\frac{kb}{N}\Big| \Big\}\le \frac{b}N.
			\end{equation*}
			Hence, it follows from \eqref{eq7:base-construction-gluing} that
			\begin{equation*}
				\max\left\{ \norm{f_{\log^+\norm{x}}(x) - f_{kb/N}(x)}, \norm{f_{\log^+\norm{z}}(z) - f_{kb/N}(z)}\right\}\le \frac{\eps}{3e^b}\norm{x}\le \frac{\eps}{3}.
			\end{equation*}
			For the remaining middle term, observe that by homogeneity, 
			\begin{equation*}
				f_{kb/N}(x) - f_{kb/N}(z) = \norm{x}\Big(f_{kb/N}\Big(\frac{x}{\norm{x}}\Big) - f_{kb/N}\Big(\frac{z}{\norm{x}}\Big)\Big).
			\end{equation*}
			Since Lemma \ref{lem:polar} guarantees that $\norm{\frac{x}{\norm{x}} - \frac{z}{\norm{x}}}\le 2\delta_2$, we can invoke \eqref{eq8:base-construction-gluing}, assuming as we may that $2\delta_2<\delta_1$, to conclude that
			\begin{equation*}
				\norm{f_{kb/N}(x) - f_{kb/N}(z)}\le \frac{\eps}{3e^b}\norm{x}\le \frac{\eps}{3}.
			\end{equation*}
			Therefore, if $\norm{z}\le a$ and $\norm{x-z}\le \delta_2$, for some $\delta_2\le \min\{\delta_1/2,b/N\}$, we have $\norm{F(x) - F(z)}<\eps$.
			Finally, if one chooses $\delta\le \min\{\delta_0,\delta_1/2,b/N\}$, then it follows from the above cases analysis that $\norm{F(x) - F(z)}<\eps$ whenever $\norm{x-z}<\delta$, ultimately showing that $F$ is uniformly continuous.
		\end{proof}

		By applying Proposition \ref{prop:base-construction-gluing} to the maps and their inverses, we obtain what will be a very useful criterion to construct coarse-Lipschitz or uniform homeomorphisms between two Banach spaces.
		
		\begin{theo}
			\label{theo:CLH-criterion}
			Let $X$ and $Y$ be Banach spaces.  In order that $X$ and
			$Y$ be coarse-Lipschitz homeomorphic, it is sufficient that there exist a
			constant $L\ge 1$ and maps $(f_t)_{t\in[0,\infty)}$ in $\cG(X,Y)$  such
			that:
			\begin{equation}
				\label{eq1:CLH-citerion}
				\sup_{t\in [0,\infty)} [[f_t]]_{e^{-2t}}\le L,
			\end{equation}
			and
			\begin{equation}
				\label{eq2:CLH-citerion}
				\Delta(f_t,f_s)\le L(\abs{t-s} + e^{-2t} + e^{-2s}), \qquad t,s\ge 0.
			\end{equation}
			If in addition the collection $(f_t)_{t\in[0,\infty)}$ is in $\mathcal{GU}(X,Y)$
			and the maps $t\mapsto f_t$ and $t\mapsto f^{-1}_t$ are
			continuous, then $X$ and $Y$ are uniformly homeomorphic.
		\end{theo}
		
		\begin{proof}  
			It is a consequence of the proof of Proposition \ref{prop:close<->ALE-spheres} that 
			$\widehat f_t$ obeys the same conditions with the
			constant $L$ replaced by $2L^2+1.$ Hence, we may assume that $f_t$
			is norm preserving for each $t$.  Then, Proposition
			\ref{prop:base-construction-gluing} allows to build a coarse-Lipschitz map
			$F\colon X\to Y$, by setting $F(x) := f_{\log^+\norm{x}}(x)$. Since the $f_t$ are norm preserving bijections, we deduce that $F$ is a bijection so that $F^{-1}(y)=f_{\log^+\norm{y}}^{-1}(y)$, for $y\in Y$. It follows from the same lemma applied to $f_t^{-1}$ that $F^{-1}$ is also  coarse-Lipschitz.  It also follows that if each $f_t\in\mathcal {GU}(X,Y)$ and the maps $t\mapsto
			f_t,t\mapsto f^{-1}_t$ are continuous, then $F$ and $F^{-1}$ are
			uniformly continuous.
		\end{proof}

		Theorem \ref{theo:CLH-criterion} is a key ingredient towards a general result, due to Kalton, about the construction of coarse-Lipschitz and uniform homeomorphisms. This result is a generalization and abstraction of Ribe's construction as can be found in \cite{BenyaminiLindenstrauss2000}.
		
		\begin{theo}
			\label{theo:Kalton-general-construction}  
			Let $X$, $Y$ and $Z$ be Banach spaces such that $Y$ and $Y\oplus Z$ are linearly isomorphic.
			\begin{enumerate}[(i)]
				\item If $X$ and $X\oplus Y$ are close, then $X^2$
				and $X^2\oplus Z$ are coarse-Lipschitz homeomorphic.
				\item If $X$ and $X\oplus Y$ are uniformly close, then $X^2$
				and $X^2\oplus Z$ are uniformly homeomorphic.
			\end{enumerate}
		\end{theo}
		
		The assumptions of Theorem \ref{theo:Kalton-general-construction} might not seem very intuitive at this point, so before we delve into the lengthy and technical proof of Theorem \ref{theo:Kalton-general-construction} we first derive several more natural corollaries. 
		
		The first one is a very clean result, but under a rather strong structural assumption on the spaces.
		
		\begin{coro}
			\label{cor:Kalton-Ribe-squares}
			Let $X$ and $Y$ be Banach spaces that are both linearly isomorphic to their squares.
			\begin{enumerate}[(i)]
				\item If $X$ and $Y$ are close, then $X$ is coarse-Lipschitz homeomorphic to $Y$.
				\item If $X$ and $Y$ are uniformly close, then $X$ is uniformly homeomorphic to $Y$.
			\end{enumerate}
		\end{coro}
		
		\begin{proof} 
			We only detail the proof of $(ii)$. We have that $X$ is uniformly close to $Y$ and hence
			$X\oplus Y$ is uniformly close to $X\oplus X$ and thus to $X$. If
			we take $Z=Y$ in Theorem \ref{theo:Kalton-general-construction} we obtain that
			$X^2$ is uniformly homeomorphic to $X^2\oplus Y$ and therefore $X$ is
			uniformly homeomorphic to $X\oplus Y$. Since the roles of $X$ and $Y$ are symmetric, we also have that $Y$ is uniformly homeomorphic to $X\oplus Y$ and therefore to $X$. 
		\end{proof}
		
		There is no additional structural assumption in the second corollary, but the conclusion is weaker.
		
		\begin{coro}
			\label{cor:Kalton-Ribe-lp-sums}
			Let $X$ and $Y$ be Banach spaces and $1\le p\le \infty$. 
			\begin{enumerate}[(i)]
				\item If $X$ and $Y$ are close, then $\ell_p(X)$ and $\ell_p(Y)$ are
				coarse-Lipschitz homeomorphic.
				\item If $X$ and $Y$ are uniformly close, then $\ell_p(X)$ and $\ell_p(Y)$ are
				uniformly homeomorphic.
			\end{enumerate}
			The case $p=\infty$ has to be interpreted as a $\co$-sum.
		\end{coro}
		
		\begin{proof} 
			We only treat the uniform case for $1\le p<\infty$. By Lemma \ref{prop:closeness-ellp-sums} $\ell_p(X)$ and $\ell_p(Y)$ are uniformly close. The result then follows from Corollary \ref{cor:Kalton-Ribe-squares}, as $\ell_p(X)$ and $\ell_p(Y)$ are linearly isometric to their squares. 
		\end{proof}
		
		In the next corollary, there is an additional structural assumption on only one of the spaces involved. 
		
		\begin{coro}
			\label{cor:Kalton-Ribe-lp-sum-square}
			Let $X$ and $Y$ be Banach spaces and suppose that $X$ is linearly isomorphic to $\ell_p(X)$ with $1\le p<\infty$, or to $c_0(X)$. 
			\begin{enumerate}[(i)]
				\item If $X$ and $Y$ are close, then $X$ is coarse-Lipschitz homeomorphic to $Y^2$.
				\item If $X$ and $Y$ are uniformly close, then $X$ uniformly homeomorphic to $Y^2$.
			\end{enumerate}
		\end{coro}
		
		\begin{proof} 
			For the uniform case with $1\le p<\infty$, by item $(ii)$ in Corollary \ref{cor:Kalton-Ribe-lp-sums}, $\ell_p(X)$ is uniformly homeomorphic to $\ell_p(Y)$. By the additional assumption on $X$, it follows that $X$ is uniformly homeomorphic to $\ell_p(Y)$ and hence to $\ell_p(Y) \oplus Y^2$, which in turn is uniformly homeomorphic to $X \oplus Y^2$. On the other hand, since $X$ is linearly isomorphic to $\ell_p(X)$, it is also isomorphic to $X^2$. So, $Y$ is uniformly close to $X^2$ and thus to $X\oplus Y$. So, we can apply Theorem \ref{theo:Kalton-general-construction}, exchanging the roles of $X$ and $Y$ and with $Z:=X$, to get that $Y^2$ is uniformly homeomorphic to $Y^2\oplus X$ and therefore to $X$. 
			
			The other cases can be treated similarly and are left to the reader.
		\end{proof}
		
		The final corollary is about the stability of uniform or coarse-Lipschitz homeomorphisms under $\ell_p$-sums.
		
		\begin{coro}
			\label{cor:lp-sums-CLH}
			Let $X$ and $Y$ be Banach spaces and $1\le p\le \infty$.
			\begin{enumerate}[(i)]
				\item If $X$ and $Y$ are coarse-Lipschitz homeomorphic, then $\ell_p(X)$ and $\ell_p(Y)$ are coarse-Lipschitz homeomorphic.
				\item If $X$ and $Y$ are uniformly homeomorphic, then $\ell_p(X)$ and $\ell_p(Y)$ are uniformly homeomorphic.
			\end{enumerate}
			The case $p=\infty$ has to be interpreted as a $\co$-sum.
		\end{coro}
		
		\begin{proof} 
			For the coarse-Lipschitz case with $1\le p<\infty$, it follows from Theorem \ref{cor:Nahum} that $X\oplus\bR$ and $Y\oplus\bR$ are close. Invoking assertion $(i)$ in Corollary \ref{cor:Kalton-Ribe-lp-sums}, we deduce that $\ell_p(X\oplus\bR)$
			and $\ell_p(Y\oplus\bR)$ are coarse-Lipschitz homeomorphic. But
			these spaces are isomorphic to $\ell_p(X)$ and $\ell_p(Y)$ respectively, and we are done. The other cases have a similar proof.
		\end{proof}
		
		The rest of this section is devoted to the proof of the uniform case of Theorem \ref{theo:Kalton-general-construction}. First note that under the assumption of Theorem \ref{theo:Kalton-general-construction} one can quickly conclude that $X$ and $X\oplus Z$ are uniformly close and hence there is a sequence of maps $(f_n)_{n\ge 0}$ in $\mathcal{GU}(X,X\oplus Z)$ such $\sup_{n\ge 0} [[f_n]]_{e^{-2n}} <\infty$. Each of these maps can only take care of a single scale and of pairs of points that are roughly at the same distance from the origins and they are not enough to construct a uniform homeomorphism between $X$ and $X\oplus Z$. In Proposition \ref{prop:base-construction-gluing}, we saw that we need an additional regularity condition to be able to glue these maps together into a uniform homeomorphism. The main point of the proof of Theorem \ref{theo:Kalton-general-construction} is to find a way to link two consecutive maps in the sequence by a continuous path of maps while maintaining condition \eqref{eq1:base-construction-gluing} and guaranteeing condition \eqref{eq2:base-construction-gluing} in Proposition \ref{prop:base-construction-gluing}. As we will see, this can be done at the expense of adding a copy of $X$. 
		
		During the linking procedure, we will need to connect certain operators to the identity. If $W$ is any Banach space, the \emph{interchange operator} is the operator $J_W\colon W\oplus W\to W\oplus W$ defined for all $w_1,w_2\in W$ by
		\begin{equation*}
			J(w_1,w_2) := (w_2,-w_1).
		\end{equation*}
		Observe that $J_W\circ J_W = - I_W$ and hence $J_W$ is invertible and with $J_W^{-1}=-J_W$. We will treat, as we may, all direct sums as $\ell_\infty$-sums and thus $J$ is an isometry. We formulate without proof an elementary proposition that we will use repeatedly for the linking procedure. 
		%For convenience, we will treat all direct sums as $\ell_\infty$ sums in the sequel of this subsection. 
		\begin{lemm}[Linking Principle]
			\label{lem:link}
			Let $W$ be any Banach space. The collection of linear maps $(\Psi_\theta)_{\theta\in [0,\pi/2]}$ on $W\oplus W$ defined by,
			\begin{equation*}
				\Psi_{\theta}(w_1,w_2) := (\cos(\theta)w_1 + \sin(\theta)w_2, \cos(\theta)w_2 - \sin(\theta)w_1),
			\end{equation*}  
			satisfies the following properties:
			\begin{enumerate}
				\item $\psi_0=I_W$ and $\psi_{\pi/2}=J_W$,
				\item For all $\theta\in[0,\pi/2]$, $\Psi_{\theta}$ is invertible and $\Psi_{\theta}^{-1}=-J_W\circ \Psi_{\frac{\pi}{2} - \theta}$,
				\item $\psi^{-1}_0=I_W$ and $\psi^{-1}_{\pi/2}=-J_W$.
				\item $\max\{\norm{\Psi_{\theta}},\norm{\Psi_{\theta}^{-1}}\}\le 2$ for all $\theta \in[0,\pi/2]$,
				\item the maps $\theta\mapsto \Psi_{\theta}$ and $\theta\mapsto\Psi_{\theta}^{-1}$
				from $[0,\pi/2] \to B(W\oplus W)$ are $2$-Lipschitz,
			\end{enumerate}
		\end{lemm}
		
		%\begin{rema}
		%There is a similar family of maps that joins $I_W$ to $-J_W=J_W^{-1}$.
		%\end{rema}
		
		When we link homogeneous maps by precomposing and postcomposing by maps as in the linking principle, we need to make sure that we can control the norm of the difference of the new homogeneous maps constructed. This technical detail is taken care of in the next lemma.
		
		\begin{lemm}
			\label{lem:diff} 
			Let $X,Y,Z$ and $W$ be Banach spaces, $f\in\mathcal H(X,Y)$ and $g\in\mathcal H(Z,W)$. Assume $h_1,h_2\in \mathcal H(Y,Z)$. Then,
			\begin{equation*}
				\norm{g\circ h_1\circ f - g\circ h_2\circ f} \le \max\{\norm{h_1},\norm{h_2}\}\norm{f}\,\omega_{g_{\restriction_{B_X}}}\Big(\frac{\norm{h_1-h_2}}{\max\{\norm{h_1},\norm{h_2}\}}\Big),
			\end{equation*}
			where $\omega_{g_{\restriction_{B_X}}}$ is the expansion modulus of $g_{\restriction_{B_X}}$.
			
			In particular, if $\eps \in (0,1]$, 
			\begin{equation*}
				\norm{g\circ h_1\circ f - g\circ h_2\circ f} \le \norm{g}_\eps\norm{f}\big(\norm{h_1-h_2} + \vep\max\{\norm{h_1},\norm{h_2}\}\big).
			\end{equation*}
		\end{lemm}
		
		\begin{proof} 
			If $\norm{x}=1$, then $\norm{h_1(f(x)) - h_2(f(x))}\le \norm{h_1-h_2}\norm{f}$.  
			Let $K := \max\{\norm{h_1},\norm{h_2}\}$ and $\omega:=\omega_{g_{\restriction_{B_X}}}$. Since $\frac{(h_1\circ f)(x)}{K\norm{f}}$ and $\frac{(h_1\circ f)(x)}{K\norm{f}}$ are in $B_X$, by homogeneity we have 
			\begin{equation*}
				\norm{g \circ h_1 \circ f(x) - g \circ h_2\circ f(x)}\le K\norm{f}\, \omega\Big(\frac{\norm{h_1-h_2}}{K}\Big).
			\end{equation*}
			Let $\eps \in (0,1]$. For $x,x' \in B_X$ with $\norm{x-x'}\le \frac{\norm{h_1-h_2}}{K}$, we have that 
			\begin{equation*}
				\norm{g(x)-g(x')}\le \norm{g}_\eps\max\Big\{\frac{\norm{h_1-h_2}}{K},\eps\Big\}\le \norm{g}_\eps\Big(\frac{\norm{h_1-h_2}}{K}+\eps\Big),
			\end{equation*}
			and the last inequality of our statement follows by homogeneity. 
		\end{proof}
		
		We will only detail the uniform case of the proof of Theorem \ref{theo:Kalton-general-construction} as the coarse-Lipschitz case is easier.
		
		\begin{proof}[Proof of Theorem \ref{theo:Kalton-general-construction}]
			Let us fix a constant $L\ge 1$ large enough so that for every $\vep>0$ we can find $g\in \mathcal {GU}(X,X\oplus Y)$ with $[[g]]_\vep\le L$ and such that there is a linear bijection $T\colon Y\to Y\oplus Z$ with $\max\{\norm{T},\norm{T^{-1}}\}\le L$.  We write $g(x) := (g_X(x),g_Y(x))$, for $x\in X$ and $T(y) :=(T_Y(y),T_Z(y))$, for $y \in Y$. Then, we define $\phi_g \colon X \to X \oplus Z$ by
			\begin{equation*}
				\phi_g(x) := \Big(g^{-1}\big(g_X(x), T_Y(g_Y(x)\big), T_Z(g_Y(x))\Big).
			\end{equation*}
			Clearly $\phi_g \in \mathcal{GU}(X,X\oplus Z)$ and we can express $\phi_g$ as the following composition:
			\begin{equation*}
				X\buildrel g\over\longrightarrow X\oplus Y \buildrel I_X \oplus T\over\longrightarrow
				X\oplus Y\oplus Z \buildrel {g^{-1}\oplus I_Z}\over \longrightarrow X\oplus Z.
			\end{equation*}
			Since all the sums are $\ell_\infty$-sums, it is not difficult to verify that $\max\{[[g^{-1}\oplus I_Z]]_\vep,[[I_X\oplus T]]_\vep\} \le L$ and it follows from Lemma \ref{lem:composition} that $[[\phi_g]]_\vep\le L^3$. This shows that $X$ and $X\oplus Z$ are uniformly close.
			
			The heart of the matter here is to find a way to link $\phi_g$ to $\phi_h$ where $g$ and $h$ are two maps for which the discussion above applies. This can be done but modulo a tensorization by the identity on $X$. 
			
			\begin{claim}
				\label{clai:linkup}  
				If $g,h\colon X\to X\oplus Y$ are such that $g,h\in\mathcal {GU}(X,X\oplus Y)$ and $\max\{[[g]]_\eps,[[h]]_\eps\}\le L$, for some $L>0$ and $\vep>0$, then there is a collection of maps $(f_t)_{t\in [0,1]}\in \mathcal {GU}(X^2,X^2\oplus Z)$ such that:
				\begin{enumerate}[(i)]
					\item $t\mapsto f_t$ and $t\mapsto f^{-1}_t$ are continuous, %for the metric $\Delta$,
					\item $[[f_t]]_\vep\le 2L^3$,
					\item $\Delta(f_t,f_s)\le 4\pi L^3\abs{t-s} + 8L^3\vep$ for all $0\le s,t\le 1$,
					%      \item $\Delta(f_t,f_s)\le 10L^3\abs{t-s} + 6L^3\vep$, for all $0\le s,t\le 1$,
					\item $f_0(x_1,x_2) = (x_1,\phi_g(x_2))$ and $f_1(x_1,x_2) = (x_1,\phi_h(x_2))$.
				\end{enumerate} 
			\end{claim}
			
			Assuming the claim for a moment, we can quickly conclude the proof of the theorem as follows. According to the initial discussion above, for each integer $n\ge 0$, we can find a map $\phi_{g_n}\in\mathcal{GU}(X,X\oplus Z)$ with
			$[[\phi_{g_n}]]_{e^{-2(n+1)}}\le L^3$. Now, with the help of Claim \ref{clai:linkup}, we can build, for each $n\ge 0$, a family $(f_t)_{t\in [n,n+1]}\subset \mathcal{GU}(X^2,X^2\oplus Z)$ linking $f_n=I_X\oplus \phi_{g_n}$ to $f_{n+1}=I_X\oplus \phi_{g_{n+1}}$ and such that for all $s,t \in [n,n+1]$, 
			\begin{equation*}
				[[f_t]]_{e^{-2t}}\le 2L^3,
			\end{equation*}
			and
			\begin{equation*}
				\Delta(f_t,f_s)\le 10L^3|t-s|+6L^3e^{-2(n+1)},
			\end{equation*}
			and so that $t\mapsto f_t$ and $t \mapsto f_t^{-1}$ are continuous on $[n,n+1]$. Collecting all these maps, one obtains a collection of maps $(f_t)_{t\in [0,\infty)}\subset \mathcal{GU}(X^2,X^2\oplus Z)$, satisfying condition \eqref{eq1:CLH-citerion} in Theorem \ref{theo:CLH-criterion} and such that the maps $t\mapsto f_t$ and $t \mapsto f_t^{-1}$ are continuous on $[0,\infty)$. To conclude that $X^2$ and $X^2\oplus Z$ are uniformly homeomorphic, one is left with verifying that condition \eqref{eq2:CLH-citerion} in Theorem \ref{theo:CLH-criterion} holds. To do so, fix $0\le s\le t<\infty$ and let $m\le n$ be integers so that $m\le s<m+1$ and $n\le t<n+1$. By triangle inequality, we have 
			\begin{equation*}
				\Delta(f_t,f_s)\le 10L^3\abs{t-s} + 6L^3\sum_{i=m}^n e^{-2(i+1)},
				%\le 10L^3\abs{t-s} + 12L^3(e^{-2t}+e^{-2s}).$$
			\end{equation*}
			but since $\sum_{i=m}^n e^{-2(i+1)} \le 2e^{-2(m+1)} \le 2e^{-2s}$, it follows that 
			\begin{equation*}
				\Delta(f_t,f_s)\le 12L^3 (\abs{t-s} + e^{-2t}+e^{-2s}).
			\end{equation*}
			
			It remains to provide the tedious details for the proof of Claim \ref{clai:linkup} which relies on some crucial observations.
			The first observation is that $f_0$ can be obtained as the composition $f_0 = F_3\circ F_2\circ F_1$ where 
			\begin{equation*}
				X^2\buildrel F_1\over\longrightarrow (X\oplus Y)^2 \buildrel F_2 \over\longrightarrow
				(X\oplus Y)^2\oplus Z \buildrel F_3\over \longrightarrow X^2\oplus Z,
			\end{equation*}
			\begin{equation*}
				F_1(x_1,x_2) := (h(x_1),g(x_2)), %=(h_X(x_1),h_Y(x_1),g_X(x_2),g_Y(x_2)),
			\end{equation*}
			\begin{equation*}
				F_2(x_1,y_1,x_2,y_2) := (x_1,y_1,x_2,Ty_2), %=(x_1,y_1,x_2,T_Yy_2,T_Zy_2),
			\end{equation*}
			and
			\begin{equation*}
				F_3(x_1,y_1,x_2,y_2,z)=(h^{-1}(x_1,y_1),g^{-1}(x_2,y_2),z).
			\end{equation*}
			The second observation is that $f_1=J_4\circ F_3\circ J_3\circ F_2\circ J_2\circ F_1\circ J_1$ where the $J_i$ are the following interchange operators, whose domains should be clear from the notation:
			\begin{equation*}
				J_1(x_1,x_2)=(x_2,-x_1),  
			\end{equation*}
			\begin{equation*}
				J_2(x_1,y_1,x_2,y_2)=(x_1,-y_2,x_2,y_1),  
			\end{equation*}
			\begin{equation*}
				J_3(x_1,y_1,x_2,y_2,z)=(x_1,y_2,x_2, -y_1,z) 
			\end{equation*}
			\begin{equation*}
				J_4(x_1,x_2,z)=(-x_2,x_1,z)
			\end{equation*}
			It is not difficult to see that we can link the interchange operators to the corresponding identities using the linking principle. Also, it is clear that $\max_{1\le i\le 3} [[F_i]]_\vep \le L$ and that $J_1,J_2,J_3$ and $J_4$ are isometries. To link $f_0$ to $f_1$ we will go one step at a time. First we will link $f_0$ with 
			\begin{equation*}
				f_{1/4} := F_3\circ F_2\circ F_1 \circ J_1,   
			\end{equation*}
			then $f_{1/4}$ with 
			\begin{equation*}
				f_{1/2} := F_3\circ F_2\circ J_2\circ F_1\circ J_1,
			\end{equation*}
			then $f_{1/2}$ with
			\begin{equation*}
				f_{3/4} := F_3\circ J_3\circ F_2\circ J_2\circ F_1\circ J_1,
			\end{equation*}
			and finally $f_{3/4}$ with $f_{1}$.
			Next, to join $f_{j-1/4}$ to $f_{j/4}$, for $j=1,\dots,4$, we just have to use the linking principle (Lemma \ref{lem:link}) to go from the appropriate identity operator to $J_j$. Composing with a shift and a rescaling from $[0,\frac{\pi}{2}]$ to $[\frac{j-1}{4},\frac{j}{4}]$, we can construct $f_t$ for $0\le t\le 1$. On the interval $[\frac{j-1}{4},\frac{j}{4}]$, $f_t$ takes the form $f_t = G_j\Phi_j(t)H_j$ where $\Phi_j(t)$ are linear maps verifying the estimates
			\begin{equation*}
				\max\big\{\norm{\Phi_j(t) - \Phi_j(s)},\norm{\Phi_j(t)^{-1} - \Phi_j(s)^{-1}}\big\}\le 4\pi\abs{t-s}% 2\sqrt{2}\pi\abs{t-s}\le 10\abs{t-s}
			\end{equation*}
			and
			\begin{equation*}
				\max\big\{\norm{\Phi_j(t)},\norm{\Phi_j(t)^{-1}}\big\}\le 2. % \sqrt 2.
			\end{equation*}
			For instance, $G_1=F_3\circ F_3 \circ F_1$, $H_1=I$ and $\Phi_1(t)=\Psi_{2\pi t}$,
			while $G_4=I$, $H_4=F_3\circ J_3\circ F_2\circ J_2\circ F_1\circ J_1$ and $\Phi_4(t)=(\Psi_{2\pi(t-\frac34)},I_Z)$. 
			
			Furthermore, we have $[[G_j]]_\vep [[H_j]]_\vep\le L^3$ and this implies that $[[f_t]]_\vep\le 2L^3$. It follows from Lemma \ref{lem:diff} that, for every $s,t\in [\frac{j-1}{4},\frac{j}{4}]$,
			\begin{align*}
				\norm{f_t - f_s} & \le \norm{G_j}_\vep\norm{H_j} \big(\norm{\Phi_j(t) - \Phi_j(s)} + \vep \max\big\{\norm{\Phi_j(t)},\norm{\Phi_j(s)}\big\}\big) \\
				& \le  L^3\big(4\pi \abs{t-s} + 2\vep\big), 
			\end{align*}
			and
			\begin{align*}
				\|f^{-1}_t - f^{-1}_s\| & \le \norm{G_j^{-1}}_\vep\norm{H_j^{-1}} \big(\norm{\Phi_j^{-1}(t) - \Phi_j^{-1}(s)} + \vep \max\big\{\norm{\Phi_j^{-1}(t)},\norm{\Phi_j^{-1}(s)}\big\}\big) \\
				& \le L^3 \big(4\pi \abs{t-s} + 2\vep\big), 
			\end{align*}
			and hence 
			\begin{equation*}
				\Delta(f_t,f_s)\le 4\pi L^3\abs{t-s} + 2L^3\vep.
			\end{equation*}  
			It then follows from the triangle inequality that for all $s,t\in[0,1]$, 
			\begin{equation*}
				\Delta(f_t,f_s)\le 4\pi L^3\abs{t-s} + 8L^3\vep.
			\end{equation*}
			Finally, since $g,h \in \mathcal {GU}(X,X\oplus Y)$, $(f_t)_{t\in [0,1]}\in \mathcal {GU}(X^2,X^2\oplus Z)$ and Lemma \ref{lem:diff} implies the continuity of the maps $t\mapsto f_t$ and $t\mapsto f_t^{-1}$.
		\end{proof}

		\subsection{Ribe's pair of uniformly equivalent but nonisomorphic separable Banach spaces}
		\label{sec:Ribe-example}
		We now explain how to derive Ribe's original example (Theorem \ref{theo:RibeExample}) of two uniformly equivalent but nonisomorphic Banach spaces from the general results of the previous section. 
		
		Let $(p_n)_{n=1}^\infty$ be a sequence in $(1,2)$ which is converging to $1$. Since it is a standard fact that $(\sum_{n=1}^\infty \ell_{p_n})_{\ell_2}$ and $(\sum_{n=1}^\infty \ell_{p_n})_{\ell_2}\oplus_2 \ell_1$ are isomorphic to their squares, in order to show that these two spaces are uniformly homeomorphic it is sufficient to show, according to Corollary \ref{cor:Kalton-Ribe-squares}, that $(\sum_{n=1}^\infty \ell_{p_n})_{\ell_2}$ and $(\sum_{n=1}^\infty \ell_{p_n})_{\ell_2}\oplus_2 \ell_1$ are uniformly close. Observing that $(\sum_{n=1}^\infty \ell_{p_n})_{\ell_2}$ is linearly isomorphic to $(\sum_{n=1}^\infty \ell_{p_n})_{\ell_2} \oplus_2 \ell_{p_{n_0}}$ with distortion 
		at most $\sqrt 2$, we simply need to show that $(\sum_{n=1}^\infty \ell_{p_n})_{\ell_2} \oplus_2 \ell_{p_{n_0}}$ and $(\sum_{n=1}^\infty \ell_{p_n})_{\ell_2}\oplus_2 \ell_1$ are uniformly close. According to Proposition \ref{prop:close<->ALE-spheres}, we will be done if we can show that $(\sum_{n=1}^\infty \ell_{p_n})_{\ell_2} \oplus_2 \ell_{p_{n_0}}$ and $(\sum_{n=1}^\infty \ell_{p_n})_{\ell_2}\oplus_2 \ell_1$ have uniformly almost Lipschitz equivalent spheres. For this last task, we need to know a little bit more about the norm-preserving and positively homogeneous Mazur maps, which we recall are the maps $\bar{M}_{p,q}\colon \ell_p \to \ell_q$ given by 
		\begin{equation*}
			\bar{M}_{p,q}(x):= \norm{x}_p^{1-p/q}\sign(x)\abs{x}^{p/q}.
		\end{equation*} 
		
		\begin{prop}
			\label{prop:Mazur-maps}
			For all $1\le q\le p<\infty$, $\bar{M}_{p,q}$ is a norm preserving, positively homogeneous map which is $\frac{p}{q}$-Lipschitz and such that the modulus of continuity $\omega$ of the restriction of its inverse $\bar{M}_{q,p}$ to the unit ball of $\ell_q$ satisfies $\omega(t)\le c2^{p/q}t^{q/p}$, for $t\in [0,2]$, where $c\ge 1$ is a universal constant.  More precisely:
			\begin{enumerate}
				\item $\norm{\bar{M}_{p,q}(x) - \bar{M}_{p,q}(x)} \le \frac{3p}{q}\norm{x-y}$ for all $x,y\in \ell_p$,
				\item $\norm{\bar{M}_{q,p}(x) - \bar{M}_{q,p}(x)} \le 4\norm{x-y}^{q/p}$ for all $x,y\in B_{\ell_p}$.
			\end{enumerate}
		\end{prop}
		We refer the reader to Theorem 9.1 in \cite{BenyaminiLindenstrauss2000} for the proof of this classical result, which is valid for any $L_p$-spaces over a measure space.
		
		So let us fix $\vep >0$ and $\eta\in(0,1)$ to be chosen small enough later. It follows from elementary calculus that there exists $\alpha_\eta \in (0,1)$ such that
		\begin{equation*}
			\forall \alpha \in [\alpha_\eta,1]\ \ \forall t \in [0,2],\ \ t^\alpha\le t+\eta,
		\end{equation*}
		and pick $n_0\in \bN$ so that $\frac{1}{p_{n_0}}\ge \alpha_\eta$. Now, consider the map $\Phi\colon (\sum_{n=1}^\infty \ell_{p_n})_{\ell_2} \oplus_2 \ell_{p_{n_0}} \to (\sum_{n=1}^\infty \ell_{p_n})_{\ell_2}\oplus_2 \ell_1$ defined by 
		\begin{equation*}
			\Phi := (\sum_{n=1}^\infty I_{\ell_{p_n}})_{\ell_2}\oplus_2 \bar{M}_{p_{n_0},1}.
		\end{equation*}
		It follows from the properties of $\bar{M}_{p_{n_0},1}$ that the map $\Phi$ is a uniform homeomorphism between the spheres. The coarse-Lipschiz character of $\Phi$ comes from our choice of $p_{n_0}$. Indeed, since $p_{n_0}\le 2$ the map $\Phi$ is a $6$-Lipschitz. As for the inverse, note that for all $w,w'$ in the unit ball of $(\sum_{n=1}^\infty \ell_{p_n})_{\ell_2}\oplus_2 \ell_1$, we have
		\begin{align*}
			\|\Phi^{-1}(w)-\Phi^{-1}(w')\|^2 & = \norm{w - w'}^2 + \norm{\bar{M}^{-1}_{p_{n_0},1}(w_{n_0}) - \bar{M}^{-1}_{p_{n_0},1}(w'_{n_0})}^2 \\
			& \le \norm{w - w'}^2 + 4^2 \norm{w_{n_0} -w'_{n_0} } ^{2/p_{n_0}} \\
			& \le \norm{w - w'}^2 + 4^2 (\norm{w-w'}+\eta)^2\\
			& \le (1+ 16)\norm{w-w'}^2 + 16(2\eta\norm{w-w'}+\eta^2)\\
			& \le 17\norm{w-w'}^2 + 80\eta\\
			& \le 17\norm{w-w'}^2 + \vep^2,
		\end{align*}
		if for the last inequality we had chosen $\eta\le \frac{\vep^2}{80}$
		Therefore, there is a constant $L>0$ such that $\|\Phi^{-1}(w)-\Phi^{-1}(w')\|\le L\norm{w-w'} + \vep$ whenever $w,w'$ are in the unit sphere of $(\sum_{n=1}^\infty \ell_{p_n})_{\ell_2}\oplus_2 \ell_1$. This completes the proof that $(\sum_{n=1}^\infty \ell_{p_n})_{\ell_2}\oplus_2 \ell_{p_{n_0}}$ and $(\sum_{n=1}^\infty \ell_{p_n})_{\ell_2}\oplus_2 \ell_1$ have uniformly almost Lipschitz isomorphic spheres. 
		
		\subsection{Kalton's counterexample to the preservation of the AUS/AUC moduli under uniform equivalences}
		\label{sec:Kalton-Kp-spaces}
		
		In this section, we detail Kalton's examples of pairs of separable and reflexive Banach spaces showing that the moduli of asymptotic smoothness and convexity are not preserved under uniform homeomorphism. In view of the results from Chapter \ref{chapter:Gorelik}, the most stunning example is the existence of a  Banach space that is uniformly equivalent to a $p$-asymptotically uniformly smooth Banach space, but does not admit any equivalent $p$-asymptotically uniformly smooth norm. 
		
		First, let us describe the spaces involved in the examples. These spaces are all unconditional sums of finite-dimensional spaces. We will use the term \emph{Banach sequence space} for a Banach space $(X,\norm{\cdot}_X)$, of real-valued sequences with a norm so that the canonical basis $(e_n)_{n=1}^{\infty}$ is a 1-unconditional basis of $X$. Given $(G_n)_{n=1}^{\infty}$ a sequence of Banach spaces and $X$ a Banach sequence space, we define $(\sum_{n=1}^{\infty}G_n)_X$ to be the space of sequences $x := (x_n)_{n=1}^{\infty}$ such that $x_n\in G_n$ and
		$(\norm{x_n})_{n=1}^{\infty}\in X$ with the norm
		\begin{equation*}
			\norm{x} := \big\|(\norm{x_n}_{G_n})_{n=1}^{\infty}\big\|_X.
		\end{equation*}
		It follows from $1$-unconditionality that $\norm{\cdot}$ is a genuine norm.
		If moreover $(G_n)_{n=1}^{\infty}$ is a sequence of finite-dimensional normed spaces dense for the Banach-Mazur distance in all finite-dimensional normed spaces, then $\Big(\sum_{n=1}^{\infty}G_n\Big)_{\ell_p}$, for $1<p<\infty$, is the space called $C_p$ already introduced in Chapter \ref{chapter:AMP_II}, Section \ref{sec:Cp}. Recall that, up to isomorphism, $C_p$ does not depend on the choice of the dense sequence $(G_n)_{n=1}^{\infty}$.
		
		To describe the other spaces showing up in the examples, we need to recall some background on $p$-convex and $q$-concave Banach lattices, and we refer to \cite{LindenstraussTzafriri1979} for more details. Banach sequence spaces can be naturally canonically equipped with a lattice structure. In particular, 
		%for $x=\sum_{n=1}^\infty x(n)e_n \in X$ and $1\le p<\infty$, we denote $|x|:=\sum_{n=1}^\infty |x(n)|e_n$ and $|x|^p:=\sum_{n=1}^\infty |x(n)|^pe_n$ (note that the last two series also converge in $X$).
		for $x=(x_n)_{n=1}^\infty\in X$, $\abs{x} = (\abs{x_n})_{n=1}^\infty\in X$ and for $p\in (1,\infty)$, expressions such as $|x|^p$ or $\Big(\sum_{j=1}^n\abs{z_j}^p\Big)^{1/p}$, for $z_1,\dots,z_n \in X$ are elements in $X$ that are defined pointwise in Banach sequence spaces and can be formally defined in general Banach lattices using functional calculus. 
		
		A Banach lattice $X$ is said to be \emph{$p$-convex} if there is a constant $M>0$ such that for all $x_1,\dots,x_n\in X$,
		\begin{equation*}
			\Big\|\Big(\sum_{j=1}^n\abs{x_j}^p\Big)^{1/p}\Big\|_X \le M\Big(\sum_{j=1}^n\norm{x_j}_X^p\Big)^{1/p},
		\end{equation*} 
		and \emph{strictly $p$-convex} in the case where $M=1$.
		
		Banach lattices and, in particular, Banach sequence spaces, can be convexified. Recall that $\Tsi$ is the asymptotic-$\ell_1$ Tsirelson space (the dual of the original space of Tsirelson, which is asymptotically-$\co$). We refer to \cite{CasazzaShura1989} for a thorough study of Tsirelson spaces. $\Tsi$ is a Banach sequence space and for $1<p<\infty$ we denote by $\Tsi_p$ its $p$-convexification, i.e., $\Tsi_p$ is the vector space $\{ x\in T \colon \norm{\abs{x}^p}_T<\infty\}$ equipped with the norm $\norm{x}_{\Tsi_p}:= \norm{\abs{x}^p}_T^{1/p}$.
		Note that $\Tsi_p$ is strictly $p$-convex. If $(G_n)_{n=1}^{\infty}$ is a sequence of finite-dimensional normed spaces dense for the Banach-Mazur distance in all finite-dimensional normed spaces, then, for $1<p<\infty$, we let 
		\begin{equation*}
			K_p := \Big(\sum_{n=1}^{\infty}G_n\Big)_{\Tsi_p}.
		\end{equation*}
		The main result of this section is that $K_p\oplus K_p$ and $C_p$ are uniformly homeomorphic. 
		
		To prove this result, we need to discuss a bit of interpolation theory and the dual notion of $p$-convexity. This dual notion, called $q$-concavity, is defined by changing the inequality sign. For $q\in (1,\infty)$, a Banach lattice $X$ is said to be \emph{$q$-concave} if there is a constant $M>0$ such that for all $x_1,\dots,x_n\in X$,
		\begin{equation*}
			\Big(\sum_{j=1}^n\norm{x_j}_X^q\Big)^{1/q}\le M\Big\|\Big(\sum_{j=1}^n\abs{x_j}^q\Big)^{1/q}\Big\|_X,
		\end{equation*}
		and \emph{strictly $q$-concave} when $M=1$.
		
		\begin{rema} 
			Let $X$ be a Banach sequence space, $p\in (1,\infty)$ and $(F_n)_{n=1}^\infty$ a sequence of finite-dimensional normed spaces. It is not difficult to see that if $X$ is strictly $p$-convex (respectively strictly $p$-concave), then $\big(\sum_{n=1}^\infty F_n\big)_X$ is $p$-AUS (respectively $p$-AUC). We leave the proof of this fact as Exercise \ref{ex:pconvex-pAUS}. 
		\end{rema}
		
		By a duality argument (see \cite[Proposition 1.d.4.]{LindenstraussTzafriri1979}), $\Tsi_p^*$ is strictly $q$-concave where $q$ is the conjugate exponent. If a Banach sequence space $X$ is strictly $p$-convex and strictly $q$-concave where $1<p\le q<\infty$, then $X$ is uniformly convex and uniformly smooth (see \cite[Theorem 1.f.1.]{LindenstraussTzafriri1979}). In this case, $X^*$ can be identified with a Banach sequence space, which is strictly $q'$-convex and strictly $p'$-concave where $p',q'$ are conjugate to $p,q$. It should also be mentioned that if $X$ is only
		$p$-convex and $q$-concave for $1< p\le q< \infty$, then it can be equivalently renormed to be a strictly $p$-convex and strictly $q$-concave Banach sequence space (see \cite[Proposition 1.d.8.]{LindenstraussTzafriri1979}).
		
		The material from interpolation theory, which is needed in the sequel, can be found in \cite{BenyaminiLindenstrauss2000}. If $X_0$ and $X_1$ are two Banach sequence spaces, then for $\theta \in (0,1)$, the Calder\'{o}n space $X_\theta := X_0^{1-\theta}X_1^{\theta}$ is described as the space of
		all $x$ so that $\abs{x}= \abs{x_0}^{1-\theta}\abs{x_1}^{\theta}$ for some $x_0\in
		X_0,\ x_1\in X_1$. The space $X_\theta$ is equipped with the norm
		\begin{equation*}
			\norm{x}_{X_{\theta}} := \inf\Big\{\max\big\{\norm{x_0}_{X_0},\norm{x_1}_{X_1}\big\}\colon \abs{x}=\abs{x_0}^{1-\theta}\abs{x_1}^{\theta}\Big\}.
		\end{equation*}
		The infimum is attained uniquely whenever $X_0$ and $X_1$ are uniformly convex.
		
		\medskip
		The following theorem is a generalization of the Mazur map for interpolation pairs.
		%The construction of the counterexample is based on the existence of uniform homeomorphisms between unit spheres of general sequence spaces, that we state now. 
		
		\begin{theo}
			\label{theo:homeo-spheres} 
			Let $X_0$ and $X_1$ be two uniformly convex sequence spaces and let $\sigma,\theta \in (0,1)$. For $x :=(x_n)_{n=1}^\infty \in X_\sigma$, let $\abs{x_n}=a_n^{1-\sigma}b_n^\sigma$ be the unique optimal factorization with
			$a :=(a_n)_{n=1}^{\infty}\in X_0,\ b :=(b_n)_{n=1}^{\infty}\in X_1$,
			and $\norm{a}_{X_0}=\norm{b}_{X_1}=\norm{x}_{X_\sigma}$. If we let
			\begin{equation*}
				f_0(x) := \big(a_n^{1-\theta}b_n^{\theta}\sign(x_n)\big)_{n=1}^{\infty},
			\end{equation*}
			then $f_0$ is a norm preserving uniform homeomorphism between $B_{X_\sigma}$ and  $B_{X_\theta}$.
		\end{theo}
		
		We refer to \cite[Theorem 9.12]{BenyaminiLindenstrauss2000} for its proof. It can also be found in \cite{Daher1995} and turns out to be the same as the original
		uniform homeomorphism given by Odell and Schlumprecht \cite{OdellSchlumprecht1994}.
		%We now detail its extension to the vector-valued setting. 
		The next result is a crucial ingredient towards the main theorem of this section.
		\begin{prop} 
			\label{prop:extrap} 
			Suppose $p\in (1,\infty)$ and $\eps>0$. Then, there exist $r\in (1,p)$ and $s\in (p,\infty)$ such that if $X$ is a Banach sequence space that is strictly $r$-convex and strictly $s$-concave and if $(Y_n)_{n=1}^{\infty}$ is any sequence of Banach spaces, then there is a bijection $f\in \mathcal
			{GU}\big((\sum_{n=1}^{\infty} Y_n)_{\ell_p}, (\sum_{n=1}^{\infty} Y_n)_X\big)$ such that $[[f]]_\eps\le 2.$
		\end{prop}
		
		\begin{proof} 
			First, we apply Corollary I.5 in \cite{BenyaminiLindenstrauss2000} to deduce that we can choose $r,s$ so that for every such $X$ there are uniformly convex Banach sequence spaces $X_0$ and $X_1$ such that $\ell_p=X_0^{1/2}X_1^{1/2}$ and $X=X_0^{1-\theta}X_1^{1-\theta}$ where $\theta \in (0,1)$ is such that $\abs{\cot\pi\theta}<\vep/8$.
			
			Then, we define $f$ as follows. If $x\in(\sum_{n=1}^{\infty}Y_n)_{\ell_p}$,
			write $x :=(x_n)_{n=1}^{\infty} = (\norm{x_n}u_n)_{n=1}^{\infty}$ where
			$\norm{u_n}=1$ and $(\norm{x_n})_{n=1}^{\infty}\in \ell_p.$  Let
			$\norm{x_n}=\sqrt{a_nb_n}$ be the unique optimal factorization with
			$a :=(a_n)_{n=1}^{\infty}\in X_0,\ b:=(b_n)_{n=1}^{\infty}\in X_1$
			and $\norm{a}_{X_0}=\norm{b}_{X_1}=\norm{x}_{(\sum_{n=1}^{\infty} Y_n)_X}$. Then, we let
			\begin{equation*}
				f(x) := (a_n^{1-\theta}b_n^{\theta}u_n)_{n=1}^{\infty}.
			\end{equation*} 
			Note that in the scalar case when $Y_n=\bR$, $f$ is the map $f_0$ given by Theorem \ref{theo:homeo-spheres}. It follows easily that $f$ is also a norm-preserving bijection in the vector case.
			
			We next prove that $f$ is uniformly continuous on the unit ball of 
			$(\sum_{n=1}^{\infty} Y_n)_{\ell_p}$.  Let $\omega_0$ be the modulus of continuity of $f_{0\restriction_{B_{\ell_p}}}$. Now suppose $x,y\in (\sum_{n=1}^{\infty} Y_n)_{\ell_p}$ with $\norm{x},\norm{y}\le 1$ and let $t :=\norm{x-y}.$ We denote by $\bA$ the set of all $n\in \bN$ for which $\norm{x_n - y_n}\ge \sqrt{t}\max\{\norm{x_n},\norm{y_n}\}$ and let $\bB :=\bN\setminus \bA$. We denote by $\car_{\bA}$ the indicator function of $\bA$ and by $\car_{\bA}x$ the sequence $(\car_{\bA}(n)x_n)_n$. We clearly have 
			\begin{equation*}
				\norm{\car_{\bA}x} \le t^{-1/2}\norm{x-y}\le
				\sqrt{t}
			\end{equation*} 
			and  similarly $\norm{\car_{\bA}y}\le \sqrt{t}$.
			
			Hence,
			\begin{equation*}
				\max\{\norm{f(x)-f(\car_{\bB}x)},\norm{f(y) - f(\car_{\bB}y)}\} \le \omega_0(\sqrt{t}).
			\end{equation*}
			Let us write $\xi_n := \car_{\bB}(n)\norm{x_n}$ and $\eta_n :=\car_{\bB}(n)\norm{y_n}$. Let $u_n,v_n$ be normalized vectors such that $x_n=\xi_nu_n$ and $y_n=\eta_nv_n$ for $n\in \bN$. Then, by Lemma \ref{lem:polar} for all $n\in\bB$,  
			\begin{equation*}
				\norm{u_n-v_n}\le \max\{\xi_n,\eta_n\}^{-1}\norm{x_n-y_n}\le 2\sqrt{t}, 
			\end{equation*}  
			and hence for all $n\in\bB$, 
			\begin{equation*}
				\norm{f_0(\xi)_nu_n -f_0(\eta)_nv_n}\le \abs{f_0(\xi)_n-f_0(\eta)_n} + 2\min\{f_0(\xi)_n,f_0(\eta)_n\}\sqrt{t}.
			\end{equation*}
			This implies that
			\begin{equation*}
				\norm{f(\car_{\mathbb B}x)-f(\car_{\mathbb B}y)}\le \omega_0(t)+ 2\sqrt{t}.
			\end{equation*}
			Combining we have
			\begin{equation*}
				\norm{f(x)-f(y)}\le \omega_0(t)+2\omega_0(\sqrt{t})+2\sqrt{t}
			\end{equation*} 
			and $f$ is uniformly continuous. The argument for $f^{-1}$ being similar, we omit it.
			
			It remains to show that $[[f]]_\eps\le 2.$  To do this, we
			first consider $\norm{f(x)-f(y)}$ when $\norm{x},\norm{y}\in B_{(\sum
				Y_n)_{\ell_p}}$, $\norm{x-y}=t$ and $x,y$ have finite supports (this
			last restriction is not really necessary, but it avoids technical
			discussions). Let us suppose we have found $a,c\in X_0,\ b,d\in
			X_1$ with $a,b,c,d\ge 0$ and such that
			$\norm{a}_{X_0} = \norm{b}_{X_1} = \norm{x}$, 
			$\norm{c}_{X_0} = \norm{d}_{X_1} = \norm{y}$ and $\norm{x_n} = \sqrt{a_nb_n},\
			\norm{y_n} = \sqrt{c_nd_n}$. Then,
			\begin{equation*}
				f(x)-f(y) = \big(a_n^{1-\theta}b_n^{\theta}u_n- c_n^{1-\theta}d_n^{\theta}v_n\big)_{n=1}^\infty
			\end{equation*}
			where $u_n$ (resp. $v_n$) is a unit vector such that $x_n = \norm{x_n}u_n$ (resp. $y_n = \norm{y_n}v_n$). 
			Let us choose $x^*\in (\sum_{n=1}^{\infty} Y_n^*)_{X^*}$ with $\norm{x^*}=1$ and
			$x^*(f(x)-f(y))=\norm{f(x)-f(y)}$. Then, we can write
			\begin{equation*}
				x^*_n = \alpha_n^{1-\theta}\beta_n^{\theta}w_n^*
			\end{equation*} 
			where $\alpha\in X_0^*$, $\beta\in X_1^*$, $\alpha,\beta\ge 0$ are such that $\norm{\alpha}_{X_0^*} = \norm{\beta}_{X_1^*}=1$.
			
			Finally, for $0<\Re z<1$, consider the analytic function
			\begin{equation*}
				F(z) := \sum_{n=1}^{\infty} \alpha_n^{1-z}\beta_n^z\big(a_n^{1-z}b_n^z w_n^*(u_n)-c_n^{1-z}d_n^zw_n^*(v_n)\big). 
			\end{equation*}
			This function is bounded on the strip $S:=\{z\in \C\colon 0<\Re z<1\}$ and extends continuously to
			the boundary. On the boundary, we can estimate
			\begin{equation*}
				\abs{F(it)}\le \norm{\alpha}_{X_0^*}(\norm{a}_{X_0} + \norm{c}_{X_0})\le 2
			\end{equation*}
			and
			\begin{equation*}
				\abs{F(1+it)}\le \norm{\beta}_{X_1^*}(\norm{b}_{X_1} + \norm{d}_{X_1})\le 2.
			\end{equation*}
			Furthermore,
			\begin{equation*}
				\abs{F(1/2)}\le \norm{x-y}=t\le 2.
			\end{equation*}
			Hence, the function $\abs{F(z)-F(1/2)}$ is bounded by $4$ on the boundary. We can write
			\begin{equation*}
				F(z)-F(1/2)=\cot(\pi z)G(z)
			\end{equation*}
			where $G$ is analytic on the strip and $\abs{G}$ is also bounded by $4$ on the boundary. Thus, by the Three Lines Theorem,  $\abs{G(z)}\le 4$ throughout the strip $S$ and in particular $\abs{G(\theta)}\le 4$. Consequently,
			\begin{align*}
				\norm{f(x)-f(y)} = \abs{F(\theta)} & \le t+4\abs{\cot\pi\theta} \\
				& \le t + \frac{\vep}{2}\\
				& \le 2\max\{t,\vep\}= 2\max\{\norm{x-y},\vep \norm{x},\vep \norm{y}\}.
			\end{align*}
			The argument for $f^{-1}$ is very similar after reversing the roles of $\theta$ and $1/2$. In this case, one should use
			\begin{equation*}
				\varphi(z) := \frac{\sin\pi(z-\theta)}{\cos\pi(z+\theta)}
			\end{equation*} 
			in place of $\cot(\pi z)$ and observe that $\abs{\varphi(1/2)} = \abs{\cot(\pi\theta)}$.
		\end{proof}
		
		We would like to apply Proposition \ref{prop:extrap} to $\Tsi_p$. Unfortunately, while $\Tsi_p$ is strictly $p$-convex and hence strictly $r$-convex for every $r<p$, it cannot be strictly $p$-concave, otherwise it would be isomorphic to an $L_p(\mu)$ space, which is impossible. This roadblock can be cleared up by passing to a tail of $\Tsi_p$.  
		For $N\in \bN$, we will denote by $\Tsi_p^{>N}$ %(respectively $\widetilde \Tsi_p^N$) 
		the closed linear span in $\Tsi_p$ %(respectively $\widetilde \Tsi_p$)
		of $(e_n)_{n>N}$, where $(e_n)_{n=1}^\infty$ is the canonical basis of $\Tsi_p$. The following fact is proved in \cite{JLS1996} (see also \cite[Proposition 10.33]{BenyaminiLindenstrauss2000}).  
		
		\begin{prop}
			\label{prop:concavity-tail-Tsirelson}  
			There is a constant $M>0$ such that for every $q>p$ there exists an integer $N := N(q)$ such that $\Tsi_p^{>N}$ has
			$q$-concavity constant at most $M$.
		\end{prop}
		
		%By duality we deduce:
		
		%\begin{prop}\label{jls2}  
		%There is a constant $M$ so that if $q<p$ there exists
		%an integer $N=N(q)$ so that $\widetilde \Tsi_p^N$ has $q$-convexity
		%constant at most $M.$\end{prop}
		
		%We are now ready for our example.
		
		%\begin{theo}
		%\label{asymconvexexample}  Let $(G_n)_{n=1}^{\infty}$ be a sequence of finite-dimensional normed spaces dense for the Banach-Mazur distance in all finite-dimensional normed spaces.  Then for $1<p<\infty$, $\big((\sum_{n=1}^{\infty}G_n)_{\Tsi_p}\big)^2$,  $\big((\sum_{n=1}^{\infty}G_n)_{\widetilde \Tsi_p}\big)^2$ and $C_p=(\sum_{n=1}^{\infty}G_n)_{\ell_p}$ are all uniformly homeomorphic.\end{theo}
		
		We are now ready to prove the uniform equivalence claimed in the introductory discussion of this section.
		
		\begin{theo}
			\label{theo:Kp-Cp}
			Let $1<p<\infty$ and $\frac1p+\frac1q=1$. Then, $K_p\oplus K_p$, $(K_q\oplus K_q)^*$ and $C_p$ are uniformly homeomorphic.
		\end{theo}
		
		%\begin{rema} Recall that $K_p := \big(\sum_{n=1}^{\infty}G_n\big)_{\Tsi_p}$ and if $(G_n)_{n=1}^{\infty}$ is chosen carefully (e.g. replaced by $(G_1,G_1,G_2,G_2,\dots)$) there is no need for the squares in this statement since in this case $K_p\oplus K_p \approx K_p$. \end{rema}
		
		\begin{proof} 
			%We only prove the statement that $\big((\sum_{n=1}^{\infty}G_n)_{\Tsi_p}\big)^2$ is uniformly homeomorphic to $C_p$, as the proof for $\big((\sum_{n=1}^{\infty}G_n)_{\widetilde \Tsi_p}\big)^2$ is similar. 
			Since $C_p$ is isomorphic to $\ell_p(C_p)$, according to Corollary \ref{cor:Kalton-Ribe-lp-sum-square}, in order to prove that $C_p$ and $K_p \oplus K_p$ are uniformly homeomorphic, it is sufficient to show that $K_p$ and $C_p$ are uniformly close.
			
			So let $\vep>0$. Then, by Proposition \ref{prop:concavity-tail-Tsirelson}, we may
			find an integer $N$ such that $\Tsi_p^{>N}$ is $M$-isomorphic to a sequence space $X$ which is strictly
			$p$-convex and strictly $q$-concave with $q$ close enough to $p$ 
			so that by Proposition \ref{prop:extrap} there is a map $f\in\mathcal
			{GU}\big((\sum_{n=N+1}^{\infty}G_n)_X, (\sum_{n=N+1}^{\infty}G_n)_{\ell_p}\big)$ with $[[f]]_{\eps}\le
			2$. It follows that as a member of $\mathcal{GU}\big((\sum_{n=N+1}^{\infty}G_n)_{\Tsi_p},(\sum_{n=N+1}^{\infty}G_n)_{\ell_p}\big)$ we have $[[f]]_{\eps}\le 2M$.
			Next, we extend $f$ to a map 
			\begin{equation*}
				g\colon V:=\Big(\sum_{n=1}^N G_n\Big)_{\Tsi_p} \oplus_\infty \Big(\sum_{n=N+1}^{\infty} G_n\Big)_{\Tsi_p}\to
				W:=\Big(\sum_{n=1}^N G_n\Big)_{\Tsi_p} \oplus_\infty \Big(\sum_{n=N+1}^{\infty}G_n\Big)_{\ell_p}
			\end{equation*}
			given by $g(x_1,x_2)=(x_1,f(x_2))$. 
			
			Then, $g\in\mathcal {GU}(V,W)$ and $[[g]]_\eps\le 2M$. Considering $g$ as a map from $(\sum_{n=1}^{\infty}G_n)_{\Tsi_p}$ to $(\sum_{n=1}^NG_n)_{\Tsi_p}\oplus_p(\sum_{n=N+1}^{\infty}G_n)_{\ell_p}$,
			we easily get a crude estimate $[[g]]_\eps\le 8M$. But the latter space is certainly $C$-isomorphic to $C_p$, for some constant $C\ge 1$. Indeed, $\{(\sum_{n=1}^N G_n)_{\Tsi_p}\} \cup \{G_n \colon n\ge 1\}$ is another dense sequence for the Banach-Mazur distance in all finite-dimensional spaces, and we can use Exercise \ref{ex:C_p}. Therefore, $[[g]]_\vep\le 8CM$ when $g$ is considered as a member of $\mathcal {GU}(K_p,C_p)$. This means that $C_p$ and $(\sum_{n=1}^{\infty}G_n)_{\Tsi_p}$ are uniformly close and concludes the proof.
			
			Since $T_q^*$ is strictly $p$-concave, the proof that $(K_q\oplus K_q)^*$ and $C_p$ are uniformly homeomorphic is similar, albeit that we need to use the fact, which can be obtained by duality from Proposition \ref{prop:concavity-tail-Tsirelson}, that there is a constant $M$ such that for all $r<p$ there exists an integer $N:=N(r)$ for which $(T_q^*)^{>N}$ is $r$-convex with constant at most $M$.
		\end{proof}
		
		The properties of the space and its dual, constructed by Kalton, show the limitations of the preservation of asymptotic uniform smoothness of convexity under uniform homeomorphisms.
		
		\begin{theo}
			\label{theo:Kalton-optimality}
			For any $1<p<\infty$, the property of having an equivalent norm that is $p$-asymptotically uniformly smooth (or $p$-asymptotic uniformly convex) is not preserved under uniform homeomorphisms. 
		\end{theo}
		
		\begin{proof}
			As an $\ell_p$-sum of finite-dimensional spaces, the space $C_p$ is both $p$-AUC and $p$-AUS. On the other hand, $K_p\oplus K_p$ cannot have an equivalent $p$-AUS norm. Indeed, due to the strict $p$-convexity of the norm of $\Tsi_p$, $K_p\oplus K_p$ does have a $p$-AUS norm. If $K_p\oplus K_p$ also had an equivalent $p$-AUC norm, it would follow from Theorem \ref{thm:subspaces-of-lp-sums} that $K_p\oplus K_p$ would linearly embed into $C_p$, but this would contradict the well-known fact that $\Tsi_p$ and hence $K_p\oplus K_p$ does not contain a linear copy of $\ell_p$. 
			
			If $\frac1p+\frac1q=1$, the space $(K_q\oplus K_q)^* :=\big((\sum_{n=1}^{\infty}G_n^*)_{T_q^*}\big)^2$ yields a similar example for the modulus of asymptotic uniform smoothness. Indeed, $(K_q\oplus K_q)^*$ admits a $p$-AUC norm because of the strict $p$-concavity of $T_q^*$ and therefore, applying again Theorem \ref{thm:subspaces-of-lp-sums}, $(K_q\oplus K_q)^*$ cannot have an equivalent $p$-AUS norm, because $T_q^*$ and thus $(K_q\oplus K_q)^*$, does not linearly contain $\ell_p$.
			
			%As an $\ell_p$-sum of finite-dimensional spaces, the space $C_p$ is both $p$-AUC and $p$-AUS. On the other hand, $\widetilde{K}_p\oplus\widetilde{K}_p := \big((\sum_{n=1}^{\infty}G_n)_{\widetilde \Tsi_p}\big)^2$ cannot have an equivalent $p$-AUS norm. Indeed, due to the strict $p$-concavity of the norm of $\widetilde{T}_p$, $\widetilde{K}_p\oplus\widetilde{K}_p$ does have a $p$-AUC norm. If $\widetilde{K}_p\oplus\widetilde{K}_p$ also had an equivalent $p$-AUS norm, it would follow from Theorem \ref{thm:subspaces-of-lp-sums} that $\widetilde{K}_p\oplus\widetilde{K}_p$ would linearly embed into $C_p$, but this would contradict the well-known fact that $\widetilde{T}_p$ and hence $\widetilde{K}_p\oplus\widetilde{K}_p$ does not contain a linear copy of $\ell_p$. 
		\end{proof}
		
		Theorem \ref{theo:Kalton-optimality} shows the optimality of Corollary \ref{cor:ausstableCL}, which states that if $X$ is uniformly homeomorphic to a $p$-AUS Banach space $Y$ and $r<p$, then $X$ admits an equivalent $r$-AUS norm. It should also be compared to the positive results in Chapter \ref{chapter:AMP_II}, like Theorem \ref{thm:CL-into-pAMUC->weak-pcoBS}.
		
		\subsection{Kalton's pair of coarsely equivalent but not uniformly equivalent separable Banach spaces}\label{subsec:CEnotUH}
		
		In this section, we prove a stunning result of Kalton, which answered in \cite{Kalton2012} the following problem of Johnson, Lindenstrauss and Schechtman raised in \cite[Section 8 (2)]{JLS1996}: Are two coarsely equivalent Banach spaces necessarily uniformly equivalent?
		The construction of a counterexample by Kalton is an astonishing breakthrough, which, due to Kalton's untimely passing, was never presented at conferences. The pair of Banach spaces in Kalton's counterexample is relatively easy to describe and comes from the lifting method. However, either showing that the spaces are coarsely equivalent or not uniformly equivalent involves a high level of sophistication where approximate midpoint, concentration inequality on the interlacing graphs and asymptotic uniform smoothness arguments are combined and mixed in an almost magical manner to produce the desired counterexample.
		
		First, we start with the description of the counterexample. Given a separable Banach space $X$, consider a quotient map $Q_X\colon \ell_1 \to X$.
		
		\begin{rema}
			It is a result of Lindenstrauss and Rosenthal (see \cite{LindenstraussRosenthal1969} or  \cite[page 108]{LindenstraussTzafriri1977}) that if $X$ is not isomorphic to $\ell_1$, then $Q_X$ is unique up to automorphisms.
		\end{rema}
		
		Now, for each $n\in \bN$, one can define an equivalent norm on $\ell_1$, denoted by $\norm{\cdot}_{X,n}$, as follows:
		\begin{equation*}
			\norm{z}_{X,n}:=\max\{2^{-n}\norm{z}, \norm{Q_X(z)}_X\}.
		\end{equation*}
		We denote by $Q_{X,n}$ the map $Q_X$ viewed as an application from $(\ell_1,\norm{\cdot}_{X,n})$ onto $X$ and 
		\begin{align*}
			\widetilde{Q}_X \colon \Big(\sum_{n=1}^\infty(\ell_1,& \norm{\cdot}_{X,n})\Big)_{\ell_1}\to X \\ 
			&(z_n)_{n=1}^\infty \mapsto  \sum_{n=1}^\infty Q_{X,n}(z_n)
		\end{align*}
		Note that $\widetilde{Q}_X$ is a quotient map.
		By taking $X=\co$ (or any separable, nonreflexive and asymptotically uniformly smooth Banach space, for that matter), we get the counterexample. We detail the simpler proof for $X=c_0$. 
		
		\begin{theo}
			\label{thm:Kalton-CE-not-UE-counterex}
			The separable Banach spaces $\big(\sum_{n=1}^\infty(\ell_1,\norm{\cdot}_{\co,n})\big)_{\ell_1}$ and $\big(\sum_{n=1}^\infty(\ell_1,\norm{\cdot}_{\co,n})\big)_{\ell_1}\oplus \co$ are coarsely equivalent but not uniformly equivalent.    
		\end{theo}
		
		The proof of Theorem \ref{thm:Kalton-CE-not-UE-counterex} relies on several rather deep results of independent interest. 
		The first one is about constructing a global coarse lifting of a quotient map from a sequence of good local liftings; this will be evidently helpful when trying to show the existence of a coarse equivalence.
		
		\begin{theo}
			\label{thm:local-global-coarse-lifting}
			Let $Q\colon Y\to X$ be a quotient map between two Banach spaces. If there exist a constant $L\ge 1$, a sequence $(\vep_n)_{n=0}^\infty$ of real numbers in $(0,1]$ with $\lim_n\vep_n=0$ and for all $n \ge 0$, a lifting $\varphi_n\colon S_X\to Y$ so that $\varphi_n$ is of CL-type $(L,\vep_n)$, then $Q$ admits a coarse-Lipschitz lifting.
		\end{theo}
		
		\begin{proof}
			By replacing $\varphi_n(x)$ with $\frac12(\varphi_n(x)-\varphi_n(-x))$, we can assume without loss of generality that $\varphi_n$ is antisymmetric on $S_X$. From this and the assumption that the $\varphi_n$ are of $CL$-type $(L,\vep_n)$, it follows easily that for all $n\in \N \cup \{0\}$ and all $x\in S_X$, $\norm{\varphi_n(x)}\le 2L$. Furthermore, the radial extension $f_n$ of $\varphi_n$ is a positively homogeneous and bounded map satisfying $Q\circ f_n = Id_X$ and, moreover, it follows from Lemma \ref{lem:local-global-homogeneous}  that $\norm{f_n}_{\vep_n} \le 8L$. By extracting a subsequence if needed, we can also assume that $\norm{f_n}_{e^{-2(n+1)}} \le 8L$ for all $n\ge 1$. After extending $n\mapsto f_n$ piecewise linearly, we obtain a collection of maps $(f_t)_{t\in[0,\infty)}$ such that $Q\circ f_t =Id_X$ and $\norm{f_t}_{e^{-2t}}\le 8L$ (we used above  that $\vep \mapsto \norm{g}_\vep$ is nonincreasing). It remains to observe that an elementary inspection of a few cases reveals that $\norm{f_s-f_t}\le 4L\abs{s-t}$ and to conclude we invoke Proposition \ref{prop:base-construction-gluing}.
		\end{proof}

		\begin{rema}
			If the local liftings in Theorem \ref{thm:local-global-coarse-lifting} are also uniformly continuous, by the second point of Proposition \ref{prop:base-construction-gluing}, the global lifting can be shown to be uniformly continuous. 
		\end{rema}
		
		The next theorem is the crucial ingredient providing an obstruction to the existence of a uniform equivalence.
		
		\begin{theo}
			\label{thm:co-not-CLU-into-l1-sum-of-Q}
			The Banach space $\co$ does not admit a coarse-Lipschitz embedding that is also uniformly continuous into a countable $\ell_1$-sum of Banach spaces whose unit balls can be uniformly embedded into a reflexive Banach space.
		\end{theo}
		
		\begin{proof}
			Let $(Y_n)_n$ be a sequence of Banach spaces such that for all $n\in\bN$ $B_{Y_n}$ can be uniformly embedded into a reflexive Banach space. Let $Y:=(\sum_{n=1}^\infty Y_n)_{\ell_1}$. Before we embark on the proof, note that if the $\ell_1$-sum were finite, then $Y$ would have property $Q$ and consequently $\co$ would not even coarsely embed into $Y$. Thus, the difficulty in the ensuing arguments comes from the fact that if the sum is truly infinite, then $Y$ might not have property $Q$ and hence the stronger faithfulness assumption for the embedding.
			
			So let $f:=(f_n)_{n=1}^\infty \colon \co \to Y$ be a coarse-Lipschitz embedding that is uniformly continuous and assume without loss of generality that there are $A>0$ such that for all $x_1,x_2\in \co$ 
			\begin{equation}
				\label{eq:lower-CL-estimate}
				\norm{f(x_1) - f(x_2)}\ge A\norm{x_1 -x_2} - 1.
			\end{equation}
			The proof of Claim \ref{clai:midpoint}, which is a particular case of Exercise \ref{ex:AUS-CLU-into-l1-sum-of-Q->reflexive} (assertion \ref{item:midpoint-claim}), follows from a now standard approximate midpoint argument.
			
			\begin{claim}
				\label{clai:midpoint}
				Given $\vep>0$, there exist $w\in \co$, $s>\frac{1}{\vep}$, $N\in\bN$ and $K\in \N$  such that for all $z\in Z_K$ with $\norm{z}\le s$,
				\begin{equation}
					\label{eq:midpoint}
					\sum_{n=N+1}^\infty \norm{f_n(w+z)}\le s\vep,
				\end{equation}
				where $Z_K$ is the closed linear span of $(e_k)_{k>K}$ in $c_0$ and $(e_k)_{k=1}^\infty$ is the canonical basis of $c_0$.   
			\end{claim}
			
			Our goal is to derive a contradiction with the nonreflexivity of $c_0$. Let $(s_n)_{n=1}^\infty$ be the summing basis of $c_0$ and $\M_1=\{n\in \N,\ n>K\}$. For all $n\neq m \in \M_1$, $s_n-s_m \in Z_K$ and for all $m_1\le n_1\le m_2\le n_2\le \dots \le m_{2k}\le n_{2k} \in \M_1$, 
			$$\Big\|\sum_{i=1}^k (s_{n_i}-s_{n_{i+k}})-(s_{m_i}-s_{m_{i+k}})\Big\|_\infty\le 1.$$ 
			For a fixed $k\ge 1$, define now $h_k\colon ([\M_1]^{2k},\dIk) \to Z_K$ by 
			\begin{equation*}
				h_k(\nbar):=f\Big(w+\frac{s}{k}\sum_{i=1}^k (s_{n_i}-s_{n_{i+k}})\Big),
			\end{equation*}
			where $w\in \co$ and $s$ are as in Claim \ref{clai:midpoint}. We clearly have that $\Lip(h_k)\le \omega_f(\frac{s}{k})$. 
			
			Now, let $P_N\colon Z \to V_N:=\big(\sum_{n=1}^N Y_n\big)_{\ell_1}$ be the canonical projection and observe that it follows from the assumptions that $B_{V_N}$ embeds uniformly into a reflexive Banach space. Therefore, $V_N$ has property $\cQ$ for some constant $Q_N$ and there is a further infinite subset of $\M_1$, denoted by $\bM$ such that for all $\nbar,\mbar \in [\M]^k$,
			\begin{equation}
				\label{eq2:co-not-CLU-into-l1-sum-of-Q}
				\norm{P_N h_k(\nbar) - P_N h_k(\mbar)} \le Q_N\, \omega_f\Big(\frac{s}{k}\Big).
			\end{equation}
			Then, it follows from the definition of $h_k$ that for all $\nbar,\mbar \in [\M]^k$,
			\begin{equation*}
				\label{eq3:co-not-CLU-into-l1-sum-of-Q}
				\norm{h_k(\nbar) - h_k(\mbar)} \stackrel{\eqref{eq:midpoint} \land \eqref{eq2:co-not-CLU-into-l1-sum-of-Q} }{\le} Q_N\,\omega_f\Big(\frac{s}{k}\Big) + 2s\vep.
			\end{equation*}
			Therefore, combining with \eqref{eq:lower-CL-estimate} we get
			\begin{equation*}
				\Big\|\frac{1}{k} \sum_{i=1}^k (s_{n_i}-s_{n_{i+k}})-(s_{m_i}-s_{m_{i+k}})\Big\|_\infty \le \frac{Q_N}{sA} \omega_f\Big(\frac{s}{k}\Big) + \frac{2\vep}{A} +\frac{1}{sA}, 
			\end{equation*}
			for all $\nbar,\mbar \in [\M]^k$. Denote by $e\in \ell_\infty$ the constant sequence equal to $1$. Fixing $n_1,\ldots,n_k,m_1,\ldots,m_{2k}$, letting $n_{2k},\ldots,n_{k+1}$ tend to $\infty$ and using the weak$^*$ lower-semicontinuity of the norm of $\ell_\infty$, we get that 
			\begin{equation*}
				\Big\|e- \frac{1}{k} \sum_{i=1}^k (s_{n_i}-s_{m_i}+s_{m_{i+k}}) \Big\|_\infty\le \frac{Q_N}{sA} \omega_f\Big(\frac{s}{k}\Big) + \frac{2\vep}{A} +\frac{1}{sA}.
			\end{equation*}
			Finally, since $s\vep>1$ and $f$ is uniformly continuous, after taking the limit when $k\to \infty$ we have
			\begin{equation*}
				d(e,\co)\le \frac{3\vep}{A},
			\end{equation*}
			but since $\vep$ can be taken arbitrarily small one gets $e\in\co$; a contradiction.
			
		\end{proof}

		\begin{rema}\label{rema:mistake}
			The argument in Theorem \ref{thm:co-not-CLU-into-l1-sum-of-Q} essentially shows that any asymptotically uniformly flat separable Banach space admitting a coarse-Lipschitz embedding that is uniformly continuous into a countable $\ell_1$-sum of Banach spaces whose unit balls can be uniformly embedded into a reflexive Banach space, must be reflexive. The same conclusion with essentially the same proof holds if the embedded space is just asymptotically uniformly smooth (cf Exercise \ref{ex:AUS-CLU-into-l1-sum-of-Q->reflexive} or the original paper \cite{Kalton2012} by Kalton). We should also warn the reader that the proof in \cite{Kalton2012} contains a small mistake. The sequence denoted by $(x_n)_n$ in \cite{Kalton2012} corresponds to $(s_n)_n$ in the case of $c_0$. We had to consider differences $s_n-s_m$ and modify the definition of the map $h_k$ to ensure that the relevant points belong to $Z_K$, which Kalton forgot to do. The same modification has to be done in Exercise \ref{ex:AUS-CLU-into-l1-sum-of-Q->reflexive}. 
		\end{rema}
		We can now proceed with the final details of the proof of the main result of this section.
		
		\begin{proof}[Proof of Theorem \ref{thm:Kalton-CE-not-UE-counterex}]
			Let $\cZ_1(\co):=\big(\sum_{n=1}^\infty(\ell_1,\norm{\cdot}_{\co,n})\big)_{\ell_1}$. It is an immediate consequence of Theorem \ref{thm:co-not-CLU-into-l1-sum-of-Q} that $\cZ_1(\co)$ is not uniformly equivalent to $\cZ_1(\co)\oplus \co$,  since any map implementing this equivalence would provide a coarse-Lipschitz embedding that is uniformly continuous of $\co$ into $\cZ_1(\co)$ which, by construction, is an $\ell_1$-sum of renormings of $\ell_1$ whose unit balls uniformly embed into, say, the reflexive space $\ell_2$ via the Mazur map.  
			
			For the coarse equivalence, it is sufficient to show that the quotient map $\widetilde{Q}_{\co}\colon \cZ_1(\co) \to \co$ admits a coarse lifting. Indeed, in this case it follows from Proposition \ref{prop:section-method} that $\cZ_1(\co)$ is coarse-Lipschitz equivalent to $\ker(\widetilde{Q}_{\co})\oplus \co$. Since $c_0$ is isomorphic to its square, we have that $\cZ_1(\co) \oplus c_0$ is coarse-Lipschitz equivalent to $\ker(\widetilde{Q}_{\co})\oplus \co^2 \simeq \ker(\widetilde{Q}_{\co})\oplus \co$ and therefore $\cZ_1(\co) \oplus c_0$ is coarse-Lipschitz equivalent to $\cZ_1(\co)$. So, according to Theorem \ref{thm:local-global-coarse-lifting}, it is enough to find a sequence of almost Lipschitz liftings of the unit sphere of $\co$. We can construct such liftings as follows. Since $Q_{\co}\colon \ell_1\to \co$ is a quotient map, there is $C>0$ such that for all $z\in \co$, we can find $\psi(z)\in \ell_1$ such that $Q_{\co}\circ \psi(z)=z$ and $\norm{\psi(z)}\le C\norm{z}_\infty$. Let $\psi_n\colon \co \to (\ell_1,\norm{\cdot}_{\co,n})$ such that $\psi_n:=\psi$ and observe that $\omega_{\psi_n}(t)\le \max\{2^{-n}\omega_{\psi}(t),t\}$.
			Defining $\varphi_n\colon S_{\co}\to \cZ_1(\co)$ by $\varphi_n:=\iota_n \circ \psi_n$ where $\iota_n\colon (\ell_1,\norm{\cdot}_{\co,n}) \to \cZ_1(\co)$ is the canonical embedding, we have that $(\widetilde{Q}_{\co}\circ \varphi_n)(z)=z$ for all $z\in S_{\co}$ and $\omega_{\varphi_n}(t)\le t +C2^{1-n}$ for all $t\in[0,2]$ and the sequence $(\varphi_n)_n$ provides the sought-after almost Lipschitz liftings on the unit sphere, thereby concluding the proof. 
		\end{proof}
		
		\section{Notes}
		
		There is another important example of nonseparable Banach spaces that are Lipschitz isomorphic but not linearly isomorphic. Indeed, Godefroy, Deville and Zizler \cite{DGZ1990} proved that if $K$ is a Hausdorff compact space, then $C(K)$  is Lipschitz equivalent to $c_0(\Gamma)$, for some set $\Gamma$ if and only if its Cantor-Bendixon index is finite. On the other hand, Ciesielski and Pol \cite{CiesielskiPol1984} proved the existence of a compact space $K$ such that $K^{(3)}=\emptyset$, but so that there is no bounded linear injection from $C(K)$ into any $c_0(\Gamma)$ space. In fact, the situation is even worse; there is no linear injection that is continuous from $(C(K),w)$ into $(c_0(\Gamma),w)$ (see  Section VI.8 in \cite{DGZ1993} for a complete presentation). 
		
		To summarize the history of the material presented in Section \ref{sec:Ribe-Kalton-machinery}, we just quote the very modest introductory sentence written by Nigel Kalton in \cite{Kalton2013}: ``The first example of two separable uniformly homeomorphic Banach spaces which are not linearly isomorphic was given in 1984 by Ribe. Ribe’s basic approach has been used subsequently to create further examples by Aharoni and Lindenstrauss and also Johnson, Lindenstrauss and Schechtman.The aim of this paper is to give further examples by exploiting the same basic technique.'' As we already said, Kalton's contributions to the construction of counterexamples contain much more, like the use of interlacing graphs or asymptotic smoothness. There is one other fundamental tool that needs to be pointed out: the use of interpolation theory. Kalton's outstanding contributions to the theory of interpolation would deserve a book on its own. We will not give here a long list of references. We refer the interested reader to the website ``Nigel Kalton Memorial'' (https://kaltonmemorial.missouri.edu), where she or he will easily find these papers. 
		
		In relation to Section \ref{section:CLstabilitySQl_p} let us mention one more example that can be found in \cite{Kalton2013}. Kalton constructed, for any  $1<p\neq 2<\infty$, a pair of subspaces of $\ell_p$ that are uniformly homeomorphic but not linearly isomorphic.

		\section{Exercises}
		
		\begin{exer}
			\label{ex:balls-to-spheres}
			Let $X$ and $Y$ be Banach spaces.
			\begin{enumerate}
				\item Show that if $f\colon S_X \to S_Y$ is a uniform homeomorphism, then the radial extension of $\bar{f}$ is a bijection between $X$ and $Y$ that is norm-preserving, positively homogeneous and such that $\bar{f}_{\restriction B_X}$ and $\bar{f}^{-1}_{\restriction B_Y}$ are uniformly continuous.
				\item Show that if $f\colon X \to Y$ is a bijection between two Banach spaces that is positively homogeneous and such that $f_{\restriction B_X}$ and $f^{-1}_{\restriction B_Y}$ are uniformly continuous, then $S_X$ and $S_Y$ are uniformly homeomorphic.
			\end{enumerate}
			
		\end{exer}
		
		\begin{exer}
			\label{ex:almost-disjoint}
			Show that there exists a family $(N_\lambda)_{\lambda \in \bR}$ of infinite subsets of $\bN$ such that $N_\lambda \cap N_\mu$  is finite for all $\lambda \neq \mu \in \bR$. 
			%Detail the proof of Lemma \ref{lem:almostdisjoint}.     
		\end{exer}
		
		\begin{exer}
			Prove the Linking Principle (Lemma \ref{lem:link}). 
		\end{exer}
		
		\begin{exer} 
			Let $p\neq q \in [1,\infty)$ and $(p_n)_{n=1}^\infty$ be a sequence in $(p,\infty)$ which is converging to $p$. Show that $(\sum_{n=1}^\infty \ell_{p_n})_{\ell_q}$ is uniformly homeomorphic but not linearly isomorphic to $(\sum_{n=1}^\infty \ell_{p_n})_{\ell_q} \oplus \ell_p$.
		\end{exer}

		\begin{exer} Let $X$ and $Y$ be two Banach spaces.
			\begin{enumerate}
				\item Assume that $X$ and $Y$ are uniformly homeomorphic. Give a very elementary proof of the fact that $\co(X)$ and $\co(Y)$ are uniformly homeomorphic. 
				\item Assume that $X$ and $Y$ are coarse-Lipschitz equivalent. Give a less but still elementary proof of the fact that $\co(X)$ and $\co(Y)$ are coarse-Lipschitz equivalent. 
			\end{enumerate} 
		\end{exer}
		
		\begin{proof}[Hint]
			For $2.$ use net equivalences.
		\end{proof}
		
		\begin{exer}
			\label{ex:pconvex-pAUS} 
			Let $X$ be a Banach sequence space and $p\in (1,\infty)$.
			\begin{enumerate}
				\item Assume that $X$ is strictly $p$-convex. Show that $X$ is $p$-AUS.
				\item Assume that $X$ is strictly $p$-concave. Show that $X$ is $p$-AUC.
				\item Let $(F_n)_{n=1}^\infty$ be a sequence of finite-dimensional normed spaces and assume that $X$ is strictly $p$-convex (resp. strictly $p$-concave). Show that $\big(\sum_{n=1}^\infty F_n\big)_X$ is $p$-AUS (resp. $p$-AUC). 
			\end{enumerate}
			
		\end{exer}

		\begin{exer}
			\label{ex:AUS-CLU-into-l1-sum-of-Q->reflexive}
			Let $X$ be an asymptotically uniformly smooth Banach space and $Y:=(\sum_{n=1}^\infty Y_n)_{\ell_1}$ where the $Y_n$ are Banach spaces. 
			\begin{enumerate}
				\item\label{item:midpoint-claim} Suppose that $f:=(f_n)_{n=1}^\infty\colon X\to Y$ is a coarse map. Show that for all $r,\vep>0$, there exist $w\in X$, $s>r$, $N\in\bN$ and a closed finite-dimensional subspace $Z$ of $X$ such that for all $z\in Z$ with $\norm{z}_X\le s$,
				\begin{equation}
					\sum_{n=N+1}^\infty \norm{f_n(w+z)}_{Y_n}\le s\vep.
				\end{equation}    
				\item Show that any separable asymptotically uniformly smooth Banach space admitting a coarse-Lipschitz embedding that is also uniformly continuous into a countable $\ell_1$-sum of Banach spaces whose unit balls can be uniformly embedded into a reflexive Banach space, must be reflexive.
			\end{enumerate}
		\end{exer}
		
		%%%%%%%%%%%%%%%%%%%%%%%%%%%%%%%%%%%%%%%%%%%%%%%%%%%%%%%%%%%%%%%%%%%%%%%%%%%%%%%%%%%%

		\chapter{The geometry of trees and applications}
		\label{chapter:trees}
		
		In this chapter, we investigate the geometry of countably branching trees. In Section \ref{sec:tree-basics}, we 
		formally define countably branching trees and prove some elementary bi-Lipschitz embedding results into the Banach spaces $\ell_1$ and $\co$. These embedding results are then refined and extended in Section \ref{sec:tree-embeddings}. The connection between coarse embeddings of trees into Banach spaces and the spreading models or asymptotic models of the latter is detailed in Section \ref{sec:tree-models}. In particular, a purely metric characterization of infinite dimensionality of a Banach space is expressed in terms of the coarse geometry of countably branching trees.
		In Section \ref{sec:tree-Szlenk}, it is shown that Banach spaces with large Szlenk index contain equi-bi-Lipschitzly the sequence of countably branching trees. Then, a Poincar\'e-type inequality for Lipschitz maps on the countably branching trees taking values into Banach spaces with property $(\beta_p)$ is obtained in Section \ref{sec:tree-umbel}, as well as distortion lower bounds and compression function upper bounds that are derived from it. From these results and the connection between Szlenk index, asymptotic renormings and property $(\beta)$ from Chapters \ref{chapter:asymptotic-moduli} and \ref{chapter:Szlenk}, several metric characterizations of asymptotic properties and rigidity results are derived in Section \ref{sec:tree-applications}.
		
		\section{Countably branching trees: basic definitions and elementary embeddings}
		\label{sec:tree-basics}
		As a graph, the countably branching (rooted) tree of infinite height, denoted by $\wT$, has vertex set $\N^{<\omega}$. We refer the reader to Section \ref{sec:summable-A} for the notation on trees such as the natural ordering, length, concatenation, or restriction of finite sequences and the immediate predecessor of a finite sequence. Two vertices $s$ and $t$ are adjacent if and only if $s\prec t$ and $|t|=|s|+1$ or $t\prec s$ and $|s|=|t|+1$. We will always equip this tree with the graph metric denoted by $\dwT$. This metric is given by the formula 
		\begin{equation}
			\dwT(s,t)=|s|+|t|-2|u|,   
		\end{equation}
		where $u$ is the greatest common predecessor of $s$ and $t$ for $\preceq$. This metric can again be realized as a metric induced by a Banach space basis. If we consider the canonical basis of $\ell_1(\bN^{<\omega})\equiv \ell_1$, denoted by $(e_{s})_{s\in \bN^{<\omega}}$, then it is immediate that 
		\begin{equation}
			\dwT(s,t)=\Big\|\sum_{v \preceq s} e_{v} - \sum_{w\preceq t} e_{w}\Big\|_1.
		\end{equation}
		In other word, the map $h_1\colon \wT\to \ell_1(\bN^{<\omega})$ given by $h_1(s)=\sum_{v \preceq s} e_{v}$ is an isometric embedding.
		
		The countably branching (rooted) tree of height $k$, denoted by $\wTk$, is defined via the same adjacency relation, albeit on the vertex set $\N^{\le k}$.
		
		\begin{rema}
			Note that if we restrict $\dwTk$ to $\bN^k$, the leaves of the tree $\wTk$, then $\dwTk$ is an ultrametric.
		\end{rema}
		
		Given $(x_{s})_{s \in \N^{\le k}}$, a normalized tree of height $k$ in a Banach space $X$, the map $h\colon \wTk \to X$ given by $h(s):= \sum_{v\preceq s} x_v$ is always $1$-Lipschitz. However, even if we require the tree $(x_{s})_{s \in \N^{\le k}}$ to be weakly null, it does not have to have some ``suppression unconditionality property" that will allow us to connect the quantity $\norm{h(s)-h(t)}$ to lower estimates on the branches of the tree $(x_{s})_{s \in \N^{\le k}}$. Trees with useful unconditionally properties are the trees generated by weakly null sequences or weakly null arrays and that appear in the theory of spreading models or asymptotic models. We will prove embedding results similar to Proposition \ref{prop:J-graph-embedding} or Proposition \ref{prop:H-graph-embedding} in the next section. 
		
		Recall that every separable metric space admits a bi-Lipschitz embedding into $\co$. However, the countably branching trees admit a particular embedding into $\co$ (different from Aharoni's embedding) that will be particularly helpful in the study of their embeddability properties. 
		Let $(e_{s})_{s\in \N^{<\omega}}$ be the canonical basis of $\co(\N^{<\omega})$ and consider the corresponding summing basis $(\sigma_{s})_{s\in \N^{<\omega}}$ given by $\sigma_{s}=  \sum_{u \preceq s} e_{u}$. 
		
		\begin{prop}
			The map $h_{\infty}\colon \wT \to \co(\N^{<\omega})$ given by $h_\infty(s)=\sum_{v \preceq s} \sigma_v$ is a bi-Lipschitz embedding with distortion $2$.
		\end{prop} 
		
		\begin{proof}
			If $u$ is the greatest common ancestor of $s$ and $t$, then, 
			
			\begin{align}
				\label{eq1:tree-co} \norm{ h_\infty(s) - h_\infty(t) }_\infty  &= \Big\|\sum_{u\prec v \preceq s} \sigma_v - \sum_{u\prec w \preceq s}\sigma_w\Big\|_\infty \\
				\label{eq2:tree-co} &=  \Big\| \sum_{i=|u| + 1}^{|s|} \sum_{v \preceq  s_{\restriction i}} e_v - \sum_{i=|u| + 1}^{|t|} \sum_{w \preceq  t_{\restriction i}} e_w \Big\|_\infty.
			\end{align}
			Since $\norm{\sigma_v}_\infty=1$ for all $v \in \N^{<\omega}$, it follows from the triangle inequality and \eqref{eq1:tree-co} that $h_\infty$ is $1$-Lipschitz. It remains to observe that the coefficients of $e_{s_{\restriction |u|+1}}$ and $e_{t_{\restriction |u|+1}}$ in \eqref{eq2:tree-co} are $\abs{s} - \abs{u}$ and $\abs{t} - \abs{u}$, respectively. Since $e_{s_{\restriction |u|+1}} \neq e_{t_{\restriction |u|+1}}$, it follows that 
			\[ \norm{ h_\infty(s) - h_\infty(t) }_\infty \ge  \max(|s|-|u|,|t|-|u|)\ge \frac12(|s|+|t|-2|u|)=\frac12\dwT(s,t).\]
		\end{proof}

		\section{Bi-Lipschitz and coarse embeddings of trees}
		\label{sec:tree-embeddings}
		
		\subsection{Embeddability of trees and asymptotic models}
		\label{sec:tree-models}
		
		The next proposition is an analog of Proposition \ref{prop:J-graph-embedding} adapted to the tree metric.
		
		\begin{prop}
			\label{prop:tree-spreading}
			Let $X$ be an infinite-dimensional Banach space admitting a spreading model $S$ generated by a normalized weakly null sequence, then for every $\nu>0$ there exists a map  $h_{\nu} \colon \wT \to X$ such that for all $s,t \in \wT$,
			\begin{equation*}
				\frac{1}{8(1+\nu)}\varphi_S\left(\dwTk(s,t)\right) \le \norm{h_{\nu}(s)-h_{\nu}(t)}_X\le \dwT(s,t).
			\end{equation*}
			In particular, $\wT$ coarsely embeds into every infinite-dimensional Banach space admitting a spreading model generated by a normalized weakly null sequence that is not isomorphic to $\co$.
		\end{prop}
		
		\begin{proof}
			Let $\vep>0$ such that $(1+\vep)^2\le (1+\nu)$. By Proposition \ref{prop:weakly-null-spreading-model}  there is a weakly null normalized basic sequence $\xn$ with basis constant  not larger than $(1+\vep)$ and thus the bi-monotonicity constant is at most $2(1+\vep)$, generating a spreading model $S$ with associated fundamental sequence $(e_i)_{i=1}^\infty$ such that for all $r\ge1$, for all $r\le n_1< n_2< \dots < n_{r}$ and for all $(\vep_i)_{i=1}^{r}\in\{-1,1\}^r$ one has
			\begin{equation}
				\label{eq:T-graph-embedding-eq2}
				\Big\|\sum_{i=1}^{r} \vep_i x_{n_i}\Big\| \ge \frac{1}{1+\vep}\Big\|\sum_{i=1}^{r} \vep_i e_i\Big\| \ge \frac{1}{2(1+\vep)}\varphi_S(r).
			\end{equation}
			Let us also recall that $(e_i)_{i=1}^\infty$ is spreading and $1$-suppression unconditional (and thus $2$-unconditional).  Let now $\Phi\colon \N^{<\omega} \to \bN$ be a bijection and define $h\colon \N^{<\omega} \to Y$ by $h(s)=\sum_{u \preceq s}x_{\Phi(u)}$, for $s \in \N^{<\omega}$.
			The map $h$ is clearly $1$-Lipschitz and note that after cancellations one has
			\begin{equation}
				\norm{h(s) - h(t)} = \Big\|\sum_{i=1}^{\sd_{\sT_\infty^\omega}(s,t)} \vep_i x_{q_i}\Big\|,
			\end{equation}
			for some $q_1<q_2<\dots<q_{\sd_{\sT_\infty^\omega}(s,t)}$ and $(\vep_i)_{i=1}^{\sd_{\sT_\infty^\omega}(s,t)}\in\{-1,1\}^{\sd_{\sT_\infty^\omega}(s,t)}$. 
			
			Consider first the case when $d:=\sd_{\sT_\infty^\omega}(s,t)$ is odd, say  $d=2r-1$ with $r\in \bN$. Note that $q_r\ge r$. So, again \eqref{eq:T-graph-embedding-eq2} and
			the assumption on the bi-monotonicity constant of $\xn$ imply that
			\begin{align*}
				\norm{h(s) - h(t)}&
				\ge \frac{1}{2(1+\vep)}\Big\|\sum_{i=r}^{2r-1}\vep_i x_{q_i}\Big\|
				\ge \frac{1}{2(1+\vep)^2}\Big\|\sum_{i=1}^{r}\vep_i e_{i}\Big\|\\
				&\ge \frac{1}{4(1+\vep)^2}\Big\|\sum_{i=1}^{r} e_{i}\Big\| \ge \frac{1}{8(1+\vep)^2}\Big\|\sum_{i=1}^{2r} e_{i}\Big\|\\
				& \ge \frac{1}{8(1+\vep)^2}\Big\|\sum_{i=1}^{2r-1} e_{i}\Big\| = \frac{1}{8(1+\vep)^2}\varphi_S(d).
			\end{align*}
			
			In the case when $d:=\sd_{\sT_\infty^\omega}(s,t)$ is even, say  $d=2r$ with $r\in \bN$. Note that again $q_{r+1}\ge r+1$. It follows from \eqref{eq:T-graph-embedding-eq2} and the assumption on the bi-monotonicity constant of $\xn$ that
			\begin{align*}
				\norm{h(s) - h(t)}&
				\ge \frac{1}{2(1+\vep)}\Big\|\sum_{i=r+1}^{2r}\vep_i x_{q_i}\Big\|
				\ge \frac{1}{2(1+\vep)^2}\Big\|\sum_{i=1}^{r}\vep_i e_{i}\Big\|\\
				&\ge \frac{1}{4(1+\vep)^2}\Big\|\sum_{i=1}^{r} e_{i}\Big\| \ge \frac{1}{8(1+\vep)^2}\Big\|\sum_{i=1}^{2r} e_{i}\Big\| = \frac{1}{8(1+\vep)^2}\varphi_S(d).
			\end{align*}
		\end{proof}

		In the next proposition, we are restricted to trees of finite height due to the limitation on the suppression unconditionally property from Proposition  \ref{prop:jointspreadingmodels}.
		
		\begin{prop}
			\label{prop:T-graph-embedding}
			Let $X$ be an infinite-dimensional Banach space admitting an asymptotic model $A=[\en]$, generated by a normalized weakly null array, then for every $\nu>0$ and $k\in \bN$ there exists a map  $h_{\nu,k} \colon \wTk \to X$ such that for all $s,t \in \wTk$,
			\begin{equation}
				\label{eq:T-graph-embedding-eq1}
				\frac{1}{1+\nu}\max\Big\{ \Big\| \sum_{i=|u|+1}^{|s|} e_i \Big\|_A, \Big\| \sum_{i=|u|+1}^{|t|} e_i \Big\|_A\Big\} \le \norm{h_{k,\nu}(s)-h_{k,\nu}(t)}_X\le \dwTk(s,t),
			\end{equation}
			where $u$ is the greatest common predecessor of $s$ and $t$.
		\end{prop}
		
		\begin{proof}
			Let $\{(x^{(i)}_j)_{j\in \bN} \colon i\in \bN\}$ be a normalized weakly null array in $X$ that generates an asymptotic model $A=[(e_j)_{j=1}^\infty]$. Fixing $k\in\bN$ and $\delta>0$ and passing to appropriate subsequences of the array, we may assume that for any $j_1<\dots<j_k$ and any $a_1,\dots,a_k$ in $[-1,1]$ we have
			\begin{equation*}
				\abs{ \Big\| \sum_{i=1}^k a_i x^{(i)}_{j_i} \Big\|_X - \Big\| \sum_{i=1}^k a_i e_i \Big\|_A} < \delta.
			\end{equation*}
			In addition,  by applying Proposition \ref{prop:jointspreadingmodels} we may also assume that for any $i_1,\dots,i_{2k}$ in $\{1,\dots,k\}$ and any pairwise different $l_1,\dots,l_{2k}$ in $\bN$ the sequence $(x^{(i_j)}_{l_j})_{j=1}^{2k}$ is $(1+\delta)$-suppression unconditional.
			
			Let $\Phi\colon \N^{\le k}\to \bN$ be a bijection and define $h:=h_{k,\nu}\colon \wTk\to X$ by $h_{k,\nu}(s):= \sum_{u\preceq s} x_{\Phi(u)}^{(|u|)}$ for all $s \in \wTk$. Observe first that the map $h$ is $1$-Lipschitz since, after cancellations one has
			\begin{equation}
				\norm{h(s) - h(t)} = \Big\| \sum_{u\prec v \preceq s} x_{\Phi(v)}^{(|v|)} - \sum_{u\prec w \preceq t} x_{\Phi(w)}^{(|w|)}\Big\|
			\end{equation}
			where $u$ is the is the greatest common predecessor of $s$ and $t$.
			Observing that the family  $(x_{\Phi(v)}^{(|v|)})_{u\prec v \preceq s} \cup 
			(x_{\Phi(w)}^{(|w|)})_{u\prec w \preceq t}$ is $(1+\delta)$-suppression unconditional, we have
			\begin{align*}
				\norm{ h(s) - h(t) }_X & \ge \frac{1}{(1+\delta)}\max\Big\{ \Big\| \sum_{u\prec v \preceq s} x_{\Phi(v)}^{(|v|)}\Big\|, \Big\| \sum_{u\prec w \preceq t} x_{\Phi(w)}^{(|w|)}\Big\| \Big\}\\
				& \ge \frac{1}{(1+\delta)}\max\Big\{\Big\| \sum_{i=\abs{u}+1}^{\abs{s}} e_i \Big\|_A - \delta,  \Big\|\sum_{i=\abs{u}+1}^{\abs{t}} e_i \Big\|_A - \delta\Big\},    
			\end{align*}
			and choosing $\delta>0$ small enough gives the result.
		\end{proof}
		The following consequence is now immediate. 
		\begin{coro}
			\label{cor:tree-model}
			Let $X$ be an infinite-dimensional Banach space admitting an $\ell_p$-asymptotic model, generated by a normalized weakly null array,  with $1\le p <\infty$. Then, there exists a constant $C\ge 1$ such that for any $k\in\bN$ there exists a map  $h_{k} \colon \wTk \to X$ so that for all $s,t \in \wTk$,
			\begin{equation*}
				\frac{1}{C}\dwTk(s,t)^{1/p} \le \norm{h_{k}(s)-h_{k}(t)}_X\le \dwTk(s,t).
			\end{equation*}
			In particular, $\sup_{k\in\bN}\cdist{X}(\wTk)<\infty$ whenever $X$ is an infinite-dimensional Banach space admitting an $\ell_1$-asymptotic model generated by a normalized weakly null array.
		\end{coro}

		The fact that for any $p\in[1,\infty)$, a $\theta$-snowflaking of $\ell_1$ coarsely embeds into $\ell_p$ but with a worse snowflake exponent as $p$ increases (namely with $\theta=1/p$) translates into the geometry of the tree when one tries to embed the tree via an embedding that is inspired by the canonical embedding of the tree metric into $\ell_1$. Indeed, note that it follows from Corollary \ref{cor:tree-model} that $\wTk$ admits a bi-Lipschitz embedding with distortion $O(k^{1-1/p})$ into every asymptotic model, generated by a normalized weakly null array that is isomorphic to $\ell_p$. Here as well, the distortion deteriorates when $p\to \infty$. It was first observed by Bourgain in \cite{Bourgain86}, in the context of binary trees, that this distortion is not tight and can be improved by using an embedding that is inspired by the canonical embedding of the tree metric into $\co$. We now detail this improvement. 
		
		\begin{prop}\, 
			\label{prop:tree-model-dist}
			\begin{enumerate}
				\item If $X$ has an $\ell_p$-asymptotic model generated by a weakly null array for some $p\in(1,\infty)$, then
				$\cdist{X}(\wTk)=O\Big((\log k)^{\frac{1}{p}}\Big).$
				\item If $X$ has a $\co$-asymptotic model generated by a weakly null array, then $\sup_{k\in \bN}\cdist{X}\big(\wTk\big)<\infty$.
			\end{enumerate}
		\end{prop}
		
		\begin{proof}
			Let $\{(x^{(i)}_j)_{j\in \bN} \colon i\in \bN\}$ be a normalized weakly null array in $X$ that generates an $\ell_p$-asymptotic model $A=[(e_j)_{j=1}^\infty]$. Let $k\ge 1$ and fix a compatible bijection $\Phi\colon \N^{\le k}\to \bN$, meaning that $\Phi(s)\le \phi(t)$ for all $s \preceq t \in \N^{\le k}$. Assume, as we may, that for any $i_1,\ldots,i_{2k}$ in $\{1,\ldots,k\}$ and any pairwise different $l_1,\ldots,l_{2k}$ in $\bN$, the sequence $(x^{(i_j)}_{l_j})_{j=1}^{2k}$ is $(1+\delta)$-suppression unconditional for a fixed $\delta>0$ (arbitrarily small). By assumption, we may also assume that there exists a constant $D \ge 1$ such that for any $j_1<\dots<j_l$ and any $a_1,\dots,a_l$ in $[-k,k]$ we have
			\begin{equation}\label{eq:l_pmodel}
				\frac1D \big(\sum_{i=1}^l |a_i|^p\big)^{\frac1p} -\delta  \le \Big\| \sum_{i=1}^l a_i x^{(i)}_{j_i} \Big\|_X \le D \big(\sum_{i=1}^l |a_i|^p\big)^{\frac1p} + \delta .
			\end{equation}
			Letting $\frac1p +\frac1q =1$, we define the map $h_p\colon \N^{\le k} \to X$ by 
			\begin{equation*}
				h_p(s)=\sum_{v\preceq s} (|s|-|v|+1)^{\frac1q}\, x_{\Phi(v)}^{(|v|)},\ \ s\in \N^{\le k}.
			\end{equation*}
			We can think of the map $h_p$ as interpolating between $h_1$ and $h_\infty$, the two canonical embeddings of the tree into $\ell_1$ and $\co$, respectively.
			Consider $s,t \in \N^{\le k}$ and let $u$ be their greatest common predecessor. Assume $s=u\smallfrown v$, $t=u \smallfrown w$ with, for instance, $|v|=j\ge h=|w|$. We have 
			\begin{align*}
				\norm{h_p(s) - h_p(t)}_X & = \Big\| \sum_{i=0}^{|u|} \underbrace{\Big((|s|-i+1)^{\frac{1}{q}}-(|t|-i+1)^{\frac{1}{q}}\Big)}_{\alpha_i}x^{(i)}_{\Phi(u_{\restriction i})}+\\
				&   \sum_{i=1}^j \underbrace{\Big(j-i+1\Big)^{\frac{1}{q}}}_{\beta_i}x^{(|u|+i)}_{\Phi(u \smallfrown v_{\restriction i} )}-\sum_{i=1}^h \underbrace{\Big(h-i+1\Big)^{\frac{1}{q}}}_{\gamma_i} x^{(|u|+i)}_{\Phi(u \smallfrown w_{\restriction i}})\Big\|_X\\
				& \le \Big\| \sum_{i=0}^{|u|}\alpha_i x^{(i)}_{\Phi(u_{\restriction i})} \Big\|_X + \Big\| \sum_{i=1}^j \beta_i x^{(|u|+i)}_{\Phi(u \smallfrown v_{\restriction i} )} \Big\|_X+\Big\| \sum_{i=1}^h \gamma_i x^{(|u|+i)}_{\Phi(u \smallfrown w_{\restriction i}}) \Big\|_X.
			\end{align*}
			Recall that for all $y>x>0$ and $a\in(0,1)$,
			\begin{equation}\label{eq:useful-ineq2}
				y^{a}-x^{a}\le \frac{y-x}{y^{1-a}}.
			\end{equation}
			Observe now that $\max\{|u|,|u|+j,|u|+h\}\le k$ and since $\Phi$ is a compatible bijection it follows that
			\begin{align*}
				\Big\| \sum_{i=0}^{|u|}\alpha_i x^{(i)}_{\Phi(u_{\restriction i})} \Big\|_X & \stackrel{\eqref{eq:l_pmodel}}{\le} D\Big(\sum_{i=0}^{|u|} \alpha_i^p\Big)^{\frac{1}{p}} +\delta \le D\Big(\sum_{i=0}^{|u|} \Big((|s|-i+1)^{\frac{1}{q}}-(|t|-i+1)^{\frac{1}{q}}\Big)^p\Big)^{\frac{1}{p}} +\delta\\
				& \stackrel{\eqref{eq:useful-ineq2}}{\le} D\Big(\sum_{i=0}^{|u|} \Big( \frac{|s|-|t|}{(|s|-i+1)^{1-\frac{1}{q}}}\Big)^p\Big)^{\frac{1}{p}} +\delta\\
				& =D(|s|-|t|) \Big( \sum_{i=0}^{|u|} \frac{1}{|s|-i+1}\Big)^{\frac{1}{p}}+\delta \le C_1(|s|-|t|)(\log k)^{\frac{1}{p}},
			\end{align*}
			for some  constant $C_1>0$. Also,
			\begin{align*}
				\Big\| \sum_{i=1}^j \beta_i x^{(|u|+i)}_{\Phi(u \smallfrown v_{\restriction i} )} \Big\|_X & \le D\Big(\sum_{i=1}^j \beta_i^p\Big)^{\frac{1}{p}} +\delta \le  D\Big(\sum_{i=1}^j (j-i+1)^{\frac{p}{q}}\Big)^{\frac{1}{p}}+\delta \\
				& = D\Big(\sum_{i=1}^{j} i^{\frac{p}{q}}\Big)^{\frac{1}{p}}+\delta = D\Big(\sum_{i=1}^{j} i^{p-1}\Big)^{\frac{1}{p}}+\delta\le C_2 j,
			\end{align*}
			for some  constant $C_2>0$. A similar computation gives 
			\begin{equation*}
				\Big\| \sum_{i=1}^h \gamma_i x^{(|u|+i)}_{\Phi(u \smallfrown w_{\restriction i}}) \Big\|_X \le C_2 h.
			\end{equation*}
			It follows that there exists $C>0$ such that for all $s,t \in \N^{\le k}$, 
			\begin{equation*}
				\norm{h_p(s)-h_p(t)}_X \le  C(\log k)^{\frac{1}{p}}\sd_{\wTk}(s,t).
			\end{equation*}
			For the lower bound, it follows from the suppression unconditionally property of the asymptotic model that 
			\begin{align*}
				\norm{h_p(s)-h_p(t)}_X & \ge   \frac{1}{1+\delta}\Big\| \sum_{i=1}^j \beta_i x^{(|u|+i)}_{\Phi(u \smallfrown v_{\restriction i} )} \Big\|_X \ge  \frac{1}{D(1+\delta)}\Big(\sum_{i=1}^j \beta_i^p\Big)^{\frac{1}{p}} -\frac{\delta}{1+\delta}\\
				&\ge  \frac{1}{D(1+\delta)}\Big(\sum_{i=1}^{j} i^{\frac{p}{q}}\Big)^{\frac{1}{p}}-\frac{\delta}{1+\delta} \\
				&\ge C_3 j\ge \frac{C_3}{2}(j+h) =  \frac{C_3}{2}\dwTk(s,t),
			\end{align*}
			for some constant $C_3>0$. The conclusion follows.
			
			For the case $p=\infty$ we define $h_\infty \colon \N^{\le k} \to X$ via the formula 
			\begin{equation*}
				h_\infty(s):=\sum_{v \preceq s} (|s|-|v|+1)x^{(|v|)}_{\Phi(v)} = \sum_{i=0}^{|s|} \sum_{j=0}^i x^{(j)}_{s_{\restriction_j}},
			\end{equation*}
			and the argument above gives a bounded distortion. 
		\end{proof}
		
		Examples of Banach spaces to which Proposition \ref{prop:tree-model-dist} applies are $\ell_p$, $(\sum_{n=1}^{\infty} F_n)_{\ell_p}$, or any asymptotic-$\ell_p$ space such as James space (with $p=2$) for $(1)$ and asymptotic-$\co$ spaces such as Tsirelson space $\Tsi^*$ for $(2)$.
		
		\medskip A quite unexpected metric characterization of infinite dimensionality can be deduced from what we have discussed so far.
		
		\begin{theo}
			Let $X$ be a Banach space. The following assertions are equivalent:
			\begin{enumerate}[(i)]
				\item $X$ is infinite-dimensional.
				%\item $\sT_\infty^\omega$ coarsely embeds into $X$.
				\item $(\sT_k^\omega)_{k\ge1}$ equi-coarsely embeds into $X$.
			\end{enumerate}
		\end{theo}
		
		\begin{proof}
			The implication $(i)\implies (ii)$ relies on Rosenthal's $\ell_1$-theorem. Indeed, if $X$ is an infinite-dimensional Banach space, then by Rosenthal's $\ell_1$-theorem, it either contains $\ell_1$ isomorphically or it has a weakly Cauchy sequence that is not norm converging. Since $\ell_1$ contains an isometric copy of $\sT_\infty^\omega$, we may assume that $X$ contains a weakly Cauchy sequence $\xn$ which is not norm converging. After eventually passing to a subsequence of $\xn$, still denoted by $\xn$, the sequence
			\begin{equation*}
				\zn:=\Big(\frac{x_{2n-1}-x_{2n}}{\norm{x_{2n-1}-x_{2n}}}\Big)_{n=1}^\infty
			\end{equation*}
			is normalized and weakly null. Therefore, either $\zn$ has a subsequence that generates a spreading model isomorphic to $\co$ and we apply Proposition \ref{prop:tree-model-dist} or $\zn$ has a subsequence that generates a spreading model not isomorphic to $\co$ and we conclude with
			Proposition \ref{prop:tree-spreading}.
			For the implication $(ii)\implies (i)$, simply observe that by compactness, a finite-dimensional Banach space cannot contain an infinite sequence that is bounded and separated. 
		\end{proof}
		
		\subsection{Embeddability of trees and Szlenk index}
		\label{sec:tree-Szlenk}
		
		The Szlenk index is naturally connected to the existence of certain countably branching trees. Indeed, if $X$ is a separable Banach space, it follows from the metrizability of
		the weak$^*$ topology on $B_{X^*}$ that if Sz$(X,\eps)>\omega$, then for all
		$k\in \bN$ there exists $(y_{s}^*)_{s\in \N^{\le k}}$ in $B_X^*$ such that for all
		$s\in \N^{\le k-1}$ and all $n\in \N$, $\norm{ y_{s\smallfrown n}^* - y_{s}^* }\ge
		\eps/2:=\eps'$ and for all $s\in \N^{\le k-1}$, $w^*-\lim_{n} y_{s\smallfrown n}^*=y_{s}^*$.
		
		Moreover, if Sz$(X)>\omega$, it follows from the submultiplicativity of the map $\eps \mapsto \Sz(X,\eps)$ that $\Sz(X,\eps)>\omega$ for any $\eps\in (0,1)$. Therefore, in the above choice of $(y_{s}^*)_{s\in \N^{\le k}}$ we can take $\eps'=\frac{1}{3}$. By further considering $z^*_{s} := y^*_{s} - y^*_{s^-}$ for $s\in \N^{\le k} \setminus  \{\emptyset\}$ (we recall that $s^-$ is the immediate predecessor of $s$), $z^*_\emptyset := y^*_\emptyset$ and rescaling, this clearly implies the
		existence for all $k\in \bN$, of a weak*-null rooted tree $(z_{s}^*)_{s\in \N^{\le k}}$ in $X^*$ so that
		\begin{itemize}
			\item $1\le \norm{ z_{s}^* }\le 6$ for all $s\in \N^{\le k} \setminus  \{\emptyset\}$,
			\item $\norm{ \sum_{v\preceq s}z^*_v }\le 3$ for all $s\in \N^{\le k}$.
		\end{itemize}
		
		We will refer to a tree satisfying the two conditions above as a \emph{semi-normalized $\ell_\infty$-tree}.
		In the next proposition, we improve the above statement by constructing an
		almost biorthogonal system associated with the $w^*$-null semi-normalized $\ell_\infty$-tree $(z_{s}^*)_{s\in \N^{\le k}}$.
		
		\begin{prop}
			\label{prop:biorthogonal}
			Let $X$ be a separable Banach space. If Sz$(X)>\omega$, then for all $k\in
			\bN$ and $\delta>0$ there exist a weak*-null semi-normalized $\ell_\infty$-tree $(x_{s}^*)_{s\in \N^{\le k}}$ in $X^*$ and $(x_{s})_{s\in \N^{\le k}}$ in $B_X$ such that
			\begin{enumerate}
				\item $x_{s}^*(x_{s})\ge\frac{1}{3}\norm{ x^*_{s} }$ for all $s \in \N^{\le k}$,
				\item[] and
				\item $\abs{x_{s}^*(x_{t})}<\delta$ for all $s\neq t \in \N^{\le k}$.
			\end{enumerate}
			
		\end{prop}

		\begin{proof}
			The proof builds upon the following elementary claim.
			
			\begin{claim}
				\label{claim:gliding-hump}
				Let $(x_n^*)_{n=0}^\infty$ be a weak*-null sequence in $X^*$ such that $\norm{
					x_n^*}\ge 1$ for all $n$ in $\bN$ and let $F$ be a finite-dimensional
				subspace of $X^*$. Then, there exists a sequence $\xn$ in $B_X$ such that
				for all $y^*\in F$, $y^*(x_n)=0$ and $\liminf x_n^*(x_n)\ge \frac{1}{2}$
			\end{claim}
			
			\begin{proof}[Proof of Claim \ref{claim:gliding-hump}]
				It is a classical consequence of Mazur's technique for
				constructing basic sequences (see, for instance, \cite{LindenstraussTzafriri1977} or Lemma \ref{lem:codim}) that $\liminf d(x_n^*,F)\ge \frac{1}{2}$.
				Let $F_\perp:=\{x\in X \colon \forall x^*\in F~x^*(x)=0\}$ be the pre-orthogonal of
				$F$. Since $F$ is finite-dimensional, we have that $F=(F_\perp)^\perp$. Therefore, for
				any $x^*\in X^*$, $d(x^*,F)= \norm{x^*}_{X^*/F}=\norm{x^*_{\restriction F_\perp}}_{F_\perp^*}$. This finishes the proof.
			\end{proof}
			
			Let $k\in \N$ and $(z_{s}^*)_{s\in \N^{\le k}}$ be a weak*-null semi-normalized $\ell_\infty$-tree in $X^*$. We pick, as we may, a bijection $g \colon \N^{\le k} \to \bN$ satisfying the following stronger compatibility conditions: 
			\begin{itemize}
				\item $g(s)< g(t)$ for all $s \prec t \in \N^{\le k}$,
				\item $g(s\smallfrown m)< g(s \smallfrown n)$ for all $s \in \N^{\le k-1}$ and all $m<n \in \N$.
			\end{itemize}
			For $i\in \N$, we denote by $s_i=g^{-1}(i)$. In particular, $\emptyset=s_1$. We now build inductively an order and length preserving  map $\Phi \colon \N^{\le k}\to
			\N^{\le k}$ and a family $(z_{\Phi(s)})_{s \in \N^{\le k}}$ in $B_{X}$ such
			that
			\begin{equation}
				\label{eq1:biortho}
				z^*_{\Phi(s)}(z_{\Phi(s)})\ge \frac{1}{3}\|z^*_{\Phi(s)}\| \textrm{ for all }s \in \N^{\le k}
			\end{equation}
			and
			\begin{equation}
				\label{eq2:biortho}
				\abs{ z^*_{\Phi(s)}(z_{\Phi(t)}) }<\delta \textrm{ for all } s \neq t \in \N^{\le k}.
			\end{equation}
			Set $\Phi(\emptyset)=\emptyset$ and pick $z_{\Phi(\emptyset)}$ in
			$B_X$ so that $z^*_{\Phi(\emptyset)}(z_{\Phi(\emptyset)})\ge
			\frac{1}{3}\| z^*_{\Phi(\emptyset)} \|$. Assume now that
			$\Phi(s_1),\dots,\Phi(s_k)$ and $z_{\Phi(s_1)},\dots,z_{\Phi(s_k)}$ have been
			constructed and satisfy \eqref{eq1:biortho} and \eqref{eq2:biortho}. Due to the compatibility properties of $g$, there exists $i\in
			\{1,\dots,k\}$ and $p\in \bN$ such that $s_{k+1}=s_i\smallfrown p$. The sequence 
			$(z^*_{\Phi(s_i)\smallfrown n})_{n\ge 1}$ is weak$^*$ null and Claim \ref{claim:gliding-hump} ensures that we can pick $n\in \bN$  and $z_{\Phi(s_i)\frown
				n}$ in $B_X$ such that 
			\begin{equation}
				\label{eq3:biortho}
				|z^*_{\Phi(s_i)\smallfrown n}(z_{\Phi(s_j)})|<\delta \textrm{ and } z^*_{\Phi(s_j)}(z_{\Phi(s_i)\smallfrown n})=0 \textrm{ for all } j\le k
			\end{equation}
			and
			\begin{equation}
				\label{eq4:biortho}
				z^*_{\Phi(s_i)\frown n}(z_{\Phi(s_i)\smallfrown n})\ge \frac{1}{3}\|z^*_{\Phi(s_i)\frown n}\|. 
			\end{equation}
			We now set $\Phi(s_{k+1})=\Phi(s_i)\smallfrown n$, where $n$ is chosen large
			enough so that all the required properties in \eqref{eq1:biortho}, \eqref{eq2:biortho}, \eqref{eq3:biortho}, \eqref{eq4:biortho} are satisfied. We conclude our construction by setting $x^*_s=z^*_{\Phi(s)}$ and $x_s=z_{\Phi(s)}$, for $s$ in $\N^{\le k}$. These families satisfy the desired properties. 
			
		\end{proof}

		Equipped with Proposition \ref{prop:biorthogonal} we can now show that if a Banach space has a large Szlenk index, then both the space and its dual will contain the sequence $(\wTk)_{k\in \bN}$ equi-bi-Lipschitzly.  
		
		\begin{theo}
			\label{thm:tree-Szlenk}
			If $X$ is a Banach space with $\Sz(X)>\omega$, then $$\max\Big\{ \sup_{k\in \bN} \cdist{X}(\wTk), \sup_{k\in \bN} \cdist{X^*}(\wTk) \Big\} <\infty.$$
			In particular, if $X$ is reflexive and $\max\{\Sz(X),\Sz(X^*)\}>\omega$, then $\sup_{k\in \bN} \cdist{X}(\wTk)<\infty$.
		\end{theo}
		
		\begin{proof} 
			We shall first embed $\wTk$ into $X$. By Theorem \ref{thm:separable-determination-Szlenk}, we can assume without loss of generality that $X$ is separable. So, let $(x^*_{s},x_{t})_{s,t\in \N^{\le k}}$ be the almost biorthogonal system given by Proposition \ref{prop:biorthogonal}. Our choice of $\delta$ will be specified later. We now mimic the
			canonical embedding of $\wTk$ into $\ell_1(\N^{\le k})$ and define $h\colon \N^{\le k}\to X$ by
			\begin{equation*}
				h(s)= \sum_{v \preceq s}x_v, \ \ s \in \N^{\le k}.
			\end{equation*}
			Since $(x_{s})_{s \in \N^{\le k}} \subset B_X$, we clearly have that $h$ is 1-Lipschitz. Let now $s \neq t$ in $\N^{\le k}$ and let $u$ be their greatest common predecessor. Set $d:=\dwTk(u,s)=|s|-|u|$ and $d':=\dwTk(u,t)=|t|-|u|$. Recall that
			$\dwTk(s,t)= d+d'$ and assume for instance that $d\ge d'$. Then,
			\begin{align*}
				\langle \sum_{v \preceq s} x^*_v, h(s) - h(t) \rangle & = \langle \sum_{v \preceq s} x^*_v,  \sum_{u\prec v \preceq s}x_v -\sum_{u\prec w \preceq t} x_w \rangle \\
				& \ge \frac{1}{3}d-\delta \abs{s}(d+d')\ge
				\frac{d}{3}- 2\delta(k+1)d\ge \frac{1}{4}d\ge \frac{1}{8}\dwTk(s,t),
			\end{align*}
			if $\delta$ was chosen less than $\frac{1}{24(k+1)}$. Since $\|\sum_{v \preceq s} x^*_v \|\le 3$, we obtain that for all $s,t$ in
			$\N^{\le k}$,
			\begin{equation*}
				\norm{ h(s) - h(t) }\ge\frac{1}{24}\dwTk(s,t).
			\end{equation*} 
			Therefore, we have proved that $\sup_{k\in \bN} \cdist{X}(\wTk)\le 24$.
			
			We now turn to the embeddability of $\wTk$ into $X^*$. There exists a subspace $Y$ of $X$ such that $X/Y$ is separable and $\Sz(X/Y)>\omega$ (see \cite{Lancien1996}). It is enough to embed $\wTk$ into $Y^\perp \equiv (X/Y)^*$. In other words, we can assume that $X$ is separable. So, let again $(x^*_{s},x_{t})_{s,t\in \N^{\le k}}$ be the almost biorthogonal system given by Proposition \ref{prop:biorthogonal}. For the embedding into the dual space, we will mimic the canonical embedding of $\wTk$ into
			$c_0(\N^{\le k})$. For $s\in \N^{\le k}$, we set $y^*_{s}:= \sum_{v \preceq s}x^*_v$. Then, we
			define $h_* \colon \N^{\le k} \to X^*$ by
			\begin{equation*}
				h_*(s) :=\sum_{v \preceq s} y^*_{v} = \sum_{i=0}^{\abs{s}} \sum_{j=0}^i x^*_{s_{\restriction j}},\ \ s\in \N^{\le k}.
			\end{equation*}
			Since $(y^*_{s})_{s \in \N^{\le k}}$ is a subset of $3B_{X^*}$, it is immediate that $h_*$ is 3-Lipschitz. Let now $s \neq t$ in $\N^{\le k}$ and consider again $u$ their greatest
			common predecessor, $d:=\dwTk(u,s)$ and $d':=\dwTk(u,t)$. Assume, for instance, that $d\ge d'$. Then,
			\begin{align}
				\norm{ h_*(s) - h_*(t) } & = \Big\| \sum_{v \preceq s} y^*_{v} - \sum_{w \preceq s} y^*_{w} \Big\|\\
				\label{eq2:tree-Szlenk}& = \Big\| \sum_{i = \abs{u} + 1}^{\abs{s}} \sum_{j=0}^{i} x^*_{s_{\restriction j}} - \sum_{i = \abs{u} + 1}^{\abs{t}} \sum_{j=0}^{i} x^*_{t_{\restriction j}} \Big\|.
			\end{align}
			The rightmost sum might be empty if $u=t$, but in any case for all $\abs{u}+1\le i\le \abs{t}$, 
			\begin{equation*}
				\Big|\langle x_{s_{\restriction |u|+1}}, \sum_{j=0}^{i} x^*_{t_{\restriction j}} \rangle \Big| \ \le (i+1)\delta \le (k+1)\delta.
			\end{equation*}
			On the other hand, for all $\abs{u}+1\le i\le \abs{s}$,
			\begin{align*}
				\Big|  \langle x_{s_{\restriction |u|+1}}, \sum_{j=0}^{i} x^*_{s_{\restriction j}} \rangle \Big|  &\ge \langle x_{s_{\restriction |u|+1}}, x^*_{s_{\restriction |u|+1}} \rangle - \sum_{\stackrel{1\le j \le i}{j\neq \abs{u}+1}} \Big| \langle x_{s_{\restriction |u|+1}}, x^*_{s_{\restriction j}} \rangle \Big|\\
				& \ge \frac13 -\delta k.
			\end{align*}
			Therefore, 
			\begin{align*}
				\norm{h_*(s) - h_*(t)} \ge \langle x_{s_{\restriction |u|+1}}, h_*(s) - h_*(t) \rangle \ge \frac{d}{3} - \delta (k+1)(d+d')  \ge \frac{\dwTk(s,t)}{8},
			\end{align*}
			if $\delta\le \frac{1}{24(k+1)}$ and we conclude that $\sup_{k\in \bN} \cdist{X^*}(\wTk)\le 24$.
		\end{proof}
		
		In Theorem \ref{thm:tree-Szlenk}, we embedded countably branching trees of finite but arbitrarily large height. With a little bit more work, we will explain how we can, in fact, embed the countably branching tree of infinite height. We will be rather sketchy, and we refer the reader to \cite{BKL2010} for the details. To prove this result, we use a variant of the barycentric gluing technique.
		
		\begin{theo}
			\label{thm:infinite-tree-Szlenk}
			If $X$ is a Banach space with Sz$(X)>\omega$, then $$\max\Big\{\cdist{X}(\sT_\infty^\omega), \cdist{X^*}(\sT_\infty^\omega)) \Big\} <\infty.$$
			In particular, if $X$ is reflexive and $\max\{\Sz(X),\Sz(X^*)\}>\omega$, then $\cdist{X}(\sT_\infty^\omega)<\infty$.
		\end{theo}
		
		\begin{proof} 
			So assume that $\Sz(X)>\omega$ and, as we may, that $X$ is separable and fix a decreasing sequence
			$(\delta_i)_{i=0}^\infty$ in $(0,1)$. By combining the technique of Proposition
			\ref{prop:biorthogonal} and a proper enumeration of $\bigcup_{i=0}^\infty
			\{i\}\times \sT^\omega_{2^i}$, one can actually build for every $i\ge 0$:
			$(x_{i,s}^*)_{s\in \sT^\omega_{2^i}}$ in $X^*$ and $(x_{i,s})_{s\in \sT^\omega_{2^i}}$ in $B_X$
			such that
			
			(i) $\forall i\ge 0,\ \forall s\in \sT^\omega_{2^i-1},\ \ x^*_{i,s\smallfrown
				n}\stackrel{w*}{\to}0$,
			
			(ii) $\forall i\ge 0,\ \forall s\in \sT^\omega_{2^i}\setminus\{\emptyset\},\ \ \Vert
			x_{i,s}^*\Vert\ge 1$
			and \ $\forall s\in \sT^\omega_{2^i},\ \Vert \sum_{v \preceq s}x_{i,v}^*\Vert\le 3$,
			
			(iii) $\forall i\ge 0,\ \forall s\in \sT^\omega_{2^i},\ \
			x_{i,s}^*(x_{i,s})\ge\frac{1}{3}\|x^*_{i,s}\|$,
			
			(iv) $\forall (i,s)\neq (j,t),\ \ |x_{i,s}^*(x_{j,t})|<\delta_i$.
			
			\noindent Let us just emphasize the fact that the whole system
			$(x_{i,s},x^*_{j,t})_{((i,s),(j,t))}$ is almost biorthogonal. We also wish to note that
			the estimate given in (iv) depends only on $i$. This last fact relies on a
			careful application of Claim \ref{claim:gliding-hump}.
			
			\noindent For $i\ge 0$, we let $\varphi_i$ be the translate of the map defined on
			$\sT^\omega_{2^{i+1}}$ in the proof of Theorem \ref{thm:tree-Szlenk}. So, let
			$$\varphi_i(\emptyset)=0\ \ {\rm and}\ \ \varphi_i(s)=\sum_{\emptyset \prec v\preceq s} x_{i+1,v},\ \ \
			s\in \sT^\omega_{2^{i+1}}\setminus\{\emptyset\}.$$ Now we adopt the gluing technique
			introduced in \cite{Baudier2007} and also used in \cite{BaudierLancien2008} and
			build our embedding as follows. For $s\in \sT^\omega_\infty\setminus\{\emptyset\}$ there
			exists $k\in \N$ such that $2^k\le |s|< 2^{k+1}$. We define
			$$h(s)=\lambda_s
			\varphi_k(s)+(1-\lambda_s)\varphi_{k+1}(s),\ \ {\rm where}\ \
			\lambda_s=\frac{2^{k+1}-|s|}{2^{k}}.$$ 
			Of course, we set $h(\emptyset)=0$. We
			clearly have that for all $s\in \sT^\omega_\infty$, $\|h(s)\|\le |s|$ and following the
			proof of Theorem 2.1 in \cite{BaudierLancien2008} that $F$ is 9-Lipschitz.
			Consider now $s\neq t\in \sT^\omega_\infty \setminus\{\emptyset\}$ and assume for
			instance that $1\le  |t|\le  |s|$. Let $2^k\le \vert t\vert\le 2^{k+1}$ and
			$2^l\le \vert s\vert\le 2^{l+1}$ with $k\le l \in \N$. Then,
			\begin{equation*}
				h(s)-h(s')=\lambda_s
				\varphi_l(s)+(1-\lambda_s)\varphi_{l+1}(s) -(\lambda_t
				\varphi_k(t)+(1-\lambda_t)\varphi_{k+1}(t))
			\end{equation*}
			Let $u$ be the greatest common predecessor of $s$ and $t$.\\ If we set 
			$(\ast)=\langle \displaystyle\sum_{u\prec v\preceq 
				s}(x_{l+1,v}^*+x_{l+2,v}^*),h(s)-h(t)\rangle$, we get
			
			\begin{equation*}
				\begin{split}
					(\ast) & \ge \lambda_s\frac{d}{3}+(1-\lambda_s)\frac{d}{3}\\
					& -\delta_{l+1}(\lambda_s d(\vert s\vert-1)+(1-\lambda_s)d\vert
					s\vert+
					\lambda_{t}d\vert t\vert+(1-\lambda_{t})d\vert t\vert)\\
					& -\delta_{l+2}((1-\lambda_s) d(\vert s\vert-1)+\lambda_s d\vert
					s\vert+
					\lambda_{t}d\vert t\vert+(1-\lambda_{t})d\vert t\vert)\\
					& \ge \frac{d}{3}-2d\vert s\vert(\delta_{l+1}+\delta_{l+2})\\
					& \ge \frac{d}{3}-2\cdotp2^{2l+2}(\delta_{l+1}+\delta_{l+2})\ge
					\frac{d}{4}\ge \frac{\sd_{\sT^\omega_\infty}(s,t)}{8},
				\end{split}
			\end{equation*}
			if the $\delta_i$ were chosen small enough. Since $\|\sum_{v\preceq s} x^*_{i,v}\|\le 3$ \ for all $i\ge 0$, we
			deduce the following lower bound
			$$\Vert h(s)-h(t)\Vert\ge \frac{\sd_{\sT^\omega_\infty}(s,t)}{96}.$$
			If $t=\emptyset\neq s$, the argument is similar but simpler. This concludes
			our proof.
			
			\medskip In order to embed $\sT^\omega_\infty$ into $X^*$, we use exactly the same
			technique. For $i\ge 0$ and $s\in \sT^\omega_{2^i}$ set $y^*_{i,s}=\sum_{v\preceq s}
			x^*_{i,v}$ and
			$$\psi_i(\emptyset)=0\ \ {\rm and}\ \ \psi_i(s)=\sum_{\emptyset \prec v \preceq s} y^*_{i+1,v},\ \
			s\in \sT^\omega_{2^{i+1}}\setminus\{\emptyset\}.$$ Then again, we set $h_*(\emptyset)=0$
			and for $s\in \sT^\omega_\infty\setminus\{\emptyset\}$:
			$$h_*(s)=\lambda_s
			\psi_k(s)+(1-\lambda_s)\psi_{k+1}(s).$$ 
			Following again the proof in
			\cite{BaudierLancien2008}, we obtain first that $G$ is 27-Lipschitz.
			
			\noindent Consider now $s\neq t\in \sT^\omega_\infty$ such that for instance $0\le 
			|t|\le  |s|$, $2^l\le \vert s\vert\le 2^{l+1}$ and  $2^k\le \vert t\vert\le
			2^{k+1}$ with $k\le l \in \N$ or $t=\emptyset$. Let $u$ be the greatest common
			predecessor of $s$ and $t$ and $v$ be the successor of $u$ such that $v\preceq s$. In
			a very similar way, by evaluating $\langle x_{l+1,v}+x_{l+2,v},
			h_*(s) - h_*(t)\rangle$, we can show that a proper choice for the $\delta_i$
			implies that
			$$\|h_*(s)- h_*(t)\|\ge \frac{\sd_{\sT^\omega_\infty}(s,t)}{16}.$$
			This concludes the proof of this proposition.
		\end{proof}

		\section{\texorpdfstring{Poincar\'e-type inequalities for $(\beta_p)$-valued maps on trees}{Poincar\'e-type inequalities for -valued maps on trees}}
		\label{sec:tree-umbel}
		
		There are several ways to show that the distortion upper bounds from Section \ref{sec:tree-embeddings} are sharp.
		One could use a Ramsey-type argument as originally done by Matousek \cite{Matousek99} for binary trees and later by Malthaner for countably branching trees \cite{Malthaner}. Or one could use what has come to be known as the ``self-improvement argument". This type of argument for metric embedding purposes was used originally by Johnson and Schechtman \cite{JohnsonSchechtman2009} for binary diamonds, then for binary trees by Kloeckner \cite{Kloeckner2014} and finally for countably branching trees in \cite{BaudierZhang16}.
		Bourgain's original argument showing that $\cdist{X}(\sB_k)\gtrsim \log(k)^{1/q}$, where $\sB_k$ is the binary tree of height $k$,  whenever $X$ is a $q$-uniformly convex space, relied implicitly on a Poincar\'e-type inequality for $X$-valued maps on the binary trees. This inequality was refined by Tessera \cite{Tessera08} in order to provide upper bounds on the compression rates for coarse embeddings of the infinite binary tree. Bourgain's circle of ideas was taken to a very sophisticated level by Lee, Mendel, Naor and Peres in \cite{LNP2009} and \cite{MendelNaor13}. The notion of Markov convexity was introduced in \cite{LNP2009} where it was shown that if a Banach space admits an equivalent norm that is $q$-uniformly convex, then it is Markov $q$-convex. In \cite{MendelNaor13}, the converse statement was proved to hold as well, thereby producing a metric characterization of the class $\langle q-{\mathbf{UC}} \rangle$ for a given fixed $q\ge 2$. In \cite{BaudierGartland24}, several Poincar\'e-type inequalities for maps on the countably branching trees were investigated. These were shown to capture the geometry of the countably branching trees in a similar way Markov convexity does for binary trees. In this section, we discuss one such Poincar\'e-type inequality.
		This metric invariant will be sufficient to provide the right order of growth for the distortion and the compression of embeddings of the countably branching trees. First we need to introduce the following quantitative version of the property $(\beta)$ of Rolewicz. 
		
		\begin{defi} Let $p \in (1,\infty)$. We say that a Banach space $X$ has property $(\beta_p)$ if there exists a constant $c>0$ such that for all $t \in (0,2]$, $\betabar_X(t)\ge ct^p$. 
		\end{defi}
		
		Similar to the approach taken by Mendel and Naor in \cite{MendelNaor13} for Markov convexity, we establish the Poincar\'e-type inequality by first proving that a curvature-like inequality for a certain point-configuration is valid in Banach spaces with property $(\beta_p)$ and then iterating this inequality to obtain the sought-after metric invariant. 
		
		\medskip We start with a simple, but crucial, lemma that is a refinement of property $(\beta_p)$.
		
		\begin{lemm}
			\label{lem:betterbeta}
			Let $X$ be a Banach space with property $(\beta_p)$ for some $p\in (1,\infty)$. Then, there exists a constant $c>0$ such that for all $x\in B_X$ and $\{z_n\}_{n\in \bN}\subseteq B_X$, 
			\begin{equation}
				\label{eq:betterbeta}
				\inf_{n\in \bN} \Big\|\frac{x-z_n}{2}\Big\|\le 1-c\inf_{i\in\bN}\liminf_{j\to \infty}\norm{z_i-z_j}^p.
			\end{equation}
		\end{lemm}
		
		\begin{proof} Let $c>0$ so that $\betabar_X(t)\ge ct^p$ for all $t\in (0,2]$. Assume, as we may, that $\inf_{i\in\bN}\liminf_{j\in \bN}\norm{z_i-z_j}=t>0$ and hence for all $i\in \bN$ we have $\liminf_{j\to\infty}\norm{z_i-z_j}\ge t$. Let $\vep \in (0,1)$ be arbitrary. A diagonal extraction argument gives a subsequence $\{z_{n_j}\}_{j\ge 1}$ such that for all $i<j\in\bN$ it holds $\norm{z_{n_i}-z_{n_j}}\ge (1-\vep)t$. Therefore, there exists an infinite subset $\bM:=\{n_1, n_2,\dots\}$ of $\bN$ such that $\inf_{i\neq j\in\bM}\norm{z_i-z_j}\ge (1-\vep)t$. Since  $\betabar_X(t)\ge c{t^p}$, it follows from the definition of the $(\beta)$-modulus that 
			there exists $m\in \bM$ such that 
			\begin{equation*}
				\Big\|\frac{x-z_m}{2}\Big\|\le 1-c(1-\vep)^pt^p=1-c(1-\vep)^p\inf_{i\in\bN}\liminf_{j\to \infty}\norm{z_i-z_j}^p.
			\end{equation*}
			Since $\vep>0$ was arbitrary, the conclusion follows.
		\end{proof}
		
		In the next proposition, we show that Banach spaces with property $(\beta_p)$ satisfy a curvature-like inequality which plays a role similar to the $q$-fork inequality for Markov $q$-convexity.
		
		\begin{prop}
			\label{pro:ban-betap-umbel}
			Let $p\in(1,\infty)$ and $X$ be a Banach space with property $(\beta_p)$. Then, there is a constant $K\in(0,\infty)$ such that for all $w,z\in X$ and $\{x_n\}_{n \in \bN} \subset X$,
			\begin{equation}
				\label{eq:ban-betap-umbel}
				\inf_{n\in \bN} \frac{\norm{w - x_n}^p }{2^p}+ \frac{1}{K^p}\inf_{i\in\bN}\liminf_{j\in \bN}\norm{ x_i - x_j}^p\le \max\{\norm{w - z}^p, \sup_{n\in\bN} \norm{x_n - z}^p\}.
			\end{equation}
		\end{prop}
		
		\begin{proof}
			Let $w,z\in X$ and $\{x_n\}_{n \in \bN} \subset X$. Since the distance in $X$ is translation invariant, we may assume $z=0$. We may assume without loss of generality that $\sup_{n\in\bN}\norm{x_n}<\infty$ and by scale-invariance of \eqref{eq:ban-betap-umbel} we can assume that $w\in B_X$ and $\{x_n\}_{n \in \bN} \subset B_X$. Thus, \eqref{eq:ban-betap-umbel} reduces to
			\begin{equation*}
				\inf_{n \in \bN} \frac{\norm{w-x_n}^p}{2^p} + \frac{1}{K^p}\inf_{i \in \bN} \liminf_{j \to \infty} \norm{x_i-x_j}^p \le  1,
			\end{equation*}
			and this inequality follows from \eqref{eq:betterbeta} with $K^p=\frac1c$ and the fact that $\frac{\norm{w-x_n}^p}{2^p} \le  \norm{\frac{w-x_n}{2}}$ whenever $w,x_n \in B_X$.
		\end{proof}
		
		Note that inequality \eqref{eq:ban-betap-umbel} only involves the norm of differences of vectors and thus it makes sense in any metric space. We now follow a time-tested process that consists of iterating a curvature-like inequality such as \eqref{eq:ban-betap-umbel} to derive a full-fledged metric invariant. This idea goes back to Enflo and its iteration of the quadrilateral inequality to derive the roundness $2$ metric invariant, and more recently by Mendel and Naor for Markov convexity. To lighten the notation we shall denote, until the end of this chapter, by $(s,t)$ the concatenation $s \smallfrown t$ of two finite sequences $s$ and $t$.

		\begin{theo}
			\label{thm:met-betap-umbel}
			Let $p\in (1,\infty)$, $K\in(0,\infty)$ and $(M,d)$ be a metric space such that for all $w,z\in M$ and $\{x_n\}_{n\in \bN}\subseteq M$ 
			\begin{equation}
				\label{eq:met-betap-umbel}
				\inf_{n\in \bN} \frac{d(w,x_n)^p}{2^p} +\frac{1}{K^p}\inf_{i\in\bN}\liminf_{j\in \bN}d(x_i,x_j)^p\le \max\{d(w,z)^p, \sup_{n\in\bN} d(x_n,z)^p\}.
			\end{equation}
			Then for all $k\ge 1$ and all $f\colon \tree^\omega_{2^k} \to M$,
			\begin{equation*}
				\sum_{l=1}^{k-1} 
				\inf_{s\in {\bN}^{\le 2^{k}-2^l}} 
				\inf_{v \in \N^{2^l}} 
				\liminf_{j\to\infty}
				\inf_{u \in \N^{2^l-1}} \frac{d(f(s,v),f(s,j,u))^p}{2^{lp}} 
				\le K^p\Lip(f)^p.
			\end{equation*}
		\end{theo}
		
		Before giving the rather technical proof of this result, let us first detail the main applications. The Poincar\'e-type inequality \eqref{eq:betap-convexity} below is a variant of similar inequalities introduced and analyzed in \cite{BaudierGartland24}. For convenience, we will use the following terminology.
		
		\begin{defi}
			Let $p\in [1,\infty)$. We say that a metric space $(M,d)$ is \emph{$\beta_p$-convex} if there is a constant $K>0$ such that for all $k\ge 1$ and all $f\colon \tree^\omega_{2^k} \to M$,
			\begin{equation}
				\label{eq:betap-convexity}
				\sum_{l=1}^{k-1} 
				\inf_{s\in{\bN}^{\le 2^{k}-2^l}} 
				\inf_{v\in{\bN}^{2^l}} \liminf_{j\to\infty}
				\inf_{u\in{\bN}^{2^l-1}} \frac{d(f(s,v),f(s,j,u))^p}{2^{lp}} 
				\le K^p\Lip(f)^p.
			\end{equation}
		\end{defi}
		
		Combining Proposition \ref{pro:ban-betap-umbel} and Theorem \ref{thm:met-betap-umbel}, we have obtained the following result.
		
		\begin{coro}\label{cor:Banachbeta_p}
			Let $p\in(1,\infty)$. If $X$ is a Banach space with property $(\beta_p)$, then $X$ is $\beta_p$-convex.
		\end{coro}
		
		The notion of $\beta_p$-convexity naturally provides distortion lower bounds on bi-Lipschitz embeddings of the countably branching trees.
		
		\begin{coro}\label{cor:distortionbeta_p}
			Let $p\in(1,\infty)$ and $X$ be a Banach space with property $(\beta_p)$. Then, $\cdist{X}\big(\tree_{k}^\omega\big) \gtrsim (\log k)^{1/p}$.
		\end{coro}
		
		\begin{proof}
			Let $f\colon \sT^\omega_{2^k} \to X$ such that for all $s,t \in \bN^{\le 2^k}$ we have 
			\begin{equation*}
				\sd_{\sT^\omega_{2^k}}(s,t) \le \norm{f(s) -f(t) }\le D \sd_{\sT^\omega_{2^k}}(s,t).
			\end{equation*}
			It suffices to observe that for all $s\in {\bN}^{\le 2^k-2^{l}}$, $i \neq j \in \bN$ and $u,v \in {\bN}^{2^{l}-1}$, it holds that
			$\dwTk((s, i,u),(s, j , v))^p= 2^{(l+1)p}$. Since $X$ satisfies \eqref{eq:betap-convexity} for the distance associated with the norm of $X$,  we have  
			\begin{equation*}
				\sum_{l=1}^{k-1}
				\inf_{s\in{\bN}^{\le 2^{k}-2^l}}
				\inf_{i \neq j \in \bN} \inf_{v,u\in {\bN}^{2^l-1}} 
				\frac{\norm{f(s,i,u) - f(s,j,v)}^p}{2^{lp}}
				\ge (k-1)2^p,
			\end{equation*}
			which in turn implies that $D\ge \frac{2(k-1)^{1/p}}{K}$. The conclusion follows. 
		\end{proof}
		
		\begin{proof}[Proof of Theorem \ref{thm:met-betap-umbel}]
			The proof is by induction on $k$. To run the induction, we will show the following stronger statement: for all maps $f\colon {\bN}^{\le 2^k}\to M$ and all $r\in \bN$ 
			
			\begin{align*}
				\inf_{s \in \N^{2^k-1}}\frac{d(f(\emptyset),f(r,s))^p}{2^{kp}} &+ 
				\frac{1}{K^p}\sum_{l=1}^{k-1}\inf_{s\in {\bN}^{\le 2^{k}-2^l}}\inf_{v\in \N^{2^l}}\liminf_{j\to \infty}\inf_{u\in \N^{2^l-1}}\frac{d(f(s,v),f(s,j,u))^p}{2^{lp}}\\		
				& \le  \max_{s\in \N^{\le 2^k}\setminus\{\emptyset\}}d(f(s),f(s^-))^p=\Lip(f)^p.
			\end{align*}
			
			For the base case $k=1$, the inequality 
			\begin{equation*}
				\inf_{n\in\bN}\frac{d(f(\emptyset),f(i,n))^p}{2^{p}} \le \max\Big\{\sup_{n\in\bN}d(f(\emptyset),f(n))^p, \sup_{(n_1,n_2)\in {\bN}^{2}}d(f(n_1),f(n_1,n_2))^p\Big\}
			\end{equation*}
			is an immediate consequence of the triangle inequality and of the fact that $x \mapsto x^p$ is increasing on $[0,\infty)$. 
			
			We now proceed to the inductive step and fix $i\in \bN$ and $f\colon {\bN}^{\le 2^{k+1}}\to M$. Given $\vep>0$ we pick $t\in {\bN}^{\le 2^k-1}$ such that 
			
			\begin{equation*}
				\frac{d(f(\emptyset),f(i,t))^p}{2^{kp}}\le \inf_{s\in{\bN}^{2^{k}-1}}\frac{d(f(\emptyset),f(i,s))^p}{2^{kp}}+\vep,
			\end{equation*}
			and for each $r\in\bN$, choose ${u}(r)\in{\bN}^{2^{k}-1}$
			
			\begin{equation*}
				\frac{d(f(i,t),f(i,t,r,{u}(r)))^p}{2^{kp}}\le \inf_{s\in{\bN}^{2^{k}-1}}\frac{d(f(i,t),f(i,t,r,s))^p}{2^{kp}}+\vep.
			\end{equation*}
			By the induction hypothesis applied to the restriction of $f$ to ${\bN}^{\le 2^{k}}$ (and with $r=i$) we get
			\begin{align*}
				&\inf_{s\in{\bN}^{2^{k}-1}}\frac{d(f(\emptyset),f(i,s))^p}{2^{kp}}\\
				&+  \frac{1}{K^p}\sum_{l=1}^{k-1}\inf_{s\in{\bN}^{\le2^k-2^l}}\inf_{v\in \bN^{2^l}}\liminf_{j\to \infty}\inf_{u\in{\bN}^{2^l-1}}\frac{d(f(s,v),f(s,j,u))^p}{2^{lp}}\\			
				& \le  \max_{s\in \N^{\le 2^k}\setminus\{\emptyset\}}d(f(s),f(s^-))^p.
			\end{align*}
			By taking the first infimum in the sum and the maximum over larger sets, we get
			\begin{align}
				\label{eq:aux1}
				\nonumber &\inf_{s\in{\bN}^{2^{k}-1}}\frac{d(f(\emptyset),f(i,s))^p}{2^{kp}}\\ 
				&+ \frac{1}{K^p}\sum_{l=1}^{k-1}\inf_{s\in{\bN}^{\le2^{k+1}-2^l}}\inf_{v\in \bN^{2^l}}\liminf_{j\to \infty}\inf_{u\in{\bN}^{2^l-1}}\frac{d(f(s,v),f(s,j,u))^p}{2^{lp}}\\			
				\nonumber & \le  \max_{s\in \N^{\le 2^{k+1}}\setminus\{\emptyset\}}d(f(s),f(s^-))^p.
			\end{align}
			On the other hand, the induction hypothesis applied to $g(s):= f(i,t,s)$ where $s\in {\bN}^{\le 2^{k}}$ gives
			\begin{align*}
				&\inf_{s\in{\bN}^{2^{k}-1}}\frac{d(g(\emptyset),g(r,s))^p}{2^{kp}}\\ 
				&+  \frac{1}{K^p}\sum_{l=1}^{k-1}\inf_{s\in{\bN}^{\le2^k-2^l}}\inf_{v\in \bN^{2^l}}\liminf_{j\to \infty}\inf_{u\in{\bN}^{2^l-1}}\frac{d(g(s,v),g(s,j,u))^p}{2^{lp}}\\			
				& \le  \max_{s\in \N^{\le 2^{k}}\setminus\{\emptyset\}}d(g(s),g(s^-))^p.
			\end{align*}
			Observe first that 
			\begin{equation*}
				\max_{s\in \N^{\le 2^{k}}\setminus\{\emptyset\}}d(g(s),g(s^-))^p 
				\le \max_{s\in \N^{\le 2^{k+1}}\setminus\{\emptyset\}}d(f(s),f(s^-))^p
			\end{equation*}
			since we are maximizing over the set of all the edges instead of a subset of it. Also,
			\begin{align*}
				&\inf_{s\in{\bN}^{\le2^k-2^l}}  \inf_{v\in \bN^{2^l}} \liminf_{j\to \infty}\inf_{u\in{\bN}^{2^l-1}}\frac{d(g(s,v),g(s,j,t))^p}{2^{lp}}\\
				&  =  \inf_{s\in{\bN}^{\le2^k-2^l}}\inf_{v\in \bN^{2^l}}\liminf_{j\to \infty}\inf_{u\in{\bN}^{2^l-1}}\frac{d(f(i,t,s,v),f(i,t,s,j,u))^p}{2^{lp}}\\
				& \ge  \inf_{s\in{\bN}^{\le2^{k+1}-2^l}}\inf_{v\in \bN^{2^l}}\liminf_{j\to \infty}\inf_{u\in{\bN}^{2^l-1}}\frac{d(f(s,v),f(s,j,u))^p}{2^{lp}},
			\end{align*}
			since $(i,t,s)\in {\bN}^{\le2^{k+1}-2^l}$ for all $s\in{\bN}^{\le2^{k}-2^l}$. Therefore, it follows from the two relaxations above that 
			\begin{align}
				\label{eq:aux2}
				\nonumber &\inf_{s\in{\bN}^{2^{k}-1}}\frac{d(f(i,t),f(i,t,r,s))^p}{2^{kp}}\\  
				&+  \frac{1}{K^p}\sum_{l=1}^{k-1}\inf_{s\in{\bN}^{\le2^{k+1}-2^l}}\inf_{v\in \bN^{2^l}}\liminf_{j\to \infty}\inf_{u\in{\bN}^{2^l-1}}\frac{d(f(s,v),f(s,j,u))^p}{2^{lp}}\\			
				\nonumber& \le \max_{s\in \N^{\le 2^{k+1}}\setminus\{\emptyset\}}d(f(s),f(s^-))^p.
			\end{align}
			Since the sum and right hand side in \eqref{eq:aux2} do not depend on $r$, it follows from \eqref{eq:aux1} and \eqref{eq:aux2} that 
			\begin{align}
				\label{eq:aux3}
				\nonumber &\frac{1}{2^{kp}}\max\Big\{\inf_{s\in{\bN}^{2^{k}-1}}d(f(\emptyset),f(i,s))^p,\ \sup_{r\in \bN} \inf_{s\in{\bN}^{2^{k}-1}}d(f(i,t),f(i,t,r,s))^p\Big\}\\
				&+ \frac{1}{K^p}\sum_{l=1}^{k-1}\inf_{s\in{\bN}^{\le2^{k+1}-2^l}}\inf_{v\in \bN^{2^l}}\liminf_{j\to \infty}\inf_{u\in{\bN}^{2^l-1}}\frac{d(f(s,v),f(s,j,u))^p}{2^{lp}}\\			
				\nonumber &\le  \max_{s\in \N^{\le 2^{k+1}}\setminus\{\emptyset\}}d(f(s),f(s^-))^p.
			\end{align}
			If we let $w:= f(\emptyset)$, $z:= f(i,t)$ and $x_r:= f(i,t,r,{u}(r))$, it follows from how $t$ and ${u}(r)$ were chosen that
			\begin{align}
				\label{eq:aux4}
				& \frac{1}{2^{kp}} \max \{d(w,z)^p,\sup_{r\in \bN}d(z,x_r)^p\} \le \\
				\nonumber  &\frac{1}{2^{kp}}\max\Big\{\inf_{s\in{\bN}^{2^{k}-1}}d(f(\emptyset),f(i,s))^p,\ \sup_{r\in \bN} \inf_{s\in{\bN}^{2^{k}-1}}d(f(i,t),f(i,t,r,s))^p\Big\} + \vep
			\end{align}
			Then, inequality \eqref{eq:met-betap-umbel} combined with \eqref{eq:aux3} and \eqref{eq:aux4} gives 
			\begin{align}
				\label{eq:aux5}
				\nonumber &\frac{1}{2^{(k+1)p}}\inf_{r\in \bN}d(w,x_r)^p+\frac{1}{K^p}\frac{1}{2^{kp}}\inf_{r\in \bN}\liminf_{\rho \to \infty}d(x_r,x_\rho)^p\\
				&+ \frac{1}{K^p}\sum_{l=1}^{k-1}\inf_{s\in{\bN}^{\le2^{k+1}-2^l}}\inf_{v\in \bN^{2^l}}\liminf_{j\to \infty}\inf_{u\in{\bN}^{2^l-1}}\frac{d(f(s,v),f(s,j,u))^p}{2^{lp}}\\			
				\nonumber &\le  \max_{s\in \N^{\le 2^{k+1}}\setminus\{\emptyset\}}d(f(s),f(s^-))^p+\vep.
			\end{align}
			Now observe that
			\begin{eqnarray*}
				\inf_{r\in \bN}d(w,x_r)^p  =  \inf_{r\in \bN}d( f(\emptyset),f(i,t,r,{u}(r)))^p  \ge \inf_{s\in {\bN}^{\le 2^{k+1}-1}}d( f(\emptyset),f(i,s))^p,
			\end{eqnarray*}
			and
			\begin{eqnarray*}
				\inf_{r\in \bN}\liminf_{\rho \to \infty}d(x_r,x_\rho)^p & 
				= & \inf_{r\in \bN}\liminf_{\rho \to \infty}d(f(i,t,r,{u}(r)),f(i,t,\rho, {u}(\rho)))^p \\
				& \ge  & \inf_{r\in \bN}\liminf_{\rho \to \infty}\inf_{u\in {\bN}^{2^k-1}}d(f(i,t,r,{u}(r)),f(i,t,\rho,u))^p\\
				& \ge & \inf_{v\in {\bN}^{2^k}}\liminf_{\rho \to \infty}\inf_{u\in {\bN}^{2^k-1}}d(f(i,t,v),f(i,t,\rho,u))^p\\
				& \ge & \inf_{s\in {\bN}^{\le 2^k}} \inf_{v\in {\bN}^{2^k}}\liminf_{\rho \to \infty}\inf_{u\in {\bN}^{2^k-1}}d(f(s,v),f(s,\rho,u))^p
			\end{eqnarray*}
			Plugging in the two relaxed inequalities above in \eqref{eq:aux5} we obtain
			\begin{align*}
				&\inf_{s\in {\bN}^{\le 2^{k+1}-1}}\frac{d( f(\emptyset),f(i,s))^p}{2^{(k+1)p}}\\
				&+\frac{1}{K^p}\inf_{s\in {\bN}^{\le 2^k}} \inf_{v\in {\bN}^{2^k}}\liminf_{\rho \to \infty}\inf_{u\in {\bN}^{2^k-1}}\frac{d(f(s,v),f(s,\rho,u))^p}{2^{kp}}\\
				&+ \frac{1}{K^p}\sum_{l=1}^{k-1}\inf_{s\in{\bN}^{\le2^{k+1}-2^l}}\inf_{v\in \bN^{2^l}}\liminf_{j\to \infty}\inf_{u\in{\bN}^{2^l-1}}\frac{d(f(s,v),f(s,j,u))^p}{2^{lp}}\\			
				&\le  \max_{s\in \N^{\le 2^{k+1}}\setminus\{\emptyset\}}d(f(s),f(s^-))^p+\vep.
			\end{align*}
			and hence 
			\begin{align*}
				&\inf_{s\in {\bN}^{\le 2^{k+1}-1}}\frac{d( f(\emptyset),f(i,s))^p}{2^{(k+1)p}}\\ 
				&+ \frac{1}{K^p}\sum_{l=1}^{k}\inf_{s\in{\bN}^{\le2^{k+1}-2^l}}\inf_{v\in \bN^{2^l}}\liminf_{j\to \infty}\inf_{u\in{\bN}^{2^l-1}}\frac{d(f(s,v),f(s,j,u))^p}{2^{lp}}\\			
				&\le  \max_{s\in \N^{\le 2^{k+1}}\setminus\{\emptyset\}}d(f(s),f(s^-))^p+\vep.
			\end{align*}
			Since $\vep>0$ is arbitrary, the induction step is completed.
			
		\end{proof}

		\begin{rema}
			\label{rem:tree-dist} 
			For all $p\in(1,\infty)$, if $K_p(k)$ is the best constant such that \eqref{eq:betap-convexity} holds, then one can easily show that $K_p(k)\ge 2(k-1)^{1/p}$ and hence $\cdist{M}\big(\tree_{k}^\omega\big)=\Omega \big( (\log k)^{1/p} \big)$ for every metric space $M$ that is $\beta_p$-convex.
		\end{rema}
		
		As cleverly observed by Tessera \cite{Tessera08}, Poincar\'e-type inequalities such as the $\beta_p$-convexity inequality \eqref{eq:betap-convexity} also provide upper bounds on the compression rate of coarse embeddings of trees.
		
		\begin{theo} 
			\label{thm:tree-compression}
			The compression rate $\rho$ of any equi-coarse embedding of $(\wTk)_{k\ge 1}$ into a Banach space with property $(\beta_p)$ satisfies  
			\begin{equation}
				\label{eq:comp-bound}
				\int_{1}^\infty \Big(\frac{\rho(t)}{t}\Big)^p\frac{dt}{t}<\infty.
			\end{equation}
		\end{theo}
		
		\begin{proof}
			Assume that $X$ has property $(\beta_p)$ with constant $K$ and that for all $k\ge 1$ there is a map $f_k\colon \sT_{2^k}^\omega\to X$ such that for all $s,t \in \sT^\omega_{2^k}$,
			\begin{equation}
				\label{eq:tree-compression}
				\rho(\sd_{\sT^\omega_{2^k}}(s,t))\le \norm{ f_k(s) - f_k(t)}\le \omega(\sd_{\sT^\omega_{2^k}}(s,t)),
			\end{equation}
			Then,
			\begin{align*}
				\sum_{l=1}^{k-1} & \inf_{s\in{\bN}^{\le 2^{k}-2^l}} \inf_{i \neq j \in \bN} \inf_{u,v\in{\bN}^{2^l-1}} \frac{\norm{f_k(s,i,u)-f_k(s,j,v)^p}}{2^{lp}} \ge\\
				\sum_{l=1}^{k-1} & \inf_{s\in{\bN}^{\le 2^{k}-2^l}} \inf_{i \neq j \in \bN} \inf_{u,v\in{\bN}^{2^l-1}} \frac{\rho(\sd_{\sT^\omega_{2^k}}((s,i,u) , (s,j,v) ) )^p}{2^{lp}} \ge
				\sum_{l=1}^{k-1}\frac{\rho(2^{(l+1)})^p}{2^{lp}},
			\end{align*}
			and hence it follows from \eqref{eq:betap-convexity}, \eqref{eq:tree-compression} and the fact that $\rho$ is nondecreasing that
			\begin{equation*}
				\sum_{l=1}^{k-1}\frac{\rho(2^{l})^p}{2^{lp}} \le K^p \omega(1)^p.
			\end{equation*} 
			But, 
			\begin{eqnarray*}
				\int_{2^{l-1}}^{2^{l}} \frac{\rho(t)^p}{t^{p}}\frac{dt}{t} & \le & \rho(2^{l})^p \int_{2^{l-1}}^{2^{l}}\frac{dt}{t^{p+1}}= \rho(2^{l})^p\frac{2^{-(l-1)p}-2^{-lp}}{p} =  \frac{2^p-1}{p}\frac{\rho(2^{l})^p}{2^{lp}},
			\end{eqnarray*}
			and hence 
			\begin{equation*}
				\int_{1}^{2^{k-1}}\frac{\rho(t)^p}{t^{p+1}}dt=\sum_{l=1}^{k-1}\int_{2^{l-1}}^{2^{l}} \frac{\rho(t)^p}{t^{p}}\frac{dt}{t}\le \frac{2^p-1}{p} \sum_{l=1}^{k-1}\frac{\rho(2^{l})^p}{2^{lp}}\le \frac{2^p-1}{p} K^p \omega(1)^p<\infty. 
			\end{equation*}
		\end{proof}
		
		\begin{rema}
			The proof of Theorem \ref{thm:tree-compression} actually shows the following purely metric statement. Let $p\in(1,\infty)$ and assume that $(M,d_M)$ is $\beta_p$-convex with constant $K$. If there are nondecreasing maps $\rho,\omega\colon [0,\infty)\to [0,\infty)$ and for all $k\ge 1$ a map $f_k\colon \wTk\to M$ such that for all $s,t\in \wTk$,
			\begin{equation*}
				\rho(\dwTk(s,t))\le d_M(f_k(s),f_k(t))\le \omega(\dwTk(s,t)),
			\end{equation*}
			then
			\begin{equation*}
				\int_{1}^\infty \Big(\frac{\rho(t)}{t}\Big)^p\frac{dt}{t} \le \frac{2^p-1}{p}K^p \omega(1)^p.
			\end{equation*}
		\end{rema}

		\section{Applications to rigidity results and metric characterizations}
		\label{sec:tree-applications}

		The following lemma can be deduced from one of James' characterizations of reflexivity and was explicitly stated and proved in \cite[Lemma 3.0.1]{DKR2016}.
		
		\begin{lemm}
			\label{lem:James}
			If $X$ is a nonreflexive Banach space, then for every $\theta\in(0,1)$, there exists a $1$-Lipschitz map $h\colon \sT_\infty^\omega\to X$ such that for all $u,s,t \in \N^{<\omega}$ with $u$ being the greatest common predecessor of $s$ and $t$ and  $s_i<t_j$ for all $|u|+1\le i \le |s|$ and all $|u|+1 \le j \le |t|$, we have 
			\begin{equation}
				\label{eq:James-sep}
				\norm{h(s) - h(t)} \ge  \frac{\theta}{3}\sd_{\sT_\infty^\omega}(s,t) 
			\end{equation}
		\end{lemm}
		
		\begin{proof}
			Let $\theta\in(0,1)$. James' characterization of reflexivity \cite{James1963/64} provides sequences $\xn\subseteq B_X$ and $\xnstar \subseteq B_{X^*}$ such that 
			$$x_n^*(x_k)= \begin{cases}
				\theta \textrm{ if }n\le  k,\\
				0 \textrm{ if } n>k.
			\end{cases}
			$$
			Define $h\colon \N^{<\omega} \to X$ by $h(\emptyset)=0$ and $h(s):=\sum_{i=1}^{|s|} x_{s_i}$ for $s\neq \emptyset$. It is immediate that $h$ is $1$-Lipschitz. Let $s,t \in \N^{<\omega}$ be as in the statement of the lemma. Then, 
			\begin{equation*}
				\norm{h(s) - h(t)} = \Big\|\sum_{i=|u|+1}^{\abs{s}} x_{s_i} - \sum_{j=|u|+1}^{\abs{t}} x_{t_j}\Big\|.
			\end{equation*}
			Therefore, if we set $d:=\abs{s}-l$, $d':=\abs{t}-l$ and  $r_1=\min\{s_i,\ |u|+1\le i \le |s|\}$ and $r_2=\min\{t_j,\ |u|+1\le j \le |t|\}$, we get 
			\begin{align*}
				\norm{h(s) - h(t)} & \ge \max\big\{|x^*_{r_1}(h(s) - h(t))|, |x^*_{r_2}(h(s) - h(t))|\big\}\\ 
				& \ge \max\{\theta|d-d'|, \theta d'\}.
			\end{align*}
			The conclusion follows from observing that $\max\{|d-d'|, d'\}\ge \frac{1}{3}(d+d')=\frac13\sd_{\sT_\infty^\omega}(s,t) $.
		\end{proof}
		
		The attentive reader will have noticed that all the results involving the Poincar\'e-type inequality \eqref{eq:betap-convexity} only needed a weaker Poincar\'e-type inequality that could be derived from \eqref{eq:met-betap-umbel}, where $\inf_i \liminf_j$ would have been replaced by $\inf_{i\neq j}$. This weaker metric invariant, called infrasup $p$-umbel convexity in \cite{BaudierGartland24}, is sufficient for quantitative embedding purposes (see Exercise \ref{ex:infrasup-umbel}). However, we need the stronger metric invariant if one wants to show that it implies reflexivity in the Banach space setting.
		
		\begin{prop}
			\label{prop:reflexivity}
			Let $X$ be a Banach space that is $\beta_p$-convex for some $p\in(1,\infty)$, then $X$ is reflexive.
		\end{prop}
		
		\begin{proof}
			Assume that there exists $p\in(1,\infty)$ such that $X$ is $\beta_p$-convex with constant $K\in(0,\infty)$ but $X$ is not reflexive and fix $k \in \N$. Consider $[\bN]^{\le 2^k}=\{\emptyset\}\cup\{\nbar=(n_1,\dots,n_i),\ i\le 2^k,\ n_1<\dots<n_i \in \N\}$ as a subset of $\sT_{\infty}^\omega$. Let $\Phi:\N^{\le 2^k} \to [\bN]^{\le 2^k}$ be a bijection, which is also a pruning of $\N^{\le 2^k}$ (we refer to Section \ref{sec:summable-A} for the definition of prunings and leave it to the readers to convince themselves of the existence of such a map $\Phi$). This induces a natural isometry from $\sT_{2^k}^\omega$ onto $([\bN]^{\le 2^k},d_{\sT_{\infty}^\omega}$). Consider now the restriction to $[\bN]^{\le 2^k}$ (still viewed as a subset of $\sT_{\infty}^\omega$) of the map $h$ from Lemma \ref{lem:James}. Then, for all $1\le l \le k-1$ and all $\nbar\in \wtree{\bN}{\le 2^k-2^{l}}$, $\bar{\delta}=(\delta_1,\dots,\delta_{2^l})\in \wtree{\bN}{2^{l}}$, $j \in \bN$ and $\bar{\eta} \in \wtree{\bN}{2^{l}-1}$ such that $(\nbar , \bar{\delta}), (\nbar , j , \bar{\eta})\in\wtree{\bN}{\le 2^k}$, it follows from \eqref{eq:James-sep} that if $j>\delta_{2^l}$, then 
			$$\|h(\nbar, \bar{\delta})-h(\nbar , j , \bar{\eta})\|^p\ge \frac{\theta^p}{3^p}2^{(l+1)p}.$$
			Therefore,  
			\begin{equation*}
				\sum_{l=1}^{k-1}\inf_{\nbar\in\wtree{\bN}{\le 2^k-2^{l}}}\inf_{\stackrel{\bar{\delta}\in\wtree{\bN}{2^l}\colon}{(\nbar,\bar{\delta})\in \wtree{\bN}{\le 2^k}}}\liminf_{j\to\infty}\inf_{\stackrel{\bar{\eta}\in\wtree{\bN}{2^l-1}\colon}{(\nbar,j,\bar{\eta})\in \wtree{\bN}{\le 2^k}}}\frac{\norm{h(\nbar,\bar{\delta})-h(\nbar,j,\bar{\eta})}^p}{2^{lp}} \ge (k-1)\Big(\frac{2\theta }{3}\Big)^p,
			\end{equation*}
			Finally, we apply the assumption of $\beta_p$-convexity with constant $K$  to the $1$-Lipschitz map $h \circ \Phi$ and get that $K^p\ge (k-1)\Big(\frac{2\theta }{3}\Big)^p$ for all $k\ge 1$; a contradiction.
		\end{proof}
		
		We have now established all the tools needed to derive purely metric characterizations of certain classes of Banach spaces. 
		The first corollary contains a metric characterization of the class $\langle \BETA \rangle$ of Banach spaces that admit an equivalent norm with property $(\beta)$.
		\begin{coro}
			\label{cor:metric-beta}
			Let $X$ be a Banach space. The following assertions are equivalent.
			\begin{enumerate}[(i)]
				\item $X$ admits an equivalent norm with property $(\beta)$.
				\item $X$ admits an equivalent norm with property $(\beta_p)$ for some $p\in(1,\infty)$.
				%\item $X$ is umbel $p$-convex for some $p\in(1,\infty)$.
				\item $X$ is $\beta_p$-convex for some $p\in(1,\infty)$.
			\end{enumerate}
		\end{coro}
		
		\begin{proof} $(i) \Rightarrow (ii)$ is ensured by Theorem \ref{thm:beta-renorming}. 
			
			$(ii) \Rightarrow (iii)$ is Corollary \ref{cor:Banachbeta_p}.
			
			$(iii) \Rightarrow (i)$. Assume that $X$ is $\beta_p$-convex for some $p\in(1,\infty)$. Then, by Proposition \ref{prop:reflexivity}, $X$ is reflexive and by Corollary \ref{cor:distortionbeta_p}, the $\sT_k^\omega$ do not equi-Lipschitz embed into $X$. Thus, it follows from Theorem \ref{thm:tree-Szlenk} that $\max\{\Sz(X),\Sz(X^*)\}\le \omega$. Finally, we apply Theorem \ref{thm:beta-renorming} to get that $X$ admits an equivalent norm with property $(\beta)$.
		\end{proof}
		
		In the next corollary, we provide a metric characterization of $\langle \BETA \rangle$ in terms of a submetric test space in Ostrovskii's terminology \cite{Ostrovskii2014}. This corollary was also independently discovered by S. Zhang \cite{Zhang2022} via a different approach and was also implicit in \cite{DKR2016}.
		
		\begin{coro}
			\label{cor:submetric-beta}
			A Banach space $X$ does not admit an equivalent norm with property $(\beta)$ if and only if there exist a constant $A>0$ and $1$-Lipschitz map  $h\colon \sT_\infty^\omega\to X$ such that for all $u,s,t \in \N^{<\omega}$ with $u$ being the greatest common predecessor of $s$ and $t$ and such that $s_i<t_j$ for all $|u|+1\le i \le |s|$ and all $|u|+1 \le j \le |t|$, we have 
			\begin{equation*}
				\norm{h(s) - h(t)} \ge  A\sd_{\sT_\infty^\omega}(s,t) 
			\end{equation*}
		\end{coro}
		
		\begin{proof}
			Assume that $X$ does not admit an equivalent norm with property $(\beta)$. If $X$ is reflexive, then, by Theorem \ref{thm:beta-renorming} it must have $\max\{\Sz(X),\Sz(X^*)\}>\omega$ and by Theorem \ref{thm:infinite-tree-Szlenk} it contains a bi-Lipschitz copy of $\sT_\infty^\omega$, the countably branching tree of infinite height and the condition is clearly satisfied.
			If $X$ is not reflexive, then we can take the map from Lemma \ref{lem:James}. Assuming now that $X$ admits an equivalent norm with property $(\beta)$, then $X$ is $\beta_p$-convex for some $p\in(1,\infty)$. It remains to observe that the proof of Proposition \ref{prop:reflexivity} shows that there cannot exist an $X$-valued map satisfying the conditions listed in the statement of Corollary \ref{cor:submetric-beta}.
		\end{proof}
		
		Historically, the corollary below, which contains a metric characterization of the class $\langle \AUC \& \AUS \rangle$ \emph{within the class of reflexive Banach spaces}, was proven in \cite{BKL2010}. This result was the first purely metric characterization of a class of Banach spaces defined in terms of asymptotic properties. It is in fact, an asymptotic analog of Bourgain's foundational result in the Ribe program. 
		\begin{coro}
			Let $X$ be a reflexive Banach space. The following assertions are equivalent:
			\begin{enumerate}[(i)]
				\item $X$ admits an equivalent norm that is asymptotically uniformly smooth and an equivalent norm that is asymptotically uniformly convex.
				\item $X$ admits an equivalent norm that is asymptotically uniformly smooth and asymptotically uniformly convex.
				\item $\max\{\Sz(X),\Sz(X^*)\}\le \omega$
				\item $X$ admits an equivalent norm with property $(\beta)$. 
				\item $\sup_{k\in \bN} \cdist{X}(\wTk)=\infty$.
				\item $\cdist{X}(\sT_\infty^\omega)=\infty$
			\end{enumerate}
		\end{coro}
		
		\begin{proof}
			The equivalence between $(i)$, $(ii)$, $(iii)$ and $(iv)$ comes from Theorem \ref{thm:beta-renorming}. $(v) \Rightarrow (vi)$ is trivial and $(vi) \Rightarrow (iii)$ follows from Theorem \ref{thm:infinite-tree-Szlenk}. Finally,  $(iv) \Rightarrow (vi)$ follows from Corollary \ref{cor:distortionbeta_p}.
		\end{proof}
		
		\begin{rema}
			Y. Perreau \cite{Perreau2020} improved the implication $(iii) \Rightarrow (v)$ above by showing that if a quasi-reflexive Banach space $X$ satisfies $\max\{\Sz(X),\Sz(X^*)\}\le \omega$, then  $\sup_{k\in \bN} \cdist{X}(\wTk)=\infty$. This applies to the quasi-reflexive (nonreflexive)  James space $\James_2$. Therefore, there is a Banach space $X\notin \langle \BETA \rangle$ such that $\sup_{k\in \bN} \cdist{X}(\wTk)=\infty$ and $\cdist{X}(\sT_\infty^\omega)=\infty$ and those metric invariants do not characterize $\langle \BETA \rangle$.
		\end{rema}
		
		It leaves open the following important question.
		
		\begin{prob} 
			Assume that $X$ is a Banach space such that  $\max\{\Sz(X),\Sz(X^*)\}\le \omega$. Does it imply that $\sup_{k\in \bN} \cdist{X}(\wTk)=\infty$ or $\cdist{X}(\sT_\infty^\omega)=\infty$? 
		\end{prob}
		
		Note that if follows from Corollary \ref{cor:metric-beta} that the class $\langle \BETA \rangle$ is stable under bi-Lipschitz embedding.
		However, this result can easily be obtained using differentiability arguments since property $(\beta)$ implies reflexivity. We shall now discuss a much more delicate stability result, namely that the class $\langle \BETA \rangle$ is stable under coarse, uniform, or Lipschitz quotients. For simplicity, we will only treat the case of uniform quotients and refer to \cite{BaudierGartland24} for more general stability results.
		
		\medskip Recall that a map $f\colon (M,d_M)\to (N,d_N)$ between metric spaces is called a \emph{uniform quotient map} and $N$ is simply said to be a \emph{uniform quotient} of $M$, if $f$ is surjective, uniformly continuous and \emph{co-uniformly continuous}, i.e. for every $r>0$ there exists $\delta(r)>0$ such that for all $x\in M$, one has
		\begin{equation*}
			\label{eq:uniquotient}
			B_{N}\left(f(x),\delta(r)\right)\subset f(B_M(x,r)).
		\end{equation*} 
		It is a standard fact that a co-uniformly continuous map into a connected space is surjective. We now state and prove the last result of this chapter, which was originally due to Dilworth, Kutzarova and Randrianarivony \cite{DKR2016}. The proof we present uses $\beta_p$-convexity and is taken from \cite{BaudierGartland24}. 
		
		\begin{theo}
			\label{thm:beta-uniform-quotient}
			Let $X$ be a Banach space and assume that $Y$ is a Banach space and a uniform quotient of $X$. 
			\begin{enumerate}[(i)]
				\item Let $p \in (1,\infty)$ and assume that $X$ is $\beta_p$-convex. Then, $Y$ is is $\beta_p$-convex.
				\item Assume that $X$ has an equivalent norm with property $(\beta)$. Then, so does $Y$. 
			\end{enumerate}
		\end{theo}
		
		\begin{proof} According to Corollary \ref{cor:metric-beta}, it is enough to prove statement $(i)$. So, assume that $X$ is $\beta_p$-convex, for some $p\in (1,\infty)$ and let $f\colon X \to Y$ be a uniform quotient map. Since $X$ is metrically convex, it is a standard fact that the map $f\colon X \to Y$ is coarse-Lipschitz, i.e. there exist $L>0$ and $A\ge 0$ such that for all $x,y\in X$ 
			\begin{equation}
				\label{eq:LLD}
				\norm{f(x) - f(y)}_Y\le L\norm{x - y}_X + A.
			\end{equation}
			Moreover, since $Y$ is also metrically convex, it is also true that $f$ is coarse co-Lipschitz and this implies that there are constants $C, K>0$ such that for all $x \in X$ and  $r>0$,
			\begin{equation}
				\label{eq:coLLD}
				B_{Y}\Big(f(x),\frac{r}{C}\Big)\subset f(B_X(x,r))+KB_Y.
			\end{equation} 
			
			We start with an elementary claim.
			
			\begin{claim}
				\label{claim:general-lifting}
				Let $f\colon X \to Y$ be a surjective map satisfying \eqref{eq:coLLD} and $g\colon {\bN}^{\le m} \to Y$ be an arbitrary map. Then, there is a map $h\colon {\bN}^{\le m} \to X$ such that 
				\begin{equation}
					\label{eq:gen-upper-lift}
					\forall s\in \N^{\le m}\setminus \{\emptyset\},\ \|h(s)-h(s^-)\|_X\le C\|g(s)-g(s^-)\|_Y+CK
				\end{equation}
				and 
				\begin{equation}
					\label{eq:gen-id-lift}
					\forall s \in \N^{\le m},\ \|f(h(s))-g(s)\|_Y \le K.
				\end{equation}
			\end{claim}
			
			\begin{proof}[Proof of Claim \ref{claim:general-lifting}]
				The proof is a simple induction on $m$. If $m=0$, by surjectivity of $f$, there is $x\in X$ such that $g(\emptyset)= f(x)$ and we let $h(\emptyset):= x$. Since $f\circ h = g$ it makes \eqref{eq:gen-id-lift} true and the condition \eqref{eq:gen-upper-lift} is vacuously true. Assume now that the result is true for $m\in \N$, let $g:\N^{\le m+1}\to Y$ and let $h:\N^{\le m}\to X$ be obtained by the induction hypothesis applied to the restriction of $g$ to $\N^{\le m}$. We will now define $h(s)$ for a given $s\in \N^{m+1}$. Set $r:=\|g(s)-g(s^-)\|_Y$. Since $\|f(h(s^-))-g(s^-)\|_Y\le K$, we have 
				$$g(s)\in B_Y(f(h(s^-)),r+K)\subset f(B_X(h(s^-),C(r+K)))+KB_Y.$$
				Let $z \in B_X(h(s^-),C(r+K))$ such that $\|f(z)-g(s)\|_Y\le K$ and set $h(s):=z$. Then, $\|f(h(s))-g(s)\|_Y\le K$ and $\|h(s)-h(s^-)\|_X\le C(r+K)=C\|g(s)-g(s^-)\|_Y+CK.$
			\end{proof}
			
			We need to show that there exists a constant $\Gamma>0$ such that for every map $g\colon  \sT_{2^k}^\omega \to Y$,
			\begin{align*}
				\sum_{l=1}^{k-1}
				\inf_{s\in{\bN}^{\le 2^{k}-2^l}} 
				\inf_{v \in{\bN}^{2^l}}\liminf_{j\to\infty}\inf_{u\in{\bN}^{2^l-1}}
				\frac{\norm{ g(s,v) - g(s,j,u)}_Y^p}{2^{lp}} 
				\le \Gamma^p \Lip(g)^p.
			\end{align*} 
			Observe that if $\Lip(g)=0$, then the left-hand side vanishes as well and there is nothing to prove. Then, by scale-invariance of the inequality above, we may assume that $\Lip(g)=1$. Let $\Pi>0$ such that $X$ is $\beta_p$-convex with constant $\Pi$. Then, for $h\colon {\bN}^{\le 2^k} \to X$, where $h$ is the lifting of $g$ with respect to $f$ as defined in Claim \ref{claim:general-lifting}, we have
			\begin{align}
				\label{eq:h}
				\sum_{l=1}^{k-1}
				\inf_{s\in{\bN}^{\le 2^{k}-2^l}} 
				\inf_{v \in{\bN}^{2^l}}\liminf_{j\to\infty}\inf_{u\in{\bN}^{2^l-1}}
				\frac{\norm{ h(s,v) - h(s,j,u)}_Y^p}{2^{lp}} 
				\le \Pi^p \Lip(h)^p.
			\end{align} 
			It follows from \eqref{eq:gen-upper-lift} and the fact that $d_{\sT_{2^k}^\omega}$ is a graph metric that $\Lip(h)\le C\Lip(g)+CK=C(K+1).$ So
			\begin{equation}
				\label{eq:h2}
				\Pi^p \Lip(h)^p \le \Pi^pD^p, \ \text{where}\ D=C(K+1).
			\end{equation} 
			Let $1\le l\le k-1$. Then, either
			\begin{equation}\label{eq:gen-case1}
				\inf_{s\in{\bN}^{\le 2^{k}-2^l}} 
				\inf_{v \in{\bN}^{2^l}}\liminf_{j\to\infty}\inf_{u\in{\bN}^{2^l-1}}
				\norm{ g(s,v) - g(s,j,u)}_Y \le  L+A+ 2K
			\end{equation}
			or
			\begin{equation}
				\label{eq:gen-case2}
				\inf_{s\in{\bN}^{\le 2^{k}-2^l}} 
				\inf_{v \in{\bN}^{2^l}}\liminf_{j\to\infty}\inf_{u\in{\bN}^{2^l-1}}
				\norm{ g(s,v) - g(s,j,u)}_Y >  L+A+ 2K.
			\end{equation}
			If \eqref{eq:gen-case2} holds, then for all $s\in{\bN}^{\le 2^{k}-2^l}$ and $v\in{\bN}^{2^l}$, we can find infinitely many $j$ such that for all $u\in{\bN}^{2^l-1}$ we have 
			\begin{equation*}
				\norm{ g(s,v) - g(s,j,u)}_Y> L + A + 2K.
			\end{equation*}
			It follows from the triangle inequality and \eqref{eq:gen-id-lift} that
			\begin{equation*}
				\norm{ f(h(s,v)) - f(h(s,j,u))}_Y \ge \norm{ g(s,v) - g(s,j,u)}_Y -2K \ge L+A.   
			\end{equation*}
			Observe now that $\norm{f(x) - f(y)}_Y< L+A$ whenever $\norm{x - y}_X < 1$ and based on the inequality above, necessarily $\norm{ h(s,v) - h(s,j,u) }_X\ge 1$. Thus, in this case,
			\begin{align*}
				\norm{ g(s,v) - g(s,j,u) }_Y & \stackrel{\eqref{eq:gen-id-lift}}{\le} \norm{f(h(s,v)) - f(h(s,j,u))}_Y + 2K\\
				& \stackrel{\eqref{eq:LLD}}{\le} L \norm{ h(s,v) - h(s,j,u) }_X + A + 2K\\
				& \le (L+ A + 2K)\norm{ h(s,v) - h(s,j,u) }_X.
			\end{align*}
			Then, letting $\gamma\eqd (L + A + 2K)$ for simplicity,
			\begin{align}
				\label{eq:gen-auxquot}
				\nonumber &\inf_{s\in{\bN}^{\le 2^{k}-2^l}} 
				\inf_{v \in{\bN}^{2^l}}\liminf_{j\to\infty}\inf_{u\in{\bN}^{2^l-1}}
				\frac{\norm{ g(s,v) - g(s,j,u)}_Y^p}{2^{lp}}\\
				& \le \gamma^p \inf_{s\in{\bN}^{\le 2^{k}-2^l}} 
				\inf_{v \in{\bN}^{2^l}}\liminf_{j\to\infty}\inf_{u\in{\bN}^{2^l-1}}
				\frac{\norm{ h(s,v) - h(s,j,u)}_Y^p}{2^{lp}}.
			\end{align}
			Therefore
			\begin{align*}
				&\sum_{l=1}^{k-1}
				\inf_{s\in{\bN}^{\le 2^{k}-2^l}} 
				\inf_{v \in{\bN}^{2^l}}\liminf_{j\to\infty}\inf_{u\in{\bN}^{2^l-1}}
				\frac{\norm{ g(s,v) - g(s,j,u)}_Y^p}{2^{lp}} \stackrel{\eqref{eq:gen-case1} \land \eqref{eq:gen-auxquot}}{\le}\\
				& \sum_{l=1}^{k-1}
				\max\Big\{\frac{\gamma^p}{2^{lp}}, \gamma^p 
				\sum_{l=1}^{k-1}
				\inf_{s\in{\bN}^{\le 2^{k}-2^l}} 
				\inf_{v \in{\bN}^{2^l}}\liminf_{j\to\infty}\inf_{u\in{\bN}^{2^l-1}}
				\frac{\norm{ h(s,v) - h(s,j,u)}_X^p}{2^{lp}}\Big\}\\
				& \le \gamma^p \sum_{l=1}^{k-1}
				\inf_{s\in{\bN}^{\le 2^{k}-2^l}} 
				\inf_{v \in{\bN}^{2^l}}\liminf_{j\to\infty}\inf_{u\in{\bN}^{2^l-1}}
				\frac{\norm{ h(s,v) - h(s,j,u)}_X^p}{2^{lp}} + \sum_{l=1}^{k-1} \frac{\gamma^p}{2^{lp}}\\
				& \stackrel{\eqref{eq:h} \land \eqref{eq:h2}}{\le}  \gamma^p\Pi^pD^p + \sum_{l=1}^{\infty} \frac{\gamma^p}{2^{lp}} < \infty,
			\end{align*} 
			which concludes the proof since $ \Lip(g)=1$ and the constant $\gamma^p\big(\Pi^pD^p + \sum_{l=1}^{\infty} \frac{1}{2^{lp}}\big)$ is independent of $k$ and $g$.
		\end{proof}
		
		\begin{rema}
			\label{rem:beta-uniform-quotient}
			The proof can be slightly modified to show that the same conclusion holds if $X$ is merely a metrically convex metric space and $Y$ a self-similar metrically convex metric space (see \cite[Section 4]{BaudierGartland24} for more details).
		\end{rema}

		\section{Notes}
		
		In a landmark paper of Bates, Johnson, Lindenstrauss, Preiss and Schechtman \cite{BJLPS1999}, it was shown using the Uniform Approximation by Affine Property (UAAP) that if a Banach space is a uniform quotient of $L_p$ for $p\in(1,\infty)$, then it is already isomorphic to a linear quotient of $L_p$. In the same paper, Bates, Johnson, Lindenstrauss, Preiss and Schechtman raised the question of whether the same can be said about $\ell_p$-spaces. Note that it follows from their work and the linear theory of Banach spaces that $\ell_q$ is not a uniform quotient of $\ell_p$ whenever $1<p<q<\infty$ and $q\neq 2$.
		The relevance of property $(\beta)$ to study uniform quotients of Banach spaces is a brilliant observation due to V. Lima and N. L. Randrianarivony. In \cite{LimaLova2012}, they use property $(\beta)$ to prove that $\ell_q$ is not a uniform quotient of $\ell_p$ whenever $1<p<q<\infty$. In particular, $\ell_2$ is not a uniform quotient of $\ell_p$ for $1<p<2$. This crucial case was the missing bit needed to conclude, using the prior work of Bates, Johnson, Lindenstrauss, Preiss and Schechtman \cite{BJLPS1999} and Johnson and Odell \cite{JohnsonOdell1974} that a uniform quotient of $\ell_p$ with $p\in (1,2)$ must be isomorphic to a linear quotient of $\ell_p$. The following problem remains open.
		
		\begin{prob}
			\label{prob:uniform-quotient-lp}
			Let $p\in(2,\infty)$. If $X$ is a uniform quotient of $\ell_p$, must $X$ be linearly isomorphic to a linear quotient of $\ell_p$?
		\end{prob}
		
		The new idea of Randrianarivony and Lima was expanded upon in \cite{DKLR2014} and \cite{DKR2016}. 
		The study of Poincar\'e-type inequalities tailored to the countably branching trees is from \cite{BaudierGartland24}, where the notions of (infrasup) umbel convexity were introduced. The notion of $(\beta)$-convexity we presented here slightly modifies those ideas but is not explicitly studied in \cite{BaudierGartland24}. 
		
		Equipped with the stability under nonlinear quotients of $(\beta)$-convexity and the fact that countably branching trees are $(\beta_p)$-convex for any $p\in(1,\infty)$, it is possible to prove, via a metric invariant approach, generalized versions of many results about the nonlinear geometry of Banach spaces with property $(\beta)$ (e.g. \cite[Theorem 4.1, Theorem 4.2, Theorem 4.3]{LimaLova2012}, \cite[Corollary 4.3, Corollary 4.5, Corollary 5.2, Corollary 5.3]{DKLR2014}, \cite[Theorem 3.0.2]{DKR2016} and \cite[Theorem 2.1, Theorem 4.6, Theorem 4.7]{BaudierZhang16}). We detail below one example illustrating the flexibility of the metric invariant approach.
		
		Corollary 4.5 in \cite{DKLR2014} states that the space $(\sum_{i=1}^\infty\ell_{p_i})_{\ell_2}$, where $\{p_i\}_{i\ge 1}$ is a decreasing sequence such that $\lim_{i\to\infty}p_i=1$, is not a uniform quotient of a Banach space that admits an equivalent norm with property $(\beta)$. The original proof uses a combination of substantial results from the nonlinear geometry of Banach spaces, which are interesting in their own right:
		\begin{itemize}
			\item Ribe's result that $(\sum_{i=1}^\infty\ell_{p_i})_{\ell_2}$ is uniformly homeomorphic to $\ell_1\oplus(\sum_{i=1}^\infty\ell_{p_i})_{\ell_2}$,
			\item the fact that $\co$ is a linear quotient of $\ell_1\oplus(\sum_{i=1}^\infty\ell_{p_i})_{\ell_2}$,
			\item a quantitative comparison of the $(\beta)$-modulus with the modulus of asymptotic uniform smoothness under uniform quotients (or the qualitative Lima-Randrianarivony theorem \cite{LimaLova2012} which states that $\co$ is not a uniform quotient of a Banach space that admits an equivalent norm with property $(\beta)$).
		\end{itemize}
		Alternatively, using the main result of \cite{DKR2016}, one could argue that the assumption implies that $(\sum_{i=1}^\infty\ell_{p_i})_{\ell_2}$ admits an equivalent norm with property $(\beta)$, hence an equivalent norm that is asymptotically uniformly smooth with power type $p$ for some $p > 1 = \lim_{i\to\infty} p_i$ and derive a contradiction using linear arguments about upper and lower tree estimates, which can be found in \cite{KOS1999} or \cite{OdellSchlumprecht2006RACSAM} for instance.
		
		The metric invariant approach helps streamline and extend the argument as follows. As we have seen, the map $\wTk\ni \nbar\mapsto \sum_{i=1}^{\abs{\nbar}} e_{(n_1,\dots,n_i)} \in \ell_{p_i}(\bN^{\le k})$ is a bi-Lipschitz embedding of $\wTk$ into $\ell_{p_i}$ with distortion at most $2$, say, if $p_i$ is chosen small enough. Therefore, $(\sum_{i=1}^\infty\ell_{p_i})_{\ell_2}$ is not $(\beta_p)$-convex for any $p\in(1,\infty)$ and $(\sum_{i=1}^\infty\ell_{p_i})_{\ell_2}$ is not a uniform quotient of a metrically convex metric space that is $(\beta_p)$-convex for some $p\in(1,\infty)$ by Remark \ref{rem:beta-uniform-quotient}.
		
		For more applications of the Poincar\'e-type inequalities to the geometry of countably branching trees, we refer to \cite{BaudierGartland24}. In particular, if $\mathbb{H}(X)$ is a Heisenberg group over $X$, then it is a result from \cite{BaudierGartland24} that $\cdist{\mathbb{H}(X)}(\wTk)\gtrsim \log(k)^{1/p}$ whenever $X$ has property $(\beta_p)$.
		
		In this chapter, we focused our attention on uniform quotients, but the arguments can be modified to handle the notion of coarse quotient that was introduced by S. Zhang in \cite{Zhang2015}. We refer to \cite[Section 4]{BaudierGartland24} for more details. 
		
		Whether the class of Banach spaces admitting an equivalent norm with property $\beta$ of power type $p$, denoted by $\langle p-\BETA\rangle$, admits a purely metric characterization remains open.
		
		\begin{prob}
			\label{pb:p-beta}
			Given $p\in(1,\infty)$, find a purely metric characterization of the class of Banach spaces admitting an equivalent norm with property $\beta$ of power type $p$.
		\end{prob}
		
		\section{Exercises}

		\begin{exer}
			Exercise to get the truncated tree into any nonreflexive Banach space.   
		\end{exer}
		
		\begin{exer}
			Let $T=(V,E)$ be a connected tree.
			\begin{enumerate}
				\item Show that $T$ embeds isometrically into $\ell_1(E)$. In particular, if $T$ is countable, then $T$ embeds isometrically into the separable Banach space $\ell_1$.
				\item Show that $T$ embeds bi-Lipschitzly into $\co(E)$ with distortion at most $2$. In particular, if $T$ is countable, then $T$ bi-Lipschitzly embeds into $\co$ with distortion at most $2$.
			\end{enumerate}
		\end{exer}
		
		\begin{exer}
			Show that for every $\vep>0$ $\sT_\infty^k$ admits a bi-Lipschitz embedding into $\sT^2_\infty$ with distortion at most $1+\vep$.
		\end{exer}
		
		\begin{exer}
			Consider a truncated binary tree of infinite height, denoted by $\bin_\infty^{trun}$ (better notation needed here) that expands only in the, say, right direction. This is a countably infinite tree that contains, for every $k\in\bN$ an isometric copy of $\bin_k$, the binary tree of height $k$. The lexicographic ordering of the vertices of $\bin_\infty^{trun}$ and consequently of its edges, is a well-order that has the property that for every fixed edge, all the edges that lie on the right and above come after in the enumeration.
			Show that if $X$ is a nonreflexive Banach space and $\vep>0$, then $\bin_\infty^{trun}$ admits a bi-Lipschitz embedding into $X$ with distortion at most $1+\vep$.
		\end{exer}
		
		\begin{exer}
			Show that $\sT_\infty^\omega$ embeds isometrically into $\co$.
		\end{exer}
		
		\begin{exer}[Infrasup Umbel convexity]
			\label{ex:infrasup-umbel}
			Let $p\in(0,\infty)$, $K\in(0,\infty)$ and $(M,d)$ be a metric space such that for all $w,z\in M$ and $\{x_n\}_{n\in \bN}\subseteq M$ 
			\begin{equation}
				\inf_{n\in \bN} \frac{d(w,x_n)^p}{2^p} +\frac{1}{K^p}\inf_{i\neq j\in\bN}d(x_i,x_j)^p\le \max\{d(w,z)^p, \sup_{n\in\bN} d(x_n,z)^p\}.
			\end{equation}
			Show that for all $k\ge 1$ and all $f\colon \tree^\omega_{2^k} \to M$,
			\begin{align*}
				\sum_{l=1}^{k-1}\inf_{s\in {\bN}^{\le 2^{k}-2^l}}\inf_{i \neq j \in \bN} \inf_{u,v \in {\bN}^{2^l-1}}\frac{d(f(s,i,v),f(s,j,u))^p}{2^{lp}}\le K^p\Lip(f)^p.
			\end{align*}
		\end{exer}
		
		\begin{exer}[Umbel convexity]\ 
			\begin{enumerate}
				\item Let $X$ be a Banach space. Show that for all $\delta,\vep>0$, $v,w \in X$ and $V,W \ge \vep$ with $\|v\| \le  V$, $\|w\| \le  W$ and $\frac{1}{2}V + \frac{1}{2}W \le 1$, if $\Big\| \frac{v}{2V} +\frac{w}{2W}\Big\| \le 1-\delta$, then $\Big\|\frac{1}{2}v+\frac{1}{2}w\Big\|\le 1-\vep\delta$.
				\item Show that if $p \in (1,\infty)$ and $X$ has property $(\beta_p)$ with constant $c$, then for all $w,z\in X$ and $\{x_n\}_{n\in \bN}\subseteq X$ 
				\begin{equation}
					\label{eq:pumbel}
					\frac{1}{2^p}\inf_{n\in \bN} \norm{w-x_n}^p+\frac{1}{K}\inf_{i\in\bN}\liminf_{j\in \bN}\norm{x_i-x_j}^p\le \frac{1}{2}\norm{w-z}^p + \frac{1}{2}\sup_{n\in\bN} \norm{x_n-z}^p
				\end{equation}
				where $K$ is the least solution in $[2c,\infty)$ to the inequality
				\begin{equation}
					\label{eq:condition}
					\frac{1}{2^p}\left(\frac{2c}{K}+\left(2-\left(\frac{2c}{K}\right)^p\right)^{1/p}\right)^p + \frac{2^{p+1}}{K} \le  1.
				\end{equation}
				\item 
				Let $p\in(0,\infty)$, $K\in(0,\infty)$ and $(M,d)$ be a metric space such that for all $w,z\in M$ and $\{x_n\}_{n\in \bN}\subseteq M$ 
				\begin{equation}
					\label{eq:pumbelmetric}
					\frac{1}{2^p}\inf_{i \in \bN} d(w,x_i)^p+\frac{1}{K^p}\inf_{i\in\bN}\liminf_{j \to \infty} d(x_i,x_j)^p\le \frac{1}{2}d(z,w)^p + \frac{1}{2}\sup_{i \in \bN} d(z,x_i)^p.
				\end{equation} 
				Show that there exists a constant $\Pi>0$ such that for all $k\ge 1$ and all $f\colon \bN^{\le 2^k}\to X$,
				\begin{align}
					\label{eq:umbel-p-convex}
					\nonumber &\sum_{l=1}^{k-1}\frac{1}{2^{k-1-l}}\sum_{t=1}^{2^{k-1-l}}\inf_{s \in {\bN}^{t2^{l+1}-2^l}}\inf_{u\in{\bN}^{2^l}}\liminf_{j\to\infty}\inf_{v\in {\bN}^{2^l-1}}\frac{d \big( f(s,u),f(s,j,v) \big)^p}{2^{lp}}\\
					&\le \Pi^p\frac{1}{2^{k}}\sum_{\ell=1}^{2^{k}}\sup_{s\in {\bN}^{\ell}}d \big( f(s),f(s^-) \big)^p.
				\end{align}
				
			\end{enumerate}
		\end{exer}
		
		\begin{exer}[Nonembeddability via the self-improvement argument]
			\label{exer:self-improvement-tree}
			Let $\sK_{\omega,1}$ be the star graph with countably many branches, i.e., the bipartite graph that has a partition into exactly two classes, one consisting of a singleton called the center, the other consisting of countably many vertices called the leaves. In this exercise, $b$ will be the center. An arbitrary leaf, denoted by $r$, is chosen and a labeling $(t_i)_{i=1}^\infty$ of the (countably many) remaining leaves is fixed. With this labeling in mind, $\sK_{\omega,1}$ can be seen as an umbel with countably many pedicels, where $r$ stands for root, $b$ for the branching point on the stem and $(t_i)_{i=1}^\infty$ is a labeling of the tips of the pedicels. As usual, $\sK_{\omega,1}$ is equipped with the shortest path metric. Let $X$ be a Banach space with property ($\beta)$.
			\begin{enumerate}
				\item Show that for every bi-Lipschitz embedding $f\colon \sK_{\omega,1}\to X$ there exists $i_0\in\bN$ such that
				\begin{equation}\label{tipbending}
					\norm{f(r)-f(t_{i_0})}_X\le 2\lip(f)\left(1-\bar{\beta}_X\left(\frac{2}{\dist(f)}\right)\right).
				\end{equation}
				\item Let $k\in \bN$ and assume that  $\sT^{\omega}_{2^k}$ bi-Lipschitzly embeds into $X$ with distortion $D$. Then, $\sT^{\omega}_{2^{k-1}}$ bi-Lipschitzly embeds into $X$ with distortion at most $D(1-\bar{\beta}_X(\frac{2}{D}))$. Deduce that if $X$ is a Banach space admitting an equivalent norm with property $(\beta)$, then $\sup_{k\ge 1}\cdist{X}(\wTk)=\infty$.
				\item If $X$ is a Banach space with property $(\beta_p)$ with $p\in(1,\infty)$ and constant $\gamma>0$, show that $\cdist{X}(\wTk)\ge 2\gamma^{\frac{1}{p}}\log(\frac{k}{2})^{\frac{1}{p}}$.
				
			\end{enumerate}
		\end{exer} 
		
		%%%%%%%%%%%%%%%%%%%%%%%%%%%%%%%%%%%%%%%%%%%%%%%%%%%%%%%%%%%%%%%%%%%%%%%%%%%%%%%%%%%%%%%%%%%%%%%

		\chapter{The geometry of diamond graphs and applications}
		\label{chapter:diamonds}
		
		The sequence of (binary) diamond graphs, denoted by $(\diak)_{k\in \bN}$, has fundamental properties. For instance, it was used to show that there is no dimension reduction in $\ell_1$ (see \cite{BrinkmanCharikar2005} and \cite{LeeNaor2004}).
		Of closer relevance to this book, Johnson and Schechtman \cite{JohnsonSchechtman2009} showed that $X\in \langle \UC \rangle$ if and only if $\sup_{k\in \bN}\cdist{X}(\diak)=\infty$. This was extended to the sequence of $r$-branching diamonds (for any $r\in \bN$) by Ostrovskii and Randrianantoanina \cite{OR2017}. Recalling Bourgain's metric characterization, this means that the sequences of binary trees and of binary diamonds both characterize the same class of Banach spaces. In Chapter \ref{chapter:trees} we identified a class of Banach spaces that is characterized by the countably branching trees. It is thus natural to try to identify a class of Banach spaces that is characterized by the countably branching diamonds. In this context, the classes are different. This can be explained by the fact that while in the local setting $\langle \UC \& \US \rangle = \langle \UC \rangle$, in the asymptotic setting, for reflexive Banach spaces, $\langle \BETA\rangle = \langle \AUC \& \AUS \rangle \neq \langle \AUC \rangle$. We have seen in Chapter \ref{chapter:trees} that the countably branching trees characterize the class $\langle \AUC \& \AUS \rangle$ within the class of reflexive spaces. The main result of this chapter states that the countably branching diamonds characterize the class $\langle \AUC \rangle$ within the class of reflexive spaces with an unconditional asymptotic structure. These two results exhibited the first discrepancy between phenomena in the Ribe and Kalton programs. The most difficult part of the metric characterization is the embeddability direction. In Section \ref{sec:diamond-basics}, we define the countably branching diamonds using the graph-theoretic slash-product construction and prove two basic embedding results. In Section \ref{sec:diamond-L1}, we give the exact distortion of the embeddings of countably branching diamond graphs into $L_1$. In Section \ref{sec:diamond-amuc}, we give lower bounds on the distortion incurred when embedding the diamond graphs into Banach spaces that are asymptotically midpoint uniformly convex. We follow the approach of Johnson and Schechtman and use a self-improvement argument. In Section \ref{sec:diamond-co} we study the distortion when embedding them into $\co$ or $\ell_p$. In the last section, we relate this problem to the Szlenk index.   
		
		\section[Basic definitions and graph metric]{Countably branching diamonds: basic definitions and graph metric}
		\label{sec:diamond-basics}
		
		Diamond graphs are constructed inductively according to the following procedure. Start with a diamond of depth $1$ that is $r$-branching with $r\in \{2,3,\dots,\omega\}$. This is a graph that has two special points $s$ and $t$ at distance $2$ and $r$ midpoints $\{m_i\}_{i=1}^w$ that are at distance $2$ from each other and at distance $1$ from $s$ and $t$. To obtain a diamond graph of depth $2$, we replace each edge of the diamond of depth $1$ with a copy of the diamond of depth $1$. We then iterate the procedure to construct a diamond of depth $k$ by replacing each edge of the diamond of depth $k-1$ with a copy of the diamond of depth $1$. To make this graph-theoretical recursive rigorous, we need to introduce the notion of $\oslash$-product (a.k.a slash-product or o-slash-product) from \cite{LeeRaghavendra2010}. 
		
		A directed $s$-$t$ graph $G=(V,E)$ is a directed graph that has two distinguished vertices $s,t\in V$. To avoid confusion, we will also write sometimes $s(G)$ and $t(G)$. Given two directed $s$-$t$ graphs $H$ and $G$, define a new directed $s$-$t$ graph $H\oslash G$ as follows:
		\begin{enumerate}[(i)]
			\item $V(H\oslash G):=V(H)\cup(E(H)\times(V(G)\backslash\{s(G),t(G)\}))$
			\item For every oriented edge $e=(u,v)\in E(H)$, there are $|E(G)|$ oriented edges in $E(H\oslash G)$,
			\begin{align*}
				&\big\{(\{e,v_1\},\{e,v_2\})| (v_1,v_2)\in E(G) \text{ and }v_1,v_2\notin \{s(G),t(G)\}\big\}\\
				\cup &\big\{(u,\{e,w\}) | (s(G),w)\in E(G)\big\}\cup \big\{(\{e,w\},u) | (w, s(G))\in E(G)\big\}\\
				\cup& \big\{(\{e,w\},v) | (w,t(G))\in E(G)\big\}\cup \big\{(v,\{e,w\}) | (t(G),w)\in E(G)\big\}
			\end{align*}
			\item $s(H\oslash G)=s(H)$ and $t(H\oslash G)=t(H)$.
		\end{enumerate}

		It is also clear that the $\oslash$-product is associative (in the sense of graph isomorphism or metric space isometry) and for a directed graph $G$ one can recursively define $G^{\oslash^{k}}$ for all $k\in\bN$ as follows:
		\begin{itemize}
			\item $G^{\oslash^{1}}:= G$.
			\item $G^{\oslash^{k+1}}:=G^{\oslash^{k}}\oslash G$, for $k\ge 1$.
		\end{itemize}
		
		Note that it is sometimes convenient, for some induction purposes, to define $G^{\oslash^{0}}$ to be the two-vertex graph with an edge connecting them. Note also that if the base graph $G$ is symmetric, then the graph $G^{\oslash^{k}}$ does not depend on the orientation of the edges.
		
		If one considers the complete bipartite infinite graph $\sK_{2,\omega}$ with two vertices on one side (such that one is $s(\sK_{2,\omega})$ and the other $t(\sK_{2,\omega})$) and countably many vertices on the other side, then the countably branching diamond graph of depth $k$ is defined as the $k$-fold $\oslash$-product of $\sK_{2,\omega}$, i.e. $\wdiak:=\sK_{2,\omega}^{\oslash^{k}}$. We also define $\mathsf{D_0^\omega}$ to be an edge with two vertices.
		%If one starts with the complete bipartite graph $K_{2,r}$ for some $r\ge 2$ instead, the graph obtained is the $r$-branching diamond graph of depth $k$. In particular $D_k^2:=K_{2,2}^{\oslash^{k}}$.
		
		We will also use a nonrecursive definition of the diamond graphs, which reveals a tree structure in the definition of the diamond graphs. We will say that a vertex in $\wdiak$ is of depth $0\le i \le k$ if it was first created in the $i$-th iteration of the $\oslash$-product So there are only two vertices of depth $0$, namely the terminal vertices of $\wdiak$ and all other vertices of depth at least $1$ will be called \emph{internal vertices}. We will think of one terminal vertex as the top vertex and the other one as the bottom vertex so that we can talk about the height of a vertex (and obviously the height of the bottom vertex will be $0$ while the height of the top vertex will be $2^{k}$. An internal vertex $x$ of depth $i$ is necessarily a midpoint of a copy of $\mathsf{D_1^\omega}$, and we can label it with some integer $n_i\in \bN$ to record its position. The vertex $x$ also comes from an edge of a copy of $\mathsf{D_1^\omega}$ attached to a midpoint of this copy of $\mathsf{D_1^\omega}$. We can also record this information using another integer $n_{i-1}$, but this would not completely characterize the vertex $x$ as we need to account that there are exactly two edges attached. Therefore, we also record a choice of sign $\vep\in\{-1,1\}$ where $-1$ will refer to the lower edge and $1$ to the upper edge. If we keep doing this, then every internal vertex in $\wdiak$ will be labeled and uniquely identified with a point in $(s,\vep)\in \bN^{i}\times \{-1,1\}^{i-1}$ with $1\le i\le k$. Note that we only need $i-1$ choices of signs for a vertex of depth $i$ because $\mathsf{D_0^\omega}$ has only one edge. Vertices of depth $1$ will only be labeled by $(n)$ with $n\in \N$. Finally, in our labeling of $\wdiak$, the top vertex will be indexed by $(\emptyset, 1)$ and the bottom vertex by $(\emptyset,-1)$.  
		%The non-recursive definition is based on the observation that to locate a point in $\wdiak$, one can record the indices of the branches from which the vertex is coming from. But of course, this is not enough as after step $1$ of the construction of $\wdiak$ we also have to record if we are going up or down in the diamond.  Thus we choose to assign $+1$ if we are going up and $-1$ if we are going down. Therefore, every internal vertex in $\wdiak$ will be labeled and uniquely identified, with a point in $(s,\vep)\in \bN^{i}\times \{-1,1\}^{i-1}$ with $1\le i\le k$. 
		A crucial parameter of a vertex $x\in \wdiak$ is its height, denoted by $\abs{x}$ and defined by 
		\begin{equation}
			\abs{x} = \begin{cases}
				0 \hskip 3cm \text{ if } x=(\emptyset,-1),\\
				2^k \hskip 3cm  \text{ if } x=(\emptyset,1),\\
				2^{k-1} + \ds\sum_{j=1}^{i-1} \vep_j 2^{k-(j+1)} \text{ if } x=(s,\vep)\in \bN^{i}\times \{-1,1\}^{i-1} \text{ with } 1\le i\le k.
			\end{cases} 
		\end{equation} 
		
		It follows from the definitions of depth and height that every vertex of depth $1$ has height $2^{k-1}$, a vertex of depth $2$ can have height $2^{k-2}$ or $2^{k-1}+2^{k-2}=3(2^{k-2})$, the height of a vertex of dept $3$ can only be $2^{k-3}$, $2^{k-2}+ 2^{k-3}=3(2^{k-3})$, $2^{k-1}+2^{k-3}=5(2^{k-3})$, or $2^{k-1}+2^{k-2}+2^{k-3}=7(2^{k-3})$, etc... Thus, for depth $1\le i \le k$, there are only $2^{i-1}$ allowable heights for a vertex of depth $i$, which are $(2l-1)2^{k-i}$ for all $1\le l\le 2^{i-1}$.
		It is worth noting that a vertex $x$ is completely determined by its height and its tree label $s$, and we should also write $x=(s,u)$ where $s\in \bN^{\le k}$ and $0\le u\le 2^k$ when it is more convenient to identify a vertex in $\wdiak$.
		For instance, note that the set of edges of $\wdiak$ is given by the set $\{\{(s,u), (t,v)\} \colon \abs{u-v}=1 \text{ and } s\prec t\}$.
		
		\begin{rema}
			In the sequel we are using the symbol $\abs{\cdot}$ in three different ways: 
			\begin{itemize}
				\item as the absolute value of a real number ($\abs{u}$ for $u\in \bR$),
				\item as the length of a point in a tree ($\abs{s}$ for $s \in \bN^{\le k}$),
				\item as the height of a vertex in $\wdiak$ ($\abs{x}$ for $x\in \wdiak$).
			\end{itemize}
			We hope this will not incur any confusion, as it should be clear from the context what is the intended meaning.
		\end{rema}
		Understanding and working with the graph metric on the diamond graphs is arguably more difficult than the tree or Hamming metric. The uniqueness of geodesics in trees and the simple closed formula for the Hamming metric usually simplify the matter. It is not straightforward to obtain a closed formula for the graph metric on the diamond graphs due to the fact that there are many geodesic paths between vertices. Nevertheless, with enough care, the situation is manageable. Some facts that are pretty obvious will be simply stated, and the reader is invited to provide a formal proof (usually a simple induction will do).
		
		\begin{defi}
			We will say that two vertices in $\wdiak$ \emph{belong to the same branch} if they belong to a simple bottom-to-top path in $\wdiak$, where a \emph{simple path} is a path with no loop. We also say that $x$ is an \emph{ancestor} of $y$ and $y$ a \emph{descendant} of $x$, if $x$ and $y$ belong to the same branch and $\abs{x}\le \abs{y}$.
		\end{defi}
		
		Obviously, any vertex and terminal vertex are always on the same branch. Also, two vertices at the same height belong to the same branch 
		if and only if they are identical.
		The distance between vertices that belong to the same branch is easy to compute, as it only involves the height of the vertices.
		If $x,y\in \wdiak$ belong to the same branch, then 
		\begin{equation}
			\dwdiak(x,y)= \big| \abs{x} -\abs{y}\big|.
		\end{equation}
		Note that if $x$ is a descendant of $y$, then $y$ is the common ancestor of $x$ and $y$ with the lowest height and vice versa $x$ is the common descendant of $x$ and $y$ with the largest height. So, the distance between vertices on the same branch could also be interpreted as the minimum of two distances (equal in this case): the distance of a path between $x$ and $y$ passing through the lowest common ancestor and the distance of a path between $x$ and $y$ passing through the highest common descendant. Let us also mention that $x$ is an ancestor of $y$ if and only we can write $x=(s,\eps)$ and $y=(t,\delta)$ in the (tree,signs) label with $s \preceq t$ and $\eps \preceq \delta$. 
		
		The distance between a pair of vertices that are not in the same branch is similar, but proving it is a bit more delicate.
		Indeed, if two vertices $x$ and $y$ are not in the same branch we can start a simple path at $x$ by going up in the diamond until we find a vertex where we can travel down until we reach $y$. This is always possible as we can always go all the way up to the top vertex and travel down from there. Note that every simple path between two vertices on the same branch has the same length, and thus it does not matter which geodesic is chosen.
		Alternatively, we could start our simple path at $x$ by going down in the diamond until we find a vertex where we can travel up until we reach $y$. These two options are most efficient when the vertex where we switch direction between up/down has the smallest height above, or the largest height below, the height of the two vertices that we are trying to connect. Whichever path is the shortest gives the graph distance. In order to formalize this intuition we will need the following lemma.
		
		\begin{lemm}
			\label{lem:diamond-ancestor}
			Let $x=(s,\vep)$ and $y=(t,\delta)$. If $v_0=(r_0,\vep_0), \dots, v_{n}=(r_n,\vep_n)$ are the vertices in a path between $x$ and $y$, then there is $0\le i\le n$ such that $r_i\preceq s \land t$, where $s\land t$ is the greatest  common ancestor of $s$ and $t$ in the tree-ordering. 
		\end{lemm}
		
		\begin{proof}
			Assume that there is a path $v_0=(r_0,\vep_0), \dots, v_{n}=(r_n,\vep_n)$ such that $v_0=x$ and $v_n=y$ and moreover that for all $1\le i\le n$, $r_i\not\preceq s\land t$. Since $r_n=t$, we would obtain a contradiction if we can show that for all $0\le i\le n$, $r_i\not\preceq t$ and $r_i\land t = s\land t$. For $i=0$, $s=r_0\not\preceq r_0\land t$ and hence $r_0\not\preceq t$, while clearly $r_0\land t =s\land t$. Assume now that $r_{i} \not\preceq t$ and $r_i\land t =s\land t$. Observe first that since we also have by assumption that $r_{i}\not\preceq s\land t$, it follows that $(r_i\land t, r_i]\cap (r_i\land t, t] =\emptyset$ with $(r_i\land t, r_i]\neq\emptyset$. We consider two cases. Either $r_{i}\preceq r_{i+1}$, in which case $r_{i+1}\not\preceq t$ since $r_{i+1}$ extends $r_i$ but $r_i\not\preceq t$ and also $r_{i+1}\land t=s\land t$ since $r_{i}\land t =s\land t$ and $r_i\not\preceq s\land t$. Or, $r_{i+1}\preceq r_i$. In this case, since  $r_{i+1}\in (s\land t, r_i]$, then we also have that $r_{i+1}\not\preceq t$ and $r_{i+1}\land t=s\land t$.
		\end{proof}
		
		Given $1\le i\le j\le k$ and two internal vertices $x=((s_1,\dots,s_i),(\vep_i,\dots,\vep_{i-1}))$ and $y=((t_1,\dots,t_j),(\delta_1,\dots,\delta_{j-1}))$, observe that $x$ and $y$ are not on the same branch if and only if exactly one of the following two situations occurs:
		\begin{enumerate}
			\item $s_1\neq t_1$,
			\item there is $2\le l \le i$ such that $\vep_{\restriction l-1}=\delta_{\restriction l-1}$, $s_{\restriction l-1}=t_{\restriction l-1}$, but $s_{l}\neq t_{l}$.
		\end{enumerate}
		If $1.$ occurs, then $s\land t =(\emptyset) $ where $s:=(s_1,\dots,s_i)$ and $t:=(t_1,\dots,t_j)$ and for any vertex $(r,w)$, written with the (tree,height)-label,  such that $r\preceq s\land t$ we have $w\in \{0,2^k\}$ and thus $w\notin [\min\{\abs{x},\abs{y}\}, \max\{\abs{x},\abs{y}\}]$.\\
		%then   path between $x$ and $y$, any vertex $z$ whose existence is provided by Lemma \ref{lem:diamond-ancestor} is either the top of bottom vertex and thus $\abs{z}\notin [\min\{\abs{x},\abs{y}\}, \max\{\abs{x},\abs{y}\}]$.
		If $2.$ occurs, then $s\land t=s_{\restriction l-1}$. Note that $z=(s\land t, \vep_{\restriction l-1})$ has depth $l-1$ and that the height of a vertex of depth $l-1$ is a multiple of $2^{k-l+1}$. On the other hand, a simple computation shows that there exists $n\in \N$ such that both $|x|$ and $|y|$ belong to $(n2^{k-l+1}, (n+1)2^{k-l+1})$. It follows that for any vertex $(r,w)$ such that $r\preceq s\land t$ we have $w\notin [\min\{\abs{x},\abs{y}\}, \max\{\abs{x},\abs{y}\}]$.
		%for every path between $x$ and $y$ the height of any vertex $z$ whose existence is provided by Lemma \ref{lem:diamond-ancestor} satisfies $\abs{z}\notin [\min\{\abs{x},\abs{y}\}, \max\{\abs{x},\abs{y}\}]$. 
		%Therefore, two vertices belong to the same branch if and only if there is a path between $x=(s,u)$ and $y=(t,v)$ (with $u\le v$) and a vertex $z=(r,w)$ such that $r\preceq s\land t$ and $w\in[u,v]$. 
		
		\medskip 
		Given $x=(s,u)$ and $y=(t,v)$ in the (tree, height) label, define 
		\begin{equation}
			h_{a}(x,y):= \max\{ 0\le w \le 2^k \colon \exists z=(r,w) \text{ s.t. } r\preceq s\land t \text{ and } w\le \min\{u,v\}\}, 
		\end{equation}
		and
		\begin{equation}
			l_{a}(x,y):= \min\{ 0\le w \le 2^k \colon \exists z=(r,w) \text{ s.t. } r\preceq s\land t \text{ and } w\ge \max\{u,v\}\}, 
		\end{equation}
		When $x$ and $y$ are not on the same branch, $h_a(x,y)$ returns the height of the highest common ancestor below $x$ and $y$ and $l_a(x,y)$ returns the height of the lowest common ancestor above $x$ and $y$. Using $h_a$ and $l_a$ we can give an explicit formula for the distance between a pair of vertices that are not on the same branch.
		
		\begin{lemm}\label{lem:formula_diamondmetric}
			Let $x=(s,u)$ and $y=(t,v)$ be two vertices in $\wdiak$, written in the (tree, height) label,  that are not on the same branch. Then,
			\begin{align*}
				\dwdiak(x,y) & = \min\{ u+v -2 h_a(x,y), 2l_a(x,y) -(u+v) \}.
			\end{align*}
		\end{lemm}
		
		\begin{proof}
			Assume without loss of generality that $u\le v$. It is easy to see that if a vertex $z=(r,w)$ is such that $r\preceq s\land t$ and $l_a(x,y)=w$, then since $z$ and $x$ are on the same branch, there is a path of length $l_a(x,y)- u$ between $x$ and $z$ and a path of length $l_a(x,y)- v$ between $y$ and $z$. By concatenating these two paths we get a path between $x$ and $y$ of length $2l_a(x,y) -(u+v)$. A similar argument shows that there is a path of length $u+ v - 2h_a(x,y)$ between $x$ and $y$. Consequently, there is always a path of length $\min\{2l_a(x,y) -(u+v), u+ v - 2h_a(x,y)\}$ between $x$ and $y$. 
			
			It remains to show that this bound is optimal. Consider a simple path $P_{xy}$ between $x$ and $y$. By Lemma \ref{lem:diamond-ancestor} and the discussion above, there is a vertex $z=(r,w)$ in this path such that $r\preceq s\land t$ and either $w\ge v\ge u$ or $w\le u\le v$. Since $r\preceq s$ and $r\preceq t$, $z$ and $x$ are on the same branch and $z$ and $y$ are on the same (but different) branch. If $w\ge v\ge u$, the length of $P_{xy}$ which is the sum of the length of $P_{xz}$ and the length of $P_{zy}$ is at least $w-u +w- v =2w -(u+v)$ (since by definition of the edge set we need at least $w-v$ edges to go from a height $v$ to $w$ and $w-u$ edges to go from a height $u$ to $w$). Thus, by definition of $l_a(x,y)$, the path between $x$ and $y$ has length at least $2l_a(x,y) -(u+v)$.
			Similarily, if $w\le u\le v$, the path between $x$ and $y$ has length at least $u+ v - 2h_a(x,y)$. This means that every path between $x$ and $y$ has length at least $\min\{2l_a(x,y) -(u+v), u+ v - 2h_a(x,y)\}$ and the conclusion follows.
		\end{proof}

		\section{\texorpdfstring{Bi-Lipschitz embeddability of diamond graphs into $L_1$}{Bi-Lipschitz embeddability of diamond graphs into}}
		\label{sec:diamond-L1}

		In a celebrated unpublished article (cf. \cite{Benyamini1985}, \cite{BenyaminiLindenstrauss2000}, or \cite{WestonThesis1993}), Enflo used an approximate midpoint argument to show that $L_1$ and $\ell_1$ are not uniformly homeomorphic. This simple, clever and extremely useful argument, is based on the fact that in $\ell_1$ the size of the approximate midpoint set between two points is rather small. In contrast, it is easy to find pairs of points in $L_1$ that have infinitely many exact midpoints that are far apart. From the purely metric point of view, Enflo's argument relies on the fact that $L_1$ contains a bi-Lipschitz copy of the countably branching diamond of depth $1$ that preserves exactly the pairwise distances of any triple $(v_s,v_i,v_t)$ where $v_s$ and $v_t$ are the two terminal vertices and $v_i$ a midpoint (see Exercise \ref{exe:diamond1-bernoulli}). Whether or not a similar embedding, for the countably branching diamond graphs of \emph{arbitrary depth} can be implemented without blowing up the distortion is not completely obvious at first sight. However, using the results from \cite{GNRS04} and an ultraproduct argument it can be shown that $\wdiak$ embeds into $L_1$ with distortion at most $2$. Below we present an embedding using Bernoulli random variables (and no ultraproduct argument) that achieves the same distortion. We also point out that the distortion is necessarily incurred by a pair of vertices that are not in the same branch.
		
		\begin{theo}
			\label{thm:diamond-L1-sharp} 
			Let $(\Omega, \Sigma, \mu)$ be an atomless probability space.
			For every $k\in\bN$ there exists $f_k \colon \wdiak\to L_1(\Omega)$ such that
			if $x$ and $y$ belong to the same branch then
			\begin{equation}
				\label{eq1:diamond-L1-sharp}
				\norm{f_k(x) - f_k(y)}_{1} = \dwdiak(x,y),
			\end{equation}
			and if $x$ and $y$ do not belong to the same branch then
			\begin{equation}
				\label{eq2:diamond-L1-sharp}
				\frac12 \dwdiak(x,y) \le \norm{f_k(x) - f_k(y)}_{1}\le  \dwdiak(x,y).
			\end{equation}
			In particular, $\sup_{k\in \bN} \cdist{L_1}(\wdiak)\le 2$.
		\end{theo}
		
		\begin{proof}
			We will prove the statement by induction on $k$. The idea is to construct measurable subsets with certain monotonicity and independence properties and whose measures capture the height of vertices. The embedding is then defined in terms of the indicator functions of these sets. Since $(\Omega, \Sigma, \mu)$
			is atomless there are countably many independent identically distributed $\{0,1\}$-valued Bernoulli random variable on $(\Omega, \Sigma, \mu)$. This leads us to prove the following statement. For every countable family $(\vep_t)_{t\in \bN^{\le k}\setminus \{\emptyset\}}$ of i.i.d. $\{0,1\}$-valued Bernoulli random variable on $(\Omega, \Sigma, \mu)$ and every $x\in \wdiak$ there exists a set $S_k(x) \in \Sigma$ such that:
			\begin{enumerate}
				\item $S_k(x)\in \sigma\{ (\vep_{(t_1,s)}):\ s\in \bN^{\le k-1}\}\}$ wherever $x=((t_1,\dots,t_j),u)$ with $j\ge 1$,
				\item $S_k(x)\subset S_k(y)$ for all $x,y$ in the same branch with $\abs{x} \le \abs{y}$,
				\item $\mu(S_k(x)) = \frac{1}{2^k}\abs{x}$ for all $x\in \wdiak$, 
				%    \item $\mu(S_k(x)\cap S_k(y))=\mu(S_k(x))\mu(S_k(y))$, whenever $x=((t_1,\dots,t_i),u)$ and $y=((s_1,\dots,s_j),v)$ with $i,j\ge1$ and $t_1\neq s_1$.
				\item The map $f_k\colon \wdiak \to L_1(\Omega)$ given by $f_k(x)=2^k\car_{S_{k}(x)}$ satisfies \eqref{eq1:diamond-L1-sharp} and \eqref{eq2:diamond-L1-sharp}.
			\end{enumerate}
			
			So let us prove statements 1 to 4 by induction on $k\in \N$. For $k=1$, let $(\vep_i)_{i\in\bN}$ be i.i.d. $\{0,1\}$-valued Bernoulli random variable on $(\Omega, \Sigma, \mu)$. We define $S_1((\emptyset,-1))=\emptyset$, $S_1((\emptyset,1))=\Omega$ and $S_1((t))= \{\vep_t =1\}$. It is completely straightforward to verify that properties $1-4$ are satisfied (note in fact that $f_1$ is an isometry). Assume now that for every countable family $(\vep_t)_{t\in \bN^{\le k}\setminus \{\emptyset\}}$ of i.i.d $\{0,1\}$-valued Bernoulli random variable on $(\Omega, \Sigma, \mu)$ and every $x\in \wdiak$ there exists a set $S_k(x) \in \Sigma$ so that properties $1-4$ hold. Fix  a countable family $(\vep_t)_{t\in \bN^{\le k+1}\setminus \{\emptyset\}}$ of i.i.d $\{0,1\}$-valued Bernoulli random variable on $(\Omega, \Sigma, \mu)$. Recall that $\mathsf{D_{k+1}^\omega}= \mathsf{D_{1}^\omega} \oslash \mathsf{D_{k}^\omega}$ and for $t_1\in \N$ and $\eta\in\{-1,1\}$ we will denote by $e_{t_1}^{\eta}\oslash \mathsf{D_{k}^\omega}$ the copy of $\wdiak$ that replaced the edge $\{(t_1), (\emptyset,1)\}$ (of $\mathsf{D_1^\omega}$) if $\eta =1$, or that replaced the edge $\{(t_1), (\emptyset,-1)\}$ if $\eta =-1$. By the induction hypothesis for every $t_1\in \bN$, $\eta \in \{-1,1\}$,  and $x \in e_{t_1}^{\eta}\oslash \mathsf{D_{k}^\omega}$ there are sets $S_{k}^{t_1,\eta}(x)\in \Sigma$ such that 
			\begin{enumerate}[(1')]
				\item $S_k^{t_1,\eta}(x)\in \sigma\{ \vep_{(t_1,s)} \colon s\in \bN^{\le k}\setminus \{\emptyset\} \}$ whenever $x$ has nonzero depth in $e_{t_1}^{\eta}\oslash \mathsf{D_{k}^\omega}$,
				\item $S_k^{t_1,\eta}(x)\subset S_k^{t_1}(y)$ for all $x,y$ in the same branch of $e_{t_1}^{\eta}\oslash \mathsf{D_{k}^\omega}$ with $\abs{x}_{k} \le \abs{y}_{k}$,
				\item $\mu(S_k^{t_1,\eta}(x)) = \frac{\abs{x}_{k}}{2^k}$ for all $x\in e_{t_1}^{\eta}\oslash \mathsf{D_{k}^\omega}$,
				%     \item $\mu(S_k^{t_1,\eta}(x_\eta)\cap S_k^{t_1,\eta}(y_\eta))=\mu(S_k^{t_1,\eta}(x_\eta))\mu(S_k^{t_1,\eta}(y_\eta))$, whenever $x_\eta=((t_2,\dots,t_i),u)$ and $y_\eta=((s_2,\dots,s_j),v)$ with $i,j\ge 2$ and $t_2\neq s_2$,
				\item The map $f_k^{t_1,\eta}\colon e_{t_1}^{\eta}\oslash \mathsf{D_{k}^\omega} \to L_1$ given by $f_k^{t_1,\eta}(x)=2^k\car_{S_k^{t_1,\eta}(x)}$ satisfies \eqref{eq1:diamond-L1-sharp} and \eqref{eq2:diamond-L1-sharp}.
			\end{enumerate}
			Here and in the sequel, $\abs{\cdot}_{k}$ is the height in $e_{t_1}^{\eta}\oslash \mathsf{D_{k}^\omega}$ and observe that 
			\begin{equation}
				\abs{x}_k=\begin{cases}
					\abs{x}_{k+1} \text{ if } \eta =-1\\
					\abs{x}_{k+1} -2^{k} \text{ if } \eta =1.
				\end{cases}
			\end{equation}
			
			For any $x\in \mathsf{D_{k+1}^\omega}$ we define the measurable sets as follows:
			$S_{k+1}((\emptyset,-1))=\emptyset$, $S_{k+1}((\emptyset,1))=\Omega$, $S_{k+1}((t_1))=\{\vep_{t_1}=1\}$ and if $x=((t_1,t),(\delta_1,\delta))$ with $t\in \N^i$, $1\le i \le k$, $\delta_1 \in \{-1,1\}$ and $\delta \in \{-1,1\}^{i-1}$:
			$$S_{k+1}(x):= \begin{cases}
				\{\vep_{t_1}=1\}\cap S_{k}^{t_1,-1}((t,\delta)) \text{ whenever } \delta_1=-1,\\
				\{\vep_{t_1}=1\}\cup( \{\vep_{t_1}=0\}\cap S_{k}^{t_1,1}((t,\delta))) \text{ whenever } \delta_1=1.
			\end{cases} 
			$$
			
			The independence property is clear since if $x=((t_1,\dots,t_i),u)$, then by definition of the sets and $(1')$ we have $S_{k+1}(x)\in  \sigma\{ \vep_{(t_1,s)} \colon s\in \bN^{\le k}\} $. 
			%while $S_{k+1}(y)\in  \sigma\{ \vep_{(s_1,s)} \colon s\in \bN^{\le k} \} $ but these two $\sigma$-algebras are independent since $t_1\neq s_1$.
			%It is clear from $(1')$ and the definition of the sets that if $x$ has non-zero depth, then $S_{k+1}(x)$ is in $\sigma\{ \vep_{s} \colon s\in \bN^{\le k+1}\setminus \{\emptyset\} \}$.
			The monotonicity property is also easy to see. If $x$ and $y$ are in the same branch and $\abs{x}_{k+1}\le \abs{y}_{k+1}$, then either $x$ and $y$ are in $e_{t_1}^{\eta}\oslash \mathsf{D_{k}^\omega}$ for some $t_1\in \bN$ and $\eta\in \{-1,1\}$ or $x\in e_{t_1}^{-1}\oslash \mathsf{D_{k}^\omega}$ and $y \in e_{t_1}^{1}\oslash \mathsf{D_{k}^\omega}$. The former case follows from $(2')$, the definition of the sets and the fact that $\abs{x}_{k}\le \abs{y}_{k}$. In the latter case, we have $S_{k+1}(x)\subset \{\vep_{t_1}=1\}\subset S_{k+1}(y)$. 
			
			For property $3$, $\mu(S_{k+1}((\emptyset,-1)))= 0$, $\mu(S_{k+1}((\emptyset,1)))= 1$ and $\mu(S_{k+1}((t_1))= \frac{1}{2} = \frac{2^{k-1}}{2^k}$. If $x\in e_{t_1}^{-1}\oslash \mathsf{D_{k}^\omega}$, then \begin{align*}
				\mu(S_k(x)) & =  \mu(\{\vep_{t_1}=1\}\cap S_{k}^{t_1,-1}((t,\delta)))\\
				& = \mu(\{\vep_{t_1}=1\})\mu(S_{k}^{t_1,-1}((t,\delta))) \text{ by independence}\\
				& \stackrel{(3')}{=} \frac12\frac{\abs{(t,\delta)}_k}{2^k}\\
				& = \frac{\abs{x}_{k+1}}{2^{k+1}}.
			\end{align*}
			If $x\in e_{t_1}^{1}\oslash \mathsf{D_{k}^\omega}$, then 
			\begin{align*}
				\mu(S_k(x)) & =  \mu(\{\vep_{t_1}=1\}\cup( \{\vep_{t_1}=0\}\cap S_{k}^{t_1,1}((t,\delta))))\\
				& = \mu(\{\vep_{t_1}=1\}) + \mu(\{\vep_{t_1}=0\})\mu(S_{k}^{t_1,1}((t,\delta))) \text{ by disjointness and independence}\\
				& \stackrel{(3')}{=} \frac12 + \frac12\frac{\abs{(t,\delta)}_k}{2^k}\\
				& = \frac12 + \frac{\abs{x}_{k+1}-2^{k}}{2^{k+1}}\\
				& = \frac{\abs{x}_{k+1}}{2^{k+1}}.
			\end{align*}
			
			It remains to verify the embeddability property. Given $x,y\in\mathsf{D_{k+1}^\omega}$ we consider a few cases. If $x$ and $y$ belong to the same branch with $\abs{x}\le \abs{y}$, then
			\begin{align*}
				\norm{f_{k+1}(x) -f_{k+1}(y)}_1 & \stackrel{(2)}{=} 2^{k+1}(\mu(S_{k+1}(y)) - \mu(S_{k+1}(x))\\
				& \stackrel{(3)}{=} 2^{k+1}\Big(\frac{\abs{y}}{2^{k+1}} - \frac{\abs{x}}{2^{k+1}}\Big)\\
				& = \mathsf{d_{D_{k+1}^\omega}}(x,y),
			\end{align*}
			and hence \eqref{eq1:diamond-L1-sharp} holds. Showing \eqref{eq2:diamond-L1-sharp} is a little bit more delicate.
			
			If $x,y\in e_{t_1}^{-1}\oslash \mathsf{D_{k}^\omega}$ for some $t_1\in \bN$, then $x=((t_1,t),(-1,\delta))$ and $y=((t_1,t'),(-1,\delta'))$ and
			\begin{align*}
				\norm{f_{k+1}(x) -f_{k+1}(y)}_1 & = 2^{k+1}\norm{ \car_{S_{k+1}(x)} -\car_{S_{k+1}(y)} }\\
				& = 2^{k+1}\norm{ \car_{\{\vep_{t_1}=1\}}( \car_{S_{k}^{t_1,-1}((t,\delta))} -\car_{S_{k}^{t_1,-1}((t',\delta'))}) }\\
				& = 2^{k} \norm{\car_{S_{k}^{t_1,-1}((t,\delta))} -\car_{S_{k}^{t_1,-1}((t',\delta'))}} \text{ (by independence)}\\
				& = \norm{f_k^{t_1,-1}((t,\delta)) - f_k^{t_1,-1}((t',\delta'))},
			\end{align*}
			and since $\dwdiak((t,\delta),(t',\delta'))= \mathsf{d_{D_{k+1}^\omega}}(x,y)$ inequalities \eqref{eq2:diamond-L1-sharp} for $f_k$ follow from inequalities \eqref{eq2:diamond-L1-sharp} for $f_k^{t_1,-1}$.
			
			If $x,y\in e_{t_1}^{1}\oslash \mathsf{D_{k}^\omega}$ for some $t_1\in \bN$, then $x=((t_1,t),(1,\delta))$ and $y=((t_1,t'),(1,\delta'))$ and
			\begin{align*}
				\norm{f_{k+1}(x) -f_{k+1}(y)}_1 & = 2^{k+1}\norm{ \car_{S_{k+1}(x)} -\car_{S_{k+1}(y)} }\\
				& = 2^{k+1}\norm{\car_{\{\vep_{t_1}=1\}} + \car_{\{\vep_{t_1}=0\}} \car_{S_{k}^{t_1,1}((t,\delta))} - (\car_{\{\vep_{t_1}=1\}} + \car_{\{\vep_{t_1}=0\}}\car_{S_{k}^{t_1,1}((t',\delta'))}) }\\
				& = 2^{k+1}\norm{ \car_{\{\vep_{t_1}=0\}}( \car_{S_{k}^{t_1,1}((t,\delta))} - \car_{S_{k}^{t_1,1}((t',\delta'))}) }\\
				& = 2^{k} \norm{\car_{S_{k}^{t_1,1}((t,\delta))} -\car_{S_{k}^{t_1,1}((t',\delta'))}} \text{ (by independence)}\\
				& = \norm{f_k^{t_1,1}((t,\delta)) - f_k^{t_1,1}((t',\delta'))},
			\end{align*}
			and since $\dwdiak((t,\delta),(t',\delta')= \mathsf{d_{D_{k+1}^\omega}}(x,y)$ inequalities \eqref{eq2:diamond-L1-sharp} for $f_k$ follow from inequalities \eqref{eq2:diamond-L1-sharp} for $f_k^{t_1,1}$.
			
			If $x\in e_{t_1}^{\eta}\oslash \mathsf{D_{k}^\omega}$ and $y \in e_{t_2}^{\eta'}\oslash \mathsf{D_{k}^\omega} $ for some $t_1\neq t_2\in \bN$ and $\eta,\eta'\in \{-1,1\}$, then 
			\begin{align*}
				\norm{f_{k+1}(x) -f_{k+1}(y)}_1 & = 2^{k+1}\mu(S_{k+1}(x)\Delta S_{k+1}(y)) \\
				& = 2^{k+1}(\mu(S_{k+1}(x)) +\mu (S_{k+1}(y)) -2\mu( S_{k+1}(x)\cap S_{k+1}(y)))\\
				& \stackrel{(1)}{=} 2^{k+1}(\mu(S_{k+1}(x)) +\mu (S_{k+1}(y)) -2\mu( S_{k+1}(x))\mu(S_{k+1}(y))).
			\end{align*}
			
			Observing that in this case 
			\begin{align*}
				\mathsf{d_{D_{k+1}^\omega}}(x,y) & = \min\{ \abs{x} +\abs{y}, 2^{k+2} -\abs{x}-\abs{y}\}\\ 
				& \stackrel{(3)}{=} 2^{k+1} \min\{ \mu(S_{k+1}(x)) + \mu(S_{k+1}(y)), 2- \mu(S_{k+1}(x)) - \mu(S_{k+1}(y)) \},
			\end{align*}
			inequalities \eqref{eq2:diamond-L1-sharp} follow immediately from the following elementary inequalities holding for all $\alpha,\beta\in [0,1]$, 
			\begin{equation}
				\label{eq2:diamond-L1}
				\alpha + \beta -2\alpha\beta \le \min\{\alpha +\beta, 2-\alpha -\beta\} \le 2(\alpha + \beta -2\alpha\beta)
			\end{equation}
			To see why \eqref{eq2:diamond-L1} holds, observe first that 
			\begin{equation*}
				\alpha +\beta -2\alpha\beta = (1-\alpha) + (1-\beta) -2(1-\alpha)(1-\beta) \le 2-\alpha -\beta,
			\end{equation*}
			since $\alpha,\beta\le 1$. The other inequality follows by applying the first inequality to $1-\alpha$ and $1-\beta$.
			
			%For the other inequality, we first show that if $\alpha + \beta \le 1$, then $ \alpha + \beta \le 2(\alpha +\beta -2\alpha\beta)$, or equivalently $\alpha (4\beta -1)\le \beta$. If $4\beta -1 \le 0$, there is nothing to do so assume that $4\beta -1>0$, then $\alpha (4\beta-1)\le (1-\beta)(4\beta -1) = 4\beta - 1 -4\beta^2+\beta = \beta - (2\beta -1)^2\le \beta$. Assume now that $\alpha +\beta >1$, then $(1-\alpha) + (1-\beta) <1$ and it follows from what we just proved that $(1-\alpha) +(1-\beta) \le 2( (1-\alpha) +(1-\beta) -2(1-\alpha)(1-\beta))$, which is the same after rearranging as $2- \alpha -\beta \le 2(\alpha +\beta -2\alpha\beta)$ and the proof is complete.
		\end{proof}
		
		\begin{rema}
			It follows from \cite{LeeRaghavendra2010} that this distortion is optimal.     
		\end{rema}

		\section[Bi-Lipschitz embeddability of diamond graphs into AMUC spaces]{Bi-Lipschitz embeddability of diamond graphs into asymptotically midpoint uniformly convex Banach spaces}
		\label{sec:diamond-amuc}

		As already mentioned in the previous section, $\ell_1$ does not contain bi-Lipschitz copies of $\mathsf{D_1^\omega}$ that are vertically isometric because its approximate midpoint sets are too small. The nonembeddability result presented in this section builds on this observation and provides sharp lower bounds for the distortion incurred by embeddings of the countably branching diamond graphs into asymptotically midpoint uniformly convex spaces. As in \cite{Baudieretal}, we follow the self-improvement approach from Johnson and Schechtman. This argument starts with the observation that if $\mathsf{D_1^\omega}$ admits a bi-Lipschitz embedding into a space that is asymptotically uniformly convex, then the images of the two terminal vertices must be slightly more contracted than expected. In fact, we only need the target space to be asymptotically midpoint uniformly convex. Recall that a Banach space $(X,\norm{\cdot})$ is asymptotically midpoint uniformly convex if $\hat{\delta}_{\norm{\cdot}}(t)>0$ for every $t\in(0,1)$, where 
		\begin{equation*}
			\hat{\delta}_{\norm{\cdot}}(t):=\inf_{x\in S_X}\sup_{Z\in\cof(X)}\inf_{z\in S_Z}\max\{\norm{x+tz}, \norm{x-tz}\}-1
		\end{equation*}
		
		Similarly to Proposition \ref{prop:moduli-nets}, one can easily prove that 
		\begin{equation}
			\label{eq:AMUCnets}
			\hat{\delta}_{\norm{\cdot}}(t):=\inf_{x\in S_X} \inf \{ \liminf_\alpha \max\{\norm{x+x_\alpha}, \norm{x-x_\alpha}\}-1 \colon x_\alpha\wtoo_\alpha 0, \norm{x_\alpha}\ge t\}.   
		\end{equation}
		
		Then, the first key lemma is the following. It was first proved in \cite{Baudieretal}. We will give a different proof with better constants, using ideas from \cite{Basset} and \cite{BLPinprepa}. The original proof from \cite{Baudieretal}, which used the Kuratowski measure of noncompactness of approximate midpoint sets, is discussed in Exercise \ref{exe:AMUC-original}.

			\begin{lemm}
				\label{lem:terminals-contraction-bis}
				Let $\mathsf{D_1^\omega}$ be the countably branching diamond graph of depth $1$ whose terminal vertices are $s$ and $t$. Let $(X,\norm{\cdot})$ be a Banach space and $f\colon \mathsf{D_1^\omega}\to X$ be a bi-Lipschitz embedding with distortion $C \ge 1$. Then,
				\begin{equation}
					\label{eq:terminals-contraction}
					\norm{f(t)-f(s)} \le \Lip(f)\Big(1-\frac12\hat{\delta}_{\norm{\cdot}}\big(\frac{1}{C}\big)\Big)\sd_{\sD_1^\omega}(t,s).
				\end{equation}
			\end{lemm}
			
			\begin{proof}
				Without loss of generality, we may assume that $\Lip(f)=1$, so that for all $u,v \in \sD_1^\omega$,
				\begin{equation*}
					\frac1C \sd_{\sD_1^\omega}(u,v)\le \norm{f(u)-f(v)}\le \sd_{\sD_1^\omega}(u,v).
				\end{equation*}
				We recall that $s$ and $t$ are the two distinguished points of $\sD_1^\omega$ at distance $2$ and $\{m_i\colon  i \in \bN\}$ is the set of exact metric midpoints of $s$ and $t$. We set  
				\begin{align*}
					f_{u,v} & :=\frac{f(u)-f(v)}{\sd_{\sD_1^\omega}(u,v)} \text{ for } u\neq v \in\sD_1^w,\\
					x & := f_{t,s},\\
					x_i & :=f_{t,m_i},\  y_i:=f_{m_i,s}, \text{ for } i \in \bN.
				\end{align*}
				Since $\Lip(f)=1$, all the above vectors belong to $B_X$. Note that for all $i \in \N$, $x=\frac12(x_i+y_i)$. Let now $V$ be a weak neighborhood of $x$ of the form $V :=\bigcap_{k=1}^n \{y\in X,\ x_k^*(y)<\alpha_k\}$ with $x_k^* \in X^*$ and $\alpha_k \in \R$. By passing to  subsequences of $(x_i)_i$ and $(y_i)_i$, we may assume that for all $1\le k \le n$, $x^*_k(x_i)\to c_k\in \R$ and $y^*_k(y_i) \to d_k \in \R$ as $i$ tends to $\infty$. Note that $x^*_k(x)=\frac12(c_k+d_k)<\alpha_k$. Thus, for $i\neq j$ large enough, $x^1_V=\frac12(x_i+y_j)$ and  $x^2_V=\frac12(x_j+y_i)$ belong to $V$. Therefore, $(x^1_V)_V$ and $(x^2_V)_V$ are nets in $B_X$ (indexed by the weak neighborhoods of $x$ of the above form) weakly converging to $x$. Now we can write $x^1_V=x+u_V$ and $x^2_V=x-u_V$ with $u_V=\frac12(y_j-y_i)=\frac12(f(m_j)-f(m_i))$, so $\|u_V\|\ge \frac1C$. Recalling that $(u_V)_V$ is a weakly null net, it follows from \eqref{eq:AMUCnets} and an homogeneity argument that 
				$\norm{x}\le \big(1+\hat{\delta}_{\norm{\cdot}}(\frac1C)\big)^{-1}\le 1-\frac12 \hat{\delta}_{\norm{\cdot}}(\frac1C)$. The last inequality is due to the fact that $\hat{\delta}_{\norm{\cdot}}(\frac1C)\le \frac1C\le 1$ and concludes the proof.
			\end{proof}
			
			The next proposition is the self-improvement argument \`a la Johnson and Schechtman adapted to our setting.
			
			\begin{prop}
				\label{pro:diamond-self-improv}
				Let $(X,\norm{\cdot})$ be a Banach space. If $k\in\bN$ and $\wdiak$ embeds bi-Lipschitzly into $X$ with distortion $C$, then $\mathsf{D_{k-1}^\omega}$ embeds bi-Lipschitzly into $X$ with distortion at most $C\big(1-\frac{1}{2}\hat{\delta}_{\norm{\cdot}}(\frac{1}{C})\big)$.
			\end{prop}
			
			\begin{proof}
				Assume $(X,\norm{\cdot})$ is a Banach space with an asymptotically midpoint uniformly convex norm. The argument is similar to the proof of Proposition $2.1$ in \cite{BaudierZhang16} and the formal (more rigorous) argument using set-theoretic representations of the graphs shall be simply sketched here. Let $f_k$ be a bi-Lipschitz embedding of $\wdiak$ into $X$ that is noncontracting and $C$-Lipschitz. Note that the subset of vertices $V(\mathsf{D_{k-1}^\omega})\subset V(\wdiak):=V(\mathsf{D_{k-1}^\omega}\oslash \mathsf{D_{1}^\omega})$ forms an isometric copy (up to a scaling factor $2$) of $\mathsf{D_{k-1}^\omega}$. Define $g_k$ to be the restriction of the embedding $f_k$ to the subset $V(\mathsf{D_{k-1}^\omega})\subset V(\wdiak):=V(\mathsf{D_{k-1}^\omega}\oslash \mathsf{D_{1}^\omega})$ rescaled by the factor $2$. Since in our setting it suffices to check the distortion on pair of adjacent vertices, it follows from Lemma \ref{lem:terminals-contraction-bis} that $g_k$ is a bi-Lipschitz embedding with distortion at most $C\big(1-\frac{1}{2}\hat{\delta}_{\norm{\cdot}}(\frac{1}{C})\big)$.
			\end{proof}
			
			Our last theorem is an asymptotic version of a result of Johnson and Schechtman in \cite{JohnsonSchechtman2009}.
			
			\begin{theo}
				\label{thm:diamond-distortion}
				If $X$ is a Banach space admitting an equivalent asymptotically midpoint uniformly convex norm, then $\sup_{k\in\bN}\cdist{X}(\wdiak)=\infty$. In particular, if the equivalent asymptotically midpoint uniformly convex norm has power type $p\in(1,\infty)$ then $\cdist{X}(\wdiak)\gtrsim k^{1/p}$.
			\end{theo}
			
			\begin{proof}
				Let $C_k:=\cdist{X}(\wdiak)$ and assume without loss of generality that $(X,\norm{\cdot})$ is a Banach space whose norm is asymptotically midpoint uniformly convex. Assume that $C:=\sup_{k\ge 1}C_k<\infty$. It follows from Proposition \ref{pro:diamond-self-improv} and the monotonicity of the modulus that
				\begin{equation*}
					C_{k-1}\le C_k\Big(1-\frac{1}{2}\hat{\delta}_{\norm{\cdot}}\big(\frac{1}{C_k}\big)\Big)\le C_k\Big(1-\frac{\delta}{2}\Big),
				\end{equation*}
				where $\delta:=\hat{\delta}_{\norm{\cdot}}(\frac{1}{C})>0$. Letting $k$ go to infinity gives a contradiction.
				
				If moreover there is a constant $\gamma>0$ such that the modulus of asymptotic midpoint uniform convexity of $\norm{\cdot}$ satisfies $\hat{\delta}_{\norm{\cdot}}(t)\ge\gamma t^p$ for some $p\in(1,\infty)$, then by Proposition \ref{pro:diamond-self-improv} we have for all $2\le j\le k$
				\begin{equation*}
					C_{j-1}\le C_j\Big(1-\frac{1}{2}\hat{\delta}_{\norm{\cdot}}\big(\frac{1}{C_j}\big)\Big)
					\le C_j\Big(1-\frac{\gamma}{2C_j^p}\Big).
				\end{equation*}
				Therefore, for all $2\le j\le k$, $C_j-C_{j-1}\ge\frac{K}{C_j^{p-1}}$, where $K:=\frac{\gamma}{2}$ and hence
				\begin{equation*}
					C_k\ge\sum_{j=2}^k\frac{K}{C_{j}^{p-1}}+C_1\ge\frac{K(k-1)}{C_{k}^{p-1}}.
				\end{equation*}
				The conclusion follows easily.
			\end{proof}

			\section{\texorpdfstring{Bi-Lipschitz embeddability of diamond graphs into $\co$ and $\ell_p$}{Bi-Lipschitz embeddability of diamond graphs into }}
			\label{sec:diamond-co}

			In the previous sections, we have seen that $(\wdiak)_{k\in \bN}$ admit equi-bi-Lipschitz embeddings into $L_1$ but not into $\ell_1$ (since $\ell_1$ is asymptotically (midpoint) uniformly convex).
			As every separable metric space, $\wdiak$ embeds bi-Lipschitzly into $\co$ with distortion at most $2$. First, we will give another embedding into $\co$ which is tailored to the countably branching diamond graphs. This embedding will have special features that are crucial in order to generalize this embedding to other Banach spaces.
			
			\begin{theo}
				\label{thm:diamond-co}
				For all $k\in \bN$, there is $f_k\colon \wdiak \to \co$ such that $\dist(f_k)\le 3$ and for all $x\in \wdiak$, $\abs{\supp(f_k(x))}\le k+1$. 
			\end{theo}
			
			Before we prove Theorem \ref{thm:diamond-co}, we derive an immediate corollary regarding the $L_p$-distortion of the countably branching diamond graphs.
			
			\begin{coro}
				\label{cor:diamond-lp}\,
				\begin{enumerate}[(i)]
					\item For $1 \le p < \infty$, $\cdist{\ell_p}(\wdiak)\approx k^{1/p}$.
					\item  For $1 < p < \infty$, $\cdist{L_p}(\wdiak)\approx \min\{k^{1/p}, \sqrt{k}\}$.
				\end{enumerate}    
			\end{coro}
			
			\begin{proof}\ \\
				$(i)$ The upper estimate simply follows from Theorem \ref{thm:diamond-co}. The lower estimate for the distortion follows from Theorem \ref{thm:diamond-distortion} and the fact that $\ell_p$ is asymptotically uniformly convex with power type $p$.\\
				$(ii)$ Since $L_p$ contains a linearly isometric copy of $\ell_2$ and an isometric copy of $\ell_p$, the upper estimate follows again from Theorem \ref{thm:diamond-co}. The lower estimate follows from Theorem \ref{thm:diamond-distortion} and the fact that $L_p$ is asymptotically uniformly convex with power type $\max\{p,2\}$.
			\end{proof}
			
			The embedding into $\co$ will be of the form $f_k((s,u)):= \sum_{t\preceq s} d(u,\abs{t}) e_{t}$, where $(e_t)_{t\in \N^{\le k}}$ is the canonical basis of $c_0(\N^{\le k})$ and for some well chosen coefficients $(d(u,\abs{t}))_{t\preceq s}$. Given a height $0\le w \le 2^{k}$ and $1\le i \le k+1$, we will define $d(w,i)$ to be the shortest distance one must travel from any vertex at height $w$ in order to reach a vertex of depth at most $i-1$. Note that the later vertex could be located above or below the vertex at height $w$. To properly define $d(w,i)$ we introduce two auxiliary parameters:
			\begin{equation}
				\underline{d}(w,i):=\begin{cases}
					\max\{0\le w' \le 2^k \colon \exists z=(r,w') \text{ s.t. } \text{depth}(r) <i \text{ and } w'\le w \},\ \text{ if } 1\le i \le k+1\\
					0, \qquad \text{ if } i=0
				\end{cases} 
			\end{equation}
			and
			\begin{equation}
				\overline{d}(w,i):=\begin{cases}
					\min\{0\le w' \le 2^k \colon \exists z=(r,w') \text{ s.t. }  \text{depth}(r)<i \text{ and } w'\ge w \},\  \text{ if } 1\le i \le k+1\\
					2^{k}, \qquad \text{ if } i=0
				\end{cases} 
			\end{equation}
			Note that $\overline{d}(w,i)$ records the first height not lesser than $w$ that has vertices of depth at most $i-1$, while $\underline{d}(w,i)$ records the last height not greater than $w$ that has vertices of depth at most $i-1$.
			Therefore, if we let 
			\begin{equation}
				d(w,i):= \begin{cases}
					\min\{ w - \underline{d}(w,i), \overline{d}(w,i) -w\}, \qquad \text{ if } 1\le i \le k+1\\
					w, \qquad \text{ if } i=0
				\end{cases} 
			\end{equation}
			then $d(w,i)$ is indeed the shortest distance one must travel to reach a vertex of depth at most $i-1$ from any vertex at height $w$.
			
			Clearly, $d(w,1) = \min\{2^{k} - w, w\}$ as the only two vertices of depth $0$ are the top and bottom vertices. Also, if $w$ is the height of a vertex of depth $j$, then $d(w,i)=0$ as soon as $i>j$.

			\begin{proof}[Proof of Theorem \ref{thm:diamond-co}]
				Let $f_k\colon \wdiak\to \co(\bN^{\le k})$ defined by 
				\begin{equation*}
					f_k((s,u)):= \sum_{t\preceq s} d(u,\abs{t}) e_{t}.
				\end{equation*}
				Let us first show that $f_k$ is $1$-Lipschitz. So, fix $x=(s,u)$ and $y=(t,v)$ such that $\{x,y\}$ is an edge in $\wdiak$. We can assume without loss of generality that $s\preceq t$ and $\abs{u-v}=1$. Then,
				\begin{align}
					\notag    \norm{ f_k(x)- f_k(y)}_\infty & = \Big\|\sum_{r\preceq s} d(u,\abs{r}) e_{r} - \sum_{r\preceq t} d(v,\abs{r}) e_{r}\Big\|_\infty\\ 
					\label{eq1:diamond-co}                           & = \Big\|\sum_{r\preceq t} (d(u,\abs{r}) - d(v,\abs{r})) e_{r}\Big\|_\infty\\                        
					\notag                           & \le \max_{0\le j\le \abs{t}} \{ \abs{d(u,j) - d(v,j)} \}\\
					\label{eq2:diamond-co}                           & \le 1,
				\end{align}
				where in $\eqref{eq1:diamond-co}$ we use the fact that $s\preceq t$ and $d(u,\abs{r})=0$ whenever $s \prec r$ and in \eqref{eq2:diamond-co} the fact $\abs{d(u,j) - d(v,j)} \le 1$ whenever $\abs{u-v}=1$.
				
				We now estimate the co-Lipschitz constant of $f_k$. Write again $x=(s,u)$ and $y=(t,v)$. If $x$ and $y$ are in the same branch, then
				\begin{equation*}
					\norm{f_k(x) - f_k(y)}_\infty \ge \norm{ (d(u,0) - d(v,0)) e_{\emptyset}}_\infty =  \abs{d(u,0) - d(v,0) } = \abs{u-v} = \dwdiak(x,y).
				\end{equation*}
				If $x$ and $y$ are not in the same branch, then 
				\begin{equation}
					\label{eq3:diamond-co}
					\norm{f_k(x) -f_k(y)}_\infty \ge \max\Big\{ d(u,\abs{s\land t} + 1), d(v,\abs{s\land t} + 1),  \abs{d(u,0) - d(v,0)} \Big\}.
				\end{equation}
				Assuming without loss of generality that $u\preceq v$, recall that 
				\begin{equation*}
					h_{a}(x,y):= \max\{ 0\le w \le 2^k \colon \exists z=(r,w) \text{ s.t. } r\preceq s\land t \text{ and } w\le u\}, 
				\end{equation*}
				and
				\begin{equation*}
					l_{a}(x,y):= \min\{ 0\le w \le 2^k \colon \exists z=(r,w) \text{ s.t. } r\preceq s\land t \text{ and } w\ge v\}. 
				\end{equation*}
				Observing that 
				\begin{equation*}
					h_a(x,y)=\underline{d}(u,\abs{s\land t} +1)=\underline{d}(v,\abs{s\land t} +1)
				\end{equation*}
				and 
				\begin{equation*}
					l_a(x,y)=\overline{d}(v,\abs{s\land t} +1)=\overline{d}(u,\abs{s\land t} +1),
				\end{equation*}
				we have that 
				\begin{equation*}
					d(u,\abs{s\land t} + 1) = \min\{ u - \underline{d}(u,\abs{s\land t}+1), \overline{d}(u,\abs{s\land t}+1) -u\} = \min \{ u - h_a(x,y), l_a(x,y) - u\}
				\end{equation*}
				and 
				\begin{equation*}
					d(v,\abs{s\land t} + 1) = \min\{ v - \underline{d}(v,\abs{s\land t}+1), \overline{d}(v,\abs{s\land t}+1) - v\} = \min \{ v - h_a(x,y), l_a(x,y) - v\}.
				\end{equation*}
				Therefore, it follows from \eqref{eq3:diamond-co} that 
				\begin{equation*}
					\norm{f_k(x) -f_k(y)}_\infty \ge \max\Big\{\min \{ u - h_a(x,y), l_a(x,y) - u\}, \min\{ v - h_a(x,y), l_a(x,y) - v\}, \abs{u-v} \Big\}.
				\end{equation*}
				
				We distinguish two cases.\\
				Either $l_a(x,y) - v \le v - h_a(x,y) $ and then (since $\abs{u-v} =v-u$),
				\begin{equation*}
					\norm{f_k(x) -f_k(y)} \ge \frac{ 2 (l_a(x,y) -v) + v-u }{3} = \frac{2l_a(x,y) - (u+v)}{3}.
				\end{equation*}
				Or $v - h_a(x,y) \le l_a(x,y) - v$, in which case $u - h_a(x,y) \le l_a(x,y) - u$ and hence 
				\begin{equation*}
					\norm{f_k(x) -f_k(y)}_\infty \ge \frac{u - h_a(x,y) +v - h_a(x,y)}{2} = \frac{ u + v - 2h_a(x,y)}{2}.
				\end{equation*} 
				Therefore, 
				\begin{equation*}
					\norm{f_k(x) -f_k(y)}_\infty \ge \frac{\min\{u + v - 2h_a(x,y), 2l_a(x,y) - (u+v)\}}{3},
				\end{equation*}
				and, by Lemma \ref{lem:formula_diamondmetric}, the proof is complete. 
			\end{proof}

			\section{Bi-Lipschitz embeddability of diamond graphs and Szlenk index}
			\label{sec:diamond-applications}
			
			A careful analysis of the $\co$-embedding of the diamond graphs reveals that it can be significantly generalized. 
			What is needed in order to implement an embedding as in the proof of Theorem \ref{thm:diamond-co} is the existence of constants $A,B\in [1,\infty)$ and for $k\in\bN$ a semi-normalized tree $(x_{t})_{t\in \bN^{\le k}}$ in a Banach $X$ together with a bijection $\phi\colon \bN^{\le k} \to \bN$ satisfying the following conditions:
			\begin{enumerate}
				\item $\forall (a_i)_{i=0}^k \in \bR^k\ \forall t\in \bN^k,\ \frac{1}{A}\max_{0\le i \le k}\abs{a_i} \le \big\|\sum_{i=0}^k a_i x_{t_{\restriction i}}\big\| \le A\max_{0\le i \le k}\abs{a_i}$,
				\item  $\phi(s)\le \phi(t)$ whenever $s\preceq t$,
				\item $\big\|\sum_{i=0}^{j_1} a_i x_{\phi^{-1}(i)}\big\| \le B \big\|\sum_{i=0}^{j_2} a_ix_{\phi^{-1}(i)} \big\|$ whenever $j_1\le j_2$ and $(a_i)_{i=0}^{j_2} \in \bR^{j_2}$.
			\end{enumerate}
			
			Condition $1.$ means that every branch of the tree is $A^2$-equivalent to the canonical basis of $\ell_\infty^{k+1}$ and condition $2.$ and $3.$ mean that there is an order-preserving enumeration of the tree according to which the vertices of the tree form a $B$-basic sequence
			
			For conciseness, we will refer to a semi-normalized tree satisfying condition $1$ to $3$ above as an \emph{$A$-$\ell_\infty$-$B$-basic tree of height $k$}. Note that in that case $\{\norm{x_t}\colon t\in \bN^{\le k}\} \subset [\frac{1}{A},A]$ and we will say that the tree is $A$-semi-normalized.  Of course, the tree given by $x_t =e_t$ for all $t\in \bN^{\le k}$, where $(e_{t})_{t\in \bN^{\le k}}$ is the canonical basis of $c_0(\bN^{\le k})$ is an example of a $1$-$\ell_\infty$-$1$-basic tree of height $k$. The proof of Theorem \ref{thm:diamond-co} can be slightly adjusted to obtain the following embedding result.
			
			\begin{theo}
				\label{thm:diamond-tree}
				Let $k\in \bN$ and $A,B\in [1,\infty)$.
				If a Banach space $X$ contains an $A$-$\ell_\infty$-$B$-basic tree of height $k$, then $\wdiak$ embeds with distortion at most $3A^2B(1+B)$ into $X$.
			\end{theo}
			
			\begin{proof}
				In this setting the natural embedding is given by $f_k\colon \wdiak\to X$ where 
				\begin{equation*}
					f_k((s,u)):= \sum_{t\preceq s} d(u,\abs{t}) x_{t}.
				\end{equation*}
				
				It follows from the left handside inequality of condition $1.$ that $f_k$ is $A$-Lipschitz.
				
				For the co-Lipschitz constant of $f_k$, we use the $B$-basicity property. Since necessarily $\phi(\emptyset)=0$, if $x$ and $y$ are in the same branch, then
				\begin{align*}
					BA\norm{f_k(x) - f_k(y)} &\ge A\norm{ (d(u,0) - d(v,0)) x_{\emptyset}}\\
					&\ge  \abs{d(u,0) - d(v,0) } = \abs{u-v} = \dwdiak(x,y).
				\end{align*}
				If $x$ and $y$ are not in the same branch, then we define $(a_i)_{i=0}^\infty\in \bR^\infty$ as follows:
				\begin{equation*}
					a_i := \begin{cases}
						d(u, |\phi^{-1}(i)| ) - d(v,  |\phi^{-1}(i)|) \qquad \text{ if } \phi^{-1}(i)\preceq s \land t\\
						d(u, |\phi^{-1}(i)| ) \qquad \text{ if } s\land t \prec \phi^{-1}(i)\preceq s,\\
						-d(v, |\phi^{-1}(i)|) \qquad \text{ if } s\land t \prec \phi^{-1}(i)\preceq t, \\
						0 \qquad \text{ otherwise.}
					\end{cases}
				\end{equation*}
				Note that condition $2.$ on the enumeration implies that there are only finitely many $a_i$ that are nonzero and 
				\begin{equation*}
					\norm{f_k(x) -f_k(y)} = \Big\|\sum_{i=0}^\infty a_i x_{\phi^{-1}(i)}\Big\|.
				\end{equation*}
				Using condition $3.$, if $r_s$ is the immediate successor of $s\land t$ in $(s\land t, s]$ and $r_t$ is the immediate successor of $s\land t$ in $(s\land t, t]$, then $\abs{r_s}=\abs{r_t}=\abs{s\land t} +1$ and we get
				\begin{align*}
					\|f_k(x) &-f_k(y)\| \ge \frac{1}{B} \max \Big\{\Big\|\sum_{i=0}^{\phi(r_s)} a_i x_{\phi^{-1}(i)} \Big\|, \Big\|\sum_{i=0}^{\phi(r_t)} a_i x_{\phi^{-1}(i)} \Big\|, \norm{a_0 x_\emptyset} \Big\}\\
					& \ge \frac{1}{B(1+B)} \max \Big\{\|a_{\phi(r_s)} x_{r_s} \|, \|a_{\phi(r_t)} x_{r_t} \|, \norm{a_0 x_\emptyset} \Big\}\\
					& \ge \frac{1}{BA(1+B)} \max \Big\{ d(u, \abs{r_s}) , d(v, \abs{r_t}) , \abs{ a_0 } \Big\}\\
					& = \frac{1}{BA(1+B)} \max\Big\{ d(u,\abs{s\land t} + 1), d(v,\abs{s\land t} + 1),  \abs{d(u,0) - d(v,0)} \Big\},
				\end{align*}
				and we conclude as in Theorem \ref{thm:diamond-co}.
			\end{proof}
			
			In the remainder of this section, we will show that, for every $\vep>0$, any Banach space with an unconditional asymptotic structure and whose Szlenk index of its dual is large contains $(1+\vep)$-basic-$(1+\vep)$-$\ell_\infty$ trees of every finite but arbitrarily large height.
			This will provide a fairly large class of Banach spaces containing equi-bi-Lipschitz copies of $(\wdiak)_{k\in\bN}$, including $\co$ and any asymptotic-$\co$ Tsirelson space with an unconditional basis.
			
			In Section \ref{sec:tree-Szlenk} of Chapter \ref{chapter:trees}, we have already seen that if $\Sz(X^*)>\omega$, then for every $k\in\bN$, $X$ contains a semi-normalized weakly null tree $(x_{t})_{t\in \bN^{\le k}}$ such that $\big\|\sum_{s\le t} x_s\big\|\le 4$. This property does not seem sufficient to prove that the embedding map in Theorem \ref{thm:diamond-tree} is Lipschitz as we need to handle coefficients that are not necessarily all equal. This issue can be remedied if we assume that the space has an unconditional asymptotic structure (see Appendix \ref{sec:asymptotic-structure}). But first, we need to prove a tree-analog of the classical fact that a weakly null sequence has a basic subsequence.
			
			\begin{prop}
				\label{prop:basic-tree}
				Let $X$ be a Banach space with separable dual and let $(x_t)_{t\in \bN^{\le k}}$ be a semi-normalized weakly null tree in $X$ and $\vep>0$. Then, there is a pruning $\psi\colon \bN^{\le k} \to \bN^{\le k}$ and an order-preserving bijection $\varphi\colon \psi(\bN^{\le k}) \to \bN$ such that $(x_{\varphi^{-1}}(n))_{n\in \bN}$ is a $(1+\vep)$-basic sequence.
			\end{prop}
			
			\begin{proof}
				Since $X^*$ is separable we choose a dense sequence $(x^*_{n})_{n\in\bN}\in S_{X^*}$. Let $\vep>0$ be fixed and $\delta>0$ be chosen small enough later. We also choose $\phi \colon \bN^{\le k} \to \bN$ an order-preserving bijection with the additional property that for all $t\in \bN^{\le k-1}$ and $m<n$ we have $\phi(t\frown m) < \phi(t\frown n)$. Let $\tau_n:=\phi^{-1}(n)$ for all $n\in \bN$. For every $i\in \bN$ we will construct $\gamma_i\in \bN^{\le k}$, a finite set $N_i\subset S_{X^*}$ and $z_i\in S_X$ satisfying the following conditions:
				\begin{enumerate}
					\item  $\gamma_1=\emptyset$ and if $i<j$ and $\tau_j=\tau_i\frown m$ for some $m\in \bN$, then $\gamma_j=\gamma_i\frown n$ for some $n\in \bN$ that is larger than every integer appearing as a coordinate of an element in $\{ \gamma_s \colon s<j\}$.
					\item $N_j\supset \{x^*_i \colon 1\le i\le j\}$,
					\item $\max_{x^*\in \N_j}\abs{x^*(x)}\ge \frac{\norm{x}}{1+\delta}$ for all $x\in span\{z_i \colon 1\le i\le j\}$,
					\item $x^*(z_j)=0$ for all $x^*\in N_i$ whenever $i<j$,
					\item $\norm{z_j - x_{\gamma_j}}<\delta$.
				\end{enumerate}
				For $j=1$ we simply let $\gamma_1=\emptyset$, $z_1=x_{\emptyset}$ and $N_1:=\{x_1^*, z_1^*\}$ where $z^*_1\in s_{X^*}$ is such that $z_1^*(z_1)=1$. Condition $4.$ is vacuously true while the others are easily verified.
				Assume now that $\gamma_1,\dots, \gamma_{i}$, $N_1,\dots, N_{i}$, $z_1,\dots,z_{i}$ have been constructed for $i\ge 1$. Observing that $\tau_{i+1}=\tau_{i_0}\frown m_0$ for some $m_0\in \bN$ and $i_0\le i$ and that $(x_{\gamma_{i_0}\frown m})_{m\in \bN}$ is weakly null, we can find an $m_1\in \bN$ that is larger than every integer appearing as a coordinate of an element in $\{ \gamma_s \colon s\le i\}$ and a vector $z_{i+1}$ in the unit sphere of the finite-codimensional space $Y:=\{x\in X\colon x^*(x)=0, \forall x^*\in N_i\}$ such that $\norm{x_{\gamma_{i_0}\frown m_1} - z_{i+1}} \le \delta$. If we let $\gamma_{i+1}:=\gamma_{i_0}\frown m_1$, then condition $1.$ and $5.$ are satisfied. Now, if $\{u_j\}_{j\in J}$ is a finite $\eta$-net of the unit sphere of $span\{z_j\colon 1\le j\le i+1\}$ and $(u^*_j)_{j\in J}\subset S_X{^*}$ such that $u^*_j(u_j)=1$ for all $j\in J$, then $N_{i+1}:=N_i\cup \{x^*_{i+1}, u_j \colon j\in J\}$ clearly satisfies $2.$ and $4.$ and it is standard to check that it satisfies $3.$ if $\eta$ is taken small enough. 
				This concludes the induction. It remains to be observed that by condition $1.$, the map $\psi \colon \bN^{\le k} \to \bN^{\le k}$ given by $\psi(t) = \gamma_{\phi(t)}$ is a pruning since $\psi(\tau_n)=\gamma_n$. If we let $\varphi:= \gamma^{-1}$, then $\varphi\colon \psi(\bN^{\le k}) \to \bN$ is an order-preserving bijection such that  $(x_{\varphi^{-1}(n)})_{n\in \bN}$ is a $\delta$-perturbation of $(z_n)_{n\in\bN}$. But conditions $3-4.$ imply that $(z_n)_{n\in\bN}$ is a $(1+\delta)$-basic sequence. It is now standard that we can choose $\delta$ small enough such that $(x_{\varphi^{-1}(n)})_{n\in \bN}$ is $(1+\vep)$-basic. 
			\end{proof}
			
			\begin{prop}
				\label{prop:Sz-linfty-trees}
				Let $X$ be a Banach space with a separable dual and an unconditional asymptotic structure. If $\Sz(X^*)>\omega$, then there is a constant $A>0$ and for every $k\in \bN$ a weakly null semi-normalized tree $(x_{t})_{t\in \bN^{\le k}}\subset X$ so that for all $(a_i)_{i=0}^k \in \bR^k$ and $t\in \bN^k$,
				\begin{equation}
					\frac{1}{A}\max_{0\le i \le k}\abs{a_i} \le \Big\|\sum_{i=0}^k a_i x_{t_{\restriction i}} \Big\| \le A\max_{0\le i \le k}\abs{a_i}.
				\end{equation}
			\end{prop}
			
			\begin{proof}
				Assume that $X^*$ is separable. Then, it follows from Chapter \ref{chapter:trees}, Section \ref{sec:tree-Szlenk} that for all $k\in \bN$ there is a weakly null rooted tree $(x_{t})_{t\in \N^{\le k}}$ in $X$ so that for all $t\in \N^{\le k}$, $1\le \norm{ x_{t} }\le 6$ and $\norm{ \sum_{i=0}^k x_{t_{\restriction i}} }\le 4$. Invoking Proposition \ref{prop:basic-tree} we can assume without loss of generality that there is an order-preserving enumeration of $(x_{t})_{t\in \N^{\le k}}$ for which it is $B$-basic sequence for $B\in(1,2]$. Since the vertices of the tree have norm at least $1$ and since the enumeration is order-preserving, it follows that for all $(a_i)_{i=0}^k\in \bR^{k+1}$ and $t\in \bN^k$,
				\begin{equation}
					\label{eq3:Sz-linfty-trees}
					\frac{1}{2B} \max_{0\le i\le k} \abs{a_i} \le \Big\|\sum_{i=0}^k a_i x_{t_{\restriction i}} \Big\|.
				\end{equation}
				For the upper bound. Let $\phi\colon \bN^{\le k} \to \bN^{\le k}$ be a pruning, then $(x_{\phi(t)})_{t\in \N^{\le k}}$ is still a weakly null semi-normalized tree in $X$ satisfying $\norm{ \sum_{i=0}^k x_{\phi(t)_{\restriction i}} }\le 4$. Thus, it follows from Proposition \ref{prop:UAS-tree} and the pruning lemma that if $X$ has a $K$-unconditional asymptotic structure, then for every $t\in \bN^{k}$ and  $\vep\in\{-1,1\}^{k+1}$,  
				\begin{equation}
					\label{eq:Sz-linfty-trees}
					\Big\|\sum_{i=0}^k \vep_i x_{\phi(t)_{\restriction i}} \Big\|\le 6K\Big\|\sum_{i=0}^k x_{\phi(t)_{\restriction i}} \Big\| \le 24K.
				\end{equation}    
				It remains to apply a classical extreme point argument. Given $t\in \bN^k$, the map $\psi\colon a:=(a_i)_{i=0}^{k+1}\in [-1,1]^{k+1} \mapsto \norm{\sum_{i=0}^k a_i x_{\phi(t)_{\restriction i}} }$ is convex. Since $[-1,1]^{k+1}$ is the convex hull of $\{-1,1\}^{k+1}$, we have that 
				$\sup_{a\in [-1,1]^{k+1} } \psi(a) = \max_{\vep\in \{-1,1\}^{k+1}} \psi(\vep)$ and we deduce from \eqref{eq:Sz-linfty-trees} and homogeneity that for all $(a_i)_{i=0}^k\in \bR^{k+1}$ and $t\in \bN^k$, 
				\begin{equation}
					\label{eq2:Sz-linfty-trees}
					\Big\|\sum_{i=0}^k a_i x_{\phi(t)_{\restriction i}} \Big\|\le 24K  \max_{0\le i\le k} \abs{a_i}.
				\end{equation}
			\end{proof}
			
			By putting together the results from Section \ref{sec:diamond-amuc} and 
			this section we obtain a metric characterization of the class $\langle \AUC \rangle$ within the class $\langle \REF \& \mathbf{UAS} \rangle$ of reflexive spaces with an unconditional asymptotic structure.
			
			\begin{theo}
				\label{thm:AUC-diamonds}
				Let $X$ be a reflexive Banach space with an unconditional asymptotic structure. Then,
				$X$ admits an equivalent norm that is asymptotically uniformly convex if and only if $\sup_{k\in \bN}\cdist{X}(\wdiak)=\infty$.   
			\end{theo}
			
			\begin{proof}
				If $X$ is reflexive but not in $\langle \AUC\rangle$, then Theorem \ref{thm:Sz-omega=AUS} and the $\AUS$-$\AUC$ duality tells us that $\Sz(X^*)>\omega$. Thus, there exists a closed subspace $Z$ of $X^*$ such that $X^*/Z$ is separable and $\Sz(X^*/Z)>\omega$ (see \cite{Lancien1996}). Moreover, if $X$ has an unconditional asymptotic structure, so does $Z^\perp \subset X$ (remember that $X$ is reflexive). In other words, we may assume that $X$ is separable. Then, it follows from Proposition \ref{prop:Sz-linfty-trees} that are constants $A,B>0$ and for every $k\in \bN$ a weakly null semi-normalized tree $(x_{t})_{t\in \bN^{\le k}}\subset X$ so that for all $(a_i)_{i=0}^k \in \bR^k$ and $t\in \bN^k$,
				\begin{equation}
					\frac{1}{A}\max_{0\le i \le k}\abs{a_i} \le \Big\|\sum_{i=0}^k a_i x_{t_{\restriction i}} \Big\| \le A\max_{0\le i \le k}\abs{a_i}.
				\end{equation}
				By Proposition \ref{prop:basic-tree}, we can assume without loss of generality that the trees are $A$-$\ell_\infty$-$2$-basic trees, and it follows from Theorem \ref{thm:diamond-tree} that $\sup_{k\in \bN}\cdist{X}(\wdiak)<\infty$.
				
				Assume now that $X\in \langle \AUC\rangle$, then $X$ admits an equivalent norm that is asymptotically midpoint uniformly convex, and it follows from Theorem \ref{thm:diamond-distortion} that $\sup_{k\in \bN}\cdist{X}(\wdiak)=\infty$
			\end{proof}
			
			\begin{rema}
				Note that for the nonembeddability direction in Theorem \ref{thm:AUC-diamonds} we do not need the reflexivity nor the unconditionally condition.    
			\end{rema}

			\section{Notes}
			\label{sec:diamond-notes}
			
			The results about the countably branching diamond graphs presented in this chapter are essentially taken from \cite{Baudieretal}. However, on many occasions, the presentation is greatly inspired by Swift's exposition in \cite{Swift2018} where those results were extended to more general bundle graphs. The difficulty in handling more general bundle graphs lies in the complexity of the coding of these graphs that is used to define the embedding map. Also, for the nonembeddability of these bundle graphs, which are iterated $\oslash$-products of a certain base graph, the base graph needs to satisfy some geometric properties so that the embeddability obstruction, which relies on the approximate midpoint set size, can be implemented. We refer to \cite{Swift2018} where solutions to overcome these intricacies are proposed. 
			We do not know if the assumption of asymptotic unconditional structure is truly needed for Theorem \ref{thm:AUC-diamonds}. We suspect it is not, but it seems to require a significant modification of the argument or a completely new approach to the embedding problem altogether. In this direction, we can mention recent results by Basset, Lancien and Proch\'azka \cite{BLPinprepa} who characterize the equi-sub-Lipschitz embeddability (in the sense of Ostrovskii \cite{Ostrovskii2014}) of the countably branching diamond graphs in terms of a countably branching analog of the classical Finite Tree Property that characterizes non super-reflexivity. 
			
			\begin{prob}
				\label{prob:AUC-diamonds}
				If $X$ is a reflexive Banach that does not admit an equivalent norm that is asymptotically uniformly convex, is it true that $\sup_{k\in \bN}\cdist{X}(\wdiak)<\infty$?
			\end{prob}
			
			In this chapter, we presented a technique to show the nonembeddability of a sequence of graphs that relies on a self-improvement argument. This technique was first used by Johnson and Schechtman \cite{JohnsonSchechtman2009} for the binary diamond graphs and later by Kloeckner \cite{Kloeckner2014} to give a simple proof of the nonembeddability of the binary trees in Bourgain's metric characterization theorem. Eventually, it was used to show the nonembeddability of the countably branching trees in \cite{BaudierZhang16} (see Exercise \ref{exer:self-improvement-tree}). Nonembeddability results for diamond graphs can also be derived from a metric invariant approach. For the binary diamond graphs, this was originally due to Lee and Naor in \cite{LeeNaor2004}. Refined Poincar\'e-type inequalities that capture the geometry of diamond graphs were further developed by Eskenazis, Mendel and Naor (see \cite{Eskenazis_PhD} and \cite{EMN}) for binary diamonds and for countably branching diamonds in \cite{AcevesBaudier}. 
			%(see Exercise \ref{exer:bicone-convexity}). 
			The problem below is a fundamental open problem in the Kalton program.
			
			\begin{prob}
				\label{pb:q-AUC}
				Given $q\in(1,\infty)$, find a purely metric characterization of the class of reflexive Banach spaces admitting an equivalent norm that is $q$-asymptotically uniformly convex.    
			\end{prob}
			
			Randrianantoanina and Ostrosvskii \cite{OR2017} have extended the Johnson-Schechtman metric characterization to finitely branching diamonds. For this purpose, they used a new embedding technique that relies on the concept of equal-signs-additive (ESA) sequences developed by Brunel and Sucheston in \cite{BrunelSucheston1974,BrunelSucheston1975,BrunelSucheston1976}. In this chapter, we have chosen to focus on diamond graphs, but there are alternative characterizations in terms of the Laakso graphs (in both the local \cite{JohnsonSchechtman2009,OR2017} and asymptotic \cite{Swift2018} setting). A significant difference between the (binary) Laakso graphs and the diamond graphs is that while the former are doubling (with a uniform bound on the doubling constant) the latter are not. Finally, Ostrovskii \cite{Ostrovskii2014} has shown that for dual Banach spaces, not having the Radon-Nikod\'ym Property is equivalent to containing a bi-Lipschitz copy of the infinite diamond graph (constructed as the inductive limit of properly rescaled binary diamond graphs). 
			
			\medskip In \cite{Basset}, E. Basset introduces dyadic and countably branching diamonds of countable ordinal height. Then she studies the dentabilty index and the weak fragmentability index of their Lipschitz-free spaces. As an application, she obtains that a separable metric space that is universal for all countable complete metric spaces and bi-Lipschitz embeddings cannot be purely-1-unrectifiable. This is pushed further in \cite{BassetLancienProchazka}, where it is proved that the same is true for a separable metric space universal for all metric spaces that are countable, discrete and complete. It is also shown in \cite{BassetLancienProchazka} that there exist  countable, discrete and complete metric spaces whose free spaces do not linearly embed into a free space of a uniformly discrete metric space. In \cite{BLPinprepa} the bi-Lipschitz embeddability into Banach spaces of these new diamond graphs is investigated in relation with the dentability and weak fragmentability index of the target space or with the Szlenk index of its predual (if any) and linear characterizations of the sub-Lipschitz embeddability of these graphs are obtained. In the spirit of Section \ref{sec:diamond-L1}, their embeddability into $L_1$ is also proved.

			\section{Exercises}
			\label{sec:diamond-exe}
			
			%\begin{exer}
			%Prove Lemma \ref{lem:amuc}. 
			%\end{exer}
			
			\begin{exer}
				\label{exe:diamond1-bernoulli}
				Show that there are functions $f_s,f_t, f_i$ with $i\in \bN$ that are in 
				$L_1$ and such that 
				\begin{equation}
					\norm{f_s - f_t}_{L_1} = 2\norm{f_s - f_i}_{L_1} = 2\norm{f_i - f_t}_{L_1} = 2\ \ \text{for all}\ i \in \N
				\end{equation}
				and 
				\begin{equation}
					\norm{f_i - f_j}_{L_1} = 2\ \ \text{for all}\ i\neq j \in \N.
				\end{equation}
			\end{exer}
			
			In the exercise, we invite the reader to discover the original proof, taken from \cite{Baudieretal}, of Lemma \ref{lem:terminals-contraction-bis}. Recall that the Kuratowski measure of noncompactness of a subset $S$ of a metric space, denoted by $\alpha(S)$, is defined as the infimum of all $\vep>0$ such that $S$ can be covered by a finite number of sets of diameter at most $\vep$.
			
			\begin{exer}
				\label{exe:AMUC-original} 
				Let $(X,\norm{\cdot})$ be a Banach space.
				\begin{enumerate}
					\item Prove that for every $t\in(0,1)$ and every $x\in S_X$ one has 
					\begin{equation*}
						\alpha\Big(\Mid\big(-x,x,\frac14\hat{\delta}_{\norm{\cdot}}(t)\big)\Big)<4t.
					\end{equation*}
					
					\item Deduce that for every $t\in(0,1)$ and every $x,y\in X$ there exists a finite subset $S$ of $X$ such that
					\begin{equation*}
						\Mid\big(x,y,\frac14\hat{\delta}_{\norm{\cdot}}(t)\big)\subset S + 2t\norm{x-y}B_X .
					\end{equation*}
					
					\item Assume now that $f\colon \mathsf{D_1^\omega}\to X$ is a bi-Lipschitz embedding with distortion $C\ge 1$ and $\Lip(f)=1$. 
					\begin{enumerate}
						\item Prove that there exists $j\in\mathbb{N}$ such that
						\begin{equation*}
							f(m_j)\notin\Mid\Big(f(v_t),f(v_b),\frac{\delta}{4}\Big).\end{equation*}
						\item deduce that 
						\begin{equation*}
							\norm{f(v_t)-f(v_b)} \le \Lip(f)\Big(1-\frac{1}{5}\hat{\delta}_{\norm{\cdot}}\big(\frac{1}{8C}\big)\Big)\sd_{\sD_1^\omega}(v_t,v_b).
						\end{equation*}
					\end{enumerate}
				\end{enumerate}

			\end{exer}
			
			\begin{exer}\, 
				\begin{enumerate}
					\item Show that for all $k\in \bN$ and $\vep>0$ there exists $p:=p(k,\vep)>0$ such that $\cdist{\ell_p}(\wdiak)\le 3+\vep$.
					\item Let $X:=(\sum_{n=1}^\infty \ell_{q_n})_{\ell_2}$ with $\lim_{n\to\infty}q_n=\infty$. Show that for all $k\in \bN$, $\cdist{X}(\wdiak)\le 3$.
				\end{enumerate}
			\end{exer}

			\begin{appendix}
				
				\chapter{Banach space theory essentials}
				\label{appendix:Banach}
				In this appendix, we recall classical definitions and results that are needed in the text and sometimes used without reference. If the reader wants to consult the proofs, they can be found in the classical textbooks on Banach space theory, such as \cite{Banach1987}, \cite{DunfordSchwartz1958}, \cite{Day1973}, \cite{DiestelUhl1977}, \cite{LindenstraussTzafriri1977}, \cite{LindenstraussTzafriri1979}, \cite{Diestel1984}, \cite{Beauzamy1985}, \cite{DiestelJarchowTonge1995}, \cite{Megginson1998}, \cite{FHHMPZ2001}, \cite{Carothers2005}, \cite{AlbiacKalton2016}, \cite{FHHMZ2011}, \cite{QueffelecLi_vol1}, and \cite{QueffelecLi_vol2}. 
				
				\section{A handful of useful functional analysis results}
				\label{appendix:FA}
				\begin{itemize}
					\item[] We first recall one of many versions of a very useful consequence of the Hahn-Banach theorem originally due to Helly (see for instance \cite[Lemma 3.3]{Brezis2011})  
					\item Let $X$ be a normed space, $x_1^*,\dots,x_n^* \in X^*$, $M>0$ and $c_1,\dots,c_n$ scalars. Then the two following conditions are equivalent:
					\begin{enumerate}[(i)]
						\item For all scalars $a_1,\dots,a_n$, $\big|\sum_{i=1}^na_ic_i\big|\le M\big\|\sum_{i=1}^n a_ix_i^*\big\|_{X^*}$. 
						\item For every $\eps>0$, there exists $x \in X$ such that $\norm{x}_X=M+\vep$ and $x_i^*(x)=c_i$, for all $i\in \{1,\dots,n\}$. 
					\end{enumerate} 
					\item[] Although the weak topology is not always metrizable, even when restricted to bounded subsets, we nevertheless have the following striking result essentially due to Rosenthal \cite{Rosenthal1978} (see also \cite{GonzalezGutierrez1992}). 
					\item Let $X$ be a Banach space that contains an isomorphic copy of $\ell_1$. Then, every bounded subset of $X$ is weakly sequentially dense in its weak closure.
					\item[] An equivalent norm on the dual $X^*$ of a Banach space $X$ is not necessarily the dual norm of an equivalent norm on $X$. The next result gives a necessary and sufficient condition guaranteeing that it is.
					\item An equivalent norm on the dual $X^*$ is the dual norm of an equivalent norm on $X$ 
					if and only if it is weak$^*$ lower-semicontinuous.
					\item[]Let $X$ and $Y$ be two Banach spaces and $B(X,Y)$ be the Banach space of bounded linear operators from $X$ into $Y$. The dual of an operator is easily seen to be weak$^*$-to-weak$^*$ continuous. The next result is about the converse.
					\item Let $X$ and $Y$ be two Banach spaces and $T \in B(Y^*,X^*)$. Then $T$ is weak$^*$-to-weak$^*$ continuous if and only if there exists $S \in B(X,Y)$ such that $T=S^*$. 
					\item The \emph{weak operator topology} is the topology on $B(X,Y)$ generated by the collection of maps $(p_{x,y^*})_{x\in X,y^*\in Y^*}$ where $p_{x,y^*}$ is the linear functional given by $p_{x,y^*}(T)=y^*(Tx)$, for all $T\in B(X,Y)$.
					\item The \emph{dual weak operator topology} is the topology on $B(X,Y)$ generated by the collection of maps $(p_{x^{**},y^*})_{x^{**}\in X^{**},y^*\in Y^*}$ where $p_{x^{**},y^*}$ is the linear functional given by $p_{x^{**},y^*}(T)=x^{**}(T^*y^*)$, for all $T\in B(X,Y)$.
					\item The weak operator topology is coarser than the dual weak operator topology, and they coincide whenever $X$ is reflexive.
					\item (Brace and Friend \cite[Theorem 8]{BraceFriend1969} and Kalton \cite[Theorem 1]{Kalton1974}). If $A\subseteq K(X,Y)$ is a subset of the closed subspace of compact operators, then $A$ is weakly compact if and only if $A$ is compact for the dual weak operator topology.
					\item (Kalton \cite[Corollary 3]{Kalton1974}) Any sequence of compact operators between two Banach spaces that converge in the dual weak operator topology to a compact operator must converge to the same limit in the weak topology. 
				\end{itemize}
				
				\section{Special bases and basic sequences}
				Standard references for basis theory are \cite{Singer1970} and \cite{LindenstraussTzafriri1977}.
				%(see also \cite{GuerreDelabriere}). 
				\begin{itemize}
					\item (Auerbach basis) A basis $(e_1,\dots,e_n)$ of a finite-dimensional normed space $E$ is called an \emph{Auerbach Basis} of $E$ if for all $1\le i \le n$, $\norm{e_i}_E=1$ and $\norm{e_i^*}_{E^*}=1$, where $(e_1^*,\dots,e_n^*)$ is the family of coordinate functionals associated with the basis $(e_1,\dots,e_n)$. Every finite-dimensional normed space admits an Auerbach basis. 
					\item (Schauder basis) A sequence $(x_n)_{n=1}^\infty$ in a Banach space $X$ is called a \emph{Schauder basis} of $X$ if for all $x\in X$, there exists a unique sequence of scalars $(a_n)_{n=1}^\infty$ such that $x=\sum_{n=1}^\infty a_n x_n$, where the convergence of this series is meant for the norm of $X$. 
					\item (Basic sequence) A sequence $(x_n)_{n=1}^\infty$ in a Banach space $X$ is a \emph{basic sequence} if it is a Schauder basis of its closed linear span. 
					\item A sequence $(x_n)_{n=1}^\infty$ of non-zero vectors in a Banach space $X$ is a basic sequence if and only if there exists a constant $K\ge 1$ such that for all $m\le n$ and all scalars $a_1,\dots,a_n$,
					\begin{equation*}
						\Big\|\sum_{i=1}^m a_ix_i\Big\| \le K \Big\|\sum_{i=1}^n a_ix_i\Big\|.
					\end{equation*}
					When $K=1$, the basic sequence is said to be \emph{monotone}. The basic sequence is called \emph{bimonotone} if for all $p\le q\le n$ and all scalars $a_1,\ldots,a_n$,
					\begin{equation*}
						\Big\|\sum_{i=p}^q a_ix_i\Big\| \le  \Big\|\sum_{i=1}^n a_ix_i\Big\|.
					\end{equation*}
					\item Every Banach space with a Schauder basis admits an equivalent norm for which this Schauder basis is bimonotone. 
					\item A basis $\xn$ of a Banach space $X$ is \emph{shrinking} if $X^*$ is the closed linear span of the biorthogonal functionals associated to $\xn$ or, equivalently, if for all $x^*\in X^*$, $\lim_{n\to\infty}\|x^*_{\restriction \overline{\linspan}\{x_k\colon k\ge n+1\}}\|=0$.
					\item A basis $\xn$ of a Banach space $X$ is \emph{boundedly complete} if for all sequence of scalars $(a_n)_n$, $\sum_n a_nx_n$ is convergent whenever $\sup_{N\ge 1} \norm{\sum_{n=1}^N a_n x_n}<\infty$.
					\item (Equivalence of basic sequences) Two basic sequences $\xn\subset X$ and $\yn\subset Y$ are \emph{equivalent} if one (and hence all) of the following three equivalent conditions holds:
					\begin{enumerate}[(i)]
						\item For every sequence of scalars $(a_n)_{n=1}^\infty$, $\sum_{n} a_n x_n$ converges if and only if $\sum_{n} a_n y_n$ converges,
						\item There is an isomorphism $T\colon X \to Y$ such that, for all $n\ge 1$, $T(x_n)=y_n$,
						\item There is a constant $C>0$, such that for all $(a_n)_{n=1}^\infty\in c_{00}$, 
						\begin{equation}
							\label{eq:equivalence-basic}
							\frac{1}{C}\Big\|\sum_{n=1}^\infty a_n y_n\Big\|_Y \le \Big\|\sum_{n=1}^\infty a_n x_n\Big\|_X \le C \Big\|\sum_{n=1}^\infty a_n y_n\Big\|_Y.
						\end{equation}
					\end{enumerate}
					\item (Permutation-unconditionality) A basic sequence $(x_n)_{n=1}^\infty$ in a Banach space $X$ is \emph{permutation-unconditional} if for any sequence of scalars $(a_n)_{n=1}^\infty$ and any permutation $\sigma$ of $\bN$, we have that $\sum_{n=1}^\infty a_{\sigma(n)}x_{\sigma(n)}$ converges if and only if $\sum_{n=1}^\infty a_n x_{n}$ converges. 
					\item (Absoluteness) A basic sequence $(x_n)_{n=1}^\infty$ in a Banach space $X$ is \emph{$c$-absolute}, for some $c\ge 1$, if  for all $(a_n)_{n=1}^\infty$ and $(b_n)_{n=1}^\infty$ in $c_{00}$ such that $\abs{a_n}\le\abs{b_n}$ for all $n\in \N$, we have 
					\begin{equation}
						\label{eq:absolute-unconditional}
						%\label{E:2.8}
						\Big\|\sum_{n=1}^\infty a_n x_{n}\Big\| \le  c\Big\|\sum_{n=1}^\infty b_n x_{n}\Big\|,
					\end{equation}
					and in this case we denote by $c_a$ the smallest number $c$ satisfying \eqref{eq:absolute-unconditional}.
					\item (Sign-unconditionality) A basic sequence $(x_n)_{n=1}^\infty$ in a Banach space $X$ is \emph{$c$-sign-unconditional} if for all $(a_n)_{n=1}^\infty \in c_{00}$ and all $(\vep_n)_{n=1}^\infty \subset \{-1,1\}$ we have
					\begin{equation}
						\label{eq:sign-unconditional}
						\Big\|\sum_{n=1}^\infty \vep_n a_n x_n\Big\| \le  c\Big\|\sum_{n=1}^\infty a_n x_n\Big\|,
					\end{equation}
					and we denote in this case the smallest number $c$ satisfying \eqref{eq:sign-unconditional} by $c_u$.
					\item (Suppression-unconditionality) A basic sequence $(x_n)_{n=1}^\infty$ in a Banach space $X$ is \emph{$c$-suppression-unconditional}, for some $c\ge 1$, if for all  $(a_n)_{n=1}^\infty \in c_{00}$  and all $A\subset \bN$
					\begin{equation}
						\label{eq:suppresion-unconditional}
						%\label{E:2.7}
						\Big\|\sum_{n\in A}a_n x_{n}\Big\| \le  c\Big\|\sum_{n=1}^\infty a_n x_{n}\Big\|,
					\end{equation}
					and in this case we denote by $c_s$ the smallest number $c$ satisfying \eqref{eq:suppresion-unconditional}.
					\item (Equivalence of unconditionality conditions, Orlicz \cite{Orlicz1929}) Let $\xn$ be a basic sequence in  a Banach space $X$. The following assertions are equivalent: 
					\begin{enumerate}[(i)]
						\item $\xn$ is sign-unconditional,
						\item $\xn$ is suppression-unconditional, 
						\item $\xn$ is absolute, 
						\item $\xn$ is permutation-unconditional,
					\end{enumerate}
					and in this case
					\begin{equation*}
						c_s\le c_u=c_a\le 2c_s.
					\end{equation*}
					A basic sequence satisfying one, and thus all of the above conditions, is simply called \emph{unconditional}.
					\item (Spreading basic sequence) A basic sequence is \emph{spreading} if it is equivalent to all its subsequences, i.e., 
					there is $C>0$ such that for all integers $1\le n_1<n_2<\dots<n_k<\dots$, and $(a_n)_{n=1}^\infty\in C_{00}$
					\begin{equation*}
						\label{eq:spreading-basic}
						\frac{1}{C}\Big\|\sum_{n=1}^\infty a_n x_n\Big\| \le \Big\|\sum_{k=1}^\infty a_{n_k} x_{n_k}\Big\| \le C \Big\|\sum_{n=1}^\infty a_n x_n\Big\|.
					\end{equation*}
					It is $1$-spreading if $C=1$. 
					\item (Symmetric basis) A basis is said to be \emph{symmetric} if it is equivalent to any of its permutations. In particular, a symmetric basis is unconditional.
					\item (Subsymmetric basis) A basis that is unconditional and spreading is called a \emph{subsymmetric basis}. Every symmetric basis is subsymmetric.
					\item (Unconditional sums of Banach spaces) If $Y$ is a Banach space with an unconditional basis $(e_n)_{n=1}^\infty$, then for every sequence of Banach spaces $(X_n)_{n=1}^\infty$, the unconditional sum with respect to $Y$, denoted by $(\sum_{n=1}^\infty X_n)_Y$, is the vector space $\Pi_{n=1}^\infty X_n$ equipped with the norm 
					\begin{equation*}
						\norm{(x_n)_{n=1}^\infty} := \Big\|\sum_{n=1}^\infty\norm{x_n}_{X_n} e_n\Big\|_Y.
					\end{equation*}
					If for all $n\ge 1$, $X_n=X$, we simply write $Y(X)$.
				\end{itemize}
				
				\section{Complemented subspaces}
				\label{appendix:complemented-subspaces}
				
				\begin{itemize}
					\item A subspace $Y$ of a Banach space $X$ is \emph{complemented} (in $X$) if $Y$ is the range of a bounded projection. In particular, $Y$ is closed and $X$ is isomorphic to $Y\oplus X/Y$.
					\item If $Y$ is a subspace of $X$, then $Y^*$ is usually not a subspace of $Y^*$ but it is always a quotient of $X^*$. In fact, $Y^*$ can be identified with $X^*/{Y^\perp}$. However, if $Y$ is a complemented subspace of $X$, then $Y^*$ is isomorphic to a complemented subspace of $X^*$.
					\item (J. Lindenstrauss \cite{Lindenstrauss1966}) Every reflexive Banach space $X$ has the \emph{separable complementation property}, i.e., for any separable subspace $Z$ of $X$, there exists a separable $1$-complemented subspace $Y$ of $X$ containing $Z$.
					\item[] A theorem of {\`{E}}del'\v{s}te{\u\i}n and Wojtaszczyk \cite{EW1976} describes the complemented subspaces of direct sums of totally incomparable infinite-dimensional spaces, i.e., spaces such that no infinite-dimensional subspace of one is isomorphic to a subspace of the other.
					\item (Complemented subspaces of direct sums of totally incomparable spaces) If $X$ and $Y$ are two totally incomparable Banach spaces, then any complemented subspace of $X\oplus Y$ is isomorphic to $E\oplus F$, where $E$ and $F$ are complemented subspaces of $X$ and $Y$, respectively.
				\end{itemize}
				
				\section{\texorpdfstring{$\cL_p$-spaces, subspaces and quotients of $L_p$}{Subspaces and quotients of}}
				\label{appendix:subspaces-Lp}
				
				In this appendix, we gather fundamental results that are needed in this book about $\cL_p$-spaces and the subspaces and quotients of $L_p$.
				
				\begin{itemize}
					\item (Lindenstrauss and Pe{\l }czy{\'n}ski \cite{LP1968}) Let $p\in[1,\infty]$. A Banach space $X$ is a \emph{$\cL_p$-space} if there exists a constant $C\ge 1$ such that for every finite-dimensional subspace $E$ of $X$, there exists a finite-dimensional subspace $F$ of $X$ such that $E\subset F$ and $F$ is $C$-isomorphic to $\ell_p^{\text{dim}(F)}$. 
					\item (Lindenstrauss and Pe{\l }czy{\'n}ski \cite{LP1968}) For $p\in[1,\infty]$, any $L_p(\mu)$-space is a $\cL_p$-space and every infinite-dimensional $\cL_p$-space has a complemented subspace isomorphic to $\ell_p$ when $p\in[1,\infty)$.
					\item (Lindenstrauss and Pe{\l }czy{\'n}ski \cite{LP1968}) For $p\in[1,\infty]$, every $\cL_p$-space $X$ is isomorphic to a subspace of an $L_p(\mu)$-space. Moreover, $X$ is isomorphic to a complemented subspace of an $L_p(\mu)$-space if and only if $X$ is complemented in its bidual. Therefore, for $p\in(1,\infty)$, every $\cL_p$-space is isomorphic to a complemented subspace of an $L_p(\mu)$-space. Moreover, if the $\cL_p$-space is separable, then it is isomorphic to a complemented subspace of  $L_p[0,1]$.
					\item (Lindenstrauss and Rosenthal \cite{LindenstraussRosenthal1969}) For $p\in(1,\infty)$, a complemented subspace of an $L_p(\mu)$-space, is either a $\cL_p$-space or isomorphic to a Hilbert space.
					\item (Kadec and Pe{\l }czy{\'n}ski \cite{KadetsPelczynski1961}) For $q\in(2,\infty)$, a subspace of $L_q$ is either isomorphic to $\ell_2$ or contains a subspace isomorphic to $\ell_q$.
					\item (Johnson and Odell \cite{JohnsonOdell1974}) For $q\in(2,\infty)$, a subspace of $L_q$ that does not contain a subspace isomorphic to $\ell_2$ is isomorphic to a subspace of $\ell_q$. 
					%Let $q\in(2,\infty)$ and $X$ be a subspace of $L_q$. Assume that $X$ does not contain any isomorphic copy of $\ell_2$. Then $X$ is isomorphic to a subspace of $\ell_q$.
					\item (Johnson and Odell \cite{JohnsonOdell1974}) Let $p\in(1,\infty)$ and $X$ be a separable $\cL_p$-space. Then, either $X$ contains a subspace isomorphic to $\ell_2$ or $X$ is isomorphic to $\ell_p$. In particular, a complemented subspace of $L_p$ with no subspace isomorphic to $\ell_2$ is isomorphic to $\ell_p$. 
					\item (Johnson \cite{Johnson1976}) Let $q\in(2,\infty)$ and $X$ be a closed subspace of $L_q$ with the approximation property. Assume that $Y$ is a quotient of $X$ with the weak-$q$-Banach-Saks property, then $Y$ is isomorphic to a quotient of $\ell_q$. In particular, a quotient of $L_q$ with the weak-$q$-Banach-Saks property is isomorphic to a quotient of $\ell_q$.
				\end{itemize} 
				
				\section{Approximation properties}
				\label{appendix:approximation}
				The survey article by Casazza in \cite{Casazzahandbook} provides a thorough exposition of approximation properties of Banach spaces. We only list below the definitions and results that are needed in this book.
				\begin{itemize}
					\item A Banach space $X$ is said to have the \emph{approximation property} (AP for short) if for every compact set $K$ in $X$ and every $\vep>0$, there is a finite rank operator $T\colon X \to X$ so that $\norm{Tx - x}<\vep$, for every $x\in K$.
					\item Let $\lambda\in [1,\infty)$. A Banach space $X$ has the
					\emph{$\lambda$-bounded approximation property} ($\lambda$-BAP for short) if for every $\vep>0$ and every compact set $K$ in X, there is a finite rank operator $T$ so that $\norm{T}\le \lambda$ and
					$\norm{Tx - x}<\vep$, for every $x\in K$. 
					\item A Banach space is said to have the \emph{bounded approximation property} (BAP for short) if it has the $\lambda$-BAP, for some $\lambda\ge 1$. 
					\item A Banach space is said to have the \emph{metric approximation property} (MAP for short) if it has the $1$-BAP.
					\item Every complemented subspace of a Banach space with the (bounded) approximation property must
					have the (bounded) approximation property.
					\item (Grothendieck \cite{Grothendieck1953}) A separable dual space with the approximation property has the metric approximation property.
					\item (Grothendieck \cite{Grothendieck1953}) Every reflexive Banach space with the approximation property also has the metric approximation property.
					\item Let $\lambda\in [1,\infty)$. A Banach space $X$ has the \emph{$\lambda$-commuting bounded approximation property} ($\lambda$-CBAP for short) if there is a net $(T_\alpha)_\alpha$ of finite rank operators on $X$ converging strongly to the identity such that $\limsup_\alpha \norm{T_\alpha}\le \lambda$ and for all $\alpha,\beta$, we have
					$T_\alpha T_\beta= T_\beta T_\alpha$. We say that X has the \emph{commuting bounded approximation property} (CBAP for short) if it has the $\lambda$-commuting bounded approximation property for some $\lambda\ge 1$.
					\item (Casazza-Kalton \cite[Proposition 2.1]{CasazzaKalton1990}) A separable Banach space has the $\lambda$-commuting bounded approximation property if and only if there is a sequence $(T_n)_n$ of finite rank operators on $X$ converging strongly to the identity such that $\limsup_n \norm{T_n}\le \lambda$ and for all $m,n\in \bN$, we have
					$T_m T_n= T_{\min\{m,n\}}$.
					\item (Casazza-Kalton \cite[Theorem 2.4]{CasazzaKalton1990}) A separable Banach space with the metric approximation property has the commuting metric approximation property.
					\item (Pe{\l }czy{\'n}ski and Rosenthal \cite{PelczynskiRosenthal1974}) A Banach space $X$ has the \emph{uniform approximation property} (UAP for
					short) if there is a constant $K$ and a function $f\colon \bN \to \bN$ so that for all $n\ge 1$ and all $n$-dimensional subspaces $E\subseteq X$, there is an operator $T\colon X \to X$ with $\mathrm{rank}(T)\le f (n)$ such that $\norm{T}\le K$ and $T_{\restriction E} = I_{\restriction E}$. 
					%If for every $0 < \vep < 1$, we can replace the $K$ above by $1+\vep$, we say that $X$ has the metric UAR And if the operators $T$ can be chosen to be projections, we say that $X$ has the uniform projection approximation property (UPAP for short).
					\item (Pe{\l }czy{\'n}ski and Rosenthal \cite{PelczynskiRosenthal1974}) For all $p\in[1,\infty]$, $L_p(\mu)$ has the uniform approximation property.
					\item(Heinrich \cite[Theorem 9.1 (ii)]{Heinrich1980}) A Banach space $X$ has the $(\lambda+\vep)$-uniform approximation property for every $\vep>0$ if and only if each ultrapower of $X$ has the $\lambda$-uniform approximation property if and only if each ultrapower of $X$ has the $\lambda$-bounded approximation property.
					\item(Heinrich \cite[Theorem 9.4]{Heinrich1980}) If a Banach space $X$ has the uniform approximation property, then so does $L_p(\mu;X)$.
					\item Implications between the various approximation properties. None of these implications is reversible in general.
					\begin{center}
						\begin{tabular}{ccccc}
							& & UAP &  &   \\
							& & $\Downarrow$ &  &   \\
							MAP & $\Rightarrow$ & BAP & $\Rightarrow$ & AP  \\
							&  & $\Uparrow$ &  &  \\
							& & CBAP  &  &\\
						\end{tabular}
					\end{center}
				\end{itemize}
				
				\section{Finite representability and super-reflexivity}
				The notion of (crude) finite representability discussed in this appendix is often attributed to R. C. James. This terminology seems indeed to have been introduced by James in 1972 in \cite{James1972b}. However, this notion was at least introduced already in 1953 by A. Grothendieck in \cite{Grothendieck1953b} under the French terminology ``$X$ a un type (m\'etrique) lin\'eaire inf\'erieur \`a $Y$'' when he conjectured what will become Dvoretzky's theorem. Note that in his 1961 paper \cite{Dvoretzky1961}, A. Dvoretzky uses the terminology ``$X$ has smaller metric type than $Y$" while in 1968, J. Lindenstrauss and A. Pe{\l }czy{\'n}ski in \cite[page 323]{LP1968} use the yet different terminology ``$X$ has finite type not exceeding $Y$" for what is now called ``$X$ is finitely representable in $Y$".
				\begin{itemize}
					\item Given a constant $\lambda\ge 1$, a Banach space $X$ is said to be \emph{$\lambda$-crudely finitely representable into} a Banach space $Y$ if for any finite-dimensional subspace $E$ of $X$, there exists a finite-dimensional subspace $F$ of $Y$ and an isomorphism $T$ from $E$ onto $F$ with $\norm{T}\,\|T^{-1}\|\le \lambda$. 
					\item A Banach space $X$ is said to be \emph{crudely finitely representable into} $Y$ if it is $\lambda$-crudely finitely representable into $Y$ for some $\lambda\ge 1$.
					\item A Banach space $X$ is said to be \emph{finitely representable into} $Y$ if it is $(1+\vep)$-crudely finitely representable into $Y$ for every $\vep>0$.
					\item (Dvoretzky's theorem \cite{Dvoretzky1961}) The separable Hilbert space $\ell_2$ is finitely representable into every infinite-dimensional Banach space.
					\item A Banach space $X$ is \emph{reflexive} if the canonical isometry from $X$ into $X^{**}$ is surjective.
					\item (James \cite{James1972b}) A Banach space is \emph{super-reflexive} if every Banach space that is finitely representable in it is reflexive. It is plain that a super-reflexive Banach space is reflexive.
					\item The bidual a Banach space $X$ is always finitely representable into $X$. However, the principle of local reflexivity gives much more precise information.
				\end{itemize}

				\section{Ultrafilters and ultrapowers}
				\label{appendix:ultraproducts}
				A standard reference on everything we need on ultrafilters is \cite{ComfortNegrepontis1974}.
				\begin{itemize}
					\item A \emph{filter} on a set $I$ is a non-empty collection of subsets of $I$ that is stable under finite intersections and by taking supersets.
					\item The collection of cofinite subsets of $I$ is a filter called the Fr\'echet filter.
					\item An \emph{ultrafilter} on a set $I$ is a maximal filter with respect to inclusion. A \emph{principal ultrafilter} on a set $I$ is an ultrafilter of the form $\cU_{i_0}:=\{S\subseteq I\colon i_0\in S\}$ for some $i_0\in I$. An ultrafilter that is not principal is called a \emph{nonprincipal ultrafilter}. 
					\item Every nonprincipal ultrafilter contains the Fr\'echet filter and is made only of infinite sets.
					\item If $f$ is a map from a set $I$ into a compact Hausdorff topological space $(X,\tau)$, then for every ultrafilter $\cU$ on $I$, $\lim_{i,\cU}f(i)$ exists and is unique, where $\ell=\lim_{i,\cU}f(i)$ if and only if for every $\tau$-neighborhood $V$ of $\ell$, $f^{-1}(V)\in \cU$. 
					\item If $(X_i)_{i\in I}$ is a collection of Banach spaces, we denote by $(\Pi_{i\in I} X_i)_{\ell_\infty}$ the space of all families $x=(x_i)_{i\in I}$ such that $\norm{x}_\infty := \sup_{i\in I}\norm{x_i}_{X_i}<\infty$. Given an ultrafilter $\cU$ on $I$, for $x=(x_i)_{i\in I} \in (\Pi_{i\in I} X_i)_{\ell_\infty}$, we define the seminorm $N_\cU(x) :=\lim_{i,\cU}\norm{x_i}_{X_i}$. The quotient of $(\Pi_{i\in I} X_i)_{\ell_\infty}$ by the closed subspace of all $x$ with $N_\cU(x)=0$ equipped with its natural quotient norm, denoted  $\norm{\cdot}_\cU$, is a Banach space called the \emph{ultraproduct} of the $X_i$ with respect to $\cU$ and is denoted by $(\prod_{i \in I} X_i)^{\cU}$. An equivalence class of $(x_i)_{i\in I}$ in $(\Pi_{i\in I} X_i)^\cU$ is usually denoted by $(x_i)^\cU$ and we have $\norm{(x_i)^\cU}_\cU=\lim_{i,\cU}\norm{x_i}_{X_i}$. If, for all $i\in I$, $X_i=X$, the ultraproduct defined above is called the \emph{ultrapower} of $X$ with respect to $\cU$ and we simply write $X^{\cU}$.
					\item The map $x\in X\mapsto (x)^\cU\in X^\cU$ is an isometric embedding of $X$ into its ultrapower $X^\cU$. 
					\item When $X$ is reflexive, the map $(x_i)^\cU \in X^{\cU} \mapsto w-\lim_{i,\cU} x_i \in X$ is well defined and is easily seen to be a projection of norm $1$. Therefore, a reflexive Banach space is always complemented in each of its ultrapowers.
					\item Any ultrapower of a Banach space $X$ is finitely representable into $X$.
					\item[] Super-reflexivity can be conveniently characterized using ultrapowers.
					\item A Banach space is super-reflexive if and only if each of its ultrapowers is reflexive.
					\item[] There is a local reflexivity principle for ultrapowers.
					\item (Heinrich \cite[Theorem 7.3]{Heinrich1980} or \cite[Proposition F.6 (v)]{BenyaminiLindenstrauss2000}) Let $X$ be a Banach space, $\cU$ an ultrafilter on $I$, $E$ be a finite-dimensional subspace of $X^\cU$, and $F$ a finite-dimensional subspace of $(X^\cU)^*$, then for every $\vep>0$  there is a subspace $G$ of $(X^*)^\cU$ and a $(1+\vep)$-isomorphism $T$ from $G$ onto $F$ such that for every $u\in E$ and $v\in G$ we have $\langle v, u\rangle = \langle Tv, u\rangle$.
					\item[] The next result is an immediate consequence of the principle of local reflexivity for ultrapowers. 
					\item If $X$ is a Banach space and $\cU$ an ultrafilter on a set $I$, then for every $\vep>0$, every finite-dimensional $\lambda$-complemented subspace of $X^\cU$ is $(1+\vep)$-isomorphic to a $(\lambda+\vep)$-complemented subspace of $X$. 
				\end{itemize}

				\section{Uniformly smooth and uniformly convex norms}
				\label{appendix:us-uc}
				
				We recall here the classical notions of uniform smoothness and uniform convexity of a Banach space, as well as classical renorming results related to these notions.
				
				\begin{itemize}
					\item (Clarkson \cite{Clarkson1936}) A Banach space $(X,\norm{\cdot})$ is \emph{uniformly convex} (UC in short) if for every $\vep>0$ there exists $\delta>0$ such that $\norm{\frac{x_1+ x_2}{2}}\le 1-\delta$ whenever $\norm{x_1}\le 1$, $\norm{x_2}\le 1$, and $\norm{x_1-x_2}\ge \delta$.
					\item (Milman \cite{Milman1938}, Kakutani \cite{Kakutani1939}, Pettis \cite{Pettis1939}) Every uniformly convex Banach space is reflexive.
					\item (Day \cite{Day1944}) The modulus of \emph{uniform convexity} of $(X,\norm{\cdot})$ is defined for $t\in [0,2]$ by 
					\begin{equation*}
						\delta_{\norm{\cdot}}(t) := \inf\Big\{1-\frac{\norm{x_1+x_2}}{2}\colon  x_1,x_2\in S_X,\ \norm{x_1-x_2}\ge t\Big\}.
					\end{equation*}
					\item It is immediate that a Banach space $(X,\norm{\cdot})$ is uniformly convex if and only if $\delta_{\norm{\cdot}}(t)>0$ for all $t\in(0,2]$.
					\item If a Banach space is uniformly convex, then it is asymptotically uniformly convex. More precisely, $\delta_{\norm{\cdot}}(t)\le \bar{\delta}_{\norm{\cdot}}(t)$ for all $t\in(0,1)$.
					%\item (Shmulyan \cite{Shmulyan1940}) Every uniformly smooth Banach space is reflexive.
					\item[] It took more time for a dual notion to uniform convexity, first investigated by Day in \cite{Day1944}, to reach its final and modern form. 
					\item A Banach space $(X,\norm{\cdot})$ is \emph{uniformly smooth} (US in short) if for every $\vep>0$ there exists $\delta>0$ such that $\norm{x_1+x_2} + \norm{x_1-x_2} \le 2 +\vep \norm{x_2}$, whenever $\norm{x_1}\le 1$ and $\norm{x_2}\le \delta$.
					\item (Lindenstrauss \cite{Lindenstrauss1963}) The modulus of \emph{uniform smoothness} of $(X,\norm{\cdot})$ is defined for $t\ge 0$ by 
					\begin{equation*}
						\rho_{\norm{\cdot}}(t) := \sup\Big\{\frac{\norm{x_1+tx_2} + \norm{x_1-tx_2}}{2} \colon x_1\in S_X,\ x_2 \in S_X\Big\}.
					\end{equation*}
					\item A Banach space $(X,\norm{\cdot})$ is uniformly smooth if and only if $\lim_{t\to 0}\frac{\rho_{\norm{\cdot}}(t)}{t}=0$.
					\item If a Banach space $(X,\norm{\cdot})$ is uniformly smooth, then it is asymptotically uniformly smooth. More precisely, $\bar{\rho}_{\norm{\cdot}}(t)\le 2\rho_{\norm{\cdot}}(t)$ for all $t\in(0,1)$.
					\item[] The duality UC-US was first elucidated by Day \cite{Day1944} and a fundamental duality formula was later discovered by Lindenstrauss \cite{Lindenstrauss1963}.
					\item (Uniformly convex-uniformly smooth duality) $X$ is uniformly convex/uniformly smooth if and only if $X^*$ is uniformly smooth/uniformly convex.
					\item For every Banach space $X$ and $t>0$, $\rho_{X^*}(t) = \sup_{\vep\in(0,2]} \{ \frac{\vep t}{2} -\delta_X(\vep)\}$.
					\item (Asplund \cite{Asplund1967}) If a Banach space admits an equivalent norm that is uniformly convex and an equivalent norm that is uniformly smooth, then it admits an equivalent norm that is simultaneously uniformly convex and uniformly smooth. 
					\item (Enflo \cite{Enflo1972}) A Banach space is super-reflexive if and only if it admits an equivalent norm that is uniformly convex (and by duality if and only if it admits an equivalent norm that is uniformly smooth).
					\item[] The quantitative renorming theory of super-reflexive Banach spaces is due to Pisier in \cite{Pisier1975} using martingale techniques and can be found in \cite{Beauzamy1985}, \cite{DGZ1993}, or \cite{Pisier_martingale}. 
					\item If $p\in(1,2]$, a Banach space $(X,\norm{\cdot})$ is said to be \emph{uniformly smooth with power-type $p$} if there exists a constant $C>0$ such that $\rho_{\norm{\cdot}}(t)\le Ct^p$, for all $t\in [0,\infty)$.
					If $p\in(1,\infty)$, $L_p(\Omega,\Sigma,\mu)$ is uniformly smooth with power-type $\min\{p,2\}$.
					\item If $q\in[2,\infty)$, a Banach space $(X,\norm{\cdot})$ is said to be \emph{uniformly convex with power-type $q$} if there exists a constant $C>0$ such that $\delta_{\norm{\cdot}}(t)\ge Ct^q$, for all $t\in [0,2]$.
					If $p\in(1,\infty)$, $L_p(\Omega,\Sigma,\mu)$ is uniformly convex with power-type $\max\{p,2\}$.
					\item If $p\in(1,2]$, a Banach space $(X,\norm{\cdot})$ is said to be \emph{$p$-uniformly smooth} if there exists a constant $S\ge 1$ such that for all $x_1,x_2\in X$,
					\begin{equation*}
						\frac{\norm{x_1+x_2}^p + \norm{x_1-x_2}^p}{2} \le \norm{x_1}^p + S^p\norm{x_2}^p.
					\end{equation*}
					\item If $q\in[2,\infty)$, a Banach space $(X,\norm{\cdot})$ is said to be \emph{$q$-uniformly convex} if there exists a constant $K\ge 1$ such that for all $x_1,x_2\in X$,
					\begin{equation*}
						\norm{x_1}^q + \frac{1}{K^q}\norm{x_2}^q \le \frac{\norm{x_1+x_2}^q + \norm{x_1-x_2}^q}{2}.
					\end{equation*}
					\item[] The equivalence between the homogeneous inequalities and the power-type estimates above can be found in \cite{BCL1994}.  
					\item If $p\in(1,2]$, a Banach space is uniformly smooth of power-type $p$ is and only if it is $p$-uniformly smooth.
					\item If $q\in[2,\infty)$, a Banach space is uniformly convex of power-type $q$ is and only if it is $q$-uniformly smooth.
					\item (Pisier \cite{Pisier1975}) A Banach space admits an equivalent uniformly convex norm (resp. uniformly smooth) if and only if it admits an equivalent norm that is $q$-uniformly convex for some $q\in[2,\infty)$ (resp. $p$-uniformly smooth for some $p\in(1,2]$).
				\end{itemize}

				\section{Classical Banach spaces from the modern era}
				\label{sec:special-spaces}
				
				\subsection{James spaces and relatives}
				\label{sec:James-space}
				
				We recall the definition and some basic properties of a space first introduced and studied by James in \cite{James1950}. We refer the reader to \cite{FetterGamboa} or \cite[Section 3.4]{AlbiacKalton2016} and references therein for details. 
				\begin{itemize}
					\item The James space $\James$ is the real Banach space of all sequences $x:=(x_n)_{n\in \bN}$ of real numbers with finite square variation and satisfying $\lim_{n \to \infty} x_n =0$. The space $\James$ is endowed with the norm
					\begin{equation*}
						\norm{x}_{\James} := \sup  \Big \{  \Big (\sum_{i=1}^{k-1} \abs{x_{p_{i+1}} - x_{p_i}}^2 \Big)^{1/2}     \; \colon \; 1 \le p_1 < p_2 < \dots < p_{k} \Big \}.
					\end{equation*}
					\item $\James$ is isomorphic to $\James^{**}$ but it is not reflexive since the canonical identification of $\James$ in its bidual is of codimension one. In fact, $\James^{**}$ can be seen as the space of all sequences of real numbers with finite square variation, which is $\James \oplus \bR \bf{1}$, where $\bf{1}$ is the constant sequence equal to $1$. In \cite{James1951}, James slightly modified the norm to ensure that there is an equivalent norm on $\James$ with which $\James$ is isometrically isomorphic to $\James^{**}$.
					\item[] A description of the complemented subspaces of $\James$ is due to P. Cassaza \cite{Cassaza1977}
					\item (Complemented subspaces of $\James$) A complemented subspace of $\James$ is either reflexive or isomorphic to $\James$.
					\item[] The next results show that in some asymptotic sense $\James$ behaves like Hilbert space. 
					\item James space $\James$ is hereditarily $\ell_2$, i.e., every infinite-dimensional subspace of $\James$ contains a further subspace that is isomorphic to $\ell_2$. In particular, $\James$ is an example of a non-reflexive Banach space that does not contain $\co$ nor $\ell_1$ (this could of course also be deduced from the separability of $\James^{**})$.
					\item (A. Andrew \cite{Andrew1981}, see also \cite{BeauzamyLapreste1984}) Every spreading model of a weakly null normalized sequence in $\James$ is isomorphic to $\ell_2$. 
					\item (G. Lancien \cite{LancienBesac} see also \cite{Netillard2018} for a more accessible reference) The norm $\norm{\cdot}_{\James}$ on $\James$ is $2$-asymptotically uniformly convex and $\James$ admits an equivalent norm which is $2$-asymptotically uniformly smooth. 
					\item[] James' construction is, of course, quite flexible and has been vastly generalized to provide Banach spaces with a wide range of behaviors. For instance, if instead of the quadratic variation we take the natural $\ell_p$-version of it we get a space that is isomorphic to $\ell_1$ when $p=1$, and a space isomorphic to $\co$ when $p=\infty$. In the reflexive range $p\in(1,\infty)$, the space obtained is usually denoted $\James_p$ (or $v_p^0$ as in \cite{Pisier1988}), and hence $\James_2=\James$. 
					\item For $p\in(1,\infty)$, the space $\James_p$ is a non-reflexive Banach space that is of codimension one in its bidual, hereditarily $\ell_p$, whose spreading models generated by weakly null normalized sequences are isomorphic to $\ell_p$. Moreover, $\James_p$ is $p$-asymptotically uniformly convex and admits an equivalent norm which is $p$-asymptotically uniformly smooth.
					\item (James tree space \cite{James1974}) Let $\bin_\infty$ be the binary tree of infinite height. The James tree space $\mathrm{JT}$ is the completion of the space of finitely supported functions $x \colon \bin_\infty \to \bR$ endowed with the norm
					\begin{equation*}
						\norm{x}_{\mathrm{JT}} := \sup  \Big \{  \Big (\sum_{i=1}^{k} \big|\sum_{t\in S_i}x(t)\big|^2 \Big)^{1/2}     \; \colon \; S_1, \dots, S_k \textrm{ are disjoint segment in } \bin_\infty\Big \},
					\end{equation*}
					where a segment in $\bin_\infty$ is a finite linearly ordered subset of consecutive elements in $\bin_\infty$ equipped with its natural tree ordering.
					\item James tree space $\mathrm{JT}$ is an example of a separable Banach space that is a dual space, has a nonseparable dual, but does not contain an isomorphic copy of $\ell_1$. Moreover, $\mathrm{JT}^{**}=\mathrm{JT}\oplus \ell_2(I)$ for some uncountable set $I$, and every infinite-dimensional subspace of $\mathrm{JT}$ contains a subspace isomorphic to $\ell_2$.
				\end{itemize} 
				
				\subsection{Tsirelson spaces and relatives}
				
				The original Tsirelson space was constructed by Tsirelson in \cite{Tsirelson1974} to answer negatively a longstanding open problem in Banach space theory which asked whether every infinite-dimensional Banach space contains an isomorphic copy of $\ell_p$, for some $p\in[1,\infty)$, or $\co$. Figiel and Johnson gave in \cite{FigielJohnson1974} an alternative approach to construct Tsirelson space. This is the one that is presented in this appendix, and we refer to the book of Casazza and Shura \cite{CasazzaShura1989} for a detailed study of Tsirelson spaces (see also \cite{BenyaminiLindenstrauss2000}, \cite{BeauzamyLapreste1984}, \cite{AlbiacKalton2016}).
				\begin{itemize}
					\item (Figiel-Johnson description of Tsirelson spaces $\Tsi$ and $\Tsi^*$) Let $(e_n)_{n=1}^\infty$ be the canonical basis of $c_{00}$.  For $x := \sum_{n=1}^\infty a_n e_n\in c_{00}$ and $E$ a finite subset of $\bN$ we write $P_Ex := \sum_{n\in E} a_n e_n$. Now, consider the following sequence of norms on $\coo$:
					\begin{itemize}
						\item $\norm{x}_0:=\norm{x}_\infty$
						\item $\norm{x}_{n}:=\max\Big\{\norm{x}_{n-1}, \frac12 \sup\sum_{j=1}^m \norm{P_{E_j}x}_{n-1}\Big\}$, for all $n \ge 1$, where the supremum is taken over all $m\in\bN$ and all admissible sequences $(E_j)_{j=1}^m$, i.e., sequences of finite subsets of $\bN$ such that $m\le \min(E_1)\le \max(E_1)<\min(E_2)\le \max(E_2)< \dots <\min(E_j)\le \max(E_j)<\dots <\min(E_m)$.
					\end{itemize}
					Since $(\norm{x}_n)_{n=1}^\infty$ is an non-decreasing sequence of positive numbers upper bounded by $\norm{x}_1$, $\lim_{n\to \infty}\norm{x}_n$ exists and we let $\norm{x}_{\Tsi}:= \lim_{n\to \infty}\norm{x}_n$. It can then be shown that $\|\cdot\|_{\Tsi}$ is the unique norm on $c_{00}$ such that for every $x\in c_{00}$
					\begin{equation}
						\label{eq:Tnorm}
						\norm{x}_{\Tsi} = \max\Big\{\norm{x}_\infty, \frac12 \sup\sum_{j=1}^n \norm{P_{E_j}x}_{\Tsi}\Big\},
					\end{equation}
					where the supremum is taken over all $n\in\bN$ and all admissible sequences $(E_j)_{j=1}^n$. The space $\Tsi$ is the completion of $c_{00}$ with this norm.
					\item The unit vector basis of $\Tsi$ is a 1-unconditional basis and since $\Tsi$ does not contain any isomorphic copy of $\co$ nor $\ell_1$, it is a reflexive space.
					\item The dual $\Tsi^*$ of $\Tsi$ is the original Tsirelson space constructed by Tsirelson in \cite{Tsirelson1974}.
					\item $\Tsi$ and $\Tsi^*$ do not contain isomorphic copies of $\ell_p$, for $p\in[1,\infty)$, $\co$, nor any spreading basic sequence or superreflexive subspace.
					\item (Casazza-Johnson-Tzafriri \cite{CJT84}) $\Tsi^*$ is minimal, i.e. every infinite-dimensional subspace of $\Tsi^*$ contains a further subspace that is isomorphic to $\Tsi^*$.
					\item (Casazza-Odell \cite{CasazzaOdell}) $\Tsi$ contains no minimal subspace at all.
					\item (Convexification of Tsirelson space) For $1<p<\infty$ we denote by $\Tsi_p$ the $p$-convexification of $\Tsi$. Since $\Tsi$ is a sequence space with a $1$-unconditional basis, $\Tsi_p$ can be simply described as the vector space $\{ x\in \Tsi \colon \norm{\abs{x}^p}_{\Tsi}<\infty\}$ equipped with the norm $\norm{x}_{\Tsi_p}:= \norm{\abs{x}^p}_{\Tsi}^{1/p}$, where $\abs{x}^p:=(\abs{x_n}^p)_{n=1}^\infty$.
					\item For every $p\in(1,\infty)$, $\Tsi_p$ is a super-reflexive Banach space that does not contain isomorphic copies of $\ell_p$, for $p\in(1,\infty)$, nor any spreading basic sequence.
					\item (Schlumprecht's space \cite{Schlumprecht91}) Schlumprecht's space $\mathrm{S}$ is the completion of $\coo$ with respect to the norm $\norm{\cdot}_\mathrm{S}$ defined for every $x\in c_{00}$ by
					\begin{equation}
						\label{eq:Snorm}
						\norm{x}_{\mathrm{S}} = \max\Big\{\norm{x}_\infty,  \sup\frac{1}{\log_2(n+1)}\sum_{j=1}^n\norm{P_{E_j}x}_{\mathrm{S}}\Big\},
					\end{equation}
					where the supremum is taken over all $n\in\bN$ and all sequences $(E_j)_{j=1}^n$ of finite subsets of $\bN$ such that $\min(E_1)\le \max(E_1)<\min(E_2)\le \max(E_2)< \dots <\min(E_j)\le \max(E_j)<\dots <\min(E_m)$.
					The difference with the admissibility condition in Tsirelson's space is that we do not require the condition $n\le \min{E_1}$. 
					Schlumprecht's space $\mathrm{S}$ was the first arbitrarily distortable Banach space.
				\end{itemize}
				
				\section{A few fundamental probabilistic inequalities and notions}
				\label{appendix:proba-Banach-spaces}
				In this appendix, we recall a few probabilistic inequalities that are instrumental in the study of the geometry of Banach spaces and their linear structure. For more on the application of probabilistic techniques to Banach space theory, we refer the readers to \cite{Schwartz}, \cite{TalagrandLedoux}, \cite{QueffelecLi_vol1}, or \cite{HNVW}. In the sequel, $(\eps_n)_{n=1}^\infty$ will be a sequence of independent identically distributed Bernoulli random variables on a probability space $(\Omega,\Sigma,\bP)$, i.e., for all $n\in \bN$,  $\bP(\vep_n=1)=\bP(\vep_n=-1)=\frac12$. These are also called Rademacher random variables by analysts.
				\begin{itemize}
					\item (Khintchine's inequalities) Let $(\eps_n)_{n=1}^\infty$ be a sequence of independent identically distributed Bernoulli random variables on a probability space $(\Omega,\Sigma,\bP)$. Then, for every $p\in(0,\infty)$ there are constants $A_p, B_p>0$ such that for every scalars $a_1,\dots, a_n$, 
					\begin{equation*}
						A_p \Big\|\sum_{i=1}^n a_i \vep_i\Big\|_{L_2(\Omega,\bP)} \le \Big\|\sum_{i=1}^n a_i \vep_i\Big\|_{L_p(\Omega,\bP)} \le B_p \Big\|\sum_{i=1}^n a_i \vep_i\Big\|_{L_2(\Omega,\bP)}.
					\end{equation*}
					\item In Banach space terms, since two i.i.d. Bernoulli random variables are orthogonal and thus $\|\sum_{i=1}^n a_i \vep_i\|_{L_2(\Omega,\bP)}=(\sum_{i=1}^n\abs{a_i}^2)^{1/2}$, Khintchine's inequalities say that a sequence of i.i.d. Bernoulli random variables is equivalent to the canonical basis of $\ell_2$. Consequently, $L_p[0,1]$ contains an isomorphic copy of $\ell_2$ for every $p\in[1,\infty)$, and this copy is complemented when $p\in(1,\infty)$.
					\item ($p$-stable random variables) For every $p\in(0,2]$, there exists a random variable $X$, called a (symmetric) $p$-stable random variable such that $\bE(e^{itX})=e^{-\abs{t}^p}$, for every $t\in \bR$. Of course, a $2$-stable random variable is a Gaussian random variable. Using these random variables, it was shown that $L_2[0,1]$ (linearly) isometrically embeds into $L_p[0,1]$ for every $p\in[1,\infty)$ and that $L_q[0,1]$ (linearly) isometrically embeds $L_p[0,1]$ for $1\le p<q\le 2$ (see for instance \cite{AlbiacKalton2016}).
					\item[] In \cite{Kahane1964}, Kahane extended the scalar-valued Khintchine's inequalities to the vector-valued case.
					\item (Kahane's inequalities) Let $(\eps_n)_{n=1}^\infty$ be a sequence of i.i.d. Bernoulli random variables on a probability space $(\Omega,\Sigma,\mu)$.  Then, for every $p\in(0,\infty)$ there are constants $A_p, B_p>0$ such that for all Banach space $X$ and vectors $x_1,\dots, x_n$ in $X$, 
					\begin{equation*}
						A_p \Big\|\sum_{i=1}^n \vep_i x_i\Big\|_{L_2(\Omega,\bP;X)} \le \Big\|\sum_{i=1}^n \vep_i x_i\Big\|_{L_p(\Omega,\bP;X)} \le B_p \Big\|\sum_{i=1}^n \vep_i x_i\Big\|_{L_2(\Omega,\bP;X)}.
					\end{equation*}
					%Then, for every $p,q\in[1,\infty)$ there is a constant $C_{p,q}>0$ such that for all Banach space $X$ and vectors $x_1,\dots, x_n$ in $X$, 
					% \begin{equation*}
					%    \Big\|\sum_{i=1}^n \vep_i x_i\Big\|_{L_p(\Omega,\bP; X)} \le C_{p,q} \Big\|\sum_{i=1}^n \vep_i x_i\Big\|_{L_q(\Omega,\bP; X)}.
					%\end{equation*}
					%\item[ We just recall this fundamental unconditionality phenomenon in vector valued $L_p$ spaces. We refer the reader, for instance, to Proposition 3.2.10 in \cite{HNVW} for the proof.]
					For the history and a discussion of the notions of type and cotype, we refer the reader to \cite{Maurey}. 
					\item (Type and cotype and linear embeddability obstructions) A Banach space $X$ has \emph{(Rademacher) type p} if there is a constant $T>0$ such that for all $x_1,\dots,x_n\in X$, 
					\begin{equation*}
						\Big\|\sum_{i=1}^n \vep_i x_i\Big\|_{L_p(\Omega,\bP;X)}\le T \Big(\sum_{i=1}^n \norm{x_i}^p\Big)^{\frac1p},
					\end{equation*}
					and $X$ has \emph{(Rademacher) cotype q} if there is a constant $C>0$ such that for all $x_1,\dots,x_n\in X$, 
					\begin{equation*}
						\Big\|\sum_{i=1}^n \vep_i x_i\Big\|_{L_q(\Omega,\bP;X)}\ge \frac{1}{C}\Big(\sum_{i=1}^n \norm{x_i}^q\Big)^{\frac1q}.
					\end{equation*}
					The basic properties of these two isomorphic invariants follow from Dvoretzky's theorem and Kahane's inequalities. In particular, every Banach space has type $p\in[1,2]$ and cotype $q\in[2,\infty]$, and type $p_1$ implies type $p_2$ whenever $p_1\ge p_2$, while cotype $q_1$ implies cotype $q_2$ whenever $q_1\le q_2$. Since $L_p(\Omega,\Sigma,\mu)$ has type $\min\{p,2\}$ and cotype $\max\{p,2\}$ (and not better), it follows that $L_p(\Omega,\Sigma,\mu)$ does not isomorphically embed into $L_q(\Omega,\Sigma,\mu)$ whenever $p>2$ and $p>q$, or $p<2$ and $p<q$. 
					\item (Contraction principle and unconditionality) Let $X$ be a real Banach space and $(\Omega,\Sigma,\bP)$ a probability space, $(X_n)_{n=1}^\infty$ sequence of random variables in $L_p(\Omega,\bP;X)$. If the $X_n$'s are independent and \emph{symmetric} (i.e. $X_n$ and $-X_n$ are identically distributed), then $(X_n)_{n=1}^\infty$ is a $1$-unconditional basic sequence. 
					If the $X_n$'s are independent and \emph{centered}, then $(X_n)_{n=1}^\infty$ is a $2$-unconditional basic sequence.
					In particular, for every Banach space $X$ and sequence of i.i.d. Bernoulli random variables $(\eps_n)_{n=1}^\infty$, for all finitely supported sequences $(a_n)_{n=1}^\infty$ in $\Rdb$ and $(x_n)_{n=1}^\infty$ in $X$ and for all $p\in[1,\infty]$ we have
					\begin{equation}
						\Big\|\sum_{n=1}^\infty a_n\eps_n x_n\Big\|_{L_p(\Omega, \bP;X)} \le \max_{n \in \N} |a_n|\,\Big\|\sum_{n=1}^\infty \eps_n x_n\Big\|_{L_p(\Omega,\bP;X)}.
					\end{equation}
					This last fact is commonly referred to as Kahane's contraction principle.
				\end{itemize}

				\chapter{Rudiments of Bochner integration}
				\label{appendix:Bochner}	
				We summarize in this appendix the basics of Bochner integration that are needed in the text. Detailed proofs and thorough treatments of Banach space-valued integration can be found in classic textbooks such as \cite{DunfordSchwartz1958} or \cite{DiestelUhl1977} and in the more recent monograph \cite{HNVW}. Throughout this section, $(\Omega,\Sigma,\mu)$ is a complete $\sigma$-finite measure space and $(X,\norm{\cdot})$ is a real Banach space. When necessary, we will write $\Omega$ as the increasing union $\Omega=\bigcup_{n=1}^\infty\Omega_n$ with $\mu(\Omega_n)$ finite for all $n\in \N$. 
				
				\begin{defi} 
					Let $f$ be a map from $\Omega$ to $X$. We say that $f$ is a \emph{simple measurable function} if there exist $x_1,\dots,x_n \in X$ and $A_1,\dots,A_n \in \Sigma$ such that $f=\sum_{i=1}^n x_i\car_{A_i}$, where $\car_A$ is the indicator function of the set $A$. Then, we say that $f$ is a \emph{measurable} function, or \emph{$X$-measurable}, if there exists a sequence $(f_n)_{n=1}^\infty$ of simple measurable functions that converges almost everywhere to $f$ with respect to the norm of $X$.
				\end{defi}
				
				When $X$ is separable, this notion of measurability coincides in a certain sense with the natural notion of measurability in terms of inverse images of open sets.
				
				\begin{prop}
					\label{prop:measurable} 
					Let $f\colon \Omega \to X$. Then, $f$ is measurable if and only if there exists $g\colon \Omega \to X$ such that
					\begin{enumerate}[(i)]
						\item $f=g$ almost everywhere.
						\item For every open subset $U$ of $X$, $g^{-1}(U)\in \Sigma$.
						\item $g(\Omega)$ is a separable subset of $X$.
					\end{enumerate}
				\end{prop}
				
				\begin{proof} 
					Assume first that $f$ is measurable. Then, there exist $A\in \Sigma$ such that $\mu(A)=0$ and a sequence $(f_n)$ of simple measurable functions converging pointwise to $f$ on $\Omega \setminus A$. Let us set $g :=f$ on $\Omega \setminus A$ and $g :=0$ on $A$. It is clear that $g$ satisfies properties $(i)$ to $(iii)$ above. 
					
					Let us now assume that $f,g\colon \Omega \to X$ satisfy conditions $(i)$ to $(iii)$.\\
					Suppose first that $\mu(\Omega)<\infty$ and fix $\eps>0$. Let $(x_n)_{n=1}^\infty$ be a dense sequence in $g(\Omega)$. By assumption, the set $A_n :=\{\omega\in \Omega,\ \norm{g(w)-x_n}\le \vep\}$ is in $\Sigma$. Set $B_1:=A_1$ and $B_n := A_n\setminus \bigcup_{k<n}A_k$ for $n>1$ and note that $\bigcup_{n=1}^\infty B_n=\Omega$. Since $\mu(\Omega)$ is finite and the $B_n$ are pairwise disjoint, we can pick $N\in \bN$ so that $\sum_{n>N}\mu(B_n)\le \vep$. If we set $E:=\bigcup_{n>N} B_n$ and $h:=\sum_{n\le N}x_n\car_{B_n}$, we have that $h$ is a simple measurable function, $\mu(E)\le \vep$ and $\norm{g-h}\le \vep$ on $\Omega \setminus E$.\\
					Assume now that $\Omega$ is the increasing union $\Omega=\bigcup_{n=1}^\infty\Omega_n$ with $\mu(\Omega_n)<\infty$. Applying the previous reasoning, we pick, for each $n$ in $\bN$, a measurable subset $E_n$ of $\Omega_n$ such that $\mu(E_n)\le 2^{-n}$ and a simple function $h_n$ such that $\norm{g-h_n}\le 2^{-n}$ on $\Omega_n\setminus E_n$. We now set $F_k:=\bigcup_{n\ge k}E_n$ and $F:=\bigcap_{k\in \bN} F_k$. It is easily checked that $F\in \Sigma$, $\mu(F)=0$ and for all $\omega \in \Omega\setminus F$, $\lim_{n\to \infty}\norm{h_n(w)-g(w)}=0$. This concludes our proof.
				\end{proof}
				
				The following characterization, due to Pettis \cite{Pettis1938}, asserts that a function is measurable if and only if it is essentially separably valued and weakly measurable.
				
				\begin{prop}
					\label{prop:Pettis} 
					Let $f\colon \Omega \to X$. Then, $f$ is measurable if and only if the two following assertions hold.
					\begin{enumerate}[(i)]
						\item For any $x^*\in X^*$, $\langle f,x^*\rangle$ is an $\bR$-valued measurable function. 
						\item There exists $A \in \Sigma$ such that $\mu(A)=0$ and $f(\Omega \setminus A)$ is separable.
					\end{enumerate}
				\end{prop}
				A function $f$ satisfying $(i)$ is said to be \emph{weakly measurable}, while a function $f$ satisfying $(ii)$ is said to be \emph{essentially separably valued}.
				
				\begin{proof} 
					It is clear from the definition that a measurable function satisfies assertions $(i)$ and $(ii)$. So, assume $f\colon \Omega \to X$ satisfies $(i)$ and $(ii)$. In particular, we may as well assume that $X$ is separable and thus we can find a sequence $(x_n^*)_{n=1}^\infty$ in $S_{X^*}$ such that
					$$\forall x\in X,\ \ \norm{x}=\sup_{n\in \N}x^*_n(x).$$
					Such a sequence is called a \emph{norming sequence} and its existence is an easy application of the Hahn-Banach theorem. Since $X$ is separable, we also have that any open set in $X$ can be written as a countable union of closed balls. In view of the proof of the previous proposition, it is therefore enough to show that for every $x\in X$ and every $\eps>0$ the set $\{\omega \in \Omega\colon \norm{f(w)-x}\le \vep\}$ is measurable. But this fact follows immediately from the assumption $(i)$ and the fact that $\norm{f(w)-x}\le \eps$ if and only if for all $n\in \bN$, $\langle f(w),x_n^*\rangle \le \langle x,x_n^*\rangle + \vep$.
				\end{proof}
				
				We now turn our attention to the integrability of vector-valued maps.
				\begin{defi}
					A map $f\colon \Omega \to X$ is said to be a \emph{simple integrable function} if there exist $x_1,\dots,x_n$ in $X$ and $A_1,\dots,A_n$ in $\Sigma$ such that $\mu(A_i)<\infty$ for all $i\le n$ and $f=\sum_{i=1}^n x_i\car_{A_i}$.
				\end{defi}
				For a simple integrable function $f=\sum_{i=1}^n x_i\car_{A_i}$, we let 
				\begin{equation*}
					\int_\Omega f\,d\mu := \sum_{i=1}^n \mu(A_i)x_i.
				\end{equation*}
				It is clear that $\int_\Omega f\,d\mu$ is well-defined since it does not depend on the representation of the simple integrable function $f$. Moreover, for every simple integrable function we immediately see that $\|\int_\Omega f\,d\mu\|\le \int_\Omega \norm{f}\,d\mu$.
				
				\begin{defi}
					A measurable function $f\colon \Omega \to X$ is called \emph{Bochner integrable} if there exists a sequence $(f_n)_{n=1}^\infty$ of simple integrable functions such that $(f_n)_{n=1}^\infty$ converges to $f$ almost everywhere and $\lim_{n\to \infty}\int_\Omega \norm{f-f_n}\,d\mu=0$. In that case, we say that $f\in L_1(\mu,X)$ and the sequence $(\int_\Omega f_n\,d\mu)_{n=1}^\infty$ is easily seen to converge in $(X,\norm{\cdot})$ to a limit that does not depend on the choice of $(f_n)_{n=1}^\infty$ and we denote this limit by $\int_\Omega f\,d\mu$. 
				\end{defi}
				If $f$ is a Bochner integrable function then $\|\int_\Omega f\,d\mu\|\le \int_\Omega \|f\|\,d\mu$. The linearity of the Bochner integral is easy to deduce from its definition. Let $X$ and $Y$ be Banach spaces. Assume that $f\colon \Omega \to X$ is Bochner integrable and let $T \in B(X,Y)$. Then, $T \circ f \colon \Omega \to Y$ is Bochner integrable and $\int_\Omega (T \circ f) \,d\mu=T(\int_\Omega f\,d\mu)$. We can now state the following fundamental and simple characterization of Bochner integrable functions.
				\begin{theo}
					\label{thm:Bochner-characterization}
					Let $f\colon \Omega \to X$ be a measurable function. Then, $f$ is Bochner integrable if and only if  $\int_\Omega \norm{f}\,d\mu<\infty.$
				\end{theo}
				
				\begin{proof} 
					The ``only if'' part is straightforward, so assume that $\int_\Omega \norm{f}\,d\mu<\infty$. It follows from the proof of Proposition \ref{prop:measurable} that we can find a sequence of simple integrable functions $(h_n)_{n=1}^\infty$ (remember that $h_n$ was supported in $\Omega_n$) such that $h_n$ tends to $f$ almost everywhere. Set now $A_n=\{\omega \in \Omega,\ \norm{h_n(\omega)}\le 2\norm{f(\omega)}\}$ and $g_n=h_n\car_{A_n}$. We still have that $g_n$ is a simple integrable function and that $g_n$ tends to $f$ almost everywhere. Moreover, $\norm{g_n-f}\le 3\norm{f}$, so it follows from the dominated convergence theorem and the integrability of $\norm{f}$ that $\lim_{n\to \infty}\int_\Omega\norm{g_n-f}\,d\mu=0$. This shows that $f\in L^1(\mu,X)$.
				\end{proof}
				
				The vector-valued version of the Lebesgue dominated convergence follows immediately from Theorem \ref{thm:Bochner-characterization}.
				\begin{theo} 
					Let $(f_n)_{n=1}^\infty$ be a sequence of Bochner integrable functions from $\Omega$ to $X$ and let $f\colon \Omega \to X$. Assume that
					\begin{enumerate}[(a)]
						\item $f_n$ tends to $f$ almost everywhere 
						\item[] and
						\item there exists $g\in L_1(\mu,\bR)$ such that for all $n\in \bN$, $\norm{f_n}\le g$ almost everywhere.
					\end{enumerate}
					Then, $f\in L_1(\mu,X)$ and $\lim_{n\to \infty}\int_\Omega \norm{f-f_n}\,d\mu=0$. In particular, 
					\begin{equation*}
						\lim_{n\to\infty}\int_\Omega f_n\,d\mu = \int_\Omega f\,d\mu.
					\end{equation*}
				\end{theo}
				We shall also need the vector-valued version of the Lebesgue differentiation theorem. So, let us denote by $\lambda_n$ the Lebesgue measure on $\bR^n$, $\Sigma_n$ the completion of the Borel $\sigma$-algebra on $\bR^n$, and $B_n(u,r)$ the closed Euclidean ball of center $u$ and radius $r$ in $\bR^n$. Let $X$ be a Banach space. A function $f\colon \bR^n\to X$ is called \emph{locally Bochner integrable} if its restriction to any compact subset of $\bR^n$ is Bochner integrable.
				
				\begin{theo}
					\label{thm:Lebesgue-point} 
					Let $X$ be a Banach space and $f\colon \bR^n\to X$ be a locally Bochner integrable function. Then, there exists $A\in \Sigma_n$ such that $\lambda_n(\bR^n\setminus A)=0$ and for all $u\in A$
					\begin{equation*}
						\lim_{r\to 0^+}\frac{1}{\lambda_n(B_n(u,r))}\int_{B_n(u,r)}\norm{f(v)-f(u)}\,d\lambda_n(v)=0.
					\end{equation*}
					In particular, for all $u\in A$
					\begin{equation*}
						\lim_{r\to 0^+}\frac{1}{\lambda_n(B_n(u,r))}\int_{B_n(u,r)}f(v)\,d\lambda_n(v)=f(u).
					\end{equation*}
				\end{theo}
				
				\begin{proof}
					As usual, we may assume that $X$ is separable. So, let $(x_k)_{k=1}^\infty$ be a dense sequence in $X$. It follows from the scalar version of the Lebesgue differentiation theorem (see \cite{Rudin1987} for instance) that there exists $A\in \Sigma_n$ such that $\lambda_n(\bR^n\setminus A)=0$ and for all $u\in A$ and all $k\in \N$,
					\begin{equation*}
						\lim_{r\to 0^+}\frac{1}{\lambda_n(B_n(u,r))}\int_{B_n(u,r)}\norm{f(v)-x_k}\,d\lambda_n(v)=\norm{f(u)-x_k}.
					\end{equation*}
					Therefore, for all $u$ in $A$ and all $k$ in $\bN$,
					\begin{equation*}
						\limsup_{r\to 0^+}\frac{1}{\lambda_n(B_n(u,r))}\int_{B_n(u,r)}\norm{f(v)-f(u)}\,d\lambda_n(v)\le 2\norm{f(u)-x_k},
					\end{equation*}
					which yields the conclusion, since $(x_k)_{k=1}^\infty$ is dense in $X$.
				\end{proof}
				
				\chapter{Ramsey Theorem for Analysts}	
				\label{appendix:Ramsey}
				Given $\bM$ an infinite subset of $\bN$, we denote by $[\bM]^\omega:=\{ S\subset \bM \colon S\text{ is infinite}\}$ the set of all \emph{infinite} subsets of $\bM$. Given $k\in \bN$ and $\bM\in [\bN]^\omega$, we denote by $[\bM]^{k}:=\{ S\subset \bM \colon |S|=k\}$, the set of subsets of $\bM$ with exactly $k$ elements. For $\bM\in [\bN]^\omega$, we always list the elements of some $\mbar$ in $[\bM]^{\omega}$, resp. in $[\bM]^{k}$, in increasing order, meaning that if we write $\mbar:=\{m_1,m_2,m_3, \dots \}$, resp. $\mbar:=\{m_1,m_2,\dots, m_k\}$, we tacitly assume that $m_1<m_2<m_3<\dots $. In other words, we identify subsets of $\bN$ with strictly increasing, finite or infinite, subsequences of $\bN$. 
				The following theorem is most commonly called the infinite Ramsey theorem. It is one of many Ramsey-type theorems that seek to find monochromatic behavior in colored objects.   
				
				\begin{theo}[Infinite Ramsey Theorem]
					\label{thm:infinite_Ramsey} 
					Let $I$ be a \emph{finite} set. Then, for every $k\in \bN$ and $f \colon [\bN]^k \to I$, there exist $i\in I$ and $\bM \in [\bN]^\omega$ such that $f(\mbar)=i$ for all $\mbar \in [\bM]^k$. In particular, if $\cA\subseteq [\bN]^k$, then there exists $\bM\in [\bM]^\omega$ such that either $[\bM]^k\subseteq \cA$ or $\cA\cap [\bM]^k=\emptyset$.
				\end{theo}  
				
				\begin{rema}
					Note that the pruning Lemma \ref{lem:pruning1} is just a variant of the classical infinite Ramsey theorem, which could be directly used in Chapter \ref{chapter:Szlenk} if $X^*$ is separable and thus $D$, the weak neighborhood basis of $0$, can be taken countable. It is worth pointing out, however, that in the countable setting, the Ramsey Theorem is formally a bit stronger than Lemma \ref{lem:pruning1}.    
				\end{rema}
				
				An elementary application of the infinite Ramsey theorem is the following concentration phenomenon. Recall that a metric space is called \emph{totally bounded} if, for every $r>0$, it can be covered by finitely many balls of radius $r$. 
				
				\begin{coro}
					\label{cor:Ramsey_concentration}
					Let $M$ be a totally bounded metric space. If $k\in \bN$ and $f\colon [\bN]^k \to M$ is a function, then for any $\vep>0$, there exist a subset $B$ of $M$ of diameter less than $\eps$ and $\bM \in [\bN]^\omega$ such that $f([\bM]^k)\subset B$.   
				\end{coro}
				
				We say that a map $f\colon [\bN]^k \to (M,d_M)$ converges to $x\in M$, and we write $\lim_{n_1<n_2<\dots<n_k \to \infty} f(n_1,\dots,n_k)=x$, or simply $\lim_{\nbar\in [\bN]^k} f(\nbar)=x$, if for every $\vep>0$ there is $N\in \bN$ such that for all $n_k>n_{k-1}>\dots> n_1\ge N$, 
				\begin{equation*}
					d_M(f(n_1,\dots, n_k),x)<\vep.
				\end{equation*} 
				
				The Infinite Ramsey Theorem is typically used by analysts in the following form. The proof is an easy combination of Corollary \ref{cor:Ramsey_concentration} and a diagonal argument. 
				
				\begin{theo}
					Let $K$ be a compact metric space and  $k\in \bN$. Then, for every map $f\colon [\bN]^k \to K$, there exists $\bM \in [\bN]^\omega$ such that $\lim_{\nbar\in [\bM]^k} f(\nbar)$ exists in $K$.  
				\end{theo}
				
				\chapter{Asymptotic notions for Banach spaces}
				\label{appendix:asymptotic}
				
				Arguably, one of the most natural ways to describe what a local property of a Banach space looks like is to say that a property is local if it is preserved under crude finite representability. When it comes to defining what an asymptotic property is, there does not seem to be a canonical way to proceed. In this appendix, we present briefly several natural approaches that lead to three notions, each of which captures some asymptotic behavior of infinite-dimensional Banach spaces. These notions are best compared when expressed in terms of properties of weakly null trees and their branches.
				
				\section{Spreading models}
				\label{sec:spreading-models}
				
				The theory of spreading models was initiated by Brunel and Sucheston \cite{BrunelSucheston1974} in the 1970s. For an excellent exposition of the theory of spreading models we refer to the monograph by Beauzamy and Laprest\'e \cite{BeauzamyLapreste1984} (see also \cite{Guerre-book92}). 
				
				\begin{defi}
					\label{defi:spreading-model}
					A sequence $\xn$ in a Banach space $X$ is said to generate a sequence $(e_n)_{n=1}^\infty$ in a semi-normed space $S$, as a \emph{spreading model}, if for every $\vep>0$ and $k\in \bN$, there is $r\in\bN$ such that for all $r\le n_1<n_2<\dots<n_k$, and any scalars $a_1,\dots ,a_k\in[-1,1]$, 
					\begin{equation*}
						\Big|\big\|\sum_{i=1}^k a_i x_{n_i}\big\|_X - \big\|\sum_{i=1}^k a_i e_{i}\big\|_S\Big|<\vep.
					\end{equation*}   
				\end{defi}
				The beautiful application of Ramsey's theorem below provides the crucial step needed for the existence of spreading models.
				\begin{prop}
					\label{prop:spreading-model}
					Every bounded sequence $\xn$ in a separable Banach space admits a subsequence $\yn$ such that for all $k\ge 1$, $(a_i)_{i=1}^k\subset \bR$,
					\begin{equation*}
						\lim_{n_1<n_2<\dots<n_k\to\infty}\norm{a_1y_{n_1}+a_2y_{n_2}+\dots+a_ky_{n_k}}
					\end{equation*}
					exists.
				\end{prop}
				
				\begin{rema}
					The subsequence extracted in Proposition \ref{prop:spreading-model} is obviously not unique
				\end{rema}
				Given a bounded sequence $\xn$ in a separable Banach space and a subsequence $\yn$ given by Proposition \ref{prop:spreading-model}, we let 
				\begin{equation*}
					N(a_1,a_2,\dots,a_k):=\lim_{n_1<n_2<\dots<n_k\to\infty}\norm{a_1y_{n_1}+a_2y_{n_2}+\dots+a_ky_{n_k}},
				\end{equation*}
				where as usual the limit means that for all $\vep>0$, there is $r\in\bN$ such that for all $r\le n_1<n_2<\dots<n_k$,
				\begin{equation*}
					\Big|\norm{a_1y_{n_1}+a_2y_{n_2}+\dots+a_ky_{n_k}} - N(a_1,a_2,\dots,a_k)\Big|<\vep.
				\end{equation*}
				It is easy to see that that if $\ei$ is the canonical basis of $\coo$, the vector space of sequences in $\bR$ which eventually vanish, the formula 
				\begin{equation*}
					\norm{a_1e_{1}+a_2e_{2}+\dots+a_ke_{k}}_S:=N(a_1,a_2,\dots,a_k)
				\end{equation*} 
				defines a semi-norm. Therefore, any bounded sequence in a Banach space has a subsequence that generates a spreading model. 
				
				It is not too difficult to show that $\xn$ does not have a converging subsequence if and only if the semi-norm is actually a norm. 
				Also, by construction, the sequence $\ei$ is $1$-spreading in the sense that for all $k\ge 1$, $(a_i)_{i=1}^k\subset \bR$, and integers $1\le n_1<n_2<\cdots<n_k$,
				\begin{equation*}
					\norm{a_1e_{n_1}+a_2e_{n_2}+\dots+a_ke_{n_k}}_S=\norm{a_1e_{1}+a_2e_{2}+\dots+a_ke_{k}}_S.
				\end{equation*}
				The completion of $\coo$ for the norm $\norm{\cdot}_S$ is a Banach space $S$ which will simply be called a \emph{spreading model of $X$ generated by the sequence $\xn$}. We will refer to the sequence $\ei$ as the \emph{fundamental sequence} of the spreading model, and we write $S=[(e_i)_i]$. The fundamental sequence $\ei$ is always a basic sequence unless $\ei$ converges weakly in $S$ towards a non-zero element of $S$. 
				
				Spreading models generated by normalized weakly null sequences play a central role due to the following proposition. 
				
				\begin{prop}(\cite[Proposition 1, p. 24]{BeauzamyLapreste1984})
					\label{prop:weakly-null-spreading-model_V0}
					If $S=[(e_i)_i]$ is a spreading model generated by a normalized and weakly null sequence, then $\ei$ is a normalized $1$-suppression unconditional basis of $S$.
				\end{prop}
				
				The function $\varphi_S\colon k \in \bN \mapsto \|\sum_{i=1}^k e_i\|_S$, where $\ei$ is the fundamental sequence of a spreading model $S$ of a Banach space, is usually called the \emph{fundamental function} of $S$.
				
				\begin{prop}
					\label{prop:weakly-null-spreading-model}
					Let $\xn$ be a normalized weakly null sequence in a Banach space $X$. Then, for every $\vep>0$ there is a normalized weakly null basic subsequence $\yn$ of $\xn$ with basis constant $(1+\vep)$, generating a spreading model $S$ with fundamental function $\varphi_S$ such that for all $k\ge1$, for all $k\le n_1< n_2< \dots < n_{k}$ and for all $(\vep_i)_{i=1}^{k}\subset\{-1,1\}$ one has
					\begin{align*}
						\frac{1}{2(1+\vep)}\varphi_S(k)&\le  \frac1{1+\vep}  \Big\|\sum_{i=1}^{k}\vep_ie_i\Big\|_S\\
						&\le \Big\|\sum_{i=1}^{k}\vep_iy_{n_i}\Big\|_X\le
						(1+\vep)\Big\|\sum_{i=1}^{k}\vep_ie_i\Big\|_S\le2 (1+\vep)\varphi_S(k).\notag
					\end{align*}
				\end{prop}
				
				We will be particularly interested in $\co$-spreading models, where $\co$ is the space of real-valued sequences converging to $0$ equipped with the sup-norm. The following proposition follows from  \cite[Lemma 1, p. 73]{BeauzamyLapreste1984} and  \cite[Lemma 4, p. 75]{BeauzamyLapreste1984}.
				
				\begin{prop}
					\label{P:3}
					A Banach space $X$ has a spreading model isomorphic to $\co$ if and only if for all $\vep>0$, $X$ has a spreading model $S$ whose fundamental sequence $\ei$ is $(1+\vep)$-equivalent to the canonical basis of $\co$, i.e., that for all $k\ge 1$, $(a_i)_{i=1}^k\subset \bR$,
					\begin{equation*}
						\frac{1}{(1+\vep)}\sup_{1\le i\le k} \abs{a_i}\le \Big\|\sum_{i=1}^k a_ie_{i}\Big\|_S\le (1+\vep)\sup_{1\le i\le k} \abs{a_i}.
					\end{equation*}
				\end{prop}
				
				We will also need the following observation which can be found in \cite[Proposition 3, p. 79]{BeauzamyLapreste1984}.
				
				\begin{prop}
					\label{P:4}
					Let $\xn$ be a normalized weakly null sequence in a Banach space $X$ that generates a spreading model $S$. Then, $S$ is not isomorphic to $\co$ if and only if there exists a subsequence $\yn$ of $\xn$ such that
					\begin{equation}
						\lim_{k\to\infty}\ \ \inf_{n_1<\cdots<n_k}\ \ \inf_{(\vep_i)_{i=1}^k\in \{-1,1\}^k}\Big\|\sum_{i=1}^k\vep_i y_{n_i}\Big\|=+\infty.
					\end{equation}
				\end{prop}
				
				The spreading models of the classical sequence spaces are well understood:
				\begin{itemize}
					\item For $p\in[1,\infty)$, every spreading model of $\ell_p$ is linearly isometric to $\ell_p$.
					\item Every spreading model of $\co$ is linearly isomorphic to $\co$ and there is a spreading model of $\co$ that is not linearly isometric to $\co$. The main reason why the situation for $\co$ is different is due to the fact that $\co$ has two non-equivalent spreading bases: the canonical basis and the summing basis.
					\item For $p,q\in[1,\infty)$, every spreading model of $\ell_p\oplus \ell_q$ is linearly isomorphic to either $\ell_p$ or $\ell_q$. 
					\item Every spreading model generated by a normalized weakly convergent sequence in James space $\James$ is linearly isomorphic to $\ell_2$.
					\item Every spreading model of Tsirelson's space $\Tsi$ is isomorphic to $\ell_1$ and every spreading model of $\Tsi^*$ is isomorphic to $\co$.
					\item If $(F_n)_{n=1}^\infty$ is a sequence of finite-dimensional spaces and $p\in(1,\infty)$, then every spreading model of $(\sum_{n=1}^\infty F_n)_{\ell_p}$ is linearly isometric to $\ell_p$. Every spreading model of $(\sum_{n=1}^\infty F_n)_{\ell_1}$, resp. $(\sum_{n=1}^\infty F_n)_{\co}$, is linearly isomorphic to $\ell_1$, resp. $\co$.
				\end{itemize}
				This last result is important to understand the nature of spreading models. It is not difficult to see that every spreading model of a Banach space $X$ is finitely representable in $X$. Thus, every spreading model of $X$ will have every local property that $X$ has, and in many cases, it will have additional local properties that are not possessed by $X$. For instance, every spreading model of $(\sum_{n=1}^\infty \ell_\infty^n)_{\ell_2}$ is isomorphic to $\ell_2$, and hence is super-reflexive, has type $2$, has cotype $2$... Due to the presence of the $\ell_\infty^n$'s, $(\sum_{n=1}^\infty \ell_\infty^n)_{\ell_2}$ has none of these properties, but the asymptotic behavior of subsequences of a bounded sequence in $(\sum_{n=1}^\infty \ell_\infty^n)_{\ell_2}$ does not detect their presence. In a certain sense, spreading models capture the asymptotic behavior of Banach spaces.
				
				Let us conclude with another point of view about spreading models. Given a sequence $\xn$ in a Banach space $X$, one can consider its tree of subsequences, $(x_{\nbar})_{\nbar\in [\bN]^{<\infty}}$ where for all $\nbar:=\{n_1,\dots,n_i\}\in [\bN]^{<\infty}$, $x_{\{n_1,\dots,n_i\}}:=x_{n_i}$. Note that if $\xn$ is a normalized weakly null sequence, then $(x_{\nbar})_{\nbar\in [\bN]^{<\infty}}$ is a normalized weakly null tree. The subsequence tree has infinite height and for every $k\in\bN$, the first $k$ levels of the subsequence tree, namely $(x_{\nbar})_{\nbar\in [\bN]^{\le k}}$, is the tree of subsequences of length at most $k$.
				If $S=[(e_i)_i]$ is a spreading model generated by a bounded sequence $\xn$, one can always extract a subtree $(x_{\mbar})_{\mbar\in [\bM]^{<\omega}}$ of the subsequence tree $(x_{\nbar})_{\nbar\in [\bN]^{<\omega}}$ of $\xn$ such that for every $k\in \bN$ and $\vep>0$, there is a finite branch in $(x_{\mbar})_{\mbar\in [\bM]^{\le k}}$ that is $(1+\vep)$-equivalent to $(e_i)_{i=1}^k$. Moreover, for every $k$ and $\vep>0$ one can extract a subtree $(x_{\mbar})_{\mbar\in [\bM]^{\le k}}$ of $(x_{\nbar})_{\nbar\in [\bN]^{\le k}}$ such that every branch in $(x_{\mbar})_{\mbar\in [\bM]^{\le k}}$ is $(1+\vep)$-equivalent to $(e_i)_{i=1}^k$. 
				
				\section{Asymptotic models}
				\label{sec:asymptotic-models}
				The notion of asymptotic models, which is a natural generalization of the concept of spreading models, was formally introduced and studied by Odell and Halbeisen in the early 2000s, and we refer to \cite{HalbeisenOdell2004} for historical comments about (implicit) predecessors of this notion. Instead of looking at the asymptotic behavior of one sequence as in the spreading model theory, Odell and Halbeisen considered the asymptotic behavior of sequences of sequences.
				We need some convenient terminology.
				\begin{itemize}
					\item A infinite sequence of infinite sequences 
					\begin{align*}
						(x^{(1)}_{n})_{n=1}^{\infty} & :=(x_{1}^{(1)}, x_{2}^{(1)},\dots,x_{n}^{(1)},\dots)\\
						(x^{(2)}_{n})_{n=1}^{\infty} & :=(x_{1}^{(2)}, x_{2}^{(2)},\dots,x_{n}^{(2)},\dots)\\
						& \quad \vdots\\
						(x^{(k)}_{n})_{n=1}^{\infty} &:=(x_{1}^{(k)}, x_{2}^{(k)},\dots,x_{n}^{(k)},\dots)\\
						& \quad \vdots
					\end{align*} 
					in a Banach space $X$, is called an \emph{infinite array} or an \emph{array of infinite height} in $X$.
					%    \item An \emph{array of infinite height} in a Banach space $X$ is an infinite sequence of sequences of vectors.
					\item For an infinite array $\{(x^{(k)}_{n})_{n=1}^{\infty} \colon k\in \bN\}\subset X$ and $i\in[1,\infty]$, we call the sequence $(x^{(i)}_{n})_{n=1}^{\infty}$ the \emph{$i$-th row of the array}. Then, we say that an infinite array is \emph{weakly null} if all rows are weakly null. Similarly, an infinite array is \emph{basic} if every row is a basic sequence.
					\item A \emph{subarray} of an infinite array $\{(x^{(k)}_{n})_{n\in \bN}) \colon k\in \bN\}\subset X$, is an array of the form $\{(x^{(k)}_{j_s})_{s=1}^{\infty} \colon k\in \bN\}$, where $(j_s)_{s\in \bN}\subset \bN$ is a subsequence. The main point here is that in the definition of a subarray, we are taking the \emph{same} subsequence in each row.
				\end{itemize}

				%The following notion was introduced by Halbeisen and Odell \cite{HalbeisenOdell2004}.
				%\begin{defi}
				%\label{def:asymptotic-models}
				%A basic sequence $(e_i)_{i\in \bN}$ is called an \emph{asymptotic model} of a Banach space $X$, if there exist a normalized infinite array $\{(x^{(k)}_{n})_{n=1}^{\infty} \colon k\in\bN\}$ in $X$ and a null-sequence $(\vep_k)_{k\in \bN}\subset(0,1)$, so that for all $k\ge 1$, all $(a_i)_{i=1}^k \subset[-1,1]$ and $k\le n_1<n_2<\dots<n_k$,
				% \begin{equation*}
				% \abs{ \Big\|\sum_{i=1}^k a_i x^{(i)}_{n_i} \Big\| - \Big\|\sum_{i=1}^k a_i e_i \Big\| }<\vep_k.
				% \end{equation*}
				%\end{defi}
				
				\begin{defi}
					\label{defi:asymptotic-model}
					A infinite array $\{(x^{(k)}_{n})_{n=1}^{\infty} \colon k\in\bN\}$ in a Banach space $X$ is said to generate a sequence $(e_n)_{n=1}^\infty$ in a semi-normed space $A$, as an  \emph{asymptotic model}, if for every $\vep>0$ and $k\in \bN$, there is $r\in\bN$ such that for all $r\le n_1<n_2<\dots<n_k$, and any scalars $a_1,\dots ,a_k\in[-1,1]$, 
					\begin{equation*}
						\Big|\big\|\sum_{i=1}^k a_i x^{(i)}_{n_i}\big\|_X - \big\|\sum_{i=1}^k a_i e_{i}\big\|_A\Big|<\vep.
					\end{equation*}   
				\end{defi}
				
				A similar use of Ramsey's theorem gives the following proposition from \cite{HalbeisenOdell2004}.
				\begin{prop}
					\label{prop:Halbeisen-Odell}
					\cite[Proposition 4.1 and Remark 4.7.5]{HalbeisenOdell2004}
					Every weakly null normalized infinite array in a Banach space $X$ admits a subarray that generates a $1$-suppression unconditional asymptotic model for $X$.
				\end{prop}
				
				Note that every spreading model of $X$ is an asymptotic model of $X$. For instance, if a normalized weakly null sequence generates a spreading model $S=[(e_i)_i]$, then by considering countably many duplicates of the sequence, the corresponding infinite array generates $S=[(e_i)_i]$ as an asymptotic model. One of the main differences with spreading models is that, in general, an asymptotic model need not be spreading and hence more basic sequences can be realized as asymptotic models. A particularly telling result from \cite{HalbeisenOdell2004} states that every normalized and bimonotone basic sequence is $1$-equivalent to an asymptotic model of $\co$. This is a striking difference with spreading models once one remembers that if $S=[(e_i)_i]$ is a spreading model of $\co$, then $\ei$ is either equivalent to the canonical basis or to the summing basis of $\co$. 
				
				It is natural to think of an infinite array $\{(x^{(k)}_{n})_{n=1}^{\infty} \colon k\in\bN\}$ in a Banach space $X$, as a tree $(x_{\nbar})_{\nbar\in [\bN]^{<\infty}}$ of infinite height where for all $\nbar:=\{n_1,\dots,n_i\}\in [\bN]^{<\infty}$, $x_{\{n_1,\dots,n_i\}}:=x_{n_i}^{(i)}$. Note that if one starts with a normalized weakly null infinite array, then $(x_{\nbar})_{\nbar\in [\bN]^{<\infty}}$ is a normalized weakly null tree.
				If $A=[(e_i)_i]$ is an asymptotic model generated by an normalized weakly null infinite array, then for every $k\in \bN$ and $\vep>0$, one can always extract a subtree $(x_{\mbar})_{\mbar\in [\bM]^{\le k}}$ of $(x_{\nbar})_{\nbar\in [\bN]^{\le k}}$ such that every branch in $(x_{\mbar})_{\mbar\in [\bM]^{\le k}}$ is $(1+\vep)$-equivalent to $(e_i)_{i=1}^k$.

				\section{Joint spreading models}
				\label{sec:joint-spreading-models}
				
				The notion of joint spreading models, introduced in 2020 by Argyros, Georgiou, Lagos, and Motakis in \cite{AGLM2020}, lies somewhere between spreading model theory and asymptotic model theory in the sense that it deals with arrays of finite height, i.e., with finite sequences of infinite sequences such as 
				\begin{align*}
					(x^{(1)}_{n})_{n=1}^{\infty} & :=(x_{1}^{(1)}, x_{2}^{(1)},\dots,x_{n}^{(1)},\dots)\\
					(x^{(2)}_{n})_{n=1}^{\infty} & :=(x_{1}^{(2)}, x_{2}^{(2)},\dots,x_{n}^{(2)},\dots)\\
					& \quad \vdots\\
					(x^{(k)}_{n})_{n=1}^{\infty} &:=(x_{1}^{(k)}, x_{2}^{(k)},\dots,x_{n}^{(k)},\dots)
				\end{align*}
				for some $k\in\bN$. The terminology introduced for infinite arrays applies to finite arrays in the obvious way.
				
				In order to define joint spreading models, we need the notion of strict plegmas.
				
				\begin{defi}[Strict plegmas]
					\label{def:plegmas}
					\cite[Definition 3]{AKT2013}
					For $\ell\in\bN$ and $k\in\bN$, a finite sequence $\{\sbar_i\}_{i=1}^\ell$ of sequences of natural numbers of length $k$ 
					\begin{align*}
						\sbar_1 & :=(s_{1}^{(1)}, s_{2}^{(1)},\dots,s_{k}^{(1)})\in \bN^k\\
						\sbar_1 & :=(s_{1}^{(2)}, s_{2}^{(2)},\dots,s_{k}^{(2)})\in \bN^k\\
						& \quad \vdots\\
						\sbar_\ell &:=(s_{1}^{(\ell)}, s_{2}^{(\ell)},\dots,s_{k}^{(\ell)})\in \bN^k
					\end{align*}
					is called a \emph{strict plegma} if
					\begin{equation*}
						s^{(1)}_1 < s_1^{(2)} < \cdots< s_1^{(k)} < s_2^{(1)}< s_2^{(2)}<\cdots< s_2^{(k)}<\cdots<s_\ell^{(1)} < s^{(2)}_\ell < \cdots < s^{(k)}_\ell.
					\end{equation*}
				\end{defi}
				
				The notion of joint spreading models generalizes the notion of spreading models.
				
				%\begin{defi}[Joint spreading models]
				%\label{def:joint-spreading-models}
				%\cite[Definition 3.1]{AGLM2020}
				%Let  $\big\{(x_j^{(i)})_{j=1}^\infty\colon  1\le  i\le  k\big\} $ and $\big\{ (e_j^{(i)})_{j=1}^\infty \colon 1\le  i\le   k \big\}$ be two normalized arrays of height $k$ in the Banach spaces $X$ and $E$, respectively, whose rows are basic. We say that a normalized and basic finite array $\big\{(x_j^{(i)})_{j=1}^\infty\colon  1\le  i\le  k\big\}$ \emph{generates the normalized and basic finite array $\big\{ (e_j^{(i)})_{j=1}^\infty \colon 1\le  i\le   k \big\}$ as a joint spreading model} if there exists a null sequence of positive real numbers $(\vep_m)_{m=1}^\infty$ so that for every $m\in\bN$, every plegma $\{\sbar_i\}_{i=1}^k:=\{(s_{1}^{(i)}, s_{2}^{(i)},\dots,s_{m}^{(i)})\}_{i=1}^k$ with $\min(\bar{s}_1)=s^{(1)}_1 \ge m$ and sequences of scalars $\{(a_j^{(i)})_{j=1}^m \colon 1\le i \le k\}$ in $[-1,1]$ we have
				%\begin{equation*}
				%    \Bigg|\Big\|\sum_{i=1}^k \sum_{j=1}^m a_j^{(i)}x^{(i)}_{s^{(i)}_j}\Big\|_X - \Big\|\sum_{i=1}^k \sum_{j=1}^m a_j^{(i)}e^{(i)}_{j}\Big\|_E\Bigg|<\vep_m.
				%\end{equation*}
				%\end{defi}
				
				\begin{defi}[Joint spreading models]
					\label{def:joint-spreading-models}
					%\cite[Definition 3.3]{AGMM2022}
					A finite array $\big\{(x_j^{(i)})_{j=1}^\infty\colon  1\le i\le \ell\big\}$ in a Banach space $X$ is said to generate a finite array $\big\{ (e_j^{(i)})_{j=1}^\infty \colon 1\le i\le  \ell \big\}$ in a semi-normed space $E$ as a \emph{joint spreading model} if for every $\vep>0$ and every $k\in \bN$, there is an $r\in\bN$, such that for every strict plegma $\{\sbar_i\}_{i=1}^\ell:=\{(s_{1}^{(i)}, s_{2}^{(i)},\dots,s_{k}^{(i)})\}_{i=1}^\ell$ with $\min(\bar{s}_1)=s^{(1)}_1 \ge r$ and sequences of scalars $\{(a_j^{(i)})_{j=1}^k \colon 1\le i \le \ell\}$ in $[-1,1]$ we have
					\begin{equation*}
						\Bigg|\Big\|\sum_{i=1}^\ell \sum_{j=1}^k a_j^{(i)}x^{(i)}_{s^{(i)}_j}\Big\|_X - \Big\|\sum_{i=1}^\ell \sum_{j=1}^k a_j^{(i)}e^{(i)}_{j}\Big\|_E\Bigg|<\vep.
					\end{equation*}
				\end{defi}
				
				The following diagram helps visualize the position of the coefficients that appear in $\Big\|\sum_{i=1}^\ell \sum_{j=1}^k a_j^{(i)}x^{(i)}_{s^{(i)}_j}\Big\|_X$:
				\begin{align*}
					a_1^{(1)} x_{s_1^{(1)}}^{(1)} &                               &         &                               &  a_2^{(1)} x_{s_2^{(1)}}^{(1)} &                               &         &                               & & & & a_k^{(1)} x_{s_k^{(1)}}^{(1)} &                               &         &                               & \\
					& a_1^{(2)} x_{s_1^{(2)}}^{(2)} &         &                               &                               &                             & a_2^{(2)} x_{s_2^{(2)}}^{(2)} &         &                               &                               & &  &a_k^{(2)} x_{s_k^{(2)}}^{(2)} &           &                           &   \\
					&                               & \ddots  &                                &                              &                               &                               &                               & \ddots  &                                &                              &                               &  \ddots    &                              & \\
					&                               &          & a_1^{(\ell)} x_{s_1^{(\ell)}}^{(\ell)} &                               &                               &          &                               &                               &          & a_2^{(\ell)} x_{s_2^{(\ell)}}^{(\ell)} &                               &                               &          &a_k^{(\ell)} x_{s_k^{(\ell)}}^{(\ell)} & \\
				\end{align*}
				Joint spreading models are naturally related to spreading models as well as asymptotic models. 
				
				On one hand, if $\ell=1$ in the definition above, then one recovers the notion of a spreading model generated by a sequence. Moreover, if $\big\{(x_j^{(i)})_{j=1}^\infty\colon  1\le i\le \ell\big\}$ generates $\big\{ (e_j^{(i)})_{j=1}^\infty \colon 1\le i\le \ell \big\}$ as a joint spreading model, then, for every $1\le i \le \ell$, $(e_j^{(i)})_{j=1}^\infty$ is the spreading model of $(x^{(i)}_j)_{j=1}^\infty$.
				
				On the other hand, if $\ell\in\bN$ and $\big\{ (x^{(i)}_j)_{j=1}^\infty \colon 1\le i \le \ell \big\}\subset X$ is a normalized weakly null array of height $\ell$, then we can extend this array to an infinite array $\big\{ (x^{(i)}_j)_{j=1}^\infty \colon i\in \bN \big\}$, by duplicating the first $\ell$ rows cyclically, i.e., we let for every $s\in\bN$ and $1\le i \le \ell$,
				\begin{equation*}
					x^{(s\ell+i)}_j := x^{(i)}_j.
				\end{equation*} 
				By Proposition \ref{prop:Halbeisen-Odell}, we can pass to a subarray $\big\{ (z^{(i)}_{j})_{j=1}^\infty \colon i\in\bN \big\}$ of $\big\{ (x^{(i)}_{j})_{j=1}^\infty \colon i\in\bN \big\}$ that generates an asymptotic model $A=[(e_j)_j]$. Now, if we let $e^{(i)}_{j} := e_{(j-1)\ell + i}$, for every $1\le i \le \ell$ and $j\in\bN$, or more visually,
				\begin{align*}
					(e^{(1)}_{j})_{j=1}^{\infty} & :=(e_{1}, e_{\ell+1}, e_{2\ell+1},\dots,e_{j\ell+1},\dots)\\
					(e^{(2)}_{j})_{j=1}^{\infty} & :=(e_{2}, e_{\ell+2}, e_{2\ell+2},\dots,e_{j\ell+2},\dots)\\
					& \quad \vdots\\
					(e^{(\ell)}_{j})_{j=1}^{\infty} &:=(e_{\ell}, e_{2\ell}, e_{3\ell},\dots,e_{(j+1)\ell},\dots),
				\end{align*}
				we observe that the array $\big\{ (e_j^{(i)})_{j=1}^\infty \colon 1\le i\le  \ell \big\}$ is the joint spreading model of $\big\{ (z^{(i)}_{j})_{j=1}^\infty \colon 1\le i \le \ell \big\}$.
				
				A joint spreading model $\big\{ (e_j^{(i)})_{j=1}^\infty \colon 1\le i\le  \ell \big\}$ is always spreading in a certain sense called \emph{plegma spreading} and under certain assumptions $\big\{ (e_j^{(i)})_{j=1}^\infty \colon 1\le i\le  \ell \big\}$ is a basic sequence when ordered properly. The following result can be found in \cite[Theorem 3.3, Proposition 3.7]{AGLM2020}. 
				
				\begin{theo}
					Every finite basic array $\big\{(x_j^{(i)})_{j=1}^\infty\colon  1\le  i\le  \ell\big\}$ generates a joint spreading model $\big\{ (e_j^{(i)})_{j=1}^\infty \colon 1\le i\le  \ell \big\}$. Moreover, if the array is weakly null, then $\big\{ (e_j^{(i)})_{j=1}^\infty \colon 1\le i\le  \ell \big\}$ is $1$-suppression unconditional and in particular for every $\vep>0$ and $k\in \bN$ there is $r\in\N$, so that for any strict plegma $(\bar{s}_i)_{i=1}^\ell$, $\bar{s}_i=(s^{(i)}_1, s^{(i)}_2,\dots, s^{(i)}_k)$, for $i=1,2,\dots, \ell$ with $r\le  s^{(1)}_1$ the family $\big\{(x^{(i)}_{s^{(i)}_j})_{j=1}^k \colon 1\le i\le \ell\big\}$ is $(1+\vep)$-suppression unconditional. 
				\end{theo}
				
				Joint spreading models generated by weakly null arrays also enjoy the following crucial suppression-unconditionality property.
				
				\begin{prop}
					\label{prop:jointspreadingmodels}
					Let $X$ be a Banach space and $\big\{(x^{(i)}_{j})_{j=1}^\infty \colon 1\le i\le \ell \big\}$ be a  normalized weakly null array of height $\ell$. Then, for every $\vep>0$ and $k\in\bN$ there exists $\bM\in[\bN]^\omega$ so that for every $i_1,\dots,i_k$ in $\{1,\dots,\ell\}$ (not necessarily different) and pairwise different $l_1,\dots,l_k\in \bM$ the sequence $(x_{l_j}^{(i_j)})_{j=1}^k$ is $(1+\vep)$-suppression unconditional.
				\end{prop}
				
				\begin{proof}
					As explained above, we may assume after passing to a subarray that $\big\{(x^{(i)}_{j})_{j=1}^\infty \colon 1\le i\le \ell \big\}$ generates a joint spreading model $\big\{ (e_j^{(i)})_{j=1}^\infty \colon 1\le i\le  \ell \big\}$ and we can find $r\in\N$, so that for any strict plegma $(\bar{s}_i)_{i=1}^\ell$, $\bar{s}_i=(s^{(i)}_1, s^{(i)}_2,\dots, s^{(i)}_k)$, for $i=1,2,\dots, \ell$ with $r\le  s^{(1)}_1$ the family $\big\{(x^{(i)}_{s^{(i)}_j})_{j=1}^k \colon 1\le i\le \ell\big\}$ is $(1+\vep)$-suppression unconditional. Let $\bL$ be the set that consists of all positive integers multiple of $2\ell$ that are greater than $r+\ell$.
					
					Let now $i_1,\dots,i_k$ in $\{1,\dots,\ell\}$ and $l_1,\dots,l_k$ be pairwise different elements of $\bL$.
					After reordering, we can assume $l_1< l_2< \dots< l_k$. Let $r_1< r_2< \dots< r_k$ be in $\bN$ so that $l_j=2\ell r_j$. We will now define a strict plegma $(\bar{s}_i)_{i=1}^\ell$, $\bar{s}_i :=(s^{(i)}_{j})_{j=1}^k$, as follows. First we define $s^{(i_j)}_j := l_j = 2\ell r_j$, for $j=1,2,\dots, k$. Then, since $l_{j+1}-l_{j}\ge 2\ell$, for every $j=1,\dots, k-1$ and $s_1^{(i_1)}> r+\ell$, we can find natural numbers $s^{(i_j)}_j< s^{(i_j+1)}_j< s^{(i_j+2)}_j<\dots s^{(\ell)}_j< s^{(1)}_{j+1}<\dots< s^{(i_{j+1})}_{j+1}$, numbers $r < s^{(1)}_1< s^{(2)}_1< \ldots < s^{(i_1-1)}_1< s^{(i_1)}_1$ and numbers $s^{(i_\ell)}_k< s^{(i_\ell+1)}_{k}<  \dots < s^{(\ell)}_k$, which means that the family $(\bar{s}_i)_{i=1}^\ell$ with $\bar{s}_i=(s_j^{(i)})_{j=1}^k$, for $i=1,2,\dots,\ell$ is a strict plegma. Thus, $\big\{(x^{(i)}_{s_j^{(i)}})_{j=1}^k \colon 1\le i \le \ell\big\}$ is $(1+\vep)$-suppression unconditional and $(x^{(i_j)}_{l_j})_{j=1}^k$ is just a subsequence of it.
				\end{proof}

				\section{Asymptotic structure}
				\label{sec:asymptotic-structure}
				
				The notion of asymptotic structure was formally introduced in the mid-1990s by B. Maurey, V. Milman, and N. Tomczak-Jaegermann \cite{MMTJ}. The asymptotic structure theory is about the collection of finite-dimensional subspaces of a Banach space that can be realized in a very specific asymptotic fashion. It is not so much about finite-dimensional spaces but rather about their normalized bases. 
				
				For a Banach space $X$, we denote by $\cof(X)$ the set of finite-codimensional subspaces of $X$. 
				\begin{defi}[Asymptotic structure]
					Let $k\in \bN$, $X$ be a Banach space, and $E$ be a $k$-dimensional Banach space with a normalized basis $(e_i)_{i=1}^k$. The space \emph{$E$ is in the $k$-th asymptotic structure of $X$} (or \emph{$X$ contains $E$ in its $k$-th asymptotic structure}), and we write $E\in \{X\}_k$, if for all $\vep>0$,
					\begin{align*}
						\label{eq:asymptotic-structure}
						\forall X_1\in \cof(X) & \quad \exists  x_1 \in S_{X_1},\\
						\forall X_2\in \cof(X) &  \quad \exists  x_2\in S_{X_2},\\
						& \vdots \\
						\forall X_k\in \cof(X) & \quad \exists  x_k\in S_{X_k},
					\end{align*}
					such that $(x_i)_{i=1}^k$ is $(1+\vep)$-equivalent to $(e_i)_{i=1}^k$.
				\end{defi}
				
				A finite-dimensional space $E$ is thus in the asymptotic structure of $X$ if it has a normalized basis that is equivalent, with arbitrarily good equivalence constant, to finite sequences of vectors in $X$ that can be found ``everywhere far away'' in $X$. The collection of all spaces in $\cup_{k\in \bN} \{X\}_k$ can be thought of as containing the information related to the asymptotic structure of $X$.
				Of course, it is much harder for a finite-dimensional subspace to be in the asymptotic structure of a Banach space rather than merely being almost isometrically isomorphic to a subspace in it. For instance, Dvoretzky's theorem asserts that one can always find arbitrarily good copies of the $\ell_2^n$'s inside any infinite-dimensional Banach space. Still, there is no reason why these copies could be found everywhere, far away. As we will see, for instance, the spaces $\ell_p$ for $p\neq 2$ do not contain the $\ell_2^n$'s in their asymptotic structure. 
				
				The definition of the $k$-th asymptotic structure of $X$ involves a finite sequence of normalized vectors that can be thought of as the outcome of a finite game between two players: a first player (Player A) chooses a finite-codimensional subspace in a Banach space to which the second player (Player B) responds with a normalized vector in this finite-codimensional subspace. The possible outcomes of this two-player game can be collected in the form of a normalized tree $(x_{\bar{Y}})_{\bar{Y}\in \cof(X)^{<\infty}}$ in $X$ indexed by finite non-empty sequences of finite-codimensional subspaces of $X$, denoted by $\bar{Y}:=(Y_1,\dots,Y_k)\in \cof(X)^{<\omega}$. If $\cof(X)$ is directed by reverse inclusion, the tree obtained is weakly null in the sense that $(x_{(\bar{Y},Z)})_{Z\in \cof(X)}$ is a weakly null net. Therefore, a normalized basis $(e_i)_{i=1}^k$ is in the $k$-th asymptotic structure of $X$ if and only if we can find for every $\vep>0$ a weakly null normalized unrooted tree $(x_{\bar{Y}})_{\bar{Y}\in \cof(X)^{\le k}}$ all of whose branches are $(1+\vep)$-equivalent to the basis $(e_i)_{i=1}^k$. Another way to put it is in terms of winning strategies. For $k\in\bN$, a $k$-dimensional space $E$ with a normalized basis $(e_i)_{i=1}^k$ is in the $k$-th asymptotic structure of $X$ if and only if for every $\vep>0$ Player B has a winning strategy to obtain a sequence $(x_i)_{i=1}^k$ that is $(1+\vep)$-equivalent to $(e_i)_{i=1}^k$. Indeed, no matter what the moves of Player A are, Player B always has a response that guarantees that the sequence of normalized vectors that the game will produce is $(1+\vep)$-equivalent to $(e_i)_{i=1}^k$. Note that, if $X$ is separable and does not contain an isomorphic copy of $\ell_1$, we can replace weakly null trees indexed by $\cof(X)^{\le k}$ by weakly null trees indexed by $[\bN]^{\le k}$.
				
				Now, for a property that makes sense for finite normalized bases, we can say that a Banach space $X$ has an asymptotic structure with property $\cP$ if, for every $k\in \bN$, every normalized basis in the $k$-th asymptotic structure has $\cP$. In terms of the game approach, this means that Player A, who picks finite-codimensional subspaces, has a winning strategy that guarantees that every outcome of the game, namely the finite sequence of normalized vectors $(x_i)_{i=1}^k$ chosen by Player B, has property $\cP$. If one prefers to think in terms of the tree approach, a Banach space $X$ has an asymptotic structure with property $\cP$ if for every $\vep>0$ and every weakly null normalized unrooted tree $(x_{\bar{Y}})_{\bar{Y}\in \cof(X)^{\le k}}$ there is a branch that is $(1+\vep)$-equivalent to a normalized basis with property $\cP$. It follows from a Ramsey-type argument (either the classical version of Ramsey Theorem from Appendix \ref{appendix:Ramsey} for trees indexed by the natural numbers or the pruning lemma \ref{lem:pruning1} for trees indexed by directed sets) that if one can find one branch with this property, then one can find a subtree all of whose branches have the same property. Here is a concrete example of a property that is needed in Chapter \ref{chapter:diamonds}.
				
				\begin{defi}
					\label{def:UAS}
					We say that a Banach space has \emph{an unconditional asymptotic structure with constant $C\ge 1$}, or a \emph{$C$-unconditional asymptotic structure}, if for all $k\in\bN$,
					\begin{align*}
						\exists X_1\in \cof(X) & \quad \forall  x_1 \in S_{X_1},\\
						\exists X_2\in \cof(X) &  \quad \forall  x_2\in S_{X_2},\\
						& \vdots \\
						\exists X_k\in \cof(X) & \quad \forall  x_k\in S_{X_k},
					\end{align*}
					one has that $(x_i)_{i=1}^k$ is $C$-unconditional.
				\end{defi}
				
				The definition of asymptotic structure naturally involves normalized unrooted trees. However, there is a bit of flexibility here as we can also work with bounded rooted trees if needed, as seen in the next proposition.
				\begin{prop}
					\label{prop:UAS-tree}
					Let $X$ be a Banach space and $C\in [1,\infty)$. The following assertions are equivalent.
					\begin{enumerate}[(i)]
						\item $X$ has a $C$-unconditional asymptotic structure.
						\item For every directed set $D$, $k\in \bN$ and weakly null tree $(x_t)_{t\in D^{\le k}}$ in $B_X$, there exists $\tau\in D^k$ such that for every $\vep\in\{-1,1\}^{k+1}$,  
						\begin{equation*}
							\Big\|\sum_{i=0}^k \vep_i x_{t_{\restriction i}} \Big\|\le C \Big\| \sum_{i=0}^k x_{t_{\restriction i}} \Big\|.
						\end{equation*} 
						\item For every directed set $D$, $k\in \bN$ and weakly null tree $(x_t)_{t\in D^{\le k}}$ in $B_X$, there exists a pruning $\phi\colon D^{\le k} \to D^{\le k}$ such that for every $\tau\in D^k$ and $\vep\in\{-1,1\}^{k+1}$,  
						\begin{equation*}
							\Big\|\sum_{i=0}^k \vep_i x_{\phi(t)_{\restriction i}} \Big\|\le C \Big\|\sum_{i=0}^k x_{\phi(t)_{\restriction i}} \Big\|. 
						\end{equation*}  
					\end{enumerate}
				\end{prop}
				
				The asymptotic properties are typically weaker than the properties from which they are derived. For instance, having an unconditional asymptotic structure for a Banach space is strictly weaker than having an unconditional basic sequence. Indeed, the Argyros-Deliyanni space \cite{ArgyrosDeliyanni1997} does not have any unconditional basic sequence but has an unconditional asymptotic structure.
				
				Nevertheless, some classical phenomena have asymptotic versions. It is well known that for $p\in[1,\infty]$, then if a Banach space contains $C$-isomorphic copies of the $\ell_p^n$'s for some $C\ge 1$, then for every $\vep>0$, $X$ contains $(1+\vep)$-isomorphic copies of the $\ell_p^n$'s. For $p\in\{1,\infty\}$, this due to James \cite{James64} and for the the other values of $p\in(1,\infty)$ this is essentially due to Krivine \cite{Krivine1976} (see \cite[Theorem 3.6]{Rosenthal1978b} for a precise quantitative statement). These local non-distortability results have asymptotic analogs. But first, we need one more definition.
				
				\begin{defi}
					For $1\le p\le\infty$, we say that \emph{$\ell_p^n$ is in the $n$-th asymptotic structure of $X$ up to a constant $C\ge 1$}, if there exists an $n$-dimensional space $E$ with a normalized basis that is $C$-equivalent to the $\ell_p^n$-unit vector basis, so that $E\in\{X\}_n$. If there is a $C\ge 1$, so that for all $n\in\N$, $\ell_p^n$ is in the $n$-th asymptotic structure of $X$ up to the constant $C$, then we say that \emph{the asymptotic structure of $X$ contains uniformly $\ell_p^n$, $n\in\bN$}.
				\end{defi}
				
				\begin{theo}
					\label{theo:asymptotic-non-distortability}
					Let $p\in[1, \infty]$. If the asymptotic structure of $X$ contains uniformly $\ell_p^n$, $n\in\bN$, then $\ell_p^n\in\{X\}_n$ for all $n\in\bN$.
				\end{theo}
				
				\begin{proof}[Sketch of proof]
					We sketch the proof for $p=\infty$.  First, it was observed in \cite[Section 1.8.3]{MMTJ} that if $E\in\{X\}_n$ with a normalized basis $(e_j)_{j=1}^n$, then any subspace $F$ of $E$ that is spanned by a normalized block basis $(f_j)_{j=1}^m$ of $(e_j)_{j=1}^n$ is in $\{X\}_m$. Secondly, James' local version of the nondistortability of $\co$ asserts that for any $m\in\bN$, $C>1$ and any $\vep>0$, there is an $n:=n(m,C,\vep)$ so that if $E$ is an $n$-dimensional space with a normalized basis $(e_j)_{j=1}^n$ that is $C$-equivalent to the $\ell_\infty^n$-unit vector basis, then $(e_j)_{j=1}^n$ admits a block basis of length $m$ that is $(1+\vep)$-equivalent to the $\ell_\infty^m$-unit vector basis. The conclusion follows. A similar argument works for $p=1$ using James' local version of the nondistortability of $\ell_1$. For the other values of $p$, it follows from the quantitative version of Krivine's Theorem proved by Rosenthal in \cite[Theorem 3.6]{Rosenthal1978b}.
				\end{proof}
				
				In \cite[Section 1.7]{MMTJ}, Maurey, Milman, and Tomczak-Jaegermann introduced the general class of asymptotic-$\ell_p$ spaces as follows.
				
				\begin{defi}
					A Banach space is called an \emph{asymptotic-$\ell_p$ space}, for $p\in[1,\infty]$, if there is a constant $C\ge 1$ such that, for every $k\in \bN$, every normalized basis $(e_i)_{i=1}^k$ in $\{X\}_k$ is $C$-equivalent to the canonical basis of $\ell_p^k$. Note that for $p=\infty$, we talk about asymptotic-$\co$ spaces rather than asymptotic-$\ell_\infty$ spaces.
				\end{defi}
				
				Typical examples of asymptotic-$\ell_p$ spaces are $\ell_p$-sums of finite-dimensional spaces or Tsirelson-type spaces. 
				
				\begin{rema}\,
					\begin{itemize}
						\item Milman and Tomczak-Jaegermann gave a definition of asymptotic-$\ell_p$ spaces for spaces with a basis in \cite{MTJ1993} prior to its coordinate-free generalization given above. It reads as follows. Let $X$ be a Banach space with a basis $(e_i)_{i=1}^\infty$ and $p\in[1,\infty)$. The basis $(e_i)_{i=1}^\infty$ of $X$ is said to be asymptotic-$\ell_p$ if there exist positive constants $C_1,C_2>0$ such that for all $n\in \bN$, there exists $r\in \bN$ with the property that whenever $r\le x_1<x_2<\dots x_n$ are vectors in $X$ with consecutive supports, then 
						\begin{equation*} 
							\frac{1}{C_2}\Big(\sum_{i=1}^n \norm{x_i}^p\Big)^{\frac{1}{p}} \le \Big\|\sum_{i=1}^n x_n\Big\| \le C_1\Big(\sum_{i=1}^n \norm{x_i}^p\Big)^{\frac{1}{p}}.
						\end{equation*}
						With this definition, it is then easy to see that $\Tsi^*$ is an asymptotic-$\co$ space while $\Tsi$ is an asymptotic-$\ell_1$ space.
						\item For $p\in[1,\infty)$, it was shown in \cite[Section 3.1]{MMTJ} that if in the definition of an asymptotic-$\ell_p$ space, we relax the assumption to simply asking that the finite-dimensional space generated by $(e_i)_{i=1}^k$ is $C$-isomorphic to $\ell_p^k$ then we get the same notion.
					\end{itemize}
					
					For simplicity of the ensuing discussion, we will always assume that the Banach spaces are separable and do not contain an isomorphic copy of $\ell_1$. For a Banach space, being an asymptotic-$\ell_p$ space is a much stronger requirement than containing uniformly the $\ell_p^n$'s in the asymptotic structure. This is particularly clear if one thinks in terms of trees. Indeed, if a Banach space is asymptotic-$\ell_p$, then, for every $k\in\bN$, \emph{every} normalized weakly null tree $(x_{\nbar})_{\nbar\in [\bN]^{\le k}}$ has a subtree $(x_{\mbar})_{\mbar\in [\bM]^{\le k}}$ such that every branch in $(x_{\mbar})_{\mbar\in [\bM]^{\le k}}$ is $(1+\vep)$-equivalent to the canonical basis of $\ell_p^k$. On the other hand, if contains uniformly the $\ell_p^n$'s in the asymptotic structure, it only means that, for every $k\in\bN$, there \emph{exists} a normalized weakly null tree $(x_{\nbar})_{\nbar\in [\bN]^{\le k}}$ whole of whose branches are $(1+\vep)$-equivalent to the canonical basis of $\ell_p^k$.  
					
					To compare the notion of asymptotic structure with spreading models and asymptotic models, it is also enlightening to think in terms of trees. Asymptotic structure is closely related to the behavior of branches of \emph{arbitrary} weakly null trees of finite but arbitrarily large height. Spreading models are related, as we have seen, with the behavior of branches of weakly null trees of finite but arbitrarily large height that are the sublevel trees of the subsequence tree of a normalized weakly null sequence, i.e., trees of the form $(x_{\nbar})_{\nbar\in [\bN]^{\le k}}$ where $x_{\{n_1,\dots,n_i\}}:= x_{n_i}$ with $\xn$ normalized and weakly null. 
					Similarly, asymptotic models are related with the behavior of branches of weakly null trees of finite but arbitrarily large height that are the sublevel trees of the natural tree induced by an infinite normalized and weakly null array, i.e., trees of the form $(x_{\nbar})_{\nbar\in [\bN]^{\le k}}$ where $x_{\{n_1,\dots,n_i\}}:= x^{(i)}_{n_i}$ with $\{(x^{(k)}_{n})_{n=1}^{\infty} \colon k\in\bN\}$ a normalized and weakly null array. From this point of view, it is immediate that if $X$ has an $\ell_p$-spreading model generated by a weakly null sequence, then $X$ has an $\ell_p$-asymptotic model generated by a weakly null array, and ultimately $X$ contains the $\ell_p^k$'s in its asymptotic structure. Also, if a Banach space is an asymptotic-$\ell_p$ space, then every asymptotic model generated by a weakly null array, and in turn every spreading model generated by a weakly null sequence, is isomorphic to $\ell_p$. Usually, none of these implications can be reversed. However, a remarkable theorem of Freeman, Odell, Sari and Zheng \cite{FOSZ2018} states that in the extreme case $p=\infty$, if every asymptotic model generated by a weakly null array of a Banach space $X$ (that is separable and does not contain an isomorphic copy of $\ell_1$) is isomorphic to $\co$, then $X$ must be asymptotic-$\co$.
					We refer the interested reader to the papers of Argyros, Georgiou, Lagos, Manoussakis and Motakis \cite{AGLM2020,AM2020,AGM2020,AGMM2022} for a rather complete panorama of the precise relationships between spreading models, asymptotic models, and asymptotic structure.
				\end{rema}
			\end{appendix}

			%%%%%%%%%%%%%%%%%%%%%%%%%%%%%%%%%%%%%%%%%%%%%%%%%%%%%%%%%%%%%%%%%%%%%%%%%%%%%%%%%%%%%%%%%%%%%%%%%%%%%%%%%

			%\addcontentsline{toc}{chapter}{\indexname}
			%\printindex
			
		\end{document}